\documentclass{amsbook}
\usepackage{amssymb}
\usepackage{mathrsfs}
\usepackage{stmaryrd} 
\usepackage{bbm}
\usepackage{oldgerm}
\usepackage[english]{babel}
\usepackage[T1]{fontenc}
\usepackage[latin1]{inputenc}
\usepackage[all]{xy}
\usepackage{hyperref}
\usepackage{upref}
\usepackage{xcolor}
\usepackage{microtype}
\usepackage{enumitem}

\usepackage[Symbol]{upgreek}

\hypersetup{
    colorlinks,
    linkcolor={red!50!black},
    citecolor={blue!50!black},
    urlcolor={blue!80!black}
}

\newtheorem{teo}[subsection]{Theorem}
\newtheorem{prop}[subsection]{Proposition}
\newtheorem{cor}[subsection]{Corollary}
\newtheorem{lem}[subsection]{Lemma}

\theoremstyle{definition}

\newtheorem{defi}[subsection]{Definition}
\newtheorem{rema}[subsection]{Remark}
\newtheorem{remas}[subsection]{Remarks}
\newtheorem{exemple}[subsection]{Example}
\newtheorem{exemples}[subsection]{Examples}

\numberwithin{equation}{subsection}

\newcommand{\gtimes}{\stackrel{\leftarrow}{\times}}

\mathchardef\mhyphen="2D

\DeclareMathSymbol{\mlq}{\mathord}{operators}{``}
\DeclareMathSymbol{\mrq}{\mathord}{operators}{`'}

\newcommand{\mB}{{\mathbb B}}

\newcommand{\mC}{{\mathbb C}}

\newcommand{\mH}{{\mathbb H}}
\newcommand{\mI}{{\mathbb I}}

\newcommand{\mQ}{{\mathbb Q}}

\newcommand{\mM}{{\mathbb M}}
\newcommand{\mN}{{\mathbb N}}

\newcommand{\mT}{{\mathbb T}}
\newcommand{\mX}{{\mathbb X}}
\newcommand{\mY}{{\mathbb Y}}
\newcommand{\mZ}{{\mathbb Z}}

\newcommand{\mG}{{\mathbb G}}
\newcommand{\mK}{{\mathbb K}}
\newcommand{\mU}{{\mathbb U}}
\newcommand{\mV}{{\mathbb V}}

\newcommand{\mi}{{\mathbbm i}}
\newcommand{\mj}{{\mathbbm j}}

\newcommand{\Indoplus}{\mlq\mlq\bigoplus \mrq\mrq}
\newcommand{\IndSym}{{\rm IS}}

\newcommand{\bA}{{\bf A}}
\newcommand{\bB}{{\bf B}}
\newcommand{\bC}{{\bf C}}
\newcommand{\bD}{{\bf D}}

\newcommand{\bIndAlg}{{\bf Ind\mhyphen Alg}}
\newcommand{\bIndMod}{{\bf Ind\mhyphen Mod}}
\newcommand{\bK}{{\bf K}}
\newcommand{\bL}{{\bf L}}

\newcommand{\bT}{{\bf T}}
\newcommand{\bP}{{\bf P}}
\newcommand{\bQ}{{\bf Q}}

\newcommand{\Et}{{\bf \acute{E}t}}
\newcommand{\Sch}{{\bf Sch}}

\newcommand{\FLS}{{\bf FLS}}
\newcommand{\Ens}{{\bf Set}}

\newcommand{\bHom}{{\bf Hom}}
\newcommand{\bCat}{{\bf Cat}}

\newcommand{\bRep}{{\bf Rep}}

\newcommand{\bAlg}{{\bf Alg}}
\newcommand{\bMod}{{\bf Mod}}
\newcommand{\Ind}{{\bf Ind}}

\newcommand{\bInd}{{\bf Ind}}
\newcommand{\red}{{\rm red}}

\newcommand{\colim}{{\underset{\longrightarrow}{\lim}}}
\newcommand{\indcolim}{{\mlq\mlq\colim \mrq\mrq}}

\newcommand{\IC}{{\rI\cC}}
\newcommand{\IcK}{{\rI\cK}}
\newcommand{\IcL}{{\rI\cL}}
\newcommand{\IrK}{{\rI\rK}}
\newcommand{\IrL}{{\rI\rL}}
\newcommand{\Idelta}{{\rI\delta}}
\newcommand{\Ifd}{{\rI\fd}}

\newcommand{\intern}{{\diamond}}
\newcommand{\lgg}{{\ttg}}

\newcommand{\et}{{\rm \acute{e}t}}
\newcommand{\fet}{{\rm f\acute{e}t}}

\newcommand{\zar}{{\rm zar}}

\newcommand{\rk}{{\rm rk}}
\newcommand{\coh}{{\rm coh}}

\newcommand{\Tot}{{\rm Tot}}

\newcommand{\Dolb}{{\rm Dolb}}
\newcommand{\sol}{{\rm sol}}
\newcommand{\sDolb}{{\rm sDolb}}
\newcommand{\ssol}{{\rm ssol}}

\newcommand{\qsol}{{\rm qsol}}

\newcommand{\rf}{{\rm f}}

\newcommand{\Dol}{{\rm Dol}}
\newcommand{\Spec}{{\rm Spec}}

\newcommand{\Spf}{{\rm Spf}}

\newcommand{\arr}{{\rm Arr}}
\newcommand{\ob}{{\rm Ob}}

\newcommand{\tor}{{\rm tor}}

\newcommand{\coker}{{\rm coker}}

\newcommand{\im}{{\rm im}}

\newcommand{\tot}{{\rm tot}}

\newcommand{\cont}{{\rm cont}}

\newcommand{\Gr}{{\rm Gr}}

\newcommand{\id}{{\rm id}}

\newcommand{\rb}{{\rm b}}
\newcommand{\Sym}{{\rm Sym}}

\newcommand{\Hom}{{\rm Hom}}

\newcommand{\End}{{\rm End}}

\newcommand{\Ext}{{\rm Ext}}

\newcommand{\bHM}{{\bf HM}}
\newcommand{\bSM}{{\bf SM}}
\newcommand{\bIMC}{{\bf IMC}}
\newcommand{\bIMIC}{{\bf IMIC}}
\newcommand{\bIndMC}{{\bf Ind\mhyphen MC}}
\newcommand{\bIndMIC}{{\bf Ind\mhyphen MIC}}
\newcommand{\bIndHM}{{\bf Ind\mhyphen HM}}
\newcommand{\bMIC}{{\bf MIC}}
\newcommand{\bIH}{{\bf IH}}
\newcommand{\bMC}{{\bf MC}}

\newcommand{\rC}{{\rm C}}

\newcommand{\rE}{{\rm E}}

\newcommand{\rH}{{\rm H}}
\newcommand{\rI}{{\rm I}}
\newcommand{\rT}{{\rm T}}
\newcommand{\rK}{{\rm K}}
\newcommand{\rL}{{\rm L}}

\newcommand{\rR}{{\rm R}}
\newcommand{\rS}{{\rm S}}

\newcommand{\rW}{{\rm W}}

\newcommand{\ro}{{\rm o}}
\newcommand{\rp}{{\rm p}}
\newcommand{\rs}{{\rm s}}
\newcommand{\rv}{{\rm v}}

\newcommand{\UB}{\mathit{UB}}
\newcommand{\UP}{\mathit{UP}}

\newcommand{\oK}{{\overline{K}}}

\newcommand{\oR}{{\overline{R}}}
\newcommand{\oS}{{\overline{S}}}

\newcommand{\oU}{{\overline{U}}}

\newcommand{\oX}{{\overline{X}}}
\newcommand{\oY}{{\overline{Y}}}
\newcommand{\oZ}{{\overline{Z}}}
\newcommand{\oa}{{\overline{a}}}

\newcommand{\of}{{\overline{f}}}
\newcommand{\ogg}{{\overline{g}}}

\newcommand{\oi}{{\overline{i}}}

\newcommand{\ou}{{\overline{u}}}

\newcommand{\ox}{{\overline{x}}}
\newcommand{\oy}{{\overline{y}}}
\newcommand{\oz}{{\overline{z}}}

\newcommand{\ocM}{{\bar{\cM}}}

\newcommand{\oeta}{{\overline{\eta}}}

\newcommand{\opsi}{{\overline{\psi}}}

\newcommand{\hGamma}{\widehat{\Gamma}}
\newcommand{\ogamma}{{\overline{\gamma}}}

\newcommand{\ovarphi}{{\overline{\varphi}}}

\newcommand{\omu}{{\overline{\mu}}}

\newcommand{\ojmath}{{\overline{\jmath}}}
\newcommand{\onabla}{{\overline{\nabla}}}
\newcommand{\oupmu}{{\overline{\upmu}}}
\newcommand{\oupnu}{{\overline{\upnu}}}

\newcommand{\ocB}{{\overline{\cB}}}

\newcommand{\ua}{{\underline{a}}}

\newcommand{\uE}{{\underline{E}}}
\newcommand{\uF}{{\underline{F}}}
\newcommand{\uG}{{\underline{G}}}
\newcommand{\uH}{{\underline{H}}}

\newcommand{\ux}{{\underline{x}}}
\newcommand{\uX}{{\underline{X}}}

\newcommand{\ue}{{\underline{e}}}

\newcommand{\ug}{{\underline{g}}}
\newcommand{\uh}{{\underline{h}}}

\newcommand{\un}{{\underline{n}}}

\newcommand{\ur}{{\underline{r}}}
\newcommand{\ut}{{\underline{t}}}
\newcommand{\uu}{{\underline{u}}}
\newcommand{\uy}{{\underline{y}}}

\newcommand{\uoX}{{\underline{\oX}}}
\newcommand{\ucD}{{\underline{\cD}}}
\newcommand{\ucE}{{\underline{\cE}}}

\newcommand{\urT}{{\underline{\rT}}}
\newcommand{\utheta}{{\underline{\theta}}}
\newcommand{\ualpha}{{\underline{\alpha}}}
\newcommand{\ubeta}{{\underline{\beta}}}

\newcommand{\udelta}{{\underline{\delta}}}

\newcommand{\upi}{{\underline{\pi}}}
\newcommand{\uPi}{{\underline{\Pi}}}

\newcommand{\ukappa}{{\underline{\kappa}}}
\newcommand{\ulambda}{{\underline{\lambda}}}
\newcommand{\uuplambda}{{\underline{\uplambda}}}
\newcommand{\uupgamma}{{\underline{\upgamma}}}
\newcommand{\uuptheta}{{\underline{\uptheta}}}

\newcommand{\uOmega}{{\underline{\Omega}}}
\newcommand{\uTheta}{{\underline{\Theta}}}

\newcommand{\uDelta}{{\underline{\Delta}}}

\newcommand{\umK}{{\underline{\mK}}}
\newcommand{\umi}{{\underline{\mi}}}
\newcommand{\utmK}{{\underline{\tmK}}}
\newcommand{\umM}{{\underline{\mM}}}

\newcommand{\uuptau}{{\underline{\uptau}}}
\newcommand{\uupphi}{{\underline{\upphi}}}

\newcommand{\uupiota}{{\underline{\upiota}}}
\newcommand{\uupnu}{{\underline{\upnu}}}
\newcommand{\uupmu}{{\underline{\upmu}}}

\newcommand{\ha}{{\widehat{a}}}

\newcommand{\hu}{{\widehat{u}}}

\newcommand{\hw}{{\widehat{w}}}

\newcommand{\hA}{{\widehat{A}}}

\newcommand{\hE}{{\widehat{E}}}

\newcommand{\hrS}{{\widehat{\rS}}}
\newcommand{\hcA}{{\widehat{\cA}}}
\newcommand{\hcF}{{\widehat{\cF}}}

\newcommand{\hbQ}{\widehat{\bQ}}

\newcommand{\hRun}{{\widehat{R_1}}}
\newcommand{\hRtau}{{\widehat{R_\uptau}}}
\newcommand{\hRunp}{{\widehat{R'_1}}}
\newcommand{\hRtaup}{{\widehat{R'_\uptau}}}

\newcommand{\hOmega}{\widehat{\Omega}}
\newcommand{\huOmega}{\widehat{\uOmega}}
\newcommand{\hvarphi}{\widehat{\varphi}}
\newcommand{\hmZ}{{\widehat{\mZ}}}

\newcommand{\halpha}{\widehat{\alpha}}

\newcommand{\hdelta}{\widehat{\delta}}

\newcommand{\htheta}{\widehat{\theta}}

\newcommand{\hupphi}{\widehat{\upphi}}

\newcommand{\hpsi}{\widehat{\psi}}

\newcommand{\hupsigma}{\widehat{\upsigma}}
\newcommand{\hupmu}{\widehat{\upmu}}

\newcommand{\huppi}{\widehat{\uppi}}
\newcommand{\hvarrho}{\widehat{\varrho}}

\newcommand{\hmj}{\widehat{\mj}}

\newcommand{\vupsigma}{\vec{\upsigma}}

\newcommand{\cA}{{\mathscr A}}
\newcommand{\cB}{{\mathscr B}}
\newcommand{\cC}{{\mathscr C}}
\newcommand{\cD}{{\mathscr D}}
\newcommand{\cE}{{\mathscr E}}
\newcommand{\cF}{{\mathscr F}}
\newcommand{\cG}{{\mathscr G}}
\newcommand{\cI}{{\mathscr I}}
\newcommand{\cJ}{{\mathscr J}}
\newcommand{\cK}{{\mathscr K}}
\newcommand{\cL}{{\mathscr L}}
\newcommand{\cP}{{\mathscr P}}
\newcommand{\co}{{\mathscr O}}
\newcommand{\cR}{{\mathscr R}}
\newcommand{\cS}{{\mathscr S}}
\newcommand{\cT}{{\mathscr T}}
\newcommand{\cH}{{\mathscr H}}
\newcommand{\cM}{{\mathscr M}}
\newcommand{\cN}{{\mathscr N}}
\newcommand{\cQ}{{\mathscr Q}}
\newcommand{\cU}{{\mathscr U}}
\newcommand{\cV}{{\mathscr V}}

\newcommand{\cZ}{{\mathscr Z}}

\newcommand{\cHom}{{\mathscr Hom}}

\newcommand{\cEnd}{{\mathscr End}}

\newcommand{\fC}{{\mathfrak C}}

\newcommand{\fF}{{\mathfrak F}}

\newcommand{\fN}{{\mathfrak N}}

\newcommand{\fS}{{\mathfrak S}}

\newcommand{\fV}{{\mathfrak V}}

\newcommand{\fX}{{\mathfrak X}}
\newcommand{\fY}{{\mathfrak Y}}
\newcommand{\fZ}{{\mathfrak Z}}
\newcommand{\fa}{{\mathfrak a}}
\newcommand{\fb}{{\mathfrak b}}
\newcommand{\fc}{{\mathfrak c}}
\newcommand{\fd}{{\mathfrak d}}
\newcommand{\fgg}{{\mathfrak g}}
\newcommand{\fh}{{\mathfrak h}}
\newcommand{\fm}{{\mathfrak m}}

\newcommand{\fp}{{\mathfrak p}}

\newcommand{\tta}{{\tt a}}
\newcommand{\ttc}{{\tt c}}

\newcommand{\ttb}{{\tt b}}

\newcommand{\ttg}{{\tt g}}
\newcommand{\tth}{{\tt h}}
\newcommand{\tti}{{\tt i}}
\newcommand{\ttt}{{\tt t}}
\newcommand{\tts}{{\tt s}}

\newcommand{\ttu}{{\tt u}}
\newcommand{\ttv}{{\tt v}}

\newcommand{\hoR}{{\widehat{\oR}}}
\newcommand{\hoRp}{{\widehat{\oR'}}}

\newcommand{\hcC}{{\widehat{\cC}}}
\newcommand{\hcD}{{\widehat{\cD}}}
\newcommand{\hcI}{{\widehat{\cI}}}
\newcommand{\hcJ}{{\widehat{\cJ}}}
\newcommand{\hcG}{{\widehat{\cG}}}
\newcommand{\hcH}{{\widehat{\cH}}}

\newcommand{\hfC}{{\widehat{\fC}}}

\newcommand{\hmX}{{\widehat{\mX}}}
\newcommand{\hmY}{{\widehat{\mY}}}

\newcommand{\hmu}{{\widehat{\mu}}}

\newcommand{\hotimes}{{\widehat{\otimes}}}

\newcommand{\httc}{{\widehat{\ttc}}}

\newcommand{\hupiota}{{\widehat{\upiota}}}

\newcommand{\bvf}{{\breve{f}}}
\newcommand{\bvu}{{\breve{u}}}

\newcommand{\bvoS}{{\breve{\oS}}}

\newcommand{\bvoX}{{\breve{\oX}}}
\newcommand{\bvocB}{{\breve{\ocB}}}

\newcommand{\bvco}{{\breve{\co}}}
\newcommand{\bvcC}{{{\breve{\cC}}}}

\newcommand{\bvcF}{{\breve{\cF}}}
\newcommand{\bvcH}{{\breve{\cH}}}

\newcommand{\bvcM}{{\breve{\cM}}}

\newcommand{\bvfX}{{\breve{\fX}}}

\newcommand{\bvpi}{{\breve{\pi}}}
\newcommand{\bvvarpi}{{\breve{\varpi}}}

\newcommand{\bvalpha}{{\breve{\alpha}}}

\newcommand{\bvsigma}{{\breve{\sigma}}}
\newcommand{\bvSigma}{{\breve{\Sigma}}}
\newcommand{\bvvarphi}{{\breve{\varphi}}}
\newcommand{\bvdelta}{{\breve{\delta}}}

\newcommand{\bvOmega}{{\breve{\Omega}}}
\newcommand{\bvupiota}{{\breve{\upiota}}}

\newcommand{\bvupphi}{{\breve{\upphi}}}
\newcommand{\bvuppsi}{{\breve{\uppsi}}}

\newcommand{\bvuptheta}{{\breve{\uptheta}}}

\newcommand{\bvtta}{{\breve{\tta}}}
\newcommand{\bvogg}{{\breve{\ogg}}}
\newcommand{\bvlgg}{{\breve{\lgg}}}

\newcommand{\ccT}{{\check{T}}}
\newcommand{\cun}{{\check{1}}}

\newcommand{\coS}{{\check{\oS}}}
\newcommand{\coU}{{\check{\oU}}}
\newcommand{\coX}{{\check{\oX}}}
\newcommand{\coY}{{\check{\oY}}}
\newcommand{\cog}{{\check{\ogg}}}

\newcommand{\cof}{{\check{\of}}}

\newcommand{\cogamma}{{\check{\ogamma}}}

\newcommand{\cphi}{{\check{\phi}}}

\newcommand{\cmu}{{\check{\mu}}}
\newcommand{\ccE}{{\check{\cE}}}
\newcommand{\ccF}{{\check{\cF}}}
\newcommand{\ctta}{{\check{\tta}}}

\newcommand{\tB}{{\widetilde{B}}}

\newcommand{\tE}{{\widetilde{E}}}

\newcommand{\tG}{{\widetilde{G}}}
\newcommand{\tuE}{{\widetilde{\uE}}}
\newcommand{\tuG}{{\widetilde{\uG}}}

\newcommand{\tN}{{\widetilde{N}}}

\newcommand{\tS}{{\widetilde{S}}}

\newcommand{\tU}{{\widetilde{U}}}
\newcommand{\tV}{{\widetilde{V}}}

\newcommand{\tX}{{\widetilde{X}}}
\newcommand{\tY}{{\widetilde{Y}}}

\newcommand{\td}{{\widetilde{d}}}

\newcommand{\tf}{{\widetilde{f}}}
\newcommand{\tg}{{\widetilde{g}}}
\newcommand{\thh}{{\widetilde{h}}}

\newcommand{\tu}{{\widetilde{u}}}
\newcommand{\tv}{{\widetilde{v}}}

\newcommand{\tcB}{{\widetilde{\cB}}}
\newcommand{\tcE}{{\widetilde{\cE}}}
\newcommand{\tcD}{{\widetilde{\cD}}}
\newcommand{\tcH}{{\widetilde{\cH}}}
\newcommand{\tcK}{{\widetilde{\cK}}}

\newcommand{\trT}{{\widetilde{\rT}}}

\newcommand{\tOmega}{{\widetilde{\Omega}}}

\newcommand{\tgamma}{{\widetilde{\gamma}}}

\newcommand{\ttheta}{{\widetilde{\theta}}}
\newcommand{\tnabla}{{\widetilde{\nabla}}}

\newcommand{\tmu}{{\widetilde{\mu}}}

\newcommand{\tpi}{{\widetilde{\pi}}}

\newcommand{\tpsi}{{\widetilde{\psi}}}
\newcommand{\tphi}{{\widetilde{\phi}}}

\newcommand{\txi}{{\widetilde{\xi}}}
\newcommand{\tvartheta}{{\widetilde{\vartheta}}}

\newcommand{\talpha}{{\widetilde{\alpha}}}
\newcommand{\tpartial}{{\widetilde{\partial}}}

\newcommand{\tupiota}{{\tilde{\upiota}}}
\newcommand{\tuupiota}{{\tilde{\uupiota}}}
\newcommand{\tcC}{{\widetilde{\cC}}}

\newcommand{\tcL}{{\widetilde{\cL}}}

\newcommand{\tmK}{{\widetilde{\mK}}}
\newcommand{\tmX}{{\widetilde{\mX}}}
\newcommand{\tmY}{{\widetilde{\mY}}}
\newcommand{\tmZ}{{\widetilde{\mZ}}}

\newcommand{\tmj}{{\widetilde{\mj}}}

\begin{document}

\title{Twisting Higgs Modules and Functorial Aspects of the $p$-adic Simpson Correspondence}
\author{Ahmed Abbes, Michel Gros and Takeshi Tsuji}
\address{A.A. Laboratoire Alexander Grothendieck, UMR 9009 du CNRS, 
Institut des Hautes \'Etudes Scientifiques, 35 route de Chartres, 91440 Bures-sur-Yvette, France}
\address{M.G. Université de Rennes, CNRS, IRMAR - UMR 6625, Campus de Beaulieu, F-35042 Rennes cedex, France}
\address{T.T. Graduate School of Mathematical Sciences, The University of Tokyo, 3-8-1, Komaba, Meguro, Tokyo 153-8914, Japan}
\email{abbes@ihes.fr}
\email{michel.gros@univ-rennes1.fr}
\email{t-tsuji@ms.u-tokyo.ac.jp}

\maketitle

\setcounter{tocdepth}{1}
\tableofcontents


\chapter*{Foreword}

The $p$-adic Simpson correspondence was introduced by Gerd Faltings in his seminal 2005 article \cite{faltings3}. 
It draws on two distinct sources: a complex-geometric origin, which provides both its name and its guiding formalism, 
and a $p$-adic origin, namely {\em Sen theory}, which shares its central objective, the description of $p$-adic representations 
of the fundamental group of an algebraic variety over a $p$-adic field, as well as the techniques of $p$-adic Hodge theory.
Since then, this program has revealed a wide array of facets and incarnations (geometric Sen theory, analytic variants, etc.), 
as well as several complementary approaches (Fontaine-type correspondences via various period rings, crystalline, prismatic, among others). 
Its role within $p$-adic Hodge theory continues to stimulate many significant developments, including contributions 
of the present authors \cite{agt,ag2,tsuji5} and, without any claim to completeness, recent works by 
Bhargav Bhatt and Mingjia Zhang, Tongmu He \cite{he}, Ben Heuer \cite{heuer}, Ben Heuer and Daxin Xu \cite{hexu}, Lue Pan \cite{pan}, 
Juan Esteban Rodríguez Camargo \cite{rod}, Daxin Xu \cite{xu}, among others.

The main goal of this book is to establish a robust framework for studying the functoriality of the $p$-adic Simpson correspondence. 
We introduce a new method for {\em twisting Higgs modules via Higgs--Tate algebras}. 
This construction builds on one of our earlier approaches to the $p$-adic Simpson correspondence \cite{agt}, which it recovers as a special case.
The resulting framework yields {\em twisted pullbacks and higher direct images of Higgs modules}, thereby enabling a systematic study of 
the functoriality of the $p$-adic Simpson correspondence under arbitrary pullbacks and proper (log)smooth direct images, 
including for morphisms that do not admit liftings to the infinitesimal deformations used in the construction of the correspondence.

The core of our approach is the introduction of a {\em twisting functor} for Higgs modules, developed in the general setting of a ringed topos 
$(X,\co_X)$, where $\co_X$ is $p$-adically complete and separated, and flat over a complete valuation ring of height $1$, 
with algebraically closed fraction field of characteristic $0$ and residue field of characteristic $p>0$.
Starting from an exact sequence of locally free $\co_X$-modules of finite type
\[
0\rightarrow \co_X\rightarrow \cF\rightarrow \cE\rightarrow 0,
\]
we associate to it a {\em Higgs--Tate $\co_X$-algebra}, defined as the direct limit
\[
\cC=\underset{\underset{n\geq 0}{\longrightarrow}}\lim\ \Sym^n_{\co_X}(\cF),
\]
which may be viewed as a twisted form of the symmetric algebra of $\cE$. 
Inspired by our previous construction of the $p$-adic Simpson correspondence \cite{ag2, agt}, we define an endofunctor
\[
\uptau\colon \bHM(\co_X,\cE) \rightarrow \bHM(\co_X,\cE),
\]
on the category of Higgs $\co_X$-modules with coefficients in $\cE$, 
together with an admissibility condition defining the class of {\em (weakly) twistable Higgs modules}. 
The study and applications of this functor rely on a Poincaré lemma; accordingly, 
the construction is carried out using a weak $p$-adic completion of $\cC$.  

\pagebreak 

This twisting functor admits several fundamental applications.

First, it provides a systematic construction of twisted pullbacks and higher direct images for Higgs modules. 
The twisted pullback for Higgs modules was introduced by Gerd Faltings in \cite{faltings3} through a different construction, 
and plays a central role in his treatment of the $p$-adic Simpson correspondence for non-small generalized representations via a descent argument. 
Twisted higher direct images for Higgs modules, on the other hand, do not appear in Faltings' work, nor elsewhere in the literature, 
and, to the best of our knowledge, are introduced here for the first time. 
They generalize, on relative Dolbeault cohomology, the Higgs-field analogue of the Katz--Oda construction 
of the Gauss--Manin connection on relative de Rham cohomology.

Second, the twisting functor yields a new construction of the $p$-adic Simpson correspondence for small Higgs bundles, 
encompassing and refining our earlier work \cite{ag2, agt}.

Finally, it provides the key tool for establishing functoriality properties of the $p$-adic Simpson correspondence, 
improving our earlier results obtained under strong liftability assumptions \cite{ag2}. 
These properties are derived from compatibility results for twisting functors, 
notably the identification of the composition of two such functors with a single one 
via a homomorphism between the corresponding Higgs--Tate algebras.

The extensions underlying the construction arise naturally from deformation theory. More precisely, given a commutative diagram of morphisms of 
(logarithmic) schemes
\[
\xymatrix{
Y \ar[r]^j \ar[d]^h \ar@/_1pc/[dd]_g & \tY \ar@/^1pc/[dd]^{\tg} \\
X \ar[r]^i \ar[d]^f & \tX \ar[d]_{\tf} \\
B \ar[r]^{\iota} & \tB,}
\]
where $\iota$ is a first-order thickening, $\tf$ is smooth, and the lower square and the exterior rectangle are Cartesian,
one considers the torsor of local liftings of $h$ over $\tX$. 
This torsor is affine under a vector bundle naturally expressed in terms of relative differentials, 
and its module of affine functions is an extension of the required form, from which the Higgs--Tate algebra, 
and hence the twisting functor, are constructed.

In this way, each twisting functor is ultimately governed by first-order thickenings. 
In the case of the $p$-adic Simpson correspondence, the morphisms $\iota$ and $j$ are given by 
{\em Fontaine's universal $p$-adic first-order thickenings}, while $i$ is an auxiliary choice on which the correspondence depends.
This dependence makes the functoriality properties of the $p$-adic Simpson correspondence more delicate and involved.

This book is organized as follows. Chapter~\ref{overview} contains a detailed overview of our results, together with a digression illustrating 
the relevance of our twisting approach in the setting of classical Sen theory \eqref{rstvt}. 

In Chapter~\ref{preliminaries}, we gather a number of technical preliminaries that are essential for the sequel. 
These go well beyond the $p$-adic framework, and most are original, including the theory of {\em quasi-fibered categories} \eqref{p2-qfc}, 
which is of independent interest. 

Chapter~\ref{twisting} is devoted to the central notion of this work, namely the twisting of Higgs modules, 
and introduces the corresponding twisted pullbacks and twisted higher direct images of Higgs modules. 

Finally, Chapter~\ref{functoriality} contains our main results: the new construction of the $p$-adic Simpson correspondence 
by twisting (for {\em small} objects) and its compatibility with twisted pullbacks and twisted proper (log-)smooth higher direct images.

\vspace{5mm}
   
\hfill Ahmed Abbes, Michel Gros, and Takeshi Tsuji

\vspace{2mm}

\hfill May 2026

\chapter{An overview}\label{overview}

\section{Introduction}\label{overview-intro}

\subsection{}\label{overview-intro1}
Initiated by Faltings \cite{faltings3} and subsequently developed through various approaches, including two due to the authors \cite{agt}, 
the {\em $p$-adic Simpson correspondence} establishes an equivalence of categories between certain 
{\em $p$-adic étale local systems} on an algebraic variety over a $p$-adic field and certain {\em Higgs bundles}.
This theory has two distinct origins: one in complex geometry and the other in the $p$-adic Hodge theory.
From the complex Simpson correspondence it inherits its name and guiding formalism, while the $p$-adic origin supplies the techniques and the main tools.

\subsection{}\label{overview-intro2}
With some partiality, yet in good conscience, one may trace the immediate origin of the $p$-adic Simpson correspondence to Sen theory \cite{sen1}.
Let $K$ be a complete discrete valuation field of characteristic $0$, with perfect residue field $k$ of characteristic $p>0$,
$\co_K$ the valuation ring of $K$, $\oK$ an algebraic closure of $K$, $\co_\oK$ the integral closure of $\co_K$ in $\oK$,
$G_K$ the Galois group of $\oK$ over $K$, $\co_C$ the $p$-adic Hausdorff completion of $\co_\oK$, $\fm_C$ its maximal ideal,
$C$ its field of fractions. 

Let $\chi\colon G_K\rightarrow \mZ_p^\times$ be the cyclotomic character of $K$, describing the action of $G_K$ on the $p$-power roots of unity in $\oK$.
Let $\log\colon \mZ_p^\times \rightarrow \mZ_p$ denote the usual $p$-adic logarithm, defined for all $x\in 1+p\mZ_p$ by  
$\log(x)=\sum_{n\ge1} \frac{(-1)^{n+1}}{n}(x-1)^n$ and in general by $\log(x)=\frac{\log(x^{p-1})}{p-1}$.
The map $\log \chi\colon G_K\rightarrow \mZ_p$ being continuous, let $H$ be its kernel and $\Gamma$ the quotient $G_K/H$.
Then $K_\infty=\oK^H$ is the cyclotomic $\mZ_p$-extension of $K$ contained in $\oK$, {\em i.e.}, the unique $\mZ_p$-extension of $K$ contained 
in the subfield of $\oK$ generated by the $p$-power order roots of the unity. We denote by $L$ the closure of $K_\infty$ in $C$. 

A $C$-representation of $G_K$ is a $C$-vector space equipped with a semilinear action of $G_K$.
We denote by $\bRep_C(G_K)$ the category of finite-dimensional continuous $C$-representations of $G_K$. 
Similarly, we define the categories $\bRep_{K_\infty}(\Gamma)$ and $\bRep_L(\Gamma)$.

\begin{teo}[Sen, \cite{sen1} Theorems 2 and 3]\label{overview-intro3} 
The functors
\begin{eqnarray}
\bRep_{K_\infty}(\Gamma)\rightarrow \bRep_L(\Gamma), && V\mapsto V\otimes_{K_\infty}L,\\
\bRep_{L}(\Gamma)\rightarrow \bRep_C(G_K), && W\mapsto W\otimes_{L}C,
\end{eqnarray}
are equivalences of categories.
\end{teo}

The functor that associates to a $C$-representation $W$ of $G_K$ the $L$-representation $W^{H}$ of 
$\Gamma$ is a quasi-inverse of the second functor. 
Sen constructed an explicit quasi-inverse of the first functor (\cite{sen1} Theorem 3). 

\begin{teo}[Sen, \cite{sen1} Theorem 4]\label{overview-intro4}
For every finite-dimensional continuous $K_\infty$-repre\-sentation $V$ of $\Gamma$, 
there exists a unique $K_\infty$-endomorphism $\sigma$ of $V$ satisfying the following property: 
for every $x\in V$, there exists an open subgroup $\Gamma_x$ of $\Gamma$ such that for every $g\in \Gamma_x$,
we have 
\begin{equation}\label{overview-intro4a}
g(x)=\exp(\log(\chi(g)) \sigma)(x).
\end{equation}
Moreover, there exists a $K_\infty$-basis of $V$ with respect to which the matrix of $\sigma$ has coefficients in $K$. 
\end{teo}

\subsection{}\label{overview-intro40}
For any field extension $F$ of $K$, 
a {\em Sen $F$-module} is a pair $(V,\sigma)$ consisting of a finite dimensional $F$-vector space $V$ equipped with  
an $F$-linear endomorphism $\sigma\colon V\rightarrow V$.
If $(V,\sigma)$ and $(V',\sigma')$ are two Sen $F$-modules, a morphism from $(V,\sigma)$ to $(V',\sigma')$ is an $F$-linear map 
$f\colon V\rightarrow V'$ such that $f\circ \sigma=\sigma'\circ f$.
Sen $F$-modules form a category, denoted by $\bSM_{F}$, which is canonically equivalent to the category of $F[X]$-modules
having finite dimension over $F$.

Theorems \ref{overview-intro3} and \ref{overview-intro4} provide an {\em exact and faithful} functor 
\begin{equation}\label{overview-intro40a}
\bRep_C(G_K) \rightarrow \bSM_{K_\infty}, \ \ \ W\mapsto (V,\sigma).
\end{equation}
Moreover, $(V,\sigma)$ determines $W$ up to an isomorphism (\cite{sen1} Theorem 7).  

By (\cite{sen1} Theorem 10), if the residue field $k$ of $\co_K$ is algebraically closed, the functor \eqref{overview-intro40a} descends into 
a non-canonical equivalence of categories  
\begin{equation}\label{overview-intro40b}
\bRep_C(G_K) \stackrel{\sim}{\rightarrow} \bSM_K. 
\end{equation}

\subsection{}\label{overview-intro5}
The $p$-adic Simpson correspondence may be viewed as an analogue of the functor \eqref{overview-intro40a}. 
It describes certain $p$-adic representations of the geometric fundamental group of an algebraic variety over $K$, 
admitting a suitable integral model, in terms of linear algebra objects, namely {\em Higgs modules}.

\subsection{}\label{overview-intro6}
In this book, we introduce a new method for {\em twisting Higgs modules} using {\em Higgs--Tate algebras}.
This construction is inspired by one of our earlier approaches to the $p$-adic Simpson correspondence, which it recovers as a special case. 
The resulting framework provides {\em twisted pullbacks} and {\em twisted higher direct images} of Higgs modules, 
allowing us to study the functoriality of the $p$-adic Simpson correspondence and 
to generalize the results of \cite{ag2}, which required a strong liftability condition on the morphism.
More precisely, we establish the functoriality of the $p$-adic Simpson correspondence under arbitrary pullbacks and 
under proper (log-)smooth direct images by morphisms that do not necessarily lift to the infinitesimal deformations of the varieties 
used in the construction of the correspondence.
We also clarify how this new twisting relates to the constructions of Heuer \cite{heuer} and Heuer--Xu \cite{hexu}, 
involving line bundles on the spectral variety \eqref{p2-rgpsc39}.

In this first chapter, we give an overview of our twisting approach to the $p$-adic Simpson correspondence, 
which we illustrate in \ref{rstvt} by providing a quasi-inverse of Sen's equivalence of categories \eqref{overview-intro40b}.

{\em The book treats the case of $\co_K$-schemes with toric singularities using logarithmic geometry; 
for simplicity, we confine this overview to the smooth case.} 

\subsection*{Acknowledgments} 
This work was supported by JSPS Grant-in-Aid for Scientific Research (Grant Numbers 20H01793 and 24K00518) 
and by a JSPS Invitational Fellowship for Research in Japan (S24045), 
which supported the first author's visit to the Graduate School of Mathematical Sciences, the University of Tokyo, in spring 2024.

\section{Twisting Higgs modules}\label{thm}

\subsection{}\label{thm0}
In this section, $V$ denotes a complete valuation ring of height $1$, with algebraically closed fraction field of characteristic $0$ 
and residue field of characteristic $p>0$. We choose a compatible system $(\beta_n)_{n>0}$ of $n$th roots of $p$ in $V$.
For any rational number $\varepsilon>0$, we set $p^\varepsilon=(\beta_n)^{\varepsilon n}$ where $n$ is an integer $>0$ such that $\varepsilon n$ is an integer.
We equip $V$ with the $p$-adic topology and set $\cS=\Spf(V)$. 

\subsection{}\label{thm1}
Let $X$ be a topos, $\co_X$ a flat $V$-algebra that is $p$-adically complete and separated, 
\begin{equation}\label{thm1a}
0\rightarrow \co_X\rightarrow \cF\rightarrow \cE \rightarrow 0
\end{equation}
an exact sequence of locally free $\co_X$-modules of finite type. 
We associate with the extension $\cF$ its {\em Higgs--Tate $\co_X$-algebra} $\cC$, defined as the direct limit 
\begin{equation}\label{thm1b}
\cC=\underset{\underset{n\geq 0}{\longrightarrow}}\lim\ \Sym^n_{\co_X}(\cF), 
\end{equation}
with transition maps $x_1\otimes\cdots \otimes x_n \mapsto 1\otimes x_1\otimes\cdots\otimes x_n$ \eqref{p1-prem1}. 
We denote by $\bHM(\co_X,\cE)$ the category of Higgs $\co_X$-modules with coefficients in $\cE$ \eqref{p1-delta-con1}. 
Inspired by our previous construction of the $p$-adic Simpson correspondence (\cite{ag2} 3.3.16, \cite{agt} II.12.15),
we associate to the extension $\cF$ a {\em twisting functor}
\begin{equation}\label{thm1c}
\uptau\colon \bHM(\co_X,\cE) \rightarrow \bHM(\co_X,\cE),
\end{equation}
and an admissibility condition, which we call {\em (weakly) twistable Higgs modules}.
The study and applications of this functor rely on a Poincaré lemma; accordingly, the construction is carried out using a weak 
$p$-adic completion of $\cC$.

\subsection{}\label{thm2}
For any rational number $r\geq 0$, we denote by $\cF^{(r)}$ the extension of $\co_X$-modules deduced from $\cF$ by pullback by the multiplication 
by $p^r$ on $\cE$; so we have an exact sequence of $\co_X$-modules
\begin{equation}\label{thm2a}
0\rightarrow \co_X\rightarrow \cF^{(r)} \rightarrow \cE \rightarrow 0. 
\end{equation} 
We denote by $\cC^{(r)}$ the associated Higgs--Tate algebra.
We also call it the {\em Higgs--Tate algebra of thickness $r$} associated with the extension $\cF$. 
Its universal $\co_X$-derivation is given by 
\begin{equation}\label{thm2b}
d_{\cC^{(r)}}\colon \cC^{(r)} \rightarrow \cE\otimes_{\co_X} \cC^{(r)}. 
\end{equation}
The algebras $\cC^{(r)}$, for $r\in \mQ_{\geq 0}$, form naturally an inductive system, and the derivations $\delta_{\cC^{(r)}}=p^rd_{\cC^{(r)}}$ are compatible.
We set 
\begin{equation}\label{thm2c}
\cC^\dagger=\underset{\underset{r\in \mQ_{>0}}{\longrightarrow}}{\lim}\ \hcC^{(r)},
\end{equation}
where $\hcC^{(r)}$ is the $p$-adic Hausdorff completion of $\cC^{(r)}$, and denote by 
\begin{equation}\label{thm2d}
\delta \colon \cC^\dagger\rightarrow \cE\otimes_{\co_X}\cC^\dagger,
\end{equation}
the $\co_X$-derivation induced by $\delta_{\cC^{(r)}}$.

\subsection{}\label{thm3}
For any Higgs $\co_X$-module $(N,\theta)$ with coefficients in $\cE$, we denote by $\uuptau(N,\theta)$ the $\co_X$-module 
\begin{equation}\label{thm3a}
\uuptau(N,\theta)=(N\otimes_{\co_X}\cC^\dagger)^{\theta_\tot=0},
\end{equation}
where $\theta_\tot=\theta\otimes \id+\id\otimes \delta$ is the total Higgs $\co_X$-field on $ N\otimes_{\co_X}\cC^\dagger$. 
The Higgs field $\theta\otimes \id$ on $N\otimes_{\co_X}\cC^\dagger$ induces a Higgs $\co_X$-field $\theta_\uptau$ on 
$\uuptau(N,\theta)$ with coefficients in $\cE$. 
We thus define a functor
\begin{equation}\label{thm3b}
\uptau\colon 
\begin{array}[t]{clcr}
\bHM(\co_X,\cE) &\rightarrow& \bHM(\co_X,\cE),\\
(N,\theta)&\mapsto&(\uuptau(N,\theta), \theta_\uptau),
\end{array}
\end{equation}
which we call the {\em twisting functor relative to the extension $\cF$ \eqref{thm1a}}; see \ref{p1-thbn9}.

We also need a ``dual'' version, namely the functor 
\begin{equation}\label{thm3c}
\uptau^\vee\colon 
\begin{array}[t]{clcr}
\bHM(\co_X,\cE) &\rightarrow& \bHM(\co_X,\cE)\\
(N,\theta)&\mapsto&(\uuptau^\vee(N,\theta), \theta_{\uptau^\vee}),
\end{array}
\end{equation}
defined similarly by replacing $\delta$ by $\delta^\vee=-\delta$.

\begin{defi}[cf. \ref{p1-thbn30}] \label{thm4}
We say that a Higgs $\co_X[\frac 1 p]$-module $(N,\theta)$ with coefficients in $\cE$ is 
{\em weakly twistable by the extension $\cF$} if the canonical $\cC^\dagger$-linear morphism
\begin{equation}\label{thm4a}
\cC^\dagger \otimes_{\co_X}\uuptau(N,\theta)\rightarrow \cC^\dagger \otimes_{\co_X}N
\end{equation}
is an isomorphism.  
\end{defi}

We reserve the adjective {\em twistable} for a slightly stronger admissibility condition \eqref{p1-thbn14}. 
The two conditions are equivalent if the $\co_X[\frac{1}{p}]$-module $N$ is flat of finite type and the topos $X$ is coherent \eqref{p1-thbn21}.

If the extension $\cF$ splits, the choice of a splitting determines for every weakly twistable Higgs $\co_X[\frac 1 p]$-module $(N,\theta)$,
a functorial isomorphism \eqref{p1-thbn31}
\begin{equation}\label{thm4b}
\uptau(N,\theta)\stackrel{\sim}{\rightarrow}(N,\theta).
\end{equation}

\begin{prop}[cf. \ref{p1-thbn41}]\label{thm10}
Let $(N,\theta)$ be a Higgs $\co_X[\frac 1 p]$-module with coefficients in $\cE$, $N^\vee$ an $\co_X$-module, 
\begin{equation}\label{thm10a}
\psi\colon N^\vee\otimes_{\co_X}\cC^\dagger\stackrel{\sim}{\rightarrow}N\otimes_{\co_X}\cC^\dagger
\end{equation}
an isomorphism of $\cC^\dagger$-modules with $\delta$-connection \eqref{p1-delta-con2}, 
where the $\delta$-connections are defined by the total Higgs fields, 
$N$ (resp.\ $N^\vee$) being endowed with the Higgs field $\theta$ (resp.\ $0$); see \ref{p1-delta-con4}. Then, 
\begin{itemize}
\item[{\rm (i)}] The isomorphism $\psi$ induces an $\co_X$-linear isomorphism $N^\vee\stackrel{\sim}{\rightarrow} \uuptau(N,\theta)$, 
where $\uuptau$ is the functor \eqref{thm3a}. In particular, $(N,\theta)$ is weakly twistable. 
We deduce a canonical Higgs $\co_X$-field $\theta^\vee$ on  $N^\vee$ with coefficients in $\cE$, so that we have an isomorphism of Higgs modules
\begin{equation}\label{thm10b}
(N^\vee,\theta^\vee)\stackrel{\sim}{\rightarrow} \uptau(N,\theta),
\end{equation}
where $\uptau$ is the functor \eqref{thm3b}. 
\item[{\rm (ii)}] The morphism $\psi$ is an isomorphism of $\cC^\dagger$-modules with $\delta^\vee$-connection, 
where the $\delta^\vee$-connections are defined by the total Higgs fields, 
$N$ (resp.\ $N^\vee$) being endowed with the Higgs field $0$ (resp.\ $\theta^\vee$); see \ref{p1-delta-con4}. 
\item[{\rm (iii)}] The isomorphism $\psi$ induces an isomorphism of Higgs $\co_X$-modules with coefficients in $\cE$, 
\begin{equation}\label{thm10c}
(N,\theta)\stackrel{\sim}{\rightarrow} \uptau^\vee(N^\vee,\theta^\vee),
\end{equation}
where $\uptau^\vee$ is the functor \eqref{thm3c}. 
\end{itemize}
\end{prop}

\begin{prop}[cf. \ref{p1-thbn44}]\label{thm5}
Let $(N,\theta)$ be a Higgs $\co_X[\frac 1 p]$-module with coefficients in $\cE$.  We set $\rS=\rS_{\co_X}(\cE^\vee)$
and denote by $\upmu\colon \rS[\frac 1 p]\rightarrow \cEnd_{\co_X}(N)$ 
the homomorphism of $\co_X[\frac 1 p]$-algebras defined by $\theta$ \eqref{p1-delta-con1j}. 
Let $B$ be a quotient $\co_X[\frac 1 p]$-algebra of $\rS[\frac 1 p]$ through which $\upmu$ factors, 
$\theta_B\colon B\rightarrow B\otimes_{\co_X}\cE$ its canonical Higgs $\co_X$-field. 
Suppose that the Higgs $\co_X[\frac 1 p]$-modules $(N,\theta)$ and $(B,\theta_B)$ are weakly twistable. 
Then, the $B$-module $\cL_B=\uptau(B,\theta_B)$ is invertible, 
and we have a canonical $B$-linear isomorphism 
\begin{equation}\label{thm5a}
N\otimes_B\cL_B\stackrel{\sim}{\rightarrow} \uptau(N,\theta).
\end{equation}
\end{prop}

We introduce in \ref{thm66} a natural condition, both stronger and easier to verify than the weak twistability conditions of the preceding proposition,
namely, the {\em CL-smallness}.

\subsection{}\label{thm6}
In the remainder of this section, let $\fX$ be a flat $\cS$-formal scheme locally of finite presentation, where $\cS = \Spf(V)$ \eqref{thm0}.
We assume that $X = \fX_{\zar}$, that $\co_X$ is an $\co_\fX$-algebra, and that $\cE = \Omega \otimes_{\co_\fX} \co_X$, 
for some locally free $\co_\fX$-module of finite type $\Omega$.

\begin{defi}[cf. \ref{p1-tshbn9}]\label{thm65}
Let $(\cN,\theta)$ be a Higgs $\co_\fX$-module with coefficients in $\Omega$, $\varepsilon\in \mQ_{> 0}$.
We say that $(\cN,\theta)$ is {\em $\varepsilon$-small} if 
\begin{equation}
\theta(\cN)\subset p^\varepsilon (\Omega\otimes_{\co_\fX} \cN).
\end{equation}
We say that $(\cN,\theta)$ is {\em small} if it is $\varepsilon'$-small for a rational number $\varepsilon'>\frac{1}{p-1}$. 
\end{defi}

\begin{defi}[cf. \ref{p1-tshbn13}]\label{thm66}
We say that a Higgs $\co_\fX[\frac 1 p]$-module $(N,\theta)$ with coefficients in $\Omega$, 
is {\em CL-small} if $N$ admits a coherent $\co_\fX$-lattice $\cN$, stable by $\theta$ such that the induced Higgs module $(\cN,\theta)$ is small.
We say that $(N,\theta)$ is {\em locally CL-small} if there exists a Zariski open covering $(U_i)_{i\in I}$ of $\fX$
such that for every $i\in I$, $(N|U_i,\theta|U_i)$ is CL-small.
\end{defi}

Observe that the rings $\co_\fX$ and $\co_\fX[\frac 1 p]$ are coherent \eqref{p1-pfs12}.

\begin{prop}[cf. \ref{p1-tshbn16}]\label{thm70}
For every locally CL-small Higgs $\co_\fX[\frac 1 p]$-module $(N,\theta)$ with coefficients in $\Omega$, 
the Higgs $\co_X[\frac 1 p]$-module $(N\otimes_{\co_\fX}\co_X,\theta\otimes \id)$ is weakly twistable. 
\end{prop}

\begin{prop}[cf. \ref{p1-tshbn27}]\label{thm7}
Let $N$ be a coherent $\co_\fX[\frac 1 p]$-module, $\theta$ a Higgs $\co_\fX$-field on $N$ with coefficients in $\Omega$, 
$\cH=\rS_{\co_\fX}(\Omega^\vee)$, $\upmu\colon \cH[\frac 1 p]\rightarrow \cEnd_{\co_\fX}(N)$ 
the homomorphism of $\co_\fX[\frac 1 p]$-algebras defined by $\theta$, $B$ its image and $\theta_B$ 
the canonical Higgs $\co_\fX$-field on $B$. 
Then, $(N,\theta)$ is CL-small if and only if $(B,\theta_B)$ is CL-small. 
\end{prop}

\subsection{}\label{thm8}
Let
\begin{equation}\label{thm8a}
0\rightarrow \cE^\vee\rightarrow \cF^\vee\rightarrow \co_X \rightarrow 0
\end{equation}
be the dual of the exact sequence \eqref{thm1a}. 
We denote by $\Lambda$ the inverse image of the section $1\in \co_X(X)$ in $\cF^\vee$.  
It is naturally equipped with a structure of a $\cE^\vee$-torsor. 

Let $\cH=\rS_{\co_\fX}(\Omega^\vee)$ and let $B$ be an $\co_\fX[\frac 1 p]$-algebra quotient of $\cH[\frac 1 p]$, 
$\theta_B$ the canonical Higgs $\co_\fX$-field on $B$.
We suppose that $(B,\theta_B)$ is locally CL-small. 
We set $(\tB,\theta_\tB)=(B,\theta)\otimes_{\co_\fX}\co_X$ and $\cL_\tB=\uptau(\tB,\theta_\tB)$, which is an invertible $\tB$-module. 
By the smallness, we define a morphism of sheaves of sets \eqref{p1-tshbn30e}
\begin{equation}\label{thm8b}
\exp \colon 
\begin{array}[t]{clcr}
\cE^\vee&\rightarrow& \tB^\times,\\
\phi&\mapsto& \exp(\phi).
\end{array}
\end{equation}
We prove that it is a group homomorphism \eqref{p1-tshbn33}. Moreover, we have the following. 

\begin{prop}[cf. \ref{p1-tshbn34}]\label{thm9}
The invertible $\tB$-module $\cL_\tB$ is canonically isomorphic to the line bundle associated with the 
$\tB^\times$-torsor $\Lambda\wedge^{\cE^\vee}\tB^\times$, 
deduced from $\Lambda$ by extension of its structural group by the homomorphism $\exp$ \eqref{thm8b}.
\end{prop}

We give in \ref{p1-tshbn41} a divided power refinement of this description of $\cL_\tB$: after replacing $\cF$ by $\cF^{(r)}$ for a suitable $r>0$, 
we construct a canonical invertible $\hGamma(\cE^\vee)$-module from $\cF^\vee$, 
where $\hGamma(\cE^\vee)$ is the completed divided power algebra of $\cE^\vee$ \eqref{p1-NC7b}, 
that induces $\cL_\tB$ by base change; see \ref{p1-tshbn35}.

\section{Torsors of liftings and Higgs--Tate algebras}\label{torlift}

\subsection{}\label{torlift1}
We will apply the twisting functor $\uptau$ \eqref{thm3b}, introduced in the previous section, 
relatively to extensions \eqref{thm1a} arising from deformation theory.
We first explain the construction of these extensions in an abstract setting of smooth schemes, 
before specializing to particular deformations that are geometrically or arithmetically meaningful. 
For simplicity, we restrict in this section to ordinary schemes; logarithmic schemes will be treated in §\ref{p1-rdt}.
Consider a commutative diagram of morphisms of schemes (without the dotted arrow)
\begin{equation}\label{torlift1a}
\xymatrix{
{Y}\ar[r]^-(0.5){j}\ar[d]^{h}\ar@/_1pc/[dd]_g&{\tY}\ar@{.>}[d]\ar@/^1pc/[dd]^\tg&\\
{X}\ar[r]^-(0.5){i}\ar[d]^f&{\tX}\ar[d]_{\tf}\\
{S}\ar[r]^-(0.5){\iota}&{\tS}}
\end{equation}
satisfying the following conditions: 
\begin{itemize}[label=--]
\item $\iota$ is a thickening of order one, defined by an invertible $\co_S$-module $I_S$; 
\item $\tf$ is smooth, and the two squares $(f,\tf,\iota,i)$ and $(g,\tg,\iota,j)$ are Cartesian; 
\item the ideal $I_Y$ defining $j$ is an invertible $\co_Y$-module; so we have $I_Y=g^*(I_S)$.    
\end{itemize}    

Let $\cL_{\tY/\tX}$ be the {\em torsor of local liftings of $h$ over $\tX$}, defined in the étale topos $Y_\et$ of $Y$, 
under the $\co_Y$-module $\cHom_{\co_Y}(h^*(\Omega^1_{X/S}),I_Y)$; see \ref{p1-rdt2}. 
Let $\cF_{\tY/\tX}$ be the $\co_Y$-module of {\em affine functions} on $\cL_{\tY/\tX}$ \eqref{p1-rdt3}.
It fits naturally into a canonical exact sequence of $\co_Y$-modules
\begin{equation}\label{torlift1b}
0\rightarrow \co_Y\rightarrow \cF_{\tY/\tX} \rightarrow \cE=h^*(I_S^{-1}\otimes_{\co_S}\Omega^1_{X/S})\rightarrow 0.
\end{equation}
We call it the {\em Higgs--Tate extension of $\cL_{\tY/\tX}$}. 
The class of the extension $I_Y\otimes_{\co_Y} \cF_{\tY/\tX}$ in $\Ext^1_{\co_Y}(Y;h^*(\tOmega^1_{X/S}),I_Y)$ is 
the obstruction to lift $g$ to $\tY$ over $\tX$. 
The torsor $\cL_{\tY/\tX}$ is representable by $\Spec(\cC_{\tY/\tX})$, where 
\begin{equation}\label{torlift1c}
\cC_{\tY/\tX}=\underset{\underset{n\geq 0}{\longrightarrow}}\lim\ \Sym^n_{\co_Y}(\cF_{\tY/\tX})
\end{equation} 
is the {\em Higgs--Tate algebra of $\cF_{\tY/\tX}$} (or of $\cL_{\tY/\tX}$).

\subsection{}\label{torlift2}
The formation of the Higgs--Tate algebra $\cC_{\tY/\tX}$ is functorial. 
Indeed, assume that diagram \eqref{torlift1a} can be enlarged into a commutative diagram 
\begin{equation}\label{torlift2a}
\xymatrix{
Y \ar[r]^j \ar[d]^{h'} \ar@/_3pc/[ddd]_g &
\tY \ar@/^3pc/[ddd]^\tg \\
X' \ar[r]^{i'} \ar[d]^{\gamma} \ar@/_1pc/[dd]_{f'} &
\tX' \ar@/^1pc/[dd]^{\tf'} \\
X \ar[r]^i \ar[d]^f &
\tX \ar[d]_{\tf} \\
S \ar[r]^\iota &\tS,}
\end{equation}
where $\tf'$ is smooth and the square $(f',\tf',\iota,i')$ is Cartesian.
By composing local liftings, we obtain an equivariant morphism of torsors 
\begin{equation}\label{torlift2b}
\nu\colon \cL_{\tY/\tX'}\times h'^+(\cL_{\tX'/\tX})\rightarrow\cL_{\tY/\tX},
\end{equation}
where $h'^+$ denotes the $h'$-affine pullback functor for torsors under quasi-coherent modules \eqref{p1-NC5}. 
The morphism $\nu$ induces a canonical homomorphism of $\co_Y$-algebras \eqref{p1-rdt6n}
\begin{equation}\label{torlift2c}
\upphi\colon \cC_{\tY/\tX}\rightarrow h'^*(\cC_{\tX'/\tX}) \otimes_{\co_Y}\cC_{\tY/\tX'},
\end{equation}
which we call the {\em composition homomorphism of Higgs--Tate algebras}. 
It is the key tool for studying of the functoriality of the $p$-adic Simpson correspondence, both under pullback and under proper higher direct images. 

\section{Twisted pullback and higher direct images of Higgs modules}\label{tphdi}

\subsection{}\label{tphdi1}
We take again the assumption and notation of \ref{thm0}. We set $S=\Spec(V)$ and consider a commutative diagram of morphisms of schemes
\begin{equation}\label{tphdi1a}
\xymatrix{
{X'}\ar[r]^-(0.5){i'}\ar[d]^{g}\ar@/_1pc/[dd]_{f'}&{\tX'}\ar@/^1pc/[dd]^{\tf'}&\\
{X}\ar[r]^-(0.5){i}\ar[d]^f&{\tX}\ar[d]_{\tf}\\
{S}\ar[r]^-(0.5){\iota}&{\tS}}
\end{equation}
satisfying the following conditions: 
\begin{itemize}[label=--]
\item $\iota$ is a thickening of order one, defined by an invertible $\co_S$-module $I_S$; 
\item the two squares $(f,\tf,\iota,i)$ and $(f',\tf',\iota,i')$ are Cartesian;
\item $\tf$ and $\tf'$ are smooth.   
\end{itemize}    
We set $\Omega=I_S^{-1}\otimes_{\co_S}\Omega^1_{X/S}$, $\Omega'=I_S^{-1}\otimes_{\co_S}\Omega^1_{X'/S}$
and $\uOmega'=I_S^{-1}\otimes_{\co_S}\Omega^1_{X'/X}$. 

We denote by $\fX$ (resp.\ $\fX'$) the $p$-adic completion of $X$ (resp.\ $X'$) and by $\fgg\colon \fX'\rightarrow \fX$ the morphism induced by $g$. 
The $p$-adic formal schemes $\fX$ and $\fX'$ are locally of finite presentation over $\cS=\Spf(V)$, and are hence idyllic (\cite{egr1} 2.6.13).  
We denote by $\hOmega$ (resp.\ $\hOmega'$, resp.\  $\huOmega'$) the $p$-adic completion of 
$\Omega$ (resp.\ $\Omega'$, resp.\ $\uOmega'$), 
which is canonically isomorphic to the module $I_S^{-1}\otimes_{\co_S}\Omega^1_{\fX/\cS}$ 
(resp.\ $I_S^{-1}\otimes_{\co_S}\Omega^1_{\fX'/\cS}$, resp.\ $I_S^{-1}\otimes_{\co_S}\Omega^1_{\fX'/\fX}$).

\subsection{}\label{tphdi2}
Following \ref{torlift1}, let $\cL_{\tX'/\tX}$ be the torsor of local liftings of $g$ over $\tX$, 
\begin{equation}
0\rightarrow \co_{X'}\rightarrow \cF_\uptau\rightarrow g^*(\Omega)\rightarrow 0
\end{equation}
the associated Higgs--Tate $\co_{X'}$-extension and $\cC_\uptau$ the associated Higgs--Tate $\co_{X'}$-algebra. 
The weak $p$-adic completion $\cC^\dagger_\uptau$ of $\cC_\uptau$ \eqref{thm2c} 
is naturally an $\co_{\fX'}$-algebra. It is equipped with a canonical $\co_{\fX'}$-derivation \eqref{thm2d}
\begin{equation}
\delta \colon \cC^\dagger_\uptau\rightarrow \fgg^*(\hOmega)\otimes_{\co_{\fX'}}\cC^\dagger_\uptau.
\end{equation}
We denote by 
\begin{equation}
\delta' \colon \cC^\dagger_\uptau\rightarrow \hOmega'\otimes_{\co_{\fX'}}\cC^\dagger_\uptau
\end{equation}
the $\co_{\fX'}$-derivation induced by $\delta$ and the canonical morphism $\fgg^*(\hOmega)\rightarrow \hOmega'$. 

\subsection{}\label{tphdi3}
We denote by $\hcF_\uptau$ the $p$-adic completion of $\cF_\uptau$; so we have an exact sequence of $\co_{\fX'}$-modules
\begin{equation}\label{tphdi3a}
0\rightarrow \co_{\fX'}\rightarrow \hcF_\uptau\rightarrow \fgg^*(\hOmega) \rightarrow 0.
\end{equation}
We defined in \eqref{thm3b} a twisting functor by the extension $\hcF_\uptau$ 
for Higgs $\co_{\fX'}$-modules with coefficients in $\fgg^*(\hOmega)$:
\begin{equation}\label{tphdi3b}
\uptau\colon \bHM(\co_{\fX'},\fgg^*(\hOmega))\rightarrow \bHM(\co_{\fX'},\fgg^*(\hOmega)).
\end{equation}
Composing with the pullback functor $\fgg^*$, we obtain a functor
\begin{equation}\label{tphdi3c}
\fgg^*_\uptau\colon 
\begin{array}[t]{clcr}
\bHM(\co_\fX,\hOmega)&\rightarrow& \bHM(\co_{\fX'},\fgg^*(\hOmega))\\
(N,\theta)&\mapsto&\uptau(\fgg^*(N),\fgg^*(\theta)),
\end{array}
\end{equation} 
that we call the {\em pullback functor by $\fgg$ twisted by the extension $\hcF_\uptau$} 
(or simply the {\em twisted pullback functor by $\fgg$} when the extension is implicit);
see \ref{p1-tphdi2}, \ref{p2-cmupiso210} and \ref{p2-fhtft4}.

\begin{prop}[cf. \ref{p2-cmupiso211}]
For every locally CL-small Higgs $\co_{\fX}[\frac 1 p]$-module $(N,\theta)$ with coefficients in $\hOmega$ \eqref{thm66}, 
the Higgs $\co_{\fX'}[\frac 1 p]$-module $\fgg^*_\uptau(N,\theta)$ \eqref{tphdi3c} is locally CL-small. 
\end{prop}

\subsection{}\label{tphdi4}
Let $(N,\theta)$ be a Higgs $\co_{\fX'}$-module with coefficients in $\hOmega'$. We equip $N\otimes_{\co_{\fX'}}\cC^\dagger_\uptau$ with 
the total Higgs $\co_{\fX'}$-field 
\begin{equation}\label{tphdi4a}
\vartheta= \theta\otimes \id+\id\otimes \delta'\colon N\otimes_{\co_{\fX'}}\cC^\dagger_\uptau \rightarrow 
\hOmega'\otimes_{\co_{\fX'}}N\otimes_{\co_{\fX'}}\cC^\dagger_\uptau,
\end{equation}
and denote by $\mK^{\bullet}(N\otimes_{\co_{\fX'}}\cC^\dagger_\uptau)$ the associated Dolbeault complex, 
that we call the {\em Dolbeault complex of $(N,\theta)$ twisted by the extension $\hcF_\uptau$} \eqref{tphdi3a}. 
The Higgs field $\id\otimes \delta$ on $N\otimes_{\co_{\fX'}}\cC^\dagger_\uptau$ 
induces a morphism of complexes of $\co_{\fX'}$-modules
\begin{equation}\label{tphdi4b}
\fd \colon \mK^\bullet(N\otimes_{\co_{\fX'}}\cC^\dagger_\uptau)\rightarrow \fgg^*(\hOmega)\otimes_{\co_{\fX'}} \mK^\bullet(N\otimes_{\co_{\fX'}}\cC^\dagger_\uptau).
\end{equation} 
By the projection formula, for any integer $q\geq 0$, $\rR^q\fgg_*(\fd)$ identifies with a Higgs $\co_\fX$-field 
\begin{equation}\label{tphdi4c}
\rR^q\fgg_*(\fd)\colon  \rR^q\fgg_*(\mK^\bullet(N\otimes_{\co_{\fX'}}\cC^\dagger_\uptau))\rightarrow \hOmega\otimes_{\co_\fX} 
\rR^q\fgg_*(\mK^\bullet(N\otimes_{\co_{\fX'}}\cC^\dagger_\uptau)). 
\end{equation}
We obtain a functor 
\begin{equation}\label{tphdi4d}
\rR^q\fgg^\uptau_*\colon 
\begin{array}[t]{clcr}
\bHM(\co_{\fX'},\hOmega')&\rightarrow&\bHM(\co_\fX,\hOmega)\\
(N,\theta)&\mapsto&(\rR^q\fgg_*(\mK^\bullet(N\otimes_{\co_{\fX'}}\cC^\dagger_\uptau)),-\rR^q\fgg_*(\fd)),
\end{array}
\end{equation}
that we call the {\em $q$th higher direct image functor by $\fgg$ twisted by the extension $\hcF_\uptau$} 
(or simply the {\em twisted $q$th higher direct image functor by $\fgg$} when the extension is implicit); 
see \ref{p1-tphdi6}, \ref{p2-cmupiso210} and \ref{p2-fhtft8}. 

\subsection{}\label{tphdi40}
Assume that $g\colon X'\rightarrow X$ is smooth. 
Therefore, we have a canonical exact sequence of locally free $\co_{X'}$-modules of finite type
\begin{equation}\label{tphdi1b}
0\rightarrow g^*(\Omega)\rightarrow \Omega' \rightarrow \uOmega' \rightarrow 0. 
\end{equation}
Let $(N,\theta)$ be a Higgs $\co_{\fX'}$-module with coefficients in $\hOmega'$,
$\utheta\colon N\rightarrow N\otimes\huOmega'$ the Higgs $\co_{\fX'}$-field induced by $\theta$, 
$\umK^\bullet(N)$ the Dolbeault complex of $(N,\utheta)$. 
The Koszul filtration associated with the $p$-adic completion of the extension \eqref{tphdi1b} induces, for every integer $q\geq 0$, 
a canonical Higgs $\co_\fX$-field on $\rR^q \fgg_*(\umK^\bullet(N))$ with coefficients in $\hOmega$, 
that we call the {\em Katz-Oda field}; see \ref{p1-tphdi4}.

\begin{prop}[cf. \ref{p2-cmupiso29} and \ref{p2-cmupiso36}]\label{tphdi5}
Assume that $g\colon X'\rightarrow X$ is smooth. 
Let $(N,\theta)$ be a locally CL-small Higgs $\co_{\fX'}[\frac 1 p]$-module with coefficients in $\hOmega'$ \eqref{thm66},
$\mK^{\bullet}(N\otimes_{\co_{\fX'}}\cC^\dagger_\uptau)$ the Dolbeault complex defined in \ref{tphdi4}, 
$\utheta\colon N\rightarrow N\otimes\huOmega'$ the Higgs $\co_{\fX'}$-field induced by $\theta$, 
$\umK^\bullet(N)$ the Dolbeault complex of $(N,\utheta)$. 
Then, a lifting of $g$ to an $\tS$-morphism $\tX'\rightarrow \tX$ determines a section 
$\cC_\uptau\rightarrow \co_{\fX'}$, which induces a quasi-isomorphism
\begin{equation}\label{tphdi5a}
\upbeta\colon \mK^\bullet(N\otimes_{\co_{\fX'}}\cC^\dagger_\uptau)\rightarrow \umK^\bullet(N).
\end{equation}
Moreover, the latter induces, for every integer $q\geq 0$, an isomorphism from 
$\rR^q \fgg^\uptau_*(N,\theta)$ to $\rR^q \fgg_*(\umK^\bullet(N))$ equipped with the Katz-Oda field \eqref{tphdi40}. 
\end{prop}

We establish in \ref{p2-cmupiso29} analogous results for Dolbeault complexes of ind-$\co_{\fX'}$-modules. 

\begin{cor}[cf. \ref{p2-cmupiso30}]\label{tphdi6}
Assume that $g\colon X'\rightarrow X$ is proper and smooth. 
Let $(N,\theta)$ be a locally CL-small Higgs $\co_{\fX'}[\frac 1 p]$-module with coefficients in $\hOmega'$ \eqref{thm66},
$\mK^{\bullet}(N\otimes_{\co_{\fX'}}\cC^\dagger_\uptau)$ the Dolbeault complex defined in \ref{tphdi4}. 
Then, for every integer $q\geq 0$, the $\co_\fX[\frac 1 p]$-module $\rR^q\fgg_*(\mK^\bullet(N\otimes_{\co_{\fX'}}\cC^\dagger_\uptau))$ is coherent.
\end{cor}

We prove in \ref{p2-cmupiso30} further cohomological properties of the Dolbeault complex $\mK^\bullet(N\otimes_{\co_{\fX'}}\cC^\dagger_\uptau)$
and its ind-$\co_{\fX'}$-modules analogue, especially, an important base change isomorphism \eqref{p2-cmupiso30c}. 

\begin{teo}[cf. \ref{p2-cmupiso47}]\label{tphdi11} 
Assume that $g\colon X'\rightarrow X$ is proper and smooth. 
Let $(N,\theta)$ be a locally CL-small Higgs $\co_{\fX'}[\frac 1 p]$-module with coefficients in $\hOmega'$ \eqref{thm66}.
Then, for every integer $q\geq 0$, the Higgs $\co_{\fX}[\frac 1 p]$-module $\rR^q\fgg^\uptau_*(N,\theta)$ \eqref{tphdi4d} is locally CL-small. 
\end{teo}

\begin{cor}[cf. \ref{p2-cmupiso49}]\label{tphdi12}
Assume that $g\colon X'\rightarrow X$ is proper and smooth. 
Let $(N,\theta)$ be a locally CL-small Higgs $\co_{\fX'}[\frac 1 p]$-module with coefficients in $\hOmega'$ \eqref{thm66},
$\utheta\colon N\rightarrow N\otimes\huOmega'$ the Higgs $\co_{\fX'}$-field induced by $\theta$, 
$\umK^\bullet(N)$ the Dolbeault complex of $(N,\utheta)$. 
Then, for every integer $q\geq 0$, the Higgs $\co_\fX[\frac 1 p]$-module $(\rR^qf_*(\umK^{\bullet}(N)),\kappa^q)$ is locally CL-small, where $\kappa^q$
is Katz-Oda Higgs $\co_\fX$-field. 
\end{cor}

Assuming that $g$ is proper and smooth and the canonical homomorphism $\co_\fX\rightarrow \fgg_*(\co_{\fX'})$ is an isomorphism,
we express in \ref{p2-cmupiso48} the difference between the isomorphisms $\rR^q\fgg_*(\upbeta)$ \eqref{tphdi5a} for two different liftings of $g$
in terms of the {\em exponential} of the Katz-Oda Higgs field $\kappa^q$ on $\rR^qf_*(\umK^{\bullet}(N))$. 

Besides \ref{tphdi5}, a key ingredient in the proofs of \ref{tphdi11} and \ref{p2-cmupiso48}  
is a {\em Taylor formula} for the Higgs--Tate algebra; see \ref{p2-cmupiso41}.

\section{\texorpdfstring{Revisiting the local $p$-adic Simpson correspondence}{Revisiting the local p-adic Simpson correspondence}}\label{rlps}

\subsection{}\label{rlps0}
As in \ref{overview-intro2}, let $K$ be a complete discrete valuation field of characteristic $0$, with perfect residue field $k$ of characteristic $p>0$,
$\co_K$ the valuation ring of $K$, $\oK$ an algebraic closure of $K$, $\co_\oK$ the integral closure of $\co_K$ in $\oK$,
$\co_C$ the $p$-adic Hausdorff completion of $\co_\oK$, $C$ its field of fractions. 
We set $S=\Spec(\co_K)$, $\oS=\Spec(\co_\oK)$ and $\coS=\Spec(\co_C)$, 
and we denote by $s$ (resp.\ $\eta$, resp.\ $\oeta$) 
the closed point of $S$ (resp.\ the generic point of $S$, resp.\ the generic point of $\oS$).

\subsection{}\label{rlps1}
Let $X=\Spec(R)$ be a {\em small} affine $S$-scheme in the sense of Faltings  (i.e., which admits
an étale $S$-morphism $X\rightarrow \mG_{m,S}^d=\Spec(\co_K[T_1^{\pm1},\dots, T_d^{\pm1}])$ for an integer $d\geq 0$) and
such that $X_s\not=\emptyset$. 
For simplicity, we assume that $X_\oeta$ is connected and fix a geometric generic point $\ox$ of $X_\oeta$. 
We set $\Delta=\pi_1(X_\oeta,\ox)$ and $R_1=R\otimes_{\co_K}\co_C$, and denote by $\hRun$ the $p$-adic Hausdorff completion of $R_1$. 
Let $(X_i)_{i\in I}$ be the universal étale cover of $X_\oeta$ at $\ox$ (\cite{ag1} 2.1.20). For any $i\in I$, we denote by $\oX_i=\Spec(R_i)$ 
the normalization of $\oX=X\times_S\oS$ in $X_i$.  Let $\oR$ be the $\co_\oK$-representation of $\Delta$ defined by 
\begin{equation}\label{rlps1b}
\oR=\underset{\underset{i\in I}{\longrightarrow}}{\lim} \ R_i.
\end{equation}
We denote by $\hoR$ its $p$-adic Hausdorff completion
and by $\bRep_{\hoR}(\Delta)$ the category of $\hoR$-representations of $\Delta$ (\cite{ag2} 2.1.2).
Consider a commutative diagram
\begin{equation}\label{rlps1a}
\xymatrix{
{\hmX=\Spec(\hoR)}\ar[r]^-(0.5)j\ar[d]_h&{\tmX=\Spec(\cA_2(\oR^\flat))}\ar@/^2pc/[dd]\\
{\coX=\Spec(R_1)}\ar[r]^-(0.5)i\ar[d]\ar@{}[rd]|\Box&{\tX}\ar[d]\\
{\coS=\Spec(\co_C)}\ar[r]^-(0.5)\iota&{\tS=\Spec(\cA_2(\co_{\oK^\flat})),}}
\end{equation}
where $\iota$ and $j$ are {\em Fontaine's universal $p$-adic first order thickenings}, 
and $\tX$ is a smooth $\tS$-deformation of $\coX$ over $\coS$, so the lower square is Cartesian.  
Recall that $\iota$ is the first infinitesimal neighborhood of the closed immersion defined by 
Fontaine's surjective homomorphism $\theta\colon \rW(\co_{\oK^\flat})\rightarrow \co_C$ (\cite{ag2} 3.1.4), 
and that the kernel of $\theta$ is generated by a non-zero divisor $\xi$; and similarly for $j$ (\cite{ag2} 3.2.8). 

Such an $\tS$-smooth deformation $\tX$ of $\coX$ exists and is unique up to non-canonical isomorphism; 
nevertheless, the $p$-adic Simpson correspondence depends on its choice.

The ideal $\xi\co_\tS$ is a free $\co_\coS$-module with basis $\xi$. We denote it by $\xi\co_C$ and denote by $\xi^{-1}\co_C$ its $\co_C$-dual. 
For any $\co_C$-algebra $A$, 
any $A$-module $M$ and any integer $i\geq 1$, we denote the $A$-modules $M\otimes_{\co_C}(\xi \co_C)^{\otimes i}$ 
and $M\otimes_{\co_C}(\xi^{-1} \co_C)^{\otimes i}$ simply by $\xi^i M$ and $\xi^{-i} M$, respectively.

\subsection{}\label{rlps100}
Following \ref{torlift1}, let $\cL_{\tmX/\tX}$ be the torsor of local liftings of $h$ over $\tX$ and 
\begin{equation}\label{rlps100a}
0\rightarrow \hoR\rightarrow \fF\rightarrow \Omega \otimes_{R_1} \hoR\rightarrow 0
\end{equation}
the associated Higgs--Tate $\hoR$-extension, where $\Omega=\xi^{-1}\Omega^1_{R/\co_K}\otimes_RR_1$, 
\begin{equation}\label{rlps100b}
\fC=\underset{\underset{n\geq 0}{\longrightarrow}}\lim\ \Sym^n_{\hoR}(\fF), 
\end{equation}
the associated Higgs--Tate $\hoR$-algebra. 
Observe that $\Delta$ acts naturally on $\fF$ in a $\hoR$-semilinear way 
such that the morphisms in \eqref{rlps100a} are $\Delta$-equivariant (\cite{ag2} 3.2.15).
The latter induces an action of $\Delta$ on $\fC$ by ring automorphisms compatible with its action on $\hoR$; see \ref{p2-rlps5}.  
These actions are continuous for the $p$-adic topology (\cite{ag2} 3.2.23). 
We consider the {\em twisting functors by the extension $\fF$}, \eqref{thm3b} and \eqref{thm3c},
\begin{eqnarray}
\uptau\colon 
\begin{array}[t]{clcr}
\bHM(\hoR, \Omega) &\rightarrow& \bHM(\hoR, \Omega),\\
(N,\theta)&\mapsto&(\uuptau(N,\theta), \theta_\uptau),
\end{array}\label{rlps1c}\\
\uptau^\vee\colon 
\begin{array}[t]{clcr}
\bHM(\hoR,\Omega) &\rightarrow& \bHM(\hoR,\Omega),\\
(N,\theta)&\mapsto&(\uuptau^\vee(N,\theta), \theta_{\uptau^\vee}).
\end{array}\label{rlps1d}
\end{eqnarray}

\subsection{}\label{rlps2}
Let $(N,\theta)$ be a Higgs $\hRun$-module with coefficients in $\Omega$,  $(\tN,\ttheta)=(N,\theta)\otimes_{\hRun}\hoR$. 
The Higgs $\hoR$-module $\uptau(\tN,\ttheta)$ is canonically equipped with an $\hoR$-semilinear action of $\Delta$ \eqref{p2-rlps17}.  
We thus define a functor 
\begin{equation}\label{rlps2a}
\mV\colon 
\begin{array}[t]{clcr}
\bHM(\hRun,\Omega)&\rightarrow& \bRep_{\hoR}(\Delta)\\
(N,\theta)&\mapsto& \uuptau(\tN,\ttheta).
\end{array}
\end{equation}
We recover the functor defined in (\cite{ag2} 3.3.7); see \ref{p2-rlps8}. 

\begin{defi}
We say that $(N,\theta)$ is {\em solvable} if $N$ is a projective $\hRun[\frac 1 p]$-module of finite type and if $(\tN,\ttheta)$ is weakly twistable 
by the extension $\fF$ in the sense of \ref{thm4}. 
\end{defi}

We recover the notion introduced in (\cite{agt} II.12.12 and \cite{ag2} 3.3.10); see \ref{p2-rlps10}. 

We prove in \ref{p2-rlps26} that under the first condition, the second condition is equivalent to $(N,\theta)$ being CL-small in the sense of \ref{thm66}.

\begin{prop}\label{rlps3} 
The restriction of the functor $\mV$ to the full subcategory of solvable Higgs $\hRun[\frac 1 p]$-modules with coefficients in $\Omega$ is fully faithful. 
We say that 
an $\hoR[\frac 1 p]$-representation of $\Delta$ is {\em Dolbeault} if it is in the essential image of this fully faithful functor, 
so that $\mV$ induces an equivalence of categories 
\begin{equation}\label{rlps3a}
\mV\colon \bHM^\sol(\hRun[\frac 1 p],\Omega)\stackrel{\sim}{\rightarrow} \bRep^\Dolb_{\hoR[\frac 1 p]}(\Delta).
\end{equation}
\end{prop}

We can explicitly describe the essential image of the functor $\mV$ in terms of a {\em smallness} condition on 
$\hoR[\frac{1}{p}]$-representations of $\Delta$ (\cite{ag2} 3.4.19). 
This description was already used implicitly by Faltings \cite{faltings3} and later proved by Tsuji (\cite{tsuji5} 13.7); see \ref{p2-rlps267}. 

A quasi-inverse functor of \eqref{rlps3a} is induced by the functor
\begin{equation}\label{rlps3b} 
\mH\colon 
\begin{array}[t]{clcr}
\bRep_{\hoR}(\Delta)&\rightarrow &\bHM(\hRun,\Omega),\\
M&\mapsto&((M\otimes_{\hoR}\fC^\dagger)^\Delta,\id\otimes \delta),
\end{array}
\end{equation}
where $\fC^\dagger$ is the weak $p$-adic completion of the Higgs--Tate algebra $\fC$ \eqref{rlps100b}, see \eqref{thm2c}, and 
$\delta\colon \fC^\dagger\rightarrow \Omega\otimes_{R_1}\fC^\dagger$ is the $\hoR$-derivation defined in \eqref{thm2d}; see \ref{p2-rlps11}.

\subsection{}\label{rlps4} 
Let $M$ be a Dolbeault $\hoR[\frac 1 p]$-representation of $\Delta$, $(N,\theta)=\mH(M)$, $(\tN,\ttheta)=(N,\theta)\otimes_{\hRun}\hoR$. 
We have a canonical isomorphism $M\stackrel{\sim}{\rightarrow} \mV(N,\theta)=\uuptau(\tN,\ttheta)$ of $\hoR[\frac 1 p]$-representations of $\Delta$. 
We deduce a canonical $\Delta$-equivariant Higgs $\hoR$-field $\theta_M$ on $M$  
with coefficients in $\Omega$, so that we have a $\Delta$-equivariant isomorphism of Higgs $\hoR$-modules
\begin{equation}\label{rlps4a}
(M,\theta_M) \stackrel{\sim}{\rightarrow} \uptau(\tN,\ttheta). 
\end{equation}
We thus define a canonical functor, depending a priori on the deformation $\tX$,
\begin{equation}\label{rlps4b}
\Theta\colon 
\begin{array}[t]{clcr}
\bRep_{\hoR[\frac 1 p]}^\Dol(\Delta)&\rightarrow&\bHM(\hoR[\frac 1 p],\Omega,\Delta)\\
M&\mapsto&(M,\theta_M),
\end{array}
\end{equation}
in the category of $\Delta$-equivariant Higgs $\hoR[\frac 1 p]$-modules with coefficients in $\Omega$, 
and a canonical functorial isomorphism
\begin{equation}\label{rlps4c}
\Theta(M) \stackrel{\sim}{\rightarrow} \uptau(\mH(M)\otimes_\hRun\hoR). 
\end{equation}
By \ref{thm10}(iii), we deduce a canonical functorial $\Delta$-equivariant isomorphism of Higgs $\hoR$-modules
\begin{equation}\label{rlps4d}
\mH(M)\otimes_\hRun\hoR \stackrel{\sim}{\rightarrow} \uptau^\vee(\Theta(M)).
\end{equation}
Taking $\Delta$-invariants, we obtain a canonical isomorphism of Higgs $\hRun$-modules
\begin{equation}\label{rlps4e}
\mH(M)\stackrel{\sim}{\rightarrow} \uptau^\vee(\Theta(M))^\Delta.
\end{equation}

\begin{prop}[cf. \ref{p2-rlps23}]\label{rlps5}
For every Dolbeault $\hoR[\frac 1 p]$-representation $M$ of $\Delta$, the canonical Higgs $\hoR$-field $\theta_M$ \eqref{rlps4b} does not depend 
on the choice of the deformation $\tX$ \eqref{rlps1a}. 
\end{prop}

Similar canonical Higgs fields were introduced by Camargo on pro-étale vector bundles (\cite{rod} 1.0.3), and by He in his Sen theory (\cite{he} 5.28). 
They were subsequently used by Heuer \cite{heuer} and Heuer-Xu \cite{hexu} in their constructions of the $p$-adic Simpson correspondence.

\section{\texorpdfstring{Functoriality of the local $p$-adic Simpson correspondence by pullback}{Functoriality of the local p-adic Simpson correspondence by pullback}}\label{nfspb}

\subsection{}\label{fspb1}
In this section, we take again the assumptions and notation of \ref{rlps0} 
and let $X=\Spec(R)$ and $X'=\Spec(R')$ be two {\em small} affine $S$-schemes, $\gamma\colon X'\rightarrow X$ an $S$-morphism. 
For simplicity, we assume that $X_\oeta$ and $X'_\oeta$ are connected and let $\ox'$ be a geometric point of $X'_\oeta$ and $\ox=\gamma(\ox')$. 
We set $\coX=X\times_S\coS$ and $\coX'=X'\times_S\coS$, and consider a commutative diagram  
\begin{equation}\label{fspb1a}
\xymatrix{
{\coX'}\ar[r]\ar[d]_{\gamma}&{\tX'}\ar@/^1pc/[dd]&\\
{\coX}\ar[r]\ar[d]&{\tX}\ar[d]\\
{\coS}\ar[r]^-(0.5)\iota&{\tS,}}
\end{equation}
where $\iota$ is Fontaine's universal $p$-adic first order thickening \eqref{rlps1a}, 
and $\tX$ (resp.\ $\tX'$) is a smooth $\tS$-deformation of $\coX$ (resp.\ $\coX'$) over $\coS$.  
We consider again the objects associated with $(X,\tX)$ introduced in §\ref{rlps}, and we associate with $(X',\tX')$ 
similar objects that we denote by the same symbols equipped with a $^\prime$ exponent.
We then have equivalences of categories quasi-inverse to each other  
\begin{eqnarray}
\xymatrix{
{\bRep^\Dolb_{\hoR[\frac 1 p]}(\Delta)}\ar@<1ex>[r]^-(0.5){\mH}&{\bHM^\sol(\hRun[\frac 1 p], \Omega),}
\ar@<1ex>[l]^-(0.5){\mV}}\label{fspb1b}\\
\xymatrix{
{\bRep^\Dolb_{\hoRp[\frac 1 p]}(\Delta')}\ar@<1ex>[r]^-(0.5){\mH'}&{\bHM^\sol(\hRunp[\frac 1 p], \Omega'),}
\ar@<1ex>[l]^-(0.5){\mV'}}\label{fspb1c}
\end{eqnarray}
where $R_1=R\otimes_{\co_K}\co_C$, $\Omega=\xi^{-1}\Omega^1_{R/\co_K}\otimes_{R}R_1$,
$R'_1=R'\otimes_{\co_K}\co_C$ and $\Omega'=\xi^{-1}\Omega^1_{R'/\co_K}\otimes_{R'}R'_1$. 

\subsection{}\label{fspb10}
We consider the torsor $\cL_{\tX'/\tX}$ of local liftings of $\gamma$ over $\tX$ \eqref{fspb1a} (see \ref{torlift1}), 
the associated Higgs--Tate $R'_1$-extension
\begin{equation}\label{fspb10a}
0\rightarrow R'_1\rightarrow \fF_\uptau\rightarrow \Omega\otimes_{R_1}R'_1 \rightarrow 0,
\end{equation}
and the associated Higgs--Tate $R'_1$-algebra
\begin{equation}\label{fspb10b}
\fC_\uptau=\underset{\underset{n\geq 0}{\longrightarrow}}\lim\ \Sym^n_{R'_1}(\fF_\uptau). 
\end{equation}
We consider the {\em pullback functor twisted by the extension $\fF_\uptau\otimes_{R'_1}\hRunp$}, \eqref{tphdi3c}, 
\begin{equation}\label{fspb10c}
\mT\colon \bHM(\hRun,\Omega) \rightarrow \bHM(\hRunp,\Omega),
\end{equation}
and the functor
\begin{equation}\label{fspb10d}
\upmu\colon \bHM(\hRunp,\Omega) \rightarrow \bHM(\hRunp,\Omega')
\end{equation} 
induced by the canonical morphism $\Omega\otimes_{R_1}R'_1\rightarrow \Omega'$. 

\begin{teo}[cf. \ref{p2-fspb9}]\label{fspb3}
Let $M$ be a Dolbeault $\hoR[\frac 1 p]$-representation of $\Delta$. 
We set $M'=M\otimes_\hoR \hoRp$ that we equip with the 
$\hoRp$-representation of $\Delta'$ induced by the canonical homomorphism $\Delta'\rightarrow \Delta$. 
Then, $M'$ is a Dolbeault $\hoRp[\frac 1 p]$-representation of $\Delta'$, and we have a canonical functorial isomorphism 
\begin{equation}\label{fspb3a}
\mH'(M')\stackrel{\sim}{\rightarrow}\upmu(\mT(\mH(M))).
\end{equation}
\end{teo}

The proof of this result relies on a twisted version \eqref{fspb2c} of the composition homomorphism of Higgs--Tate algebras \eqref{torlift2c}. 
Consider the commutative diagram 
\begin{equation}\label{fspb2a}
\xymatrix{
{\hmX'=\Spec(\hoRp)}\ar[r]^-(0.5){j'}\ar[d]_{h'}&{\tmX'=\Spec(\cA_2(\oR'^\flat))}\ar@/^2pc/[ddd]\\
{\coX'}\ar[r]\ar[d]_{\gamma}&{\tX'}\ar@/^1pc/[dd]&\\
{\coX}\ar[r]\ar[d]&{\tX}\ar[d]\\
{\coS}\ar[r]^-(0.5)\iota&{\tS,}}
\end{equation}
where $j'$ is Fontaine's universal $p$-adic first order thickening \eqref{rlps1a}.  
We prove that $\fC_{\hoRp}=\fC\otimes_{\hoR}\hoRp$ \eqref{rlps100b} is the Higgs--Tate algebra associated 
with the torsor $\cL_{\tmX'/\tX}$ \eqref{torlift1}; see \ref{p1-rdt4}. 
The composition homomorphism of Higgs--Tate algebras \eqref{torlift2c} relative to \eqref{fspb2a} is a homomorphism of $\hoRp$-algebras
\begin{equation}\label{fspb2b}
\upphi\colon \fC_{\hoRp}\rightarrow \fC'\otimes_{R'_1}\fC_\uptau. 
\end{equation}

For any rational numbers $r \geq r' \geq 0$, $\upphi$ can be refined into a canonical homomorphism of $\hoRp$-algebras
\begin{equation}\label{fspb2c}
\upphi^{(r,r')}\colon \fC^{(r)}_{\hoRp}\rightarrow \fC'^{(r)} \otimes_{R'_1}\fC^{(r')}_\uptau,
\end{equation}
which is the key ingredient for the functoriality of the $p$-adic Simpson correspondence under pullback \eqref{fspb3}. 
However, the study of functoriality with respect to higher direct images requires a further refinement, depending on three radii that we introduce below.

\subsection{}\label{fspb4}
Assume that $\gamma\colon X'\rightarrow X$ is smooth. 
We have an exact sequence of free $R'_1$-modules of finite type 
\begin{equation}\label{fspb4a}
0\rightarrow \Omega\otimes_{R_1}R'_1 \stackrel{u}{\rightarrow} \Omega' \rightarrow \uOmega'\rightarrow 0,
\end{equation}
where $\uOmega'=\xi^{-1}\Omega^1_{R'/R}\otimes_{R'}R'_1$. 
For any rational number $r\geq 0$, we set
\begin{equation}\label{fspb4b}
\Omega^{(r)}=p^r\Omega,\ \ \ \Omega'^{(r)}=p^r\Omega',\ \ \ \uOmega'^{(r)}=p^r\uOmega'.
\end{equation}
For any rational numbers $r\geq r'\geq 0$, we define the $R'_1$-module $\Omega'^{(r,r')}$ by the commutative diagram
\begin{equation}\label{fspb4c}
\xymatrix{
0\ar[r]&{\Omega^{(r')}\otimes_{R_1}R'_1}\ar[r]^-(0.5){u^{(r')}}&{\Omega'^{(r')}}\ar[r]&{\uOmega'^{(r')}}\ar[r]&0\\
0\ar[r]&{\Omega^{(r')}\otimes_{R_1}R'_1}\ar[r]^-(0.5){u^{(r,r')}}\ar@{=}[u]&{\Omega'^{(r,r')}}\ar[r]\ar@{^(->}[u]\ar@{}[ru]|{(2)}&{\uOmega'^{(r)}}\ar[r]\ar[u]&0\\
0\ar[r]&{\Omega^{(r)}\otimes_{R_1}R'_1}\ar[r]^-(0.5){u^{(r)}}\ar[u]\ar@{}[ru]|{(1)}&{\Omega'^{(r)}}\ar[r]\ar@{^(->}[u]&{\uOmega'^{(r)}}\ar[r]\ar@{=}[u]&0,}
\end{equation}
where square $(1)$ is co-Cartesian, square $(2)$ is Cartesian and the morphisms $u^{(r)}$ and $u^{(r,r')}$ are induced by $u$. 

Using a different, but equivalent, normalization from \ref{thm2}, we consider the Higgs--Tate extension 
$\fF^{(r)}$ deduced from the Higgs--Tate extension $\fF$ \eqref{rlps100a}
by pullback by the canonical morphism $\Omega^{(r)}\otimes_{R_1}\hoR \rightarrow \Omega\otimes_{R_1}\hoR$;
so we have an exact sequence of $\hoR$-modules  
\begin{equation}
0\rightarrow \hoR\rightarrow \fF^{(r)}\rightarrow \Omega^{(r)} \otimes_{R_1} \hoR\rightarrow 0.
\end{equation}
We consider similarly the Higgs--Tate extension $\fF'^{(r')}$, relative to $(X',\tX')$;
so we have an exact sequence of $\hoRp$-modules  
\begin{equation}
0\rightarrow \hoRp\rightarrow \fF'^{(r')}\rightarrow \Omega'^{(r')} \otimes_{R'_1} \hoRp\rightarrow 0.
\end{equation}
We denote by $\fF'^{(r,r')}$ the pullback of this extension 
by the canonical morphism $\Omega'^{(r,r')}\otimes_{R'_1}\hoRp \rightarrow \Omega'^{(r')}\otimes_{R'_1}\hoRp$; 
so we have an exact sequence of $\hoRp$-modules  
\begin{equation}\label{fspb4d}
0\rightarrow \hoRp\rightarrow \fF'^{(r,r')}\rightarrow \Omega'^{(r,r')} \otimes_{R'_1} \hoRp\rightarrow 0.
\end{equation}
We consider then the $\hoRp$-algebra
\begin{equation}\label{fspb4e}
\fC'^{(r,r')}=\underset{\underset{n\geq 0}{\longrightarrow}}\lim\ \Sym^n_{\hoRp}(\fF'^{(r,r')}).
\end{equation}

We denote by $I$ the set of triples of rational numbers $\ur=(r_1,r_2,r_3)$ such that $r_1\geq r_2\geq r_3\geq 0$. 
For any such a triple, we can refine the homomorphism \eqref{fspb2c} into a canonical homomorphism of $\hoRp$-algebras 
\begin{equation}\label{fspb4f}
\upphi^\ur\colon \fC^{(r_2)}_{\hoRp}\rightarrow \fC'^{(r_1,r_2)} \otimes_{R'_1}\fC^{(r_3)}_\uptau.
\end{equation}

\section{\texorpdfstring{Revisiting the global $p$-adic Simpson correspondence}{Revisiting the global p-adic Simpson correspondence}}
\label{ft}

\subsection{}\label{ft0}
In this section, we take again the assumptions and notation of \ref{rlps0} and let $X$ be a smooth $S$-scheme of finite type. 
The $p$-adic Simpson correspondence constructed in §\ref{rlps} does not readily glue to yield a global correspondence for schemes 
that are not small affine. To remedy this, we sheafify the Higgs--Tate algebra in the Faltings topos of the canonical morphism 
$h\colon X_\oeta\rightarrow X$ and use it to construct 
a global $p$-adic Simpson correspondence, parallel to the local construction.
For small affine schemes, this global correspondence is equivalent to the local one (\cite{ag2} 4.8.31 and 4.8.32).

Let $\iota\colon \coS\rightarrow \tS$ be Fontaine's universal $p$-adic first order thickening \eqref{rlps1a}.
{\em We suppose that there exists a smooth $\tS$-deformation $\tX$ of $\coX=X\times_S\coS$ that we fix:}
\begin{equation}\label{ft0a}
\xymatrix{
{\coX}\ar[r]\ar[d]\ar@{}[rd]|\Box&{\tX}\ar[d]\\
{\coS}\ar[r]^-(0.5){\iota}&{\tS.}}
\end{equation} 
The existence of such an $\tS$-smooth deformation imposes a strong condition on $X$, and the $p$-adic Simpson correspondence depends on its choice.

\subsection{}\label{ft1}
We denote by $E$ the category of morphisms $(V\rightarrow U)$ over $h$, i.e., commutative diagrams
\begin{equation}\label{ft1a}
\xymatrix{V\ar[r]\ar[d]&U\ar[d]\\
{X_\oeta}\ar[r]^h&X,}
\end{equation}
such that $U$ is étale over $X$ and the canonical morphism $V\rightarrow U_\oeta$ is {\em finite étale}.
We consider it as fibered by the functor
\begin{equation}\label{ft1b}
\pi\colon E\rightarrow \Et_{/X}, \ \ \ (V\rightarrow U)\mapsto U,
\end{equation}
over the category of étale schemes over $X$ (\cite{ag2} 4.3.2). We refer to \ref{p2-ncgt6} for the notation related to the étale topology. 
The fiber category over any $U\in \ob(\Et_{/X})$ is canonically equivalent to
the category $\Et_{\rf/U_\oeta}$ of finite étale schemes over $U_\oeta$. 
We consider $\pi$ as a {\em covanishing fibered site} \eqref{p2-cmt46} by endowing each fiber category with the étale topology. 
We call it the {\em Faltings (covanishing) fibered site} associated with $h$; see \ref{p2-htaft3}.

We endow $E$ with the {\em covanishing topology} associated with $\pi$, that is, with
the topology generated by the coverings $\{(V_i\rightarrow U_i)\rightarrow (V\rightarrow U)\}_{i\in I}$
of the following two types:
\begin{itemize}
\item[(v)] $U_i=U$ for every $i\in I$, and $(V_i\rightarrow V)_{i\in I}$ is an étale covering.
\item[(c)] $(U_i\rightarrow U)_{i\in I}$ is an étale covering and $V_i=U_i\times_UV$ for every $i\in I$.
\end{itemize}
The resulting covanishing site $E$  is also called the {\em Faltings  site} associated with $h$.
We denote by $\tE$ the topos of sheaves of sets on $E$, called the {\em Faltings  topos} associated with $h$ (\cite{agt} VI.10.1, \cite{ag2} 4.3.2);
see \ref{p2-htaft4}.
By (\cite{agt} VI.5.11), $\tE$ is canonically equivalent to the category of collections $F=\{U \mapsto F_U\}$ such that:
\begin{itemize}
\item[(i)] for every object $U$ of $\Et_{/X}$, $F_U$ is a sheaf of $U_{\oeta,\fet}$;
\item[(ii)] for every morphism $f\colon U'\rightarrow U$ of $\Et_{/X}$, we are given a morphism $\gamma_f\colon F_U\rightarrow f_{\oeta*}(F_{U '})$;
\end{itemize}
and these data satisfy the usual cocycle condition (with respect to composition of morphisms) and the usual gluing condition 
(with respect to coverings in $\Et_{/X}$).

The functor
\begin{equation}\label{ft1c}
\sigma^+\colon 
\begin{array}[t]{clcr}
\Et_{/X}&\rightarrow& E \\
U&\mapsto &(U_\oeta\rightarrow U)
\end{array}
\end{equation}
is continuous and left exact (\cite{agt} VI.10.6). Hence it defines a morphism of topos
\begin{equation}\label{ft1d}
\sigma\colon \tE \rightarrow X_\et.
\end{equation}

\subsection{}\label{ft5}
For any object $(V\rightarrow U)$ of $E$, we denote by $\oU^V$ the integral closure of $\oU=U\times_S\oS$ in $V$.
We denote by $\ocB$ the presheaf on $E$ defined for any $(V\rightarrow U)\in \ob(E)$, by
\begin{equation}\label{ft5a}
\ocB((V\rightarrow U))=\Gamma(\oU^V,\co_{\oU^V}).
\end{equation}
It is a sheaf for the covanishing topology (\cite{agt} III.8.16).
Writing $\ocB=\{U\mapsto \ocB_U\}$, for every étale $X$-scheme which is affine, $U=\Spec(R)$, and every geometric point $\oy$ of $U_\oeta$, 
the stalk $\ocB_{U,\oy}$ is isomorphic to the $\co_\oK$-representation $\oR_U$ 
of $\pi_1(U_\oeta,\oy)$ defined in \eqref{rlps1b}. 

\subsection{}\label{ft6}
We denote by $\fX$ the formal scheme $p$-adic completion of $\oX=X\times_S\oS$, and by
$\Omega^1_{\fX/\cS}$ the sheaf of relative differentials.  
Similarly, we consider the $p$-adic completion of the ringed topos $(\tE,\ocB)$ defined as follows. 
We denote by $\tE_s$ the closed subtopos of $\tE$ complement of the open object $\sigma^*(X_\eta)$ 
(i.e., the full subcategory of $\tE$ consisting of those sheaves $F$ whose pullback to $\tE_{/\sigma^*(X_\eta)}$ is a final object; see \cite{ sga4} IV, 8.3),
by $\delta\colon \tE_s\rightarrow \tE$ the canonical embedding and by
\begin{equation}\label{ft6a}
\sigma_s\colon \tE_s\rightarrow X_{s,\et}
\end{equation}
the morphism of topos induced by $\sigma$; see \ref{p2-htaft5}. 
Observe that for every integer $n\geq 0$, the ring $\ocB/p^n\ocB$ is naturally an object of $\tE_s$.
We consider the topos $\tE^{\mN^\circ}_s$ of the projective systems of $\tE_s$ indexed by the ordered set $\mN$ \eqref{p2-ncgt5}, 
ringed by $\bvocB=(\ocB/p^n\ocB)_ {n\geq 0}$.
The morphism $\sigma$ induces a morphism of ringed topos 
\begin{equation}\label{ft6b}
\hupsigma \colon (\tE^{\mN^\circ}_s,\bvocB)\rightarrow (X_{s,\zar},\co_\fX).
\end{equation}

\subsection{}\label{ft7}
We now introduce the main categories involved in the $p$-adic Simpson correspondence.
We denote by $\bMod_{\mQ}(\bvocB)$ the category of {\em $\bvocB$-modules up to isogeny} (\cite{ag2} 2.9.1), which is the global analogue
of the category of $\hoR[\frac 1 p]$-representations of $\Delta$ considered in \ref{rlps1}. 
Since this category does not admit filtered direct limits,
we embed it in the abelian category $\bIndMod(\bvocB)$ of {\em ind-$\bvocB$-modules}, which does not have this defect; see \ref{p2-htaft27}. 
We will use the notation $\indcolim$ to denote direct limits in $\bIndMod(\bvocB)$ and we will keep the notation
$\underset{\longrightarrow}{\lim}$ to denote direct limits in $\bMod(\bvocB)$.
In the same way, we can naturally embed the category of {\em coherent} $\co_\fX[\frac 1 p]$-modules into the category $\bIndMod(\co_\fX)$ 
of ind-$\co_\fX$-modules; see \ref{p2-cmupiso2}.

The morphism $\hupsigma$ induces two adjoint functors 
\begin{equation}\label{ft7a}
\xymatrix{
{\bIndMod(\bvocB)}\ar@<1ex>[r]^-(0.5){\rI \hupsigma_*}&{\bIndMod(\co_\fX).}\ar@<1ex>[l]^-(0.5){\rI \hupsigma^*}}
\end{equation}
which extend the adjoint functors $\hupsigma^*$ and $\hupsigma_*$.
Let $\vupsigma_*$ be the composed functor
\begin{equation}\label{ft7b}
\vupsigma_*\colon 
\xymatrix{\bIndMod(\bvocB)\ar[r]^-(0.5){\rI \hupsigma_*}&{\bIndMod(\co_\fX)}\ar[r]^-(0.5){\kappa}&{\bMod(\co_\fX),}}
\end{equation}
where $\kappa$ is the functor defined by 
\begin{equation}\label{ft7c}
\kappa(\indcolim\alpha)= \underset{\longrightarrow}{\lim}\ \alpha.
\end{equation}

For $A=\co_\fX$ or $\co_\fX[\frac 1 p]$, we consider Higgs $A$-modules with coefficients in $\xi^{-1}\Omega^1_{\fX/\cS}$ 
\eqref{p1-delta-con1}. Such a Higgs module is said to be {\em coherent} if the underlying $A$-module is coherent.
The categories of these modules will be denoted by $\bHM(A,\xi^{-1}\Omega^1_{\fX/\cS})$ and $\bHM^\coh(A, \xi^{-1}\Omega^1_{\fX/\cS})$, respectively. 
We denote by $\bIndHM(\co_\fX,\xi^{-1}\Omega^1_{\fX/\cS})$ the category of Higgs ind-$\co_\fX$-modules with coefficients in 
$\xi^{-1}\Omega^1_{\fX/\cS}$ \eqref{p1-indmal20}. 
The category of coherent Higgs $\co_\fX[\frac{1}{p}]$-modules naturally embeds into the category of Higgs ind-$\co_\fX$-modules \eqref{p2-cmupiso31h}.
We will therefore consider any coherent Higgs $\co_\fX[\frac 1 p]$-module also as a Higgs ind-$\co_\fX$-module.

We denote by $\bIndHM(\bvocB,\hupsigma^*(\xi^{-1}\Omega^1_{\fX/\cS}))$ the category of Higgs ind-$\bvocB$-modules
with coefficients in $\hupsigma^*(\xi^{-1}\Omega^1_{\fX/\cS} )$. We can extend the functors \eqref{ft7a} to Higgs ind-modules; see \ref{p2-htaft23}.

\subsection{}\label{dolb1}
We sheafify the Higgs--Tate extension $\fF$ introduced in \eqref{rlps100a}, equipped with its natural $\hoR$-semilinear action of $\Delta$,  
into a canonical exact sequence of $\bvocB$-modules
\begin{equation}\label{dolb1a}
0\rightarrow \bvocB\rightarrow \bvcF\rightarrow \hupsigma^*(\xi^{-1}\Omega^1_{\fX/\cS})\rightarrow 0;
\end{equation}
see \ref{p2-htaft16} and \ref{p2-htaft18}. 
For any rational number $r\geq 0$, we denote by 
\begin{equation}\label{dolb1b}
0\rightarrow \bvocB\rightarrow \bvcF^{(r)}\rightarrow \hupsigma^*(\xi^{-1}\Omega^1_{\fX/\cS})\rightarrow 0. 
\end{equation}
the extension of $\bvocB$-modules deduced from $\bvcF$ by pullback by the morphism of multiplication by $p^r$ on 
$\hupsigma^*(\xi^{-1}\Omega^1_{\fX/\cS})$,
and by $\hfC^{(r)}$ the associated Higgs--Tate $\bvocB$-algebra:  
\begin{equation}\label{dolb1c}
\bvcC^{(r)}=\underset{\underset{m\geq 0}{\longrightarrow}}\lim\ \Sym^m_{\bvocB}(\bvcF^{ (r)}).
\end{equation}
We provide in \ref{p2-htaft18} alternative, but essentially equivalent, definitions of $\bvcF^{(r)}$ and $\bvcC^{(r)}$.
We denote by
\begin{equation}\label{dolb1d}
\delta_{\bvcC^{(r)}}\colon \bvcC^{(r)}\rightarrow \hupsigma^*(\xi^{-1}\Omega^1_{\fX/\cS})\otimes_{\bvocB}\bvcC^{(r)}
\end{equation}
the $p^r$ multiple of the universal $\bvocB$-derivation of $\bvcC^{(r)}$ \eqref{p2-htaft19d}. 
The algebras $\bvcC^{(r)}$, for $r\in \mQ_{\geq 0}$, form naturally a direct system, for which the derivations $\delta_{\bvcC^{(r)}}$ are compatible.  

We consider the ind-$\bvocB$-algebra \eqref{p1-indmal2}
\begin{equation}\label{dolb1e}
\IC^\dagger=\underset{\underset{r\in \mQ_{>0}}{\longrightarrow}}{\mlq\mlq\lim \mrq\mrq}\ \bvcC^{(r)},
\end{equation}
that we equip with the derivation \eqref{p1-indmal21}
\begin{equation}\label{dolb1f}
\delta_{\IC^\dagger}=\underset{\underset{r\in \mQ_{>0}}{\longrightarrow}}{\mlq\mlq\lim \mrq\mrq}\ \delta_{\bvcC^{(r)}}\colon 
\IC^\dagger \rightarrow \hupsigma^*(\xi^{-1}\Omega^1_{\fX/\cS})\otimes_{\bvocB}\IC^\dagger.
\end{equation}
The latter is a Higgs $\bvocB$-field \eqref{p1-indmal20}.

\subsection{}\label{dolb10}
Generalizing the construction of \eqref{thm3b}, we define 
the {\em twisting functor by the extension $\bvcF$ for Higgs ind-$\bvocB$-modules} \eqref{p2-rgpsc1c}
\begin{equation}\label{dolb10a}
\upnu\colon 
\begin{array}[t]{clcr}
\bIndHM(\bvocB, \hupsigma^*(\xi^{-1}\Omega^1_{\fX/\cS}))&\rightarrow &\bIndHM(\bvocB, \hupsigma^*(\xi^{-1}\Omega^1_{\fX/\cS})),\\
(\cN,\theta)&\mapsto& (\uupnu(\cN,\theta),\theta_\upnu). 
\end{array}
\end{equation}
We reserve the symbol $\uptau$ for a different twisting functor, which will be introduced later \eqref{p2-fhtft4a}. 
Similar to \ref{thm5}, under a suitable admissibility condition for Higgs ind-modules, 
we describe in \ref{p2-rgpsc31} the functor $\upnu$ as a twist by a line bundle.
We also consider the functor \eqref{p2-rgpsc1d}
\begin{equation}\label{dolb10b}
\fh\colon 
\begin{array}[t]{clcr}
\bIndHM(\bvocB, \hupsigma^*(\xi^{-1}\Omega^1_{\fX/\cS}))&\rightarrow &\bHM(\co_\fX,\xi^{-1}\Omega^1_{\fX/\cS}),\\
(\cN,\theta)&\mapsto& \vupsigma_*(\cN\otimes_{\bvocB} \IC^\dagger,\theta_\tot),
\end{array}
\end{equation}
where $\theta_\tot$ denotes the total Higgs field and the functor $\vupsigma_*$ is defined in \eqref{ft7b}.

\subsection{}\label{dolb11}
We consider the functor \eqref{p2-rgpsc14a}
\begin{equation}\label{dolb11a}
\cH\colon 
\begin{array}[t]{clcr}
\bIndMod(\bvocB)&\rightarrow &\bHM(\co_\fX,\txi^{-1}\tOmega^1_{\fX/\cS}),\\
\cM&\mapsto& \fh(\cM,0),
\end{array}
\end{equation}
where $\fh$ is defined in \eqref{dolb10b}, and the functor \eqref{p2-rgpsc14b}
\begin{equation}\label{dolb11b}
\cV\colon 
\begin{array}[t]{clcr}
\bHM^\coh(\co_\fX[\frac 1 p],\xi^{-1}\Omega^1_{\fX/\cS})&\rightarrow &\bIndMod(\bvocB),\\
(N,\theta)&\mapsto& \uupnu(\rI\hupsigma^*(N,\theta)),
\end{array}
\end{equation}
where $\uupnu$ is introduced in \eqref{dolb10a}.

\begin{defi}[cf. \ref{p2-rgpsc11}]\label{dolb2}
Let $\cM$ be an ind-$\bvocB$-module,
$(N,\theta)$ a Higgs $\co_\fX[\frac 1 p]$-{\em bundle} with coefficients in $\xi^{-1}\Omega^1_{\fX/\cS}$
(i.e., the $\co_\fX[\frac 1 p]$-module $N$ is locally projective of finite type \eqref{p1-delta-con6}), 
that we also consider as a Higgs ind-$\co_\fX$-module with coefficients in $\xi^{-1}\Omega^1_{\fX/\cS}$ \eqref{p2-cmupiso31h}.  
We say that $\cM$ and $(N,\theta)$ are {\em associated} if there exists a rational number $r>0$  and 
an isomorphism of ind-$\bvcC^{(r)}$-modules with integrable $\delta_{\bvcC^{(r)}}$-connection, in the sense of \ref{p1-indmal22},
\begin{equation}\label{dolb2a}
\bvcC^{(r)}\otimes_{\bvocB}\cM\stackrel{\sim}{\rightarrow} \bvcC^{(r)}\otimes_{\bvocB}\rI\hupsigma^*(N,\theta),
\end{equation}
where the $\delta_{\bvcC^{(r)}}$-connections are defined by the total Higgs fields, 
$N$ (resp.\ $\cM$) being endowed with the Higgs field $\theta$ (resp.\ $0$); see \ref{p1-delta-con9}. 
We then say that $(N,\theta)$ is {\em solvable} and that $\cM$ is {\em Dolbeault}.  
\end{defi}

\begin{teo}[cf. \ref{p2-rgpsc15}]\label{dolb12}
The functors $\cH$ \eqref{dolb11a} and $\cV$ \eqref{dolb11b} induce equivalences of categories quasi-inverse to each other
between the category of Dolbeault ind-$\bvocB$-modules and the category of solvable $\co_\fX[\frac 1 p]$-modules 
with coefficients in $\xi^{-1}\Omega^1_{\fX/\cS}$, 
\begin{equation}\label{dolb12a}
\xymatrix{
{\bIndMod^\Dolb(\bvocB)}\ar@<1ex>[r]^-(0.5){\cH}&{\bHM^\sol(\co_\fX[\frac 1 p], \xi^{-1}\Omega^1_{\fX/\cS}).}
\ar@<1ex>[l]^-(0.5){\cV}}
\end{equation}
\end{teo}

\begin{teo}[cf. \ref{p2-rgpsc16}]\label{dolb3}
For every Dolbeault ind-$\bvocB$-module $\cM$, we have a canonical functorial isomorphism of $\bD^+(\bMod(\co_\fX))$
\begin{equation}
\rR\vupsigma_*(\cM)\stackrel{\sim}{\rightarrow}\mK^\bullet(\cH(\cM)),
\end{equation}
where $\rR\vupsigma_*$ is the right derived functor of \eqref{ft7b} and $\mK^\bullet(\cH(\cM))$ is the Dolbeault complex of $\cH(\cM)$. 
\end{teo}

\begin{defi}[cf. \ref{p2-rgpsc20}]\
\begin{itemize}
\item[(i)] A $\bvocB_\mQ$-module is said to be {\em Dolbeault} if it is so as an ind-$\bvocB$-module \eqref{ft7}, i.e., if it is associated with
a Higgs $\co_\fX[\frac 1 p]$-bundle with coefficients in $\xi^{-1}\Omega^1_{\fX/\cS}$ \eqref{dolb2}.
\item[(ii)] A $\bvocB_\mQ$-module is said to be {\em strongly Dolbeault} if it is Dolbeault and adic of finite type (\cite{ag2} 4.3.14).
\item[(iii)] A Higgs $\co_\fX[\frac 1 p]$-bundle with coefficients in $\xi^{-1}\Omega^1_{\fX/\cS}$ \eqref{p1-delta-con6}
is said to be {\em rationally solvable} (resp.\ {\em strongly solvable})
if it is associated with a $\bvocB_\mQ$-module (resp.\ an adic $\bvocB_\mQ$-module of finite type).
\end{itemize}
\end{defi}

\begin{prop}[cf. \ref{p2-rgpsc240}]
Every strongly solvable Higgs $\co_\fX[\frac 1 p]$-bundle with coefficients in $\xi^{-1}\Omega^1_{\fX/\cS} $ is locally CL-small \eqref{thm66}.
\end{prop}

We give in \ref{p2-htaft58} and \ref{p2-rgpsc24} other relations between the various admissibility conditions and the CL-smallness. 

\subsection{}\label{dolb4}
Let $\cM$ be a Dolbeault ind-$\bvocB$-module, $(N,\theta)=\cH(\cM)$, 
that we also consider as a Higgs ind-$\co_\fX$-module with coefficients in $\xi^{-1}\Omega^1_{\fX/\cS}$.
By \ref{dolb12} and the definition of the functor $\cV$ \eqref{dolb11b}, we have a canonical isomorphism of ind-$\bvocB$-modules
\begin{equation}\label{dolb4a}
\cM\stackrel{\sim}{\rightarrow} \uupnu(\rI\hupsigma^*(N,\theta)),
\end{equation}
where the functor $\uupnu$ is defined in \eqref{dolb10a}. We deduce a Higgs $\bvocB$-field $\theta_{\cM}$ on $\cM$ 
with coefficients in $\hupsigma^*(\xi^{-1}\Omega^1_{\fX/\cS})$, such that 
\eqref{dolb4a} is underlying an isomorphism of Higgs ind-$\bvocB$-modules
\begin{equation}\label{dolb4b}
(\cM,\theta_{\cM})\stackrel{\sim}{\rightarrow} \upnu(\rI\hupsigma^*(N,\theta)). 
\end{equation}
We thus define a functor 
\begin{equation}\label{dolb4c}
\begin{array}[t]{clcr}
\bIndMod^\Dolb(\bvocB)&\rightarrow &\bIndHM(\bvocB,\hupsigma^*(\xi^{-1}\Omega^1_{\fX/\cS})),\\
\cM&\mapsto & (\cM,\theta_{\cM}).
\end{array}
\end{equation}

\begin{prop}[cf. \ref{p2-dpscc6}]
For every Dolbeault ind-$\bvocB$-module $\cM$ such that the $\co_{\fX}[\frac 1 p]$-bundle $\cH(\cM)$ is locally CL-small \eqref{thm66}, 
the canonical Higgs $\bvocB$-field $\theta_{\cM}$ \eqref{dolb4c} does not depend on the choice of the deformation $\tX$ \eqref{ft0a}. 
\end{prop}

\section[Functoriality of the $p$-adic Simpson correspondence]
{\texorpdfstring{Functoriality of the $p$-adic Simpson correspondence by pullback and proper higher direct images}
{Functoriality of the p-adic Simpson correspondence by pullback and proper higher direct images}}\label{fpscpdi}

\subsection{}\label{fpscpdi1} 
In this section, we again adopt the assumptions and notation of \ref{rlps0}. Let $\iota\colon \coS \to \tS$ 
be Fontaine's universal $p$-adic first-order thickening \eqref{rlps1a}, defined by the ideal $\xi \co_\tS$, 
which is free of rank one over $\co_\coS$ with basis $\xi$.
Let $g\colon X'\rightarrow X$ be a morphism of smooth $S$-schemes of finite type. 
We denote by $\fX$ (resp.\ $\fX'$) the formal scheme $p$-adic completion of $\coX=X\times_S\coS$ (resp.\ $\coX'=X'\times_S\coS$),
and by $\cog\colon \coX'\rightarrow \coX$ and $\fgg\colon \fX'\rightarrow \fX$ the morphisms induced by $g$. 
We set $\Omega=\xi^{-1}\Omega^1_{\coX/\coS}$ (resp.\ $\Omega'=\xi^{-1}\Omega^1_{\coX'/\coS}$,
resp.\ $\uOmega'=\xi^{-1}\Omega^1_{\coX'/\coX}$) and denote by $\hOmega$ (resp.\ $\hOmega'$, resp.\  $\huOmega'$) its $p$-adic completion, 
which is canonically isomorphic to the module $\xi^{-1}\Omega^1_{\fX/\cS}$ (resp.\ $\xi^{-1}\Omega^1_{\fX'/\cS}$, 
resp.\ $\xi^{-1}\Omega^1_{\fX'/\fX}$). 

By functoriality of the Faltings topos, $g$ induces a canonical morphism of topos $\Theta\colon \tE'\rightarrow \tE$, 
from the Faltings topos associated with $X'_\oeta\rightarrow X'$ to the Faltings topos associated with $X_\oeta\rightarrow X$ \eqref{ft1}, 
that fits into a commutative diagram
\begin{equation}\label{fpscpdi1a}
\xymatrix{
\tE'\ar[d]_{\sigma'}\ar[r]^{\Theta}&{\tE}\ar[d]^-(0.5){\sigma}\\
{X'_\et}\ar[r]^g&{X_\et,}}
\end{equation}
where $\sigma$ and $\sigma'$ are the canonical morphisms of topos \eqref{ft1d}. 
We have a canonical homomorphism of rings $\ocB\rightarrow \Theta_*(\ocB')$. 
Diagram \eqref{fpscpdi1a} induces by $p$-adic completion \eqref{ft6} a commutative diagram of morphisms of ringed topos 
\begin{equation}\label{fpscpdi1c}
\xymatrix{
{(\tE'^{\mN^\circ}_s,\bvocB')}\ar[d]_{\hupsigma'}\ar[r]^{\bvuptheta}&{(\tE^{\mN^\circ}_s,\bvocB)}\ar[d]^-(0.5){\hupsigma}\\
{(X'_{s,\zar},\co_{\fX'})}\ar[r]^{\fgg}&{(X_{s,\zar},\co_{\fX}).}}
\end{equation}
The morphism $\bvuptheta$ induces two adjoint functors 
\begin{equation}\label{fpscpdi1d}
\xymatrix{
{\bIndMod(\bvocB')}\ar@<1ex>[r]^-(0.5){\rI \bvuptheta_*}&{\bIndMod(\bvocB),}\ar@<1ex>[l]^-(0.5){\rI \bvuptheta^*}}
\end{equation}
which extend the adjoint functors $\bvuptheta^*$ and $\bvuptheta_*$, and similarly for $\hupsigma$ \eqref{ft7a}, $\hupsigma'$ and $\fgg$.

\subsection{}\label{fpscpdi2}
We suppose that there exists a commutative diagram
\begin{equation}\label{fpscpdi2a}
\xymatrix{
{\coX'}\ar[d]_{\cog}\ar[r]&{\tX'}\ar@/^1pc/[dd]\\
{\coX}\ar[r]\ar[d]\ar@{}[rd]|\Box&{\tX}\ar[d]\\
{\coS}\ar[r]^-(0.5)\iota&{\tS,}}
\end{equation}
where $\tX$ (resp.\ $\tX'$) are smooth $\tS$-deformations of $\coX$ (resp.\ $\coX'$) over $\coS$. 
We consider the torsor $\cL_{\tX'/\tX}$ of local liftings of $\cog$ over $\tX$ \eqref{torlift1},
the associated Higgs--Tate $\co_{\coX'}$-extension
\begin{equation}\label{fpscpdi2b}
0\rightarrow \co_{\coX'}\rightarrow \cF_\uptau\rightarrow \cog^*(\Omega)\rightarrow 0,
\end{equation}
and the associated Higgs--Tate $\co_{\coX'}$-algebra $\cC_\uptau$; see \ref{p2-fhtft3}. For any rational number $r\geq 0$, 
we denote by $\cC^{(r)}_\uptau$ the Higgs--Tate $\co_{\coX'}$-algebra of thickness $r$, equipped with 
the $p^r$-multiple of its universal $\co_{\coX'}$-derivation
\begin{equation}\label{fpscpdi2i}
\delta_{\cC^{(r)}_\uptau}\colon \cC^{(r)}_\uptau \rightarrow
\cog^*(\Omega)\otimes_{\co_{\coX'}}\cC^{(r)}_\uptau.
\end{equation}
We also consider the $\co_{\fX'}$-derivation
\begin{equation}\label{fpscpdi2e}
\delta'_{\cC^{(r)}_\uptau}=(u\otimes \id)\circ \delta_{\cC^{(r)}_\uptau}\colon \cC^{(r)}_\uptau \rightarrow \Omega'\otimes_{\co_{\coX'}}\cC^{(r)}_\uptau
\end{equation}
where $u\colon \cog^*(\Omega)\rightarrow \Omega'$ is the canonical morphism; see \ref{p2-fhtft6}.

We consider the ind-$\co_{\fX'}$-algebra \eqref{p2-fhtft6k}
\begin{equation}\label{fpscpdi2f}
\IC^\dagger_{\uptau}=\underset{\underset{r\in \mQ_{>0}}{\longrightarrow}}{\mlq\mlq\lim \mrq\mrq}\ \hcC^{(r)}_{\uptau}, 
\end{equation}
where $\hcC^{(r)}_{\uptau}$ is the $p$-adic completion of $\cC^{(r)}_{\uptau}$, 
equipped with the $\co_{\fX'}$-derivation \eqref{p2-fhtft6i}
\begin{equation}\label{fpscpdi2j}
\Idelta_{\uptau}=\underset{\underset{r\in \mQ_{>0}}{\longrightarrow}}{\mlq\mlq\lim \mrq\mrq}\ \delta_{\hcC^{(r)}_\uptau}\colon \IC^\dagger_\uptau \rightarrow 
\fgg^*(\hOmega)\otimes_{\co_{\fX'}}\IC^\dagger_\uptau,
\end{equation}
where $\delta_{\hcC^{(r)}_\uptau}$ is the extension of $\delta_{\cC^{(r)}_\uptau}$ \eqref{fpscpdi2i} to the $p$-adic completions. 
We also consider the $\co_{\fX'}$-derivation
\begin{equation}\label{fpscpdi2g}
\Idelta'_{\uptau}=(\hu\otimes \id)\circ \Idelta_{\uptau}\colon \IC^\dagger_\uptau \rightarrow 
\hOmega'\otimes_{\co_{\fX'}}\IC^\dagger_\uptau,
\end{equation}
where $\hu\colon \fgg^*(\hOmega)\rightarrow \hOmega'$ is the canonical morphism.

We denote by $\hcF_\uptau$ the $p$-adic completion of $\cF_\uptau$, which fits into an exact sequence of $\co_{\fX'}$-modules
\begin{equation}\label{fpscpdi2bb}
0\rightarrow \co_{\fX'}\rightarrow \hcF_\uptau\rightarrow \fgg^*(\hOmega)\rightarrow 0.
\end{equation}
We consider for Higgs modules the pullback functor by $\fgg$ twisted by the extension $\hcF_\uptau$, \eqref{tphdi3c},
\begin{equation}\label{pscpdi2d}
\fgg^*_\uptau\colon 
\bHM(\co_\fX,\hOmega)\rightarrow \bHM(\co_{\fX'},\fgg^*(\hOmega)),
\end{equation} 
and for any integer $q\geq 0$, 
the $q$th higher direct image functor by $\fgg$ twisted by the extension $\hcF_\uptau$, \eqref{tphdi4d},
\begin{equation}\label{pscpdi2c}
\rR^q\fgg^\uptau_*\colon 
\bHM(\co_{\fX'},\hOmega')\rightarrow\bHM(\co_\fX,\hOmega).
\end{equation}
We denote by
\begin{equation}\label{pscpdi2h}
\upmu\colon \bHM(\co_{\fX'},\fgg^*(\hOmega)) \rightarrow \bHM(\co_{\fX'},\hOmega')
\end{equation} 
the functor induced by the canonical morphism $\hu\colon \fgg^*(\hOmega)\rightarrow \hOmega'$. 

\subsection{}\label{fpscpdi30}
We consider again the objects associated with $(X,\tX)$ introduced in §\ref{ft}, and we associate with $(X',\tX')$ 
similar objects that we denote by the same symbols equipped with a $^\prime$ exponent.
In particular, we denote by $\bIndMod^\Dolb(\bvocB')$ the category of Dolbeault ind-$\bvocB'$-modules 
and by $\bHM^\sol(\co_{\fX'}[\frac 1 p], \hOmega')$
the category of solvable Higgs $\co_{\fX'}[\frac 1 p]$-modules with coefficients in $\hOmega'$ \eqref{dolb2}. We denote by
\begin{equation}\label{fpscpdi30a}
\cH'\colon \bIndMod(\bvocB')\rightarrow \bHM(\co_{\fX'}, \hOmega')
\end{equation}
the functor defined in \eqref{dolb11a}, associated with $(X',\tX')$.
By \ref{dolb12}, this induces an equivalence of categories that we denote again by
\begin{equation}\label{fpscpdi30b}
\cH'\colon \bIndMod^\Dolb(\bvocB')\stackrel{\sim}{\rightarrow} \bHM^\sol(\co_{\fX'}[\frac 1 p], \hOmega').
\end{equation}
 
\begin{teo}[cf. \ref{p2-fgpscp2}]\label{fpscpdi3}
For every Dolbeault ind-$\bvocB$-module $\cM$ \eqref{dolb2} such that the associated 
$\co_{\fX}[\frac 1 p]$-module $\cH(\cM)$ \eqref{dolb11a} is locally CL-small \eqref{thm66}, 
the ind-$\bvocB'$-module $\rI\bvuptheta^*(\cM)$ \eqref{fpscpdi1d} is Dolbeault, 
its associated $\co_{\fX'}[\frac 1 p]$-bundle $\cH'(\rI\bvuptheta^*(\cM))$ \eqref{fpscpdi30a} is locally CL-small, 
and we have a canonical functorial isomorphism 
\begin{equation}\label{p2-fgpscp3a}
\cH'(\rI\bvuptheta^*(\cM))\stackrel{\sim}{\rightarrow}\upmu(\fgg^*_\uptau(\cH(\cM))),
\end{equation}
where the functors $\fgg^*_\uptau$ and $\upmu$ are  defined in \eqref{pscpdi2d} and \eqref{pscpdi2h}, respectively. 
\end{teo}

The proof relies on a composition homomorphism between twisted Higgs--Tate algebras, introduced below in \eqref{fpscpdi6a}.

\begin{teo}[cf. \ref{p2-rcdim10}]\label{fpscpdi4}
Let $\cM$ be an ind-$\bvocB'$-module, $(N,\theta)$ a Higgs $\co_{\fX'}[\frac 1 p]$-bundle with coefficients in $\hOmega'$ \eqref{p1-delta-con6},
that we consider as a Higgs ind-$\co_{\fX'}$-module \eqref{ft7}, $q$ an integer. 
We assume that the following conditions are satisfied: 
\begin{itemize}
\item[{\rm (a)}] $g\colon X'\rightarrow X$ is proper and smooth; 
\item[{\rm (b)}] $\cM$ and $(N,\theta)$ are associated in the sense of \ref{dolb2};
\item[{\rm (c)}] $(N,\theta)$ is locally CL-small in the sense of \ref{thm66}. 
\end{itemize}
Then, 
\begin{itemize}
\item[{\rm (i)}] The Higgs $\co_\fX[\frac 1 p]$-module $\rR^q\fgg^\uptau_*(N,\theta)$ \eqref{pscpdi2c} 
is (locally) CL-small and in particular coherent \eqref{ft7}. 
We consider it as a Higgs ind-$\co_{\fX}$-module. 
\item[{\rm (ii)}] There exist a rational number $r>0$, independent of $q$, and an isomorphism of ind-$\bvcC^{(r)}$-modules with 
$\delta_{\bvcC^{(r)}}$-connection, in the sense of \ref{p1-indmal22}, 
\begin{equation}\label{fpscpdi4a}
\bvcC^{(r)}\otimes_{\bvocB}\rR^q\rI\bvuptheta_*(\cM)\stackrel{\sim}{\rightarrow}
\bvcC^{(r)}\otimes_{\bvocB}\rI\hupsigma^*(\rR^q\fgg^\uptau_*(N,\theta)),
\end{equation}
where the $\delta_{\bvcC^{(r)}}$-connections are defined by the total Higgs fields, 
$\rR^q\rI\bvuptheta_*(\cM)$ \eqref{fpscpdi1d} being endowed with the Higgs field $0$.  
\end{itemize}
\end{teo}

\begin{cor}[cf. \ref{p2-rcdim14}]\label{fpscpdi5}
Let $\cM$ be a Dolbeault ind-$\bvocB'$-module \eqref{dolb2} such that the Higgs $\co_{\fX'}[\frac 1 p]$-module $\cH'(\cM)$ \eqref{fpscpdi30a}
is locally CL-small \eqref{thm66}, $q$ an integer.
Suppose that the morphism $g\colon X'\rightarrow X$ is proper and smooth and the $\co_\fX[\frac 1 p]$-module underlying the 
Higgs $\co_\fX[\frac 1 p]$-module $\rR^q\fgg^\uptau_*(\cH'(\cM))$ is locally projective of finite type.
Then, the ind-$\bvocB$-module $\rR^q\rI\bvuptheta_*(\cM)$ \eqref{fpscpdi1d} is Dolbeault, the Higgs $\co_\fX[\frac 1 p]$-module 
$\rR^q\fgg^\uptau_*(\cH'(\cM))$ \eqref{pscpdi2c} is solvable and (locally) CL-small, 
and we have a canonical functorial isomorphism of Higgs $\co_\fX[\frac 1 p]$-bundles
\begin{equation}\label{fpscpdi5a}
\cH(\rR^q\rI\bvuptheta_*(\cM))\stackrel{\sim}{\rightarrow} \rR^q\fgg^\uptau_*(\cH'(\cM)),
\end{equation}
where $\cH$ is the functor \eqref{dolb11a}. 
\end{cor}

\subsection{}\label{fpscpdi40}
For the remainder of this section, we assume that $g\colon X'\rightarrow X$ is smooth.
To prove theorem \ref{fpscpdi4}, we need to introduce a {\em relative Faltings topos} $\tG$. The quickest definition is as
the fiber product of $\sigma$ and $g$ in diagram \eqref{fpscpdi1a}, so we have a commutative diagram of morphisms of topos with a Cartesian square
\begin{equation}\label{fpscpdi40a}
\xymatrix{
\tE'\ar[dr]_{\sigma'}\ar[r]^{\varpi}\ar@/^2pc/[rr]^{\Theta}&{\tG}\ar[r]^-(0.4){\pi}\ar[d]_{\lgg}\ar@{}[rd]|\Box&{\tE}\ar[d]^-(0.5){\sigma}\\
&{X'_\et}\ar[r]^g&{X_\et,}}
\end{equation}
where $\varpi$, $\pi$ and $\lgg$ are the canonical morphisms. 
However, $\tG$ can also be defined via an explicit site similar to Faltings's site, which justifies the terminology {\em relative Faltings topos}; see \ref{p2-rftchta1} and (\cite{ag1}, 3.4.19). It is this description that we use to study this topos.
We set
\begin{equation}\label{fpscpdi40b}
\ocB^!=\varpi_*(\ocB').
\end{equation}

We denote by $(\tG^{\mN^\circ}_s,\bvocB^!)$ the $p$-adic completion of the ringed topos $(\tG,\ocB^!)$; see \ref{ft6} and \ref{p2-rftchta4}. 
Diagram \eqref{fpscpdi40a} induces a commutative diagram of morphisms of ringed topos 
\begin{equation}\label{fpscpdi40c}
\xymatrix{
{(\tE'^{\mN^\circ}_s,\bvocB')}\ar[r]^-(0.5){\bvvarpi}\ar[rd]_{\hupsigma'}\ar@/^2pc/[rr]^{\bvuptheta}&{(\tG^{\mN^\circ}_s,\bvocB^!)}
\ar[d]^-(0.5){\huppi}\ar[r]^-(0.5){\bvlgg}&{(\tE^{\mN^\circ}_s,\bvocB)}\ar[d]^-(0.5){\hupsigma}\\
&{(X'_{s,\zar},\co_{\fX'})}\ar[r]^{\fgg}&{(X_{s,\zar},\co_{\fX}).}}
\end{equation}
The morphism $\huppi$ induces two adjoint functors 
\begin{equation}\label{fpscpdi40d}
\xymatrix{
{\bIndMod(\bvocB^!)}\ar@<1ex>[r]^-(0.5){\rI \huppi_*}&{\bIndMod(\co_{\fX'}),}\ar@<1ex>[l]^-(0.5){\rI \huppi^*}}
\end{equation}
which extend the adjoint functors $\huppi^*$ and $\huppi_*$, and similarly for $\bvvarpi$ and $\bvlgg$.  

The proof of theorem \ref{fpscpdi4} relies on two main ingredients:
\begin{itemize}
\item[(i)] A computation of the Dolbeault cohomology twisted by the extension $\hcF_\uptau$ \eqref{fpscpdi2bb} 
and of the higher direct images under $\bvvarpi$ 
of the Higgs--Tate algebras $\bvcC'^{(r')}$ relative to $(X',\tX')$, expressed in terms of the Higgs--Tate algebras $\bvcC^{(r)}$ relative to $(X,\tX)$.
This computation uses the composition homomorphism between twisted Higgs--Tate algebras introduced below in~\eqref{fpscpdi6a}; see \ref{fpscpdi7}
and \ref{fpscpdi8}. 
\item[(ii)] A base change theorem for twisted Dolbeault complexes with respect to the square in diagram~\eqref{fpscpdi40c}; see \ref{fpscpdi10}
\end{itemize}

\subsection{}\label{fpscpdi6}
Let $r,r'$ be rational numbers such that $r\geq r'\geq 0$.  
The sheaves introduced in the affine case in \eqref{fspb4c} can be extended naturally to the general case. 
We can sheafify the $\hoRp$-algebra $\fC'^{(r,r')}$ introduced in \eqref{fspb4e} into a $\bvocB'$-algebra $\bvcC'^{(r,r')}$ of $\tE'^{\mN^\circ}_s$; 
see \ref{p2-fhtft28}. It is equipped with a normalized $\bvocB'$-derivation \eqref{p2-fhtft29b}
\begin{equation}\label{fpscpdi6b}
\delta_{\bvcC'^{(r,r')}}\colon \bvcC'^{(r,r')}\rightarrow \hupsigma'^*(\hOmega')\otimes_{\bvocB'}\bvcC'^{(r,r')},
\end{equation}
which is also a Higgs $\bvocB'$-field on $\bvcC'^{(r,r')}$. 

Let $I$ be the set of triples of rational numbers $\ur=(r_1,r_2,r_3)$ such that $r_1\geq r_2\geq r_3\geq 0$. 
For such a triple $\ur$, we set 
\begin{equation}
\bvcC'^\ur=\bvcC'^{(r_1,r_2)} \otimes_{\bvocB'}\hupsigma'^*(\hcC^{(r_3)}_\uptau).
\end{equation}
The homomorphism \eqref{fspb4f} can be sheafified into a homomorphism of $\bvocB'$-algebras
\begin{equation}\label{fpscpdi6a}
\bvupphi^\ur\colon \bvuptheta^*(\bvcC^{(r_2)})\rightarrow \bvcC'^\ur;
\end{equation}
see \ref{p2-fhtft38}. 
The construction of this homomorphism of Faltings topos is quite delicate; see \ref{p2-fhtft35} and \ref{p2-fhtft362}. 
It relies on a careful study of the functorial properties of 
categories of presheaves of sets on a category, together with a new notion of {\em quasi-fibered categories}, 
refining Grothendieck's notion of fibered categories; see §~\ref{p2-cmt} and §~\ref{p2-qfc}. 

We equip $\bvcC'^\ur$ with the $\bvocB'$-derivation
\begin{equation}\label{fpscpdi6c}
\delta^\ur\colon \bvcC'^{(r_1,r_2)}\otimes_{\bvocB'}\hupsigma'^*(\hcC^{(r_3)}_\uptau)
\rightarrow \hupsigma'^*(\hOmega')\otimes_{\bvocB'}\bvcC'^\ur
\end{equation}
defined by 
\begin{equation}\label{fpscpdi6d}
\delta^\ur=\delta_{\bvcC'^{(r_1,r_2)}}\otimes \id- \id \otimes \hupsigma'^*(\delta'_{\hcC^{(r_3)}_\uptau}),
\end{equation}
where $\delta_{\bvcC'^{(r_1,r_2)}}$ is defined in \eqref{fpscpdi6b} and $\delta'_{\hcC^{(r_3)}_\uptau}$ is induced by \eqref{fpscpdi2e}. 
We denote by $\mK^\bullet(\bvcC'^{\ur})$ the Dolbeault complex of $(\bvcC'^{\ur},-\bvdelta^\ur)$, and by $\mK^\bullet_\mQ(\bvcC'^{\ur})$ its image 
in the category $\bMod_{\mQ}(\bvocB')$ of $\bvocB'$-modules up to isogeny; see \ref{p2-rftchta17}.  

\begin{prop}[cf. \ref{p2-rftchta20}]\label{fpscpdi7}
Let $\ut=(t_1,t_2,t_3)$ be a triple of rational numbers such that 
$t_1\geq t_2 > t_3 \geq 0$, $I_\ut$ the subset of elements $\ur=(r_1,r_2,r_3)$ of $I$ \eqref{fpscpdi6} 
such that $r_1>t_1$, $r_2\geq t_2$  and $r_3>t_3$. Then, the morphism of complexes of ind-$\bvocB'$-modules
\begin{equation}\label{fpscpdi7a}
\bvuptheta^*(\bvcC^{(t_2)})_\mQ[0] \rightarrow  
\underset{\underset{\ur\in I_\ut}{\longrightarrow}}{\mlq\mlq\lim \mrq\mrq} \ 
\mK^\bullet_\mQ(\bvcC'^{\ur})
\end{equation}
induced by $\bvupphi^{\ur}$ \eqref{fpscpdi6a}, is a quasi-isomorphism.
\end{prop}

We deduce this proposition from analogous local statements; see \ref{p2-fhtal16} and \ref{p2-fhtal17}. 

\begin{cor}[cf. \ref{p2-rcdim7}]\label{fpscpdi77}
Let $\cM$ be a rational and flat ind-$\bvocB'$-module, $\ut=(t_1,t_2,t_3)$ a triple of rational numbers such that $t_1\geq t_2 > t_3 \geq 0$, 
$I_\ut$ the subset of elements $\ur=(r_1,r_2,r_3)$ of $I$ \eqref{fpscpdi6} 
such that $r_1>t_1$, $r_2\geq t_2$  and $r_3>t_3$. Then, the morphism of complexes of ind-$\bvocB'$-modules
\begin{equation}\label{fpscpdi77a}
\cM\otimes_{\bvocB'}\bvuptheta^*(\bvcC^{(t_2)})[0] \rightarrow  
\underset{\underset{\ur\in I_\ut}{\longrightarrow}}{\mlq\mlq\lim \mrq\mrq} \ 
\mK^\bullet(\cM\otimes_{\bvocB'}\bvcC'^{\ur})
\end{equation}
induced by $\bvupphi^{\ur}$ \eqref{fpscpdi6a}, is a quasi-isomorphism.
\end{cor}

An ind-$\bvocB'$-module $\cM$ is said to be {\em rational} if the multiplication by $p$ on $\cM$ is an isomorphism \eqref{p2-htaft26}. 
It is said to be {\em flat} if the tensor product $-\otimes_{\bvocB'}\cM$ is an exact endofunctor of the category of ind-$\bvocB'$-modules 
(\cite{ag2} 2.7.9). 

\begin{prop}[cf. \ref{p2-rftchta8}]\label{fpscpdi8}
Let $t_2,t_3$ be two rational numbers such that $t_2\geq t_3 \geq 0$, $I_{t_2,t_3}$ the subset of elements $\ur=(r_1,r_2,r_3)$ of $I$ \eqref{fpscpdi6} 
such that $r_1>r_2\geq t_2$ and $r_3\geq t_3$. Then, 
\begin{itemize}
\item[{\rm (i)}] With the notation of \eqref{fpscpdi40c}, the homomorphism \eqref{fpscpdi6a} induces an isomorphism of $\bvocB^!_\mQ$-algebras
\begin{equation}\label{fpscpdi8a}
(\bvlgg^*(\bvcC^{(t_2)})\otimes_{\bvocB^!}\huppi^*(\hcC^{(t_3)}_\uptau))_\mQ\stackrel{\sim}{\rightarrow} 
\underset{\underset{\ur\in I_{t_2,t_3}}{\longrightarrow}}{\lim}\ \bvvarpi_*(\bvcC'^{\ur})_\mQ; 
\end{equation}
and for every integer $q\geq 1$, we have 
\begin{equation}\label{fpscpdi8b}
\underset{\underset{\ur\in I_{t_2,t_3}}{\longrightarrow}}{\lim}\ \rR^q\bvvarpi_*(\bvcC'^{\ur})_\mQ= 0,
\end{equation}
where the limits are taken in the category $\bMod_\mQ(\bvocB^!)$ of $\bvocB^!$-modules up to isogeny. These are, in particular, representable.

\item[{\rm (ii)}] The homomorphism \eqref{fpscpdi6a} induces an isomorphism of ind-$\bvocB^!$-algebras
\begin{equation}\label{fpscpdi8c}
\bvlgg^*(\bvcC^{(t_2)})\otimes_{\bvocB^!}\rI\huppi^*(\hcC^{(t_3)}_{\uptau,\mQ})\rightarrow 
\underset{\underset{\ur\in I_{t_2,t_3}}{\longrightarrow}}{\mlq\mlq\lim\mrq\mrq}\ \rI \bvvarpi_*(\bvcC'^{\ur}_\mQ),
\end{equation}
where $\hcC^{(t_3)}_{\uptau,\mQ}$ and $\bvcC'^{\ur}_\mQ$ are considered naturally as ind-algebras
(see \ref{p1-bcim5}); and for every integer $q\geq 1$, we have 
\begin{equation}\label{fpscpdi8d}
\underset{\underset{\ur\in I_{t_2,t_3}}{\longrightarrow}}{\mlq\mlq\lim\mrq\mrq}\ \rR^q\bvvarpi_*(\bvcC'^{\ur})_\mQ= 0.
\end{equation} 
\end{itemize}
\end{prop}

We deduce this proposition from the analogous local statement \ref{p2-fhtal13}. 

\begin{cor}[cf. \ref{p2-rcdim9}]\label{fpscpdi9}
Let $(N,\theta)$ be a Higgs $\co_{\fX'}[\frac 1 p]$-bundle with coefficients in $\hOmega'$, 
that we consider as a Higgs ind-$\co_{\fX'}$-module, 
$t_2$ a rational number $>0$, $I_{t_2}$ the subset of elements $\ur=(r_1,r_2,r_3)$ of $I$ \eqref{fpscpdi6} 
such that $r_1>r_2\geq t_2$ and $r_3>0$. Then, we have a canonical isomorphism of $\bD^+(\bIndMod(\bvlgg^*(\bvcC^{(t_2)})))$
\begin{equation}\label{fpscpdi9a}
\bvlgg^*(\bvcC^{(t_2)})\otimes_{\bvocB^!}\rI\huppi^*(\mK^\bullet(N\otimes_{\co_{\fX'}}\IC^\dagger_{\uptau}))\stackrel{\sim}{\rightarrow}
\rR\rI\bvvarpi_*(\underset{\underset{\ur\in I_{t_2}}{\longrightarrow}} {\mlq\mlq\lim\mrq\mrq} \ \mK^\bullet(\rI\hupsigma'^*(N)\otimes_{\bvocB'}\bvcC'^{\ur})),
\end{equation}
where $\mK^\bullet(\rI\hupsigma'^*(N)\otimes_{\co_{\fX'}}\bvcC'^{\ur})$ is 
the Dolbeault complex of $(\rI\hupsigma'^*(N)\otimes_{\co_{\fX'}}\bvcC'^{\ur},\rI\hupsigma'^*(\theta)\otimes\id-\id\otimes \delta^{\ur})$ \eqref{fpscpdi6c}
and $\mK^\bullet(N\otimes_{\co_{\fX'}}\IC^\dagger_{\uptau})$ is 
the Dolbeault complex of $(N\otimes_{\co_{\fX'}}\IC^\dagger_{\uptau},\theta\otimes\id+\id\otimes \Idelta'_{\uptau})$ \eqref{fpscpdi2g}.
\end{cor}

\subsection{}
Let $\cM$ be an ind-$\bvocB'$-module and $(N,\theta)$ a Higgs $\co_{\fX'}[\frac 1 p]$-bundle with coefficients in $\hOmega'$,
which are associated in the sense of \ref{dolb2}. Starting from an admissibility isomorphism,
\ref{fpscpdi77} and \ref{fpscpdi9} imply the existence of a rational number $r>0$ together with an isomorphism of 
$\bD^+(\bIndMod(\bvlgg^*(\bvcC^{(r)})))$ 
\begin{equation}
\rR\rI\bvvarpi_*(\bvuptheta^*(\bvcC^{(r)})\otimes_{\bvocB'}\cM)
\stackrel{\sim}{\rightarrow}
\bvlgg^*(\bvcC^{(r)})\otimes_{\bvocB^!}\rI\huppi^*(\mK^{\bullet}(\cN\otimes_{\co_{\fX'}}\IC^\dagger_\uptau)),
\end{equation}
where on the right, the tensor product and $\rI\huppi^*$ are defined term by term; see \ref{p2-rcdim8}. 
We prove that this isomorphism is compatible with the Higgs $\bvocB^!$-fields with coefficients in $\huppi^*(\fgg^*(\hOmega))$ \eqref{fpscpdi40c},
induced by the derivations \eqref{dolb1d} of $\bvcC^{(r)}$ and \eqref{fpscpdi2j} of $\IC^\dagger_\uptau$. 
Applying the functor $\rR^q\rI\bvlgg_*$ to the above isomorphism, 
and using the projection formula for $\rR\rI\bvuptheta_*$ and $\rR\rI\bvlgg_*$ \eqref{p2-rcdim11},
the final step reduces to a base change theorem, which we explain below.

\subsection{}\label{fpscpdi11}
Let $q$ be an integer. For every bounded from below complex of $\co_{\fX'}$-modules $\cM^\bullet$, 
there exists a canonical functorial base change morphism of $\bvocB$-modules, with respect to the square in \eqref{fpscpdi40c},
\begin{equation}\label{fpscpdi11a}
\hupsigma^*(\rR^q\fgg_*(\cM^\bullet))\rightarrow \rR^q\bvlgg_*(\huppi^*(\cM^\bullet)),
\end{equation}
where the pullback $\huppi^*(\cM^\bullet)$ is defined term by term (not derived); see \eqref{p1-bcim4d} and \eqref{p2-bch7b}. 
For every $\co_{\fX'}$-module $\cM$, the base change morphism \eqref{fpscpdi11a} 
for $\cM[0]$ coincides with the classical base change morphism (\cite{egr1} 1.2.3).

For every bounded from below complex of ind-$\co_{\fX'}$-modules $\cF^\bullet$, 
there exists a canonical functorial base change morphism of ind-$\bvocB$-modules, with respect to the square in \eqref{fpscpdi40c},
\begin{equation}\label{fpscpdi11c}
\rI\hupsigma^*(\rR^q\rI \fgg_*(\cF^\bullet))\rightarrow \rR^q\rI\bvlgg_*(\rI\huppi^*(\cF^\bullet)),
\end{equation}
where the pullback $\rI\huppi^*(\cF^\bullet)$ \eqref{fpscpdi40d} is defined term by term (not derived); see \eqref{p1-bcim4f} and \eqref{p2-bch7d}.

\begin{teo}[cf. \ref{p2-bch12}]\label{fpscpdi10}
We assume that the morphism $g\colon X'\rightarrow X$ is proper and smooth \eqref{fpscpdi40}. 
Let $(N,\theta)$ be a {\em locally CL-small} Higgs $\co_{\fX'}[\frac 1 p]$-module with coefficients in $\hOmega'$, 
such that the $\co_{\fX'}[\frac 1 p]$-module $N$ is flat, $q$ an integer. 
We consider $(N,\theta)$ as a Higgs ind-$\co_{\fX'}$-module 
and denote by $\mK^{\bullet}(N\otimes_{\co_{\fX'}}\IC^\dagger_\uptau)$ the Dolbeault complex of $N\otimes_{\co_{\fX'}}\IC^\dagger_\uptau$ 
equipped with the total Higgs field $\theta\otimes \id+\id\otimes \Idelta'_\uptau$ \eqref{fpscpdi2g}. 
Then, the base change morphism \eqref{fpscpdi11c}
\begin{equation}\label{fpscpdi10a}
\rI\hupsigma^*(\rR^q\rI \fgg_*(\mK^\bullet(N\otimes_{\co_{\fX'}}\IC^\dagger_\uptau)))\rightarrow \rR^q\rI\bvlgg_*(\rI\huppi^*(\mK^\bullet(N\otimes_{\co_{\fX'}}\IC^\dagger_\uptau)))
\end{equation}
is an isomorphism of ind-$\bvocB$-modules. 
\end{teo}

The proof of this base change theorem relies on a fine study of the cohomology of ind-modules; see \ref{p2-bch11} and \ref{p2-cmupiso28}. 
The square in diagram \eqref{fpscpdi40c} naturally decomposes into 3 commutative squares 
\begin{equation}
\xymatrix{
{(\tG^{\mN^\circ}_s,\bvocB^!)}\ar[r]^-(0.5){\bvpi}\ar[d]_-(0.5){\bvlgg}&{(X'^{\mN^\circ}_{s,\et},\co_{\bvoX'})}\ar[r]\ar[d]_ {\bvogg}
&{(X'^{\mN^\circ}_{s,\zar},\co_{\bvoX'})}\ar[r]^-(0.5){\uplambda'}\ar[d] _{\bvogg}&{(X'_{s,\zar},\co_{\fX'})}\ar[d]^{\fgg}\\
{(\tE^{\mN^\circ}_s,\bvocB)}\ar[r]^-(0.5){\bvsigma}&{(X^{\mN^\circ}_{s,\et},\co_{\bvoX})}\ar[r]&
{(X^{\mN^\circ}_{s,\zar},\co_{\bvoX})}\ar[r]^-(0.5)\uplambda&{(X_{s,\zar},\co_{\fX}),}}
\end{equation}
where $\co_{\bvoX}$ (resp.\ $\co_{\bvoX'}$) is the ring $(\co_{\oX_{n}})_{n\in \mN}$ (resp.\ $(\co_{\oX'_{n}})_{n\in \mN}$), 
$\uplambda$ (resp.\ $\uplambda'$) is the morphism of ringed topos whose 
direct image functor associates with any inverse system its inverse limit \eqref{p2-ncgt5a}.
The most delicate study is devoted to the right square; see \ref{p2-cmupiso17}. 
For non-twisted Dolbeault cohomology, the analogous result is much simpler; see (\cite{ag1} 6.5.38) and (\cite{ag2} 6.2.14).

\section{Revisiting Sen theory via twisting}\label{rstvt}

\subsection{}\label{rstvt1}
We take again the assumptions and notation of \ref{overview-intro2}. 
We give in this section a new construction of Sen endomorphism \eqref{overview-intro4} by the twisting approach developed in this book. 

Let $\cF_{\mZ_p}=\mZ_p\oplus \mZ_p$ be the free $\mZ_p$-module of rank two equipped with the $p$-adically continuous representation of $G_K$ 
\begin{equation}\label{rstvt1a}
\begin{pmatrix}1&\log(\chi)\\ 0&1\end{pmatrix}
\end{equation}
We consider the exact sequence of $\mZ_p[G_K]$-modules 
\begin{equation}\label{rstvt1b}
0\rightarrow \mZ_p\rightarrow \cF_{\mZ_p}\rightarrow \mZ_p\rightarrow 0,
\end{equation}
where the second (resp.\ third) morphism is defined by $a\mapsto (a,0)$ (resp.\ $(a,b)\mapsto b$), 
that represents the character $\log(\chi)$ in $\rH^1_\cont(G_K,\mZ_p)=\Hom_\cont(G_K,\mZ_p)$.

We set $\cF=\cF_{\mZ_p}\otimes_{\mZ_p}\co_C$, equipped with the $\co_C$-semilinear action of $G_K$ induced by the action of $G_K$ on $\cF_{\mZ_p}$. 
We consider the exact sequence 
\begin{equation}\label{rstvt1c}
0\rightarrow \co_C\rightarrow \cF\rightarrow \co_C \rightarrow 0
\end{equation}
induced by \eqref{rstvt1b}. We denote by $\cC$ the associated Higgs--Tate algebra 
\begin{equation}\label{rstvt1d}
\cC=\underset{\underset{n\geq 0}{\longrightarrow}}\lim\ \Sym^n_{\co_C}(\cF),
\end{equation}
equipped with the action of $G_K$ compatible with its ring structure and induced by the action on $\cF$.  
We set $T=(0,1)\in \cF\subset \cC$. The splitting of \eqref{rstvt1c} defined by $b\mapsto bT$ induces 
a $G_K$-equivariant isomorphism of $\co_C$-algebras
\begin{equation}\label{rstvt1e}
\co_C[T]\stackrel{\sim}{\rightarrow} \cC, 
\end{equation}
where for every $\sigma\in G_K$, $\sigma(T)=\log(\chi(\sigma))+T$. 
The same splitting defines a homomorphism of $\co_C$-algebras $\alpha_T\colon \cC\rightarrow \co_C$; it corresponds to the section of $\co_C[T]$
defined by $T=0$. 
We equip $\cC$ with the $\co_C$-derivation 
\begin{equation}\label{rstvt1f}
\delta\colon \cC\rightarrow \cC
\end{equation}
defined by $\delta(T)=-1$. It is clearly $G_K$-equivariant. 

We denote by $\hcC$ the $p$-adic Hausdorff completion of $\cC$ and by 
\begin{eqnarray}
\hdelta\colon \hcC&\rightarrow &\hcC,\label{rstvt1g}\\
\halpha_T\colon \hcC&\rightarrow &\co_C, \label{rstvt1h}
\end{eqnarray}
the extensions of $\delta$ and $\alpha_T$, respectively, to the completions. 

\subsection{}\label{rstvt2}
With the terminology and notation introduced in \ref{overview-intro40}, 
let $(V,\sigma)$ be a Sen $C$-module. We denote by $\uuptau(V,\sigma)$ the kernel of the $C$-linear map 
\begin{equation}\label{rstvt2a}
\sigma_\tot=\hdelta\otimes \id+\id\otimes \sigma\colon  \hcC\otimes_{\co_C}V\rightarrow \hcC\otimes_{\co_C}V.
\end{equation}
The endomorphism $\id\otimes \sigma$ on $\hcC\otimes_{\co_C}V$ induces a $C$-linear endomorphism $\sigma_\uptau$ of $\uuptau(V,\sigma)$. 
We thus define an endo-functor 
\begin{equation}\label{rstvt2b}
\uptau\colon 
\begin{array}[t]{clcr}
\bSM_C&\rightarrow &\bSM_C\\
(V,\sigma)&\mapsto & (\uuptau(V,\sigma),\sigma_\uptau).
\end{array}
\end{equation}

For any Sen $C$-module $(V,\sigma)$, we denote by 
\begin{equation}\label{rstvt2c}
u_T\colon \xymatrix{
{\uuptau(V,\sigma)}\ar[r]&{\hcC\otimes_{\co_C}V}\ar[r]^-(0.5){\halpha_T\otimes\id}&V}
\end{equation}
the $C$-linear composed morphism, where the first map is the canonical injection and $\halpha_T$ is the section of $\hcC$ defined in \eqref{rstvt1h}.

\begin{defi}\label{rstvt3}
A Sen $C$-module $(V,\sigma)$ is said to be {\em convergent} if for every $x\in V$, 
the sequence $\frac{1}{n!}\sigma^n(x)$ converges to $0$.  
\end{defi}

\begin{lem}\label{rstvt4}
For every Sen $C$-module $(V,\sigma)$, the canonical $C$-linear map $u_T\colon \uuptau(V,\sigma)\rightarrow V$ \eqref{rstvt2c} is injective. 
Its image is the sub-$C$-vector space of $V$ made of the elements $x\in V$ such that the sequence $\frac{1}{n!}\sigma^n(x)$ converges to $0$.
\end{lem}

Observe first that every $\hcC$-module of finite type is $p$-adically complete and separated (\cite{egr1} 1.10.2). 
In particular, $\hcC\otimes_{\co_C}V$ is complete and separated for the $p$-adic topology defined by any $\hcC$-lattice of finite type 
(or equivalently coherent $\hcC$-lattice). 
By \eqref{rstvt1e}, any element $x\in \hcC\otimes_{\co_C}V$ can be uniquely written as $x=\sum_{n\geq 0}T^n\otimes x_n$, where $x_n\in V$ 
and $x_n$ converges to $0$. We have 
\begin{equation}\label{rstvt4b}
\sigma_\tot(x)=-\sum_{n\geq 1}nT^{n-1}\otimes x_n+ \sum_{n\geq 0}T^n\otimes \sigma(x_n),
\end{equation}
where $\sigma_\tot$ is defined in \eqref{rstvt2a}. Then, $\sigma_\tot(x)=0$ if and only if for every integer $n\geq 0$, 
\begin{equation}\label{rstvt4c}
\sigma(x_n)=(n+1) x_{n+1}, 
\end{equation}
which is equivalent to $x_n=\frac{1}{n!}\sigma^n(x_0)$. The proposition follows. 

\begin{defi}\label{rstvt5}
A Sen $C$-module $(V,\sigma)$ is said to be {\em $\hcC$-admissible} if the canonical $\hcC$-linear morphism 
\begin{equation}\label{rstvt5a}
\hcC\otimes_{\co_C}\uuptau(V,\sigma) \rightarrow \hcC\otimes_{\co_C}V
\end{equation}
is an isomorphism. 
\end{defi}

\begin{prop}\label{rstvt6}
A Sen $C$-module $(V,\sigma)$ is convergent in the sense of \ref{rstvt3}
if and only if it is $\hcC$-admissible in the sense of \ref{rstvt5}. 
\end{prop}

Indeed, assume first that $(V,\sigma)$ is $\hcC$-admissible. The canonical $\hcC$-linear isomorphism
\begin{equation}
\hcC\otimes_{\co_C}\uuptau(V,\sigma) \stackrel{\sim}{\rightarrow} \hcC\otimes_{\co_C}V
\end{equation}
induces by the base change by $\halpha_T$ \eqref{rstvt1h},
a $C$-linear isomorphism $u\colon \uuptau(V,\sigma)\stackrel{\sim}{\rightarrow} V$. Since $u$ is clearly equal to the canonical morphism  
$u_T\colon \uuptau(V,\sigma)\rightarrow V$ \eqref{rstvt2c}, we deduce by \ref{rstvt4} that $(V,\sigma)$ is convergent. 

Conversely, assume that $(V,\sigma)$ is convergent. The canonical $C$-linear morphism 
$u_T\colon \uuptau(V,\sigma)\rightarrow V$ \eqref{rstvt2c} is an isomorphism by \ref{rstvt4}.
Moreover, the proof of \ref{rstvt4} shows that the composed morphism 
\begin{equation}
\xymatrix{
V\ar[r]^-(0.5){u_T^{-1}}&{\uuptau(V,\sigma)}\ar[r]&{\hcC\otimes_{\co_C}V,}}
\end{equation}
where the second morphism is the canonical injection, is defined, for any $x\in V$, by 
\begin{equation}
\exp(T\otimes \sigma)(1\otimes x)=\sum_{n\geq 0}\frac{1}{n!}T^n\otimes \sigma^n(x),
\end{equation}
which is well defined. We immediately see that the $\hcC$-linear morphism 
\begin{equation}
\exp(T\otimes \sigma)=\sum_{n\geq 0}\frac{1}{n!}T^n\otimes \sigma^n\colon \hcC\otimes_{\co_C}V\rightarrow \hcC\otimes_{\co_C}V
\end{equation}
is a well defined isomorphism, with inverse $\exp(-T\otimes \sigma)$. Hence, $(V,\sigma)$ is $\hcC$-admissible. 

\begin{cor}\label{rstvt7}
For every $\hcC$-admissible Sen $C$-module $(V,\sigma)$, 
the canonical $C$-linear morphism $u_T\colon \uuptau(V,\sigma)\rightarrow V$ \eqref{rstvt2c} is an isomorphism;
it coincides with the $C$-linear isomorphism deduced from  
the $\hcC$-linear admissibility isomorphism \eqref{rstvt5a} by base change by $\halpha_T$ \eqref{rstvt1h}. 
\end{cor}

\subsection{}\label{rstvt8}
For any field extension $F$ of $K$, we consider the base change functor 
\begin{equation}\label{rstvt8a}
\begin{array}[t]{clcr}
\bSM_K&\rightarrow &\bSM_F\\
(W,\sigma)&\mapsto&(W_F,\sigma_F)=(F\otimes_KW,\id\otimes\sigma).
\end{array}
\end{equation}

A Sen $K$-module $(W,\sigma)$ is said to be {\em convergent} (resp.\ {\em $\hcC$-admissible}) if so is $(W_C,\sigma_C)$. 

Let $(W,\sigma)$ be a Sen $K$-module. Since the derivation $\delta$ \eqref{rstvt1f} is $G_K$-equivariant,
the $\hcC$-semi-linear action of $G_K$ on $\hcC\otimes_{\co_K}W=\hcC\otimes_{\co_C}W_C$ 
that fixes the factor $W$, induces a $C$-semi-linear action of $G_K$ on $\uuptau(W_C,\sigma_C)$. 
Moreover, the $C$-linear endomorphism $\sigma_\uptau$ of $\uuptau(W_C,\sigma_C)$ is $G_K$-equivariant \eqref{rstvt2b}. 
Therefore, the Sen $C$-module $\uptau(W_C,\sigma_C)$ is canonically equipped with a $G_K$-equivariant structure.
We denote by  $\bSM_C(G_K)$ the category of $G_K$-equivariant Sen $C$-modules. 
The functor $(W,\sigma)\mapsto \uptau(W_C,\sigma_C)$ can therefore be naturally upgraded into a functor
\begin{equation}\label{rstvt8c}
\upphi\colon \bSM_K\rightarrow \bSM_C(G_K).
\end{equation}
Forgetting the Sen endomorphism $\sigma_\uptau$, we get a functor 
\begin{equation}\label{rstvt8d}
\uupphi\colon \bSM_K\rightarrow \bRep_C(G_K).
\end{equation}

\begin{prop}\label{rstvt9}
For every convergent Sen $K$-module $(W,\sigma)$, the canonical $C$-linear isomorphism \eqref{rstvt2c}
\begin{equation}\label{rstvt9a}
u_T\colon \uuptau(W_C,\sigma_C)\rightarrow W_C
\end{equation}
is $G_K$-equivariant for the $C$-semi-linear action of $G_K$ on $\uuptau(W_C,\sigma_C)$ defined in \eqref{rstvt8d} and 
that on $W_C=C\otimes_KW$ defined for any $g\in G_K$ and $x\in W$, by 
\begin{equation}\label{rstvt9b}
g(1\otimes x)=1\otimes\exp(\log(\chi(g))\sigma)(x).
\end{equation}
\end{prop}

Indeed, by the proof of \ref{rstvt6}, for all $g\in G_K$ and $x\in W$, we have 
\begin{eqnarray*}
g(u_T^{-1}(1\otimes x))&=& g(\exp(T\otimes \sigma)(1\otimes x))\\
&=&\exp((T+\log(\chi(g)))\otimes \sigma)(1\otimes x)\\
&=&\exp(T\otimes \sigma)(\exp(\log(\chi(g))\otimes \sigma)(1\otimes x))\\
&=&u_T^{-1}(1\otimes \exp(\log(\chi(g))\sigma)(x)).
\end{eqnarray*}

\begin{cor}\label{rstvt90}
For every convergent Sen $K$-module $(W,\sigma)$, 
$(W_{K_\infty},\sigma_{K_\infty})$ is the Sen $K_\infty$-module associated by Sen theory \eqref{overview-intro40a} 
with the continuous $C$-representation $\uuptau(W_C,\sigma_C)$ of $G_K$ defined in \eqref{rstvt8d}. 
\end{cor}

It follows immediately from \ref{overview-intro4} and \ref{rstvt9}.

\subsection{}\label{rstvt10}
For any nonzero polynomial $f\in K[X]$, 
we set $B_f=K[X]/(f)$ and $B_{f,C}=C\otimes_KB_f$ equipped with the $C$-semi-linear action of $G_K$ that fixes the factor $B_f$. 
We denote by $\UB_{f,C}$ the subset of {\em principal units} of the finite dimensional $C$-algebra $B_{f,C}$,
i.e., the units $x\in  B^\times_{f,C}$ such that the sequence $(x-1)^n$ converges to $0$. 

\begin{lem}\label{rstvt11}
For every nonzero polynomial $f\in K[X]$, $\UB_{f,C}$ is a subgroup of $B^\times_{f,C}$.
\end{lem}

Indeed, let $x_1,\dots,x_d$ be the zeros of $f$ in $C$ and $n_1,\dots,n_d$ be their multiplicities. We have an isomorphism of $C$-algebras
$B_{f,C}\stackrel{\sim}{\rightarrow}\prod_{1\leq i\leq d}C[X]/(X-x_i)^{n_i}$. It is therefore enough to consider the case where $f=X^n$. 
For $a\in C^\times$ and $g\in C[X]$, setting $x=a+g(X)X$, we have 
\begin{equation}\label{rstvt11a}
(x-1)^m=(a-1+g(X)X)^m=\sum_{j=0}^{n-1}\binom{m}{j}(a-1)^{m-j}g(X)^jX^j.
\end{equation}
The sequence $(x-1)^m$ converges to $0$ if and only if $a-1\in \fm_C$. We deduce a canonical isomorphism 
\begin{equation}\label{rstvt11b}
\UB_{X^n,C}\stackrel{\sim}{\rightarrow}(1+\fm_C)\times (1+X C[X]/(X^n)).
\end{equation}

\begin{lem}\label{rstvt12}
For every nonzero polynomial $f\in K[X]$, the map 
\begin{equation}\label{rstvt12a}
\log \colon \UB_{f,C}\rightarrow B_{f,C}
\end{equation}
given, for any $x\in \UB_{f,C}$, by 
\begin{equation}\label{rstvt12b}
\log(x)=\sum_{n\geq 1}(-1)^{n-1}\frac{(x-1)^n}{n},
\end{equation}
is well defined. It is a $G_K$-equivariant surjective homomorphism with kernel  
the group $\mu_{p^\infty}(B_{f,C})$ of $p$-power roots of unity in $B_{f,C}$. 
\end{lem}

By the proof of \ref{rstvt11}, we may reduce to the case where $f=X^n$. Let $x\in \fm_C$, $g\in C[X]/(X^n)$. 
Taking into account \eqref{rstvt11b}, we have 
\begin{equation}
\log(1+x,1+X g(X)) = \log_C(1+x)+\log_{X}(1+X g(X))\in C\oplus X C[X]/(X^n)=B_{X^n,C},
\end{equation}
where the maps $\log_C\colon 1+\fm_C\rightarrow C$ and $\log_{X}\colon 1+X C[X]/(X^n) \rightarrow X C[X]/(X^n)$ 
are given by the same power series \eqref{rstvt12b}.  The proposition follows from the facts that $\log_C$ is a well defined surjective homomorphism
with kernel $\mu_{p^\infty}(C)$, and $\log_X$ is an isomorphism of groups, with inverse the exponential. Observe that 
the canonical homomorphism $\mu_{p^\infty}(C)\rightarrow \mu_{p^\infty}(B_{X^n,C})$ is an isomorphism. 

\begin{lem}\label{rstvt13}
Let $f$ be a nonzero polynomial in $K[X]$ and let $g\in K[X]$ be a polynomial such that $B_g=(B_f)_{\red}$. 
We set $B_{f,\oK}=B_f\otimes_K\oK$ and $B_{g,\oK}=B_g\otimes_K\oK$. Then all the homomorphisms of the commutative diagram 
\begin{equation}
\xymatrix{
{\mu_{p^\infty}(B_{f,\oK})}\ar[r]\ar[d]&{\mu_{p^\infty}(B_{g,\oK})}\ar[d]\\
{\mu_{p^\infty}(B_{f,C})}\ar[r]&{\mu_{p^\infty}(B_{g,C})}}
\end{equation}
are isomorphisms. 
\end{lem}

Indeed, if $x_1,\dots,x_d$ are the zeros of $f$ in $\oK$ with multiplicities  $n_1,\dots,n_d$, we have an isomorphism of $\oK$-algebras
$B_{f,\oK}\stackrel{\sim}{\rightarrow}\prod_{1\leq i\leq d}\oK[X]/(X-x_i)^{n_i}$, and hence 
$B_{g,\oK}\stackrel{\sim}{\rightarrow}\prod_{1\leq i\leq d}\oK [X]/(X-x_i)$. The proposition then follows from the fact 
that the canonical homomorphism $\mu_{p^\infty}(\oK)\rightarrow \mu_{p^\infty}(C)$ is an isomorphism.

\begin{prop}\label{rstvt14}
Let $f$ be a nonzero polynomial in $K[X]$. Then, 
\begin{itemize}
\item[{\rm (i)}] The canonical morphism 
\begin{equation}
\rH^1(G_K,\mu_{p^\infty}(B_{f,C}))\rightarrow \rH^1_\cont(G_K,B^\times_{f,C})
\end{equation}
is trivial. 
\item[{\rm (ii)}] If the residue field $k$ of $\co_K$ is algebraically closed, then we have $\rH^2(G_K,\mu_{p^\infty}(B_{f,C}))=0$.
\end{itemize}
\end{prop}

Let $g\in K[X]$ be a polynomial such that $B_g=(B_f)_{\red}$. 

(i) Since the homomorphism $K\rightarrow B_g$ is étale, 
there exists a $K$-homomorphism $\varphi\colon B_g\rightarrow B_f$ which is a section of the canonical surjection $B_f\rightarrow B_g$. 
By \ref{rstvt13}, using the section $\varphi$, it is enough to prove that the canonical morphism 
\begin{equation}
\rH^1(G_K,\mu_{p^\infty}(B_{g,\oK}))\rightarrow \rH^1(G_K,B^\times_{g,\oK})
\end{equation}
is trivial; observe that $B_{g,\oK}$ is a discrete $G_K$-module. Let $\phi\colon \Spec(B_g)\rightarrow \Spec(K)$ be the canonical morphism. 
Since $B_g$ is a finite product of finite separable field extensions of $K$, we have 
\[
\rH^1(G_K,B^\times_{g,\oK})=\rH^1_\et(\Spec(K),\phi_*(\mG_m))=\rH^1_\et(\Spec(B_g),\mG_m)=0,
\]
which implies the proposition. 

(ii) By \ref{rstvt13}, it is enough to prove that $\rH^2(G_K,\mu_{p^\infty}(B_{g,\oK}))=0$. 
Since $B_g$ is étale over $K$, 
and the residue field of $\co_K$ is algebraically closed, we have 
\[
\rH^2(G_K,B^\times_{g,\oK})=\rH^2_\et(\Spec(B_g),\mG_m)=0. 
\]
For every integer $n\geq 1$, we have the Kummer exact sequence 
\begin{equation}
0\longrightarrow \mu_{p^n}(B_{g,\oK}) \longrightarrow B^\times_{g,\oK} \stackrel{p^n}{\longrightarrow} B^\times_{g,\oK} \longrightarrow 0. 
\end{equation}
We deduce that $\rH^2(G_K,\mu_{p^n}(B_{g,\oK}))=0$, which implies the required assertion. 

\begin{lem}\label{rstvt18}
Let $G$ be a profinite group, 
\begin{equation}
0 \rightarrow L \rightarrow M \stackrel{\sigma}{\rightarrow} N \rightarrow 0
\end{equation}
a short exact sequence of topological $G$-modules such that the topology of $L$ is induced by that of $M$ and  
$\sigma$ admits locally continuous sections (these sections being just maps of topological sets), i.e., for every $y \in N$, 
there exists an open neighborhood $U$ of $y$ in $N$ and a continuous map $s\colon U \rightarrow M$ such that $\sigma \circ s = \id_U$.
Then, the induced sequence of complexes of continuous cochains
\begin{equation}
0 \rightarrow \rC^\bullet_{\cont}(G,L) \rightarrow \rC^\bullet_{\cont}(G,M) \rightarrow \rC^\bullet_{\cont}(G,N) \rightarrow 0
\end{equation}
is exact. In particular, there is a natural long exact sequence of continuous cohomology
\begin{equation}
\rH^n_{\cont}(G,L)\rightarrow \rH^n_{\cont}(G,M)\rightarrow \rH^n_{\cont}(G,N) \rightarrow \rH^{n+1}_{\cont}(G,L). 
\end{equation}
\end{lem}

The proposition is stated in (\cite{tate1} §2) for a general topological group $G$ under the stronger assumption that $\sigma$ admits 
a continuous section. 

We only need to prove that for every integer $n\geq 1$, the morphism 
\begin{equation}
\begin{array}[t]{clcr}
\rC^n_{\cont}(G,M)&\rightarrow&\rC^n_{\cont}(G,N),\\
f&\mapsto&\sigma\circ f
\end{array}
\end{equation}
is surjective. Let $\varphi\colon G^n\rightarrow N$ be a continuous map. 
Since $G$ is profinite, $G^n$ is compact and admits a basis of open sets which are also closed.  
Therefore, $\varphi(G^n)$ is compact. It follows that there exist finitely many open subsets $U_1,\dots,U_\ell$ of $N$ covering $\varphi(G^n)$, 
together with continuous maps $s_i\colon U_i \to M$ for $1 \leq i \leq \ell$, such that $\sigma \circ s_i = \mathrm{id}_{U_i}$.
We set $V_i=\varphi^{-1}(U_i)$ so we have $G^n=\cup_{1\leq i\leq \ell}V_i$. There exists a finite partition $G^n=\sqcup_{1\leq j\leq r}W_i$
of $G^n$ into closed and open subsets $W_j$ such that for every $1\leq j\leq r$, $W_j\subset V_{i(j)}$ for an integer $1\leq i(j)\leq \ell$. 
For every $\gamma\in W_j$, we have $\varphi(\gamma)\in U_{i(j)}$. We set 
\begin{equation}
f(\gamma)=s_{i(j)}(\varphi(\gamma)) \in M.
\end{equation}
Then, $f\colon G^n\rightarrow M$ is a well defined continuous map and we have $\sigma\circ f=\varphi$, which proves the required statement.

\begin{prop}\label{rstvt15}
Assume that the residue field $k$ of $\co_K$ is algebraically closed. Let $f$ be a nonzero polynomial in $K[X]$,
$P$ a continuous $G_K$-equivariant torsor under the group $B_{f,C}$. Then, there exists a continuous $G_K$-equivariant torsor $\UP$ 
under the group $\UB_{f,C}$ such that the torsor obtained from $\UP$ by extending the structural group by the homomorphism 
$\log\colon \UB_{f,C}\rightarrow B_{f,C}$ \eqref{rstvt12a} is isomorphic to $P$ \eqref{p1-NC4}. 
Moreover, the torsor obtained from $\UP$ by extending the structural group by the canonical injection 
$\UB_{f,C}\rightarrow B^\times_{f,C}$ is unique up to isomorphism. 
\end{prop}

Indeed, by \ref{rstvt12}, we have a $G_K$-equivariant exact sequence of abelian representations of $G_K$,  
\begin{equation}\label{rstvt15a}
\xymatrix{
0\ar[r]&{\mu_{p^\infty}(B_{f,C})}\ar[r]&{\UB_{f,C}}\ar[r]^-(0.5){\log}&{B_{f,C}}\ar[r]&0.}
\end{equation}
Observe that the $p$-adic topology on $\UB_{f,C}$ induces the discrete topology on $\mu_{p^\infty}(B_{f,C})$, 
and $\log$ is continuous for the $p$-adic topology and admits locally continuous sections (induced by the exponential). 
Hence, by \ref{rstvt14}(ii) and \ref{rstvt18}, the exact sequence \eqref{rstvt15a} induces an exact sequence
\begin{equation}
\rH^1(G_K,\mu_{p^\infty}(B_{f,C}))\rightarrow \rH^1_{\cont}(G_K,\UB_{f,C})\rightarrow \rH^1_\cont(G_K,B_{f,C})\rightarrow 0.
\end{equation}
We deduce the existence of the required continuous $G_K$-equivariant torsor $\UP$ under the group $\UB_{f,C}$. 
By \ref{rstvt14}(i), the image of the class of such a torsor $\UP$ by the canonical homomorphism 
\begin{equation}
\rH^1_{\cont}(G_K,\UB_{f,C})\rightarrow \rH^1_{\cont}(G_K,B^\times_{f,C})
\end{equation} 
is uniquely determined, which finishes the proof of the proposition. 

\subsection{}\label{rstvt16}
We set $\ccF_{\mZ_p}=\Hom_{\mZ_p}(\cF_{\mZ_p},\mZ_p)$ equipped with the action of $G_K$ induced by that on $\cF_{\mZ_p}$ \eqref{rstvt1}. 
The exact sequence \eqref{rstvt1b} induces an exact sequence of $\mZ_p[G_K]$-modules 
\begin{equation}\label{rstvt16a}
0\rightarrow \mZ_p\rightarrow \ccF_{\mZ_p}\rightarrow \mZ_p\rightarrow 0,
\end{equation}
We denote by $(\cun,\ccT)$ the  $\mZ_p$-basis of $\ccF_{\mZ_p}$ dual of the $\mZ_p$-basis $(1,T)$ of $\cF_{\mZ_p}$; 
so $\ccT$ is the image of $1\in \mZ_p$ by the second map and $\cun$ determines the splitting of \eqref{rstvt16a} dual of the second map of \eqref{rstvt1b}.  
For every $g\in G_K$, we have 
\begin{eqnarray}
g(\cun)&=&\cun-\log(\chi(g))\ccT,\label{rstvt16b1}\\
g(\ccT)&=&\ccT.\label{rstvt16b2}
\end{eqnarray}

Let $f$ be a nonzero polynomial in $K[X]$. We take again the notation of \ref{rstvt10}.
We denote by $\ccF_{f,C}$ the pushout of the extension $\ccF_{\mZ_p}$ by the $\mZ_p$-linear morphism 
$\mZ_p\ccT\rightarrow B_{f,C}$ sending $\ccT$ to the class of $-X$. 
So $\ccF_{f,C}$ is a continuous abelian representation of $G_K$ that fits into a $G_K$-equivariant exact sequence 
\begin{equation}\label{rstvt16c}
0\rightarrow B_{f,C}\rightarrow \ccF_{f,C}\rightarrow \mZ_p\rightarrow 0. 
\end{equation}

Assume that the residue field $k$ of $\co_K$ is algebraically closed. 
Let $P_{f,C}\subset \ccF_{f,C}$ be the inverse image of $1\in \mZ_p$ by the third map of  \eqref{rstvt16c}, 
which is a continuous $G_K$-equivariant torsor under the group $B_{f,C}$.
We choose a continuous $G_K$-equivariant torsor $\UP_{f,C}$ under the group $\UB_{f,C}$ 
such that the torsor obtained from $\UP_{f,C}$ by extending the structural group by the homomorphism 
$\log\colon \UB_{f,C}\rightarrow B_{f,C}$ \eqref{rstvt12a} is isomorphic to $P_{f,C}$; such a torsor exists by \ref{rstvt15}. 
We denote by
\begin{equation}\label{rstvt16d}
P^\times_{f,C}=\UP_{f,C}\wedge^{\UB_{f,C}}B^\times_{f,C}
\end{equation}
the torsor obtained from $\UP_{f,C}$ by extending the structural group by the canonical injection 
$\UB_{f,C}\rightarrow B^\times_{f,C}$, which by \ref{rstvt15} does not depend on the choice of $\UP_{f,C}$ up to isomorphism. 
Let $L_{f,C}$ be the associated continuous $G_K$-equivariant invertible $B_{f,C}$-module. 

\begin{prop}\label{rstvt17}
Assume that the residue field $k$ of $\co_K$ is algebraically closed. 
Let $(W,\sigma)$ be a Sen $K$-module, $f\in K[X]$ a nonzero polynomial such that $f(\sigma)=0$. We take again the notation of \ref{rstvt10} 
and consider $W$ as a $B_f$-module by letting $X$ acts by $\sigma$. 
Let $L_{f,C}$ be the continuous $G_K$-equivariant invertible $B_{f,C}$-module defined in \ref{rstvt16}. 
Then, $(W,\sigma)$ is the Sen $K$-module associated by Sen theory \eqref{overview-intro40b} 
with the continuous $C$-representation $V=L_{f,C}\otimes_{B_f}W$ of $G_K$ that fixes the factor $W$. 
\end{prop}

Consider $\cun\in P_{f,C}\subset \ccF_{f,C}$ and let $g\in G_K$; see \ref{rstvt16}. By \eqref{rstvt16b1}, we have $g(\cun)=\cun+\log(\chi(g))X$. 
Choose a continuous $G_K$-equivariant torsor $\UP_{f,C}$ under the group $\UB_{f,C}$ as in \ref{rstvt16}, 
and let $e\in \UP_{f,C}$ be such that the image of $(e,1)\in \UP_{f,C}\times B_{f,C}$ is $\cun\in P_{f,C}$.  
Then, there exists a continuous cocycle $c\colon G_K\rightarrow UB_{f,C}$ such that for every $g\in G_K$, we have 
\begin{equation}
g(e)= c(g)\cdot e.
\end{equation}
Then, $\log(c(g))=\log(\chi(g))X$ for every $g\in G_K$. By \ref{rstvt12}, 
there exists an open subgroup $H$ of $G_K$ such that for every $g\in H$, $\exp(\log(\chi(g))X)$ is well defined in $UB_{f,C}$,
and we have $c(g)=\eta(g) \exp(\log(\chi(g))X)$, for an element $\eta(g)\in \mu_{p^\infty}(B_{f,C})$. 
Since the cocycle $c$ is continuous, the map $\eta\colon H\rightarrow  \mu_{p^\infty}(B_{f,C})$ is continuous. Hence, after shrinking $H$,
we may assume that $\eta(g)=1$ for every $g\in H$. 
Setting  $F=\oK^H$, $B^\times_{f,F}\cdot e \subset P^\times_{f,C}$ \eqref{rstvt16d}
is then a torsor under the group $B^\times_{f,F}$ which is stable by the action of $G_F=H$. 
Let $L_{f,F}$ be the associated continuous $G_F$-equivariant 
invertible $B_{f,F}$-module. We have a canonical $G_F$-equivariant $C$-linear isomorphism 
\begin{equation}
L_{f,F}\otimes_FC\stackrel{\sim}{\rightarrow}L_{f,C}.
\end{equation}
The element $e$ induces a $B_{f,F}$-linear isomorphism 
\begin{equation}
W_F=F\otimes_KW \stackrel{\sim}{\rightarrow} L_{f,F}\otimes_{B_f} W=V_F.
\end{equation}
The latter is $G_F$-equivariant when we equip $V_F=L_{f,F}\otimes_{B_f} W$ with the action of $G_F$ that fixes the factor $W$, 
and $W_F$ with the action of $G_F$ defined for any $g\in G_F$ by $\exp(\log(\chi(g))\sigma)$.  
We deduce $G_F$-equivariant $C$-linear isomorphisms 
\begin{equation}
C\otimes_FW_F \stackrel{\sim}{\rightarrow} C\otimes_FV_F\stackrel{\sim}{\rightarrow}V,
\end{equation}
where $V=L_{f,C}\otimes_{B_f}W$ is equipped with the action of $G_K$ that fixes the factor $W$. 
The proposition follows then by (\cite{sen1} Proposition 8').

\chapter{Preliminaries}\label{preliminaries}

\section{Notation and conventions}\label{p1-NC}

\subsection{}\label{p1-NC0}
Throughout this book, we fix a universe $\mU$ with an element of infinite cardinality.
A set is said to be {\em $\mU$-small} (or, when no confusion arises, {\em small}) if it is isomorphic to
an element of $\mU$. We also use the terminology: {\em small group}, {\em small ring}, {\em small category}...
We say that a category $\cC$ is a {\em $\mU$-category} if for all objects $X,Y$ of $\cC$, the set $\Hom_\cC(X,Y)$
is $\mU$-small (\cite{sga4} I 1.1).
The {\em category of $\mU$-sets}, denoted by $\Ens$, is the category of sets that are in $\mU$.
It is a punctual $\mU$-topos (\cite{sga4} IV 2.2).
We denote by $\Sch$ the category of  schemes that are elements of $\mU$.
Unless explicitly stated otherwise, it will be understood that the rings and
the logarithmic schemes (and in particular the schemes) considered in this book are elements of the
universe~$\mU$.

\subsection{}\label{p2-ncgt4}
Let $\cC$ be a $\mU$-category (\cite{sga4} I 1.1). We denote by $\hcC$ the category
of presheaves of $\mU$-sets on $\cC$, that is, the category of contravariant  functors
on $\cC$ with values in $\Ens$ (\cite{sga4} I 1.2), and by 
\begin{equation}\label{p2-ncgt4a}
\tth_{\cC}\colon \cC\rightarrow\hcC, \ \ \ 
X\mapsto (\tth_\cC(X)\colon Y\mapsto \Hom_{\cC}(Y,X)),
\end{equation}
the canonical functor (\cite{sga4} I 1.3), which is fully faithful (\cite{sga4} I 1.4).
If $\cC$ is endowed with a topology (\cite{sga4} II 1.1), we denote by $\tcC$ the topos of
sheaves of $\mU$-sets on $\cC$ (\cite{sga4} II 2.1).

For an object $F$  of $\hcC$, we denote by $\cC_{/F}$ the following category
(\cite{sga4} I 3.4.0). The objects of $\cC_{/F}$ are the pairs consisting of an object $X$ of $\cC$
and a morphism $u$ from $X$  to $F$. If $(X,u)$ and $(Y,v)$ are two objects,
a morphism from $(X,u)$ to $(Y,v)$ is a morphism $g\colon X\rightarrow Y$
such that $u=v\circ g$.

\subsection{}\label{p2-ncgt5}
Let $X$ be a $\mU$-topos (\cite{sga4} IV 1.1.2). Inverse systems of objects of $X$ indexed by 
the ordered set of natural numbers $\mN$, form a topos that we denote by $X^{\mN^\circ}$. We refer to
(\cite{agt} III.7) for useful facts on this type of topos. We recall, in particular, that we have a morphism of topos
\begin{equation}\label{p2-ncgt5a}
\uplambda\colon X^{\mN^\circ}\rightarrow X,
\end{equation}
whose pullback functor $\uplambda^*$ associates with any object $F$ of $X$ the constant inverse system of 
value $F$, and whose direct image functor $\uplambda_*$ associates with any inverse system its inverse limit (\cite{agt} III.7.4).

\subsection{}\label{p1-NC4}
Let $T$ be a $\mU$-topos, $u\colon F\rightarrow G$ a group homomorphism of $T$, $P$ a right $F$-torsor of $T$.  
The {\em $u$-twist of $P$} is the contracted product $P\wedge^FG$, where $F$ acts on the left on $G$ by left translations,
i.e., the quotient of the product $P\times G$ by the group $F$ acting diagonally: 
\begin{equation}
P\times G \times F\rightarrow P\times G, \ \ \ (p, g,f)\mapsto (pf, u(f^{-1})g);
\end{equation}
it is naturally a right $G$-torsor (\cite{giraud2} III 1.4.6). We say that $P\wedge^FG$ {\em is obtained from $P$ by extension of its structural group by $u$} or simply 
that $P\wedge^FG$ {\em is the $u$-twist of $P$}.

\subsection{}\label{p1-NC7}
When we consider a ringed $\mU$-topos $(X,A)$ (\cite{sga4} IV 11.1.1), the ring $A$ is
assumed to be commutative unless stated otherwise. 
For such a ringed $\mU$-topos, we denote by $\bMod(A)$ or $\bMod(A,X)$ the category of $A$-modules of $X$.

If $M$ is an $A$-module, we denote by $\rS_A(M)$
(resp.\ $\wedge_A(M)$, resp.\ $\Gamma_A(M)$) the symmetric  (resp.\ exterior, resp.\ divided powers) algebra
of $M$ and for any integer $n\geq 0$, by $\rS_A^n(M)$ (resp.\ $\wedge_A^n(M)$,
resp.\ $\Gamma_A^n(M)$) its homogeneous part of degree $n$.
The formations of these algebras commute with localization over an object of $X$.
We will omit the ring $A$ from the notation when there is no risk of ambiguity.

For any integer $n\geq 0$, we set $\Gamma^{\geq n}(M)=\oplus_{i\geq n}\Gamma^i(M)$; 
so $\Gamma^{\geq 1}(M)$ is the divided power ideal of $\Gamma(M)$, and $\Gamma^{\geq n}(M)$ 
is its $n$th divided power (\cite{bo} 3.24). We set 
\begin{eqnarray}
\Gamma^{\leq n}(M)&=&\Gamma(M)/\Gamma^{\geq n+1}(M),\label{p1-NC7a}\\
\hGamma(M)&=&\underset{\underset{n\geq 0}{\longleftarrow}}\lim\ \Gamma^{\leq n}(M),\label{p1-NC7b}
\end{eqnarray}
which are naturally $A$-algebras with divided powers. 
We denote by
\begin{equation}\label{p1-NC7c}
\exp_{M}\colon \begin{array}[t]{clcr}
M&\rightarrow& \hGamma(M)^\times,\\
x&\mapsto&\sum_{n\geq 0} x^{[n]},
\end{array}
\end{equation}
the exponential homomorphism. 

\subsection{}\label{p1-NC8}
For any abelian category $\bA$, we denote by $\bC(\bA)$ the category complexes of $\bA$ and by $\bC^-(\bA)$, $\bC^+(\bA)$, and $\bC^\rb(\bA)$ the full subcategories of $\bC(\bA)$ of complexes bounded from above, from below, and from both sides, respectively.
We denote by $\bD(\bA)$ the derived category of $\bA$ and by $\bD^-(\bA)$, $\bD^+(\bA)$, and $\bD^\rb(\bA)$ the full
subcategories of $\bD(\bA)$ of complexes with cohomology
bounded from above, from below, and from both sides, respectively.
Unless mentioned otherwise, complexes in $\bA$ have a differential
of degree $+1$, the degree being written as an exponent.

Following (\cite{hodge2} 1.4.6) and (\cite{sp} \href{https://stacks.math.columbia.edu/tag/0118}{0118}),
for any $F\in \ob(\bD(\bA))$ and any integer $q$, we define the canonical truncation functors by 
\begin{equation}\label{p1-NC8a}
\tau_{\leq q}(F)^n=
\left\{
\begin{array}{clcr}
0 & {\rm for} \ n > q,\\
\ker(d) & {\rm for} \ n=q,\\
F^n &  {\rm for} \ n < q,
\end{array}
\right.
\ \ \ {\rm and}\ \ \
\tau_{\geq q}(F)^n=
\left\{
\begin{array}{clcr}
F^n & {\rm for} \ n > q,\\
\coker(d) & {\rm for} \ n=q,\\
0 &  {\rm for} \ n < q.
\end{array}
\right.
\end{equation}
We have canonical distinguished triangles
\begin{eqnarray}
\tau_{\leq q-1}(F)\rightarrow \tau_{\leq q}(F)\rightarrow \cH^q(F)[-q]\rightarrow \tau_{\leq q-1}(F)[1],\label{p1-NC8b1}\\
\cH^q(F)[-q]\rightarrow \tau_{\geq q}(F)\rightarrow \tau_{\geq q+1}(F) \rightarrow \cH^q(F)[-q+1],\label{p1-NC8b2}\\
\tau_{\leq q}(F)\rightarrow F\rightarrow \tau_{\geq q+1}(F)\rightarrow \tau_{\leq q}(F)[1],\label{p1-NC8b3}\\
\end{eqnarray}

For any complex $F^\bullet$ of $\bA$ and any integer $q$, we define the stupid truncation functors by 
\begin{equation}\label{p1-NC8c}
\sigma_{\leq q}(F^\bullet)^n=
\left\{
\begin{array}{clcr}
0 & {\rm for} \ n > q,\\
F^n &  {\rm for} \ n \leq q,
\end{array}
\right.
\ \ \ {\rm and}\ \ \
\sigma_{\geq q}(F^\bullet)^n=
\left\{
\begin{array}{clcr}
F^n & {\rm for} \ n \geq q,\\
0 &  {\rm for} \ n < q.
\end{array}
\right.
\end{equation}
We have canonical exact sequences 
\begin{eqnarray}
0\rightarrow F^q[-q]\rightarrow \sigma_{\leq q}(F^\bullet)\rightarrow \sigma_{\leq q-1}(F^\bullet)\rightarrow 0,\label{p1-NC8d1}\\
0\rightarrow \sigma_{\geq q+1}(F^\bullet) \rightarrow \sigma_{\geq q}(F^\bullet)\rightarrow F^q[-q] \rightarrow 0,\label{p1-NC8d2}\\
0\rightarrow \sigma_{\geq q+1}(F^\bullet) \rightarrow F^\bullet \rightarrow \sigma_{\leq q}(F^\bullet) \rightarrow 0.\label{p1-NC8d3}
\end{eqnarray}

\subsection{}\label{p2-ncgt6}
Let $X$ be a scheme. We denote by $\Et_{/X}$ the {\em étale site} of $X$,
i.e., the full subcategory of $\Sch_{/X}$ made up  of the étale schemes 
over $X$, endowed with the étale topology; it is a $\mU$-site (\cite{sga4} II 3.0.2).
We denote by $X_\et$ the {\em étale topos} of $X$, that is, the topos of sheaves of 
$\mU$-sets on $\Et_{/X}$.

We denote by $\Et_{\coh/X}$ the full subcategory of $\Et_{/X}$ made up of étale schemes of finite 
presentation on $X$, endowed with the topology induced by that of $\Et_{/X}$; it is a $\mU$-small site.
If $X$ is quasi-separated, the restriction functor from $X_\et$ to the topos of sheaves of $\mU$-sets on
$\Et_{\coh/X}$ is an equivalence of categories (\cite{sga4} VII 3.1 and 3.2).

We denote by $\Et_{\rf/X}$ the full subcategory of $\Et_{/X}$ made up of finite étale schemes on $X$,
endowed with the topology induced by that of $\Et_{/X}$; it is a $\mU$-small site.
The {\em finite étale topos} of $X$, denote by $X_\fet$, is the topos of sheaves of $\mU$-sets 
on $\Et_{\rf/X}$ (see \cite{agt} VI.9.2). The canonical injection $\Et_{\rf/X}\rightarrow \Et_{/X}$ induces a 
morphism of topos
\begin{equation}\label{p2-ncgt6a}
\rho_X\colon X_\et\rightarrow X_\fet.
\end{equation}

We denote by $X_\zar$ the Zariski topos of $X$, and by
\begin{equation}\label{p1-NC6a}
u_X\colon X_{\et}\rightarrow X_{\zar}
\end{equation}
the canonical morphism (\cite{sga4} VII 4.2.2). If $F$ is a quasi-coherent $\co_X$-module of $X_\zar$,
we denote by $\iota(F)$ the sheaf of $X_\et$ defined for any  étale $X$-scheme $U$ by 
\begin{equation}\label{p1-NC6b}
\iota(F)(U)=\Gamma(U,F\otimes_{\co_X}\co_U).
\end{equation}
It is convenient, when there is no risk of confusion, to abusively denote $\iota(F)$ by $F$.
Note that $\iota(\co_X)$ is a ring of $X_\et$ and that $\iota(F)$ is a $\iota(\co_X)$-module (see \cite{ag2} 2.1.13).

For every quasi-coherent $\co_X$-module $F$ of $X_\zar$, we have a canonical isomorphism
\begin{equation}\label{p1-NC6d}
F\stackrel{\sim}{\rightarrow}u_{X*}(\iota(F)).
\end{equation}
We therefore consider $u_X$ as a morphism of ringed topos
\begin{equation}\label{p1-NC6e}
u_X\colon (X_{\et},\co_X)\rightarrow (X_{\zar},\co_X).
\end{equation}
For modules, we use the notation $u^{-1}_X$ to denote the pullback in the sense of abelian sheaves,
and we keep the notation $u^*_X$ for the pullback in the sense of modules.
By (\cite{ag1} 2.1.18), the isomorphism \eqref{p1-NC6d} induces by adjunction an isomorphism
\begin{equation}\label{p1-NC6f}
u^*_X(F)\stackrel{\sim}{\rightarrow}\iota(F).
\end{equation}

Any morphism of schemes  $f\colon Y\rightarrow X$ induces morphisms of topos that we often denote (abusively) 
also by  $f\colon Y_\zar\rightarrow X_\zar$, $f\colon Y_\et\rightarrow X_\et$ and $f\colon Y_\fet\rightarrow X_\fet$. 
For $\co_X$-modules of $X_\zar$ (resp.\ $X_\et$), we will use the notation $f^{-1}$ to denote the
pullback in the sense of abelian sheaves and will keep the
notation $f^*$ for the pullback in the sense of modules. 

\subsection{}\label{p1-NC5}
Let $f\colon X'\rightarrow X$ be a morphism of schemes, 
$\rT$ a quasi-coherent $\co_X$-module, and $\cL$ a $\rT$-torsor of $X_\zar$ (resp.\ $X_\et$).  
The {\em affine pullback} of $\cL$ under $f$, denoted by $f^+(\cL)$,
is the $f^*(\rT)$-torsor of $X'_\zar$ (resp.\ $X'_\et$) deduced from the
$f^{-1}(\rT)$-torsor $f^*(\cL)$ by extending its structural group by
the canonical homomorphism $f^{-1}(\rT)\rightarrow f^*(\rT)$ \eqref{p1-NC4}. 

\subsection{}\label{p1-NC1}
We refer to (\cite{agt} § II.5) for a lexicon of logarithmic geometry and we adopt the same notation. 
We denote by $\FLS$ the category of fine logarithmic schemes (for the étale topology). 
For any object $(X,\cM_X)$ of $\FLS$, we denote by $\FLS_{/(X,\cM_{X})}$ 
the category of fine logarithmic schemes above it. 

\subsection{}\label{p1-NC3}
Let $i\colon (S,\cM_{S})\rightarrow (\tS,\cM_{\tS})$ be a thickening of fine logarithmic schemes (\cite{ogus} IV 2.1.1), 
$I$ the ideal of $\co_\tS$ defining the closed immersion of schemes $S\rightarrow \tS$. 
Then, a diagram of $\FLS$
\begin{equation}\label{p1-NC3b}
\xymatrix{
{(X,\cM_{X})}\ar[r]^-(0.5){j}\ar[d]&{(\tX,\cM_{\tX})}\ar[d]\\
{(S,\cM_{S})}\ar[r]^-(0.5){i}&{(\tS,\cM_{\tS})}}
\end{equation}
is Cartesian if and only if it is Cartesian in the category of logarithmic schemes; in this case, $j$ is a thickening and the underlying 
closed immersion of schemes $X\rightarrow \tX$ is defined by the ideal $I\co_\tX$. 
Indeed, let $(X',\cM_{X'})$ be the logarithmic scheme defined by the Cartesian diagram in the category of logarithmic schemes
\begin{equation}\label{p1-NCEc}
\xymatrix{
{(X',\cM_{X'})}\ar[r]^-(0.5){j'}\ar[d]\ar@{}[rd]|{\Box}&{(\tX,\cM_{\tX})}\ar[d]\\
{(S,\cM_{S})}\ar[r]^-(0.5){i}&{(\tS,\cM_{\tS}).}}
\end{equation}
It is clear that $j'$ is a strict closed immersion. Since the underlying closed immersion of schemes $X'\rightarrow \tX$ 
is defined by the ideal $I\co_\tX$, $j'$ is a thickening. 
As the logarithmic scheme $(\tX,\cM_\tX)$ is fine, so is $(X',\cM_{X'})$ by (\cite{ogus} IV 2.1.3). Hence $(X,\cM_X)=(X',\cM_{X'})$.

\section{Formal schemes}\label{p1-pfs}

\begin{defi}\label{p1-pfs1}
Let $\fX$ be an adic formal scheme (\cite{egr1} 2.1.24), $\cA$ an $\co_\fX$-algebra. 
\begin{itemize}
\item[(i)] We say that $\cA$ is {\em adic} if for every ideal of definition of finite type $\cJ$ of $\fX$, 
the $(\co_{\fX}/\cJ)$-algebra $\cA/\cJ\cA$ is quasi-coherent and the canonical morphism 
\begin{equation}\label{p1-pfs1a}
\cA\rightarrow \underset{\underset{n\geq 1}{\longleftarrow}}\lim\ \cA/\cJ^{n}\cA
\end{equation}
is an isomorphism. 
\item[(ii)] We say that $\cA$ is {\em topologically of finite type} (resp.\ {\em topologically of finite presentation}) 
if it is adic and if for every ideal of definition of finite type $\cJ$ of $\fX$, the $(\co_{\fX}/\cJ)$-algebra $\cA/\cJ\cA$ is of finite type 
(resp.\ of finite presentation).
\end{itemize}
\end{defi}

\begin{lem}\label{p1-pfs2}
Let $\fX$ be an adic formal scheme, $\cJ$ an ideal of definition of finite type of $\fX$. 
For an $\co_\fX$-algebra $\cA$ to be adic, it is necessary and sufficient that it satisfies the following conditions:
\begin{itemize}
\item[{\rm (i)}] for every integer $n\geq 1$, the $(\co_{\fX}/\cJ^n)$-algebra $\cA/\cJ^n\cA$ is quasi-coherent; 
\item[{\rm (ii)}] the canonical morphism 
\begin{equation}\label{p1-pfs2a}
\cA\rightarrow \underset{\underset{n\geq 1}{\longleftarrow}}\lim\ \cA/\cJ^{n}\cA
\end{equation}
is an isomorphism. 
\end{itemize}
\end{lem}

Obvious. 

\begin{lem}\label{p1-pfs3}
Let $\fX$ be an adic formal scheme, $\cJ$ an ideal of definition of finite type of $\fX$. 
For an $\co_\fX$-algebra $\cA$ to be topologically of finite type (resp.\ topologically of finite presentation), 
it is necessary and sufficient that it is adic and that for every integer $n\geq 1$, 
the $(\co_{\fX}/\cJ^n)$-algebra $\cA/\cJ^n\cA$ is of finite type (resp.\ finite presentation). 
\end{lem}

Obvious.

\subsection{}\label{p1-pfs4}
Let $\fX$ be an adic formal scheme, $\cJ$ an ideal of definition of finite type of $\fX$, $\cA$ an adic $\co_\fX$-algebra. 
For any integer $n\geq 1$, we set $\fX_n=(\fX,\co_\fX/\cJ^n)$. We associate with $\cA$ an adic formal scheme, 
denoted by $\fY=\Spf(\cA)$, equipped with an affine adic morphism $f\colon \fY\rightarrow \fX$ (\cite{egr1} 2.3.9 and 2.3.4). 
Then, $\cK=f^*(\cJ)\co_\fY$ is an ideal of definition of $\fY$ and for every integer $n\geq 1$, there exists a canonical 
$\fX_n$-isomorphism $(\fY,\co_\fY/\cK^n)\stackrel{\sim}{\rightarrow}\Spec(\cA/\cJ^n\cA)$ (\cite{egr1} 2.2.14).  
We obviously have $f_*(\co_\fY)=\cA$ and $f_*(\co_\fY/\cK^n)=\cA/\cJ^n\cA$. We deduce that $f_*(\cK^n)=\cJ^n\cA$. 

The contravariant functor which associates with $\cA$ the pair $(\fY,f)$ is an equivalence between the category of adic $\co_\fX$-algebras and 
that of affine adic $\fX$-schemes (\cite{egr1} 2.3.10). 

If the $\co_\fX$-algebra $\cA$ is topologically of finite type (resp.\ topologically of finite presentation),  
then the morphism $f$ is of finite type (resp.\ of finite presentation) (\cite{egr1} 2.3.13, resp.\ 2.3.15).

\begin{lem}\label{p1-pfs5}
We keep the assumptions and notation of \ref{p1-pfs4}. We assume, moreover, that  $\fX$ is an affine formal scheme {\rm (\cite{egr1} 2.1.34)}. Then, 
\begin{itemize} 
\item[{\rm (i)}] We have $\Gamma(\fX,\cJ\cA)=\cJ(\fX)\cA(\fX)$ and $\Gamma(\fX,\cA/\cJ\cA)=\cA(\fX)/\cJ(\fX)\cA(\fX)$. 
\item[{\rm (ii)}] The $\co_\fX(\fX)$-algebra $\cA(\fX)$ is adic. 
\item[{\rm (iii)}] The $\co_\fX$-algebra $\cA$ is topologically of finite type,  
if and only if the $\co_\fX(\fX)$-algebra $\cA(\fX)$ is topologically of finite type.
\item[{\rm (iv)}] Suppose that $\fX$ is globally idyllic {\rm (\cite{egr1} 2.6.1)}. 
Then, the $\co_\fX$-algebra $\cA$ is topologically of finite presentation,  
if and only if the $\co_\fX(\fX)$-algebra $\cA(\fX)$ is topologically of finite presentation.
\end{itemize}
\end{lem}

(i) It follows from \ref{p1-pfs4} and (\cite{egr1} 2.2.8). 

(ii) It follows from (i) and definition \ref{p1-pfs1}(i) since $\Gamma(\fX,-)$ commutes with inverse limits.  

(iii) It follows from (i), (ii) and (\cite{egr1} 1.8.19). 

(iv) It follows from (i), (ii) and (\cite{egr1} 1.10.4).

\subsection{}\label{p1-pfs6}
We keep the assumptions and notation of \ref{p1-pfs4}. We assume, moreover, that the formal scheme $\fX$ is idyllic (\cite{egr1} 2.6.2) 
and that the $\co_\fX$-algebra $\cA$ is topologically of finite presentation. Since $f$ is of finite presentation, 
the formal scheme $\fY$ is also idyllic (\cite{egr1} 2.6.13). If $\fX$ is affine formal globally idyllic, then so is $\fY$ (\cite{egr1} 2.6.12). 
Recall that the sheaves of rings $\co_\fX$ and $\co_\fY$ are coherent (\cite{egr1} 2.8.1). 

The morphism $f$ induces a morphism of ringed topos 
\begin{equation}\label{p1-pfs6a}
\varphi\colon (\fY_\zar,\co_\fY)\rightarrow (\fX_\zar, \cA),
\end{equation}
where the subscript $_\zar$ denotes the Zariski topos; so we have $\varphi_*=f_*$. 
We denote by $\bMod^\coh(\co_\fY)$ the category of coherent $\co_\fY$-modules of $\fY_\zar$, and by $\bMod^\coh(\cA)$ the category of coherent $\cA$-modules of $\fX_\zar$.

\begin{prop}\label{p1-pfs7}
We keep the assumptions and notation of \ref{p1-pfs6} and let $\cF$ be a coherent $\co_\fY$-module. Then, 
\begin{itemize}
\item[{\rm (i)}] We have $\rR^if_*(\cF)=0$ for every $i\geq 1$. 
\item[{\rm (ii)}] The functor $f_*$ transforms any short exact sequence of $\co_\fY$-modules
$0\rightarrow \cM'\rightarrow \cM\rightarrow \cM''\rightarrow 0$ where $\cM''$ is coherent, into a short exact sequence. 
\item[{\rm (iii)}] The $\cA$-module $f_*(\cF)$ is of finite presentation. 
\item[{\rm (iv)}] The canonical morphisms 
\begin{eqnarray}
f_*(\cF)/\cJ f_*(\cF) &\rightarrow& f_*(\cF/\cK\cF),\label{p1-pfs8a}\\
\cJ f_*(\cF) &\rightarrow& f_*(\cK\cF),\label{p1-pfs8b}
\end{eqnarray}
are isomorphisms. 
\end{itemize}
\end{prop}

(i) Indeed, for any affine formal globally idyllic open subscheme $U$ of $\fX$, we have $\rH^i(f^{-1}(U),\cF)=0$ (\cite{egr1} 2.11.1). 
The proposition follows from (\cite{egr1} 2.6.15) and (\cite{sga4} V 5.1). 

(ii) It follows from (i). 

(iii) We may assume that $\fX$ is affine formal globally idyllic, and hence so is $\fY$. By (\cite{egr1} 2.7.2), there exists an exact sequence 
$\co_\fY^n\rightarrow \co_\fY^m\rightarrow \cF\rightarrow 0$. Since the category $\bMod^\coh(\co_\fY)$ is abelian, we deduce from (ii) by applying the functor $f_*$
that the sequence $\cA^n\rightarrow \cA^m\rightarrow f_*(\cF)\rightarrow 0$ is exact, which proves the proposition. 

(iv) By (ii), it is enough to prove that the first morphism in an isomorphism. 
We may assume that $\fX$ is affine formal globally idyllic, and hence so is $\fY$. By (\cite{egr1} 2.7.2), there exists an exact sequence 
$\co_\fY^n\rightarrow \co_\fY^m\rightarrow \cF\rightarrow 0$. By (ii), we are reduced to the case where $\cF=\co_\fY$, which is obvious (see \ref{p1-pfs4}). 

\begin{cor}\label{p1-pfs8}
Let $\fX$ be an idyllic formal scheme, $\cA$ an $\co_\fX$-algebra topologically of finite presentation. 
Then, $\cA$ is a coherent sheaf of rings. 
\end{cor}
It is enough to prove that if $\fX$ is affine formal globally idyllic, for every $\cA$-linear morphism $u\colon \cA^n\rightarrow \cA$, the $\cA$-module $\ker(u)$ is of finite type. 
Let $u_1,\dots,u_n$ be the sections of $\cA(\fX)$ defining $u$. We take the notation of \ref{p1-pfs6}.
Since $\co_\fY(\fY)=\cA(\fX)$, the $u_i$'s define a morphism $v\colon \co_\fY^n\rightarrow \co_\fY$. 
It is clear that $u=f_*(v)$. Since the $\co_\fY$-module $\ker(v)$ is of finite presentation, the $\cA$-module $\ker(u)=f_*(\ker(v))$ is of finite presentation by \ref{p1-pfs7}(iii), 
which gives the required assertion.

\begin{cor}\label{p1-pfs9}
Under the assumptions of \ref{p1-pfs6} and with the same notation, the morphism $\varphi$ \eqref{p1-pfs6a} 
induces equivalences of abelian categories, quasi-inverse to each other
\begin{equation}\label{p1-pfs9a}
\xymatrix{
{\bMod^\coh(\co_\fY)} \ar@<1ex>[r]^-(0.5){\varphi_*}& {\bMod^\coh(\cA)} \ar@<1ex>[l]^-(0.5){\varphi^*}}.
\end{equation}
\end{cor}

Observe first that the two functors are well defined by \ref{p1-pfs7}(iii) and \ref{p1-pfs8}. 
We claim that for every coherent $\co_\fY$-module $\cF$, the adjunction morphism $\varphi^*(\varphi_*(\cF))\rightarrow\cF$ is an isomorphism. 
Indeed, we may assume that $\fX$ is affine formal globally idyllic, and hence so is $\fY$. By (\cite{egr1} 2.7.2), there exists an exact sequence 
$\co_\fY^n\rightarrow \co_\fY^m\rightarrow \cF\rightarrow 0$. By \ref{p1-pfs7}(ii), we are reduced to the case where $\cF=\co_\fY$, which is obvious. 

We claim that for every coherent $\cA$-module $\cM$, the adjunction morphism $\cM\rightarrow \varphi_*(\varphi^*(\cM))$ is an isomorphism. 
Indeed, we may assume that there exists an exact sequence $\cA^n\rightarrow \cA^m\rightarrow \cM\rightarrow 0$. By \ref{p1-pfs7}(ii), 
we are reduced to the case where $\cM=\cA$, which is obvious.

\begin{cor}\label{p1-pfs15}
We keep the assumptions and notation of \ref{p1-pfs6}. Let $\cF$ be a coherent $\cA$-module, 
$\cF_\tor$ its torsion submodule defined in {\rm (\cite{egr1} 2.10.1.4)}. 
Then, $\cF_\tor$ is a coherent $\cA$-module. 
\end{cor}

Observe first that, although $\cF_\tor$ is defined in (\cite{egr1} 2.10.1.4) as an sub-$\co_\fX$-module of $\cF$, it is a sub-$\cA$-module. 
Indeed, the property being local, it is enough to check that for $\fX$ affine globally idyllic, $\cF_\tor(\fX)$ is a sub-$\cA(\fX)$-module of $\cF(\fX)$,
which follows from (\cite{egr1} 2.10.5) and the definition in (\cite{egr1} 1.8.30). 

The canonical morphism $\cF_\tor\rightarrow \varphi_*(\varphi^*(\cF)_\tor)$ is an isomorphism by (\cite{egr1} 2.10.5). 
The proposition follows then from \ref{p1-pfs9} and (\cite{egr1} 2.10.14).

\begin{cor}\label{p1-pfs13}
We keep the assumptions and notation of \ref{p1-pfs6}. Then, 
\begin{itemize}
\item[{\rm (i)}] For every coherent $\co_\fX$-modules $\cE$, the canonical morphism
\begin{equation}\label{p1-pfs13a}
\cE\otimes_{\co_\fX}\cA\rightarrow f_*(f^*(\cE))
\end{equation}
is an isomorphism.
\item[{\rm (ii)}] For every coherent $\co_\fY$-modules $\cF$ and $\cG$, the canonical morphism
\begin{equation}\label{p1-pfs13b}
f_*(\cF)\otimes_\cA f_*(\cG)\rightarrow f_*(\cF\otimes_{\co_\fY}\cG)
\end{equation}
is an isomorphism.
\end{itemize}
\end{cor}

(i) Indeed, the morphism \eqref{p1-pfs13a} can be identified with the adjunction morphism 
$\cE\otimes_{\co_\fX}\cA\rightarrow \varphi_*(\varphi^*(\cE\otimes_{\co_\fX}\cA))$, which is an isomorphism by \ref{p1-pfs9}. 
 
(ii) The question being local on $\fX$, we may assume that it is affine formal globally idyllic, and hence so is $\fY$.
Then, there exists an exact sequence of $\co_\fY$-modules $\co_\fY^n\rightarrow \co_\fY^m\rightarrow \cF\rightarrow 0$.  
By \ref{p1-pfs7}(ii), it is therefore enough to consider the case where $\cF=\co_\fY$, in which case the statement follows from \ref{p1-pfs7}(ii).  

\begin{cor}\label{p1-pfs10}
Let $\fX$ be an idyllic formal scheme, $\cA$ an $\co_\fX$-algebra topologically of finite presentation, $\cM$ a coherent $\cA$-module. 
Then $\cM$ is adic, i.e., for every ideal of definition of finite type $\cJ$ of $\fX$, the canonical morphism 
\begin{equation}\label{p1-pfs10a}
\cM\rightarrow \underset{\underset{n\geq 1}{\longleftarrow}}\lim\ \cM/\cJ^{n}\cM
\end{equation}
is an isomorphism. 
\end{cor}

Indeed, by \ref{p1-pfs9}, with the notation of \ref{p1-pfs6}, there exists a coherent $\co_\fY$-module $\cF$ such that $\cM=f_*(\cF)$. We set $\cK=f^*(\cJ)\co_\fY$. 
By (\cite{egr1} 2.8.5), the canonical morphism of $\co_\fY$-modules
\begin{equation}\label{p1-pfs9b}
\cF\rightarrow \underset{\underset{n\geq 1}{\longleftarrow}}\lim\ \cF/\cK^{n}\cF
\end{equation}
is an isomorphism. The proposition follows from \ref{p1-pfs7}(iv) and the fact that the functor $f_*$ commutes with inverse limits.

\subsection{}\label{p1-pfs140}
We take again the assumptions and notation of \ref{p1-pfs6}. Let $\cB$ be an $\co_\fX$-algebra topologically of finite presentation, 
$\rho\colon \cA\rightarrow \cB$ a homomorphism of $\co_\fX$-algebras. We set $\fZ=\Spf(\cB)$ and denote by $g\colon \fZ\rightarrow \fX$ the canonical morphism
and by $h\colon \fZ\rightarrow \fY$ the $\fX$-morphism induced by $\rho$. 
The morphism $g$ induces a morphism of ringed topos 
\begin{equation}\label{p1-pfs140a}
\psi\colon (\fZ_\zar,\co_\fZ)\rightarrow (\fX_\zar, \cB),
\end{equation}
where $\fZ_\zar$ denotes the Zariski topos of $\fZ$; so we have $\psi_*=g_*$. 
We immediately check that the following diagram
\begin{equation}\label{p1-pfs140b}
\xymatrix{
{\bMod(\co_\fY)}\ar[rr]^-(0.5){h^*}&&{\bMod(\co_\fZ)}\\
{\bMod(\cA)}\ar[u]^{\varphi^*}\ar[rr]^-(0.5){-\otimes_\cA\cB}&&{\bMod(\cB)}\ar[u]_{\psi^*}}
\end{equation}
is commutative. 

For every coherent $\cA$-modules $\cE$, the canonical morphism
\begin{equation}\label{p1-pfs140c}
\cE\otimes_{\cA}\cB\rightarrow g_*(h^*(\varphi^*(\cE)))
\end{equation}
is an isomorphism. Indeed, it can be identified with the adjunction morphism 
$\cE\otimes_{\cA}\cB\rightarrow \psi_*(\psi^*(\cE\otimes_{\cA}\cB))$, which is an isomorphism by \ref{p1-pfs9}.

\begin{prop}\label{p1-pfs14}
We keep the assumptions and notation of \ref{p1-pfs140}. Let $\cM$ be a coherent $\cB$-module, $\cJ$ an ideal of definition of finite type of $\fX$. 
Then, the following conditions are equivalent:
\begin{itemize} 
\item[{\rm (i)}] The $\cB$-module $\cM$ is $\cA$-flat. 
\item[{\rm (ii)}] The $\co_\fZ$-module $\psi^*(\cM)$ is $h$-flat.
\item[{\rm (iii)}] For every $n\geq 1$, the $(\cB/\cJ^n\cB)$-module $\cM/\cJ^n\cM$ is $(\cA/\cJ^n\cA)$-flat.
\end{itemize}
\end{prop}

Indeed, for every coherent $\cA$-module $\cN$, we have a canonical functorial isomorphism 
\begin{equation}
h^*(\varphi^*(\cN))\otimes_{\co_\fZ}\psi^*(\cM)\stackrel{\sim}{\rightarrow}\psi^*(\cN\otimes_\cA\cM).
\end{equation}
The equivalence of (i) and (ii) follows then from \ref{p1-pfs9} and (\cite{egr1} 1.3.17). The equivalence of (ii) and (iii) is a consequence of (\cite{egr1} 5.1.2).

\subsection{}\label{p1-pfs11}
We take again the assumptions and notation of \ref{p1-pfs6}. 
We assume, moreover, that $\fX$ admits a monogenic ideal of definition generated by $\mu\in \co_\fX(\fX)$. 
It follows from (\cite{egr1} 2.10.31.2) that for every $\co_\fY$-module $\cF$, the canonical morphism
\begin{equation}\label{p1-pfs11a}
f_*(\cF)[\frac 1 \mu]\rightarrow f_*(\cF[\frac 1 \mu])
\end{equation}
is an isomorphism. In particular, $f$ induces a morphism of ringed topos 
\begin{equation}\label{p1-pfs11b}
\phi\colon (\fY_\zar,\co_\fY[\frac 1 \mu])\rightarrow (\fX_\zar, \cA[\frac 1 \mu]),
\end{equation}
where the subscript $_\zar$ denotes the Zariski topos; so we have $\phi_*=f_*$. 
We denote by $\bMod^\coh(\co_\fY[\frac 1 \mu])$ the category of coherent $\co_\fY[\frac 1 \mu]$-modules of $\fY_\zar$, 
and by $\bMod^\coh(\cA[\frac 1 \mu])$ the full subcategory of coherent $\cA[\frac 1 \mu]$-modules of $\fX_\zar$.

\begin{prop}\label{p1-pfs12}
We keep the assumptions and notation of \ref{p1-pfs11}. Then, 
\begin{itemize}
\item[{\rm (i)}] For every coherent $\cA$-module $\cM$, $\cM[\frac 1 \mu]$ is a coherent $\cA[\frac 1 \mu]$-module; in particular, $\cA[\frac 1 \mu]$ is a coherent sheaf of rings. 
\item[{\rm (ii)}] If $\fX$ is quasi-compact, any coherent $\cA[\frac 1 \mu]$-module is of the form $\cM[\frac 1 \mu]$ for a coherent $\cA$-module $\cM$. 
\item[{\rm (iii)}] If $\fX$ is quasi-compact, any exact sequence of coherent $\cA[\frac 1 \mu]$-modules is the image of an exact sequence of coherent $\cA$-modules by localization by $\mu$. 
\item[{\rm (iv)}] The functor $f_*$ transforms short exact sequences of coherent $\co_\fY[\frac 1 \mu]$-modules into short exact sequences. 
\item[{\rm (v)}] The morphism $\phi$ \eqref{p1-pfs11b} 
induces equivalences of abelian categories, quasi-inverse to each other
\begin{equation}\label{p1-pfs12a}
\xymatrix{
{\bMod^\coh(\co_\fY[\frac 1 \mu])} \ar@<1ex>[r]^-(0.5){\phi_*}& {\bMod^\coh(\cA[\frac 1 \mu]).} \ar@<1ex>[l]^-(0.5){\phi^*}}
\end{equation}
\item[{\rm (vi)}] For every coherent $\co_\fX[\frac 1 \mu]$-module $\cE$, the canonical morphism
\begin{equation}\label{p1-pfs12b}
\cE\otimes_{\co_\fX}\cA\rightarrow f_*(f^*(\cE))
\end{equation}
is an isomorphism.
\item[{\rm (vii)}] For every coherent $\co_\fY[\frac 1 \mu]$-modules $\cF$ and $\cG$, the canonical morphism
\begin{equation}\label{p1-pfs12c}
f_*(\cF)\otimes_\cA f_*(\cG)\rightarrow f_*(\cF\otimes_{\co_\fY}\cG)
\end{equation}
is an isomorphism.
\end{itemize}
\end{prop}

The proofs of (i), (ii) and (iii) are identical to those of (\cite{egr1} 2.10.24). 

(iv) This follows from (iii) and \ref{p1-pfs7}(ii). 

(v) Recall that we have the functorial canonical isomorphism \eqref{p1-pfs11a}.  
Observe that the two functors in \eqref{p1-pfs12a} are well defined by (ii), \ref{p1-pfs9} and (\cite{egr1} 2.10.24(ii)). By the proof of \ref{p1-pfs9}, 
for every coherent $\co_\fY$-module $\cF$, the adjunction morphism $\varphi^*(\varphi_*(\cF))\rightarrow\cF$ is an isomorphism.  
We deduce by (\cite{egr1} 2.10.24(ii)) that for every coherent $\co_\fY[\frac 1 \mu]$-module $\cF$, the adjunction morphism $\phi^*(\phi_*(\cF))\rightarrow\cF$ is an isomorphism.  
Likewise, by the proof of \ref{p1-pfs9}, for every coherent $\cA$-module $\cM$, the adjunction morphism $\cM\rightarrow \varphi_*(\varphi^*(\cM))$ is an isomorphism. 
We deduce by (ii) that for every coherent $\cA[\frac 1 \mu]$-module $\cM$, the adjunction morphism $\cM\rightarrow \phi_*(\phi^*(\cM))$ is an isomorphism. 

The proofs of (vi) and (vii) are identical to those of \ref{p1-pfs13}.

\begin{prop}\label{p1-pfs16}
Let $V$ be a complete non-discrete valuation ring of height $1$, 
$\pi$ a non zero element of the maximal ideal of $V$. We equip $V$ with the $\pi$-adic topology and set $\cS=\Spf(V)$. 
Let $\fX$ be an $\cS$-formal scheme locally of finite presentation \eqref{p1-thbn0} (\cite{egr1} 2.3.15); 
so it is an idyllic formal scheme (\cite{egr1} 2.6.13).  
Let $\cA$ an $\co_\fX$-algebra topologically of finite presentation, $\cE$ a locally free $\cA$-module of finite type, 
$\cH=\rS_\cA(\cE)$ the symmetric $\cA$-algebra of $\cE$ \eqref{p1-NC7}, $\hcH=\hrS_\cA(\cE)$ its $\pi$-adic Hausdorff completion, 
$\cB$ an $\cA[\frac 1 \pi]$-algebra which is a coherent $\cA[\frac 1 \pi]$-module, 
$\psi\colon \cH \rightarrow \cB$ a homomorphism of $\cA$-algebras. Then, the following conditions are equivalent:
\begin{itemize}
\item[{\rm (i)}] The image of $\psi$ is a coherent $\cA$-module;
\item[{\rm (ii)}] The morphism $\psi$ extends to a homomorphism of $\cA$-algebras $\hpsi\colon \hcH\rightarrow \cB$.
\end{itemize}
\end{prop}

Indeed, we have (i)$\Rightarrow$(ii) since every coherent $\cA$-module is $\cJ$-adically complete and separated by \ref{p1-pfs10}.
To prove that (ii)$\Rightarrow$(i), we may assume that $\fX=\Spf(R)$, where $R$ is a $V$-algebra topologically of finite presentation, 
and $\cE$ is a free $\cA$-module of rank $d\geq 1$. We set $A=\cA(\fX)$, 
which is an $R$-algebra topologically of finite presentation by \ref{p1-pfs5}, and $E=\cE(\fX)$ which is a free $A$-module of rank $d$. 
Let $e_1,\dots,e_d$ be an $A$-basis of $E$.
By (\cite{egr1} 2.10.25), there exists an $\cA$-algebra $\cC$ which is a coherent $\cA$-module such that $\cB=\cC[\frac 1 \pi]$.  
We put $C=\cC(\fX)$, which is an $A[\frac 1 \pi]$-algebra and a coherent $A$-module, and $B=C[\frac 1 \pi]=\cB(\fX)$. 
We may assume that $C$ is $\pi$-torsion free by \ref{p1-pfs15}. 
By assumption, we have a commutative diagram 
\begin{equation}
\xymatrix{
{A[T_1,\dots,T_d]}\ar[r]^-(0.5){\varphi}\ar[d]_\iota&B\\
{A\{T_1,\dots,T_d\}}\ar[ru]_-(0.5){\hvarphi}&}
\end{equation}
where $\varphi$ (resp.\ $\hvarphi$) is the evaluation of $\psi$ (resp.\ $\hpsi$) on $\fX$ and $\iota$ is the canonical homomorphism. 

We equip $A$, $A\{T_1,\dots,T_d\}$ and $C$ with the $\pi$-adic topology and 
$B$ with the $\pi$-adic topology defined by $C$ (i.e., the unique topology compatible with its additive group structure, for which
the $A$-modules $\pi^nC$, for $n\in \mN$, form a fundamental system of neighborhoods of $0$). 
Therefore, the homomorphism $\hvarphi$ is continuous by (\cite{bgr} 6.1.3/1), and for every $1\leq i\leq d$, $\varphi(T_i)$ is power bounded in $B$. 
Since $B$ is finite over $A[\frac 1 p]$, we deduce by (\cite{bgr} 6.3.4/1) that $\varphi(T_i)$ is integral over $A$. 
Hence, the $A$-module $\im(\varphi)$ is of finite type. There exists an integer $n\geq 1$ such that $\im(\varphi)\subset \pi^{-m} C$. 
As the $A$-module $\pi^{-m} C$ is coherent, we deduce that the $A$-module $\im(\varphi)$ is coherent. 
The coherent $\cA$-module associated with $\im(\varphi)$ by \ref{p1-pfs9} and (\cite{egr1} 2.7.2), is clearly the image of $\psi$, 
which proves the required implication.

\section{\texorpdfstring{Higgs modules and $\delta$-connections}{Higgs modules and delta-connections}}\label{p1-delta-con}

In this section, we fix a ringed $\mU$-topos $(X,A)$ \eqref{p1-NC7}. 

\begin{defi}\label{p1-delta-con1}
Let $\Omega$ be an $A$-module.
\begin{itemize}
\item[(i)]
A {\em Higgs $A$-module with coefficients in $\Omega$} is 
a pair $(M,\theta)$ consisting of an  $A$-module $M$ and an  $A$-linear morphism
\begin{equation}\label{p1-MH1a}
\theta\colon M\rightarrow M\otimes_A\Omega
\end{equation}
such that $\theta\wedge \theta=0$. We then say that $\theta$ is a {\em  Higgs $A$-field} on $M$
with coefficients in $\Omega$.
\item[(ii)] If $(M,\theta)$ and $(M',\theta')$ are two Higgs $A$-modules with coefficients in $\Omega$,
a morphism from $(M,\theta)$ to $(M',\theta')$ is an  $A$-linear morphism
$u\colon M\rightarrow M'$ such that $(u\otimes\id_\Omega)\circ \theta=\theta'\circ u$.
\end{itemize}
\end{defi}

Higgs $A$-modules with coefficients in $\Omega$ form a category that we denote by $\bHM(A,\Omega)$. 
We can complete the terminology and make the following remarks.

\addtocounter{subsubsection}{1}
\addtocounter{equation}{1}

\subsubsection{}\label{p1-delta-con1b}
Let $(M,\theta)$ be a  Higgs $A$-module with coefficients in $\Omega$. For each $i\geq 1$, we denote by
\begin{equation}\label{p1-delta-con1c}
\theta^i\colon M\otimes_A \wedge^i\Omega \rightarrow M\otimes_A \wedge^{i+1}\Omega
\end{equation}
the $A$-linear morphism defined for all local sections
$m$ of $M$ and $\omega$ of $\wedge^i\Omega$ by $\theta^i(m\otimes \omega)=\theta(m)\wedge \omega$ \eqref{p1-NC7}.
We have $\theta^{i+1}\circ \theta^i=0$. The {\em Dolbeault} complex of $(M,\theta)$,
denoted by $\mK^\bullet(M,\theta)$, is the complex of cochains of $A$-modules
\begin{equation}\label{p1-delta-con1d}
M\stackrel{\theta}{\longrightarrow}M\otimes_A\Omega\stackrel{\theta^1}{\longrightarrow} M\otimes_A\wedge^2\Omega \dots,
\end{equation}
where $M$ is placed in degree $0$ and the differentials are of degree $1$.

\addtocounter{subsubsection}{2}
\addtocounter{equation}{1}

\subsubsection{}\label{p1-delta-con1g}
Let $(M,\theta),(M',\theta')$ be two  Higgs $A$-modules with coefficients in $\Omega$.
The {\em total} Higgs field on $M\otimes_AM'$ is the $A$-linear morphism
\begin{equation}\label{p1-delta-con1h}
\theta_\tot\colon M\otimes_AM'\rightarrow M\otimes_AM'\otimes_A\Omega
\end{equation}
defined by
\begin{equation}\label{p1-delta-con1i}
\theta_\tot=\theta\otimes \id_{M'}+\id_{M}\otimes \theta'.
\end{equation}
We say that $(M\otimes_AM',\theta_\tot)$ is the {\em tensor product} of $(M,\theta)$ and $(M',\theta')$.

\addtocounter{subsubsection}{2}
\addtocounter{equation}{1}

\subsubsection{}\label{p1-delta-con1e}
Let $(M,\theta)$ be a Higgs $A$-module with coefficients in $\Omega$ with $M$ a locally projective $A$-module
of finite type (\cite{ag2} 2.1.11). Consider, for an integer $i\geq 1$, the composed morphism
\begin{equation}\label{p1-delta-con1f}
\xymatrix{
{\wedge^iM}\ar[r]^-(0.5){\wedge^i\theta}&{\wedge^i(M\otimes_A\Omega)}\ar[r]&
{\wedge^iM\otimes_A\rS^i\Omega,}}
\end{equation}
where the second arrow is the canonical morphism (\cite{illusie1} V 4.5).
Observe that for any locally projective $A$-module of finite type $P$ and any $A$-module $Q$, the canonical morphism
\begin{equation}\label{p1-delta-con1ef}
\cHom_A(P,A)\otimes_AQ\rightarrow \cHom_A(P,Q)
\end{equation}
is an isomorphism. We can then define the trace of the morphism \eqref{p1-delta-con1f} as a section of $\Gamma(X,\rS^i\Omega)$. 
We call it the {\em $i$th characteristic invariant} of $(M,\theta)$ (or of $\theta$), and denote it by $\lambda_i(\theta)$.

\addtocounter{subsubsection}{2}
\addtocounter{equation}{1}

\subsubsection{}\label{p1-delta-con1j}
Suppose that $\Omega$ is locally projective of finite type over $A$ (\cite{ag2} 2.1.11).
Let $\rT=\cHom_A(\Omega,A)$ be its dual, and $\rS_A(\rT)$ the symmetric $A$-algebra of $\rT$  \eqref{p1-NC7}.
Observe that the $A$-module $T$ is locally projective of finite type and the canonical morphism 
$\Omega\rightarrow \cHom_A(\rT,A)$ is an isomorphism. 
For every $A$-module $M$, the canonical morphism
\begin{equation}\label{p1-delta-con1k}
\cEnd_A(M)\otimes_A\Omega\rightarrow \cHom_A(M,M\otimes_A\Omega)
\end{equation}
is an isomorphism. Moreover, by \eqref{p1-delta-con1ef}, the canonical morphism 
\begin{equation}\label{p1-delta-con1l}
\cHom_A(\rT,A) \otimes_A\cEnd_A(M)\rightarrow \cHom_A(\rT,\cEnd_A(M))
\end{equation}
is an isomorphism. Hence, giving an $A$-linear morphism 
$\theta\colon M\rightarrow M\otimes_A\Omega$ is equivalent to giving an $A$-linear morphism $\mu\colon \rT\rightarrow \cEnd_A(M)$. 
Working locally in $X$ and embedding $\Omega$ as a direct factor of a free $A$-module of finite type, 
we see that $\theta$ is a Higgs $A$-field 
if and only if the image of $\mu$ consists of endomorphisms of $M$ that commute to each other, i.e., 
giving a Higgs $A$-field $\theta$ on $M$ is equivalent to giving an $\rS_A(\rT)$-module structure on $M$ that is
compatible with its $A$-module structure.

\begin{lem}\label{p1-delta-con7}
Let $\Omega$ be a locally projective $A$-module of finite type, $\rT=\cHom_A(\Omega,A)$,
$(M,\theta)$, $(M',\theta')$ two Higgs $A$-modules with coefficients in $\Omega$. Then, 
\begin{itemize}
\item[{\rm (i)}] The Higgs field
\begin{equation}\label{p1-delta-con7a}
\theta\colon M\rightarrow M\otimes_A\Omega 
\end{equation}
is $\rS_A(\rT)$-linear \eqref{p1-delta-con1j}. 
\item[{\rm (ii)}] The $A$-module underlying $M\otimes_{\rS_A(\rT)}M'$ is the maximal quotient of 
$M\otimes_AM'$ on which the Higgs $A$-fields $\theta\otimes\id_{M'}$ and $\id_M\otimes\theta'$ induce the same Higgs $A$-field 
with coefficients in $\Omega$, namely its canonical Higgs field.
\end{itemize}
\end{lem}

(i) Indeed, the Higgs $A$-field $\theta$ is the composition of the morphisms 
\begin{equation}\label{p1-delta-con7b}
M\rightarrow M\otimes_A\rT\otimes_A \Omega\rightarrow M\otimes_A\rS_A(\rT)\otimes_A \Omega
\rightarrow M\otimes_A \Omega,
\end{equation}
the first being induced by the section $\id_\Omega$ of $\rT\otimes_A \Omega$ identified with $\cHom_A(\Omega,\Omega)$ \eqref{p1-delta-con1ef},
the second being induced by the canonical injection $\rT\rightarrow \rS_A(\rT)$ and the third being induced by $\upmu$. 
Each of these morphisms is obviously $\rS_A(\rT)$-linear. 

(ii) The question being local, we may assume that $\Omega$ is a direct factor of a free $A$-module of finite type $\Omega'$. 
Then $\theta$ (resp.\ $\theta'$) is determined by its coordinates $\theta_i\in \End_A(M)$ (resp.\  $\theta'_i\in \End_A(M')$) relatively to a 
basis of $\Omega'$, which commute to each other. The $A$-module underlying $M\otimes_{\rS_A(\rT)}M'$ 
is the quotient of $M\otimes_AM'$ by the images of $\theta_i\otimes\id_{M'}-\id_M\otimes\theta_i$, for all $i$, which implies the proposition. 

\begin{defi}\label{p1-delta-con6}
Let $\Omega$ be an $A$-module.
A {\em Higgs $A$-bundle with coefficients in $\Omega$}
is a Higgs $A$-module $(M,\theta)$ with coefficients in $\Omega$
such that the $A$-module $M$ is locally projective of finite type (\cite{ag2} 2.1.11).
\end{defi}

\begin{defi}\label{p1-delta-con2} 
Let $B$ be a commutative $A$-algebra, $\Omega$ an $A$-module, $\delta\colon B\rightarrow \Omega\otimes_AB$ an $A$-derivation,
which is also a Higgs $A$-field on $B$ with coefficients in $\Omega$.
\begin{itemize}
\item[(i)] A {\em $B$-module with $\delta$-connection} is a pair $(M,\nabla)$ consisting of a $B$-module $M$ and an $A$-linear morphism
\begin{equation}\label{p1-delta-con2a}
\nabla\colon M\rightarrow \Omega\otimes_AM
\end{equation}
such that for every local sections $\xi$ of $B$ and $m$ of $M$, we have
\begin{equation}\label{p1-delta-con2b}
\nabla(\xi m)=\delta(\xi)\otimes_Bm+\xi\nabla(m).
\end{equation} 
We then say that $\nabla$ is a {\em $\delta$-connection} on $M$. 
The $\delta$-connection $\nabla$ is said to be {\em integrable} if it is a Higgs $A$-field on $M$ with coefficients in $\Omega$. 
\item[(ii)] Let $(M,\nabla)$, $(M',\nabla')$ be two $B$-modules with $\delta$-connection.
A morphism from $(M,\nabla)$ to $(M',\nabla')$ is a $B$-linear morphism $u\colon M\rightarrow M'$
such that $(\id \otimes u)\circ \nabla=\nabla'\circ u$.
\end{itemize}
\end{defi}

The $B$-modules with $\delta$-connection form a category that we denote by $\bMC_\delta(B/A)$. We denote by $\bMIC_\delta(B/A)$ the full subcategory 
made of $B$-modules with integrable $\delta$-connection.

\subsection{}\label{p1-delta-con4}
We take again the assumption and notation of \ref{p1-delta-con2}, and 
let $(M,\nabla)$ be a $B$-module with $\delta$-connection, $(N,\theta)$ a Higgs $A$-module with coefficients in $\Omega$. 
Then, the $A$-linear morphism
\begin{equation}\label{p1-delta-con4a}
\nabla'\colon M\otimes_A N\rightarrow \Omega\otimes_A M\otimes_A N
\end{equation}
defined by
\begin{equation}\label{p1-delta-con4b}
\nabla'=\nabla\otimes_A \id_N+ \id_M\otimes_A \theta,
\end{equation}
is a $\delta$-connection on $M\otimes_AN$. If $\nabla$ is integrable, so is $\nabla'$.

\subsection{}\label{p1-delta-con5} 
We take again the assumption and notation of \ref{p1-delta-con2}, and let $B'$ be a commutative $A$-algebra, $\Omega'$ an $A$-module, 
$\delta'\colon B'\rightarrow \Omega'\otimes_AB'$ an $A$-derivation, which is also a Higgs $A$-field on $B'$ with coefficients in $\Omega'$.
Let $\iota\colon B\rightarrow B'$ be a homomorphism of $A$-algebras, 
$u\colon \Omega\rightarrow \Omega'$ an $A$-linear morphism such that $\delta'\circ \iota=(u\otimes \iota) \circ \delta$. 
Let $(M,\nabla)$ be a $B$-module with $\delta$-connection. Then, there exists a unique $A$-linear morphism
\begin{equation}\label{p1-delta-con5a}
\nabla'\colon M\otimes_BB'\rightarrow \Omega'\otimes_AM\otimes_BB'
\end{equation}
such that for all local sections $m$ of $M$ and $\xi'$ of $B'$, we have
\begin{equation}\label{p1-delta-con5b}
\nabla'(m\otimes_B \xi')=(u\otimes \id)(\nabla(m))\otimes_B \xi'+ m\otimes_B \delta'(\xi').
\end{equation}
It is a $\delta'$-connection on $M\otimes_BB'$. If $\nabla$ is integrable, so is $\nabla'$.

\subsection{}\label{p1-delta-con8}
We take again the assumption and notation of \ref{p1-delta-con2}. A {\em $\delta$-isoconnection} with respect to the extension $B/A$
(or simply a {\em $\delta$-isoconnection} when there is no risk of confusion) is a quadruple
\begin{equation}\label{p1-delta-con8a}
(M,N,u\colon M\rightarrow N,\nabla\colon M\rightarrow \Omega\otimes_AN),
\end{equation}
where $M$ and $N$ are $B$-modules, $u$ an isogeny of $B$-modules (\cite{ag2} 2.9.1), 
and $\nabla$ is an $A$-linear morphism such that for all local sections
$x$ of $B$ and $t$ of $M$, we have
\begin{equation}\label{p1-delta-con8b}
\nabla(xt)=\delta(x) \otimes u(t)+x\nabla(t).
\end{equation}
For every $B$-linear morphism $v\colon N\rightarrow M$ for which there exists an integer $n$
such that $u\circ v=n\cdot \id_N$ and $v\circ u=n\cdot \id_M$,
the pairs $(M,(\id\otimes v)\circ \nabla)$ and $(N,\nabla\circ v)$ are $B$-modules with $(n\delta)$-connections \eqref{p1-delta-con2},
and $u$ is a morphism from $(M,(\id\otimes v)\circ \nabla)$  to $(N,\nabla\circ v)$. 

The $\delta$-isoconnection $(M,N,u,\nabla)$ is said to be {\em integrable} if 
there exist a $B$-linear morphism $v\colon N\rightarrow M$ and an integer $n\not= 0$
such that $u\circ v=n\cdot \id_N$, $v\circ u=n\cdot \id_M$ and the $(n\delta)$-connections $(\id\otimes v)\circ \nabla$ on $M$ and $\nabla\circ v$
on $N$ are integrable in the sense of \ref{p1-delta-con2}(i).

Let $(M,N,u,\nabla)$, $(M',N',u',\nabla')$ be two $\delta$-isoconnections.
A morphism from $(M,N,u,\nabla)$  to $(M',N',u',\nabla')$ consists of
two $B$-linear morphisms $\alpha\colon M\rightarrow M'$ and $\beta\colon N\rightarrow N'$
such that $\beta\circ u=u'\circ \alpha$ and $(\id \otimes \beta)\circ \nabla=\nabla'\circ \alpha$.

We denote by $\bIMC_\delta(B/A)$ the category of $\delta$-isoconnections and by 
$\bIMIC_\delta(B/A)$ the subcategory of integrable $\delta$-isoconnections.

\begin{exemple}\label{p1-delta-con3}
Let $B$ be a commutative $A$-algebra, $d\colon B\rightarrow \Omega^1_{B/A}$ the universal $A$-derivation, $\Omega_{B/A}=\wedge_B\Omega^1_{B/A}$.  
We denote also by $d\colon \Omega_{B/A}\rightarrow \Omega_{B/A}$ the unique $A$-anti-derivation of degree $1$ and square zero that extends $d$. 
We suppose that there exists an $A$-module $\Omega$ and a $B$-isomorphism $\upgamma\colon \Omega\otimes_AB\stackrel{\sim}{\rightarrow}\Omega^1_{B/A}$
such that for every local section $\omega$ of $\Omega$, we have $d(\upgamma(\omega\otimes 1))=0$. 
Let $\lambda\in A(X)$. Consider the $A$-derivation $\delta=\lambda \upgamma^{-1}\circ  d\colon B\rightarrow \Omega\otimes_AB$, which is a Higgs $A$-field with coefficients in $\Omega$. 
Let $M$ be a $B$-module. An $A$-linear morphism 
\begin{equation}\label{p1-delta-con3a}
\nabla\colon M\rightarrow \Omega\otimes_AM
\end{equation}
is a $\delta$-connection (resp.\ an integrable $\delta$-connection) if and only if 
\begin{equation}\label{p1-delta-con3b}
\tnabla=(\upgamma\otimes \id_M)\circ \nabla\colon M\rightarrow \Omega^1_{B/A}\otimes_BM
\end{equation}
is a $\lambda$-connection (resp.\ an integrable $\lambda$-connection) with respect to the extension $B/A$ in the sense of (\cite{ag2} 2.5.24). 
Indeed, the diagram 
\begin{equation}
\xymatrix{
{\Omega\otimes_AM}\ar[r]^-(0.5){-\id\wedge \nabla}\ar[d]&{(\wedge^2\Omega)\otimes_AM}\ar[d]\\
{\Omega^1_{B/A}\otimes_BM}\ar[r]^-(0.4){\tnabla}&{\Omega^2_{B/A}\otimes_BM,}}
\end{equation}
where the vertical maps are the isomorphisms induced by $\upgamma$, is clearly commutative. 
\end{exemple}

\section{Abelian categories up to isogeny and indization}\label{p1-abisoind}

\subsection{}\label{p1-abisoind1}
In this section, we fix an abelian $\mU$-category $\cA$ \eqref{p1-NC0}. 
We denote by $\cA_\mQ$ the category of  objects of $\cA$ up to isogeny (\cite{ag2} 2.9.1) and by 
\begin{equation}\label{p1-abisoind1a}
Q\colon 
\cA\rightarrow\cA_\mQ,\ \ \ X\mapsto X_\mQ,
\end{equation}
the localization functor; the category $\cA_\mQ$ is abelian and the functor $Q$ is exact (\cite{ag2} 2.9.2). 
We denote by $\bInd(\cA)$ the category of ind-objects of $\cA$ (\cite{ag2} 2.6.3) and by 
\begin{equation}\label{p1-abisoind1b}
\iota_\cA\colon \cA\rightarrow \bInd(\cA)
\end{equation}
the canonical functor; the category $\bInd(\cA)$ is abelian and the functor $\iota_\cA$ is fully faithful and exact (\cite{ag2} 2.6.7). 

For any object $X$ of $\cA$, we denote by $I^X$ the category of isogenies of source $X$ and by
\begin{equation}\label{p1-abisoind1c}
\alpha^X\colon I^X\rightarrow \cA, \ \ \ (X\rightarrow X')\mapsto X',
\end{equation}
the ``target'' functor (see \cite{ag2} 2.6.6 and 2.9.2). The category $I^X$ being filtered and essentially small, 
we can therefore consider the fully faithful functor (\cite{ag2} (2.9.2.3))
\begin{equation}\label{p1-abisoind1d}
\upalpha_\cA\colon 
\cA_\mQ\rightarrow \bInd(\cA),\ \ \ X\mapsto \indcolim \alpha^X,
\end{equation}
which is exact by (\cite{ag2} 2.9.3(ii)). We have a canonical morphism  (\cite{ag2} (2.9.2.4))
\begin{equation}\label{p1-abisoind1e}
\iota_\cA\rightarrow \upalpha_\cA\circ Q,
\end{equation}
which is not an isomorphism in general.

\begin{prop}\label{p1-abisoind2}
The functor $\upalpha_\cA$ \eqref{p1-abisoind1d} makes $\cA_\mQ$ a thick subcategory of $\bInd(\cA)$. 
\end{prop}

It is enough to prove that $\cA_\mQ$ is stable by extension in $\bInd(\cA)$. Let $M,N$ be objects of $\cA$ and 
\begin{equation}\label{p1-abisoind2a}
0\rightarrow \upalpha_{\cA}(M_\mQ)\rightarrow E \rightarrow \upalpha_{\cA}(N_\mQ)\rightarrow 0
\end{equation}
an exact sequence of $\bInd(\cA)$. 
By (\cite{ks2} 8.6.9), there exists an epimorphism $\varphi\colon L\rightarrow N$ of $\cA$ and a commutative diagram 
\begin{equation}
\xymatrix{
{\iota_\cA(L)}\ar[r]^{\iota_\cA(\varphi)}\ar[d]_v&{\iota_\cA(N)}\ar[d]\\
E\ar[r]&{\upalpha_{\cA}(N_\mQ),}}
\end{equation} 
where the right vertical arrow is the canonical morphism \eqref{p1-abisoind1e}.
Since $n\cdot \id_E$ is an isomorphism of $E$ \eqref{p1-abisoind2a}, by (\cite{ag2} 2.9.3(i)), 
$v$ induces a morphism $\nu\colon \upalpha_\cA(L_\mQ)\rightarrow E$ that fits into a commutative diagram 
\begin{equation}
\xymatrix{
{\upalpha_\cA(L_\mQ)}\ar[rd]^{\upalpha_\cA(\varphi_\mQ)}\ar[d]_\nu&\\
E\ar[r]&{\upalpha_{\cA}(N_\mQ).}}
\end{equation}
Since $\upalpha_\cA$ is exact and fully faithful, setting $K=\ker(\varphi)$, 
$\nu$ induces a morphism $u\colon K_\mQ \rightarrow M_\mQ$ such that the diagram 
\begin{equation}\label{p1-abisoind2b}
\xymatrix{
0\ar[r]&{\upalpha_\cA(K_\mQ)}\ar[r]\ar[d]_{\upalpha_\cA(u)}&{\upalpha_\cA(L_\mQ)}\ar[r]\ar[d]_\nu&{\upalpha_\cA(N_\mQ)}\ar[r]\ar@{=}[d]&0\\
0\ar[r]&{\upalpha_\cA(M_\mQ)}\ar[r]&E\ar[r]&{\upalpha_\cA(N_\mQ)}\ar[r]&0}
\end{equation}
is commutative. 
Let $\cL$ be the pushout in $\cA_\mQ$ of the extension $0\rightarrow K_\mQ\rightarrow L_\mQ\rightarrow N_\mQ\rightarrow 0$ 
by the the morphism $u$. 
Then \eqref{p1-abisoind2b} induces an isomorphism $\upalpha_\cA(\cL)\stackrel{\sim}{\rightarrow} E$, which proves the required assertion.

\begin{lem}\label{p1-abisoind5}
Let $\rho\colon I\rightarrow \cC$ be a functor, $\varphi\colon \cC\rightarrow \cC'$ a fully faithful functor. 
Assume that the direct limit of $\varphi \circ \rho$ is representable in $\cC'$ and that there exist an
object $X$ of $\cC$ and an isomorphism of $\cC'$ 
\begin{equation}\label{p1-abisoind5a}
\nu\colon \varphi(X)\stackrel{\sim}{\rightarrow} \underset{\underset{I}{\longrightarrow}}{\lim} \ \varphi\circ \rho. 
\end{equation}
Then, the direct limit of $\rho$ is representable in $\cC$, and $\nu$ induces 
an isomorphism of $\cC$
\begin{equation}\label{p1-abisoind5b}
\nu^\flat \colon  \underset{\underset{I}{\longrightarrow}}{\lim} \ \rho  \stackrel{\sim}{\rightarrow} X. 
\end{equation}
In particular, we have an isomorphism of $\cC'$
\begin{equation}\label{p1-abisoind5c}
\nu \circ \varphi(\nu^\flat) \colon  \varphi(\underset{\underset{I}{\longrightarrow}}{\lim} \ \rho)  \stackrel{\sim}{\rightarrow} 
\underset{\underset{I}{\longrightarrow}}{\lim} \ \varphi\circ \rho. 
\end{equation}
\end{lem}

Indeed, $\nu$ induces a morphism of functors $\lambda\colon \varphi \circ \rho\rightarrow \varrho_{\varphi(X)}$, where
$\varrho_{\varphi(X)}\colon I\rightarrow \cC'$ is the constant functor with value $\varphi(X)$. The direct limit of $\varphi \circ \rho$
is representable by $\varphi(X)$ via $\lambda$. 
By the full faithfulness of $\varphi$, there exists a unique morphism of functors $\mu\colon \rho\rightarrow \varrho_{X}$, where
$\varrho_{X}\colon I\rightarrow \cC$ is the constant functor with value $X$, such that $\varphi\circ \mu =\lambda$. 
We see that $X$ represents the direct limit of $\rho$ via $ \mu$. Indeed, let $Y\in \ob(\cC)$, $u\colon \rho\rightarrow  \varrho_Y$ a morphism of functors, 
where $\varrho_Y\colon I\rightarrow \cC$ is the constant functor with value $Y$. By the full faithfulness of $\varphi$,
giving a morphism $v\colon X\rightarrow Y$ of $\cC$
such that $u=\varrho_{v}\circ \mu$ is equivalent to giving a morphism $w\colon \varphi(X)\rightarrow \varphi(Y)$ of $\cC'$
such that $\varphi\circ u=\varrho_{w}\circ \lambda$, which proves the assertion.

\begin{cor}\label{p1-abisoind6}
Let $\rho\colon I\rightarrow \cC$ be a functor, where $I$ is a small filtered category, $X\in \ob(\cC)$,
\begin{equation}\label{p1-abisoind6a}
\nu\colon \iota_\cC(X)\stackrel{\sim}{\rightarrow} \underset{\underset{I}{\longrightarrow}}{\mlq\mlq\lim \mrq\mrq} \iota_\cC\circ \rho, 
\end{equation}
an isomorphism of $\bInd(\cC)$, where $\iota_\cC\colon \cC\rightarrow \bInd(\cC)$ is the canonical functor {\rm (\cite{ag2} 2.6.3)}. 
Then, the direct limit of $\rho$ is representable in $\cC$, and $\nu$ induces 
an isomorphism of $\cC$
\begin{equation}\label{p1-abisoind6b}
\nu^\flat \colon  \underset{\underset{I}{\longrightarrow}}{\lim} \ \rho  \stackrel{\sim}{\rightarrow} X. 
\end{equation}
In particular, we have an isomorphism of $\bInd(\cC)$
\begin{equation}\label{p1-abisoind6c}
\nu \circ \iota_\cC(\nu^\flat) \colon  \iota_\cC(\underset{\underset{I}{\longrightarrow}}{\lim} \ \rho)  \stackrel{\sim}{\rightarrow} 
\underset{\underset{I}{\longrightarrow}}{\mlq\mlq\lim \mrq\mrq} \iota_\cC\circ \rho. 
\end{equation}
\end{cor}

\begin{rema}\label{p1-abisoind9}
Let $\rho\colon I\rightarrow \cC$ be a functor, where $I$ is a small filtered category. Then, the following conditions are equivalent: 
\begin{itemize}
\item[(i)] There exist $X\in \ob(\cC)$ and an isomorphism of $\bInd(\cC)$
\begin{equation}\label{p1-abisoind8a}
\nu\colon \iota_\cC(X)\stackrel{\sim}{\rightarrow} \underset{\underset{I}{\longrightarrow}}{\mlq\mlq\lim \mrq\mrq} \iota_\cC\circ \rho. 
\end{equation}
\item[(ii)] There exist $X\in \ob(\cC)$, a morphism $\nu^\flat_\bullet\colon \rho\rightarrow \varrho_X$ of functors from $I$ to $\cC$, 
where $\varrho_X$ is the constant functor with value $X$, $i\in \ob(I)$ and a morphism $\nu_i\colon X\rightarrow \rho(i)$ of $\cC$, 
such that the following conditions are satisfied:
\begin{itemize}
\item[(a)] The composed morphism $\nu^\flat_i\circ \nu_i\colon X\rightarrow \rho(i)\rightarrow X$ is the identity of $X$.
\item[(b)] For every $j\in \ob(I)$, there exist $k\in \ob(I)$ and two morphisms $u\colon i\rightarrow k$ and $v\colon j\rightarrow k$ of $I$
such that $\rho(v)$ is the composed morphism 
\begin{equation}
\xymatrix{
{\rho(j)}\ar[r]^-(0.5){\nu^\flat_j}&X\ar[r]^-(0.5){\nu_i}&{\rho(i)}\ar[r]^-(0.5){\rho(u)}&{\rho(k).}}
\end{equation}
\end{itemize}
\end{itemize}
The equivalence follows immediately from (\cite{ag2} (2.6.1.2)).
When these conditions are satisfied, we say that $\rho$ is {\em essentially constant, as an ind-object, with value $X$} (see \cite{sp}  \href{https://stacks.math.columbia.edu/tag/05PU}{05PU}).
\end{rema}

\begin{cor}\label{p1-abisoind7}
Let $I$ be a small filtered category, $\rho\colon I\rightarrow \cA_\mQ$ a functor, $X\in \ob(\cA_\mQ)$ 
\begin{equation}\label{p1-abisoind7a}
\nu\colon \upalpha_\cA(X)\stackrel{\sim}{\rightarrow} \underset{\underset{I}{\longrightarrow}}{\mlq\mlq\lim \mrq\mrq} \upalpha_\cA\circ \rho 
\end{equation}
an isomorphism of $\bInd(\cA)$. 
Then, the direct limit of $\rho$ is representable in $\cA_\mQ$, and $\nu$ induces 
an isomorphism of $\cA_\mQ$
\begin{equation}\label{p1-abisoind7b}
\nu^\flat \colon  \underset{\underset{I}{\longrightarrow}}{\lim} \ \rho  \stackrel{\sim}{\rightarrow} X. 
\end{equation}
In particular, we have an isomorphism of $\bInd(\cA)$
\begin{equation}\label{p1-abisoind7c}
\nu \circ \upalpha_\cA(\nu^\flat) \colon  \upalpha_\cC(\underset{\underset{I}{\longrightarrow}}{\lim} \ \rho)  \stackrel{\sim}{\rightarrow} 
\underset{\underset{I}{\longrightarrow}}{\mlq\mlq\lim \mrq\mrq} \upalpha_\cA\circ \rho. 
\end{equation}
\end{cor}

\subsection{}\label{p1-abisoind3}
We assume that {\em small filtered direct limits are representable in $\cA$}.
By (\cite{ks2} 6.3.1), the functor $\iota_\cA$ \eqref{p1-abisoind1b} admits a left adjoint 
\begin{equation}\label{p1-abisoind3a}
\kappa_\cA \colon \bInd(\cA)\rightarrow \cA. 
\end{equation}
For every small filtered  category $J$ and every functor $\alpha\colon J\rightarrow \cA$, we have an isomorphism
\begin{equation}\label{p1-abisoind3b}
\kappa_\cA(\underset{\underset{J}{\longrightarrow}}{\mlq\mlq\lim \mrq\mrq} \alpha)\stackrel{\sim}{\rightarrow} \underset{ \underset{J}{\longrightarrow}}{\lim}\ \alpha.
\end{equation} 
Moreover, the canonical morphism $\kappa_\cA\circ \iota_\cA \rightarrow \id_{\cA}$ is an isomorphism. 

We consider the functor
\begin{equation}\label{p1-abisoind3c}
-\otimes \mQ\colon \cA\rightarrow \cA,\ \ \ X\mapsto \underset{\underset{I^X}{\longrightarrow}}{\lim}\ \alpha^X,
\end{equation}
where $\alpha^X$ is the functor \eqref{p1-abisoind1c}. For every $X\in \ob(\cA)$, we have a canonical functorial isomorphism 
\begin{equation}\label{p1-abisoind3d}
X\otimes \mQ \stackrel{\sim}{\rightarrow} \kappa_\cA(\upalpha_\cA(X_\mQ)),
\end{equation}
and a canonical functorial morphism $X\rightarrow X\otimes \mQ$.

\begin{lem}\label{p1-abisoind4}
If small filtered direct limits are representable in $\cA$, and are exact, 
then the functor $\kappa_\cA$ \eqref{p1-abisoind3a} is exact and it commutes with small direct limits.
\end{lem}
Indeed, $\kappa_\cA$ commutes with small direct limits since it admits a right adjoint.
Let
\begin{equation}\label{p1-abisoind4a}
0\rightarrow F'\rightarrow F\rightarrow F''\rightarrow 0
\end{equation}
be an exact sequence of $\bInd(\cA)$. By (\cite{ks2} 8.6.6), there exists a small filtered  category $J$ and an exact sequence
of functors from $J$ to $\cA$
\begin{equation}\label{p1-abisoind4b}
0\rightarrow \varphi'\rightarrow \varphi\rightarrow \varphi''\rightarrow 0
\end{equation}
which induces the sequence \eqref{p1-abisoind4a} by taking the direct limit in $\bInd(\cA)$ (\cite{ag2} 2.6.7). By assumption, the sequence
\begin{equation}\label{p1-abisoind4c}
0\rightarrow \underset{\underset{J}{\longrightarrow}}{\lim}\ \varphi'\rightarrow
\underset{\underset{J}{\longrightarrow}}{\lim}\ \varphi\rightarrow
\underset{\underset{J}{\longrightarrow}}{\lim}\ \varphi''\rightarrow 0
\end{equation}
obtained by taking the direct limit in $\cA$ is exact. Hence, $\kappa_\cA$ is exact.

\begin{prop}\label{p1-abisoind40}
Let $V$ be a complete valuation ring of height $1$ and mixed characteristic $(0,p)$, equipped with the $p$-adic topology, $\cS=\Spf(V)$. 
Let $\fX$ be an $\cS$-formal scheme of finite presentation, $\cA$ an $\co_\fX$-algebra topologically of finite presentation \eqref{p1-pfs1}. 
For $A=\cA$ or $\cA[\frac 1 p]$, we denote by $\bMod^\coh(A)$ the category of coherent $A$-modules of $\fX_\zar$ and by 
$\bMod^{\coh}_\mQ(\cA)$ the category of coherent $\cA$-modules up to isogeny. Then, the canonical functor 
\begin{equation}
\bMod^{\coh}(\cA)\rightarrow \bMod^{\coh}(\cA[\frac 1p]),\ \ \ \cF\mapsto \cF[\frac 1 p],
\end{equation}
induces an equivalence of abelian categories
\begin{equation}
\bMod^{\coh}_\mQ(\cA)\stackrel{\sim}{\rightarrow} \bMod^{\coh}(\cA[\frac 1 p]). 
\end{equation} 
\end{prop}

The proof is identical to that of (\cite{agt} III.6.16) using \ref{p1-pfs15} and \ref{p1-pfs12}(i)-(ii).

\section{Ind-modules and ind-algebras}\label{p1-indma}

\subsection{}\label{p1-indmal1} 
In this section, we fix a ringed $\mU$-topos $(X,A)$ \eqref{p1-NC7}. 
We denote by $\bMod(A)$ the category of $A$-modules, by $\bIndMod(A)$ the category of ind-$A$-modules (\cite{ag2} 2.7.1) and by 
\begin{equation}\label{p1-indmal1a}
\iota_A\colon \bMod(A)\rightarrow \bIndMod(A)
\end{equation}
the canonical functor (\cite{ag2} 2.6.3) which is fully faithful. 
When there is no risk of ambiguity, we  identify $\bMod(A)$ with a full subcategory of $\bIndMod(A)$
by the functor $\iota_A$ that we omit from the notation.

Small inverse limits are representable in $\bIndMod(A)$ and $\iota_A$ commutes with these limits (\cite{ks2} 6.1.17(ii)). 
Small {\em filtered} direct limits are representable in $\bIndMod(A)$ (\cite{ks2} 6.1.8), however, $\iota_A$ does not commute in general with direct limits.
We will use the notation $\indcolim$ to denote direct limits in $\bIndMod(A)$ and we will keep the notation
$\underset{\longrightarrow}{\lim}$ to denote direct limits in $\bMod(A)$.

The category $\bIndMod(A)$ is abelian and the functor $\iota_A$ is exact (\cite{ks2} 8.6.5). 
The functor $\iota_A$ admits a left adjoint (\cite{ks2} 6.3.1)
\begin{equation}\label{p1-indmal1b}
\kappa_A \colon \bIndMod(A)\rightarrow \bMod(A). 
\end{equation}
For every small filtered  category $J$ and every functor $\alpha\colon J\rightarrow \bMod(A)$, we have an isomorphism
\begin{equation}\label{p1-indmal1c}
\kappa_A(\underset{\underset{J}{\longrightarrow}}{\mlq\mlq\lim \mrq\mrq} \alpha)\stackrel{\sim}{\rightarrow} \underset{ \underset{J}{\longrightarrow}}{\lim}\ \alpha.
\end{equation} 
The canonical morphism $\kappa_A\circ \iota_A \rightarrow \id_{\bMod(A)}$ is an isomorphism. The functor $\kappa_A$ is exact (\cite{ag2} 2.7.2). 

There is a canonical tensor product (\cite{ag2} (2.7.1.4))
\begin{equation}\label{p1-indmal1d}
\otimes_A\colon \bIndMod(A)\times \bIndMod(A)\rightarrow \bIndMod(A),
\end{equation}
that extends the usual tensor product of $\bMod(A)$. It makes $\bIndMod(A)$ into a symmetric monoidal category, with $A$ as the unit object. 

\begin{defi}\label{p1-indmal2} 
An {\em ind-$A$-algebra} is an algebra of the monoidal category $(\bIndMod(A),\otimes_A)$. 
If $B$ is an ind-$A$-algebra, a {\em left ind-$B$-module} is a left module over $B$ in this monoidal category. 
\end{defi}

Ind-$A$-algebras form naturally a category that we denote by $\bIndAlg(A)$. We denote by $\bAlg(A)$ the category of $A$-algebras. The functor $\iota_A$ induces a 
fully faithful functor that we denote also by
\begin{equation}\label{p1-indmal2b} 
\iota_A\colon \bAlg(A)\rightarrow \bIndAlg(A). 
\end{equation}
When there is no risk of ambiguity, we identify $\bAlg(A)$ with a full subcategory of $\bIndAlg(A)$
by the functor $\iota_A$ that we omit from the notation.

For every ind-$A$-algebra $B$, left ind-$B$-modules form naturally a category that we denote by $\bIndMod(B)$. It is an abelian category, and the forgetful functor 
$\bIndMod(B)\rightarrow \bIndMod(A)$ is exact. 

\subsection{}\label{p1-indmal3}
Let $B$ be a commutative $A$-algebra, $E$ an ind-$B$-module. By (\cite{ag2} 2.6.3.2 and 2.6.3.4), with the notation of (\cite{ag2} 2.1.8), 
the category $\bMod(B)_{/E}$ is filtered and essentially small, and we have 
\begin{equation}\label{p1-indmal3a}
E=\underset{\underset{(M\rightarrow E)\in\bMod(B)_{/E}}{\longrightarrow}}{\mlq\mlq\lim \mrq\mrq} M. 
\end{equation}
The canonical forgetful functor $\bMod(B)\rightarrow \bMod(A)$ induces an additive functor
\begin{equation}\label{p1-indmal3b}
\varphi_{B/A}\colon\bIndMod(B)\rightarrow \bIndMod(A).
\end{equation}
We deduce from this a functor $\jmath_E\colon \bMod(B)_{/E}\rightarrow \bMod(A)_{/\varphi_{B/A}(E)}$, which is cofinal by (\cite{sga4} I 8.1.3(b)). 
The {\em $B$-multiplication morphism of $E$}, denoted by $\mu_E\colon B\otimes_AE\rightarrow E$,  
is the morphism of ind-$A$-modules, direct limit of the $B$-multiplication morphisms $\mu_M\colon B\otimes_AM\rightarrow M$, for $M\in \ob(\bMod(B)_{/E})$.  
We thus define a functor 
\begin{equation}\label{p1-indmal3d}
\phi_{B/A}\colon 
\begin{array}[t]{clcr}
\bIndMod(B)&\rightarrow &\bIndMod(\iota_A(B))\\
E&\mapsto& (\varphi_{B/A}(E),\mu_E).
\end{array}
\end{equation}
It is an equivalence of categories by \ref{p1-indmal4} below. When there is no risk of ambiguity, we will identify the categories $\bIndMod(B)$ and $\bIndMod(\iota_A(B))$ 
by this functor that we omit from the notation.

\begin{prop}\label{p1-indmal4}
For every commutative $A$-algebra $B$, the functor $\phi_{B/A}$ \eqref{p1-indmal3d} is an equivalence of categories. 
\end{prop}

We prove first that the functor $\phi_{B/A}$ is fully faithful. 
It is enough to prove that for every $B$-module $M$ and every ind-$B$-module $E$, the morphism 
\begin{equation}\label{p1-indmal4a}
\Hom_{\bIndMod(B)}(M,E)\rightarrow \Hom_{\bIndMod(\iota_A(B))}(\phi_{B/A}(M),\phi_{B/A}(E))
\end{equation}
induced by $\phi_{B/A}$,  is an isomorphism. It is clearly injective. 
Consider an element of the right hand side of \eqref{p1-indmal4a}, namely, a morphism of ind-$A$-modules $u\colon M\rightarrow E$ such that the diagram 
\begin{equation}
\xymatrix{
{B\otimes_AM}\ar[r]^-(0.5){\id\otimes u}\ar[d]_{\mu_M}&{B\otimes_AE}\ar[d]^{\mu_E}\\
M\ar[r]^-(0.5)u&E,}
\end{equation}
where the vertical arrows are the $B$-multiplication morphisms, is commutative. 
Since the canonical (forgetful) functor $\jmath_E\colon \bMod(B)_{/E}\rightarrow \bMod(A)_{/\varphi_{B/A}(E)}$ is cofinal \eqref{p1-indmal3b}, the morphism $u$ factors into 
\begin{equation}
M\stackrel{v}{\rightarrow}N  \stackrel{w}{\rightarrow} E,
\end{equation}
where $N$ is a $B$-module, $v$ is $A$-linear and $w$ is a morphism of ind-$B$-modules. Moreover, by (\cite{sga4} I 8.1.3(b)), we can choose $N$ such that the diagram 
\begin{equation}
\xymatrix{
{B\otimes_AM}\ar[r]^-(0.5){\id\otimes v}\ar[d]_{\mu_M}&{B\otimes_AN}\ar[d]^{\mu_N}\\
M\ar[r]^-(0.5)v&N}
\end{equation}
is commutative, i.e., such that $v$ is $B$-linear. We deduce that \eqref{p1-indmal4a} is surjective and hence an isomorphism. 

We prove next that the functor $\phi_{B/A}$ is essentially surjective. Let $E$ be an ind-$A$-module, 
$\mu_E\colon B\otimes_AE\rightarrow E$ a morphism ind-$A$-modules 
such that $(E,\mu_E)$ is an ind-$\iota_A(B)$-module. We denote by $I$ the category made of pairs $(M,u)$, 
where $M$ is a $B$-module and $u\colon M\rightarrow E$ is a morphism of ind-$A$-modules such that the diagram 
\begin{equation}
\xymatrix{
{B\otimes_AM}\ar[r]^-(0.5){\id\otimes u}\ar[d]_{\mu_M}&{B\otimes_AE}\ar[d]^{\mu_E}\\
M\ar[r]^-(0.5)u&E,}
\end{equation}
where $\mu_M$ is the $B$-multiplication morphism of $M$, is commutative. 
The category $I$ is filtered and essentially small, since  
the sum of two objects of $I$ and the coequalizer of two morphisms of $I$ (with the same source and target) are representable in $I$. 
It is enough to show that the canonical forgetful functor $\rho\colon I\rightarrow \bMod(A)_{/E}$ is cofinal.
Let $P$ be an $A$-module, $u\colon P\rightarrow E$ a morphism of ind-$A$-modules. Then, $u$ factors into 
\begin{equation}
P\stackrel{v}{\rightarrow}B\otimes_AP  \stackrel{w}{\rightarrow} E,
\end{equation}
where $v$ is the canonical morphism and $w=\mu_E\circ (\id_B\otimes_A u)$. The pair $(B\otimes_AP,w)$ is clearly an object of $I$. 
We deduce, by (\cite{sga4} I 8.1.3(b)), that the functor $\rho$ is cofinal, which finishes the proof of the proposition.

\subsection{}\label{p1-indmal5}
Let $B$ be a commutative ind-$A$-algebra, $E$ and $F$ two ind-$B$-modules, $\mu_E\colon B\otimes_AE\rightarrow E$ and $\mu_F\colon B\otimes_AF\rightarrow F$ 
the $B$-multiplications of $E$ and $F$, respectively. 
We denote by $E\otimes_BF$ the cokernel of the morphism of ind-$A$-modules 
\begin{equation}\label{p1-indmal5a}
\mu_E\otimes_A\id_F-\id_E\otimes_A\mu_F\colon B\otimes_AE\otimes_AF\rightarrow E\otimes_AF,
\end{equation}
where we identified $B\otimes_AE\otimes_AF$ and $E\otimes_AB\otimes_AF$. 
We make $E\otimes_BF$ into an ind-$B$-module 
by the morphism $\mu_{E\otimes_BF}\colon B\otimes_A(E\otimes_BF)\rightarrow E\otimes_BF$ 
induced by $\mu_E\otimes_A\id_F$ (or equivalently by $\id_E\otimes_A\mu_F$). 
We thus define a bifunctor
\begin{equation}\label{p1-indmal5b}
\begin{array}[t]{clcr}
\bIndMod(B)\times \bIndMod(B)&\rightarrow& \bIndMod(B)\\
(E,F)&\mapsto& E\otimes_BF.
\end{array}
\end{equation}

Let $B\rightarrow C$ be a homomorphism of commutative ind-$A$-algebras. 
If $E$ is an ind-$B$-module and $F$ is an ind-$C$-module, then $E\otimes_BF$ is naturally an ind-$C$-module. 
In particular, the tensor product \eqref{p1-indmal5b} induces a functor
\begin{equation}\label{p1-indmal5c}
\begin{array}[t]{clcr}
\bIndMod(B)&\rightarrow& \bIndMod(C)\\
E&\mapsto& E\otimes_BC.
\end{array}
\end{equation}
It is a left adjoint of the forgetful functor 
\begin{equation}\label{p1-indmal5d}
\bIndMod(C) \rightarrow \bIndMod(B).
\end{equation}

Applied to the identity of $B$, this adjunction shows that for all ind-$B$-modules $E$ and $F$, we have a canonical isomorphism 
\begin{equation}\label{p1-indmal5e}
\Hom_{\bIndMod(B)}(E,F)\stackrel{\sim}{\rightarrow} \Hom_{\bIndMod(B)}(E\otimes_BB,F).
\end{equation}
Therefore, we have a canonical functorial isomorphism of ind-$B$-modules
\begin{equation}\label{p1-indmal5f}
E\otimes_BB\stackrel{\sim}{\rightarrow} E.
\end{equation}

For all ind-$B$-modules $E$ and $E'$, and all ind-$C$-modules $F$ and $F'$, we have canonical isomorphisms of ind-$C$-modules
\begin{eqnarray}
E\otimes_B (F\otimes_CF')\stackrel{\sim}{\rightarrow} (E\otimes_B F)\otimes_CF',\label{p1-indmal5g}\\
(E\otimes_B E')\otimes_BF\stackrel{\sim}{\rightarrow} E\otimes_B (E'\otimes_BF). \label{p1-indmal5h}
\end{eqnarray}

\begin{defi}\label{p1-indmal6}
Let $E$ be an ind-$A$-module, $n$ an integer $\geq 0$. 

The {\em $n$th symmetric power} of $E$, denoted by $\rS^n_A(E)$, is the direct limit of $\iota_A(\rS^n_A(M))$ \eqref{p1-indmal1a}, 
for $M\in \ob(\bMod(A)_{/E})$ (\cite{ag2} 2.6.3.4), where $\rS^n_A(M)$ is the $n$th symmetric power of $M$ \eqref{p1-NC7}. 
We have $\rS^0_A(E)=A$ and $\rS^1(E)=E$. 
For all integers $n,m\geq 0$, there is a canonical morphism $\mu^{m,n}\colon \rS^m_A(E)\otimes_A\rS^n_A(E)\rightarrow \rS^{n+m}_A(E)$. 

The {\em symmetric ind-$A$-algebra} of $E$, denoted by $\IndSym_A(E)$, is the ind-$A$-algebra defined by 
\begin{equation}\label{p1-indmal6a}
\IndSym_A(E)=\underset{n\geq 0}{\Indoplus}\ \rS^n_A(E),
\end{equation}
where $\Indoplus$ denotes the direct sum in $\bIndMod(A)$ (\cite{ag2} 2.6.7), equipped with the multiplication induced by the morphisms $\mu^{m,n}$. 
\end{defi}

For every $A$-module $M$, we have a canonical morphism of ind-$A$-algebras 
\begin{equation}\label{p1-indmal6b}
\IndSym_A(M)\rightarrow \iota_A(\rS_A(M)),
\end{equation}
where $\iota_A$ is the functor \eqref{p1-indmal2b}, which is not in general an isomorphism. However, it induces an isomorphism 
of $A$-algebras $\kappa_A(\IndSym_A(M))\stackrel{\sim}{\rightarrow} \rS_A(M)$.

\begin{lem}\label{p1-indmal7}
Let $E$ be an ind-$A$-module, $B$ an ind-$A$-algebra, $\mu_B\colon B\otimes_AB\rightarrow B$ the multiplication of $B$, 
$\varphi\colon E\rightarrow B$ a morphism of ind-$A$-modules such that the diagram 
\begin{equation}
\xymatrix{
{E\otimes_AE}\ar[d]_{\sigma}\ar[rr]^-(0.5){\mu_B\circ (\varphi\otimes \varphi)}&&B\\
{E\otimes_AE}\ar[rru]_{\mu_B\circ (\varphi\otimes \varphi)}&&}
\end{equation}
where $\sigma$ is the morphism that switches the factors, is commutative. Then, there exists a unique homomorphism of ind-$A$-algebras 
$\phi\colon \IndSym_A(E)\rightarrow B$ extending $\varphi$. 
\end{lem}

Indeed, by considering the essentially small filtered category $\bMod(A)_{/E}$, we are reduced to the case where $E$ is an $A$-module. The assertion is then immediate.

\subsection{}\label{p1-indmal8}
Let $F, G$ be two ind-$A$-modules. We consider the presheaf $\Phi(F,G)$ of sets on $\bMod(A)$ defined for any $A$-module $M$, by 
\begin{equation}\label{p1-indmal8b}
\Phi(F,G)(M)=\Hom_{\bIndMod(A)}(M\otimes_AF,G).
\end{equation}
We first prove that it is representable by an ind-$A$-module. 
Let $I$ and $J$ two small filtered categories, $\alpha\colon I\rightarrow \bMod(A)$ and $\beta\colon J\rightarrow \bMod(A)$ two functors such that 
\begin{equation}\label{p1-indmal8a}
F=\underset{\underset{I}{\longrightarrow}}{\mlq\mlq\lim \mrq\mrq} \alpha \ \ \ {\rm and}\ \ \ 
G=\underset{\underset{J}{\longrightarrow}}{\mlq\mlq\lim \mrq\mrq} \beta.
\end{equation} 
By (\cite{ag2} 2.7.3), we have a canonical isomorphism of ind-$A$-modules 
\begin{equation}\label{p1-indmal8c}
M\otimes_AF\stackrel{\sim}{\rightarrow} \underset{\underset{i\in I}{\longrightarrow}}{\mlq\mlq\lim \mrq\mrq} M\otimes_A\alpha(i).
\end{equation}
We deduce canonical isomorphisms 
\begin{eqnarray}\label{p1-indmal8d}
\Phi(F,G)(M)&\stackrel{\sim}{\rightarrow}&\underset{\underset{i\in I}{\longleftarrow}}{\lim}\ \underset{\underset{j\in J}{\longrightarrow}}{\lim} \ \Hom_{\bMod(A)}(M\otimes_A\alpha(i),\beta(j))\\
&\stackrel{\sim}{\rightarrow}&\underset{\underset{i\in I}{\longleftarrow}}{\lim}\ \underset{\underset{j\in J}{\longrightarrow}}{\lim} \ \Hom_{\bMod(A)}(M,\cHom_A(\alpha(i),\beta(j)))\nonumber\\
&\stackrel{\sim}{\rightarrow}&\underset{\underset{i\in I}{\longleftarrow}}{\lim} \ \Hom_{\bIndMod(A)}(M, \underset{\underset{j\in J}{\longrightarrow}}{\mlq\mlq\lim \mrq\mrq} 
\cHom_A(\alpha(i),\beta(j)))\nonumber\\
&\stackrel{\sim}{\rightarrow}&\Hom_{\bIndMod(A)}(M, \underset{\underset{i\in I}{\longleftarrow}}{\lim} \ \underset{\underset{j\in J}{\longrightarrow}}{\mlq\mlq\lim \mrq\mrq}  
\cHom_A(\alpha(i),\beta(j))),\nonumber
\end{eqnarray}
where $\cHom_A(\alpha(i),\beta(j))$ is the internal Hom of $\bMod(A)$; the first, third and fourth isomorphisms follow from 
(\cite{sga4} I 3.1, \cite{ks2} 6.1.8 and 6.1.17(ii)), 
and the second one is induced by the Cartan isomorphism. 
Therefore, $\Phi(F,G)$ is representable by an ind-$A$-module that we denote by $\cHom_A(F,G)$. 
We thus define a bifunctor
\begin{equation}\label{p1-indmal8e}
\cHom_A\colon \bIndMod(A)^\circ\times \bIndMod(A)\rightarrow \bIndMod(A),
\end{equation}
which clearly extends the internal Hom of $\bMod(A)$. The notation hence does not lead to any ambiguity. 

\begin{prop}\label{p1-indmal9}
Let $I,J$ be two small categories such that $I$ is filtered, 
$\alpha\colon I\rightarrow \bIndMod(A)$, $\beta\colon J^\circ \rightarrow \bIndMod(A)$ two functors. 
\begin{itemize}
\item[{\rm (i)}] For every $A$-module $E$, we have a canonical isomorphism
\begin{equation}\label{p1-indmal9a}
\cHom_A(E,\underset{\underset{I}{\longrightarrow}}{\mlq\mlq\lim \mrq\mrq} \alpha) \stackrel{\sim}{\rightarrow}
\underset{\underset{i\in I}{\longrightarrow}}{\mlq\mlq\lim \mrq\mrq} \cHom_A(E,\alpha(i)).
\end{equation}
\item[{\rm (ii)}] For every ind-$A$-module $F$, we have a canonical isomorphism
\begin{equation}\label{p1-indmal9b}
\cHom_A(\underset{\underset{I}{\longrightarrow}}{\mlq\mlq\lim \mrq\mrq} \alpha,F) \stackrel{\sim}{\rightarrow}
\underset{\underset{i\in I}{\longleftarrow}}{\lim} \ \cHom_A(\alpha(i),F).
\end{equation}
\item[{\rm (iii)}] For every ind-$A$-module $F$, we have a canonical isomorphism
\begin{equation}\label{p1-indmal9c}
\cHom_A(F,\underset{\underset{J}{\longleftarrow}}{\lim} \ \beta) \stackrel{\sim}{\rightarrow}
\underset{\underset{j\in J}{\longleftarrow}}{\lim} \ \cHom_A(F,\beta(j)).
\end{equation}
\end{itemize}
\end{prop}

(i) It follows from the composed isomorphism \eqref{p1-indmal8d}. 

(ii) Indeed, for every $A$-module $M$, we have canonical isomorphisms 
\begin{eqnarray*}
\Hom_{\bIndMod(A)}(M,\cHom_A(\underset{\underset{I}{\longrightarrow}}{\mlq\mlq\lim \mrq\mrq} \alpha,F))
&\stackrel{\sim}{\rightarrow}&
\Hom_{\bIndMod(A)}(\underset{\underset{i\in I}{\longrightarrow}}{\mlq\mlq\lim \mrq\mrq} M\otimes_A \alpha(i),F)\\
&\stackrel{\sim}{\rightarrow}&
\underset{\underset{i\in I}{\longleftarrow}}{\lim}\ \Hom_{\bIndMod(A)}(M\otimes_A \alpha(i),F)\\
&\stackrel{\sim}{\rightarrow}&
\underset{\underset{i\in I}{\longleftarrow}}{\lim}\ \Hom_{\bIndMod(A)}(M,\cHom_A(\alpha(i),F))\\
&\stackrel{\sim}{\rightarrow}&
\Hom_{\bIndMod(A)}(M,\underset{\underset{i\in I}{\longleftarrow}}{\lim}\ \cHom_A(\alpha(i),F)). 
\end{eqnarray*}

(iii) Indeed, for every $A$-module $M$, we have canonical isomorphisms 
\begin{eqnarray*}
\Hom_{\bIndMod(A)}(M,\cHom_A(F,\underset{\underset{J}{\longleftarrow}}{\lim} \ \beta))
&\stackrel{\sim}{\rightarrow}&
\Hom_{\bIndMod(A)}(M\otimes_AF, \underset{\underset{J}{\longleftarrow}}{\lim} \ \beta)\\
&\stackrel{\sim}{\rightarrow}&
\underset{\underset{j\in J}{\longleftarrow}}{\lim} \ \Hom_{\bIndMod(A)}(M\otimes_AF, \beta(j))\\
&\stackrel{\sim}{\rightarrow}&
\underset{\underset{j\in J}{\longleftarrow}}{\lim} \ \Hom_{\bIndMod(A)}(M,\cHom_A(F, \beta(j)))\\
&\stackrel{\sim}{\rightarrow}&
\Hom_{\bIndMod(A)}(M,\underset{\underset{j\in J}{\longleftarrow}}{\lim} \ \cHom_A(F, \beta(j))).
\end{eqnarray*}

\begin{cor}\label{p1-indmal10}
The bifunctor \eqref{p1-indmal8e} is left exact.
\end{cor}

This follows from \ref{p1-indmal9}(i)-(ii) and (\cite{ks2} 8.6.5(iii) and 8.6.6(a)).

\begin{prop}\label{p1-indmal11}
Let $E,F,G,H$ be four ind-$A$-modules. Then, 
\begin{itemize}
\item[{\rm (i)}] We have a canonical functorial isomorphism of ind-$A$-modules
\begin{equation}\label{p1-indmal11h}
E\stackrel{\sim}{\rightarrow}\cHom_A(A,E).
\end{equation}
\item[{\rm (ii)}] We have a canonical functorial morphism of ind-$A$-modules
\begin{equation}\label{p1-indmal11i}
\cHom_A(E,F)\rightarrow\cHom_A(E\otimes_AG,F\otimes_AG),
\end{equation}
induced by the tensor product with $G$. 
\item[{\rm (iii)}] The composition of the morphisms 
\begin{equation}
\cHom_A(E,F)\rightarrow\cHom_A(E\otimes_AG,F\otimes_AG)\rightarrow\cHom_A(E\otimes_AG\otimes_AH,F\otimes_AG\otimes_AH),
\end{equation}
induced by the tensor products with $G$ and $H$, respectively, is the morphism induced by the tensor product with $G\otimes_AH$. 
\item[{\rm (iv)}] We have a canonical functorial isomorphism of abelian groups
\begin{equation}\label{p1-indmal11a}
\Hom_{\bIndMod(A)}(E\otimes_AF,G)\stackrel{\sim}{\rightarrow} \Hom_{\bIndMod(A)}(E,\cHom_A(F,G)).
\end{equation}
In other terms, the symmetric monoidal category $\bIndMod(A)$ is closed 
with internal Hom $\cHom_A$. 
\item[{\rm (v)}] The isomorphism \eqref{p1-indmal11a}
\begin{equation}\label{p1-indmal11m}
\Hom_{\bIndMod(A)}(E\otimes_AF,E\otimes_AF)\stackrel{\sim}{\rightarrow} \Hom_{\bIndMod(A)}(E,\cHom_A(F,E\otimes_AF))
\end{equation}
maps the identity of $E\otimes_AF$ to the morphism $\mu\colon E\rightarrow \cHom_A(F,E\otimes_AF)$ induced by the tensor product with $F$ \eqref{p1-indmal11i}. 
\item[{\rm (vi)}] We have a canonical functorial isomorphism of ind-$A$-modules
\begin{equation}\label{p1-indmal11b}
\cHom_A(E\otimes_AF,G)\stackrel{\sim}{\rightarrow} \cHom_A(E,\cHom_A(F,G)).
\end{equation}
\item[{\rm (vii)}] The diagram 
\begin{equation}\label{p1-indmal11c}
\xymatrix{
{\cHom_A(E\otimes_AF\otimes_AG,H)}\ar[r]\ar[d]&{\cHom_A(E\otimes_AF,\cHom_A(G,H))}\ar[d]\\
{\cHom_A(E,\cHom_A(F\otimes_AG,H))}\ar[r]&{\cHom_A(E,\cHom_A(F,\cHom_A(G,H))),}}
\end{equation}
where all morphisms are induced by \eqref{p1-indmal11b}, is commutative. 
\item[{\rm (viii)}] The diagram 
\begin{equation}\label{p1-indmal11k}
\xymatrix{
{\cHom_A(E\otimes_AG,F)}\ar[r]\ar[d]&{\cHom_A(E\otimes_AG\otimes_AH,F\otimes_AH)}\ar[d]\\
{\cHom_A(G,\cHom_A(E,F))}\ar[r]&{\cHom_A(G,\cHom_A(E\otimes_AH,F\otimes_AH)),}}
\end{equation}
where the vertical morphisms are \eqref{p1-indmal11b} and the horizontal morphisms are induced by the tensor product with $H$, is commutative. 
\end{itemize}
\end{prop}

(i) The required morphism corresponds for any $A$-module $M$, to the identity of 
\[
\Hom_{\bIndMod(A)}(M,E).
\] 

(ii) The required morphism corresponds for any $A$-module $M$, to the morphism 
\begin{equation}\label{p1-indmal11j}
\Hom_{\bIndMod(A)}(M\otimes_AE,F)\rightarrow\Hom_{\bIndMod(A)}(M\otimes_AE\otimes_AG,F\otimes_AG),
\end{equation}
induced by the tensor product with $G$. 

(iii) Clear from the definition.

(iv) Let $I$ be a small filtered category, $\alpha\colon I\rightarrow \bMod(A)$ a functor such that 
\begin{equation}\label{p1-indmal11d}
E=\underset{\underset{I}{\longrightarrow}}{\mlq\mlq\lim \mrq\mrq} \alpha.
\end{equation}
For every $i\in \ob(I)$, we have a canonical functorial isomorphism 
\begin{equation}\label{p1-indmal11e}
\Hom_{\bIndMod(A)}(\alpha(i)\otimes_AF,G) \stackrel{\sim}{\rightarrow}
\Hom_{\bIndMod(A)}(\alpha(i),\cHom_A(F,G)).
\end{equation}
The inverse limit of these isomorphisms gives the required isomorphism \eqref{p1-indmal11a}, which does not depend on $\alpha$ by (\cite{ag2} 2.6.3.4). 

(v) Indeed, by definition, for every $A$-module $M$, the diagram 
\begin{equation}
\xymatrix{
{\Hom_{\bIndMod(A)}(M,E)}\ar[r]\ar[rd]&{\Hom_{\bIndMod(A)}(M\otimes_AF,E\otimes_AF)}\ar[d]\\
&{\Hom_{\bIndMod(A)}(M,\cHom_A(F,E\otimes_AF)),}}
\end{equation}
where the horizontal morphism is induced by the tensor product with $F$, 
the vertical one is the canonical isomorphism and the slanted one in induced by $\mu$, is commutative. 
By taking an inverse limit as in the proof of (iv), we deduce that the 
\begin{equation}
\xymatrix{
{\Hom_{\bIndMod(A)}(E,E)}\ar[r]\ar[rd]&{\Hom_{\bIndMod(A)}(E\otimes_AF,E\otimes_AF)}\ar[d]\\
&{\Hom_{\bIndMod(A)}(E,\cHom_A(F,E\otimes_AF)),}}
\end{equation}
where the horizontal morphism is induced by the tensor product with $F$, 
the vertical one is \eqref{p1-indmal11a} and the slanted one in induced by $\mu$, is commutative. 
We deduce the proposition by following the image of $\id_E$.

(vi) The required morphism is defined, for every $A$-module $M$, by the lower horizontal morphism of the commutative diagram 
\begin{equation}\label{p1-indmal11f}
\xymatrix{
{\Hom_{\star}(M\otimes_AE\otimes_AF,G)}\ar[r]\ar[d]&{\Hom_{\star}(M\otimes_AE,\cHom_A(F,G))}\ar[d]\\
{\Hom_{\star}(M,\cHom_A(E\otimes_AF,G))}\ar[r]&{\Hom_{\star}(M,\cHom_A(E,\cHom_A(F,G))),}}
\end{equation}
where $\star=\bIndMod(A)$, the vertical and the upper horizontal arrows are the isomorphisms \eqref{p1-indmal11a}, 
taking into account (\cite{ag2} (2.7.5.1)). 

(vii) By taking an inverse limit of \eqref{p1-indmal11f} as in the proof of (iv), we see that the diagram
\begin{equation}\label{p1-indmal11g}
\xymatrix{
{\Hom_{\star}(E\otimes_AF\otimes_AG,H)}\ar[r]\ar[d]&{\Hom_{\star}(E\otimes_AF,\cHom_A(G,H))}\ar[d]\\
{\Hom_{\star}(E,\cHom_A(F\otimes_AG,H))}\ar[r]&{\Hom_{\star}(E,\cHom_A(F,\cHom_A(G,H))),}}
\end{equation}
where $\star=\bIndMod(A)$ and the morphisms are induced either by \eqref{p1-indmal11a} or \eqref{p1-indmal11b}, is commutative. 
The proposition follows by replacing $E$ by $M\otimes_AE$ for a variable $A$-module $M$. 

(viii) It follows from the definition of \eqref{p1-indmal11i}, by taking an inverse limit, that the diagram
\begin{equation}\label{p1-indmal11l}
\xymatrix{
{\Hom_{\star}(E\otimes_AG,F)}\ar[r]\ar[d]&{\Hom_{\star}(E\otimes_AG\otimes_AH,F\otimes_AH)}\ar[d]\\
{\Hom_{\star}(G,\cHom_A(E,F))}\ar[r]&{\Hom_{\star}(G,\cHom_A(E\otimes_AH,F\otimes_AH)),}}
\end{equation}
where $\star=\bIndMod(A)$, the vertical morphisms are \eqref{p1-indmal11a}, 
the horizontal morphisms are induced by the tensor product with $H$, is commutative. 
The proposition follows by replacing $G$ by $M\otimes_AG$ for a variable $A$-module $M$.

\begin{lem}\label{p1-indmal12}
Let $E,F,G,H$ be four ind-$A$-modules. Then, 
\begin{itemize}
\item[{\rm (i)}] The isomorphism \eqref{p1-indmal11b} induces by \ref{p1-indmal11}{\rm (iv)} a canonical functorial morphism of ind-$A$-modules
\begin{equation}\label{p1-indmal12a}
E\otimes_A\cHom_A(E\otimes_AF,G)\rightarrow \cHom_A(F,G).
\end{equation}
In particular, we have  a canonical functorial morphism of ind-$A$-modules
\begin{equation}\label{p1-indmal12b}
E\otimes_A\cHom_A(E,F)\rightarrow F.
\end{equation} 
\item[{\rm (ii)}] The composition of the morphisms induced by \eqref{p1-indmal12a}
\begin{equation}
E\otimes_AF\otimes_A\cHom_A(E\otimes_AF\otimes_AG,H)\rightarrow E\otimes_A\cHom_A(E\otimes_AG,H)\rightarrow \cHom_A(G,H)
\end{equation}
is also of type \eqref{p1-indmal12a}. 
\item[{\rm (iii)}] The inverse image of a morphism of ind-$A$-modules $u\colon E\rightarrow \cHom_A(F,G)$ by the isomorphism \eqref{p1-indmal11a} is the composed morphism
\begin{equation}\label{p1-indmal12c}
\xymatrix{
{F\otimes_AE}\ar[r]^-(0.5){\id\otimes u}&{F\otimes_A\cHom_A(F,G)}\ar[r]&G},
\end{equation}
where the second one is \eqref{p1-indmal12b}.
\item[{\rm (iv)}] We have a canonical functorial morphism of ind-$A$-modules
\begin{equation}\label{p1-indmal12d}
G\otimes_A\cHom_A(E,F)\rightarrow \cHom_A(E,F\otimes_AG),
\end{equation}
composed of the morphisms 
\begin{equation}
G\otimes_A\cHom_A(E,F)\rightarrow G\otimes_A\cHom_A(E\otimes_AG,F\otimes_AG) \rightarrow \cHom_A(E,F\otimes_AG),
\end{equation}
where the first one is induced by the tensor product with $G$ and the second one is \eqref{p1-indmal12a}. 
\item[{\rm (v)}] The composition of the morphisms induced by \eqref{p1-indmal12d}
\begin{equation}\label{p1-indmal12e}
\cHom_A(E,F)\otimes_AG\otimes_AH \rightarrow \cHom_A(E,F\otimes_AG)\otimes_AH \rightarrow \cHom_A(E,F\otimes_AG\otimes_AH)
\end{equation}
is also of type \eqref{p1-indmal12d}. 
\item[{\rm (vi)}] The morphism \eqref{p1-indmal12d} corresponds by \ref{p1-indmal11}{\rm (iv)} to the morphism 
\begin{equation}
G\otimes_AE\otimes_A\cHom_A(E,F)\rightarrow F\otimes_AG,
\end{equation}
induced by \eqref{p1-indmal12b}. 
\end{itemize}
\end{lem}

(i) Obvious.  

(ii) It follows from \ref{p1-indmal11}(vii). 

(iii) It follows immediately from the functoriality of \eqref{p1-indmal11a} and the definition of \eqref{p1-indmal12b}. 

(iv) Obvious. 

(v) It follows from (ii) and \ref{p1-indmal11}(iii)-(viii).  

(vi) By (iii), we need to show that the diagram 
\begin{equation}\label{p1-indmal12f}
\xymatrix{
{G\otimes_AE\otimes_A\cHom_A(E,F)}\ar[r]\ar[rdd]&{G\otimes_AE\otimes_A\cHom_A(E\otimes_AG,F\otimes_AG)}\ar[d]\\
&{E\otimes_A\cHom_A(E,F\otimes_AG)}\ar[d]\\
&{F\otimes G,}}
\end{equation}
where the horizontal morphism is induced by the tensor product with $G$ and the others by \eqref{p1-indmal12a}, is commutative. 
By \ref{p1-indmal11}(viii), the diagram 
\begin{equation}\label{p1-indmal12g}
\xymatrix{
{E\otimes_A\cHom_A(E,F)}\ar[r]\ar[d]&{E\otimes_A\cHom_A(E\otimes_AG,F\otimes_AG)}\ar[d]\\
{\cHom_A(A,F)}\ar[r]&{\cHom_A(G,F\otimes_AG),}}
\end{equation}
where the vertical morphisms are induced by \eqref{p1-indmal12a} and the horizontal ones by the tensor product with $G$, is commutative. 
Then by (ii), we are further reduced to proving that the diagram 
\begin{equation}\label{p1-indmal12h}
\xymatrix{
{G\otimes_A\cHom_A(A,F)}\ar[r]\ar[rd]_{\id}&{G\otimes_A\cHom_A(G,F\otimes_AG)}\ar[d]\\
&{F\otimes_AG,}}
\end{equation}
where the vertical morphism is induced by \eqref{p1-indmal12b} and the horizontal one by the tensor product with $G$, is commutative. 
It is a consequence of (iii) and \ref{p1-indmal11}(v).

\begin{prop}\label{p1-indmal13}
For all ind-$A$-modules $E,F,G$, there exists a functorial canonical morphism of ind-$A$-modules
\begin{equation}\label{p1-indmal13a}
\cHom_A(E,F)\otimes_A\cHom_A(F,G)\rightarrow \cHom_A(E,G).
\end{equation}
Moreover, it satisfies the following properties:
\begin{itemize}
\item[{\rm (i)}] The diagram 
\begin{equation}\label{p1-indmal13c}
\xymatrix{
{E\otimes_A\cHom_A(E,F)\otimes_A\cHom_A(F,G)}\ar[r]\ar[d]&{E\otimes_A\cHom_A(E,G)}\ar[d]\\
{F\otimes_A\cHom_A(F,G)}\ar[r]&G,}
\end{equation}
where the upper horizontal morphism is induced by \eqref{p1-indmal13a} and the other morphisms are induced by \eqref{p1-indmal12b}, is commutative. 
\item[{\rm (ii)}] For every ind-$A$-module $H$, the diagram
\begin{equation}\label{p1-indmal130b}
\xymatrix{
{\cHom_A(E,F)\otimes_A\cHom_A(F,G)\otimes_A\cHom_A(G,H)}\ar[d]\ar[r]&{\cHom_A(E,F)\otimes_A\cHom_A(F,H)}\ar[d]\\
{\cHom_A(E,G)\otimes_A\cHom_A(G,H)}\ar[r]&{\cHom_A(E,H),}}
\end{equation}
where all morphisms are induced by \eqref{p1-indmal13a}, is commutative.
\end{itemize}
\end{prop}

Indeed, the required morphism \eqref{p1-indmal13a} is composed of 
\[
\cHom_A(E,F)\otimes_A\cHom_A(F,G)\rightarrow\cHom_A(E,F\otimes_A\cHom_A(F,G))\rightarrow\cHom_A(E,G),
\]
where the first arrow is the morphism \eqref{p1-indmal12d} and the second arrow is induced by the morphism \eqref{p1-indmal12b}.  

(i) It follows from \ref{p1-indmal12}(ii)-(iii)  and \ref{p1-indmal11}(v)-(viii); see \eqref{p1-indmal12g} and \eqref{p1-indmal12h}. 

(ii) It follows from (i) and \ref{p1-indmal12}(v).

\begin{cor}\label{p1-indmal14}
For every ind-$A$-module $E$, the canonical morphism \eqref{p1-indmal13a} 
\begin{equation}\label{p1-indmal14a}
\cEnd_A(E)\otimes_A\cEnd_A(E)\rightarrow \cEnd_A(E)
\end{equation}
makes the ind-$A$-module $\cEnd_A(E)$ into an ind-$A$-algebra \eqref{p1-indmal2}.
\end{cor}
It follows from \ref{p1-indmal13}. 

\begin{rema}\label{p1-indmal130}
Let $u\colon E\rightarrow F$, $v\colon F\rightarrow G$ be two morphisms of ind-$A$-modules that we consider as morphisms 
$u\colon A\rightarrow \cHom_A(E,F)$ and $v\colon A\rightarrow \cHom_A(F,G)$. Then, the composition 
\begin{equation}
\xymatrix{
A\ar[r]^-(0.5){u\otimes v}&{\cHom_A(E,F)\otimes_A \cHom_A(F,G)}\ar[r]&{\cHom_A(E,G),}}
\end{equation}
where the second morphism is \eqref{p1-indmal13a}, corresponds to $v\circ u\colon E\rightarrow G$. 
This follows easily from \ref{p1-indmal12}(iii) and \ref{p1-indmal13}(i). 
\end{rema}

\begin{prop}\label{p1-indmal15}
Let $E$ be an ind-$A$-module, $B$ an ind-$A$-algebra. Giving a left ind-$B$-module structure on $E$ 
is equivalent to giving a homomorphism of ind-$A$-algebras $B\rightarrow \cEnd_A(E)$. 
\end{prop}

Indeed, giving a morphism of ind-$A$-modules $\mu_E\colon B\otimes_AE\rightarrow E$ is equivalent to giving a morphism of ind-$A$-modules 
$\varphi\colon B\rightarrow \cEnd_A(E)$ by \ref{p1-indmal11}(iv). 
The fact that $\mu_E$ defines on $E$ an ind-$B$-module structure is equivalent to the commutativity of the following diagram
\begin{equation}
\xymatrix{
{B\otimes_AB\otimes_AE}\ar[r]^-(0.5){\mu_B\otimes \id}\ar[d]_{\id\otimes \varphi\otimes \id}\ar[rrd]|{\id\otimes \mu_E}&{B\otimes_AE}\ar[r]^-(0.5){\mu_E}&{E}\\
{B\otimes_A\cEnd_A(E)\otimes_AE}\ar[rr]&&{B\otimes_AE,}\ar[u]_{\mu_E}}
\end{equation}
where $\mu_B$ is the multiplication of $B$ and the unlabelled arrow is induced by \eqref{p1-indmal12b}. Observe that the left triangle is commutative by \ref{p1-indmal12}(iii).
It is further equivalent, by \ref{p1-indmal12}(vi), to the commutativity of the diagram
\begin{equation}
\xymatrix{
{B\otimes_AB}\ar[r]^-(0.5){\mu_B}\ar[d]_{\id\otimes \varphi}&{B}\ar[r]^\varphi&{\cEnd_A(E)}\\
{B\otimes_A\cEnd_A(E)}\ar[rr]&&{\cHom_A(E,B\otimes_AE),}\ar[u]_{\cHom_A(E,\mu_E)}}
\end{equation}
where the unlabelled morphism is \eqref{p1-indmal12d}. 
On the other hand, by \ref{p1-indmal12}(iii), the diagram 
\begin{equation}
\xymatrix{
{B\otimes_A\cEnd_A(E)}\ar[r]\ar[d]_{\varphi\otimes \id}&{\cHom_A(E,B\otimes_AE)}\ar[d]^{\cHom_A(E,\mu_E)}\\
{\cEnd_A(E)\otimes_A\cEnd(E)}\ar[r]&{\cEnd_A(E),}}
\end{equation}
where the upper horizontal morphism is \eqref{p1-indmal12d} and the lower one is \eqref{p1-indmal13a}, is commutative.   
The proposition follows.

\subsection{}\label{p1-indmal16}
Let $U$ an object of $X$.
We denote by $j_U\colon X_{/U}\rightarrow X$ the localization morphism of $X$ over $U$.
For any $F\in \ob(X)$, the sheaf $j_U^*(F)$ will also be denoted by $F|U$. The topos $X_{/U}$ will be ringed by $A|U$.
By (\cite{ag2} 2.7.10), $j_U$ induces two functors  
\begin{eqnarray}
\rI j_U^*\colon \bIndMod(A)&\rightarrow& \bIndMod(A|U),\label{p1-indmal16a}\\
\rI j_{U*}\colon \bIndMod(A|U)&\rightarrow& \bIndMod(A).\label{p1-indmal16b} 
\end{eqnarray}
The functor $\rI j_U^*$ is a left adjoint of the functor $\rI j_{U*}$.
By (\cite{ks2} 8.6.8), the functor $\rI j_U^*$ (resp.\ $\rI  j_{U*}$) is right exact (resp.\ left exact). 
For any ind-$A$-module $F$, the ind-$(A|U)$-module $\rI j_U^*(F)$ will also be denoted by $F|U$.

The extension by zero functor $j_{U!}\colon \bMod(A|U)\rightarrow \bMod(A)$ is exact and faithful (\cite{sga4} IV 11.3.1).
By (\cite{ag2} 2.6.4.4) and (\cite{ks2} 8.6.8), it induces an exact and faithful functor
\begin{equation}\label{p1-indmal16c}
\rI j_{U!}\colon \bIndMod(A|U)\rightarrow \bIndMod(A).
\end{equation}
This is a left adjoint of the functor $\rI j_U^*$. Therefore, the functor $\rI j_U^*$ commutes with representable direct and inverse limits and is in particular exact.

\begin{prop}\label{p1-indmal17}
Let $U$ be an object of $X$, $L$ an ind-$(A|U)$-module, $E,F,G$ three ind-$A$-modules. Then,
\begin{itemize}
\item[{\rm (i)}] We have a canonical functorial isomorphism 
\begin{eqnarray}
(E\otimes_AF)|U &\stackrel{\sim}{\rightarrow}& (E|U)\otimes_{A|U}(F|U),\label{p1-indmal17a}\\
\rI j_{U!}(L\otimes_{A|U}(E|U)) &\stackrel{\sim}{\rightarrow}& \rI j_{U!}(L)\otimes_AE.\label{p1-indmal17b}
\end{eqnarray}
\item[{\rm (ii)}] We have a canonical functorial isomorphism 
\begin{equation}\label{p1-indmal17c}
\cHom_A(E,F)|U\stackrel{\sim}{\rightarrow}\cHom_{A|U}(E|U,F|U).
\end{equation}
\item[{\rm (iii)}] The diagram 
\begin{equation}\label{p1-indmal17d}
\xymatrix{
{\Hom_{\star}(\rI j_{U!}(L)\otimes_AE,F)}\ar[r]\ar[d]_a&{\Hom_{\star}(\rI j_{U!}(L),\cHom_A(E,F))}\ar[d]^c\\
{\Hom_{\star}(\rI j_{U!}(L\otimes_{A|U}(E|U)),F)}\ar[d]_b&{\Hom_{\star|U}(L,\cHom_A(E,F)|U)}\ar[d]^d\\
{\Hom_{\star|U}(L\otimes_{A|U}(E|U),F|U)}\ar[r]&{\Hom_{\star|U}(L,\cHom_{A|U}(E|U,F|U)),}}
\end{equation}
where $\star=\bIndMod(A)$, $\star|U=\bIndMod(A|U)$, the horizontal morphisms are \eqref{p1-indmal11a}, 
$a$ is induced by \eqref{p1-indmal17b}, $b$ and $c$
are the adjunction isomorphisms and $d$ is induced by \eqref{p1-indmal17c}, is commutative. 
\item[{\rm (iv)}] The diagram 
\begin{equation}\label{p1-indmal17e}
\xymatrix{
{\Hom_{\star}(E\otimes_AF,G)}\ar[r]\ar[d]&{\Hom_{\star}(E,\cHom_A(F,G))}\ar[d]\\
{\Hom_{\star|U}((E|U)\otimes_A(F|U),G|U)}\ar[r]&{\Hom_{\star|U}(E|U,\cHom_{A|U}(F|U,G|U)),}}
\end{equation}
where $\star=\bIndMod(A)$, $\star|U=\bIndMod(A|U)$, the vertical morphisms are induced by $\rI j_U^*$, 
\eqref{p1-indmal17a} and \eqref{p1-indmal17c}, 
and the horizontal morphisms are \eqref{p1-indmal11a}, is commutative. 
\item[{\rm (v)}] The diagram 
\begin{equation}\label{p1-indmal17f}
\xymatrix{
{\cHom_A(E\otimes_AF,G)|U}\ar[r]\ar[d]&{\cHom_A(E,\cHom_A(F,G))|U}\ar[d]\\
{\cHom_{A|U}((E|U)\otimes_{A|U}(F|U),G|U)}\ar[r]&{\cHom_{A|U}(E|U,\cHom_{A|U}(F|U,G|U)),}}
\end{equation}
where the horizontal morphisms are induced by \eqref{p1-indmal11b} and the vertical ones by \eqref{p1-indmal17c}, is commutative. 
\item[{\rm (vi)}] The diagram 
\begin{equation}\label{p1-indmal17g}
\xymatrix{
{(E\otimes_A\cHom_A(E\otimes_AF,G))|U}\ar[r]\ar[d]&{\cHom_A(F,G)|U}\ar[d]\\
{(E|U)\otimes_{A|U}\cHom_{A|U}((E|U)\otimes_{A|U}(F|U),G|U)}\ar[r]&{\cHom_{A|U}(F|U,G|U),}}
\end{equation}
where the horizontal morphisms are induced by \eqref{p1-indmal12a} 
and the vertical ones by \eqref{p1-indmal17a} and \eqref{p1-indmal17c}, is commutative. 
\item[{\rm (vii)}] The diagram 
\begin{equation}\label{p1-indmal17h}
\xymatrix{{(G\otimes_A\cHom_A(E,F))|U}\ar[r]\ar[d]&{\cHom_A(E,F\otimes_AG)|U}\ar[d]\\
{(G|U)\otimes_{A|U}\cHom_{A|U}(E|U,F|U)}\ar[r]&{\cHom_{A|U}(E|U,(F|U)\otimes_{A|U}(G|U)),}}
\end{equation}
where the horizontal morphisms are induced by \eqref{p1-indmal12d} and the vertical ones by \eqref{p1-indmal17c}, is commutative. 
\item[{\rm (viii)}] The diagram 
\begin{equation}
\xymatrix{
{(\cHom_A(E,F)\otimes_A\cHom_A(F,G))|U}\ar[r]\ar[d]&{\cHom_A(E,G)|U}\ar[d]\\
{\cHom_{A|U}(E,F)\otimes_{A|U}\cHom_{A|U}(F|U,G|U)}\ar[r]&{\cHom_{A|U}(E|U,G|U),}}
\end{equation}
where the horizontal morphisms are induced by \eqref{p1-indmal13a} and the vertical ones by \eqref{p1-indmal17c}, is commutative. 
\end{itemize}
\end{prop}

(i) It follows from (\cite{sga4} 12.11), (\cite{ag2} 2.7.3) and the fact that $\rI j_U^*$ and $\rI j_{U!}$ commute with small filtered direct limits.  

(ii)  It follows from \ref{p1-indmal9}, (\cite{sga4} 12.3(b)) and the fact that $\rI j_U^*$ commutes with representable direct and inverse limits.  

(iii) Indeed, by taking direct and inverse limits, we may assume that $L$ is an $(A|U)$-module and $E$ and $F$ are $A$-modules, in which case the assertion follows from the definition 
of the morphism \eqref{p1-indmal17b}, see the proof of (\cite{sga4} 12.11(b)). 

(iv) It follows from (iv) by taking for $L$ (resp.\ $E$, resp.\  $F$) $E|U$ (resp.\ $F$, resp.\ $G$), and combining with the adjunction morphism $\rI j_{U!}\rI j_U^*\rightarrow \id$. 

(v) Indeed, for every $(A|U)$-module $M$, the diagram 
\begin{equation}
{\small 
\xymatrix{
{\Hom_{\star|U}(M,\cHom_A(E\otimes_AF,G)|U)}\ar[r]\ar[d]\ar@{}[rd]|{(1)}&{\Hom_{\star|U}(M,\cHom_A(E,\cHom_A(F,G))|U)}\ar[d]\\
{\Hom_{\star}(j_{U!}(M),\cHom_A(E\otimes_AF,G))}\ar@{}[rd]|{(2)}\ar[r]\ar[d]&{\Hom_{\star}(j_{U!}(M),\cHom_A(E,\cHom_A(F,G)))}\ar[d]\\
{\Hom_{\star}(j_{U!}(M)\otimes_AE\otimes_AF,G)}\ar@{}[rd]|{(3)}\ar[r]\ar[d]&{\Hom_{\star}(j_{U!}(M)\otimes_AE,\cHom_A(F,G))}\ar[d]\\
{\Hom_{\star|U}(M\otimes_{A|U}(E|U)\otimes_{A|U}(F|U),G|U)}\ar@{}[rd]|{(4)}\ar[r]\ar[d]&{\Hom_{\star|U}(M\otimes_{A|U}(E|U),\cHom_{A|U}(F|U,G|U))}\ar[d]\\
{\Hom_{\star|U}(M,\cHom_{A|U}((E|U)\otimes_A(F|U),G|U))}\ar[r]&{\Hom_{\star|U}(M,\cHom_{A|U}(E|U,\cHom_{A|U}(F|U,G|U))),}}}
\end{equation}
where $\star=\bIndMod(A)$, $\star|U=\bIndMod(A|U)$, $(1)$ is the functoriality of the adjunction for the morphism \eqref{p1-indmal12d}, 
$(2)$ and $(4)$ correspond to the commutative diagram 
\eqref{p1-indmal11f} and $(3)$ corresponds to the external rectangle of \eqref{p1-indmal17d}, is commutative. 
It also follows from \eqref{p1-indmal17d} that the composition of the vertical morphisms are induced by \eqref{p1-indmal17c}. The proposition follows. 

(vi) If follows from \eqref{p1-indmal17e} and \eqref{p1-indmal17f}, by computing in two different ways the morphism 
\begin{equation}
(E|U)\otimes_{A|U}(\cHom_A(E\otimes_AF,G)|U)\rightarrow \cHom_{A|U}(E|U,F|U)
\end{equation}
corresponding to the diagonal morphism of \eqref{p1-indmal17f} by \ref{p1-indmal11}(iv).

(vii) Indeed, for every $A$-module $M$, the diagram
\begin{equation}
\xymatrix{
{\Hom_{\star}(j_{U!}(M)\otimes_AE,F)}\ar[r]\ar[d]&{\Hom_{\star}(j_{U!}(M)\otimes_AE\otimes_AG,F\otimes_AG)}\ar[d]\\
{\Hom_{\star|U}(M\otimes_{A|U}(E|U),F|U)}\ar[r]&{\Hom_{\star|U}(M\otimes_{A|U}(E|U)\otimes_{A|U}(G|U),(F|U)\otimes_{A|U}(G|U)),}}
\end{equation}
where $\star=\bIndMod(A)$, $\star|U=\bIndMod(A|U)$, the horizontal morphisms are induced by the tensor products with $G$ and $G|U$ \eqref{p1-indmal11i}, 
and the vertical ones by \eqref{p1-indmal17b} and the adjunction, is commutative 
We deduce by (iii) that the diagram 
\begin{equation}
\xymatrix{
{\cHom_A(E,F)|U}\ar[r]\ar[d]&{\cHom_A(E\otimes_AG,F\otimes_AG)|U}\ar[d]\\
{\cHom_{A|U}(E|U,F|U)}\ar[r]&{\cHom_{A|U}((E|U)\otimes_{A|U}(G|U),(F|U)\otimes_{A|U}(G|U)),}}
\end{equation}
where the horizontal morphisms are induced by the tensor products with $G$ and $G|U$ \eqref{p1-indmal11i}, and the vertical ones by \eqref{p1-indmal17a} and \eqref{p1-indmal17c}, is commutative. 
The proposition follows from the commutativity of the diagram above and (vi).

(viii) It follows from (vi) and (vii).

\begin{lem}\label{p1-indmal18}
Let $E,F,G$ be three ind-$A$-modules. Suppose that there exists a finite covering $(U_p)_{1\leq p\leq n}$ of the final object of $X$, such that for every $1\leq p\leq n$, 
either $E|U_i$ or $G|U_i$ is isomorphic to a direct factor of a free $A|U_i$-module of finite type \eqref{p1-indmal1a}. Then, the canonical morphism of ind-$A$-modules \eqref{p1-indmal12d}
\begin{equation}\label{p1-indmal18a}
G\otimes_A\cHom_A(E,F)\rightarrow \cHom_A(E,F\otimes_AG)
\end{equation}
is an isomorphism. 
\end{lem}

Indeed, by \ref{p1-indmal17}(vi) and (\cite{ag2} 2.7.17(iii)), we may assume that either
$E$ or $G$ is isomorphic to a direct factor of a free $A$-module of finite type. 
By functoriality of \eqref{p1-indmal18a}, we may further reduce to the case where 
$E$ or $G$ is isomorphic to a free $A$-module of finite type. The proposition is then obvious.

\begin{defi}\label{p1-indmal20}
Let $\Omega$ be an $A$-module.
\begin{itemize}
\item[(i)]
A {\em Higgs ind-$A$-module with coefficients in $\Omega$}
is a pair $(M,\theta)$ consisting of an ind-$A$-module $M$ and a morphism of ind-$A$-modules
\begin{equation}\label{p1-indmal20a}
\theta\colon M\rightarrow M\otimes_A\Omega
\end{equation}
such that the composition
\begin{equation}\label{p1-indmal20b}
\xymatrix{M\ar[r]^-(0.5)\theta&{M\otimes_A\Omega}\ar[rr]^-(0.5){\theta\otimes_A\id_\Omega}&&
{M\otimes_A\Omega\otimes_A\Omega}\ar[r]&{M\otimes_A\wedge^2\Omega,}}
\end{equation}
where the third morphism is induced by the exterior product $\Omega\otimes_A\Omega\rightarrow \wedge^2\Omega$ \eqref{p1-NC7}, vanishes.
We then say that $\theta$ is a {\em  Higgs $A$-field} on $M$ with coefficients in $\Omega$.
\item[(ii)] If $(M,\theta)$ and $(M',\theta')$ are two Higgs ind-$A$-modules,
a morphism from $(M,\theta)$ to $(M',\theta')$ is a morphism of ind-$A$-modules
$u\colon M\rightarrow M'$ such that $(u\otimes\id_\Omega)\circ \theta=\theta'\circ u$.
\end{itemize}
\end{defi}

Higgs ind-$A$-modules with coefficients in $\Omega$ form a category denoted by $\bIndHM(A,\Omega)$.

\addtocounter{subsubsection}{2}
\addtocounter{equation}{1}

\subsubsection{}\label{p1-indmal20c}
Let $(M,\theta)$ be a Higgs ind-$A$-module with coefficients in $\Omega$. For any $i\geq 1$, we denote by
\begin{equation}\label{p1-indmal20d}
\theta_i\colon M\otimes_A \wedge^i\Omega \rightarrow M\otimes_A \wedge^{i+1}\Omega
\end{equation}
the composition 
\begin{equation}\label{p1-indmal20e}
\xymatrix{{M\otimes_A \wedge^i\Omega}\ar[rr]^-(0.5){\theta\otimes_A\id_{\wedge^i\Omega}}&&
{M\otimes_A\Omega\otimes_A\wedge^i \Omega}\ar[r]&{M\otimes_A\wedge^{i+1}\Omega,}}
\end{equation}
where the second morphism is induced by the exterior product $\Omega\otimes_A\wedge^i \Omega\rightarrow \wedge^{i+1}\Omega$.
We have $\theta_{i+1}\circ \theta_i=0$. The {\em Dolbeault complex of $(M,\theta)$},
denoted by $\mK^\bullet(M,\theta)$, is the complex of cochains of ind-$A$-modules
\begin{equation}\label{p1-indmal20f}
M\stackrel{\theta}{\longrightarrow}M\otimes_A\Omega\stackrel{\theta_1}{\longrightarrow} M\otimes_A\wedge^2\Omega \dots,
\end{equation}
where $M$ is placed in degree $0$ and the differentials are of degree $1$.

\addtocounter{subsubsection}{3}
\addtocounter{equation}{1}

\subsubsection{}\label{p1-indmal20g}
Let $(M,\theta),(M',\theta')$ be two Higgs ind-$A$-modules with coefficients in $\Omega$.
The {\em total Higgs field} on $M\otimes_AM'$ is the $A$-linear morphism
\begin{equation}\label{p1-indmal20h}
\theta_\tot\colon M\otimes_AM'\rightarrow M\otimes_AM'\otimes_A\Omega
\end{equation}
defined by
\begin{equation}\label{p1-indmal20i}
\theta_\tot=\theta\otimes \id_{M'}+\id_{M}\otimes \theta'.
\end{equation}
We say that $(M\otimes_AM',\theta_\tot)$ is {\em the tensor product} of $(M,\theta)$ and $(M',\theta')$.

\addtocounter{subsubsection}{2}
\addtocounter{equation}{1}

\subsubsection{}\label{p1-indmal20j}
Suppose that there exists a finite covering $(U_p)_{1\leq p\leq n}$ of the final object of $X$, such that for every $1\leq p\leq n$, 
$\Omega|U_i$ is a direct factor of a free $A|U_i$-module of finite type.
Let $\rT=\cHom_A(\Omega,A)$ be the dual of $\Omega$, and let $\IndSym_A(\rT)$ the {\em symmetric ind-$A$-algebra of $\rT$} \eqref{p1-indmal6}. 
Observe that the canonical morphism $\Omega\rightarrow \cHom_A(\rT,A)$ is an isomorphism. 
Let $M$ be an ind-$A$-module. The canonical morphism \eqref{p1-indmal12d}
\begin{equation}\label{p1-indmal20j1}
\cEnd_A(M)\otimes_A\Omega\rightarrow \cHom_A(M,M\otimes_A\Omega)
\end{equation}
is an isomorphism by \ref{p1-indmal18}. Similarly, the canonical morphism 
\begin{equation}\label{p1-indmal20j2}
\cHom_A(\rT,A) \otimes_A\cEnd_A(M)\rightarrow \cHom_A(\rT,\cEnd_A(M))
\end{equation}
is an isomorphism. We deduce an isomorphism 
\begin{equation}\label{p1-indmal20j3}
\cHom_A(M,M\otimes_A\Omega)\stackrel{\sim}{\rightarrow}\cHom_A(\rT,\cEnd_A(M)).
\end{equation}
Applied to $A$, the latter gives an isomorphism 
\begin{equation}\label{p1-indmal20j4}
\begin{array}[t]{clcr}
\Hom_{\bIndMod(A)}(M,M\otimes_A\Omega)&\stackrel{\sim}{\rightarrow}&\Hom_{\bIndMod(A)}(\rT,\cEnd_A(M)),\\
\theta&\mapsto&\varphi.
\end{array}
\end{equation}
We claim that $\theta$ is a Higgs field if and only if the diagram 
\begin{equation}\label{p1-indmal20j5}
\xymatrix{
{\rT\otimes_A\rT}\ar[d]_{\sigma}\ar[rr]^-(0.5){\mu\circ (\varphi\otimes \varphi)}&&{\cEnd_A(M),}\\
{\rT\otimes_A\rT}\ar[rru]_{\mu\circ (\varphi\otimes \varphi)}&&}
\end{equation}
where $\mu$ is the multiplication of $\cEnd_A(M)$ \eqref{p1-indmal14} and $\sigma$ is the morphism that switches the factors, is commutative. 
Indeed, we may assume that $\Omega$ is a direct factor of a free $A$-module of finite type by (\cite{ag2} 2.7.17(i)), and even that $\Omega$ is a free $A$-module of finite type, in which case 
the claim is obvious \eqref{p1-indmal130}.  We deduce by \ref{p1-indmal7} that {\em giving a Higgs $A$-field $\theta$ on $M$ with coefficients in $\Omega$ 
is equivalent to giving an ind-$\IndSym_A(\rT)$-module structure on the ind-$A$-module $M$}. 

\begin{lem}\label{p1-indmal23}
Let $\Omega$ be an $A$-module, $\rT=\cHom_A(\Omega,A)$, $\IndSym_A(\rT)$ the {\em symmetric ind-$A$-algebra of $\rT$} \eqref{p1-indmal6}, 
$(M,\theta)$, $(M',\theta')$ two Higgs ind-$A$-modules with coefficients in $\Omega$. 
Suppose that there exists a finite covering $(U_p)_{1\leq p\leq n}$ of the final object of $X$, such that for every $1\leq p\leq n$, 
$\Omega|U_i$ is a direct factor of a free $A|U_i$-module of finite type. Then, 
\begin{itemize}
\item[{\rm (i)}] The Higgs field
\begin{equation}\label{p1-indmal23a}
\theta\colon M\rightarrow M\otimes_A\Omega 
\end{equation}
is a morphism of ind-$\IndSym_A(\rT)$-modules \eqref{p1-indmal20j}. 
\item[{\rm (ii)}] The Higgs $A$-fields $\theta\otimes\id_{M'}$ and $\id_M\otimes\theta'$ on $M\otimes_AM'$ induce
on the quotient $M\otimes_{\IndSym_A(\rT)}M'$ \eqref{p1-indmal5b} the same Higgs $A$-field 
with coefficients in $\Omega$, namely its canonical Higgs field. 
\end{itemize}
\end{lem}

(i) Indeed, the Higgs $A$-field $\theta$ is the composition of the morphisms 
\begin{equation}\label{p1-indmal23b}
M\rightarrow M\otimes_A\rT\otimes_A \Omega\rightarrow M\otimes_A\IndSym_A(\rT)\otimes_A \Omega
\rightarrow M\otimes_A \Omega,
\end{equation}
the first being induced by the section $\id_\Omega$ of $\rT\otimes_A \Omega$ identified with $\cHom_A(\Omega,\Omega)$ \eqref{p1-delta-con1ef},
the second being induced by the canonical injection $\rT\rightarrow \IndSym_A(\rT)$ and the third being induced by the multiplication 
$\upmu\colon M\otimes_A\IndSym_A(\rT)\rightarrow M$ defined by $\theta$ \eqref{p1-indmal20j}. 
Each of these morphisms is obviously a morphism of ind-$\IndSym_A(\rT)$-modules. 

(ii) It follows directly from the definitions.

\begin{defi}\label{p1-indmal21}
Let $B$ be a commutative ind-$A$-algebra, $\Omega$ an ind-$B$-module \eqref{p1-indmal2}. A morphism of ind-$A$-modules $\delta\colon B\rightarrow \Omega$ is an {\em $A$-derivation} 
if the diagram 
\begin{equation}\label{p1-indmal21a}
\xymatrix{
{B\otimes_AB}\ar[r]^-(0.5){\mu_B}\ar[d]_{\delta\otimes \id_B+\id_B\otimes \delta}&B\ar[d]^{\delta}\\
{B\otimes_A\Omega}\ar[r]^-(0.5){\mu_\Omega}&{\Omega,}}
\end{equation}
where $\mu_B$ (resp.\ $\mu_\Omega$) is the multiplication of $B$ (resp.\ $\Omega$), is commutative.
\end{defi}

\begin{defi}\label{p1-indmal22} 
Let $B$ be a commutative ind-$A$-algebra \eqref{p1-indmal2}, $\Omega$ an $A$-module, $\delta\colon B\rightarrow \Omega\otimes_AB$ an $A$-derivation \eqref{p1-indmal21},
which is also a Higgs $A$-field on $B$ with coefficients in $\Omega$ \eqref{p1-indmal20}. 
\begin{itemize}
\item[(i)] An {\em ind-$B$-module with $\delta$-connection} is a pair $(M,\nabla)$ consisting of an ind-$B$-module $M$  and 
a morphism of ind-$A$-modules 
\begin{equation}\label{p1-indmal22a}
\nabla\colon M\rightarrow \Omega\otimes_AM
\end{equation}
such that the diagram of morphisms of ind-$A$-modules 
\begin{equation}
\xymatrix{{B\otimes_AM}\ar[rrr]^-(0.5){\id_B\otimes\nabla+\delta\otimes\id_M}\ar[d]_{\mu_M}&&&{B\otimes_A\Omega\otimes_AM}\ar[d]^{\id_\Omega\otimes\mu_{M}}\\
M\ar[rrr]^-(0.5)\nabla&&&{\Omega\otimes_AM,}}
\end{equation}
where $\mu_M$ is the multiplication of $M$ (\cite{ag2} (2.7.5.1)), is commutative. 
We also say that $\nabla$ is a {\em $\delta$-connection on $M$}. 
We say that $\nabla$ is {\em integrable} if it is a Higgs $A$-field on $M$ with coefficients in $\Omega$. 
\item[(ii)] If $(M,\nabla)$ and $(M',\nabla')$ are two ind-$B$-modules with $\delta$-connection,
a morphism from $(M,\nabla)$ to $(M',\nabla')$ is a morphism of ind-$B$-modules
$u\colon M\rightarrow M'$ such that $(u\otimes_A\id_\Omega)\circ \nabla=\nabla'\circ u$.
\end{itemize}
\end{defi}

The ind-$B$-modules with $\delta$-connection form a category that we denote by $\bIndMC_\delta(B/A)$. We denote by $\bIndMIC_\delta(B/A)$ the full subcategory 
made of ind-$B$-modules with integrable $\delta$-connection.

\subsection{}\label{p1-delta-con9}
We take again the assumption and notation of \ref{p1-indmal22}, and 
let $(M,\nabla)$ be an ind-$B$-module with $\delta$-connection, $(N,\theta)$ a Higgs ind-$A$-module with coefficients in $\Omega$. 
Then, the morphism of ind-$A$-modules
\begin{equation}\label{p1-delta-con9a}
\nabla'\colon M\otimes_A N\rightarrow \Omega\otimes_A M\otimes_A N
\end{equation}
defined by
\begin{equation}\label{p1-delta-con9b}
\nabla'=\nabla\otimes_A \id_N+ \id_M\otimes_A \theta,
\end{equation}
is a $\delta$-connection on $M\otimes_AN$. If $\nabla$ is integrable, so is $\nabla'$ \eqref{p1-indmal20g}.

\subsection{}\label{p1-delta-con10} 
We take again the assumption and notation of \ref{p1-indmal22}, and let $B'$ be a commutative ind-$A$-algebra, $\Omega'$ an $A$-module, 
$\delta'\colon B'\rightarrow \Omega'\otimes_AB'$ an $A$-derivation, which is also a Higgs $A$-field on $B'$ with coefficients in $\Omega'$.
Let $\iota\colon B\rightarrow B'$ be a homomorphism of ind-$A$-algebras, 
$u\colon \Omega\rightarrow \Omega'$ an $A$-linear morphism such that $\delta'\circ \iota=(u\otimes \iota) \circ \delta$. 
Let $(M,\nabla)$ be an ind-$B$-module with $\delta$-connection. Then, there exists a unique morphism of ind-$A$-modules
\begin{equation}\label{p1-delta-con10a}
\nabla'\colon M\otimes_BB'\rightarrow \Omega'\otimes_AM\otimes_BB'
\end{equation}
that fits into a commutative diagram 
\begin{equation}
\xymatrix{
{M\otimes_AB'}\ar[rrrr]^-(0.5){((u\otimes_A\id_M)\circ \nabla)\otimes\id_{B'}+\id_M\otimes \delta'}\ar[d]_\gamma&&&&{\Omega'\otimes_AM\otimes_AB'}\ar[d]^{\id_{\Omega'}\otimes \gamma}\\
{M\otimes_BB'}\ar[rrrr]^-(0.5){\nabla'}&&&&{\Omega'\otimes_AM\otimes_BB',}}
\end{equation}
where $\gamma$ is the canonical morphism \eqref{p1-indmal5}. 
It is a $\delta'$-connection on $M\otimes_BB'$, which is integrable, if so is $\nabla$ \eqref{p1-indmal20g}.

\begin{exemples}
Let $B$ a commutative $A$-algebra, $d\colon B\rightarrow \Omega^1_{B/A}$ the universal $A$-derivation, $\Omega_{B/A}=\wedge_B\Omega^1_{B/A}$.  
We denote also by $d\colon \Omega_{B/A}\rightarrow \Omega_{B/A}$ the unique $A$-anti-derivation of degree $1$ and square zero that extends $d$. 
We suppose that there exists an $A$-module $\Omega$ and a $B$-isomorphism $\upgamma\colon \Omega\otimes_AB\stackrel{\sim}{\rightarrow}\Omega^1_{B/A}$
such that for every local section $\omega$ of $\Omega$, we have $d(\upgamma(\omega\otimes 1))=0$. 
Let $\lambda\in A(X)$. Consider the $A$-derivation $\delta=\lambda \upgamma^{-1}\circ  d\colon B\rightarrow \Omega\otimes_AB$, which is a Higgs $A$-field with coefficients in $\Omega$. 
Let $M$ be an ind-$B$-module. A morphism of ind-$A$-modules 
\begin{equation}
\nabla\colon M\rightarrow \Omega\otimes_AM 
\end{equation}
is a $\delta$-connection (resp.\ an integrable $\delta$-connection) if and only if $\tnabla=(\upgamma\otimes \id_M)\circ \nabla\colon M\rightarrow \Omega^1_{B/A}\otimes_BM$  
(\cite{ag2} (2.7.5.1)) is a $\lambda$-connection (resp.\ an integrable $\lambda$-connection) with respect to the extension $B/A$ in the sense of (\cite{ag2} 2.8.4) (resp.\ \cite{ag2} 2.8.8).  
Indeed, the diagram 
\begin{equation}
\xymatrix{
{\Omega\otimes_AM}\ar[r]^-(0.5){-\id\otimes_A\nabla}\ar[d]&{\Omega\otimes_A\Omega\otimes_AM}\ar[rr]^-(0.5){w\otimes_A\id_M}&&{(\wedge^2\Omega)\otimes_AM}\ar[d]\\
{\Omega^1_{B/A}\otimes_BM}\ar[rrr]^-(0.5){\tnabla}&&&{\Omega^2_{B/A}\otimes_BM,}}
\end{equation}
where $w\colon \Omega\otimes_A\Omega\rightarrow \wedge^2\Omega$ is the exterior product, 
the lower horizontal map $\tnabla$ is defined in (\cite{ag2} 2.8.6) and 
the vertical maps are the isomorphisms induced by $\upgamma$, is clearly commutative. 
\end{exemples}

\section{Base change}\label{p1-bcim}

\subsection{}\label{p1-bcim15}
Let $f_*\colon \cA\rightarrow \cB$ be a left exact functor of abelian $\mU$-categories \eqref{p1-NC0} having a right derived functor 
$\rR f_*\colon \bD^+(\cA)\rightarrow \bD^+(\cB)$ (\cite{ks2} 13.3.1) which is a triangulated functor.  
For every integer $q$, we set
\begin{equation}\label{p1-bcim15a}
\rR^q f_*=\cH^q\circ \rR f_* \colon \bD^+(\cA)\rightarrow \cB. 
\end{equation}
For every $X\in \ob(\cA)$ and every quasi-isomorphism $X[0]\rightarrow F$ of $\bC^+(\cA)$, 
the composed morphism $X[0]\rightarrow F\rightarrow \tau_{\geq 0}(F)$ is a quasi-isomorphism \eqref{p1-NC8b3}. 
We deduce a canonical isomorphism of $\bInd(\bD^+(\cB))$
\begin{equation}\label{p1-bcim15e}
\nu\colon \rR f_*(X[0])\stackrel{\sim}{\rightarrow} \underset{\underset{X[0]\rightarrow F\in \cQ_X}{\longrightarrow}}{\mlq\mlq\lim \mrq\mrq} f_*(F),
\end{equation}
where $\cQ_X$ is the category of quasi-isomorphisms $X[0]\rightarrow F$ of $\bC^+(\cA)$ such that $F^i=0$ for every $i<0$. 
By \ref{p1-abisoind9}, the isomorphism $\nu$ determines the following data:
\begin{itemize}
\item[(i)] for every $F\in \ob(\cQ_X)$, a functorial morphism $\nu^\flat_F\colon f_*(F) \rightarrow \rR f_*(X[0])$ of $\bD^+(\cB)$,
\item[(ii)] a morphism $u\colon G\rightarrow G'$ of $\cQ_X$, and 
\item[(iii)] a morphism $\nu_G\colon \rR f_*(X[0])\rightarrow f_*(G)$ of $\bD^+(\cB)$,
\end{itemize} 
such that the following conditions are satisfied: 
\begin{itemize}
\item[(a)] The composed morphism $\nu^\flat_G\circ \nu_G\colon \rR f_*(X[0])\rightarrow f_*(G)\rightarrow \rR f_*(X[0])$ is the identity of $\rR f_*(X[0])$.
\item[(b)] The composed morphism 
\begin{equation}
\xymatrix{
f_*(X)[0]\ar[r]^-(0.5){\nu^\flat_{X[0]}}&{\rR f_*(X[0])}\ar[r]^-(0.5){\nu_G}&{f_*(G)}\ar[r]^-(0.5){f_*(u)}&{f_*(G')}}
\end{equation}
is the canonical morphism. 
\end{itemize}

We deduce from (a) that $\rR^if_*(X[0])=0$ for every $i<0$.  Moreover, for every $F\in \ob(\cQ_X)$, the canonical morphism $X[0]\rightarrow F$ induces an isomorphism 
$f_*(X)\stackrel{\sim}{\rightarrow} \cH^0(f_*(F))$. Hence (a) and (b) imply that the canonical morphism $f_*(X)\rightarrow \rR^0f_*(X[0])$ 
is an isomorphism with inverse induced by $\cH^0(\nu_G)$. 
 
Similarly, for every $F\in \ob(\bD^+(\cA))$ and every integer $q$ such that $\cH^i(F)=0$ for all $i< q$, 
we have $\rR^i f_*(F)=0$ for every $i< q$; moreover, the canonical distinguished triangle \eqref{p1-NC8b2}
\begin{equation}\label{p1-bcim15b}
\cH^q(F)[-q]\rightarrow F\rightarrow \tau_{\geq q+1} (F)\rightarrow \cH^q(F)[-q+1]
\end{equation}
induces an isomorphism 
\begin{equation}\label{p1-bcim15c}
f_*(\cH^q(F))\stackrel{\sim}{\rightarrow} \rR^qf_*(F).
\end{equation}
Therefore, for every $G\in \ob(\bD^+(\cA))$, the canonical morphism $G\rightarrow \tau_{\geq q}(G)$ induces a canonical morphism 
\begin{equation}\label{p1-bcim15d}
\rR^qf_*(G)\rightarrow f_*(\cH^q(G)). 
\end{equation}

\begin{lem}\label{p1-bcim16}
Let $f_*\colon \cA\rightarrow \cB$ be a left exact functor of abelian $\mU$-categories \eqref{p1-NC0} having a right derived functor 
$\rR f^*\colon \bD^+(\cA)\rightarrow \bD^+(\cB)$, $\cB_0$ a thick subcategory of $\cB$, i.e. a full subcategory of $\cB$ 
closed by kernels, cokernels and extensions (see {\rm \cite{ks2} 8.3.21}).
\begin{itemize}
\item[{\rm (i)}]  Let  $F^\bullet$ be a bounded from below complex of $\cA$ such that 
for all integers $q$ and $r$, $\rR^qf_*(F^r)$ is an object of $\cB_0$, Then, for every integer $q$, $\rR^qf_*(F^\bullet)$ is an object of $\cB_0$. 
\item[{\rm (ii)}]  Let  $F$ be an object of $\bD^+(\cA)$ such that for all integers $q$ and $r$, $\rR^qf_*(\cH^r(F))$ is an object of $\cB_0$.
Then, for every integer $q$, $\rR^qf_*(F)$ is an object of $\cB_0$. 
\end{itemize}
\end{lem}

(i) Indeed, for every integer $r$, we have an exact sequence \eqref{p1-NC8d1}
\begin{equation}\label{p1-bcim16a}
0\rightarrow F^r[-r]\rightarrow \sigma_{\leq r}(F^\bullet)\rightarrow \sigma_{\leq r-1}(F^\bullet)\rightarrow 0.
\end{equation}
Since $\sigma_{\leq r}(F^\bullet)=0$ for $r\ll 0$, we deduce, by induction, that for all integers $r$ and $q$,
$\rR^qf_*(\sigma_{\leq r}(F^\bullet))$ is an object of $\cB_0$. 
For every integer $r>q$, the canonical exact sequence \eqref{p1-NC8d3}
\begin{equation}
0\rightarrow \sigma_{\geq r+1}(F^\bullet) \rightarrow F^\bullet \rightarrow \sigma_{\leq r}(F^\bullet) \rightarrow 0
\end{equation}
induces an isomorphism $\rR^qf_*(F)\stackrel{\sim}{\rightarrow} \rR^qf_*(\sigma_{\leq r}(F))$, hence the proposition. 

(ii) Indeed, for every integer $r$, we have a distinguished triangle \eqref{p1-NC8b1}
\begin{equation}\label{p1-bcim16b}
\tau_{\leq r-1}(F)\rightarrow \tau_{\leq r}(F)\rightarrow \cH^r(F)[-r]\rightarrow \tau_{\leq r-1}(F)[1]. 
\end{equation}
Since $\tau_{\leq r}(F)=0$ for $r\ll 0$, we deduce, by induction, that for all integers $r$ and $q$,
$\rR^qf_*(\tau_{\leq r}(F))$ is an object of $\cB_0$. 
For every integer $r\geq q$, the canonical distinguished triangle \eqref{p1-NC8b3}
\begin{equation}
\tau_{\leq r}(F) \rightarrow F \rightarrow \tau_{\geq r+1}(F) \rightarrow \tau_{\leq r}(F) [+1]
\end{equation}
induces an isomorphism $\rR^qf_*(\tau_{\leq r}(F))\stackrel{\sim}{\rightarrow} \rR^qf_*(F)$, hence the proposition.

\subsection{}\label{p1-bcim13}
Let $f^*\colon \cB\rightarrow \cA$ be a right exact functor of abelian $\mU$-categories \eqref{p1-NC0}. 
Then, for every complex $G^\bullet$ of $\cB$, and every integer $q$, 
we have a canonical functorial morphism of $\cA$ 
\begin{equation}\label{p1-bcim13a}
f^*(\cH^q(G^\bullet))\rightarrow \cH^q(f^*(G^\bullet)),
\end{equation}
where the pullback $f^*(G^\bullet)$ is defined term by term (not derived). 
Indeed, let $\cZ^q(G^\bullet)$ be the kernel of the differential $d^q_{G^\bullet}\colon G^q\rightarrow G^{q+1}$ and  
\begin{equation}\label{p1-bcim13b}
G^{q-1}\stackrel{\delta^{q-1}_{G^\bullet}}{\longrightarrow} \cZ^q(G^\bullet)\stackrel{i^{q}_{G^\bullet}}{\longrightarrow} G^q
\end{equation} 
the canonical factorization of $d^{q-1}_{G^\bullet}$. Then, $f^*(i^q_{G^\bullet})$ induces a morphism $\varepsilon^q\colon f^*(\cZ^q(G^\bullet))\rightarrow \cZ^q(f^*(G^\bullet))$ that makes the diagram
\begin{equation}\label{p1-bcim13c}
\xymatrix{
{f^*(G^{q-1})}\ar[rr]^-(0.5){f^*(\delta^{q-1}_{G^\bullet})}\ar[rrd]_{\delta^{q-1}_{f^*(G^\bullet)}}&&{f^*(\cZ^q(G^\bullet))}\ar[d]^{\varepsilon^q}
\ar[rrd]^{f^*(i^{q}_{G^\bullet})}&&\\
&&{\cZ^q(f^*(G^\bullet))}\ar[rr]^-(0.4){i^{q}_{f^*(G^\bullet)}}&&{f^*(G^q)}}
\end{equation}
commutative. The morphism \eqref{p1-bcim13a} follows since $f^*$ is right exact.

Let $g^*\colon \cC\rightarrow \cB$ be another right exact functor of abelian $\mU$-categories. 
It immediately follows from \eqref{p1-bcim13c} that for every complex $H^\bullet$ of $\cC$, and every integer $q$, the diagram of morphisms of $\cA$ 
\begin{equation}\label{p1-bcim13d}
\xymatrix{
{f^*(g^*(\cH^q(H^\bullet))}\ar[r]\ar[rd]&{f^*(\cH^q(g^*(H^\bullet)))}\ar[d]\\
&{f^*(g^*(\cH^q(H^\bullet))),}}
\end{equation}
where all arrows are induced by the morphisms \eqref{p1-bcim13a} for the functors $f^*$, $g^*$ and $f^*\circ g^*$, is commutative.

\subsection{}\label{p1-bcim10}
We consider a diagram of additive functors of abelian $\mU$-categories \eqref{p1-NC0}
\begin{equation}\label{p1-bcim10a}
\xymatrix{
{\cA'}\ar[r]^{g'_*}\ar[d]_{f'_*}&{\cA}\ar[d]^{f_*}\\
{\cB'}\ar[r]^{g_*}&{\cB}}
\end{equation}
and an isomorphism   
\begin{equation}\label{p1-bcim10b}
f_*\circ g'_*\stackrel{\sim}{\rightarrow}g_*\circ f'_*.
\end{equation}
We assume that the following conditions hold:
\begin{itemize}
\item[(i)] The functors $f_*,f'_*,g_*$ and $g'_*$ admit left adjoint functors $f^*,f'^*,g^*$ and $g'^*$, respectively. 
In particular, the functors $f_*,f'_*,g_*$ and $g'_*$ are left exact. 
\item[(ii)] The functors $f_*,f'_*,g_*$ and $g'_*$ are right derivable for $\bD^+(\cA)$, $\bD^+(\cA')$, $\bD^+(\cB')$ and $\bD^+(\cA')$, 
respectively (\cite{ks2} 13.3.1), which are triangulated functors.
\item[(iii)] The functor $\rR f_*\circ \rR g'_*$ (resp.\ $\rR g_* \circ \rR f'_*$) is a right derived functor of $f_*\circ g'_*$ (resp.\ $g_* \circ f'_*$). 
Hence, the isomorphism \eqref{p1-bcim10b} induces an isomorphism 
\begin{equation}\label{p1-bcim10c}
(\rR f_*)\circ (\rR g'_*)\stackrel{\sim}{\rightarrow}(\rR g_*)\circ (\rR f'_*).
\end{equation}
\end{itemize}

For every bounded from below complex $F^\bullet$ of $\cA$ and every integer $q$, 
there exists a canonical functorial morphism of $\cB'$,
\begin{equation}\label{p1-bcim10d}
g^*(\rR^q f_*(F^\bullet))\rightarrow \rR^q f'_*(g'^*(F^\bullet)),
\end{equation}
where the pullback $g'^*(F^\bullet)$ is defined term by term (not derived), called {\em base change morphism}. 
Indeed, this amounts to giving a morphism
\begin{equation}\label{p1-bcim10e}
\rR^q f_*(F^\bullet)\rightarrow g_*(\rR^q f'_*(g'^*(F^\bullet))),
\end{equation}
and we take the composed morphism
\begin{eqnarray}\label{p1-bcim10f}
\lefteqn{\rR^q f_*(F^\bullet)\rightarrow \rR^q f_*(g'_*(g'^*(F^\bullet)))\rightarrow} \\
&&\rR^qf_*(\rR g'_*(g'^*(F^\bullet)))\stackrel{\sim}{\rightarrow} \rR^qg_*(\rR f'_*(g'^*(F^\bullet)))
\rightarrow g_*(\rR^qf'_*(g'^*(F^\bullet))),\nonumber
\end{eqnarray}
where the first arrow is induced by the adjunction morphism
$\id\rightarrow g'_*g'^*$, defined term by term, the second arrow by the canonical morphism $g'_*\rightarrow \rR g'_*$, 
where $g'_*$ is defined for complexes term by term, 
the third arrow by the isomorphism \eqref{p1-bcim10c} and the fourth arrow by the morphism \eqref{p1-bcim15d}. 

Note that the morphism \eqref{p1-bcim10d} is not defined for objects of $\bD^+(\cA)$. 

\begin{rema}\label{p1-bcim11}
We keep the assumptions and notation of \ref{p1-bcim10}. For every $F'\in \ob(\bD^+(\cA'))$, we have a canonical morphism \eqref{p1-bcim15d}
\begin{equation}\label{p1-bcim11b}
\rR^q g_*(\rR f'_*(F'))\rightarrow g_*(\rR^qf_*(F')).
\end{equation}

If the second hypercohomology spectral sequence of the functor $g_*$ 
with respect to the complex $\rR f'_*(F')$ exists (\cite{ega3} 0.11.4.3), the morphism \eqref{p1-bcim11b} 
coincides with the edge-homomorphism of this spectral sequence.

If $F'\in \ob(\cA')$ and if the Cartan-Leray spectral sequence for the composed functor $g_*\circ f_*$ exists, the morphism
\eqref{p1-bcim11b} coincides with the edge-homomorphism of this spectral sequence.
\end{rema}

\begin{prop}\label{p1-bcim14}
We take again the assumptions and notation of \ref{p1-bcim10}, moreover, let $\cB_0$ be a thick subcategory of $\cB$ 
such that the restriction $g_0^*\colon \cB_0\rightarrow \cB'$ of the functor $g^*$ is exact, and let 
$F^\bullet$ be a bounded from below complex of $\cA$. Assume that one of the following conditions is satisfied: 
\begin{itemize}
\item[{\rm (i)}] For all integers $q$ and $r$, $\rR^q f_*(F^r)$ is an object of $\cB_0$ and the base change morphism \eqref{p1-bcim10d}
\begin{equation}\label{p1-bcim14a}
g^*(\rR^q f_*(F^r))\rightarrow \rR^q f'_*(g'^*(F^r))
\end{equation}
is an isomorphism. 
\item[{\rm (ii)}] For all integers $q$ and $r$, $\rR^q f_*(\cH^r(F^\bullet))$ is an object of $\cB_0$ and the morphisms 
\begin{eqnarray}
g'^*(\cH^r(F^\bullet))&\rightarrow& \cH^r(g'^*(F^\bullet)),\label{p1-bcim14b}\\
g^*(\rR^q f_*(\cH^r(F^\bullet)))&\rightarrow& \rR^q f'_*(g'^*(\cH^r(F^\bullet))),\label{p1-bcim14c}
\end{eqnarray}
defined in \eqref{p1-bcim13a} and \eqref{p1-bcim10d} respectively, are isomorphisms. 
\end{itemize}
Then, the base change morphism \eqref{p1-bcim10d}
\begin{equation}\label{p1-bcim14d}
g^*(\rR^q f_*(F^\bullet))\rightarrow \rR^q f'_*(g'^*(F^\bullet))
\end{equation}
is an isomorphism. 
\end{prop}

(i) Let $r,q$ be two integers. The canonical exact sequence \eqref{p1-NC8d1}
\begin{equation}\label{p1-bcim14f1}
0\longrightarrow F^r[-r]\longrightarrow \sigma_{\leq r}(F^\bullet)\stackrel{\pi_r}{\longrightarrow} \sigma_{\leq r-1}(F^\bullet)\longrightarrow 0
\end{equation}
induces a long exact sequence of cohomology 
\begin{equation}\label{p1-bcim14f2}
\rightarrow \rR^{q-r}f_*(F^r)\rightarrow \rR^qf_*(\sigma_{\leq r}(F^\bullet))\rightarrow \rR^qf_*(\sigma_{\leq r-1}(F^\bullet))\rightarrow 
\end{equation}
We need to describe explicitly the connecting morphism. 
We denote by $C_r$ the mapping cone of $\pi_r$ \eqref{p1-bcim14f1}. We have a canonical distinguished triangle of complexes of $\cA$,
\begin{equation}\label{p1-bcim14e1}
\sigma_{\leq r}(F^\bullet)\stackrel{\pi_r}{\longrightarrow} \sigma_{\leq r-1}(F^\bullet)\longrightarrow C_r\longrightarrow \sigma_{\leq r}(F^\bullet) [1], 
\end{equation}
and a canonical quasi-isomorphism 
\begin{equation}\label{p1-bcim14e2}
\iota_r\colon F^r[-r+1]\rightarrow C_r. 
\end{equation}
We deduce a long exact sequence 
\begin{equation}\label{p1-bcim14e3}
\rR^qf_*(\sigma_{\leq r}(F^\bullet))\rightarrow \rR^qf_*(\sigma_{\leq r-1}(F^\bullet))\rightarrow \rR^qf_*(C_r)\rightarrow \rR^{q+1}f_*(\sigma_{\leq r}(F^\bullet)), 
\end{equation}
and an isomorphism 
\begin{equation}\label{p1-bcim14e4}
\rR^{q-r+1}f_*(F^r)\stackrel{\sim}{\rightarrow} \rR^qf_*(C_r),
\end{equation}
that induce the exact sequence \eqref{p1-bcim14f2}.
For all integers $s,t$, $\rR^tf_*(\sigma_{\leq s}(F^\bullet))$ and $\rR^tf_*(C^\bullet_s)$ are objects of $\cB_0$, by \ref{p1-bcim16}(i).
Since the restriction $g_0^*\colon \cB_0\rightarrow \cB'$ of the functor $g^*$ is exact, the long sequence 
\begin{equation}\label{p1-bcim14g}
g^*(\rR^qf_*(\sigma_{\leq r}(F^\bullet)))\rightarrow g^*(\rR^qf_*(\sigma_{\leq r-1}(F^\bullet)))\rightarrow g^*(\rR^qf_*(C_r))\rightarrow g^*(\rR^{q+1}f_*(\sigma_{\leq r}(F^\bullet)))
\end{equation}
is exact.

Similarly, we denote by $\pi'_r\colon \sigma_{\leq r}(g'^*(F^\bullet))\rightarrow \sigma_{\leq r-1}(g'^*(F^\bullet))$ the canonical morphism 
and by $C'_r$ its mapping cone. We have a canonical distinguished triangle of complexes of $\cA'$,
\begin{equation}\label{p1-bcim14h1}
\sigma_{\leq r}(g'^*(F^\bullet))\stackrel{\pi'_r}{\longrightarrow} \sigma_{\leq r-1}(g'^*(F^\bullet))\longrightarrow C'_r\longrightarrow \sigma_{\leq r}(g'^*(F^\bullet)) [1], 
\end{equation}
and a canonical quasi-isomorphism 
\begin{equation}\label{p1-bcim14h2}
\iota'_r\colon g'^*(F^r)[-r+1]\rightarrow C'_r. 
\end{equation}
Since $\pi'_r=g'^*(\pi_r)$, we see that \eqref{p1-bcim14h1} is canonically identified with the image of \eqref{p1-bcim14e1} by $g'^*$. We deduce a commutative diagram 
\begin{equation}\label{p1-bcim14h3}
\xymatrix{
{g^*(\rR^qf_*(\sigma_{\leq r}(F^\bullet)))}\ar[r]\ar[d]&{g^*(\rR^qf_*(\sigma_{\leq r-1}(F^\bullet)))}\ar[r]\ar[d]&{g^*(\rR^qf_*(C_r))}\ar[r]\ar[d]&{g^*(\rR^{q+1}f_*(\sigma_{\leq r}(F^\bullet)))}\ar[d]\\
{\rR^qf'_*(\sigma_{\leq r}(g'^*(F^\bullet)))}\ar[r]&{\rR^qf'_*(\sigma_{\leq r-1}(g'^*(F^\bullet)))}\ar[r]&{\rR^qf'_*(C'_r)}\ar[r]&{\rR^{q+1}f'_*(\sigma_{\leq r}(g'^*(F^\bullet)),}}
\end{equation}
where the horizontal lines are induced by \eqref{p1-bcim14g} and \eqref{p1-bcim14h1} and 
the vertical arrows are the base change morphisms \eqref{p1-bcim10d}. 
Moreover, we have a commutative diagram
\begin{equation}\label{p1-bcim14h4}
\xymatrix{
{g^*(\rR^{q-r+1}f_*(F^r))}\ar[r]\ar[d]&{g^*(\rR^qf_*(C_r))}\ar[d]\\
{\rR^{q-r+1}f'_*(g'^*(F^r))}\ar[r]&{\rR^{q-r+1}f'_*(C'_r),}}
\end{equation}
where the horizontal arrows are the isomorphisms induced by $\iota_r$ and $\iota'_r$ and the vertical arrows are the base change morphisms \eqref{p1-bcim10d}. 
By assumption, the left vertical arrow of \eqref{p1-bcim14h4} is an isomorphism. 
As $\sigma_{\leq r}(F^\bullet)=0$ for $r\ll 0$, we deduce by induction on $r$ that the base change morphism 
\begin{equation}\label{p1-bcim14k}
g^*(\rR^q f_*(\sigma_{\leq r}(F^\bullet)))\rightarrow \rR^q f'_*(g'^*(\sigma_{\leq r}(F^\bullet)))
\end{equation}
is an isomorphism. Consider the commutative diagram 
\begin{equation}\label{p1-bcim14l}
\xymatrix{
{g^*(\rR^q f_*(F^\bullet))}\ar[r]\ar[d]&{g^*(\rR^q f_*(\sigma_{\leq r}(F^\bullet)))}\ar[d]\\
{\rR^q f'_*(g'^*(F^\bullet))}\ar[r]&{\rR^q f'_*(g'^*(\sigma_{\leq r}(F^\bullet))),}}
\end{equation}
where the horizontal arrows are induced by the canonical morphisms and the vertical arrows are the base change morphisms. 
The proposition follows since for every $r>q$, the horizontal arrows are isomorphisms. 

(ii) Let $r,q$ be two integers. We denote by $u^\bullet_r\colon \tau_{\leq r-1}(F^\bullet)\rightarrow \tau_{\leq r}(F^\bullet)$
the canonical morphism \eqref{p1-NC8a} and by $C^\bullet_r$ its mapping cone. We have a canonical distinguished triangle of complexes of $\cA$,
\begin{equation}\label{p1-bcim14m}
\tau_{\leq r-1}(F^\bullet)\stackrel{u^\bullet_r}{\longrightarrow}
\tau_{\leq r}(F^\bullet) \longrightarrow C^\bullet_r \longrightarrow \tau_{\leq {r-1}}(F^\bullet) [1].
\end{equation}
We denote by $u'^\bullet_r\colon \tau_{\leq r-1}(g'^*(F^\bullet))\rightarrow \tau_{\leq r}(g'^*(F^\bullet))$
the canonical morphism and by $C'^\bullet_r$ its mapping cone. We have a commutative diagram of complexes of $\cA'$
\begin{equation}\label{p1-bcim14n}
\xymatrix{
{g'^*(\tau_{\leq r-1}(F^\bullet))}\ar[rr]^-(0.5){g'^*(u^\bullet_r)}\ar[d]_{t^\bullet_{r-1}}&&{g'^*(\tau_{\leq r}(F^\bullet))}\ar[d]^{t^\bullet_{r}}\\
{\tau_{\leq r-1}(g'^*(F^\bullet))}\ar[rr]^-(0.5){u'^\bullet_r}&&{\tau_{\leq r}(g'^*(F^\bullet)),}}
\end{equation}
where the vertical arrows are the canonical morphisms; see \eqref{p1-bcim13c}. 
Since the mapping cone of $g'^*(u^\bullet_r)$ identifies with $g'^*(C^\bullet_r)$, 
we deduce a canonical morphism of complexes $v_r^\bullet\colon g'^*(C^\bullet_r) \rightarrow C'^\bullet_r$ 
that fits into a commutative diagram 
\begin{equation}\label{p1-bcim14o}
\xymatrix{
{g'^*(\tau_{\leq r-1}(F^\bullet))}\ar[r]\ar[d]&{g'^*(\tau_{\leq r}(F^\bullet))}\ar[d]\ar[r]&{g'^*(C^\bullet_r)}\ar[r]\ar[d]^{v_r^\bullet}&
{g'^*(\tau_{\leq r-1}(F^\bullet)[1])}\ar[d]\\
{\tau_{\leq r-1}(g'^*(F^\bullet))}\ar[r]&{\tau_{\leq r}(g'^*(F^\bullet))}\ar[r]&{C'^\bullet_r}\ar[r]&{\tau_{\leq r-1}(g'^*(F^\bullet))[1],}}
\end{equation}
where the upper (resp.\ lower) line is the image of the sequence \eqref{p1-bcim14m} by $g'^*$ (resp.\ the canonical sequence). 
By construction \eqref{p1-bcim13c}, we have a commutative diagram 
\begin{equation}\label{p1-bcim14p}
\xymatrix{
{g'^*(C^\bullet_r)}\ar[r]\ar[d]_{v_r^\bullet}&{g'^*(\cH^r(F^\bullet))[-r]}\ar[d]\\
{C'^\bullet_r}\ar[r]&{\cH^r(g'^*(F^\bullet))[-r],}}
\end{equation}
where the horizontal arrows are the canonical morphisms and 
the right vertical arrow is induced by the morphism \eqref{p1-bcim14b}. 

We have a commutative diagram of $\cB'$
\begin{equation}\label{p1-bcim14q}
{\small
\xymatrix{
{g^*(\rR^qf_*(\tau_{\leq r-1}(F^\bullet)))}\ar[r]\ar[d]&{g^*(\rR^qf_*(\tau_{\leq r}(F^\bullet)))}\ar[r]\ar[d]&
{g^*(\rR^qf_*(C^\bullet_r))}\ar[r]\ar[d]&{g^*(\rR^{q+1}f_*(\tau_{\leq r-1}(F^\bullet)))}\ar[d]\\
{\rR^qf'_*(g'^*(\tau_{\leq r-1}(F^\bullet)))}\ar[r]\ar[d]&{\rR^qf'_*(g'^*(\tau_{\leq r}(F^\bullet)))}\ar[r]\ar[d]&
{\rR^qf'_*(g'^*(C^\bullet_r))}\ar[r]\ar[d]^{\rR^qf'_*(v^\bullet_r)}&{\rR^{q+1}f'_*(g'^*(\tau_{\leq r-1}(F^\bullet)))}\ar[d]\\
{\rR^qf'_*(\tau_{\leq r-1}(g'^*(F^\bullet)))}\ar[r]&{\rR^qf'_*(\tau_{\leq r}(g'^*(F^\bullet)))}\ar[r]&
{\rR^qf'_*(C'^\bullet_r)}\ar[r]&{\rR^{q+1}f'_*(\tau_{\leq r-1}(g'^*(F^\bullet))),}}}
\end{equation}
where the upper vertical arrows are the base change morphisms \eqref{p1-bcim10d} and the lower ones are induced by \eqref{p1-bcim14o}. 
For all integers $s,t$, $\rR^tf_*(\tau_{\leq s}(F^\bullet))$ and $\rR^tf_*(C^\bullet_s)$ are objects of $\cB_0$, by \ref{p1-bcim16}(ii). Since the restriction $g_0^*\colon \cB_0\rightarrow \cB'$ of $g^*$ is exact, we deduce that the upper line of \eqref{p1-bcim14q} is exact. The lower line is obviously exact.

We have commutative diagrams 
\begin{equation}
\xymatrix{
{g^*(\rR^qf_*(C^\bullet_r))}\ar[r]\ar[d]&{g^*(\rR^{q-r}f_*(\cH^r(F^\bullet)))}\ar[d]\\
{\rR^qf'_*(g'^*(C^\bullet_r))}\ar[r]\ar[d]_{\rR^qf'_*(v_r^\bullet)}&{\rR^{q-r}f'_*(g'^*(\cH^r(F^\bullet)))}\ar[d]\\
{\rR^qf'_*(C'^\bullet_r)}\ar[r]&{\rR^{q-r}f'_*(\cH^r(g'^*(F^\bullet))),}}
\end{equation}
where the horizontal arrows are the canonical morphisms, 
the upper vertical arrows are the base change morphisms \eqref{p1-bcim10d}, the right one being an isomorphism by assumption,
and the right lower vertical arrow is induced by the morphism \eqref{p1-bcim14b}, which is an isomorphism by assumption. 
Since the upper and lower horizontal arrows are isomorphisms, we deduce that the composition of the left vertical arrows is an isomorphism. 

As $\tau_{\leq r}(F^\bullet)=0$ for $r\ll 0$, we deduce from \eqref{p1-bcim14q}, by induction on $r$, that the composed morphism 
\begin{equation}
g^*(\rR^qf_*(\tau_{\leq r}(F^\bullet)))\rightarrow \rR^qf_*(g'^*(\tau_{\leq r}(F^\bullet)))\rightarrow \rR^qf'_*(\tau_{\leq r}(g'^*(F^\bullet))),
\end{equation}
where the first arrow is the base change morphism and the second arrow is induced by ${t^\bullet_{r}}$ \eqref{p1-bcim14n}, is an isomorphism.
On the other hand, we have a commutative diagram 
\begin{equation}
\xymatrix{
{g^*(\rR^qf_*(\tau_{\leq r}(F^\bullet)))}\ar[r]\ar[d]&{\rR^qf_*(g'^*(\tau_{\leq r}(F^\bullet)))}\ar[r]&{\rR^qf'_*(\tau_{\leq r}(g'^*(F^\bullet)))}\ar[d]\\
{g^*(\rR^qf_*(F^\bullet))}\ar[rr]&&{\rR^qf'_*(g'^*(F^\bullet)),}}
\end{equation}
where the vertical arrows are induced by the canonical morphisms and the lower horizontal arrow is the base change morphism  \eqref{p1-bcim10d}. 
The proposition follows since for $r\geq q$, the vertical arrows are isomorphisms.

\subsection{}\label{p1-bcim17}
Let $f_*\colon \cA\rightarrow \cB$ be a left exact functor between abelian $\mU$-categories \eqref{p1-NC0}. We take again the notation 
and conventions of §~\ref{p1-abisoind}. By (\cite{ag2} 2.6.4), there exists essentially unique additive functors
\begin{equation}\label{p1-bcim17a}
\rI f_*\colon \bInd(\cA) \rightarrow \bInd(\cB)
\end{equation}
that commute with filtered direct limits and that fit into commutative diagrams, up to canonical isomorphism,  
\begin{equation}\label{p1-bcim17b}
\xymatrix{
{\cA}\ar[r]^-(0.5){f_*}\ar[d]_{\iota_\cA}&{\cB}\ar[d]^{\iota_{\cB}}\\
{\bInd(\cA)}\ar[r]^-(0.5){\rI f_*}&{\bInd(\cB),}}
\end{equation}
where the vertical arrows are the canonical functors \eqref{p1-abisoind1b}. 
The functor $\rI f_*$ is left exact by (\cite{ks2} 8.6.8).
Suppose that there exists an {\em $f_*$-injective} subcategory $\cJ$ of $\cA$ in the sense of (\cite{ks2} 13.3.4).
By virtue of (\cite{ks2} 13.3.5), the functor $f_*$ admits a right derived functor
\begin{equation}\label{p1-bcim17c}
\rR f_*\colon \bD^+(\cA)\rightarrow \bD^+(\cB).
\end{equation}
It is a triangulated functor by (\cite{ks2} 10.3.3). 
For any integer $q$, we denote by $\rR^qf_*\colon \cA\rightarrow \cB$ the $q$th right derived functor of $f_*$ and by
\begin{equation}\label{p1-bcim17d}
\rI(\rR^if_*)\colon \Ind(\cA)\rightarrow \Ind(\cB)
\end{equation}
the associated functor (\cite{ag2} 2.6.4). By (\cite{ks2} 15.3.2), the category $\Ind(\cJ)$ is $\rI f_*$-injective and the functor $\rI f_*$
admits a right derived functor
\begin{equation}\label{p1-bcim17e}
\rR\rI f_*\colon \bD^+(\Ind(\cA))\rightarrow \bD^+(\Ind(\cB)).
\end{equation}
Moreover, for every integer $q$, we have a canonical isomorphism
\begin{equation}\label{p1-bcim17f}
\rR^q\rI f_*\stackrel{\sim}{\rightarrow}\rI(\rR^q f_*).
\end{equation}
In particular, $\rR^q\rI f_*$ commutes with small filtered direct limits.
Moreover, since the canonical functor $\iota_\cA$ is exact (\cite{ag2} 2.6.3) and $\iota_\cA(\cJ)\subset \bInd(\cJ)$, the diagram
\begin{equation}\label{p1-bcim17g}
\xymatrix{
{\bD^+(\cA)}\ar[r]^-(0.5){\rR f_*}\ar[d]_{\iota_\cA}&{\bD^+(\cB)} \ar[d]^{\iota_{\cB}}\\
{\bD^+(\Ind(\cA))}\ar[r]^-(0.5){\rR\rI f_*}&{\bD^+(\Ind(\cB)),} }
\end{equation}
is commutative up to canonical isomorphism.

\subsection{}\label{p1-bcim170}
Let $f_*\colon \cA\rightarrow \cB$ be a left exact functor between abelian $\mU$-categories \eqref{p1-NC0}. 
We suppose that the following conditions hold true: 
\begin{itemize}
\item[(i)] There exists an {\em $f_*$-injective} subcategory $\cJ$ of $\cA$ in the sense of (\cite{ks2} 13.3.4).
\item[(ii)] Filtered direct limits are representable in $\cA$ and $\cB$, and are exact. 
\end{itemize}
By (i) and \ref{p1-bcim17}, the functors $f_*$ and $\rI f_*$ \eqref{p1-bcim17b} admit right derived functors $\rR f_*$ \eqref{p1-bcim17c} and $\rR(\rI f_*)$ \eqref{p1-bcim17e}. 
By (ii) and (\cite{ag2} 2.6.5), the canonical functors $\iota_\cA$ and $\iota_\cB$ \eqref{p1-abisoind1b} admit left adjoints 
\begin{equation}\label{p1-bcim170a}
\kappa_\cA\colon \bInd(\cA)\rightarrow \cA, \ \ \ \kappa_\cB\colon \bInd(\cB)\rightarrow \cB.
\end{equation}
The latters are exact; the proof is identical to that of (\cite{ag2} 2.7.2) using (ii) and (\cite{ks2} 8.6.6(a)). 
In particular, $\bInd(\cA)$ is $\kappa_\cA$-injective, $\kappa_\cA$ admits a right derived functor
\begin{equation}\label{p1-bcim170b}
\rR \kappa_\cA\colon \bD^+(\bInd(\cA))\rightarrow \bD^+(\cA),
\end{equation}
and similarly for $\kappa_\cB$. 

By (\cite{ks2} 13.3.13(ii)), since the category $\bInd(\cJ)$ is $\rI f_*$-injective \eqref{p1-bcim17}, 
it is $\kappa_\cB \circ \rI f_*$-injective and we have a canonical isomorphism
\begin{equation}\label{p1-bcim170c}
\rR(\kappa_\cB \circ \rI f_*)\stackrel{\sim}{\rightarrow} \rR\kappa_\cB \circ \rR \rI f_*.
\end{equation}

Consider the composed morphism 
\begin{equation}\label{p1-bcim170d}
\rI f_*\rightarrow \rI f_*\circ \iota_\cA \circ \kappa_\cA \rightarrow \iota_\cB  \circ  f_*\circ \kappa_\cA, 
\end{equation}
where the first arrow is induced by the adjunction morphism $\id\rightarrow \iota_\cA \circ \kappa_\cA$ and the second by \eqref{p1-bcim17b}. 
We deduce by adjunction a morphism 
\begin{equation}\label{p1-bcim170e}
\kappa_\cB  \circ  \rI f_*\rightarrow f_*\circ \kappa_\cA. 
\end{equation}

\begin{prop}\label{p1-bcim171}
We keep the assumptions and notation of \ref{p1-bcim170}. 
\begin{itemize}
\item[{\rm (i)}] If $f_*$ commutes with filtered direct limits, then the morphism \eqref{p1-bcim170e} is an isomorphism. 
In particular, the functor $f_*\circ \kappa_\cA$ admits a right derived functor and we have a canonical morphism 
\begin{equation}\label{p1-bcim171a}
\rR(f_*\circ \kappa_\cA)\rightarrow \rR f_*\circ \rR \kappa_\cA.
\end{equation}
\item[{\rm (ii)}] If for every integer $q$, the functor that we denote abusively by
\begin{equation}\label{p1-bcim171b}
\rR^qf_*\colon
\begin{array}[t]{clcr}
\cA&\rightarrow &\cB\\
X&\mapsto&\rR^qf_*(X[0]),
\end{array} 
\end{equation}
commutes with filtered direct limits, 
then the morphism \eqref{p1-bcim171a} is an isomorphism.  
\end{itemize}
\end{prop}

(i) Indeed, for every small filtered category $J$ and every functor $\alpha\colon J\rightarrow \cA$, the composition of the canonical morphisms
\begin{equation}\label{p1-bcim171c}
\kappa_\cB  (\rI f_*(\underset{\underset{J}{\longrightarrow}}{\mlq\mlq\lim \mrq\mrq} \alpha))\stackrel{\sim}{\rightarrow} 
\kappa_\cB  (\underset{\underset{J}{\longrightarrow}}{\mlq\mlq\lim \mrq\mrq} f_*\circ \alpha) 
\stackrel{\sim}{\rightarrow} \underset{\underset{J}{\longrightarrow}}{\lim} \ f_*\circ \alpha
\rightarrow f_*(\kappa_\cA (\underset{\underset{J}{\longrightarrow}}{\mlq\mlq\lim \mrq\mrq} \alpha))
\end{equation}
is none other than the morphism \eqref{p1-bcim170e}. The first (resp.\ second) morphism is an isomorphism by the definition of $\rI f_*$ (resp.\ $\kappa_\cB$), and the last morphism is an isomorphism by assumption, hence the first assertion. The second assertion follows by (\cite{ks2} 13.3.13(i)). 

(ii) Let $\cJ'$ be the full subcategory of right $f_*$-acyclic objects of $\cA$ (\cite{ks2} 13.3.6(ii)). By (\cite{ks2} 13.3.12), we have $\cJ\subset \cJ'$
and $\cJ'$ is $f_*$-injective. On the one hand, since $\bInd(\cJ)$ is $\rI f_*$-injective \eqref{p1-bcim17}, for every object $X$ of $\bInd(\cA)$,
there exists a quasi-isomorphism $X\rightarrow Y$ with $Y\in \bK^+(\bInd(\cJ))$. Since $\kappa_\cA$ is exact \eqref{p1-bcim170}, we deduce that 
$\bInd(\cJ)$ is $\kappa_\cA$-injective. On the other hand, we have $\kappa_\cA(\bInd(\cJ))\subset \cJ'$. 
Indeed, for every small filtered category $J$, every functor $\alpha\colon J\rightarrow \cJ$ and every integer $q\not=0$, we have 
$\rR^q f_*\circ \alpha=0$ (\cite{ks2} 13.3.5(i)). Since the canonical morphism
\begin{equation}\label{p1-bcim171d}
\underset{\underset{J}{\longrightarrow}}\lim \ \rR^q f_*\circ \alpha \rightarrow 
\rR^q f_*(\underset{\underset{J}{\longrightarrow}}\lim \ \alpha)=
\rR^q f_*(\kappa_\cA(\underset{\underset{J}{\longrightarrow}}{\mlq\mlq\lim \mrq\mrq} \alpha))
\end{equation}
is an isomorphism by assumption, we deduce that the target vanishes, which proves the claim. 
The proposition follows then from (\cite{ks2} 13.3.13(ii)). 

\begin{cor}\label{p1-bcim172}
We keep the assumptions and notation of \ref{p1-bcim170}, and assume, 
moreover, that for every integer $q$, the functor $\rR^qf_*$ \eqref{p1-bcim171b} commutes with filtered direct limits,
then we have a canonical functorial isomorphism 
\begin{equation}\label{p1-bcim172a}
\rR \kappa_\cB\circ \rR \rI f_* \stackrel{\sim}{\rightarrow} \rR f_* \circ \rR \kappa_\cA. 
\end{equation}
In particular, for every $X\in \bD^+(\bInd(\cA))$ and every integer $q$, we have a canonical functorial isomorphism 
\begin{equation}\label{p1-bcim172b}
\kappa_\cB(\rR^q \rI f_* (X))\stackrel{\sim}{\rightarrow}\rR^q f_* (\rR \kappa_\cA(X)).
\end{equation}
\end{cor}

Observe first that the functor $f_*$ commutes with filtered direct limits since for every $X\in \ob(\cA)$, 
the canonical morphism $f_*(X)\rightarrow \rR^0f_*(X[0])$ is an isomorphism by \ref{p1-bcim15}. 
By \ref{p1-bcim171}, we have the morphisms (without the dotted arrow)
\begin{equation}
\xymatrix{
{\rR (\kappa_\cB\circ \rI f_*)}\ar[r]\ar[d]&{\rR (f_* \circ \kappa_\cA)}\ar[d]\\
{\rR \kappa_\cB\circ \rR \rI f_*}\ar@{.>}[r]&{\rR f_* \circ \rR \kappa_\cA,}}
\end{equation}
where the upper horizontal arrow is the isomorphism induced by \eqref{p1-bcim170e} 
and the vertical arrows are the isomorphisms \eqref{p1-bcim170c} and \eqref{p1-bcim171a}.  We take for the isomorphism 
\eqref{p1-bcim172a} the unique dotted arrow that makes the above diagram commutative. 

The second assertion follows from the fact that the functors $\kappa_\cA$ and $\kappa_\cB$ are exact. 

\begin{lem}\label{p1-bcim18}
Let $f_*\colon \cA\rightarrow \cB$ be a left exact functor between abelian $\mU$-categories, 
$\cJ$ an $f_*$-injective subcategory of $\cA$, $I$ a small filtered category, $\varphi\colon I\rightarrow \bC^+(\cA)$ a functor \eqref{p1-NC8}. 
For any $i\in \ob(I)$, we set $\cF^\bullet_i=\varphi(i)$. 
We assume that there exists an integer $s$ such that for every $r\leq s$ and every $i\in \ob(I)$, $\cF^r_i=0$. 
\begin{itemize}
\item[{\rm (i)}] For every integer $q$, the canonical morphism 
\begin{equation}\label{p1-bcim18a}
\underset{\underset{i\in I}{\longrightarrow}}{\mlq\mlq\lim \mrq\mrq} \ \rR^q f_* (\cF^\bullet_i)\rightarrow 
\rR^q \rI f_* (\underset{\underset{i\in I}{\longrightarrow}}{\mlq\mlq\lim \mrq\mrq} \ \cF^\bullet_i),
\end{equation}
where the direct limit on the right hand side is computed term by term, is an isomorphism. 
\item[{\rm (ii)}] Assume that filtered direct limits are representable in $\cA$ and $\cB$ and are exact and that, for every integer $q$, 
$\rR^q f_*$ commutes with these limits. Then, for every integer $q$, the canonical morphism 
\begin{equation}\label{p1-bcim18b}
\underset{\underset{i\in I}{\longrightarrow}}{\lim} \ \rR^q f_* (\cF^\bullet_i)\rightarrow 
\rR^q f_* (\underset{\underset{i\in I}{\longrightarrow}}{\lim} \ \cF^\bullet_i),
\end{equation}
where the direct limit on the right hand side is computed term by term, is an isomorphism. 
\item[{\rm (iii)}] Under the assumptions of {\rm (ii)}, for every integer $q$, we have a canonical functorial isomorphism 
\begin{equation}\label{p1-bcim18c}
\kappa_\cB(\rR^q \rI f_* (\underset{\underset{i\in I}{\longrightarrow}}{\mlq\mlq\lim \mrq\mrq} \ \cF^\bullet_i))
\stackrel{\sim}{\rightarrow}
\rR^q f_* (\kappa_\cA(\underset{\underset{i\in I}{\longrightarrow}}{\mlq\mlq\lim \mrq\mrq} \ \cF^\bullet_i)),
\end{equation}
where $\kappa_\cA$ is computed term by term. 
\end{itemize}
\end{lem}

(i) Indeed, by (\cite{ag2} 2.6.7.5), for every integer $r$ and every $i\in \ob(I)$, we have exact sequences 
\begin{eqnarray}
0\rightarrow \tau_{\leq r-1}(\cF^\bullet_i)\rightarrow \tau_{\leq r}(\cF^\bullet_i)\rightarrow \cL_{i,r}^\bullet \rightarrow 0,\\
0\rightarrow \underset{\underset{i\in I}{\longrightarrow}}{\mlq\mlq\lim \mrq\mrq}\ \tau_{\leq r-1}(\cF^\bullet_i)\rightarrow 
\underset{\underset{i\in I}{\longrightarrow}}{\mlq\mlq\lim \mrq\mrq}\ \tau_{\leq r}(\cF^\bullet_i)\rightarrow 
\underset{\underset{i\in I}{\longrightarrow}}{\mlq\mlq\lim \mrq\mrq}\ \cL_{i,r}^\bullet \rightarrow 0,
\end{eqnarray}
and quasi-isomorphisms
\begin{eqnarray}
\cL_{i,r}^\bullet \rightarrow \cH^r(\cF^\bullet_i)[-r],\\
\underset{\underset{i\in I}{\longrightarrow}}{\mlq\mlq\lim \mrq\mrq}\ \cL_{i,r}^\bullet \rightarrow 
\underset{\underset{i\in I}{\longrightarrow}}{\mlq\mlq\lim \mrq\mrq}\ (\cH^r(\cF^\bullet_i))[-r]. 
\end{eqnarray}
The formers induce a canonical morphism of long exact sequences of cohomology 
\begin{equation}
\xymatrix{
\ar[r]&{\underset{\underset{i\in I}{\longrightarrow}}{\mlq\mlq\lim \mrq\mrq}\ \rR^q f_* (\tau_{\leq r-1}(\cF^\bullet_i))}\ar[r]\ar[d]
&{\underset{\underset{i\in I}{\longrightarrow}}{\mlq\mlq\lim \mrq\mrq}\ \rR^q f_* (\tau_{\leq r}(\cF^\bullet_i))}\ar[r]\ar[d]&
{\underset{\underset{i\in I}{\longrightarrow}}{\mlq\mlq\lim \mrq\mrq}\ \rR^q f_* (\cL_{i,r}^\bullet)}\ar[r]\ar[d]&\\
\ar[r]&{\rR^q \rI f_* (\underset{\underset{i\in I}{\longrightarrow}}{\mlq\mlq\lim \mrq\mrq}\ \tau_{\leq r-1}(\cF^\bullet_i))}\ar[r]&
{\rR^q \rI f_* (\underset{\underset{i\in I}{\longrightarrow}}{\mlq\mlq\lim \mrq\mrq}\ \tau_{\leq r}(\cF^\bullet_i))}\ar[r]& 
{\rR^q \rI f_* (\underset{\underset{i\in I}{\longrightarrow}}{\mlq\mlq\lim \mrq\mrq}\ \cL_{i,r}^\bullet)}\ar[r]&}
\end{equation}
The latter induce a commutative diagram
\begin{equation}
\xymatrix{
{\underset{\underset{i\in I}{\longrightarrow}}{\mlq\mlq\lim \mrq\mrq}\ \rR^q f_* (\cL_{i,r}^\bullet)}\ar[r]\ar[d]&
{\underset{\underset{i\in I}{\longrightarrow}}{\mlq\mlq\lim \mrq\mrq}\ \rR^{q-r} f_* (\cH^r(\cF^\bullet_i))}\ar[d]\\
{\rR^q \rI f_* (\underset{\underset{i\in I}{\longrightarrow}}{\mlq\mlq\lim \mrq\mrq}\ \cL_{i,r}^\bullet)}\ar[r]&
{\rR^{q-r} \rI f_* (\underset{\underset{i\in I}{\longrightarrow}}{\mlq\mlq\lim \mrq\mrq}\ (\cH^r(\cF^\bullet_i))),}} 
\end{equation}
where the horizontal arrows are isomorphisms, and the vertical arrows are the canonical morphisms, the right one being an isomorphism by 
\eqref{p1-bcim17f}. Since $\tau_{\leq t}(\cF^\bullet_i)=0$ for $t\ll 0$ and every $i\in \ob(I)$, we deduce by induction that for every integers $r$ and $q$, 
the canonical morphism 
\begin{equation}
\underset{\underset{i\in I}{\longrightarrow}}{\mlq\mlq\lim \mrq\mrq} \rR^q f_* (\tau_{\leq r}(\cF^\bullet_i))\rightarrow 
\rR^q \rI f_* (\underset{\underset{i\in I}{\longrightarrow}}{\mlq\mlq\lim \mrq\mrq} \ \tau_{\leq r}(\cF^\bullet_i))
\end{equation}
is an isomorphism. 

For all integers $t\leq r$ and every $i\in \ob(I)$, since $\cH^t(\cF^\bullet_i/\tau_{\leq r}(\cF^\bullet_i))=0$, we have
\begin{equation}
\cH^t((\underset{\underset{i\in I}{\longrightarrow}}{\mlq\mlq\lim \mrq\mrq} \ \cF^\bullet_i)/
(\underset{\underset{i\in I}{\longrightarrow}}{\mlq\mlq\lim \mrq\mrq} \ \tau_{\leq r}(\cF^\bullet_i))=0. 
\end{equation}
Hence, for all integers $q\leq r$, in the commutative canonical diagram 
\begin{equation}
\xymatrix{
{\underset{\underset{i\in I}{\longrightarrow}}{\mlq\mlq\lim \mrq\mrq}\ \rR^q f_* (\tau_{\leq r}(\cF^\bullet_i))}\ar[r]\ar[d]&
{\underset{\underset{i\in I}{\longrightarrow}}{\mlq\mlq\lim \mrq\mrq}\ \rR^q f_* (\cF^\bullet_i)}\ar[d]\\
{\rR^q \rI f_* (\underset{\underset{i\in I}{\longrightarrow}}{\mlq\mlq\lim \mrq\mrq}\ \tau_{\leq r}(\cF^\bullet_i))}\ar[r]&
{\rR^q \rI f_* (\underset{\underset{i\in I}{\longrightarrow}}{\mlq\mlq\lim \mrq\mrq}\ \cF^\bullet_i),}} 
\end{equation}
the horizontal arrows are isomorphisms \eqref{p1-bcim15}. Since the left vertical arrow is an isomorphism, so is the right vertical arrow.

(ii) It is enough to copy the proof of (i). 

(iii) Indeed, \eqref{p1-bcim18a} induces an isomorphism
\begin{equation}
\kappa_\cB(\underset{\underset{i\in I}{\longrightarrow}}{\mlq\mlq\lim \mrq\mrq} \ \rR^q f_* (\cF^\bullet_i))\stackrel{\sim}{\rightarrow}
\kappa_\cB(\rR^q \rI f_* (\underset{\underset{i\in I}{\longrightarrow}}{\mlq\mlq\lim \mrq\mrq} \ \cF^\bullet_i)).
\end{equation}
On the other hand, we have canonical isomorphisms
\begin{equation}
\kappa_\cB(\underset{\underset{i\in I}{\longrightarrow}}{\mlq\mlq\lim \mrq\mrq} \ \rR^q f_* (\cF^\bullet_i))\stackrel{\sim}{\rightarrow}
\underset{\underset{i\in I}{\longrightarrow}}{\lim} \ \rR^q f_* (\cF^\bullet_i)\stackrel{\sim}{\rightarrow}
\rR^q f_* (\underset{\underset{i\in I}{\longrightarrow}}{\lim} \ \cF^\bullet_i)\stackrel{\sim}{\rightarrow}
\rR^q f_* (\kappa_\cA(\underset{\underset{i\in I}{\longrightarrow}}{\mlq\mlq\lim \mlq\mlq} \ \cF^\bullet_i)),
\end{equation}
the second one being \eqref{p1-bcim18b}. 

\begin{rema}
We can prove that the isomorphism \eqref{p1-bcim18c} is induced by \eqref{p1-bcim172b}.
\end{rema}

\subsection{}\label{p1-bcim5}
Let $(X,A)$ be a ringed $\mU$-topos.
We denote by $\bMod_{\mQ}(A)$ the category of $A$-modules of $X$ up to isogeny and by
\begin{equation}\label{p1-bcim5a}
Q_A\colon \bMod(A)\rightarrow \bMod_\mQ(A),\ \ \ M\mapsto M_\mQ,
\end{equation}
the canonical functor (see \ref{p1-abisoind1}). Objects of $\bMod_\mQ(A)$ are also called {\em $A_\mQ$-modules}.
The tensor product of $\bMod(A)$ induces a right exact bifunctor
\begin{equation}\label{p1-bcim5b}
\begin{array}[t]{clcr}
\bMod_{\mQ}(A)\times \bMod_{\mQ}(A)&\rightarrow &\bMod_{\mQ}(A),\\ 
(M,N)&\mapsto &M\otimes_{A_\mQ} N,
\end{array}
\end{equation}
making $\bMod_{\mQ}(A)$ into a symmetric monoidal category, having $A_\mQ$ as unit object.

We denote by $\bIndMod(A)$ the category of ind-$A$-modules (see \ref{p1-indmal1}). 
We have a canonical fully faithful and exact functor \eqref{p1-abisoind1d}
\begin{equation}\label{p1-bcim5c}
\upalpha_{A}\colon \bMod_\mQ(A)\rightarrow \bIndMod(A).
\end{equation}
For all $A$-modules $M$ and $N$, we have canonical isomorphisms
\begin{equation}\label{p1-bcim5d}
\upalpha_{A}(M_\mQ)\otimes_A N\stackrel{\sim}{\rightarrow} \upalpha_{A}(M_\mQ)\otimes_A\upalpha_{A}(N_\mQ)
\stackrel{\sim}{\rightarrow} \upalpha_{A}(M_\mQ\otimes_{A_\mQ}N_\mQ),
\end{equation}
where the tensor product of ind-$A$-modules is defined in \eqref{p1-indmal1d}. 

\subsection{}\label{p1-bcim6}
Let $f\colon (X,A)\rightarrow (Y,B)$ be a morphism of ringed $\mU$-topos.
The adjoint functors $f^*$ and $f_*$ between the categories $\bMod(A)$ and $\bMod(B)$ induce adjoint additive functors
\begin{eqnarray}
f^*_\mQ\colon \bMod_\mQ(B) \rightarrow \bMod_{\mQ}(A),\label{p1-bcim6a}\\
f_{\mQ*}\colon \bMod_{\mQ}(A) \rightarrow \bMod_\mQ(B),\label{p1-bcim6b}
\end{eqnarray}
whose former (resp.\ latter) is right (resp.\ left) exact. By (\cite{ag2} 2.9.3(iii)), the canonical functor 
$Q_A$ \eqref{p1-bcim5a} transforms injective $A$-modules into injective $A_\mQ$-modules.
We denote by $\cI_A$ the full subcategory of $\bMod(A)$ made up of injective $A$-modules and by 
$\cI_{\cA,\mQ}$ the full subcategory of $\bMod_{\mQ}(A)$
made up of the images by the functor $Q_A$ of the injective $A$-modules.
Then, $\cI_A$ (resp.\ $\cI_{A,\mQ}$) is cogenerating in $\bMod(A)$ (resp.\ $\bMod_{\mQ}(A)$), 
in particular $\bMod_{\mQ}(A)$ has enough injectives.
The functor $f_{\mQ*}$ admits therefore a right derived functor
\begin{equation}\label{p1-bcim6c}
\rR f_{\mQ*}\colon \bD^+(\bMod_\mQ(A))\rightarrow \bD^+(\bMod_\mQ(B)).
\end{equation}
The diagram
\begin{equation}\label{p1-bcim6d}
\xymatrix{
{\bD^+(\bMod(A))}\ar[r]^{\rR f_*}\ar[d]_{Q_A}&{\bD^+(\bMod(B))}\ar [d]^{Q_B}\\
{\bD^+(\bMod_\mQ(A))}\ar[r]^{\rR f_{\mQ*}}&{\bD^+(\bMod_\mQ(B))}}
\end{equation}
is commutative up to canonical isomorphism.

\subsection{}\label{p1-bcim1}
Let $f\colon (X,A)\rightarrow (Y,B)$ be a morphism of ringed $\mU$-topos. 
By (\cite{ag2} 2.6.4), there exist essentially unique additive functors
\begin{eqnarray}
\rI f^*\colon \bIndMod(B) \rightarrow \bIndMod(A),\label{p1-bcim1g}\\
\rI f_*\colon \bIndMod(A) \rightarrow \bIndMod(B),\label{p1-bcim1h}
\end{eqnarray}
that commute with filtered direct limits and that fit into commutative diagrams, up to canonical isomorphism,  
\begin{equation}\label{p1-bcim1a}
\xymatrix{
{\bMod(B)}\ar[r]^-(0.5){f^*}\ar[d]_{\iota_B}&{\bMod(A)}\ar[d]^{\iota_{A}}\\
{\bIndMod(B)}\ar[r]^-(0.5){\rI f^*}&{\bIndMod(A),}} 
\ \ \
\xymatrix{
{\bMod(A)}\ar[r]^-(0.5){f_*}\ar[d]_{\iota_A}&{\bMod(B)}\ar[d]^{\iota_{B}}\\
{\bIndMod(A)}\ar[r]^-(0.5){\rI f_*}&{\bIndMod(B),}}
\end{equation}
where the vertical arrows are the canonical functors \eqref{p1-indmal1a}. 
The functor $\rI f^*$ is a left adjoint of the functor $\rI f_*$.
By (\cite{ks2} 8.6.8), the functor $\rI f^*$ (resp.\ $\rI f_*$) is right (resp.\ left) exact. 

By (\cite{ag2} 2.6.10), the functor $\rI f_*$ admits a right derived functor
\begin{equation}\label{p1-bcim1b}
\rR\rI f_*\colon \bD^+(\bIndMod(A))\rightarrow \bD^+(\bIndMod(B)).
\end{equation}
By \eqref{p1-bcim17g}, the diagram 
\begin{equation}\label{p1-bcim1d}
\xymatrix{
{\bD^+(\bMod(A))}\ar[r]^-(0.5){\rR f_*}\ar[d]_{\iota_A}&{\bD^+(\bMod(B))}\ar[d]^{\iota_{B}}\\
{\bD^+(\bIndMod(A))}\ar[r]^-(0.5){\rR\rI f_*}&{\bD^+(\bIndMod(B))}}
\end{equation}
is commutative up to canonical isomorphism. 

By (\cite{ag2} 2.9.5(ii)), the diagram
\begin{equation}\label{p1-bcim1e}
\xymatrix{
{\bMod_\mQ(B)}\ar[r]^{f^*_{\mQ}}\ar[d]_{\upalpha_{B}}&{\bMod_\mQ(A)}\ar [d]^{\upalpha_{A}}\\
{\bIndMod(B)}\ar[r]^{\rI f^*}&{\bIndMod(A),}}
\end{equation}
where the vertical arrows are the canonical functors \eqref{p1-bcim5c}, 
is commutative up to canonical isomorphism.
By (\cite{ag2} (2.9.6.7)), the diagram
\begin{equation}\label{p1-bcim1f}
\xymatrix{
{\bD^+(\bMod_\mQ(A))}\ar[r]^-(0.5){\rR f_{\mQ*}}\ar[d]_{\upalpha_{A}} &{\bD^+(\bMod_\mQ(B))}\ar[d]^{\upalpha_{B}}\\
{\bD^+(\bIndMod(A))}\ar[r]^-(0.5){\rR\rI f_*}&{\bD^+(\bIndMod(B))}}
\end{equation}
is commutative up to canonical isomorphism. 

\begin{lem}\label{p1-bcim7}
Let $f\colon (X,A)\rightarrow (Y,B)$, $g\colon (Y,B)\rightarrow (Z,C)$ be two morphisms of ringed $\mU$-topos, $h=g\circ f\colon (X,A)\rightarrow (Z,C)$.
Then, 
\begin{itemize}
\item[{\rm (i)}] There is a canonical isomorphism
\begin{equation}\label{p1-bcim7a}
\rR h_*\stackrel{\sim}{\rightarrow} (\rR g_*)\circ (\rR f_*).
\end{equation}
\item[{\rm (ii)}] There is a canonical isomorphism
\begin{equation}\label{p1-bcim7b}
\rR h_{\mQ *}\stackrel{\sim}{\rightarrow} (\rR g_{\mQ *})\circ (\rR f_{\mQ *}).
\end{equation}
Moreover, the diagram 
\begin{equation}\label{p1-bcim7c}
\xymatrix{
{Q_C\circ (\rR h_*)}\ar[r]\ar[d]&{Q_C\circ(\rR g_*)\circ (\rR f_*)}\ar[d]\\
{(\rR h_{\mQ *})\circ Q_A}\ar[r]&{(\rR \rI g_{\mQ *})\circ (\rR f_{\mQ *})\circ Q_A,}}
\end{equation}
where the vertical arrows are induced by the isomorphism underlying \eqref{p1-bcim6d}, is commutative up to canonical isomorphism. 
\item[{\rm (iii)}] There is a canonical isomorphism
\begin{equation}\label{p1-bcim7d}
\rR \rI h_*\stackrel{\sim}{\rightarrow} (\rR \rI g_*)\circ (\rR\rI f_*).
\end{equation}
Moreover, the diagrams 
\begin{equation}\label{p1-bcim7e}
\xymatrix{
{\iota_C\circ (\rR h_*)}\ar[r]\ar[d]&{\iota_C\circ(\rR g_*)\circ (\rR f_*)}\ar[d]\\
{(\rR \rI h_*)\circ \iota_A}\ar[r]&{(\rR \rI g_*)\circ (\rR\rI f_*)\circ \iota_A,}}
\end{equation}
\begin{equation}\label{p1-bcim7f}
\xymatrix{
{\upalpha_C\circ (\rR h_{\mQ *})}\ar[r]\ar[d]&{\upalpha_C\circ(\rR g_{\mQ *})\circ (\rR f_{\mQ *})}\ar[d]\\
{(\rR \rI h_*)\circ \upalpha_A}\ar[r]&{(\rR \rI g_*)\circ (\rR\rI f_*)\circ \upalpha_A,}}
\end{equation}
where  the vertical arrows are induced by the isomorphism underlying \eqref{p1-bcim1d} and \eqref{p1-bcim1f}, 
are commutative up to canonical isomorphism.
\end{itemize}
\end{lem}  
 
(i) We denote by $\cF_A$ (resp.\ $\cF_B$) the full subcategory of $\bMod(A)$ (resp.\ $\bMod(B)$) made up of flabby modules (\cite{sga4} V 4.1).
Injective $A$-modules being flabby, the category $\cF_A$ is cogenerating in $\bMod(A)$.
By (\cite{sga4} V 5.2), the objects of $\cF_A$ are acyclic for $f_*$.
Hence, the category $\cF_A$ is $f_*$-injective by (\cite{ks2} 13.3.8).
Since $f_*(\cF_A)\subset \cF_B$ by (\cite{sga4} V 4.9), the proposition follows by (\cite{ks2} 13.3.13). 

(ii) We denote by $\cF_{A,\mQ}$ (resp.\ $\cF_{B,\mQ}$) the full subcategory of $\bMod_\mQ(A)$ (resp.\ $\bMod_\mQ(B)$) 
made up of the images by $Q_A$ (resp.\ $Q_B$) of flabby modules.
The category $\cF_{A,\mQ}$ is cogenerating in $\bMod_\mQ(A)$ and the objects of $\cF_{A,\mQ}$ are acyclic for $f_{\mQ*}$.
Hence, the category $\cF_{A,\mQ}$ is $f_{\mQ*}$-injective by (\cite{ks2} 13.3.8).
Moreover, we have $f_{\mQ*}(\cF_{A,\mQ})\subset \cF_{B,\mQ}$ by (\cite{sga4} V 4.9). 
The proposition follows then by (\cite{ks2} 13.3.13). 

(iii) The category $\Ind(\cF_A)$ is $\rI f_*$-injective by (\cite{ks2} 15.3.2). 
We have $\upalpha_A(\cF_{A,\mQ})\subset \bInd(\cF_A)$ and $\rI f_*(\Ind(\cF_A))\subset \Ind(\cF_B)$.
The proposition follows then by (\cite{ks2} 13.3.13). 

\begin{lem}\label{p1-bcim22}
Let $f\colon (X,A)\rightarrow (Y,B)$ be a morphism of ringed $\mU$-topos,
$E$ a locally free $B$-module of finite type. 
Then, 
\begin{itemize}
\item[{\rm (i)}] For every complex $F$ of $\bD^+(\bMod(A))$, 
there is a canonical functorial isomorphism of $\bD^+(\bMod(B))$ 
\begin{equation}\label{p1-bcim22a}
E\otimes_B \rR f_* (F)\stackrel{\sim}{\rightarrow} \rR f_* (f^*(E)\otimes_A F).
\end{equation}
\item[{\rm (ii)}] For every complex $F$ of $\bD^+(\bMod_\mQ(A))$, 
there is a canonical functorial isomorphism of $\bD^+(\bMod_\mQ(B))$ 
\begin{equation}\label{p1-bcim22b}
E\otimes_B \rR f_{\mQ*} (F)\stackrel{\sim}{\rightarrow} \rR f_{\mQ*} (f^*(E)\otimes_A F).
\end{equation}
\item[{\rm (iii)}] For every complex $F$ of $\bD^+(\bIndMod(A))$, 
there is a canonical functorial isomorphism of $\bD^+(\bIndMod(B))$ 
\begin{equation}\label{p1-bcim22c}
E\otimes_B \rR\rI f_* (F)\stackrel{\sim}{\rightarrow} \rR\rI f_* (f^*(E)\otimes_A F).
\end{equation}
\end{itemize}
\end{lem}

We prove only (iii); the proofs of the two other propositions are similar.
Recall (\cite{ks2} 13.3.1) that for every complex $C$ of $\bK^+(\bIndMod(A))$, there is a canonical functorial isomorphism 
\begin{equation}\label{p1-bcim22d}
\rR \rI f_*(C) \stackrel{\sim}{\rightarrow} \underset{\underset{C\rightarrow C'}{\longrightarrow}}{\mlq\mlq\lim \mrq\mrq} \ \rI f_*(C'),
\end{equation}
where the direct limit is taken in the category of ind-objects of $\bD^+(\bIndMod(B))$, over the category of quasi-isomorphisms of $\bK^+(\bIndMod(A))$. 

For every ind-$A$-module $M$, the adjunction morphism $\rI f^*(\rI f_*(M))\rightarrow M$ induces, by adjunction, a morphism of ind-$B$-modules
\begin{equation}\label{p1-bcim22e}
E\otimes_B \rI f_* (M)\rightarrow \rI f_* (f^*(E)\otimes_A M).
\end{equation}
It is an isomorphism by localization and (\cite{ag2} 2.7.17(iii)). We deduce, for every complex $C$ of $\bK^+(\bIndMod(A))$, an isomorphism of complexes of 
ind-$B$-modules 
\begin{equation}\label{p1-bcim22f}
E\otimes_B \rI f_* (C)\stackrel{\sim}{\rightarrow} \rI f_* (f^*(E)\otimes_A C).
\end{equation}
We obtain from this and \eqref{p1-bcim22d}, for every complex $C$ of $\bK^+(\bIndMod(A))$, a canonical functorial morphism of $\bD^+(\bIndMod(B))$ 
\begin{equation}\label{p1-bcim22g}
E\otimes_B \rR\rI f_* (C)\rightarrow \rR\rI f_* (f^*(E)\otimes_A C).
\end{equation}
It is obviously compatible with localization. 
In order to prove that it is an isomorphism, we take cohomology. Then, by (\cite{ag2} 2.7.17(iii)), we may reduce to the case where $E$ is a free $B$-module 
of finite type, in which case the assertion is obvious. 

\subsection{}\label{p1-bcim4}
We consider a diagram of morphisms of ringed $\mU$-topos
\begin{equation}\label{p1-bcim4a}
\xymatrix{
{(X',A')}\ar[r]^-(0.5){g'}\ar[d]_{f'}&{(X,A)}\ar[d]^f\\
{(Y',B')}\ar[r]^-(0.5)g&{(Y,B)}}
\end{equation}
commutative up to canonical isomorphism; in other words, we have an isomorphism
\begin{equation}\label{p1-bcim4b}
f_*g'_*\stackrel{\sim}{\rightarrow}g_*f'_*
\end{equation}
and the diagram
\begin{equation}\label{p1-bcim4c}
\xymatrix{
{f_*(A)}\ar[d]_{f_*(g'^\#)}&B\ar[r]^-(0.5){g^\#}\ar[l]_-(0.5 ){f^\#}&{g_*(B')}\ar[d]^{g_*(f'^\#)}\\
{f_*(g'_*(A'))}\ar[rr]&&{g_*(f'_*(A')),}}
\end{equation}
where the unlabeled arrow is the isomorphism \eqref{p1-bcim4b}, is commutative.

Let $q$ be an integer. By \ref{p1-bcim10} and \ref{p1-bcim7}, for every bounded from below complex of $A$-modules $F^\bullet$, 
there exists a canonical functorial morphism of $B'$-modules \eqref{p1-bcim10d}
\begin{equation}\label{p1-bcim4d}
g^*(\rR^q f_*(F^\bullet))\rightarrow \rR^q f'_*(g'^*(F^\bullet)),
\end{equation}
where the pullback $g'^*(F^\bullet)$ is defined term by term (not derived), called {\em base change morphism}. 

Similarly, for every bounded from below complex of $A_\mQ$-modules $F^\bullet$, 
there exists a canonical functorial morphism of $B'_\mQ$-modules
\begin{equation}\label{p1-bcim4e}
g^*_\mQ(\rR^q f_{\mQ *}(F^\bullet))\rightarrow \rR^q f'_{\mQ *}(g'^*_\mQ(F^\bullet)),
\end{equation}
where the pullback $g'^*_\mQ(F^\bullet)$ is defined term by term (not derived), 
called {\em base change morphism}.

Similarly, for every bounded from below complex of ind-$A$-modules $F^\bullet$, 
there exists a canonical functorial morphism of ind-$B'$-modules
\begin{equation}\label{p1-bcim4f}
\rI g^*(\rR^q\rI f_*(F^\bullet))\rightarrow \rR^q\rI f'_*(\rI g'^*(F^\bullet)),
\end{equation}
where the pullback $\rI g'^*(F^\bullet)$ is defined term by term (not derived), called {\em base change morphism}. 

\begin{rema}\label{p1-bcim12}
We keep the assumptions and notation of \ref{p1-bcim4}.
For every $A$-module $F$ and every integer $q\geq 0$, it follows from \ref{p1-bcim11} that the base change morphism \eqref{p1-bcim4d} 
for $F[0]$ coincides with the classical base change morphism (\cite{egr1} 1.2.3).
\end{rema}

\subsection{}\label{p1-bcim9}
We keep the assumptions and notation of \ref{p1-bcim4} and let $q$ be an integer. 
For every bounded from below complex of $A$-modules $F^\bullet$, 
the base change morphism \eqref{p1-bcim4e}
\begin{equation}\label{p1-bcim9a}
g^*_\mQ(\rR^q f_{\mQ*}(F^\bullet_\mQ))\rightarrow \rR^q f'_{\mQ*}(g'^*_\mQ(F^\bullet_\mQ))
\end{equation}
is induced by \eqref{p1-bcim4d}. More precisely, the diagram 
\begin{equation}\label{p1-bcim9b}
\xymatrix{
{Q_{B'}(g^*(\rR^q f_*(F^\bullet)))}\ar[r]\ar[d]&{Q_{B'}(\rR^q f'_*(g'^*(F^\bullet)))}\ar[d]\\
{g^*_\mQ(\rR^q f_{\mQ*}(F^\bullet_\mQ))}\ar[r]&{\rR^q f'_{\mQ*}(g'^*_\mQ(F^\bullet_\mQ)),}}
\end{equation}
where the vertical arrows are the isomorphisms induced by \eqref{p1-bcim6d}, 
is commutative up to canonical isomorphism. 

The base change morphism \eqref{p1-bcim4f}
\begin{equation}\label{p1-bcim9c}
\rI g^*(\rR^q\rI f_*(\iota_A(F^\bullet)))\rightarrow \rR^q\rI f'_*(\rI g'^*(\iota_A(F^\bullet)))
\end{equation}
is induced by \eqref{p1-bcim4d}. More precisely, the diagram 
\begin{equation}\label{p1-bcim9d}
\xymatrix{
{\iota_{B'}(g^*(\rR^q f_*(F^\bullet)))}\ar[r]\ar[d]&{\iota_{B'}(\rR^q f'_*(g'^*(F^\bullet)))}\ar[d]\\
{\rI g^*(\rR^q\rI f_*(\iota_A(F^\bullet)))}\ar[r]&{\rR^q\rI f'_*(\rI g'^*(\iota_A(F^\bullet))),}}
\end{equation}
where the vertical arrows are induced by the isomorphisms underlying \eqref{p1-bcim1a} and \eqref{p1-bcim1d}, is commutative up to canonical isomorphism.  

For every bounded from below complex of $A_\mQ$-modules $F^\bullet$, the base change morphism \eqref{p1-bcim4f}
\begin{equation}\label{p1-bcim9e}
\rI g^*(\rR^q\rI f_*(\upalpha_A(F^\bullet)))\rightarrow \rR^q\rI f'_*(\rI g'^*(\upalpha_A(F^\bullet)))
\end{equation}
is induced by \eqref{p1-bcim4e}. More precisely, the diagram 
\begin{equation}\label{p1-bcim9f}
\xymatrix{
{\upalpha_{B'}(g^*_\mQ(\rR^q f_{\mQ *}(F^\bullet))}\ar[r]\ar[d]&{\upalpha_{B'}(\rR^q f'_{\mQ *}(g'^*_\mQ(F^\bullet))}\ar[d]\\
{\rI g^*(\rR^q\rI f_*(\upalpha_A(F^\bullet)))}\ar[r]&{\rR^q\rI f'_*(\rI g'^*(\upalpha_A(F^\bullet))),}}
\end{equation}
where the vertical arrows are the isomorphisms induced by \eqref{p1-bcim1e} and \eqref{p1-bcim1f}, 
is commutative up to canonical isomorphism.

\begin{lem}\label{p1-bcim20}
Let $\cA$ be an abelian category, 
\begin{equation}\label{p1-bcim20a}
f^{\bullet,\bullet},g^{\bullet,\bullet}\colon C^{\bullet,\bullet}\rightrightarrows D^{\bullet,\bullet}
\end{equation}
two morphisms of double complexes of $\cA$, bounded from below in both variables, 
$K^\bullet=\Tot(C^{\bullet,\bullet})$ (resp.\ $L^\bullet=\Tot(D^{\bullet,\bullet})$) the total complex associated with $C^{\bullet,\bullet}$
(resp.\ $D^{\bullet,\bullet}$) (with the sign conventions of {\rm \cite{sp} \href{https://stacks.math.columbia.edu/tag/0FNB}{0FNB}}),  
\begin{equation}\label{p1-bcim20b}
\varphi^{\bullet},\psi^{\bullet}\colon K^{\bullet}\rightrightarrows L^{\bullet}
\end{equation}
the morphisms associated with $f^{\bullet,\bullet}$ and $g^{\bullet,\bullet}$, respectively. 
For any integer $j$, let 
\begin{equation}\label{p1-bcim20c}
k^{i,j}\colon C^{i,j}\rightarrow D^{i-1,j}, \ \ \ i\in \mZ,
\end{equation}
be a homotopy between $f^{\bullet,j}=(f^{i,j})_{i\in \mZ}$ and $g^{\bullet,j}=(g^{i,j})_{i\in \mZ}$, i.e., for every integer $i$, we have 
\begin{equation}\label{p1-bcim20d}
f^{i,j}-g^{i,j}=d_{D,1}^{i-1,j}\circ k^{i,j} + k^{i+1,j} \circ d_{C,1}^{i,j},
\end{equation}
where $d_{C,1}^{i,j}$ and $d_{D,1}^{i,j}$ denote the horizontal differentials of $C^{\bullet,\bullet}$ and $D^{\bullet,\bullet}$, respectively. 
We assume, furthermore, that for all integers $i,j$, the diagram 
\begin{equation}\label{p1-bcim20e}
\xymatrix{
{C^{i,j+1}}\ar[r]^-(0.5){k^{i,j+1}}&{D^{i-1,j+1}}\\
{C^{i,j}}\ar[r]^-(0.5){k^{i,j}}\ar[u]^{d_{C,2}^{i,j}}&{D^{i-1,j},}\ar[u]_{d_{D,2}^{i-1,j}}}
\end{equation}
where $d_{C,2}^{i,j}$ and $d_{D,2}^{i,j}$ denote the vertical differentials of $C^{\bullet,\bullet}$ and $D^{\bullet,\bullet}$, respectively, 
is commutative.
For any integer $q$, we define a morphism 
\begin{equation}\label{p1-bcim20f}
\lambda^q\colon K^q\rightarrow L^{q-1}
\end{equation}
by setting $\lambda^q|C^{i,j}=k^{i,j}$ for any $(i,j)\in \mZ^2$ such that $i+j=q$. Then, $(\lambda^q)_{q\in \mZ}$ gives a homotopy between 
$\varphi^{\bullet}$ and $\psi^{\bullet}$, i.e., 
\begin{equation}\label{p1-bcim20g}
\varphi^q-\psi^q=d_L^{q-1}\circ\lambda^q+\lambda^{q+1}\circ d_K^q. 
\end{equation}
\end{lem}

Indeed, for all integers $i,j,q$ such that $q=i+j$, we have 
\begin{eqnarray*}
d_L^{q-1}\circ\lambda^q|C^{i,j}&=&d_L^{q-1}\circ k^{i,j}=d^{i-1,j}_{D,1}\circ k^{i,j} + (-1)^{i-1}d^{i-1,j}_{D,2}\circ k^{i,j},\\
\lambda^{q+1}\circ d_K^{q}|C^{i,j}&=&\lambda^{q+1}\circ  d^{i,j}_{C,1} + (-1)^i\lambda^{q+1}\circ  d^{i,j}_{C,2}\\
&=&k^{i+1,j}\circ  d^{i,j}_{C,1} + (-1)^i k^{i,j+1}\circ  d^{i,j}_{C,2}.
\end{eqnarray*}

\begin{lem}\label{p1-bcim21}
Let $\cA$ be an abelian category, $C^{\bullet,\bullet}$ a double complex of $\cA$, bounded from below in both variables, such that $C^{i,j}=0$ 
for all $i<0$, $d_{C,1}^{\bullet,\bullet}$  (resp.\ $d_{C,2}^{\bullet,\bullet}$) its horizontal (resp.\ vertical) differential, 
\begin{equation}
\varepsilon^\bullet \colon (M^\bullet,d_M^\bullet)\rightarrow (C^{0,\bullet},d_{C,2}^{0,\bullet})
\end{equation}
a morphism of complexes of $\cA$ whose composition with the morphism 
\begin{equation}
d_{C,1}^{0,\bullet}\colon (C^{0,\bullet},d_{C,2}^{0,\bullet})\rightarrow (C^{1,\bullet},d_{C,2}^{1,\bullet})
\end{equation}
vanishes. We set $K^\bullet=\Tot(C^{\bullet,\bullet})$ the total complex associated with $C^{\bullet,\bullet}$ and denote by 
$\rho^\bullet \colon M^\bullet \rightarrow K^\bullet$ the morphism induced by $\varepsilon^\bullet$. Then, the mapping cone of $\rho^\bullet$
is canonically isomorphic to the total complex associated with the double complex $(D^{\bullet,\bullet},d_{D,1}^{\bullet,\bullet},d_{D,2}^{\bullet,\bullet})$
defined by 
\begin{equation}
\begin{array}{clcr}
D^{i,j}=C^{i,j},&d_{D,1}^{i,j}=d_{C,1}^{i,j},&d_{D,2}^{i,j}=d_{C,2}^{i,j},& \forall i\in \mZ_{\geq 0},  \forall j\in \mZ,\\
D^{-1,j}=M^i,&d_{D,1}^{-1,j}=\varepsilon^j,&d_{D,2}^{-1,j}=d^j_M,& \forall j\in \mZ,\\
D^{i,j}=0,&&&\forall i\in \mZ_{\leq -2},  \forall j\in \mZ.
\end{array}
\end{equation}
\end{lem}

Indeed, if $(L^\bullet,d_L^\bullet)$ denotes the mapping cone of $\rho^\bullet$, we have for every integer $q$, $L^q=M^{q+1}\oplus K^q$ and 
\begin{eqnarray*}
d_L^q|_{M^{q+1}}&=&-d_M^{q+1}+\varepsilon^{q+1},\\
d_L^q|_{K^q}&=&d_K^q. 
\end{eqnarray*}

\section{Affine bundles}

\subsection{}\label{p1-prem1}
Let $X$ be a scheme, $\rT$ a locally free $\co_X$-module of finite type, $\Omega=\cHom_{\co_X}(\rT,\co_X)$ its dual, $\cL$ a $\rT$-torsor of $X_\zar$. 
We denote by $\cF$ the {\em sheaf of affine functions on $\cL$} (\cite{agt} II.4.9). It is an $\co_X$-module that naturally fits into a canonical exact sequence 
\begin{equation}\label{p1-prem1a}
0\longrightarrow \co_X\stackrel{c}{\longrightarrow}\cF\stackrel{\nu}{\longrightarrow} \Omega\longrightarrow 0.
\end{equation}
The latter induces for every integer $n\geq 1$, an exact sequence \eqref{p1-NC7}
\begin{equation}\label{p1-prem1b}
0\rightarrow \rS^{n-1}_{\co_X}(\cF)\rightarrow \rS^{n}_{\co_X}(\cF)\rightarrow \rS^n_{\co_X}(\Omega)\rightarrow 0.
\end{equation}
The $\co_X$-modules $(\rS^{n}_{\co_X}(\cF))_{n\in \mN}$ form thus a filtered direct system, whose direct limit
\begin{equation}\label{p1-prem1c}
\cC=\underset{\underset{n\geq 0}{\longrightarrow}}\lim\ \rS^n_{\co_X}(\cF)
\end{equation}
is naturally equipped with an $\co_X$-algebra structure.  We call it the {\em Higgs--Tate algebra associated with the extension $\cF$} \eqref{p1-prem1a}. 
We set $\bT=\Spec(\rS_{\co_X}(\Omega))$ the $X$-vector bundle associated with $\Omega$. 
Then, $\bL=\Spec(\cC)$ is canonically a principal homogeneous $\bT$-bundle over $X$ (\cite{agt} II.4.10). 

For every $s\in \cL(X)$, the morphism $\rho_s\colon \cF\rightarrow \co_X$ that associates with any local section $f$ of $\cF$, $f(s)$, 
is a splitting of the exact sequence \eqref{p1-prem1a}. The latter extends to a unique homomorphism of $\co_X$-algebras
$\varrho_s\colon\cC\rightarrow \co_X$, which induces a section $\pi_s\in \bL(X)$. 
By (\cite{agt} II.4.10), the correspondence $s\mapsto \pi_s$ defines an isomorphism of $\rT$-torsors of $X_\zar$ 
\begin{equation}\label{p1-prem1d}
\cL\stackrel{\sim}{\rightarrow} \cHom_X(-,\bL).
\end{equation}
The endomorphism $\id_\cF-c\circ \rho_s$ of $\cF$ induces an $\co_X$-linear morphism $\sigma_s\colon \Omega\rightarrow \cF$, which is a section of $\nu$, 
and hence an isomorphism of $\co_X$-algebras 
\begin{equation}\label{p1-prem1g}
\iota_s\colon \rS_{\co_X}(\Omega)\stackrel{\sim}{\rightarrow} \cC.
\end{equation}
The latter is also induced by the trivialization of the principal homogenous $\bT$-bundle $\bL$ defined by the image of $s$ by \eqref{p1-prem1d} (see \cite{agt} II.4.10).

\subsection{}\label{p1-prem24}
Let $X$ be a scheme, and for $i=1,2$, let $\rT_i$ be a locally free $\co_X$-module of finite type, 
$\Omega_i=\cHom_{\co_X}(\rT_i,\co_X)$ its dual, $\cL_i$ a $\rT_i$-torsor of $X_\zar$. 
We denote by $\cF_i$ the $\co_X$-module of affine functions on $\cL_i$, by 
\begin{equation}\label{p1-prem24a}
0\longrightarrow \co_X\stackrel{c_i}{\longrightarrow}\cF_i\stackrel{\nu_i}{\longrightarrow} \Omega_i\longrightarrow 0
\end{equation}
the associated canonical extension and by 
\begin{equation}\label{p1-prem24b}
\cC_i=\underset{\underset{n\geq 0}{\longrightarrow}}\lim\ \rS^n_{\co_X}(\cF_i)
\end{equation}
the associated $\co_X$-algebra.  

The product $\cL_1\times \cL_2$ is naturally a torsor of $X_\zar$ under the $\co_X$-module $\rT_1\oplus \rT_2$. 
We denote by $\cF_{1,2}$ the $\co_X$-module of affine functions on $\cL_1\times \cL_2$, by 
\begin{equation}\label{p1-prem24c}
0\longrightarrow \co_X\stackrel{c_{12}}{\longrightarrow}\cF_{1,2}\stackrel{\nu_{12}}{\longrightarrow} \Omega_1\oplus \Omega_2\longrightarrow 0
\end{equation}
the associated canonical extension and by 
\begin{equation}\label{p1-prem24d}
\cC_{1,2}=\underset{\underset{n\geq 0}{\longrightarrow}}\lim\ \rS^n_{\co_X}(\cF_{1,2})
\end{equation}
the associated $\co_X$-algebra. 

For $i=1,2$, we denote by $\pi_i\colon \cL_1\times \cL_2\rightarrow \cL_i$ the canonical projection. 
It is equivariant with respect to the canonical projection $\rT_1\oplus \rT_2\rightarrow \rT_i$. 
We deduce by pullback (\cite{agt} II.4.12) a commutative diagram 
\begin{equation}\label{p1-prem24e}
\xymatrix{
0\ar[r]&{\co_X\oplus\co_X}\ar[r]^-(0.5){c_1\oplus c_2}\ar[d]_{\sigma}&{\cF_1\oplus\cF_2}\ar[r]^-(0.5){\nu_1\oplus\nu_2}\ar[d]_{\pi_1^*+\pi_2^*}&{\Omega_1\oplus\Omega_2}\ar[r]\ar@{=}[d]&0\\
0\ar[r]&{\co_X}\ar[r]^-(0.5){c_{1,2}}&{\cF_{1,2}}\ar[r]^-(0.5){\nu_{1,2}}&{\Omega_1\oplus\Omega_2}\ar[r]&0,}
\end{equation}
where $\sigma(a\oplus b)=a+b$. The morphisms $\pi_1^*$ and $\pi_2^*$ induce a homomorphism of $\co_X$-algebras
\begin{equation}\label{p1-prem24f}
\psi\colon \cC_1\otimes_{\co_X}\cC_2\rightarrow \cC_{1,2}. 
\end{equation}
Moreover, the $\co_X$-linear morphism 
\begin{equation}\label{p1-prem24g}
\cF_1\oplus \cF_2\rightarrow  \cC_1\otimes_{\co_X}\cC_2, \ \ \ f_1\oplus f_2\mapsto f_1\otimes 1+1\otimes f_2
\end{equation}
factors naturally through $\pi_1^*+\pi_2^*$ \eqref{p1-prem24e} and induces an $\co_X$-linear morphism 
$\cF_{1,2}\rightarrow  \cC_1\otimes_{\co_X}\cC_2$ whose composition with $c_{1,2}$ is the canonical
homomorphism. We deduce a homomorphism of $\co_X$-algebras
\begin{equation}\label{p1-prem24h}
\phi\colon \cC_{1,2} \rightarrow \cC_1\otimes_{\co_X}\cC_2.
\end{equation}
We easily see that $\psi$ and $\phi$ are isomorphisms inverse to each other. By \eqref{p1-prem24g}, $\phi$ induces an $\co_X$-linear morphism 
\begin{equation}\label{p1-prem24i}
\varphi\colon \cF_{1,2} \rightarrow \cF_1\otimes_{\co_X}\cF_2
\end{equation}
that fits into a commutative diagram 
\begin{equation}\label{p1-prem24j}
\xymatrix{
\cF_{1,2}\ar[rr]^-(0.5){\nu_{1,2}}\ar[d]_{\varphi}&&{\Omega_1\oplus \Omega_2}\ar[d]^{(\id_{\Omega_1}\otimes c_2)\oplus (c_1\otimes \id_{\Omega_2})}\\
{\cF_1\otimes_{\co_X}\cF_2}\ar[rr]^-(0.5){(\nu_1\otimes \id)\oplus(\id\otimes \nu_2)}&&{(\Omega_1\otimes \cF_2)\oplus (\cF_1\otimes_{\co_X}\Omega_2).}} 
\end{equation}

\subsection{}\label{p1-prem26}
We keep the assumptions and notation of \ref{p1-prem24}, moreover, let $\rT$ be a locally free $\co_X$-module of finite type, 
$\Omega=\cHom_{\co_X}(\rT,\co_X)$ its dual, $\cL$ a $\rT$-torsor of $X_\zar$. We denote by $\cF$ the $\co_X$-module of affine functions on $\cL$, by 
\begin{equation}\label{p1-prem26a}
0\longrightarrow \co_X\stackrel{c}{\longrightarrow}\cF\stackrel{\nu}{\longrightarrow} \Omega\longrightarrow 0,
\end{equation}
the associated canonical extension and by 
\begin{equation}\label{p1-prem26b}
\cC=\underset{\underset{n\geq 0}{\longrightarrow}}\lim\ \rS^n_{\co_X}(\cF)
\end{equation}
the associated $\co_X$-algebra. 

Let $v\colon \rT_1\oplus \rT_2\rightarrow \rT$ be an $\co_X$-linear morphism, $\varpi\colon \cL_1\times\cL_2 \rightarrow \cL$ a $v$-equivariant morphism. 
The latter induces by pullback (\cite{agt} II.4.12) a canonical $\co_X$-linear morphism $\varpi^*\colon \cF\rightarrow \cF_{1,2}$ that fits into a commutative diagram  
\begin{equation}\label{p1-prem26c}
\xymatrix{
0\ar[r]&{\co_X}\ar[r]\ar@{=}[d]&{\cF}\ar[r]\ar[d]^{\varpi^*}&{\Omega}\ar[r]\ar[d]^{u}&0\\
0\ar[r]&{\co_X}\ar[r]&{\cF_{1,2}}\ar[r]&{\Omega_1\oplus \Omega_2}\ar[r]&0,}
\end{equation}
where $u$ is the dual of $v$. In view of \ref{p1-prem24}, we deduce a canonical $\co_X$-linear morphism
\begin{equation}\label{p1-prem26d}
\upvarphi\colon \cF \rightarrow \cF_1\otimes_{\co_X}\cF_2,
\end{equation}
and a canonical homomorphism of $\co_X$-algebras
\begin{equation}\label{p1-prem26e}
\upphi\colon \cC\rightarrow \cC_1\otimes_{\co_X}\cC_2,
\end{equation}
compatible with $\upvarphi$. It follows from \eqref{p1-prem24j} that the diagram 
\begin{equation}\label{p1-prem26f}
\xymatrix{
\cF\ar[rr]^-(0.5)\nu\ar[d]_{\upvarphi}&&{\Omega}\ar[d]^{(u_1\otimes c_2)\oplus (c_1\otimes u_2)}\\
{\cF_1\otimes_{\co_X}\cF_2}\ar[rr]^-(0.5){(\nu_1\otimes \id)\oplus (\id\otimes \nu_2)}&&
{(\Omega_1\otimes_{\co_X}\cF_2)\oplus (\cF_1\otimes_{\co_X}\Omega_2),}} 
\end{equation}
where, for $i=1,2$, $u_i\colon \Omega\rightarrow \Omega_i$ is the composition of $u$ with the canonical projection, is commutative.

\begin{lem}\label{p1-prem18}
We keep the assumptions and notation of \ref{p1-prem26}, moreover, we set $\bT=\Spec(\rS_{\co_X}(\Omega))$, $\bL=\Spec(\cC)$, and for $i=1,2$, 
$\bT_i=\Spec(\rS_{\co_X}(\Omega_i))$ and $\bL_i=\Spec(\cC_i)$. Then, the morphism 
\begin{equation}
\bL_1\times_X\bL_2\rightarrow \bL,
\end{equation}
induced by $\upphi$ \eqref{p1-prem26e}, represents the morphism $\varpi\colon \cL_1\times\cL_2 \rightarrow \cL$ via the isomorphism \eqref{p1-prem1d}. 
In particular, for any $(s_1,s_2)\in \cL_1(X)\times \cL_2(X)$ and $s=\varpi(s_1,s_2)\in \cL(X)$, the diagram 
\begin{equation}
\xymatrix{
{\rS_{\co_X}(\Omega)}\ar[r]^-(0.5)\Delta\ar[d]_{\iota_s}&
{\rS_{\co_X}(\Omega_1)\otimes_{\co_X}\rS_{\co_X}(\Omega_2)}\ar[d]^{\iota_{s_1}\otimes \iota_{s_2}}\\
{\cC}\ar[r]^-(0.5){\upphi}&{\cC_1\otimes_{\co_X}\cC_2,}}
\end{equation}
where $\iota_{s_i}$ (resp.\ $\iota_s$) is the isomorphism \eqref{p1-prem1g} defined by $s_i$ (resp.\ $s$)
and $\Delta$ is the homomorphism of $\co_X$-algebras defined, for any local section $\omega$ of $\Omega$, 
by $\Delta(\omega)=u_1(\omega)\otimes 1 + 1 \otimes u_2(\omega)$, is commutative. 
\end{lem}

It follows immediately from the definitions and (\cite{agt} II.4.10).

\subsection{}\label{p1-prem12}
We take again the assumptions and notation of \ref{p1-prem26}, moreover, let $s\in \cL_1(X)$, $\tau_s\colon \rT_1\rightarrow \cL_1$ be the associated trivialization, 
$\rho_s\colon \cF_1\rightarrow \co_X$ the associated splitting of the extension $\cF_1$ \eqref{p1-prem24a}, 
$\varrho_s\colon\cC_1\rightarrow \co_X$ the induced homomorphism of $\co_X$-algebras (see \ref{p1-prem1}). 
We denote by $v_2\colon \rT_2\rightarrow \rT$ the $\co_X$-linear morphism induced by $v$, by $\varpi_s$ the $v_2$-equivariant morphism
\begin{equation}\label{p1-prem12a}
\varpi_s=\varpi\circ (s\times \id_{\cL_2})\colon \cL_2\rightarrow \cL,
\end{equation}
and by $\varpi_s^*\colon \cF \rightarrow  \cF_2$ the $\co_X$-linear morphism induced by pullback by $\varpi_s$. Then, we have
\begin{equation}\label{p1-prem12c}
\varpi_s^*=(\rho_s\otimes \id_{\cF_2})\circ \upvarphi. 
\end{equation}
Indeed, we are immediately reduced to the case where $\cL=\cL_1\times \cL_2$, in which case $\varpi_s=s\times \id_{\cL_2}\colon \cL_2\rightarrow \cL_1\times\cL_2$. 
In the following diagram
\begin{equation}\label{p1-prem12b}
\xymatrix{
{\cF_1\oplus\cF_2}\ar[r]^-(0.5){\pi_1^*+\pi_2^*}\ar[rrd]_{c_2\circ \rho_s+\id_{\cF_2}}&
{\cF_{1,2}}\ar[rd]^{\varpi_s^*}\ar[r]^-(0.5)\varphi&{\cF_1\otimes_{\co_X}\cF_2}\ar[d]^{\rho_s\otimes \id_{\cF_2}}\\
&&{\cF_2,}}
\end{equation}
the left triangle and the exterior triangle are obviously commutative \eqref{p1-prem24g}. Since $\pi_1^*+\pi_2^*$ is surjective \eqref{p1-prem24e}, so is the right triangle, which gives the required assertion 
\eqref{p1-prem12c}. 

It follows from \eqref{p1-prem12c} that the homomorphism of $\co_X$-algebras 
\begin{equation}\label{p1-prem12d}
(\varrho_s\otimes \id_{\cC_2})\circ \upphi \colon \cC\rightarrow \cC_2
\end{equation}
is induced by $\varpi_s^*$.

\subsection{}\label{p1-prem25}
Let $X$ be a scheme, and for $i=1,2,3,4$, let $\rT_i$ be a locally free $\co_X$-module of finite type, 
$\Omega_i=\cHom_{\co_X}(\rT_i,\co_X)$ its dual, $\cL_i$ a $\rT_i$-torsor of $X_\zar$. 
We denote by $\cF_i$ the $\co_X$-module of affine functions on $\cL_i$, by 
\begin{equation}\label{p1-prem25a}
0\rightarrow \co_X\rightarrow\cF_i\rightarrow \Omega_i\rightarrow 0
\end{equation}
the associated canonical extension and by 
\begin{equation}\label{p1-prem25b}
\cC_i=\underset{\underset{n\geq 0}{\longrightarrow}}\lim\ \rS^n_{\co_X}(\cF_i)
\end{equation}
the associated $\co_X$-algebra.  

For $i=1,3$, we denote by $\cF_{i,i+1}$ the $\co_X$-module of affine functions on $\cL_i\times \cL_{i+1}$, by 
\begin{equation}\label{p1-prem25c}
0\rightarrow \co_X\rightarrow\cF_{i,i+1}\rightarrow \Omega_i\oplus\Omega_{i+1}\rightarrow 0,
\end{equation}
the associated canonical extension and by 
\begin{equation}\label{p1-prem25d}
\cC_{i,i+1}=\underset{\underset{n\geq 0}{\longrightarrow}}\lim\ \rS^n_{\co_X}(\cF_{i,i+1})
\end{equation}
the associated $\co_X$-algebra. 

Let $v\colon \rT_1\oplus \rT_2\rightarrow \rT_3\oplus \rT_4$ be an $\co_X$-linear morphism, 
$\varpi\colon \cL_1\times \cL_2\rightarrow \cL_3\times\cL_4$ a $v$-equivariant morphism. 
By pullback (\cite{agt} II.4.12), $\varpi$ induces an $\co_X$-linear morphism $\varpi^*\colon \cF_{3,4} \rightarrow \cF_{1,2}$ 
that fits into a commutative diagram  
\begin{equation}\label{p1-prem25e}
\xymatrix{
0\ar[r]&{\co_X}\ar[r]\ar@{=}[d]&{\cF_{3,4}}\ar[r]\ar[d]^{\varpi^*}&{\Omega_3\oplus\Omega_4}\ar[r]\ar[d]^{u}&0\\
0\ar[r]&{\co_X}\ar[r]&{\cF_{1,2}}\ar[r]&{\Omega_1\oplus \Omega_2}\ar[r]&0,}
\end{equation}
where the lines are the exact sequences \eqref{p1-prem25c} and $u$ is the dual of $v$. 
By \ref{p1-prem24}, for $i=1,3$, we have a canonical isomorphism of $\co_X$-algebras 
\begin{equation}
\cC_{i,i+1} \stackrel{\sim}{\rightarrow} \cC_i\otimes_{\co_X}\cC_{i+1},
\end{equation}
and $\varpi^*$ induces a morphism of $\co_X$-algebras 
\begin{equation}\label{p1-prem25f}
\cC_3\otimes_{\co_X}\cC_4 \rightarrow \cC_1\otimes_{\co_X}\cC_2.
\end{equation}

\subsection{}\label{p1-prem9}
Let $X$ be a scheme, and for $i=1,2$, let $\rT_i$ be a locally free $\co_X$-module of finite type, 
$\Omega_i=\cHom_{\co_X}(\rT_i,\co_X)$ its dual, $\cL_i$ a $\rT_i$-torsor of $X_\zar$. 
We denote by $\cF_i$ the $\co_X$-module of affine functions on $\cL_i$, by 
\begin{equation}\label{p1-prem9a}
0\longrightarrow \co_X\stackrel{c_i}{\longrightarrow}\cF_i\stackrel{\nu_i}{\longrightarrow} \Omega_i\longrightarrow 0
\end{equation}
the associated canonical extension and by 
\begin{equation}\label{p1-prem9b}
\cC_i=\underset{\underset{n\geq 0}{\longrightarrow}}\lim\ \rS^n_{\co_X}(\cF_i)
\end{equation}
the associated $\co_X$-algebra.

Let $v\colon \rT_1\rightarrow \rT_2$ be an $\co_X$-linear morphism. 
We denote by $\cHom_v(\cL_1,\cL_2)$ the {\em sheaf of $v$-equivariant morphisms from $\cL_1$ to $\cL_2$}, 
i.e., the sheaf whose sections over any open subscheme $U$ of $X$ are the $v|U$-equivariant morphisms $\cL_1|U\rightarrow \cL_2|U$ (\cite{sga4} IV 10.2). 
The canonical action of $\rT_2$ on $\cL_2$ induces an action on $\cHom_v(\cL_1,\cL_2)$ making it into a $\rT_2$-torsor of $X_\zar$. 
We denote by $\cF$ the $\co_X$-module of affine functions on $\cHom_v(\cL_1,\cL_2)$, by 
\begin{equation}\label{p1-prem9c}
0\longrightarrow \co_X\stackrel{c}{\longrightarrow}\cF\stackrel{\nu}{\longrightarrow} \Omega_2\longrightarrow 0,
\end{equation}
the associated canonical extension and by 
\begin{equation}\label{p1-prem9d}
\cC=\underset{\underset{n\geq 0}{\longrightarrow}}\lim\ \rS^n_{\co_X}(\cF)
\end{equation}
the associated $\co_X$-algebra. 

We set $w=v+\id_{\rT_2}\colon \rT_1\oplus \rT_2\rightarrow \rT_2$ and denote by 
\begin{equation}\label{p1-prem9e}
\cL_1 \times \cHom_{\rT}(\cL_1,\cL_2) \rightarrow \cL_2
\end{equation}
the ``evaluation'', which is a $w$-equivariant morphism. 
By \eqref{p1-prem26d}, the latter induces a canonical $\co_X$-linear morphism 
\begin{equation}\label{p1-prem9f}
\upvarphi\colon \cF_2\rightarrow \cF_1\otimes_{\co_X} \cF.
\end{equation}
By \eqref{p1-prem26e}, it also induces a canonical homomorphism of $\co_X$-algebras 
\begin{equation}\label{p1-prem9g}
\upphi\colon \cC_2\rightarrow \cC_1\otimes_{\co_X} \cC,
\end{equation}
compatible with $\upvarphi$.

\subsection{}\label{p1-prem11}
We keep the assumptions and notation of \ref{p1-prem9}. 
Let $s\in \cL_1(X)$, $\rho_s\colon \cF_1\rightarrow \co_X$ be the associated splitting of the extension $\cF_1$ \eqref{p1-prem9a}. 
The evaluation at $s$ defines an isomorphism of $\rT_2$-torsors
\begin{equation}\label{p1-prem11a}
\delta_s\colon \cHom_v(\cL_1,\cL_2)\stackrel{\sim}{\rightarrow} \cL_2.
\end{equation}
By \eqref{p1-prem12c}, the associated pullback isomorphism $\delta_s^*\colon \cF_2 \stackrel{\sim}{\rightarrow}  \cF$ is given by 
\begin{equation}\label{p1-prem11b}
\delta_s^*=(\rho_s\otimes \id_\cF)\circ \upvarphi,
\end{equation}
where $\upvarphi$ is the morphism \eqref{p1-prem9f}. 

\subsection{}\label{p1-prem13} 
We keep the assumptions and notation of \ref{p1-prem9}. 
Let $\gamma\colon \cL_1\rightarrow \cL_2$ be a $v$-equivariant morphism, $\gamma^*\colon \cF_2 \rightarrow  \cF_1$ the associated pullback morphism, 
$\rho_\gamma\colon \cF\rightarrow \co_X$ the associated splitting of the extension \eqref{p1-prem9c}.
It follows from \eqref{p1-prem12c} that we have 
\begin{equation}\label{p1-prem13b}
\gamma^*=(\id_{\cF_1}\otimes \rho_\gamma)\circ \upvarphi,
\end{equation}
where $\upvarphi$ is the morphism \eqref{p1-prem9f}.

\subsection{}\label{p1-prem10}
We keep the assumptions and notation of \ref{p1-prem9}.
We have a canonical $\co_X$-linear morphism 
\begin{equation}\label{p1-prem10a}
\upmu\colon \cHom_{\co_X}(\cHom_{\co_X}(\cF_2,\cF_1),\co_X)\rightarrow \cF
\end{equation}
defined, for every open subscheme $U$ of $X$, by mapping an $\co_U$-linear map 
\begin{equation}\label{p1-prem10b}
g\colon \cHom_{\co_U}(\cF_2|U,\cF_1|U)\rightarrow \co_U 
\end{equation}
to the affine function $f\colon \cHom_v(\cL_1|U,\cL_2|U)\rightarrow \co_U$ defined, for every open subscheme $V$ of $U$,  by 
sending a $v$-equivariant morphism $\gamma\colon \cL_1|V\rightarrow \cL_2|V$ to $f(\gamma)=g(\gamma^*)$, 
where $\gamma^*\colon \cF_2|V\rightarrow \cF_1|V$ is the pullback morphism by $\gamma$. 
Indeed, for every $t\in \rT_2(V)$, that we consider as an $\co_V$-linear map $\Omega_2|V\rightarrow \co_V$, we have 
\begin{equation}\label{p1-prem10d}
f(\gamma+t)-f(\gamma)=g(\gamma^*+c_1\circ t\circ \nu_2)-g(\gamma^*)=g(c_1\circ t\circ \nu_2),
\end{equation}
which is clearly $\co_X(V)$-linear as a function of $t$. 

The diagram 
\begin{equation}\label{p1-prem10e}
\xymatrix{
{\cF_2}\ar[d]\ar[r]^-(0.5){\upvarphi}&{\cF\otimes_{\co_X}\cF_1}\\
{\cHom_{\co_X}(\cHom_{\co_X}(\cF_2,\cF_1),\cF_1)}\ar[r]^-(0.5)\sim&
{\cHom_{\co_X}(\cHom_{\co_X}(\cF_2,\cF_1),\co_X)\otimes_{\co_X}\cF_1,}\ar[u]_{\upmu\otimes\id}}
\end{equation} 
where the lower horizontal arrow is the canonical isomorphism and the left vertical arrow is the canonical morphism
defined by mapping a local section $f_2$ of $\cF_2$ to the morphism defined  
by $(u\colon \cF_2\rightarrow \cF_1)\mapsto u(f_2)$, is commutative. Indeed, since $\cF_1$ is a locally free $\co_X$-module of finite type, we may identify 
$\cF\otimes_{\co_X}\cF_1$ with (a submodule of) the $\co_X$-module of ``affine'' functions on $\cHom_v(\cL_1,\cL_2)$ with values in $\cF_1$. 
By \eqref{p1-prem13b}, for every open subscheme $U$ of $X$, every $f_2\in \cF_2(U)$ and every $\gamma\in \cHom_v(\cL_1,\cL_2)(U)$, we have
\begin{equation}
\upvarphi(f_2)(\gamma)= \gamma^*(f_2),
\end{equation}
where $\gamma^*\colon \cF_2|U\rightarrow \cF_1|U$ is the pullback by $\gamma$, which implies the claim in view of the definition of $\upmu$.  

\begin{rema}\label{p1-prem100}
We can reformulate the definition of $\upmu$ \eqref{p1-prem10a} by saying that the morphism 
\begin{equation}\label{p1-prem100a}
\cHom_v(\cL_1,\cL_2)\rightarrow \cHom_{\co_X}(\cF_2,\cF_1), \ \ \ \gamma\mapsto \gamma^*,
\end{equation}
is an affine function on the $\rT_2$-torsor $\cHom_v(\cL_1,\cL_2)$ with values in $\cHom_{\co_X}(\cF_2,\cF_1)$; 
this notion is an obvious generalization of that of an affine function with values in $\co_X$ defined in (\cite{agt} II.4.9). 
Its linear form is the image of the section $\nu_2\otimes c_1(1)$ by the canonical morphism
\begin{equation}\label{p1-prem100b}
\cHom_{\co_X}(\cF_2,\Omega_2)\otimes_{\co_X}\cF_1\rightarrow \cHom_{\co_X}(\cF_2,\cF_1)\otimes_{\co_X}\Omega_2.
\end{equation}
The morphism $\upmu$ \eqref{p1-prem10a} is then defined by evaluation of \eqref{p1-prem100a}, namely by sending a local linear form $g$ on $\cHom_{\co_X}(\cF_2,\cF_1)$ to the 
local affine function $\gamma\mapsto g(\gamma^*)$ on $\cHom_v(\cL_1,\cL_2)$. 
\end{rema}

\section{Koszul filtrations on divided power polynomial algebras}\label{p1-imdpa}

\subsection{}\label{p1-imdpa1}
We consider in this section a ringed $\mU$-topos  $(X,A)$ (\cite{sga4} IV 1.1) and an exact sequence of locally free $A$-modules of finite type 
\begin{equation}\label{p1-imdpa1a}
0\longrightarrow A\stackrel{c}{\longrightarrow} \cF\stackrel{\nu}{\longrightarrow} \cE \longrightarrow 0.
\end{equation}
By (\cite{illusie1} I 4.3.1.7), this sequence induces, for every integer $n\geq 0$, an exact sequence \eqref{p1-NC7}
\begin{equation}\label{p1-imdpa1b}
0\rightarrow \rS^{n}(\cF)\rightarrow \rS^{n+1}(\cF)\rightarrow \rS^{n+1}(\cE)\rightarrow 0.
\end{equation}
The $A$-modules $(\rS^{n}(\cF))_{n\in \mN}$ therefore form a filtered direct system, whose direct limit
\begin{equation}\label{p1-imdpa1c}
\cC=\underset{\underset{n\geq 0}{\longrightarrow}}\lim\ \rS^n(\cF)
\end{equation}
is naturally endowed with an $A$-algebra structure, called the {\em Higgs--Tate algebra associated with the extension $\cF$} \eqref{p1-imdpa1a}. 

The canonical injection $\cF\rightarrow \cC$ induces a surjective homomorphism of $A$-algebras 
\begin{equation}\label{p1-imdpa1d}
\phi\colon \rS(\cF)\rightarrow \cC,
\end{equation}
whose kernel is the ideal generated by $c(1)-1$.  
It induces for every integer $n\geq 0$, a surjective $A$-linear morphism
\begin{equation}\label{p1-imdpa1e}
\phi^n\colon \begin{array}[t]{clcr}
\oplus_{0\leq i\leq n}\rS^i(\cF)&\rightarrow &\rS^n(\cF),\\
\oplus_{0\leq i\leq n}x_i&\mapsto&\sum_{i=0}^{n} x_i \cdot c(1)^{n-i}.
\end{array}
\end{equation}

There is a unique homomorphism of $A$-algebras 
\begin{equation}\label{p1-imdpa1f}
\mu\colon \cC\rightarrow \rS(\cE)\otimes_A\cC,
\end{equation}
such that for every local section $x$ of $\cF$, we have 
\begin{equation}\label{p1-imdpa1g}
\mu(x)=\nu(x)\otimes 1+1\otimes x.
\end{equation}
It induces, for every integer $n\geq 0$, a morphism of $A$-modules
\begin{equation}\label{p1-imdpa1h}
\mu^n\colon \rS^n(\cF)\rightarrow \oplus_{0\leq i\leq n}\rS^i(\cE)\otimes_A\rS^{n}(\cF). 
\end{equation}

\begin{rema}\label{p1-imdpa100}
Assume given a splitting $\sigma\colon \cE\rightarrow \cF$ of  \eqref{p1-imdpa1a} and let $\rho\colon \cF\rightarrow A$ 
be the morphism induced by $\id_\cF-\sigma\circ \nu$. 
Then $\rho$ extends uniquely to a homomorphism of $A$-algebras $\varrho\colon \cC\rightarrow \co_X$ 
and $\sigma$ induces an isomorphism of $A$-algebras
\begin{equation}\label{p1-imdpa100a}
\iota_\sigma\colon \rS(\cE)\stackrel{\sim}{\rightarrow} \cC.
\end{equation}
\end{rema}

\subsection{}\label{p1-imdpa19}
There is a unique $A$-derivation of $\cC$ \eqref{p1-imdpa1c}
\begin{equation}\label{p1-imdpa19a}
d_\cC\colon \cC \rightarrow \cE\otimes_A \cC
\end{equation}
extending $\nu$ \eqref{p1-imdpa1a}. It canonically identifies by \ref{p1-imdpa100} with the universal $A$-derivation of $\cC$. 
We check immediately that it is also a Higgs $A$-field on $\cC$ with coefficients in $\cE$ \eqref{p1-delta-con1}. For every integer $n\geq 0$, 
$d_\cC$ induces a Higgs $A$-field 
\begin{equation}\label{p1-imdpa19b}
\theta^n\colon \rS^n(\cF)\rightarrow \cE\otimes_A\rS^n(\cF).
\end{equation}
We set $\cE^\vee=\cHom_A(\cE,A)$ and equip $\cC$ (resp.\ $\rS^n(\cF)$) with the $\rS(\cE^\vee)$-module structure defined by $d_\cC$ (resp.\ $\theta^n$) 
\eqref{p1-delta-con1j}. The canonical injection $\rS^n(\cF)\rightarrow \cC$ is then $\rS(\cE^\vee)$-linear.

\subsection{}\label{p1-imdpa3}
For any $A$-module $\cM$, we set $\cM^\vee=\cHom_A(\cM,A)$. By (\cite{bo} A.10), for every integer $n\geq 0$, we have a 
canonical $A$-linear morphism \eqref{p1-NC7}
\begin{equation}\label{p1-imdpa3a}
 \rS^n(\cM^\vee)\rightarrow \Gamma^n(\cM)^\vee.
\end{equation}
It induces compatible $A$-linear morphisms 
\begin{eqnarray}
\rS(\cM^\vee)&\rightarrow &\cHom_A(\Gamma(\cM),A), \label{p1-imdpa3b}\\
\oplus_{0\leq i\leq n}\rS^i(\cM^\vee)&\rightarrow &\cHom_A(\Gamma^{\leq n}(\cM),A), \label{p1-imdpa3c}
\end{eqnarray}
where $\Gamma^{\leq n}(\cM)$ is the $A$-algebra defined in \eqref{p1-NC7a}. 

If $\cM$ is locally free of finite type, \eqref{p1-imdpa3a}, \eqref{p1-imdpa3b} and \eqref{p1-imdpa3c} are isomorphisms. 

If the $A$-module $\cM$ is free with basis $x_1,\dots,x_d$,
we denote by $y_1,\dots,y_d$ the dual basis of $\cM^\vee$, and we set,  
for any $\un=(n_1,\dots,n_d)\in \mN^d$, $|\un|=\sum_{i=1}^dn_i$, and 
\begin{eqnarray}
\uy^\un&=&\prod_{i=1}^dy_i^{n_i},\label{p1-imdpa3d}\\
\ux^{[\un]}&=&x_1^{[n_1]}\dots x_d^{[n_d]}.\label{p1-imdpa3e}
\end{eqnarray}
Then, \eqref{p1-imdpa3a} transforms the basis $(\uy^\un)_{|\un|=n}$ of $\rS^n(\cM^\vee)$ 
into the dual basis of $(\ux^{[\un]})_{|\un|=n}$ of $\Gamma^n(\cM)^\vee$.

\subsection{}\label{p1-imdpa4}
The exact sequence \eqref{p1-imdpa1a} being locally split, it induces by duality an exact sequence of $A$-modules
\begin{equation}\label{p1-imdpa4a}
0\longrightarrow \cE^\vee\stackrel{\nu^\vee}{\longrightarrow} \cF^\vee\stackrel{c^\vee}{\longrightarrow} A \longrightarrow 0.
\end{equation}
By (\cite{illusie1} I 4.3.1.7), this sequence induces, for every integer $n\geq 0$, an exact sequence 
\begin{equation}\label{p1-imdpa4b}
0\rightarrow \Gamma^{n+1}(\cE^\vee)\rightarrow \Gamma^{n+1}(\cF^\vee)\rightarrow \Gamma^{n}(\cF^\vee)\rightarrow 0,
\end{equation}
which identifies by \eqref{p1-imdpa3a} with the dual of the sequence \eqref{p1-imdpa1b} (\cite{illusie1} I 4.3.1.4). 
The $A$-modules $(\Gamma^n(\cF^\vee))_{n\in \mN}$ therefore form a cofiltered inverse system. We set
\begin{equation}\label{p1-imdpa4c}
\cV=\underset{\underset{n\geq 0}{\longleftarrow}}\lim\ \Gamma^n(\cF^\vee). 
\end{equation}

\begin{prop}\label{p1-imdpa5}
For every integer $n\geq 0$, the morphisms $\phi^n$ \eqref{p1-imdpa1e} and $\mu^n$ \eqref{p1-imdpa1h}
induce by duality \eqref{p1-imdpa3} two $A$-linear morphisms 
\begin{eqnarray}
\cphi^n\colon \Gamma^n(\cF^\vee)\rightarrow \Gamma^{\leq n}(\cF^\vee),\label{p1-imdpa5a}\\
\cmu^n\colon \Gamma^{\leq n}(\cE^\vee)\otimes_A\Gamma^n(\cF^\vee)\rightarrow \Gamma^n(\cF^\vee),\label{p1-imdpa5b}
\end{eqnarray}
where $\Gamma^{\leq n}$ is the functor defined in \eqref{p1-NC7a}, satisfying the following properties:
\begin{itemize}
\item[{\rm (i)}] The morphism $\cmu^n$ equips $\Gamma^n(\cF^\vee)$ with a $\Gamma^{\leq n}(\cE^\vee)$-module structure. 
\item[{\rm (ii)}] The morphism $\cphi^n$ is $\Gamma^{\leq n}(\cE^\vee)$-linear and injective, 
where $\Gamma^{\leq n}(\cF^\vee)$ is considered as a $\Gamma^{\leq n}(\cE^\vee)$-algebra by $\nu^\vee$ \eqref{p1-imdpa4a}. 
\end{itemize}
Moreover, the diagram 
\begin{equation}\label{p1-imdpa5c}
\xymatrix{
{\Gamma^{n+1}(\cF^\vee)}\ar[r]^{\cphi^{n+1}}\ar[d]&{\Gamma^{\leq n+1}(\cF^\vee)}\ar[d]\\
{\Gamma^n(\cF^\vee)}\ar[r]^{\cphi^n}&{\Gamma^{\leq n}(\cF^\vee),}}
\end{equation}
where the left (resp.\ right) vertical arrow is the morphism defined in \eqref{p1-imdpa4b} (resp.\ the canonical morphism), is commutative. 
\end{prop}

Indeed, $\cphi^n$ is clearly injective. The diagram
\begin{equation}\label{p1-imdpa5d}
\xymatrix{
{\oplus_{0\leq i\leq n}\rS^i(\cF)}\ar[rrr]^-(0.5){\phi^n}\ar[d]&&&{\rS^n(\cF)}\ar[d]^{\mu^n}\\
{(\oplus_{0\leq i\leq n}\rS^i(\cF))\otimes_A(\oplus_{0\leq j\leq n}\rS^j(\cF))}\ar[rrr]^-(0.5){(\oplus_{0\leq i\leq n}\rS^i(\nu))\otimes\phi^n}&&&
{\oplus_{0\leq i\leq n}\rS^i(\cE)\otimes_A\rS^n(\cF),}}
\end{equation}
where the unlabeled arrow is induced by the canonical comultiplication of $\rS(\cF)$, is commutative. 
Recall that the multiplication of $\Gamma(\cF^\vee)$ and the comultiplication of $\rS(\cF)$ are dual (\cite{ov} 5.19(ii)). 
Therefore, diagram \eqref{p1-imdpa5d} induces by duality a commutative diagram 
\begin{equation}\label{p1-imdpa5e}
\xymatrix{
{\Gamma^{\leq n}(\cE^\vee)\otimes_A\Gamma^n(\cF^\vee)}\ar[rr]^-(0.5){\Gamma^{\leq n}(\nu^\vee)\otimes\cphi^n}\ar[d]_{\cmu^n}&&
{\Gamma^{\leq n}(\cF^\vee)\otimes_A\Gamma^{\leq n}(\cF^\vee)}\ar[d]\\
{\Gamma^n(\cF^\vee)}\ar[rr]^-(0.5){\cphi^n}&&{\Gamma^{\leq n}(\cF^\vee),}}
\end{equation}
where the unlabeled arrow is the multiplication of the $A$-algebra $\Gamma^{\leq n}(\cF^\vee)$. This proves propositions (i) and (ii). 

On the other hand, the diagram 
\begin{equation}
\xymatrix{
{\oplus_{0\leq i\leq n}\rS^i(\cF)}\ar[rr]^-(0.5){\phi^n}\ar@{^(->}[d]&&{\rS^n(\cF)}\ar[d]^{\cdot c(1)}\\
{\oplus_{0\leq i\leq n+1}\rS^i(\cF)}\ar[rr]^-(0.5){\phi^{n+1}}&&{\rS^{n+1}(\cF),}}
\end{equation}
where the unlabeled arrow is the canonical injection, is obviously commutative. 
It follows by duality that the diagram \eqref{p1-imdpa5c} is commutative.

\begin{cor}\label{p1-imdpa6}
The morphisms $\cphi^n$ \eqref{p1-imdpa5a} and $\cmu^n$ \eqref{p1-imdpa5b} define two $A$-linear morphisms 
\begin{eqnarray}
\cphi\colon \cV\rightarrow \hGamma(\cF^\vee),\label{p1-imdpa6a}\\
\cmu\colon \hGamma(\cE^\vee)\otimes_A\cV\rightarrow \cV,\label{p1-imdpa6b}
\end{eqnarray}
where $\hGamma$ is the functor defined in \eqref{p1-NC7b}
and $\cV$ is the $A$-module defined in \eqref{p1-imdpa4c}, satisfying the following properties:
\begin{itemize}
\item[{\rm (i)}] The morphism $\cmu$ equips $\cV$ with a $\hGamma(\cE^\vee)$-module structure. 
\item[{\rm (ii)}] The morphism $\cphi$ is $\hGamma(\cE^\vee)$-linear and injective, 
where  $\hGamma(\cF^\vee)$ is considered as a $\hGamma(\cE^\vee)$-algebra by $\nu^\vee$ \eqref{p1-imdpa4a}.  
\end{itemize}
\end{cor}

\begin{prop}\label{p1-imdpa8}
Let $\sigma\colon \cE\rightarrow \cF$ be an $A$-linear section of the extension \eqref{p1-imdpa1a}, 
$\sigma^\vee\colon \cF^\vee\rightarrow \cE^\vee$ the dual morphism. 
\begin{itemize}
\item[{\rm (i)}] For every integer $n\geq 0$, the composition of the $\Gamma^{\leq n}(\cE^\vee)$-linear morphisms 
\begin{equation}\label{p1-imdpa8a}
\xymatrix{
{\Gamma^n(\cF^\vee)}\ar[r]^{\cphi^n}&{\Gamma^{\leq n}(\cF^\vee)}\ar[rr]^{\Gamma^{\leq n}(\sigma^\vee)}&&
{\Gamma^{\leq n}(\cE^\vee),}}
\end{equation}
where $\cphi^n$ is defined in \eqref{p1-imdpa5a}, is an isomorphism. 
\item[{\rm (ii)}] The composition of the $\hGamma(\cE^\vee)$-linear morphisms  
\begin{equation}\label{p1-imdpa8b}
\xymatrix{
{\cV}\ar[r]^-(0.5){\cphi}&{\hGamma(\cF^\vee)}\ar[rr]^{\hGamma(\sigma^\vee)}&&
{\hGamma(\cE^\vee),}}
\end{equation}
where $\cphi$ is defined in \eqref{p1-imdpa6a}, is an isomorphism. 
\end{itemize}
\end{prop}

(i) Indeed, \eqref{p1-imdpa8a} is the dual of the isomorphism
\begin{equation}\label{p1-imdpa8c}
\phi^n_\sigma\colon \begin{array}[t]{clcr}
\oplus_{0\leq i\leq n}\rS^i(\cE)&\stackrel{\sim}{\rightarrow}& \rS^n(\cF),\\
\oplus_{0\leq i\leq n}x_i&\mapsto&\sum_{i=1}^n\rS^i(\sigma)(x_i)\cdot c(1)^{n-i}.
\end{array}
\end{equation}

(ii) Follows immediately from (i). 

\begin{cor}\label{p1-imdpa9}
Let $n$ be an integer $\geq 0$. 
\begin{itemize}
\item[{\rm (i)}] The $\Gamma^{\leq n}(\cE^\vee)$-module $\Gamma^n(\cF^\vee)$ is invertible.
\item[{\rm (ii)}] The $\hGamma(\cE^\vee)$-module $\cV$ \eqref{p1-imdpa6b} is invertible. 
\end{itemize}
\end{cor}

\subsection{}\label{p1-imdpa10}
We denote by $\Lambda$ the sheaf of (local) $A$-linear forms $\rho\colon \cF\rightarrow A$ such that $\rho\circ c=\id_A$. 
The canonical injection
\begin{equation}\label{p1-imdpa10a}
\omega\colon \Lambda\rightarrow \cF^\vee
\end{equation}
identifies $\Lambda$ with the inverse image of $1$ by the morphism $c^\vee$ \eqref{p1-imdpa4a}. It follows that $\Lambda$ is naturally 
equipped with the structure of a $\cE^\vee$-torsor, i.e. a torsor under the $A$-module $\cE^\vee$.

\begin{prop}\label{p1-imdpa11}
Let $\rho\colon \cF\rightarrow A$ be an $A$-linear form such that $\rho\circ c=\id_A$ \eqref{p1-imdpa1a},
that we consider as a section of $\cF^\vee$, $n$ an integer $\geq 0$. Then, 
\begin{itemize}
\item[{\rm (i)}] The canonical morphism $\Gamma^{n+1}(\cF^\vee)\rightarrow \Gamma^n(\cF^\vee)$ \eqref{p1-imdpa4b} maps 
$\rho^{[n+1]}$ to $\rho^{[n]}$. 
\item[{\rm (ii)}] The section $\rho^{[n]}$ is a basis of the $\Gamma^{\leq n}(\cE^\vee)$-module $\Gamma^n(\cF^\vee)$.
\item[{\rm (iii)}] The section $\cphi^n(\rho^{[n]})$ of $\Gamma^{\leq n}(\cF^\vee)$ \eqref{p1-imdpa5a} is the class of 
the unit $\exp_{\cF^\vee}(\rho)$ \eqref{p1-NC7c}. 
\end{itemize}
\end{prop}

Indeed, the questions being local, we may assume that the $A$-module $\cE$ is free. 
Let $y'_2,\dots,y'_d$ be an $A$-basis of $\cE$, $x'_2,\dots,x'_d$ the dual $A$-basis of $\cE^\vee$. 
We denote by $\sigma\colon \cE\rightarrow \cF$ the section of the extension \eqref{p1-imdpa1a} induced by $\id_\cF-c\circ \rho$. 
We set $y_1=c(1),y_2=\sigma(y'_2),\dots,y_d=\sigma(y'_d)$  \eqref{p1-imdpa1a}, 
and $x_1=\rho, x_2=\nu^\vee(x'_2),\dots,x_d=\nu^\vee(x'_d)$  \eqref{p1-imdpa4a},
which are dual $A$-bases of $\cF$ and $\cF^\vee$ respectively.  
For any $\un=(n_1,\dots,n_d)\in \mN^d$,
we define $\uy^\un$ (resp.\ $\ux^{[\un]}$) by the formula \eqref{p1-imdpa3d} (resp.\ \eqref{p1-imdpa3e}). 
For $\ualpha=(\alpha_1,\dots,\alpha_n)\in \Gamma(X,A)^d$, we set $\ualpha^\un=\prod_{i=1}^d\alpha_i^{n_i}$.

Let 
\begin{eqnarray}
\varphi^n\colon \Gamma^n(\cF^\vee)&\stackrel{\sim}{\rightarrow} &\rS^n(\cF)^\vee,\label{p1-imdpa11a}\\
\varphi^{\leq n}\colon \Gamma^{\leq n}(\cF^\vee)&\stackrel{\sim}{\rightarrow} &\cHom_A(\oplus_{0\leq i\leq n}\rS^i(\cF),A), \label{p1-imdpa11b}
\end{eqnarray}
be the canonical duality isomorphisms \eqref{p1-imdpa3}. 
Then, $\varphi^n$ maps the $A$-basis $(\ux^{[\un]})_{|\un|=n}$ of $\Gamma^n(\cF^\vee)$ 
to the dual $A$-basis of $(\uy^\un)_{|\un|=n}$.
Similarly, $\varphi^{\leq n}$ maps the $A$-basis of $\Gamma^{\leq n}(\cF^\vee)$ 
made of the classes of $(\ux^{[\un]})_{|\un|\leq n}$ to the dual $A$-basis of $(\uy^\un)_{|\un|\leq n}$.

For any  $\un\in \mN^d$ with $|\un|=n$, we have 
\begin{equation}\label{p1-imdpa11c}
\varphi^n(\rho^{[n]})(\uy^\un)=
\left\{
\begin{array}{clcr}
1 &{\rm if}\ \un=(n,0,\dots,0), \\
0 &{\rm otherwise}.
\end{array}
\right.
\end{equation}
Proposition (i) follows immediately since \eqref{p1-imdpa4b} is the dual of \eqref{p1-imdpa1b}. 

For every $x\in \Gamma(X,\cF^\vee)$, with $x=\sum_{i=1}^d\alpha_ix_i$ and $\alpha_i\in \Gamma(X,A)$ for all $1\leq i\leq d$, we have 
\begin{equation}\label{p1-imdpa11d}
\exp_{\cF^\vee}(x)=\sum_{\un\in \mN^d}\ualpha^\un\ux^{[\un]}. 
\end{equation}
Therefore, denoting the class of $\exp_{\cF^\vee}(x)$ in $\Gamma^{\leq n}(\cF^\vee)$ by the same notation, 
we have, for every $\un\in \mN^d$ with $|\un|\leq n$, 
\begin{equation}\label{p1-imdpa11e}
\varphi^{\leq n}(\exp_{\cF^\vee}(x))(\uy^\un)=\ualpha^\un.
\end{equation}
In particular, for any $\un\in \mN^d$ with $|\un|\leq n$, we have 
\begin{equation}\label{p1-imdpa11f}
\varphi^{\leq n}(\exp_{\cF^\vee}(\rho))(\uy^\un)=
\left\{
\begin{array}{clcr}
1 &{\rm if}\ \un=(n_1,0,\dots,0), \ \forall n_1\in \mN,\\
0 &{\rm otherwise}.
\end{array}
\right.
\end{equation}
Proposition (iii) follows immediately from the definition of $\cphi_n$ \eqref{p1-imdpa5a}. 

Let $\sigma\colon \cE\rightarrow \cF$ be the section of the extension \eqref{p1-imdpa1a} induced by $\id_\cF-c\circ \rho$,
$\sigma^\vee\colon \cF^\vee\rightarrow \cE^\vee$ the dual morphism. 
By (iii) and the functoriality of the exponential morphism \eqref{p1-NC7c}, we have 
\begin{equation}
\Gamma^{\leq n}(\sigma^\vee)(\cphi^n(\rho^{[n]}))=\Gamma^{\leq n}(\sigma^\vee)(\exp_{\cF^\vee}(\rho))=
\exp_{\cE^\vee}(\sigma^\vee(\rho))=1,
\end{equation}
which is a basis of the $\Gamma^{\leq n}(\cE^\vee)$-module $\Gamma^{\leq n}(\cE^\vee)$. 
Proposition (ii) follows since $\Gamma^{\leq n}(\sigma^\vee)\circ \cphi^n \colon \Gamma^n(\cF^\vee)
\rightarrow \Gamma^{\leq n}(\cE^\vee)$ is a $\Gamma^{\leq n}(\cE^\vee)$-linear isomorphism by \ref{p1-imdpa8}(i).

\begin{cor}\label{p1-imdpa12}
Let $\rho\colon \cF\rightarrow A$ be an $A$-linear form such that $\rho\circ c=\id_A$ \eqref{p1-imdpa1a},
that we consider as a section of $\cF^\vee$. Then, 
\begin{itemize}
\item[{\rm (i)}] There exists a canonical section $\rho^{[\infty]}$ of $\Gamma(X,\cV)$ \eqref{p1-imdpa4c} 
whose image in $\Gamma^n(\cF^\vee)$ is $\rho^{[n]}$, for every $n\geq 0$. 
\item[{\rm (ii)}] The section $\rho^{[\infty]}$ is a basis of the $\hGamma(\cE^\vee)$-module $\cV$.
\item[{\rm (iii)}] We have $\cphi(\rho^{[\infty]})=\exp_{\cF^\vee}(\rho)$ \eqref{p1-NC7c}.
\end{itemize} 
\end{cor}

It follows immediately from \ref{p1-imdpa11}.

\begin{cor}\label{p1-imdpa13}
There is a canonical morphism 
\begin{equation}\label{p1-imdpa13a}
\exp_\Lambda\colon \Lambda\rightarrow \cV,
\end{equation} 
where $\Lambda$ is the $\cE^\vee$-torsor defined in  \ref{p1-imdpa10}, 
fitting into the commutative diagram
\begin{equation}\label{p1-imdpa13b}
\xymatrix{
{\Lambda}\ar[rr]^-(0.5){\exp_\Lambda}\ar[d]_{\omega}&&{\cV}\ar[d]^-(0.5){\cphi}\\
{\cF^\vee}\ar[rr]^-(0.5){\exp_{\cF^\vee}}&&{\hGamma(\cF^\vee),}}
\end{equation}
where $\exp_{\cF^\vee}$, $\cphi$ and $\omega$ are defined in \eqref{p1-NC7c}, \eqref{p1-imdpa6a} and \eqref{p1-imdpa10a}, respectively. 
Moreover, $\exp_\Lambda$ is equivariant with respect to the homomorphism 
$\exp_{\cE^\vee}\colon \cE^\vee\rightarrow \hGamma(\cE^\vee)^\times$, and its image is contained in 
the subsheaf of local bases of the invertible $\hGamma(\cE^\vee)$-module $\cV$, i.e., the $\hGamma(\cE^\vee)^\times$-torsor 
associated with  $\cV$. 
\end{cor}

Indeed, $\exp_\Lambda$ is defined by mapping a local section $\rho$ of $\Lambda$, considered as a local section of $\cF^\vee$, to the 
local section $\rho^{[\infty]}$ of $\cV$ defined in \ref{p1-imdpa12}(i). 
The commutativity of the diagram \eqref{p1-imdpa13b} follows from \ref{p1-imdpa12}(iii). 
The equivariance of  $\exp_\Lambda$ with respect to the homomorphism $\exp_{\cE^\vee}$ follows from \eqref{p1-imdpa13b} and the fact 
that for any $A$-module $\cM$, $\exp_\cM\colon \cM\rightarrow \hGamma(\cM)^\times$ is a group homomorphism, functorial in $\cM$.  
The last assertion follows from \ref{p1-imdpa12}(ii).

\begin{cor}\label{p1-imdpa14}
The invertible $\hGamma(\cE^\vee)$-module $\cV$ \eqref{p1-imdpa6b} is canonically isomorphic to the line bundle associated with the 
$\hGamma(\cE^\vee)^\times$-torsor $\Lambda\wedge^{\cE^\vee}\hGamma(\cE^\vee)^\times$ \eqref{p1-NC4}, 
deduced from $\Lambda$ by extension of its structural group by the homomorphism $\exp_{\cE^\vee}$ \eqref{p1-NC7c}.
\end{cor}

\subsection{}\label{p1-imdpa15}
We consider $\cF$ \eqref{p1-imdpa1a} as an extension of $\cE$ by $A$ by the morphisms $\nu$ and $-c$ and denote it for clarity by 
adding a prime exponent:
\begin{equation}\label{p1-imdpa15a}
0\longrightarrow A\stackrel{c'}{\longrightarrow} \cF'\stackrel{\nu'}{\longrightarrow} \cE \longrightarrow 0,
\end{equation}
so $\cF'=\cF$, $\nu'=\nu$ and $c'=-c$. We denote by
\begin{equation}\label{p1-imdpa15b}
0\longrightarrow \cE^\vee\stackrel{\nu'^\vee}{\longrightarrow} \cF'^\vee\stackrel{c'^\vee}{\longrightarrow} A \longrightarrow 0.
\end{equation}
the dual exact sequence, so $\cF'^\vee=\cF^\vee$ \eqref{p1-imdpa4a}. 

The extension \eqref{p1-imdpa15a} defines an $A$-algebra 
\begin{equation}\label{p1-imdpa15e}
\cC'=\underset{\underset{n\geq 0}{\longrightarrow}}\lim\ \rS^n(\cF'),
\end{equation}
similar to $\cC$ \eqref{p1-imdpa1c}. The multiplication by $-1$ on $\cF$ induces an isomorphism of $A$-algebras 
$\cC\stackrel{\sim}{\rightarrow}\cC'$. 

We denote by $\Lambda'$ 
the sheaf of (local) $A$-linear forms $\rho'\colon \cF'\rightarrow A$ such that $\rho'\circ c'=\id_A$, i.e., 
the $\cE^\vee$-torsor $(c'^\vee)^{-1}(1)$, and by $\omega'\colon \Lambda'\rightarrow \cF'$ the canonical injection. 
The multiplication by $-1$ on $\cF^\vee$ induces an isomorphism $\Lambda\stackrel{\sim}{\rightarrow}\Lambda'$ that fits into
a commutative diagram
\begin{equation}\label{p1-imdpa15c}
\xymatrix{
\Lambda\ar[r]^-(0.5){\omega}\ar[d]&\cF^\vee\ar[d]\\
\Lambda'\ar[r]^-(0.5){\omega'}&\cF'^\vee.}
\end{equation}
We deduce a canonical isomorphism of $\cE^\vee$-torsors
\begin{equation}\label{p1-imdpa15d}
\Lambda\wedge^{\cE^\vee,-1}\cE^\vee\stackrel{\sim}{\rightarrow} \Lambda',
\end{equation}
where the source is the $\cE^\vee$-torsor deduced from $\Lambda$ by extension of its structural group 
by the multiplication by $-1$ on $\cE^\vee$.  

\subsection{}\label{p1-imdpa16}
By \ref{p1-imdpa5}, for every integer $n\geq 0$, $\Gamma^n(\cF'^\vee)$ is canonically equipped with 
a $\Gamma^{\leq n}(\cE^\vee)$-module structure defined by the extension \eqref{p1-imdpa15a}, 
and we have a canonical injective $\Gamma^{\leq n}(\cE^\vee)$-linear morphism
\begin{equation}\label{p1-imdpa16a}
\cphi'^n\colon \Gamma^n(\cF'^\vee)\rightarrow \Gamma^{\leq n}(\cF'^\vee). 
\end{equation}
Observe that $\Gamma^{\leq n}(\cF'^\vee)=\Gamma^{\leq n}(\cF^\vee)$ as $\Gamma^{\leq n}(\cE^\vee)$-modules since $\nu'=\nu$. 
However, the $\Gamma^{\leq n}(\cE^\vee)$-module structures on $\Gamma^n(\cF'^\vee)$ and $\Gamma^n(\cF^\vee)$
and the injections $\cphi'^n$ and $\cphi^n$ are different, 
although the sheaves of sets underlying  these modules are the same.

The $A$-modules $(\Gamma^n(\cF'))_{n\in \mN}$ form naturally a cofiltered inverse system defined by \eqref{p1-imdpa15b}. We set
\begin{equation}\label{p1-imdpa16b}
\cV'=\underset{\underset{n\geq 0}{\longleftarrow}}\lim\ \Gamma^n(\cF'^\vee),
\end{equation}
which we take care not to confuse with $\cV$ \eqref{p1-imdpa4c} defined by the exact sequence \eqref{p1-imdpa4a}. 
By \ref{p1-imdpa6}, $\cV'$ is canonically equipped with a $\hGamma(\cE^\vee)$-module structure, and we have a canonical 
injective $\hGamma(\cE^\vee)$-linear morphism 
\begin{equation}\label{p1-imdpa16c}
\cphi'\colon \cV'\rightarrow \hGamma(\cF'^\vee). 
\end{equation}
Observe that $\hGamma(\cF'^\vee)=\hGamma(\cF^\vee)$ as $\hGamma(\cE^\vee)$-modules since $\nu'=\nu$.

\begin{prop}\label{p1-imdpa17}
Let $n$ be an integer $\geq 0$. For all local sections $x$ of $\Gamma^n(\cF^\vee)$ and $x'$ of $\Gamma^n(\cF'^\vee)$, 
the product $\cphi^n(x) \cphi'^n(x')$ is contained in 
$\Gamma^{\leq n}(\cE^\vee) \subset \Gamma^{\leq n}(\cF^\vee)=\Gamma^{\leq n}(\cF'^\vee)$.   
The induced $\Gamma^{\leq n}(\cE^\vee)$-bilinear pairing
\begin{equation}\label{p1-imdpa17a}
\langle \ ,\ \rangle \colon 
\begin{array}[t]{clcr}
\Gamma^n(\cF^\vee)\times\Gamma^n(\cF'^\vee)&\rightarrow& \Gamma^{\leq n}(\cE^\vee),\\
(x,x')&\mapsto& \cphi^n(x) \cphi'^n(x'),
\end{array}
\end{equation}
is perfect, i.e., it induces an isomorphism 
\begin{equation}\label{p1-imdpa17b}
\Gamma^n(\cF^\vee)\otimes_{\Gamma^{\leq n}(\cE^\vee)}\Gamma^n(\cF'^\vee)\stackrel{\sim}{\rightarrow}\Gamma^{\leq n}(\cE^\vee).
\end{equation}
\end{prop}

Observe first that the pairing $(x,x')\mapsto \cphi^n(x) \cphi'^n(x')$ with values in $\Gamma^{\leq n}(\cF^\vee)$ is clearly 
$\Gamma^{\leq n}(\cE^\vee)$-bilinear. The statements of the proposition being local, we may assume 
that the $A$-module $\cE$ is free. Let $\sigma\colon \cE\rightarrow \cF$, $\rho\colon \cF\rightarrow A$ be 
two $A$-linear morphisms such that $\id_\cF=\sigma\circ \nu+c\circ\rho$. 
We set $\sigma'=\sigma\colon \cE\rightarrow \cF'$ and $\rho'=-\rho\colon \cF'\rightarrow A$, so that 
$\id_{\cF'}=\sigma'\circ \nu'+c'\circ\rho'$. 
By \ref{p1-imdpa11}(ii), $\Gamma^n(\cF^\vee)$ (resp.\ $\Gamma^n(\cF'^\vee)$) is a free $\Gamma^{\leq n}(\cE^\vee)$-module 
generated by $\rho^{[n]}$ (resp.\ $\rho'^{[n]}$). By \ref{p1-imdpa11}(iii) and with the same notation, 
we have 
\begin{equation}\label{p1-imdpa17c}
\cphi^n(\rho^{[n]}) \cphi'^n(\rho'^{[n]})= \exp_{\cF^\vee}(\rho) \exp_{\cF^\vee}(\rho') = \exp_{\cF^\vee}(\rho+\rho')=1,
\end{equation}
which proves the proposition. 

\begin{cor}\label{p1-imdpa18}
For all local sections $x$ of $\cV$ and $x'$ of $\cV'$ \eqref{p1-imdpa16b}, 
the product $\cphi(x) \cphi'(x')$ is contained in 
$\hGamma(\cE^\vee) \subset \hGamma(\cF^\vee)=\hGamma(\cF'^\vee)$.  
The induced $\hGamma(\cE^\vee)$-bilinear pairing
\begin{equation}\label{p1-imdpa18a}
\langle \ ,\ \rangle \colon 
\begin{array}[t]{clcr}
\cV\times\cV'&\rightarrow& \hGamma(\cE^\vee),\\
(x,x')&\mapsto& \cphi(x) \cphi'(x'),
\end{array}
\end{equation}
is perfect, i.e., it induces an isomorphism 
\begin{equation}\label{p1-imdpa18b}
\cV\otimes_{\hGamma(\cE^\vee)}\cV'\stackrel{\sim}{\rightarrow}\hGamma(\cE^\vee).
\end{equation}
\end{cor}

\begin{prop}\label{p1-imdpa20}
For every integer $n\geq 0$, the duality isomorphism \eqref{p1-imdpa3a}
\begin{equation}\label{p1-imdpa20a}
\rS^n(\cF)\stackrel{\sim}{\rightarrow} \cHom_A(\Gamma^n(\cF^\vee),A),
\end{equation}
is $\rS(\cE^\vee)$-linear when we equip $\rS^n(\cF)$ with the action of $\rS(\cE^\vee)$ defined by the Higgs $A$-field 
$\theta^n$ \eqref{p1-imdpa19b}, see \ref{p1-delta-con1j}, and 
$\Gamma^n(\cF^\vee)$ with  the action of $\rS(\cE^\vee)$ induced by $\cmu^n$ \eqref{p1-imdpa5b} and the canonical morphism
$\rS(\cE^\vee)\rightarrow \Gamma^{\leq n}(\cE^\vee)$. 
\end{prop}

Indeed, let $\pi_0\colon \rS(\cE)\rightarrow A$ and $\pi_1\colon \rS(\cE)\rightarrow \cE$ be the canonical projections. 
We denote by $\delta$ the composition
\begin{equation}
\delta\colon 
\xymatrix{
{\cC}\ar[r]^-(0.5){\mu}&{\rS(\cE)\otimes_A\cC}\ar[r]^-(0.5){\pi_1\otimes \id}&{\cE\otimes_A\cC,}}
\end{equation}
where $\mu$ is the morphism defined in \eqref{p1-imdpa1f}.  
Since $\pi_0\circ \mu=\id_\cC$ and $\mu$ is a ring homomorphism, for all local sections $x,y$ of $\cC$, we have 
\begin{equation}
\delta(xy)= x \cdot \delta(y)+y \cdot \delta(x). 
\end{equation}
Therefore, $\delta$ is an $A$-derivation. By \eqref{p1-imdpa1g}, the restriction 
of $\delta$ to $\cF$ coincides with $\nu$ \eqref{p1-imdpa1a}. Hence $\delta=d_\cC$. 
We deduce that $\theta^n$ \eqref{p1-imdpa19b} is equal to the composition 
\begin{equation}
\xymatrix{
{\rS^n(\cF)}\ar[r]^-(0.5){\mu^n}&{\oplus_{0\leq i\leq n}\rS^i(\cE)\otimes_A\rS^n(\cF)}\ar[r]^-(0.5){\pi^n_1\otimes \id}&{\cE\otimes_A\rS^n(\cF),}}
\end{equation}
where $\mu^n$ is defined in \eqref{p1-imdpa1h} and 
$\pi^n_1\colon \oplus_{0\leq i\leq n}\rS^i(\cE)\rightarrow \cE$ is the canonical projection. 
The proposition follows in view of the definition of $\cmu^n$ \eqref{p1-imdpa5b}.

\begin{prop}\label{p1-imdpa21}
Let $(N,\theta)$ be a Higgs $A$-module with coefficients in $\cE$ \eqref{p1-delta-con1}, $n$ an integer $\geq 0$. 
We equip $\rS^n(\cF)$ with the Higgs $A$-field $\theta^n$ \eqref{p1-imdpa19b},
and $\rS^n(\cF)\otimes_AN$ with the total Higgs $A$-field
\begin{equation}\label{p1-imdpa21b}
\theta^n_\tot=\theta^n\otimes \id+\id\otimes \theta. 
\end{equation}
We equip $N$ with the action of $\rS(\cE^\vee)$ induced by $\theta$ \eqref{p1-delta-con1j}, 
and $\Gamma^n(\cF'^\vee)$ \eqref{p1-imdpa15a} with the action of $\rS(\cE^\vee)$ defined by the composition
\begin{equation}
\xymatrix{
{\rS(\cE^\vee)\otimes_A\Gamma^n(\cF'^\vee)}\ar[r]&
{\Gamma^{\leq n}(\cE^\vee)\otimes_A\Gamma^n(\cF'^\vee)}\ar[r]^-(0.5){\cmu'^n}&{\Gamma^n(\cF'^\vee),}}
\end{equation}
where the first arrow is the canonical morphism and the second is defined in \eqref{p1-imdpa5b} relatively to the extension \eqref{p1-imdpa15a}. 
Then,
\begin{itemize}
\item[{\rm (i)}] The morphism 
\begin{equation}\label{p1-imdpa21c}
\begin{array}{clcr}
\cE^\vee\otimes_A\cHom_A(\Gamma^n(\cF'^\vee),N)&\rightarrow& \cHom_A(\Gamma^n(\cF'^\vee),N),\\
\xi\otimes \varphi&\mapsto &(x\mapsto \xi\cdot\varphi(x)-\varphi(\xi\cdot x)),
\end{array}
\end{equation}
defines a Higgs $A$-field on $\cHom_A(\Gamma^n(\cF'^\vee),N)$ with coefficients in $\cE$. 
\item[{\rm (ii)}]  The composition of the isomorphisms 
\begin{equation}\label{p1-imdpa21a}
\rS^n(\cF)\otimes_AN\stackrel{\sim}{\rightarrow}\rS^n(\cF')\otimes_AN\stackrel{\sim}{\rightarrow} \cHom_A(\Gamma^n(\cF'^\vee),N),
\end{equation}
where the first one is induced by the multiplication by $-1$ on $\cF$ and the second by the duality isomorphism \eqref{p1-imdpa3a}, 
is an isomorphism of Higgs $A$-modules for the Higgs fields defined in \eqref{p1-imdpa21b} and {\rm (i)}.  
\item[{\rm (iii)}] The isomorphism \eqref{p1-imdpa21a} induces an isomorphism 
\begin{equation}\label{p1-imdpa21d}
(\rS^n(\cF)\otimes_AN)^{\theta^n_\tot=0}\stackrel{\sim}{\rightarrow} \cHom_{\rS(\cE^\vee)}(\Gamma^n(\cF'^\vee),N). 
\end{equation}
\item[{\rm (iv)}] The diagram
\begin{equation}\label{p1-imdpa21e}
\xymatrix{
{(\rS^n(\cF)\otimes_AN)^{\theta^n_\tot=0}}\ar[r]\ar[d]&{\cHom_{\rS(\cE^\vee)}(\Gamma^n(\cF'^\vee),N)}\ar[d]\\
{(\rS^{n+1}(\cF)\otimes_AN)^{\theta^{n+1}_\tot=0}}\ar[r]&{\cHom_{\rS(\cE^\vee)}(\Gamma^{n+1}(\cF'^\vee),N),}}
\end{equation}
where the horizontal arrows are the isomorphisms \eqref{p1-imdpa21d}, the left (resp.\ right) vertical arrow is induced by the morphism 
defined in \eqref{p1-imdpa1b} (resp.\ \eqref{p1-imdpa4b}), is commutative. 
\end{itemize}
\end{prop}

(i) It is obvious, see \ref{p1-delta-con1j}.

(ii) We denote by $\iota\colon \cF\rightarrow \cF'$ the isomorphism induced by the multiplication by $-1$ on $\cF$
and by $\theta'^n\colon \rS^n(\cF')\rightarrow \cE\otimes_A\rS^n(\cF')$ the Higgs $A$-field induced by 
$d_{\cC'}$ \eqref{p1-imdpa15e}, similar to \eqref{p1-imdpa19b}. The diagrams
\begin{equation}\label{p1-imdpa21f}
\xymatrix{
0\ar[r]&A\ar[r]^{c}\ar@{=}[d]&\cF\ar[r]^{\nu}\ar[d]^{\iota}&\cE\ar[r]\ar[d]^{-1}&0\\
0\ar[r]&A\ar[r]^{c'}&\cF'\ar[r]^{\nu'}&\cE\ar[r]&0}
\end{equation}
\begin{equation}\label{p1-imdpa21g}
\xymatrix{
{\rS^n(\cF)}\ar[r]^-(0.5){\rS^n(\iota)}\ar[d]_{\theta^n}&{\rS^n(\cF')}\ar[d]^{-\theta'^n}\\
{\cE\otimes_A\rS^n(\cF)}\ar[r]^-(0.5){\id\otimes \iota}&{\cE\otimes_A\rS^n(\cF')}}
\end{equation}
are commutative. The proposition follows then immediately from \ref{p1-imdpa20}. 

(iii) It follows from (ii). 

(iv) It follows immediately from the duality of the exact sequences \eqref{p1-imdpa1b} and \eqref{p1-imdpa4b} (\cite{illusie1} I 4.3.1.4)
and the commutativity of \eqref{p1-imdpa21f}. 

\subsection{}\label{p1-imdpa22}
Let $n$ an integer $\geq 0$. We defined the Higgs $A$-field 
$\theta^n\colon \rS^n(\cF)\rightarrow \cE\otimes_A\rS^n(\cF)$ in \eqref{p1-imdpa19b}. 
On the other hand, the canonical homomorphism 
$\varsigma \colon \rS(\cE^\vee)\rightarrow \Gamma(\cE^\vee)$ defines a Higgs $A$-field
\begin{equation}\label{p1-imdpa22a}
\gamma\colon \Gamma(\cE^\vee)\rightarrow \cE\otimes_A \Gamma(\cE^\vee).
\end{equation}
For every ideal $I$ of $\Gamma(\cE^\vee)$, we have $\gamma(I)\subset \cE\otimes_A I$. 
Therefore, $\gamma$ induces a Higgs $A$-field on $\Gamma^{\leq n}(\cE^\vee)$  with coefficients in $\cE$, that we denote by $\gamma^{n}$. 
We equip $\rS^n(\cF)\otimes_A\Gamma^{\leq n}(\cE^\vee)$ with the total Higgs $A$-field
\begin{equation}\label{p1-imdpa22b}
\gamma^{n}_\tot=\theta^n\otimes \id+\id\otimes \gamma^{n}.
\end{equation}
Then, by \ref{p1-imdpa17} and \ref{p1-imdpa21}, 
we deduce a canonical $A$-linear morphism 
\begin{equation}\label{p1-imdpa22c}
\Gamma^n(\cF^\vee)\rightarrow (\rS^n(\cF)\otimes_A\Gamma^{\leq n}(\cE^\vee))^{\gamma^n_\tot=0}
\end{equation}
that fits into a commutative diagram 
\begin{equation}\label{p1-imdpa22d}
\xymatrix{
{\Gamma^n(\cF^\vee)}\ar[r]\ar[d]&{(\rS^n(\cF)\otimes_A\Gamma^{\leq n}(\cE^\vee))^{\gamma^n_\tot=0}}\ar[d]\\
{\cHom_{\Gamma^{\leq n}(\cE^\vee)}(\Gamma^n(\cF'^\vee),\Gamma^{\leq n}(\cE^\vee))}\ar[r]&
{\cHom_{\rS(\cE^\vee)}(\Gamma^n(\cF'^\vee),\Gamma^{\leq n}(\cE^\vee)),}}
\end{equation}
where the left vertical arrow is the isomorphism induced by \eqref{p1-imdpa17b}, 
the right vertical arrow is the isomorphism \eqref{p1-imdpa21d} and the lower horizontal arrow is the injection 
induced by the canonical homomorphism 
$\varsigma^{\leq n} \colon \rS(\cE^\vee)\rightarrow \Gamma^{\leq n}(\cE^\vee)$.

Embedding $\rS^n(\cF)$ in $\cC$ \eqref{p1-imdpa1c}, the morphism \eqref{p1-imdpa22c} induces a $\cC$-linear morphism  
\begin{equation}\label{p1-imdpa22e}
\cC\otimes_A\Gamma^n(\cF^\vee)\rightarrow \cC\otimes_A\Gamma^{\leq n}(\cE^\vee).
\end{equation}
It is in fact a morphism of $\cC$-modules with $d_{\cC}$-connection \eqref{p1-imdpa19a}, 
where the $d_{\cC}$-connections are defined as in \ref{p1-delta-con4}, 
$\Gamma^n(\cF^\vee)$ (resp.\ $\Gamma^{\leq n}(\cE^\vee)$) being endowed with the Higgs field $0$ (resp.\ $\gamma^n$).

\begin{prop}\label{p1-imdpa23}
Let $n$ be an integer $\geq 0$. Then, 
\begin{itemize}
\item[{\rm (i)}] The Higgs fields $\gamma^n$ \eqref{p1-imdpa22a} and $\gamma^{n}_\tot$ \eqref{p1-imdpa22b} 
are $\Gamma^{\leq n}(\cE^\vee)$-linear, for the natural actions of $\Gamma^{\leq n}(\cE^\vee)$.  
In particular, $\ker(\gamma^{n}_\tot)$ is a sub-$\Gamma^{\leq n}(\cE^\vee)$-module of $\rS^n(\cF)\otimes_A\Gamma^{\leq n}(\cE^\vee)$.
\item[{\rm (ii)}] The morphism \eqref{p1-imdpa22c} is $\Gamma^{\leq n}(\cE^\vee)$-linear when the left-hand (resp.\ right-hand) side is equipped
with the action of $\Gamma^{\leq n}(\cE^\vee)$ defined in \eqref{p1-imdpa5b} (resp.\ in {\rm (i)}). 
\item[{\rm (iii)}] If $p$ is not a zero divisor in $A$, the morphism \eqref{p1-imdpa22c} is an isomorphism.
\item[{\rm (iv)}] The morphism \eqref{p1-imdpa22e} is a $\cC\otimes_A\Gamma^{\leq n}(\cE^\vee)$-linear isomorphism. 
\item[{\rm (v)}] The diagram
\begin{equation}
\xymatrix{
{\cC\otimes_A\Gamma^{n+1}(\cF^\vee)}\ar[r]\ar[d]&{\cC\otimes_A\Gamma^{\leq n+1}(\cE^\vee)}\ar[d]\\
{\cC\otimes_A\Gamma^n(\cF^\vee)}\ar[r]&{\cC\otimes_A\Gamma^{\leq n}(\cE^\vee),}}
\end{equation}
where the horizontal arrows are the isomorphisms \eqref{p1-imdpa22e}, the left (resp.\ right) vertical arrow is 
induced by \eqref{p1-imdpa4b} (resp.\ is the canonical morphism), is commutative. 
\end{itemize} 
\end{prop}

(i) It follows immediately from the definitions since the $A$-algebra $\Gamma^{\leq n}(\cE^\vee)$ is commutative. 

(ii) It follows from the commutativity of \eqref{p1-imdpa22d} 
since the vertical arrows and the lower horizontal arrow are $\Gamma^{\leq n}(\cE^\vee)$-linear, where 
$\Gamma^{\leq n}(\cE^\vee)$ acts on $\cHom_{\Gamma^{\leq n}(\cE^\vee)}(\Gamma^n(\cF'^\vee),\Gamma^{\leq n}(\cE^\vee))$
and $\cHom_{\rS(\cE^\vee)}(\Gamma^n(\cF'^\vee),\Gamma^{\leq n}(\cE^\vee))$ by acting only on the targets. 

(iii) Indeed, the lower horizontal arrow of \eqref{p1-imdpa22d} is an isomorphism since 
$p$ is not a zero divisor in the locally free $A$-module $\Gamma^{\leq n}(\cE^\vee)$,
and every local section of the cokernel of the canonical homomorphism 
$\varsigma^{\leq n} \colon \rS(\cE^\vee)\rightarrow \Gamma^{\leq n}(\cE^\vee)$ is annihilated by a power of $p$. 

(iv) Indeed, the morphism \eqref{p1-imdpa22e} is $\cC\otimes_A\Gamma^{\leq n}(\cE^\vee)$-linear by (ii). 
Observe that the kernel of the canonical homomorphism 
$\varpi^n\colon \Gamma^{\leq n}(\cE^\vee)\rightarrow \Gamma^{\leq 0}(\cE^\vee)=A$ is nilpotent.  
By (v) below, the $\varpi^n$-scalar extension of \eqref{p1-imdpa22e}, 
\begin{equation}
(\cC\otimes_A\Gamma^n(\cF^\vee))\otimes_{\Gamma^{\leq n}(\cE^\vee)}A\rightarrow 
(\cC\otimes_A\Gamma^{\leq n}(\cE^\vee))\otimes_{\Gamma^{\leq n}(\cE^\vee)}A,
\end{equation}
coincides with the morphism \eqref{p1-imdpa22e} for $n=0$, which is clearly an isomorphism. 
Hence, the morphism \eqref{p1-imdpa22e} is surjective. 
Since the left and right-hand sides are invertible $\cC\otimes_A\Gamma^{\leq n}(\cE^\vee)$-modules by \ref{p1-imdpa9}(i),
\eqref{p1-imdpa22e} is an isomorphism. 

(v) The diagram 
\begin{equation}
\xymatrix{
{\Gamma^{n+1}(\cF^\vee)}\ar[dd]\ar[r]&{\cHom_{\Gamma^{\leq n+1}(\cE^\vee)}(\Gamma^{n+1}(\cF'^\vee),\Gamma^{\leq n+1}(\cE^\vee))}\ar[d]\\
&{\cHom_{\Gamma^{\leq n+1}(\cE^\vee)}(\Gamma^{n+1}(\cF'^\vee),\Gamma^{\leq n}(\cE^\vee))}\\
{\Gamma^{n}(\cF^\vee)}\ar[r]&{\cHom_{\Gamma^{\leq n}(\cE^\vee)}(\Gamma^{n}(\cF'^\vee),\Gamma^{\leq n}(\cE^\vee)),}\ar[u]}
\end{equation}
where the horizontal arrows are the isomorphisms induced by \eqref{p1-imdpa17b}, is commutative. On the other hand, by \ref{p1-imdpa21}(iv),
the diagram 
\begin{equation}
\xymatrix{
{(\rS^{n+1}(\cF)\otimes_A\Gamma^{\leq n+1}(\cE^\vee))^{\gamma^{n+1}_\tot=0}}\ar[r]\ar[d]&
{\cHom_{\rS(\cE^\vee)}(\Gamma^{n+1}(\cF'^\vee),\Gamma^{\leq n+1}(\cE^\vee))}\ar[d]\\
{(\rS^{n+1}(\cF)\otimes_A\Gamma^{\leq n}(\cE^\vee))^{\gamma'^n_\tot=0}}\ar[r]&
{\cHom_{\rS(\cE^\vee)}(\Gamma^{n+1}(\cF'^\vee),\Gamma^{\leq n}(\cE^\vee))}\\
{(\rS^n(\cF)\otimes_A\Gamma^{\leq n}(\cE^\vee))^{\gamma^n_\tot=0}}\ar[r]\ar[u]&
{\cHom_{\rS(\cE^\vee)}(\Gamma^n(\cF'^\vee),\Gamma^{\leq n}(\cE^\vee)),}\ar[u]}
\end{equation}
where $\gamma'^n_\tot=\theta^{n+1}\otimes \id+\id\otimes \gamma^n$ and
the horizontal arrows are the isomorphisms \eqref{p1-imdpa21d}, is commutative. 
We deduce that the diagram 
\begin{equation}
\xymatrix{
{\Gamma^{n+1}(\cF^\vee)}\ar[dd]\ar[r]&{(\rS^{n+1}(\cF)\otimes_A\Gamma^{\leq n+1}(\cE^\vee))^{\gamma^{n+1}_\tot=0}}\ar[d]\\
&{(\rS^{n+1}(\cF)\otimes_A\Gamma^{\leq n}(\cE^\vee))^{\gamma'^n_\tot=0}}\\
{\Gamma^{n}(\cF^\vee)}\ar[r]&{(\rS^{n}(\cF)\otimes_A\Gamma^{\leq n}(\cE^\vee))^{\gamma^n_\tot=0},}\ar[u]}
\end{equation}
where the horizontal arrows are the morphisms \eqref{p1-imdpa22c}, is commutative. The proposition follows. 

\subsection{}\label{p1-imdpa24}
We set 
\begin{eqnarray}
\cC\hotimes_A\hGamma(\cE^\vee)&=&\underset{\underset{n\geq 0}{\longleftarrow}}\lim\ \cC\otimes_A\Gamma^{\leq n}(\cE^\vee),\\
\cC\hotimes_A\cV&=&\underset{\underset{n\geq 0}{\longleftarrow}}\lim\ \cC\otimes_A\Gamma^{n}(\cF^\vee),
\end{eqnarray}
where the transition morphisms of the second inverse limit are induced by \eqref{p1-imdpa4b}. 
By \ref{p1-imdpa6} and \ref{p1-imdpa9}, $\cC\hotimes_A\cV$ has a canonical structure of an invertible $\cC\otimes_A\hGamma(\cE^\vee)$-module. 
By \ref{p1-imdpa23}, the isomorphisms \eqref{p1-imdpa22e} induce a $\cC\hotimes_A\hGamma(\cE^\vee)$-linear isomorphism 
\begin{equation}
\cC\hotimes_A\cV\stackrel{\sim}{\rightarrow} \cC\hotimes_A\hGamma(\cE^\vee).
\end{equation}

\subsection{}\label{p1-imdpa25}
Let $\rho\colon \cF\rightarrow A$ be an $A$-linear form such that 
$\rho\circ c=\id_A$, i.e., a section of $\Gamma(X,\Lambda)$ \eqref{p1-imdpa10}. 
Then $\rho$ extends uniquely to a homomorphism of $A$-algebras $\varrho\colon \cC\rightarrow A$ \eqref{p1-imdpa1d}. 
By \ref{p1-imdpa23}(iv), for every integer $n\geq 0$, the $\cC\otimes_A\Gamma^{\leq n}(\cE^\vee)$-linear isomorphism \eqref{p1-imdpa22e} 
induces by extension of scalars by $\varrho$ a $\Gamma^{\leq n}(\cE^\vee)$-linear isomorphism
\begin{equation}\label{p1-imdpa25a}
\ttt_{\rho,n}\colon \Gamma^n(\cF^\vee)\stackrel{\sim}{\rightarrow} \Gamma^{\leq n}(\cE^\vee).
\end{equation}
By \ref{p1-imdpa23}(v), the latter induces a $\hGamma(\cE^\vee)$-linear isomorphism
\begin{equation}\label{p1-imdpa25b}
\ttt_{\rho}\colon \cV\stackrel{\sim}{\rightarrow} \hGamma(\cE^\vee).
\end{equation}

\begin{prop}\label{p1-imdpa26}
Under the assumption of \ref{p1-imdpa25} and with the same notation, we have 
\begin{equation}
\ttt_{\rho}(\exp_\Lambda(\rho))=1,
\end{equation}
where $\exp_\Lambda$ is the morphism defined in \eqref{p1-imdpa13a}. 
\end{prop}

Indeed, for every integer $n\geq 0$, the $A$-linear form $\varrho|\rS^n(\cF) \colon \rS^n(\cF) \rightarrow A$ corresponds by the duality 
isomorphism \eqref{p1-imdpa11a} to the section $\rho^{[n]}$ of $\Gamma^n(\cF^\vee)$. Let $\rho'=-\rho\colon \cF'\rightarrow A$ \eqref{p1-imdpa15a}. 
In view of the definition \eqref{p1-imdpa22d} of the morphism \eqref{p1-imdpa22c}, we have by \eqref{p1-imdpa17c}
\begin{equation}
\ttt_{\rho,n}(\rho^{[n]})= \langle \rho^{[n]}, \rho'^{[n]}\rangle = \cphi^n(\rho^{[n]})\cphi'^n(\rho'^{[n]})=1. 
\end{equation}
The proposition follows since $\exp_\Lambda(\rho)=\rho^{[\infty]}$, see the proof of \ref{p1-imdpa13}.

\section{Revisiting deformation theory}\label{p1-rdt}

\subsection{}\label{p1-rdt1}
In this section, we consider a commutative diagram of the category of fine logarithmic schemes $\FLS$ \eqref{p1-NC1}
\begin{equation}\label{p1-rdt1a}
\xymatrix{
{(Y,\cM_{Y})}\ar[r]^-(0.5){j}\ar[d]^{h}\ar@/_2pc/[dd]_{g}&{(\tY,\cM_{\tY})}\ar@/^2pc/[dd]^{\tg}&\\
{(X,\cM_{X})}\ar[r]^-(0.5){i}\ar[d]^f\ar@{}[rd]|{\Box}&{(\tX,\cM_{\tX})}\ar[d]_{\tf}\\
{(S,\cM_{S})}\ar[r]^-(0.5){\iota}&{(\tS,\cM_{\tS})}}
\end{equation}
satisfying the following conditions:
\begin{itemize}
\item[(i)] The morphism $\iota$ is a thickening of order one, i.e., a strict closed immersion  of fine logarithmic schemes 
such that the underlying closed immersion of schemes $S\rightarrow \tS$ is defined by an ideal $I$ of $\co_\tS$ of square zero.
\item[(ii)] The morphism $\tf$ is smooth.
\item[(iii)] The lower square and the external rectangle of diagram \eqref{p1-rdt1a} are Cartesian, and $g=f\circ h$. 
Then, by \ref{p1-NC3}, $i$ (resp.\ $j$) is a thickening of order one and the underlying closed immersion of schemes $X\rightarrow \tX$
(resp.\ $Y\rightarrow \tY$) is defined by the ideal $I_X=I\co_\tX$ (resp.\ $I_Y=I\co_{\tY}$).
\item[(iv)] The $\co_{Y}$-module $I_Y$ is invertible. 
\end{itemize}

To lighten the notation, we set 
\begin{equation}\label{p1-rdt1c}
\tOmega^1_{X/S}=\Omega^1_{(X,\cM_X)/(S,\cM_S)}. 
\end{equation}

\subsection{}\label{p1-rdt2}
We denote by $\cL_{\tY/\tX}$ the étale sheaf on $Y$ defined for any étale $Y$-scheme $U$ by setting 
$\cL_{\tY/\tX}(U)$ to be the set of liftings of $h_{|U}$ \eqref{p1-rdt1a}, i.e., the set of $(\tS,\cM_{\tS})$-morphisms 
\begin{equation}\label{p1-rdt2a}
\gamma\colon (\tU,\cM_{\tY}|\tU)\rightarrow (\tX,\cM_{\tX})
\end{equation}
such that $\gamma\circ \jmath_U=i\circ h_{|U}$,
where $\tU$ is the essentially unique étale $\tY$-scheme equipped with 
a $Y$-isomorphism $U\stackrel{\sim}{\rightarrow} Y\times_{\tY}\tU$ (\cite{ega4} 18.1.2), and $\jmath_U\colon U\rightarrow \tU$ 
is the morphism induced by the canonical projection. 
By (\cite{kato1} 3.9 or \cite{ogus} IV 2.2.2), $\cL_{\tY/\tX}$ is a torsor of $Y_\et$ under the $\co_Y$-module 
\begin{equation}\label{p1-rdt2b}
\cHom_{\co_Y}(h^*(\tOmega^1_{X/S}),I_Y).
\end{equation}
We call it the {\em torsor of liftings of $h$ to $\tY$ over $\tX$}. 

By the local-global spectral sequence of $\Ext$, we have a canonical isomorphism
\begin{equation}\label{p1-rdt2c}
\rH^1(Y,\cHom_{\co_Y}(h^*(\tOmega^1_{X/S}),I_Y)) \stackrel{\sim}{\rightarrow} \Ext^1_{\co_Y}(Y;h^*(\tOmega^1_{X/S}),I_Y),
\end{equation}
where the cohomology groups are indifferently computed for the étale or the Zariski topologies. 

The {\em obstruction to lift $g$ to $\tY$ over $\tX$} \eqref{p1-rdt1a}, denoted by
\begin{equation}\label{p1-rdt2d}
\ro_{\tY/\tX}(h)\in \Ext^1_{\co_Y}(Y;h^*(\tOmega^1_{X/S}),I_Y),
\end{equation}
is the image of the class of the torsor $\cL_{\tY/\tX}$ by the isomorphism \eqref{p1-rdt2c}. 

\subsection{}\label{p1-rdt3}
We identify the dual $\co_Y$-module of \eqref{p1-rdt2b} with $I_Y^{-1}\otimes h^*(\tOmega^1_{X/S})$. 
We denote by $\cF_{\tY/\tX}$ the sheaf of affine functions on the restriction of the torsor $\cL_{\tY/\tX}$ to $Y_\zar$ \eqref{p1-prem1}. 
It is an $\co_Y$-module that fits into a canonical exact sequence
\begin{equation}\label{p1-rdt3a}
0\rightarrow \co_Y\rightarrow \cF_{\tY/\tX}\rightarrow I_Y^{-1}\otimes_{\co_Y} h^*(\tOmega^1_{X/S}) \rightarrow 0.
\end{equation}
We denote by 
\begin{equation}\label{p1-rdt3c}
\cC_{\tY/\tX}=\underset{\underset{n\geq 0}{\longrightarrow}}\lim\ \rS^n_{\co_Y}(\cF_{\tY/\tX})
\end{equation}
the associated $\co_Y$-algebra \eqref{p1-prem1c}. Since $\Spec(\cC_{\tY/\tX})$ canonically represents the restriction of the torsor $\cL_{\tY/\tX}$ to $Y_\zar$ \eqref{p1-prem1d}, 
a section of $\cL_{\tY/\tX}$ over a Zariski open $U$ of $Y$ corresponds to a splitting of the restriction of the exact sequence \eqref{p1-rdt3a} over $U$. 
We easily deduce that the 
class of the extension $I_Y\otimes_{\co_Y} \cF_{\tY/\tX}$ in $\Ext^1_{\co_Y}(Y;h^*(\tOmega^1_{X/S}),I_Y)$ is $\ro_{\tY/\tX}(h)$ (see the proof of \cite{agt} II.4.10).

\begin{defi}\label{p1-rdt46}
The torsor $\cL_{\tY/\tX}$ is called the {\em torsor of local liftings of the morphism $h$ to $\tY$ over $\tX$}.
The $\co_Y$-module (resp.\ $\co_Y$-algebra) $\cF_{\tY/\tX}$ (resp.\ $\cC_{\tY/\tX}$)
is called the {\em Higgs--Tate extension} (resp.\ {\em Higgs--Tate algebra}) of $\cL_{\tY/\tX}$ or of $\tY$ over $\tX$ (resp.\ or associated with $\cF_{\tY/\tX}$). 
\end{defi}

\begin{lem}\label{p1-rdt5}
If the morphisms $\tg$ and $h$ are smooth \eqref{p1-rdt1a}, then any lifting
\begin{equation}\label{p1-rdt5a}
\thh\colon (\tY,\cM_{\tY})\rightarrow (\tX,\cM_{\tX})
\end{equation}
of $h$ is smooth. 
\end{lem}

Indeed, by (\cite{kato1} 3.12 or \cite{ogus} IV 3.2.3), it suffices to prove that the canonical morphism
\begin{equation}\label{p1-rdt5b}
\thh^*(\Omega^1_{(\tX,\cM_\tX)/(\tS,\cM_{\tS})})\rightarrow \Omega^1_{(\tY,\cM_{\tY})/(\tS,\cM_{\tS})}
\end{equation}
is locally left invertible, i.e., that it induces locally an isomorphism from the source into a direct factor of the target. 
It amounts to say that for every point $y$ of $\tY$ (or equivalently of $Y$), of image $x$ in $\tX$, the canonical morphism 
\begin{equation}\label{p1-rdt5c}
\thh^*(\Omega^1_{(\tX,\cM_\tX)/(\tS,\cM_{\tS})})_y\rightarrow \Omega^1_{(\tY,\cM_{\tY})/(\tS,\cM_{\tS}),y}
\end{equation}
is left invertible. Since the source and the target are free $\co_{\tY,y}$-modules of finite type, 
the last condition is equivalent to the fact that the canonical morphism 
\begin{equation}\label{p1-rdt5d}
\Omega^1_{(\tX,\cM_\tX)/(\tS,\cM_{\tS}),x}\otimes_{\co_{\tX,x}}\kappa(y)\rightarrow \Omega^1_{(\tY,\cM_{\tY})/(\tS,\cM_{\tS}),y}\otimes_{\co_{\tY,y}}\kappa(y)
\end{equation}
is injective by (\cite{ega4} 0.19.1.12). By \ref{p1-rdt1}(iii), the latter is identified with the morphism
\begin{equation}\label{p1-rdt5e}
\Omega^1_{(X,\cM_X)/(S,\cM_{S}),x}\otimes_{\co_{X,x}}\kappa(y)\rightarrow \Omega^1_{(Y,\cM_Y)/(S,\cM_{S}),y}\otimes_{\co_{Y,y}}\kappa(y)
\end{equation}
induced by the smooth morphism $h$. It is therefore injective; the proposition follows.

\subsection{}\label{p1-rdt4}
Consider a Cartesian diagram of $\FLS$
\begin{equation}\label{p1-rdt4a}
\xymatrix{
{(Y',\cM_{Y'})}\ar[r]^-(0.5){j'}\ar[d]_\phi\ar@{}[rd]|\Box&{(\tY',\cM_{\tY'})}\ar[d]^{\tphi}\\
{(Y,\cM_{Y})}\ar[r]^-(0.5){j}&{(\tY,\cM_{\tY}).}}
\end{equation}
We set $h'=h\circ \phi\colon (Y',\cM_{Y'})\rightarrow (X,\cM_X)$. 
Then, $j'$ is a thickening of order one and the underlying closed immersion of schemes $Y'\rightarrow \tY'$ is defined by the ideal $I_{Y'}=I\co_{\tY'}$ of $\co_{\tY'}$ \eqref{p1-NC3}. 
We assume further that the $\co_{Y'}$-module $I_{Y'}$ is invertible. Thus, the canonical surjective morphism $\lambda\colon \phi^*(I_Y)\rightarrow I_{Y'}$ is an isomorphism.
We denote by   
\begin{equation}\label{p1-rdt4b}
\uu\colon \cHom_{\co_Y}(h^*(\tOmega^1_{X/S}),I_Y)\rightarrow \phi_*(\cHom_{\co_{Y'}}(h'^*(\tOmega^1_{X/S}),I_{Y'}))
\end{equation}
the $\co_Y$-linear morphism induced by $\lambda$, by 
\begin{equation}\label{p1-rdt4f}
\ou\colon \phi^{-1}(\cHom_{\co_Y}(h^*(\tOmega^1_{X/S}),I_Y))\rightarrow \cHom_{\co_{Y'}}(h'^*(\tOmega^1_{X/S}),I_{Y'})
\end{equation}
its adjoint as a morphism of abelian sheaves, and by
\begin{equation}\label{p1-rdt4g}
u\colon \phi^*(\cHom_{\co_Y}(h^*(\tOmega^1_{X/S}),I_Y))\rightarrow \cHom_{\co_{Y'}}(h'^*(\tOmega^1_{X/S}),I_{Y'})
\end{equation}
the $\co_{Y'}$-linear adjoint of $\uu$ as a morphism of $\co_Y$-modules, which is an isomorphism.

Let  $\cL_{\tY'/\tX}$ be the torsor of liftings of $h'$ to $\tY'$ over $\tX$, which is a torsor of $Y'_\et$ under the $\co_{Y'}$-module $\cHom_{\co_{Y'}}(h'^*(\tOmega^1_{X/S}),I_{Y'})$
and let $\phi^+(\cL_{\tY/\tX})$ be the affine pullback of $\cL_{\tY/\tX}$ by $\phi$ \eqref{p1-NC5}, which is a torsor of $Y'_\et$ under the $\co_{Y'}$-module 
$\phi^*(\cHom_{\co_{Y'}}(h^*(\tOmega^1_{X/S}),I_Y))$. 

We have a canonical morphism of sheaves 
\begin{equation}\label{p1-rdt4c}
\uupmu\colon \cL_{\tY/\tX}\rightarrow \phi_*(\cL_{\tY'/\tX})
\end{equation}
defined for any $U\in \ob(\Et_{/Y})$, by sending a section $\sigma\in \cL_{\tY/\tX}(U)$ to the section of $\cL_{\tY'/\tX}(U\times_YY')$ 
corresponding to the composed morphism 
\begin{equation}\label{p1-rdt4d}
\xymatrix{(\tU',\cM_{\tY'}|\tU')\ar[r]^{\tphi_{\tU}}&(\tU,\cM_{\tY}|\tU)\ar[r]^\sigma& (\tX,\cM_{\tX}),}
\end{equation}
where $\tU'=\tU\times_{\tY}\tY'$ and $\tphi_{\tU}$ is the canonical projection. We see immediately that $\uupmu$ is $\uu$-equivariant \eqref{p1-rdt4b}. The adjoint 
\begin{equation}\label{p1-rdt4e}
\oupmu\colon \phi^*(\cL_{\tY/\tX})\rightarrow \cL_{\tY'/\tX}
\end{equation}
of $\uupmu$ is $\ou$-equivariant \eqref{p1-rdt4f}. 
By (\cite{agt} II.4.14), it induces a $u$-equivariant isomorphism \eqref{p1-rdt4g}
\begin{equation}\label{p1-rdt4h}
\upmu\colon \phi^+(\cL_{\tY/\tX})\stackrel{\sim}{\rightarrow} \cL_{\tY'/\tX}. 
\end{equation}
We deduce an $\co_{Y'}$-linear isomorphism \eqref{p1-rdt3a}
\begin{equation}\label{p1-rdt4i}
\cF_{\tY'/\tX}\stackrel{\sim}{\rightarrow} \phi^*(\cF_{\tY/\tX})
\end{equation}
that fits into a commutative diagram
\begin{equation}\label{p1-rdt4j}
\xymatrix{
0\ar[r]&{\co_{Y'}}\ar[r]\ar@{=}[d]&{\cF_{\tY'/\tX}}\ar[r]\ar[d]&{I_{Y'}^{-1}\otimes_{\co_Y} h'^*(\tOmega^1_{X/S})}\ar[r]\ar[d]&0\\
0\ar[r]&{\co_{Y'}}\ar[r]&{\phi^*(\cF_{\tY/\tX})}\ar[r]&{\phi^*(I_Y^{-1})\otimes_{\co_{Y'}} h'^*(\tOmega^1_{X/S})}\ar[r]&0,}
\end{equation}
where the right vertical arrow is induced by the canonical isomorphism $\lambda\colon \phi^*(I_Y)\stackrel{\sim}{\rightarrow} I_{Y'}$. 
We thus obtain an isomorphism of $\co_Y$-algebras \eqref{p1-rdt3c}
\begin{equation}\label{p1-rdt4k}
\cC_{\tY'/\tX}\stackrel{\sim}{\rightarrow} \phi^*(\cC_{\tY/\tX}).
\end{equation}

\subsection{}\label{p1-rdt40}
Consider a commutative diagram of $\FLS$ with Cartesian squares
\begin{equation}\label{p1-rdt40a}
\xymatrix{
{(Y'',\cM_{Y''})}\ar[r]^-(0.5){j''}\ar[d]_{\phi'}\ar@{}[rd]|\Box\ar@/_2.5pc/[dd]_{\phi''}&{(\tY'',\cM_{\tY''})}\ar[d]^{\tphi'}\ar@/^2.5pc/[dd]^{\tphi''}\\
{(Y',\cM_{Y'})}\ar[r]^-(0.5){j'}\ar[d]_\phi\ar@{}[rd]|\Box&{(\tY',\cM_{\tY'})}\ar[d]^{\tphi}\\
{(Y,\cM_{Y})}\ar[r]^-(0.5){j}&{(\tY,\cM_{\tY}).}}
\end{equation}
Then, the diagram 
\begin{equation}\label{p1-rdt40b}
\xymatrix{
{\phi'^+(\phi^+(\cL_{\tY/\tX}))}\ar[rr]^-(0.5){\phi'^+(\upmu)}\ar[d]_c&&{\phi'^+(\cL_{\tY'/\tX})}\ar[d]^-(0.5){\upmu'}\\
{\phi''^+(\cL_{\tY/\tX})}\ar[rr]^-(0.5){\upmu''}&&{\cL_{\tY''/\tX},}}
\end{equation}
where $\upmu$, $\upmu'$ and $\upmu''$ are the isomorphisms \eqref{p1-rdt4h} and $c$ is the canonical isomorphism, is commutative. 

\subsection{}\label{p1-rdt6}
We consider a commutative diagram of $\FLS$
\begin{equation}\label{p1-rdt6a}
\xymatrix{
{(Y,\cM_{Y})}\ar[r]^-(0.5){j}\ar[d]^{h'}\ar@/_4pc/[ddd]_{g}&{(\tY,\cM_{\tY})}\ar@/^4pc/[ddd]^\tg&\\
{(X',\cM_{X'})}\ar[r]^-(0.5){i'}\ar[d]^{\gamma}\ar@/_2pc/[dd]_{f'}&{(\tX',\cM_{\tX'})}\ar@/^2pc/[dd]^{\tf'}&\\
{(X,\cM_{X})}\ar[r]^-(0.5){i}\ar[d]^f&{(\tX,\cM_{\tX})}\ar[d]_{\tf}\\
{(S,\cM_{S})}\ar[r]^-(0.5){\iota}&{(\tS,\cM_{\tS})}}
\end{equation}
extending diagram \eqref{p1-rdt1a} and satisfying the following conditions:
\begin{itemize}
\item[(v)] The morphism $\tf'$ is smooth. 
\item[(vi)] The morphisms $(f',\tf',\iota,i')$ define a Cartesian diagram, and we have $f'=f\circ \gamma$ and $h=\gamma\circ h'$.
Then, $i'$ is a thickening of order one and the underlying closed immersion of schemes $X'\rightarrow \tX'$ is defined by the ideal $I_{X'}=I\co_{\tX'}$ \eqref{p1-NC3}.
\item[(vii)] The $\co_{X'}$-module $I_{X'}$ is invertible. 
\end{itemize}
To lighten the notation, we set 
\begin{equation}
\tOmega^1_{X'/S}=\Omega^1_{(X',\cM_{X'})/(S,\cM_S)}. 
\end{equation}

Let  $\cL_{\tY/\tX'}$ be the torsor of liftings of $h'$ to $\tY$ over $\tX'$, 
which is a torsor of $Y_\et$ under the $\co_Y$-module $\cHom_{\co_Y}(h'^*(\tOmega^1_{X'/S}),I_Y)$  \eqref{p1-rdt2}.
Let  $\cL_{\tX'/\tX}$ be the torsor of liftings of $\gamma$ to $\tX'$ over $\tX$, 
which is a torsor of $X'_\et$ under the $\co_{X'}$-module $\cHom_{\co_{X'}}(\gamma^*(\tOmega^1_{X/S}),I_{X'})$,
and let $h'^+(\cL_{\tX'/\tX})$ be its affine pullback by $h'$ \eqref{p1-NC5}, which is a torsor of $Y_\et$ under the $\co_Y$-module 
$h'^*(\cHom_{\co_{X'}}(\gamma^*(\tOmega^1_{X/S}),I_{X'}))$.

Any section $\thh'\in \cL_{\tY/\tX'}(Y)$ determines a surjective, and hence bijective $\co_Y$-linear morphism $\lambda_{\thh'}\colon h'^*(I_{X'})\rightarrow I_Y$. 
The latter does not actually depend on $\thh'$ because $I_{X'}=I\co_{\tX'}$ and $II_Y=0$. In fact, $\lambda_{\thh'}$ is characterized by the commutative diagram 
\begin{equation}
\xymatrix{
{h'^*(f'^*(I))}\ar[r]^-(0.5)\sim\ar@{->>}[d]&{g^*(I)}\ar@{->>}[d]\\
{h'^*(I_{X'})}\ar[r]^-(0.5){\lambda_{\thh'}}&{I_Y.}}
\end{equation}

Since the morphism $\tf'$ is smooth, for every affine open covering 
$(V_\ell)_{1\leq \ell\leq n}$ of $Y$, there exist sections $\thh'_\ell\in \cL_{\tY/\tX'}(V_\ell)$ for all $1\leq \ell\leq n$. 
The induced isomorphisms $\lambda_{\ell}\colon h'^*(I_{X'})|V_\ell\stackrel{\sim}{\rightarrow} I_Y|V_\ell$, for $1\leq \ell\leq n$, 
glue to define an $\co_Y$-linear isomorphism 
\begin{equation}\label{p1-rdt6b}
\lambda\colon h'^*(I_{X'})\stackrel{\sim}{\rightarrow} I_Y.
\end{equation}

We denote by  
\begin{equation}\label{p1-rdt6d}
v\colon \cHom_{\co_Y}(h'^*(\tOmega^1_{X'/S}),I_Y)\rightarrow  \cHom_{\co_Y}(h^*(\tOmega^1_{X/S}),I_Y)
\end{equation}
the canonical morphism, by 
\begin{equation}\label{p1-rdt6c}
\uu\colon \cHom_{\co_{X'}}(\gamma^*(\tOmega^1_{X/S}),I_{X'})\rightarrow h'_*(\cHom_{\co_Y}(h^*(\tOmega^1_{X/S}),I_Y))
\end{equation}
the $\co_{X'}$-linear morphism induced by $\lambda$, by 
\begin{equation}\label{p1-rdt6e}
\ou\colon h'^{-1}(\cHom_{\co_{X'}}(\gamma^*(\tOmega^1_{X/S}),I_{X'}))\rightarrow \cHom_{\co_Y}(h^*(\tOmega^1_{X/S}),I_Y)
\end{equation}
its adjoint as a morphism of abelian sheaves, and by
\begin{equation}\label{p1-rdt6f}
u\colon h'^*(\cHom_{\co_{X'}}(\gamma^*(\tOmega^1_{X/S}),I_{X'}))\rightarrow \cHom_{\co_Y}(h^*(\tOmega^1_{X/S}),I_Y)
\end{equation}
the $\co_Y$-linear adjoint of $\uu$ as a morphism of $\co_{X'}$-modules, which is an isomorphism.

We denote by $\cHom_{v}(\cL_{\tY/\tX'},\cL_{\tY/\tX})$ the sheaf of $v$-equivariant morphisms from $\cL_{\tY/\tX'}$ to $\cL_{\tY/\tX}$, 
which is canonically a torsor of $Y_\et$ under the $\co_Y$-module $\cHom_{\co_Y}(h^*(\tOmega^1_{X/S}),I_Y)$ \eqref{p1-prem9}. We have a canonical morphism 
\begin{equation}\label{p1-rdt6g}
\uupnu\colon \cL_{\tX'/\tX}\rightarrow h'_*(\cHom_{v}(\cL_{\tY/\tX'},\cL_{\tY/\tX}))
\end{equation}
defined for any $U'\in \ob(\Et_{/X'})$ by sending a section $\sigma\in \cL_{\tX'/\tX}(U')$ to the $v$-equivariant morphism 
\begin{equation}\label{p1-rdt6h}
\uupnu(\sigma)\colon \cL_{\tY/\tX'}|U'\times_{X'}Y\rightarrow \cL_{\tY/\tX}|U'\times_{X'}Y,
\end{equation}
defined for any $V\in \ob(\Et_{/U'\times_{X'}Y})$ by sending a section $\tau\in \cL_{\tY/\tX'}(V)$ 
to the section of $\cL_{\tY/\tX}(V)$ corresponding to the composed morphism 
\begin{equation}\label{p1-rdt6i}
\xymatrix{(\tV,\cM_{\tY}|\tV)\ar[r]^-(0.5){\tau'}&(\tU',\cM_{\tX'}|\tU')\ar[r]^-(0.5){\sigma}&(\tX,\cM_{\tX}),}
\end{equation}
where $\tV$ is the essentially unique étale $\tY$-scheme equipped with 
a $Y$-isomorphism $V\stackrel{\sim}{\rightarrow} Y\times_{\tY}\tV$ and 
$\tau'$ is the unique lifting of $\tau$ that fits into the commutative diagram
\begin{equation}\label{p1-rdt6j}
\xymatrix{
{(V,\cM_Y|V)}\ar[r]\ar[d]&{(\tV,\cM_\tY|\tV)}\ar[d]_{\tau'}\ar@/^3pc/[dd]^{\tau}\\
{(U',\cM_{X'}|U')}\ar[r]\ar[d]&{(\tU',\cM_{\tX'}|\tU')}\ar[d]\\
{(X',\cM_{X'})}\ar[r]&{(\tX',\cM_{\tX'}).}}
\end{equation} 
We immediately see that $\uupnu(\sigma)$ is $v$-equivariant \eqref{p1-rdt6d} and that $\uupnu$ is $\uu$-equivariant \eqref{p1-rdt6c}. 

The adjoint 
\begin{equation}\label{p1-rdt6k}
\oupnu\colon h'^*(\cL_{\tX'/\tX})\rightarrow \cHom_{v}(\cL_{\tY/\tX'},\cL_{\tY/\tX})
\end{equation}
of $\uupnu$ is $\ou$-equivariant \eqref{p1-rdt6e}. By (\cite{agt} II.4.14), it induces a $u$-equivariant isomorphism \eqref{p1-rdt6f}
\begin{equation}\label{p1-rdt6l}
\upnu\colon h'^+(\cL_{\tX'/\tX})\stackrel{\sim}{\rightarrow} \cHom_{v}(\cL_{\tY/\tX'},\cL_{\tY/\tX}).
\end{equation}
By \eqref{p1-prem9f}, we deduce a canonical $\co_Y$-linear morphism
\begin{equation}\label{p1-rdt6m}
\upvarphi\colon \cF_{\tY/\tX}\rightarrow h'^*(\cF_{\tX'/\tX}) \otimes_{\co_Y}\cF_{\tY/\tX'}.
\end{equation}
By \eqref{p1-prem9g}, we also deduce a canonical morphism of $\co_Y$-algebras 
\begin{equation}\label{p1-rdt6n}
\upphi\colon \cC_{\tY/\tX}\rightarrow h'^*(\cC_{\tX'/\tX}) \otimes_{\co_Y}\cC_{\tY/\tX'},
\end{equation}
compatible with $\upvarphi$, which we call the {\em composition morphism of Higgs--Tate algebras}.

\begin{remas}\label{p1-rdt8}
We keep the notation and assumptions of \ref{p1-rdt6}.
\begin{itemize}
\item[(i)] For every $\tgamma\in \cL_{\tX'/\tX}(X')$, the $v$-equivariant morphism \eqref{p1-rdt6g}
\begin{equation}\label{p1-rdt8a}
\uupnu(\tgamma)\colon \cL_{\tY/\tX'}\rightarrow \cL_{\tY/\tX}
\end{equation}
induces by pullback an $\co_Y$-linear morphism
\begin{equation}\label{p1-rdt8d}
\uupnu(\tgamma)^*\colon \cF_{\tY/\tX}\rightarrow \cF_{\tY/\tX'}
\end{equation}
that fits into a commutative diagram
\begin{equation}\label{p1-rdt8e}
\xymatrix{
0\ar[r]&{\co_Y}\ar[r]\ar@{=}[d]&{\cF_{\tY/\tX}}\ar[r]\ar[d]^{\uupnu(\tgamma)^*}&{I^{-1}_Y\otimes_{\co_Y} h^*(\tOmega^1_{X/S})}\ar[r]\ar[d]&0\\
0\ar[r]&{\co_Y}\ar[r]&{\cF_{\tY/\tX'}}\ar[r]&{I^{-1}_Y\otimes_{\co_Y} h'^*(\tOmega^1_{X'/S})}\ar[r]&0,}
\end{equation}
where the right vertical arrow is the canonical morphism. 
By \ref{p1-prem13}, if we denote by $\rho_\tgamma \colon \cF_{\tX'/\tX}\rightarrow \co_{X'}$ the splitting of the canonical extension associated with $\cF_{\tX'/\tX}$, defined by $\tgamma$, 
we have 
\begin{equation}
\uupnu(\tgamma)^*=(h'^*(\rho_\tgamma)\otimes \id_{\cF_{\tY/\tX'}})\circ \upvarphi.
\end{equation}
\item[(ii)] Given $\thh' \in \cL_{\tY/\tX'}(Y)$, the morphisms 
\begin{eqnarray}
\uupmu\colon \cL_{\tX'/\tX}\rightarrow h'_*(\cL_{\tY/\tX}),\\
\upmu\colon h'^+(\cL_{\tX'/\tX})\rightarrow \cL_{\tY/\tX},
\end{eqnarray}
defined in \eqref{p1-rdt4c} and \eqref{p1-rdt4h} are nothing but the evaluation of the morphisms $\uupnu$ \eqref{p1-rdt6g} and $\upnu$ \eqref{p1-rdt6l} at $\thh'$, respectively. 
Therefore, by \ref{p1-prem11}, if we denote by $\rho_{\thh'} \colon \cF_{\tY/\tX'}\rightarrow \co_Y$ the splitting of the canonical extension associated with $\cF_{\tY/\tX'}$, 
defined by $\thh'$, the composed morphism 
\begin{equation}
(\id_{h'^*(\cF_{\tX'/\tX})}\otimes \rho_{\thh'})\circ \upvarphi\colon \cF_{\tY/\tX}\rightarrow h'^*(\cF_{\tX'/\tX})
\end{equation}
is induced by $\upmu$. Observe that in \ref{p1-rdt4}, we do not assume $\tf'$ smooth, unlike in \ref{p1-rdt6}. 
\end{itemize}
\end{remas}

\subsection{}\label{p1-rdt9}
We take again the notation and assumptions of \ref{p1-rdt6}, moreover, we consider a Cartesian diagram of $\FLS$
\begin{equation}\label{p1-rdt9a}
\xymatrix{
{(Y',\cM_{Y'})}\ar[r]^-(0.5){j'}\ar[d]_\phi\ar@{}[rd]|\Box&{(\tY',\cM_{\tY'})}\ar[d]^{\tphi}\\
{(Y,\cM_{Y})}\ar[r]^-(0.5){j}&{(\tY,\cM_{\tY}).}}
\end{equation}
We set $\varphi'=h'\circ \phi\colon (Y',\cM_{Y'})\rightarrow (X',\cM_{X'})$ and $\varphi=\gamma\circ \varphi' \colon (Y',\cM_{Y'})\rightarrow (X,\cM_{X})$, and denote by  
\begin{equation}\label{p1-rdt9b}
v'\colon \cHom_{\co_{Y'}}(\varphi'^*(\tOmega^1_{X'/S}),I\co_{Y'})\rightarrow  \cHom_{\co_{Y'}}(\varphi^*(\tOmega^1_{X/S}),I\co_{Y'})
\end{equation}
the canonical morphism. Then, the diagram 
\begin{equation}\label{p1-rdt9c}
\xymatrix{
{\phi^+(h'^+(\cL_{\tX'/\tX}))}\ar[r]^-(0.5){\phi^+(\upnu)}\ar[dd]_a&{\phi^+(\cHom_{v}(\cL_{\tY/\tX'},\cL_{\tY/\tX}))}\ar[d]^b\\
&{\cHom_{v'}(\phi^+(\cL_{\tY/\tX'}),\phi^+(\cL_{\tY/\tX}))}\ar[d]^c\\
{\varphi'^+(\cL_{\tX'/\tX})}\ar[r]^-(0.5){\upnu'}&{\cHom_{v'}(\cL_{\tY'/\tX'},\cL_{\tY'/\tX}),}}
\end{equation}
where $\upnu$ and $\upnu'$ are the isomorphisms \eqref{p1-rdt6l}, $a$ and $b$ are the canonical isomorphisms and $c$ is the isomorphism induced by \eqref{p1-rdt4h},
is commutative.

\section{Revisiting continuous functors}\label{p2-cmt}

\begin{defi}\label{p2-cmt0}
Let $\cC, \cD$ be two categories, $\varphi \colon \cC\rightarrow \cD$ a functor. 
We say that $\varphi$ is a {\em $\mU$-functor} if for every $X\in \ob(\cC)$ and $Y\in \ob(\cD)$, $\Hom_{\cD}(\varphi(X),Y)$ is $\mU$-small \eqref{p2-ncgt4}. 
\end{defi}

If $\cD$ is a $\mU$-category (\cite{sga4} I 1.1), then every functor $\varphi \colon \cC\rightarrow \cD$ is a $\mU$-functor. 

Let $\psi\colon \cD\rightarrow \cD'$ be a fully faithful functor. Then $\varphi \colon \cC\rightarrow \cD$ is a $\mU$-functor if and only if so is 
$\psi\circ \varphi \colon \cC\rightarrow \cD'$. 

\subsection{}\label{p2-cmt1}
Let $\cC_1, \cC_2$ be two $\mU$-categories, $\hcC_i$ the category of presheaves of $\mU$-sets on $\cC_i$ \eqref{p2-ncgt4}, for $i=1,2$, 
$\varphi \colon \cC_1\rightarrow \hcC_2$ a $\mU$-functor \eqref{p2-cmt0}. 
In view of (\cite{sga4} I 1.3), we associate with $\varphi$ the functor 
\begin{equation}\label{p2-cmt1a}
\varphi_{\rp}\colon \hcC_2\rightarrow \hcC_1
\end{equation}
defined for any $G\in \ob(\hcC_2)$ and $X\in \ob(\cC_1)$ by 
\begin{equation}\label{p2-cmt1b}
\varphi_\rp(G)(X)=\Hom_{\hcC_2}(\varphi(X),G).
\end{equation} 

\begin{rema}\label{p2-cmt101}
In the situation of \ref{p2-cmt1}, there are two important cases where $\varphi \colon \cC_1\rightarrow \hcC_2$ is a $\mU$-functor:
\begin{itemize}
\item[{\rm (a)}] If the category $\cC_2$ is $\mU$-small, then $\hcC_2$ is a $\mU$-category (\cite{sga4} I 1.2) and hence $\varphi$ is $\mU$-functor. 
\item[{\rm (b)}] If $\varphi$ factors as $\varphi=\tth_2\circ \phi$, for a functor $\phi\colon \cC_1\rightarrow \cC_2$, 
where $\tth_2\colon \cC_2\rightarrow \hcC_2$ is the canonical functor \eqref{p2-ncgt4a}, then $\varphi$ is a $\mU$-functor. 
Indeed, $\tth_2$ is obviously a $\mU$-functor. 
In this case, the functor $\varphi_{\rp}$ \eqref{p2-cmt1a} is defined by composition with $\phi$, 
and is denoted by $\phi^*$ in (\cite{sga4} I 5.0). For coherence of notation, we will use the notation $\phi_\rp$ instead of $\phi^*$; so we have 
$\phi_{\rp}=\varphi_{\rp}$. 
\end{itemize}
\end{rema}

\subsection{}\label{p2-cmt100}
We keep the assumptions and notation of \ref{p2-cmt1}. We denote by $\mI_\varphi$ the category defined as follows. 
Objects of $\mI_\varphi$ are triples $(X,Y,f)$, where $X\in \ob(\cC_1)$, $Y\in \ob(\cC_2)$ and $f\in \varphi(X)(Y)$; 
such an object will also be denoted by $f\colon Y\rightarrow \varphi(X)$. 
Let $f\colon Y\rightarrow \varphi(X)$ and $f'\colon Y'\rightarrow \varphi(X')$ be two objects of $\mI_\varphi$. 
A morphism from $f\colon Y\rightarrow \varphi(X)$ to $f'\colon Y'\rightarrow \varphi(X')$ is a couple of morphisms 
$x\colon X\rightarrow X'$ of $\cC_1$ and $y\colon Y\rightarrow Y'$ of $\cC_2$ such that the diagram
\begin{equation}\label{p2-cmt100a}
\xymatrix{
Y\ar[r]^-(0.5)f\ar[d]_{y}&{\varphi(X)}\ar[d]^{\varphi(x)}\\
Y'\ar[r]^-(0.5){f'}&{\varphi(X')}}
\end{equation}
is commutative. We consider the functors 
\begin{eqnarray}
\ttb\colon 
\begin{array}[t]{clcr}
\mI_\varphi&\rightarrow& \cC_1,\\
(f\colon Y\rightarrow \varphi(X)) &\mapsto& X,
\end{array}\label{p2-cmt100b}\\
\tts\colon 
\begin{array}[t]{clcr}
\mI_\varphi&\rightarrow& \cC_2,\\
(f\colon Y\rightarrow \varphi(X)) &\mapsto& Y.
\end{array}\label{p2-cmt100c}
\end{eqnarray}

For any object $Y$ of $\cC_2$, we denote by $\mI_\varphi^Y$ the fiber category of the functor $\tts$ above $Y$. Objects of $\mI_\varphi^Y$ are pairs $(X,f)$, 
where $X$ is an object of $\cC_1$ and $f\in \varphi(X)(Y)$. Let $(X,f)$ and $(X',f')$ be two objects of $\mI_\varphi^Y$.
A morphism from $(X,f)$ to $(X',f')$ is a morphism $u\colon X\rightarrow X'$ of $\cC_1$ such that $\varphi(u)(f)=f'$. 
For every morphism $y\colon Y\rightarrow Y'$ of $\cC_2$, we have a functor 
\begin{equation}\label{p2-cmt100d}
\mI_\varphi^y\colon 
\begin{array}[t]{clcr}
\mI_\varphi^{Y'}&\rightarrow &\mI_\varphi^Y,\\ 
(X,f)&\mapsto& (X,f\circ y).
\end{array}
\end{equation}

\begin{defi}\label{p2-cmt10}
Under the assumptions of \ref{p2-cmt100}, for any $F_1\in \ob(\hcC_1)$ and $F_2\in \ob(\hcC_2)$, 
an {\em $\mI_\varphi$-system of morphisms from $F_1$ to $F_2$} is the data for any object 
$f\colon Y\rightarrow \varphi(X)$ of $\mI_\varphi$ of a {\em functorial} morphism
\begin{equation}\label{p2-cmt10a}
\xi_f\colon F_1(X)\rightarrow F_2(Y).
\end{equation}
\end{defi}

The functoriality is expressed as follows. For every morphism $(x\colon X\rightarrow X',y\colon Y\rightarrow Y')$ 
from $f\colon Y\rightarrow \varphi(X)$ to $f'\colon Y'\rightarrow \varphi(X')$ in $\mI_\varphi$, the diagram 
\begin{equation}\label{p2-cmt10b}
\xymatrix{
{F_1(X')}\ar[r]^{\xi_{f'}}\ar[d]_{x^*}&{F_2(Y')}\ar[d]^{y^*}\\
{F_1(X)}\ar[r]^{\xi_f}&{F_2(Y)}}
\end{equation}
is commutative.

\begin{prop}\label{p2-cmt2}
We keep the assumptions and notation of \ref{p2-cmt1}.  
\begin{itemize}
\item[{\rm (i)}] If for every $Y\in \ob(\cC_2)$, the category $\mI_\varphi^Y$ defined in \ref{p2-cmt100} admits a $\mU$-small cofinal subcategory 
{\rm (\cite{sga4} I 8.1.1)}, then the functor $\varphi_\rp$ \eqref{p2-cmt1a} admits a left adjoint 
\begin{equation}\label{p2-cmt2a}
\varphi^\rp\colon \hcC_1\rightarrow \hcC_2.
\end{equation} 
The latter, which is defined only up to isomorphism, can be chosen so that $\varphi=\varphi^\rp\circ \tth_1$, 
where $\tth_1\colon \cC_1\rightarrow \hcC_1$ is the canonical functor \eqref{p2-ncgt4a}. 
\item[{\rm (ii)}] If for every $Y\in \ob(\cC_2)$, the category $\mI_\varphi^Y$ admits a $\mU$-small cofinal cofiltered 
subcategory, then the functor $\varphi^\rp$ \eqref{p2-cmt2a} is left exact. 
\end{itemize}
\end{prop}

(i) The proof is similar to those of (\cite{sga4} I 5.1, 5.2 and 5.4). For any $F\in \ob(\hcC_1)$
and $Y\in \ob(\cC_2)$, we set 
\begin{equation}\label{p2-cmt2b}
\varphi^\rp(F)(Y)=\underset{\underset{(X,f)\in \mI_\varphi^Y}{\longrightarrow}}{\lim}\ F(X).
\end{equation}
By \eqref{p2-cmt100d}, $\varphi^\rp(F)$ defines naturally a presheaf on $\cC_2$. 

Let $F_1\in \ob(\hcC_1)$. 
For any $X\in \ob(\cC_1)$ and $s\in F_1(X)$, let $\xi_{F_1}(X)(s)\colon \varphi(X)\rightarrow \varphi^\rp(F_1)$ 
be the morphism defined for any $Y\in \ob(\cC_2)$, by the map 
\begin{equation}\label{p2-cmt2c}
\varphi(X)(Y)\rightarrow \varphi^\rp(F_1)(Y)=\underset{\underset{(Z,g)\in \mI_\varphi^Y}{\longrightarrow}}{\lim}\ F_1(Z)
\end{equation}
sending a section $f\in \varphi(X)(Y)$ to the image of $s$ by the morphism $F_1(X)\rightarrow \varphi^\rp(F_1)(Y)$ corresponding to $(X,f) \in \ob(\mI_\varphi^Y)$. 
We obtain a morphism 
\begin{equation}\label{p2-cmt2d}
\xi_{F_1}(X)\colon 
\begin{array}[t]{clcr}
F_1(X)&\rightarrow&\varphi_\rp(\varphi^\rp(F_1))(X)=\Hom_{\hcC_2}(\varphi(X),\varphi^\rp(F_1))\\
s&\mapsto& \xi_{F_1}(X)(s). 
\end{array}
\end{equation}
We thus define a morphism $\xi_{F_1}\colon F_1\rightarrow \varphi_\rp(\varphi^\rp(F_1))$ of $\hcC_1$.

Let $F_2\in \ob(\hcC_2)$. 
For every $Y\in \ob(\cC_2)$ and $(X,f)\in \ob(\mI_\varphi^Y)$, we consider the 
morphism $\varphi_\rp(F_2)(X)=\Hom_{\hcC_2}(\varphi(X),F_2)\rightarrow F_2(Y)$ defined by composition with $f$. 
These morphisms define a direct system indexed by $\mI_\varphi^Y$. We deduce  
a morphism $\zeta_{F_2}(Y)\colon \varphi^\rp(\varphi_\rp(F_2))(Y)\rightarrow F_2(Y)$. 
We thus obtain a morphism $\zeta_{F_2}\colon \varphi^\rp(\varphi_\rp(F_2))\rightarrow F_2$ of $\hcC_2$. 

We easily check that $\xi_{F_1}$ and $\zeta_{F_2}$ induce a canonical bifunctorial isomorphism 
\begin{equation}\label{p2-cmt2f}
\Hom_{\hcC_2}(\varphi^\rp(F_1),F_2)\stackrel{\sim}{\rightarrow}\Hom_{\hcC_1}(F_1,\varphi_\rp(F_2)).
\end{equation}

It is obvious from the definition \eqref{p2-cmt2b} that we have $\varphi=\varphi^\rp\circ \tth_1$. 

(ii) It follows from the definition \eqref{p2-cmt2b} and the fact that finite inverse limits commute with filtered direct limits. 

\begin{rema}\label{p2-cmt20}
We keep the assumptions and notation of \ref{p2-cmt1}.  If the category $\cC_1$ is $\mU$-small, then the assumption of \ref{p2-cmt2}(i) is satisfied. 
Indeed, for every $Y\in \ob(\cC_2)$, the category $\mI_\varphi^Y$ is $\mU$-small. 
If, moreover, finite inverse limits are representable in $\cC_1$ and $\varphi$ is left exact, then the assumption of \ref{p2-cmt2}(ii) is satisfied. 
Indeed, for every $Y\in \ob(\cC_2)$, the category $(\mI_\varphi^Y)^\circ$ is filtered. 
\end{rema}

\begin{cor}\label{p2-cmt9}
We keep the assumptions and notation of \ref{p2-cmt2}{\rm (i)}, and let $F_1\in \ob(\hcC_1)$, $F_2\in \ob(\hcC_2)$. 
The following data are equivalent:
\begin{itemize}
\item[(i)] a morphism $\xi\colon F_1\rightarrow \varphi_\rp(F_2)$ of $\hcC_1$ \eqref{p2-cmt1a} ;
\item[(ii)] a morphism $\zeta\colon \varphi^\rp(F_1)\rightarrow F_2$ of $\hcC_2$ \eqref{p2-cmt2a};
\item[(iii)] an $\mI_\varphi$-system of morphisms from $F_1$ to $F_2$ in the sense of \ref{p2-cmt10}. 
\end{itemize}
 \end{cor}

\begin{defi}\label{p2-cmt3}
Let $\cC_1$, $\cC_2$ be two $\mU$-sites (\cite{sga4} II 3.0.2). 
We say that a $\mU$-functor $\varphi \colon \cC_1\rightarrow \hcC_2$ is {\em continuous} 
if for every sheaf of $\mU$-sets $G$ on $\cC_2$, the presheaf $\varphi_\rp(G)$ \eqref{p2-cmt1b} is a sheaf on $\cC_1$. 
\end{defi} 

For $i=1,2$, we denote by $\tcC_i$ the topos of sheaves of $\mU$-sets on $\cC_i$. 
To say that the functor $\varphi\colon \cC_1\rightarrow \hcC_2$ is continuous amounts to saying that there exists 
a unique functor $\varphi_\rs\colon \tcC_2\rightarrow \tcC_1$ 
making strictly commutative the following diagram
\begin{equation}\label{p2-cmt3a}
\xymatrix{
{\tcC_2}\ar[r]^-(0.5){\varphi_\rs}\ar[d]_{\tti_2}&{\tcC_1}\ar[d]^{\tti_1}\\
{\hcC_2}\ar[r]^-(0.5){\varphi_\rp}&{\hcC_1,}}
\end{equation}
where $\tti_1$ and $\tti_2$ are the canonical functors. 

\subsection{}\label{p2-cmt4}
Let $\cC_1$, $\cC_2$ be two $\mU$-categories \eqref{p2-ncgt4}, $\phi\colon \cC_1\rightarrow \cC_2$ a functor.
We set $\varphi=\tth_2\circ \phi\colon \cC_1\rightarrow \hcC_2$, where $\tth_2\colon \cC_2\rightarrow \hcC_2$ is the canonical functor \eqref{p2-ncgt4a}.
We associate with $\phi$ the functor 
\begin{equation}\label{p2-cmt4a}
\phi_\rp\colon \hcC_2\rightarrow \hcC_1
\end{equation}
defined by composition with $\phi$. We have $\varphi_\rp=\phi_\rp$, see \ref{p2-cmt101}(b). 
The category $\mI_\varphi$ defined in \ref{p2-cmt100}, will also be denoted by $\mI_\phi$. 
For any object $Y$ of $\cC_2$, the category $\mI_\varphi^Y$ will also be denoted by $\mI_\phi^Y$. 
These two categories can be explicitly defined using only the functor $\phi$. 

By \ref{p2-cmt2} or (\cite{sga4} I 5.1), if the category $\cC_1$ is $\mU$-small, the functor $\phi_\rp$ admits a left adjoint 
\begin{equation}\label{p2-cmt4b}
\phi^\rp\colon \hcC_1\rightarrow \hcC_2,
\end{equation}
denoted by $\phi_!$ in (\cite{sga4} I 5.1). We clearly have $\phi^\rp=\varphi^\rp$.
It is defined only up to isomorphism, but it can be chosen so that the diagram 
\begin{equation}\label{p2-cmt4d}
\xymatrix{
{\cC_1}\ar[d]_{\tth_1}\ar[r]^\phi&{\cC_2}\ar[d]^{\tth_2}\\
{\hcC_1}\ar[r]^{\phi^\rp}&{\hcC_2,}}
\end{equation}
where the vertical arrows are the canonical functors, is commutative. 

The functor $\varphi$ is continuous \eqref{p2-cmt3} if and only if so is $\phi$ (\cite{sga4} III 3.1).  
The terminology introduced in \ref{p2-cmt3} is therefore coherent with that of (\cite{sga4} III 3.1).

\begin{prop}\label{p2-cmt5}
Let $\cC_1$ be a $\mU$-small site, $\cC_2$ a $\mU$-site, $\varphi \colon \cC_1\rightarrow \hcC_2$ a continuous $\mU$-functor \eqref{p2-cmt3}. 
Then, the composed functor 
\begin{equation}\label{p2-cmt5a}
\varphi^\rs\colon 
\xymatrix{ 
{\tcC_1}\ar[r]^-(0.5){\tti_1}&{\hcC_1}\ar[r]^-(0.5){\varphi^\rp}&{\hcC_2}\ar[r]^-(0.5){\tta_2}&{\tcC_2,}}
\end{equation}
where $\varphi^\rp$ is defined in \ref{p2-cmt2}{\rm (i)} and 
$\tta_2$ is the ``associated sheaf'' functor, is a left adjoint of $\varphi_\rs\colon \tcC_2\rightarrow \tcC_1$ \eqref{p2-cmt3a}. Moreover, the diagram
\begin{equation}\label{p2-cmt5b}
\xymatrix{
{\hcC_1}\ar[r]^-(0.5){\varphi^\rp}\ar[d]_{\tta_1}&{\hcC_2}\ar[d]^{\tta_2}\\
{\tcC_1}\ar[r]^-(0.5){\varphi^\rs}&{\tcC_2}}
\end{equation}
is commutative up to an isomorphism. 
\end{prop}

Indeed, for every $F_1\in \ob(\tcC_1)$ and $F_2\in \ob(\tcC_2)$, we have canonical isomorphisms 
\begin{eqnarray}
\Hom_{\tcC_2}(\varphi^\rs(F_1),F_2)&\stackrel{\sim}{\rightarrow}&\Hom_{\hcC_2}(\varphi^\rp(\tti_1(F_1)),\tti_2(F_2))\\
&\stackrel{\sim}{\rightarrow}&\Hom_{\hcC_1}(\tti_1(F_1),\varphi_\rp(\tti_2(F_2)))\\
&\stackrel{\sim}{\rightarrow}&\Hom_{\hcC_1}(\tti_1(F_1),\tti_1(\varphi_\rs(F_2)))\\
&\stackrel{\sim}{\rightarrow}&\Hom_{\tcC_1}(F_1,\varphi_\rs(F_2)),
\end{eqnarray}
since $\tti_2$ (resp.\ $\varphi_\rp$) is a right adjoint of $\tta_2$ (resp.\ $\varphi^\rp$) and $\tti_1$ is fully faithful. The first assertion follows. 
We then deduce from \eqref{p2-cmt3a} that the diagram \eqref{p2-cmt5b} is commutative up to isomorphism.

\subsection{}\label{p2-cmt7}
Let $\cC_1$, $\cC_2$ be two $\mU$-sites, $\varphi \colon \cC_1\rightarrow \cC_2$ a continuous functor, $\cC'_1$ (resp.\ $\cC'_2$)
a full subcategory of $\cC_1$ (resp.\ $\cC_2$) whose objects form a topologically generating $\mU$-small family (\cite{sga4} II 3.0.1). 
We denote by $\varphi_\rp \colon \hcC_2\rightarrow \hcC_1$ the functor defined by composition with $\varphi$ \eqref{p2-cmt4a} and by  
$\varphi_\rs\colon \tcC_2\rightarrow \tcC_1$ the functor defined by the relation $\tti_1\circ \varphi_\rs=\varphi_\rp\circ \tti_2$. 
By (\cite{sga4} III 1.3), the functor $\varphi_\rs$ admits a left adjoint $\varphi^\rs\colon \tcC_1\rightarrow \tcC_2$. 
For $i=1,2$, we denote by $u_i\colon \cC'_i\rightarrow \cC_i$ the canonical functor and by $u_{i,\rp}\colon \hcC_i\rightarrow \hcC'_i$ the  
functor defined by composition with $u_i$ \eqref{p2-cmt4a}. 
We equip $\cC'_i$ with the topology induced by $u_i$ (\cite{sga4} III 3.1). By (\cite{sga4} III 4.1 and its proof), 
the functor $u_i$ is continuous and cocontinuous,
and the functor $u_{i,\rp}$ induces an equivalence of categories $u_{i,\rs}\colon \tcC_i\stackrel{\sim}{\rightarrow} \tcC'_i$. 
Let $\varphi'$ be the composed functor 
\begin{equation}\label{p2-cmt7b}
\varphi'\colon
\xymatrix{
{\cC'_1}\ar[r]^-(0.4){u_1}&{\cC_1}\ar[r]^-(0.4){\varphi}&{\cC_2}\ar[r]^-(0.4){\tth_2}&{\hcC_2}\ar[r]^-(0.4){u_{2,\rp}}&{\hcC'_2.}}
\end{equation}
It is a $\mU$-functor by \ref{p2-cmt101}(a). We denote by $\varphi'_\rp\colon \hcC'_2\rightarrow \hcC'_1$ the associated functor defined in \eqref{p2-cmt1a} and  
$\varphi'^\rp\colon \hcC'_1\rightarrow \hcC'_2$ a left adjoint of $\varphi'_\rp$, which exists by \ref{p2-cmt2}.  

\begin{prop}\label{p2-cmt77}
We keep the assumptions and notation of \ref{p2-cmt7}. 
\begin{itemize}
\item[{\rm (i)}] The functor $\varphi'$ \eqref{p2-cmt7b} is continuous in the sense of \ref{p2-cmt3}. 
Let $\varphi'_\rs\colon \tcC'_2\rightarrow \tcC'_1$ be the functor 
making strictly commutative the following diagram
\begin{equation}\label{p2-cmt77a}
\xymatrix{
{\tcC'_2}\ar[r]^-(0.5){\varphi'_\rs}\ar[d]_{\tti'_2}&{\tcC'_1}\ar[d]^{\tti'_1}\\
{\hcC'_2}\ar[r]^-(0.5){\varphi'_\rp}&{\hcC'_1.}}
\end{equation}
Then, the diagram 
\begin{equation}\label{p2-cmt77b}
\xymatrix{
{\tcC_2}\ar[r]^{\varphi_\rs}\ar[d]_{u_{2,\rs}}&{\tcC_1}\ar[d]^{u_{1,\rs}}\\
{\tcC'_2}\ar[r]^{\varphi'_\rs}&{\tcC'_1}}
\end{equation}
is commutative up to isomorphism.
\item[{\rm (ii)}] The composed functor 
\begin{equation}\label{p2-cmt77c}
\varphi'^\rs\colon 
\xymatrix{ 
{\tcC'_1}\ar[r]^-(0.5){\tti'_1}&{\hcC'_1}\ar[r]^-(0.5){\varphi'^\rp}&{\hcC'_2}\ar[r]^-(0.5){\tta'_2}&{\tcC'_2}}
\end{equation}
is a left adjoint of $\varphi'_\rs$. Moreover, the diagrams  
\begin{equation}\label{p2-cmt77d}
\xymatrix{
{\hcC'_1}\ar[r]^-(0.5){\varphi'^\rp}\ar[d]_{\tta'_1}&{\hcC'_2}\ar[d]^{\tta'_2}\\
{\tcC'_1}\ar[r]^-(0.5){\varphi'^\rs}&{\tcC'_2,}}
\ \ \
\xymatrix{
{\tcC'_1}\ar[r]^{\varphi'^\rs}\ar[d]_{u^\rs_1}&{\tcC'_2}\ar[d]^{u^\rs_2}\\
{\tcC_1}\ar[r]^{\varphi^\rs}&{\tcC_2,}}
\end{equation}
where $u_i^\rs$ is a left adjoint of $u_{i,\rs}$, are commutative up to isomorphism. 
\end{itemize}
\end{prop}

(i) We give two proofs, the second one is available in \ref{p2-cmt78}.
We set $v_2=u_{2,\rp}\circ \tth_2\colon \cC_2\rightarrow \hcC'_2$ and let 
$v_{2,\rp}\colon \hcC'_2\rightarrow \hcC_2$ be the associated functor defined in \eqref{p2-cmt1a}. 
Since we have $\varphi'=v_2\circ \varphi \circ u_1$, the functor $\varphi'_\fp$ is equal to the composition 
\begin{equation}
\xymatrix{
{\hcC'_2}\ar[r]^-(0.5){v_{2,\rp}}&{\hcC_2}\ar[r]^-(0.5){\varphi_\rp}&{\hcC_1}\ar[r]^-(0.5){u_{1,\rp}}&{\hcC'_1.}}
\end{equation}
It is therefore enough to prove that $v_2$ is continuous and that the functor $v_{2,\rs}\colon \tcC'_2\rightarrow \tcC_2$ 
making strictly commutative the following diagram
\begin{equation}
\xymatrix{
{\tcC'_2}\ar[r]^-(0.5){v_{2,\rs}}\ar[d]_{\tti'_2}&{\tcC_2}\ar[d]^{\tti_2}\\
{\hcC'_2}\ar[r]^-(0.5){v_{2,\rp}}&{\hcC_2,}}
\end{equation}
is a quasi-inverse of $u_{2,\rs}$. By (\cite{sga4} I 5.1), the functor $u_{2,\rp}\colon \hcC_2\rightarrow \hcC'_2$ 
admits a right adjoint $u_2^!\colon \hcC'_2\rightarrow \hcC_2$. For every $X\in \ob(\cC_2)$ and $G\in \ob(\tcC'_2)$,
we have canonical bifunctorial isomorphisms 
\begin{eqnarray}
u_2^!(G)(X)&\stackrel{\sim}{\rightarrow}&\Hom_{\hcC_2}(\tth_2(X),u_2^!(G))\\
&\stackrel{\sim}{\rightarrow}&\Hom_{\hcC'_2}(u_{2,\rp}(\tth_2(X)),G)\nonumber\\
&\stackrel{\sim}{\rightarrow}&\Hom_{\hcC'_2}(v_2(X),G).\nonumber\
\end{eqnarray}
We have therefore an isomorphism $v_{2,\rp}\stackrel{\sim}{\rightarrow}u_2^!$. Since $u_{2}$ is continuous and cocontinuous, 
we deduce that $v_{2}$ is continuous in the sense of \ref{p2-cmt3}, and $v_{2,\rs}$ is a right adjoint of $u_{2,\rs}$ by (\cite{sga4} III 2.3(3)). 
Hence $v_{2,\rs}$ is a quasi-inverse of $u_{2,\rs}$, since the latter is an equivalence of categories. 

(ii) The first assertion follows from \ref{p2-cmt5} and the second assertion from (i).

\begin{rema}\label{p2-cmt78}
We give another proof of \ref{p2-cmt77}(i). 
It is enough to prove that the diagram 
\begin{equation}
\xymatrix{
{\tcC_2}\ar[r]^{\varphi_\rs}\ar[d]_{u_{2,\rs}}&{\tcC_1}\ar[d]^{u_{1,\rs}}\\
{\tcC'_2}\ar[d]_{\tti'_2}&{\tcC'_1}\ar[d]^{\tti'_1}\\
{\hcC'_2}\ar[r]^{\varphi'_\rp}&{\hcC'_1}}
\end{equation}
is commutative up to isomorphism. By (\cite{sga4} III 2.3(2)), the diagram 
\begin{equation}\label{p2-cmt7a}
\xymatrix{
{\hcC_i}\ar[r]^-(0.5){u_{i,\rp}}\ar[d]_{\tta_i}&{\hcC'_i}\ar[d]^{\tta'_i}\\
{\tcC_i}\ar[r]^-(0.5){u_{i,\rs}}&{\tcC'_i}}
\end{equation}
is commutative up to an isomorphism. 
For every $X\in \ob(\cC'_1)$ and $F\in \ob(\tcC_2)$, we have canonical bifunctorial isomorphisms 
\begin{eqnarray}
\varphi'_\rp(\tti'_2(u_{2,\rs}(F)))(X)&\stackrel{\sim}{\rightarrow}&\Hom_{\hcC'_2}(\varphi'(X),\tti'_2(u_{2,\rs}(F)))\\
&\stackrel{\sim}{\rightarrow}&\Hom_{\tcC'_2}(\tta'_2(u_{2,\rp}(\tth_2(\varphi(u_1(X))))),u_{2,\rs}(F))\nonumber\\
&\stackrel{\sim}{\rightarrow}&\Hom_{\tcC'_2}(u_{2,\rs}(\tta_2(\tth_2(\varphi(u_1(X))))),u_{2,\rs}(F))\nonumber\\
&\stackrel{\sim}{\rightarrow}&\Hom_{\tcC_2}(\tta_2(\tth_2(\varphi(u_1(X)))),F)\nonumber\\
&\stackrel{\sim}{\rightarrow}&F(\varphi(u_1(X))),\nonumber
\end{eqnarray}
\begin{eqnarray}
\tti'_1(u_{1,\rs}(\varphi_\rs(F)))(X)&\stackrel{\sim}{\rightarrow}&u_{1,\rs}(\varphi_\rs(F))(X)\\
&\stackrel{\sim}{\rightarrow}&\varphi_\rs(F)(u_1(X))\nonumber\\
&\stackrel{\sim}{\rightarrow}&F(\varphi(u_1(X))).\nonumber
\end{eqnarray}
\end{rema}

\begin{rema}\label{p2-cmt8}
We keep the assumptions and notation of \ref{p2-cmt7}, and 
let $F_1\in \ob(\hcC'_1)$ (resp.\ $F_2\in \ob(\hcC'_2)$), $F_1^\tta$ (resp.\ $F_2^\tta$) be the associated sheaf.
By \ref{p2-cmt9}, an $\mI_{\varphi'}$-system of morphisms from $F_1$ to $F_2$ \eqref{p2-cmt10} determines a morphism 
$\varphi'^\rp(F_1)\rightarrow F_2$ of $\hcC_2$ \eqref{p2-cmt2a}, and hence, by \eqref{p2-cmt77d}, a morphism 
$\varphi'^\rs(F^\tta_1)\rightarrow F^\tta_2$ of $\tcC'_2$.
\end{rema}

\section{Quasi-fibered categories}\label{p2-qfc}

\subsection{}\label{p2-qfc1}
We consider in this section a functor $\pi\colon \cC\rightarrow \cI$ between two $\mU$-categories. 
For any $i\in \ob(\cI)$, we denote by $\cC_i=\pi^{-1}(i)$ the fiber category above $i$. 
Recall that the objects of $\cC_i$ are the objects of $\cC$ above $i$ and the morphisms of $\cC_i$ 
are the morphisms of $\cC$ above the identity of $i$. 
We denote by 
\begin{equation}\label{p2-qfc1a}
\alpha_i\colon \cC_i\rightarrow \cC
\end{equation} 
the canonical functor and by 
\begin{equation}\label{p2-qfc1ab}
\alpha_{i,\rp}\colon \hcC\rightarrow \hcC_i
\end{equation} 
the functor defined by composition with $\alpha_i$,
where $\hcC$ and $\hcC_i$ are the categories of presheaves of $\mU$-sets on $\cC$ and $\cC_i$ respectively. 
Observe that in terms of $\pi$, the category $\mI_{\alpha_i}$, defined in \ref{p2-cmt4}, 
identifies canonically with the category of morphisms of $\cC$ with target above $i$. 
For every $Y\in \ob(\cC)$,  the category $\mI_{\alpha_i}^Y$ \eqref{p2-cmt4} identifies canonically with the category of morphisms of $\cC$ 
with source $Y$ and target above $i$.

For any objects $X$ and $Y$ of $\cC$ and every morphism $f\colon \pi(Y)\rightarrow \pi(X)$ of $\cI$, we denote by 
$\Hom_{\cC/f}(Y,X)$ the set of morphisms $u\colon Y\rightarrow X$ of $\cC$ above $f$ (i.e., such that $\pi(u)=f$). 
Since $\cC$ is a $\mU$-category, $\Hom_{\cC/f}(Y,X)$ is $\mU$-small. 
We can then associate with any morphism $f\colon j\rightarrow i$ of $\cI$, the functor 
\begin{equation}\label{p2-qfc1c}
f^*\colon \cC_i\rightarrow \hcC_j
\end{equation}
defined, for any $X\in \ob(\cC_i)$ and $Y\in \ob(\cC_j)$, by the relation 
\begin{equation}\label{p2-qfc1d}
f^*(X)(Y)=\Hom_{\cC/f}(Y,X).
\end{equation}
We call $f^*$ the {\em generalized inverse image functor by $f$ in $\cC$}. The terminology is justified in \ref{p2-qfc2}. 

Observe that in terms of $\pi$, the category $\mI_{f^*}$, defined in \ref{p2-cmt100}, identifies canonically with the category of morphisms of $\cC$ above $f$. 
For every object $Y$ above $j$,  the category $\mI_{f^*}^Y$ \eqref{p2-cmt100} identifies canonically with the category of morphisms of $\cC$ with source $Y$, above $f$.

\begin{rema}\label{p2-qfc5}
Let $i,j\in \ob(\cI)$. For every morphism $f\colon j\rightarrow i$ of $\cI$, the diagram
\begin{equation}\label{p2-qfc5a}
\xymatrix{
{\cC_i}\ar[d]_{\alpha_i}\ar[r]^{f^*}&{\hcC_j}\\
{\cC}\ar[r]^\tth&{\hcC}\ar[u]_{\alpha_{j,\rp}}}
\end{equation}
may not be commutative, where $\tth$ is the canonical functor \eqref{p2-ncgt4a}. But there is a canonical morphism 
\begin{equation}\label{p2-qfc5b}
\upgamma_f\colon f^*\rightarrow \alpha_{j,\rp}\circ \tth\circ \alpha_i.
\end{equation}
Moreover, for every $X\in \ob(\cC_i)$, the induced morphism 
\begin{equation}\label{p2-qfc5c}
\sqcup_{f\in \Hom_\cI(j,i)}f^*(X)\rightarrow \alpha_{j,\rp}(\tth(\alpha_i(X)))
\end{equation}
is an isomorphism. 
\end{rema}

\subsection{}\label{p2-qfc2}
Let $f\colon j\rightarrow i$ be a morphism of $\cI$, $f^*\colon \cC_i\rightarrow \hcC_j$ the generalized inverse image functor by $f$ in $\cC$ \eqref{p2-qfc1c}, 
$X\in \ob(\cC_i)$. Then, the presheaf $f^*(X)$ is representable by an object $Y$ of $\cC_j$ if and only if there exists a 
Cartesian morphism $F\colon Y\rightarrow X$ of $\cC$ above $f$, in the sense of (\cite{sga1} VI 5.1). 
In particular, {\em the inverse image functor by $f$ in $\cC$ exists} in the terminology of (\cite{sga1} VI §5) if and only if there exists 
a functor $f^+\colon \cC_i\rightarrow \cC_j$ making the diagram 
\begin{equation}\label{p2-qfc2a}
\xymatrix{
{\cC_i}\ar[r]^-(0.5){f^*}\ar[rd]_-(0.5){f^+}&{\hcC_j}\\
&{\cC_j}\ar[u]_{\tth_j}}
\end{equation}
commutative up to a unique isomorphism, where $\tth_j$ is the canonical functor \eqref{p2-ncgt4a}. 
When such a functor $f^+$ exists, it is unique up to a unique isomorphism. 

We see immediately that for every $i\in \ob(\cI)$, the inverse image functor by $\id_i$ in $\cC$ exists, and we have $(\id_i)^+=\id_{\cC_i}$,
and hence $(\id_i)^*=\tth_i$.

\subsection{}\label{p2-qfc16}
We assume in the remaining part of this section that the functor $\pi\colon \cC\rightarrow \cI$ satisfies the following conditions:

\addtocounter{equation}{2}

\subsubsection{}\label{p2-qfc16a}
For every morphism $f\colon j\rightarrow i$ of $\cI$, the functor $f^*\colon \cC_i\rightarrow \hcC_j$ \eqref{p2-qfc1c} is a $\mU$-functor \eqref{p2-cmt0}. 

\subsubsection{}\label{p2-qfc16b}
For every morphism $f\colon j\rightarrow i$ of $\cI$ and every $Y\in \ob(\cC_j)$, the category $\mI_{f^*}^Y$ \eqref{p2-qfc1}
admits a $\mU$-small cofinal subcategory. 

\vspace{2mm}

By \ref{p2-cmt1} and \ref{p2-qfc16a}, the functor $f^*\colon \cC_i\rightarrow \hcC_j$ \eqref{p2-qfc1c} determines a functor \eqref{p2-cmt1a}
\begin{equation}\label{p2-qfc16c}
f_\rp\colon \hcC_j\rightarrow \hcC_i.
\end{equation}
By \ref{p2-cmt2}(i) and \ref{p2-qfc16b}, the functor $f_\rp$ admits a left adjoint 
\begin{equation}\label{p2-qfc16d}
f^\rp\colon \hcC_i\rightarrow \hcC_j,
\end{equation}
defined by \eqref{p2-cmt2b}. In particular, we have $f^*=f^\rp\circ \tth_i$.

\begin{remas}\label{p2-qfc3}\
\begin{itemize}
\item[(i)] By \ref{p2-cmt101}, condition \ref{p2-qfc16a} is satisfied if for every morphism $f\colon j\rightarrow i$ of $\cI$,
one of the following conditions is satisfied:
\begin{itemize}
\item[(a)] either the category $\cC_j$ is $\mU$-small; 
\item[(b)] or the inverse image functor $f^+\colon \cC_i\rightarrow \cC_j$ by $f$ in $\cC$ exists \eqref{p2-qfc2}.
\end{itemize}
\item[(ii)] 
Condition \ref{p2-qfc16b} is satisfied if for every $i\in\ob(\cI)$ which is the target of a morphism of $\cI$ not equal to $\id_i$,  
the category $\cC_i$ is $\mU$-small. Indeed, if $f\not=\id_i$, then $\cC_i$ is $\mU$-small and hence so is $\mI_{f^*}^Y$. 
If $f=\id_i$, then $\mI_{\id_i}^Y$ admits $\id_Y$ as an initial object. 
\end{itemize}
\end{remas}

\subsection{}\label{p2-qfc6}
For all $i\in \ob(\cI)$ and $Y\in \ob(\cC)$, the category $\mI_{\alpha_i}^Y$ \eqref{p2-qfc1} admits a $\mU$-small cofinal subcategory.
Indeed, setting $j=\pi(Y)$, we have 
\begin{equation}\label{p2-qfc6c}
\mI_{\alpha_i}^Y =\bigsqcup_{f\in \Hom_\cI(j,i)}\mI_{f^*}^Y,
\end{equation}
where $\mI_{\alpha_i}^Y$ is the category defined in \ref{p2-qfc1}, i.e. the objects  of $\mI_{\alpha_i}^Y$ are the disjoint union of the objects of 
$\mI_{f^*}^Y$ for all $f\in \Hom_\cI(j,i)$, the canonical functors $\mI_{f^*}^Y\rightarrow \mI_{\alpha_i}^Y$ are fully faithful and 
for $f\not=f'$, there is no morphism from any object of $\mI_{f'^*}^Y$ to any object of $\mI_{f^*}^Y$. 
The required property follows then from \ref{p2-qfc16b}. 

By \ref{p2-cmt2}(i), we deduce that the functor $\alpha_{i,\rp}\colon \hcC\rightarrow \hcC_i$ \eqref{p2-qfc1ab} admits a left adjoint
\begin{equation}\label{p2-qfc6a}
\alpha_i^\rp\colon \hcC_i\rightarrow \hcC,
\end{equation}
defined by \eqref{p2-cmt2b}. In particular, the diagram 
\begin{equation}\label{p2-qfc6b}
\xymatrix{
{\cC_i}\ar[r]^{\alpha_i}\ar[d]_{\tth_i}&{\cC}\ar[d]^{\tth}\\
{\hcC_i}\ar[r]^{\alpha^\rp_i}&{\hcC}}
\end{equation}
is commutative. 

\begin{exemple}\label{p2-qfc7}
We denote by $[1]$ the category with only two distincts objects $0$ and $1$, such that 
$\Hom_{[1]}(i,i)=\{\id_i\}$ for $i=0,1$ and $[1]$ has a single morphism $\lambda\colon 1\rightarrow 0$ and no morphisms from $0$ to $1$.  
There is a one-to-one correspondence between the following data:
\begin{itemize}
\item[(i)] a functor $\phi\colon \cD\rightarrow [1]$ such that $\cD$ is a $\mU$-category; 
\item[(ii)] two $\mU$-categories $\cD_0$ and $\cD_1$ and a functor $\varphi\colon \cD_0\rightarrow \hcD_1$.
\end{itemize}
We associate with a functor $\phi\colon \cD\rightarrow [1]$ as in (i) the triple made of the fiber categories $\cD_0$ and $\cD_1$
and the generalized inverse image functor by $\lambda$ in $\cD$, $\varphi=\lambda^*\colon \cD_0\rightarrow \hcD_1$  \eqref{p2-qfc1c}.  
For $i=0,1\in \ob([1])$, we denote by $\alpha_i\colon \cD_i\rightarrow \cD$ the canonical functor, by 
$\alpha_{i,\rp}\colon \hcD\rightarrow \hcD_i$ the functor defined by composition with $\alpha_i$, 
and by $\tth_\cD\colon \cD\rightarrow \hcD$ the canonical functor \eqref{p2-ncgt4a}. 
By \ref{p2-qfc5} and \eqref{p2-qfc5c}, the diagram 
\begin{equation}\label{p2-qfc7d}
\xymatrix{
{\cD_0}\ar[r]^\varphi \ar[d]_{\alpha_0}&{\hcD_1}\\
{\cD}\ar[r]^{\tth_\cD}&{\hcD}\ar[u]_{\alpha_{1,\rp}}}
\end{equation}
is strictly commutative.

Conversely, we associate with a triple $(\cD_0,\cD_1,\varphi)$ as in (ii), the category $\cD$ whose objects 
are the objects of $\cD_1$ and $\cD_0$ and for $X,Y\in \ob(\cD)$, 
\begin{equation}\label{p2-qfc7a}
\Hom_\cD(Y,X)=\left\{
\begin{array}{clcr}
\Hom_{\cD_1}(Y,X)& {\rm if}\ X,Y\in \ob(\cD_1),\\
\Hom_{\cD_0}(Y,X)& {\rm if}\ X,Y\in \ob(\cD_0),\\
\varphi(X)(Y)& {\rm if}\ X \in \ob(\cD_0) \ {\rm and} \ Y\in \ob(\cD_1),\\
\emptyset & {\rm otherwise}.
\end{array}
\right.
\end{equation}
The composition in $\cD$ is defined using the bifunctor $\cD_0\times\cD_1^\circ\rightarrow \Ens$, $(X,Y)\mapsto \varphi(X)(Y)$. 
We immediately see that  there is a canonical functor $\phi\colon \cD\rightarrow [1]$ and $\cD$ is a $\mU$-category. We denote $\cD$ by
$\cD_1\rhd^\varphi\cD_0$. Observe that $\phi$ satisfies condition \ref{p2-qfc16a} if and only $\varphi$ is a $\mU$-functor. 
It satisfies condition \ref{p2-qfc16b} if and only if for every $Y\in \ob(\cD_1)$, the category $\mI_\varphi^Y$ admits a $\mU$-small 
cofinal subcategory \eqref{p2-cmt100}. 

It is clear that these two correspondences between data of type (i) and (ii) are inverse to each other. 

We denote by $\bCat_{/[1]}$ the category of functors  
$\phi\colon \cD\rightarrow [1]$ such that $\cD$ is a $\mU$-category (belonging to a given universe $\mV$ containing $\mU$), 
and by $\cT$ the category of triples $(\cD_0,\cD_1,\varphi)$ where $\cD_0$ and $\cD_1$ are 
two $\mU$-categories (belonging to $\mV$) and $\varphi\colon \cD_0\rightarrow \hcD_1$ is a functor. 
Given two objects $(\cD_0,\cD_1,\varphi)$ and $(\cD'_0,\cD'_1,\varphi')$ of $\cT$, a morphism from
$(\cD_0,\cD_1,\varphi)$ to $(\cD'_0,\cD'_1,\varphi')$ in $\cT$ is given by two functors $F_0\colon \cD_0\rightarrow \cD'_0$ and 
$F_1\colon \cD_1\rightarrow \cD'_1$ and a morphism of functors 
\begin{equation}\label{p2-qfc7b}
\varphi \rightarrow F_{1,\rp}\circ \varphi'\circ F_0,
\end{equation}
where $F_{1,\rp}\colon \hcD'_1\rightarrow \hcD_1$ is defined by composition with $F_1$. 
We immediately see that the two correspondences between data of type (i) and (ii) defined above can be upgraded to two
isomorphisms of categories inverse to each other 
\begin{equation}\label{p2-qfc7c}
\xymatrix{
{\bCat_{/[1]}}\ar@<1ex>[r]&{\cT.}\ar@<1ex>[l]}
\end{equation}
\end{exemple}

\subsection{}\label{p2-qfc17}
We denote by $[\bullet]$ the final category, i.e., the category having a unique object $\bullet$ and a unique morphism $\id_\bullet$. 
The category of presheaves  of $\mU$-sets on $[\bullet]$ identifies with the category $\Ens$ of $\mU$-sets, and 
the canonical functor $\tth_\bullet\colon [\bullet]\rightarrow \Ens$ sends $\bullet$ to the set $\{\id_\bullet\}$. 
We denote by $\varphi_\bullet\colon \cI\rightarrow \Ens$ the constant functor with value $\{\id_\bullet\}$ that we identify with 
the composition of the canonical functors $\cI\rightarrow [\bullet]\rightarrow \Ens$. We can then form the $\mU$-category 
$[\bullet]\rhd^{\varphi_\bullet} \cI$ \eqref{p2-qfc7} that we denote simply by $[\bullet]\rhd \cI$. 
It is naturally equipped with a functor $[\bullet]\rhd \cI\rightarrow [1]$, 
where the fiber category above $1$ is $[\bullet]$ and the fiber category above $0$ is $\cI$. 
For all $i,j\in \ob(\cI)$, that we consider as objects of $[\bullet]\rhd \cI$, and similarly for $\bullet$,  
we have $\Hom_{[\bullet]\rhd \cI}(i,j)=\Hom_\cI(i,j)$, $\Hom_{[\bullet]\rhd \cI}(\bullet,j)=\{\id_\bullet\}$  
and $\Hom_{[\bullet]\rhd \cI}(i,\bullet)=\emptyset$. We denote the unique morphism $\bullet\rightarrow j$ by $\beta_j$. 
So for every morphism $f\colon j\rightarrow i$ of $\cI$, we have $\beta_i=f\circ \beta_j$. 

We denote simply by $\cC\rhd \cC$ the $\mU$-category $\cC\rhd^\tth \cC$ defined by the canonical functor $\tth\colon \cC\rightarrow \hcC$ \eqref{p2-qfc7}. 
It is naturally equipped with a functor $\pi_\cC\colon \cC\rhd \cC\rightarrow [1]$, where the fiber categories above $0$ and $1$ are $\cC$. 
For all $X,Y\in \ob(\cC)$, considered as objects of $\cC\rhd \cC$, we have 
\begin{equation}\label{p2-qfc17a}
\Hom_{\cC\rhd \cC}(Y,X)=\left\{
\begin{array}{clcr}
\emptyset & {\rm if}\ \pi_\cC(X)=1 \ {\rm and} \ \pi_\cC(Y)=0,\\
\Hom_\cC(Y,X) & {\rm otherwise.}
\end{array}
\right.
\end{equation}

By \eqref{p2-qfc7c}, the canonical functor $\tau\colon \cC\rightarrow [\bullet]$, the functor $\pi\colon \cC\rightarrow \cI$ and the canonical morphism 
\begin{equation}\label{p2-qfc17b}
\tth \rightarrow \tau_\rp\circ \varphi_\bullet\circ \pi,
\end{equation}
where $\tau_\rp\colon \Ens\rightarrow \hcC$ is the ``constant presheaf functor'', define a $[1]$-functor 
\begin{equation}\label{p2-qfc17c}
\tpi\colon \cC\rhd \cC\rightarrow  [\bullet]\rhd \cI. 
\end{equation}
We have $\tpi_1=\tau\colon \cC\rightarrow [\bullet]$ and $\tpi_0=\pi\colon \cC\rightarrow \cI$. For all $X,Y\in \ob(\cC)$, 
considered as objects of $\cC\rhd \cC$ such that $\pi_\cC(X)=0$ and $\pi_\cC(Y)=1$, $\tpi$ maps 
$\Hom_\cC(Y,X)$ to $\{\beta_{\pi_\cC(X)}\}$. We deduce the following properties:
\begin{itemize}
\item[(i)] For every $i\in \ob(\cI)$ and every $X\in \ob(\cC)$ above $i$, $\id_X$ considered as a morphism of $\cC\rhd \cC$ above $\beta_i\colon \bullet\rightarrow i$ is Cartesian for $\tpi$. Therefore, the inverse image functor by $\beta_i$ in $\cC\rhd \cC$ exists and identifies with
the canonical functor $\alpha_i\colon \cC_i\rightarrow \cC$ \eqref{p2-qfc1a}, see \ref{p2-qfc2}.
\item[(ii)] For every morphism $f\colon j\rightarrow i$ of $\cI$, considered as a morphism of $[\bullet]\rhd \cI$, the generalized inverse image functor by $f$ 
in $\cC\rhd \cC$ identifies with the generalized inverse image functor by $f$ in $\cC$ \eqref{p2-qfc1c}. 
\end{itemize}
We deduce by \ref{p2-qfc3}(i) and \ref{p2-qfc6} that the functor $\tpi$ \eqref{p2-qfc17c} satisfies conditions \ref{p2-qfc16a} and \ref{p2-qfc16b}.

\subsection{}\label{p2-qfc4}
Let $g\colon \ell\rightarrow j$ and $f\colon j\rightarrow i$ be two composable morphisms of $\cI$. 
For any morphism $u\colon Y\rightarrow X$ of $\cC$ with $\pi(u)=f$, we associate a morphism 
\begin{equation}\label{p2-qfc4a}
g^\star(u)\colon g^*(Y)\rightarrow (f\circ g)^*(X)
\end{equation}
defined for any $Z\in \ob(\cC_\ell)$ by the map
\begin{equation}\label{p2-qfc4b}
\begin{array}[t]{clcr}
\Hom_{\cC/g}(Z,Y)&\rightarrow&\Hom_{\cC/f\circ g}(Z,X),\\
v&\mapsto &u\circ v.
\end{array}
\end{equation}
We chose the notation $g^\star$ to avoid confusion with the functor $g^*$ \eqref{p2-qfc1c}. 
The correspondence $u\mapsto g^\star(u)$ is clearly bifunctorial in $X$ and $Y$. We deduce a canonical functorial
morphism $f^*(X)\rightarrow g_\rp((f\circ g)^*(X))$ of $\hcC_j$, and hence a canonical morphism of functors
\begin{equation}\label{p2-qfc4g}
\delta_{f,g}\colon f^* \rightarrow g_\rp\circ (f\circ g)^*. 
\end{equation}

Let $X\in \ob(\cC_i)$, $H\in \ob(\hcC_\ell)$. The composition
\begin{equation}\label{p2-qfc4d}
\Hom_{\hcC_\ell}((f\circ g)^*(X),H)\rightarrow \Hom_{\hcC_j}(g_\rp((f\circ g)^*(X)),g_\rp(H))\rightarrow 
\Hom_{\hcC_j}(f^*(X),g_\rp(H)),
\end{equation}
where the first morphism is induced by the functor $g_\rp$ and the second one is defined by composing with 
$\delta_{f,g}(X)$ \eqref{p2-qfc4g}, defines a functorial morphism $(f\circ g)_\rp(H)(X)\rightarrow 
f_\rp(g_\rp(H))(X)$. We thus obtain a canonical morphism 
\begin{equation}\label{p2-qfc4e}
c_{\rp}^{f,g}\colon (f\circ g)_\rp\rightarrow f_\rp\circ g_\rp. 
\end{equation}
The latter induces by adjunction a canonical morphism 
\begin{equation}\label{p2-qfc4f}
c^{\rp}_{f,g}\colon g^\rp\circ f^\rp \rightarrow (f\circ g)^\rp.
\end{equation}

The morphism \eqref{p2-qfc4g} induces by adjunction a canonical morphism 
\begin{equation}\label{p2-qfc4c}
c^*_{f,g}\colon g^\rp\circ f^* \rightarrow (f\circ g)^*. 
\end{equation}

\begin{rema}\label{p2-qfc400}
The morphism $c^*_{f,g}$ \eqref{p2-qfc4c} can be explicitly described as follows. 
For all $X\in \ob(\cC_i)$, $Z\in \ob(\cC_\ell)$ and $(Y,v)\in \mI_{g^*}^Z$, we have 
a canonical morphism 
\begin{equation}\label{p2-qfc400a}
\begin{array}[t]{clcr}
\Hom_{\cC/f}(Y,X) &\rightarrow& \Hom_{\cC/f\circ g}(Z,X)\\
u&\mapsto& u\circ v.
\end{array}
\end{equation}
We deduce, for all $X\in \ob(\cC_i)$ and $Z\in \ob(\cC_\ell)$, a morphism 
\begin{equation}\label{p2-qfc400b}
\underset{\underset{(Y,v)\in \mI_{g^*}^Z}{\longrightarrow}}{\lim}\ \Hom_{\cC/f}(Y,X) \rightarrow \Hom_{\cC/f\circ g}(Z,X),
\end{equation}
which is $c^*_{f,g}(X)(Z)$ in view of \eqref{p2-cmt2b}.
\end{rema}

\begin{lem}\label{p2-qfc401}
Let $h\colon m\rightarrow \ell$, $g\colon \ell\rightarrow j$, $f\colon j\rightarrow i$ be composable morphisms of $\cI$, $v\colon Z\rightarrow Y$, 
$u\colon Y\rightarrow X$ composable morphisms of $\cC$ such that $\pi(v)=g$ and $\pi(u)=f$. Then, with the notation of \ref{p2-qfc4}, the diagrams 
\begin{equation}\label{p2-qfc401a}
\xymatrix{
{h^*(Z)}\ar[r]^-(0.5){h^\star(u\circ v)}\ar[d]_-(0.5){h^\star(v)}&{(f\circ g\circ h)^*(X)}\\
{(g\circ h)^*(Y)}\ar[ru]_-(0.5){(g\circ h)^\star(u)}&}
\end{equation}
\begin{equation}\label{p2-qfc401b}
\xymatrix{
{g^*(Y)}\ar[rr]^-(0.5){g^\star(u)}\ar[d]_{\delta_{g,h}(Y)}&&{(f\circ g)^*(X)}\ar[d]^{\delta_{f\circ g,h}(X)}\\
{h_\rp((g\circ h)^*(Y))}\ar[rr]^-(0.5){h_\rp((g\circ h)^\star(u))}&&{h_\rp((f\circ g\circ h)^*(X))}}
\end{equation}
are commutative. 
\end{lem}

Indeed, for every $T\in\ob(\cC_m)$, the evaluation of the diagram \eqref{p2-qfc401a} at $T$ coincides with the commutative diagram 
\begin{equation}\label{p2-qfc401c}
\xymatrix{
{\Hom_{\cC/h}(T,Z)}\ar[r]^-(0.5){u\circ v\circ -}\ar[d]_-(0.5){v\circ -}&{\Hom_{\cC/f\circ g\circ h}(T,X).}\\
{\Hom_{\cC/g\circ h}(T,Y)}\ar[ru]_-(0.5){u\circ -}&}
\end{equation}
Hence \eqref{p2-qfc401a} is commutative. 

Considering $v$ as a section of $g^*(Y)(Z)$, we have $g^\star(u)(v)=u\circ v\in (f\circ g)^*(X)(Z)$ and 
$\delta_{g,h}(Y)(v)=h^\star(v)\colon h^*(Z)\rightarrow (g\circ h)^*(Y)$. The commutativity of \eqref{p2-qfc401b} follows then from that of \eqref{p2-qfc401a}
by varying $v$. 

\begin{prop}\label{p2-qfc8}
With the notation of \ref{p2-qfc4}, we have the following properties. 
\begin{itemize}
\item[{\rm (i)}] For every $i\in \ob(\cI)$, we have $(\id_i)_\rp=\id_{\hcC_i}$ \eqref{p2-qfc16c}. 
\item[{\rm (ii)}] Let $f\colon j\rightarrow i$ be a morphism of $\cI$. 
We have $\delta_{f,\id_j}=\id_{f^*}$, $c^{f,\id_j}_\rp=\id_{f_\rp}$ and $c^{\id_i,f}_\rp=\id_{f_\rp}$.
The morphism $\delta_{\id_i,f}\colon \tth_i\rightarrow f_\rp\circ f^*$ is defined for all $X,Y\in \ob(\cC_i)$ by the morphism 
\begin{equation}
\Hom_{\cC_i}(Y,X)\rightarrow   \Hom_{\hcC_j}(f^*(Y),f^*(X))
\end{equation}
defined by the functor $f^*$. 
\item[{\rm (iii)}] Let $g\colon \ell\rightarrow j$, $f\colon j\rightarrow i$ be two composable morphisms of $\cI$.  
For all $X\in \ob(\cC_i)$ and $H\in \ob(\hcC_\ell)$, the map $c_{\rp}^{f,g}(H)(X)$ \eqref{p2-qfc4e} is the composition
\begin{equation}\label{p2-qfc8f}
\Hom_{\hcC_\ell}((f\circ g)^*(X),H)\rightarrow \Hom_{\hcC_\ell}(g^\rp(f^*(X)),H)\stackrel{\sim}{\rightarrow} 
\Hom_{\hcC_j}(f^*(X),g_\rp(H)),
\end{equation}
where the first map is induced by $c^*_{f,g}$ \eqref{p2-qfc4c} and the second one is the adjunction isomorphism.
Moreover, the diagram 
\begin{equation}\label{p2-qfc8a}
\xymatrix{
{g^\rp\circ f^*}\ar[rr]^-(0.5){c^*_{f,g}}\ar[d]&&{(f\circ g)^*}\ar[d]\\
{g^\rp\circ f^\rp\circ \tth_i}\ar[rr]^-(0.5){c^\rp_{f,g}\circ \tth_i}&&{(f\circ g)^\rp\circ \tth_i,}}
\end{equation}
where the vertical arrows are the canonical isomorphisms, is commutative. 
\item[{\rm (iv)}] For all composable morphisms $h\colon m\rightarrow \ell$, $g\colon \ell\rightarrow j$ and $f\colon j\rightarrow i$ of $\cI$,
the diagrams
\begin{equation}\label{p2-qfc8b}
\xymatrix{
f^*\ar[rr]^-(0.5){\delta_{f,g}}\ar[d]_-(0.5){\delta_{f,g\circ h}}&&{g_\rp\circ(f\circ g)^*}\ar[d]^-(0.5){g_\rp\circ \delta_{f\circ g,h}}\\
{(g\circ h)_\rp\circ (f\circ g\circ h)^*}\ar[rr]^-(0.5){c_\rp^{g,h}\circ(f\circ g\circ h)^*}&&{g_\rp\circ h_\rp \circ (f\circ g\circ h)^*}}
\end{equation}
\begin{equation}\label{p2-qfc8c}
\xymatrix{
{(f\circ g\circ h)_\rp}\ar[rr]^-(0.5){c^{f\circ g,h}_\rp}\ar[d]_{c^{f,g\circ h}_\rp}&&{(f\circ g)_\rp\circ h_\rp}\ar[d]^{c^{f,g}_\rp\circ h_\rp}\\
{f_\rp\circ (g\circ h)_\rp}\ar[rr]^-(0.5){f_\rp\circ c^{g,h}_\rp}&&{f_\rp\circ g_\rp \circ h_\rp}}
\end{equation}
\begin{equation}\label{p2-qfc8d}
\xymatrix{
{h^\rp\circ g^\rp \circ f^\rp}\ar[rr]^-(0.5){c^\rp_{g,h}\circ f^\rp}\ar[d]_{h^\rp\circ c^\rp_{f,g}}&&{(g\circ h)^\rp\circ f^\rp}\ar[d]^{c^\rp_{f,g\circ h}}\\
{h^\rp\circ (f\circ g)^\rp}\ar[rr]^-(0.5){c^\rp_{f\circ g,h}}&&{(f\circ g\circ h)^\rp}}
\end{equation}
\begin{equation}\label{p2-qfc8e}
\xymatrix{
{h^\rp\circ g^\rp \circ f^*}\ar[rr]^-(0.5){c^\rp_{g,h}\circ f^*}\ar[d]_{h^\rp\circ c^*_{f,g}}&&{(g\circ h)^\rp\circ f^*}\ar[d]^{c^*_{f,g\circ h}}\\
{h^\rp\circ (f\circ g)^*}\ar[rr]^-(0.5){c^*_{f\circ g,h}}&&{(f\circ g\circ h)^*}}
\end{equation}
are commutative. 
\end{itemize}
\end{prop}

(i) Indeed, we have $(\id_i)^*=\tth_i\colon \cC_i\rightarrow \hcC_i$ is the canonical functor. 

(ii) This is obvious from the definitions. 

(iii) Indeed, for all $X\in \ob(\cC_i)$ and $H\in \ob(\hcC_\ell)$, the diagram 
\begin{equation}
\xymatrix{
{\Hom_{\hcC_\ell}((f\circ g)^*(X),H)}\ar[r]\ar[d]&{\Hom_{\hcC_j}(g_\rp((f\circ g)^*(X)),g_\rp(H))}\ar[d]\\
{\Hom_{\hcC_\ell}(g^\rp(f^*(X)),H)}\ar[r]\ar[rd]&{\Hom_{\hcC_\ell}(g_\rp(g^\rp(f^*(X))),g_\rp(H))}\ar[d]\\
&{\Hom_{\hcC_\ell}(f^*(X),g_\rp(H)),}}
\end{equation}
where the horizontal (resp.\ upper vertical) arrows are induced by the functor $g_\rp$ (resp.\ the morphism $c^*_{f,g}$ \eqref{p2-qfc4c}), 
the lower vertical arrow is induced by the adjunction morphism $\id\rightarrow g_\rp\circ g^\rp$ and the slant arrow is 
the canonical adjunction isomorphism, is commutative. The composed right vertical arrows is then induced by 
the morphism $\delta_{f,g}$ \eqref{p2-qfc4g}. The first assertion \eqref{p2-qfc8f} follows. 

On the other hand, for all $G\in \ob(\hcC_i)$ and $H\in \ob(\hcC_\ell)$, the diagram 
\begin{equation}
\xymatrix{
{\Hom_{\hcC_\ell}((f\circ g)^\rp(G),H)}\ar[r]\ar[d]&{\Hom_{\hcC_\ell}(g^\rp(f^\rp(G)),H)}\ar[d]\\
{\Hom_{\hcC_i}(G,(f\circ g)_\rp(H))}\ar[r]&{\Hom_{\hcC_\ell}(G,f_\rp(g_\rp(H))),}}
\end{equation}
where the upper (resp.\ lower) horizontal arrow is induced by the composition with $c^\rp_{f,g}$ (resp.\ $c_\rp^{f,g}$) 
and the vertical arrows are the canonical adjunction isomorphisms, is commutative. Taking $G=\tth_i(X)$ for $X\in \ob(\cC_i)$, we deduce 
that the map $c_{\rp}^{f,g}(H)(X)$ is the composition
\begin{equation}
\Hom_{\hcC_\ell}((f\circ g)^\rp(\tth_i(X)),H)\rightarrow \Hom_{\hcC_\ell}(g^\rp(f^\rp(\tth_i(X))),H)\stackrel{\sim}{\rightarrow} 
\Hom_{\hcC_\ell}(f^\rp(\tth_i(X)),g_\rp(H)),
\end{equation}
where the first map is induced by $c^\rp_{f,g}$ and the second one is the adjunction isomorphism.
The second assertion \eqref{p2-qfc8a} follows then from the first one.

(iv) The commutativity of the diagram \eqref{p2-qfc8b} evaluated at $X\in \ob(\cC_i)$ 
is equivalent to the commutativity of the diagram \eqref{p2-qfc401b} for all morphisms $u\colon Y\rightarrow X$ of $\cC$ above $f$. 
The commutativity of \eqref{p2-qfc8b} implies the commutativity of the three other diagrams.

\begin{defi}\label{p2-qfc9}
We say that the functor $\pi\colon \cC\rightarrow \cI$ is {\em quasi-fibering} 
and that the category $\cC$ is {\em quasi-fibered} over $\cI$ if for
all composable morphisms $g\colon \ell\rightarrow j$ and $f\colon j\rightarrow i$ of $\cI$, the morphism $c_{f,g}^*$ \eqref{p2-qfc4c} is an isomorphism. 
\end{defi}

It follows then from \ref{p2-qfc8}(iii) and by adjunction that for all composable morphisms $g\colon \ell\rightarrow j$ and $f\colon j\rightarrow i$ of $\cI$, 
the morphisms $c^{f,g}_\rp$ \eqref{p2-qfc4e} and $c_{f,g}^\rp$ \eqref{p2-qfc4f} are isomorphisms.  

Whenever we consider a quasi-fibered category $\cC$ over $\cI$, we imply that $\cC$ and $\cI$ are $\mU$-categories and $\pi$ 
satisfies conditions \eqref{p2-qfc16a} and \eqref{p2-qfc16b}.

\begin{rema}\label{p2-qfc166}
We may require that the functor $\pi$ satisfies condition \ref{p2-qfc16a} and a condition weaker than \ref{p2-qfc16b}, namely: 

\addtocounter{equation}{1}

\subsubsection{}\label{p2-qfc166b}
For every morphism $f\colon j\rightarrow i$ of $\cI$, the functor $f_\rp\colon \hcC_j\rightarrow \hcC_i$ 
admits a left adjoint 
\begin{equation}\label{p2-qfc166d}
f^\rp\colon \hcC_i\rightarrow \hcC_j, 
\end{equation}
and the morphism 
\begin{equation}\label{p2-qfc166e}
f^\rp\circ \tth_i\rightarrow f^*
\end{equation}
deduced by adjunction from the morphism $\tth_i\rightarrow f_\rp \circ f^*$ defined by $\id_{f^*}\colon f^*\rightarrow f^*$, is an isomorphism.

\vspace{2mm}
Indeed, by \ref{p2-cmt2}(i), condition \ref{p2-qfc16b} implies \ref{p2-qfc166b}. 

If we assume that the functor $\pi$ satisfies conditions \ref{p2-qfc16a} and \ref{p2-qfc166b} rather than \ref{p2-qfc16b}, then Proposition~\ref{p2-qfc8} still holds, and we may therefore define the property of $\pi$ being \emph{quasi-fibering} in \ref{p2-qfc9}. However, the conclusions of \ref{p2-qfc6} no longer hold under this weaker assumption, and these are required in the proof of \ref{p2-qfc19} below.
\end{rema}

\subsection{}\label{p2-qfc10}
Let $g\colon \ell\rightarrow j$, $f\colon j\rightarrow i$ be two composable morphisms of $\cI$ such that 
the inverse image functors by $f$, $g$ and $g\circ f$ in $\cC$ exist \eqref{p2-qfc2}. Consider the diagram 
\begin{equation}\label{p2-qfc10b}
\xymatrix{
{\cC_i}\ar[r]^-(0.5){f^*}\ar[rd]_-(0.5){f^+}&{\hcC_j}\ar[rd]^{g^\rp}&\\
&{\cC_j}\ar[u]_-(0.5){\tth_j}\ar[r]^{g^*}\ar[rd]_{g^+}&{\hcC_\ell}\\
&&{\cC_\ell}\ar[u]_-(0.5){\tth_\ell}}
\end{equation}
whose triangles are commutative up to canonical isomorphisms by \eqref{p2-qfc2a} and \ref{p2-cmt2}(i). Then the canonical 
morphism $c_{f,g}^*\colon g^\rp\circ f^*\rightarrow (f\circ g)^*$ \eqref{p2-qfc4c} induces a canonical morphism 
\begin{equation}\label{p2-qfc10c}
c_{f,g}^+\colon g^+\circ f^+\rightarrow (f\circ g)^+. 
\end{equation}
We see immediately that $c_{f,g}^+$ is an isomorphism if and only if so is $c_{f,g}^*$. We deduce the following:

\begin{prop}\label{p2-qfc11}
Assume that the functor $\pi\colon \cC\rightarrow \cI$ is prefibering in the sense of {\em (\cite{sga1} VI 6.1)}, 
i.e., that for every morphism of $\cI$, the inverse image functor by $f$ in $\cC$ exists \eqref{p2-qfc10}. 
Then, $\pi$ is quasi-fibering  in the sense of \ref{p2-qfc9} if and only if it is fibering in the sense of {\em (\cite{sga1} VI 6.1)}, i.e.,
for all composable morphisms $g\colon \ell\rightarrow j$ and $f\colon j\rightarrow i$ of $\cI$, the morphism $c_{f,g}^+$ \eqref{p2-qfc10c} 
is an isomorphism.
\end{prop}

\begin{lem}\label{p2-qfc12}
Let $g\colon \ell\rightarrow j$ and $f\colon j\rightarrow i$ be two composable morphisms of $\cI$. 
Then, the morphism $c_{f,g}^*$ \eqref{p2-qfc4c} is an isomorphism if the following conditions are satisfied:
\begin{itemize}
\item[{\rm (i)}] For every morphism $w\colon Z\rightarrow X$ of $\cC$ above $f\circ g$, there exist an object $Y$ of $\cC_j$,
a morphism $v\colon Z\rightarrow Y$ of $\cC$ above $g$ and a morphism $u\colon Y\rightarrow X$ of $\cC$ above $f$ such that $w=u\circ v$. 
\item[{\rm (ii)}] For all $X\in \ob(\cC_i)$, $Y_1,Y_2\in \ob(\cC_j)$, $Z\in \ob(\cC_\ell)$ and morphisms 
\begin{equation}\label{p2-qfc12a}
\xymatrix{
Z\ar[rd]_-(0.5){v_2}\ar[r]^-(0.5){v_1}&Y_1\ar[rd]^-(0.6){u_1}&\\
&Y_2\ar[r]^-(0.4){u_2}&X}
\end{equation}
such that $u_1$ and $u_2$ (resp.\ $v_1$ and $v_2$) are above $f$ (resp.\ $g$) and $u_1\circ v_1=u_2\circ v_2$, 
there exist $Y\in \ob(\cC_j)$,  two morphisms $v\colon Z\rightarrow Y$ and $u\colon Y\rightarrow X$ of $\cC$ and two morphisms 
$y_1\colon Y\rightarrow Y_1$ and $y_2\colon Y\rightarrow Y_2$ of $\cC_j$ such that $v_1=y_1\circ v$, $v_2=y_2\circ v$ and $u=u_1\circ y_1=u_2\circ y_2$. 
\end{itemize}
\end{lem}

It follows directly from the definition \eqref{p2-qfc400b} of $c_{f,g}^*$.  Observe that for all $X\in \ob(\cC_i)$ and $Z\in \ob(\cC_\ell)$, the canonical morphism 
\begin{equation}
\underset{(Y,v)\in\ob(\mI_{g^*}^Z)}{\bigcup} \Hom_{\cC/f}(Y,X)  \rightarrow  \underset{\underset{(Y,v)\in \mI_{g^*}^Z}{\longrightarrow}}{\lim}\ \Hom_{\cC/f}(Y,X)
\end{equation}
is surjective. 

\begin{lem}\label{p2-qfc13}
We keep the assumptions and notation of \ref{p2-qfc12}. Then, the condition \ref{p2-qfc12}{\rm (ii)} is satisfied whenever the following two conditions are satisfied:
\begin{itemize}
\item[{\rm (a)}] The products of two objects are representable in $\cC_i$ and in $\cC_j$ and 
the generalized inverse image functors $f^*\colon \cC_i\rightarrow \hcC_j$, $g^*\colon \cC_j\rightarrow \hcC_\ell$ 
and $(f\circ g)^*\colon \cC_i\rightarrow \hcC_\ell$
commute with these products. 
\item[{\rm (b)}] For every morphism $Y\rightarrow X$ of $\cC$ above $f$ and every morphism $X'\rightarrow X$ of $\cC_i$, 
the fiber product $Y'=X'\times_XY$ is representable in $\cC$ and the canonical projection $Y'\rightarrow Y$ is a morphism of $\cC_j$. 
\end{itemize}
\end{lem} 

Let us consider a commutative diagram \eqref{p2-qfc12a} satisfying the assumptions of \ref{p2-qfc12}(ii). 
The product $X\times X$ in $\cC_i$ (resp.\ $Y_1\times Y_2$ in $\cC_j$) is well defined. 
Since $f^*(X\times X)=f^*(X)\times f^*(X)$, there exists a unique morphism 
$u_1\times u_2\colon Y_1\times Y_2\rightarrow X\times X$ of $\cC$ over $f$ that fits into commutative diagrams, for $i=1,2$,  
\begin{equation}
\xymatrix{
{Y_1\times Y_2}\ar[r]^-(0.5){u_1\times u_2}\ar[d]&{X\times X}\ar[d]\\
Y_i\ar[r]^{u_i}&X,}
\end{equation} 
where the vertical arrows are the projections on the $i$th factor. 
Let $Y$ be the fibered product in $\cC$ of $u_1\times u_2$ and the diagonal morphism $\Delta\colon X\rightarrow X\times X$, 
and let $y$ and $u$ be the canonical projections:
\begin{equation}\label{p2-qfc13a}
\xymatrix{
Y\ar[rr]^u\ar[d]_y\ar@{}[rrd]|\Box&&X\ar[d]^\Delta\\
{Y_1\times Y_2}\ar[rr]^{u_1\times u_2}&&{X\times X.}}
\end{equation}
Then, $y$ is a morphism of $\cC_j$ and $u$ is above $f$. Let $y_1\colon Y\rightarrow Y_1$ and $y_2\colon Y\rightarrow Y_2$ be 
the morphisms induced by $y$. We set $w=(v_1,v_2)\in g^*(Y_1\times Y_2)(Z)=g^*(Y_1)(Z)\times g^*(Y_2)(Z)$, 
that we consider as a morphism $w\colon Z\rightarrow Y_1\times Y_2$ of $\cC$ above $g$. 
Since $(f\circ g)^*(X\times X)=(f\circ g)^*(X)\times (f\circ g)^*(X)$, we deduce that $(u_1\times u_2)\circ w=\Delta\circ u_1\circ v_1=\Delta\circ u_2\circ v_2$. 
Therefore, $w$ and $u_1\circ v_1=u_2\circ v_2$ induce a morphism $v\colon Z\rightarrow Y$. 
We obviously have $v_1=y_1\circ v$, $v_2=y_2\circ v$ and $u=u_1\circ y_1=u_2\circ y_2$.

\begin{prop}\label{p2-qfc14}
Assume that the category $\cC$ is quasi-fibered over $\cI$ \eqref{p2-qfc9}. 
Then, there exists a canonical cleaved and normalized fibered category {\rm (\cite{sga1} VI 7.1)}
\begin{equation}\label{p2-qfc14a}
\cP^\vee\rightarrow \cI^\circ
\end{equation}
such that for every $i\in \ob(\cI)$, the fiber category $\cP^\vee_i$ of $\cP^\vee$ over $i$ is $\hcC_i$, 
and for every morphism $f\colon j\rightarrow i$ of $\cI$, the inverse image functor associated with $f^\circ$ is $f_\rp\colon \hcC_j\rightarrow \hcC_i$ \eqref{p2-qfc16c}.
\end{prop}

It follows from \ref{p2-qfc8}, (\cite{sga1} VI 7.2 and 8).

\subsection{}\label{p2-qfc15}
For every morphism $f\colon j\rightarrow i$ of $\cI$, we have the canonical morphism \eqref{p2-qfc5b}
\begin{equation}\label{p2-qfc15a}
\upgamma_f\colon f^*\rightarrow \alpha_{j,\rp}\circ \tth\circ \alpha_i
\end{equation}
defined for any $X\in \ob(\cC_i)$ and $Y\in \ob(\cC_j)$, by the canonical injection 
\begin{equation}\label{p2-qfc15c}
\Hom_{\cC/f}(Y,X)\rightarrow \Hom_\cC(Y,X).
\end{equation}

Let $F\in \ob(\hcC)$, $X\in \ob(\cC_i)$. The composed morphism 
\begin{equation}\label{p2-qfc15e}
\Hom_{\hcC}(\tth\circ \alpha_i(X),F)\rightarrow \Hom_{\hcC_j}(\alpha_{j,\rp}\circ \tth\circ \alpha_i(X),\alpha_{j,\rp}(F))
\rightarrow \Hom_{\hcC_j}(f^*(X),\alpha_{j,\rp}(F)),
\end{equation}
where the first morphism is induced by the functor $\alpha_{j,\rp}$ and the second one is defined by composition with  $\upgamma_f$ \eqref{p2-qfc15a},
defines a functorial morphism $\alpha_{i,\rp}(F)(X)\rightarrow f_\rp(\alpha_{j,\rp}(F))(X)$. 
We thus obtain a morphism of functors 
\begin{equation}\label{p2-qfc15d}
\upgamma_{f,\rp}\colon \alpha_{i,\rp}\rightarrow f_\rp\circ \alpha_{j,\rp}. 
\end{equation}
The latter induces by adjunction a morphism 
\begin{equation}\label{p2-qfc15f}
\upgamma^\rp_f\colon  \alpha_j^\rp\circ f^\rp\rightarrow \alpha_i^\rp. 
\end{equation}

The morphism \eqref{p2-qfc15a} induces by adjunction a canonical morphism 
\begin{equation}\label{p2-qfc15b}
\upgamma_f^*\colon \alpha_j^\rp\circ f^*\rightarrow \tth\circ \alpha_i. 
\end{equation}

\begin{rema}\label{p2-qfc18}
For every morphism $f\colon j\rightarrow i$ of $\cI$,
we have $\beta_i=f\circ \beta_j \colon \bullet\rightarrow i$ in $[\bullet]\rhd \cI$. In view of \ref{p2-qfc17}(i)-(ii), we have
\begin{equation}\label{p2-qfc18a}
\upgamma_f=\delta_{f,\beta_j}\colon f^*\rightarrow \alpha_{j,\rp}\circ \tth\circ \alpha_i,
\end{equation}
where $\upgamma_f$ is defined in \eqref{p2-qfc15a} and $\delta_{f,\beta_j}$ is defined in \eqref{p2-qfc4g} relatively to 
the functor $\tpi\colon \cC\rhd \cC\rightarrow  [\bullet]\rhd \cI$ \eqref{p2-qfc17c}.
\end{rema}

\begin{prop}\label{p2-qfc19}
With the notation of \ref{p2-qfc4} and \ref{p2-qfc15}, we have the following properties:
\begin{itemize}
\item[{\rm (i)}] For every $i\in \ob(\cI)$, $\upgamma^*_{\id_i}\colon \alpha_i^\rp\circ \tth_i\rightarrow \tth\circ \alpha_i$
is the canonical isomorphism making commutative the diagram \eqref{p2-qfc6b}, see \ref{p2-qfc2}. 
In particular, $\upgamma_{\id_i,\rp}=\id_{\alpha_{i,\rp}}$. 
\item[{\rm (ii)}] For all composable morphisms $g\colon \ell\rightarrow j$ and $f\colon j\rightarrow i$ of $\cI$, the diagrams 
\begin{equation}\label{p2-qfc19d}
\xymatrix{
f^*\ar[rr]^-(0.5){\upgamma_{f}}\ar[d]_-(0.5){\delta_{f,g}}&&{\alpha_{j,\rp}\circ \tth\circ \alpha_i}\ar[d]^-(0.5){\upgamma_{g,\rp}\circ \tth\circ \alpha_i}\\
{g_\rp\circ (f\circ g)^*}\ar[rr]^-(0.5){g_\rp\circ\upgamma_{f\circ g}}&&{g_\rp\circ \alpha_{\ell,\rp} \circ \tth\circ \alpha_i}}
\end{equation}
\begin{equation}\label{p2-qfc19a}
\xymatrix{
{\alpha_{i,\rp}}\ar[rr]^-(0.5){\upgamma_{f,\rp}}\ar[d]_{\upgamma_{f\circ g,\rp}}&&{f_\rp\circ \alpha_{j,\rp}}\ar[d]^{f_\rp\circ \upgamma_{g,\rp}}\\
{(f\circ g)_\rp\circ \alpha_{\ell,\rp}}\ar[rr]^-(0.5){c_\rp^{f,g}\circ \alpha_{\ell,\rp}}&&{f_\rp\circ g_\rp\circ \alpha_{\ell,\rp}}}
\end{equation}
\begin{equation}\label{p2-qfc19c}
\xymatrix{
{\alpha_{\ell}^{\rp}\circ g^\rp\circ f^\rp}\ar[r]^-(0.5){\upgamma_g^\rp\circ f^\rp}\ar[d]_{\alpha_{\ell}^{\rp}\circ c^\rp_{f,g}}&{\alpha_{j}^\rp\circ f^\rp}\ar[d]^{\upgamma^\rp_f}\\
{\alpha_{\ell}^\rp\circ (f\circ g)^\rp}\ar[r]^-(0.5){\upgamma^\rp_{f\circ g}}&{\alpha^\rp_i}}
\end{equation}
\begin{equation}\label{p2-qfc19b}
\xymatrix{
{\alpha_{\ell}^{\rp}\circ g^\rp\circ f^*}\ar[r]^-(0.5){\upgamma_g^\rp\circ f^*}\ar[d]_{\alpha_{\ell}^{\rp}\circ c^*_{f,g}}&{\alpha_{j}^\rp\circ f^*}\ar[d]^{\upgamma^*_f}\\
{\alpha_{\ell}^\rp\circ (f\circ g)^*}\ar[r]^-(0.5){\upgamma^*_{f\circ g}}&{\tth\circ \alpha_{i}}}
\end{equation}
are commutative. 
\end{itemize}
\end{prop}

It follows from \ref{p2-qfc8} in view of \ref{p2-qfc18}. 

\subsection{}\label{p2-qfc20}
Condition \ref{p2-qfc16b} is not necessary for this subsection and the next proposition \ref{p2-qfc21}.
We define the category $\cP(\cC/\cI)$ as follows. Objects of $\cP(\cC/\cI)$ are collections $(F_i,\gamma_f)_{i,f}$, where for each $i\in \ob(\cI)$, 
$F_i$ is an object of $\hcC_i$ and for each morphism $f\colon j\rightarrow i$ of $\cI$, 
$\gamma_f\colon F_i\rightarrow f_\rp(F_j)$ is a morphism of $\hcC_i$, where $f_\rp$ is the functor \eqref{p2-qfc16c}, 
satisfying the following conditions: 
\begin{itemize}
\item[(i)] for every $i\in \ob(\cI)$, $\gamma_{\id_i}=\id_{F_i}$, see \ref{p2-qfc8}(i);
\item[(ii)] for all composable morphisms $g\colon \ell\rightarrow j$ and $f\colon j\rightarrow i$ of $\cI$, the diagram 
\begin{equation}\label{p2-qfc20a}
\xymatrix{
F_i\ar[r]^-(0.5){\gamma_{f}}\ar[d]_{\gamma_{f\circ g}}&{f_\rp(F_j)}\ar[d]^-(0.5){f_\rp(\gamma_{g})}\\
{(f\circ g)_\rp(F_\ell)}\ar[r]^{c^{f,g}_\rp(F_\ell)}&{f_\rp(g_\rp(F_\ell)),}}
\end{equation}
where $c^{f,g}_\rp$ is the morphism \eqref{p2-qfc4e}, is commutative. 
\end{itemize}
Such an object will also be denoted by $\{i\mapsto F_i\}$, omitting the $\gamma_{f}$, when 
there is no risk of ambiguity. Let $(F_i,\gamma_{f})_{i,f}$, $(F'_i,\gamma'_{f})_{i,f}$ be two objects of 
$\cP(\cC/\cI)$. A morphism from $(F_i,\gamma_{f})_{i,f}$
to $(F'_i,\gamma'_{f})_{i,f}$ is a collection of morphisms $u_i\colon F_i\rightarrow F'_i$ of $\hcC_i$, for each $i\in \ob(\cI)$, 
such that for every morphism 
$f\colon j\rightarrow i$ of $\cI$, the diagram 
\begin{equation}\label{p2-qfc20b}
\xymatrix{
{F_i}\ar[r]^-(0.5){\gamma_{f}}\ar[d]_{u_i}&{f_\rp(F_j)}\ar[d]^{f_*(u_j)}\\
{F'_i}\ar[r]^-(0.5){\gamma'_{f}}&{f_\rp(F'_j)}}
\end{equation}
is commutative. 

Let $F\in \ob(\hcC)$. For any $i\in \ob(\cI)$, we set $F_i=\alpha_{i,\rp}(F)\in \ob(\hcC_i)$, where $\alpha_{i,\rp}$ is the morphism \eqref{p2-qfc1ab}. 
For any morphism $f\colon j\rightarrow i$ of $\cI$, we set 
\begin{equation}\label{p2-qfc20c}
\gamma_f=\upgamma_{f,\rp}(F)\colon F_i\rightarrow f_\rp(F_j),
\end{equation}
where $\upgamma_{f,\rp}$ is the morphism \eqref{p2-qfc15d}. 
It follows from the strict commutativity of \eqref{p2-qfc19a} that $(F_i,\gamma_f)_{i,f}$ is naturally an object of $\cP(\cC/\cI)$. 
We thus define a functor 
\begin{equation}\label{p2-qfc20d}
\iota\colon \hcC\rightarrow \cP(\cC/\cI).
\end{equation}

Conversely, let $(F_i,\gamma_f)_{i,f}$ be an object of $\cP(\cC/\cI)$. For any $X\in \ob(\cC)$ with $\pi(X)=i$, we set 
\begin{equation}\label{p2-qfc20e}
F(X)=F_i(X). 
\end{equation}
Let $u\colon Y\rightarrow X$ be a morphism of $\cC$ above a morphism $f\colon j\rightarrow i$ of $\cI$. 
We define a morphism $u^\flat\colon F(X)\rightarrow F(Y)$ by the composition
\begin{equation}\label{p2-qfc20f}
\xymatrix{
{F_i(X)}\ar[r]^-(0.35){\gamma_f(X)}&{\Hom_{\hcC_j}(f^*(X),F_j)}\ar[r]^-(0.5){F_j(u)}&{F_j(Y),}}
\end{equation}
where $F_j(u)$ is induced by the composition with $u$ seen as a morphism $\tth_j(Y)\rightarrow f^*(X)$ of $\hcC_j$ \eqref{p2-qfc1d}. 
Let $v\colon Z\rightarrow Y$ be a second morphism of $\cC$ above a morphism $g\colon \ell\rightarrow j$ of $\cI$. 
Then, we have $(u\circ v)^\flat=v^\flat\circ u^\flat$. Indeed, the diagram 
\begin{equation}\label{p2-qfc20g}
\xymatrix{
{F_i(X)}\ar[rr]^-(0.5){\gamma_f(X)}\ar[d]_{\gamma_{f\circ g}(X)}&&{f_\rp(F_j)(X)}\ar[rr]^-(0.5){F_j(u)}
\ar[d]^{\gamma_g(f^*(X))}&&{F_j(Y)}\ar[d]^{\gamma_g(Y)}\\
{(f\circ g)_\rp(F_\ell)(X)}\ar[drrrr]_{F_\ell(u\circ v)}\ar[rr]^-(0.5){c_\rp^{f,g}(F_\ell)(X)}&&
{f_\rp(g_\rp(F_\ell))(X)}\ar[rr]^-(0.5){g_\rp(F_\ell)(u)}&&{g_\rp(F_\ell)(Y)}\ar[d]^{F_\ell(v)}\\
&&&&{F_\ell(Z)}} 
\end{equation}
is commutative: the left (resp.\ right) square in commutative by assumption \eqref{p2-qfc20a} (resp.\ functoriality) and the lower triangle is commutative 
because $g_\rp(F_\ell)(u) \circ c_\rp^{f,g}(F_\ell)(X)$ is induced by the composition with $g^\star(u)\colon g^*(Y)\rightarrow (f\circ g)^*(X)$
by the definition \eqref{p2-qfc4d} of  $c_\rp^{f,g}$. 
Hence, $F$ is naturally a presheaf on $\cC$. 
We thus define a functor 
\begin{equation}\label{p2-qfc20h}
\jmath\colon \cP(\cC/\cI)\rightarrow \hcC.
\end{equation}

\begin{prop}\label{p2-qfc21}
The functors $\iota$ \eqref{p2-qfc20d} and $\jmath$ \eqref{p2-qfc20h} are isomorphisms of categories inverse to each other. 
\end{prop}

Condition \ref{p2-qfc16b} is not necessary for this proposition.

\begin{rema}\label{p2-qfc22}
Assume that the category $\cC$ is quasi-fibered over $\cI$ \eqref{p2-qfc9}.
Let $\cP^\vee\rightarrow \cI^\circ$ be the fibered category defined in \ref{p2-qfc14}. 
Then, $\cP(\cC/\cI)$ \eqref{p2-qfc20} identifies with the category $\bHom_{\cI^\circ}(\cI^\circ,\cP^\vee)$
of $\cI^\circ$-functors $\cI^\circ \rightarrow \cP^\vee$, see (\cite{ag2} 2.1.5). The functor $\iota$ \eqref{p2-qfc20d} 
induces an equivalence of categories 
\begin{equation}\label{p2-qfc22a}
\hcC\stackrel{\sim}{\rightarrow} \bHom_{\cI^\circ}(\cI^\circ,\cP^\vee),
\end{equation}
generalizing the one defined in (\cite{ag2} 2.10.2). 
\end{rema} 

\begin{exemple}\label{p2-qfc220}
Let $\cD$ be a $\mU$-category, $\phi\colon \cD\rightarrow [1]$ a functor, $\varphi\colon \cD_0\rightarrow \hcD_1$ the associated functor defined in \ref{p2-qfc7}. 
We assume that $\varphi$ is a $\mU$-functor \eqref{p2-cmt0}; so $\phi$ satisfies condition \ref{p2-qfc16a}. 
We denote by $\varphi_\rp\colon \hcD_1\rightarrow \hcD_0$ the associated functor \eqref{p2-cmt1a}. 

Recall that we have a canonical isomorphism \eqref{p2-qfc7d}
\begin{equation}\label{p2-qfc220a}
\upgamma\colon \varphi\stackrel{\sim}{\rightarrow} \alpha_{1,\rp}\circ \tth_\cD\circ \alpha_0,
\end{equation}
which is none other that the morphism \eqref{p2-qfc15a} for $\lambda\colon 1\rightarrow 0$ relative to $\phi$. 
As in \eqref{p2-qfc15d}, we deduce a morphism of functors 
\begin{equation}\label{p2-qfc220b}
\upgamma_{\rp}\colon \alpha_{0,\rp}\rightarrow \varphi_\rp\circ \alpha_{1,\rp}. 
\end{equation}

We consider the category $\cP(\cD/[1])$ defined in \ref{p2-qfc20}, 
of triples $(F_0,F_1,\gamma\colon F_0\rightarrow \varphi_\rp(F_1))$, where $F_i\in \ob(\hcD_i)$, for $i=0,1\in \ob([1])$. 
By \ref{p2-qfc21}, the functor 
\begin{equation}\label{p2-qfc220c}
\iota_\phi\colon 
\begin{array}[t]{clcr}
\hcD&\rightarrow&\cP(\cD/[1])\\
F&\mapsto& (\alpha_{0,\rp}(F),\alpha_{1,\rp}(F),\upgamma_{\rp}(F))
\end{array}
\end{equation}
is an equivalence of categories. 
\end{exemple}

\begin{exemple}\label{p2-qfc221}
Let $\cD,\cD'$ be two $\mU$-categories, $\kappa\colon \cD'\rightarrow \cD$, $\phi\colon \cD\rightarrow [1]$ two functors, 
$\phi'=\phi\circ \kappa\colon \cD'\rightarrow [1]$. 
We assume that the functor $\varphi\colon \cD_0\rightarrow \hcD_1$ (resp.\ $\varphi'\colon \cD'_0\rightarrow \hcD'_1$)
associated with $\phi$ (resp.\ $\phi'$) \eqref{p2-qfc7} is a $\mU$-functor \eqref{p2-cmt0}. 
We take again the notation of \ref{p2-qfc220} for $\phi$ and consider the similar constructions associated with $\phi'$, 
that we denote by the same notation equipped with a $\prime$ exponent. 
For $i=0,1\in \ob([1])$, $\kappa$ induces a functor $\kappa_{i}\colon \cD'_i\rightarrow \cD_i$ that fits into a strictly commutative diagram 
\begin{equation}\label{p2-qfc221a}
\xymatrix{
{\cD'_i}\ar[r]^{\alpha'_i}\ar[d]_{\kappa_i}&{\cD'}\ar[d]^\kappa\\
{\cD_i}\ar[r]^{\alpha_i}&{\cD,}}
\end{equation}
where the horizontal arrows are the canonical functors. We denote by $\kappa_{\rp}\colon \hcD\rightarrow \hcD'$  
(resp.\ $\kappa_{i,\rp}\colon \hcD_i\rightarrow \hcD'_i$) the functor defined by composition with $\kappa$ (resp.\ $\kappa_i$). 

The functor $\kappa$ defines a morphism \eqref{p2-qfc7b}
\begin{equation}\label{p2-qfc221b}
\upeta\colon \varphi' \rightarrow \kappa_{1,\rp}\circ \varphi\circ \kappa_0.
\end{equation}
The functors $\phi$ and $\phi'$ define two isomorphisms \eqref{p2-qfc220a}
\begin{eqnarray}
\upgamma\colon \varphi &\stackrel{\sim}{\rightarrow}& \alpha_{1,\rp}\circ \tth_\cD\circ \alpha_0,\label{p2-qfc221c1}\\
\upgamma'\colon \varphi' &\stackrel{\sim}{\rightarrow}& \alpha'_{1,\rp}\circ \tth_{\cD'}\circ \alpha'_0.\label{p2-qfc221c2} 
\end{eqnarray}
We have a morphism 
\begin{equation}\label{p2-qfc221d}
\uprho\colon \tth_{\cD'}\rightarrow \kappa_\rp\circ \tth_\cD\circ \kappa
\end{equation}
defined for $X',Y'\in \ob(\cD')$ by the map induced by $\kappa$
\begin{equation}\label{p2-qfc221e}
\Hom_{\cD'}(X',Y')\rightarrow \Hom_\cD(\kappa(X'),\kappa(Y')). 
\end{equation}
We check easily that the diagram 
\begin{equation}\label{p2-qfc221f}
\xymatrix{
{\varphi'}\ar[r]^-(0.5){\upeta}\ar[d]_{\upgamma'}&{\kappa_{1,\rp}\circ\varphi\circ\kappa_0}\ar[rr]^-(0.5){\kappa_{1,\rp}\circ\upgamma\circ\kappa_0}&&
{\kappa_{1,\rp}\circ\alpha_{1,\rp}\circ\tth_\cD\circ \alpha_0 \circ\kappa_0}\ar[d]\\
{\alpha'_{1,\rp}\circ \tth_{\cD'}\circ \alpha'_0}\ar[rrr]^-(0.5){\alpha'_{1,\rp}\circ \uprho\circ \alpha'_0}&&&
{\alpha'_{1,\rp}\circ \kappa_\rp\circ \tth_\cD\circ \kappa\circ \alpha'_0,}}
\end{equation}
where the unlabelled arrow is the isomorphism induced by the strictly commutative diagrams \eqref{p2-qfc221a}, is commutative. 

Let $F\in \ob(\hcD_1)$, $X\in \ob(\cD'_0)$. The composed morphism 
\begin{equation}\label{p2-qfc221g}
\Hom_{\hcD_1}(\varphi(\kappa_0(X)),F)\rightarrow \Hom_{\hcD'_1}(\kappa_{1,\rp}(\varphi(\kappa_0(X))),\kappa_{1,\rp}(F))
\rightarrow \Hom_{\hcD'_1}(\varphi'(X),\kappa_{1,\rp}(F)),
\end{equation}
where the first morphism is induced by the functor $\kappa_{1,\rp}$ and the second one is defined by composition with $\upeta$ \eqref{p2-qfc221b},
defines a functorial morphism $\kappa_{0,\rp}(\varphi_\rp(F))(X)\rightarrow \varphi'_\rp(\kappa_{1,\rp}(F))(X)$. 
We thus obtain a morphism of functors 
\begin{equation}\label{p2-qfc221h}
\upeta_{\rp}\colon \kappa_{0,\rp}\circ \varphi_\rp\rightarrow \varphi'_\rp\circ \kappa_{1,\rp}. 
\end{equation}
It follows from \ref{p2-qfc221f} that the diagram 
\begin{equation}\label{p2-qfc221i}
\xymatrix{
{\kappa_{0,\rp}\circ \alpha_{0,\rp}}\ar[r]^-(0.5){\kappa_{0,\rp}\circ \upgamma_\rp}\ar[d]
&{\kappa_{0,\rp}\circ \varphi_\rp\circ \alpha_{1,\rp}}\ar[r]^-(0.5){\upeta_\rp\circ \alpha_{1,\rp}}&
{\varphi'_\rp\circ \kappa_{1,\rp}\circ \alpha_{1,\rp}}\ar[d]\\
{\alpha'_{0,\rp}\circ \kappa_\rp}\ar[rr]^{\upgamma'_\rp\circ \kappa_\rp}&&{\varphi'_\rp\circ \alpha'_{1,\rp}\circ \kappa_\rp,}}
\end{equation}
where the vertical arrows are the isomorphisms induced by \eqref{p2-qfc221a}, is commutative. 
We deduce a strictly commutative diagram 
\begin{equation}\label{p2-qfc221j}
\xymatrix{
{\hcD}\ar[r]^-(0.5){\kappa_\rp}\ar[d]_-(0.4){\iota_{\phi}}&{\hcD'}\ar[d]^-(0.4){\iota_{\phi'}}\\
{\cP(\cD/[1])}\ar[r]^-(0.5){\kappa_\cP}&{\cP(\cD'/[1]),}}
\end{equation}
where $\kappa_\cP$ is the functor defined by 
\begin{equation}\label{p2-qfc221k}
\kappa_\cP(F_0,F_1,\gamma)=(\kappa_{0,\rp}(F_0),\kappa_{1,\rp}(F_1),\gamma'),
\end{equation}
and $\gamma'$ is the composed morphism 
\begin{equation}\label{p2-qfc221l}
\gamma'\colon \xymatrix{
{\kappa_{0,\rp}(F_0)}\ar[r]^-(0.5){\kappa_{0,\rp}(\gamma)}&
{\kappa_{0,\rp}(\varphi_\rp(F_1))}\ar[r]^{\upeta_\rp(F_1)}&{\varphi'_\rp(\kappa_{1,\rp}(F_1)).}}
\end{equation}
\end{exemple}

\subsection{}\label{p2-qfc23}
Let $\psi\colon \cI'\rightarrow \cI$ be a functor. We denote by $\cC'$ the fiber product of $\cC$ and $\cI'$ over $\cI$ 
and by $\pi'$ and $\Psi$ the canonical projections:
\begin{equation}\label{p2-qfc23a}
\xymatrix{
{\cC'}\ar[r]^\Psi\ar[d]_{\pi'}&{\cC}\ar[d]^\pi\\
{\cI'}\ar[r]^\psi&{\cI.}}
\end{equation}
Following (\cite{sga1} VI § 3), we say that {\em $\pi'$ is deduced from $\pi$ by base change by $\psi$}. Recall that we have 
\begin{eqnarray}
\ob(\cC')&=&\ob(\cC)\times_{\ob(\cI)}\ob(\cI'),\label{p2-qfc23b}\\
\arr(\cC')&=&\arr(\cC)\times_{\arr(\cI)}\arr(\cI'),\label{p2-qfc23c}
\end{eqnarray}
where $\arr$ stands for the set of arrows. For every $i'\in \ob(\cI')$, 
the fiber category $\cC'_{i'}$ is canonically equivalent to the fiber category $\cC_{\psi(i')}$; they will be identified in the following.
Let $f'\colon j'\rightarrow i'$ be a morphism of $\cI'$, $X\in \ob(\cC_{\psi(i')})$, $Y\in \ob(\cC_{\psi(j')})$. Considering $X$ and $Y$ as objects of $\cC'$, 
we have 
\begin{equation}\label{p2-qfc23d}
\Hom_{\cC'/f'}(Y,X)=\Hom_{\cC/\psi(f')}(Y,X).
\end{equation}
Therefore, the generalized inverse image functor by $f'$ in $\cC'$, 
\begin{equation}\label{p2-qfc23e}
f'^*\colon \cC'_{i'}\rightarrow \hcC'_{j'},
\end{equation}
identifies with the generalized inverse image functor by $\psi(f')$ in $\cC$, $\psi(f')^*\colon \cC_{\psi(i')}\rightarrow \hcC_{\psi(j')}$. 
In particular, the functor $\pi'$ satisfies the conditions \ref{p2-qfc16a} and \ref{p2-qfc16b}. We can therefore define the adjoint functors
$f'_\rp\colon \hcC'_{j'}\rightarrow \hcC'_{i'}$ \eqref{p2-qfc16c} and $f'^\rp\colon \hcC'_{i'}\rightarrow \hcC'_{j'}$ \eqref{p2-qfc16d}. 
For all composable morphisms $g'\colon \ell'\rightarrow j'$ and $f'\colon j'\rightarrow i'$ of $\cI'$, we immediately check that 
the morphism $\delta'_{f',g'}\colon f'^* \rightarrow g'_\rp\circ (f'\circ g')^*$ defined relatively to the functor $\pi'$ \eqref{p2-qfc4g}
identifies with $\delta_{\psi(f'),\psi(g')}$ defined relatively to the functor $\pi$. Hence, $c'^{f',g'}_\rp=c^{\psi(f'),\psi(g')}_\rp$ \eqref{p2-qfc4e}, 
$c'^\rp_{f',g'}=c^\rp_{\psi(f'),\psi(g')}$ \eqref{p2-qfc4f} and $c'^*_{f',g'}=c^*_{\psi(f'),\psi(g')}$ \eqref{p2-qfc4c}. 
In particular, if the functor $\pi$ is quasi-fibering, then so is $\pi'$.

We denote by $\cP(\cC'/\cI')$ the category associated with the functor $\pi'$ in \ref{p2-qfc20}. 
The functor $\Psi\colon \cC'\rightarrow \cC$ induces by composition a functor $\Psi_\rp\colon \hcC\rightarrow \hcC'$ \eqref{p2-cmt4a}.  
We immediately check that we have a strictly commutative diagram 
\begin{equation}\label{p2-qfc23f}
\xymatrix{
{\hcC}\ar[r]^-(0.5){\Psi_\rp}\ar[d]_-(0.4){\iota}&{\hcC'}\ar[d]^-(0.4){\iota'}\\
{\cP(\cC/\cI)}\ar[r]^-(0.5){\Psi_\cP}&{\cP(\cC'/\cI'),}}
\end{equation}
where $\iota$ and $\iota'$ are the equivalences of categories \eqref{p2-qfc20d} and $\Psi_\cP$ is the functor defined by 
\begin{equation}\label{p2-qfc23g}
\Psi_\cP((F_i,\gamma_f)_{i,f})=(F_{\psi(i')},\gamma_{\psi(f')})_{i',f'}.
\end{equation}

\subsection{}\label{p2-qfc240}
Let $\fa\colon \cI\rightarrow [1]$ be a functor, where $[1]$ is the category defined in \ref{p2-qfc7}, 
$\fb=\fa\circ \pi \colon \cC\rightarrow [1]$. We consider $\pi\colon \cC\rightarrow \cI$ as a $[1]$-functor. 
For $o=1,0\in \ob([1])$, we denote by $\pi_{(o)}\colon \cC_{(o)}\rightarrow \cI_{(o)}$ the fiber of $\pi$ above $o$ and by 
$\pi_{(o),\rp}\colon \hcI_{(o)}\rightarrow \hcC_{(o)}$ the functor defined by composition with $\pi_{(o)}$. 
We use parentheses $()$ to distinguish the fibers of $\fb$ from those of $\pi$. 
Let $\Phi\colon \cC_{(0)}\rightarrow \hcC_{(1)}$ and $\varphi\colon \cI_{(0)}\rightarrow \hcI_{(1)}$ be the inverse images by the morphism 
$\lambda\colon 1\rightarrow 0$ of $[1]$ in $\cC$ and $\cI$, respectively. 
By \eqref{p2-qfc7c}, the functor $\pi$ defines a morphism of functors
\begin{equation}\label{p2-qfc240a}
\Phi\rightarrow \pi_{(1),\rp}\circ \varphi \circ \pi_{(0)}.
\end{equation}

We assume that $\Phi\colon \cC_{(0)}\rightarrow \hcC_{(1)}$ is a $\mU$-functor \eqref{p2-cmt0}. 
We denote by $\Phi_\rp\colon \hcC_{(1)}\rightarrow \hcC_{(0)}$ the functor associated with $\Phi$ \eqref{p2-cmt1a}. 

For all $i\in \ob(\cI_0)$ and $j\in \ob(\cI_1)$, we have a strictly commutative diagram \eqref{p2-qfc7d}
\begin{equation}\label{p2-qfc240b}
\xymatrix{
{\cC_{i}}\ar[r]^-(0.5){\alpha_i^{(0)}}\ar[rd]_{\alpha_i}&
{\cC_{(0)}}\ar[r]^-(0.5){\Phi}\ar[d]|{\alpha_{(0)}}&{\hcC_{(1)}}\ar[r]^-(0.5){\alpha_{j,\rp}^{(1)}}&{\hcC_j}\\
&{\cC}\ar[r]^-(0.5)\tth&{\hcC}\ar[u]|{\alpha_{(1),\rp}}\ar[ru]_{\alpha_{j,\rp}}&}
\end{equation}
where the arrows of the left (resp.\ right) triangle are the canonical functors (resp.\ the functors defined by composition with the canonical functors). 
Hence, for every morphism $f\colon j\rightarrow i$ of $\cI$ above $\lambda\colon 1\rightarrow 0$
(i.e., $\fa(f)=\lambda$), the morphisms $\upgamma_f$ \eqref{p2-qfc15a} and $\upgamma_{f,\rp}$ \eqref{p2-qfc15d} induce 
morphisms that we denote again (abusively) by 
\begin{eqnarray}
\upgamma_f\colon f^*&\rightarrow &\alpha_{j,\rp}^{(1)}\circ \Phi\circ \alpha_i^{(0)},\label{p2-qfc240c}\\
\upgamma_{f,\rp}\colon  \alpha_{i,\rp}^{(0)}\circ \alpha_{(0),\rp}&\rightarrow &f_\rp\circ \alpha_{j,\rp}^{(1)}\circ \alpha_{(1),\rp}.\label{p2-qfc240d}
\end{eqnarray}

By \eqref{p2-qfc19a}, for every commutative diagram of $\cI$
\begin{equation}\label{p2-qfc240e}
\xymatrix{
{j'}\ar[r]^{f'}\ar[d]_v\ar[rd]^g&{i'}\ar[d]^u\\
{j}\ar[r]^f&{i,}}
\end{equation}
where $u$ (resp.\ $v$) is above $\id_0$ (resp.\ $\id_1$), we have strictly commutative diagrams 
\begin{equation}\label{p2-qfc240f}
\xymatrix{
{\alpha_{i,\rp}^{(0)}\circ \alpha_{(0),\rp}}\ar[rr]^-(0.5){\upgamma_{f,\rp}}\ar[d]_{\upgamma_{g,\rp}}&&
{f_\rp\circ \alpha_{j,\rp}^{(1)}\circ \alpha_{(1),\rp}}\ar[d]^{f_\rp\circ \upgamma_{v,\rp}}\\
{g_\rp\circ \alpha_{j',\rp}^{(1)}\circ \alpha_{(1),\rp}}\ar[rr]^-(0.5){c_\rp^{f,v}\circ \alpha_{j',\rp}^{(1)}\circ \alpha_{(1),\rp}}&&
{f_\rp\circ v_\rp\circ \alpha_{j',\rp}^{(1)}\circ \alpha_{(1),\rp},}}
\end{equation}
\begin{equation}\label{p2-qfc240g}
\xymatrix{
{\alpha_{i,\rp}^{(0)}\circ \alpha_{(0),\rp}}\ar[rr]^-(0.5){\upgamma_{u,\rp}}\ar[d]_{\upgamma_{g,\rp}}&&
{u_\rp\circ \alpha_{i',\rp}^{(0)}\circ \alpha_{(0),\rp}}\ar[d]^{u_\rp\circ \upgamma_{f',\rp}}\\
{g_\rp\circ \alpha_{j',\rp}^{(1)}\circ \alpha_{(1),\rp}}\ar[rr]^-(0.5){c_\rp^{u,f'}\circ \alpha_{j',\rp}^{(1)}\circ \alpha_{(1),\rp}}&&
{u_\rp\circ f'_\rp\circ \alpha_{j',\rp}^{(1)}\circ \alpha_{(1),\rp}.}}
\end{equation}

\subsection{}\label{p2-qfc241}
We keep the assumptions and notation of \ref{p2-qfc240}. 
For $o=0,1\in \ob([1])$, we see immediately that the functor
$\pi_{(o)}\colon \cC_{(o)}\rightarrow \cI_{(o)}$ is deduced from $\pi$ by base change by the canonical functor
$\psi_{(o)}\colon \cI_{(o)}\rightarrow \cI$ \eqref{p2-qfc23}, and $\alpha_{(o)}\colon \cC_{(o)}\rightarrow \cC$ is the canonical projection. 
Hence, by \ref{p2-qfc23}, $\pi_{(o)}$ satisfies conditions \ref{p2-qfc16a} and \ref{p2-qfc16b}.
We can then consider the category $\cP(\cC_{(o)}/\cI_{(o)})$ defined in \ref{p2-qfc20}. We have a strictly commutative diagram \eqref{p2-qfc23f}
\begin{equation}\label{p2-qfc241b}
\xymatrix{
{\hcC}\ar[r]^-(0.5){\alpha_{(o),\rp}}\ar[d]_-(0.4){\iota}&{\hcC_{(o)}}\ar[d]^-(0.4){\iota_{(o)}}\\
{\cP(\cC/\cI)}\ar[r]^-(0.5){\alpha_{(o),\cP}}&{\cP(\cC_{(o)}/\cI_{(o)}),}}
\end{equation}
where the vertical arrows are the functors \eqref{p2-qfc20d}, which are equivalences of categories by \ref{p2-qfc21}, and 
$\alpha_{(o),\cP}$ is the functor defined by restriction to $\cI_{(o)}$ \eqref{p2-qfc23g}. 

Since the functor $\fb\colon \cC\rightarrow [1]$ satisfies condition \ref{p2-qfc16a}, we can also consider the category $\cP(\cC/[1])$. 
By \ref{p2-qfc21}, we have an equivalence of categories 
\begin{equation}\label{p2-qfc241a}
\upiota\colon \hcC\stackrel{\sim}{\rightarrow} \cP(\cC/[1]). 
\end{equation}

For any $F\in \ob(\hcC)$, we set $F_{(o)}=\alpha_{(o),\rp}(F)$ and $\uF_{(o)}=\iota_{(o)}(F_{(o)})=(F_{(o),i},\gamma_{u,F_{(o)}})_{i,u\in \cI_{(o)}}$. 
We have a canonical morphism $\lambda_F\colon F_{(o)}\rightarrow \Phi_\rp(F_{(1)})$, and  
\begin{equation}\label{p2-qfc241c}
\upiota(F)=(F_{(o)},F_{(1)},\lambda_F). 
\end{equation}
For every morphism $f\colon j\rightarrow i$ of $\cI$ above $\lambda$, the morphism $\gamma_{f,\rp}$ \eqref{p2-qfc240d} induces a morphism of $\hcC_i$
\begin{equation}\label{p2-qfc241d}
\gamma_{f,F}\colon F_{(0),i}\rightarrow f_\rp(F_{(1),j}).
\end{equation}
By \eqref{p2-qfc240f} and \eqref{p2-qfc240g}, for every commutative diagram of $\cI$
\begin{equation}\label{p2-qfc241e}
\xymatrix{
{j'}\ar[r]^{f'}\ar[d]_v\ar[rd]^g&{i'}\ar[d]^u\\
{j}\ar[r]^f&{i,}}
\end{equation}
where $f'$, $g$ and $f$ are above $\lambda$, the diagrams 
\begin{equation}\label{p2-qfc241f}
\xymatrix{
{F_{(0),i}}\ar[r]^-(0.5){\gamma_{f,F}}\ar[d]_{\gamma_{g,F}}&{f_\rp(F_{(1),j})}\ar[d]^{f_\rp(\gamma_{v,F_{(1)}})}\\
{g_\rp(F_{(1),j'})}\ar[r]^-(0.5){c_\rp^{f,v}}&{f_\rp(v_\rp(F_{(1),j'}))}}
\end{equation}
\begin{equation}\label{p2-qfc241g}
\xymatrix{
{F_{(0),i}}\ar[r]^-(0.5){\gamma_{u,F_{(0)}}}\ar[d]_{\gamma_{g,F}}&{u_\rp(F_{(0),i'})}\ar[d]^{u_\rp(\gamma_{f',F})}\\
{g_\rp(F_{(1),j'})}\ar[r]^-(0.5){c_\rp^{u,f'}}&{u_\rp(f'_\rp(F_{(1),j'}))}}
\end{equation}
are commutative. We have 
\begin{equation}\label{p2-qfc241h}
\iota(F)=(\uF_{(0)},\uF_{(1)},\gamma_{f,F})_{f},
\end{equation}
where $f$ describes the morphisms of $\cI$ above $\lambda$. These considerations motivate the following definition:

\begin{defi}\label{p2-qfc242}
Under the assumptions of \ref{p2-qfc240} and with the notation of \ref{p2-qfc241}, 
for any objects $\uF=(F_i,\gamma_u)_{i,u\in \cI_{(0)}}$ of $\cP(\cC_{(0)}/\cI_{(0)})$ and $\uG=(G_j,\gamma_v)_{j,v\in \cI_{(1)}}$ of $\cP(\cC_{(1)}/\cI_{(1)})$, 
an {\em $\mI_\varphi$-system of morphisms from $\uF$ to $\uG$} \eqref{p2-cmt100}
is the data for any morphism $f\colon j\rightarrow i$ of $\cI$ above $\lambda\colon 1\rightarrow 0$ 
of a {\em bifunctorial} morphism
\begin{equation}\label{p2-qfc242a}
\gamma_f\colon F_i\rightarrow f_\rp(G_j).
\end{equation}
The bifunctoriality is expressed as follows. For every commutative diagram \eqref{p2-qfc240e} of $\cI$ with $f,f'$ and $g$ above $\lambda$, 
the diagrams 
\begin{equation}\label{p2-qfc242b}
\xymatrix{
{F_{i}}\ar[r]^-(0.5){\gamma_{f}}\ar[d]_{\gamma_{g}}&{f_\rp(G_j)}\ar[d]^{f_\rp(\gamma_{v})}\\
{g_\rp(G_{j'})}\ar[r]^-(0.5){c_\rp^{f,v}}&{f_\rp(v_\rp(G_{j'}))}}
\ \ \ 
\xymatrix{
{F_{i}}\ar[r]^-(0.5){\gamma_{u}}\ar[d]_{\gamma_{g}}&{u_\rp(F_{i'})}\ar[d]^{u_\rp(\gamma_{f'})}\\
{g_\rp(G_{j'})}\ar[r]^-(0.5){c_\rp^{u,f'}}&{u_\rp(f'_\rp(G_{j'}))}}
\end{equation}
are commutative.  
\end{defi}

Recall that the category $\mI_{\varphi}$, defined in \ref{p2-cmt100}, identifies canonically with the category of morphisms of $\cI$ above $\lambda$.

\begin{prop}\label{p2-qfc243}
We keep the assumptions and notation of \ref{p2-qfc240} and \ref{p2-qfc241}. 
Let $F\in \ob(\hcC_{(0)})$, $G \in \ob(\hcC_{(1)})$, $\uF=\iota_{(0)}(F)$ and $\uG=\iota_{(1)}(G)$ \eqref{p2-qfc241b}. 
The following data are then equivalent:
\begin{itemize}
\item[{\rm (i)}] an object $H$ of $\hcC$ such that $\alpha_{(0),\rp}(H)=F$ and $\alpha_{(1),\rp}(H)=G$;
\item[{\rm (ii)}] an object $\uH$ of $\cP(\cC/\cI)$ such that $\alpha_{(0),\cP}(\uH)=\uF$ and $\alpha_{(1),\cP}(H)=\uG$;
\item[{\rm (iii)}] a morphism $F\rightarrow \Phi_\rp(G)$ of $\hcC_{(0)}$; 
\item[{\rm (iv)}] an $\mI_\varphi$-system of morphisms from $\uF$ to $\uG$ in the sense of \ref{p2-qfc242}.
\end{itemize}
Moreover, the equivalence is bifunctorial in $F$ and $G$, and we have $\iota(H)=\uH$ \eqref{p2-qfc241b}. 
\end{prop}

Indeed, data (i) and (iii) are equivalent because $\upiota$ \eqref{p2-qfc241a} is an equivalence of categories. 
Data (ii) and (iv) are obviously equivalent.
Data (i) and (ii) are equivalent because $\iota$, 
$\iota_{(0)}$ and $\iota_{(1)}$ are equivalences of categories by \ref{p2-qfc21}, and \eqref{p2-qfc241b} is strictly commutative.

\begin{rema}\label{p2-qfc244}
Condition \ref{p2-qfc16b} is not used in definition \ref{p2-qfc242} and proposition \ref{p2-qfc243}; only condition \ref{p2-qfc16a} is necessary. 
However, when $\pi$ satisfies condition \ref{p2-qfc16b}, we can reformulate definition \ref{p2-qfc242} in an equivalent way: 
giving an $\mI_\varphi$-system of morphisms from $\uF$ to $\uG$ amounts to giving for any 
morphism $f\colon j\rightarrow i$ of $\cI$ above $\lambda\colon 1\rightarrow 0$, a {\em bifunctorial} morphism
\begin{equation}\label{p2-qfc244a}
\mu_f\colon f^\rp(F_i)\rightarrow G_j.
\end{equation}
The bifunctoriality is expressed as follows. For every commutative diagram \eqref{p2-qfc240e} of $\cI$ with $f,f'$ and $g$ above $\lambda$, 
the diagrams 
\begin{equation}\label{p2-qfc244b}
\xymatrix{
{v^\rp(f^\rp(F_{i}))}\ar[r]^-(0.5){v^\rp(\mu_{f})}\ar[d]_{c^\rp_{f,v}}&{v^\rp(G_j)}\ar[d]^{\mu_{v}}\\
{g^\rp(F_i)}\ar[r]^-(0.5){\mu_{g}}&{G_{j'},}}
\ \ \ 
\xymatrix{
{f'^\rp(u^\rp(F_{i}))}\ar[r]^-(0.5){f'^\rp(\mu_{u})}\ar[d]_{c^\rp_{u,f'}}&{f'^\rp(F_{i'})}\ar[d]^{\mu_{f'}}\\
{g^\rp(F_i)}\ar[r]^-(0.5){\mu_{g}}&{G_{j'},}}
\end{equation}
where $\mu_u\colon u^\rp(F_i)\rightarrow F_{i'}$ (resp.\ $\mu_v$) is the adjoint morphism of $\gamma_u$ (resp.\ $\gamma_v$), are commutative.  
\end{rema}

\begin{rema}\label{p2-qfc245}
If the functor $\pi$ is quasi-fibering \eqref{p2-qfc9}, the bifunctoriality in \ref{p2-qfc242} 
can also be expressed as follows. For every commutative diagram \eqref{p2-qfc240e} of $\cI$ with $f,f'$ and $g$ above $\lambda$, 
the diagram
\begin{equation}\label{p2-qfc245a}
\xymatrix{
{F_{i}}\ar[r]^-(0.5){\gamma_{f}}\ar[d]_{\gamma_{u}}&{f_\rp(G_j)}\ar[r]^-(0.5){f_\rp(\gamma_{v})}&{f_\rp(v_\rp(G_{j'}))}\ar[d]^{(c_\rp^{f,v})^{-1}}\\
{u_\rp(F_{i'})}\ar[r]^-(0.5){u_\rp(\gamma_{f'})}&{u_\rp(f'_\rp(G_{j'}))}\ar[r]^-(0.5){(c_\rp^{u,f'})^{-1}}&{g_\rp(G_{j'})}}
\end{equation}
is commutative. Indeed, since $c_\rp^{f,v}$ et $c_\rp^{u,f'}$ are isomorphisms, \eqref{p2-qfc245a} follows from \eqref{p2-qfc242b}. Conversely,  
\eqref{p2-qfc245a} applied with $u=\id_i$ (resp.\ $v=\id_j$), implies \eqref{p2-qfc242b}.
\end{rema}

\subsection{}\label{p2-qfc249}
We take again the assumptions and notation of \ref{p2-qfc240} and \ref{p2-qfc241}. 
Let $\psi\colon \cI'\rightarrow \cI$ be a functor. We denote by $\cC'$ the fiber product of $\cC$ and $\cI'$ over $\cI$ 
and by $\pi'$ and $\Psi$ the canonical projections:
\begin{equation}\label{p2-qfc249a}
\xymatrix{
{\cC'}\ar[r]^\Psi\ar[d]_{\pi'}&{\cC}\ar[d]^\pi\\
{\cI'}\ar[r]^\psi&{\cI.}}
\end{equation}
We set $\fa'=\fa\circ \psi\colon \cI'\rightarrow [1]$ and $\fb'=\fa'\circ \pi'=\fb\circ \Psi \colon \cC'\rightarrow [1]$. 
We consider for the functors $(\fa',\pi')$ the constructions similar to those for the functors $(\fa,\pi)$ introduced in \ref{p2-qfc240} and \ref{p2-qfc241},  
that we denote by the same notation equipped with a $\prime$ exponent. For $o=0,1\in \ob([1])$, diagram \eqref{p2-qfc249a} induces a commutative diagram
\begin{equation}\label{p2-qfc249b}
\xymatrix{
{\cC'_{(o)}}\ar[r]^{\Psi_{(o)}}\ar[d]_{\pi'_{(o)}}&{\cC_{(o)}}\ar[d]^{\pi_{(o)}}\\
{\cI'_{(o)}}\ar[r]^{\psi_{(o)}}&{\cI_{(o)}.}}
\end{equation}
We denote by $\Psi_{(o),\rp}\colon \hcC_{(o)}\rightarrow \hcC'_{(o)}$ (resp.\ $\Psi_\rp\colon \hcC\rightarrow \hcC'$) the functor defined by composition 
with $\Psi_{(o)}$ (resp.\ $\Psi$). 
Recall that we have a canonical morphism \eqref{p2-qfc221h}
\begin{equation}\label{p2-qfc249c}
\upeta_{\rp}\colon \Psi_{0,\rp}\circ \Phi_\rp\rightarrow \Phi'_\rp\circ \Psi_{1,\rp}. 
\end{equation}

Let $H\in \ob(\hcC)$, $F=\alpha_{(0),\rp}(H)\in \ob(\hcC_{(0)})$ and $G=\alpha_{(1),\rp}(H) \in \ob(\hcC_{(1)})$,
$\uH=\iota(H)$, $\uF=\iota_{(0)}(F)=(F_i,\gamma_u)_{i,u\in \cI_{(0)}}$ and $\uG=\iota_{(1)}(G)=(G_j,\gamma_v)_{j,v\in \cI_{(1)}}$ \eqref{p2-qfc241b}. 
By \eqref{p2-qfc241b}, we have $\uF=\alpha_{(0),\cP}(\uH)$ and $\uG=\alpha_{(1),\cP}(H)$. 
By \ref{p2-qfc243}, $H$ determines (and is determined by any of the following data): 
\begin{itemize}
\item[(i)] a morphism $\gamma\colon F\rightarrow \Phi_\rp(G)$ of $\hcC_{(0)}$; 
\item[(ii)] an $\mI_\varphi$-system of morphisms from $\uF$ to $\uG$ in the sense of \ref{p2-qfc242}, i.e., 
the data for any morphism $f\colon j\rightarrow i$ of $\cI$ above $\lambda\colon 1\rightarrow 0$, 
of a morphism $\gamma_f\colon F_i\rightarrow f_\rp(G_j)$ of $\hcC_i$. 
\end{itemize}

We set $H'=\Psi_\rp(H)$, $F'=\alpha'_{(0),\rp}(H')=\Psi_{(0),\rp}(F)$, $G'=\alpha'_{(1),\rp}(H')=\Psi_{(1),\rp}(G)$, $\uH'=\iota'(H')$, 
$\uF'=\iota'_{(0)}(F')$ and $\uG'=\iota'_{(1)}(G')$. 
By \eqref{p2-qfc23f}, we have $\uH'=\Psi_\cP(\uH)$, $\uF'=\alpha'_{(0),\cP}(\uH')=\Psi_{(0),\cP}(\uF)=(F_{\psi(i)},\gamma_{\psi(u)})_{i,u\in \cI'_{(0)}}$ 
and $\uG'=\alpha'_{(1),\cP}(\uH)=\Psi_{(1),\cP}(\uG)=(G_{\psi(j)},\gamma_{\psi(v)})_{j,v\in \cI'_{(1)}}$.
By \ref{p2-qfc243}, $H'$ determines (and is determined by any of the following data): 
\begin{itemize}
\item[(i')] a morphism $\gamma'\colon F'\rightarrow \Phi'_\rp(G')$ of $\hcC'_{(0)}$, which by \eqref{p2-qfc221j}, is the composition \eqref{p2-qfc221l}
\begin{equation}\label{p2-qfc249d}
\gamma'\colon \xymatrix{
{\Psi_{(0),\rp}(F)}\ar[r]^-(0.5){\Psi_{(0),\rp}(\gamma)}&
{\Psi_{(0),\rp}(\Phi_\rp(G))}\ar[r]^{\upeta_\rp(G)}&{\Phi'_\rp(\Psi_{(1),\rp}(G)).}}
\end{equation}
\item[(ii')] an $\mI_{\varphi'}$-system of morphisms from $\uF'$ to $\uG'$, which by \eqref{p2-qfc23f}, is given 
by the morphisms $\gamma_{\psi(f)}\colon F_{\psi(i)}\rightarrow \psi(f)_\rp(G_{\psi(j)})$ of $\hcC_{\psi(i)}$ 
for all morphism $f\colon j\rightarrow i$ of $\cI'$ above $\lambda\colon 1\rightarrow 0$. 
\end{itemize}

\subsection{}\label{p2-qfc30}
We assume that the following conditions are satisfied: 
\begin{itemize}
\item[(i)] for every $i\in \ob(\cI)$, the fiber category $\cC_i$ is $\mU$-small and is equipped with a topology $\cT_i$; 
\item[(ii)] for every morphism $f\colon j\rightarrow i$ of $\cI$, the functor $f^*\colon \cC_i\rightarrow \hcC_j$ \eqref{p2-qfc1c} is continuous in the sense of \ref{p2-cmt3}.
\end{itemize}

For every morphism $f\colon j\rightarrow i$ of $\cI$, there exists a unique functor $f_\rs\colon \tcC_j\rightarrow \tcC_i$ 
making strictly commutative the following diagram
\begin{equation}\label{p2-qfc30a}
\xymatrix{
{\tcC_j}\ar[r]^-(0.5){f_\rs}\ar[d]_{\tti_j}&{\tcC_i}\ar[d]^{\tti_i}\\
{\hcC_j}\ar[r]^-(0.5){f_\rp}&{\hcC_i.}}
\end{equation}
By \ref{p2-cmt5}, the functor $f_\rs$ admits a left adjoint 
\begin{equation}\label{p2-qfc30b}
f^\rs\colon \tcC_j\rightarrow \tcC_i. 
\end{equation}

For all composable morphisms $g\colon \ell\rightarrow j$ and $f\colon j\rightarrow i$ of $\cI$, the morphism $c^{f,g}_\rp$ \eqref{p2-qfc4e} 
induces a canonical morphism 
\begin{equation}\label{p2-qfc30c}
d^{f,g}_\rs\colon (f\circ g)_\rs \rightarrow f_\rs\circ g_\rs.
\end{equation}
It satisfies properties analogous to those satisfied by $c^{f,g}_\rp$ in \ref{p2-qfc8}. 
The morphism $d^{f,g}_\rs$ induces by adjunction a canonical morphism 
\begin{equation}\label{p2-qfc30d}
d_{f,g}^\rs\colon g^\rs \circ f^\rs \rightarrow (f\circ g)^\rs.
\end{equation}

\begin{defi}\label{p2-cmt45}
Under the assumptions of \ref{p2-qfc30}, a {\em v-presheaf} of $\mU$-sets on $\cC$ is any presheaf $F$ of $\mU$-sets on $\cC$
such that for every $i\in \ob(\cI)$, $\alpha_{i,\rp}(F)$ is a sheaf on $\cC_i$ \eqref{p2-qfc1ab}.
\end{defi} 

We denote by $\hcC_{\rv}$ the category of v-presheaves,
i.e., the full subcategory of $\hcC$ made up  of the v-presheaves, 
and by $\cP_\rv(\cC/\cI)$ the full subcategory of $\cP(\cC/\cI)$ \eqref{p2-cmt20} made up of the collections $(F_i,\gamma_f)_{i,f}$
such that for each $i\in \ob(\cI)$, $F_i$  is a sheaf on $\cC_i$. 
Then, the functor \eqref{p2-qfc20d} induces an equivalence of categories
\begin{equation}\label{p2-cmt45a}
\begin{array}[t]{clcr}
\hcC_\rv&\stackrel{\sim}{\rightarrow} & \cP_\rv(\cC/\cI),\\
F&\mapsto& \{i\mapsto \alpha_{i,\rp}(F)\}.
\end{array}
\end{equation}

\begin{rema}\label{p2-cmt450}
We keep the assumptions and notation of \ref{p2-qfc30} and let $\uF=(F_i,\gamma_u)_{i,u\in \cI_{(0)}}$ be an object 
of $\cP_\rv(\cC_{(0)}/\cI_{(0)})$ \eqref{p2-cmt45} and $\uG=(G_j,\gamma_v)_{j,v\in \cI_{(1)}}$ an object of $\cP_\rv(\cC_{(1)}/\cI_{(1)})$.  
Giving an $\mI_\varphi$-system of morphisms from $\uF$ to $\uG$ \eqref{p2-qfc242} amounts to giving for any 
morphism $f\colon j\rightarrow i$ of $\cI$ above $\lambda\colon 1\rightarrow 0$, a {\em bifunctorial} morphism
\begin{equation}\label{p2-cmt450a}
\mu_f\colon f^\rs(F_i)\rightarrow G_j.
\end{equation}
The bifunctoriality is expressed as follows. For every commutative diagram \eqref{p2-qfc240e} of $\cI$ with $f,f'$ and $g$ above $\lambda$, 
the diagrams 
\begin{equation}\label{p2-cmt450b}
\xymatrix{
{v^\rs(f^\rs(F_{i}))}\ar[r]^-(0.5){v^\rs(\mu_{f})}\ar[d]_{d^\rs_{f,v}}&{v^\rs(G_j)}\ar[d]^{\mu_{v}}\\
{g^\rs(F_i)}\ar[r]^-(0.5){\mu_{g}}&{G_{j'},}}
\ \ \ 
\xymatrix{
{f'^\rs(u^\rs(F_{i}))}\ar[r]^-(0.5){f'^\rs(\mu_{u})}\ar[d]_{d^\rs_{u,f'}}&{f'^\rs(F_{i'})}\ar[d]^{\mu_{f'}}\\
{g^\rs(F_i)}\ar[r]^-(0.5){\mu_{g}}&{G_{j'},}}
\end{equation}
where $\mu_u\colon u^\rs(F_i)\rightarrow F_{i'}$ (resp.\ $\mu_v$) is the adjoint morphism of $\gamma_u\colon F_i\rightarrow u_\rs(F_{i'})$ 
(resp.\ $\gamma_v\colon G_j\rightarrow v_\rs(G_{j'})$), are commutative.  
\end{rema}

\begin{defi}\label{p2-qfc31}
Under the assumptions of \ref{p2-qfc30} and with the same notation, we say that $\cC$ is a {\em quasi-fibered site} over $\cI$ if for
every composable morphisms $g\colon \ell\rightarrow j$ and $f\colon j\rightarrow i$ of $\cI$, the morphism $d^{f,g}_\rs$ \eqref{p2-qfc30c} is an isomorphism. 
It follows then that the morphism $d_{f,g}^\rs$ \eqref{p2-qfc30d}  is an isomorphism. 
\end{defi}

Whenever we consider a quasi-fibered site, we imply that conditions \ref{p2-qfc30}(i)-(ii) are satisfied.
The first one implies conditions \ref{p2-qfc16a} and \ref{p2-qfc16b}.

If the functor $\pi\colon \cC\rightarrow \cI$ is quasi-fibering in the sense of \ref{p2-qfc9} and satisfies conditions \ref{p2-qfc30}(i)-(ii), 
then $\cC$ is a quasi-fibered site over $\cI$.

\begin{lem}\label{p2-cmt42}
We keep the assumptions and notation of \ref{p2-qfc30}.  Let $g\colon \ell\rightarrow j$ and $f\colon j\rightarrow i$ be two composable morphisms of $\cI$. 
Then, the morphisms $d^{f,g}_\rs$ \eqref{p2-qfc30c} and $d_{f,g}^\rs$ \eqref{p2-qfc30d} are isomorphisms if the following conditions are satisfied:
\begin{itemize}
\item[{\rm (i)}] For every morphism $w\colon Z\rightarrow X$ of $\cC$ above $f\circ g$, there exists a covering
 $(z_\alpha\colon Z_\alpha\rightarrow Z)_{\alpha\in \Sigma}$ of $\cC_\ell$ 
and for every $\alpha\in \Sigma$, an object $Y_\alpha$ of $\cC_j$,
a morphism $v_\alpha\colon Z_\alpha\rightarrow Y_\alpha$ of $\cC$ above $g$ and a morphism $u_\alpha\colon Y_\alpha\rightarrow X$ of 
$\cC$ above $f$ such that $w\circ z_\alpha=u_\alpha\circ v_\alpha$. 
\item[{\rm (ii)}] For all $X\in \ob(\cC_i)$, $Y_1,Y_2\in \ob(\cC_j)$, $Z\in \ob(\cC_\ell)$ and morphisms 
\begin{equation}\label{p2-cmt42a}
\xymatrix{
Z\ar[rd]_{v_2}\ar[r]^-(0.4){v_1}&Y_1\ar[rd]^{u_1}&\\
&Y_2\ar[r]^-(0.4){u_2}&X}
\end{equation}
such that $u_1$ and $u_2$ (resp.\ $v_1$ and $v_2$) are above $f$ (resp.\ $g$) and $u_1\circ v_1=u_2\circ v_2$, 
there exist $Y\in \ob(\cC_j)$,  two morphisms $v\colon Z\rightarrow Y$ and $u\colon Y\rightarrow X$ of $\cC$ and two morphisms 
$y_1\colon Y\rightarrow Y_1$ and $y_2\colon Y\rightarrow Y_2$ of $\cC_j$ such that $v_1=y_1\circ v$, $v_2=y_2\circ v$ and $u=u_1\circ y_1=u_2\circ y_2$. 
\end{itemize}
\end{lem}

Indeed, it follows from \ref{p2-qfc400} that the morphism 
\begin{equation}
\tta_\ell\circ c_{f,g}^*\colon \tta_\ell\circ g^\rp\circ f^*\rightarrow \tta_\ell\circ (f\circ g)^*,
\end{equation}
where $\tta_\ell\colon \hcC_\ell\rightarrow \tcC_\ell$ is the ``associated sheaf'' functor, is an isomorphism. 
By \ref{p2-qfc8}(iii), we deduce that $d^{f,g}_\rs$ \eqref{p2-qfc30c} is an isomorphism, and hence by adjunction so is $d_{f,g}^\rs$ \eqref{p2-qfc30d}. 

\begin{prop}\label{p2-cmt44}
We assume that $\cC$ is a quasi-fibered site over $\cI$ in the sense of \ref{p2-qfc31}. Then,
\begin{itemize}
\item[{\rm (i)}] There exists a canonical cleaved and normalized fibered category {\rm (\cite{sga1} VI 7.1)}
\begin{equation}\label{p2-cmt44a}
\cF^\vee\rightarrow \cI^\circ
\end{equation}
such that for every $i\in \ob(\cI)$, the fiber category $\cF^\vee_i$ of $\cF^\vee$ over $i$ is the topos $\tcC_i$, 
and for every morphism $f\colon j\rightarrow i$ of $\cI$, the inverse image functor associated with $f^\circ$ is $f_\rs\colon \tcC_j\rightarrow \tcC_i$
\eqref{p2-qfc30a}. 
\item[{\rm (ii)}] There exists a canonical cleaved and normalized fibered category 
\begin{equation}\label{p2-cmt44b}
\cF\rightarrow \cI
\end{equation}
such that for every $i\in \cI$, the fiber category $\cF_i$ of $\cF$ over $i$ is the topos $\tcC_i$ and for every morphism $f\colon j\rightarrow i$ of $\cI$, 
the associated inverse image functor is $f^\rs\colon \tcC_i\rightarrow \tcC_j$ \eqref{p2-qfc30b}. 
\end{itemize}
\end{prop}

It follows from (\cite{sga1} VI 7.2 and 8) since the isomorphisms $d^{f,g}_\rs$ \eqref{p2-qfc30c} (resp. $d_{f,g}^\rs$ \eqref{p2-qfc30d})
satisfy properties analogous to those in \ref{p2-qfc8}(ii) and (iv).

\begin{defi}[\cite{agt} §VI.5 and \cite{ag2} 2.10.1]\label{p2-cmt46}
A {\em covanishing fibered site} is a cleaved and normalized fibered category $\kappa\colon \cD\rightarrow \cJ$ (\cite{sga1} VI 7.1)
equipped with a topology on $\cJ$ and for every $j\in \ob(\cJ)$, a topology on the fiber category $\cD_j$, 
satisfying the following conditions:
\begin{itemize}
\item[(i)] The base $\cJ$ is a $\mU$-site and fiber products are representable in $\cJ$. 
\item[(ii)] For every $j\in \ob(\cJ)$, $\cD_j$ is $\mU$-small and admits finite inverse limits. 
\item[(iii)] For every morphism $f\colon j'\rightarrow j$ of $\cJ$, the inverse image functor $f^+\colon \cC_j\rightarrow \cC_{j'}$
is continuous and left exact. It therefore defines a morphism of topos that we (abusively) denote also by
$f\colon \tcC_{j'}\rightarrow \tcC_j$ (\cite{sga4} IV 4.9.2); we denote by $f^\rs$ (resp.\ $f_\rs$) the associated pullback (resp.\ direct image) functor. 
\end{itemize}
\end{defi}

Every covanishing fibered site is clearly a quasi-fibered category \eqref{p2-qfc11} and a quasi-fibered site \eqref{p2-qfc31}. 
We can then use the notation of \ref{p2-qfc1} and \ref{p2-qfc30}, which are compatible with the notation above. 

\subsection{}\label{p2-cmt47}
Let $\kappa\colon \cD\rightarrow \cJ$ be a covanishing fibered site \eqref{p2-cmt46}. 
Following (\cite{agt} VI.5.3), The {\em covanishing} topology of $\cD$ is 
the topology generated by the families of coverings $(X_n\rightarrow X)_{n\in \Sigma}$
of the following two types:
\begin{itemize}
\item[(v)] There exists $j\in \ob(\cJ)$ such that $(X_n\rightarrow X)_{n\in \Sigma}$ is
a covering family of $\cD_j$.
\item[(c)] There exists a covering family of morphisms $(f_n\colon j_n\rightarrow j)_{n\in \Sigma}$ of $\cJ$
such that $\kappa(X)=j$ and for every $n\in \Sigma$, $X_n$ is isomorphic to $f_n^+(X)$.
\end{itemize}
The coverings of type (v) are called {\em vertical}, and those of type (c) are called {\em Cartesian}.
The resulting site is called the {\em covanishing site} associated with the covanishing fibered site $\kappa$; it is a $\mU$-site.
The associated topos of sheaves of $\mU$-sets is called
the {\em covanishing topos} associated with the covanishing fibered site $\kappa$.

\subsection{}\label{p2-qfc246}
Let $\kappa\colon \cD\rightarrow \cJ$ be a functor of $\mU$-categories satisfying condition \ref{p2-qfc16a}, 
fitting into a strictly commutative diagram of functors:
\begin{equation}\label{p2-qfc246a}
\xymatrix{
{\cC}\ar[r]^\Phi\ar[d]_{\pi}&{\cD}\ar[d]^{\kappa}\\
{\cI}\ar[r]^\varphi&{\cJ}}
\end{equation}
i.e., $\kappa\circ \Phi=\varphi\circ \pi$. For every $i\in \ob(\cI)$, $\Phi$ induces a functor $\Phi_i\colon \cC_i\rightarrow \cD_{\varphi(i)}$. We denote by 
\begin{equation}\label{p2-qfc246j}
\Phi_{i,\rp}\colon \hcD_{\varphi(i)}\rightarrow \hcC_i
\end{equation}
the functor defined by composition with $\Phi_{i}$ \eqref{p2-cmt4a}. 
Let $\cJ'$ be a full subcategory of $\cJ$, 
$\sigma\colon \cJ'\rightarrow \cJ$ the canonical functor,
that we omit from the notation when there is no risk of confusion. 
We denote by $\kappa'\colon \cD'\rightarrow \cJ'$
the functor deduced from $\kappa$
by base change by $\sigma$, and by $\nu\colon \cD'\rightarrow \cD$ the canonical projection \eqref{p2-qfc23}, 
that we also omit from the notation when there is no risk of confusion. 
We consider the composed functor 
\begin{equation}\label{p2-qfc246b}
\Phi'\colon 
\xymatrix{
{\cC}\ar[r]^{\Phi}&{\cD}\ar[r]^{\tth_\cD}&{\hcD}\ar[r]^{\nu_\rp}&{\hcD',}}
\end{equation}
where $\tth_\cD$ is the canonical functor and $\nu_\rp$ is defined by composition with $\nu$ \eqref{p2-cmt4a}. 
We denote simply by $\cD'\rhd\cC$ the $\mU$-category over the category $[1]$ defined in \ref{p2-qfc7} by the functor $\Phi'$. 

We consider similarly the composed functor 
\begin{equation}\label{p2-qfc246c}
\varphi'\colon 
\xymatrix{
{\cI}\ar[r]^{\varphi}&{\cJ}\ar[r]^{\tth_\cJ}&{\hcJ}\ar[r]^{\sigma_\rp}&{\hcJ',}}
\end{equation}
and denote simply by $\cJ'\rhd\cI$ the $\mU$-category over the category $[1]$ defined by the functor $\varphi'$.

We have a canonical morphism of functors 
\begin{equation}\label{p2-qfc246d}
\Phi'\rightarrow \kappa'_\rp\circ\varphi'\circ \pi,
\end{equation}
defined for all $X\in \ob(\cC)$ and $Y\in \ob(\cD')$, with $i=\pi(X)$ and $j=\kappa'(Y)$, by the bifunctorial map 
\begin{equation}\label{p2-qfc246e}
\Hom_\cD(\nu(Y),\Phi(X))\rightarrow \Hom_\cJ(\sigma(j),\varphi(i))
\end{equation}
induced by $\kappa$ and the relation $\kappa\circ \Phi=\varphi\circ \pi$. By \eqref{p2-qfc7c}, the functors $\pi$ and $\kappa'$ and the 
morphism \eqref{p2-qfc246d} define a $[1]$-functor 
\begin{equation}\label{p2-qfc246f}
\varpi\colon \cD'\rhd\cC\rightarrow \cJ'\rhd \cI.
\end{equation}
We have $\varpi_1=\kappa'\colon \cD'\rightarrow \cJ'$, $\varpi_0=\pi\colon \cC\rightarrow \cI$, and for $X\in \ob(\cC)$ and $Y\in \ob(\cD')$ with $i=\pi(X)$ and $j=\kappa'(Y)$, the diagram 
\begin{equation}\label{p2-qfc246g}
\xymatrix{
{\Hom_{\cD'\rhd\cC}(Y,X)}\ar[r]^\varpi\ar@{=}[d]&{\Hom_{\cJ'\rhd\cI}(j,i)}\ar@{=}[d]\\
{\Hom_{\cD}(Y,\Phi(X))}\ar[r]^\kappa&{\Hom_{\cJ}(j,\varphi(i))}}
\end{equation}
is commutative; we omit here and below the functors $\sigma,\nu$ from the notation. We deduce the following properties:
\begin{itemize}
\item[(i)] For every morphism $u\colon i'\rightarrow i$ of $\cI$, considered as a morphism of $\cJ'\rhd \cI$, the generalized inverse image functor by $u$ 
in $\cD'\rhd \cC$ identifies with the generalized inverse image functor by $u$ in $\cC$, $u^*\colon \cC_i\rightarrow \hcC_{i'}$ \eqref{p2-qfc1c}. 
\item[(ii)] For every morphism $v\colon j'\rightarrow j$ of $\cJ'$, considered as a morphism of $\cJ'\rhd \cI$, the generalized inverse image functor by $v$ 
in $\cD'\rhd \cC$ identifies with the generalized inverse image functor by $v$ in $\cD'$, $v^*\colon \cD_{j}\rightarrow \hcD_{j'}$. 
\item[(iii)] For every $i\in \ob(\cI)$ and $j\in \ob(\cJ')$ and every morphism $f\colon j\rightarrow \varphi(i)$ of $\cJ$, 
considered as a morphism $\tf\colon j\rightarrow i$ of $\cJ'\rhd \cI$, 
the generalized inverse image functor by $\tf$ in $\cD'\rhd \cC$ identifies with the composed functor
\begin{equation}\label{p2-qfc246h}
\xymatrix{
\Phi_f^*\colon {\cC_{i}}\ar[r]^-(0.5){\Phi_{i}}&{\cD_{\varphi(i)}}\ar[r]^-(0.5){f^*}&{\hcD_j.}}
\end{equation}
\end{itemize}

The functor $\varpi$ \eqref{p2-qfc246f} clearly satisfies condition \ref{p2-qfc16a}. 
By \ref{p2-cmt1}, for every $i\in \ob(\cI)$ and $j\in \ob(\cJ')$ and every morphism $f\colon j\rightarrow \varphi(i)$ of $\cJ$,
we associate with $\Phi_f^*$ \eqref{p2-qfc246h} a functor $\Phi_{f,\rp}\colon \hcD_j\rightarrow \hcC_{i}$ which is none other than the composed functor 
\begin{equation}\label{p2-qfc246i}
\xymatrix{
\Phi_{f,\rp}\colon {\hcD_j}\ar[r]^-(0.5){f_\rp}&{\hcD_{\varphi(i)}}\ar[r]^-(0.4){\Phi_{i,\rp}}&{\hcC_{i}.}}
\end{equation}

Let $u\colon i'\rightarrow i$ (resp.\ $v\colon j'\rightarrow j$) be a morphism of $\cI$ (resp.\ $\cJ'$), 
\begin{equation}\label{p2-qfc246k}
\xymatrix{
j'\ar[r]^-(0.5){f'}\ar[d]_{v}\ar[rd]^g&{\varphi(i')}\ar[d]^{\varphi(u)}\\
j\ar[r]^-(0.5)f&{\varphi(i)}}
\end{equation}
a commutative diagram of $\cJ$. We have $\tg=u\circ \tf'=\tf\circ v$ in $\cJ'\rhd \cI$. 
Applying the construction \eqref{p2-qfc4g} to the functor $\varpi$, we deduce two morphisms 
\begin{eqnarray}\label{p2-qfc246l}
\fd_{u,f'}\colon u^*\rightarrow  \Phi_{f',\rp}\circ \Phi_g^*,&&
\fd_{f,v}\colon \Phi_f^*\rightarrow v_\rp\circ \Phi_g^*.
\end{eqnarray}
By \eqref{p2-qfc4e}, they induce morphisms 
\begin{eqnarray}\label{p2-qfc246m}
\fc^{u,f'}_{\rp}\colon \Phi_{g,\rp} \rightarrow u_\rp\circ \Phi_{f',\rp},&&
\fc^{f,v}_{\rp}\colon \Phi_{g,\rp} \rightarrow  \Phi_{f,\rp} \circ v_\rp.
\end{eqnarray}

Applying \ref{p2-qfc8}(iv) to the functor $\varpi$ \eqref{p2-qfc246f}, we deduce that these morphisms satisfy compatibility relations 
that we do not spell out in full; we only give one example:
For a morphism $w\colon j''\rightarrow j'$ of $\cJ'$, setting $h=g\circ w\colon j''\rightarrow \varphi(i)$, we have a commutative diagram \eqref{p2-qfc8c}
\begin{equation}\label{p2-qfc246o}
\xymatrix{
{\Phi_{h,\rp}}\ar[rr]^-(0.5){\fc^{g,w}_\rp}\ar[d]_{\fc^{f,v\circ w}_\rp}&&{\Phi_{g,\rp}\circ w_\rp}\ar[d]^{\fc^{f,v}_\rp\circ w_\rp}\\
{\Phi_{f,\rp}\circ (v\circ w)_\rp}\ar[rr]^-(0.5){\Phi_{f,\rp}\circ c^{v,w}_\rp}&&{\Phi_{f,\rp}\circ v_\rp \circ w_\rp.}}
\end{equation}
 
Writing $\tf=\widetilde{\id}_{\varphi(i)}\circ f$ in $\cJ'\rhd \cI$, we easily check that we have
\begin{equation}\label{p2-qfc246p}
\fd_{\id_{\varphi(i)},f}=\delta_{\id_{\varphi(i)},f}\circ \Phi_i\colon \tth_{\cD_{\varphi(i)}}\circ \Phi_i\rightarrow f_\rp\circ f^* \circ \Phi_i=f_\rp\circ\Phi_f^*,
\end{equation}
where $\tth_{\varphi(i)}\colon \cD_{\varphi(i)}\rightarrow \hcD_{\varphi(i)}$ is the canonical functor and 
$\delta_{\id_{\varphi(i)},f}\colon  \tth_{\cD_{\varphi(i)}}\rightarrow f_\rp\circ f^*$ is the morphism \eqref{p2-qfc4g} relative to the functor $\kappa$. 
We deduce that 
\begin{equation}\label{p2-qfc246q}
\fc_\rp^{\id_{\varphi(i)},f}\colon \Phi_{f,\rp}\rightarrow \Phi_{i,\rp}\circ f_\rp
\end{equation}
is the canonical isomorphism \eqref{p2-qfc246i}. 

Since $\kappa'$ satisfies condition \ref{p2-qfc16a}, we can consider the category $\cP(\cD'/\cJ')$
defined in \ref{p2-qfc20}. We denote by 
\begin{equation}\label{p2-qfc246n}
\iota_{\cD'/\cJ'}\colon \hcD'\rightarrow \cP(\cD'/\cJ')
\end{equation}
the canonical functor \eqref{p2-qfc20d}, which is an equivalence of categories by \ref{p2-qfc21}. 

\begin{defi}\label{p2-qfc29}
We keep the assumptions of \ref{p2-qfc246} and let $\uF=(F_i,\gamma_u)_{i,u\in \cI}$ be an object of $\cP(\cC/\cI)$ \eqref{p2-qfc20}, 
$\uG=(G_j,\gamma_v)_{j,v\in \cJ'}$ an object of $\cP(\cD'/\cJ')$. 
An {\em $\mI_{\varphi'}$-system of morphisms from $\uF$ to $\uG$} \eqref{p2-cmt100} is
the data for any $i\in \ob(\cI)$ and $j\in \ob(\cJ')$ and any morphism $f\colon j\rightarrow \varphi(i)$ of $\cJ$ 
of a {\em bifunctorial} morphism
\begin{equation}\label{p2-qfc29a}
\tgamma_f\colon F_i\rightarrow \Phi_{f,\rp}(G_j).
\end{equation}
The bifunctoriality is expressed as follows. For every commutative diagram \eqref{p2-qfc246k} of $\cJ$ where 
$u\colon i'\rightarrow i$ (resp.\ $v\colon j'\rightarrow j$) is a morphism of $\cI$ (resp.\ $\cJ'$), the diagrams 
\begin{equation}\label{p2-qfc29b}
\xymatrix{
{F_{i}}\ar[r]^-(0.5){\tgamma_{f}}\ar[d]_{\tgamma_{g}}&{\Phi_{f,\rp}(G_j)}\ar[d]^{\Phi_{f,\rp}(\gamma_{v})}\\
{\Phi_{g,\rp}(G_{j'})}\ar[r]^-(0.5){\fc_\rp^{f,v}}&{\Phi_{f,\rp}(v_\rp(G_{j'}))}}
\ \ \ 
\xymatrix{
{F_{i}}\ar[r]^-(0.5){\gamma_{u}}\ar[d]_{\tgamma_{g}}&{u_\rp(F_{i'})}\ar[d]^{u_\rp(\tgamma_{f'})}\\
{\Phi_{g,\rp}(G_{j'})}\ar[r]^-(0.5){\fc_\rp^{u,f'}}&{u_\rp(\Phi_{f',\rp}(G_{j'}))}}
\end{equation}
are commutative.  
\end{defi}

This notion corresponds to the notion of $\mI_{\varphi'}$-system of morphisms from $\uF$ to $\uG$ relatively 
to the functor $\varpi$ \eqref{p2-qfc246f}, in the sense of \ref{p2-qfc242}. 

\begin{prop}\label{p2-qfc251}
We keep the assumptions and notation of \ref{p2-qfc246}, and assume moreover that $\Phi'$ \eqref{p2-qfc246b} is a $\mU$-functor \eqref{p2-cmt0}. 
We denote by $\Phi'_\rp\colon \hcD'\rightarrow \hcC$ the associated functor defined in \eqref{p2-cmt1a}. 
Let $F\in \ob(\hcC)$, $G \in \ob(\hcD')$, $\uF=\iota(F)$ \eqref{p2-qfc20d} and $\uG=\iota_{\cD'/\cJ'}(G)$ \eqref{p2-qfc246n}. 
The following data are then equivalent:
\begin{itemize}
\item[{\rm (i)}] a morphism $F\rightarrow \Phi'_\rp(G)$ of $\hcC$; 
\item[{\rm (ii)}] an $\mI_{\varphi'}$-system of morphisms from $\uF$ to $\uG$ in the sense of \ref{p2-qfc29}.
\end{itemize}
Moreover, the equivalence is bifunctorial in $F$ and $G$. 
\end{prop}

It is a special case of \ref{p2-qfc243}.

\subsection{}\label{p2-qfc255}
We keep the assumptions and notation of \ref{p2-qfc251}. 
Let $\Phi_\rp\colon \hcD\rightarrow \hcC$ be the functor defined by composition with $\Phi\colon \cC\rightarrow \cD$ \eqref{p2-qfc246a} and 
$\mI_\varphi$ the category associated with the functor $\varphi\colon \cI\rightarrow \cJ$, defined in \ref{p2-cmt4}. 
We have the notion of {\em $\mI_{\varphi}$-system of morphisms from an object of $\cP(\cC/\cI)$ to an object of $\cP(\cD/\cJ)$} 
(defined in \ref{p2-qfc29} by taking $\cD'=\cD$ and $\nu=\id_\cD$). 
As for the definition of \eqref{p2-qfc221h}, the relation $\Phi'=\nu_\rp\circ \tth_\cD\circ \Phi$ \eqref{p2-qfc246b} induces a morphism 
\begin{equation}\label{p2-qfc255a}
\upeta_\rp\colon \Phi_\rp\rightarrow \Phi'_\rp\circ \nu_\rp.
\end{equation}
Since $\kappa$ satisfies condition \ref{p2-qfc16a}, we can consider the category $\cP(\cD/\cJ)$ defined in \ref{p2-qfc20} and 
the canonical functor \eqref{p2-qfc20d} 
\begin{equation}\label{p2-qfc255b}
\iota_{\cD/\cJ}\colon \hcD\rightarrow \cP(\cD/\cJ),
\end{equation}
which is an equivalence of categories by \ref{p2-qfc21}. 

Let $F\in \ob(\hcC)$, $G \in \ob(\hcD)$, $\uF=\iota(F)=(F_i,\gamma_u)_{i,u\in \cI}$ \eqref{p2-qfc20d} and 
$\uG=\iota_{\cD/\cJ}(G)=(G_j,\gamma_v)_{j,v\in \cJ}$. 
Let $\gamma\colon F\rightarrow \Phi_\rp(G)$ be a morphism of $\hcC$, and for any morphism $f\colon j\rightarrow \varphi(i)$ of $\cJ$, 
let $\tgamma_f\colon F_i\rightarrow \Phi_{f,\rp}(G_j)$ be the morphism of $\cC_i$, such that the collection of $\tgamma_f$'s
form the $\mI_{\varphi}$-system of morphisms from $\uF$ to $\uG$ associated with $\gamma$ by \ref{p2-qfc251}. 

We set $G'=\nu_\rp(G)$, $\uG'=\iota_{\cD'/\cJ'}(G')$ and 
\begin{equation}\label{p2-qfc255c}
\gamma'=\eta_\rp(G)\circ \gamma \colon F\rightarrow \Phi'_\rp(G'),
\end{equation}
where $\upeta_\rp$ is defined in \eqref{p2-qfc255a}. 
By \eqref{p2-qfc23f}, we have $\uG'=\nu_\cP(\uG)=(G_{j},\gamma_{v})_{j,v\in \cJ'}$. By \ref{p2-qfc249}, the collection of morphisms 
$\tgamma_f\colon F_i\rightarrow \Phi_{f,\rp}(G_j)$ of $\cC_i$, for all $i\in \ob(\cI)$, $j\in \ob(\cJ')$ and morphisms $f\colon j\rightarrow \varphi(i)$ of $\cJ$, 
is the $\mI_{\varphi'}$-system of morphisms from $\uF$ to $\uG'$ associated with $\gamma'$ by \ref{p2-qfc251}.

\subsection{}\label{p2-qfc254} 
We keep the assumptions and notation of \ref{p2-qfc251}. 
For any object $j\in \cJ$ (resp.\ $j\in \cJ'$), we denote by $\alpha_j\colon \cD_j\rightarrow \cD$ (resp.\ $\alpha'_j\colon \cD_j\rightarrow \cD'$) 
the canonical functor and by $\alpha_{j,\rp}\colon \hcD\rightarrow \hcD_j$ (resp.\ $\alpha'_{j,\rp}\colon \hcD'\rightarrow \hcD_j$) the functor defined 
by composition with $\alpha_j$ (resp.\ $\alpha'_j$). 
Let $i\in \ob(\cI)$, $j\in \ob(\cJ')$, $f\colon j\rightarrow \varphi(i)$ a morphism of $\cJ$. 
Consider the diagram 
\begin{equation}\label{p2-qfc254a}
\xymatrix{
{\cC_i}\ar[d]_{\alpha_i}\ar[r]^{\Phi_i}&{\cD_{\varphi(i)}}\ar[r]^{f^*}\ar[d]_{\alpha_{\varphi(i)}}&{\hcD_j}&\\
{\cC}\ar[r]^\Phi&{\cD}\ar[r]^{\tth_\cD}&{\hcD}\ar[r]^-(0.5){\nu_\rp}\ar[u]^{\alpha_{j,\rp}}&{\hcD'}\ar[ul]_{\alpha'_{j,\rp}}}
\end{equation}
whose left (resp.\ central) square is (resp.\ in general not) commutative, and right triangle is commutative. 
The morphism $\upgamma_f$ for $\kappa$ \eqref{p2-qfc15a} induces a morphism 
\begin{equation}\label{p2-qfc254b}
\upmu_f\colon \Phi_f^*\rightarrow  \alpha'_{j,\rp}\circ \Phi'\circ \alpha_i.
\end{equation}
Let $F\in \ob(\hcD')$, $X\in \ob(\cC_i)$. The composed morphism 
\[
\Hom_{\hcD'}(\Phi'\circ \alpha_i(X),F)\rightarrow 
\Hom_{\hcD_j}(\alpha'_{j,\rp}\circ \Phi'\circ \alpha_i(X),\alpha'_{j,\rp}(F)) \rightarrow \Hom_{\hcD_j}(\Phi_f^*(X),\alpha'_{j,\rp}(F)),
\]
where the first morphism is induced by the functor $\alpha'_{j,\rp}$ and the second one is defined by composition with $\upmu_f$ \eqref{p2-qfc254b},
defines a functorial morphism $\alpha_{i,\rp}(\Phi'_{\rp}(F))(X)\rightarrow \Phi_{f,\rp}(\alpha'_{j,\rp}(F))(X)$. 
We thus obtain a morphism of functors 
\begin{equation}\label{p2-qfc254c}
\upmu_{f,\rp}\colon \alpha_{i,\rp} \circ \Phi'_{\rp}\rightarrow \Phi_{f,\rp}\circ \alpha'_{j,\rp}. 
\end{equation}

\begin{prop}\label{p2-qfc256}
We keep the assumptions and notation of \ref{p2-qfc254} and let $F\in \ob(\hcC)$, $G \in \ob(\hcD')$, 
$\uF=\iota(F)=(F_i,\gamma_u)_{i,u\in \cI}$ \eqref{p2-qfc20d}, 
$\uG=\iota_{\cD'/\cJ'}(G)=(G_j,\gamma_v)_{j,v\in \cJ'}$ \eqref{p2-qfc246n}.
Let $\gamma\colon F\rightarrow \Phi'_\rp(G)$ be a morphism of $\hcC$.  
For any $i\in \ob(\cI)$ and $j\in \ob(\cJ')$ and any morphism $f\colon j\rightarrow \varphi(i)$ of $\cJ$, we denote by $\tgamma_f$ the composed morphism
\begin{equation}\label{p2-qfc256a}
\tgamma_f=\upmu_{f,\rp}(G)\circ \alpha_{i,\rp}(\gamma)\colon F_i\rightarrow \alpha_{i,\rp}(\Phi'_\rp(G))\rightarrow \Phi_{f,\rp}(G_j),
\end{equation}
where $\upmu_{f,\rp}$ is the morphism \eqref{p2-qfc254c}. 
Then, the $\tgamma_f$'s form the $\mI_{\varphi'}$-system of morphisms from $\uF$ to $\uG$ associated with $\gamma$ by \ref{p2-qfc251}. 
\end{prop}

Indeed, let $i\in \ob(\cI)$, $j\in \ob(\cJ')$, $f\colon j\rightarrow \varphi(i)$ be a morphism of $\cJ$. We denote by $\cD_j\rhd\cC_i$ the $\mU$-category 
associated with the functor $\Phi_f^*\colon \cC_i\rightarrow \hcD_j$ \eqref{p2-qfc246h}
and by $\varpi_f\colon \cD_j\rhd\cC_i\rightarrow [1]$ the associated functor defined in \ref{p2-qfc7}. Let
\begin{equation}
\iota_f\colon (\cD_j\rhd\cC_i)^\wedge \rightarrow \cP(\cD_j\rhd\cC_i/[1])
\end{equation}
be the functor \eqref{p2-qfc220c}, which is an equivalence of categories by \ref{p2-qfc21}. 

Let $\psi_f\colon [1]\rightarrow \cJ'\rhd \cI$ be the functor sending $1$ to $j$, $0$ to $i$ and the morphism $\lambda\colon 1\rightarrow 0$ of $[1]$ 
to the morphism $\tf\colon j\rightarrow i$ of $\cJ'\rhd \cI$ defined by the morphism $f\colon j\rightarrow \varphi(i)$ of $\cJ$, see \ref{p2-qfc246}. 
By \eqref{p2-qfc7c}, the functors $\alpha'_j$ and $\alpha_i$ and the canonical morphism 
$\upmu_f\colon \Phi_f^*\rightarrow \alpha'_{j,\rp}\circ \Phi'\circ \alpha_i$ \eqref{p2-qfc254b} define a $[1]$-functor 
$\Psi_f\colon \cD_j\rhd\cC_i\rightarrow \cD'\rhd\cC$. We immediately check that the diagram 
\begin{equation}
\xymatrix{
{\cD_j\rhd\cC_i}\ar[r]^-(0.5){\Psi_f}\ar[d]_{\varpi_f}&{\cD'\rhd\cC}\ar[d]^{\varpi}\\
{[1]}\ar[r]^-(0.5){\psi_f}&{\cJ'\rhd \cI}}
\end{equation}
is strictly commutative and Cartesian; so $\varpi_f$ identifies with the functor deduced from $\varpi$ by base change by $\psi_f$, see \ref{p2-qfc23}. 
Observe that the morphism $\upmu_{f,\rp}$ \eqref{p2-qfc254c} coincides with the morphism $\upeta_\rp$ \eqref{p2-qfc249c} associated 
with the $[1]$-functor $\Psi_f$.

By \ref{p2-qfc243}, the morphism $\gamma$ determines a presheaf of $\mU$-sets $H$ on $\cD'\rhd\cC$, which determines the component 
$\tgamma_f\colon F_i\rightarrow \Phi_{f,\rp}(G_j)$ of the $\mI_{\varphi'}$-system of morphisms from $\uF$ to $\uG$ associated with $\gamma$ by \ref{p2-qfc251}. 
Let $\Psi_{f,\rp}(H)$ be the presheaf of $\mU$-sets on $\cD_j\rhd\cC_i$ deduced from $H$ by composition with $\Psi_f$. We have 
\begin{equation}
\iota_f(\Psi_{f,\rp}(H))=(F_i,G_j,\upmu_{f,\rp}(G)\circ \alpha_{i,\rp}(\gamma))
\end{equation}
Indeed, the morphism $\upmu_{f,\rp}(G)\circ \alpha_{i,\rp}(\gamma)$ is none other than the morphism associated with $\gamma$ in \eqref{p2-qfc249d}
relatively to the $[1]$-functor $\Psi_f$. It follows then from \ref{p2-qfc249} that $\tgamma_f=\upmu_{f,\rp}(G)\circ \alpha_{i,\rp}(\gamma)$.

\subsection{}\label{p2-qfc247}
We take again the assumptions and notation of \ref{p2-qfc246}. We assume, moreover, that the functor $\kappa$ satisfies condition \ref{p2-qfc16b} and that for 
every $i\in \ob(\cI)$, the functor $\Phi_{i,\rp}$ \eqref{p2-qfc246j} admits a left adjoint $\Phi^\rp_i\colon \hcC_{i}\rightarrow \hcD_{\varphi(i)}$ such that the diagram 
\begin{equation}
\xymatrix{
{\cC_i}\ar[r]^-(0.5){\Phi_i}\ar[d]_{\tth_i}&{\cD_{\varphi(i)}}\ar[d]^{\tth_{\varphi(i)}}\\
{\hcC_i}\ar[r]^-(0.5){\Phi_i^\rp}&{\hcD_{\varphi(i)}}}
\end{equation}
is commutative up to a canonical isomorphism. 

For every $i\in \ob(\cI)$ and $j\in \ob(\cJ')$ and every morphism $f\colon j\rightarrow \varphi(i)$ of $\cJ$,
the morphism $f_\rp \colon \hcD_j\rightarrow \hcD_{\varphi(i)}$ admits a left adjoint $f^\rp$. 
We denote by $\Phi_f^\rp$ the composed functor
\begin{equation}\label{p2-qfc247a}
\xymatrix{
\Phi_f^\rp\colon {\hcC_i} \ar[r]^-(0.5){\Phi_i^\rp}&{\hcD_{\varphi(i)}}\ar[r]^-(0.5){f^\rp}&{\hcD_j,}}
\end{equation}
which is a left adjoint of $\Phi_{f,\rp}$ \eqref{p2-qfc246i}. The functor $\varpi\colon \cD'\rhd\cC\rightarrow \cJ'\rhd \cI$ \eqref{p2-qfc246f} 
clearly satisfies condition \ref{p2-qfc166b}.

We consider again the commutative diagram \eqref{p2-qfc246k} of $\cJ$ where 
$u\colon i'\rightarrow i$ (resp.\ $v\colon j'\rightarrow j$) is a morphism of $\cI$ (resp.\ $\cJ'$). 
The morphisms \eqref{p2-qfc246l} and \eqref{p2-qfc246m} induce by adjunction 
the morphisms 
\begin{eqnarray}
\fc^*_{u,f'}\colon \Phi^\rp_{f'}\circ u^* \rightarrow \Phi_g^*,&&\fc^*_{f,v}\colon v^\rp\circ \Phi_f^*\rightarrow\Phi_g^*,\label{p2-qfc247b2}\\
\fc^{\rp}_{u,f'}\colon \Phi^\rp_{f'}\circ u^\rp \rightarrow \Phi_g^\rp,&&\fc^{\rp}_{f,v}\colon v^\rp \circ \Phi^\rp_f \rightarrow \Phi_g^\rp.\label{p2-qfc247b1}
\end{eqnarray}

We deduce from \eqref{p2-qfc246p} that we have
\begin{equation}\label{p2-qfc247d}
\fc^*_{\id_{\varphi(i)},f}=c^*_{\id_{\varphi(i)},f}\circ \Phi_i \colon f^\rp\circ \tth_{\cD_{\varphi(i)}}\circ \Phi_i\rightarrow f^*\circ \Phi_i=\Phi^*_f, 
\end{equation} 
where $c^*_{\id_{\varphi(i)},f} \colon f^\rp\circ \tth_{\cD_{\varphi(i)}}\rightarrow f^*$ is the morphism \eqref{p2-qfc4c} relative to the functor $\kappa$,
which is an isomorphism. We deduce from \eqref{p2-qfc246q} that
\begin{equation}\label{p2-qfc247e}
\fc^\rp_{\id_{\varphi(i)},f}\colon f^\rp\circ \Phi_i^\rp \rightarrow \Phi_f^\rp
\end{equation}
is the canonical isomorphism \eqref{p2-qfc247a}.

Taking $f'=\id_{\varphi(i')}$, in view of \eqref{p2-qfc246h}, \eqref{p2-qfc246i} and \eqref{p2-qfc247a}, we deduce canonical morphisms 
\begin{eqnarray}
\fd_{u,\id_{\varphi(i')}}\colon u^*&\rightarrow&\Phi_{i',\rp}\circ \varphi(u)^*\circ \Phi_{i},\label{p2-qfc247c1}\\
\fc^{u,\id_{\varphi(i')}}_{\rp}\colon \Phi_{i,\rp}\circ \varphi(u)_\rp&\rightarrow& u_\rp \circ \Phi_{i',\rp},\label{p2-qfc247c2}\\
\fc_{u,\id_{\varphi(i')}}^{\rp}\colon \Phi_{i'}^\rp \circ u^\rp&\rightarrow& \varphi(u)^\rp\circ \Phi_i^\rp,\label{p2-qfc247c3}\\
\fc_{u,\id_{\varphi(i')}}^*\colon \Phi_{i'}^\rp\circ u^*&\rightarrow& \varphi(u)^*\circ \Phi_i.\label{p2-qfc247c4}
\end{eqnarray}

\begin{prop}\label{p2-qfc26}
Under the assumptions of \ref{p2-qfc247}, the functor $\varpi\colon \cD'\rhd\cC\rightarrow \cJ'\rhd \cI$ \eqref{p2-qfc246f} is quasi-fibering, 
see \ref{p2-qfc9} and \ref{p2-qfc166}, if and only if the functors $\pi\colon \cC\rightarrow \cI$ and $\kappa'\colon \cD'\rightarrow \cJ'$ 
are quasi-fibering and for every morphism $u\colon i'\rightarrow i$ of $\cI$, the morphism \eqref{p2-qfc247c4}
\begin{equation}\label{p2-qfc26a}
\fc_{u,\id_{\varphi(i')}}^*\colon \Phi_{i'}^\rp\circ u^*\rightarrow \varphi(u)^*\circ \Phi_i
\end{equation}
is an isomorphism. 
\end{prop}

Indeed, the functor $\varpi$ is quasi-fibering if and only if the functors $\pi$ and $\kappa'$ are quasi-fibering 
and for every commutative diagram \eqref{p2-qfc246k} of $\cJ$ where 
$u\colon i'\rightarrow i$ (resp.\ $v\colon j'\rightarrow j$) is a morphism of $\cI$ (resp.\ $\cJ'$), the morphisms \eqref{p2-qfc247b2}
\begin{eqnarray}\label{p2-qfc26b}
\fc^*_{u,f'}\colon \Phi^\rp_{f'}\circ u^* \rightarrow \Phi_g^*,&&\fc^*_{f,v}\colon v^\rp\circ \Phi_f^*\rightarrow\Phi_g^*
\end{eqnarray}
are isomorphisms. Taking for $f'=\id_{\varphi(i')}$, we deduce that the condition of the proposition is necessary. 
Conversely, we assume the functors $\pi$ and $\kappa'$ are quasi-fibering and for every morphism $u\colon i'\rightarrow i$ of $\cI$, the morphism 
$\fc_{u,\id_{\varphi(i')}}^*$ \eqref{p2-qfc247c4} is an isomorphism. Consider a commutative diagram \eqref{p2-qfc246k} of $\cJ$ where 
$u\colon i'\rightarrow i$ (resp.\ $v\colon j'\rightarrow j$) is a morphism of $\cI$ (resp.\ $\cJ'$). 
Applying \eqref{p2-qfc8e} to the functor $\varpi$ and using \eqref{p2-qfc247d}, we deduce that the diagram
\begin{equation}\label{p2-qfc26c}
\xymatrix{
{v^\rp\circ f^\rp\circ \tth_{\cD_{\varphi(i)}}\circ \Phi_{i}}\ar[rr]^-(0.5){c^\rp_{f,v}\circ\Phi_{i}}\ar[d]_{v^\rp\circ c^*_{\id_{\varphi(i)},f}\circ \Phi_i}
&&{g^\rp\circ \tth_{\cD_{\varphi(i)}}\circ \Phi_{i}}\ar[d]^{c^*_{\id_{\varphi(i)},g}\circ \Phi_i}\\
{v^\rp\circ f^*\circ \Phi_{i}}\ar[rr]^-(0.5){c^*_{f,v}\circ\Phi_{i}}\ar@{=}[d]&&{g^*\circ \Phi_{i}}\ar@{=}[d]\\
{v^\rp\circ \Phi_f^*}\ar[rr]^-(0.5){\fc^*_{f,v}}&&{\Phi_g^*}}
\end{equation}
is commutative and the vertical arrows are isomorphisms. In particular, we have $\fc^*_{\varphi(u),f'}=c^*_{\varphi(u),f'}\circ \Phi_{i}$. 
Applying again \eqref{p2-qfc8e} to the functor $\varpi$ and using \eqref{p2-qfc247e}, we obtain that the diagram 
\begin{equation}\label{p2-qfc26d}
\xymatrix{
{f'^\rp\circ \Phi_{i'}^\rp\circ u^*}\ar[d]_{f'^\rp\circ \fc_{u,\id_{\varphi(i')}}^*}\ar@{=}[rrr]&&&{\Phi_{f'}^\rp\circ u^*}\ar[d]^{\fc^*_{u,f'}}\\
{f'^\rp \circ \varphi(u)^*\circ \Phi_{i}}\ar[rr]^-(0.5){c^*_{\varphi(u),f'}\circ \Phi_{i}}&&{g^*\circ \Phi_{i}}\ar@{=}[r]&{\Phi_{g}^*}}
\end{equation}
is commutative. 
We deduce that $\fc^*_{u,f'}$ and $\fc^*_{f,v}$ are isomorphisms, and hence $\varpi$ is quasi-fibering. 

\begin{rema}\label{p2-qfc248}
We keep the assumptions of \ref{p2-qfc247} and let $\uF=(F_i,\gamma_u)_{i,u\in \cI}$ be an object of $\cP(\cC/\cI)$, 
$\uG=(G_j,\gamma_v)_{j,v\in \cJ'}$ an object of $\cP(\cD'/\cJ')$. 
Giving an $\mI_{\varphi'}$-system of morphisms from $\uF$ to $\uG$ \eqref{p2-qfc29} amounts to giving 
for any $i\in \ob(\cI)$ and $j\in \ob(\cJ')$ and any morphism $f\colon j\rightarrow \varphi(i)$ of $\cJ$ 
a {\em bifunctorial} morphism
\begin{equation}\label{p2-qfc248a}
\tmu_f\colon \Phi_f^\rp(F_i)\rightarrow G_j.
\end{equation}
The bifunctoriality is expressed as follows. For every commutative diagram \eqref{p2-qfc246k} of $\cJ$ where 
$u\colon i'\rightarrow i$ (resp.\ $v\colon j'\rightarrow j$) is a morphism of $\cI$ (resp.\ $\cJ'$), the diagrams 
\begin{equation}\label{p2-qfc248b}
\xymatrix{
{v^\rp(\Phi_f^\rp(F_{i}))}\ar[r]^-(0.5){v^\rp(\tmu_{f})}\ar[d]_{\fc^\rp_{f,v}}&{v^\rp(G_j)}\ar[d]^{\mu_{v}}\\
{\Phi_g^\rp(F_i)}\ar[r]^-(0.5){\tmu_{g}}&{G_{j'},}}
\ \ \ 
\xymatrix{
{\Phi_{f'}^\rp(u^\rp(F_{i}))}\ar[r]^-(0.5){\Phi_{f'}^\rp(\mu_{u})}\ar[d]_{\fc^\rp_{u,f'}}&{\Phi_{f'}^\rp(F_{i'})}\ar[d]^{\tmu_{f'}}\\
{\Phi_g^\rp(F_i)}\ar[r]^-(0.5){\tmu_{g}}&{G_{j'},}}
\end{equation}
where $\mu_u\colon u^\rp(F_i)\rightarrow F_{i'}$ (resp.\ $\mu_v$) is the adjoint morphism of $\gamma_u$ (resp.\ $\gamma_v$), are commutative.  

As in \ref{p2-qfc245}, if the functor $\varpi$ \eqref{p2-qfc246f} is quasi-fibering \eqref{p2-qfc9}, the bifunctoriality condition above can be replaced by the commutativity of the diagram
\begin{equation}\label{p2-qfc248c}
\xymatrix{
{\Phi_g^\rp(F_i)}\ar[r]^-(0.5){(\fc^\rp_{f,v})^{-1}}\ar[d]_{(\fc^\rp_{u,f'})^{-1}}&{v^\rp(\Phi_f^\rp(F_{i}))}\ar[r]^-(0.5){v^\rp(\tmu_{f})}&
{v^\rp(G_j)}\ar[d]^{\mu_v}\\
{\Phi_{f'}^\rp(u^\rp(F_{i}))}\ar[r]^-(0.5){\Phi_{f'}^\rp(\mu_{u})}&{\Phi_{f'}^\rp(F_{i'})}\ar[r]^-(0.5){\tmu_{f'}}&{G_{j'}.}}
\end{equation}
\end{rema}

\begin{rema}\label{p2-qfc28}
We keep the assumptions and notation of \ref{p2-qfc247}, and assume, moreover, 
that the functors $\pi$ and $\kappa$ are fibering in the sense of (\cite{sga1} VI 6.1). 
For every morphism $u\colon i'\rightarrow i$ of $\cI$ (resp.\ $v\colon j'\rightarrow j$ of $\cJ$), we denote by 
$u^+\colon \cC_i\rightarrow \cC_{i'}$ (resp.\ $v^+\colon \cD_{j}\rightarrow \cD_{j'}$) the inverse image functor by $u$ in $\cC$ (resp.\ $v$ in $\cD$), 
see \ref{p2-qfc2}. 

Let $i\in \ob(\cI)$, $j\in \ob(\cJ')$, $f\colon j\rightarrow \varphi(i)$ a morphism of $\cJ$. 
Considering $f$ as a morphism $\tf\colon j\rightarrow i$ 
of $\cJ'\rhd \cI$, we see that the inverse image functor by $\tf$ in $\cD'\rhd \cC$ exists and identifies with the composed functor
\begin{equation}\label{p2-qfc28a}
\xymatrix{
\Phi_f^+\colon {\cC_{i}}\ar[r]^-(0.5){\Phi_{i}}&{\cD_{\varphi(i)}}\ar[r]^-(0.5){f^+}&{\cD_j.}}
\end{equation}

Let $u\colon i'\rightarrow i$ (resp.\ $v\colon j'\rightarrow j$) be a morphism of $\cI$ (resp.\ $\cJ'$),
\begin{equation}\label{p2-qfc28b}
\xymatrix{
j'\ar[r]^-(0.5){f'}\ar[d]_{v}\ar[rd]^g&{\varphi(i')}\ar[d]^{\varphi(u)}\\
j\ar[r]^-(0.5)f&{\varphi(i)}}
\end{equation}
a commutative diagram of $\cJ$. By \ref{p2-cmt2}(i), the diagram 
\begin{equation}\label{p2-qfc28c}
\xymatrix{
{\cC_{i}}\ar[r]^{u^+}\ar[d]_{\tth_{i}}&{\cC_{i'}}\ar[r]^-(0.5){\Phi_{i'}}\ar[d]^{\tth_{i'}}&
{\cD_{\varphi(i')}}\ar[r]^-(0.5){f'^+}\ar[d]^{\tth_{\cD,\varphi(i')}}&{\cD_{j'}}\ar[d]^{\tth_{\cD,j'}}\\
{\hcC_{i}}\ar[r]^{u^\rp}&{\hcC_{i'}}\ar[r]^-(0.5){\Phi^\rp_{i'}}&{\hcD_{\varphi(i')}}\ar[r]^-(0.5){f'^\rp}&{\hcD_{j'}}}
\end{equation}
is commutative. Therefore, the morphism $\fc^*_{u,f'}$ \eqref{p2-qfc247b2} induces a morphism 
\begin{equation}\label{p2-qfc28d}
\fc^+_{u,f'}\colon \Phi^+_{f'}\circ u^+ \rightarrow \Phi^+_g. 
\end{equation}
Similarly,  the morphism $\fc^*_{f,v}$ \eqref{p2-qfc247b2} induces a morphism 
\begin{equation}\label{p2-qfc28e}
\fc^+_{f,v}\colon v^+\circ \Phi^+_f\rightarrow\Phi^+_g. 
\end{equation}

Taking $f'=\id_{\varphi(i')}$, in view of \eqref{p2-qfc28a}, we deduce a canonical morphism
\begin{equation}\label{p2-qfc28f}
\fc_{u,\id_{\varphi(i')}}^+\colon \Phi_{i'}\circ u^+\rightarrow \varphi(u)^+\circ \Phi_{i}.
\end{equation}

It follows from \ref{p2-qfc11} and \ref{p2-qfc26} that the functor $\varpi \colon \cD'\rhd\cC\rightarrow \cJ'\rhd \cI$ \eqref{p2-qfc246f} is quasi-fibering
if and only if for every morphism $u\colon i'\rightarrow i$ of $\cI$,
the morphism $\fc_{u,\id_{\varphi(i')}}^+$ \eqref{p2-qfc28f} is an isomorphism. 

Assume, moreover, that $\cJ'=\cJ$. We denote by $\ukappa\colon \ucD\rightarrow \cI$
the fibered category deduced from $\kappa$ by base change by $\varphi$ \eqref{p2-qfc23}, and by
\begin{equation}\label{p2-qfc28g}
\Psi\colon \cC\rightarrow \ucD
\end{equation}
the $\cI$-functor induced by \eqref{p2-qfc246a}. We deduce by (\cite{sga1} VI end of §12) that 
the functor $\varpi$ is quasi-fibering if and only if the functor $\Psi$ is Cartesian in the sense of (\cite{sga1} VI 5.2). 
\end{rema}

\subsection{}\label{p2-qfc260}
We take again the assumptions and notation of \ref{p2-qfc246} and assume, moreover, that the functor $\pi$ is quasi-fibering \eqref{p2-qfc9},
the functor $\kappa$ satisfies condition \ref{p2-qfc16b} and is quasi-fibring, 
and for every morphism $u\colon i'\rightarrow i$ of $\cI$ (resp.\ $v\colon j'\rightarrow j$ of $\cJ'$) and every commutative diagram of $\cJ$
\begin{equation}\label{p2-qfc260a}
\xymatrix{
j'\ar[r]^-(0.5){f'}\ar[d]_{v}\ar[rd]^g&{\varphi(i')}\ar[d]^{\varphi(u)}\\
j\ar[r]^-(0.5)f&{\varphi(i),}}
\end{equation}
the morphisms \eqref{p2-qfc246m}
\begin{eqnarray}\label{p2-qfc260b}
\fc^{u,f'}_{\rp}\colon \Phi_{g,\rp} \rightarrow u_\rp\circ \Phi_{f',\rp},&&
\fc^{f,v}_{\rp}\colon \Phi_{g,\rp} \rightarrow  \Phi_{f,\rp} \circ v_\rp
\end{eqnarray}
are isomorphisms. 

For any object $j$ of $\cJ$, we denote by $\cJ'_{/j}$ the full subcategory of $\cJ_{/j}$ made of the morphisms $j'\rightarrow j$ of $\cJ$
such that $j'$ is an object of $\cJ'$. Objects of $\cJ'_{/j}$ will be denoted by $(j',v)$ or $v\colon j'\rightarrow j$. 

Let $(G_j,\gamma_v)_{j,v\in \cJ'}$ be an object of $\cP(\cD'/\cJ')$ \eqref{p2-qfc20}.
For any $i\in \ob(\cI)$ and any morphism $v\colon (j',g)\rightarrow (j,f)$ of $\cJ'_{/\varphi(i)}$,
we denote by $\tau_{v,f}$ the composed morphism
\begin{equation}\label{p2-qfc260c}
\tau_{v,f}\colon \xymatrix{
{\Phi_{f,\rp}(G_j)}\ar[rr]^-(0.5){\Phi_{f,\rp}(\gamma_v)}&&{\Phi_{f,\rp}(v_\rp(G_{j'}))}\ar[rr]^-(0.5){(\fc^{f,v}_{\rp})^{-1}}&&{\Phi_{g,\rp}(G_{j'}).}}
\end{equation}
Applying \eqref{p2-qfc8c} to the functor $\varpi$ \eqref{p2-qfc246f}, we deduce the following compatibility relations:
\begin{itemize}
\item[(i)] For every composable morphisms $w\colon (j'',h)\rightarrow (j',g)$ and $v\colon (j',g)\rightarrow (j,f)$ of $\cJ'_{/\varphi(i)}$, we have
\begin{equation}\label{p2-qfc260d}
\tau_{v\circ w,f}=\tau_{w,g}\circ \tau_{v,f}.
\end{equation}
\item[(ii)] For every commutative diagram of $\cJ$
\begin{equation}\label{p2-qfc260h}
\xymatrix{
j'\ar[rr]^-(0.5){g'}\ar[dd]_{v}\ar[rrdd]^(0.3)g&&{\varphi(i')}\ar[dd]^{\varphi(u)}\\
&&\\
j\ar[rr]^-(0.5)f\ar[rruu]|!{[uu];[rr]}\hole^-(0.3){f'}&&{\varphi(i),}}
\end{equation}
where $u\colon i'\rightarrow i$ (resp.\ $v\colon j'\rightarrow j$) is a morphism of $\cI$ (resp.\ of $\cJ'$), the diagram 
\begin{equation}\label{p2-qfc260i}
\xymatrix{
{\Phi_{f,\rp}(G_j)}\ar[r]^-(0.5){\fc^{u,f'}_{\rp}}\ar[d]_{\tau_{v,f}}&{u_\rp(\Phi_{f',\rp}(G_j))}\ar[d]^{u_\rp(\tau_{v,f'})}\\
{\Phi_{g,\rp}(G_{j'})}\ar[r]^-(0.5){\fc^{u,g'}_{\rp}}&{u_\rp(\Phi_{g',\rp}(G_{j'}))}}
\end{equation}
is commutative. 
\end{itemize} 

Therefore, for any $i\in \ob(\cI)$, we can set 
\begin{equation}\label{p2-qfc260e}
F_i=\underset{\underset{(j,f)\in (\cJ'_{/\varphi(i)})^\circ}{\longleftarrow}}{\lim}\ \Phi_{f,\rp}(G_j), 
\end{equation}
where the transition morphisms are the morphisms $\tau_{v,f}$ \eqref{p2-qfc260c}. 
For any morphism $u\colon i'\rightarrow i$ of $\cI$, since $u_\rp$ commutes with inverse limits, 
the morphisms $\fc^{u,f'}_{\rp}$ \eqref{p2-qfc260b} induce a morphism 
\begin{equation}\label{p2-qfc260f}
\gamma_u\colon F_i\rightarrow u_\rp(F_{i'}).
\end{equation}
It follows from \eqref{p2-qfc8c} applied to the functor $\varpi$ \eqref{p2-qfc246f} that $(F_i,\gamma_u)_{i,u\in \cI}$ form an object of $\cP(\cC/\cI)$. 
We thus define a functor 
\begin{equation}\label{p2-qfc260g}
\Phi'_\cP\colon
\begin{array}[t]{clcr}  
\cP(\cD'/\cJ')&\rightarrow &\cP(\cC/\cI),\\
(G_j,\gamma_v)_{j,v\in \cJ'}&\mapsto &(F_i,\gamma_u)_{i,u\in \cI}.
\end{array}
\end{equation}

\begin{prop}\label{p2-qfc261}
We keep the assumptions and notation of \ref{p2-qfc260}, and assume moreover that $\Phi'$ \eqref{p2-qfc246b} is a $\mU$-functor \eqref{p2-cmt0}. 
We denote by $\Phi'_\rp\colon \hcD'\rightarrow \hcC$ the associated functor defined in \eqref{p2-cmt1a}. 
Then, the diagram 
\begin{equation}\label{p2-qfc261a}
\xymatrix{
{\hcD'}\ar[r]^-(0.5){\Phi'_\rp}\ar[d]_{\iota_{\cD'/\cJ'}}&{\hcC}\ar[d]^{\iota}\\
{\cP(\cD'/\cJ')}\ar[r]^-(0.5){\Phi'_\cP}&{\cP(\cC/\cI),}}
\end{equation}
where the vertical arrows are the equivalences of categories \eqref{p2-qfc20d} and $\Phi'_\cP$ is the functor \eqref{p2-qfc260g}, 
is commutative up to canonical isomorphism. 
\end{prop}

Indeed, let $F\in \ob(\hcC)$, $G \in \ob(\hcD')$, $\uF=\iota(F)$ and $\uG=\iota_{\cD'/\cJ'}(G)$.  
Giving an $\mI_{\varphi'}$-system of morphisms from $\uF$ to $\uG$ in the sense of \ref{p2-qfc29} amounts to giving a morphism 
$\uF\rightarrow \Phi'_\cP(\uG)$ of $\cP(\cC/\cI)$. By \ref{p2-qfc251}, we have a canonical bifunctorial isomorphism 
\begin{equation}
\Hom_{\hcC}(F, \Phi'_{\rp}(G))\stackrel{\sim}{\rightarrow}
\Hom_{\cP(\cC/\cI)}(\uF, \Phi'_\cP(\uG)).
\end{equation}
Since $\iota$ is an equivalence of categories, we deduce an isomorphism 
\begin{equation}
\iota\circ \Phi'_\rp\stackrel{\sim}{\rightarrow} \Phi'_\cP\circ \iota_{\cD'/\cJ'}.
\end{equation}

\subsection{}\label{p2-qfc262}
We take again the assumptions and notation of \ref{p2-qfc247}, and assume, moreover, that the functor $\varpi$ \eqref{p2-qfc246f} is quasi-fibering \eqref{p2-qfc9},
so the morphisms \eqref{p2-qfc246m}, \eqref{p2-qfc247b2} and \eqref{p2-qfc247b1} are isomorphisms, see \ref{p2-qfc26}. 
We consider the category $\mI_{\varphi'}$ defined in \ref{p2-cmt100} and the functor 
\begin{equation}\label{p2-qfc262a}
\tts\colon 
\begin{array}[t]{clcr}
\mI_{\varphi'}&\rightarrow& \cJ',\\
(f\colon j\rightarrow \varphi(i)) &\mapsto& j.
\end{array}
\end{equation}
For any $j\in \ob(\cJ')$, we denote by $\mI_{\varphi'}^j$ the fiber category of the functor $\tts$ above $j$. Objects of $\mI_{\varphi'}^j$ are pairs $(i,f)$, 
where $i\in \ob(\cI)$ and $f\colon j\rightarrow \varphi(i)$ is a morphism of $\cJ$. Let $(i,f)$ and $(i',f')$ be two objects of $\mI_{\varphi'}^j$.
A morphism from $(i',f')$ to $(i,f)$ is a morphism $u\colon i'\rightarrow i$ of $\cI$ such that $\varphi(u)\circ f'=f$. 
For every morphism $v\colon j'\rightarrow j$ of $\cJ'$, we have a functor 
\begin{equation}\label{p2-qfc262b}
\mI_{\varphi'}^v\colon 
\begin{array}[t]{clcr}
\mI_{\varphi'}^{j}&\rightarrow& \mI_{\varphi'}^{j'},\\ 
(i,f)&\mapsto& (i,f\circ v).
\end{array}
\end{equation}

Let $(F_i,\gamma_u)_{i,u\in \cI}$ be an object of $\cP(\cC/\cI)$. 
For any $j\in \ob(\cJ')$ and any morphism $u\colon (i',f')\rightarrow (i,f)$ of $\mI_{\varphi'}^j$, we denote by $\rho_{u,f'}$ the 
composed morphism 
\begin{equation}\label{p2-qfc262c}
\rho_{u,f'}\colon \xymatrix{
{\Phi_{f}^\rp(F_i)}\ar[rr]^-(0.5){(\fc_{u,f'}^\rp)^{-1}}&&{\Phi_{f'}^\rp(u^\rp(F_i))}\ar[rr]^-(0.5){\Phi_{f'}^\rp(\mu_u)}&&{\Phi_{f'}^\rp(F_{i'}),}}
\end{equation}
where $\mu_u\colon u^\rp(F_i)\rightarrow F_{i'}$ is the adjoint of $\gamma_u$. 
Applying \eqref{p2-qfc8d} to the functor $\varpi$ \eqref{p2-qfc246f}, we deduce the following compatibility relations:
\begin{itemize}
\item[(i)] For every composable morphisms $t\colon (i'',f'')\rightarrow (i',f')$ and $u\colon (i',f')\rightarrow (i,f)$ of $\mI_{\varphi'}^j$, we have
\begin{equation}\label{p2-qfc262d}
\rho_{f,u\circ t}=\rho_{f',t}\circ \tau_{f,u}.
\end{equation}
\item[(ii)] For every commutative diagram of $\cJ$
\begin{equation}\label{p2-qfc262e}
\xymatrix{
j'\ar[rr]^-(0.5){g'}\ar[dd]_{v}\ar[rrdd]^(0.3)g&&{\varphi(i')}\ar[dd]^{\varphi(u)}\\
&&\\
j\ar[rr]^-(0.5)f\ar[rruu]|!{[uu];[rr]}\hole^-(0.3){f'}&&{\varphi(i),}}
\end{equation}
where $u\colon i'\rightarrow i$ (resp.\ $v\colon j'\rightarrow j$) is a morphism of $\cI$ (resp.\ of $\cJ'$), the diagram 
\begin{equation}\label{p2-qfc262f}
\xymatrix{
{v^\rp(\Phi_f^\rp(F_i))}\ar[rr]^-(0.5){v^\rp(\rho_{u,f'})}\ar[d]_-(0.5){\fc_{f,v}^{\rp}}&&{v^\rp(\Phi_{f'}^\rp(F_{i'}))}\ar[d]^-(0.5){\fc_{f',v}^{\rp}}\\
{\Phi_{g}^\rp(F_i)}\ar[rr]^-(0.5){\rho_{u,g'}}&&{\Phi_{g'}^\rp(F_{i'})}}
\end{equation}
is commutative. 
\end{itemize} 

Therefore, for any $j\in \ob(\cJ')$, we can set 
\begin{equation}\label{p2-qfc262g}
G_j=\underset{\underset{(i,f)\in (\mI_{\varphi'}^j)^\circ}{\longrightarrow}}{\lim}\ \Phi_f^\rp(F_i), 
\end{equation}
where the transition morphisms are the morphisms $\rho_{u,f'}$ \eqref{p2-qfc262c}. 
For any morphism $v\colon j'\rightarrow j$ of $\cJ'$, since $v^\rp$ commutes with direct limits, 
the morphisms $\fc_{f,v}^{\rp}$ \eqref{p2-qfc247b1} induce a morphism 
\begin{equation}\label{p2-qfc262h}
\mu_v\colon v^\rp(G_j)\rightarrow G_{j'}.
\end{equation}
We denote by $\gamma_v\colon G_j\rightarrow v_\rp(G_{j'})$ the adjoint. 
It follows from \eqref{p2-qfc8d} applied to the functor $\varpi$ \eqref{p2-qfc246f} that $(G_j,\gamma_v)_{j,v\in \cJ'}$ form an object of $\cP(\cD'/\cJ')$. 
We thus define a functor 
\begin{equation}\label{p2-qfc262i}
\Phi'^\cP\colon
\begin{array}[t]{clcr}  
\cP(\cC/\cI)&\rightarrow &\cP(\cD'/\cJ'),\\
(F_i,\gamma_u)_{i,u\in \cI}&\mapsto &(G_j,\gamma_v)_{j,v\in \cJ'}.
\end{array}
\end{equation}

\begin{prop}\label{p2-qfc263}
We keep the assumptions and notation of \ref{p2-qfc262}, and assume moreover that $\Phi'$ \eqref{p2-qfc246b} is a $\mU$-functor \eqref{p2-cmt0}. 
We denote by $\Phi'_\rp\colon \hcD'\rightarrow \hcC$ the associated functor defined in \eqref{p2-cmt1a}. 
We assume further that  $\Phi'_\rp$ admits a left adjoint $\Phi'^\rp\colon \hcC\rightarrow \hcD'$. Then, the diagram 
\begin{equation}\label{p2-qfc263a}
\xymatrix{
{\hcC}\ar[r]^-(0.5){\Phi'^\rp}\ar[d]_{\iota}&{\hcD'}\ar[d]^{\iota_{\cD'/\cJ'}}\\
{\cP(\cC/\cI)}\ar[r]^-(0.5){\Phi'^\cP}&{\cP(\cD'/\cJ'),}}
\end{equation}
where the vertical arrows are the equivalences of categories \eqref{p2-qfc20d} and $\Phi'^\cP$ is the functor \eqref{p2-qfc262i}, 
is commutative up to canonical isomorphism. 
\end{prop}

It follows from \ref{p2-qfc261}. 

\subsection{}\label{p2-qfc250}
We keep the assumptions of \ref{p2-qfc246}, we assume moreover that the following conditions are satisfied:
\begin{itemize}
\item[(i)] the category $\cC$ (resp.\ $\cD$) is equipped with a topology, and that $\cD'$ is a topologically generating subcategory of $\cD$. 
We equip $\cD'$ with the topology induced by that of $\cD$ and denote by $\tcD'$ (resp.\ $\tcC$) the topos of sheaves of $\mU$-sets on 
$\cD'$ (resp.\ $\cC$); 
\item[(ii)] for every $i\in \ob(\cI)$ (resp.\ $j\in \ob(\cJ)$), the category $\cC_i$ (resp.\ $\cD_j$) is equipped with a topology 
having a $\mU$-small topologically generating family; 
\item[(iii)] the categories $\cC$, $\cD'$ and $\cD_j$, for every $j\in \ob(\cJ)$, are $\mU$-small;
\item[(iv)] the functor $\Phi\colon \cC\rightarrow \cD$ is continuous \eqref{p2-qfc246a}; 
\item[(v)] for every morphism $u\colon i'\rightarrow i$ of $\cI$, the functor $u^*\colon \cC_i\rightarrow \hcC_{i'}$ is continuous in the sense of \ref{p2-cmt3};
\item[(vi)] for every morphism $v\colon j'\rightarrow j$ of $\cJ$, the functor $v^*\colon \cD_j\rightarrow \hcD_{j'}$ is continuous;
\item[(vii)] for every $i\in \ob(\cI)$, the functor $\Phi_i\colon \cC_i\rightarrow \cD_{\varphi(i)}$ is continuous. 
\end{itemize}
The functor $\kappa$ satisfies clearly condition \ref{p2-qfc16b}. 
For every $i\in \ob(\cI)$, the category $\cC_i$ is $\mU$-small. Then by \ref{p2-cmt2}(i), the functor $\Phi_{i,\rp}$ \eqref{p2-qfc246j} 
admits a left adjoint $\Phi^\rp_i\colon \hcC_{i}\rightarrow \hcD_{\varphi(i)}$. We take again the notation of \ref{p2-qfc247}. 
 
For every morphism $u\colon i'\rightarrow i$ of $\cI$, there exists a unique functor $u_\rs\colon \tcC_{i'}\rightarrow \tcC_i$ 
making strictly commutative the following diagram
\begin{equation}
\xymatrix{
{\tcC_{i'}}\ar[r]^-(0.5){u_\rs}\ar[d]&{\tcC_i}\ar[d]\\
{\hcC_{i'}}\ar[r]^-(0.5){u_\rp}&{\hcC_i,}}
\end{equation}
where the vertical arrows are the canonical functors. 
By \ref{p2-cmt5}, the functor $u_\rs$ admits a left adjoint 
\begin{equation}
u^\rs\colon \tcC_{i'}\rightarrow \tcC_i. 
\end{equation} 

For every morphism $v\colon j'\rightarrow j$ of $\cJ$, there exists a unique functor $v_\rs\colon \tcC_{j'}\rightarrow \tcC_j$ 
making strictly commutative the following diagram
\begin{equation}
\xymatrix{
{\tcD_{j'}}\ar[r]^-(0.5){v_\rs}\ar[d]&{\tcD_j}\ar[d]\\
{\hcD_{j'}}\ar[r]^-(0.5){v_\rp}&{\hcD_j,}}
\end{equation}
where the vertical arrows are the canonical functors. 
By \ref{p2-cmt5}, the functor $v_\rs$ admits a left adjoint 
\begin{equation}
v^\rs\colon \tcD_{j'}\rightarrow \tcD_j. 
\end{equation} 
 
For every $i\in \ob(\cI)$, there exists a unique functor $\Phi_{i,\rs}\colon \tcD_{\varphi(i)}\rightarrow \tcC_i$ 
making strictly commutative the following diagram
\begin{equation}\label{p2-qfc250a}
\xymatrix{
{\tcD_{\varphi(i)}}\ar[r]^-(0.4){\Phi_{i,\rs}}\ar[d]&{\tcC_i}\ar[d]\\
{\hcD_{\varphi(i)}}\ar[r]^-(0.4){\Phi_{i,\rp}}&{\hcC_i.}}
\end{equation}
By \ref{p2-cmt5}, the functor $\Phi_{i,\rs}$ admits a left adjoint 
\begin{equation}\label{p2-qfc250b}
\Phi_i^\rs\colon\tcC_i \rightarrow  \tcD_{\varphi(i)}. 
\end{equation}

For every $i\in \ob(\cI)$ and $j\in \ob(\cJ')$ and every morphism $f\colon j\rightarrow \varphi(i)$ of $\cJ$, we set 
$\Phi_{f,\rs}$ to be the composed functor 
\begin{equation}
\xymatrix{
\Phi_{f,\rs}\colon {\tcD_j}\ar[r]^-(0.5){f_\rs}&{\tcD_{\varphi(i)}}\ar[r]^-(0.4){\Phi_{i,\rs}}&{\tcC_{i}.}}
\end{equation}
It is the unique functor that makes strictly commutative the following diagram
\begin{equation}
\xymatrix{
{\tcD_{j}}\ar[r]^-(0.4){\Phi_{f,\rs}}\ar[d]&{\tcC_i}\ar[d]\\
{\hcD_{j}}\ar[r]^-(0.4){\Phi_{f,\rp}}&{\hcC_i.}}
\end{equation}
By \ref{p2-cmt5}, the functor $\Phi_{f,\rs}$ admits a left adjoint 
\begin{equation}
\Phi_f^\rs\colon\tcC_i \rightarrow  \tcD_j. 
\end{equation}

We denote by $\Phi_\rp \colon \hcD\rightarrow \hcC$ the functor defined by composition with the continuous functor 
$\Phi\colon \cC\rightarrow \cD$ \eqref{p2-qfc246a} and by  
$\Phi_\rs\colon \tcD\rightarrow \tcC$ the functor 
making strictly commutative the following diagram
\begin{equation}
\xymatrix{
{\tcD}\ar[r]^-(0.5){\Phi_\rs}\ar[d]&{\tcC}\ar[d]\\
{\hcD}\ar[r]^-(0.5){\Phi_\rp}&{\hcC,}}
\end{equation}
where the vertical arrows are the canonical functors. By (\cite{sga4} III 1.3), the functor $\Phi_\rs$ admits a left adjoint 
\begin{equation}
\Phi^\rs\colon \tcC\rightarrow \tcD.
\end{equation} 

By \ref{p2-cmt101}(a), $\Phi'$ \eqref{p2-qfc246b} is a $\mU$-functor. 
Let $\Phi'_\rp\colon \hcD'\rightarrow \hcC$ be the associated functor \eqref{p2-cmt1a} and  
$\Phi'^\rp\colon \hcC\rightarrow \hcD'$ a left adjoint of $\Phi'_\rp$, which exists by \ref{p2-cmt2}. 
By \ref{p2-cmt77}(i), the functor $\Phi'$ is continuous in the sense of \ref{p2-cmt3}. 
Let $\Phi'_\rs\colon \tcD'\rightarrow \tcC$ be the functor 
making strictly commutative the following diagram
\begin{equation}
\xymatrix{
{\tcD'}\ar[r]^-(0.5){\Phi'_\rs}\ar[d]&{\tcC}\ar[d]\\
{\hcD'}\ar[r]^-(0.5){\Phi'_\rp}&{\hcC.}}
\end{equation}
By \ref{p2-cmt5}, the functor $\Phi'_{\rs}$ admits a left adjoint 
\begin{equation}
\Phi'^\rs\colon\tcC \rightarrow  \tcD'. 
\end{equation}

By (\cite{sga4} III 4.1 and its proof), 
the canonical functor $\nu\colon \cD'\rightarrow \cD$ is continuous and cocontinuous,
and the functor $\nu_\rp\colon \hcD\rightarrow \hcD'$, defined by composition with $\nu$, 
induces an equivalence of categories $\nu_\rs\colon \tcD\stackrel{\sim}{\rightarrow} \tcD'$.
By \eqref{p2-cmt77b} and \eqref{p2-cmt77d}, the diagrams 
\begin{equation}
\xymatrix{
{\tcD}\ar[r]^{\Phi_\rs}\ar[d]_-(0.5){\nu_\rs}&{\tcC}\\
{\tcD'}\ar[ru]_{\Phi'_\rs}&}
\ \ \
\xymatrix{
{\tcC}\ar[r]^{\Phi'^\rs}\ar[rd]_{\Phi^\rs}&{\tcD'}\ar[d]^-(0.5){\nu^\rs}\\
&{\tcD}}
\end{equation}
where $\nu^\rs$ is a felt adjoint of $\nu_\rs$,  are commutative up to isomorphisms. 

\begin{rema}\label{p2-qfc252}
We keep the assumptions of \ref{p2-qfc250} and let $\uF=(F_i,\gamma_u)_{i,u\in \cI}$ be an object of $\cP_\rv(\cC/\cI)$ \eqref{p2-cmt45}, 
$\uG=(G_j,\gamma_v)_{j,v\in \cJ'}$ an object of $\cP_\rv(\cD'/\cJ')$. 
Giving an $\mI_{\varphi'}$-system of morphisms from $\uF$ to $\uG$ \eqref{p2-qfc29} amounts to giving 
for any $i\in \ob(\cI)$ and $j\in \ob(\cJ')$ and any morphism $f\colon j\rightarrow \varphi(i)$ of $\cJ$ 
a {\em bifunctorial} morphism
\begin{equation}\label{p2-qfc252a}
\tmu_f\colon \Phi_f^\rs(F_i)\rightarrow G_j.
\end{equation}
The bifunctoriality is expressed as follows. For every commutative diagram \eqref{p2-qfc246k} of $\cJ$ where 
$u\colon i'\rightarrow i$ (resp.\ $v\colon j'\rightarrow j$) is a morphism of $\cI$ (resp.\ $\cJ'$), 
the diagrams 
\begin{equation}\label{p2-qfc252b}
\xymatrix{
{v^\rs(\Phi_f^\rs(F_{i}))}\ar[r]^-(0.5){v^\rs(\tmu_{f})}\ar[d]_{\fd^\rs_{f,v}}&{v^\rs(G_j)}\ar[d]^{\mu_{v}}\\
{\Phi_g^\rs(F_i)}\ar[r]^-(0.5){\tmu_{g}}&{G_{j'},}}
\ \ \ 
\xymatrix{
{\Phi_{f'}^\rs(u^\rs(F_{i}))}\ar[r]^-(0.5){\Phi_{f'}^\rs(\mu_{u})}\ar[d]_{\fd^\rs_{u,f'}}&{\Phi_{f'}^\rs(F_{i'})}\ar[d]^{\tmu_{f'}}\\
{\Phi_g^\rs(F_i)}\ar[r]^-(0.5){\tmu_{g}}&{G_{j'},}}
\end{equation}
where $\fd^\rs_{f,v}$ (resp.\ $\fd^\rs_{u,f'}$) is the morphism induced by $\fc_\rp^{f,v}$ (resp.\ $\fc_\rp^{u,f'}$) \eqref{p2-qfc246m}, 
and $\mu_u\colon u^\rs(F_i)\rightarrow F_{i'}$ (resp.\ $\mu_v$) is the adjoint morphism of $\gamma_u$ (resp.\ $\gamma_v$), are commutative.  

As in \ref{p2-qfc245}, if the functor $\varpi$ \eqref{p2-qfc246f} is quasi-fibering \eqref{p2-qfc9}, the bifunctoriality 
condition above can be replaced by the commutativity of the diagram
\begin{equation}\label{p2-qfc250c}
\xymatrix{
{\Phi_g^\rs(F_i)}\ar[r]^-(0.5){(\fd^\rs_{f,v})^{-1}}\ar[d]_{(\fd^\rs_{u,f'})^{-1}}&{v^\rs(\Phi_f^\rs(F_{i}))}\ar[r]^-(0.5){v^\rs(\tmu_{f})}&
{v^\rs(G_j)}\ar[d]^{\mu_v}\\
{\Phi_{f'}^\rs(u^\rp(F_{i}))}\ar[r]^-(0.5){\Phi_{f'}^\rs(\mu_{u})}&{\Phi_{f'}^\rs(F_{i'})}\ar[r]^-(0.5){\tmu_{f'}}&{G_{j'}.}}
\end{equation}
\end{rema}

\begin{prop}\label{p2-qfc253}
We keep the assumptions and notation of \ref{p2-qfc250}. 
Let $F\in \ob(\hcC)$, $G\in \ob(\hcD')$, $\uF=\iota(F)$ \eqref{p2-qfc20d}, $\uG=\iota_{\cD'/\cJ'}(G)$ \eqref{p2-qfc246n}, $F^a$ (resp.\ $G^a$) 
the sheaf associated with $F$ (resp.\ $G$). Then, an $\mI_{\varphi'}$-system of morphisms from $\uF$ to $\uG$ in the sense of \ref{p2-qfc29}
determines a morphism $\Phi'^\rs(F^a)\rightarrow G^a$ of $\tcC$. 
\end{prop}

Indeed, an $\mI_{\varphi'}$-system of morphisms from $\uF$ to $\uG$ is equivalent to a morphism $\Phi'^\rp(F)\rightarrow G$ of $\hcC$ by \ref{p2-qfc251},  
which induces a morphism $\Phi'^\rs(F^a)\rightarrow G^a$ of $\tcC$ by \ref{p2-cmt77}(ii).


\chapter{Twisting Higgs modules}\label{twisting}

\section{Twistable Higgs modules}\label{p1-thbn}

\subsection{}\label{p1-thbn0}
In this chapter, $V$ denotes a complete valuation ring of height $1$, with algebraically closed fraction field of characteristic $0$ 
and residue field of characteristic $p>0$. We choose a compatible system $(\beta_n)_{n>0}$ of $n$th roots of $p$ in $V$.
For any rational number $\varepsilon>0$, we set $p^\varepsilon=(\beta_n)^{\varepsilon n}$ where $n$ is an integer $>0$ 
such that $\varepsilon n$ is an integer. We equip $V$ with the $p$-adic topology and set $\cS=\Spf(V)$.

\subsection{}\label{p1-thbn1}
We consider in this section a $\mU$-topos  $X$ (\cite{sga4} IV 1.1), a $V$-algebra $\co_X$ of $X$ which is $p$-adic 
(i.e., $p$-adically complete and separated) and $V$-flat 
(i.e., $p$-torsion free), and an exact sequence of locally free $\co_X$-modules of finite type
\begin{equation}\label{p1-thbn1a}
0\longrightarrow \co_X\stackrel{c}{\longrightarrow} \cF\stackrel{\nu}{\longrightarrow} \cE \longrightarrow 0.
\end{equation}
By (\cite{illusie1} I 4.3.1.7), this sequence induces for every integer $n\geq 1$ an exact sequence \eqref{p1-NC7}
\begin{equation}\label{p1-thbn1b}
0\rightarrow \rS^{n-1}_{\co_X}(\cF)\rightarrow \rS^{n}_{\co_X}(\cF)\rightarrow \rS^n_{\co_X}(\cE)\rightarrow 0.
\end{equation}
The $\co_X$-modules $(\rS^{n}_{\co_X}(\cF))_{n\in \mN}$ therefore form a filtered direct system, whose direct limit
\begin{equation}\label{p1-thbn1c}
\cC=\underset{\underset{n\geq 0}{\longrightarrow}}\lim\ \rS^n_{\co_X}(\cF)
\end{equation}
is naturally endowed with an $\co_X$-algebra structure. We call it the {\em Higgs--Tate algebra} 
associated with the extension $\cF$ \eqref{p1-thbn1a}. There is a unique $\co_X$-derivation 
\begin{equation}\label{p1-thbn1d}
d_\cC\colon \cC \rightarrow \cE\otimes_{\co_X} \cC
\end{equation}
extending $\nu$ \eqref{p1-thbn1a}. It identifies canonically with the universal $\co_X$-derivation of $\cC$ \eqref{p1-imdpa100}.
We check immediately that it is also a Higgs $\co_X$-field on $\cC$ with coefficients in $\cE$ \eqref{p1-delta-con1}.

\begin{rema}\label{p1-thbn18}
Let $U$ be an object of $X$, $j_U\colon X_{/U}\rightarrow X$ the localization of $X$ at $U$. We equip $X_{/U}$ with the ring $j_U^*(\co_X)$.  
For any sheaf $F$ of $X$, we may also denote $j_U^*(F)$ by $F|U$.  
Since the functor $j_U^*$ commutes with direct limits, the formation of the algebra $\cC$ commutes with $j_U^*$. 
On the other hand, since $j_U^*$ admits a left adjoint, it commutes with inverse limits. 
Hence, the property for an $\co_X$-module of being $p$-adic (i.e. $p$-adically complete and Hausdorff) is local in $X$.
In particular, $\cE$ and $\cF$ are $p$-adic.
\end{rema}

\subsection{}\label{p1-thbn3}
For any rational number $r\geq 0$, we denote by $\cF^{(r)}$ the extension of $\co_X$-modules deduced from $\cF$ by pullback by the multiplication 
by $p^r$ on $\cE$; so we have an exact sequence of $\co_X$-modules
\begin{equation}\label{p1-thbn3a}
0\longrightarrow \co_X\stackrel{c^{(r)}}{\longrightarrow} \cF^{(r)}\stackrel{\nu^{(r)}}{\longrightarrow} \cE \longrightarrow 0. 
\end{equation} 
We denote by 
\begin{equation}\label{p1-thbn3b}
\cC^{(r)}=\underset{\underset{n\geq 0}{\longrightarrow}}\lim\ \rS^n_{\co_X}(\cF^{(r)})
\end{equation}
the associated Higgs--Tate algebra, see \eqref{p1-thbn1c}. We also call it the {\em Higgs--Tate algebra of thickness $r$} 
associated with the extension $\cF$ \eqref{p1-thbn1a}.
There exists a unique $\co_X$-derivation
\begin{equation}\label{p1-thbn3c}
d_{\cC^{(r)}}\colon \cC^{(r)} \rightarrow \cE\otimes_{\co_X} \cC^{(r)}
\end{equation}
extending $\nu^{(r)}$. It is the universal $\co_X$-derivation of $\cC^{(r)}$, and is also a Higgs $\co_X$-field. We set 
\begin{equation}\label{p1-thbn3d}
\delta_{\cC^{(r)}}=p^rd_{\cC^{(r)}}\colon \cC^{(r)}\rightarrow \cE\otimes_{\co_X}\cC^{(r)}. 
\end{equation}

\subsection{}\label{p1-thbn4}
For all rational numbers $r\geq r'\geq 0$, we have a canonical $\co_X$-linear morphism
\begin{equation}\label{p1-thbn4a}
\tta^{r,r'}\colon \cF^{(r)}\rightarrow \cF^{(r')}
\end{equation}
which lifts the multiplication by $p^{r-r'}$ on $\cE$ and which extends the identity on $\co_X$ \eqref{p1-thbn3a}.
It induces a homomorphism of $\co_X$-algebras
\begin{equation}\label{p1-thbn4b}
\alpha^{r,r'}\colon \cC^{(r)}\rightarrow \cC^{(r')}.
\end{equation}
We clearly have
\begin{equation}\label{p1-thbn4c}
(\id \otimes \alpha^{r,r'}) \circ \delta_{\cC^{(r)}}=\delta_{\cC^{(r')}} \circ \alpha^{r,r'}.
\end{equation}

For every rational numbers $r\geq r'\geq r''\geq 0$, we have
\begin{equation}\label{p1-thbn4d}
\tta^{r,r''}=\tta^{r',r''} \circ \tta^{r,r'} \ \ \ {\rm and}\ \ \ \alpha^{r,r''}=\alpha^{r',r''} \circ \alpha^{r,r'}.
\end{equation}

\subsection{}\label{p1-thbn6}
For any rational number $r\geq 0$, we denote by $\hcC^{(r)}$ the $p$-adic completion of $\cC^{(r)}$:
\begin{equation}\label{p1-thbn6aa}
\hcC^{(r)}=\underset{\underset{n\geq 0}{\longleftarrow}}\lim\ \cC^{(r)}/p^n \cC^{(r)}.
\end{equation}
Observe that $\cE\otimes_{\co_X}\hcC^{(r)}$ is the $p$-adic completion of $\cE\otimes_{\co_X}\cC^{(r)}$ \eqref{p1-thbn18}. We denote by 
\begin{eqnarray}
d_{\hcC^{(r)}} \colon \hcC^{(r)} \rightarrow \cE\otimes_{\co_X} \hcC^{(r)},\label{p1-thbn6a}\\
\delta_{\hcC^{(r)}}\colon \hcC^{(r)} \rightarrow \cE\otimes_{\co_X} \hcC^{(r)},\label{p1-thbn6b}
\end{eqnarray}
the extensions of $d_{\cC^{(r)}}$ and $\delta_{\cC^{(r)}}$ to the $p$-adic completions; so $d_{\hcC^{(r)}}$ is a Higgs $\co_X$-field 
on $\hcC^{(r)}$ with coefficients in $\cE$, and we have $\delta_{\hcC^{(r)}}=p^rd_{\hcC^{(r)}}$. 

For all rational numbers $r\geq r'\geq 0$, we denote by 
\begin{equation}\label{p1-thbn6c}
\halpha^{r,r'}\colon \hcC^{(r)}\rightarrow \hcC^{(r')}
\end{equation}
the extension of $\alpha^{r,r'}$  to the $p$-adic completions. By \eqref{p1-thbn4c}, we have 
\begin{equation}\label{p1-thbn6d}
(\id \otimes \halpha^{r,r'}) \circ \delta_{\hcC^{(r)}}=\delta_{\hcC^{(r')}} \circ \halpha^{r,r'}.
\end{equation}

We define the $\co_X$-algebra
\begin{equation}\label{p1-thbn6e}
\hcC^{(r+)}=\underset{\underset{t\in \mQ_{>r}}{\longrightarrow}}{\lim}\ \hcC^{(t)}. 
\end{equation}
By \eqref{p1-thbn6d}, the derivations $\delta_{\hcC^{(t)}}$ \eqref{p1-thbn6b} induce an $\co_X$-derivation
\begin{equation}\label{p1-thbn6f}
\delta_{\hcC^{(r+)}}\colon \hcC^{(r+)}\rightarrow \cE\otimes_{\co_X}\hcC^{(r+)},
\end{equation}
which is also a Higgs $\co_X$-field.
For simplicity, we set $\cC^\dagger=\hcC^{(0+)}$ and $\delta_{\cC^\dagger}=\delta_{\hcC^{(0+)}}$.

\subsection{}\label{p1-thbn16}
Let $\sigma\colon \cE\rightarrow \cF$, $\rho\colon \cF\rightarrow \co_X$ be associated splittings of the extension \eqref{p1-thbn1a}, 
i.e., $\sigma\circ \nu=\id_\cF-c\circ \rho$. 
For every rational number $r\geq 0$, they induce associated splittings 
$\sigma^{(r)}\colon \cE\rightarrow \cF^{(r)}$ and $\rho^{(r)}\colon \cF^{(r)}\rightarrow \co_X$ 
of the extension \eqref{p1-thbn3a}. The latter induce homomorphisms of $\co_X$-algebras 
\begin{eqnarray}
\varsigma^{(r)}\colon \rS_{\co_X}(\cE)\stackrel{\sim}{\rightarrow} \cC^{(r)},\label{p1-thbn16a}\\
\varrho^{(r)}\colon \cC^{(r)}\rightarrow \co_X,\label{p1-thbn16b}
\end{eqnarray}
the first of which is an isomorphism \eqref{p1-imdpa100}.
For all rational numbers $r\geq r'\geq 0$, we immediately check that the diagram 
\begin{equation}\label{p1-thbn16c}
\xymatrix{
{\rS_{\co_X}(\cE)}\ar[r]^-(0.5){\varsigma^{(r)}}\ar[d]_{\mu^{(r-r')}}&{\cC^{(r)}}\ar[d]_-(0.5){\alpha^{r,r'}}\ar[rd]^{\varrho^{(r)}}&\\
{\rS_{\co_X}(\cE)}\ar[r]^-(0.5){\varsigma^{(r')}}&{\cC^{(r')}}\ar[r]^{\varrho^{(r')}}&{\co_X,}}
\end{equation}
where $\mu^{(r-r')}$ is the homomorphism of $\co_X$-algebras induced by the multiplication by $p^{r-r'}$ on $\cE$, is commutative. 

\begin{rema}\label{p1-thbn220}
Set $\cG=\rS_{\co_X}(\cE)$ and denote by $\hcG$ its $p$-adic completion.  
Then, for every object $U$ of $X$, $\hcG(U)$ is canonically isomorphic to the 
$\co_X(U)$-algebra of series $(s_i)_{i\in \mN}\in \prod_{i\in \mN}\rS^i_{\co_X}(\cE)(U)$
such that for every integer $m\geq 0$, there exists $n\in \mN$ such that for every $i\geq n$, the class of $s_i$ in 
$\rS^i_{\co_X}(\cE/p^m\cE)(U)$ vanishes. Indeed, for every $i\geq 0$, the $\co_X$-module $\rS^i_{\co_X}(\cE)$ is locally free of finite type 
and hence $p$-adic \eqref{p1-thbn18}. 
\end{rema}

\begin{lem}\label{p1-thbn22}
Let $r,r'$ be rational numbers such that $r\geq r'\geq 0$. Then, 
\begin{itemize}
\item[{\rm (i)}] The $V$-algebra $\hcC^{(r)}$ is flat. 
\item[{\rm (ii)}] The homomorphism $\halpha^{r,r'}$ \eqref{p1-thbn6c} is injective. 
\item[{\rm (iii)}] We have $\ker(d_{\hcC^{(r)}})=\co_X$. 
\item[{\rm (iv)}] We have $\ker(\delta_{\hcC^{(r+)}})=\ker(\delta_{\hcC^{(r)}})=\co_X$. 
\end{itemize}
\end{lem}

(i), (ii)  \& (iii)  Indeed, the propositions being local on $X$, 
we may assume that the $\co_X$-module $\cE$ is free of finite type, so the extension \eqref{p1-thbn1a} is split. 
Then by \ref{p1-thbn16}, we are reduced to proving the analogous propositions for the algebra $\cG=\rS_{\co_X}(\cE)$,
which follow from \ref{p1-thbn220}. 

(iv) It follows from (i) and (iii). 

\subsection{}\label{p1-thbn35}
For every rational number $r\geq 0$, we denote by $\mK^\bullet(\hcC^{(r)})$
the Dolbeault complex of $(\hcC^{(r)},\delta_{\hcC^{(r)}})$ \eqref{p1-thbn6b} and by $\tmK^\bullet(\hcC^{( r)})$
the augmented Dolbeault complex 
\begin{equation}\label{p1-thbn35a}
\co_X\rightarrow \mK^0(\hcC^{(r)})\rightarrow \mK^1(\hcC^{(r)})\rightarrow \dots
\rightarrow \mK^n(\hcC^{(r)})\rightarrow \dots,
\end{equation}
where $\co_X$ is placed in degree $-1$ and the differential $\co_X\rightarrow\hcC^{(r)}$ is the canonical homomorphism.
We denote by $\mK^\bullet(\hcC^{(r+)})$ the Dolbeault complex of $(\hcC^{(r+)},\delta_{\hcC^{(r+)}})$ \eqref{p1-thbn6f}. For 
$r=0$, we denote it by $\mK^\bullet(\cC^\dagger)$. 

For all rational numbers $r\geq r'\geq 0$, by \eqref{p1-thbn6d}, $\halpha^{r,r'}$ induces a morphism of complexes
\begin{equation}\label{p1-thbn35b}
\upiota^{r,r'}\colon \tmK^\bullet(\hcC^{(r)})\rightarrow \tmK^\bullet(\hcC^{(r')}).
\end{equation}

\begin{prop}\label{p1-thbn36}
Suppose that the extension \eqref{p1-thbn1a} is split and that the $\co_X$-module $\cE$ is free. 
Then, for all rational numbers such that $r>r'>0$, there exists a rational number $\alpha\geq 0$ such that
\begin{equation}\label{p1-thbn36a}
p^\alpha\upiota^{r,r'}\colon \tmK^\bullet(\hcC^{(r)})\rightarrow \tmK^\bullet(\hcC^{(r')}),
\end{equation}
where $\upiota^{r,r'}$ is the morphism \eqref{p1-thbn35b}, is homotopic to $0$ by an $\co_X$-linear homotopy.
\end{prop}

Indeed, we are reduced by \ref{p1-thbn16} to proving the analogous proposition for the algebra $\rS_{\co_X}(\cE)$,
in which case the homotopy is explicitly given in the proof of (\cite{agt} II.11.2), taking into account \ref{p1-thbn220} and the fact 
that the $\co_X$-modules $\rS^i_{\co_X}(\cE)$ are $V$-flat for all $i\geq 0$.  

\begin{cor}\label{p1-thbn37}
For every rational number $r\geq 0$ and every $\co_X[\frac 1 p]$-module $M$, 
the complex $M\otimes_{\co_X}\mK^\bullet(\hcC^{(r+)})$ is a resolution of $M$, 
where the tensor product is defined term by term (non derived).
\end{cor}

\subsection{}\label{p1-thbn9}
For any Higgs $\co_X$-module $(N,\theta)$ with coefficients in $\cE$, we denote by $\uuptau(N,\theta)$ the $\co_X$-module 
\begin{equation}\label{p1-thbn9a}
\uuptau(N,\theta)=(N\otimes_{\co_X}\cC^\dagger)^{\theta_\tot=0},
\end{equation}
where $\theta_\tot=\theta\otimes \id+\id\otimes \delta_{\cC^\dagger}$ is the total Higgs $\co_X$-field on $ N\otimes_{\co_X}\cC^\dagger$. 
It is clear that $\uuptau(N,\theta)$ is functorial in $(N,\theta)$.  
The Higgs field $\theta\otimes \id$ on $N\otimes_{\co_X}\cC^\dagger$ induces a Higgs $\co_X$-field $\theta_\uptau$ on $\uuptau(N,\theta)$ with coefficients in $\cE$. 
Indeed, the question being local on $X$, we may assume that the $\co_X$-module $\cE$ is free of finite type. 
Then $\theta$ (resp.\ $\delta_{\cC^\dagger}$) is determined by its components on a basis of $\cE$, 
which are $\co_X$-linear endomorphisms of $N$ (resp.\ $\cC^\dagger$) 
that commute to each other. Denoting by $\bHM(\co_X,\cE)$ 
the category of Higgs $\co_X$-modules with coefficients in $\cE$,  
we thus define a functor, \begin{equation}\label{p1-thbn9b}
\uptau\colon 
\begin{array}[t]{clcr}
\bHM(\co_X,\cE) &\rightarrow& \bHM(\co_X,\cE),\\
(N,\theta)&\mapsto&(\uuptau(N,\theta), \theta_\uptau),
\end{array}
\end{equation}
which we call the {\em twisting functor relative to the extension $\cF$}.

\begin{rema}\label{p1-thbn34}
Set $\cE^\vee=\cHom_{\co_X}(\cE,\co_X)$. 
For every Higgs $\co_X$-module $(N,\theta)$ with coefficients in $\cE$, the canonical morphism 
\begin{equation}\label{p1-thbn34a}
\uuptau(N,\theta)\rightarrow \cC^\dagger\otimes_{\co_X}N
\end{equation}
is $\rS_{\co_X}(\cE^\vee)$-linear, where $\rS_{\co_X}(\cE^\vee)$ acts on $N$ (resp.\  $\uuptau(N,\theta)$) via the 
Higgs field $\theta$ (resp.\ $\theta_\uptau$) \eqref{p1-delta-con1j} 
and does not act on $\cC^\dagger$. Since the actions of  
$\cC^\dagger$ and $\rS_{\co_X}(\cE^\vee)$ on $\cC^\dagger \otimes_{\co_X}N$ 
(resp.\ $\cC^\dagger \otimes_{\co_X}\uuptau(N,\theta)$) commute, we deduce that 
the canonical $\cC^\dagger$-linear morphism
\begin{equation}\label{p1-thbn34b}
\cC^\dagger \otimes_{\co_X}\uuptau(N,\theta)\rightarrow \cC^\dagger \otimes_{\co_X}N
\end{equation}
is $\rS_{\co_X}(\cE^\vee)$-linear. 
\end{rema}

\begin{defi}\label{p1-thbn30}
A Higgs $\co_X[\frac 1 p]$-module $(N,\theta)$ with coefficients in $\cE$ is said to be 
{\em weakly twistable by the extension \eqref{p1-thbn1a}} (or simply {\em weakly twistable} if there is no risk of ambiguity)
if the canonical $\cC^\dagger$-linear morphism
\begin{equation}\label{p1-thbn30a}
\cC^\dagger \otimes_{\co_X}\uuptau(N,\theta)\rightarrow \cC^\dagger \otimes_{\co_X}N
\end{equation}
is an isomorphism \eqref{p1-thbn9a}.  
\end{defi}

We keep the adjective {\em twistable} for a stronger notion \eqref{p1-thbn14}, see \ref{p1-thbn21}.

\begin{prop}\label{p1-thbn31}
Suppose that the extension \eqref{p1-thbn1a} splits. 
Let $(N,\theta)$ be a weakly twistable Higgs $\co_X[\frac 1 p]$-module with coefficients in $\cE$. 
Then, the choice of a splitting of the extension \eqref{p1-thbn1a} determines a functorial isomorphism of Higgs $\co_X$-modules
\begin{equation}\label{p1-thbn31a}
\uptau(N,\theta)\stackrel{\sim}{\rightarrow}(N,\theta).
\end{equation}
\end{prop}

Indeed, the choice of a splitting of the extension \eqref{p1-thbn1a} determines a homomorphism $\varrho^\dagger \colon \cC^\dagger \rightarrow \co_X$, by \ref{p1-thbn16}. 
Consider the commutative diagram 
\begin{equation}
\xymatrix{
{\uuptau(N,\theta)}\ar[r]\ar[rd]_\id&{\cC^\dagger \otimes_{\co_X}\uuptau(N,\theta)}\ar[r]\ar[d]^{\varrho^\dagger\otimes \id}&{\cC^\dagger \otimes_{\co_X}N}\ar[d]^{\varrho^\dagger\otimes\id}\\
&{\uuptau(N,\theta)}\ar[r]^u&N,}
\end{equation}
where the upper horizontal arrows are the canonical morphisms, and $u$ is defined by base change by $\varrho^\dagger$. 
It shows that $u$ is an isomorphism, and is compatible with the Higgs fields $\theta_\uptau$ and $\theta$. 
We then take for \eqref{p1-thbn31a} the isomorphism $u$. 

\begin{cor}\label{p1-thbn32}
For every weakly twistable Higgs $\co_X[\frac 1 p]$-module $(N,\theta)$ with coefficients in $\cE$, the $\co_X[\frac 1 p]$-module 
$\uuptau(N,\theta)$ is locally isomorphic to $N$. 
\end{cor}

\begin{cor}\label{p1-thbn33}
For every weakly twistable Higgs $\co_X[\frac 1 p]$-bundle $(N,\theta)$ with coefficients in $\cE$ \eqref{p1-delta-con6}, 
$(N,\theta)$ and $\uptau(N,\theta)$ have the same characteristic invariants \eqref{p1-delta-con1e}.
\end{cor}

Indeed, the question being local, we may assume that the extension \eqref{p1-thbn1a} splits, in which case the proposition 
follows immediately from \ref{p1-thbn31}.

\subsection{}\label{p1-thbn40}
For any rational number $r\geq 0$, we set 
\begin{equation}\label{p1-thbn40a}
\delta^\vee_{\hcC^{(r)}}=-\delta_{\hcC^{(r)}}, \ \ \ \delta^\vee_{\hcC^{(r+)}}=-\delta_{\hcC^{(r+)}}\ \ \ {\rm and} \ \ \ \delta^\vee_{\cC^\dagger}=-\delta_{\cC^\dagger},
\end{equation}
where $\delta_{\hcC^{(r)}}$ is defined in \eqref{p1-thbn6b} and $\delta_{\hcC^{(r+)}}$ in \eqref{p1-thbn6f}. 
For any Higgs $\co_X$-module $(N,\theta)$ with coefficients in $\cE$, we denote by $\uuptau^\vee(N,\theta)$ the $\co_X$-module 
\begin{equation}\label{p1-thbn40b}
\uuptau^\vee(N,\theta)=(N\otimes_{\co_X}\cC^\dagger)^{\theta^\vee_\tot=0},
\end{equation}
where $\theta^\vee_\tot=\theta\otimes \id+\id\otimes \delta^\vee_{\cC^\dagger}$ is the total Higgs $\co_X$-field on $N\otimes_{\co_X}\cC^\dagger$. 
It is clear that $\uuptau^\vee(N,\theta)$ is functorial in $(N,\theta)$.  
The Higgs field $\theta\otimes \id$ on $N\otimes_{\co_X}\cC^\dagger$ induces a Higgs $\co_X$-field $\theta_{\uptau^\vee}$ on $\uuptau^\vee(N,\theta)$ with coefficients in $\cE$,
by the same proof as in \ref{p1-thbn9}. We thus define a functor 
\begin{equation}\label{p1-thbn40c}
\uptau^\vee\colon 
\begin{array}[t]{clcr}
\bHM(\co_X,\cE) &\rightarrow& \bHM(\co_X,\cE),\\
(N,\theta)&\mapsto&(\uuptau^\vee(N,\theta), \theta_{\uptau^\vee}).
\end{array}
\end{equation}

\begin{prop}\label{p1-thbn41}
Let $(N,\theta)$ be a Higgs $\co_X[\frac 1 p]$-module with coefficients in $\cE$, $N^\vee$ an $\co_X$-module, 
\begin{equation}\label{p1-thbn41a}
\psi\colon N^\vee\otimes_{\co_X}\cC^\dagger\stackrel{\sim}{\rightarrow}N\otimes_{\co_X}\cC^\dagger
\end{equation}
an isomorphism of $\cC^\dagger$-modules with $\delta_{\cC^\dagger}$-connection, 
where the $\delta_{\cC^\dagger}$-connections are defined as in \ref{p1-delta-con4}, $N$ (resp.\ $N^\vee$) being endowed with the Higgs field $\theta$ (resp.\ $0$).  Then, 
\begin{itemize}
\item[{\rm (i)}] The isomorphism $\psi$ induces an $\co_X$-linear isomorphism $N^\vee\stackrel{\sim}{\rightarrow} \uuptau(N,\theta)$, 
where $\uuptau$ is the functor \eqref{p1-thbn9a}. In particular, $(N,\theta)$ is weakly twistable. 
We deduce a canonical Higgs $\co_X$-field $\theta^\vee$ on  $N^\vee$ with coefficients in $\cE$, so that we have an isomorphism of Higgs modules
\begin{equation}\label{p1-thbn41b}
(N^\vee,\theta^\vee)\stackrel{\sim}{\rightarrow} \uptau(N,\theta),
\end{equation}
where $\uptau$ is the functor \eqref{p1-thbn9b}. 
\item[{\rm (ii)}] The morphism $\psi$ is an isomorphism of $\cC^\dagger$-modules with $\delta^\vee_{\cC^\dagger}$-connection \eqref{p1-thbn40a}, 
where the $\delta^\vee_{\cC^\dagger}$-connections are defined as in \ref{p1-delta-con4}, $N$ (resp.\ $N^\vee$) being endowed with the Higgs field $0$ (resp.\ $\theta^\vee$). 
\item[{\rm (iii)}] The isomorphism $\psi$ induces an isomorphism of Higgs $\co_X$-modules with coefficients in $\cE$, 
\begin{equation}\label{p1-thbn41c}
(N,\theta)\stackrel{\sim}{\rightarrow} \uptau^\vee(N^\vee,\theta^\vee),
\end{equation}
where $\uptau^\vee$ is the functor \eqref{p1-thbn40c}. 
\end{itemize}
\end{prop}

(i) It follows immediately from \ref{p1-thbn16} that the $\co_X$-modules $N$ and $N^\vee$ are locally isomorphic (see the proof of \ref{p1-thbn31}). 
In particular, $N^\vee$ is an $\co_X[\frac 1 p]$-module.
The proposition follows then from \ref{p1-thbn37}. 

(ii) It follows immediately from the fact that the Higgs field $\theta^\vee$ on $N^\vee$ 
is induced by $\id\otimes\delta^\vee_{\cC^\dagger}$ on $N\otimes_{\co_\fX}\cC^\dagger$, which is a $\delta^\vee_{\cC^\dagger}$-connection.

(iii) It  follows from (ii) and \ref{p1-thbn37} that $\psi^{-1}$ induces an isomorphism of the $\co_X$-modules underlying \eqref{p1-thbn41c}.  
 Let $x$ be a local section of $N$. Since $\theta^\vee_\tot(\psi^{-1}(x\otimes 1))=0$, we have 
\begin{eqnarray}
(\theta^\vee\otimes \id)(\psi^{-1}(x\otimes 1))&=&(\id\otimes\delta_{\cC^\dagger})(\psi^{-1}(x\otimes 1))\nonumber \\
&=&(\theta\otimes \id+\id\otimes\delta_{\cC^\dagger})(x\otimes 1)\nonumber \\
&=&\theta(x)\otimes 1,
\end{eqnarray}
the second equality being a consequence of the assumptions. Therefore, \eqref{p1-thbn41c} is an isomorphism of Higgs $\co_X$-modules.

\subsection{}\label{p1-thbn10}
Let $r$ be a rational number $\geq 0$. 
We denote by $\bMIC_r(\hcC^{(r)}/\co_X)$ the category of $\hcC^{(r)}$-modules with integrable 
$\delta_{\hcC^{(r)}}$-connection \eqref{p1-delta-con2}.
By \ref{p1-delta-con4}, we have a functor 
\begin{equation}\label{p1-thbn10a}
\uppi^{(r)}\colon 
\begin{array}[t]{clcr}
\bHM(\co_X,\cE)&\rightarrow&\bMIC_r(\hcC^{(r)}/\co_X)\\
(N,\theta)&\mapsto&(N\otimes_{\co_X}\hcC^{(r)},\theta\otimes \id+\id\otimes\delta_{\hcC^{(r)}}). 
\end{array}
\end{equation}
Restricting this functor to Higgs $\co_X$-modules with vanishing Higgs field, we obtain the functor
\begin{equation}\label{p1-thbn10b}
\fS^{(r)}\colon 
\begin{array}[t]{clcr}
\bMod(\co_X)&\rightarrow&\bMIC_r(\hcC^{(r)}/\co_X)\\
M&\mapsto&(M\otimes_{\co_X}\hcC^{(r)},\id\otimes\delta_{\hcC^{(r)}}).
\end{array}
\end{equation}

\subsection{}\label{p1-thbn11}
Let $r\geq r'\geq 0$ be rational numbers, $(M,\nabla)$ a $\hcC^{(r)}$-module with integrable $\delta_{\hcC^{(r)}}$-connection. It follows from \ref{p1-delta-con5} 
that there exists a unique integrable $\delta_{\hcC^{(r')}}$-connection
\begin{equation}\label{p1-thbn11a}
\nabla'\colon M\otimes_{\hcC^{(r)}}\hcC^{(r')}\rightarrow \cE\otimes_{\co_X}M\otimes_{\hcC^{(r)}}\hcC^{(r')}
\end{equation}
such that for any local sections $m$ of $M$ and $\xi'$ of $\hcC^{(r')}$, we have
\begin{equation}\label{p1-thbn11b}
\nabla'(m\otimes_{\hcC^{(r)}}\xi')=\nabla(m)\otimes_{\hcC^{(r)}} \xi' + m\otimes_{\hcC^{(r)}}\delta_{\hcC^{(r')}}(\xi'). 
\end{equation} 
We thus define a functor 
\begin{equation}\label{p1-thbn11c}
\varepsilon^{r,r'}\colon 
\begin{array}[t]{clcr}
\bMIC_r(\hcC^{(r)}/\co_X)&\rightarrow& \bMIC_{r'}(\hcC^{(r')}/\co_X)\\
(M,\nabla)&\mapsto&(M\otimes_{\hcC^{(r)}}\hcC^{(r')},\nabla').
\end{array}
\end{equation} 
We have a canonical isomorphism of functors from $\bMod(\co_X)$ to $\bMIC_{r'}(\hcC^{(r')}/\co_X)$ 
\begin{equation}\label{p1-thbn11d}
\varepsilon^{r,r'}\circ \fS^{(r)}\stackrel{\sim}{\longrightarrow} \fS^{(r')}.
\end{equation}
We have a canonical isomorphism of functors from $\bHM(\co_X,\cE)$ to $\bMIC_{r'}(\hcC^{(r')}/\co_X)$ 
\begin{equation}\label{p1-thbn11e}
\varepsilon^{r,r'}\circ \uppi^{(r)}\stackrel{\sim}{\longrightarrow} \uppi^{(r')}.
\end{equation}

\begin{defi}\label{p1-thbn12}
Let $M$ be an $\co_X$-module, $(N,\theta)$ a Higgs $\co_X$-module with coefficients in $\cE$. 
\begin{itemize}
\item[(i)] Let $r$ be a rational number $>0$.  We say that $M$ and  $(N,\theta)$ are {\em $r$-associated} 
if there exists an isomorphism of $\bMIC_r(\hcC^{(r)}/\co_X)$ 
\begin{equation}\label{p1-thbn12a}
\alpha\colon \fS^{(r)}(M)\stackrel{\sim}{\rightarrow}\uppi^{(r)}(N,\theta).
\end{equation}
We say then that the triple $(M,(N,\theta),\alpha)$ is {\em $r$-admissible}. 
\item[(ii)] We say that $M$ and $(N,\theta)$ are {\em associated} if there exists a rational number $r>0$ such that they are $r$-associated.
\end{itemize}
\end{defi}

Observe that for every rational numbers $r\geq r'>0$,
if $M$ and $(N,\theta)$ are $r$-associated, they are $r'$-associated, in view of \eqref{p1-thbn11d} and \eqref{p1-thbn11e}.

\begin{defi}\label{p1-thbn14}
We say that a Higgs $\co_X[\frac 1 p]$-module with coefficients in $\cE$ is 
{\em twistable by the extension \eqref{p1-thbn1a}} (or simply {\em twistable} if there is no risk of ambiguity) 
if it is associated with an $\co_X$-module \eqref{p1-thbn12}. 
\end{defi}

\begin{rema}\label{p1-thbn26}
The notions introduced in \ref{p1-thbn12}(ii) and in \ref{p1-thbn14} are clearly stable by restriction to any object of $X$ \eqref{p1-thbn18}. 
We prove in \ref{p1-thbn24} that the second one is in fact local under extra conditions. 
\end{rema}

\begin{lem}\label{p1-thbn15}
Let $M$ be an $\co_X$-module, $(N,\theta)$ a Higgs $\co_X[\frac 1 p]$-module with coefficients in $\cE$
such that $M$ and $(N,\theta)$ are associated relatively to the extension \eqref{p1-thbn1a}. 
Then, the $\co_X$-modules $M$ and $N$ are locally isomorphic. In particular, $M$ is an $\co_X[\frac 1 p]$-module.
\end{lem}

It follows immediately from \ref{p1-thbn16} (see the proof of \ref{p1-thbn31}).

\begin{prop}\label{p1-thbn21}
Let $(N,\theta)$ be a Higgs $\co_X[\frac 1 p]$-module with coefficients in $\cE$. 
If $(N,\theta)$ is twistable by the extension \eqref{p1-thbn1a} in the sense of \ref{p1-thbn14}, then it is weakly twistable by the extension \eqref{p1-thbn1a} in the sense of \ref{p1-thbn30}. 
The converse is true if the $\co_X[\frac 1 p]$-module $N$ is flat of finite type and 
the topos $X$ is coherent {\rm (\cite{sga4} VI 2.3)}.
\end{prop}

We assume first that the Higgs module $(N,\theta)$ is twistable. Then, there exist an $\co_X$-module $M$, 
a rational number $r>0$ and an isomorphism of $\bMIC_{r}(\hcC^{(r)}/\co_X)$
\begin{equation}\label{p1-thbn21a}
\fS^{(r)}(M)\stackrel{\sim}{\rightarrow}\uppi^{(r)}(N,\theta). 
\end{equation}
By \ref{p1-delta-con5}, the latter induces by extension of scalars by the canonical injection $\hcC^{(r)}\rightarrow \cC^\dagger$,  
an isomorphism of $\cC^\dagger$-modules with $\delta_{\cC^\dagger}$-connection \eqref{p1-thbn6f},
\begin{equation}\label{p1-thbn21b}
\cC^\dagger \otimes_{\co_X}M\stackrel{\sim}{\rightarrow}\cC^\dagger \otimes_{\co_X} N, 
\end{equation}
where the $\delta_{\cC^\dagger}$-connections are defined as in \ref{p1-delta-con4}, $N$ (resp.\ $M$) being endowed with the Higgs field $\theta$ (resp.\ $0$). 
By \ref{p1-thbn41}(i), \eqref{p1-thbn21b} induces an isomorphism 
\begin{equation}\label{p1-thbn21c}
M \stackrel{\sim}{\rightarrow} \uuptau(N,\theta).
\end{equation}
Moreover, the isomorphism \eqref{p1-thbn21b} identifies with the canonical morphism \eqref{p1-thbn30a}. 
Therefore, $(N,\theta)$ is weakly twistable.

Conversely, we assume that the $\co_X[\frac 1 p]$-module $N$ is flat of finite type, the topos $X$ is coherent, 
and the Higgs bundle $(N,\theta)$ is weakly twistable,
i.e., the canonical $\cC^\dagger$-linear morphism 
\begin{equation}\label{p1-thbn21g}
\cC^\dagger \otimes_{\co_X}\uuptau(N,\theta)\rightarrow \cC^\dagger \otimes_{\co_X}N
\end{equation}
is an isomorphism. By \ref{p1-thbn32}, the $\co_X[\frac 1 p]$-module $\uuptau(N,\theta)$ is flat of finite type. 
By (\cite{sga4} VI 1.23 and 2.2), in view of \ref{p1-thbn22}(ii), there exists a rational number $r>0$ such that 
the image of $\uuptau(\cN,\theta)$ by \eqref{p1-thbn21g} is contained in $\hcC^{(r)} \otimes_{\co_X}N$.  We deduce a 
$\hcC^{(r)}$-linear morphism 
\begin{equation}\label{p1-thbn21f}
\hcC^{(r)}\otimes_{\co_X}\uuptau(N,\theta)\rightarrow \hcC^{(r)}\otimes_{\co_X}N.
\end{equation}
It is injective since the morphism \eqref{p1-thbn21g} is injective and the $\co_X[\frac 1 p]$-module $\uuptau(N,\theta)$ is flat. 
After replacing $r$ by a smaller rational number $>0$, we may assume that it is surjective and hence bijective. 
Hence $(N,\theta)$ is twistable. 

\begin{cor}\label{p1-thbn23}
Any Higgs $\co_X[\frac 1 p]$-module $(N,\theta)$ with coefficients in $\cE$ which is twistable 
is associated with the $\co_X$-module $\uuptau(N,\theta)$ \eqref{p1-thbn9a}.
\end{cor}

\begin{cor}\label{p1-thbn24}
Suppose that the topos $X$ is coherent. Let $N$ be a flat $\co_X[\frac 1 p]$-module of finite type, $\theta$ 
a Higgs $\co_X[\frac 1 p]$-field with coefficients in $\cE$. 
Then, $(N,\theta)$ is twistable by the extension \eqref{p1-thbn1a} if and only if 
there exists a covering family $(U_i)_{i\in I}$ of the final object of $X$ such that for every $i\in I$, the Higgs module
$(N|U_i,\theta|U_i)$ of the ringed topos $(X_{/U_i},\co_X|U_i)$ \eqref{p1-thbn18} 
is twistable by the restriction of the extension \eqref{p1-thbn1a} to $U_i$.
\end{cor}

\begin{rema}\label{p1-thbn25}
Let $(N,\theta)$ be a Higgs $\co_X[\frac 1 p]$-module with coefficients in $\cE$, $N^\vee$ an $\co_X$-module, $r$ a rational number $>0$,
\begin{equation}\label{p1-thbn25a}
N^\vee\otimes_{\co_X}\hcC^{(r)}\stackrel{\sim}{\rightarrow}N\otimes_{\co_X}\hcC^{(r)}
\end{equation} 
an isomorphism of $\hcC^{(r)}$-modules with $\delta_{\hcC^{(r)}}$-connection, 
where the $\delta_{\hcC^{(r)}}$-connections are defined as in \ref{p1-delta-con4}, $N$ (resp.\ $N^\vee$) being endowed with the Higgs field $\theta$ (resp.\ $0$).  
The following properties are slight variations on those proved in \ref{p1-thbn41}:
\begin{itemize}
\item[{\rm (i)}] The isomorphism \eqref{p1-thbn25a} induces an injective $\co_X$-linear morphism
\begin{equation}\label{p1-thbn25b}
N^\vee \rightarrow (N\otimes_{\co_X}\hcC^{(r)})^{\theta^{(r)}_\tot=0},
\end{equation}
where $\theta^{(r)}_\tot=\theta\otimes \id+\id\otimes\delta_{\hcC^{(r)}}$ is the total Higgs $\co_X$-field on $N\otimes_{\co_X}\hcC^{(r)}$. 
By \ref{p1-thbn22}(iii), it is an isomorphism if $N^\vee$ is $\co_X$-flat. 
\item[{\rm (ii)}]  The Higgs field $\theta\otimes \id$ on $N\otimes_{\co_X}\hcC^{(r)}$ induces a Higgs $\co_X$-field $\theta^\vee$
on $N^\vee$ with coefficients in $\cE$.
\item[{\rm (iii)}] The morphism \eqref{p1-thbn25a} is an isomorphism of $\hcC^{(r)}$-modules with $\delta^\vee_{\hcC^{(r)}}$-connection 
\emergencystretch=1em 
\eqref{p1-thbn40a}, where the $\delta^\vee_{\hcC^{(r)}}$-connections are defined as in \ref{p1-delta-con4}, $N$ 
(resp.\ $N^\vee$) being endowed with the Higgs field $0$ (resp.\ $\theta^\vee$). 
\item[{\rm (iv)}] By \ref{p1-delta-con5}, \eqref{p1-thbn25a}  induces by extension of scalars by the canonical injection $\hcC^{(r)}\rightarrow \cC^\dagger$,  
an isomorphism of $\cC^\dagger$-modules with $\delta_{\cC^\dagger}$-connection \eqref{p1-thbn6f},
\begin{equation}
N^\vee\otimes_{\co_X}\cC^\dagger\stackrel{\sim}{\rightarrow}N\otimes_{\co_X}\cC^\dagger,
\end{equation}
where the $\delta_{\cC^\dagger}$-connections are defined as in \ref{p1-delta-con4}, $N$ (resp.\ $N^\vee$) being endowed with the Higgs field $\theta$ (resp.\ $0$). 
By \ref{p1-thbn41}, the latter induces isomorphisms of Higgs $\co_X$-modules with coefficients in $\cE$, 
\begin{eqnarray}
(N,\theta)&\stackrel{\sim}{\rightarrow}& \uptau^\vee(N^\vee,\theta^\vee),\label{p1-thb28c}\\
(N^\vee,\theta^\vee)&\stackrel{\sim}{\rightarrow}& \uptau(N,\theta).\label{p1-thb28d}
\end{eqnarray}
\end{itemize}
\end{rema}

\subsection{}\label{p1-thbn42}
We set $\cE^\vee=\cHom_{\co_X}(\cE,\co_X)$ and $\rS=\rS_{\co_X}(\cE^\vee)$ \eqref{p1-NC7}. 
Let $B$ be an $\co_X[\frac 1 p]$-algebra quotient of $\rS[\frac 1 p]$.
By \ref{p1-delta-con1j}, the composed homomorphism 
\begin{equation}\label{p1-thbn42a}
\upmu_B\colon \rS[\frac 1 p]\rightarrow B\rightarrow \cEnd_{\co_X}(B),
\end{equation}
where the first homomorphism is the canonical surjection and the second one is defined by sending 
a local section $P$ of $B$ to the endomorphism $(Q\mapsto PQ)$,
defines a Higgs $\co_X$-field
\begin{equation}\label{p1-thbn42b}
\theta_B\colon B\rightarrow B \otimes_{\co_X} \cE,
\end{equation}
that we qualify by {\em canonical}. 
Equivalently, $\theta_B$ is the composed morphism 
\begin{equation}\label{p1-thbn42c}
B\rightarrow B\otimes_{\co_X}\cE^\vee\otimes_{\co_X}\cE\rightarrow B\otimes_{\co_X}\rS\otimes_{\co_X}\cE\rightarrow B\otimes_{\co_X}\cE,
\end{equation}
where the first morphism is induced by the section $\id_\cE$ of $\cE^\vee\otimes_{\co_X}\cE$, canonically identified with $\cEnd_{\co_X}(\cE)$ 
\eqref{p1-delta-con1ef}, the second morphism by the embedding of $\cE^\vee$ in $\rS$ and the third morphism by the multiplication in $B$. 
If $U$ is an object of $X$, $e_1,\dots,e_d$ is an $\co_X|U$-basis of $\cE|U$ 
and $e^\vee_1,\dots,e^\vee_d$ the dual basis of $\cE^\vee|U$, then 
\begin{equation}\label{p1-thbn42d}
\theta_B|U=\sum_{i=1}^d\theta_{B,i}\otimes e_i,
\end{equation}
where $\theta_{B,i}$ is the endomorphism of $B|U$ defined by the multiplication by the class of $e^\vee_i$. 

We set 
\begin{equation}\label{p1-thbn42e}
\cL_B=\uptau(B,\theta_B),
\end{equation} 
where $\uptau$ is the functor \eqref{p1-thbn9b}. We will also set (abusively) $\cL_B=\uuptau(B,\theta_B)$, which is the $\co_X$-module
underlying $\uptau(B,\theta_B)$. 
We denote by $\theta^\vee_B\colon \cL_B\rightarrow \cL_B\otimes_{\co_X}\cE$ the Higgs $\co_X$-field associated with $\cL_B$
and by $\upmu^\vee_B\colon \rS[\frac 1 p]\rightarrow \cEnd_{\co_X}(\cL_B)$ the associated homomorphism of $\co_X$-algebras. 
Then, $\upmu^\vee_B$ factors through $B$ \eqref{p1-thbn34}, so that $\cL_B$ is naturally a $B$-module.

If the Higgs $\co_X[\frac 1 p]$-module $(B,\theta_B)$ is weakly twistable by the extension \eqref{p1-thbn1a}, the $B$-module $\cL_B$ is invertible by \ref{p1-thbn31}.

\subsection{}\label{p1-thbn43}
Let $(N,\theta)$ be a Higgs $\co_X[\frac 1 p]$-module with coefficients in $\cE$, $\upmu\colon \rS[\frac 1 p]\rightarrow \cEnd_{\co_X}(N)$ 
the homomorphism of $\co_X$-algebras defined by $\theta$ \eqref{p1-delta-con1j}. We assume that $\upmu$ factors 
through a quotient $\co_X[\frac 1 p]$-algebra $B$ of $\rS[\frac 1 p]$, and let $\theta_B\colon B\rightarrow B\otimes_{\co_X}\cE$ 
be the canonical Higgs $\co_X$-field \eqref{p1-thbn42}. Since $\theta_B$ is $B$-linear by \ref{p1-delta-con7}(i), we may consider the morphism 
\begin{equation}\label{p1-thbn43a}
\id_N\otimes_B\theta_B\colon N\otimes_BB\rightarrow N\otimes_BB\otimes_{\co_X}\cE.
\end{equation}
It is the canonical Higgs $\co_X$-field on the $B$-module $N\otimes_BB$ by \ref{p1-delta-con7}(ii). 
It fits into a commutative diagram
\begin{equation}\label{p1-thbn43b}
\xymatrix{
{N\otimes_BB}\ar[rr]^-(0.5){\id_N\otimes\theta_B}\ar[d]&&{N\otimes_BB\otimes_{\co_X}\cE}\ar[d]\\
N\ar[rr]^-(0.5)\theta&&{N\otimes_{\co_X}\cE,}}
\end{equation}
where the vertical arrows are the isomorphisms induced by $\upmu$.  Indeed, if $U$ is an object of $X$, 
$e_1,\dots,e_d$ is an $\co_X$-basis of $\cE|U$ and $e^\vee_1,\dots,e^\vee_d$ is the dual basis of $\cE^\vee|U$, then we have 
\begin{eqnarray}
\theta&=&\sum_{i=1}^d\theta_i\otimes e_i,\label{p1-thbn43c}\\
\theta_B&=&\sum_{i=1}^d\theta_{B,i}\otimes e_i,\label{p1-thbn43d}
\end{eqnarray}
where $\theta_i$ (resp.\ $\theta_{B,i}$) is the morphism of $N$ (resp.\ $B$) induced by the multiplication by the class of $e^\vee_i$. 

The homomorphism $\upmu$ induces an isomorphism 
\begin{equation}\label{p1-thbn43e}
N\otimes_BB\otimes_{\co_X}\cC^\dagger \rightarrow N\otimes_{\co_X}\cC^\dagger.
\end{equation} 
By \eqref{p1-thbn43b}, it is compatible with the Higgs $\co_X$-fields indicated on the same line of the table
\begin{equation}\label{p1-thbn43f}
\begin{tabular}{|c|c|}
\hline
left hand side&  right  hand  side\\
\hline
$\id\otimes \theta_B\otimes \id$ &  $\theta\otimes \id$\\
\hline
$\id\otimes \id\otimes \delta_{\cC^\dagger}$ & $\id\otimes \delta_{\cC^\dagger}$\\
\hline
\end{tabular}
\end{equation}
We deduce an $\co_X$-linear morphism 
\begin{equation}\label{p1-thbn43g}
N\otimes_B\cL_B\rightarrow \uptau(N,\theta),
\end{equation}
where $\cL_B=\uptau(B,\theta_B)$ \eqref{p1-thbn42e}. Since \eqref{p1-thbn43e} is $B$-linear, 
\eqref{p1-thbn43g} is also $B$-linear by \ref{p1-thbn34}. It is therefore a morphism of Higgs $\co_X$-modules,
where $N\otimes_B\cL_B$ is equipped with its canonical Higgs field, see \ref{p1-delta-con7}(ii).

\begin{prop}\label{p1-thbn44}
Let $(N,\theta)$ be a Higgs $\co_X[\frac 1 p]$-module with coefficients in $\cE$.  
We denote by $\upmu\colon \rS[\frac 1 p]\rightarrow \cEnd_{\co_X}(N)$ 
the homomorphism of $\co_X[\frac 1 p]$-algebras defined by $\theta$. 
Let $B$ be a quotient $\co_X[\frac 1 p]$-algebra of $\rS[\frac 1 p]$ through which $\upmu$ factors, 
$\theta_B\colon B\rightarrow B\otimes_{\co_X}\cE$ its canonical Higgs $\co_X$-field. 
Suppose that the Higgs $\co_X[\frac 1 p]$-modules $(N,\theta)$ and $(B,\theta_B)$ are weakly twistable. 
Then, the $B$-module $\cL_B=\uptau(B,\theta_B)$ is invertible \eqref{p1-thbn42e}, 
and the canonical $B$-linear morphism \eqref{p1-thbn43g}
\begin{equation}\label{p1-thbn44a}
N\otimes_B\cL_B\rightarrow \uptau(N,\theta)
\end{equation}
is an isomorphism. 
\end{prop}

Indeed, the first assertion is mentioned as a reminder \eqref{p1-thbn42}. To prove that the morphism \eqref{p1-thbn44a} is an isomorphism, we may assume that the extension \eqref{p1-thbn1a} splits. 
The choice of a splitting determines a homomorphism 
$\varrho^\dagger \colon \cC^\dagger \rightarrow \co_X$, by \ref{p1-thbn16}. 
Since $(N,\theta)$ and $(B,\theta_B)$ are weakly twistable, $\varrho^\dagger$ induces two $B$-linear isomorphisms $\ttt_N\colon \uptau(N,\theta)\rightarrow N$ and 
$\ttt_B\colon \cL_B\rightarrow B$ by \ref{p1-thbn31}. The morphism \eqref{p1-thbn44a} then fits into a commutative diagram 
\begin{equation}\label{p1-thbn44b}
\xymatrix{
{N\otimes_B\uptau(B,\theta_B)}\ar[r]\ar[d]\ar@/_4pc/[dd]_{\id\otimes_B\ttt_{B}}&{\uptau(N,\theta)}\ar[d]\ar@/^4pc/[dd]^{\ttt_N}\\
{N\otimes_BB\otimes_{\co_X}\cC^\dagger}\ar[d]^{\id\otimes\varrho^\dagger}\ar[r]&
{N\otimes_{\co_X}\cC^\dagger}\ar[d]_{\id\otimes\varrho^\dagger}\\
{N\otimes_BB}\ar[r]&{N}}
\end{equation}
where the higher vertical arrows are induced by the canonical morphisms and the two lower horizontal arrows 
are induced by the multiplication by $B$ in $N$. We deduce that \eqref{p1-thbn44a} is an isomorphism.

\subsection{}\label{p1-thbn17}
Consider a second exact sequence of locally free $\co_X$-modules of finite type 
\begin{equation}\label{p1-thbn17a}
0\rightarrow \co_X\rightarrow \cF'\rightarrow \cE' \rightarrow 0,
\end{equation}
that fits into a commutative digram
\begin{equation}\label{p1-thbn17b}
\xymatrix{
0\ar[r]&\co_X\ar[r]\ar@{=}[d]&\cF\ar[r]\ar[d]^v&\cE\ar[r]\ar[d]^u&0\\
0\ar[r]&\co_X\ar[r]&\cF'\ar[r]&\cE'\ar[r]&0.}
\end{equation}
We associate with \eqref{p1-thbn17a} objects analogous to those associated with the extension \eqref{p1-thbn1a} 
that we denote by the same symbols equipped with a prime exponent: 
$\cC'$, $\cF'^{(r)}$ and $\cC'^{(r)}$, for rational numbers $r\geq 0$. 
We immediately see that for every integer $n\geq 0$, the induced diagram 
\begin{equation}\label{p1-thbn17c}
\xymatrix{
0\ar[r]&{\rS^{n-1}_{\co_X}(\cF)}\ar[r]\ar[d]_{\rS^{n-1}_{\co_X}(v)}&{\rS^{n}_{\co_X}(\cF)}\ar[r]\ar[d]^{\rS^n_{\co_X}(v)}&{\rS^{n}_{\co_X}(\cE)}\ar[r]\ar[d]^{\rS^n_{\co_X}(u)}&0\\
0\ar[r]&{\rS^{n-1}_{\co_X}(\cF')}\ar[r]&{\rS^n_{\co_X}(\cF')}\ar[r]&{\rS^n_{\co_X}(\cE')}\ar[r]&0}
\end{equation}
is commutative. Therefore, $v$ induces a homomorphism of $\co_X$-algebras 
\begin{equation}\label{p1-thbn17d}
w\colon \cC\rightarrow \cC'.
\end{equation}

For every rational number $r\geq 0$, $v$ induces an $\co_X$-linear morphism $v^{(r)}\colon \cF^{(r)}\rightarrow \cF'^{(r)}$ 
that fits into a commutative digram
\begin{equation}\label{p1-thbn17e}
\xymatrix{
0\ar[r]&\co_X\ar[r]\ar@{=}[d]&\cF^{(r)}\ar[r]\ar[d]^-(0.5){v^{(r)}}&\cE\ar[r]\ar[d]^u&0\\
0\ar[r]&\co_X\ar[r]&\cF'^{(r)}\ar[r]&\cE'\ar[r]&0.}
\end{equation}
We deduce a homomorphism of $\co_X$-algebras 
\begin{equation}\label{p1-thbn17f}
w^{(r)}\colon \cC^{(r)}\rightarrow \cC'^{(r)}.
\end{equation}
We denote by $\hw^{(r)}\colon \hcC^{(r)}\rightarrow \hcC'^{(r)}$ its extension to the $p$-adic completions. 
The following diagram 
\begin{equation}\label{p1-thbn17h}
\xymatrix{
{\cC^{(r)}}\ar[r]^-(0.5){d_{\cC^{(r)}}}\ar[d]_-(0.5){w^{(r)}}&{\cE\otimes_{\co_X}\cC^{(r)}}\ar[d]^-(0.5){u\otimes w^{(r)}}\\
{\cC'^{(r)}}\ar[r]^-(0.5){d_{\cC'^{(r)}}}&{\cE'\otimes_{\co_X}\cC'^{(r)}}}
\end{equation}
is clearly commutative.

For all rational numbers $r\geq r'\geq 0$, the following diagram 
\begin{equation}\label{p1-thbn17g}
\xymatrix{
\cC^{(r)}\ar[r]^{\alpha^{r,r'}}\ar[d]_-(0.5){w^{(r)}}&\cC^{(r')}\ar[d]^-(0.5){w^{(r')}}\\
\cC'^{(r)}\ar[r]^{\alpha'^{r,r'}}&\cC'^{(r')}}
\end{equation}
is clearly commutative. Therefore, the $\hw^{(t)}$'s, for $t\in \mQ_{>r}$, induce a homomorphism of $\co_X$-algebras 
\begin{equation}\label{p1-thbn17i}
\hw^{(r+)}\colon \hcC^{(r+)}\rightarrow \hcC'^{(r+)}.
\end{equation}
For simplicity, we set $w^\dagger=\hw^{(0+)}\colon \cC^\dagger\rightarrow \hcC'^\dagger$.

\subsection{}\label{p1-thbn19}
We keep the notation and assumptions of \ref{p1-thbn17}. 
Composing with $\id\otimes u$ defines a functor that we denote by
\begin{equation}\label{p1-thbn19a}
\upmu\colon \bHM(\co_X,\cE)\rightarrow \bHM(\co_X,\cE').
\end{equation}
By \ref{p1-delta-con5} and \eqref{p1-thbn17h}, for every rational number $r\geq 0$, the extension of scalars by $\hw^{(r)}$ \eqref{p1-thbn17f} defines a functor
\begin{equation}\label{p1-thbn19b}
\upomega^{(r)}\colon \bMIC_r(\hcC^{(r)}/\co_X)\rightarrow \bMIC_r(\hcC'^{(r)}/\co_X),
\end{equation}
where the target is the category of $\hcC'^{(r)}$-modules with integrable $\delta_{\hcC'^{(r)}}$-connection \eqref{p1-delta-con2}, 
see \eqref{p1-delta-con5b}. 

We have a canonical isomorphism of functors from $\bHM(\co_X,\cE)$ to $\bMIC_r(\hcC'^{(r)}/\co_X)$ 
\begin{equation}\label{p1-thbn19d}
\upomega^{(r)}\circ \uppi^{(r)}\stackrel{\sim}{\longrightarrow} \uppi'^{(r)}\circ \upmu,
\end{equation}
where the functor $\uppi'^{(r)}$ is the analogue of $\uppi^{(r)}$ \eqref{p1-thbn10a} for the extension \eqref{p1-thbn17a}.

Let $(N,\theta)$ be an object of $\bHM(\co_X,\cE)$, $(N,\theta')=\upmu(N,\theta)$. The diagram 
\begin{equation}\label{p1-thbn19f}
\xymatrix{
{\cC^\dagger\otimes_{\co_X}N}\ar[r]^-(0.5){\theta_\tot}\ar[d]&{\cE\otimes_{\co_X}\cC^\dagger\otimes_{\co_X}N}\ar[d]\\
{\cC'^\dagger\otimes_{\co_X}N}\ar[r]^-(0.5){\theta'_\tot}&{\cE'\otimes_{\co_X}\cC'^\dagger\otimes_{\co_X}N,}}
\end{equation}
where $\theta_\tot$ is defined in \eqref{p1-thbn9a}, $\theta'_\tot$ is defined likewise for $\theta'$ instead of $\theta$ and the vertical morphisms are induced by $u$ and $w^\dagger$ \eqref{p1-thbn17i}, 
is commutative. We deduce a 
canonical morphism of functors from $\bHM(\co_X,\cE)$ to $\bHM(\co_X,\cE')$
\begin{equation}\label{p1-thbn19g}
\upmu\circ \uptau\rightarrow \uptau'\circ \upmu, 
\end{equation}
where the functor $\uptau'$  is the analogue of $\uptau$ \eqref{p1-thbn9b} for the extension \eqref{p1-thbn17a}.

\begin{lem}\label{p1-thbn20}
We keep the assumptions and notation of \ref{p1-thbn19}. 
Let $M$ be an $\co_X$-module, $(N,\theta)$ a Higgs $\co_X$-module with coefficients in $\cE$, $r$ a rational number $>0$.
If $M$ and $(N,\theta)$ are $r$-associated relatively to the extension \eqref{p1-thbn1a}, in the sense of  \ref{p1-thbn12}, 
then they are also $r$-associated relatively to the extension \eqref{p1-thbn17a}. 
\end{lem}
Indeed, an isomorphism of $\bMIC_r(\hcC^{(r)}/\co_X)$ 
\begin{equation}
\alpha\colon \fS^{(r)}(M)\stackrel{\sim}{\rightarrow}\uppi^{(r)}(N,\theta)
\end{equation}
induces, by \eqref{p1-thbn19d}, an isomorphism of $\bMIC_r(\hcC'^{(r)}/\co_X)$ 
\begin{equation}
\upomega^{(r)}(\alpha)\colon \fS'^{(r)}(M)\stackrel{\sim}{\rightarrow}\uppi'^{(r)}(\upmu(N,\theta)).
\end{equation}

\begin{prop}\label{p1-thbn27}
We keep the assumptions and notation of \ref{p1-thbn19}. Let $(N,\theta)$ be a Higgs $\co_X[\frac 1 p]$-module with coefficients in $\cE$, 
twistable by the extension \eqref{p1-thbn1a}. Then, the Higgs module $\upmu(N,\theta)$ is twistable by the extension \eqref{p1-thbn17a}, and the canonical morphism \eqref{p1-thbn19g}
\begin{equation}\label{p1-thbn27a}
\upmu(\uptau(N,\theta))\rightarrow \uptau'(\upmu(N,\theta))
\end{equation}
is an isomorphism.
\end{prop}

Indeed, there exist an $\co_X$-module $M$, a rational number $r>0$, and an isomorphism of $\bMIC_{r}(\hcC^{(r)}/\co_X)$ \eqref{p1-thbn12}
\begin{equation}\label{p1-thbn27b}
\alpha\colon \fS^{(r)}(M)\stackrel{\sim}{\rightarrow}\uppi^{(r)}(N,\theta). 
\end{equation}
As in the proof of \ref{p1-thbn21}, we deduce an isomorphism \eqref{p1-thbn21c}
\begin{equation}\label{p1-thbn27c}
\beta\colon M \stackrel{\sim}{\rightarrow} \uuptau(N,\theta).
\end{equation}
In view of \eqref{p1-thbn19d}, $\upomega^{(r)}(\alpha)$ induces an isomorphism of $\bMIC_r(\hcC'^{(r)}/\co_X)$
\begin{equation}\label{p1-thbn27d}
\alpha'\colon \fS'^{(r)}(M)\stackrel{\sim}{\rightarrow}\uppi'^{(r)}(\upmu(N,\theta)),
\end{equation}
where the functor $\fS'^{(r)}$ is defined as in \eqref{p1-thbn10b} for the extension \eqref{p1-thbn17a}. 
Therefore, the Higgs module $\upmu(N,\theta)$ is twistable by the extension \eqref{p1-thbn17a}. 
Similarly to $\alpha$, we deduce from $\alpha'$ an isomorphism 
\begin{equation}\label{p1-thbn27e}
\beta'\colon M \stackrel{\sim}{\rightarrow} \uuptau'(\upmu(N,\theta)).
\end{equation}
It is clear that the diagram of $\co_X$-linear morphisms 
\begin{equation}\label{p1-thbn27f}
\xymatrix{
{M}\ar[r]^-(0.5)\beta\ar[rd]_{\beta'}&{\uuptau(N,\theta)}\ar[d]\\
&{\uuptau'(\upmu(N,\theta)),}}
\end{equation}
where the vertical arrow is the morphism underlying \eqref{p1-thbn19g}, is commutative. The proposition follows.

\subsection{}\label{p1-thbn28}
Let $f\colon (X',\co_{X'})\rightarrow (X,\co_X)$ be a morphism of ringed topos, such that $\co_{X'}$ is $V$-flat and $p$-adic. 
The exact sequence \eqref{p1-thbn1a} induces by pullback an exact sequence of $\co_{X'}$-modules
\begin{equation}\label{p1-thbn28a}
0\rightarrow \co_{X'}\rightarrow f^*(\cF)\rightarrow f^*(\cE) \rightarrow 0
\end{equation}
for which we associate objects analogous to those associated with the extension \eqref{p1-thbn1a} that we equip with a $'$ exponent:
$\cF'^{(r)}$ and $\cC'^{(r)}$, $\delta_{\cC'^{(r)}}$ for $r\in \mQ_{\geq 0}$. The morphism $f$ defines by pullback a functor 
\begin{equation}\label{p1-thbn28b}
f^*\colon \bHM(\co_{X},\cE)\rightarrow \bHM(\co_{X'},f^*(\cE)),
\end{equation}
where the target is the category of Higgs $\co_{X'}$-modules with coefficients in $f^*(\cE)$

Let $r$ be a rational number $\geq 0$. We have a canonical isomorphism of $\co_{X'}$-algebras 
\begin{equation}\label{p1-thbn28c}
w^{(r)}\colon f^*(\cC^{(r)})\stackrel{\sim}{\rightarrow} \cC'^{(r)}.
\end{equation} 
It induces a homomorphism of $\co_{X'}$-algebras $\omega^{(r)}\colon f^*(\hcC^{(r)})\rightarrow \hcC'^{(r)}$.
The following diagram 
\begin{equation}\label{p1-thbn28d}
\xymatrix{
{f^*(\hcC^{(r)})}\ar[rr]^-(0.5){f^*(d_{\hcC^{(r)}})}\ar[d]_-(0.5){\omega^{(r)}}&&{f^*(\cE\otimes_{\co_{X}}\hcC^{(r)})}\ar[d]^-(0.5){\id\otimes \omega^{(r)}}\\
{\hcC'^{(r)}}\ar[rr]^-(0.5){d_{\hcC'^{(r)}}}&&{f^*(\cE)\otimes_{\co_{X'}}\hcC'^{(r)}}}
\end{equation}
is clearly commutative

For all rational numbers $r\geq r'\geq 0$, the following diagram 
\begin{equation}\label{p1-thbn28e}
\xymatrix{
{f^*(\hcC^{(r)})}\ar[rr]^{f^*(\halpha^{r,r'})}\ar[d]_-(0.5){\omega^{(r)}}&&{f^*(\hcC^{(r')})}\ar[d]^-(0.4){\omega^{(r')}}\\
\hcC'^{(r)}\ar[rr]^{\halpha'^{r,r'}}&&\hcC'^{(r')}}
\end{equation}
is clearly commutative. Therefore, the $\omega^{(t)}$'s, for $t\in \mQ_{>r}$, induce a homomorphism of $\co_{X'}$-algebras 
\begin{equation}\label{p1-thbn28f}
\omega^{(r+)}\colon f^*(\hcC^{(r+)})\rightarrow \hcC'^{(r+)}.
\end{equation}
For simplicity, we set $\omega^\dagger=\omega^{(0+)}\colon f^*(\cC^\dagger)\rightarrow \cC'^\dagger$.

The pullback by $f$ and the extension of scalars by $\omega^{(r)}$ define a functor 
\begin{equation}\label{p1-thbn28g}
f^{(r)+}\colon \bMIC_r(\hcC^{(r)}/\co_X)\rightarrow \bMIC_r(\hcC'^{(r)}/\co_{X'}),
\end{equation}
where the target is the category of $\hcC'^{(r)}$-modules with integrable $\delta_{\hcC'^{(r)}}$-connection, see \eqref{p1-delta-con5b}.

There is a canonical isomorphism of functors from $\bHM(\co_X,\cE)$ to $\bMIC_r(\hcC'^{(r)}/\co_{X'})$ 
\begin{equation}\label{p1-thbn28i}
f^{(r)+}\circ \uppi^{(r)}\stackrel{\sim}{\longrightarrow} \uppi'^{(r)}\circ f^*,
\end{equation}
where the functor $\uppi'^{(r)}$ is the analogue of $\uppi^{(r)}$ \eqref{p1-thbn10a} for the extension \eqref{p1-thbn28a}.

Let $(N,\theta)$ be an object of $\bHM(\co_{X},\cE)$, $(N',\theta')=f^*(N,\theta)$. The diagram 
\begin{equation}\label{p1-thbn28k}
\xymatrix{
{f^*(\cC^\dagger\otimes_{\co_X}N)}\ar[r]^-(0.5){f^*(\theta_\tot)}\ar[d]&{f^*(\cE\otimes_{\co_X}\cC^\dagger\otimes_{\co_X}N)}\ar[d]\\
{\cC'^\dagger\otimes_{\co_{X'}}N'}\ar[r]^-(0.5){\theta'_\tot}&{f^*(\cE)\otimes_{\co_{X'}}\cC'^\dagger\otimes_{\co_{X'}}N',}}
\end{equation}
where $\theta_\tot$ is defined in \eqref{p1-thbn9a}, $\theta'_\tot$ is defined likewise and the vertical arrows are induced by $\omega^\dagger$ \eqref{p1-thbn28f}, is commutative. 
We deduce a canonical morphism of functors from $\bHM(\co_X,\cE)$ to $\bHM(\co_{X'},\cE')$
\begin{equation}\label{p1-thbn28j}
f^*\circ \uptau\rightarrow \uptau'\circ f^*, 
\end{equation}
where the functor $\uptau'$  is the analogue of $\uptau$ \eqref{p1-thbn9b} for the extension \eqref{p1-thbn28a}.

\begin{prop}\label{p1-thbn29}
We keep the assumptions and notation of \ref{p1-thbn28}. Let $(N,\theta)$ be a Higgs $\co_X[\frac 1 p]$-module with coefficients in $\cE$ 
twistable by the extension \eqref{p1-thbn1a}. Then, the Higgs module $f^*(N,\theta)$ is twistable by the extension \eqref{p1-thbn28a}, 
and the canonical morphism \eqref{p1-thbn28j}
\begin{equation}\label{p1-thbn29a}
f^*(\uptau(N,\theta))\rightarrow \uptau'(f^*(N,\theta))
\end{equation}
is an isomorphism.
\end{prop}

Indeed, there exist an $\co_X$-module $M$, a rational number $r>0$, and an isomorphism of $\bMIC_{r}(\hcC^{(r)}/\co_X)$ \eqref{p1-thbn12}
\begin{equation}\label{p1-thbn29b}
\alpha\colon \fS^{(r)}(M)\stackrel{\sim}{\rightarrow}\uppi^{(r)}(N,\theta). 
\end{equation}
As in the proof of \ref{p1-thbn21}, we deduce an isomorphism \eqref{p1-thbn21c}
\begin{equation}\label{p1-thbn29c}
\beta\colon M \stackrel{\sim}{\rightarrow} \uuptau(N,\theta).
\end{equation}
In view of \eqref{p1-thbn28i}, $f^{(r)+}(\alpha)$ induces an isomorphism of $\bMIC_r(\hcC'^{(r)}/\co_{X'})$
\begin{equation}\label{p1-thbn29d}
\alpha'\colon \fS'^{(r)}(f^*(M))\stackrel{\sim}{\rightarrow}\uppi'^{(r)}(f^*(N,\theta)),
\end{equation}
where the functor $\fS'^{(r)}$ is the analogue of $\fS^{(r)}$ \eqref{p1-thbn10b} for the extension \eqref{p1-thbn28a}.
Therefore, the Higgs module $f^*(N,\theta)$ is twistable by the extension \eqref{p1-thbn28a}. 
Similarly to $\alpha$, we deduced from $\alpha'$ an isomorphism 
\begin{equation}
\beta'\colon f^*(M) \stackrel{\sim}{\rightarrow} \uuptau'(f^*(N,\theta)).
\end{equation}
It is clear that the diagram of $\co_{X'}$-linear morphisms 
\begin{equation}
\xymatrix{
{f^*(M)}\ar[r]^-(0.5){f^*(\beta)}\ar[rd]_{\beta'}&{f^*(\uuptau(\cN,\theta))}\ar[d]\\
&{\uuptau'(f^*(\cN,\theta)),}}
\end{equation}
where the vertical arrow is the morphism underlying \eqref{p1-thbn28j}, is commutative. The proposition follows.

\section{Twisted pullback and higher direct images of Higgs modules}\label{p1-tphdi}

\subsection{}\label{p1-tphdi1}
In this section, $f\colon (X',\co_{X'})\rightarrow (X,\co_X)$ deotes a morphism of ringed $\mU$-topos, 
where $\co_X$ and $\co_{X'}$ are $p$-adic $V$-flat $V$-algebras \eqref{p1-thbn0}, 
and the canonical homomorphism $f^*(\co_X)\rightarrow \co_{X'}$ is a morphism of $V$-algebras. 
Let $\cE$ be a locally free $\co_X$-module of finite type, 
\begin{equation}\label{p1-tphdi1a}
0\rightarrow \co_{X'}\rightarrow \cF\rightarrow f^*(\cE) \rightarrow 0
\end{equation}
an exact sequence of $\co_{X'}$-modules, $\cE'$ a locally free $\co_{X'}$-module of finite type, and 
\begin{equation}\label{p1-tphdi1b}
u\colon f^*(\cE)\rightarrow \cE'
\end{equation} 
an $\co_{X'}$-linear morphism. We will add further assumptions in \ref{p1-tphdi4}. 

Recall that for any rational number $r\geq 0$, we associated with the extension \eqref{p1-tphdi1a} the $\co_{X'}$-algebra $\cC^{(r)}$ \eqref{p1-thbn3b}, 
its $p$-adic completion $\hcC^{(r)}$, the $\co_{X'}$-algebra $\hcC^{(r+)}$ \eqref{p1-thbn6e} and the $\co_{X'}$-derivations 
\begin{eqnarray}\label{p1-tphdi5a}
d_{\hcC^{(r)}}\colon \hcC^{(r)}\rightarrow f^*(\cE)\otimes_{\co_{X'}}\hcC^{(r)},\\
\delta_{\hcC^{(r+)}}\colon \hcC^{(r+)}\rightarrow f^*(\cE)\otimes_{\co_{X'}}\hcC^{(r+)},
\end{eqnarray}
defined in \eqref{p1-thbn6a} and \eqref{p1-thbn6f}, respectively, which are also Higgs $\co_{X'}$-fields. 
For simplicity, we set $\cC^\dagger=\hcC^{(0+)}$. By  \eqref{p1-thbn6d}, $\delta_{\hcC^{(0+)}}$ and $d_{\hcC^{(0)}}$ are compatible via 
the canonical homomorphism $\cC^\dagger\rightarrow \hcC=\hcC^{(0)}$ \eqref{p1-thbn6c}. 
Hence, there is no risk of ambiguity to denote $\delta_{\hcC^{(0+)}}$ by 
\begin{equation}\label{p1-tphdi5b}
d\colon \cC^\dagger\rightarrow f^*(\cE)\otimes_{\co_{X'}}\cC^\dagger.
\end{equation}
We set 
\begin{equation}\label{p1-tphdi5c}
d'=(u\otimes \id)\circ d\colon \cC^\dagger\rightarrow \cE'\otimes_{\co_{X'}}\cC^\dagger.
\end{equation}

\subsection{}\label{p1-tphdi2}
Recall that we defined in \eqref{p1-thbn9b} a twisting functor by the extension \eqref{p1-tphdi1a} 
for Higgs $\co_{X'}$-modules with coefficients in $f^*(\cE)$:
\begin{equation}\label{p1-tphdi2a}
\uptau\colon \bHM(\co_{X'},f^*(\cE))\rightarrow \bHM(\co_{X'},f^*(\cE)).
\end{equation}
However, it only behaves well for {\em weakly twistable} Higgs bundles \eqref{p1-thbn30}. 
Composing with the pullback functor $f^*$, we obtain a functor
\begin{equation}\label{p1-tphdi2b}
f^*_\uptau\colon 
\begin{array}[t]{clcr}
\bHM(\co_X,\cE)&\rightarrow& \bHM(\co_{X'},f^*(\cE))\\
(N,\theta)&\mapsto&\uptau(f^*(N),f^*(\theta)),
\end{array}
\end{equation} 
that we call the {\em pullback by $f$ twisted by the extension \eqref{p1-tphdi1a}} (or simply the {\em twisted pullback by $f$} when the extension is implicit).

\begin{defi}\label{p1-tphdi3}
We say that a Higgs $\co_X[\frac 1 p]$-module $(N,\theta)$ with coefficients in $\cE$ is 
{\em weakly twistable} (resp.\ {\em twistable}) by the extension \eqref{p1-tphdi1a} 
if the Higgs $\co_{X'}[\frac 1 p]$-module $(f^*(N),f^*(\theta))$ with coefficients 
in $f^*(\cE)$ is weakly twistable (resp.\ twistable) by the extension \eqref{p1-tphdi1a} 
in the sense of \ref{p1-thbn30} (resp.\ \ref{p1-thbn14}). 
\end{defi}

\subsection{}\label{p1-tphdi6}
Let $(N,\theta)$ be a Higgs $\co_{X'}$-module with coefficients in $\cE'$. 
We equip $N\otimes_{\co_{X'}}\cC^\dagger$ with the Higgs $\co_{X'}$-field 
\begin{equation}\label{p1-tphdi6a}
\vartheta= \theta\otimes \id+\id\otimes d'\colon N\otimes_{\co_{X'}}\cC^\dagger \rightarrow \cE'\otimes_{\co_{X'}}N\otimes_{\co_{X'}}\cC^\dagger,
\end{equation}
and denote by $\mK^{\bullet}(N\otimes_{\co_{X'}}\cC^\dagger)$ the associated Dolbeault complex, and by $\vartheta^\bullet$ its differentials.  
The Higgs field $\id\otimes d$ on $N\otimes_{\co_{X'}}\cC^\dagger$ induces a morphism of complexes of $\co_{X'}$-modules
\begin{equation}\label{p1-tphdi6b}
\fd \colon \mK^\bullet(N\otimes_{\co_{X'}}\cC^\dagger)\rightarrow f^*(\cE)\otimes_{\co_{X'}} \mK^\bullet(N\otimes_{\co_{X'}}\cC^\dagger).
\end{equation} 
Indeed, the question being local on $X'$, we may assume that the $\co_X$-module $\cE$ is free of finite type. 
The assertion then follows by considering  the components, with respect to a basis of $f^*(\cE)$, 
of the derivation $d$ of $\cC^\dagger$, which are $\co_{X'}$-linear endomorphisms of $\cC^\dagger$ 
commuting with each other since $d$ is a Higgs field. 

By the projection formula, for any integer $q\geq 0$, $\rR^qf_*(\fd)$ identifies with an $\co_X$-linear morphism
\begin{equation}\label{p1-tphdi6c}
\rR^qf_*(\fd)\colon  \rR^qf_*(\mK^\bullet(N\otimes_{\co_{X'}}\cC^\dagger))\rightarrow \cE\otimes_{\co_X} \rR^qf_*(\mK^\bullet(N\otimes_{\co_{X'}}\cC^\dagger)). 
\end{equation}
By choosing locally a basis of $\cE$, we easily check that it is a Higgs $\co_X$-field. We thus obtain a functor that we denote by 
\begin{equation}\label{p1-tphdi6d}
\rR^qf^\uptau_*\colon 
\begin{array}[t]{clcr}
\bHM(\co_{X'},\cE')&\rightarrow&\bHM(\co_X,\cE)\\
(N,\theta)&\mapsto&(\rR^qf_*(\mK^\bullet(N\otimes_{\co_{X'}}\cC^\dagger)),-\rR^qf_*(\fd)),
\end{array}
\end{equation}
and call it the {\em $q$th higher direct image functor by $f$ twisted by the extension $\cF$ \eqref{p1-tphdi1a}} 
(or simply the {\em twisted $q$th higher direct image functor by $f$} when the extension is implicit). 
Although omitted from the notation, this notion depends also on the morphism $u$ \eqref{p1-tphdi1b}.

\subsection{}\label{p1-tphdi4}
For the remainder of this section, we assume that the morphism $u$ \eqref{p1-tphdi1b} 
fits into an exact sequence of locally free $\co_{X'}$-modules of finite type
\begin{equation}\label{p1-tphdi4a}
\xymatrix{
0\ar[r]&{f^*(\cE)}\ar[r]^u&{\cE'}\ar[r]^{u'}&{\ucE'}\ar[r]&0.}
\end{equation}
For any integers $i,n\geq 0$, we denote by $\rW^i\wedge^n\cE'$ the image of the morphism 
\begin{equation}\label{p1-tphdi4i}
\wedge^i f^*(\cE)\otimes_{\co_{X'}}\wedge^{n-i}\cE'\rightarrow \wedge^n\cE'
\end{equation}
induced by $u$. 
We thus define an exhaustive decreasing filtration $(\rW^i \wedge^\bullet \cE')_{i\geq 0}$ of the exterior algebra $\wedge^\bullet \cE'$  by ideals,
called the {\em Koszul filtration associated with the extension \eqref{p1-tphdi4a}}. 
We denote by $\Gr_\rW^\bullet\wedge^\bullet \cE'$ the associated graded algebra, which is a bigraded algebra.  
The canonical injection $\ucE' \oplus f^*(\cE) \rightarrow \Gr^\bullet_\rW\wedge^\bullet \cE'$
extends canonically into a homomorphism of bigraded algebras
\begin{equation}\label{p1-tphdi4j}
\wedge^\bullet(\ucE' \oplus f^*(\cE))\rightarrow \Gr^\bullet_\rW\wedge^\bullet \cE'.
\end{equation}
It is an isomorphism according to (\cite{illusie1} V 4.1.6).
Furthermore, by virtue of (\cite{alg1-3} III § 7.7 prop.~10), we have a canonical isomorphism of bigraded algebras
\begin{equation}\label{p1-tphdi4k}
\wedge^\bullet \ucE'\ {^g\otimes} \wedge^\bullet f^*(\cE) \stackrel{\sim}{\rightarrow}\wedge^\bullet(\ucE' \oplus f^*(\cE)),
\end{equation}
where the symbol ${^g\otimes}$ designates the left tensor product (see \cite{alg1-3} III § 4.7 remarks page 49).
In particular, for any integer $i\geq 0$, we have a canonical isomorphism
\begin{equation}\label{p1-tphdi4l}
\Gr^i_\rW\wedge^\bullet \cE'\stackrel{\sim}{\rightarrow}\wedge^i f^*(\cE)\otimes_A\wedge^{\bullet -i} \ucE'.
\end{equation}

Let $(N,\theta)$ be a  Higgs $\co_{X'}$-module with coefficients in $\cE'$.
We denote by $\utheta\colon N\rightarrow N\otimes_{\co_{X'}}\ucE'$ the  Higgs $\co_{X'}$-field induced by $\theta$,
and by $\mK^\bullet(N)$ (resp.\ $\umK^\bullet(N)$) the Dolbeault complex of $(N,\theta)$ (resp.\ $(N,\utheta)$) \eqref{p1-delta-con1b}.
For any integers $i,n\geq 0$, we set
\begin{equation}\label{p1-tphdi4b}
\rW^i\mK^n(N)=N\otimes_{\co_{X'}}(\rW^i\wedge^n\cE').
\end{equation}
We thus define an exhaustive decreasing filtration of the $\wedge^\bullet \cE'$-graded module $\mK^\bullet(N)$
by graded sub-$\wedge^\bullet \cE'$-modules $(\rW^i\mK^\bullet(N))_{i\geq 0}$, stable by the differential of $\mK^\bullet(N)$,
called the {\em Koszul filtration of $\mK^\bullet(N)$ associated with the exact sequence \eqref{p1-tphdi4a}};
in particular, it is a filtration of the complex $\mK^\bullet(N)$ by subcomplexes. We denote by $\Gr_\rW^\bullet\mK^\bullet(N)$ the 
associated graded complex.

By \eqref{p1-tphdi4j} and \eqref{p1-tphdi4k}, for every integer $i\geq 0$, we have a canonical isomorphism of complexes
\begin{equation}\label{p1-tphdi4c}
\Gr_\rW^i\mK^\bullet(N)\stackrel{\sim}{\rightarrow}\wedge^i f^*(\cE)\otimes_{\co_{X'}} \umK^{\bullet}(N)[-i],
\end{equation}
where the differentials of $\umK^{\bullet}(N)[-i]$ are those of $\umK^\bullet(N)$ multiplied by $(-1)^i$.
We therefore have an exact sequence of complexes of $\co_{X'}$-modules
\begin{equation}\label{p1-tphdi4d}
0\rightarrow \wedge^{i+1} f^*(\cE)\otimes_{\co_{X'}} \umK^{\bullet}(N)[-i-1]\rightarrow \rW^i/\rW^{i+2}(\mK^\bullet(N))
\rightarrow \wedge^i f^*(\cE)\otimes_{\co_{X'}} \umK^{\bullet}(N)[-i]\rightarrow 0.
\end{equation}
This induces a morphism in the derived category $\bD^+(\bMod(\co_{X'}))$,
\begin{equation}\label{p1-tphdi4e}
\partial^i\colon \wedge^i f^*(\cE)\otimes_{\co_{X'}} \umK^{\bullet}(N)\rightarrow \wedge^{i+1}f^*(\cE)\otimes_{\co_{X'}} \umK^{\bullet}(N),
\end{equation}
which we will call the {\em boundary associated with the Koszul filtration of $\mK^\bullet(N)$}.

By virtue of (\cite{sp} \href{https://stacks.math.columbia.edu/tag/015W}{015W}), there exists a canonical convergent spectral sequence
\begin{equation}\label{p1-tphdi4f}
\rE_1^{i,j}=\rR^{i+j}f_*(\Gr_\rW^i(\mK^\bullet(N)))\Rightarrow \rR^{i+j}f_*(\mK^\bullet(N)),
\end{equation}
called the hypercohomology spectral sequence of the filtered complex $\mK^\bullet(N)$.
By \eqref{p1-tphdi4c} and the projection formula (\cite{sp} \href{https://stacks.math.columbia.edu/tag/0B54}{0B54}), we have 
a canonical isomorphism
\begin{equation}\label{p1-tphdi4g}
\rE_1^{i,j}\stackrel{\sim}{\rightarrow} \rR^jf_*(\umK^\bullet(N)) \otimes_{\co_\fX}\wedge^i\cE.
\end{equation}
For every integer $j\geq 0$, we obtain a complex of $\co_X$-modules
\begin{equation}\label{p1-tphdi4h}
\rR^jf_*(\umK^\bullet(N))\stackrel{d_1^{0,j}}{\longrightarrow} \rR^jf_*(\umK^\bullet(\cN))\otimes_{\co_X}\cE\stackrel{d_1^{1,j }}{\longrightarrow}
\rR^jf_*(\umK^\bullet(N))\otimes_{\co_X}\wedge^2\cE\stackrel{d_1^{2,j}}{\longrightarrow} \dots
\end{equation}
whose differentials are none other than the higher direct images of the boundary morphisms \eqref{p1-tphdi4e}.
The morphism $d_1^{0,j}$ is a Higgs $\co_X$-field on $\rR^jf_*(\umK^\bullet(N))$ with coefficients in $\cE$, 
that we call the {\em Katz-Oda field}, whose Dolbeault complex is none other than the complex $(\rE_1^{i,j})_{i\geq 0}$ (\cite{ag2} 2.5.16).

\subsection{}\label{p1-tphdi5}
We set 
\begin{equation}\label{p1-tphdi5d}
\tcE=f^*(\cE)\oplus \cE'
\end{equation}
and consider the exact sequences 
\begin{eqnarray}
\xymatrix{
0\ar[r]&{f^*(\cE)}\ar[r]^{\iota}&{\tcE}\ar[r]^{\iota'}&{\cE'}\ar[r]&0,}\label{p1-tphdi5e}\\
\xymatrix{
0\ar[r]&{f^*(\cE)}\ar[r]^{\gamma}&{\tcE}\ar[r]^{\gamma'}&{\cE'}\ar[r]&0,}\label{p1-tphdi5f}
\end{eqnarray}
where for all local sections $x$ of $f^*(\cE)$ and $x'$ of $\cE'$, $\iota(x)=(x,0)$, 
$\iota'(x,x')=x'$, $\gamma(x)=(x,u(x))$ and $\gamma'(x,x')=x'-u(x)$ \eqref{p1-tphdi4a}.

We equip the exterior algebra $\wedge^\bullet \tcE$ with the Koszul filtration $(\rW^i \wedge^\bullet \tcE)_{i\geq 0}$ associated with the extension \eqref{p1-tphdi5e}. 
By virtue of (\cite{alg1-3} III § 7.7 prop.~10), we have a canonical isomorphism of bigraded algebras
\begin{equation}\label{p1-tphdi5h}
\wedge^\bullet f^*(\cE)\ {^g\otimes} \wedge^\bullet \cE' \stackrel{\sim}{\rightarrow}\wedge^\bullet\tcE,
\end{equation}
where the symbol ${^g\otimes}$ denotes the left tensor product (see \cite{alg1-3} III § 4.7 remarks page 49). 
It induces, for all integers $i,n\geq 0$, an $\co_{X'}$-linear isomorphism
\begin{equation}\label{p1-tphdi5i}
\oplus_{j\geq i} \wedge^jf^*(\cE)\otimes_{\co_{X'}}\wedge^{n-j}\cE' \stackrel{\sim}{\rightarrow}\rW^i \wedge^n\tcE.
\end{equation}

We denote abusively by $\gamma'\colon \wedge^\bullet \tcE\rightarrow \wedge^\bullet \cE'$ the morphism induced by $\gamma'$ \eqref{p1-tphdi5f}. We see easily that for all integers 
$i,n\geq 0$, we have 
\begin{equation}\label{p1-tphdi5j}
\gamma'(\rW^i \wedge^n\tcE)=\rW^i \wedge^n\cE',
\end{equation}
where $\rW^i \wedge^\bullet\cE'$ is the Koszul filtration of the exterior algebra $\wedge^\bullet\cE'$ associated with the extension \eqref{p1-tphdi4a}. 

We denote by
\begin{equation}\label{p1-tphdi5g}
\td\colon \cC^\dagger\rightarrow \tcE\otimes_{\co_{X'}}\cC^\dagger
\end{equation}
the $\co_{X'}$-derivation $\td=(\gamma\otimes \id)\circ d=d\oplus d'$, where $d$ and $d'$ are the derivations 
defined in \eqref{p1-tphdi5b} and \eqref{p1-tphdi5c}.

\subsection{}\label{p1-tphdi7}
We take again the assumption and notation of \ref{p1-tphdi6} and consider the $\co_{X'}$-linear morphism 
\begin{equation}\label{p1-tphdi7a}
\tvartheta=(\id\otimes d)\oplus \vartheta\colon N\otimes_{\co_{X'}}\cC^\dagger \rightarrow \tcE\otimes_{\co_{X'}}N\otimes_{\co_{X'}}\cC^\dagger,
\end{equation}
where the direct sum corresponds to the decomposition \eqref{p1-tphdi5d} of $\tcE$. Since $d$ \eqref{p1-tphdi5b} is a Higgs field, 
$\tvartheta$ is a Higgs $\co_{X'}$-field on $N\otimes_{\co_{X'}}\cC^\dagger$ with coefficients in  $\tcE$;
see the proof of the existence of $\fd$ \eqref{p1-tphdi6b}. 
We denote by $\tmK^\bullet(N\otimes_{\co_{X'}}\cC^\dagger)$ the Dolbeault complex of the Higgs module $(N\otimes_{\co_{X'}}\cC^\dagger,\tvartheta)$,
and by $\tvartheta^\bullet$ its differentials.
We equip it with the Koszul filtration $\rW^\bullet\tmK^\bullet(N\otimes_{\co_{X'}}\cC^\dagger)$ associated with the extension \eqref{p1-tphdi5e}, see \ref{p1-tphdi4}. 
For all integers $i,n\geq 0$, we have
\begin{equation}\label{p1-tphdi7g}
\rW^i\tmK^n(N\otimes_{\co_{X'}}\cC^\dagger)=N\otimes_{\co_{X'}}\cC^\dagger\otimes_{\co_{X'}}(\rW^i\wedge^n\tcE),
\end{equation}
where $\rW^i\wedge^\bullet\tcE$ is the Koszul filtration of the exterior algebra $\wedge^\bullet\tcE$  associated with the extension \eqref{p1-tphdi5e}. 
Hence, by \eqref{p1-tphdi5i}, we have a canonical $\co_{X'}$-linear isomorphism
\begin{equation}\label{p1-tphdi7b}
\rW^i\tmK^n(N\otimes_{\co_{X'}}\cC^\dagger) \stackrel{\sim}{\rightarrow} \oplus_{j\geq i} \wedge^jf^*(\cE)\otimes_{\co_{X'}}\mK^{n-j}(N\otimes_{\co_{X'}}\cC^\dagger),  
\end{equation}
where $\mK^\bullet(N\otimes_{\co_{X'}}\cC^\dagger)$ is the Dolbeault complex of the Higgs module \eqref{p1-tphdi6a}. 
Moreover, the differentials $\tvartheta^\bullet$ are the sum of the morphisms 
\begin{equation}\label{p1-tphdi7c}
(-1)^i\id \otimes \vartheta^{n}\colon
\wedge^if^*(\cE)\otimes_{\co_{X'}}\mK^{n}(N\otimes_{\co_{X'}}\cC^\dagger)\rightarrow \wedge^if^*(\cE)\otimes_{\co_{X'}}\mK^{n+1}(N\otimes_{\co_{X'}}\cC^\dagger), 
\end{equation}
\begin{equation}\label{p1-tphdi7d}
\wedge^if^*(\cE)\otimes_{\co_{\fX'}}\mK^{n}(N\otimes_{\co_{X'}}\cC^\dagger)\rightarrow \wedge^{i+1}f^*(\cE)\otimes_{\co_{X'}}\mK^{n}(N\otimes_{\co_{X'}}\cC^\dagger), 
\end{equation} 
the second being induced by $d$ \eqref{p1-tphdi5a}. In particular, we have a canonical isomorphism \eqref{p1-tphdi4c}
\begin{equation}\label{p1-tphdi7e}
\Gr_\rW^i\tmK^\bullet(N\otimes_{\co_{X'}}\cC^\dagger) \stackrel{\sim}{\rightarrow} \wedge^if^*(\cE)\otimes_{\co_{X'}} \mK^\bullet(N\otimes_{\co_{X'}}\cC^\dagger)[-i].
\end{equation}
We denote by
\begin{equation}\label{p1-tphdi7f}
\tpartial \colon \mK^\bullet(N\otimes_{\co_{X'}}\cC^\dagger)\rightarrow f^*(\cE)\otimes_{\co_{X'}}\mK^\bullet(N\otimes_{\co_{X'}}\cC^\dagger)
\end{equation}
the morphism of $\bD^+(\bMod(\co_{X'}))$, boundary of the Koszul filtration on $\tmK^\bullet(N\otimes_{\co_{X'}}\cC^\dagger)$, see \ref{p1-tphdi4}. 

\begin{prop}\label{p1-tphdi8}
Under the assumptions of \ref{p1-tphdi6} and \ref{p1-tphdi7} and with the same notation, 
the morphisms $\fd$ \eqref{p1-tphdi6b} and $\tpartial$ \eqref{p1-tphdi7f} are equal in $\bD^+(\bMod(\co_{X'}))$. 
\end{prop}

For simplicity, we set $\mK^\bullet=\mK^\bullet(N\otimes_{\co_{X'}}\cC^\dagger)$ \eqref{p1-tphdi6a}, 
$\tmK^\bullet=\tmK^\bullet(N\otimes_{\co_{X'}}\cC^\dagger)$ \eqref{p1-tphdi7a} and $\mG^\bullet=\rW^0/\rW^2(\tmK^\bullet)$ \eqref{p1-tphdi7b}.
We (abusively) denote the differentials of $\mG^\bullet$ also by $\tvartheta^\bullet$ which does not induce any ambiguity.
For every integer $n\geq 0$, we have 
\begin{equation}\label{p1-tphdi8a}
\mG^n=\mK^n\oplus  (f^*(\cE)\otimes_{\co_{X'}}\mK^{n-1}),
\end{equation}
the differential $\tvartheta^n\colon \mG^n\rightarrow \mG^{n+1}$ being given by the matrix 
\begin{equation}\label{p1-tphdi8b}
\begin{pmatrix}
\vartheta^n &0 \\
d & -\id\otimes\vartheta^{n-1}
\end{pmatrix},
\end{equation}
where we denoted abusively by $d\colon \mK^n\rightarrow f^*(\cE)\otimes_{\co_{X'}}\mK^n$ the morphism induced by $d$ \eqref{p1-tphdi5a}. 

We have a canonical exact sequence of complexes 
\begin{equation}\label{p1-tphdi8c}
0\longrightarrow f^*(\cE)\otimes_{\co_{X'}}\mK^\bullet[-1]\stackrel{a}{\longrightarrow} \mG^\bullet \longrightarrow \mK^\bullet \longrightarrow 0,
\end{equation}
whose associated boundary in $\bD^+(\bMod(\co_{X'}))$ is none other that $\tpartial$ \eqref{p1-tphdi7f}. 
We denote by $\rC^\bullet$ the mapping cone of the morphism $a$ and by $c^\bullet$ its differentials.
For every integer $n$, we therefore have 
\begin{equation}\label{p1-tphdi8d}
\rC^n=(f^*(\cE)\otimes_{\co_{X'}}\mK^n)\oplus \mG^n=(f^*(\cE)\otimes_{\co_{X'}}\mK^n)\oplus \mK^n\oplus  (f^*(\cE)\otimes_{\co_{X'}}\mK^{n-1}),
\end{equation}
the differential $c^n\colon C^n\rightarrow C^{n+1}$ being given by the matrix  
\begin{equation}\label{p1-tphdi8e}
\begin{pmatrix}
\id\otimes \vartheta^n &0 & 0\\
0 &\vartheta^n &0\\
\id& d & -\id\otimes\vartheta^{n-1}
\end{pmatrix}.
\end{equation}
We denote by $\pi_1$ (resp.\ $\pi_2$) the projection of $\rC^\bullet$ onto
$f^*(\cE)\otimes_{\co_{X'}}\mK^\bullet$ (resp.\ $\mK^\bullet$) \eqref{p1-tphdi8d}. 
Then, $\pi_2$ is a quasi-isomorphism, and $-\tpartial$ \eqref{p1-tphdi7f} is the composition in $\bD^+(\bMod(\co_{X'}))$ of the inverse of 
$\pi_2$ and $\pi_1$ (\cite{sp} \href{https://stacks.math.columbia.edu/tag/09KF}{09KF}).

We immediately see that the morphism
\begin{equation}
\pi_1+\fd\circ \pi_2\colon \rC^\bullet\rightarrow f^*(\cE) \otimes_{\co_{X'}}\mK^\bullet,
\end{equation}
where $\fd$ is the morphism \eqref{p1-tphdi6b}, is homotopic to zero by the homotopy defined by the projection on the third factor \eqref{p1-tphdi8d}
\begin{equation}
\rC^n \rightarrow f^*(\cE)\otimes_{\co_{X'}}\mK^{n-1},
\end{equation}
which implies the proposition.

\subsection{}\label{p1-tphdi9}
We take again the assumptions and notation of \ref{p1-tphdi6} and \ref{p1-tphdi4} and let $\rho\colon \cF\rightarrow \co_{X'}$ be a splitting 
of the extension \eqref{p1-tphdi1a}. By \ref{p1-thbn16}, the latter induces a homomorphism of $\co_{X'}$-algebras 
\begin{equation}\label{p1-tphdi9a}
\varrho^\dagger\colon \cC^\dagger\rightarrow \co_{X'}.
\end{equation}
We immediately see that the diagram 
\begin{equation}\label{p1-tphdi9f}
\xymatrix{
{N\otimes_{\co_{X'}}\cC^\dagger}\ar[r]^-(0.5){\vartheta}\ar[d]_{\id\otimes\varrho^\dagger}&
{\cE'\otimes_{\co_{X'}}N\otimes_{\co_{X'}}\cC^\dagger}\ar[d]^{u'\otimes\id\otimes \varrho^\dagger}\\
{N}\ar[r]^-(0.5){\utheta}&{\ucE'\otimes_{\co_{X'}}N,}}
\end{equation}
where $u'$ is the morphism defined in \eqref{p1-tphdi4a} and $\utheta$ is the Higgs field induced by $\theta$, is commutative. 
We deduce a morphism of Dolbeault complexes 
\begin{equation}\label{p1-tphdi9g}
\upbeta\colon \mK^\bullet(N\otimes_{\co_{X'}}\cC^\dagger) \rightarrow \umK^\bullet(N). 
\end{equation}

With the notation of \ref{p1-tphdi5}, we have $(\gamma'\otimes \id)\circ \td=0$, 
where $\gamma'$ and $\td$ are defined in \eqref{p1-tphdi5f} and \eqref{p1-tphdi5g}. 
Hence, the diagram
\begin{equation}\label{p1-tphdi9c}
\xymatrix{
{N\otimes_{\co_{X'}}\cC^\dagger}\ar[r]^-(0.5){\tvartheta}\ar[d]_{\id\otimes\varrho^\dagger}&
{\tcE\otimes_{\co_{X'}}N\otimes_{\co_{X'}}\cC^\dagger}\ar[d]^{\gamma'\otimes\id\otimes \varrho^\dagger}\\
{N}\ar[r]^-(0.5){\theta}&{\cE'\otimes_{\co_{X'}}N}}
\end{equation}
is commutative. We deduce a morphism of Dolbeault complexes 
\begin{equation}\label{p1-tphdi9d}
\upalpha\colon \tmK^\bullet(N\otimes_{\co_{X'}}\cC^\dagger) \rightarrow \mK^\bullet(N).
\end{equation}
By \eqref{p1-tphdi5j}, for every $i\geq 0$, we have 
\begin{equation}\label{p1-tphdi9e}
\upalpha (\rW^i \tmK^\bullet(N\otimes_{\co_{X'}}\cC^\dagger))\subset \rW^i\mK^\bullet(N),
\end{equation}
where $\rW^i$ denotes the Koszul filtrations defined in \eqref{p1-tphdi4b} and \eqref{p1-tphdi7g}. 
Moreover, since $\gamma'\circ \iota=-u$, the diagram
\begin{equation}\label{p1-tphdi9h}
\xymatrix{
{\Gr_\rW^i \tmK^\bullet(N\otimes_{\co_{X'}}\cC^\dagger)}\ar[d]\ar[r]^-(0.5){(-1)^i\Gr_\rW^i(\upalpha)}&{\Gr_\rW^i\mK^\bullet(N)}\ar[d]\\
{\wedge^if^*(\cE)\otimes_{\co_{X'}}\mK^\bullet(N\otimes_{\co_{X'}}\cC^\dagger)[-i]}\ar[r]^-(0.4){\id\otimes \upbeta}&{\wedge^if^*(\cE)\otimes_{\co_{X'}}\umK^\bullet(N)[-i],}}
\end{equation}
where the vertical arrows are the isomorphisms \eqref{p1-tphdi4c} and \eqref{p1-tphdi7e}, is commutative. 

The diagram
\begin{equation}\label{p1-tphdi9i}
{\scriptsize
\xymatrix{
0\ar[r]&{f^*(\cE)\otimes_{\co_{X'}}\mK^\bullet(N\otimes_{\co_{X'}}\cC^\dagger)[-1]}\ar[r]\ar[d]_{-\id\otimes \upbeta}&
{\rW^0/\rW^2(\tmK^\bullet(N\otimes_{\co_{X'}}\cC^\dagger))}\ar[r]\ar[d]^{\upalpha}&
{\mK^\bullet(N\otimes_{\co_{X'}}\cC^\dagger)}\ar[r]\ar[d]^{\upbeta}&0\\
0\ar[r]&{f^*(\cE)\otimes_{\co_{X'}}\umK^\bullet(N)[-1]}\ar[r]&{\rW^0/\rW^2(\mK^\bullet(N\otimes_{\co_{X'}}\cC^\dagger))}\ar[r]&
{\umK^\bullet(N)}\ar[r]&0,}}
\end{equation}
where the horizontal lines are the canonical exact sequences, is commutative. We deduce the following.

\begin{prop}\label{p1-tphdi10}
Under the assumptions of \ref{p1-tphdi9} and with the same notation, 
for every $q\geq 0$, the morphism $\upbeta$ \eqref{p1-tphdi9g} induces a morphism of Higgs $\co_{X}$-modules
\begin{equation}\label{p1-tphdi10a}
\rR^qf^\uptau_*(N,\theta) \rightarrow (\rR^qf_*(\umK^\bullet(N)),\kappa^q), 
\end{equation}
where the source is defined in \eqref{p1-tphdi6d} and the target is  
equipped with the Katz-Oda field \eqref{p1-tphdi4h}.
\end{prop}

\begin{rema}
We give in \ref{p2-cmupiso29} an important case where \eqref{p1-tphdi10a} is an isomorphism. 
\end{rema}

\section{Small Higgs modules}\label{p1-tshbn}

\subsection{}\label{p1-tshbn1}
In this section, we adopt the assumptions and notation of §\ref{p1-thbn},
moreover, let $\fX$ be an $\cS$-formal scheme locally of finite presentation \eqref{p1-thbn0} (\cite{egr1} 2.3.15); 
so it is an idyllic formal scheme (\cite{egr1} 2.6.13).  
We suppose that $X$ \eqref{p1-thbn1} is the Zariski topos of $\fX$ and that $\co_X$ is an $\co_\fX$-algebra. 
Let $\cR$ be an $\co_\fX$-algebra, which is topologically of finite presentation \eqref{p1-pfs1} and $V$-flat, i.e., $p$ is not a zero-divisor in $\cR$. 
We suppose given a homomorphism of $\co_\fX$-algebras $\cR\rightarrow \co_X$ and a {\em locally free $\cR$-module of finite type} $\Omega$ 
such that $\cE=\Omega\otimes_{\cR}\co_X$. 

For every rational number $r\leq 0$, we consider the $\co_X$-module $\cF^{(r)}$ \eqref{p1-thbn3a} as an extension 
\begin{equation}\label{p1-tshbn1a}
0\rightarrow \co_X\rightarrow \cF^{(r)}\rightarrow p^r\Omega\otimes_\cR\co_X \rightarrow 0, 
\end{equation} 
so that for all rational number $r'$ with $0\leq r'\leq r$, the morphism $\tta^{r,r'}$ \eqref{p1-thbn4a}
is compatible with the inclusion $p^r\Omega\otimes_\cR\co_X\rightarrow p^{r'}\Omega\otimes_\cR\co_X$.

\begin{prop}\label{p1-tshbn2}
Let $\cM$ be a locally free $\cR$-module of finite type, $\cA=\rS_{\cR}(\cM)$, $\hcA$ its $p$-adic completion. Then, 
\begin{itemize}
\item[{\rm (i)}] The $V$-algebra $\hcA$ is flat, and for every integer $n\geq 1$, the canonical morphism 
\begin{equation}\label{p1-tshbn2a}
\cA/p^n\cA\rightarrow \hcA/p^n\hcA
\end{equation}
is an isomorphism. 
\item[{\rm (ii)}] The $\co_{\fX}$-algebra $\hcA$ is topologically of finite presentation.
\item[{\rm (iii)}] The sheaves of rings $\hcA$ and $\hcA[\frac 1 p]$ are coherent. 
\item[{\rm (iv)}] The $\cR$-algebra $\hcA$ is flat. 
\item[{\rm (v)}] Any coherent $\hcA$-module is $p$-adic, i.e., $p$-adically complete and Hausdorff. 
\end{itemize}
\end{prop}

(i) Indeed, since $\cR$ is $V$-flat, $\cA$ is $V$-flat. It is enough to prove that the canonical exact sequence 
\begin{equation}\label{p1-tshbn2b}
0\longrightarrow \cA \stackrel{p^n}{\longrightarrow} \cA \longrightarrow \cA/p^n\cA\longrightarrow 0
\end{equation}
induces by $p$-adic completion an exact sequence
\begin{equation}\label{p1-tshbn2c}
0\longrightarrow \hcA \stackrel{p^n}{\longrightarrow} \hcA \longrightarrow \cA/p^n\cA\longrightarrow 0.
\end{equation}
Let $U$ be an affine formal open subscheme of $\fX$. It follows from \ref{p1-pfs9} and \ref{p1-pfs7}(i) that for any integer $i\geq 0$, the canonical morphism 
\begin{equation}\label{p1-tshbn2d}
\Gamma(U,\rS^i_{\cR}(\cM))/p^n \Gamma(U,\rS^i_{\cR}(\cM)) \rightarrow \Gamma(U,\rS^i_{\cR}(\cM/p^n\cM)) 
\end{equation}
is an isomorphism. By (\cite{egr1} 2.6.10), we deduce that the canonical morphism 
\begin{equation}\label{p1-tshbn2e}
\Gamma(U,\cA)/p^n \Gamma(U,\cA)  \rightarrow \Gamma(U,\cA/p^n\cA) 
\end{equation}
is an isomorphism. The ring $\Gamma(U,\hcA)$ is therefore the $p$-adic completion of the ring 
$\Gamma(U,\cA)$. Hence by (\cite{ac} chap.~III §2.11 prop.~14 and cor.~1), we have an isomorphism
\begin{equation}\label{p1-tshbn2f}
\Gamma(U,\hcA)/p^n \Gamma(U,\hcA) \stackrel{\sim}{\rightarrow} \Gamma(U,\cA)/p^n \Gamma(U,\cA). 
\end{equation}
Since $\Gamma(U,\cA)$ is $V$-flat, we deduce that $\Gamma(U,\hcA)$ is also $V$-flat. 
Therefore, the evaluation of the sequence \eqref{p1-tshbn2c} at $U$ is exact, which implies the required assertion. 

(ii) This follows from \eqref{p1-tshbn2a}.

(iii) This follows from (ii), \ref{p1-pfs8} and \ref{p1-pfs12}(i). 

(iv) This follows from (i)--(iii) and  \ref{p1-pfs14}.

(v) This follows from (ii) and \ref{p1-pfs10}.  

\subsection{}\label{p1-tshbn3}
Recall that $\Omega$ is a locally free $\cR$-module of finite type \eqref{p1-tshbn1}.
For any rational number $r\geq 0$, we denote by $\cG^{(r)}$ the $\cR$-algebra
\begin{equation}\label{p1-tshbn3a}
\cG^{(r)}=\rS_{\cR}(p^r\Omega),
\end{equation}
and by $\hcG^{(r)}$ its $p$-adic completion. 
For an affine open formal subscheme $U$ of $\fX$, $\hcG^{(r)}(U)$ is canonically isomorphic to the 
$\cR(U)$-algebra of series $(s_i)_{i\in \mN}\in \prod_{i\in \mN}\rS^i_{\cR}(p^r\Omega)(U)$
such that for every integer $m\geq 0$, there exists $n\in \mN$ such that for every $i\geq n$, $s_i\in p^m  \rS^i_{\cR}(p^r\Omega)(U)$. 

We denote by $\delta_{\cG^{(r)}}$ the $\cR$-derivation of $\cG^{(r)}$ composed of 
\begin{equation}\label{p1-tshbn3c}
\delta_{\cG^{(r)}}\colon \cG^{(r)}\rightarrow p^r\Omega\otimes_{\cR}\cG^{(r)}\rightarrow \Omega\otimes_{\cR}\cG^{(r)},
\end{equation}
where the first map is the universal $\cR$-derivation of $\cG^{(r)}$ and the second map is the canonical injection, and by 
\begin{equation}\label{p1-tshbn3e}
\delta_{\hcG^{(r)}}\colon \hcG^{(r)} \rightarrow \Omega\otimes_{\cR} \hcG^{(r)}
\end{equation}
its extension to the $p$-adic completions. 
Observe that $\Omega\otimes_{\cR}\hcG^{(r)}$ is the $p$-adic completion of $\Omega\otimes_{\cR}\cG^{(r)}$ \eqref{p1-thbn18}.
Then $\delta_{\cG^{(r)}}$ and $\delta_{\hcG^{(r)}}$ are Higgs $\cR$-fields. 

For all rational numbers $r\geq r'\geq 0$, we have an injective homomorphism
$a^{r,r'}\colon \cG^{(r)}\rightarrow \cG^{(r')}$ induced by the canonical injection $p^r\Omega\rightarrow p^{r'}\Omega$. We immediately check that the induced homomorphism
$\ha^{r,r'}\colon\hcG^{(r)}\rightarrow \hcG^{(r')}$ is injective. We have
\begin{equation}\label{p1-tshbn3f}
(\id \otimes a^{r,r'}) \circ \delta_{\cG^{(r)}}=\delta_{\cG^{(r')}}\circ a^{r,r'}.
\end{equation}

We define the $\cR$-algebra
\begin{equation}\label{p1-tshbn3g}
\hcG^{(r+)}=\underset{\underset{t\in \mQ_{>r}}{\longrightarrow}}{\lim}\ \hcG^{(t)}. 
\end{equation}
By \eqref{p1-tshbn3f}, the derivations $\delta_{\hcG^{(t)}}$ \eqref{p1-tshbn3e} induce an $\cR$-derivation
\begin{equation}\label{p1-tshbn3h}
\delta_{\hcG^{(r+)}}\colon \hcG^{(r+)}\rightarrow \Omega\otimes_{\cR}\hcG^{(r+)},
\end{equation}
which is also a Higgs $\cR$-field.
For simplicity, we set $\cG^\dagger=\hcG^{(0+)}$ and $\delta_{\cG^\dagger}=\delta_{\hcG^{(0+)}}$. 

\subsection{}\label{p1-tshbn4}
We denote by $\fm$ the maximal ideal of $V$ \eqref{p1-thbn0} and by $v$ its valuation normalized by $v(p)=1$. 
Recall that on $V$, the logarithmic function $\log(x)$ converges when $x\in 1+\fm$, and the exponential function $\exp(x)$ converges when $v(x)>\frac{1}{p-1}$. 
For every $x\in V$ such that $v(x)>\frac{1}{p-1}$, we have  
\begin{eqnarray*}
\exp(x)&\equiv& 1+x \mod(x\fm),\\
\log(1+x)&\equiv& x \mod(x\fm),
\end{eqnarray*}
$\exp(\log(1+x))=1+x$ and $\log(\exp(x))=x$.

\begin{defi}\label{p1-tshbn5}
Let $\cN$ be an $\cR$-module, 
$u$ an $\cR$-endomorphism of $\cN$, $\varepsilon$ a rational number $>0$.
\begin{itemize}
\item[{\rm (i)}] We say that $u$ is {\em $\varepsilon$-small} if $u(\cN)\subset p^\varepsilon \cN$. 
\item[{\rm (ii)}] We say that $u$ is {\em small} if it is $\varepsilon'$-small for a rational number $\varepsilon'>\frac{1}{p-1}$. 
\end{itemize}
\end{defi}

\subsection{}\label{p1-tshbn6}
Let $\cN$ be a coherent $\cR$-module which is $V$-flat, $u$ a small $\cR$-endomorphism of $\cN$.
The $\cR$-module $\cEnd_{\cR}(\cN)$ is $p$-adic by \ref{p1-pfs10}. It is also $V$-flat.
We identify $u$ with with an element of $p^\varepsilon \cEnd_{\cR}(\cN)$ for a rational number $\varepsilon>\frac{1}{p-1}$. 
Then, the series 
\begin{equation}\label{p1-tshbn6a}
\exp(u)= \sum_{n\geq 0} \frac{1}{n!}u^n
\end{equation}
defines an $\cR$-automorphism of $\cN$, with inverse $\exp(-u)$.  

\begin{defi}\label{p1-tshbn7}
Let $N$ be an $\cR[\frac 1 p]$-module, $u$ an $\cR$-endomorphism of $N$. 
\begin{itemize}
\item[{\rm (i)}] 
We say that $u$ is {\em CL-small} if $N$ admits a {\em coherent $\cR$-lattice}, stable by $u$ such that the induced endomorphism is small, i.e., 
if there exist a coherent sub-$\cR$-module $\cN$ of $N$
which generates it over $\cR[\frac 1 p]$ and a rational number $\varepsilon>\frac{1}{p-1}$ such that $u(\cN)\subset p^\varepsilon \cN$. 
\item[{\rm (ii)}] We say that $u$ is {\em locally CL-small} if there exists a Zariski open covering $(U_i)_{i\in I}$ of $\fX$
such that for every $i\in I$, $u|U_i$ is CL-small.
\end{itemize}
\end{defi}

By \ref{p1-pfs12}(ii), if $\fX$ is quasi-compact, any coherent $\cR[\frac 1 p]$-module $N$ admits a coherent lattice, i.e., there exists a coherent sub-$\cR$-module $\cN$ of $N$
which generates it over $\cR[\frac 1 p]$. 

\subsection{}\label{p1-tshbn8}
Let $N$ be an $\cR[\frac 1 p]$-module, $u$ a CL-small $\cR$-endomorphism of $N$, $\cN$ a coherent $\cR$-lattice of $N$, 
stable by $u$ such that the induced endomorphism $v$ of $\cN$ is small. 
Then, the $\cR[\frac 1 p]$-automorphism $\exp(u)\colon N\rightarrow N$ induced by the $\cR$-automorphism $\exp(v)\colon \cN\rightarrow \cN$, 
does not depend on the choice of $\cN$. Indeed, for every quasi-compact open formal subscheme $U$ of $\fX$ and every $z\in \Gamma(U,N)$, 
there exists $m\geq 0$ such that $p^mz\in \Gamma(U,\cN)$. We have
\begin{equation}\label{p1-tshbn8a}
\exp(u)(z)=\frac{1}{p^m}  \sum_{n\geq 0} \frac{1}{n!}u^n(p^mz).
\end{equation}

\subsection{}\label{p1-tshbn80}
Let $N$ be an $\cR[\frac 1 p]$-module, $u$ a locally CL-small $\cR$-endomorphism of $N$. By gluing the exponentials defined locally in \eqref{p1-tshbn8a}, 
we define a canonical $\cR[\frac 1 p]$-automorphism of $N$ that we denote by $\exp(u)\colon N\rightarrow N$.

\begin{lem}\label{p1-tshbn81}
Let $N$ be an $\cR[\frac 1 p]$-module, $u$ a (resp.\ locally) CL-small $\cR$-endomorphism of $N$, $M$ a sub-$\cR[\frac 1 p]$-module of finite type of $N$
stable by $u$. Then, the restriction of $u$ to $M$ is (resp.\ locally) CL-small. 
\end{lem}

Indeed, we may assume that $N$ is CL-small. Let $\cN$ a coherent $\cR$-lattice of $N$, 
stable by $u$ such that the induced endomorphism of $\cN$ is small. We set $\cM=M\cap \cN$. 
Since the quotient $\cN/\cM$ is $p$-torsion free, the $\cR$-module $\cM$ is of finite type by (\cite{egr1} 1.9.14). 
It is therefore a coherent $cR$-lattice of $M$ and it is obviously stable by $u$.

\begin{defi}\label{p1-tshbn9}
Let $(\cN,\theta)$ be a Higgs $\cR$-module with coefficients in $\Omega$ \eqref{p1-delta-con6}, $\varepsilon$ a rational number~$> 0$.
\begin{itemize}
\item[(i)] We say that $(\cN,\theta)$ is {\em $\varepsilon$-small} if 
\begin{equation}\label{p1-tshbn9a}
\theta(\cN)\subset p^\varepsilon (\Omega\otimes_{\cR} \cN).
\end{equation} 
We say then that the Higgs $\cR$-field $\theta$ is {\em $\varepsilon$-small}. 
\item[{\rm (ii)}] We say that $(\cN,\theta)$ is {\em small} if it is $\varepsilon'$-small for a rational number $\varepsilon'>\frac{1}{p-1}$. 
We say then that the Higgs $\cR$-field $\theta$ is {\em small}. 
\end{itemize}
\end{defi}

\subsection{}\label{p1-tshbn10}
Let $r$, $\varepsilon$ be two rational numbers such that $r\geq 0$ and $\varepsilon>r+\frac{1}{p-1}$, 
$(\cN,\theta)$ be an $\varepsilon$-small Higgs $\cR$-module with coefficients in $\Omega$ such that the $\cR$-module $\cN$ is coherent and $V$-flat.  
Consider the composed $\cR$-linear morphism 
\begin{equation}\label{p1-tshbn10a}
\cN\rightarrow p^r\Omega\otimes_{\cR}\cN\rightarrow \cG^{(r)}\otimes_{\cR}\cN,
\end{equation}
where the first morphism is induced by $\theta$ and the second one by the canonical inclusion $p^r\Omega\rightarrow \cG^{(r)}$. 
It induces, by extension of scalars, a $\hcG^{(r)}$-linear morphism
\begin{equation}\label{p1-tshbn10b}
\vartheta^{(r)}\colon \hcG^{(r)}\otimes_{\cR}\cN \rightarrow  \hcG^{(r)}\otimes_{\cR}\cN.
\end{equation}
By \ref{p1-tshbn2}, the $\co_\fX$-algebra $\hcG^{(r)}$ is topologically of finite presentation, and the $\hcG^{(r)}$-module $\hcG^{(r)}\otimes_{\cR}\cN$ is coherent.
Moreover, the $\cR$-algebra $\hcG^{(r)}$ being flat by \ref{p1-tshbn2}(iv),
the $\hcG^{(r)}$-module $\hcG^{(r)}\otimes_{\cR}\cN$ is $V$-flat. 
Since the endomorphism $\vartheta^{(r)}$ is $(\varepsilon-r)$-small, we can define the $\hcG^{(r)}$-linear automorphism
\begin{equation}\label{p1-tshbn10c}
\exp(\vartheta^{(r)})=\sum_{n\geq 0}\frac{1}{n!}(\vartheta^{(r)})^n\colon \hcG^{(r)}\otimes_{\cR}\cN\rightarrow \hcG^{(r)}\otimes_{\cR}\cN,
\end{equation}
with inverse $\exp(-\vartheta^{(r)})$. 
For every rational number $r'$ such that $0\leq r'\leq r$, the diagram
\begin{equation}\label{p1-tshbn10d}
\xymatrix{
{\hcG^{(r)}\otimes_{\cR}\cN}\ar[rr]^-(0.5){\exp(\vartheta^{(r)})}\ar[d]_-(0.4){\ha^{r,r'}\otimes \id}&&
{\hcG^{(r)}\otimes_{\cR}\cN}\ar[d]^-(0.4){\ha^{r,r'}\otimes \id}\\
{\hcG^{(r')}\otimes_{\cR}\cN}\ar[rr]^-(0.5){\exp(\vartheta^{(r')})}&&{\hcG^{(r')}\otimes_{\cR}\cN}}
\end{equation}
is commutative. We can therefore denote $\exp(\vartheta^{(r)})$ simply by $\exp(\theta)$ without risk of ambiguity.

\begin{rema}\label{p1-tshbn11}
We keep the assumptions and notation of \ref{p1-tshbn10} and suppose moreover that the $\cR$-module $\Omega$ is free 
with basis $e_1,\dots,e_d$.    
We can uniquely write 
\begin{equation}\label{p1-tshbn11a}
\theta=\sum_{i=1}^d\theta_i \otimes e_i,
\end{equation}
where the $\theta_i$ are endomorphisms of $\cN$ belonging to $p^{\varepsilon}\End_{\cR}(\cN)$ and commuting to each other.
For any $\un=(n_1,\dots,n_d)\in \mN^d$, we set $|\un|=\sum_{i=1}^d n_i$, $\un!=\prod_{i=1}^d n_i!$, $\utheta^{\un}=\prod_{i =1}^d\theta_i^{n_i} \in \End_{\cR}(\cN)$
and $\ue^{\un}=\prod_{i=1}^de_i^{n_i}\in \cG$. Then, we have
\begin{equation}\label{p1-tshbn11b}
\exp(\theta) = \sum_{\un\in \mN^d} \frac{1}{\un!}\utheta^{\un}\otimes\ue^{\un}.
\end{equation}
\end{rema}

\begin{prop}\label{p1-tshbn12}
Under the assumptions of \ref{p1-tshbn10}, the $\hcG^{(r)}$-linear automorphism $\exp(\vartheta^{(r)})$ \eqref{p1-tshbn10c}, that we denoted simply by
\begin{equation}\label{p1-tshbn12a}
\exp(\theta)\colon \hcG^{(r)}\otimes_{\cR}\cN\rightarrow \hcG^{(r)}\otimes_{\cR}\cN,
\end{equation}
is underlying an isomorphism of $\hcG^{(r)}$-modules with $\delta_{\hcG^{(r)}}$-connection, 
where the $\delta_{\hcG^{(r)}}$-connections are defined as in \ref{p1-delta-con4}, the module $\cN$ of the source (resp.\ target) being endowed with the Higgs field $\theta$ (resp.\ $0$). 
\end{prop}

Indeed, it immediately follows from the definition \eqref{p1-tshbn10c} that the diagram
\begin{equation}
\xymatrix{
{\hcG^{(r)}\otimes_{\cR}\cN}\ar[rr]^-(0.5){\exp(\theta)}
\ar[d]_{\nabla^{(r)}}&&{\hcG^{(r)}\otimes_{\cR}\cN}\ar[d]^{\delta_{\hcG^{(r)}}\otimes\id}\\
{\Omega\otimes_{\cR}\hcG^{(r)}\otimes_{\cR}\cN}\ar[rr]^-(0.5){\id \otimes \exp(\theta)}&&
{\Omega\otimes_{\cR}\hcG^{(r)}\otimes_{\cR}\cN,}}
\end{equation}
where $\nabla^{(r)}= \delta_{\hcG^{(r)}}\otimes \id+\id\otimes \theta$, 
is commutative modulo $p^n$ for all $n\geq 1$. It is therefore commutative, which proves the proposition.

\begin{defi}\label{p1-tshbn13}
Let $(N,\theta)$ be a Higgs $\cR[\frac 1 p]$-module with coefficients in $\Omega$.
\begin{itemize}
\item[(i)] We say that $(N,\theta)$ is {\em CL-small} if $N$ admits a coherent $\cR$-lattice, stable by $\theta$ such that the induced Higgs field is small, 
i.e., if there exists a coherent sub-$\cR$-module $\cN$ of $N$
which generates it over $\cR[\frac 1 p]$ and a rational number $\varepsilon>\frac{1}{p-1}$ such that
\begin{equation}\label{p1-tshbn13a}
\theta(\cN)\subset p^\varepsilon \Omega\otimes_{\cR} \cN.
\end{equation}
\item[(ii)] We say that $(N,\theta)$ is {\em locally CL-small} if there exists a Zariski open covering $(U_i)_{i\in I}$ of $\fX$
such that for every $i\in I$, $(N|U_i,\theta|U_i)$ is CL-small.
\end{itemize}
\end{defi}

\begin{rema}
(Locally) CL-small Higgs modules were called (locally) small in (\cite{ag2} 4.2.11); see also \ref{p2-rlps25}. 
\end{rema}

\begin{lem}\label{p1-tshbn130}
Let  $(N,\theta)$ be a (resp.\ locally) CL-small Higgs $\cR[\frac 1 p]$-module with coefficients in $\Omega$, 
$M$ a sub-$\cR[\frac 1 p]$-module of finite type of $N$ such that $\theta(M)\subset \Omega\otimes_\cR M$. 
We denote by $\theta_M\colon M\rightarrow \Omega\otimes_\cR M$ the Higgs $\cR[\frac 1 p]$-field induced by $\theta$. 
Then, the Higgs $\cR[\frac 1 p]$-module $(M,\theta_M)$ is (resp.\ locally) CL-small. 
\end{lem}

The proof is similar to that of \ref{p1-tshbn81}.

\subsection{}\label{p1-tshbn14}
Let $(N,\theta)$ be a CL-small Higgs $\cR[\frac 1 p]$-module with coefficients in $\Omega$,
$r$, $\varepsilon$ two rational numbers such that $r\geq 0$ and $\varepsilon>r+\frac{1}{p-1}$, 
$\cN$ a coherent $\cR$-lattice of $N$ such that $\theta(\cN)\subset p^\varepsilon \Omega\otimes_{\cR} \cN$. 
By \ref{p1-tshbn2}, $\hcG^{(r)}\otimes_{\cR} \cN$ is a coherent $\hcG^{(r)}$-lattice of $\hcG^{(r)}\otimes_{\cR} N$. 
We denote by $\kappa$ the Higgs $\cR$-field induced by $\theta$ on $\cN$. 
Then, the $\hcG^{(r)}[\frac 1 p]$-automorphism 
\begin{equation}\label{p1-tshbn14a}
\exp(\theta)\colon \hcG^{(r)}\otimes_{\cR} N\rightarrow \hcG^{(r)}\otimes_{\cR} N,
\end{equation}
induced by the $\hcG^{(r)}$-automorphism $\exp(\kappa)$ of $\hcG^{(r)}\otimes_{\cR} \cN$ \eqref{p1-tshbn10}, does not depend on the choice of $\cN$ \eqref{p1-tshbn8}. 
By \ref{p1-tshbn12}, it is an isomorphism of $\hcG^{(r)}$-modules with $\delta_{\hcG^{(r)}}$-connection, 
where the $\delta_{\hcG^{(r)}}$-connections are defined as in \ref{p1-delta-con4}, the module $N$ of the source (resp.\ target) being endowed with the Higgs field $\theta$ (resp.\ $0$).

\begin{prop}\label{p1-tshbn15}
Suppose that the extension \eqref{p1-thbn1a} splits. Let $(N,\theta)$ be a CL-small Higgs 
$\cR[\frac 1 p]$-module with coefficients in $\Omega$. 
Set $(\tN,\ttheta)=(N\otimes_{\cR}\co_X,\theta\otimes \id)$, which is a Higgs $\co_X[\frac 1 p]$-module with coefficients in 
$\cE=\Omega\otimes_\cR\co_X$.
Then, $(\tN,\ttheta)$ is twistable by the extension \eqref{p1-thbn1a}, in the sense of \ref{p1-thbn14}. 
\end{prop}

Indeed, by \ref{p1-tshbn14}, there exists a rational number $r>0$ such that the automorphism \eqref{p1-tshbn14a}
\begin{equation}\label{p1-tshbn15b}
\exp(\theta)\colon \hcG^{(r)}\otimes_{\cR} N\stackrel{\sim}{\rightarrow} \hcG^{(r)}\otimes_{\cR} N,
\end{equation}
is an isomorphism of $\hcG^{(r)}$-modules with $\delta_{\hcG^{(r)}}$-connection,
where the $\delta_{\hcG^{(r)}}$-connections are defined as in \ref{p1-delta-con4}, the module $N$ of the source (resp.\ target) being endowed with the Higgs field $\theta$ (resp.\ $0$).

By \ref{p1-thbn16}, a splitting of the extension \eqref{p1-thbn1a} induces an isomorphism of $\co_X$-algebras
\begin{equation}
\cG^{(r)}\otimes_{\cR}\co_X\stackrel{\sim}{\rightarrow} \cC^{(r)},
\end{equation}
compatible with the derivations $\delta_{\cG^{(r)}}$ \eqref{p1-tshbn3c} and $\delta_{\cC^{(r)}}$ \eqref{p1-thbn3d}.
Hence, we deduce from \eqref{p1-tshbn15b}, by extension of scalars from $\hcG^{(r)}$ to $\hcC^{(r)}$,
that the Higgs $\co_X[\frac 1 p]$-module $(\tN,\ttheta)$ is twistable by the extension \eqref{p1-thbn1a}.

\begin{cor}\label{p1-tshbn16}
For every locally CL-small Higgs $\cR[\frac 1 p]$-module $(N,\theta)$ with coefficients in $\Omega$, 
the Higgs $\co_X[\frac 1 p]$-module $(N\otimes_{\cR}\co_X,\theta\otimes \id)$ is weakly twistable by the extension \eqref{p1-thbn1a}, 
in the sense of \ref{p1-thbn30}. 
\end{cor}

It follows from \ref{p1-tshbn15} and \ref{p1-thbn21}. 

\subsection{}\label{p1-tshbn31}
Let $(N,\theta)$ be a locally CL-small Higgs $\cR[\frac 1 p]$-module with coefficients in $\Omega$ \eqref{p1-tshbn13}. 
We set $\cG=\cG^{(0)}$ and $\hcG=\hcG^{(0)}$ \eqref{p1-tshbn3a}, and denote by 
\begin{equation}\label{p1-tshbn31b}
\vartheta\colon \hcG\otimes_{\cR}N \rightarrow  \hcG \otimes_{\cR}N
\end{equation}
the $\hcG$-linear morphism induced by the composed $\cR$-linear morphism 
\begin{equation}\label{p1-tshbn31c}
N\stackrel{\theta}{\longrightarrow} \Omega\otimes_{\cR}N\longrightarrow \hcG\otimes_{\cR}N,
\end{equation}
where the second morphism is induced by the canonical inclusion $\Omega\rightarrow \cG$. 
By \ref{p1-tshbn2}, the $\co_\fX$-algebra $\hcG$ is topologically of finite presentation and is $\cR$-flat. 
Therefore, by \ref{p1-tshbn2}(v), the $\hcG$-linear morphism $\vartheta$ is locally CL-small \eqref{p1-tshbn7}. 
By \ref{p1-tshbn80}, we can then define the $\hcG[\frac 1 p]$-automorphism 
\begin{equation}\label{p1-tshbn31d}
\exp(\vartheta)\colon \hcG\otimes_{\cR}N \rightarrow  \hcG \otimes_{\cR}N.
\end{equation}
The integrability of $\theta$ implies that the diagram
\begin{equation}\label{p1-tshbn31e}
\xymatrix{
{\hcG\otimes_{\cR} N}\ar[rr]^-(0.5){\vartheta}
\ar[d]_{\id\otimes \theta}&&{\hcG\otimes_{\cR} N}\ar[d]^{\id\otimes \theta}\\
{\Omega\otimes_{\cR}\hcG\otimes_{\cR} N}\ar[rr]^-(0.5){\id \otimes\vartheta}&&
{\Omega\otimes_{\cR}\hcG\otimes_{\cR} N}}
\end{equation}
is commutative, see \ref{p1-tshbn11}. Hence, so is the diagram 
\begin{equation}\label{p1-tshbn31f}
\xymatrix{
{\hcG\otimes_{\cR} N}\ar[rr]^-(0.5){\exp(\vartheta)}
\ar[d]_{\id\otimes \theta}&&{\hcG\otimes_{\cR} N}\ar[d]^{\id\otimes \theta}\\
{\Omega\otimes_{\cR}\hcG\otimes_{\cR} N}\ar[rr]^-(0.5){\id \otimes\exp(\vartheta)}&&
{\Omega\otimes_{\cR}\hcG\otimes_{\cR} N.}}
\end{equation}
It follows that $\exp(\vartheta)$ is underlying an isomorphism of Higgs $\hcG[\frac 1 p]$-modules
\begin{equation}\label{p1-tshbn31i}
\exp(\vartheta)\colon (\hcG\otimes_{\cR}N, \id\otimes \theta) \rightarrow  (\hcG \otimes_{\cR}N,\id\otimes \theta).
\end{equation}

We set $(\tN,\ttheta)=(N\otimes_{\cR}\co_X,\theta\otimes \id)$, which is a Higgs $\co_X[\frac 1 p]$-module with coefficients in $\cE$.  
Let $\phi\colon \cE\rightarrow \co_X$ be an $\co_X$-linear form. Since $\co_X$ is $p$-adic, $\phi$ induces a homomorphism 
of $\cR$-algebras
\begin{equation}\label{p1-tshbn31j}
\varphi\colon \hcG\rightarrow \co_X. 
\end{equation}
We set 
\begin{equation}\label{p1-tshbn31k}
\ttheta_\phi=(\phi\otimes \id_\tN)\circ \ttheta \colon \tN\rightarrow \tN, 
\end{equation}
which is canonically identified with $\vartheta\otimes_{\hcG,\varphi}\co_X$. 
We cannot define ``$\exp(\ttheta_\phi)$'' in general, but we can replace it by the isomorphism of Higgs $\co_X[\frac 1 p]$-modules
\begin{equation}\label{p1-tshbn31a}
\exp_\phi(\theta)=\exp(\vartheta)\otimes_{\hcG,\varphi}\co_X \colon (\tN,\ttheta)\rightarrow (\tN,\ttheta),
\end{equation}
deduced from \eqref{p1-tshbn31i} by base change by $\varphi$. 

By \ref{p1-tshbn16}, the Higgs $\co_X[\frac 1 p]$-module $(\tN,\ttheta)$ is weakly twistable by the extension \eqref{p1-thbn1a}. 
Therefore, by \ref{p1-thbn31}, every splitting $\rho\colon \cF\rightarrow \co_X$ of the extension \eqref{p1-thbn1a} 
induces an isomorphism of Higgs $\co_X$-modules 
\begin{equation}\label{p1-tshbn31l}
\ttt_\rho\colon \uptau(\tN,\ttheta)\stackrel{\sim}{\rightarrow}(\tN,\ttheta),
\end{equation}
where $\uptau$ is the functor \eqref{p1-thbn9b}. 

\begin{prop}\label{p1-tshbn32}
Let $(N,\theta)$ be a locally CL-small Higgs $\cR[\frac 1 p]$-module with coefficients in $\Omega$, 
$(\tN,\ttheta)=(N,\theta)\otimes_{\cR}\co_X$, 
$\rho_i\colon \cF\rightarrow \co_X$, for $i=1,2$, two splittings of the extension \eqref{p1-thbn1a}.  
We denote by $\phi\colon \cE\rightarrow \co_X$ the $\co_X$-linear form induced by $\rho_1-\rho_2$. 
Then, the diagram
\begin{equation}\label{p1-tshbn32a}
\xymatrix{
{\uptau(\tN,\ttheta)}\ar[r]^-(0.5){\ttt_{\rho_1}}\ar[rd]_-(0.5){\ttt_{\rho_2}}&{(\tN,\ttheta)}\ar[d]^{\exp_\phi(\theta)}\\
&{(\tN,\ttheta),}}
\end{equation}
where $\ttt_{\rho_i}$ and $\exp_\phi(\theta)$ are defined in \eqref{p1-tshbn31l} and \eqref{p1-tshbn31a}, is commutative.   
\end{prop}

Indeed, the question being local, we may assume the Higgs $\cR[\frac 1 p]$-module $(N,\theta)$ CL-small.
For $i=1,2$, let $\sigma_i\colon \cE\rightarrow \cF$ be the splitting of the extension \eqref{p1-thbn1a} associated with $\rho_i$,
i.e., $\id_{\cF}=\sigma_i\circ \nu+c\circ \rho_i$. 
For every rational number $r\geq 0$, $\sigma_i$ and $\rho_i$ induce two associated splittings 
$\sigma^{(r)}_i\colon \cE\rightarrow \cF^{(r)}$ and $\rho^{(r)}_i\colon \cF^{(r)}\rightarrow \co_X$ of the extension \eqref{p1-thbn3a}. 
We immediately see that the diagram 
\begin{equation}\label{p1-tshbn32b}
\xymatrix{
\cE\ar[r]^{\sigma^{(r)}_2}\ar[rd]_{p^r\phi}&{\cF^{(r)}}\ar[d]^{\rho^{(r)}_1}\\
&{\co_X}}
\end{equation}
is commutative.  The morphism $\sigma^{(r)}_i$ (resp.\ $\rho^{(r)}_i$) induces a homomorphism of $\cR$-algebras
$\varsigma^{(r)}_i\colon \hcG^{(r)}\rightarrow \hcC^{(r)}$ \eqref{p1-tshbn3a}
(resp.\ $\co_X$-algebras $\varrho^{(r)}_i\colon \hcC^{(r)}\rightarrow \co_X$). 
We deduce a homomorphism 
of $\cR$-algebras $\lambda^{(r)}_i\colon \hcG^{(r)}\rightarrow \cC^\dagger$ 
(resp.\ $\co_X$-algebras $\varrho^\dagger_i\colon \cC^\dagger\rightarrow \co_X$ \eqref{p1-thbn16}).
Observe that the composition 
\begin{equation}\label{p1-tshbn32c}
\zeta^{(r)}\colon \xymatrix{
{\hcG^{(r)}}\ar[r]^-(0.5){\varsigma^{(r)}_i}&{\hcC^{(r)}}\ar[r]^{\varrho^{(r)}_i}&{\co_X}}
\end{equation}
is induced by the zero morphism $p^r\Omega\rightarrow \co_X$. 
It follows from \eqref{p1-tshbn32b} that the composition 
\begin{equation}\label{p1-tshbn32d}
\varphi^{(r)}\colon 
\xymatrix{
{\hcG^{(r)}}\ar[r]^-(0.5){\varsigma^{(r)}_2}&{\hcC^{(r)}}\ar[r]^{\varrho^{(r)}_1}&{\co_X}}
\end{equation}
is induced by the composed morphism
\begin{equation}\label{p1-tshbn32e}
p^r\Omega\longrightarrow \Omega \longrightarrow \cE \stackrel{\phi}{\longrightarrow} \co_X,
\end{equation}
where the first two arrows are the canonical morphisms. 

The Higgs field $\theta$ induces a $\hcG$-linear morphism 
$\vartheta\colon \hcG\otimes_{\cR}N\rightarrow \hcG\otimes_{\cR}N$ 
defined in \eqref{p1-tshbn31b}, and its exponential 
$\exp(\vartheta)\colon \hcG\otimes_{\cR}N\stackrel{\sim}{\rightarrow} \hcG\otimes_{\cR}N$ defined in \eqref{p1-tshbn31d}. 
On the other hand, by \ref{p1-tshbn14}, 
there exist a rational number $r>0$ and an isomorphism of $\hcG^{(r)}$-modules with $\delta_{\hcG^{(r)}}$-connection
\begin{equation}\label{p1-tshbn32f}
\exp(\theta)\colon \hcG^{(r)}\otimes_{\cR} N\stackrel{\sim}{\rightarrow} \hcG^{(r)}\otimes_{\cR} N,
\end{equation}
where the $\delta_{\hcG^{(r)}}$-connections are defined as in \ref{p1-delta-con4}, 
the module $N$ of the source (resp.\ target) being endowed with the Higgs field $\theta$ (resp.\ $0$).
By \eqref{p1-tshbn10d}, the diagram 
\begin{equation}\label{p1-tshbn32g}
\xymatrix{
{\hcG^{(r)}\otimes_{\cR}N}\ar[rr]^-(0.5){\exp(\theta)}\ar[d]_-(0.4){\ha^{r,0}\otimes \id}&&
{\hcG^{(r)}\otimes_{\cR}N}\ar[d]^-(0.4){\ha^{r,0}\otimes \id}\\
{\hcG\otimes_{\cR}N}\ar[rr]^-(0.5){\exp(\vartheta)}&&{\hcG\otimes_{\cR}N}}
\end{equation}
is commutative.
We deduce by the description \eqref{p1-tshbn32e} of $\varphi^{(r)}$ that the diagram
\begin{equation}\label{p1-tshbn32h}
\xymatrix{
{\hcG^{(r)}\otimes_{\cR}N}\ar[rr]^{\exp(\theta)}\ar[d]_{\varphi^{(r)}\otimes \id}&&
{\hcG^{(r)}\otimes_{\cR}N}\ar[d]^{\varphi^{(r)}\otimes \id}\\
{\tN}\ar[rr]^{\exp_\phi(\theta)}&&{\tN}}
\end{equation}
is commutative.  On the other hand, the diagram
\begin{equation}\label{p1-tshbn32i}
\xymatrix{
{\hcG^{(r)}\otimes_{\cR}N}\ar[rr]^{\exp(\theta)}\ar[d]_{\zeta^{(r)}\otimes \id}&&
{\hcG^{(r)}\otimes_{\cR}N}\ar[d]^{\zeta^{(r)}\otimes \id}\\
{\tN}\ar[rr]^{\id_\tN}&&{\tN,}}
\end{equation}
where $\zeta^{(r)}$ is the homomorphism \eqref{p1-tshbn32c}, is commutative.

Extending scalars by $\lambda_2^{(r)}\colon \hcG^{(r)}\rightarrow \cC^\dagger$, 
the isomorphism $\exp(\theta)$ \eqref{p1-tshbn32f} induces an isomorphism of $\cG^\dagger$-modules with 
$\delta_{\cC^\dagger}$-connection
\begin{equation}\label{p1-tshbn32k}
\upepsilon\colon \cC^\dagger\otimes_{\co_X} \tN\stackrel{\sim}{\rightarrow} \cC^\dagger\otimes_{\cR} \tN,
\end{equation}
where the $\delta_{\cC^\dagger}$-connections are defined as in \ref{p1-delta-con4}, 
the module $\tN$ of the source (resp.\ target) being endowed with the Higgs field $\ttheta$ (resp.\ $0$).
By \ref{p1-thbn37}, the latter induces an isomorphism $\varepsilon_1$ 
that fits into a commutative diagram
\begin{equation}\label{p1-tshbn32l}
\xymatrix{
{\uuptau(\tN,\ttheta)}\ar[r]_-(0.5)\sim^-(0.5){\varepsilon_2}\ar[d]_-(0.5)\iota&{\tN}\ar[d]^-(0.5)\jmath\\
{\cC^\dagger\otimes_{\co_X} \tN}\ar[r]_\sim^{\upepsilon}&{\cC^\dagger\otimes_{\cR} \tN,}}
\end{equation}
where $\iota$ is the canonical injection \eqref{p1-thbn9a} and $\jmath$ is induced by the canonical homomorphism $\co_X\rightarrow \cC^\dagger$. 

In view of the definition \eqref{p1-tshbn32c} of $\zeta^{(r)}$, 
diagrams \eqref{p1-tshbn32i} and \eqref{p1-tshbn32l} induce a diagram with commutative squares
\begin{equation}\label{p1-tshbn32m}
\xymatrix{
{\uuptau(\tN,\ttheta)}\ar[r]^-(0.5){\varepsilon_2}\ar[d]_-(0.5)\iota&{\tN}\ar[d]^-(0.5)\jmath\\
{\cC^\dagger\otimes_{\co_X} \tN}\ar[r]^{\upepsilon}\ar[d]_{\varrho_2^\dagger\otimes \id}&
{\cC^\dagger\otimes_{\cR} \tN}\ar[d]^{\varrho_2^\dagger\otimes \id}\\
{\tN}\ar[r]^{\id_\tN}&{\tN.}}
\end{equation}
Since $(\varrho_2^\dagger\otimes \id)\circ \jmath=\id_\tN$ and $(\varrho_2^\dagger\otimes \id)\circ \iota=\ttt_{\rho_2}$, 
we deduce that $\ttt_{\rho_2}=\varepsilon_2$.    

On the other hand, in view of the definition \eqref{p1-tshbn32d} of $\varphi^{(r)}$, 
the diagrams \eqref{p1-tshbn32h} and \eqref{p1-tshbn32l} induce a diagram with commutative squares
\begin{equation}\label{p1-tshbn32n}
\xymatrix{
{\uuptau(\tN,\ttheta)}\ar[r]^-(0.5){\varepsilon_2}\ar[d]_-(0.5)\iota&{\tN}\ar[d]^-(0.5)\jmath\\
{\cC^\dagger\otimes_{\co_X} \tN}\ar[r]^{\upepsilon}\ar[d]_{\varrho_1^\dagger\otimes \id}&
{\cC^\dagger\otimes_{\cR} \tN}\ar[d]^{\varrho_1^\dagger\otimes \id}\\
{\tN}\ar[r]^{\exp_\phi(\theta)}&{\tN.}}
\end{equation}
Since $(\varrho_1^\dagger\otimes \id)\circ \jmath=\id_\tN$ and $(\varrho_1^\dagger\otimes \id)\circ \iota=\ttt_{\rho_1}$, 
the proposition follows.

\begin{prop}\label{p1-tshbn18}
Let $(N,\theta)$ be a locally CL-small Higgs $\cR[\frac 1 p]$-module with coefficients in $\Omega$, 
$(\tN,\ttheta)=(N,\theta)\otimes_{\cR}\co_X$, 
$\phi_i\colon \cE\rightarrow \co_X$, for $i=1,2$, two $\co_X$-linear forms.  Then, with the notation of \ref{p1-tshbn31},
we have the equalities of the automorphisms of the Higgs module $(\tN,\ttheta)$ \eqref{p1-tshbn31a},
\begin{equation}\label{p1-tshbn18a}
\exp_{\phi_1+\phi_2}(\theta) =\exp_{\phi_1}(\theta)\circ \exp_{\phi_2}(\theta)
=\exp_{\phi_2}(\theta)\circ \exp_{\phi_1}(\theta).
\end{equation}
\end{prop}

We identify $\cG\otimes_\cR\cG$ with the symmetric algebra $\rS_\cR(\Omega\oplus \Omega)$ and denote by 
$\cG\hotimes_\cR\cG$ its $p$-adic completion. 
Let $\delta\colon \Omega\rightarrow \Omega\oplus\Omega$ be the $\cR$-linear morphism defined by $\delta(x)=(x,x)$,
and $\Delta\colon \hcG\rightarrow \cG\hotimes_\cR\cG$ the induced homomorphism.   
Let $\iota_1,\iota_2\colon \Omega\rightarrow \Omega\oplus\Omega$ be the $\cR$-linear morphisms defined
by $\iota_1(x)=(x,0)$ and $\iota_2(x)=(0,x)$, 
and $\pi_1,\pi_2\colon \hcG\rightarrow \cG\hotimes_\cR\cG$ the respective induced homomorphisms.   

The Higgs field $\theta$ induces a $\hcG$-linear morphism  \eqref{p1-tshbn31b}
\begin{equation}
\vartheta\colon \hcG\otimes_\cR N\rightarrow \hcG\otimes_\cR N.
\end{equation}
Similarly, the Higgs $\cR[\frac 1 p]$-field 
\begin{equation}
\theta\oplus \theta=(\id\otimes \delta)\circ \theta \colon N\rightarrow N\otimes_{\cR}(\Omega\oplus\Omega)
\end{equation}
induces a $\cG\hotimes_\cR\cG$-linear morphism 
\begin{equation}
\varTheta\colon (\cG\hotimes_\cR\cG)\otimes_\cR N\rightarrow (\cG\hotimes_\cR\cG)\otimes_\cR N. 
\end{equation}
Observe that the Higgs $\cR[\frac 1 p]$-field $\theta\oplus \theta$ is locally CL-small. 

The diagram 
\begin{equation}
\xymatrix{
{\hcG\otimes_\cR N}\ar[r]^-(0.5){\vartheta}\ar[d]_{\Delta\otimes \id}&{\hcG\otimes_\cR N}\ar[d]^{\Delta\otimes \id}\\
{(\cG\hotimes_\cR\cG)\otimes_\cR N}\ar[r]^-(0.5){\varTheta}&{(\cG\hotimes_\cR\cG)\otimes_\cR N}}
\end{equation}
is clearly commutative. Hence, so is the diagram 
\begin{equation}\label{p1-tshbn18b}
\xymatrix{
{\hcG\otimes_\cR N}\ar[r]^-(0.5){\exp(\vartheta)}\ar[d]_{\Delta\otimes \id}&{\hcG\otimes_\cR N}\ar[d]^{\Delta\otimes \id}\\
{(\cG\hotimes_\cR\cG)\otimes_\cR N}\ar[r]^-(0.5){\exp(\varTheta)}&{(\cG\hotimes_\cR\cG)\otimes_\cR N,}}
\end{equation}
where the $\exp(\vartheta)$ and $\exp(\varTheta)$ are the automorphisms defined in \eqref{p1-tshbn31d}.

We have 
\begin{equation}
\varTheta=\vartheta \otimes_{\hcG,\pi_1}(\cG\hotimes_\cR\cG)+\vartheta \otimes_{\hcG,\pi_2}(\cG\hotimes_\cR\cG).
\end{equation}
The endomorphisms $\vartheta \otimes_{\hcG,\pi_1}(\cG\hotimes_\cR\cG)$ and $\vartheta \otimes_{\hcG,\pi_2}(\cG\hotimes_\cR\cG)$
of $(\cG\hotimes_\cR\cG)\otimes_\cR N$ commute. Indeed, the question being local, we may assume that the $\cR$-module $\Omega$ is free, 
in which case the required property follows from the commutation of the components of $\theta$ relatively to a basis of $\Omega$.  
We deduce that 
\begin{eqnarray}\label{p1-tshbn18c}
\exp(\varTheta)&=&(\exp(\vartheta) \otimes_{\hcG,\pi_1}(\cG\hotimes_\cR\cG))\circ (\exp(\vartheta) \otimes_{\hcG,\pi_2}(\cG\hotimes_\cR\cG))\\
&=&(\exp(\vartheta) \otimes_{\hcG,\pi_2}(\cG\hotimes_\cR\cG)) \circ (\exp(\vartheta) \otimes_{\hcG,\pi_1}(\cG\hotimes_\cR\cG)).\nonumber
\end{eqnarray}

Let $\varphi_i\colon \hcG\rightarrow \co_X$, for $i=1,2$, (resp.\ $\psi\colon \hcG\rightarrow \co_X$) 
be the homomorphism induced by $\phi_i$ (resp.\ $\phi_1+\phi_2$) \eqref{p1-tshbn31j}, 
$\phi\colon \Omega\oplus\Omega\rightarrow \co_X$ the $\co_X$-linear form defined by $\phi(x,y)=\phi_1(x)+\phi_2(y)$,
$\varphi\colon \cG\hotimes_\cR\cG\rightarrow \co_X$ the homomorphism induced by $\phi$. 
Then, we have 
\begin{equation}
\psi= \varphi  \circ \Delta\ \ \ {\rm and} \ \ \ \varphi_i= \varphi \circ \pi_i \ \ \ {\rm for} \ \ \ i=1,2. 
\end{equation}  
The proposition follows from \eqref{p1-tshbn18b} and \eqref{p1-tshbn18c}, by scalar extension by $\varphi$,  
in view of the definition \eqref{p1-tshbn31a}.

\subsection{}\label{p1-tshbn17}
Let $(U_i)_{i\in I}$ be an open covering of $\fX$, and for each $i\in I$, 
$\rho_i\colon \cF|U_i\rightarrow \co_X|U_i$ a splitting of the extension \eqref{p1-thbn1a}. 
For any $(i,j)\in I^2$, we set $U_{ij}=U_i\times_\fX U_j$. 
The morphism $\rho_i|U_{ij}-\rho_j|U_{ij}\colon \cF|U_{ij}\rightarrow \co_X|U_{ij}$ induces a morphism
\begin{equation}\label{p1-tshbn17a}
\phi_{ij}\colon \cE|U_{ij}\rightarrow \co_X|U_{ij}.
\end{equation}

Let $(N,\theta)$ be a locally CL-small Higgs $\cR[\frac 1 p]$-module with coefficients in $\Omega$. 
We set $(\tN,\ttheta)=(N\otimes_{\cR}\co_X,\theta\otimes \id)$,
which is a Higgs $\co_X[\frac 1 p]$-module with coefficients in $\cE$. 
By \ref{p1-tshbn16},  $(\tN,\ttheta)$ is weakly twistable by the extension \eqref{p1-thbn1a}. 
Then by \ref{p1-thbn31}, for every $i\in I$, $\rho_i$ induces an isomorphism of Higgs $(\co_X|U_i)$-modules 
\begin{equation}\label{p1-tshbn17b}
\ttt_i\colon \uptau(\tN,\ttheta)|U_i\stackrel{\sim}{\rightarrow}(\tN|U_i,\ttheta|U_i). 
\end{equation}
For any $(i,j)\in I^2$, applying the construction of \ref{p1-tshbn31} over $U_{ij}$, we define a canonical isomorphism of Higgs $\co_X[\frac 1 p]|U_{ij}$-modules
\eqref{p1-tshbn31a}
\begin{equation}\label{p1-tshbn17d}
\exp_{\phi_{ij}}(\theta)\colon (\tN|U_{ij},\ttheta|U_{ij})\rightarrow (\tN|U_{ij},\ttheta|U_{ij}).
\end{equation}
Then by \ref{p1-tshbn32}, the diagram
\begin{equation}\label{p1-tshbn17e}
\xymatrix{
{\uptau(\tN,\theta)|U_{ij}}\ar[r]^-(0.5){\ttt_i|U_{ij}}\ar[rd]_-(0.5){\ttt_j|U_{ij}}&{(\tN|U_{ij},\ttheta|U_{ij})}\ar[d]^{\exp_{\phi_{ij}}(\theta)}\\
&{(\tN|U_{ij},\ttheta|U_{ij})}}
\end{equation}
is commutative.  We proved the following:
 
\begin{prop}\label{p1-tshbn21}
Under the assumptions of \ref{p1-tshbn17} and with the same notation, 
$\uptau(\tN,\ttheta)$ is canonically isomorphic to the Higgs $\co_X[\frac 1 p]$-module 
obtained by gluing the Higgs modules $(\tN|U_i,\ttheta|U_i)$, for $i\in I$, by the isomorphisms $\exp_{\phi_{ij}}(\theta)$, for $(i,j)\in I^2$. 
\end{prop}

\subsection{}\label{p1-tshbn22}
We set $\rT=\cHom_\cR(\Omega,\cR)$, which is a locally free $\cR$-module of finite type, and for any rational number $\varepsilon\geq 0$, 
\begin{equation}\label{p1-tshbn22a}
\rT^{(\varepsilon)}=\cHom_\cR(p^\varepsilon\Omega,\cR),
\end{equation}
that we identify with $p^{-\varepsilon}\rT$, since $\rT$ is $V$-flat. We also set $\cH=\rS_\cR(\rT)$ \eqref{p1-NC7} and 
\begin{equation}\label{p1-tshbn22b}
\cH^{(\varepsilon)}=\rS_\cR(\rT^{(\varepsilon)}).
\end{equation}
We denote by $\hcH^{(\varepsilon)}$ its $p$-adic completion. 

Let $\cB$ be an $\cR$-algebra, $\phi\colon \cH\rightarrow \cB$ a homomorphism of $\cR$-algebras. 
By \eqref{p1-delta-con1j}, the homomorphism $\cH\rightarrow \cEnd_{\cR}(\cB)$, induced by $\phi$ and the multiplication in $\cB$, 
defines a Higgs $\cR$-field 
\begin{equation}\label{p1-tshbn22d}
\theta_\phi\colon \cB\rightarrow \cB \otimes_\cR \Omega.
\end{equation}

Let $(\cN,\theta)$ be a Higgs $\cR$-module with coefficients in $\Omega$, such that $\cN$ is $V$-flat, 
$\upmu\colon \cH\rightarrow \cEnd_\cR(\cN)$
the associated homomorphism of $\cR$-algebras \eqref{p1-delta-con1j}. 
Then, $(\cN,\theta)$ is $\varepsilon$-small for a rational number $\varepsilon >0$ \eqref{p1-tshbn9}, if and only if $\mu$ can be extended 
to a homomorphism $\upmu^{(\varepsilon)}$ that fits into a commutative diagram 
\begin{equation}\label{p1-tshbn22c}
\xymatrix{
{\cH}\ar[r]^-(0.5){\upmu}\ar[d]&{\cEnd_\cR(\cN)}\\
{\cH^{(\varepsilon)}}\ar[ru]_-(0.5){\upmu^{(\varepsilon)}}&}
\end{equation}
where the vertical arrow is induced by the canonical injection $\rT\rightarrow \rT^{(\varepsilon)}$.

\begin{lem}\label{p1-tshbn270}
Let $\cN$ be a coherent  $\cR$-module which is $V$-flat,  $\varepsilon$ a rational number $\geq 0$, 
$\theta$ an $\varepsilon$-small Higgs $\cR$-field on $\cN$ with coefficients in $\Omega$ \eqref{p1-tshbn9}. 
We denote by $\mu^{(\varepsilon)}\colon \cH^{(\varepsilon)}\rightarrow \cEnd_\cR(\cN)$
the homomorphism of $\cR$-algebras induced by the Higgs $\cR$-field $\theta$ \eqref{p1-tshbn22c}
and by $\cG^{(\varepsilon)}$ its image. Then, the $\cR$-module $\cG^{(\varepsilon)}$ is coherent. 
\end{lem}

Indeed, the question being local, we may assume that the $\cR$-module $\Omega$ is free with basis $\xi_1,\dots,\xi_d$. 
Let $t_1,\dots,t_d$ be the dual basis of $\rT$. 
The $\cR$-algebra $\cH^{(\varepsilon)}$ is generated by $p^{-\varepsilon}t_1,\dots,p^{-\varepsilon}t_d$. 
We may further assume that $\fX$ is affine formal and the $\cR$-module $\cN$ is generated by $x_1,\dots,x_n\in \Gamma(\fX,\cN)$. 
Let $1\leq \ell\leq d$, $u_\ell=\upmu^{(\varepsilon)}(p^{-\varepsilon}t_\ell)$. There exists an $n\times n$ matrix 
$A_\ell=(a^{(\ell)}_{i,j})_{1\leq i,j\leq n}$ with coefficients in $\Gamma(\fX,\cR)$ such that for every $1\leq j\leq n$, we have 
\begin{equation}\label{p1-tshbn270a}
u_\ell(x_j)=\sum_{1\leq i\leq n} a^{(\ell)}_{ij}x_i.
\end{equation}
By the Cayley-Hamilton theorem, $A_\ell$ is annihilated by a monic polynomial with coefficients in $\Gamma(\fX,\cR)$, and hence so 
is $u_\ell$. Therefore, the $\cR$-module $\cG^{(\varepsilon)}$ is of finite type. 
Since the $\cR$-module $\cEnd_\cR(\cN)$ is coherent, we deduce that the $\cR$-module $\cG^{(\varepsilon)}$ is coherent.

\begin{prop}\label{p1-tshbn27}
Let $N$ be a coherent $\cR[\frac 1 p]$-module, $\theta$ a Higgs $\cR$-field on $N$ with coefficients in $\Omega$, 
$\upmu\colon \cH[\frac 1 p]\rightarrow \cEnd_\cR(N)$ the homomorphism of $\cR[\frac 1 p]$-algebras defined by $\theta$ \eqref{p1-delta-con1j}. 
We denote by $B$ the image of $\upmu$ and by $\theta_B\colon B\rightarrow B\otimes_\cR\Omega$ 
its canonical Higgs $\cR$-field \eqref{p1-tshbn22d}. Then, the following conditions are equivalent:
\begin{itemize}
\item[{\rm (i)}] The Higgs $\cR[\frac 1 p]$-module $(N,\theta)$ is CL-small in the sense of \ref{p1-tshbn13}.
\item[{\rm (ii)}] There exists a rational number $\varepsilon >\frac{1}{p-1}$ such that 
the image $\cB^{(\varepsilon)}$ of $\cH^{(\varepsilon)}$ in $B$ \eqref{p1-tshbn22b}, 
is a coherent $\cR$-module.
\item[{\rm (ii')}] There exists a rational number $\varepsilon >\frac{1}{p-1}$ 
such that the canonical homomorphism $\cH^{(\varepsilon)}\rightarrow B$ extends 
to a homomorphism of $\cR$-algebras $\hcH^{(\varepsilon)}\rightarrow B$ \eqref{p1-tshbn22b}.
\item[{\rm (ii'')}] There exists a rational number $\varepsilon >\frac{1}{p-1}$ 
such that the homomorphism $\cH^{(\varepsilon)}\rightarrow \cEnd_\cR(N)$  induced by $\upmu$ extends 
to a homomorphism of $\cR$-algebras $\hcH^{(\varepsilon)}\rightarrow \cEnd_\cR(N)$ \eqref{p1-tshbn22b}.
\item[{\rm (iii)}] The Higgs $\cR[\frac 1 p]$-module $(B,\theta_B)$ is CL-small. 
\end{itemize}
When condition {\rm (ii)} is satisfied, for every rational number $\varepsilon'$ such that $0\leq \varepsilon'<\varepsilon$,  
the image $\cB^{(\varepsilon')}$ of $\cH^{(\varepsilon')}$ in $B$ is a coherent $\cR$-module. 

If, moreover, the formal scheme $\fX$ is coherent, i.e., of finite presentation over $\cS$, 
then the three conditions are equivalent to the following condition:
\begin{itemize}
\item[{\rm (i')}] The Higgs $\cR[\frac 1 p]$-module $(N,\theta)$ is locally CL-small.
\end{itemize}
\end{prop}

Observe first that the $\cR[\frac 1 p]$-module $\cEnd_\cR(N)$ is coherent. 
Hence, the $\cR[\frac 1 p]$-module $B$ is of finite type (by the same proof of \ref{p1-tshbn270})
and hence coherent. Conditions (ii) and (ii') (resp.\ (ii) and (ii'')) are therefore equivalent by \ref{p1-pfs16}.

We prove that (i)$\Rightarrow$(ii). Let $\cN$ be a coherent $\cR$-lattice of $N$, $\varepsilon$ a rational number $>\frac{1}{p-1}$ such that  
$\theta(\cN)\subset p^\varepsilon \cN\otimes_\cR\Omega$. 
We denote by $\mu^{(\varepsilon)}\colon \cH^{(\varepsilon)}\rightarrow \cEnd_\cR(\cN)$
the homomorphism of $\cR$-algebras induced by the Higgs $\cR$-field $\theta$ \eqref{p1-tshbn22c}.
\begin{equation}\label{p1-tshbn27a}
\xymatrix{
{\cH^{(\varepsilon)}}\ar@/^1pc/[rr]^{\upmu^{(\varepsilon)}}
\ar@{->>}[r]\ar@{^(->}[d]&{\cB^{(\varepsilon)}}\ar@{^(->}[d]\ar@{^(->}[r]&{\cEnd_\cR(\cN)}\ar@{^(->}[d]\\
{\cH[\frac 1 p]}\ar@/_1pc/[rr]_{\upmu}\ar@{->>}[r]&B\ar@{^(->}[r]&{\cEnd_{\cR[\frac 1 p]}(N).}}
\end{equation}
Therefore, $\cB^{(\varepsilon)}$ is the image of $\upmu^{(\varepsilon)}$.
Hence, the $\cR$-module $\cB^{(\varepsilon)}$ is coherent by \ref{p1-tshbn270}.  

We have (ii)$\Rightarrow$(iii) since 
\begin{equation}\label{p1-tshbn27b}
\theta_B(\cB^{(\varepsilon)})\subset p^\varepsilon \Omega \otimes_\cR \cB^{(\varepsilon)}. 
\end{equation}

The implication (iii)$\Rightarrow$(ii) is a special case of the already proven implication (i)$\Rightarrow$(ii), 
since the homomorphism $B\rightarrow \cEnd_\cR(B)$ defined by the multiplication in $B$, is injective

We prove that (ii)$\Rightarrow$(i). 
The $\cR$-algebra $\cB^{(\varepsilon)}$ is topologically of finite presentation by (\cite{egr1} 1.10.4 and 6.2.10). 
It is therefore topologically of finite presentation over $\co_\fX$ by (\cite{egr1} 1.10.7). 
As a $B$-module, $N$ is coherent. 
By \ref{p1-pfs12}(ii), there exists a coherent $\cB^{(\varepsilon)}$-module $\cN$ such that $N=\cN[\frac 1 p]$. 
We may replace $\cN$ by its image in $N$ by \ref{p1-pfs15}. Then, $\cN$ is a coherent $\cR$-lattice of $N$, 
and it is obviously $\varepsilon$-small \eqref{p1-tshbn22c}. 

Under condition (ii), for every rational number $\varepsilon'$ such that $0\leq \varepsilon'<\varepsilon$,  
$\cB^{(\varepsilon')}$ is a sub-$\cR$-module of finite type of $\cB^{(\varepsilon)}$, and hence a coherent $\cR$-module.

Finally, we deduce from the equivalence of conditions (i) and (ii) and the remark above 
that if the formal scheme $\fX$ is coherent, then (i')$\Rightarrow$(ii).

\begin{cor}\label{p1-tshbn24}
Let $B$ be an $\cR[\frac 1 p]$-algebra quotient of $\cH[\frac 1 p]$ which is a coherent $\cR[\frac 1 p]$-module, 
$\theta_B\colon B\rightarrow B\otimes_\cR\Omega$ its canonical Higgs $\cR$-field \eqref{p1-tshbn22d}.
Then, the following conditions are equivalent:
\begin{itemize}
\item[{\rm (i)}] The Higgs module $(B,\theta_B)$ is CL-small.
\item[{\rm (ii)}] There exists a rational number $\varepsilon >\frac{1}{p-1}$ such that 
the image $\cB^{(\varepsilon)}$ of $\cH^{(\varepsilon)}$ in $B$ \eqref{p1-tshbn22b} is a coherent $\cR$-module.
\end{itemize}
\end{cor}

\subsection{}\label{p1-tshbn30}
We set $\trT=\rT\otimes_{\cR}\co_X$ and $\ccF=\cHom_{\co_X}(\cF,\co_X)$ and consider the dual  exact sequence of \eqref{p1-thbn1a}, 
\begin{equation}\label{p1-tshbn30a}
0\rightarrow \trT\rightarrow \ccF\rightarrow \co_X \rightarrow 0.
\end{equation}
We denote by $\Lambda$ the inverse image of the section $1\in \co_X(X)$ in $\ccF$.  
It is naturally equipped with a structure of a $\trT$-torsor. 

Let $B$ be an $\cR[\frac 1 p]$-algebra quotient of $\cH[\frac 1 p]$, 
$\theta_B\colon B\rightarrow B\otimes_\cR\Omega$ the canonical Higgs $\cR$-field.
{\em We suppose that $(B,\theta_B)$ is locally CL-small}. 
We set $(\tB,\theta_\tB)=(B,\theta)\otimes_\cR\co_X$ and $\cL_\tB=\uptau(\tB,\theta_\tB)$, 
where $\uptau$ is the functor \eqref{p1-thbn9b}. It is naturally a $\tB$-module \eqref{p1-thbn42e}.
By \ref{p1-tshbn16}, the Higgs $\co_X[\frac 1 p]$-module $(\tB,\theta_\tB)$ is weakly twistable by the extension \eqref{p1-thbn1a}.
Hence, by \ref{p1-thbn31}, the $\tB$-module $\cL_\tB$ is invertible.

By \ref{p1-thbn31}, for every open formal subscheme $U$ of $\fX$, and every $\rho\in \Lambda(U)$, 
we have a canonical isomorphism of Higgs $(\co_X|U)$-modules 
\begin{equation}\label{p1-shbn30b}
\ttt_\rho\colon \cL_\tB|U\stackrel{\sim}{\rightarrow}(\tB|U,\theta_\tB|U),
\end{equation}
that we consider as an isomorphism of $\tB|U$-modules $\ttt_\rho\colon \cL_\tB|U\stackrel{\sim}{\rightarrow}\tB|U$. 
We hence define a morphism of sheaves of sets
\begin{equation}\label{p1-tshbn30c}
\uplambda \colon 
\begin{array}[t]{clcr}
\Lambda&\rightarrow& \cL_\tB,\\
\rho&\mapsto&  \ttt_\rho^{-1}(1). 
\end{array}
\end{equation}

For every open formal subscheme $U$ of $\fX$, and every $\phi\in \trT(U)$, the $\co_X$-linear morphism 
$\ttheta_{B,\phi}\colon \tB|U\rightarrow \tB|U$ defined in \eqref{p1-tshbn31k}
is none other than the multiplication by the canonical image of $\phi$ in $\tB(U)$. 
We can define as in \eqref{p1-tshbn31a} the automorphism $\exp_\phi(\theta_B)$ of the Higgs $(\co_X|U)$-module $(\tB|U,\theta_\tB|U)$,
that we denote abusively by $\exp(\phi)$ and consider as an automorphism of the $\tB|U$-module $\tB|U$, 
\begin{equation}\label{p1-shbn30d}
\exp(\phi)\colon \tB|U\stackrel{\sim}{\rightarrow}\tB|U,
\end{equation}
and hence as a unit $\exp(\phi)\in \Gamma(U,\tB^\times)$. 
We thus define a morphism of sheaves of sets
\begin{equation}\label{p1-tshbn30e}
\exp \colon 
\begin{array}[t]{clcr}
\trT&\rightarrow& \tB^\times,\\
\phi&\mapsto& \exp(\phi).
\end{array}
\end{equation}

By \ref{p1-tshbn32} and \ref{p1-tshbn18}, for every local sections $\rho$ of $\Lambda$ and $\phi,\phi'$ of $\trT$, we have 
\begin{eqnarray}
\uplambda(\rho+\phi)&=&\exp(\phi)\uplambda(\rho),\label{p1-tshbn30f}\\
\exp(\phi+\phi')&=&\exp(\phi)\exp(\phi'). \label{p1-tshbn30g}
\end{eqnarray}
We proved the following. 

\begin{prop}\label{p1-tshbn33}
Under the assumption of \ref{p1-tshbn30} and with the same notation, the morphism $\exp\colon \trT\rightarrow \tB^\times$ \eqref{p1-tshbn30e} 
is a group homomorphism, and the morphism $\uplambda\colon \Lambda\rightarrow \cL_\tB$ \eqref{p1-tshbn30c} is $\exp$-equivariant. 
\end{prop}

\begin{cor}\label{p1-tshbn34}
Under the assumption of \ref{p1-tshbn30}, 
the invertible $\tB$-module $\cL_\tB$ is canonically isomorphic to the line bundle associated with the 
$\tB^\times$-torsor $\Lambda\wedge^{\trT}\tB^\times$ \eqref{p1-NC4}, 
deduced from $\Lambda$ by extension of its structural group by the homomorphism $\exp$ \eqref{p1-tshbn30e}.
\end{cor}

Indeed,  since the image of $\uplambda$ \eqref{p1-tshbn30c} is clearly contained in 
the subsheaf of local bases of the invertible $\tB$-module $\cL_\tB$ (i.e., the $\tB^\times$-torsor 
associated with  $\cL_\tB$), the proposition follows from \ref{p1-tshbn33}.

\subsection{}\label{p1-tshbn35}
Let $r$ be a rational number $\geq 0$. We set $\trT^{(r)}=\rT^{(r)}\otimes_{\cR}\co_X$ \eqref{p1-tshbn22a} and
$\ccF^{(r)}=\cHom_{\co_X}(\cF^{(r)},\co_X)$ and consider the dual  exact sequence of \eqref{p1-tshbn1a}, 
\begin{equation}\label{p1-tshbn35a}
0\rightarrow \trT^{(r)}\rightarrow \ccF^{(r)}\rightarrow \co_X \rightarrow 0.
\end{equation}
We denote by $\Lambda^{(r)}$ the inverse image of the section $1\in \co_X(X)$ in $\ccF^{(r)}$.  
It is naturally equipped with a structure of a $\trT^{(r)}$-torsor. 

For every integer $n\geq 0$, the sequence \eqref{p1-tshbn35a} induces an exact sequence 
\begin{equation}\label{p1-tshbn35b}
0\rightarrow \Gamma^{n+1}_{\co_X}(\trT^{(r)})\rightarrow \Gamma^{n+1}_{\co_X}(\ccF^{(r)})\rightarrow \Gamma^{n}_{\co_X}(\ccF^{(r)})\rightarrow 0,
\end{equation}
where $\Gamma_{\co_X}$ denotes the divided powers algebra \eqref{p1-NC7}.  
The $\co_X$-modules $(\Gamma^n_{\co_X}(\ccF^{(r)}))_{n\in \mN}$ therefore form a cofiltered inverse system. We set
\begin{equation}\label{p1-tshbn35c}
\cV^{(r)}=\underset{\underset{n\geq 0}{\longleftarrow}}\lim\ \Gamma^n_{\co_X}(\ccF^{(r)}). 
\end{equation}
In the remaining part of this section, when there is no risk of ambiguity, we will omit the subscript $\co_X$ 
from the notation of the divided powers algebra $\Gamma_{\co_X}$ and the symmetric algebra $\rS_{\co_X}$ \eqref{p1-NC7}. 
Recall that in \ref{p1-imdpa5} (resp.\ \ref{p1-imdpa6}), we equipped $\Gamma^{n}(\ccF^{(r)})$ (resp.\ $\cV^{(r)}$) with a canonical 
$\Gamma^{\leq n}(\trT^{(r)})$-module (resp.\ $\hGamma(\trT^{(r)})$-module) structure, depending on the extension \eqref{p1-tshbn35a}. 
It is invertible by \ref{p1-imdpa9}.

The composition of the canonical homomorphisms 
\begin{equation}\label{p1-tshbn35d}
\rS(\trT)\rightarrow \rS(\trT^{(r)})\rightarrow \Gamma(\trT^{(r)}),
\end{equation} 
where the first one is induced by the injection $p^r\Omega\rightarrow \Omega$, defines a Higgs $\co_X$-field
\begin{equation}\label{p1-tshbn35e}
\gamma^{(r)}\colon \Gamma(\trT^{(r)})\rightarrow \cE\otimes_{\co_X}\Gamma(\trT^{(r)}).
\end{equation}
For every ideal $I$ of $\Gamma(\trT^{(r)})$, we have $\gamma^{(r)}(I)\subset \cE\otimes_{\co_X} I$. 
Therefore, for every integer $n\geq 0$,  
$\gamma^{(r)}$ induces a Higgs $\co_X$-field on $\Gamma^{\leq n}(\trT^{(r)})$, that we also denote (abusively) by $\gamma^{(r)}$. 
We equip $\cC^{(r)}\otimes_{\co_X}\Gamma(\trT^{(r)})$ (resp.\ $\cC^{(r)}\otimes_{\co_X}\Gamma^{\leq n}(\trT^{(r)})$) 
with the total Higgs $\co_X$-field
\begin{equation}\label{p1-tshbn35g}
\gamma^{(r)}_\tot=\delta_{\cC^{(r)}}\otimes \id+\id\otimes \gamma^{(r)},
\end{equation}
where $\delta_{\cC^{(r)}}$ is the $\co_X$-derivation of $\cC^{(r)}$ defined in \eqref{p1-thbn3d}. 

By \ref{p1-imdpa23}(iv), for every integer $n\geq 0$, 
we have a canonical $\cC^{(r)}\otimes_{\co_X}\Gamma^{\leq n}(\trT^{(r)})$-linear isomorphism \eqref{p1-imdpa22e}
\begin{equation}\label{p1-tshbn35f}
\cC^{(r)}\otimes_{\co_X}\Gamma^n(\ccF^{(r)})\stackrel{\sim}{\rightarrow} \cC^{(r)}\otimes_{\co_X}\Gamma^{\leq n}(\trT^{(r)}).
\end{equation}
It is in fact an isomorphism of $\cC^{(r)}$-modules with $\delta_{\cC^{(r)}}$-connection, 
where the $\delta_{\cC^{(r)}}$-connections are defined as in \ref{p1-delta-con4}, 
$\Gamma^n(\ccF^{(r)})$ (resp.\ $\Gamma^{\leq n}(\trT^{(r)})$) being endowed with the Higgs field $0$ (resp.\ $\gamma^{(r)}$). 

By \ref{p1-imdpa23}(i), the morphism $\gamma^{(r)}_\tot$ is $\Gamma^{\leq n}(\trT^{(r)})$-linear,
and by \ref{p1-imdpa23}(iii), \eqref{p1-tshbn35f} induces a $\Gamma^{\leq n}(\trT^{(r)})$-linear isomorphism 
\begin{equation}\label{p1-tshbn35h}
\Gamma^n(\ccF^{(r)})\stackrel{\sim}{\rightarrow} (\cC^{(r)}\otimes_{\co_X}\Gamma^{\leq n}(\trT^{(r)}))^{\gamma^{(r)}_\tot=0}.
\end{equation}

\subsection{}\label{p1-tshbn39}
Let $r,r'$ be rational numbers such that $r\geq r'\geq 0$. 
The dual $\ctta^{r,r'}\colon \ccF^{(r')}\rightarrow \ccF^{(r)}$ 
of the canonical morphism $\tta^{r,r'}\colon \cF^{(r)}\rightarrow \cF^{(r')}$ \eqref{p1-thbn4a} fits into a commutative diagram 
\begin{equation}\label{p1-tshbn39a}
\xymatrix{
0\ar[r]&{\trT^{(r')}}\ar[r]\ar[d]&{\ccF^{(r')}}\ar[r]\ar[d]^-(0.4){\ctta^{r,r'}}&{\co_X}\ar@{=}[d]\ar[r]&0\\
0\ar[r]&{\trT^{(r)}}\ar[r]&{\ccF^{(r)}}\ar[r]&{\co_X}\ar[r]&0,}
\end{equation}
where the horizontal arrows are the exact sequences \eqref{p1-tshbn35a}, and 
the left vertical arrow is induced by the dual of the canonical injection $p^{r}\Omega\rightarrow p^{r'}\Omega$. 
We deduce a canonical isomorphism of $\trT^{(r)}$-torsors
\begin{equation}\label{p1-tshbn39d}
\Lambda^{(r')}\wedge^{\trT^{(r')}}\trT^{(r)}\stackrel{\sim}{\rightarrow} \Lambda^{(r)},
\end{equation}
where $\Lambda^{(r)}$ is the $\trT^{(r)}$-torsor defined in \ref{p1-tshbn35} and 
the source is the $\trT^{(r)}$-torsor deduced from $\Lambda^{(r')}$ by extension of its structural group 
by the canonical morphism $\trT^{(r')}\rightarrow \trT^{(r)}$.  
We also deduce from \eqref{p1-tshbn39a} a homomorphism of $\co_X$-algebras
\begin{equation}\label{p1-tshbn39b}
\hGamma(\trT^{(r')})\rightarrow \hGamma(\trT^{(r)}),
\end{equation}
and a morphism of invertible modules over it \eqref{p1-tshbn35c}
\begin{equation}\label{p1-tshbn39c}
\cV^{(r')}\rightarrow \cV^{(r)}.
\end{equation}
The diagram 
\begin{equation}\label{p1-tshbn39e}
\xymatrix{
{\Lambda^{(r')}}\ar[rr]^-(0.5){\exp_{\Lambda^{(r')}}}\ar[d]&&{\cV^{(r')}}\ar[d]\\
{\Lambda^{(r)}}\ar[rr]^-(0.5){\exp_{\Lambda^{(r)}}}&&{\cV^{(r)},}}
\end{equation}
where the horizontal arrows are the morphisms defined in \eqref{p1-imdpa13a}, is commutative. 
This is a consequence of \ref{p1-imdpa13}, taking into account the functoriality and the injectivity of $\cphi$ \eqref{p1-imdpa6a} by \ref{p1-imdpa6}(ii), 
and the functoriality of the exponential homomorphism \eqref{p1-NC7c}.

\begin{prop}\label{p1-tshbn40}
For all rational numbers $r\geq r'\geq 0$, the canonical morphism \eqref{p1-tshbn39c}
\begin{equation}\label{p1-tshbn40a}
\cV^{(r')}\otimes_{\hGamma(\trT^{(r')})}\hGamma(\trT^{(r)})\rightarrow \cV^{(r)}
\end{equation}
is an isomorphism. 
\end{prop}

Indeed, it follows from \eqref{p1-tshbn39e} and \ref{p1-imdpa13} that 
the morphism \eqref{p1-tshbn39c} maps a local basis of $\cV^{(r')}$ to a local basis of $\cV^{(r)}$.

\subsection{}\label{p1-tshbn36}
Let $B$ be an $\cR[\frac 1 p]$-algebra quotient of $\cH[\frac 1 p]$ \eqref{p1-tshbn22}, $\theta_B\colon B\rightarrow B \otimes_\cR \Omega$ 
the canonical Higgs $\cR$-field \eqref{p1-tshbn22d}.
We set  $(\tB,\ttheta_B)=(B,\theta_B)\otimes_\cR\co_X$ and $\cL_\tB=\uptau(\tB,\ttheta_B)$ \eqref{p1-thbn9b}. 
Let $\varepsilon$ be a rational number $>\frac{1}{p-1}$,  
$\cB^{(\varepsilon)}$ the image of $\cH^{(\varepsilon)}$ in $B$ \eqref{p1-tshbn22b}. 
{\em We assume that the $\cR$-module $\cB^{(\varepsilon)}$ is coherent}; so $(B,\theta_B)$ is CL-small by \ref{p1-tshbn24}, 
$(\tB,\ttheta_B)$ is weakly twistable by the extension \eqref{p1-thbn1a} by \ref{p1-tshbn16}, 
and the $\tB$-module $\cL_\tB$ is invertible by \ref{p1-thbn31}. 

The composition of the canonical homomorphisms $\cH\rightarrow \cH^{(\varepsilon)}\rightarrow \cB^{(\varepsilon)}$ 
defines a Higgs $\cR$-field  
\begin{equation}\label{p1-tshbn36a}
\theta_{\cB^{(\varepsilon)}}\colon  \cB^{(\varepsilon)}\rightarrow \Omega\otimes_\cR \cB^{(\varepsilon)},
\end{equation}
which is none other than the Higgs field induced by $\theta_B$ on $\cB^{(\varepsilon)}$ in view of \eqref{p1-tshbn27b}. 

We de note by $\tcB^{(\varepsilon)}$ the image of the canonical homomorphism $\cB^{(\varepsilon)}\otimes_\cR\co_X\rightarrow \tB$,
which is an $\co_X$-lattice of finite type of $\tB$. We clearly have $\ttheta_B(\tcB^{(\varepsilon)})\subset \cE\otimes_{\co_X}\tcB^{(\varepsilon)}$. 
We denote by $\mB^{(\varepsilon)}$ the $p$-adic completion of $\tcB^{(\varepsilon)}$, and by
\begin{equation}
\ttheta_{\cB^{(\varepsilon)}}\colon \tcB^{(\varepsilon)}\rightarrow \cE\otimes_{\co_X}\tcB^{(\varepsilon)}
\end{equation}
the Higgs $\co_X$-field induced by $\ttheta_B$.  

Let $r$ be a rational number such that $0\leq r < \varepsilon-\frac{1}{p-1}$.  
Since $\cB^{(\varepsilon)}$ is $V$-flat, the composed homomorphism of $\cR$-algebras 
$\cH^{(r)}\rightarrow \cH^{(\varepsilon)}\rightarrow \cB^{(\varepsilon)}$
factors through the canonical homomorphism $\cH^{(r)}=\rS_\cR(\rT^{(r)})\rightarrow \Gamma_\cR(\rT^{(r)})$ \eqref{p1-NC7}, 
and induces a homomorphism of $\cR$-algebras
\begin{equation}\label{p1-tshbn36bb}
\psi^{r,\varepsilon}\colon \Gamma_\cR(\rT^{(r)})\rightarrow \cB^{(\varepsilon)}.
\end{equation}
We deduce a homomorphism of $\co_X$-algebras 
\begin{equation}\label{p1-tshbn36b}
\tpsi^{r,\varepsilon}\colon \Gamma_{\co_X}(\trT^{(r)})\rightarrow \tcB^{(\varepsilon)}.
\end{equation}

We will omit $\co_X$ from the notation of the divided power algebra $\Gamma_{\co_X}$. 
For every integer $n\geq 0$, we have 
$\tpsi^{r,\varepsilon}(\Gamma^{\geq n}(\trT^{(r)}))\subset p^{\lfloor \frac{n}{d}\rfloor\alpha}\tcB^{(\varepsilon)}$, where $d$ is an integer 
$\geq {\rm rank}_{\co_X}(\trT)$ and $\alpha=\varepsilon-r-\frac{1}{p-1}>0$. We set 
\begin{equation}\label{p1-tshbn36c}
\tcB^{(\varepsilon)}_n=\tcB^{(\varepsilon)}/p^{\lfloor \frac{n+1}{d}\rfloor\alpha}\tcB^{(\varepsilon)},
\end{equation} 
and equip it with the Higgs $\co_X$-field induced by $\ttheta_{\cB^{(\varepsilon)}}$, 
that we denote also (abusively) by $\ttheta_{\cB^{(\varepsilon)}}$. 
We deduce from $\tpsi^{r,\varepsilon}$ homomorphisms of $\co_X$-algebras 
\begin{eqnarray}
\tpsi^{r,\varepsilon}_n\colon \Gamma^{\leq n}(\trT^{(r)})&\rightarrow& \tcB^{(\varepsilon)}_n,\label{p1-tshbn36d}\\
\Psi^{r,\varepsilon}\colon \hGamma(\trT^{(r)})&\rightarrow& \mB^{(\varepsilon)}.\label{p1-tshbn36dd}
\end{eqnarray}

We equip $\cC^{(r)}\otimes_{\co_X}\tcB^{(\varepsilon)}$ with the total Higgs $\co_X$-field
\begin{equation}\label{p1-tshbn36e}
\ttheta^{(r,\varepsilon)}_{\tot}=\delta_{\cC^{(r)}}\otimes \id+\id\otimes \ttheta_{\cB^{(\varepsilon)}},
\end{equation}
where $\delta_{\cC^{(r)}}$ is the $\co_X$-derivation of $\cC^{(r)}$ defined in \eqref{p1-thbn3d}. 
The Higgs fields $\ttheta_{\cB^{(\varepsilon)}}$ and $\ttheta^{(r,\varepsilon)}_{\tot}$ are $\tcB^{(\varepsilon)}$-linear,
and the homomorphism $\tpsi^{r,\varepsilon}$ \eqref{p1-tshbn36b} is compatible with the Higgs
fields $\gamma^{(r)}$ \eqref{p1-tshbn35e} and $\ttheta_{\cB^{(\varepsilon)}}$. 
Therefore, for every integer $n\geq 0$, the isomorphism \eqref{p1-tshbn35h} induces, by extension of scalars by $\tpsi^{r,\varepsilon}_n$, 
a $\tcB^{(\varepsilon)}_n$-linear morphism 
\begin{equation}\label{p1-tshbn36f}
\Gamma^n(\ccF^{(r)})\otimes_{\Gamma^{\leq n}(\trT^{(r)})}\tcB^{(\varepsilon)}_n
\rightarrow (\cC^{(r)}\otimes_{\co_X}\tcB^{(\varepsilon)}_n)^{\ttheta^{(r,\varepsilon)}_{\tot}=0}.
\end{equation}
On the other hand, the isomorphism \eqref{p1-tshbn35f} induces, by extension of scalars by $\tpsi^{r,\varepsilon}_n$, 
a $\cC^{(r)}\otimes_{\co_X}\tcB^{(\varepsilon)}_n$-linear isomorphism
\begin{equation}\label{p1-tshbn36g}
\cC^{(r)}\otimes_{\co_X}(\Gamma^n(\ccF^{(r)})\otimes_{\Gamma^{\leq n}(\trT^{(r)})}\tcB^{(\varepsilon)}_n)
\stackrel{\sim}{\rightarrow} \cC^{(r)}\otimes_{\co_X}\tcB^{(\varepsilon)}_n.
\end{equation}
Since this morphism is compatible with \eqref{p1-tshbn36f}, 
it is in fact an isomorphism of $\cC^{(r)}$-modules with $\delta_{\cC^{(r)}}$-connection, 
where the $\delta_{\cC^{(r)}}$-connections are defined as in \ref{p1-delta-con4}, 
$\Gamma^n(\ccF^{(r)})\otimes_{\Gamma^{\leq n}(\trT^{(r)})}\tcB^{(\varepsilon)}_n$ (resp.\ $\tcB^{(\varepsilon)}_n$) 
being endowed with the Higgs field $0$ (resp.\ $\ttheta_{\cB^{(\varepsilon)}}$).

\subsection{}\label{p1-tshbn37}
We keep the assumptions and notation of \ref{p1-tshbn36}. We set 
\begin{eqnarray}
\cV^{(r)}\hotimes_{\Gamma(\trT^{(r)})}\tcB^{(\varepsilon)}&=&
\underset{\underset{n\geq 0}{\longleftarrow}}\lim\ \Gamma^n(\ccF^{(r)})\otimes_{\Gamma^{\leq n}(\trT^{(r)})}\tcB^{(\varepsilon)}_n,
\label{p1-tshbn37a1}\\
\cC^{(r)}\hotimes_{\co_X}\tcB^{(\varepsilon)}&=&
\underset{\underset{n\geq 0}{\longleftarrow}}\lim\ \cC^{(r)}\otimes_{\co_X}\tcB^{(\varepsilon)}_n,\label{p1-tshbn37a2}\\
\cC^{(r)}\hotimes_{\co_X}(\cV^{(r)}\hotimes_{\Gamma(\trT^{(r)})}\tcB^{(\varepsilon)})&=&
\underset{\underset{n\geq 0}{\longleftarrow}}\lim\ 
\cC^{(r)}\otimes_{\co_X}(\Gamma^n(\ccF^{(r)})\otimes_{\Gamma^{\leq n}(\trT^{(r)})}\tcB^{(\varepsilon)}_n).\label{p1-tshbn37a3}
\end{eqnarray}
It follows from \ref{p1-imdpa8}(i) that $\cV^{(r)}\hotimes_{\Gamma(\trT^{(r)})}\tcB^{(\varepsilon)}$ is an invertible $\mB^{(\varepsilon)}$-module. 
By \ref{p1-imdpa23}(v), the isomorphisms \eqref{p1-tshbn36g} induce a $\cC^{(r)}\hotimes_{\co_X}\tcB^{(\varepsilon)}$-linear isomorphism
\begin{equation}\label{p1-tshbn37b}
\cC^{(r)}\hotimes_{\co_X}(\cV^{(r)}\hotimes_{\Gamma(\trT^{(r)})}\tcB^{(\varepsilon)})
\stackrel{\sim}{\rightarrow} \cC^{(r)}\hotimes_{\co_X}\tcB^{(\varepsilon)}.
\end{equation}
It is in fact an isomorphism of $\hcC^{(r)}$-modules with $\delta_{\hcC^{(r)}}$-connection \eqref{p1-thbn6b}, 
where the $\delta_{\hcC^{(r)}}$-connections are the inverse limits of the $\delta_{\cC^{(r)}}$-connections defined on 
the left and right hand sides of \eqref{p1-tshbn36g}. We denote by 
\begin{equation}\label{p1-tshbn37c}
\htheta^{(r,\varepsilon)}_{\tot}=\delta_{\cC^{(r)}}\hotimes \id+\id\hotimes \ttheta_{\cB^{(\varepsilon)}}
\colon \cC^{(r)}\hotimes_{\co_X}\tcB^{(\varepsilon)}\rightarrow \cE\otimes_{\co_X}
\cC^{(r)}\hotimes_{\co_X}\tcB^{(\varepsilon)}
\end{equation}
the extension of $\ttheta^{(r,\varepsilon)}_{\tot}$ \eqref{p1-tshbn36e} to the $p$-adic completions \eqref{p1-thbn18}. 
It is $\mB^{(\varepsilon)}$-linear. 
The isomorphism \eqref{p1-tshbn37b} induces then a $\mB^{(\varepsilon)}$-linear morphism
\begin{equation}\label{p1-tshbn37d}
\cV^{(r)}\hotimes_{\Gamma(\trT^{(r)})}\tcB^{(\varepsilon)}
\rightarrow (\cC^{(r)}\hotimes_{\co_X}\tcB^{(\varepsilon)})^{\htheta^{(r,\varepsilon)}_{\tot}=0}.
\end{equation}

\begin{prop}\label{p1-tshbn38}
Under the assumptions of \ref{p1-tshbn37} and with the same notation, the morphism \eqref{p1-tshbn37d} is an isomorphism. 
\end{prop}

Indeed, in view of the isomorphism \eqref{p1-tshbn37b}, it is enough to prove that 
the canonical morphism 
\begin{equation}\label{p1-tshbn38a}
\cV^{(r)}\hotimes_{\Gamma(\trT^{(r)})}\tcB^{(\varepsilon)} \rightarrow \ker(\delta_{\cC^{(r)}}\hotimes \id)
\end{equation}
is an isomorphism, where 
\begin{equation}\label{p1-tshbn38d}
\delta_{\cC^{(r)}}\hotimes \id\colon \cC^{(r)}\hotimes_{\co_X}(\cV^{(r)}\hotimes_{\Gamma(\trT^{(r)})}\tcB^{(\varepsilon)})
\rightarrow \cE\otimes_{\co_X} \cC^{(r)}\hotimes_{\co_X}(\cV^{(r)}\hotimes_{\Gamma(\trT^{(r)})}\tcB^{(\varepsilon)})
\end{equation}
denotes the inverse limit of the Higgs $\co_X$-fields $\delta_{\cC^{(r)}}\otimes \id$ \eqref{p1-tshbn37a3}. 
The question being local, denoting also (abusively) by 
\begin{equation}\label{p1-tshbn38e}
\delta_{\cC^{(r)}}\hotimes \id\colon \cC^{(r)}\hotimes_{\co_X}\tcB^{(\varepsilon)}
\rightarrow \cE\otimes_{\co_X} \cC^{(r)}\hotimes_{\co_X}\tcB^{(\varepsilon)}
\end{equation}
the inverse limit of the Higgs $\co_X$-fields $\delta_{\cC^{(r)}}\otimes \id$ \eqref{p1-tshbn37a2},
we are reduced, by \ref{p1-imdpa8}(i), to proving that the canonical morphism 
\begin{equation}\label{p1-tshbn38b}
\mB^{(\varepsilon)} \rightarrow \ker(\delta_{\cC^{(r)}}\hotimes \id)
\end{equation}
is an isomorphism. 

Observe that $\cC^{(r)}\hotimes_{\co_X}\tcB^{(\varepsilon)}$ is the $p$-adic completion of 
$\cC^{(r)}\otimes_{\co_X}\mB^{(\varepsilon)}$. 
The latter is the Higgs--Tate $\mB^{(\varepsilon)}$-algebra \eqref{p1-thbn3b} associated with the 
exact sequence of $\mB^{(\varepsilon)}$-modules 
\begin{equation}\label{p1-tshbn38c}
0\rightarrow \mB^{(\varepsilon)}\rightarrow \cF^{(r)}\otimes_{\co_X}\mB^{(\varepsilon)}\rightarrow 
\cE\otimes_{\co_X}\mB^{(\varepsilon)} \rightarrow 0, 
\end{equation} 
deduced from \eqref{p1-thbn3a}. Since $\tcB^{(\varepsilon)}$ is $V$-flat by definition \eqref{p1-tshbn36}, $\mB^{(\varepsilon)}$ is $V$-flat. 
Moreover, $\delta_{\cC^{(r)}}\hotimes \id$ \eqref{p1-tshbn38e} is the $\mB^{(\varepsilon)}$-analogue of \eqref{p1-thbn6b}. 
Therefore, by \ref{p1-thbn22}(iii), applied to the ringed topos $(\fX_\zar,\mB^{(\varepsilon)})$ and the extension \eqref{p1-tshbn38c}, 
we obtain that \eqref{p1-tshbn38b} is an isomorphism. 

\subsection{}\label{p1-tshbn43}
We keep the assumption and notation of \ref{p1-tshbn36} and \ref{p1-tshbn37} and let $\rho\colon \cF^{(r)}\rightarrow \co_X$
be an $\co_X$-linear form such that $\rho\circ c^{(r)}=\id_{\co_X}$ \eqref{p1-thbn3a}, i.e., $\rho\in \Lambda^{(r)}(X)$  \eqref{p1-tshbn35}. 
Then $\rho$ extends uniquely to a homomorphism of $\co_X$-algebras  $\varrho\colon \cC^{(r)}\rightarrow \co_X$ \eqref{p1-imdpa1d}. 
For every integer $n\geq 0$, the $\cC^{(r)}\otimes_{\co_X}\tcB^{(\varepsilon)}_n$-linear isomorphism \eqref{p1-tshbn36g} 
induces by extension of scalars by $\varrho$ a $\tcB^{(\varepsilon)}_n$-linear isomorphism
\begin{equation}\label{p1-tshbn43a}
\ttt_{\rho,n}\colon \Gamma^n(\ccF^{(r)})\otimes_{\Gamma^{\leq n}(\trT^{(r)})}\tcB^{(\varepsilon)}_n
\stackrel{\sim}{\rightarrow} \tcB^{(\varepsilon)}_n. 
\end{equation}
The latter induces a $\mB^{(\varepsilon)}$-linear isomorphism
\begin{equation}\label{p1-tshbn43b}
\ttt_{\rho}\colon \cV^{(r)}\hotimes_{\Gamma(\trT^{(r)})}\tcB^{(\varepsilon)}
\stackrel{\sim}{\rightarrow} \mB^{(\varepsilon)}.
\end{equation}

\begin{prop}\label{p1-tshbn44}
Under the assumption of \ref{p1-tshbn43}, and with the same notation, we have
\begin{equation}
\ttt_{\rho}(\exp_{\Lambda^{(r)}}(\rho)\hotimes 1)=1,
\end{equation}
where $\exp_{\Lambda^{(r)}}$ is the morphism \eqref{p1-imdpa13a}. 
\end{prop}

It follows immediately from \ref{p1-imdpa26}.

\subsection{}\label{p1-tshbn42}
We take again the assumptions and notation of \ref{p1-tshbn36},  
and let $r'$ be a rational number such that $r\geq r'\geq 0$. It follows from the definition \eqref{p1-imdpa22d} of the morphism \eqref{p1-imdpa22c}
that the diagram 
\begin{equation}\label{p1-tshbn42a}
\xymatrix{
{\Gamma^n(\ccF^{(r')})}\ar[r]\ar[d]&{\cC^{(r')}\otimes_{\co_X}\Gamma^{\leq n}(\trT^{(r')})}\ar[r]&
{\cC^{(r')}\otimes_{\co_X}\Gamma^{\leq n}(\trT^{(r)})}\\
{\Gamma^n(\ccF^{(r)})}\ar[r]&{\cC^{(r)}\otimes_{\co_X}\Gamma^{\leq n}(\trT^{(r)})}\ar[ru]&}
\end{equation}
where the left horizontal arrows are induced by \eqref{p1-tshbn35h}, is commutative. Since the diagram 
\begin{equation}\label{p1-tshbn42b}
\xymatrix{
{\Gamma^{\leq n}(\trT^{(r')})}\ar[r]^-(0.5){\tpsi^{r',\varepsilon}_n}\ar[d]&{\tcB^{(\varepsilon)}_n}\\
{\Gamma^{\leq n}(\trT^{(r)})}\ar[ru]_{\tpsi^{r,\varepsilon}_n}&}
\end{equation}
where $\tpsi^{r,\varepsilon}_n$ is defined in \eqref{p1-tshbn36d}, is commutative, we deduce that the diagram
\begin{equation}\label{p1-tshbn42c}
\xymatrix{
{\Gamma^n(\ccF^{(r')})\otimes_{\Gamma^{\leq n}(\trT^{(r')})}\tcB^{(\varepsilon)}_n}\ar[r]\ar[d]&{\cC^{(r')}\otimes_{\co_X}\tcB^{(\varepsilon)}_n}\\
{\Gamma^n(\ccF^{(r)})\otimes_{\Gamma^{\leq n}(\trT^{(r)})}\tcB^{(\varepsilon)}_n}\ar[r]&{\cC^{(r)}\otimes_{\co_X}\tcB^{(\varepsilon)}_n,}\ar[u]}
\end{equation}
where the horizontal arrows are induced by \eqref{p1-tshbn36f}, is commutative.

\begin{prop}\label{p1-tshbn41}
We take again the notation of \ref{p1-tshbn36} and \ref{p1-tshbn37} and assume, moreover, that $r>0$. Then, 
\begin{itemize}
\item[{\rm (i)}] For every rational number $r'$ such that $0< r'\leq r$, 
the diagram 
\begin{equation}\label{p1-tshbn41a}
\xymatrix{
{\cV^{(r')}\hotimes_{\Gamma(\trT^{(r')})}\tcB^{(\varepsilon)}}\ar[r]\ar[d]&
{(\cC^{(r')}\hotimes_{\co_X}\tcB^{(\varepsilon)})^{\htheta^{(r',\varepsilon)}_{\tot}=0}}\\
{\cV^{(r)}\hotimes_{\Gamma(\trT^{(r)})}\tcB^{(\varepsilon)}}\ar[r]&
{(\cC^{(r)}\hotimes_{\co_X}\tcB^{(\varepsilon)})^{\htheta^{(r,\varepsilon)}_{\tot}=0},}\ar[u]}
\end{equation}
where the horizontal arrows are induced by \eqref{p1-tshbn37d}, is commutative. Moreover, all arrows are isomorphisms. 
\item[{\rm (ii)}] There exists a canonical $\tcB^{(\varepsilon)}$-linear morphism 
\begin{equation}\label{p1-tshbn41b}
\upnu^{(r,\varepsilon)}_B\colon \cL_\tB\rightarrow (\cV^{(r)}\hotimes_{\Gamma(\trT^{(r)})}\tcB^{(\varepsilon)})\otimes_{\mZ_p}\mQ_p
\end{equation}
such that for all rational numbers $0< r'\leq r$, the diagram 
\begin{equation}\label{p1-tshbn41c}
\xymatrix{
&{(\cV^{(r')}\hotimes_{\Gamma(\trT^{(r')})}\tcB^{(\varepsilon)})\otimes_{\mZ_p}\mQ_p}\ar[d]\\
{\cL_\tB}\ar[r]_-(0.5){\upnu^{(r,\varepsilon)}_B}\ar[ru]^-(0.5){\upnu^{(r',\varepsilon)}_B}&
{(\cV^{(r)}\hotimes_{\Gamma(\trT^{(r)})}\tcB^{(\varepsilon)})\otimes_{\mZ_p}\mQ_p}\ar[d]\\
{(\hcC^{(r)}\otimes_{\cR}B)^{\theta^{(r)}_{\tot}=0}}\ar[r]\ar[u]&
{(\cC^{(r)}\hotimes_{\co_X}\tcB^{(\varepsilon)})^{\htheta^{(r,\varepsilon)}_{\tot}=0}\otimes_{\mZ_p}\mQ_p,}}
\end{equation}
where the lower horizontal arrow is induced by the morphism 
$\hcC^{(r)}\otimes_{\cR}\cB^{(\varepsilon)} \rightarrow \cC^{(r)}\hotimes_{\co_X}\tcB^{(\varepsilon)}$
and $\theta^{(r)}_{\tot}=\delta_{\hcC^{(r)}}\otimes \id+\id\otimes \theta_B$, is commutative. 
\item[{\rm (iii)}] The morphism \eqref{p1-tshbn41b} induces a $\mB^{(\varepsilon)}$-linear isomorphism 
\begin{equation}\label{p1-tshbn41d}
\cL_\tB \otimes_{\tcB^{(\varepsilon)}}\mB^{(\varepsilon)}
\stackrel{\sim}{\rightarrow} (\cV^{(r)}\hotimes_{\Gamma(\trT^{(r)})}\tcB^{(\varepsilon)})\otimes_{\mZ_p}\mQ_p. 
\end{equation}
\item[{\rm (iv)}] The diagram
\begin{equation}\label{p1-tshbn41e}
\xymatrix{
\Lambda\ar[d]_{\uplambda}\ar[r]^-(0.5){\exp_\Lambda}&{\cV}\ar[r]&{\cV^{(r)}}\ar[d]\\
{\cL_\tB}\ar[rr]^-(0.5){\upnu^{(r,\varepsilon)}_B}&&{\cV^{(r)}\hotimes_{\Gamma(\trT^{(r)})}\tcB^{(\varepsilon)},}}
\end{equation}
where $\Lambda=\Lambda^{(0)}$, $\cV=\cV^{(0)}$ \eqref{p1-tshbn35}, 
$\exp_\Lambda$ is defined in \eqref{p1-imdpa13a} and $\uplambda$ in \eqref{p1-tshbn30c},  is commutative.
\end{itemize}
\end{prop}

(i) Indeed, the commutativity of the diagram \eqref{p1-tshbn41a} follows from that of \eqref{p1-tshbn42c}. 
The last assertion follows then from \ref{p1-tshbn38} and \ref{p1-tshbn40}. 

(ii) It follows from (i) and the definition of the functor $\uptau$ \eqref{p1-thbn9b}.

(iii) The question being local, we may assume that there exists a splitting  
$\rho\colon \cF\rightarrow \co_X$ of the extension \eqref{p1-thbn1a}. 
By \ref{p1-tshbn16}, $(\tB,\theta_\tB)$ is weakly twistable by the extension \eqref{p1-thbn1a}. 
Then by \ref{p1-thbn31}, $\rho$ induces a $\tB$-linear isomorphism 
\begin{equation}\label{p1-tshbn41f}
\ttt_\rho\colon \cL_\tB\stackrel{\sim}{\rightarrow}\tB. 
\end{equation}
The splitting $\rho$ induces a splitting $\rho^{(r)}\colon \cF^{(r)}\rightarrow \co_X$ 
of the extension \eqref{p1-thbn3a}, and hence by \ref{p1-tshbn43}, a $\mB^{(\varepsilon)}$-linear isomorphism
\begin{equation}\label{p1-tshbn41g}
\ttt_{\rho^{(r)}}\colon \cV^{(r)}\hotimes_{\Gamma(\trT^{(r)})}\tcB^{(\varepsilon)}
\stackrel{\sim}{\rightarrow} \mB^{(\varepsilon)}.
\end{equation}
It follows from the commutativity of the lower square of \eqref{p1-tshbn41c} that the diagram 
\begin{equation}\label{p1-tshbn41h}
\xymatrix{
{\cL_\tB}\ar[r]^-(0.5){\upnu^{(r,\varepsilon)}_B}\ar[d]_{\ttt_\rho}&
{(\cV^{(r)}\hotimes_{\Gamma(\trT^{(r)})}\tcB^{(\varepsilon)})\otimes_{\mZ_p}\mQ_p}\ar[d]^{\ttt_{\rho^{(r)}}}\\
{\tB}\ar[r]&{\mB^{(\varepsilon)}[\frac 1 p],}}
\end{equation}
where the lower horizontal morphism is induced by the canonical homomorphism $\tB^{(\varepsilon)}\rightarrow \mB^{(\varepsilon)}$, is commutative, which proves the proposition. 

(iv) It follows from \ref{p1-tshbn44}, \eqref{p1-tshbn41h} and \eqref{p1-tshbn39e}.

\begin{lem}\label{p1-tshbn50}
Let $d$ be an integer such that $\rk_{\co_X}(\cE)\leq d$ \eqref{p1-thbn1a}, 
and let $\beta,s$ be rational numbers such that $0\leq \beta<s+\frac{1}{p-1}$. We set 
\begin{equation}\label{p1-tshbn50a}
\gamma=\inf\{ v_p(n!)-(\beta-s)n; n\in \mN\},
\end{equation}
which is a rational number $\leq 0$. Then, we have a commutative diagram of $\co_X$-linear morphisms
\begin{equation}\label{p1-tshbn50b}
\xymatrix{
{\rS(\trT)\otimes_{\co_X}\cC^{(s)}}\ar[r]^-(0.5){\mu^s}\ar[d]&{\cC^{(s)}}\ar[d]\\
{\rS(\trT^{(\beta)})\otimes_{\co_X}\cC^{(s)}}\ar[r]^-(0.5){\mu^s_\beta}&{p^{\gamma d}\cC^{(s)},}}
\end{equation}
where $\trT=\rT\otimes_{\cR}\co_X$ \eqref{p1-tshbn22}, $\trT^{(\beta)}=\rT^{(\beta)}\otimes_{\cR}\co_X$ \eqref{p1-tshbn35},
$\rS=\rS_{\co_X}$ denotes the symmetric $\co_X$-algebra \eqref{p1-NC7}, $\mu^s$ is defined by the Higgs $\co_X$-field 
$\delta_{\cC^{(s)}}$ on $\cC^{(s)}$ \eqref{p1-thbn3d} and the vertical arrows are the canonical morphisms. 
Moreover, $\mu^s_\beta$ is uniquely determined by \eqref{p1-tshbn50b}. 
\end{lem}

Observe first that $\lim_{n\mapsto +\infty} v_p(n!)-(\beta-s)n =+\infty$, which implies that $\gamma$ is a rational number. 
Since $\cC^{(s)}$ and $\trT$ are $p$-torsion free, $\mu^s_\beta$ is uniquely determined by \eqref{p1-tshbn50b}. 
Hence, we may assume that $\cE$ is a free $\co_X$-module of rank $d$ and the exact sequence \eqref{p1-thbn1a} splits.
A splitting of \eqref{p1-thbn1a} induces a splitting of the extension \eqref{p1-tshbn1a}
\begin{equation}
0\rightarrow \co_X\rightarrow \cF^{(s)}\rightarrow p^s\cE \rightarrow 0,
\end{equation}
and hence an isomorphism $\rS(p^s\cE)\stackrel{\sim}{\rightarrow} \cC^{(s)}$. 
We choose an $\co_X$-basis $X_1,\dots,X_d$ of $\cE$ and denote by $\partial_1,\dots,\partial_d$ the dual $\co_X$-basis of $\trT$. 
We deduce isomorphisms of $\co_X$-algebras $\cC^{(s)}\stackrel{\sim}{\rightarrow} \co_X[p^sX_1,\dots,p^sX_d]$,  
$\rS(\trT)\stackrel{\sim}{\rightarrow} \co_X[\partial_1,\dots,\partial_d]$ and 
$\rS(\trT^{(\beta)})\stackrel{\sim}{\rightarrow} \co_X[p^{-\beta}\partial_1,\dots,p^{-\beta}\partial_d]$. The morphism $\mu^s$ \eqref{p1-tshbn50b} 
is defined, for any $1\leq i\leq d$, by letting $\partial_i$ acts on 
$\co_X[p^sX_1,\dots,p^sX_d]$ as the partial derivative relatively to $X_i$. The proposition follows immediately.

\section{Coherent modules up to isogeny on formal schemes}\label{p2-cmupiso}

\subsection{}\label{p2-cmupiso1}
In this section,  $\fX$ denotes a formal $\cS$-scheme of finite presentation \eqref{p1-thbn0} (\cite{egr1} 2.3.15);
so it is a {\em quasi-compact} idyllic formal scheme (\cite{egr1} 2.6.13). For any $\co_\fX$-algebra $\cA$ of $\fX_\zar$, 
we denote by $\bMod(\cA)$ the category of $\cA$-modules of $\fX_{\zar}$, 
by $\bIndMod(\cA)$ the category of ind-$\cA$-modules  and by
\begin{equation}\label{p2-cmupiso1a}
\iota_{\cA}\colon \bMod(\cA)\rightarrow \bIndMod(\cA)
\end{equation}
the canonical functor, which is exact and fully faithful (see \ref{p1-indmal1}).
When there is no risk of ambiguity, we will identify $\bMod(\cA)$ with a full subcategory of $\bIndMod(\cA)$
by the functor $\iota_{\cA}$, which we will omit from the notation.

The functor $\iota_{\cA}$ admits a left adjoint 
\begin{equation}\label{p2-cmupiso1b}
\kappa_{\cA}\colon \bIndMod(\cA) \rightarrow \bMod(\cA).
\end{equation}
It is exact and the adjunction morphism $\kappa_{\cA}\circ \iota_{\cA} \rightarrow \id_{\bMod(\cA)}$ is an isomorphism. 

We denote by $\bMod_\mQ(\cA)$ the category of $\cA$-modules up to isogeny (i.e. $\cA_{\mQ}$-modules), by 
\begin{equation}\label{p2-cmupiso1d}
Q_\cA\colon \bMod(\cA)\rightarrow \bMod_\mQ(\cA),\ \ \ \cM\mapsto \cM_\mQ,
\end{equation}
the localisation functor and by
\begin{equation}\label{p2-cmupiso1c}
\upalpha_{\cA}\colon \bMod_\mQ(\cA)\rightarrow \bIndMod(\cA)
\end{equation}
the canonical exact and fully faithful functor (see \ref{p1-bcim5}).

\subsection{}\label{p2-cmupiso2}
Let $\cA$ an $\co_\fX$-algebra topologically of finite presentation \eqref{p1-pfs1}. 
For $A=\cA$ or $\cA[\frac 1 p]$, we denote by $\bMod^\coh(A)$ the category of {\em coherent} $A$-modules,
by $\bMod^{\coh}_\mQ(\cA)$ the category of coherent $\cA$-modules up to isogeny (see \ref{p1-abisoind1})
and by $\jmath_\cA\colon \bMod^\coh(\cA)\rightarrow \bMod(\cA)$ the canonical injection functor. 
By \ref{p1-abisoind40}, the functor
\begin{equation}\label{p2-cmupiso2a}
\bMod^{\coh}_\mQ(\cA)\stackrel{\sim}{\rightarrow} \bMod^{\coh}(\cA[\frac 1 p]),\ \ \ 
\cF \mapsto \cF[\frac  1 p],
\end{equation}
is an equivalence of abelian categories. 
We denote by $\ojmath_\cA$ the composition of a quasi-inverse of \eqref{p2-cmupiso2a} with the fully faithful functor $\jmath_{\cA,\mQ}$,
\begin{equation}\label{p2-cmupiso2b}
\ojmath_\cA\colon \xymatrix{
{\bMod^{\coh}(\cA[\frac 1 p])}\ar[r]&{\bMod^{\coh}_\mQ(\cA)}\ar[r]^{\jmath_{\cA,\mQ}}&{\bMod_\mQ(\cA).}} 
\end{equation}
Let $\upalpha^\coh_{\cA}$ be the composed functor
\begin{equation}\label{p2-cmupiso2c}
\upalpha^\coh_{\cA}=\upalpha_{\cA}\circ \ojmath_\cA\colon\bMod^{\coh}(\cA[\frac 1 p])\rightarrow \bIndMod(\cA),
\end{equation}
which is fully faithful and exact. 
When there is no risk of ambiguity, we will identify $\bMod^{\coh}(\cA[\frac 1 p])$ with a full subcategory of $\bIndMod(\cA)$ 
by $\upalpha^\coh_{\cA}$. 
We will therefore consider any coherent $\cA[\frac 1 p]$-module also as an ind-$\cA$-module. 

The functor $\upalpha^\coh_{\cA}$ is compatible with the natural monoidal structures (\cite{ag2} (2.9.7.4)).  
It therefore transforms an $\cA[\frac 1 p]$-algebra which is a coherent 
$\cA[\frac 1 p]$-module into an ind-$\cA$-algebra \ref{p1-indmal2}.  

By (\cite{ag2} (2.6.5.2) and 2.9.3(i)), the adjunction morphism $\kappa_{\cA}\circ \iota_{\cA} \rightarrow \id_{\bMod(\cA)}$,
which is an isomorphism, induces for every coherent $\cA[\frac 1 p]$-module $N$, a canonical functorial isomorphism
\begin{equation}\label{p2-cmupiso2d}
\kappa_{\cA}(\upalpha^\coh_{\cA}(N))\stackrel{\sim}{\rightarrow} N.
\end{equation}

We denote by $\bIndMod^{\mQ\coh}(\cA)$ the essential image of the functor \eqref{p2-cmupiso2c} and (abusively) by 
\begin{equation}\label{p2-cmupiso2h}
\upalpha^\coh_{\cA}\colon\bMod^{\coh}(\cA[\frac 1 p])\rightarrow \bIndMod^{\mQ\coh}(\cA)
\end{equation}
the induced equivalence of categories. By \eqref{p2-cmupiso2d}, the functor $\kappa_{\cA}$
induces a functor 
\begin{equation}\label{p2-cmupiso2e}
\kappa^\coh_{\cA}\colon \bIndMod^{\mQ\coh}(\cA) )\rightarrow \bMod^{\coh}(\cA[\frac 1 p]), 
\end{equation}
which is in fact a quasi-inverse of $\upalpha^\coh_{\cA}$. 

\begin{prop}\label{p2-cmupiso3}
Under the assumptions of \ref{p2-cmupiso2} and with the same notation, 
the functor $\upalpha^\coh_{\cA}$ \eqref{p2-cmupiso2c} makes the category 
$\bMod^{\coh}(\cA[\frac 1 p])$ a thick subcategory of $\bIndMod(\cA)$. 
\end{prop}

Indeed, the diagram of functors 
\begin{equation}
\xymatrix{
{\bMod^{\coh}_\mQ(\cA)}\ar[r]^-(0.5){\upalpha}\ar[d]_{\jmath_{\cA,\mQ}}&{\bInd(\bMod^\coh(\cA))}\ar[d]^{\rI \jmath_\cA}\\
{\bMod_\mQ(\cA)}\ar[r]^-(0.5){\upalpha_{\cA}}&{\bIndMod(\cA),}}
\end{equation}
where $\upalpha$ is the functor defined in \eqref{p1-abisoind1d}, 
is commutative up to isomorphism. By \ref{p1-abisoind2}, $\upalpha$ makes $\bMod^{\coh}_\mQ(\cA)$ a thick subcategory of
$\bInd(\bMod^\coh(\cA))$. Since $\jmath_\cA$ makes $\bMod^\coh(\cA)$ a thick subcategory of $\bMod(\cA)$, 
$\rI\jmath_\cA$ makes $\bInd(\bMod^\coh(\cA))$ a thick subcategory of $\bIndMod(\cA)$ by (\cite{ks2} 6.1.10, 8.6.8 and 8.6.12).  
Therefore, $\upalpha_{\cA}\circ \jmath_{\cA,\mQ}$ makes $\bMod^{\coh}_\mQ(\cA)$ a thick subcategory of $\bIndMod(\cA)$, 
which proves the proposition.

\subsection{}\label{p2-cmupiso31}
Let $\cE$ be a coherent $\co_\fX$-module. For $A=\co_\fX$ or $\co_\fX[\frac 1 p]$, 
we consider Higgs $A$-modules with coefficients in $\cE\otimes_{\co_\fX}A$ \eqref{p1-delta-con1}. 
We say abusively that they have coefficients in $\cE$.
Such a Higgs module is said to be {\em coherent} if the underlying $A$-module is coherent.
The categories of these modules will be denoted by $\bHM(A,\cE)$ and $\bHM^\coh(A, \cE)$, respectively. 
We denote by $\bIndHM(\co_\fX,\cE)$ the category of Higgs ind-$\co_\fX$-modules with coefficients in $\cE$ \eqref{p1-indmal20}. 
By (\cite{ag2} (2.7.1.5)), the functors $\iota_{\co_\fX}$ and $\kappa_{\co_\fX}$ induce naturally functors that we denote (abusively) also by
\begin{eqnarray}
\iota_{\co_\fX}\colon \bHM(\co_\fX,\cE)\rightarrow \bIndHM(\co_\fX,\cE),\label{p2-cmupiso31a}\\
\kappa_{\co_\fX}\colon \bIndHM(\co_\fX,\cE) \rightarrow \bHM(\co_\fX,\cE),\label{p2-cmupiso31b}
\end{eqnarray}
the second one being a left adjoint of the first one. 

We denote by $\bIH(\co_\fX,\cE)$ the category of Higgs  $\co_\fX$-isogenies
with coefficients in $\cE$ (\cite{ag2} 2.9.9) and by
$\bIH^\coh(\co_\fX,\cE)$ the full subcategory
made up  of the quadruples $(\cM,\cN,u,\theta)$ such that $\cM$ and $\cN$ are coherent $\co_\fX$-modules. These are additive categories.
We denote by $\varphi\colon \bIH^\coh(\co_\fX,\cE)\rightarrow \bIH(\co_\fX,\cE)$ the canonical injection functor and 
by $\bIH_\mQ(\co_\fX,\cE)$ (resp.\ $\bIH^\coh_\mQ(\co_\fX, \cE)$) the category of objects of $\bIH(\co_\fX,\cE)$
(resp.\ $\bIH^\coh(\co_\fX,\cE)$) up to isogeny (\cite{agt} III.6.1). The functor
\begin{equation}\label{p2-cmupiso31c}
\begin{array}[t]{clcr}
\bIH(\co_\fX,\cE)&\rightarrow& \bHM(\co_\fX[\frac 1 p], \cE),\\
(\cM,\cN,u,\theta)&\mapsto& (\cM_{\mQ_p}, (\id \otimes u_{\mQ_p}^{-1})\circ\theta_{\mQ_p}),
\end{array}
\end{equation}
induces a functor
\begin{equation}\label{p2-cmupiso31d}
\bIH_\mQ(\co_\fX,\cE)\rightarrow \bHM(\co_\fX[\frac 1 p], \cE).
\end{equation}
By (\cite{agt} III.6.21), the latter induces an equivalence of categories
\begin{equation}\label{p2-cmupiso31e}
\bIH^\coh_\mQ(\co_\fX,\cE)\stackrel{\sim}{\rightarrow}
\bHM^\coh(\co_\fX[\frac 1 p], \cE).
\end{equation}
We denote by $\ovarphi$ the composition of a quasi-inverse of \eqref{p2-cmupiso31e} with the fully faithful functor $\varphi_\mQ$,
\begin{equation}\label{p2-cmupiso31f}
\ovarphi\colon 
\xymatrix{
{\bHM^\coh(\co_\fX[\frac 1 p], \cE)}\ar[r]&{\bIH^\coh_\mQ(\co_\fX,\cE)}\ar[r]^-(0.5){\varphi_\mQ}&{\bIH_\mQ(\co_\fX,\cE).}}
\end{equation}

In view of (\cite{ag2} (2.9.7.4)), the functor $\upalpha_{\co_\fX}$ induces a functor that we denote (abusively) also by 
\begin{equation}\label{p2-cmupiso31g}
\upalpha_{\co_\fX}\colon 
\begin{array}[t]{clcr}
\bIH_\mQ(\co_\fX,\cE)&\rightarrow& \bIndHM(\co_\fX,\cE)\\
(\cM,\cN,u,\theta)&\mapsto& (\upalpha_{\co_\fX}(\cM_\mQ),(\id \otimes \upalpha_{\co_\fX}(u_\mQ)^ {-1})\circ\upalpha_{\co_\fX}(\theta_\mQ)).
\end{array}
\end{equation}
It is fully faithful by (\cite{ag2} (2.6.6.4) and (2.9.10.3)).
We denote (abusively) by $\upalpha^\coh_{\co_\fX}$ the composed functor
\begin{equation}\label{p2-cmupiso31h}
\upalpha^\coh_{\co_\fX}=\upalpha_{\co_\fX}\circ \ovarphi\colon 
\bHM^\coh(\co_\fX[\frac 1 p], \cE)\rightarrow \bIndHM(\co_\fX,\cE).
\end{equation}
We will identify $\bHM^\coh(\co_\fX[\frac 1 p], \cE)$ with a full subcategory of $\bIndHM (\co_\fX,\cE)$
by this functor, which we will often omit from the notation. We will therefore consider any
coherent Higgs $\co_\fX[\frac 1 p]$-module with coefficients in $\cE$ as
a Higgs ind-$\co_\fX$-module with coefficients in $\cE$.

Since the equivalences of categories \eqref{p2-cmupiso31e} and \eqref{p2-cmupiso2a} are compatible in an obvious way, 
so are their quasi-inverses and hence also \eqref{p2-cmupiso31h} and \eqref{p2-cmupiso2c}, denoted both by $\upalpha^\coh_{\co_\fX}$;
this abuse of notation does thus not lead to any ambiguity.

\subsection{}\label{p2-cmupiso310}
Let $\cA$ be an $\co_\fX$-algebra topologically of finite presentation, 
$\cE$ a coherent $\co_\fX$-module, $\delta\colon \cA\rightarrow \cE\otimes_{\co_\fX}\cA$ an $\co_\fX$-derivation,
which is also a Higgs $\co_\fX$-field on $\cA$ with coefficients in $\cE$. 
We denote by $\bMIC_\delta(\cA[\frac 1 p]/\co_\fX[\frac 1 p])$ the category of $\cA[\frac 1 p]$-modules with integrable $\delta[\frac 1 p]$-connection 
\eqref{p1-delta-con2} and by $\bMIC^\coh_\delta(\cA[\frac 1 p]/\co_\fX[\frac 1 p])$ the full subcategory of objects whose underlying 
$\cA[\frac 1 p]$-module is coherent. 
We denote by $\bIMIC_\delta(\cA/\co_\fX)$ the category of {\em $\delta$-isoconnections} with respect to the extension $\cA/\co_\fX$ 
\eqref{p1-delta-con8} and by $\bIMIC^\coh_\delta(\cA/\co_\fX)$ the full subcategory
made up  of the quadruples 
\begin{equation}\label{p2-cmupiso310a}
(\cM,\cN,u\colon \cM\rightarrow \cN,\nabla\colon \cM\rightarrow \cE\otimes_{\co_\fX}\cN)
\end{equation} 
such that the $\cA$-modules $\cM$ and $\cN$ are coherent. These are additive categories.
We denote by $\psi_\delta\colon \bIMIC^\coh_\delta(\cA/\co_\fX)\rightarrow \bIMIC_\delta(\cA/\co_\fX)$ the canonical injection functor and 
by $\bIMIC_{\delta,\mQ}(\cA/\co_\fX)$ (resp.\ $\bIMIC^\coh_{\delta,\mQ}(\cA/\co_\fX)$) the category of objects of 
$\bIMIC_\delta(\cA/\co_\fX)$ (resp.\ $\bIMIC^\coh_\delta(\cA/\co_\fX)$) up to isogeny (\cite{agt} III.6.1).

The functor
\begin{equation}\label{p2-cmupiso310b}
\begin{array}[t]{clcr}
\bIMIC_\delta(\cA/\co_\fX)&\rightarrow& \bMIC_\delta(\cA[\frac 1 p]/\co_\fX[\frac 1 p]),\\
(\cM,\cN,u,\nabla)&\mapsto& (\cM[\frac 1 p], (\id \otimes (u[\frac 1 p])^{-1})\circ\nabla[\frac 1 p]),
\end{array}
\end{equation}
induces a functor
\begin{equation}\label{p2-cmupiso310c}
\bIMIC_{\delta,\mQ}(\cA/\co_\fX)\rightarrow \bMIC_\delta(\cA[\frac 1 p]/\co_\fX[\frac 1 p]).
\end{equation}
By \ref{p2-cmupiso311}(ii) below, the latter induces an equivalence of categories
\begin{equation}\label{p2-cmupiso310d}
\bIMIC^\coh_{\delta,\mQ}(\cA/\co_\fX)\stackrel{\sim}{\rightarrow}
\bMIC^\coh_\delta(\cA[\frac 1 p]/\co_\fX[\frac 1 p]).
\end{equation}
We denote by $\opsi_\delta$ the composition of a quasi-inverse of \eqref{p2-cmupiso310d} with the fully faithful functor $\psi_{\delta,\mQ}$,
\begin{equation}\label{p2-cmupiso310e}
\opsi_\delta\colon 
\xymatrix{
{\bMIC^\coh_\delta(\cA[\frac 1 p]/\co_\fX[\frac 1 p])}\ar[r]&{\bIMIC^\coh_{\delta,\mQ}(\cA/\co_\fX)}\ar[r]^-(0.5){\psi_{\delta,\mQ}}&
{\bIMIC_{\delta,\mQ}(\cA/\co_\fX).}}
\end{equation}

We denote by $\bIndMIC_\delta(\cA/\co_\fX)$ the category of ind-$\cA$-modules with integrable $\delta$-connection \eqref{p1-indmal22}.
The functor $\upalpha_{\cA}$ \eqref{p2-cmupiso1c} induces a functor that we denote by 
\begin{equation}\label{p2-cmupiso310f}
\upalpha_{\delta}\colon 
\begin{array}[t]{clcr}
\bIMIC_{\delta,\mQ}(\cA/\co_\fX)&\rightarrow& \bIndMIC_\delta(\cA/\co_\fX)\\
(\cM,\cN,u,\nabla)&\mapsto& (\upalpha_{\cA}(\cM_\mQ),(\id \otimes \upalpha_{\cA}(u_\mQ)^ {-1})\circ\upalpha_{\co_{\fX}}(\nabla_\mQ)).
\end{array}
\end{equation}
We denote by $\upalpha^\coh_{\delta}$ the composed functor
\begin{equation}\label{p2-cmupiso310g}
\upalpha^\coh_{\delta}=\upalpha_{\delta}\circ \opsi_\delta\colon 
\bMIC^\coh_\delta(\cA[\frac 1 p]/\co_\fX[\frac 1 p])\rightarrow \bIndMIC_\delta(\cA/\co_\fX).
\end{equation}

\begin{prop}\label{p2-cmupiso311}
We keep the assumptions and notation of \ref{p2-cmupiso310}. Then, 
\begin{itemize}
\item[{\rm (i)}] The functor $\upalpha_{\delta}$ \eqref{p2-cmupiso310f} is fully faithful.
\item[{\rm (ii)}] The functor \eqref{p2-cmupiso310c} induces an equivalence of categories
\begin{equation}\label{p2-cmupiso311a}
\bIMIC^\coh_{\delta,\mQ}(\cA/\co_\fX)\stackrel{\sim}{\rightarrow}
\bMIC^\coh_\delta(\cA[\frac 1 p]/\co_\fX[\frac 1 p]).
\end{equation}
\end{itemize}
\end{prop}

(i) Let $(\cM,\cN,u,\nabla)$, $(\cM',\cN',u',\nabla')$ be two $\delta$-isoconnections with respect to the extension $\cA/\co_\fX$. Since 
$u_\mQ$ is an isomorphism of $\bMod_\mQ(\cA)$, we can consider the composed morphism of $\bMod_\mQ(\co_\fX)$
\begin{equation}\label{p2-cmupiso311b}
\onabla=(u^{-1}_\mQ\otimes \id_{\cE})\circ \nabla_\mQ\colon \cM_\mQ\rightarrow \cM_\mQ\otimes_{\cA_\mQ}\cE_\mQ.
\end{equation}
We define in the same way the morphism $\onabla'\colon \cM'_\mQ\rightarrow \cM'_\mQ\otimes_{\cA_\mQ}\cE_\mQ$.
The morphism
\begin{equation}\label{p2-cmupiso311c}
\begin{array}[t]{clcr}
\Hom_{\bIMIC_\delta(\cA/\co_\fX)}((\cM,\cN,u,\nabla),(\cM',\cN',u',\nabla'))&\rightarrow &\Hom_\cA(\cM,\cM') \\
(m,n)&\mapsto &m
\end{array}
\end{equation}
then induces an isomorphism
\begin{eqnarray}\label{p2-cmupiso311d}
\lefteqn{\Hom_{\bIMIC_{\delta,\mQ}(\cA/\co_\fX)}((\cM,\cN,u,\nabla),(\cM',\cN',u',\nabla'))\stackrel{\sim}{ \rightarrow}}\\
&& \{m\in \Hom_\cA(\cM,\cM')\otimes_\mZ\mQ, \ (m\otimes \id_\cE)\circ \onabla=\onabla'\circ m\}.\nonumber
\end{eqnarray}
Indeed, since $u_\mQ$ and $u'_\mQ$ are isomorphisms of $\bMod_\mQ(\cA)$, the morphism \eqref{p2-cmupiso311c} induces
an injective morphism
\begin{equation}\label{p2-cmupiso311e}
\Hom_{\bIMIC_{\delta,\mQ}(\cA/\co_\fX)}((\cM,\cN,u,\nabla),(\cM',\cN',u',\nabla'))\rightarrow \Hom_\cA(\cM,\cM')\otimes_\mZ\mQ.
\end{equation}

We claim that every $m\in \Hom_\cA(\cM,\cM')\otimes_\mZ\mQ$ such that $(m\otimes \id_\cE)\circ \onabla=\onabla'\circ m$,
is in the image of the morphism \eqref{p2-cmupiso311d}. Indeed, since the latter is $\mQ$-linear,
we may assume that $m \in \Hom_\cA(\cM,\cM')$. There exist two nonzero integers $n$ and $n'$ and two $\cA$-linear morphisms
$v\colon \cN\rightarrow \cM$ and $v'\colon \cN'\rightarrow \cM'$ such that $v\circ u=n\cdot \id_\cM$, $u\circ v=n\cdot \id_\cN$, 
$v'\circ u'=n'\cdot \id_{\cM'}$, $u'\circ v'=n'\cdot \id_{\cN'}$, and that 
the $\cA$-module with $(n\delta)$-connection $(\cM,(v\otimes \id_\cE)\circ \nabla)$ (resp. $(n'\delta)$-connection $(\cM',(v'\otimes \ id_\cE)\circ \nabla')$)
is integrable in the sense of \ref{p1-delta-con2}(i). 

After multiplying $n$ and $n'$ if necessary, we may assume that $n=n'$ and that the diagram
\begin{equation}\label{p2-cmupiso311f}
\xymatrix{
\cM\ar[r]^-(0.5)\nabla\ar[d]_m&{\cN\otimes_\cA\cE}\ar[r]^-(0.5){v\otimes \id_\cE}&{\cM\otimes_\cA\cE}\ar[d]^{m\otimes\id_\cE}\\
\cM'\ar[r]^-(0.5){\nabla'}&{\cN'\otimes_\cA\cE}\ar[r]^-(0.5){v'\otimes \id_\cE}&{\cM'\otimes_\cA\cE}}
\end{equation}
is commutative. Hence, $(m,m)$ defines a morphism of $\bIMIC_{n\delta}(\cA/\co_\fX)$ from $(\cM,\cM,\id,(v\otimes\id_\cE)\circ \nabla)$ into
$(\cM',\cM',\id,(v'\otimes\id_\cE)\circ \nabla')$. Since the morphisms
\begin{eqnarray}
(\id_\cM,v)\colon (\cM,\cN,u,\nabla)&\rightarrow& (\cM,\cM,n\id,(v\otimes\id_\cE)\circ \nabla),\label{p2-cmupiso311g}\\
(\id_{\cM'},v')\colon (\cM',\cN',u',\nabla')&\rightarrow& (\cM',\cM',n'\id,(v'\otimes\id_\cE)\circ \nabla'),\label{p2-cmupiso311h}
\end{eqnarray}
are isogenies of $\bMIC_\delta(\cA/\co_\fX)$, we deduce that the morphism \eqref{p2-cmupiso311d} is surjective and therefore bijective.
Observe that the inverses of \eqref{p2-cmupiso311g} (resp. \eqref{p2-cmupiso311h}) up to multiplication
by $n$ (resp. $n'$) is given by $(n\id_\cM,u)$ (resp. $(n'\id_{\cM'},u')$).

It then follows from \eqref{p2-cmupiso311d} and (\cite{ag2} (2.6.6.4)) that the functor $\upalpha_{\delta}$ \eqref{p2-cmupiso310f} is fully faithful.

(ii) The proof is similar to that of (\cite{agt} III.6.21) using \ref{p1-pfs15}, \ref{p1-pfs12} and \ref{p1-abisoind40}.

\subsection{}\label{p2-cmupiso18}
Let $f\colon \fX'\rightarrow \fX$ be a morphism of finite presentation \eqref{p2-cmupiso1}. 
We also denote by $f\colon (\fX'_\zar,\co_{\fX'})\rightarrow (\fX_\zar,\co_\fX)$ the morphism of ringed topos induced by $f$, for which 
we will use the notation/conventions of \ref{p1-bcim6} and \ref{p1-bcim1}. It follows from \eqref{p1-bcim6d} that the diagram 
\begin{equation}\label{p2-cmupiso18a}
\xymatrix{
{\bD^+(\bMod_\mQ(\co_{\fX'}))}\ar[d]\ar[r]^{\rR f_{\mQ*}}&{\bD^+(\bMod_\mQ(\co_{\fX}))}\ar[d]\\
{\bD^+(\bMod(\co_{\fX'}[\frac 1 p]))}\ar[r]^{\rR f_*}&{\bD^+(\bMod(\co_{\fX}[\frac 1 p])),}}
\end{equation}
where the vertical arrows are the canonical functors, is commutative up to canonical isomorphism. 
Let $M$ be a coherent $\co_{\fX'}[\frac 1 p]$-module, $q$ an integer $\geq 0$. If the $\co_\fX[\frac 1 p]$-module $\rR^qf_*(M)$ is coherent, 
diagram \eqref{p2-cmupiso18a} induces a canonical functorial isomorphism 
\begin{equation}\label{p2-cmupiso18b}
\ojmath_{\co_{\fX}}(\rR^qf_*(M))\stackrel{\sim}{\rightarrow} \rR^qf_{\mQ*}(\ojmath'_{\co_{\fX'}}(M)),
\end{equation}
where $\ojmath_{\co_{\fX}}$ is the functor \eqref{p2-cmupiso2b} and $\ojmath'_{\co_{\fX'}}$ is the analogous functor for $\co_{\fX'}$.

\subsection{}\label{p2-cmupiso10}
Let $\cU$ be a Zariski open subscheme of $\fX$, $j\colon \cU\rightarrow \fX$ the canonical open immersion. 
By \ref{p1-indmal16}, $j$ induces three functors  
\begin{eqnarray}
\rI j_!\colon \bIndMod(\co_\cU)&\rightarrow& \bIndMod(\co_\fX),\label{p2-cmupiso10a}\\
\rI j^*\colon \bIndMod(\co_\fX)&\rightarrow& \bIndMod(\co_\cU),\label{p2-cmupiso10b}\\
\rI j_*\colon \bIndMod(\co_\cU)&\rightarrow& \bIndMod(\co_\fX).\label{p2-cmupiso10c} 
\end{eqnarray}
The functor $\rI j_*$ (resp.\ $\rI j_!$) is a right (resp.\ left) adjoint of the functor $\rI j^*$.
Therefore, the functor $\rI j^*$ commutes with representable direct and inverse limits and is in particular exact.
For any ind-$\co_\fX$-module $M$, the ind-$\co_\cU$-module $\rI j^*(M)$ will also be denoted by $M|U$. 

The diagrams of functors 
\begin{equation}\label{p2-cmupiso10d}
\xymatrix{
{\bMod(\co_{\fX})}\ar[r]^-(0.5){\iota_{\co_\fX}}\ar[d]_{j^*}&{\bIndMod(\co_{\fX})}\ar[d]^{\rI j^*}\ar[r]^-(0.5){\kappa_{\co_\fX}}&
{\bMod(\co_{\fX})}\ar[d]_{j^*}\ar[r]^-(0.5){Q_{\co_\fX}}&{\bMod_\mQ(\co_{\fX})}\ar[d]^{j^*_\mQ}\\
{\bMod(\co_{\cU})}\ar[r]^-(0.5){\iota_{\co_\cU}}&{\bIndMod(\co_{\cU})}\ar[r]^-(0.5){\kappa_{\co_\cU}}&{\bMod(\co_{\cU})}
\ar[r]^-(0.5){Q_{\co_\cU}}&{\bMod_\mQ(\co_{\cU}),}}
\end{equation}
\begin{equation}\label{p2-cmupiso10e}
\xymatrix{
{\bMod_\mQ(\co_{\fX})}\ar[r]^-(0.5){\upalpha_{\co_\fX}}\ar[d]_{j^*_\mQ}&{\bIndMod(\co_{\fX})}\ar[d]^{\rI j^*}\\
{\bMod_\mQ(\co_{\cU})}\ar[r]^-(0.5){\upalpha_{\co_\cU}}&{\bIndMod(\co_{\cU}),}}
\ \ \
\xymatrix{
{\bMod^\coh(\co_{\fX}[\frac{1}{p}])}\ar[r]^-(0.5){\upalpha^\coh_{\co_\fX}}\ar[d]_{j^*}&{\bIndMod(\co_{\fX})}\ar[d]^{\rI j^*}\\
{\bMod^\coh(\co_{\cU}[\frac{1}{p}])}\ar[r]^-(0.5){\upalpha^\coh_{\co_\cU}}&{\bIndMod(\co_{\cU}),}}
\end{equation}
are commutative up to canonical isomorphisms.

\subsection{}\label{p2-cmupiso8}
Let $I$ be a small filtered category, $\rho\colon I\rightarrow \bMod(\co_\fX)$ a functor. 
We set, for any $i\in \ob(I)$, $\cM_i=\rho(i)$,
\begin{equation}\label{p2-cmupiso8a}
M=\underset{\underset{I}{\longrightarrow}}{\mlq\mlq\lim \mrq\mrq} \ \upalpha_{\co_\fX}(\cM_{i,\mQ}) \ \ \ {\rm and}\ \ \ 
M_\infty=\underset{\underset{i\in I}{\longrightarrow}}{\lim}\ \cM_{i}[\frac 1 p].
\end{equation}
We denote by $\psi_i\colon \upalpha_{\co_\fX}(\cM_{i,\mQ})\rightarrow M$ and $\varphi_i\colon \cM_{i}[\frac 1 p]\rightarrow M_\infty$ the canonical morphisms. 
By \eqref{p1-abisoind3d}, we have a canonical isomorphism 
\begin{equation}\label{p2-cmupiso8b}
\kappa_{\co_\fX}(M)\stackrel{\sim}{\rightarrow} M_\infty.
\end{equation}

We assume that there exists a Zariski open covering $(\cU_\lambda)_{\lambda\in \Lambda}$ of $\fX$ such that for every $\lambda\in \Lambda$, 
$M|\cU_\lambda$ \eqref{p2-cmupiso10} belongs to the essential image of the functor 
\begin{equation}\label{p2-cmupiso8c}
\upalpha^\coh_{\co_{\cU_\lambda}}\colon\bMod^{\coh}(\co_{\cU_\lambda}[\frac 1 p])\rightarrow \bIndMod(\co_{\cU_\lambda}).
\end{equation}
We deduce by \eqref{p2-cmupiso2d}, \eqref{p2-cmupiso10d} and \eqref{p2-cmupiso8b} that the $\co_\fX[\frac 1 p]$-module $M_\infty$ is coherent. 
Let $\cM_\infty$ be a coherent $\co_\fX$-module, $\mu_\infty \colon M_\infty\stackrel{\sim}{\rightarrow}\cM_\infty[\frac 1 p]$ an $\co_\fX$-linear isomorphism.
We denote by $\varrho_\infty$ the constant functor from $I$ to $\bMod_\mQ(\co_\fX)$ with value $\cM_{\infty,\mQ}$. 

\begin{prop}\label{p2-cmupiso9}
We keep the assumptions and notation of \ref{p2-cmupiso8}. Then, 
\begin{itemize}
\item[{\rm (i)}] For every $i\in \ob(I)$, the canonical morphism
\begin{equation}\label{p2-cmupiso9a}
\Hom_{\co_\fX}(\cM_i,\cM_\infty)\otimes_{\mZ}\mQ\rightarrow \Hom_{\co_\fX}(\cM_i[\frac 1 p],M_\infty) 
\end{equation}
is injective, and its image contains the canonical morphism $\varphi_i\colon \cM_{i}[\frac 1 p]\rightarrow M_\infty$ \eqref{p2-cmupiso8a}.  
We denote by $\phi_i \colon  \cM_{i,\mQ}\rightarrow \cM_{\infty,\mQ}$ the inverse image of $\varphi_i$. 
The collection of these morphisms define a morphism of functors $\phi_\bullet\colon Q_{\co_\fX}\circ \rho\rightarrow Q_{\co_\fX}\circ \varrho_\infty$, 
where $Q_{\co_\fX}$ is the functor \eqref{p2-cmupiso1d}.
\item[{\rm (ii)}] The morphism $\upalpha_{\co_\fX}(\phi_\bullet)$ induces an isomorphism 
\begin{equation}\label{p2-cmupiso9b}
M \stackrel{\sim}{\rightarrow} \upalpha_{\co_\fX}(\cM_{\infty,\mQ}). 
\end{equation}
\end{itemize}
\end{prop}

(i) The kernel of the canonical morphism $\cM_\infty\rightarrow M_\infty$ is a coherent $\co_\fX$-module 
(\cite{egr1} 2.10.14), and is hence annihilated by a non-zero integer. Therefore, the morphism \eqref{p2-cmupiso9a} is injective. 
To prove the proposition, we may replace $\cM_\infty$ by its image in $M_\infty$. We are then reduced to proving that there exists an integer $n_i\geq 1$ 
such that $n_i \gamma_i(\cM_i)\subset \cM_\infty$, where $\gamma_i\colon \cM_{i}\rightarrow M_\infty$ is the morphism induced by $\varphi_i$. 

This assertion being local on $\fX$, we may assume that there exists a coherent $\co_\fX$-module $\cM$ 
and an isomorphism $\lambda\colon \upalpha_{\co_\fX}(\cM_{\mQ})\stackrel{\sim}{\rightarrow} M$. 
By construction, the isomorphism \eqref{p2-cmupiso8b} fits into a commutative diagram 
\begin{equation}
\xymatrix{
&{\cM_i[\frac 1 p]}\ar[r]^-(0.5)\sim\ar[d]_{\varphi_i}&{\kappa_{\co_\fX}(\upalpha_{\co_\fX}(\cM_{i,\mQ}))}\ar[d]^{\kappa_{\co_\fX(\psi_i)}}\\
{\cM_\infty[\frac 1 p]}&{M_\infty}\ar[l]_-(0.5){\mu_\infty}^-(0.5)\sim&{\kappa_{\co_\fX}(M),}\ar[l]_-(0.5)\sim}
\end{equation}
where the upper horizontal arrow is defined in \eqref{p1-abisoind3d}. 
By the full faithfulness of $\upalpha_{\co_\fX}$, there exists a unique morphism $\nu_i \colon  \cM_{i,\mQ}\rightarrow \cM_{\mQ}$ such that 
$\psi_i=\lambda\circ \upalpha_{\co_\fX}(\nu_i)$. The previous diagram implies then that the diagram 
\begin{equation}
\xymatrix{
&{\cM_i[\frac 1 p]}\ar[d]_{\varphi_i}\ar[r]^-(0.5){\nu_i[\frac 1 p]}&{\cM[\frac 1 p]}\ar[r]^-(0.5)\sim&{\kappa_{\co_\fX}(\upalpha_{\co_\fX}(\cM_{\mQ}))}
\ar[d]^{\kappa_{\co_\fX}(\lambda)}\\
{\cM_\infty[\frac 1 p]}&{M_\infty}\ar[l]_-(0.5){\mu_\infty}^-(0.5)\sim&&{\kappa_{\co_\fX}(M)}\ar[ll]_-(0.5)\sim}
\end{equation}
is commutative. 
Since $\cM$ and $\cM_\infty$ are coherent, the composite isomorphism $\cM[\frac 1 p]\stackrel{\sim}{\rightarrow}\cM_\infty[\frac 1 p]$ appearing 
above is of the form $\mu[\frac 1 p]$ for a unique isomorphism $\mu\colon \cM_\mQ\rightarrow\cM_{\infty,\mQ}$. This completes the proof. 

(ii) The question being local on $\fX$ by (\cite{ag2} 2.7.17(iii)), we may assume that there exists a coherent $\co_\fX$-module $\cM$ 
and an isomorphism $\lambda\colon \upalpha_{\co_\fX}(\cM_{\mQ})\stackrel{\sim}{\rightarrow} M$. 
Then the morphism $\phi_\infty\colon M\rightarrow \upalpha_{\co_\fX}(\cM_{\infty,\mQ})$ induced by $\upalpha_{\co_\fX}(\phi_\bullet)$
is the unique morphism satisfying $\phi_\infty\circ \psi_i=\upalpha_{\co_\fX}(\phi_i)$ for every $i\in \ob(I)$. 
Since $\psi_i=\lambda\circ \upalpha_{\co_\fX}(\nu_i)$, we see that 
$\phi_\infty \circ \lambda$ is the unique morphism satisfying $(\phi_\infty\circ \lambda)\circ \upalpha_{\co_\fX}(\nu_i)= \upalpha_{\co_\fX}(\phi_i)$
for every $i\in \in \ob(I)$. 
Since $\phi_i=\mu\circ \nu_i$ by the proof of (i), we deduce that $\phi_\infty\circ \lambda=\upalpha_{\co_\fX}(\mu)$. 
Since $\lambda$ and $\mu$ are isomorphisms, so is $\phi_\infty$, which proves the proposition.

\begin{cor}\label{p2-cmupiso11}
Let $M$ be an ind-$\co_\fX$-module, $(\cU_\lambda)_{\lambda\in \Lambda}$ a Zariski open covering of $\fX$, 
such that for every $\lambda\in \Lambda$, $M|\cU_\lambda$ \eqref{p2-cmupiso10} belongs to the essential image of the functor 
\begin{equation}\label{p2-cmupiso11a}
\upalpha^\coh_{\co_{\cU_\lambda}}\colon\bMod^{\coh}(\co_{\cU_\lambda}[\frac 1 p])\rightarrow \bIndMod(\co_{\cU_\lambda}).
\end{equation}
Then, $M$ belongs to the essential image of the functor $\upalpha^\coh_{\co_{\fX}}$ \eqref{p2-cmupiso2c}. 
\end{cor}

Indeed, let $I$ be a small filtered category, $\rho\colon I\rightarrow \bMod(\co_\fX)$ a functor, such that 
\begin{equation}
M=\underset{\underset{I}{\longrightarrow}}{\mlq\mlq\lim \mrq\mrq} \rho.
\end{equation}
By (\cite{ag2} 2.7.17(iii)), for every integer $n\geq 1$, the multiplication by $n$ on $M$ is an isomorphism. We deduce an isomorphism 
\begin{equation}
M\stackrel{\sim}{\rightarrow}\underset{\underset{I}{\longrightarrow}}{\mlq\mlq\lim \mrq\mrq} \upalpha_{\co_\fX} \circ \rho.
\end{equation}
The proposition follows then from \ref{p2-cmupiso9}.

\begin{prop}\label{p2-cmupiso12}
Let $f\colon \fX'\rightarrow \fX$ be a proper morphism of finite presentation. Then, 
\begin{itemize}
\item[{\rm (i)}] For every coherent $\co_{\fX'}[\frac 1 p]$-module $M$ and every integer $q$, the $\co_\fX[\frac 1 p]$-module $\rR^qf_*(M)$ is coherent, 
and we have a canonical functorial isomorphism 
\begin{equation}\label{p2-cmupiso12a}
\upalpha_{\co_\fX}^\coh(\rR^qf_*(M))\stackrel{\sim}{\rightarrow} \rR^q\rI f_*(\upalpha_{\co_{\fX'}}^\coh(M)).
\end{equation} 
\item[{\rm (ii)}] For every bounded from below complex $M^\bullet$ of ind-$\co_{\fX'}$-modules such that for every integer $i$, 
$\cH^i(M^\bullet)$ is in the essential image of the functor $\upalpha_{\co_{\fX'}}^\coh$,
$\rR^q\rI f_*(M^\bullet)$ is in the essential image of $\upalpha_{\co_\fX}^\coh$ \eqref{p2-cmupiso2c}, for every integer $q$. 
\end{itemize}
\end{prop}

(i) Let $\cM$ be a coherent $\co_{\fX'}$-module such that $M=\cM[\frac 1 p]$. By (\cite{egr1} 2.11.5), the $\co_\fX$-module $f_*(\cM)$ is coherent. 
Therefore, the $\co_\fX[\frac 1 p]$-module $\rR^qf_*(M)=\rR^qf_*(\cM)[\frac 1 p]$ is coherent. 
We denote by $\Delta$ the multiplicative monoid $\mZ-\{0\}$ and by $\uDelta$
the filtered category whose objects are the elements of $\Delta$ and the morphisms are determined by the divisibility relation in $\Delta$.
Let $\upmu\colon \uDelta \rightarrow \bMod(\co_\fX)$ be the functor which sends a nonzero integer $n$ to $\cM$
and a morphism $m|n$ of $\uDelta$ to the isogeny of multiplication by $n/m$ on $\cM$. 
By (\cite{ag2} 2.9.3(i)), we have canonical isomorphisms
\begin{eqnarray}
\rR^q\rI f_*(\upalpha_{\co_{\fX'}}^\coh(M))&\stackrel{\sim}{\rightarrow} & 
\rR^q\rI f_*(\underset{\underset{n\in \uDelta}{\longrightarrow}}{\mlq\mlq\lim \mrq\mrq} \cM)\\
&\stackrel{\sim}{\rightarrow} & \underset{\underset{n\in \uDelta}{\longrightarrow}}{\mlq\mlq\lim \mrq\mrq} \rR^q f_*(\cM)\\
&\stackrel{\sim}{\rightarrow} & \upalpha_{\co_\fX}^\coh(\rR^q f_*(M)). 
\end{eqnarray}

(ii) It follows from (i), \ref{p2-cmupiso3} and \ref{p1-bcim16}.

\begin{prop}\label{p2-cmupiso13}
Let $f\colon \fX'\rightarrow \fX$ be a morphism of finite presentation, 
$I$ a small filtered category, $\varphi\colon I\rightarrow \bC^+(\bMod(\co_{\fX'}))$ a functor. 
For any $i\in \ob(I)$, we set $\cM^\bullet_i=\varphi(i)$. 
We assume that there exists an integer $s$ such that for every $r\leq s$ and every $i\in \ob(I)$, $\cM^r_i=0$. 
We set $\cM^\bullet=\underset{\underset{i\in I}{\longrightarrow}}{\mlq\mlq\lim \mrq\mrq} \ \cM^\bullet_i \in \bC^+(\bIndMod(\co_{\fX'}))$, 
where the direct limit is computed term by term. Then, 
\begin{itemize}
\item[{\rm (i)}] For every integer $q$, we have a canonical isomorphism 
\begin{equation}
\rR^qf_*(\kappa_{\co_{\fX'}}(\cM^\bullet))\stackrel{\sim}{\rightarrow} \kappa_{\co_{\fX}} (\rR^q\rI f_*(\cM^\bullet)).
\end{equation}
\item[{\rm (ii)}] Assume that $f$ is proper and that for every integer $q$, $\cH^q(\cM^\bullet)$ is in the essential image of the functor
$\upalpha_{\co_{\fX'}}^\coh$. Then, for every integer $q$,  $\rR^qf_*(\kappa_{\co_{\fX'}}(\cM^\bullet))$ is a coherent $\co_\fX[\frac 1 p]$-module 
and we have a canonical isomorphism
\begin{equation}
\rR^q\rI f_*(\cM^\bullet)\stackrel{\sim}{\rightarrow} \upalpha^\coh_{\co_{\fX}} (\rR^qf_*(\kappa_{\co_{\fX'}}(\cM^\bullet))).
\end{equation}
\end{itemize}
\end{prop}

(i) Indeed, the morphism of topos $f\colon \fX'_\zar\rightarrow \fX_\zar$ being coherent (\cite{sga4} VI 3.3), $\rR^qf_*$ 
commutes with filtered direct limits of abelian sheaves by (\cite{sga4} VI 5.1). 
The proposition follows then from \ref{p1-bcim18}(iii) and (\cite{sga4} II 4.1 and 4.3).

(ii) It follows from (i), \ref{p2-cmupiso12}(ii) and \eqref{p2-cmupiso2d}.

\subsection{}\label{p2-cmupiso5}
We denote by $\fX_\zar^{\mN^\circ}$ the topos of inverse systems of $\fX_\zar$,
indexed by the ordered set $\mN$ of natural numbers \eqref{p2-ncgt5}, and by $\co_{\bvfX}$ the ring $(\co_{\fX_{n+1}})_{n\in \mN}$ of $\fX_\zar^{\mN^\circ}$,
where $\co_{\fX_n}=\co_{\fX}/p^n\co_{\fX}$. We have a morphism of topos
\begin{equation}\label{p2-cmupiso5a}
\uplambda\colon \fX_\zar^{\mN^\circ}\rightarrow \fX_\zar,
\end{equation}
whose direct image functor $\uplambda_*$ associates with any inverse system its inverse limit \eqref{p2-ncgt5a}. It is underlying a canonical 
morphism of ringed topos that we also (abusively) denote by 
\begin{equation}\label{p2-cmupiso5b}
\uplambda\colon (\fX_\zar^{\mN^\circ},\co_\bvfX)\rightarrow (\fX_\zar,\co_\fX). 
\end{equation}
For modules, we use the notation $\uplambda^{-1}$ 
to denote the  pullback in the sense of abelian sheaves, and we keep the notation
$\uplambda^*$ for the pullback in the sense of  modules.

For any $\co_\fX$-module $\cF$, we set $\bvcF=(\cF/p^{n+1}\cF)_{n\in \mN}$. 
By (\cite{agt} III.7.12), we have a canonical functorial isomorphism 
\begin{equation}\label{p2-cmupiso5c}
\uplambda^*(\cF)\stackrel{\sim}{\rightarrow} \bvcF.
\end{equation}

We denote by $\bMod(\co_{\bvfX})$ the category of $\co_{\bvfX}$-modules of $\fX^{\mN^\circ}_{\zar}$, 
by $\bIndMod(\co_\bvfX)$ the category of ind-$\co_\bvfX$-modules  and by
\begin{equation}
\iota_{\co_\bvfX}\colon \bMod(\co_\bvfX)\rightarrow \bIndMod(\co_\bvfX)
\end{equation}
the canonical functor (see \ref{p1-indmal1}).
Following \ref{p1-bcim1}, we denote by 
\begin{eqnarray}
\rI\uplambda_*\colon \bIndMod(\co_\bvfX)\rightarrow \bIndMod(\co_\fX),\label{p2-cmupiso5d1}\\
\rI\uplambda^*\colon \bIndMod(\co_\fX)\rightarrow \bIndMod(\co_\bvfX),\label{p2-cmupiso5d2}
\end{eqnarray}
the adjoint functors induced by $\uplambda$.

\begin{prop}\label{p2-cmupiso6}
The composed functor 
\begin{equation}\label{p2-cmupiso6a}
\xymatrix{
{\bMod^{\coh}(\co_{\fX}[\frac 1 p])}\ar[r]^-(0.5){\ojmath_{\co_\fX}}&{\bMod_\mQ(\co_\fX)}\ar[r]^-(0.5){\uplambda^*_\mQ}&{\bMod_\mQ(\co_\bvfX),}} 
\end{equation}
is exact.
\end{prop}

By (\cite{egr1} 2.10.24(iii)), since $\fX$ is quasi-compact, it is enough to prove that for any exact sequence of coherent $\co_\fX$-modules 
\begin{equation}\label{p2-cmupiso6b}
\xymatrix{
0\ar[r]&{\cM'}\ar[r]&{\cM}\ar[r]^\varphi&{\cM''}\ar[r]&0,}
\end{equation}
the sequence of $\co_{\bvfX,\mQ}$-modules 
\begin{equation}\label{p2-cmupiso6c}
\xymatrix{
0\ar[r]&{\bvcM'_\mQ}\ar[r]&{\bvcM_\mQ}\ar[r]^\varphi&{\bvcM''_\mQ}\ar[r]&0}
\end{equation}
is exact \eqref{p2-cmupiso5c}. Let $\cM''_\tor$ (resp.\ $\ocM''$) be the kernel (resp.\ image) of the canonical morphism $\cM''\rightarrow \cM''[\frac 1 p]$. 
Consider the commutative diagram 
\begin{equation}\label{p2-cmupiso6d}
\xymatrix{
0\ar[r]&{\cM'}\ar[r]\ar[d]&{\cM}\ar[r]^\varphi\ar@{=}[d]&{\cM''}\ar[r]\ar[d]&0\\
0\ar[r]&{\cM'_1}\ar[r]&{\cM}\ar[r]^\ovarphi&{\ocM''}\ar[r]&0,}
\end{equation}
where the vertical arrows are the canonical morphisms and $\cM'_1$ is the kernel of $\ovarphi$. 
By (\cite{egr1} 2.10.14), the $\co_\fX$-module $\cM''_\tor$ is coherent and is hence annihilated by a non-zero integer. 
Therefore, the $\co_\fX$-modules $\ocM''$ and $\cM'_1$ are coherent, and
the vertical morphisms of \eqref{p2-cmupiso6d} induce isomorphisms of $\bMod_\mQ(\co_\fX)$. 
We are hence reduced to the case where 
$\cM''$ is $p$-torsion free, in which case the sequence of $\co_\bvfX$-modules
\begin{equation}
\xymatrix{
0\ar[r]&{\bvcM'}\ar[r]&{\bvcM}\ar[r]^\bvvarphi&{\bvcM''}\ar[r]&0}
\end{equation}
is exact, which implies the proposition. 

\begin{cor}\label{p2-cmupiso7}
The composed functor 
\begin{equation}
\xymatrix{
{\bMod^{\coh}(\co_{\fX}[\frac 1 p])}\ar[r]^-(0.5){\upalpha_{\co_\fX}^\coh}&{\bIndMod(\co_\fX)}\ar[r]^-(0.5){\rI\uplambda^*}&{\bIndMod(\co_\bvfX)}} 
\end{equation}
is exact.
\end{cor}

Indeed, the diagram 
\begin{equation}
\xymatrix{
{\bMod_\mQ(\co_{\fX})}\ar[r]^-(0.5){\upalpha_{\co_\fX}}\ar[d]_{\uplambda^*_{\mQ}}&{\bIndMod(\co_\fX)}\ar[d]^-(0.5){\rI\uplambda^*}\\
{\bMod_\mQ(\co_{\bvfX})}\ar[r]^-(0.5){\upalpha_{\bvco_\fX}}&{\bIndMod(\co_\bvfX)}} 
\end{equation}
being commutative up to isomorphism, the proposition follows from \ref{p2-cmupiso6} and the exactness of the functor $\upalpha_{\co_\bvfX}$
(\cite{ag2} 2.9.3(ii)).

\begin{prop}\label{p2-cmupiso20}
Let $M^\bullet$ be a bounded from below complex of ind-$\co_\fX$-modules, 
$(\cU_\lambda)_{\lambda\in \Lambda}$ a Zariski open covering of $\fX$. 
We assume that for every $\lambda\in \Lambda$, 
there exists a bounded from below complex of coherent $\co_{\cU_\lambda}[\frac 1 p]$-modules $N_\lambda^\bullet$ and a morphism
\begin{equation}\label{p2-cmupiso20a}
u_\lambda\colon \upalpha^\coh_{\co_{\cU_\lambda}}(N_\lambda^\bullet)\rightarrow M^\bullet|\cU_\lambda,
\end{equation}
where both the source and the target are defined term by term, whose mapping cone is homotopically equivalent to zero. Then, 
for every integer $q$, $\cH^q(M^\bullet)$ is in the essential image of $\upalpha^\coh_{\co_{\fX}}$ and the canonical morphism \eqref{p1-bcim13a}
\begin{equation}\label{p2-cmupiso20b}
\rI \uplambda^*(\cH^q(M^\bullet))\rightarrow \cH^q(\rI \uplambda^*M^\bullet)
\end{equation}
is an isomorphism.  
\end{prop}

By \ref{p2-cmupiso11} and (\cite{ag2} 2.7.17(iii)), 
we may assume that the covering $(\cU_\lambda)_{\lambda\in \Lambda}$ is reduced to $\fX$; 
so there exists a bounded from below complex of coherent $\co_{\fX}[\frac 1 p]$-modules $N^\bullet$ and a morphism
\begin{equation}\label{p2-cmupiso20c}
u\colon \upalpha^\coh_{\co_{\fX}}(N^\bullet)\rightarrow M^\bullet
\end{equation}
whose mapping cone $C^\bullet$ is homotopically equivalent to zero. 
Since $\upalpha^\coh_{\co_{\fX}}$ is exact, we have a canonical isomorphism 
\begin{equation}\label{p2-cmupiso20d}
\cH^q(\upalpha^\coh_{\co_{\fX}}(N^\bullet))\rightarrow \upalpha^\coh_{\co_{\fX}}(\cH^q(N^\bullet)),
\end{equation}
which implies the first assertion. We have canonical distinguished triangles 
\begin{eqnarray}
\upalpha^\coh_{\co_{\fX}}(N^\bullet)\rightarrow M^\bullet\rightarrow C^\bullet \rightarrow \upalpha^\coh_{\co_{\fX}}(N^\bullet)[1],\\
\rI \uplambda^*(\upalpha^\coh_{\co_{\fX}}(N^\bullet))\rightarrow \rI \uplambda^*(M^\bullet)\rightarrow 
\rI \uplambda^*(C^\bullet) \rightarrow \rI \uplambda^*(\upalpha^\coh_{\co_{\fX}}(N^\bullet))[1],
\end{eqnarray}
the second being the image of the first one by the functor $\rI \uplambda^*$. 
They induce a commutative diagram 
\begin{equation}
\xymatrix{
{\rI \uplambda^*(\cH^q(\upalpha^\coh_{\co_{\fX}}(N^\bullet)))}\ar[r]\ar[d]&{\rI \uplambda^*(\cH^q(M^\bullet))}\ar[r]\ar[d]&
{\rI \uplambda^*(\cH^q(C^\bullet))}\ar[d]\ar[r]&{\rI \uplambda^*(\cH^q(\upalpha^\coh_{\co_{\fX}}(N^\bullet)))[1]}\ar[d]\\
{\cH^q(\rI \uplambda^*(\upalpha^\coh_{\co_{\fX}}(N^\bullet)))}\ar[r]&{\cH^q(\rI \uplambda^*(M^\bullet))}\ar[r]&
{\cH^q(\rI \uplambda^*(C^\bullet))}\ar[r]&{\cH^q(\rI \uplambda^*(\upalpha^\coh_{\co_{\fX}}(N^\bullet)))[1],}}
\end{equation}
where the vertical arrows are the morphisms \eqref{p1-bcim13a}. By assumption, $\cH^i(C^\bullet)$ and $\cH^i(\rI \uplambda^*(C^\bullet))$
vanish for all $i$. We are therefore reduced to the case $M^\bullet=\upalpha^\coh_{\co_{\fX}}(N^\bullet)$. 

By \eqref{p1-bcim13d}, the diagram 
\begin{equation}
\xymatrix{
{\rI \uplambda^*(\upalpha^\coh_{\co_{\fX}}(\cH^q(N^\bullet)))}\ar[r]\ar[rd]&
{\rI \uplambda^*(\cH^q(\upalpha^\coh_{\co_{\fX}}(N^\bullet)))}\ar[d]\\
&{\cH^q(\rI \uplambda^*(\upalpha^\coh_{\co_{\fX}}(N^\bullet))),}}
\end{equation}
where all arrows are induced by the morphisms \eqref{p1-bcim13a}, is commutative. Since $\upalpha^\coh_{\co_{\fX}}$ and 
$\rI \uplambda^*\circ \upalpha^\coh_{\co_{\fX}}$ are exact by \ref{p2-cmupiso7}, the horizontal and the slanting arrows are isomorphisms. 
Hence so is the vertical arrow, which proves the second assertion.

\subsection{}\label{p2-cmupiso14}
Let $f\colon \fX'\rightarrow \fX$ be a {\em proper} morphism of finite presentation \eqref{p2-cmupiso1}. The diagram of morphisms of ringed topos 
\begin{equation}\label{p2-cmupiso14a}
\xymatrix{
{(\fX'^{\mN^\circ}_\zar,\co_{\bvfX'})}\ar[r]^-(0.5){\uplambda'}\ar[d]_{\bvf}&{(\fX'_\zar,\co_{\fX'})}\ar[d]^f\\
{(\fX^{\mN^\circ}_\zar,\co_\bvfX)}\ar[r]^-(0.5){\uplambda}&{(\fX_\zar,\co_\fX),}}
\end{equation}
where $\uplambda$ and $\uplambda'$ are the canonical morphisms \eqref{p2-cmupiso5b} and the morphism $\bvf$ is induced by $f$, 
is commutative up to canonical isomorphism. We use the notation/conventions of \ref{p1-bcim6} and \ref{p1-bcim1} for these morphisms of ringed topos. 

Let $q$ be an integer. For every bounded from below complex of $\co_{\fX'}$-modules $\cM^\bullet$, 
there exists a canonical functorial base change morphism of $\co_{\bvfX}$-modules \eqref{p1-bcim4d}, with respect to \eqref{p2-cmupiso14a},
\begin{equation}\label{p2-cmupiso14b}
\uplambda^*(\rR^qf_*(\cM^\bullet))\rightarrow \rR^q\bvf_*(\uplambda'^*(\cM^\bullet)),
\end{equation}
where the pullback $\uplambda'^*(\cM^\bullet)$ is defined term by term (not derived). 
By \ref{p1-bcim12}, for every $\co_{\fX'}$-module $\cM$, the base change morphism \eqref{p2-cmupiso14b} 
for $\cM[0]$ coincides with the classical base change morphism (\cite{egr1} 1.2.3).

For every bounded from below complex of $\co_{\fX',\mQ}$-modules $M^\bullet$, 
there exists a canonical functorial base change morphism of $\co_{\bvfX,\mQ}$-modules \eqref{p1-bcim4e}, with respect to \eqref{p2-cmupiso14a},
\begin{equation}\label{p2-cmupiso14c}
\uplambda^*_\mQ(\rR^qf_{\mQ *}(M^\bullet))\rightarrow \rR^q\bvf_{\mQ*}(\uplambda'^*_\mQ(M^\bullet)),
\end{equation}
where the pullback $\uplambda'^*_\mQ(M^\bullet)$ is defined term by term (not derived). 

For every bounded from below complex of ind-$\co_{\fX'}$-modules $\cF^\bullet$, 
there exists a canonical functorial base change morphism of ind-$\co_{\bvfX}$-modules \eqref{p1-bcim4f}, with respect to \eqref{p2-cmupiso14a},
\begin{equation}\label{p2-cmupiso14d}
\rI\uplambda^*(\rR^q\rI f_*(\cF^\bullet))\rightarrow \rR^q\rI\bvf_*(\rI\uplambda'^*(\cF^\bullet)),
\end{equation}
where the pullback $\rI\uplambda'^*(\cF^\bullet)$ is defined term by term (not derived).

\begin{prop}\label{p2-cmupiso15}
We keep the assumptions and notation of \ref{p2-cmupiso14} and let $\cM$ be a coherent $\co_{\fX'}$-module, $q$ an integer.
Then, 
\begin{itemize}
\item[{\rm (i)}]  The $\co_\fX$-module $\rR^qf_*(\cM)$ is coherent. 
\item[{\rm (ii)}] There exists an integer $N\geq 0$ such that the kernel and the cokernel of the base change morphism \eqref{p2-cmupiso14b}
\begin{equation}\label{p2-cmupiso15a}
\uplambda^*(\rR^qf_*(\cM))\rightarrow \rR^q\bvf_*(\uplambda'^*(\cM))
\end{equation}
are annihilated by $p^N$.
\end{itemize}
\end{prop}

(i) This is the statement of (\cite{egr1} 2.11.5). 

(ii) This is the statement of (\cite{ag1}, 6.5.38). Observe that the latter is stated under more restrictive conditions, 
specific to the context of {\em loc. cit.}, but the proof applies mutatis mutandis to our setting.

\begin{cor}\label{p2-cmupiso16}
We keep the assumptions and notation of \ref{p2-cmupiso14} and let $M$ be a coherent $\co_{\fX'}[\frac 1 p]$-module, $q$ an integer.
Then, 
\begin{itemize}
\item[{\rm (i)}]  The $\co_\fX[\frac 1 p]$-module $\rR^qf_*(M)$, is coherent. 
We consider $M$ (resp.\ $\rR^qf_*(M)$) as an $\co_{\fX',\mQ}$-module (resp.\ $\co_{\fX,\mQ}$-module) 
via $\ojmath_{\co_{\fX}}$ (resp.\ $\ojmath'_{\co_{\fX'}}$) \eqref{p2-cmupiso2b}. 
\item[{\rm (ii)}] The base change morphism \eqref{p2-cmupiso14c}
\begin{equation}\label{p2-cmupiso16a}
\uplambda^*_\mQ(\rR^qf_*(M))\rightarrow \rR^q\bvf_{\mQ*}(\uplambda'^*_\mQ(M))
\end{equation}
is an isomorphism.
\end{itemize}
\end{cor}

It follows from \ref{p2-cmupiso15} and \eqref{p1-bcim9b}.

\begin{prop}\label{p2-cmupiso17}
We keep the assumptions and notation of \ref{p2-cmupiso14} and 
let $\cM^\bullet$ be a bounded from below complex of ind-$\co_{\fX'}$-modules satisfying the following conditions:
\begin{itemize}
\item[{\rm (i)}] For every integer $q$, $\cH^q(\cM^\bullet)$ is in the essential image of $\upalpha_{\co_{\fX'}}^\coh$ \eqref{p2-cmupiso2c}.
\item[{\rm (ii)}] For every integer $q$, the canonical morphism \eqref{p1-bcim13a}
\begin{equation}\label{p2-cmupiso17a}
\rI \uplambda'^*(\cH^q(\cM^\bullet))\rightarrow \cH^q(\rI \uplambda'^*\cM^\bullet)
\end{equation}
is an isomorphism.  
\end{itemize}
Then, for every integer $q$, the base change morphism \eqref{p2-cmupiso14d}
\begin{equation}
\rI\uplambda^*(\rR^q\rI f_*(\cM^\bullet))\rightarrow \rR^q\rI \bvf_*(\rI\uplambda'^*(\cM^\bullet))
\end{equation}
is an isomorphism. 
\end{prop}

It follows from \ref{p1-bcim14}(ii) applied to \eqref{p2-cmupiso14a} and the thick subcategory $\bMod^{\coh}(\co_\fX[\frac 1 p])$ of 
$\bIndMod(\co_\fX)$ \eqref{p2-cmupiso3}. Recall \ref{p2-cmupiso7} and observe that for all integers $q$ and $r$, 
$\rR^q\rI f_*(\cH^r (\cM^\bullet))$ is in the essential image of $\upalpha_{\co_{\fX}}^\coh$
by \ref{p2-cmupiso16}(i), \eqref{p2-cmupiso18b} and \eqref{p1-bcim1f}; and the base change morphism \eqref{p2-cmupiso14d}
\begin{equation}
\rI\uplambda^*(\rR^q\rI f_*(\cH^r(\cM^\bullet)))\rightarrow \rR^q\rI \bvf_*(\uplambda'^*(\cH^r(\cM^\bullet)))
\end{equation}
is an isomorphism by \ref{p2-cmupiso16}(ii) and \eqref{p1-bcim9f}.

\subsection{}\label{p2-cmupiso32}
We assume that $\fX$ is $\cS$-flat \eqref{p2-cmupiso1} and consider an exact sequence of locally free $\co_{\fX}$-modules of finite type
\begin{equation}\label{p2-cmupiso32a}
0\rightarrow \co_{\fX}\rightarrow \cF\rightarrow  \cE \rightarrow 0.
\end{equation}
Recall that for any rational number $r\geq 0$, we associate with the extension \eqref{p2-cmupiso32a} an $\co_{\fX}$-algebra $\cC^{(r)}$
\eqref{p1-thbn3b}, its $p$-adic completion $\hcC^{(r)}$, the $\co_{\fX}$-algebra $\hcC^{(r+)}$ \eqref{p1-thbn6e}
and the $\co_{\fX}$-derivations 
\begin{eqnarray}
\delta_{\hcC^{(r)}}\colon \hcC^{(r)}\rightarrow \cE\otimes_{\co_{\fX}}\hcC^{(r)},\label{p2-cmupiso32b}\\
\delta_{\hcC^{(r+)}}\colon \hcC^{(r+)}\rightarrow \cE\otimes_{\co_{\fX}}\hcC^{(r+)},\label{p2-cmupiso32c}
\end{eqnarray}
defined in \eqref{p1-thbn6b} and \eqref{p1-thbn6f}, respectively, 
which are also Higgs $\co_{\fX}$-fields. For simplicity, we set $\cC^\dagger=\hcC^{(0+)}$ and 
\begin{equation}\label{p2-cmupiso32d}
\delta=\delta_{\hcC^{(0+)}} \colon \cC^\dagger\rightarrow \cE\otimes_{\co_{\fX}}\cC^\dagger.
\end{equation}

We consider the ind-$\co_{\fX}$-algebra \eqref{p1-indmal2}
\begin{equation}\label{p2-cmupiso32e}
\IC^\dagger=\underset{\underset{r\in \mQ_{>0}}{\longrightarrow}}{\mlq\mlq\lim \mrq\mrq}\ \hcC^{(r)}. 
\end{equation}
In view of \eqref{p1-thbn6d}, we equip it with the $\co_{\fX}$-derivation \eqref{p1-indmal21}, 
which is also a Higgs $\co_{\fX}$-field \eqref{p1-indmal20},
\begin{equation}\label{p2-cmupiso32f}
\Idelta=\underset{\underset{r\in \mQ_{>0}}{\longrightarrow}}{\mlq\mlq\lim \mrq\mrq}\ \delta_{\hcC^{(r)}}\colon 
\IC^\dagger \rightarrow \cE\otimes_{\co_{\fX}}\IC^\dagger.
\end{equation}

We consider the ind-$\co_{\fX}$-algebra 
\begin{equation}\label{p2-cmupiso32g}
\IC^\dagger_\mQ=\underset{\underset{r\in \mQ_{>0}}{\longrightarrow}}{\mlq\mlq\lim \mrq\mrq}\ \upalpha_{\co_{\fX}}(\hcC^{(r)}_\mQ),
\end{equation}
that we equip with the $\co_{\fX}$-derivation, which is also a Higgs $\co_{\fX}$-field, 
\begin{equation}\label{p2-cmupiso32h}
\Idelta_\mQ=\underset{\underset{r\in \mQ_{>0}}{\longrightarrow}}{\mlq\mlq\lim \mrq\mrq}\ \upalpha_{\co_{\fX}}(\delta_{\hcC^{(r)},\mQ})\colon 
\IC^\dagger_\mQ \rightarrow \cE\otimes_{\co_{\fX}}\IC^\dagger_\mQ.
\end{equation}

\begin{prop}\label{p2-cmupiso33}
We keep the assumptions and notation of \ref{p2-cmupiso32} and assume, moreover, 
that the extension \eqref{p2-cmupiso32a} is split and that the $\co_{\fX}$-module $\cE$ is free. 
We denote by $(\cC^\dagger\otimes_{\co_{\fX}}\wedge^\bullet \cE,\delta^\bullet)$ the Dolbeault complex of $(\cC^\dagger,\delta)$
and by $(\IC^\dagger_\mQ\otimes_{\co_{\fX}}\wedge^\bullet \cE, \rI \delta^\bullet_\mQ)$ the Dolbeault complex of $(\IC^\dagger_\mQ,\Idelta_\mQ)$. 
Then, 
\begin{itemize}
\item[{\rm (i)}] The augmented Dolbeault complex 
\begin{equation}\label{p2-cmupiso33a}
\xymatrix{
{\co_{\fX}[\frac 1 p]}\ar[r]^-(0.5){\varepsilon}&{\cC^\dagger[\frac 1 p]}\ar[r]^-(0.5){\delta^0[\frac 1 p]}&
{\cC^\dagger[\frac 1 p]\otimes_{\co_{\fX}}\cE}\ar[r]^-(0.5){\delta^1[\frac 1 p]}&
{\cdots}\ar[r]& {\cC^\dagger[\frac 1 p]\otimes_{\co_{\fX}}\wedge^n \cE}\ar[r]^-(0.5){\delta^n[\frac 1 p]}&\cdots,}
\end{equation}
where $\co_{\fX}[\frac 1 p]$ is placed in degree $-1$ and $\varepsilon$ is the canonical homomorphism, 
is homotopically equivalent to zero, i.e., there exist $\co_{\fX}$-linear morphisms 
$h^0\colon \cC^\dagger[\frac 1 p]\rightarrow \co_{\fX}[\frac 1 p]$ and, for any $i\geq 1$, 
\begin{equation}\label{p2-cmupiso33b}
h^i\colon \cC^\dagger[\frac 1 p]\otimes_{\co_{\fX}}\wedge^i \cE\rightarrow \cC^\dagger[\frac 1 p]\otimes_{\co_{\fX}}\wedge^{i-1} \cE,
\end{equation}
such that $h^0\circ \varepsilon =\id$, $\varepsilon \circ h^0+h^1\circ \delta^0[\frac 1 p]=\id$ and for every $i\geq 1$, 
$\delta^{i-1}[\frac 1 p] \circ h^i+h^{i+1}\circ \delta^i[\frac 1 p]=\id$. 
\item[{\rm (ii)}] The augmented Dolbeault complex 
\begin{equation}\label{p2-cmupiso33c}
\xymatrix{
{\co_{\fX,\mQ}}\ar[r]^-(0.5){\epsilon}&{\IC^\dagger_\mQ}\ar[r]^-(0.5){\Idelta^0_\mQ}&
{\IC^\dagger_\mQ\otimes_{\co_{\fX}}\cE}\ar[r]^-(0.5){\Idelta^1_\mQ}&{\cdots}\ar[r]&{\IC^\dagger_\mQ\otimes_{\co_{\fX}}\wedge^n \cE}
\ar[r]^-(0.5){\Idelta^n_\mQ}&{\cdots,}}
\end{equation}
where $\co_{\fX,\mQ}$ is placed in degree $-1$ and $\epsilon$ is the canonical morphism,
is homotopically equivalent to zero, i.e., there exist morphisms of ind-$\co_{\fX}$-modules 
$\rI h^0\colon \IC^\dagger_\mQ\rightarrow \co_{\fX,\mQ}$ and, for any $i\geq 1$, 
\begin{equation}\label{p2-cmupiso33d}
\rI h^i\colon \IC^\dagger_\mQ\otimes_{\co_{\fX}}\wedge^i \cE\rightarrow \IC^\dagger_\mQ\otimes_{\co_{\fX}}\wedge^{i-1} \cE,
\end{equation}
such that $\rI  h^0\circ \epsilon=\id$, $\epsilon \circ \rI h^0+\rI  h^1\circ \rI \delta^0_\mQ=\id$ and for every $i\geq 1$, 
$\rI \delta^{i-1}_\mQ \circ \rI h^i+\rI h^{i+1}\circ \rI \delta^i_\mQ=\id$. 
\end{itemize}
\end{prop}

(i) Indeed, by \ref{p1-thbn36}, for every rational number $r>0$, there is a homotopy 
between the canonical morphism 
\begin{equation}
\tmK^\bullet(\hcC^{(2r)})[\frac 1 p]\rightarrow \tmK^\bullet(\hcC^{(r)})[\frac 1 p]
\end{equation}
and the zero morphism, where $\tmK^\bullet(\hcC^{(r)})$ is the augmented Dolbeault complex of $(\hcC^{(r)},\delta_{\hcC^{(r)}})$
defined in \eqref{p1-thbn35a}. 
The homotopy is explicitly given in the proof of (\cite{agt} II.11.2)
by choosing a splitting of the extension \eqref{p2-cmupiso32a} and a basis of $\cE$ \eqref{p1-thbn16}.  
We see, in particular, that we can choose these homotopies to be compatible for all $r\in \mQ_{>0}$. 
Their direct limit yields the desired homotopy. 

(ii) Similarly, by the proof of \eqref{p1-thbn36}, for every rational number $r>0$, there is a homotopy 
between the canonical morphism 
\begin{equation}
\tmK^\bullet_\mQ(\hcC^{(2r)})\rightarrow \tmK^\bullet_\mQ(\hcC^{(r)})
\end{equation}
and the zero morphism, where $\tmK^\bullet_\mQ(\hcC^{(r)})$ is the class of $\tmK^\bullet(\hcC^{(r)})$ in $\bMod_\mQ(\co_\fX)$. 
We can choose these homotopies to be compatible for all $r\in \mQ_{>0}$. The direct limit of their images by $\upalpha_{\co_\fX}$ 
yields the desired homotopy in $\bIndMod(\co_\fX)$.

\subsection{}\label{p2-cmupiso210}
Let $f\colon \fX'\rightarrow \fX$ be a morphism of finite presentation \eqref{p2-cmupiso1}, 
$\cE$ a locally free $\co_\fX$-module of finite type, 
\begin{equation}\label{p2-cmupiso210a}
0\rightarrow \co_{\fX'}\rightarrow \cF\rightarrow f^*(\cE)\rightarrow 0
\end{equation}
an exact sequence of locally free $\co_{\fX'}$-modules of finite type, $\cE'$ a locally free $\co_{\fX'}$-module of finite type, 
\begin{equation}\label{p2-cmupiso210d}
u\colon f^*(\cE)\rightarrow \cE'
\end{equation}
an $\co_{\fX'}$-linear morphism. We assume that $\fX$ and $\fX'$ are $\cS$-flat. 

Following \ref{p2-cmupiso32}, we associate with the extension \eqref{p2-cmupiso210a} an $\co_{\fX'}$-algebra $\cC^\dagger=\hcC^{(0+)}$ 
and an $\co_{\fX'}$-derivation \eqref{p2-cmupiso32d}
\begin{equation}\label{p2-cmupiso21d}
\delta=\delta_{\hcC^{(0+)}} \colon \cC^\dagger\rightarrow f^*(\cE)\otimes_{\co_{\fX'}}\cC^\dagger.
\end{equation}
It is a Higgs $\co_{\fX'}$-field. We set
\begin{equation}\label{p2-cmupiso21e}
\delta'=(u\otimes \id)\circ \delta\colon \cC^\dagger\rightarrow \cE'\otimes_{\co_{\fX'}}\cC^\dagger.
\end{equation}

We also associate with the extension \eqref{p2-cmupiso210a} an ind-$\co_{\fX'}$-algebra $\IC^\dagger$ \eqref{p2-cmupiso32e}
and an $\co_{\fX'}$-derivation \eqref{p2-cmupiso32f}
\begin{equation}\label{p2-cmupiso21g}
\Idelta\colon \IC^\dagger \rightarrow f^*(\cE)\otimes_{\co_{\fX'}}\IC^\dagger.
\end{equation}
It is a Higgs $\co_{\fX'}$-field \eqref{p1-indmal20}. We set
\begin{equation}\label{p2-cmupiso21h}
\Idelta'=(u\otimes \id)\circ \Idelta\colon \IC^\dagger\rightarrow \cE'\otimes_{\co_{\fX'}}\IC^\dagger.
\end{equation}

We defined in \eqref{p1-thbn9b} a twisting functor by the extension $\cF$ 
for Higgs $\co_{\fX'}$-modules with coefficients in $f^*(\cE)$:
\begin{equation}\label{p2-cmupiso210b}
\uptau\colon \bHM(\co_{\fX'},f^*(\cE))\rightarrow \bHM(\co_{\fX'},f^*(\cE)).
\end{equation}
Composing with the pullback functor $f^*$, we obtain a functor
\begin{equation}\label{p2-cmupiso210c}
f^*_\uptau\colon 
\begin{array}[t]{clcr}
\bHM(\co_\fX,\cE)&\rightarrow& \bHM(\co_{\fX'},f^*(\cE))\\
(N,\theta)&\mapsto&\uptau(f^*(N),f^*(\theta)). 
\end{array}
\end{equation} 
called the {\em pullback functor by $f$ twisted by the extension $\cF$} (or the {\em twisted pullback functor by $f$} when the extension is implicit); 
see \ref{p1-tphdi2}. 

Let $(N,\theta)$ be a Higgs $\co_{\fX'}$-module with coefficients in $\cE'$. We equip $N\otimes_{\co_{\fX'}}\cC^\dagger$ with the Higgs 
$\co_{\fX'}$-field 
\begin{equation}\label{p2-cmupiso22a}
\vartheta= \theta\otimes \id+\id\otimes \delta'\colon N\otimes_{\co_{\fX'}}\cC^\dagger \rightarrow \cE'\otimes_{\co_{\fX'}}N\otimes_{\co_{\fX'}}\cC^\dagger,
\end{equation}
and denote by $\mK^{\bullet}(N\otimes_{\co_{\fX'}}\cC^\dagger)$ the associated Dolbeault complex,
that we call the {\em Dolbeault complex of $(N,\theta)$ twisted by the extension $\cF$} \eqref{p2-cmupiso21a}. 
By \ref{p1-tphdi6}, the Higgs field $\id\otimes \delta$ on $N\otimes_{\co_{\fX'}}\cC^\dagger$ 
induces a morphism of complexes of $\co_{\fX'}$-modules
\begin{equation}\label{p2-cmupiso22b}
\fd \colon \mK^\bullet(N\otimes_{\co_{\fX'}}\cC^\dagger)\rightarrow f^*(\cE)\otimes_{\co_{\fX'}} \mK^\bullet(N\otimes_{\co_{\fX'}}\cC^\dagger).
\end{equation} 
By the projection formula, for any integer $q\geq 0$, $\rR^qf_*(\fd)$ identifies with a Higgs $\co_\fX$-field 
\begin{equation}\label{p2-cmupiso22c}
\rR^qf_*(\fd)\colon  \rR^qf_*(\mK^\bullet(N\otimes_{\co_{\fX'}}\cC^\dagger))\rightarrow \cE\otimes_{\co_\fX} \rR^qf_*(\mK^\bullet(N\otimes_{\co_{\fX'}}\cC^\dagger)). 
\end{equation}
We thus obtain a functor that we denote by 
\begin{equation}\label{p2-cmupiso22d}
\rR^qf^\uptau_*\colon 
\begin{array}[t]{clcr}
\bHM(\co_{\fX'},\cE')&\rightarrow&\bHM(\co_\fX,\cE)\\
(N,\theta)&\mapsto&(\rR^qf_*(\mK^\bullet(N\otimes_{\co_{\fX'}}\cC^\dagger)),-\rR^qf_*(\fd)),
\end{array}
\end{equation}
and call the {\em $q$th higher direct image functor by $f$ twisted by the extension $\cF$ \eqref{p2-cmupiso210a}} 
(or the twisted $q$th higher direct image functor by $f$ when the extension is implicit); see \eqref{p1-tphdi6d}. 

Let $(\cN,\uptheta)$ be a Higgs ind-$\co_{\fX'}$-module with coefficients in $\cE'$ \eqref{p1-indmal20}. 
We equip $\cN\otimes_{\co_{\fX'}}\IC^\dagger$ with the Higgs $\co_{\fX'}$-field 
\begin{equation}\label{p2-cmupiso23a}
\rI\vartheta=\uptheta\otimes \id+\id\otimes \Idelta'\colon \cN\otimes_{\co_{\fX'}}\IC^\dagger \rightarrow 
\cE'\otimes_{\co_{\fX'}}\cN\otimes_{\co_{\fX'}}\IC^\dagger,
\end{equation}
and denote by $\mK^{\bullet}(\cN\otimes_{\co_{\fX'}}\IC^\dagger)$ the associated Dolbeault complex. 
The Higgs field $\id\otimes \Idelta$ on $\cN\otimes_{\co_{\fX'}}\IC^\dagger$ 
induces a morphism of complexes of ind-$\co_{\fX'}$-modules
\begin{equation}\label{p2-cmupiso23b}
\Ifd \colon \mK^\bullet(\cN\otimes_{\co_{\fX'}}\IC^\dagger)\rightarrow f^*(\cE)\otimes_{\co_{\fX'}} \mK^\bullet(\cN\otimes_{\co_{\fX'}}\IC^\dagger).
\end{equation} 
Indeed, the question being local on $\fX'$ by (\cite{ag2} 2.7.17(i)), we may assume that the $\co_\fX$-module $\cE$ is free of finite type. 
The assertion then follows by considering  the components, with respect to a basis of $f^*(\cE)$, 
of the derivation $\Idelta$ of $\IC^\dagger$, which are $\co_{\fX'}$-linear endomorphisms of $\IC^\dagger$ 
commuting with each other, since the $\delta_{\hcC^{(r)}}$'s are Higgs fields \eqref{p2-cmupiso32b}. 
By the projection formula (\cite{ag2} 2.7.14), for any integer $q\geq 0$, $\rR^q\rI f_*(\Ifd)$ identifies with a Higgs $\co_\fX$-field 
\begin{equation}\label{p2-cmupiso23c}
\rR^q\rI f_*(\fd)\colon  \rR^q\rI f_*(\mK^\bullet(\cN\otimes_{\co_{\fX'}}\IC^\dagger))\rightarrow 
\cE\otimes_{\co_\fX} \rR^qf_*(\mK^\bullet(\cN\otimes_{\co_{\fX'}}\IC^\dagger)). 
\end{equation}
We thus obtain a functor that we denote by 
\begin{equation}\label{p2-cmupiso23d}
\rR^q\rI f^\uptau_*\colon 
\begin{array}[t]{clcr}
\bIndHM(\co_{\fX'},\cE')&\rightarrow&\bIndHM(\co_\fX,\cE)\\
(\cN,\uptheta)&\mapsto&(\rR^q\rI f_*(\mK^\bullet(\cN\otimes_{\co_{\fX'}}\IC^\dagger)),-\rR^q\rI f_*(\Ifd)),
\end{array}
\end{equation} 
and call the {\em $q$th higher direct image functor by $\rI f$ twisted by the extension \eqref{p2-cmupiso210a}}
(or  the {\em  twisted $q$th higher direct image functor by $\rI f$} when the extension is implicit).

\begin{prop}\label{p2-cmupiso211}
Under the assumptions of \ref{p2-cmupiso210} and with the same notation, for every 
locally CL-small Higgs $\co_{\fX}[\frac 1 p]$-module $(N,\theta)$ with coefficients in $\cE$ \eqref{p1-tshbn13}, 
the Higgs $\co_{\fX'}[\frac 1 p]$-module $f^*_\uptau(N,\theta)$ \eqref{p2-cmupiso210c} is locally CL-small. 
\end{prop}

Indeed, the Higgs $\co_{\fX'}[\frac 1 p]$-module $(f^*(N),f^*(\theta))$ is locally CL-small by (\cite{egr1} 2.10.14). 
The proposition follows by \ref{p1-thbn31}.

\begin{prop}\label{p2-cmupiso24}
Under the assumptions of \ref{p2-cmupiso210} and with the same notation,
for every Higgs ind-$\co_{\fX'}$-module $(\cN,\uptheta)$ with coefficients in $\cE'$
and every integer $q\geq 0$, we have a canonical functorial isomorphism
\begin{equation}\label{p2-cmupiso24a}
\rR^qf^\uptau_*(\kappa_{\co_{\fX'}}(\cN,\uptheta))\stackrel{\sim}{\rightarrow} \kappa_{\co_\fX}(\rR^q\rI f^\uptau_*(\cN,\uptheta)),
\end{equation}
where $\kappa_{\co_\fX}$ denotes the functor \eqref{p2-cmupiso31b} and $\kappa_{\co_{\fX'}}$ the analogous functor for $\fX'$. 
\end{prop}

We set $(N,\theta)= \kappa_{\co_{\fX'}}(\cN,\uptheta)$. By (\cite{ag2} (2.7.1.5)), we have a canonical isomorphism 
\begin{equation}\label{p2-cmupiso24b}
\mK^\bullet(N\otimes_{\co_{\fX'}}\cC^\dagger) \stackrel{\sim}{\rightarrow}
\kappa_{\co_{\fX'}}(\mK^\bullet(\cN\otimes_{\co_{\fX'}}\IC^\dagger)). 
\end{equation}
We deduce, by \ref{p2-cmupiso13}(i), a canonical isomorphism 
\begin{equation}\label{p2-cmupiso24c}
\rR^q f_*(\mK^\bullet(N\otimes_{\co_{\fX'}}\cC^\dagger))\stackrel{\sim}{\rightarrow}
\kappa_{\co_{\fX}}(\rR^q\rI f_*(\mK^\bullet(\cN\otimes_{\co_{\fX'}}\IC^\dagger))).
\end{equation}
The diagram
\begin{equation}\label{p2-cmupiso24d}
\xymatrix{
{\mK^\bullet(N\otimes_{\co_{\fX'}}\cC^\dagger)}\ar[rr]^-(0.5){\fd}\ar[d]&&{f^*(\cE)\otimes_{\co_{\fX'}} \mK^\bullet(N\otimes_{\co_{\fX'}}\cC^\dagger)}\ar[d]\\
{\kappa_{\co_{\fX'}}(\mK^\bullet(\cN\otimes_{\co_{\fX'}}\IC^\dagger))}\ar[rr]^-(0.5){\kappa_{\co_{\fX'}}(\rI\fd)}&&
{\kappa_{\co_{\fX'}}(f^*(\cE)\otimes_{\co_{\fX'}}\mK^\bullet(\cN\otimes_{\co_{\fX'}}\IC^\dagger)),}}
\end{equation}
where the horizontal arrows are defined in \eqref{p2-cmupiso22b} and \eqref{p2-cmupiso23b},
and the vertical arrows are induced by \eqref{p2-cmupiso24b}, is commutative. 
Moreover, the diagram 
\begin{equation}
\xymatrix{
{\rR^q f_*(f^*(\cE)\otimes_{\co_{\fX'}} \mK^\bullet(N\otimes_{\co_{\fX'}}\cC^\dagger))}\ar[r]\ar[d]&
{\kappa_{\co_{\fX}}(\rR^q\rI f_*(f^*(\cE)\otimes_{\co_{\fX'}} \mK^\bullet(\cN\otimes_{\co_{\fX'}}\IC^\dagger)))}\ar[d]\\
{\cE\otimes_{\co_{\fX}}\rR^q f_*(\mK^\bullet(N\otimes_{\co_{\fX'}}\cC^\dagger)))}\ar[r]&
{\kappa_{\co_{\fX}}(\cE\otimes_{\co_{\fX}}\rR^q\rI f_*(\mK^\bullet(\cN\otimes_{\co_{\fX'}}\IC^\dagger))),}}
\end{equation}
where the upper (resp.\ lower) horizontal arrow is the isomorphism 
defined as in (resp.\ induced by) \eqref{p2-cmupiso24c} and the vertical arrows are the projection 
isomorphisms (\cite{ag2} 2.7.14), is commutative. Indeed, to check it, we may assume that $\cE$ is free, in which case it is obvious. 
The proposition follows.

\subsection{}\label{p2-cmupiso21}
In the remainder of this section, we fix a morphism of finite presentation $f\colon \fX'\rightarrow \fX$, 
a locally free $\co_\fX$-module of finite type $\cE$ and two exact sequences of locally free $\co_{\fX'}$-modules of finite type
\begin{eqnarray}
0\rightarrow \co_{\fX'}\rightarrow \cF\rightarrow f^*(\cE)\rightarrow 0,\label{p2-cmupiso21a}\\
0\rightarrow f^*(\cE)\stackrel{u}{\rightarrow}\cE'\stackrel{u'}{\rightarrow}\ucE'\rightarrow 0.\label{p2-cmupiso21b}
\end{eqnarray}
We assume that $\fX$ and $\fX'$ are $\cS$-flat. In this context, we use the notation introduced in \ref{p2-cmupiso210}. 

\subsection{}\label{p2-cmupiso26}
We keep the assumptions and notation of \ref{p2-cmupiso21}. 
Let $\rho\colon \cF\rightarrow \co_{\fX'}$ be a splitting of the extension \eqref{p2-cmupiso21a}. 
By \ref{p1-thbn16}, $\rho$ induces a homomorphism of $\co_{\fX'}$-algebras 
\begin{equation}\label{p2-cmupiso26a}
\varrho^\dagger\colon \cC^\dagger\rightarrow \co_{\fX'},
\end{equation}
and a homomorphism of ind-$\co_{\fX'}$-algebras 
\begin{equation}\label{p2-cmupiso26b}
\rI\varrho^\dagger\colon \IC^\dagger\rightarrow \co_{\fX'}.
\end{equation}

For every Higgs $\co_{\fX'}$-module  $(N,\theta)$ with coefficients in $\cE'$, the diagram 
\begin{equation}\label{p2-cmupiso26c}
\xymatrix{
{N\otimes_{\co_{\fX'}}\cC^\dagger}\ar[r]^-(0.5){\vartheta}\ar[d]_{\id\otimes\varrho^\dagger}&
{\cE'\otimes_{\co_{\fX'}}N\otimes_{\co_{\fX'}}\cC^\dagger}\ar[d]^{u'\otimes\id\otimes \varrho^\dagger}\\
{N}\ar[r]^-(0.5){\utheta}&{\ucE'\otimes_{\co_{\fX'}}N,}}
\end{equation}
where $\vartheta$ is the total Higgs field defined in \eqref{p2-cmupiso22a}, 
$u'$ is the morphism defined in \eqref{p2-cmupiso21b} and $\utheta=(u'\otimes \id)\circ \theta$, is commutative. 
We deduce a morphism of Dolbeault complexes 
\begin{equation}\label{p2-cmupiso26d}
\upbeta\colon \mK^\bullet(N\otimes_{\co_{\fX'}}\cC^\dagger) \rightarrow \umK^\bullet(N). 
\end{equation}

For every Higgs ind-$\co_{\fX'}$-module  $(\cN,\uptheta)$ with coefficients in $\cE'$, the diagram 
\begin{equation}\label{p2-cmupiso26e}
\xymatrix{
{\cN\otimes_{\co_{\fX'}}\IC^\dagger}\ar[r]^-(0.5){\rI\vartheta}\ar[d]_{\id\otimes\rI\varrho^\dagger}&
{\cE'\otimes_{\co_{\fX'}}\cN\otimes_{\co_{\fX'}}\IC^\dagger}\ar[d]^{u'\otimes\id\otimes \rI\varrho^\dagger}\\
{\cN}\ar[r]^-(0.5){\uuptheta}&{\ucE'\otimes_{\co_{\fX'}}\cN,}}
\end{equation}
where $\rI\vartheta$ is the total Higgs field defined in \eqref{p2-cmupiso23a} and $\uuptheta=(u'\otimes \id)\circ \uptheta$, is commutative. 
We deduce a morphism of Dolbeault complexes 
\begin{equation}\label{p2-cmupiso26f}
\rI\upbeta\colon \mK^\bullet(\cN\otimes_{\co_{\fX'}}\IC^\dagger) \rightarrow \umK^\bullet(\cN). 
\end{equation}

\begin{prop}\label{p2-cmupiso28}
We keep the assumptions and notation of \ref{p2-cmupiso21}, and assume, moreover, 
that the $\co_{\fX'}$-module $f^*(\cE)$ is free and that we are given a splitting  $\rho\colon \cF\rightarrow \co_{\fX'}$ 
of the extension \eqref{p2-cmupiso21a} and a splitting $v'\colon \ucE'\rightarrow \cE'$ of the extension \eqref{p2-cmupiso21b}. 
Let $(N,\theta)$ be a {\em locally CL-small} Higgs $\co_{\fX'}[\frac 1 p]$-module with 
coefficients in $\cE'$  \eqref{p1-tshbn13}; it is, in particular, coherent \eqref{p2-cmupiso31}. 
We set $\utheta=(\id \otimes u')\circ \theta \colon N\rightarrow N\otimes_{\co_{\fX'}}\ucE'$ \eqref{p2-cmupiso21b} and denote by 
$\umK(N)$ the Dolbeault complex of $(N,\utheta)$. 
\begin{itemize}
\item[{\rm (i)}] Let $\mK^{\bullet}(N\otimes_{\co_{\fX'}}\cC^\dagger)$ be the Dolbeault complex of $N\otimes_{\co_{\fX'}}\cC^\dagger$ equipped 
with the total Higgs field $\vartheta = \theta\otimes \id+\id\otimes \delta'$ \eqref{p2-cmupiso22a}. 
Then, $\rho$ and $v'$ induce a morphism of complexes of $\co_{\fX'}$-modules
\begin{equation}\label{p2-cmupiso28a}
\umK^\bullet(N)\rightarrow \mK^{\bullet}(N\otimes_{\co_{\fX'}}\cC^\dagger),
\end{equation}
whose cone is homotopically equivalent to zero. The latter fits into a commutative diagram 
\begin{equation}\label{p2-cmupiso28aa}
\xymatrix{
{\umK^\bullet(N)}\ar[r]\ar[rd]_\id&{\mK^{\bullet}(N\otimes_{\co_{\fX'}}\cC^\dagger)}\ar[d]^{\upbeta}\\
&{\umK^\bullet(N),}}
\end{equation}
where $\upbeta$ is defined in \eqref{p2-cmupiso26d} by the choice of $\rho$.
\item[{\rm (ii)}] Assume that the $\co_{\fX'}[\frac 1 p]$-module $N$ is flat. 
Let $(\cN,\uptheta)=\upalpha^\coh_{\co_{\fX'}}(N,\theta)$ \eqref{p2-cmupiso31h},
$\mK^{\bullet}(\cN\otimes_{\co_{\fX'}}\IC^\dagger)$ be the Dolbeault complex of $\cN\otimes_{\co_{\fX'}}\IC^\dagger$ equipped 
with the total Higgs field $\rI\vartheta= \uptheta\otimes \id+\id\otimes \Idelta'$ \eqref{p2-cmupiso23a}. 
Then, $\rho$ and $v'$ induce a morphism of complexes of ind-$\co_{\fX'}$-modules
\begin{equation}\label{p2-cmupiso28b}
\upalpha^\coh_{\co_{\fX'}}(\umK^\bullet(N))\rightarrow \mK^{\bullet}(\cN\otimes_{\co_{\fX'}}\IC^\dagger),
\end{equation}
whose cone is homotopically equivalent to zero. 
The latter fits into a commutative diagram 
\begin{equation}\label{p2-cmupiso28bb}
\xymatrix{
{\upalpha^\coh_{\co_{\fX'}}(\umK^\bullet(N))}\ar[r]\ar[rd]_\id&{\mK^{\bullet}(\cN\otimes_{\co_{\fX'}}\IC^\dagger)}\ar[d]^{\rI\upbeta}\\
&{\upalpha^\coh_{\co_{\fX'}}(\umK^\bullet(N)),}}
\end{equation}
where $\rI\upbeta$ is defined in \eqref{p2-cmupiso26f} by the choice of $\rho$.
\end{itemize}
\end{prop}

We denote by $\cE'_0$ (resp.\ $\cE'_1$) 
the image of $u$ (resp.\ $v'$) \eqref{p2-cmupiso21b}, so we have a decomposition $\cE'=\cE'_0\oplus\cE'_1$. 
The Higgs field $\theta$ induces two Higgs $\co_{\fX'}$-fields
\begin{equation}\label{p2-cmupiso28ic}
\theta_i\colon N\rightarrow N\otimes_{\co_{\fX'}}\cE'_i, \ \ \ {\rm for}\ \ \ i=0,1.
\end{equation}
Identifying $f^*(\cE)$ with $\cE'_0$ by the isomorphism induced by $u$, 
$\delta$ \eqref{p2-cmupiso21d} induces an $\co_{\fX'}$-derivation 
\begin{equation}\label{p2-cmupiso28d}
\delta'_0\colon \cC^\dagger\rightarrow \cC^\dagger \otimes_{\co_{\fX'}}\cE'_0,
\end{equation}
which in turn induces the derivation $\delta'$ \eqref{p2-cmupiso21e} via the inclusion $\cE'_0\rightarrow \cE'$. 
Similarly, $\Idelta$ \eqref{p2-cmupiso21g} induces an $\co_{\fX'}$-derivation 
\begin{equation}\label{p2-cmupiso28dd}
\Idelta'_0\colon \IC^\dagger\rightarrow \IC^\dagger \otimes_{\co_{\fX'}}\cE'_0,
\end{equation}
which in turn induces the derivation $\Idelta'$ \eqref{p2-cmupiso21h} via the inclusion $\cE'_0\rightarrow \cE'$. 

(i) We consider the total Higgs $\co_{\fX'}$-field
\begin{equation}\label{p2-cmupiso28i1}
\vartheta_0=\theta_0\otimes \id+\id\otimes \delta'_0\colon N\otimes_{\co_{\fX'}}\cC^\dagger \rightarrow 
N\otimes_{\co_{\fX'}}\cC^\dagger\otimes_{\co_{\fX'}} \cE'_0,
\end{equation}
and the double complex (with the sign conventions of \cite{sp} \href{https://stacks.math.columbia.edu/tag/0FNB}{0FNB})
\begin{equation}\label{p2-cmupiso28i2}
\cK^{i,j}=N\otimes_{\co_{\fX'}}\cC^\dagger\otimes_{\co_{\fX'}}\wedge^i\cE'_0\otimes_{\co_{\fX'}}\wedge^j\cE'_1,
\end{equation}
where the differentials are defined for any local sections $\nu$ of $N\otimes_{\co_{\fX'}}\cC^\dagger$, $\omega_0$ 
of $\wedge^i\cE'_0$ and  $\omega_1$ of $\wedge^j\cE'_1$ by
\begin{eqnarray}
d_{\cK,1}^{i,j}(\nu\otimes \omega_0 \otimes \omega_1)&=&\vartheta_0(\nu)\wedge \omega_0 \otimes \omega_1,\label{p2-cmupiso28i3}\\
d_{\cK,2}^{i,j}(\nu\otimes \omega_0 \otimes \omega_1)&=&(-1)^i(\theta_1\otimes \id)(\nu)\wedge \omega_0 \otimes \omega_1,\label{p2-cmupiso28i4}
\end{eqnarray}
where the $\wedge$ product is taken in the left tensor $\co_{\fX'}$-algebra 
$\wedge^\bullet \cE'_0\ {^g\otimes}_{\co_{\fX'}} \wedge^\bullet \cE'_1$ (cf. \cite{alg1-3} III § 4.7 remarks page 49).  
We denote by $\rK^\bullet$ the total complex associated with $\cK^{\bullet,\bullet}$ (with the sign conventions of \cite{sp} \href{https://stacks.math.columbia.edu/tag/012Z}{012Z}).

By (\cite{alg1-3} III § 7.7 prop.~10), we have a canonical isomorphism of bigraded algebras
\begin{equation}
\wedge^\bullet \cE'_0\ {^g\otimes} \wedge^\bullet \cE'_1 \stackrel{\sim}{\rightarrow}\wedge^\bullet\cE'.
\end{equation}
Since 
\begin{equation}\label{p2-cmupiso28i6}
\vartheta_0\oplus(\theta_1\otimes \id)\colon N\otimes_{\co_{\fX'}}\cC^\dagger \rightarrow N\otimes_{\co_{\fX'}}\cC^\dagger 
\otimes_{\co_{\fX'}}(\cE'_0\oplus \cE'_1)
\end{equation}
identifies with the total Higgs field $\vartheta=\theta\otimes \id+\id\otimes \delta'$, we deduce an isomorphism 
\begin{equation}\label{p2-cmupiso28i17}
\rK^{\bullet}\stackrel{\sim}{\rightarrow}\mK^\bullet(N\otimes_{\co_{X'}}\cC^\dagger).
\end{equation}

We set 
\begin{equation}\label{p2-cmupiso28i5}
M=(N\otimes_{\co_{\fX'}}\cC^\dagger)^{\vartheta_0=0}.
\end{equation} 
Since $\vartheta_0\oplus(\theta_1\otimes \id)$ is a Higgs $\co_{\fX'}$-field \eqref{p2-cmupiso28i6}, 
$\theta_1\otimes \id$ induces a Higgs $\co_{\fX'}$-field on $M$ that we denote by 
\begin{equation}\label{p2-cmupiso28i7}
\lambda_1\colon M\rightarrow M\otimes_{\co_{\fX'}}\cE'_1.
\end{equation} 
By \ref{p1-tshbn16}, since $(N,\theta_0)$ is a locally CL-small Higgs $\co_{\fX'}[\frac 1 p]$-module,
it is weakly twistable by the extension \eqref{p2-cmupiso21a} in the sense of \eqref{p1-thbn30}. 
Therefore, the canonical $\cC^\dagger$-linear morphism 
\begin{equation}\label{p2-cmupiso28i8}
M\otimes_{\co_{\fX'}}\cC^\dagger\rightarrow N\otimes_{\co_{\fX'}}\cC^\dagger
\end{equation}
is an isomorphism, compatible with the Higgs $\co_{\fX'}$-fields indicated on the same line of the table
\begin{equation}
\begin{tabular}{|c|c|c|}
\hline
left hand side&  right  hand  side & coefficients\\
\hline
$\id\otimes \delta'_0$ &  $\vartheta_0=\theta_0\otimes \id+\id\otimes \delta'_0$ & $\cE'_0$\\
\hline
$\lambda_1\otimes \id$ & $\theta_1\otimes \id$ & $\cE'_1$\\
\hline
\end{tabular}
\end{equation}

By \ref{p1-thbn16}, the splitting $\rho$ of \eqref{p2-cmupiso21a} induces a homomorphism of $\co_{\fX'}$-algebras 
$\varrho^\dagger\colon \cC^\dagger\rightarrow \co_{\fX'}$. 
The $\cC^\dagger$-linear isomorphism \eqref{p2-cmupiso28i8} induces by scalar extension by $\varrho^\dagger$, 
an $\co_{\fX'}$-linear isomorphism $\iota\colon M\stackrel{\sim}{\rightarrow} N$.
We denote by $\ulambda_1\colon M\rightarrow M\otimes_{\co_{\fX'}}\ucE'$ (resp.\ $\utheta_1\colon N\rightarrow N\otimes_{\co_{\fX'}}\ucE'$) 
the Higgs $\co_{\fX'}$-field induced by $\lambda_1$ (resp.\ $\theta_1$) and the isomorphism $\cE'_1\stackrel{\sim}{\rightarrow} \ucE'$ 
deduced from $u'$ \eqref{p2-cmupiso21b}. Observe that $\utheta_1=\utheta$. 
Then, $\iota$ is underlying an isomorphism of Higgs $\co_{\fX'}$-modules with coefficients in $\ucE'$ 
\begin{equation}\label{p2-cmupiso28i18}
\iota\colon (M,\ulambda_1) \stackrel{\sim}{\rightarrow} (N,\utheta),
\end{equation}
that fits into a commutative diagram of Higgs $\co_{\fX'}$-modules 
\begin{equation}\label{p2-cmupiso28i180}
\xymatrix{
{(M,\ulambda_1)}\ar[r]\ar[d]_\iota&{(M\otimes_{\co_{\fX'}}\cC^\dagger,\ulambda_1\otimes \id)}\ar[d]\\
{(N,\utheta)}&{(N\otimes_{\co_{\fX'}}\cC^\dagger,\utheta_1\otimes \id),}\ar[l]_-(0.45){\id\otimes\varrho^\dagger}}
\end{equation}
where the right vertical arrow is the isomorphism \eqref{p2-cmupiso28i8}. 

We consider the double complex 
\begin{equation}\label{p2-cmupiso28i9}
\cL^{i,j}=M\otimes_{\co_{\fX'}}\cC^\dagger\otimes_{\co_{\fX'}}\wedge^i\cE'_0\otimes_{\co_{\fX'}}\wedge^j\cE'_1,
\end{equation}
where the differentials are defined for any local sections $\mu$ of $M\otimes_{\co_{\fX'}}\cC^\dagger$, $\omega_0$ 
of $\wedge^i\cE'_0$ and  $\omega_1$ of $\wedge^j\cE'_1$ by
\begin{eqnarray}
d_{\cL,1}^{i,j}(\mu\otimes \omega_0 \otimes \omega_1)&=&(\id\otimes \delta'_0)(\mu)\wedge \omega_0 \otimes \omega_1,\label{p2-cmupiso28i10}\\
d_{\cL,2}^{i,j}(\mu\otimes \omega_0 \otimes \omega_1)&=&(-1)^i(\lambda_1\otimes \id)(\mu)\wedge \omega_0 \otimes \omega_1,\label{p2-cmupiso28i100}
\end{eqnarray}
where the $\wedge$ product is taken in the left tensor $\co_{\fX'}$-algebra 
$\wedge^\bullet \cE'_0\ {^g\otimes}_{\co_{\fX'}} \wedge^\bullet \cE'_1$.  
The isomorphism \eqref{p2-cmupiso28i8} induces an isomorphism of bicomplexes 
\begin{equation}\label{p2-cmupiso28i12}
\cL^{\bullet,\bullet}\stackrel{\sim}{\rightarrow}\cK^{\bullet,\bullet}. 
\end{equation}
Denoting by $\rL^\bullet$ the total complex associated with $\cL^{\bullet,\bullet}$, we deduce an isomorphism 
\begin{equation}\label{p2-cmupiso28i13}
\rL^{\bullet}\stackrel{\sim}{\rightarrow}\rK^{\bullet}. 
\end{equation}

We denote by $(\mK^\bullet(M),\lambda_1^\bullet)$ the Dolbeault complex of $(M,\lambda_1)$. We have a canonical morphism of complexes 
\begin{equation}\label{p2-cmupiso28i14}
\epsilon^\bullet\colon (\mK^\bullet(M),\lambda_1^\bullet)\rightarrow (\cL^{0,\bullet},d_{\cL,2}^{0,\bullet}).
\end{equation}
It follows from the definition \eqref{p2-cmupiso28i10} that the composition of $\epsilon^\bullet$ with the morphism 
\begin{equation}\label{p2-cmupiso28i15}
d_{\cL,1}^{0,\bullet}\colon (\cL^{0,\bullet},d_{\cL,2}^{0,\bullet})\rightarrow (\cL^{1,\bullet},d_{\cL,2}^{1,\bullet})
\end{equation}
vanishes. We deduce a morphism 
\begin{equation}\label{p2-cmupiso28i16}
\varepsilon^\bullet\colon \mK^\bullet(M)\rightarrow \rL^{\bullet}.
\end{equation}
By \ref{p2-cmupiso34} below, the mapping cone of $\varepsilon^\bullet$ is homotopically equivalent to zero. 
We take then for \eqref{p2-cmupiso28a} the morphism induced by \eqref{p2-cmupiso28i16} and the isomorphisms \eqref{p2-cmupiso28i17}, \eqref{p2-cmupiso28i18} 
and \eqref{p2-cmupiso28i13}. 
An explicit computation shows that the composition of morphisms of complexes
\begin{equation}
\mK^\bullet(M)\stackrel{\varepsilon^\bullet}{\longrightarrow} \rL^\bullet \longrightarrow \rK^\bullet \longrightarrow 
\mK^\bullet(N\otimes_{\co_{X'}}\cC^\dagger)\stackrel{\upbeta}{\longrightarrow} \umK^\bullet(N),
\end{equation}
where the second (resp.\ third) arrow is the ismorphism \eqref{p2-cmupiso28i13} (resp.\ \eqref{p2-cmupiso28i17}), is the morphism of Dolbeault 
complexes induced by $\iota$ \eqref{p2-cmupiso28i18}, which proves the commutativity of \eqref{p2-cmupiso28aa}.

(ii) The Higgs $\co_{\fX'}$-field $\uptheta$ on $\cN$ induces two Higgs $\co_{\fX'}$-fields
\begin{equation}\label{p2-cmupiso28iic}
\uptheta_i\colon \cN\rightarrow \cN\otimes_{\co_{\fX'}}\cE'_i, \ \ \ {\rm for}\ \ \ i=0,1.
\end{equation}
We consider the total Higgs $\co_{\fX'}$-field
\begin{equation}\label{p2-cmupiso28ii1}
\rI\vartheta_0=\uptheta_0\otimes \id+\id\otimes \Idelta'_0\colon \cN\otimes_{\co_{\fX'}}\IC^\dagger \rightarrow 
\cN\otimes_{\co_{\fX'}}\IC^\dagger\otimes_{\co_{\fX'}} \cE'_0,
\end{equation}
and denote by $(\mK^\bullet(\cN\otimes_{\co_{\fX'}}\IC^\dagger),\rI\vartheta_0^\bullet)$ the Dolbeault complex of 
$(\cN\otimes_{\co_{\fX'}}\IC^\dagger,\rI\vartheta_0)$. 
We consider the double complex 
\begin{equation}\label{p2-cmupiso28ii2}
\IcK^{i,j}=\mK^i(\cN\otimes_{\co_{\fX'}}\IC^\dagger)\otimes_{\co_{\fX'}}\wedge^j\cE'_1=
\cN\otimes_{\co_{\fX'}}\IC^\dagger\otimes_{\co_{\fX'}}\wedge^i\cE'_0\otimes_{\co_{\fX'}}\wedge^j\cE'_1,
\end{equation}
where the differentials are defined by
\begin{eqnarray}
d_{\IcK,1}^{i,j}&=&\rI\vartheta_0^i \otimes \id_{\wedge^j\cE'_1},\label{p2-cmupiso28ii3}\\
d_{\IcK,2}^{i,j}&=&\id_{\wedge^i\cE'_0} \otimes \uptheta^j_1.\label{p2-cmupiso28ii4}
\end{eqnarray}
We denote by $\IrK^\bullet$ the total complex associated with $\IcK^{\bullet,\bullet}$. Since 
\begin{equation}
\rI\vartheta_0\oplus(\uptheta_1\otimes \id)\colon \cN\otimes_{\co_{\fX'}}\IC^\dagger \rightarrow \cN\otimes_{\co_{\fX'}}\IC^\dagger 
\otimes_{\co_{\fX'}}(\cE'_0\oplus \cE'_1)
\end{equation}
identifies with the total Higgs field $\rI\vartheta=\uptheta\otimes \id+\id\otimes \Idelta'$, we deduce an isomorphism 
\begin{equation}\label{p2-cmupiso28ii20}
\IrK^{\bullet}\stackrel{\sim}{\rightarrow}\mK^\bullet(\cN\otimes_{\co_{X'}}\IC^\dagger).
\end{equation}

For any rational number $r\geq 0$, we consider the $\co_{\fX'}$-algebra $\hcC^{(r)}$ associated with the extension \eqref{p2-cmupiso21a} 
(see \ref{p2-cmupiso32}) and its $\co_{\fX'}$-derivation \eqref{p2-cmupiso32b}
\begin{equation}\label{p2-cmupiso28ii200}
\delta_{\hcC^{(r)}}\colon \hcC^{(r)}\rightarrow f^*(\cE)\otimes_{\co_{\fX'}}\hcC^{(r)}.
\end{equation}
Identifying $f^*(\cE)$ with $\cE'_0$ by the isomorphism induced by $u$, 
the latter induces an $\co_{\fX'}$-derivation 
\begin{equation}\label{p2-cmupiso28ii201}
\delta'_{0,\hcC^{(r)}}\colon \hcC^{(r)}\rightarrow \hcC^{(r)}\otimes_{\co_{\fX'}}\cE'_0. 
\end{equation}

By \ref{p1-tshbn16} and \ref{p1-thbn21}, since  $N$ is $\co_{\fX'}[\frac 1 p]$-flat,
$(N,\theta_0)$ is twistable by the extension \eqref{p2-cmupiso21a} in the sense of \ref{p1-thbn14}. 
Therefore, there exist a rational number $r>0$, an $\co_{\fX'}$-module $M$ and an isomorphism 
of $\hcC^{(r)}$-modules with $\delta'_{0,\hcC^{(r)}}$-connection 
\begin{equation}\label{p2-cmupiso28ii7}
M\otimes_{\co_{\fX'}}\hcC^{(r)}\stackrel{\sim}{\rightarrow}N\otimes_{\co_{\fX'}}\hcC^{(r)},
\end{equation} 
where the $\delta'_{0,\hcC^{(r)}}$-connections are defined as in \ref{p1-delta-con4}, 
$N$ (resp.\ $M$) being endowed with the Higgs field $\theta_0$ (resp.\ $0$).  
Observe that $M$ is in fact a coherent $\co_{\fX'}[\frac 1 p]$-module \eqref{p1-thbn15}. 
By \ref{p1-delta-con5}, for every rational number $r'$ such that $r\geq r' >0$, \eqref{p2-cmupiso28ii7} induces 
an isomorphism of $\hcC^{(r')}$-modules with $\delta'_{0,\hcC^{(r')}}$-connection, 
\begin{equation}\label{p2-cmupiso28ii8}
M\otimes_{\co_{\fX'}}\hcC^{(r')}\stackrel{\sim}{\rightarrow}N\otimes_{\co_{\fX'}}\hcC^{(r')}. 
\end{equation}  
By \ref{p1-thbn15} and \ref{p1-thbn25}(i), since  $N$ is $\co_{\fX'}[\frac 1 p]$-flat, the latter induces an isomorphism of $\co_{\fX'}$-modules
\begin{equation}\label{p2-cmupiso28ii80}
M\stackrel{\sim}{\rightarrow}(N\otimes_{\co_{\fX'}}\hcC^{(r')})^{\vartheta^{(r')}_0=0},
\end{equation}
where $\vartheta^{(r')}_0=\theta_0\otimes \id+\id\otimes \delta'_{0,\hcC^{(r')}}$ is the total Higgs field on $N\otimes_{\co_{\fX'}}\hcC^{(r')}$.
Moreover, $M$ coincides with the module defined in \eqref{p2-cmupiso28i5}. We equip it with the Higgs field $\lambda_1$
defined in \eqref{p2-cmupiso28i7}. We set $(\cM,\uplambda_1)=\upalpha^\coh_{\co_{\fX'}}(M,\lambda_1)$ \eqref{p2-cmupiso31h}, 
and denote by $(\mK^\bullet(\cM),\uplambda_1^\bullet)$ its Dolbeault complex. 

By \eqref{p2-cmupiso310g}, the isomorphisms \eqref{p2-cmupiso28ii8} induce an isomorphism of ind-$\IC^\dagger$-modules
\begin{equation}\label{p2-cmupiso28ii9}
\cM\otimes_{\co_{\fX'}}\IC^\dagger\stackrel{\sim}{\rightarrow}\cN\otimes_{\co_{\fX'}}\IC^\dagger, 
\end{equation} 
that is compatible with the Higgs $\co_{\fX'}$-fields indicated on the same line of the table
\begin{equation}
\begin{tabular}{|c|c|c|}
\hline
left hand side&  right  hand  side & coefficients\\
\hline
$\id\otimes \Idelta'_{0}$ &  $\rI\vartheta_{0}=\uptheta_0\otimes \id+\id\otimes \Idelta'_{0}$ & $\cE'_0$\\
\hline
$\uplambda_1\otimes \id$ & $\uptheta_1\otimes \id$ & $\cE'_1$\\
\hline
\end{tabular}
\end{equation}

The splitting $\rho$ of the extension \eqref{p2-cmupiso21a}
induces a morphism of ind-$\co_{\fX'}$-algebras $\rI\varrho^\dagger\colon \IC^\dagger\rightarrow \co_{\fX'}$ \eqref{p2-cmupiso26b}.
We denote by $\uuplambda_1\colon \cM\rightarrow \cM\otimes_{\co_{\fX'}}\ucE'$ (resp.\ $\uuptheta_1\colon \cN\rightarrow \cN\otimes_{\co_{\fX'}}\ucE'$)  
the Higgs $\co_{\fX'}$-field induced by $\uplambda_1$ (resp.\ $\uptheta_1$) and the isomorphism $\cE'_1\stackrel{\sim}{\rightarrow} \ucE'$ 
deduced from $u'$ \eqref{p2-cmupiso21b}. Observe that $(\cN,\uuptheta_1)=\upalpha^\coh_{\co_{\fX'}}(N,\utheta)$ \eqref{p2-cmupiso31h}.
Then, the diagram of Higgs ind-$\co_{\fX'}$-modules with coefficients in $\ucE'$ 
\begin{equation}\label{p2-cmupiso28ii90}
\xymatrix{
{(\cM,\uuplambda_1)}\ar[r]\ar[d]&{(\cM\otimes_{\co_{\fX'}}\IC^\dagger,\uuplambda_1\otimes \id)}\ar[d]\\
{\upalpha^\coh_{\co_{\fX'}}(N,\utheta)}&{(\cN\otimes_{\co_{\fX'}}\IC^\dagger,\uuptheta_1\otimes \id),}\ar[l]_-(0.45){\id\otimes\rI\varrho^\dagger}}
\end{equation}
where the left (resp.\ right) vertical arrow is induced by the isomorphism \eqref{p2-cmupiso28i18} (resp. 
\emergencystretch=1em
\eqref{p2-cmupiso28ii9}), is commutative.

We consider the double complex 
\begin{equation}\label{p2-cmupiso28ii10}
\IcL^{i,j}=\IC^\dagger\otimes_{\co_{\fX'}}\wedge^i\cE'_0\otimes_{\co_{\fX'}}\mK^j(\cM)=
\IC^\dagger\otimes_{\co_{\fX'}}\wedge^i\cE'_0\otimes_{\co_{\fX'}}\cM\otimes_{\co_{\fX'}}\wedge^j\cE'_1,
\end{equation}
where the differentials are defined by
\begin{eqnarray}
d_{\IcL,1}^{i,j}&=&\Idelta'^i_{0}\otimes \id_{\cM\otimes \wedge^j\cE'_1},\label{p2-cmupiso28ii11}\\
d_{\IcL,2}^{i,j}&=&\id_{\IC^\dagger\otimes_{\co_{\fX'}}\wedge^i\cE'_0}\otimes \uplambda^j_1.\label{p2-cmupiso28ii110}
\end{eqnarray}
The isomorphism \eqref{p2-cmupiso28ii9} induces an isomorphism of bicomplexes 
\begin{equation}\label{p2-cmupiso28ii12}
\IcL^{\bullet,\bullet}\stackrel{\sim}{\rightarrow}\IcK^{\bullet,\bullet}. 
\end{equation}
Denoting by $\IrL^\bullet$ the total complex associated with $\IcL^{\bullet,\bullet}$, we deduce an isomorphism 
\begin{equation}\label{p2-cmupiso28ii13}
\IrL^{\bullet}\stackrel{\sim}{\rightarrow}\IrK^{\bullet}. 
\end{equation}

We denote by $(\mK^\bullet(\cM),\uplambda_1^\bullet)$ the Dolbeault complex of $(\cM,\uplambda_1)$. We have a canonical morphism of complexes 
\begin{equation}\label{p2-cmupiso28ii14}
\rI\epsilon^\bullet\colon (\mK^\bullet(\cM),\uplambda_1^\bullet)\rightarrow (\IcL^{0,\bullet},d_{\IcL,2}^{0,\bullet}).
\end{equation}
It follows from \eqref{p2-cmupiso28ii11} and the fact that $\Idelta'_{0}(\co_{\fX'})=0$, 
that the composition of $\rI\epsilon^\bullet$ with the morphism 
\begin{equation}\label{p2-cmupiso28ii15}
d_{\IcL,1}^{0,\bullet}\colon (\IcL^{0,\bullet},d_{\IcL,2}^{0,\bullet})\rightarrow (\IcL^{1,\bullet},d_{\IcL,2}^{1,\bullet})
\end{equation}
vanishes. We deduce a morphism 
\begin{equation}\label{p2-cmupiso28ii16}
\rI\varepsilon^\bullet\colon \mK^\bullet(\cM)\rightarrow \IrL^{\bullet}.
\end{equation}
By \ref{p2-cmupiso34} below, the mapping cone of $\rI\epsilon^\bullet$ is homotopically equivalent to zero. 
We take then for \eqref{p2-cmupiso28b} the morphism induced by \eqref{p2-cmupiso28ii16} and the isomorphisms \eqref{p2-cmupiso28ii20}, \eqref{p2-cmupiso28i18} 
and \eqref{p2-cmupiso28ii13}. An explicit computation \eqref{p2-cmupiso28ii90} shows that the composition of morphisms of complexes
\begin{equation}
\mK^\bullet(\cM)\stackrel{\rI\varepsilon^\bullet}{\longrightarrow} \IrL^\bullet \longrightarrow \IrK^\bullet \longrightarrow 
\mK^\bullet(\cN\otimes_{\co_{X'}}\IC^\dagger)\stackrel{\rI\upbeta}{\longrightarrow} \upalpha^\coh_{\co_{\fX'}}(\umK^\bullet(N)),
\end{equation}
where the second (resp.\ third) arrow is the ismorphism \eqref{p2-cmupiso28ii13} (resp.\ \eqref{p2-cmupiso28ii20}), is the morphism of Dolbeault 
complexes induced by $\iota$ \eqref{p2-cmupiso28i18}, which proves the commutativity of \eqref{p2-cmupiso28bb}.

\begin{lem}\label{p2-cmupiso34}
We keep the assumption of \ref{p2-cmupiso28} and the notation of its proof. 
\begin{itemize}
\item[{\rm (i)}] The mapping cone of the morphism $\varepsilon^\bullet\colon \mK^\bullet(M)\rightarrow \rL^{\bullet}$ \eqref{p2-cmupiso28i16} is homotopically equivalent to zero. 
\item[{\rm (ii)}] The mapping cone of the morphism $\rI\varepsilon^\bullet\colon \mK^\bullet(\cM)\rightarrow \IrL^{\bullet}$ \eqref{p2-cmupiso28ii16} is homotopically equivalent to zero. 
\end{itemize}
\end{lem}

(i) By \ref{p1-bcim21}, the mapping cone of $\varepsilon^\bullet$ is canonically isomorphic to the total complex associated with the double complex 
$(\tcL^{\bullet,\bullet},d_{\tcL,1}^{\bullet,\bullet},d_{\tcL,2}^{\bullet,\bullet})$ defined by 
\begin{equation}
\begin{array}{clcr}
\tcL^{i,j}=\cL^{i,j},&d_{\tcL,1}^{i,j}=d_{\cL,1}^{i,j},&d_{\tcL,2}^{i,j}=d_{\cL,2}^{i,j},& \forall i\in \mZ_{\geq 0},  \forall j\in \mZ,\\
\tcL^{-1,j}=\mK^i(M),&d_{\tcL,1}^{-1,j}=\varepsilon^j,&d_{\tcL,2}^{-1,j}=\lambda^j_1,& \forall j\in \mZ,\\
\tcL^{i,j}=0,&&&\forall i\in \mZ_{\leq -2},  \forall j\in \mZ.
\end{array}
\end{equation}
For any integers $i,j$, we consider the morphism 
\begin{equation}
h^{i,j}=h^i\otimes \id_{M\otimes \wedge^j\cE'_1}\colon \tcL^{i,j}\rightarrow \tcL^{i-1,j}, 
\end{equation}
where $(h^\bullet)$ is the homotopy given in \ref{p2-cmupiso33}(i). We immediately see that the diagram 
\begin{equation}
\xymatrix{
{\tcL^{i,j+1}}\ar[r]^-(0.5){h^{i,j+1}}&{\tcL^{i-1,j+1}}\\
{\tcL^{i,j}}\ar[r]^-(0.5){h^{i,j}}\ar[u]^{d_{\tcL,2}^{i,j}}&{\tcL^{i-1,j}}\ar[u]_{d_{\tcL,2}^{i-1,j}}}
\end{equation}
is commutative.
The proposition follows then from \ref{p1-bcim20}.  

(ii) We consider the ind-$\co_{\fX'}$-algebra
\begin{equation}\label{p2-cmupiso28ii5}
\IC^\dagger_\mQ=\underset{\underset{r\in \mQ_{>0}}{\longrightarrow}}{\mlq\mlq\lim \mrq\mrq}\ \upalpha_{\co_{\fX'}}(\hcC^{(r)}_\mQ),
\end{equation}
that we equip with the $\co_{\fX'}$-derivation
\begin{equation}\label{p2-cmupiso28ii6}
\Idelta'_{0,\mQ}=\underset{\underset{r\in \mQ_{>0}}{\longrightarrow}}{\mlq\mlq\lim \mrq\mrq}\ \upalpha_{\co_{\fX'}}(\delta'_{0,\hcC^{(r)},\mQ})\colon 
\IC^\dagger_\mQ \rightarrow \cE'_0\otimes_{\co_{\fX'}}\IC^\dagger_\mQ,
\end{equation}
where $\delta'_{0,\hcC^{(r)}}$ is the derivation defined in \eqref{p2-cmupiso28ii201}. 
We denote by $(\mK^\bullet(\IC^\dagger_\mQ),\Idelta'^\bullet_{0,\mQ})$ the Dolbeault complex of 
$(\IC^\dagger_\mQ,\Idelta'_{0,\mQ})$. 
By \eqref{p1-bcim5d}, we can replace in the definition \eqref{p2-cmupiso28ii10} of $\IcL^{\bullet,\bullet}$, 
$\IC^\dagger$ by $\IC^\dagger_\mQ$ and $\Idelta'_0$ \eqref{p2-cmupiso28dd} by $\Idelta'_{0,\mQ}$. 
The proposition follows then from \ref{p1-bcim20}, \ref{p1-bcim21} and \ref{p2-cmupiso33}(ii), as in the proof of (i).

\begin{cor}\label{p2-cmupiso29}
We take the assumptions and notation of \ref{p2-cmupiso21} and assume, moreover, that we are given
a splitting $\rho\colon \cF\rightarrow \co_{\fX'}$ of the extension \eqref{p2-cmupiso21a}. 
Let $(N,\theta)$ be a locally CL-small Higgs $\co_{\fX'}[\frac 1 p]$-module with coefficients in $\cE'$ \eqref{p1-tshbn13}.  Then, 
\begin{itemize}
\item[{\rm (i)}] The morphism of complexes of $\co_{\fX'}$-modules \eqref{p2-cmupiso26d}
\begin{equation}\label{p2-cmupiso29a}
\upbeta\colon \mK^\bullet(N\otimes_{\co_{\fX'}}\cC^\dagger) \rightarrow \umK^\bullet(N)
\end{equation}
is a quasi-isomorphism. 
\item[{\rm (ii)}]  Assume that the $\co_{\fX'}[\frac 1 p]$-module $N$ is flat. Then, the morphism of complexes of ind-$\co_{\fX'}$-modules \eqref{p2-cmupiso26f}
\begin{equation}\label{p2-cmupiso29b}
\rI\upbeta\colon \mK^\bullet(\cN\otimes_{\co_{\fX'}}\IC^\dagger) \rightarrow \upalpha^\coh_{\co_{\fX'}}(\umK^\bullet(N)) 
\end{equation}
is a quasi-isomorphism. 
\end{itemize}
\end{cor}

Indeed, the questions being local, we may assume that the extension \eqref{p2-cmupiso21b} is also split and 
that the $\co_{\fX'}$-module $f^*(\cE)$ is free. The proposition follows then from the commutative diagrams \eqref{p2-cmupiso28aa} and \eqref{p2-cmupiso28bb},
whose horizontal arrows are homotopically equivalent to zero by \ref{p2-cmupiso28}.  

\begin{cor}\label{p2-cmupiso36}
Under the assumptions of \ref{p2-cmupiso29}, for every $q\geq 0$, 
the morphism $\upbeta$ \eqref{p2-cmupiso26d} induces an isomorphism of Higgs $\co_{\fX}$-modules
\begin{equation}\label{p2-cmupiso36a}
\rR^qf^\uptau_*(N,\theta) \stackrel{\sim}{\rightarrow} (\rR^qf_*(\umK^\bullet(N)),\kappa^q), 
\end{equation}
where the source is defined in \eqref{p2-cmupiso22d} and the target is   
equipped with the Katz-Oda field \eqref{p1-tphdi4h}.
\end{cor}

It follows from \ref{p2-cmupiso29} and \ref{p1-tphdi10}.

\begin{teo}\label{p2-cmupiso30}
We keep the assumptions and notation of \ref{p2-cmupiso21} and assume moreover that the morphism $f\colon \fX'\rightarrow \fX$ is proper. 
Let $(N,\theta)$ be a {\em locally CL-small} Higgs $\co_{\fX'}[\frac 1 p]$-module with coefficients in $\cE'$ \eqref{p1-tshbn13}, 
$q$ an integer. We set $(\cN,\uptheta)=\upalpha^\coh_{\co_{\fX'}}(N,\theta)$ \eqref{p2-cmupiso31h}. 
We denote by $\mK^{\bullet}(N\otimes_{\co_{\fX'}}\cC^\dagger)$ the Dolbeault complex of $N\otimes_{\co_{\fX'}}\cC^\dagger$ equipped 
with the total Higgs field $\vartheta = \theta\otimes \id+\id\otimes \delta'$ \eqref{p2-cmupiso22a}
and by $\mK^{\bullet}(\cN\otimes_{\co_{\fX'}}\IC^\dagger)$ the Dolbeault complex of $\cN\otimes_{\co_{\fX'}}\IC^\dagger$ equipped 
with the total Higgs field $\rI\vartheta= \uptheta\otimes \id+\id\otimes \Idelta'$ \eqref{p2-cmupiso23a}. 
We assume, for {\rm (ii)}, {\rm (iii)} and {\rm (iv)}, that the $\co_{\fX'}[\frac 1 p]$-module $N$ is flat. Then, 
\begin{itemize}
\item[{\rm (i)}] The $\co_\fX[\frac 1 p]$-module $\rR^qf_*(\mK^\bullet(N\otimes_{\co_{\fX'}}\cC^\dagger))$ is coherent. 
\item[{\rm (ii)}] We have a canonical functorial isomorphism of ind-$\co_{\fX}$-modules
\begin{equation}\label{p2-cmupiso30a}
\upalpha^\coh_{\co_{\fX}}(\rR^qf_*(\mK^\bullet(N\otimes_{\co_{\fX'}}\cC^\dagger)))\stackrel{\sim}{\rightarrow}
\rR^q\rI f_*(\mK^\bullet(\cN\otimes_{\co_{\fX'}}\IC^\dagger)),
\end{equation}
where the functor $\upalpha^\coh_{\co_{\fX}}$ is defined in \eqref{p2-cmupiso2c}.
\item[{\rm (iii)}] The isomorphism \eqref{p2-cmupiso30a} is underlying an isomorphism of Higgs ind-$\co_{\fX}$-modules with coefficients in $\cE$, 
\begin{equation}\label{p2-cmupiso30b}
\upalpha^\coh_{\co_{\fX}}(\rR^qf^\uptau_*(N,\theta))\stackrel{\sim}{\rightarrow}
\rR^q\rI f^\uptau_*(\cN,\uptheta),
\end{equation}
where the functors $\upalpha^\coh_{\co_{\fX}}$, $\rR^qf^\uptau_*$ and $\rR^q\rI f^\uptau_*$ are defined in \eqref{p2-cmupiso31h}, 
\eqref{p2-cmupiso22d} and 
\emergencystretch=1em
\eqref{p2-cmupiso23d}, respectively. 
\item[{\rm (iv)}] The base change morphism \eqref{p2-cmupiso14d}
\begin{equation}\label{p2-cmupiso30c}
\rI\uplambda^*(\rR^q\rI f_*(\mK^\bullet(\cN\otimes_{\co_{\fX'}}\IC^\dagger)))\rightarrow \rR^q\rI\bvf_*(\rI\uplambda'^*(\mK^\bullet(\cN\otimes_{\co_{\fX'}}\IC^\dagger)))
\end{equation}
is an isomorphism of ind-$\co_{\bvfX}$-modules. 
\end{itemize}
\end{teo}

(i) It is a consequence of \ref{p2-cmupiso28}(i), \ref{p2-cmupiso12}(i) and \ref{p1-bcim16}(ii)
applied to the thick subcategory 
$\bMod^\coh(\co_\fX[\frac 1 p])$ of $\bMod(\co_\fX)$.  

(ii) Let $\fN$ be a coherent $\co_{\fX'}$-lattice of $N$, i.e. a coherent sub-$\co_{\fX'}$-module of $N$ which generates it over $\co_{\fX'}[\frac 1 p]$, 
$t$ a rational number $\geq 0$ such that $\theta(\fN)\subset p^{-t}\fN\otimes_{\co_{\fX'}}\cE'$. 
For any rational number $r>0$, we consider the complex of $\co_{\fX'}$-modules 
\begin{equation}
p^{-t\bullet} \fN\otimes_{\co_{\fX'}} \hcC^{(r)}\otimes_{\co_{\fX'}} \wedge^\bullet \cE',
\end{equation}
where the differentials are induced by the morphism 
\begin{equation}
\theta\otimes \id+\id\otimes \delta_{\hcC^{(r)}}\colon  \fN\otimes_{\co_{\fX'}} \hcC^{(r)}\rightarrow 
p^{-t} \fN\otimes_{\co_{\fX'}}\hcC^{(r)}\otimes_{\co_{\fX'}}\cE'. 
\end{equation}
These complexes form naturally an inductive system indexed by the ordered set $\mQ_{>0}$ \eqref{p1-thbn6}. 
We denote by $\Delta$ the multiplicative monoid $\mZ-\{0\}$ and by $\uDelta$
the filtered category whose objects are the elements of $\Delta$ and the morphisms are determined by the divisibility relation in $\Delta$.
Let $\upmu\colon \uDelta \rightarrow \bMod(\co_\fX)$ be the functor which sends a nonzero integer $n$ to $\fN$
and a morphism $m|n$ of $\uDelta$ to the isogeny of multiplication by $n/m$ on $\fN$. 
By (\cite{ag2} 2.9.3(i)), we have a canonical isomorphism of ind-$\co_{\fX'}$-modules
\begin{equation}
\cN\stackrel{\sim}{\rightarrow}
\underset{\underset{n\in \uDelta}{\longrightarrow}}{\mlq\mlq\lim \mrq\mrq}\fN.
\end{equation}
By (\cite{ag2} 2.7.3), we deduce an isomorphism 
\begin{equation}
\mK^\bullet(\cN\otimes_{\co_{\fX'}}\IC^\dagger)\stackrel{\sim}{\rightarrow}
\underset{\underset{r\in \mQ_{>0}; n\in \uDelta}{\longrightarrow}}{\mlq\mlq\lim \mrq\mrq}
p^{-t\bullet} \fN \otimes_{\co_{\fX'}} \hcC^{(r)}\otimes_{\co_{\fX'}} \wedge^\bullet \cE'.
\end{equation}

By \ref{p2-cmupiso11} and \ref{p2-cmupiso28}(ii), for every integer $q$, $\cH^q(\mK^\bullet(\cN\otimes_{\co_{\fX'}}\IC^\dagger))$ is in the essential image of the functor
$\upalpha_{\co_{\fX'}}^\coh$. The proposition follows then from \ref{p2-cmupiso13}(ii) and \eqref{p2-cmupiso24b}. 

(iii) It follows from (ii), the functoriality of \eqref{p2-cmupiso30a} and \ref{p2-cmupiso24}. 

(iv) It follows from \ref{p2-cmupiso17}, applied to $\mK^\bullet(\cN\otimes_{\co_{\fX'}}\IC^\dagger)$, whose conditions are satisfied by 
\ref{p2-cmupiso6}, \ref{p2-cmupiso20} and \ref{p2-cmupiso28}(ii).

\subsection{}\label{p2-cmupiso40}
We take again the assumptions and notation of \ref{p2-cmupiso21}, and set $\ccE=\cHom_{\co_\fX}(\cE,\co_\fX)$. 
Let $\beta,s$ be rational numbers such that $0\leq \beta<s+\frac{1}{p-1}$ and $s\geq 0$. 
We set 
\begin{equation}\label{p2-cmupiso40a}
\gamma=\inf\{ v_p(n!)-(\beta-s)n; n\in \mN\},
\end{equation}
which is a rational number $\leq 0$. By \ref{p1-tshbn50}, we have a commutative diagram of $\co_{\fX'}$-linear morphisms
\begin{equation}\label{p2-cmupiso40b}
\xymatrix{
{\rS(f^*(\ccE))\otimes_{\co_{\fX'}}\cC^{(s)}}\ar[r]^-(0.5){\mu^s}\ar[d]&{\cC^{(s)}}\ar[d]\\
{\rS(p^{-\beta}f^*(\ccE))\otimes_{\co_{\fX'}}\cC^{(s)}}\ar[r]^-(0.5){\mu^s_\beta}&{p^{\gamma d}\cC^{(s)},}}
\end{equation}
where $\rS=\rS_{\co_{\fX'}}$ denotes the symmetric $\co_{\fX'}$-algebra \eqref{p1-NC7}, $\cC^{(s)}$ is the Higgs--Tate 
$\co_{\fX'}$-algebra of thickness $s$ associated with the extension \eqref{p2-cmupiso21a},  
$\mu^s$ is defined by the Higgs $\co_{\fX'}$-field 
$\delta_{\cC^{(s)}}$ on $\cC^{(s)}$ \eqref{p1-thbn3d} and the vertical arrows are the canonical morphisms. 
Moreover, $\mu^s_\beta$ is uniquely determined by \eqref{p2-cmupiso40b}. 
It follows immediately from \ref{p1-delta-con7}(i) that the diagram 
\begin{equation}\label{p2-cmupiso40k}
\xymatrix{
{\rS(p^{-\beta}f^*(\ccE))\otimes_{\co_{\fX'}}\cC^{(s)}}\ar[r]^-(0.5){\mu^s_\beta}\ar[d]_{\id\otimes \delta_{\cC^{(s)}}}&
{p^{\gamma d}\cC^{(s)}}\ar[d]^{\delta_{\cC^{(s)}}}\\
{\rS(p^{-\beta}f^*(\ccE))\otimes_{\co_{\fX'}}\cC^{(s)}\otimes_{\co_{\fX'}}f^*(\ccE)}\ar[r]^-(0.5){\mu^s_\beta}&
{p^{\gamma d}\cC^{(s)}\otimes_{\co_{\fX'}}f^*(\ccE)}}
\end{equation}
is commutative. 

We denote by $\hrS(f^*(\ccE))$ (resp.\ $\hrS(p^{-\beta}f^*(\ccE))$) the $p$-adic completion of the $\co_{\fX'}$-algebra 
$\rS(f^*(\ccE))$ (resp.\ $\rS(p^{-\beta}f^*(\ccE))$). 
Since $\cC^{(s)}$ is $p$-torsion free, extending $\mu^s$ and $\mu^s_\beta$ to the $p$-adic completions, restricting to 
the usual tensor products and inverting $p$, we obtain a commutative diagram of $\co_{\fX'}$-linear morphisms
\begin{equation}\label{p2-cmupiso40c}
\xymatrix{
{\hrS(f^*(\ccE))\otimes_{\co_{\fX'}}\hcC^{(s)}}\ar[r]^-(0.5){\hmu^s}\ar[d]&{\hcC^{(s)}}\ar[d]\\
{\hrS(p^{-\beta}f^*(\ccE))\otimes_{\co_{\fX'}}\hcC^{(s)}[\frac 1 p]}\ar[r]^-(0.5){\hmu^s_\beta}&{\hcC^{(s)}[\frac 1 p],}}
\end{equation}
where the vertical arrows are the canonical morphisms. We immediately check that the morphisms $\hmu^s$ and $\hmu^s_\beta$ are compatible when 
$s$ and $\beta$ vary. It follows from \eqref{p2-cmupiso40k} that the Higgs field $\delta_{\hcC^{(s)}}\otimes_{\mZ_p}\mQ_p$ \eqref{p1-thbn6b} 
is $\hrS(p^{-\beta}f^*(\ccE))$-linear.

Let $r$ be a rational number such that $r\geq 0$ and $\beta\leq r+ \frac{1}{p-1}$. We consider the $\co_{\fX'}$-algebra \eqref{p1-thbn6e}
\begin{equation}\label{p2-cmupiso40i}
\hcC^{(r+)}=\underset{\underset{s\in \mQ_{>r}}{\longrightarrow}}{\lim}\ \hcC^{(s)},
\end{equation}
equipped with the $\co_{\fX'}$-derivation \eqref{p1-thbn6f}
\begin{equation}\label{p2-cmupiso40j}
\delta_{\hcC^{(r+)}}\colon \hcC^{(r+)}\rightarrow f^*(\cE)\otimes_{\co_{\fX'}}\hcC^{(r+)}.
\end{equation}

We deduce from \eqref{p2-cmupiso40c} a commutative diagram of $\co_{\fX'}$-linear morphisms 
\begin{equation}\label{p2-cmupiso40d}
\xymatrix{
{\hrS(f^*(\ccE))\otimes_{\co_{\fX'}}\hcC^{(r+)}}\ar[r]^-(0.5){\hmu^{r+}}\ar[d]&{\hcC^{(r+)}}\ar[d]\\
{\hrS(p^{-\beta}f^*(\ccE))\otimes_{\co_{\fX'}}\hcC^{(r+)}[\frac 1 p]}\ar[r]^-(0.5){\hmu^{r+}_\beta}&{\hcC^{(r+)}[\frac 1 p],}}
\end{equation}
where the vertical arrows are the canonical morphisms. Observe that $\hmu^{r+}$ is compatible with the Higgs field $\delta_{\hcC^{(r+)}}$. 
It follows from above that the Higgs field $\delta_{\hcC^{(r+)}}\otimes_{\mZ_p}\mQ_p$ \eqref{p2-cmupiso40j} 
is $\hrS(p^{-\beta}f^*(\ccE))$-linear.

Assume that $\beta>\frac{1}{p-1}$. The canonical homomorphism $\rS(f^*(\ccE))\rightarrow \rS(p^{-\beta}f^*(\ccE))$ factors through 
the canonical homomorphism $\rS(f^*(\ccE))\rightarrow \Gamma(f^*(\ccE))$ \eqref{p1-NC7}, and induces a homomorphism of $\co_{\fX'}$-algebras
\begin{equation}\label{p2-cmupiso40e}
\psi^{\beta}\colon \Gamma(f^*(\ccE))\rightarrow \rS(p^{-\beta}f^*(\ccE)). 
\end{equation}
For every integer $n\geq 0$, we have $\psi^{\beta}(\Gamma^{\geq n}(f^*(\ccE))\subset 
p^{\lfloor \frac{n}{d}\rfloor r} \rS(p^{-\beta}f^*(\ccE))$ \eqref{p1-NC7}, where $d$ is an integer such that $\rk_{\co_{\fX}}(\cE)\leq d$. 
Therefore, $\psi^{\beta}$ induces a homomorphism of $\co_{\fX'}$-algebras 
\begin{equation}\label{p2-cmupiso40f}
\Psi^{\beta}\colon \hGamma(f^*(\ccE))\rightarrow \hrS(p^{-\beta}f^*(\ccE)).
\end{equation}

For any $\phi\in\Gamma(\fX',f^*(\ccE))$, we denote by $\exp(\phi)\in \hGamma(\fX',\hGamma(f^*(\ccE)))$ the image of $\phi$
by the exponential homomorphism \eqref{p1-NC7c}, and also (abusively) by 
\begin{equation}\label{p2-cmupiso40g}
\exp(\phi)\in \Gamma(\fX',\hrS(p^{-\beta}f^*(\ccE)))
\end{equation}
its image by the homomorphism $\Psi^{\beta}$ \eqref{p2-cmupiso40f}. 

\begin{lem}\label{p2-cmupiso45} 
Let $f\colon \fX'\rightarrow \fX$ be a morphism of finite presentation,  $\cM$ a locally free $\co_\fX$-module of finite type. Then, 
\begin{itemize}
\item[{\rm (i)}] The diagram 
\begin{equation}\label{p2-cmupiso45a} 
\xymatrix{
{f^{-1}(\cM)}\ar[rr]^-(0.5){f^{-1}(\exp_{\cM})}\ar[d]&&{f^{-1}(\hGamma(\cM))}\ar[d]\\
{f^*(\cM)}\ar[rr]^-(0.5){\exp_{f^*(\cM)}}&&{\hGamma(f^*(\cM)),}}
\end{equation}
where the vertical arrows are the canonical morphisms and $\exp$ is the exponential homomorphism \eqref{p1-NC7c}, is commutative. 
\item[{\rm (ii)}] For every rational number $\beta>\frac{1}{p-1}$, the diagram of canonical morphisms
\begin{equation}\label{p2-cmupiso45b} 
\xymatrix{
{f^{-1}(\hGamma(\cM))}\ar[d]\ar[r]&{f^{-1}(\hrS(p^{-\beta}\cM))}\ar[d]\\
{\hGamma(f^*(\cM))}\ar[r]&{\hrS(p^{-\beta}f^*(\cM))}}
\end{equation}
is commutative; see \eqref{p2-cmupiso40f} for the definition of the horizontal arrows. 
\end{itemize}
\end{lem}

Obvious. 

\begin{prop}\label{p2-cmupiso41}
We take again the assumptions and notation of \ref{p2-cmupiso21}, and set $\ccE=\cHom_{\co_\fX}(\cE,\co_\fX)$. 
Let $\beta,s$ be rational numbers such that $0<\beta-\frac{1}{p-1}<s$,
and let $\rho_1,\rho_2\colon \cF\rightarrow \co_{\fX'}$ be two splittings of the extension \eqref{p2-cmupiso21a}. 
We denote by $\phi\colon f^*(\cE)\rightarrow \co_{\fX'}$ the $\co_{\fX'}$-linear form induced by $\rho_1-\rho_2$,
and by $\hvarrho^{(s)}_i\colon \hcC^{(s)}\rightarrow \co_{\fX'}$ the homomorphism of $\co_{\fX'}$-algebras 
induced by $\rho_i$, for $i=1,2$ \eqref{p1-thbn16}. Then, the diagram 
\begin{equation}\label{p2-cmupiso41a}
\xymatrix{
{\hcC^{(s)}[\frac 1 p]}\ar[r]^-(0.5){\hvarrho^{(s)}_1}\ar[d]_{\exp(\phi)}&{\co_{\fX'}[\frac 1 p]}\\
{\hcC^{(s)}[\frac 1 p]}\ar[ru]_-(0.5){\hvarrho^{(s)}_2}&}
\end{equation}
where $\exp(\phi)\in \Gamma(\fX',\hrS(p^{-\beta}f^*(\ccE)))$ is defined in \eqref{p2-cmupiso40g} and acts on $\hcC^{(s)}[\frac 1 p]$
via $\hmu^s_\beta$ \eqref{p2-cmupiso40c}, is commutative. 
\end{prop}

Observe first that $\exp(\phi)(\hcC^{(s)})\subset p^{\gamma d}\hcC^{(s)}$, where $\gamma \in \mQ_{\leq 0}$ 
and $d\in \mZ_{\geq 1}$ are as in \ref{p2-cmupiso40}. Moreover, we have $\exp(\phi)(\cC^{(s)})\subset p^{\gamma d}\cC^{(s)}$, so 
$\exp(\phi)\colon \hcC^{(s)}\rightarrow p^{\gamma d}\hcC^{(s)}$ is the $p$-adic extension of its restriction that we also denote by 
$\exp(\phi)\colon \cC^{(s)}\rightarrow p^{\gamma d}\cC^{(s)}$. 
It is therefore enough to prove that the diagram 
\begin{equation}\label{p2-cmupiso41b}
\xymatrix{
{\cC^{(s)}[\frac 1 p]}\ar[r]^-(0.5){\varrho^{(s)}_1}\ar[d]_{\exp(\phi)}&{\co_{\fX'}[\frac 1 p]}\\
{\cC^{(s)}[\frac 1 p]}\ar[ru]_-(0.5){\varrho^{(s)}_2}&}
\end{equation}
where $\varrho^{(s)}_i$ is the homomorphism of $\co_{\fX'}$-algebras induced by $\rho_i$, for $i=1,2$, is commutative. 
The canonical homomorphism $\alpha^{s,0}$ \eqref{p1-thbn4b} induces an isomorphism $\cC^{(s)}[\frac 1 p]\stackrel{\sim}{\rightarrow}\cC[\frac 1 p]$. 
We set $\partial_\phi=(\phi\otimes \id_\cC)\circ d_\cC$ as an $\co_{\fX'}$-linear endomorphism $\cC$, 
where $d_\cC$ is the universal $\co_{\fX'}$-derivation of $\cC$ \eqref{p1-thbn3c}. 
Then, the $\co_{\fX'}$-linear morphism 
\begin{equation}\label{p2-cmupiso41c}
\exp(\phi)=\sum_{n\geq 0} \frac{1}{n!} \partial_\phi^n\colon \cC[\frac 1 p]\rightarrow \cC[\frac 1 p],
\end{equation} 
is well defined and it coincides with the endomorphism of $\cC^{(s)}[\frac 1 p]$ denoted by $\exp(\phi)$ defined above. 
Hence, we are reduced to proving a Taylor formula, namely, that the diagram 
\begin{equation}\label{p2-cmupiso41d}
\xymatrix{
{\cC[\frac 1 p]}\ar[r]^-(0.5){\varrho_1}\ar[d]_{\exp(\phi)}&{\co_{\fX'}[\frac 1 p]}\\
{\cC[\frac 1 p]}\ar[ru]_-(0.5){\varrho_2}&}
\end{equation}
where $\varrho_i$ is the homomorphism of $\co_{\fX'}$-algebras induced by $\rho_i$, for $i=1,2$, is commutative.

The question being local on $\fX'$, we may assume that the $\co_{\fX'}$-module $f^*(\cE)$ is free with basis $X_1,\dots,X_d$. 
We denote by $\partial_1,\dots,\partial_d$ the dual basis of the dual $\co_{\fX'}$-module $f^*(\ccE)$. 
The splitting $\rho_2$ of \eqref{p1-thbn1a} induces an isomorphism of $\co_{\fX'}$-algebras $\rS(f^*(\cE))\stackrel{\sim}{\rightarrow} \cC$. 
We deduce an isomorphism of $\co_{\fX'}$-algebras $\cC\stackrel{\sim}{\rightarrow} \co_{\fX'}[X_1,\dots,X_d]$. The homomorphism  
$\varrho_2\colon \cC\rightarrow \co_{\fX'}$ identifies with the homomorphism 
$\co_{\fX'}[X_1,\dots,X_d]\rightarrow \co_{\fX'}$ that maps $F$ to $F(0)$. 
We write $\phi=\sum_{1\leq i\leq d}a_i\partial_i \in \Gamma(\fX',f^*(\ccE))$, where $\ua=(a_1,\dots,a_d)\in \co_{\fX'}(\fX')^d$. 
The homomorphism $\varrho_1\colon \cC\rightarrow \co_{\fX'}$ identifies with the homomorphism 
$\co_{\fX'}[X_1,\dots,X_d]\rightarrow \co_{\fX'}$ that maps $F$ to $F(\ua)$. 
Moreover, 
\begin{equation}\label{p2-cmupiso41e}
\partial_\phi=(\phi\otimes \id_\cC)\circ d_\cC\colon \co_{\fX'}[X_1,\dots,X_d] \rightarrow \co_{\fX'}[X_1,\dots,X_d]
\end{equation}
is defined for any $F\in \co_{\fX'}[X_1,\dots,X_d]$ by 
\begin{equation}\label{p2-cmupiso41f}
\partial_\phi(F)=\sum_{1\leq i\leq d} a_i\frac{\partial }{\partial X_i}(F). 
\end{equation}
Setting $\uX=(X_1,\dots,X_d)$, we have the Taylor formula
\begin{equation}\label{p2-cmupiso41g}
F(\uX+\ua)=\sum_{n\geq 0} \frac{1}{n!} \partial_\phi^n(F)(\uX)=\exp(\phi)(F),
\end{equation}
the sum being finite, which implies, by evaluation at $\uX=0$, the commutativity of \eqref{p2-cmupiso41d}.

\begin{cor}\label{p2-cmupiso42}
We take again the assumptions and notation of \ref{p2-cmupiso21}, and set $\ccE=\cHom_{\co_\fX}(\cE,\co_\fX)$. 
Let $\beta,r$ be rational numbers such that $0<\beta-\frac{1}{p-1}\leq r$,
and let $\rho_1,\rho_2\colon \cF\rightarrow \co_{\fX'}$ be two splittings of the extension \eqref{p2-cmupiso21a}. 
We denote by $\phi\colon f^*(\cE)\rightarrow \co_{\fX'}$ the $\co_{\fX'}$-linear form induced by $\rho_1-\rho_2$,
and by $\hvarrho^{(r+)}_i\colon \hcC^{(r+)}\rightarrow \co_{\fX'}$ the homomorphism of $\co_{\fX'}$-algebras 
induced by $\rho_i$, for $i=1,2$ \eqref{p1-thbn16}. Then, the diagram 
\begin{equation}\label{p2-cmupiso42a}
\xymatrix{
{\hcC^{(r+)}[\frac 1 p]}\ar[r]^-(0.5){\hvarrho^{(r+)}_1}\ar[d]_{\exp(\phi)}&{\co_{\fX'}[\frac 1 p]}\\
{\hcC^{(r+)}[\frac 1 p]}\ar[ru]_-(0.5){\hvarrho^{(r+)}_2}&}
\end{equation}
where $\exp(\phi)\in \Gamma(\fX',\hrS(p^{-\beta}f^*(\ccE)))$ is defined in \eqref{p2-cmupiso40g} and acts on $\hcC^{(r+)}[\frac 1 p]$
via $\hmu^{r+}_\beta$ \eqref{p2-cmupiso40d}, is commutative. 
\end{cor}

\subsection{}\label{p2-cmupiso43}
We take again the assumptions and notation of \ref{p2-cmupiso21}. 
Let $(N,\theta)$ be a Higgs $\co_{\fX'}[\frac 1 p]$-module with coefficients in $\cE'$, $r$ a rational number $\geq 0$. 
We consider the $\co_{\fX'}$-algebra $\hcC^{(r+)}$ \eqref{p2-cmupiso40i}, equipped with the $\co_{\fX'}$-derivation 
\begin{equation}\label{p2-cmupiso43a}
\delta'_{\hcC^{(r+)}}=(u\otimes \id)\circ \delta_{\hcC^{(r+)}}\colon \hcC^{(r+)}\rightarrow \cE'\otimes_{\co_{\fX'}}\hcC^{(r+)},
\end{equation}
where $\delta_{\hcC^{(r+)}}$ is the $\co_{\fX'}$-derivation defined in \eqref{p2-cmupiso40j} and $u\colon f^*(\cE)\rightarrow \cE'$ 
is the morphism given in \eqref{p2-cmupiso21b}. 

We equip $N\otimes_{\co_{\fX'}}\hcC^{(r+)}$ with the Higgs 
$\co_{\fX'}$-field 
\begin{equation}\label{p2-cmupiso43b}
\vartheta^{(r+)}= \theta\otimes \id+\id\otimes \delta'_{\hcC^{(r+)}}\colon N\otimes_{\co_{\fX'}}\hcC^{(r+)} \rightarrow 
\cE'\otimes_{\co_{\fX'}}N\otimes_{\co_{\fX'}}\hcC^{(r+)},
\end{equation}
and denote by $\mK^{\bullet}(N\otimes_{\co_{\fX'}}\hcC^{(r+)})$ the associated Dolbeault complex. For $r=0$, we recover the Dolbeault complex 
$\mK^{\bullet}(N\otimes_{\co_{\fX'}}\cC^\dagger)$ defined in \ref{p2-cmupiso210}. 

For any rational number $r'$ such that $r\geq r'\geq 0$, we have a canonical homomorphism of $\co_{\fX'}$-algebras 
\begin{equation}\label{p2-cmupiso43c}
\talpha^{r,r'}\colon \hcC^{(r+)}\rightarrow \hcC^{(r'+)}
\end{equation} 
induced by the homomorphisms $\halpha^{s,s'}\colon \hcC^{(s)}\rightarrow \hcC^{(s')}$ \eqref{p1-thbn6c}, for rational numbers 
$s,s'$ such that $s>r$, $s'>r'$ and $s\geq s'$. It follows from \eqref{p1-thbn6d} that we have 
\begin{equation}\label{p2-cmupiso43d}
(\id \otimes \talpha^{r,r'}) \circ \delta'_{\hcC^{(r+)}}=\delta'_{\hcC^{(r'+)}} \circ \talpha^{r,r'}.
\end{equation}
We deduce a morphism of Dolbeault complexes 
\begin{equation}\label{p2-cmupiso43e}
\upiota^{r,r'}\colon \mK^{\bullet}(N\otimes_{\co_{\fX'}}\hcC^{(r+)})\rightarrow \mK^{\bullet}(N\otimes_{\co_{\fX'}}\hcC^{(r'+)}). 
\end{equation}

Let $\rho\colon \cF\rightarrow \co_{\fX'}$ be a splitting of the extension \eqref{p2-cmupiso21a}, 
$\hvarrho^{(r+)}\colon \hcC^{(r+)}\rightarrow \co_{\fX'}$ the induced homomorphism of $\co_{\fX'}$-algebras \eqref{p1-thbn16}. 
The diagram 
\begin{equation}\label{p2-cmupiso43f}
\xymatrix{
{N\otimes_{\co_{\fX'}}\hcC^{(r+)}}\ar[r]^-(0.5){\vartheta^{(r+)}}\ar[d]_{\id\otimes\hvarrho^{(r+)}}&
{\cE'\otimes_{\co_{\fX'}}N\otimes_{\co_{\fX'}}\hcC^{(r+)}}\ar[d]^{u'\otimes\id\otimes \hvarrho^{(r+)}}\\
{N}\ar[r]^-(0.5){\utheta}&{\ucE'\otimes_{\co_{\fX'}}N,}}
\end{equation}
where $\utheta$ is the Higgs field induced by $\theta$ and 
$u'$ is the morphism defined in \eqref{p2-cmupiso21b}, is commutative. 
We deduce a morphism of Dolbeault complexes 
\begin{equation}\label{p2-cmupiso43g}
\upbeta^{(r)}\colon \mK^\bullet(N\otimes_{\co_{\fX'}}\hcC^{(r+)}) \rightarrow \umK^\bullet(N). 
\end{equation}
We immediately check that we have 
\begin{equation}\label{p2-cmupiso43h}
\upbeta^{(r)}=\upbeta^{(r')}\circ \upiota^{r,r'}.
\end{equation}

\subsection{}\label{p2-cmupiso46}
We keep the assumptions and notation of \ref{p2-cmupiso43}, and set $\ccE=\cHom_{\co_\fX}(\cE,\co_\fX)$. 
By the same argument as in \ref{p1-tphdi6}, the Higgs field $\id\otimes \delta_{\hcC^{(r+)}}$ on $N\otimes_{\co_{\fX'}}\hcC^{(r+)}$ 
induces a morphism of complexes of $\co_{\fX'}$-modules
\begin{equation}\label{p2-cmupiso46a}
\fd^{(r+)} \colon \mK^\bullet(N\otimes_{\co_{\fX'}}\hcC^{(r+)})\rightarrow f^*(\cE)\otimes_{\co_{\fX'}} \mK^\bullet(N\otimes_{\co_{\fX'}}\hcC^{(r+)}).
\end{equation} 
By the projection formula, for any integer $q\geq 0$, $\rR^qf_*(\fd^{(r+)})$ identifies with an $\co_\fX$-linear morphism
\begin{equation}\label{p2-cmupiso46b}
\rR^qf_*(\fd^{(r+)})\colon  \rR^qf_*(\mK^\bullet(N\otimes_{\co_{\fX'}}\hcC^{(r+)}))\rightarrow \cE\otimes_{\co_\fX} 
\rR^qf_*(\mK^\bullet(N\otimes_{\co_{\fX'}}\hcC^{(r+)})). 
\end{equation}
By choosing locally a basis of $\cE$, we easily check that it is a Higgs $\co_{\fX}$-field. It hence defines an action of the $\co_\fX$-algebra 
$\rS(\ccE)=\rS_{\co_\fX}(\ccE)$ on $\rR^qf_*(\mK^{\bullet}(N\otimes_{\co_{\fX'}}\hcC^{(r+)}))$. 

Alternatively, this action can be defined as follows.
We consider the action of $\rS(f^*(\ccE))$ on $\hcC^{(r+)}$ defined by the Higgs field $\delta_{\hcC^{(r+)}}$ \eqref{p2-cmupiso40j}. 
By \ref{p1-delta-con7}(i), $\delta_{\hcC^{(r+)}}$ is $\rS(f^*(\ccE))$-linear, and hence so is $\vartheta^{(r+)}$ \eqref{p2-cmupiso43b}. 
Therefore, $\mK^{\bullet}(N\otimes_{\co_{\fX'}}\hcC^{(r+)})$ is naturally a complex of $\rS(f^*(\ccE))$-modules. 
We deduce an action of $f_*(\rS(f^*(\ccE)))$ on $\rR^qf_*(\mK^{\bullet}(N\otimes_{\co_{\fX'}}\hcC^{(r+)}))$. 
The two actions on $\rR^qf_*(\mK^{\bullet}(N\otimes_{\co_{\fX'}}\hcC^{(r+)}))$ defined above are compatible via the adjunction morphism 
$\rS(\ccE)\rightarrow f_*(\rS(f^*(\ccE)))$. 

For every rational number $r'$ such that $r\geq r'\geq 0$, we immediately see that we have
\begin{equation}\label{p2-cmupiso46c}
(\id\otimes \upiota^{r,r'})\circ \fd^{(r+)}= \fd^{(r'+)} \circ \upiota^{r,r'}. 
\end{equation}

Let $\beta$ be a rational number such that $0\leq \beta\leq r+\frac{1}{p-1}$. 
We consider the action of $\hrS(p^{-\beta}f^*(\ccE))$ on $\hcC^{(r+)}$ defined by $\hmu^{r+}_\beta$ \eqref{p2-cmupiso40d}. 
We have seen in \ref{p2-cmupiso40} that the Higgs field $\delta_{\hcC^{(r+)}}\otimes_{\mZ_p}\mQ_p$ \eqref{p2-cmupiso40j} is 
$\hrS(p^{-\beta}f^*(\ccE))$-linear, and hence so is $\vartheta^{(r+)}$ \eqref{p2-cmupiso43b}. Therefore, $\mK^{\bullet}(N\otimes_{\co_{\fX'}}\hcC^{(r+)})$
is naturally a complex of $\hrS(p^{-\beta}f^*(\ccE))$-modules. We deduce an action of $f_*(\hrS(p^{-\beta}f^*(\ccE)))$ on 
$\rR^qf_*(\mK^{\bullet}(N\otimes_{\co_{\fX'}}\hcC^{(r+)}))$ for any integer $q\geq 0$. 
It follows from \eqref{p2-cmupiso40d} that this action is compatible with the actions defined above via the composition of the canonical homomorphisms
\begin{equation}\label{p2-cmupiso46d}
\rS(\ccE)\rightarrow f_*(\rS(f^*(\ccE)))\rightarrow f_*(\hrS(p^{-\beta}f^*(\ccE))).
\end{equation}

We denote by $\umK^\bullet(N)$ the Dolbeault complex of the Higgs $\co_{\fX'}$-module $(N,\utheta)$, where  
$\utheta =(u'\otimes \id)\circ \theta \colon N\rightarrow \ucE'\otimes_{\co_{\fX'}}N$ \eqref{p2-cmupiso21b}. 
Let $\rho_1,\rho_2\colon \cF\rightarrow \co_{\fX'}$ be two splittings of the extension \eqref{p2-cmupiso21a},
$\phi\colon f^*(\cE)\rightarrow \co_{\fX'}$ the $\co_{\fX'}$-linear form induced by $\rho_1-\rho_2$. 
If $\beta>\frac{1}{p-1}$, which implies that $r>0$, it follows from \ref{p2-cmupiso42} that the diagram 
\begin{equation}\label{p2-cmupiso46e}
\xymatrix{
{\mK^\bullet(N\otimes_{\co_{\fX'}}\hcC^{(r+)})}\ar[r]^-(0.5){\upbeta^{r}_1}\ar[d]_{\exp(\phi)}&{\umK^\bullet(N)}\\
{\mK^\bullet(N\otimes_{\co_{\fX'}}\hcC^{(r+)})}\ar[ru]_-(0.5){\upbeta^{r}_2}&}
\end{equation}
where $\upbeta^r_i$ ($i=1,2$) is the morphism of complexes induced by $\rho_i$ defined in \eqref{p2-cmupiso43g},
and $\exp(\phi)\in \Gamma(\fX',\hrS(p^{-\beta}f^*(\ccE)))$ is defined in \eqref{p2-cmupiso40g}, is commutative.

\begin{prop}\label{p2-cmupiso44} 
We take again the assumptions and notation of \ref{p2-cmupiso21}. 
Let $(N,\theta)$ be a Higgs $\co_{\fX'}[\frac 1 p]$-module with coefficients in $\cE'$ \eqref{p1-tshbn9}, 
and let $r,\varepsilon$ be two rational numbers such that $r>0$ and $\varepsilon > r+\frac{1}{p-1}$. Assume that Zariski locally on $\fX'$, 
there exists a coherent $\co_{\fX'}$-lattice $\cN$ of $N$, stable by $\theta$ such that the induced Higgs field is $\varepsilon$-small \eqref{p1-tshbn9}. 
Then, the canonical morphism \eqref{p2-cmupiso43e}
\begin{equation}\label{p2-cmupiso44a}
\upiota^{r,0}\colon \mK^{\bullet}(N\otimes_{\co_{\fX'}}\hcC^{(r+)})\rightarrow \mK^{\bullet}(N\otimes_{\co_{\fX'}}\cC^\dagger)
\end{equation}
is a quasi-isomorphism. 
\end{prop}

Indeed, the question being local on $\fX'$, we may assume that there exist 
a splitting $\rho\colon \cF\rightarrow \co_{\fX'}$ of the extension \eqref{p2-cmupiso21a},
and a coherent $\co_{\fX'}$-lattice $\cN$ of $N$, stable by $\theta$ such that the induced Higgs field is $\varepsilon$-small.
It induces a homomorphism of $\co_{\fX'}$-algebras $\varrho^\dagger\colon \cC^\dagger\rightarrow \co_{\fX'}$ \eqref{p2-cmupiso26a}. 
We set $\utheta=(\id\otimes u')\circ \theta \colon N\rightarrow N\otimes_{\co_{\fX'}}\ucE'$ \eqref{p2-cmupiso21b}, 
and denote by $\umK(N)$ the Dolbeault complex of the Higgs $\co_{\fX'}$-module $(N,\utheta)$ with coefficients in $\ucE'$. 
By \ref{p2-cmupiso29}(i), the morphism of complexes of $\co_{\fX'}$-modules induced by $\varrho^\dagger$ \eqref{p2-cmupiso26d}
\begin{equation}\label{p2-cmupiso44b} 
\upbeta\colon \mK^\bullet(N\otimes_{\co_{\fX'}}\cC^\dagger) \rightarrow \umK^\bullet(N)
\end{equation}
is a quasi-isomorphism. 

We consider the extension 
\begin{equation}\label{p2-cmupiso44c} 
0\rightarrow \co_{\fX'}\rightarrow \cF^{(r)}\rightarrow p^rf^*(\cE) \rightarrow 0
\end{equation}
deduced from $\cF$ \eqref{p2-cmupiso21a} by pullback by the canonical morphism injection $p^rf^*(\cE)\rightarrow f^*(\cE)$, 
which defines the Higgs--Tate algebra $\cC^{(r)}$. Replacing the extension $\cF$ by the extension $\cF^{(r)}$ 
and the extension \eqref{p2-cmupiso21b} by the extension 
\begin{equation}\label{p2-cmupiso44d} 
0\rightarrow p^rf^*(\cE)\rightarrow p^r\cE'\rightarrow p^r\ucE'\rightarrow 0,
\end{equation}
we see that the $\co_{\fX'}$-algebra $\cC^\dagger=\hcC^{(0+)}$ is replaced by $\hcC^{(r+)}$, and the $\co_{\fX'}$-derivation 
$\delta'_{\cC^\dagger}\otimes_{\mZ_p}\mQ_p$ \eqref{p2-cmupiso21e} by $\delta'_{\hcC^{(r+)}}\otimes_{\mZ_p}\mQ_p$ \eqref{p2-cmupiso43a}. 
Moreover, the Higgs $\co_{\fX'}$-field $\colon \cN\rightarrow p^r\cE'\otimes_{\co_{\fX'}}\cN$ induced by $\theta$, is $(\varepsilon-r)$-small \eqref{p1-tshbn9}. 
The splitting $\rho$ of the extension \eqref{p2-cmupiso21a} induces a splitting $\rho^{(r)} \colon \cF^{(r)} \rightarrow \co_{\fX'}$ of the extension 
\eqref{p2-cmupiso44c}, and hence a homomorphism of $\co_{\fX'}$-algebras $\hvarrho^{(r+)}\colon \hcC^{(r+)}\rightarrow \co_{\fX'}$, 
which is in fact none other than the restriction of $\varrho^\dagger$ to $\hcC^{(r+)}$. 
Hence, by \eqref{p2-cmupiso26c}, $\varrho^\dagger$ induces a morphism of complexes of $\co_{\fX'}$-modules 
\begin{equation}\label{p2-cmupiso44e} 
\upbeta^{(r)}\colon \mK^\bullet(N\otimes_{\co_{\fX'}}\hcC^{(r+)}) \rightarrow \umK^\bullet(N),
\end{equation}
which is a quasi-isomorphism by \ref{p2-cmupiso29}(i). This morphism is none other than \eqref{p2-cmupiso43g}.
Since we clearly have $\upbeta^{(r)}=\upbeta \circ \upiota^{r,0}$, 
we deduce that $\upiota^{r,0}$ is a quasi-isomorphism. 

\begin{teo}\label{p2-cmupiso47} 
We take again the assumptions and notation of \ref{p2-cmupiso21} and assume moreover that the morphism $f\colon \fX'\rightarrow \fX$ is proper. 
Let $(N,\theta)$ be a locally CL-small Higgs $\co_{\fX'}[\frac 1 p]$-module with coefficients in $\cE'$ \eqref{p1-tshbn13}, 
$q$ an integer. Then, the Higgs $\co_{\fX}[\frac 1 p]$-module $\rR^qf^\uptau_*(N,\theta)$ \eqref{p2-cmupiso22d} is locally CL-small. 
\end{teo}

Indeed, we may assume that the formal scheme $\fX$ is coherent, and hence so is $\fX'$.  
By \ref{p1-tshbn27}, there exist a rational number $\varepsilon > \frac{1}{p-1}$ and a coherent $\co_{\fX'}$-lattice $\cN$ of $N$, 
stable by $\theta$ such that the induced Higgs field is $\varepsilon$-small. 
Let $r$ be a rational number such that $0<r<\varepsilon-\frac{1}{p-1}$, $\beta=r+\frac{1}{p-1}$. 
We take again the notation of \ref{p2-cmupiso43}. By \ref{p2-cmupiso44}, the morphism \eqref{p2-cmupiso43e}
\begin{equation}\label{p2-cmupiso47a} 
\rR^qf_*(\upiota^{r,0})\colon \rR^qf_*(\mK^{\bullet}(N\otimes_{\co_{\fX'}}\hcC^{(r+)}))\rightarrow 
\rR^qf_*(\mK^{\bullet}(N\otimes_{\co_{\fX'}}\cC^\dagger))
\end{equation}
is an isomorphism. We set $\ccE=\cHom_{\co_{\fX}}(\cE,\co_\fX)$.  
By \eqref{p2-cmupiso46c}, the morphism \eqref{p2-cmupiso47a} is $\rS(\ccE)$-linear for the actions of $\rS(\ccE)$ defined in \ref{p2-cmupiso46}. 
Moreover, we have seen in \ref{p2-cmupiso46} that $f_*(\hrS(p^{-\beta}f^*(\ccE)))$ acts naturally on 
$\rR^qf_*(\mK^{\bullet}(N\otimes_{\co_{\fX'}}\hcC^{(r+)}))$ and this action is compatible with the action of $\rS(\ccE)$ 
via the composition of the canonical homomorphisms
\begin{equation}\label{p2-cmupiso47b} 
\rS(\ccE)\rightarrow f_*(\rS(f^*(\ccE)))\rightarrow f_*(\hrS(p^{-\beta}f^*(\ccE))).
\end{equation}
The latter factors through the canonical homomorphisms 
\begin{equation}\label{p2-cmupiso47c} 
\rS(\ccE)\rightarrow \hrS(p^{-\beta}\ccE)\rightarrow f_*(f^*(\hrS(p^{-\beta}\ccE)))\stackrel{\sim}{\rightarrow}
f_*(\hrS(f^*(p^{-\beta}\ccE)))\stackrel{\sim}{\rightarrow} f_*(\hrS(p^{-\beta}f^*(\ccE))).
\end{equation}
Hence, $\hrS(p^{-\beta}\ccE)$ acts on $\rR^qf_*(\mK^{\bullet}(N\otimes_{\co_{\fX'}}\hcC^{(r+)}))$ and 
this action is compatible with the action of $\rS(\ccE)$. 
Since the $\co_\fX[\frac 1 p]$-module $\rR^qf_*(\mK^{\bullet}(N\otimes_{\co_{\fX'}}\cC^\dagger))$ is coherent by \ref{p2-cmupiso30}(i), 
the proposition follows by \ref{p1-tshbn27}. 

\begin{cor}\label{p2-cmupiso49} 
Let $f\colon \fX'\rightarrow \fX$ be a proper morphism of finite presentation such that $\fX$ and $\fX'$ are $\cS$-flat,
$\cE$ a locally free $\co_\fX$-module of finite type, 
\begin{equation}
0\rightarrow f^*(\cE)\rightarrow \cE'\rightarrow \ucE'\rightarrow 0
\end{equation} 
an exact sequence of locally free $\co_{\fX'}$-modules of finite type. 
Let $(N,\theta)$ a locally CL-small Higgs $\co_{\fX'}[\frac 1 p]$-module with coefficients in $\cE'$ \eqref{p1-tshbn13}, 
$q$ an integer $\geq 0$. We denote by $\utheta \colon N\rightarrow N\otimes_{\co_{\fX'}}\ucE'$ the Higgs field induced by $\theta$, by 
$\umK^{\bullet}(N)$ the associated Dolbeault complex, and by 
\begin{equation}
\kappa^q\colon \rR^qf_*(\umK^{\bullet}(N)) \rightarrow \cE\otimes_{\co_{\fX}}\rR^qf_*(\umK^{\bullet}(N))
\end{equation}
the Katz-Oda Higgs $\co_\fX$-field \eqref{p1-tphdi4h}. 
Then, the Higgs $\co_\fX[\frac 1 p]$-module $(\rR^qf_*(\umK^{\bullet}(N)),\kappa^q)$ is locally CL-small \eqref{p1-tphdi4h}. 
\end{cor}

It follows immediately from \ref{p2-cmupiso36} and \ref{p2-cmupiso47} applied with  the trivial extension $\cF=\co_{\fX'}\oplus f^*(\cE)$ 
of $f^*(\cE)$ by $\co_{\fX'}$.

\begin{prop}\label{p2-cmupiso48} 
We take again the assumptions and notation of \ref{p2-cmupiso21} and assume moreover that the morphism $f\colon \fX'\rightarrow \fX$ is proper
and the canonical homomorphism $\co_\fX\rightarrow f_*(\co_{\fX'})$ is an isomorphism. 
Let $\rho_1,\rho_2\colon \cF\rightarrow \co_{\fX'}$ be two splittings of the extension \eqref{p2-cmupiso21a}. 
We set $\ccE=\cHom_{\co_\fX}(\cE,\co_\fX)$, and denote by $\phi\colon f^*(\cE)\rightarrow \co_{\fX'}$ the $\co_{\fX'}$-linear 
form induced by $\rho_1-\rho_2$, that we consider as a section $\phi\in \Gamma(\fX,\ccE)$ 
since the adjunction morphism $\ccE\rightarrow f_*(f^*(\ccE))$ is an isomorphism. 

Let $(N,\theta)$ be a locally CL-small Higgs $\co_{\fX'}[\frac 1 p]$-module with coefficients in $\cE'$ \eqref{p1-tshbn13}. 
We denote by $\utheta \colon N\rightarrow N\otimes_{\co_{\fX'}}\ucE'$ the Higgs field induced by $\theta$, by 
$\umK^{\bullet}(N)$ the associated Dolbeault complex, and by $\mK^{\bullet}(N\otimes_{\co_{\fX'}}\cC^\dagger)$ 
the Dolbeault complex of $N\otimes_{\co_{\fX'}}\cC^\dagger$ equipped 
with the total Higgs field $\vartheta = \theta\otimes \id+\id\otimes \delta'$ \eqref{p2-cmupiso22a}. 
Let $q$ an integer $\geq 0$, 
\begin{equation}
\kappa^q\colon \rR^qf_*(\umK^{\bullet}(N)) \rightarrow \cE\otimes_{\co_{\fX}}\rR^qf_*(\umK^{\bullet}(N))
\end{equation}
the Katz-Oda Higgs $\co_\fX[\frac 1 p]$-field \eqref{p1-tphdi4h}, which is CL-small by \ref{p2-cmupiso49}. 
The $\co_\fX$-linear endomorphism $\kappa^q_\phi=(\phi\otimes \id)\circ \kappa^q$ of $\rR^qf_*(\umK^{\bullet}(N))$ being therefore 
locally CL-small \eqref{p1-tshbn7}, we can \eqref{p1-tshbn80} define the $\co_{\fX}$-linear isomorphism 
\begin{equation}\label{p2-cmupiso48a} 
\exp(\kappa^q_\phi)\colon \rR^qf_*(\umK^{\bullet}(N))\stackrel{\sim}{\rightarrow} \rR^qf_*(\umK^{\bullet}(N)). 
\end{equation}
Then, the diagram 
\begin{equation}\label{p2-cmupiso48b} 
\xymatrix{
{\rR^qf_*(\mK^\bullet(N\otimes_{\co_{\fX'}}\cC^\dagger))}\ar[r]^-(0.5){\upbeta^q_2}\ar[rd]_-(0.5){\upbeta^q_1}&{\rR^qf_*(\umK^\bullet(N))}\ar[d]^{\exp(\kappa^q_\phi)}\\
&{\rR^qf_*(\umK^\bullet(N)),}}
\end{equation}
where $\upbeta^q_i$ is the (iso)morphism induced by the morphism \eqref{p2-cmupiso26d} defined by $\rho_i$, is commutative. 
\end{prop}

Indeed, we may assume that the formal scheme $\fX$ is coherent.  
By \ref{p1-tshbn27}, \ref{p2-cmupiso47} and \ref{p2-cmupiso49}, there exists a rational number $\beta>\frac{1}{p-1}$ such that 
the canonical actions of $\rS(\ccE)$ on $\rR^qf_*(\mK^\bullet(N\otimes_{\co_{\fX'}}\cC^\dagger))$ and $\rR^qf_*(\umK^\bullet(N)$
extend to actions of the $\co_\fX$-algebra $\hrS(p^{-\beta}\ccE)$. 
Recall that the homomorphism $\rS(\ccE)\rightarrow \rS(p^{-\beta}\ccE)$ factors through 
$\rS(\ccE)\rightarrow \Gamma(\ccE)$ \eqref{p1-NC7} and induces a homomorphism of $\co_{\fX}$-algebras
\begin{equation}\label{p2-cmupiso48c} 
\hGamma(\ccE)\rightarrow \hrS(p^{-\beta}\ccE). 
\end{equation}
We can therefore consider $\exp(\phi)\in \Gamma(\fX, \hrS(p^{-\beta}\ccE))$. 
The morphism $\upbeta^q_2$ is $\rS(\ccE)$-linear by \ref{p1-tphdi10}. 
It is therefore $\hrS(p^{-\beta}\ccE)$-linear by \ref{p1-tshbn27}. 
Hence, the diagram 
\begin{equation}\label{p2-cmupiso48d} 
\xymatrix{
{\rR^qf_*(\mK^\bullet(N\otimes_{\co_{\fX'}}\cC^\dagger))}\ar[r]^-(0.5){\upbeta^q_2}\ar[d]_{\exp(\phi)}&
{\rR^qf_*(\umK^\bullet(N))}\ar[d]^{\exp(\kappa^q_\phi)}\\
{\rR^qf_*(\mK^\bullet(N\otimes_{\co_{\fX'}}\cC^\dagger))}\ar[r]^-(0.5){\upbeta^q_2}&{\rR^qf_*(\umK^\bullet(N))}}
\end{equation}
is commutative. We take again the notation of \ref{p2-cmupiso43} and \ref{p2-cmupiso46} for $r=\beta-\frac{1}{p-1}$. 
We can then naturally consider $\mK^\bullet(N\otimes_{\co_{\fX'}}\hcC^{(r+)})$ as a complex of $\hrS(p^{-\beta}f^*(\ccE))$-modules.  
After eventually shrinking $\beta$, the diagram 
\begin{equation}\label{p2-cmupiso48e} 
\xymatrix{
{\rR^qf_*(\mK^\bullet(N\otimes_{\co_{\fX'}}\hcC^{(r+)}))}\ar[rr]^-(0.5){\rR^qf_*(\upiota^{r,0})}\ar[d]_{\rR^qf_*(\exp(\phi))}&&
{\rR^qf_*(\mK^\bullet(N\otimes_{\co_{\fX'}}\cC^\dagger))}\ar[d]^{\exp(\phi)}\\
{\rR^qf_*(\mK^\bullet(N\otimes_{\co_{\fX'}}\hcC^{(r+)}))}\ar[rr]^-(0.5){\rR^qf_*(\upiota^{r,0})}&&
{\rR^qf_*(\mK^\bullet(N\otimes_{\co_{\fX'}}\cC^\dagger)),}}
\end{equation}
where $\upiota^{r,0}$ is defined in \eqref{p2-cmupiso43e} and 
in the left vertical arrow, $\exp(\phi)\in \Gamma(\fX',\hrS(p^{-\beta}f^*(\ccE)))$, is commutative. 
Indeed, $\upiota^{r,0}$ is $\rS(f^*(\ccE))$-linear by \eqref{p2-cmupiso46c}. 
After eventually shrinking $\beta$, $\upiota^{r,0}$ is a quasi-isomorphism by \ref{p2-cmupiso44}.
Since the $\co_\fX[\frac 1 p]$-module $\rR^qf_*(\mK^\bullet(N\otimes_{\co_{\fX'}}\cC^\dagger))$ is coherent by \ref{p2-cmupiso30}(i), 
we deduce that after eventually shrinking $\beta$, $\rR^qf_*(\upiota^{r,0})$ is $\hrS(p^{-\beta}\ccE)$-linear by \ref{p1-tshbn27}. 
The commutativity of \eqref{p2-cmupiso48e} follows by \ref{p2-cmupiso45}. 
The proposition follows then from \eqref{p2-cmupiso48d}, \eqref{p2-cmupiso48e}, \eqref{p2-cmupiso46e} and \eqref{p2-cmupiso43h}.


\chapter{\texorpdfstring{Functorial aspects of the $p$-adic Simpson correspondence}{Functorial aspects of the p-adic Simpson correspondence}}
\label{functoriality}

\section{Assumptions and notation}\label{p2-NC}

The conventions and notation of §\ref{p1-NC} are in force in this chapter, to which we add the following.

\subsection{}\label{p2-ncgt1}
In this chapter, $K$ denotes a complete discrete valuation field of characteristic $0$, 
with algebraically closed residue field $k$ of characteristic $p>0$,  
$\co_K$ the valuation ring of $K$, $\oK$ an algebraic closure of $K$, $\co_\oK$ the integral closure of $\co_K$ in $\oK$,
$\fm_\oK$ the maximal ideal of $\co_\oK$ and $G_K$ the Galois group of $\oK$ over $K$.  
We denote by $\co_C$ the $p$-adic completion of $\co_\oK$, $\fm_C$ its maximal ideal,
$C$ its field of fractions and $v$ its valuation, normalized by $v(p)=1$. 
For any abelian group $A$, we denote by $\hA$ its $p$-adic completion. 

We set $S=\Spec(\co_K)$, $\oS=\Spec(\co_\oK)$, $\coS=\Spec(\co_C)$ and $\cS=\Spf(\co_C)$.  
We denote by $s$ (resp.\ $\eta$, resp.\ $\oeta$) the closed point of $S$ (resp.\  generic point of $S$, resp.\ generic point of $\oS$).
We equip $S$ with the logarithmic structure $\cM_S$ defined by its closed point, 
and $\oS$ and $\coS$ with the logarithmic structures $\cM_\oS$ and $\cM_\coS$ pullbacks of $\cM_S$. 

For any integer $n\geq 1$, we set $S_n=\Spec(\co_K/p^n\co_K)$. For any $S$-scheme $X$, we set 
\begin{equation}\label{p2-ncgt1a}
\oX=X\times_S\oS,  \ \ \ \coX=X\times_S\coS \ \ \ {\rm and}\ \ \  X_n=X\times_SS_n.
\end{equation} 

\subsection{}\label{p2-ncgt2}
Since $\co_\oK$ is a non-discrete valuation ring of height $1$,
it is possible to develop the $\alpha$-algebra (or almost-algebra) on this ring (\cite{agt} §~V, \cite{ag1} 2.6-2.10).
We choose a compatible system $(\beta_n)_{n>0}$ of $n$th roots of $p$ in $\co_\oK$. For any rational number $\varepsilon>0$,
we set $p^\varepsilon=(\beta_n)^{\varepsilon n}$ where $n$ is an integer $>0$ such that $\varepsilon n$ is an integer.

\subsection{}\label{p2-ncgt3}
With the notation of (\cite{ag2} 3.1), we denote by $(\tS,\cM_\tS)$ one of the following two schemes  
\begin{equation}\label{p2-ncgt3a}
(\cA_2(\oS),\cM_{\cA_2(\oS)})\ \ \ {\rm or} \ \ \ (\cA^{\ast}_2(\oS/S),\cM_{\cA^{\ast}_2(\oS/S)});
\end{equation}
the first case will be called {\em absolute} and the second one will be called {\em relative} (see \cite{ag2} 3.1.14). 
Observe that the relative setting depends on the choice of a uniformizer $\pi$ of $\co_K$ and a sequence $(\pi_n)_{n\geq 0}$
of elements of $\co_\oK$ such that $\pi_0=\pi$ and $\pi_{n+1}^p=\pi_n$ for all $n\geq 0$. 
We have a canonical strict closed immersion
\begin{equation}\label{p2-ncgt3b}
\iota\colon (\coS,\cM_\coS)\rightarrow (\tS,\cM_{\tS}),
\end{equation} 
whose underlying closed immersion of schemes $\coS\rightarrow \tS$ is defined by the square-zero ideal $\txi\co_\tS$ generated by the canonical element $\txi=\xi$ or $\txi=\xi^*_\pi$ 
depending on whether we are in the absolute or relative case (see \cite{ag2} 3.1.4 and 3.1.6). 
The ideal $\txi\co_\tS$ is a free $\co_\coS$-module with basis $\txi$. We denote it by $\txi\co_\coS$ and denote by $\txi^{-1}\co_\coS$ its $\co_\coS$-dual. 

For any ringed topos $(X,\co_X)$ over $(\coS,\co_\coS)$, 
any $\co_X$-module $M$ and any integer $i\geq 1$, we denote the $\co_X$-modules $M\otimes_{\co_\coS}(\txi \co_\coS)^{\otimes i}$ 
and $M\otimes_{\co_\coS}(\txi^{-1} \co_\coS)^{\otimes i}$ simply by $\txi^i M$ and $\txi^{-i} M$, respectively.

\section{The Higgs--Tate algebras}

\subsection{}\label{p2-hta1}
We consider in this section a commutative diagram of the category $\FLS$ of fine logarithmic schemes \eqref{p1-NC1}
\begin{equation}\label{p2-hta1a}
\xymatrix{
{(Y,\cM_{Y})}\ar[r]^-(0.5){j}\ar[d]^{h}\ar@/_2pc/[dd]_{g}&{(\tY,\cM_{\tY})}\ar@/^2pc/[dd]^{\tg}&\\
{(\coX,\cM_{\coX})}\ar[r]^-(0.5){i}\ar[d]^{\cof}\ar@{}[rd]|{\Box}&{(\tX,\cM_{\tX})}\ar[d]_{\tf}\\
{(\coS,\cM_{\coS})}\ar[r]^-(0.5){\iota}&{(\tS,\cM_{\tS}),}}
\end{equation}
where $\iota$ is the strict closed immersion defined in \eqref{p2-ncgt3b}, such that the following conditions are satisfied:
\begin{itemize}
\item[(i)] The lower square $(\cof,\tf,\iota,i)$ of \eqref{p2-hta1a} is Cartesian in $\FLS$, or equivalently in the category of logarithmic schemes by \ref{p1-NC3}. 
\item[(ii)] The morphism $\tf$ is smooth. 
\item[(iii)] The external rectangle $(g,\tg,\iota,j)$ of \eqref{p2-hta1a} is Cartesian and $g=\cof\circ h$. 
Therefore, $j$ is a thickening of order one and the underlying closed immersion of schemes $Y\rightarrow \tY$ is defined by the ideal 
$\txi\co_{\tY}$ by \ref{p1-NC3}. 
\item[(iv)] The $\co_{Y}$-module $\txi\co_{\tY}$ is invertible. 
Hence, with the convention of \ref{p2-ncgt3}, the canonical $\co_Y$-linear morphism $\txi\co_{Y}\rightarrow \txi\co_{\tY}$ is an isomorphism.   
We will identify $\txi\co_{\tY}$ with $\txi\co_{Y}$ by this isomorphism. 
\end{itemize}

\begin{lem}\label{p2-hta2}
The morphism of logarithmic schemes $\tf$ is integral. In particular, the morphism of schemes $\tX\rightarrow \tS$ is flat. 
\end{lem}

Since $\tS$ has a chart of the form $\mN\rightarrow \Gamma(\tS,\cM_{\tS})$ (\cite{ag2} 3.1.12), $\tf$ is integral by (\cite{kato1} 4.4). 
The second assertion follows from (\cite{kato1} 4.5). 

\subsection{}\label{p2-hta4}
To lighten the notation, we set, with the conventions of \ref{p2-ncgt3}, 
\begin{eqnarray}
\tOmega^1_{\coX/\coS}&=&\Omega^1_{(\coX,\cM_\coX)/(\coS,\cM_\coS)}, \label{p2-hta4b}\\
\Omega&=&\txi^{-1}\tOmega^1_{\coX/\coS}, \label{p2-hta4c}\\
\rT&=&\cHom_{\co_{\coX}}(\Omega,\co_{\coX}).\label{p2-hta4cc}
\end{eqnarray} 
Observe that $\rT$ is canonically isomorphic to $\cHom_{\co_{\coX}}(\tOmega^1_{\coX/\coS},\txi\co_{\coX})$ and that 
$\Omega$ and $\rT$ are $\coS$-flat. 
For any rational number $r\geq 0$, we set 
\begin{eqnarray}
\Omega^{(r)}&=&p^r\Omega,\label{p2-hta4d}\\ 
\rT^{(r)}&=&\cHom_{\co_{\coX}}(\Omega^{(r)},\co_{\coX}).\label{p2-hta4dd}
\end{eqnarray}
The canonical injection $\pi^{(r)}\colon \Omega^{(r)}\rightarrow \Omega$ induces an injection $\pi^{(r)\vee}\colon \rT\rightarrow \rT^{(r)}$ 
that identifies $\rT^{(r)}$ with the submodule $p^{-r}\rT$ of $\rT\otimes_{\mZ_p}\mQ_p$. 

For any rational numbers $r\geq r'\geq 0$, we denote by 
\begin{equation}\label{p2-hta4e}
\pi^{(r,r')}\colon \Omega^{(r)}\rightarrow \Omega^{(r')}
\end{equation} 
the canonical injection, so that $\pi^{(r)}=\pi^{(r,0)}$. We denote by $\pi^{(r,r')\vee}\colon \rT^{(r')}\rightarrow \rT^{(r)}$ its dual. 

\subsection{}\label{p2-hta5}
The conditions of \ref{p1-rdt1} being satisfied by diagram \eqref{p2-hta1a}, we denote by 
$\cL_{\tY/\tX}$ the torsor of liftings of $h$ to $\tY$ over $\tX$ \ref{p1-rdt2}, which is a torsor of $Y_\et$ 
under the $\co_Y$-module $h^*(\rT)$ \eqref{p2-hta4cc}, by $\cF_{\tY/\tX}$ the $\co_Y$-module of affine functions on $\cL_{\tY/\tX}$ 
\ref{p1-rdt3}, by
\begin{equation}\label{p2-hta5b}
0\rightarrow \co_{Y}\rightarrow \cF_{\tY/\tX}\rightarrow h^*(\Omega) \rightarrow 0
\end{equation}
the associated exact sequence \eqref{p1-rdt3a} and by
\begin{equation}\label{p2-hta5c}
\cC_{\tY/\tX}=\underset{\underset{n\geq 0}{\longrightarrow}}\lim\ \rS^n_{\co_{Y}}(\cF_{\tY/\tX})
\end{equation}
the associated $\co_{Y}$-algebra \eqref{p1-rdt3c}. 

For any rational number $r\geq 0$, we denote by $\cL^{(r)}_{\tY/\tX}$ the $h^*(\pi^{(r)\vee})$-twist of $\cL_{\tY/\tX}$ \eqref{p2-hta4dd} (see \ref{p1-NC4}), 
which is a $h^*(\rT^{(r)})$-torsor of $Y_\et$, by $\cF^{(r)}_{\tY/\tX}$  the $\co_Y$-module of affine functions on $\cL^{(r)}_{\tY/\tX}$ \eqref{p1-prem1}, by
\begin{equation}\label{p2-hta5d}
0\rightarrow \co_{Y}\rightarrow \cF^{(r)}_{\tY/\tX}\rightarrow h^*(\Omega^{(r)}) \rightarrow 0
\end{equation} 
the associated exact sequence \eqref{p1-prem1a} and by
\begin{equation}\label{p2-hta5e}
\cC^{(r)}_{\tY/\tX}=\underset{\underset{n\geq 0}{\longrightarrow}}\lim\ \rS^{n}_{\co_{Y}}(\cF^{(r)}_{\tY/\tX})
\end{equation}
the associated $\co_{Y}$-algebra \eqref{p1-prem1c}. 

The canonical morphism $\varpi^{(r)}_{\tY/\tX}\colon \cL_{\tY/\tX}\rightarrow \cL^{(r)}_{\tY/\tX}$ being $h^*(\pi^{(r)\vee})$-equivariant, it induces by pullback (\cite{agt} II.4.12) 
a canonical $\co_Y$-linear morphism
\begin{equation}
\varpi^{(r)*}_{\tY/\tX}\colon \cF^{(r)}_{\tY/\tX} \rightarrow \cF_{\tY/\tX}
\end{equation}
that fits into a commutative diagram
\begin{equation}\label{p2-hta5f}
\xymatrix{
0\ar[r]&\co_Y\ar@{=}[d]\ar[r]&{\cF^{(r)}_{\tY/\tX}}\ar[r]\ar[d]^{\varpi^{(r)*}_{\tY/\tX}}&{h^*(\Omega^{(r)})}\ar[r]\ar[d]^{h^*(\pi^{(r)})}&0\\
0\ar[r]&\co_Y\ar[r]&{\cF_{\tY/\tX}}\ar[r]&{h^*(\Omega)}\ar[r]&0.}
\end{equation}
The extension $\cF^{(r)}_{\tY/\tX}$ is therefore canonically isomorphic to the pullback of $\cF_{\tY/\tX}$ by $h^*(\pi^{(r)})$ \eqref{p2-hta4d}.

Recall \eqref{p1-thbn1d} that there exists a unique $\co_Y$-derivation 
\begin{equation}\label{p2-hta5g}
d_{\cC^{(r)}_{\tY/\tX}}\colon \cC^{(r)}_{\tY/\tX} \rightarrow h^*(\Omega^{(r)})\otimes_{\co_Y}\cC^{(r)}_{\tY/\tX}
\end{equation}
extending the canonical morphism $\cF^{(r)}_{\tY/\tX}\rightarrow h^*(\Omega^{(r)})$. 
It identifies canonically with the universal $\co_Y$-derivation of $\cC^{(r)}_{\tY/\tX}$.
It is also a Higgs $\co_Y$-field on $\cC^{(r)}_{\tY/\tX}$ with coefficients in $h^*(\Omega^{(r)})$ \eqref{p1-delta-con1}. We denote by
\begin{equation}\label{p2-hta5h}
\delta_{\cC^{(r)}_{\tY/\tX}}\colon \cC^{(r)}_{\tY/\tX} \rightarrow h^*(\Omega)\otimes_{\co_Y}\cC^{(r)}_{\tY/\tX}
\end{equation}
the $\co_Y$-derivation $(h^*(\pi^{(r)})\otimes \id)\circ d_{\cC^{(r)}_{\tY/\tX}}$. 

For any rational numbers $r\geq r'\geq 0$, the morphism $\varpi^{(r)}_{\tY/\tX}$ factors uniquely into
\begin{equation}
\xymatrix{
{\cL_{\tY/\tX}}\ar[r]^{\varpi^{(r')}_{\tY/\tX}}&{\cL^{(r')}_{\tY/\tX}}\ar[r]^{\varpi^{(r,r')}_{\tY/\tX}}&{\cL^{(r')}_{\tY/\tX},}}
\end{equation}
where $\varpi^{(r,r')}_{\tY/\tX}$ is $\pi^{(r,r')\vee}$-equivariant \eqref{p2-hta4e}. The latter induces by pullback a canonical $\co_{Y}$-linear morphism
\begin{equation}\label{p2-hta5i}
\varpi^{(r,r')*}_{\tY/\tX}\colon\cF_{\tY/\tX}^{(r)}\rightarrow \cF_{\tY/\tX}^{(r')}
\end{equation}
that fits into a commutative diagram
\begin{equation}\label{p2-hta5m}
\xymatrix{
0\ar[r]&\co_Y\ar@{=}[d]\ar[r]&{\cF^{(r)}_{\tY/\tX}}\ar[r]\ar[d]^{\varpi^{(r,r')*}_{\tY/\tX}}&{h^*(\Omega^{(r)})}\ar[r]\ar[d]^{h^*(\pi^{(r,r')})}&0\\
0\ar[r]&\co_Y\ar[r]&{\cF^{(r')}_{\tY/\tX}}\ar[r]&{h^*(\Omega^{(r')})}\ar[r]&0.}
\end{equation}
We deduce a canonical homomorphism of $\co_{Y}$-algebras
\begin{equation}\label{p2-hta5j}
\alpha_{\tY/\tX}^{r,r'}\colon \cC_{\tY/\tX}^{(r)}\rightarrow \cC_{\tY/\tX}^{(r')}.
\end{equation}

We have
\begin{equation}\label{p2-hta5k}
(h^*(\pi^{(r,r')})\otimes \alpha^{r,r'}_{\tY/\tX}) \circ d_{\cC^{(r)}_{\tY/\tX}}=d_{\cC^{(r')}_{\tY/\tX}} \circ \alpha^{r,r'}_{\tY/\tX},
\end{equation}
and hence 
\begin{equation}\label{p2-hta5l}
(\id\otimes \alpha^{r,r'}_{\tY/\tX}) \circ \delta_{\cC^{(r)}_{\tY/\tX}}=\delta_{\cC^{(r')}_{\tY/\tX}} \circ \alpha^{r,r'}_{\tY/\tX}.
\end{equation}

\begin{defi}\label{p2-hta7}
The torsor $\cL_{\tY/\tX}$ (resp.\ $\cL^{(r)}_{\tY/\tX}$) is called the (resp.\ $(r)$-twisted) 
{\em torsor of local liftings of the morphism $h$ to $\tY$ over $\tX$} \eqref{p2-hta5}.
The $\co_Y$-algebra $\cC_{\tY/\tX}$ (resp.\ $\cC^{(r)}_{\tY/\tX}$) is called the {\em Higgs--Tate algebra of $\tY$ over $\tX$} (resp.\ {\em of thickness $r$}).
The $\co_Y$-module $\cF_{\tY/\tX}$ (resp.\ $\cF^{(r)}_{\tY/\tX}$) is called the {\em Higgs--Tate extension of $\tY$ over $\tX$} (resp.\ {\em of thickness $r$}).
\end{defi}

The constructions in \ref{p2-hta5} are motivated and inspired by the construction of 
the period rings of the $p$-adic Simpson correspondence in (\cite{agt} II.10.5 and II.12.1), 
see \ref{p2-rlps5}, which explains the terminology in \ref{p2-hta7}. We adopt here a slightly different, although equivalent, normalization, see \ref{p2-hta71}. 

\begin{rema}\label{p2-hta71}
Since the $\co_X$-module $\Omega$ is $\co_C$-flat, for every rational number $r\geq 0$, we can identify $\pi^{(r)}\colon \Omega^{(r)}\rightarrow \Omega$ \eqref{p2-hta4d} 
with the morphism of multiplication by $p^r$ of $\Omega$. Hence by \eqref{p2-hta5f}, the extension $\cF^{(r)}_{\tY/\tX}$ \eqref{p2-hta5d} is canonically isomorphic to the extension obtained 
by pullback of the extension $\cF_{\tY/\tX}$ \eqref{p2-hta5b} by 
the morphism of multiplication by $p^r$ on $h^*(\Omega)$, which corresponds to the construction introduced in \ref{p1-thbn3}. 
\end{rema}

\subsection{}\label{p2-hta6}
Any section $\tgamma \in \cL_{\tY/\tX}(Y)$ determines associated splittings $\rho_\tgamma\colon \cF_{\tY/\tX}\rightarrow \co_Y$ and 
$\sigma_\tgamma\colon h^*(\Omega)\rightarrow \cF_{\tY/\tX}$ of the extension \eqref{p2-hta5b} (see \ref{p1-thbn16}). 
For every rational number $r\geq 0$, $\tgamma^{(r)}=\varpi^{(r)}(\tgamma)\in \cL^{(r)}_{\tY/\tX}(Y)$ 
determines associated splittings $\rho_{\tgamma^{(r)}}\colon \cF^{(r)}_{\tY/\tX}\rightarrow \co_Y$ and 
$\sigma_{\tgamma^{(r)}}\colon h^*(\Omega^{(r)})\rightarrow \cF^{(r)}_{\tY/\tX}$ of the extension \eqref{p2-hta5d}. 
We immediately see that $\rho_\tgamma$ and $\rho_{\tgamma^{(r)}}$ (resp.\ $\sigma_\tgamma$ and $\sigma_{\tgamma^{(r)}}$) are compatible with \eqref{p2-hta5f}. 
The morphism $\rho_{\tgamma^{(r)}}$ extends to a homomorphism of $\co_Y$-algebras $\varrho_{\tgamma^{(r)}}\colon \cC^{(r)}_{\tY/\tX}\rightarrow \co_Y$, 
and the morphism $\sigma_{\tgamma^{(r)}}$ induces an isomorphism of $\co_Y$-algebras 
\begin{equation}\label{p2-hta6a}
\varsigma_{\tgamma^{(r)}}\colon \rS_{\co_Y}(h^*(\Omega^{(r)}))\stackrel{\sim}{\rightarrow} \cC^{(r)}_{\tY/\tX}.
\end{equation}

\subsection{}\label{p2-hta120}
We also consider in this section a commutative diagram of $\FLS$ 
\begin{equation}\label{p2-hta12a}
\xymatrix{
{(Y,\cM_{Y})}\ar[r]^-(0.5){j}\ar[d]^{h'}\ar@/_2.2pc/[dd]_{h}\ar@/_4pc/[ddd]_{g}&{(\tY,\cM_{\tY})}\ar@/^4pc/[ddd]^{\tg}&\\
{(\coX',\cM_{\coX'})}\ar[r]^-(0.5){i'}\ar[d]^{\cogamma}\ar@/^2pc/[dd]|\hole^-(0.3){\cof'}&{(\tX',\cM_{\tX'})}\ar@/^2.2pc/[dd]^{\tf'}&\\
{(\coX,\cM_{\coX})}\ar[r]^-(0.5){i}\ar[d]^{\cof}&{(\tX,\cM_{\tX})}\ar[d]_{\tf}\\
{(\coS,\cM_{\coS})}\ar[r]^-(0.5){\iota}&{(\tS,\cM_{\tS})}}
\end{equation}
extending diagram \eqref{p2-hta1a} and satisfying the following conditions:
\begin{itemize}
\item[(v)] We have $h=\cogamma\circ h'$,  $\cof'=\cof\circ \cogamma$ and $g=\cof\circ h=\cof'\circ h'$. 
\item[(vi)] The rectangle $(\cof',\tf',\iota,i')$ of \eqref{p2-hta12a} is Cartesian. 
Therefore, $i'$ is a thickening of order one and the underlying closed immersion of schemes 
$\coX'\rightarrow \tX'$ is defined by the ideal $\txi\co_{\tX'}$ \eqref{p1-NC3}. 
\item[(vii)] The morphism $\tf'$ is smooth. So the morphism of schemes $\tX'\rightarrow \tS$ is flat by \ref{p2-hta2}. 
Therefore, with the convention of \ref{p2-ncgt3}, the canonical $\co_{\coX'}$-linear morphism 
$\txi\co_{\coX'}\rightarrow \txi\co_{\tX'}$ is an isomorphism. We will identify $\txi\co_{\tX'}$ with $\txi\co_{\coX'}$ by this isomorphism.
\end{itemize}
We will add an extra-condition in \ref{p2-hta12}. 

To lighten the notation, we set, with the convention of \ref{p2-ncgt3},   
\begin{eqnarray}
\tOmega^1_{\coX'/\coS}&=&\Omega^1_{(\coX',\cM_{\coX'})/(\coS,\cM_\coS)},\label{p2-hta120a}\\ 
\Omega'&=&\txi^{-1}\tOmega^1_{\coX'/\coS},\label{p2-hta120b}\\
\rT'&=&\cHom_{\co_{\coX'}}(\Omega',\co_{\coX'}).\label{p2-hta120c}
\end{eqnarray}
The canonical morphism $\cogamma^*(\tOmega^1_{\coX/\coS})\rightarrow \tOmega^1_{\coX'/\coS}$ induces dual morphisms 
of locally free $\co_{\coX'}$-modules of finite type
\begin{eqnarray}
u\colon \cogamma^*(\Omega)\rightarrow \Omega',\label{p2-hta120d}\\ 
v\colon \rT' \rightarrow \cogamma^*(\rT).\label{p2-hta120e}
\end{eqnarray}
For any rational number $r\geq 0$, extending the notation of \ref{p2-hta4}, we set 
\begin{eqnarray}
\Omega'^{(r)}&=&p^r\Omega',\label{p2-hta13a}\\
\rT'^{(r)}&=&\cHom_{\co_{\coX'}}(\Omega'^{(r)},\co_{\coX'}).\label{p2-hta13aa}
\end{eqnarray} 
The canonical injection $\pi'^{(r)}\colon \Omega'^{(r)}\rightarrow \Omega'$ induces an injection $\pi'^{(r)\vee}\colon \rT'\rightarrow \rT'^{(r)}$ 
that identifies $\rT'^{(r)}$ with the submodule $p^{-r}\rT'$ of $\rT'\otimes_{\mZ_p}\mQ_p$. 
The morphisms $u$ and $v$ induce dual morphisms  
\begin{eqnarray}
u^{(r)}\colon \cogamma^*(\Omega^{(r)})\rightarrow \Omega'^{(r)},\label{p2-hta120dd}\\ 
v^{(r)}\colon \rT'^{(r)} \rightarrow \cogamma^*(\rT^{(r)}).\label{p2-hta120ee}
\end{eqnarray}

\subsection{}\label{p2-hta121}
The conditions \ref{p2-hta1}(i)-(ii) being satisfied by the rectangle $(\cof',\tf',\iota,i')$ of \eqref{p2-hta12a},
we denote by $\cL_{\tY/\tX'}$ the torsor of liftings of $h'$ to $\tY$ over $\tX'$ \eqref{p2-hta7}, which is a torsor of $Y_\et$ 
under the $\co_Y$-module $h'^*(\rT')$, and by $\cF_{\tY/\tX'}$ (resp.\ $\cC_{\tY/\tX'}$) the Higgs--Tate extension (resp.\ algebra) 
of $\tY$ over $\tX'$; so we have a canonical exact sequence of $\co_Y$-modules
\begin{equation}\label{p2-hta12e}
0\rightarrow \co_{Y}\rightarrow \cF_{\tY/\tX'}\rightarrow h'^*(\Omega') \rightarrow 0. 
\end{equation}

For any rational number $r\geq 0$, we denote by $\cL^{(r)}_{\tY/\tX'}$ the $(r)$-twisted torsor of liftings of the morphism $h'$ to $\tY$ over $\tX'$, 
which is a torsor of $Y_\et$ under the $\co_Y$-module $h'^*(\rT'^{(r)})$,
and by $\cF^{(r)}_{\tY/\tX'}$ (resp.\ $\cC^{(r)}_{\tY/\tX'}$) the Higgs--Tate extension (resp.\ algebra) of $\tY$ over $\tX'$ of thickness $r$ \eqref{p2-hta7};
so we have a canonical exact sequence of $\co_Y$-modules
\begin{equation}\label{p2-hta13e}
0\rightarrow \co_{Y}\rightarrow \cF^{(r)}_{\tY/\tX'}\rightarrow h'^*(\Omega'^{(r)}) \rightarrow 0. 
\end{equation}

\subsection{}\label{p2-hta122}
The conditions \ref{p2-hta1}(iii)-(iv) being satisfied by the rectangle $(\cof',\tf',\iota,i')$ of \eqref{p2-hta12a}, 
we denote by $\cL_{\tX'/\tX}$ the torsor of liftings of $\cogamma$ to $\tX'$ over $\tX$, which is a torsor of $\coX'_\et$ 
under the $\co_{\coX'}$-module $\cogamma^*(\rT)$, 
and by $\cF_{\tX'/\tX}$ (resp.\ $\cC_{\tX'/\tX}$) the Higgs--Tate extension (resp.\ algebra) of $\tX'$ over $\tX$ \eqref{p2-hta7}; 
so we have a canonical exact sequence of $\co_{\coX'}$-modules
\begin{equation}\label{p2-hta8b}
0\rightarrow \co_{\coX'}\rightarrow \cF_{\tX'/\tX}\rightarrow \cogamma^*(\Omega) \rightarrow 0.
\end{equation} 

For any rational number $r\geq 0$, we denote by $\cL^{(r)}_{\tX'/\tX}$ the $(r)$-twisted torsor of liftings of the morphism $\cogamma$ to $\tX'$ over $\tX$,
which is a torsor of $\coX'_\et$ under the $\co_{\coX'}$-module $\cogamma^*(\rT^{(r)})$,
and by $\cF^{(r)}_{\tX'/\tX}$ (resp.\ $\cC^{(r)}_{\tX'/\tX}$) the Higgs--Tate extension (resp.\ algebra) of $\tX'$ over $\tX$ of thickness $r$  \eqref{p2-hta7};
so we have a canonical exact sequence of $\co_{\coX'}$-modules
\begin{equation}\label{p2-hta13f}
0\rightarrow \co_{\coX'}\rightarrow \cF^{(r)}_{\tX'/\tX}\rightarrow \cogamma^*(\Omega^{(r)}) \rightarrow 0. 
\end{equation}

\subsection{}\label{p2-hta10}
We denote by $\cHom_{v}(\cL_{\tY/\tX'},\cL_{\tY/\tX})$ the sheaf of $h'^*(v)$-equivariant morphisms from $\cL_{\tY/\tX'}$ to $\cL_{\tY/\tX}$, 
which is canonically a torsor of $Y_\et$ under the $\co_Y$-module $h^*(\rT)$.
By \eqref{p1-rdt6l}, we have a canonical $h^*(\rT)$-equivariant isomorphism
\begin{equation}\label{p2-hta10a}
\upnu\colon h'^+(\cL_{\tX'/\tX})\stackrel{\sim}{\rightarrow} \cHom_{v}(\cL_{\tY/\tX'},\cL_{\tY/\tX}),
\end{equation}
where $h'^+(\cL_{\tX'/\tX})$ is the affine pullback of $\cL_{\tX'/\tX}$ by $h'$ \eqref{p1-NC5}, 
which is a torsor of $Y_\et$ under the $\co_Y$-module $h^*(\rT)=h'^*(\cogamma^*(\rT))$. 
We deduce by ``evaluation'' a canonical morphism 
\begin{equation}\label{p2-hta10b}
\upphi_\cL\colon \cL_{\tY/\tX'}\times h'^+(\cL_{\tX'/\tX}) \rightarrow \cL_{\tY/\tX},
\end{equation}
which is equivariant with respect to the $\co_Y$-linear morphism $h'^*(v)+\id_{h^*(\rT)}\colon h'^*(\rT')\oplus h^*(\rT)\rightarrow h^*(\rT)$ 
\eqref{p2-hta120e}. By \ref{p1-prem26}, $\upphi_\cL$ induces a canonical $\co_{Y}$-linear morphism \eqref{p1-prem26d}
\begin{equation}\label{p2-hta10c}
\upvarphi\colon \cF_{\tY/\tX}\rightarrow \cF_{\tY/\tX'}\otimes_{\co_{Y}} h'^*(\cF_{\tX'/\tX}),
\end{equation}
and a canonical morphism of $\co_{Y}$-algebras \eqref{p1-prem26e}
\begin{equation}\label{p2-hta10d}
\upphi\colon \cC_{\tY/\tX}\rightarrow \cC_{\tY/\tX'} \otimes_{\co_{Y}} h'^*(\cC_{\tX'/\tX}),
\end{equation}
compatible with $\upvarphi$.

It is also convenient to consider the morphism 
\begin{equation}\label{p2-hta10e}
\uppsi_\cL=\upphi_\cL\times \uppi_2\colon \cL_{\tY/\tX'}\times h'^+(\cL_{\tX'/\tX}) \rightarrow \cL_{\tY/\tX} \times h'^+(\cL_{\tX'/\tX}), 
\end{equation}
where $\uppi_2$ is the projection on the second factor of $\cL_{\tY/\tX'}\times h'^+(\cL_{\tX'/\tX})$, which is equivariant with respect to the $\co_Y$-linear morphism 
$w\colon h'^*(\rT')\oplus h^*(\rT)\rightarrow h^*(\rT)\oplus h^*(\rT)$ defined by the matrix 
\begin{equation}\label{p2-hta10f}
\begin{pmatrix}
h'^*(v)&\id_{h^*(\rT)}\\0&\id_{h^*(\rT)}
\end{pmatrix}.
\end{equation}
By \ref{p1-prem25},  $\uppsi_\cL$ induces a canonical morphism of $\co_{Y}$-algebras
\begin{equation}\label{p2-hta10g}
\uppsi\colon \cC_{\tY/\tX}\otimes_{\co_{Y}} h'^*(\cC_{\tX'/\tX})\rightarrow \cC_{\tY/\tX'} \otimes_{\co_{Y}} h'^*(\cC_{\tX'/\tX}),
\end{equation}
which is none other than the $h'^*(\cC_{\tX'/\tX})$-linearization of $\upphi$.

\subsection{}\label{p2-hta150} 
For any rational numbers $r\geq r'\geq 0$, we denote by 
\begin{equation}\label{p2-hta150a} 
w^{(r,r')}\colon h'^*(\rT'^{(r)})\oplus h^*(\rT^{(r')})\rightarrow h^*(\rT^{(r)}) \oplus h^*(\rT^{(r')})
\end{equation}
the $\co_Y$-linear morphism defined by the matrix 
\begin{equation}\label{p2-hta150b}
\begin{pmatrix}
h'^*(v^{(r)})&h^*(\pi^{(r,r')\vee})\\
0&\id_{h^*(\rT^{(r')})}
\end{pmatrix}.
\end{equation}
where $v^{(r)}$ is defined in \eqref{p2-hta120ee} and $\pi^{(r,r')\vee}$ is the dual of \eqref{p2-hta4e}. In view of \eqref{p2-hta14e}, the diagram 
\begin{equation}\label{p2-hta150c} 
\xymatrix{
{h'^*(\rT')\oplus h^*(\rT)}\ar[r]^-(0.5)w\ar[d]_{h'^*(\pi'^{(r)\vee})\oplus h^*(\pi^{(r')\vee})}&{h^*(\rT)\oplus h^*(\rT)}\ar[d]^{h^*(\pi^{(r)\vee})\oplus 
h^*(\pi^{(r')\vee})}\\
{h'^*(\rT'^{(r)})\oplus h^*(\rT^{(r')})}\ar[r]^-(0.5){w^{(r,r')}}&{h^*(\rT^{(r)}) \oplus h^*(\rT^{(r')}),}}
\end{equation}
where $w$ is the morphism \eqref{p2-hta10f}, is commutative. 
Therefore, the morphism $\uppsi_\cL$ \eqref{p2-hta10e} induces by twisting, a $w^{(r,r')}$-equivariant morphism 
\begin{equation}\label{p2-hta150d}
\uppsi^{(r,r')}_\cL\colon \cL^{(r)}_{\tY/\tX'}\times h'^+(\cL^{(r')}_{\tX'/\tX}) \rightarrow \cL^{(r)}_{\tY/\tX}\times h'^+(\cL^{(r')}_{\tX'/\tX}). 
\end{equation}
By \ref{p1-prem25}, the latter induces a canonical morphism of $\co_{Y}$-algebras \eqref{p1-prem25f} 
\begin{equation}\label{p2-hta150e}
\uppsi^{(r,r')}\colon \cC^{(r)}_{\tY/\tX}\otimes_{\co_{Y}} h'^*(\cC^{(r')}_{\tX'/\tX})\rightarrow 
\cC^{(r)}_{\tY/\tX'} \otimes_{\co_{Y}} h'^*(\cC^{(r')}_{\tX'/\tX}).
\end{equation}

Likewise, the morphism $\upphi_\cL$ \eqref{p2-hta10b} induces by twisting, a canonical morphism 
\begin{equation}\label{p2-hta150f}
\upphi^{(r,r')}_\cL\colon \cL^{(r)}_{\tY/\tX'}\times h'^+(\cL^{(r')}_{\tX'/\tX}) \rightarrow \cL^{(r)}_{\tY/\tX},
\end{equation}
which is equivariant with respect to the $\co_Y$-linear morphism
\begin{equation}\label{p2-hta150g}
h'^*(v^{(r)})+h^*(\pi^{(r,r')\vee})\colon h'^*(\rT'^{(r)})\oplus h^*(\rT^{(r')})\rightarrow h^*(\rT^{(r)}). 
\end{equation} 
By \ref{p1-prem26}, $\upphi^{(r,r')}_\cL$ induces a canonical $\co_{Y}$-linear morphism \eqref{p1-prem26d}
\begin{equation}\label{p2-hta150h}
\upvarphi^{(r,r')}\colon \cF^{(r)}_{\tY/\tX}\rightarrow \cF^{(r)}_{\tY/\tX'}\otimes_{\co_{Y}} h'^*(\cF^{(r')}_{\tX'/\tX}),
\end{equation}
and a canonical morphism of $\co_{Y}$-algebras \eqref{p1-prem26e}
\begin{equation}\label{p2-hta150i}
\upphi^{(r,r')}\colon \cC^{(r)}_{\tY/\tX}\rightarrow \cC^{(r)}_{\tY/\tX'} \otimes_{\co_{Y}} h'^*(\cC^{(r')}_{\tX'/\tX}),
\end{equation}
compatible with $\upvarphi^\ur$.  
Observe that $\uppsi^{(r,r')}$ is none other than the $h'^*(\cC^{(r')}_{\tX'/\tX})$-linearization of $\upphi^{(r,r')}$. 

Let $t\geq t'\geq 0$ be two rational numbers such that $r\geq t$ and $r'\geq t'$. 
We immediately check that the diagram 
\begin{equation}\label{p2-hta150k} 
\xymatrix{
{\cC^{(r)}_{\tY/\tX}\otimes_{\co_{Y}} h'^*(\cC^{(r')}_{\tX'/\tX})}\ar[r]^-(0.5){\uppsi^{(r,r')}}\ar[d]_{\alpha_{\tY/\tX}^{r,t}\otimes h'^*(\alpha_{\tX'/\tX}^{r',t'})}
&{\cC^{(r)}_{\tY/\tX'} \otimes_{\co_{Y}} h'^*(\cC^{(r')}_{\tX'/\tX})}
\ar[d]^{\alpha_{\tY/\tX'}^{r,t}\otimes h'^*(\alpha_{\tX'/\tX}^{r',t'})}\\
{\cC^{(t)}_{\tY/\tX}\otimes_{\co_{Y}} h'^*(\cC^{(t')}_{\tX'/\tX})}\ar[r]^-(0.5){\uppsi^{(t,t')}}&{\cC^{(t)}_{\tY/\tX'} \otimes_{\co_{Y}} h'^*(\cC^{(t')}_{\tX'/\tX}),}}
\end{equation}
where the vertical arrows are defined in \eqref{p2-hta5j}, is commutative. We have a similar commutative diagram for $\upphi^{(r,r')}$.

\subsection{}\label{p2-hta18} 
Let $r$ be a rational number $\geq 0$. We denote by 
\begin{eqnarray}
d_{\cC^{(r)}_{\tY/\tX}}\colon \cC^{(r)}_{\tY/\tX} &\rightarrow& h^*(\Omega^{(r)})\otimes_{\co_Y}\cC^{(r)}_{\tY/\tX},\label{p2-hta18a}\\
d_{\cC^{(r)}_{\tY/\tX'}}\colon \cC^{(r)}_{\tY/\tX'} &\rightarrow& h'^*(\Omega'^{(r)})\otimes_{\co_Y}\cC^{(r)}_{\tY/\tX'},\label{p2-hta18b}
\end{eqnarray}
the universal $\co_Y$-derivations of $\cC^{(r)}_{\tY/\tX}$ and $\cC^{(r)}_{\tY/\tX'}$, respectively, \eqref{p2-hta5g} and by
\begin{equation}\label{p2-hta18c}
d_{\cC^{(r)}_{\tX'/\tX}}\colon \cC^{(r)}_{\tX'/\tX} \rightarrow \cogamma^*(\Omega^{(r)})\otimes_{\co_{\coX'}}\cC^{(r)}_{\tX'/\tX}
\end{equation}
the universal $\co_{\coX'}$-derivation of $\cC^{(r)}_{\tX'/\tX}$. We set 
\begin{eqnarray}
\delta_{\cC^{(r)}_{\tY/\tX}}=(h^*(\pi^{(r)})\otimes \id)\circ d_{\cC^{(r)}_{\tY/\tX}}\colon \cC^{(r)}_{\tY/\tX} &\rightarrow& h^*(\Omega)\otimes_{\co_Y}\cC^{(r)}_{\tY/\tX},\label{p2-hta18d}\\
\delta_{\cC^{(r)}_{\tY/\tX'}}=(h'^*(\pi'^{(r)})\otimes \id)\circ d_{\cC^{(r)}_{\tY/\tX}}\colon \cC^{(r)}_{\tY/\tX'} &\rightarrow& h'^*(\Omega')\otimes_{\co_Y}\cC^{(r)}_{\tY/\tX'},\label{p2-hta18e}\\
\delta_{\cC^{(r)}_{\tX'/\tX}}=(\cogamma^*(\pi^{(r)})\otimes \id)\circ d_{\cC^{(r)}_{\tX'/\tX}}\colon \cC^{(r)}_{\tX'/\tX} &\rightarrow& \cogamma^*(\Omega)\otimes_{\co_{\coX'}}\cC^{(r)}_{\tX'/\tX}.\label{p2-hta18f}
\end{eqnarray}
We also consider the derivations 
\begin{eqnarray}
\delta'_{\cC^{(r)}_{\tY/\tX}}=(h'^*(u)\otimes \id)\circ \delta_{\cC^{(r)}_{\tY/\tX}}\colon \cC^{(r)}_{\tY/\tX} &\rightarrow& h'^*(\Omega')\otimes_{\co_Y}\cC^{(r)}_{\tY/\tX},\label{p2-hta18g}\\
\delta'_{\cC^{(r)}_{\tX'/\tX}}=(u\otimes \id)\circ \delta_{\cC^{(r)}_{\tX'/\tX}}\colon \cC^{(r)}_{\tX'/\tX} &\rightarrow&\Omega'\otimes_{\co_{\coX'}}\cC^{(r)}_{\tX'/\tX},\label{p2-hta18h}
\end{eqnarray}
where $u$ is the morphism \eqref{p2-hta120d}.

For any rational numbers $r\geq r'\geq 0$, we denote by
\begin{equation}\label{p2-hta180k}
\delta^{(r,r')}\colon \cC^{(r)}_{\tY/\tX'}\otimes_{\co_Y}h'^*(\cC^{(r')}_{\tX'/\tX})\rightarrow 
h'^*(\Omega')\otimes_{\co_Y}\cC^{(r)}_{\tY/\tX'}\otimes_{\co_Y}h'^*(\cC^{(r')}_{\tX'/\tX})
\end{equation}
the $\co_Y$-derivation defined by 
\begin{equation}\label{p2-hta180l}
\delta^{(r,r')}=\delta_{\cC^{(r)}_{\tY/\tX'}}\otimes \id- \id \otimes h'^*(\delta'_{\cC^{(r')}_{\tX'/\tX}}).
\end{equation}

\begin{lem}\label{p2-hta190}
For all rational numbers $r\geq r'\geq 0$, the following diagram is commutative 
\begin{equation}\label{p2-hta190a}
\xymatrix{
{\cC^{(r)}_{\tY/\tX}\otimes_{\co_{Y}} h'^*(\cC^{(r')}_{\tX'/\tX})}\ar[r]^-(0.5){\partial^{(r,r')}}\ar[d]_{\uppsi^{(r,r')}}&
{h^*(\Omega\oplus \Omega)\otimes_{\co_Y}\cC^{(r)}_{\tY/\tX}\otimes_{\co_{Y}} h'^*(\cC^{(r')}_{\tX'/\tX})}\ar[d]^{h'^*(\sigma)\otimes \uppsi^{(r,r')}}\\
{\cC^{(r)}_{\tY/\tX'}\otimes_{\co_Y}h'^*(\cC^{(r')}_{\tX'/\tX})}\ar[r]^-(0.5){\updelta^{(r,r')}}&
{(h'^*(\Omega')\oplus h^*(\Omega))\otimes_{\co_Y}\cC^{(r)}_{\tY/\tX'}\otimes_{\co_Y}h'^*(\cC^{(r')}_{\tX'/\tX}),}}
\end{equation}
where $\uppsi^{(r,r')}$ is defined in \eqref{p2-hta150e}, 
\begin{eqnarray}
\partial^{(r,r')}&=&(\delta_{\cC^{(r)}_{\tY/\tX}}\otimes \id) \oplus (\id \otimes h'^*(\delta_{\cC^{(r')}_{\tX'/\tX}})),\\
\updelta^{(r,r')}&=&(\delta_{\cC^{(r)}_{\tY/\tX'}}\otimes \id) \oplus (\id \otimes h'^*(\delta_{\cC^{(r')}_{\tX'/\tX}})),
\end{eqnarray}
the derivations $\delta_{\cC^{(r)}_{\tY/\tX'}}$ and $\delta_{\cC^{(r')}_{\tX'/\tX}}$ being defined in \eqref{p2-hta18e} and \eqref{p2-hta18f}, and 
$\sigma\colon \cogamma^*(\Omega)\oplus \cogamma^*(\Omega)\rightarrow \Omega'\oplus \cogamma^*(\Omega)$ is the $\co_{Y}$-linear morphism represented by the matrix 
\begin{equation}\label{p2-hta190d} 
\begin{pmatrix}
u&0\\
\id_{\cogamma^*(\Omega)}&\id_{\cogamma^*(\Omega)}
\end{pmatrix},
\end{equation}
with $u$ defined in \eqref{p2-hta120d}. 

In particular, we have 
\begin{equation}\label{p2-hta190c}
\delta^{(r,r')}\circ \upphi^{(r,r')}=0,
\end{equation}
where $\delta^{(r,r')}$ is defined in \eqref{p2-hta180l} and $\upphi^{(r,r')}$ in \eqref{p2-hta150i}. 
\end{lem}

The commutativity of diagram \eqref{p2-hta190a} follows from the definition of $\uppsi^{(r,r')}$, especially \eqref{p2-hta150a}, and  \eqref{p1-prem24j}.
The second assertion is an immediate consequence of the first one. 

\subsection{}\label{p2-hta12}
We keep the assumptions and notation of \ref{p2-hta120}, and assume moreover that the following condition is satisfied:
\begin{itemize} 
\item[(viii)] the morphism $\cogamma\colon (\coX',\cM_{\coX'})\rightarrow (\coX,\cM_\coX)$ is smooth.
\end{itemize}
To lighten the notation, we set, with the convention of \ref{p2-ncgt3},   
\begin{eqnarray}
\tOmega^1_{\coX'/\coX}&=&\Omega^1_{(\coX',\cM_{\coX'})/(\coX,\cM_\coX)},\\
\uOmega'&=&\txi^{-1}\tOmega^1_{\coX'/\coX},\label{p2-hta12b}\\
\urT'&=&\cHom_{\co_{\coX'}}(\uOmega',\co_{\coX'}).\label{p2-hta12c}
\end{eqnarray}
We have a canonical exact sequence of locally free $\co_{\coX'}$-modules of finite type
\begin{equation}\label{p2-hta11c}
0\rightarrow \cogamma^*(\tOmega^1_{\coX/\coS})\rightarrow \tOmega^1_{\coX'/\coS} \rightarrow \tOmega^1_{\coX'/\coX} \rightarrow 0. 
\end{equation}
It induces dual exact sequences of locally free $\co_{\coX'}$-modules of finite type 
\begin{eqnarray}
0\rightarrow \cogamma^*(\Omega)\stackrel{u}{\rightarrow} \Omega' \rightarrow \uOmega' \rightarrow 0, \label{p2-hta11g}\\
0\rightarrow \urT'\rightarrow \rT' \stackrel{v}{\rightarrow} \cogamma^*(\rT) \rightarrow 0, \label{p2-hta12d}
\end{eqnarray}
where $u$ (resp.\ $v$) is the morphism \eqref{p2-hta120d} (resp.\ \eqref{p2-hta120e}). 

For any rational number $r\geq 0$, extending the notation of \ref{p2-hta4} and \ref{p2-hta120}, we set 
\begin{eqnarray}
\uOmega'^{(r)}&=&p^r\uOmega',\\
\urT'^{(r)}&=&\cHom_{\co_{\coX'}}(\uOmega'^{(r)},\co_{\coX'}).\label{p2-hta13b}
\end{eqnarray}
The canonical injection $\upi'^{(r)}\colon \uOmega'^{(r)}\rightarrow \uOmega'$ induces an injection $\upi'^{(r)\vee}\colon \urT'\rightarrow \urT'^{(r)}$ 
that identifies $\urT'^{(r)}$ with the submodule $p^{-r}\urT'$ of $\urT'\otimes_{\mZ_p}\mQ_p$.

The exact sequence \eqref{p2-hta11g} induces an exact sequence of locally free $\co_{\coX'}$-modules of finite type
\begin{equation}\label{p2-hta13c}
0\longrightarrow \cogamma^*(\Omega^{(r)})\stackrel{u^{(r)}}{\longrightarrow} \Omega'^{(r)} \longrightarrow \uOmega'^{(r)} \longrightarrow 0. 
\end{equation}
The latter induces by duality an exact sequence of locally free $\co_{\coX'}$-modules of finite type
\begin{equation}\label{p2-hta13d}
0\longrightarrow \urT'^{(r)}\longrightarrow \rT'^{(r)} \stackrel{v^{(r)}}{\longrightarrow} \cogamma^*(\rT^{(r)}) \longrightarrow 0. 
\end{equation}

\subsection{}\label{p2-hta14} 
Let $r,r'$ be rational numbers such that $r\geq r'\geq 0$. We denote by 
\begin{eqnarray}
\pi'^{(r,r')}\colon \Omega'^{(r)}&\rightarrow& \Omega'^{(r')},\label{p2-hta14a}\\
\upi'^{(r,r')}\colon \uOmega'^{(r)}&\rightarrow& \uOmega'^{(r')},\label{p2-hta14aa}
\end{eqnarray} 
the canonical injections and by $\Omega'^{(r,r')}$ the inverse image of $\upi'^{(r,r')}(\uOmega'^{(r)})$ in $\Omega'^{(r')}$, 
so that we have a commutative diagram 
\begin{equation}\label{p2-hta14b}
\xymatrix{
0\ar[r]&{\cogamma^*(\Omega^{(r')})}\ar[r]^{u^{(r')}}&{\Omega'^{(r')}}\ar[r]&{\uOmega'^{(r')}}\ar[r]&0\\
0\ar[r]&{\cogamma^*(\Omega^{(r')})}\ar[r]^{u^{(r,r')}}\ar@{=}[u]&{\Omega'^{(r,r')}}\ar[r]\ar@{^(->}[u]\ar@{}[ru]|{(2)}&{\uOmega'^{(r)}}\ar[r]\ar[u]_{\upi'^{(r,r')}}&0\\
0\ar[r]&{\cogamma^*(\Omega^{(r)})}\ar[r]^{u^{(r)}}\ar[u]^{\cogamma^*(\pi^{(r,r')})}\ar@{}[ru]|{(1)}&{\Omega'^{(r)}}\ar[r]\ar@{^(->}[u]&{\uOmega'^{(r)}}\ar[r]\ar@{=}[u]&0,}
\end{equation}
where square $(1)$ is co-Cartesian, square $(2)$ is Cartesian, the morphism $u^{(r,r')}$ is induced by $u^{(r')}$ 
and the composition of the middle vertical arrows is $\pi'^{(r,r')}$.  We set 
\begin{equation}\label{p2-hta13i}
\rT'^{(r,r')}= \cHom_{\co_{\coX'}}(\Omega'^{(r,r')},\co_{\coX'}).
\end{equation}
The dual of the middle vertical injections of \eqref{p2-hta14b} induce injections 
\begin{equation}\label{p2-hta14c}
\rT'^{(r')}\hookrightarrow \rT'^{(r,r')}\hookrightarrow \rT'^{(r)}.
\end{equation}
The middle horizontal exact sequence of \eqref{p2-hta14b} induces an exact sequence of locally free $\co_{\coX'}$-modules of finite type
\begin{equation}\label{p2-hta14d}
\xymatrix{
0\ar[r]&\urT'^{(r)} \ar[r]& \rT'^{(r,r')}\ar[r]^-(0.5){v^{(r,r')}}&\cogamma^*(\rT^{(r')})\ar[r]& 0.}
\end{equation}
The diagram
\begin{equation}\label{p2-hta14e}
\xymatrix{
{\rT'^{(r,r')}}\ar[r]^-(0.4){v^{(r,r')}}&{\cogamma^*(\rT^{(r')})}\\
{\rT'}\ar[r]^-(0.5)v\ar[u]^{\pi'^{(r,r')\vee}}&{\cogamma^*(\rT),}\ar[u]_{\cogamma^*(\pi^{(r')\vee})}}
\end{equation}
where $\pi'^{(r,r')\vee}$ is the composition of $\pi'^{(r')\vee}$ \eqref{p2-hta13aa} with the first injection in \eqref{p2-hta14c}, is commutative.
Observe that the composition of $\pi'^{(r,r')\vee}$ with the second injection in \eqref{p2-hta14c} is $\pi'^{(r)\vee}$. 

We denote by $\cL^{(r,r')}_{\tY/\tX'}$ the $h'^*(\pi'^{(r,r')\vee})$-twist of $\cL_{\tY/\tX'}$ \eqref{p1-NC4}, 
which is a torsor of $Y_\et$ under the $\co_Y$-module $h'^*(\rT'^{(r,r')})$,
by $\cF^{(r,r')}_{\tY/\tX'}$ the $\co_Y$-module of affine functions on $\cL^{(r,r')}_{\tY/\tX'}$, by 
\begin{equation}\label{p2-hta14f}
0\rightarrow \co_{Y}\rightarrow \cF^{(r,r')}_{\tY/\tX'}\rightarrow h'^*(\Omega'^{(r,r')}) \rightarrow 0
\end{equation}
the associated exact sequence and by
\begin{equation}\label{p2-hta14g}
\cC^{(r,r')}_{\tY/\tX'}=\underset{\underset{n\geq 0}{\longrightarrow}}\lim\ \rS^{n}_{\co_{Y}}(\cF^{(r,r')}_{\tY/\tX'})
\end{equation}
the associated $\co_{Y}$-algebra. 

Let $t,t'$ be rational numbers such that $t\geq t'\geq 0$, $r\geq t$ and $r'\geq t'$. We have a commutative diagram of canonical $\co_{\coX'}$-linear morphisms
\begin{equation}\label{p2-hta14h}
\xymatrix{
0\ar[r]&{\cogamma^*(\Omega^{(r')})}\ar[r]\ar[d]_{\cogamma^*(\pi^{(r',t')})}&{\Omega'^{(r,r')}}\ar[r]\ar[d]^-(0.5){\pi'^{(r,r',t,t')}}&
{\uOmega'^{(r)}}\ar[r]\ar[d]^{\upi'^{(r,t)}}&0\\
0\ar[r]&{\cogamma^*(\Omega^{(t')})}\ar[r]&{\Omega'^{(t,t')}}\ar[r]&{\uOmega'^{(t)}}\ar[r]&0.}
\end{equation}
Denoting by $\pi'^{(r,r',t,t')\vee}\colon \rT'^{(t,t')}\rightarrow \rT'^{(r,r')}$ the dual of the morphism $\pi'^{(r,r',t,t')}$, 
we have a canonical $h'^*(\pi'^{(r,r',t,t')\vee})$-equivariant morphism $\varpi^{(r,r',t,t')}_{\tY/\tX'}\colon \cL^{(t,t')}_{\tY/\tX'}\rightarrow \cL^{(r,r')}_{\tY/\tX'}$.  
It induces by pullback a canonical $\co_{Y}$-linear morphism
\begin{equation}\label{p2-hta14i}
\varpi^{(r,r',t,t')*}_{\tY/\tX'}\colon\cF^{(r,r')}_{\tY/\tX'}\rightarrow \cF^{(t,t')}_{\tY/\tX'}
\end{equation}
that fits into a commutative diagram
\begin{equation}\label{p2-hta14j}
\xymatrix{
0\ar[r]&\co_Y\ar@{=}[d]\ar[r]&{\cF^{(r,r')}_{\tY/\tX'}}\ar[r]\ar[d]^{\varpi^{(r,r',t,t')*}_{\tY/\tX'}}&{h'^*(\Omega'^{(r,r')})}\ar[r]\ar[d]^{h'^*(\pi'^{(r,r',t,t')})}&0\\
0\ar[r]&\co_Y\ar[r]&{\cF^{(t,t')}_{\tY/\tX'}}\ar[r]&{h'^*(\Omega'^{(t,t')})}\ar[r]&0.}
\end{equation}
We deduce a canonical homomorphism of $\co_{Y}$-algebras
\begin{equation}\label{p2-hta14k}
\alpha_{\tY/\tX'}^{r,r',t,t'}\colon \cC_{\tY/\tX'}^{(r,r')}\rightarrow \cC_{\tY/\tX'}^{(t,t')}.
\end{equation}

\subsection{}\label{p2-hta17}
Let $r\geq r'\geq 0$ be rational numbers. 
We denote by $\cHom_{v^{(r,r')}}(\cL^{(r,r')}_{\tY/\tX'},\cL^{(r')}_{\tY/\tX})$ the sheaf of $h'^*(v^{(r,r')})$-equivariant morphisms from 
$\cL^{(r,r')}_{\tY/\tX'}$ to $\cL^{(r')}_{\tY/\tX}$ \eqref{p2-hta14d}, 
which is canonically a torsor of $Y_\et$ under the $\co_Y$-module $h^*(\rT^{(r')})$.
In view of \eqref{p2-hta14e}, this torsor is canonically isomorphic to the $h^*(\pi^{(r')\vee})$-twist of $\cHom_{v}(\cL_{\tY/\tX'},\cL_{\tY/\tX})$ \eqref{p2-hta4dd}. 
Hence, the isomorphism $\upnu$ \eqref{p2-hta10a} induces a canonical $h^*(\rT^{(r')})$-equivariant isomorphism
\begin{equation}\label{p2-hta17a}
\upnu^{(r,r')}\colon h'^+(\cL^{(r')}_{\tX'/\tX})\stackrel{\sim}{\rightarrow} \cHom_{v^{(r,r')}}(\cL^{(r,r')}_{\tY/\tX'},\cL^{(r')}_{\tY/\tX}). 
\end{equation}

Let $\tgamma\in \cL_{\tX'/\tX}(\coX')$ and $\tgamma^{(r')}$ its canonical image in $\cL^{(r')}_{\tX'/\tX}(\coX')$. 
The diagram 
\begin{equation}\label{p2-hta17b}
\xymatrix{
{\cL_{\tY/\tX'}}\ar[rrr]^-(0.5){\upnu(h'^+(\tgamma))}\ar[d]&&&{\cL_{\tY/\tX}}\ar[d]\\
{\cL^{(r,r')}_{\tY/\tX'}}\ar[rrr]^-(0.5){\upnu^{(r,r')}(h'^+(\tgamma^{(r')}))}&&&{\cL^{(r')}_{\tY/\tX},}}
\end{equation}
where $\upnu$ is the morphism \eqref{p2-hta10a}, is commutative. By pullback, $\upnu^{(r,r')}(h'^+(\tgamma^{(r')}))$ induces an $\co_Y$-linear morphism 
\begin{equation}\label{p2-hta17c}
\upvarphi_\tgamma^{(r,r')}\colon \cF^{(r')}_{\tY/\tX}\rightarrow \cF^{(r,r')}_{\tY/\tX'}
\end{equation}
that fits into a commutative diagram 
\begin{equation}\label{p2-hta17d}
\xymatrix{
0\ar[r]&{\co_Y}\ar@{=}[d]\ar[r]&{\cF^{(r')}_{\tY/\tX}}\ar[r]\ar[d]^{\upvarphi_\tgamma^{(r,r')}}&{h^*(\Omega^{(r')})}\ar[r]\ar[d]&0\\
0\ar[r]&{\co_Y}\ar[r]&{\cF^{(r,r')}_{\tY/\tX'}}\ar[r]&{h'^*(\Omega'^{(r,r')})}\ar[r]&0,}
\end{equation}
where the right vertical arrow is the canonical morphism \eqref{p2-hta14b}. We deduce a canonical morphism of $\co_Y$-algebras 
\begin{equation}\label{p2-hta17e}
\upphi_\tgamma^{(r,r')}\colon \cC^{(r')}_{\tY/\tX}\rightarrow \cC^{(r,r')}_{\tY/\tX'}. 
\end{equation}

Let $t,t'$ be rational numbers such that $t\geq t'\geq 0$, $r\geq t$ and $r'\geq t'$. We immediately check that the diagram 
\begin{equation}\label{p2-hta17f} 
\xymatrix{
{\cC^{(r')}_{\tY/\tX}}\ar[r]^-(0.5){\upphi_\tgamma^{(r,r')}}\ar[d]_{\alpha^{t,t'}_{\tY/\tX}}&{\cC^{(r,r')}_{\tY/\tX'}}\ar[d]^{\alpha_{\tY/\tX'}^{r,r',t,t'}}\\
{\cC^{(t')}_{\tY/\tX}}\ar[r]^-(0.5){\upphi_\tgamma^{(t,t')}}&{\cC^{(t,t')}_{\tY/\tX'},}}
\end{equation}
where the vertical arrows are defined in \eqref{p2-hta5j} and \eqref{p2-hta14k}, is commutative.

\begin{lem}\label{p2-hta170}
For all rational numbers $r\geq r'\geq 0$ and every $\tgamma\in \cL_{\tX'/\tX}(\coX')$, the canonical diagram \eqref{p2-hta17f} 
\begin{equation}\label{p2-hta170a}
\xymatrix{
{\cC^{(r)}_{\tY/\tX}}\ar[r]^-(0.5){\upphi_\tgamma^{(r,r)}}\ar[d]_{\alpha^{r,r'}_{\tY/\tX}}&{\cC^{(r)}_{\tY/\tX'}}\ar[d]^{\alpha_{\tY/\tX'}^{r,r,r,r'}}\\
{\cC^{(r')}_{\tY/\tX}}\ar[r]^-(0.5){\upphi_\tgamma^{(r,r')}}&{\cC^{(r,r')}_{\tY/\tX'}}}
\end{equation}
is co-Cartesian. 
\end{lem}

Indeed, we may assume $Y$ affine, in which case there exists $\thh'\in \cL_{\tY/\tX'}(Y)$ \eqref{p2-hta121}. 
We set $\upnu(h'^+(\tgamma))(\thh')=\thh$ \eqref{p2-hta10a}, so that we have $\thh=\tgamma\circ \thh'$. We denote by $\thh^{(r')}$ (resp.\ $\thh'^{(r,r')}$) the 
canonical image of $\thh$ in $\cL^{(r)}_{\tY/\tX}(Y)$ (resp.\ of $\thh'$ in $\cL^{(r,r')}_{\tY/\tX'}(Y)$).  Then, we have a commutative diagram 
\begin{equation}
\xymatrix{
{\rS_{\co_Y}(h^*(\Omega^{(r')}))}\ar[rrr]^-(0.5){\rS_{\co_Y}(h'^*(u^{(r,r')}))}\ar[d]_{\varsigma_{\thh^{(r')}}}&&&{\rS_{\co_Y}(h'^*(\Omega'^{(r,r')}))}\ar[d]^{\varsigma_{\thh'^{(r,r')}}}\\
{\cC^{(r')}_{\tY/\tX}}\ar[rrr]^-(0.5){\upphi_\tgamma^{(r,r')}}&&&{\cC^{(r,r')}_{\tY/\tX'},}}
\end{equation}
where $u^{(r,r')}$ is the morphism defined in \eqref{p2-hta14b}, and 
$\varsigma_{\thh^{(r')}}$ and $\varsigma_{\thh'^{(r,r')}}$ are the isomorphisms \eqref{p2-hta6a} associated with the sections $\thh^{(r')}$ and $\thh'^{(r,r')}$, respectively. 
This follows from the fact that the morphism $\upnu^{(r,r')}(h'^+(\tgamma^{(r')}))$ \eqref{p2-hta17a}  is $h^*(v^{(r,r')})$-equivariant \eqref{p2-hta14d} and the relation \eqref{p2-hta17b}
\begin{equation}
\upnu^{(r,r')}(h'^+(\tgamma^{(r)}))(\thh'^{(r,r')})=\thh^{(r')}. 
\end{equation}
On the other hand, the diagram
\begin{equation}
\xymatrix{
{\rS_{\co_Y}(h^*(\Omega^{(r)}))}\ar[rrr]^-(0.5){\rS_{\co_Y}(h'^*(u^{(r)}))}\ar[d]&&&{\rS_{\co_Y}(h'^*(\Omega'^{(r)}))}\ar[d]\\
{\rS_{\co_Y}(h^*(\Omega^{(r')}))}\ar[rrr]^-(0.5){\rS_{\co_Y}(h'^*(u^{(r,r')}))}&&&{\rS_{\co_Y}(h'^*(\Omega'^{(r,r')})),}}
\end{equation}
where the vertical arrows are the canonical homomorphisms, is co-Cartesian. This is a consequence of the fact that square $(1)$ of diagram \eqref{p2-hta14b} is co-Cartesian. The proposition follows.

\subsection{}\label{p2-hta15} 
We denote by $I$ the set of triples of rational numbers $\ur=(r_1,r_2,r_3)$ such that $r_1\geq r_2\geq r_3\geq 0$. 
For such a triple $\ur$, we denote by 
\begin{equation}\label{p2-hta15a} 
w^\ur\colon h'^*(\rT'^{(r_1,r_2)})\oplus h^*(\rT^{(r_3)})\rightarrow h^*(\rT^{(r_2)}) \oplus h^*(\rT^{(r_3)})
\end{equation}
the $\co_Y$-linear morphism defined by the matrix 
\begin{equation}\label{p2-hta15b}
\begin{pmatrix}
h'^*(v^{(r_1,r_2)})&h^*(\pi^{(r_2,r_3)\vee})\\
0&\id_{h^*(\rT^{(r_3)})}
\end{pmatrix},
\end{equation}
where $v^{(r_1,r_2)}$ is defined in \eqref{p2-hta14d} and $\pi^{(r_2,r_3)\vee}$ is the dual of \eqref{p2-hta4e}. In view of \eqref{p2-hta14e}, the diagram 
\begin{equation}\label{p2-hta15c} 
\xymatrix{
{h'^*(\rT')\oplus h^*(\rT)}\ar[r]^-(0.5)w\ar[d]_{h'^*(\pi'^{(r_1,r_2)\vee})\oplus h^*(\pi^{(r_3)\vee})}&{h^*(\rT)\oplus h^*(\rT)}\ar[d]^{h^*(\pi^{(r_2)\vee})\oplus h^*(\pi^{(r_3)\vee})}\\
{h'^*(\rT'^{(r_1,r_2)})\oplus h^*(\rT^{(r_3)})}\ar[r]^-(0.5){w^\ur}&{h^*(\rT^{(r_2)}) \oplus h^*(\rT^{(r_3)}),}}
\end{equation}
where $w$ is the morphism \eqref{p2-hta10f}, is commutative. Therefore, the morphism $\uppsi_\cL$ \eqref{p2-hta10e} induces by twisting, 
a $w^\ur$-equivariant morphism 
\begin{equation}\label{p2-hta15d}
\uppsi^\ur_\cL\colon \cL^{(r_1,r_2)}_{\tY/\tX'}\times h'^+(\cL^{(r_3)}_{\tX'/\tX}) \rightarrow \cL^{(r_2)}_{\tY/\tX}\times h'^+(\cL^{(r_3)}_{\tX'/\tX}). 
\end{equation}
By \ref{p1-prem25}, the latter induces a canonical morphism of $\co_{Y}$-algebras \eqref{p1-prem25f} 
\begin{equation}\label{p2-hta15e}
\uppsi^\ur\colon \cC^{(r_2)}_{\tY/\tX}\otimes_{\co_{Y}} h'^*(\cC^{(r_3)}_{\tX'/\tX})\rightarrow \cC^{(r_1,r_2)}_{\tY/\tX'} \otimes_{\co_{Y}} h'^*(\cC^{(r_3)}_{\tX'/\tX}).
\end{equation}

Likewise, the morphism $\upphi_\cL$ \eqref{p2-hta10b} induces by twisting, a canonical morphism 
\begin{equation}\label{p2-hta15f}
\upphi^\ur_\cL\colon \cL^{(r_1,r_2)}_{\tY/\tX'}\times h'^+(\cL^{(r_3)}_{\tX'/\tX}) \rightarrow \cL^{(r_2)}_{\tY/\tX},
\end{equation}
which is equivariant with respect to the $\co_Y$-linear morphism
\begin{equation}\label{p2-hta15g}
h'^*(v^{(r_1,r_2)})+h^*(\pi^{(r_2,r_3)\vee})\colon h'^*(\rT'^{(r_1,r_2)})\oplus h^*(\rT^{(r_3)})\rightarrow h^*(\rT^{(r_2)}). 
\end{equation} 
By \ref{p1-prem26}, $\upphi^\ur_\cL$ induces a canonical $\co_{Y}$-linear morphism \eqref{p1-prem26d}
\begin{equation}\label{p2-hta15h}
\upvarphi^\ur\colon \cF^{(r_2)}_{\tY/\tX}\rightarrow \cF^{(r_1,r_2)}_{\tY/\tX'}\otimes_{\co_{Y}} h'^*(\cF^{(r_3)}_{\tX'/\tX}),
\end{equation}
and a canonical morphism of $\co_{Y}$-algebras \eqref{p1-prem26e}
\begin{equation}\label{p2-hta15i}
\upphi^\ur\colon \cC^{(r_2)}_{\tY/\tX}\rightarrow \cC^{(r_1,r_2)}_{\tY/\tX'} \otimes_{\co_{Y}} h'^*(\cC^{(r_3)}_{\tX'/\tX}),
\end{equation}
compatible with $\upvarphi^\ur$.  
Observe that $\uppsi^\ur$ is none other than the $h'^*(\cC^{(r_3)}_{\tX'/\tX})$-linearization of $\upphi^\ur$.

The canonical $\pi^{(r_2)\vee}$-equivariant morphism $\cL_{\tX'/\tX} \rightarrow \cL^{(r_2)}_{\tX'/\tX}$ \eqref{p2-hta4dd} induces a
canonical $\pi^{(r_2,r_3)\vee}$-equivariant morphism $\cL^{(r_3)}_{\tX'/\tX} \rightarrow \cL^{(r_2)}_{\tX'/\tX}$ \eqref{p2-hta4e}. Composing with $\upnu^{(r_1,r_2)}$ \eqref{p2-hta17a},
we get a canonical $h^*(\pi^{(r_2,r_3)\vee})$-equivariant morphism
\begin{equation}\label{p2-hta15j}
\upnu^\ur\colon h'^+(\cL^{(r_3)}_{\tX'/\tX})\rightarrow \cHom_{v^{(r_1,r_2)}}(\cL^{(r_1,r_2)}_{\tY/\tX'},\cL^{(r_2)}_{\tY/\tX}). 
\end{equation}
We immediately see that $\upphi^\ur_\cL$ \eqref{p2-hta15f} and $\uppsi^\ur_\cL$ \eqref{p2-hta15d} are deduced from $\upnu^\ur$ by ``evaluation'' in the same way as
$\upphi_\cL$ \eqref{p2-hta10b} and $\uppsi_\cL$ \eqref{p2-hta10e} were deduced from $\upnu$.

Let $\ur=(r_1,r_2,r_3)$, $\ur'=(r'_1,r'_2,r'_3)$ be two elements of $I$ such that $r_i\geq r'_i$ for all $1\leq i\leq 3$. 
We immediately check that the diagram 
\begin{equation}\label{p2-hta15k} 
\xymatrix{
{\cC^{(r_2)}_{\tY/\tX}\otimes_{\co_{Y}} h'^*(\cC^{(r_3)}_{\tX'/\tX})}\ar[r]^-(0.5){\uppsi^\ur}\ar[d]_{\alpha_{\tY/\tX}^{r_2,r'_2}\otimes h'^*(\alpha_{\tX'/\tX}^{r_3,r'_3})}
&{\cC^{(r_1,r_2)}_{\tY/\tX'} \otimes_{\co_{Y}} h'^*(\cC^{(r_3)}_{\tX'/\tX})}
\ar[d]^{\alpha_{\tY/\tX'}^{r_1,r_2,r'_1,r'_2}\otimes h'^*(\alpha_{\tX'/\tX}^{r_3,r'_3})}\\
{\cC^{(r'_2)}_{\tY/\tX}\otimes_{\co_{Y}} h'^*(\cC^{(r'_3)}_{\tX'/\tX})}\ar[r]^-(0.5){\uppsi^{\ur'}}&{\cC^{(r'_1,r'_2)}_{\tY/\tX'} \otimes_{\co_{Y}} h'^*(\cC^{(r'_3)}_{\tX'/\tX}),}}
\end{equation}
where the vertical arrows are defined in \eqref{p2-hta5j} and \eqref{p2-hta14k}, is commutative. We have a similar commutative diagram for $\upphi^\ur$. 

\begin{rema}\label{p2-hta151}
For all rational numbers $r\geq r'\geq 0$, setting $\ur=(r,r,r') \in I$ \eqref{p2-hta15}, we have 
$\uppsi^{(r,r')}_\cL=\uppsi^\ur_\cL$, see \eqref{p2-hta150d} and \eqref{p2-hta15d},  so 
$\uppsi^{(r,r')}=\uppsi^\ur$, see \eqref{p2-hta150e} and \eqref{p2-hta15e}; and 
$\upphi^{(r,r')}_\cL=\upphi^\ur_\cL$, see \eqref{p2-hta150f} and \eqref{p2-hta15f},  so 
$\upphi^{(r,r')}=\upphi^\ur$, see \eqref{p2-hta150i} and \eqref{p2-hta15i}. 
\end{rema}

\subsection{}\label{p2-hta16} 
Let $r,r'$ be two rational numbers such that $r\geq r'\geq 0$, $\tgamma\in \cL_{\tX'/\tX}(\coX')$. 
Consider the $\co_{\coX}$-linear isomorphism 
\begin{equation}\label{p2-hta16a} 
\lambda^{(r,r')}\colon \rT^{(r)}\oplus \rT^{(r')}\rightarrow \rT^{(r)} \oplus \rT^{(r')}
\end{equation}
defined by the matrix 
\begin{equation}\label{p2-hta16b} 
\begin{pmatrix}
\id_{\rT^{(r)}}&\pi^{(r,r')\vee}\\
0&\id_{\rT^{(r')}}
\end{pmatrix},
\end{equation}
where $\pi^{(r,r')\vee}$ is the dual of \eqref{p2-hta4e}.
There is a canonical $h^*(\lambda^{(r,r')})$-equivariant isomorphism
\begin{equation}\label{p2-hta16c} 
\Lambda^{(r,r')}_{\cL,\tgamma}\colon \cL^{(r)}_{\tY/\tX} \times h'^+(\cL^{(r')}_{\tX'/\tX}) \stackrel{\sim}{\rightarrow} \cL^{(r)}_{\tY/\tX} \times h'^+(\cL^{(r')}_{\tX'/\tX})
\end{equation}
that maps a local section $(s,s')$ to $(s+\pi^{(r,r')\vee}(s'-h'^+(\tgamma^{(r')})),s')$, where $\tgamma^{(r')}$ is the canonical image of $\tgamma$ in $\cL^{(r')}_{\tX'/\tX}(\coX')$. 
By \ref{p1-prem25}, $\Lambda^{(r,r')}_{\cL,\tgamma}$ induces a canonical isomorphism of $\co_{Y}$-algebras
\begin{equation}\label{p2-hta16d}
\Lambda^{(r,r')}_\tgamma\colon \cC^{(r)}_{\tY/\tX}\otimes_{\co_{Y}} h'^*(\cC^{(r')}_{\tX'/\tX})
\stackrel{\sim}{\rightarrow} \cC^{(r)}_{\tY/\tX}\otimes_{\co_{Y}} h'^*(\cC^{(r')}_{\tX'/\tX}).
\end{equation}

For all rational numbers $t,t'$ such that $t\geq t'\geq 0$, $r\geq t$ and $r'\geq t'$, we immediately see that the diagram 
\begin{equation}\label{p2-hta16i} 
\xymatrix{
{\cC^{(r)}_{\tY/\tX}\otimes_{\co_{Y}} h'^*(\cC^{(r')}_{\tX'/\tX})}\ar[r]^-(0.5){\Lambda^{(r,r')}_\tgamma}\ar[d]_{\alpha^{r,t}_{\tY/\tX}\otimes_{\co_{Y}} h'^*(\alpha^{r',t'}_{\tX'/\tX})}&
{\cC^{(r)}_{\tY/\tX}\otimes_{\co_{Y}} h'^*(\cC^{(r')}_{\tX'/\tX})}\ar[d]^{\alpha^{r,t}_{\tY/\tX}\otimes_{\co_{Y}} h'^*(\alpha^{r',t'}_{\tX'/\tX})}\\
{\cC^{(t)}_{\tY/\tX}\otimes_{\co_{Y}} h'^*(\cC^{(t')}_{\tX'/\tX})}\ar[r]^-(0.5){\Lambda^{(t,t')}_\tgamma}&{\cC^{(t)}_{\tY/\tX}\otimes_{\co_{Y}} h'^*(\cC^{(t')}_{\tX'/\tX}),}}
\end{equation}
where the vertical arrows are defined in \eqref{p2-hta5j}, is commutative. 

\begin{lem}\label{p2-hta160} 
Let $\ur=(r_1,r_2,r_3)\in I$ \eqref{p2-hta15}, $\tgamma\in \cL_{\tX'/\tX}(\coX')$. Then, 
\begin{itemize}
\item[{\rm (i)}] We have 
\begin{eqnarray}
\upvarphi_\tgamma^{(r_1,r_2)}&=&(\id \otimes_{\co_Y} h'^*(\rho_{\tgamma^{(r_3)}}))\circ \upvarphi^\ur,\label{p2-hta160c}\\
\upphi_\tgamma^{(r_1,r_2)}&=&(\id \otimes_{\co_Y} h'^*(\varrho_{\tgamma^{(r_3)}}))\circ \upphi^\ur,\label{p2-hta160d} 
\end{eqnarray}
where $\upvarphi^\ur$, $\upphi^\ur$, $\upvarphi_\tgamma^{(r_1,r_2)}$ and $\upphi_\tgamma^{(r_1,r_2)}$ are defined in 
\eqref{p2-hta15h}, \eqref{p2-hta15i}, \eqref{p2-hta17c} and \eqref{p2-hta17e}, $\tgamma^{(r_3)}$ is the canonical image of $\tgamma$ in $\cL^{(r_3)}_{\tX'/\tX}(\coX')$,
$\rho_{\tgamma^{(r_3)}}\colon \cF_{\tX'/\tX}^{(r_3)}\rightarrow \co_{\coX'}$ is the associated splitting of the extension \eqref{p2-hta13f} and 
$\varrho_{\tgamma^{(r_3)}}\colon \cC_{\tX'/\tX}^{(r_3)}\rightarrow \co_{\coX'}$ is the associated homomorphism of $\co_{\coX'}$-algebras \eqref{p1-prem1}. 
\item[{\rm (ii)}] The diagram 
\begin{equation}\label{p2-hta160a} 
\xymatrix{
{\cL^{(r_1,r_2)}_{\tY/\tX'}\times h'^+(\cL^{(r_3)}_{\tX'/\tX})}\ar[rd]^{\uppsi^\ur_\cL}\ar[d]_{\upnu^{\ur}(h'^+(\tgamma^{(r_3)}))\times \id}&\\
{\cL^{(r_2)}_{\tY/\tX}\times h'^+(\cL^{(r_3)}_{\tX'/\tX})}\ar[r]^{\Lambda^{(r_2,r_3)}_{\cL,\tgamma}}&{\cL^{(r_2)}_{\tY/\tX}\times h'^+(\cL^{(r_3)}_{\tX'/\tX})}}
\end{equation}
is commutative. 
\item[{\rm (iii)}] We have 
\begin{equation}\label{p2-hta160b} 
\uppsi^\ur=(\upphi_\tgamma^{(r_1,r_2)}\otimes \id)\circ \Lambda^{(r_2,r_3)}_\tgamma.
\end{equation} 
\end{itemize}
\end{lem}

(i) This follows from \ref{p1-prem13}. 

(ii) Indeed, for every local section $s$ of $\cL^{(r_1,r_2)}_{\tY/\tX'}$, we have 
\begin{eqnarray}
\Lambda^{(r_2,r_3)}_{\cL,\tgamma}(\upnu^{\ur}(h'^+(\tgamma^{(r_3)}))(s),h'^+(\tgamma^{(r_3)}))&=&(\upnu^{\ur}(h'^+(\tgamma^{(r_3)}))(s),h'^+(\tgamma^{(r_3)}))\\
&=&\uppsi^\ur_\cL(s,h'^+(\tgamma^{(r_3)})).\nonumber\
\end{eqnarray}
Moreover, $\upnu^{\ur}(h'^+(\tgamma^{(r_3)}))$ is $h'^*(v^{(r_1,r_2)})$-equivariant \eqref{p2-hta15j}, $\Lambda^{(r_2,r_3)}_{\cL,\tgamma}$ is $h^*(\lambda^{(r_2,r_3)})$-equivariant, 
$\uppsi^\ur_\cL$ is $w^\ur$-equivariant \eqref{p2-hta15a} and we have 
\begin{equation}
w^\ur=h^*(\lambda^{(r_2,r_3)}) \circ (h'^*(v^{(r_1,r_2)})\oplus \id_{h^*(\rT^{(r_3)})}).
\end{equation}
Hence, the diagram \eqref{p2-hta160a} is commutative. 

(iii) It follows immediately from (ii). 

\subsection{}\label{p2-hta180} 
Let $r,r'$ be rational numbers such that $r\geq r'\geq 0$. Recall \eqref{p1-thbn1d} that there is a unique $\co_Y$-derivation 
\begin{equation}\label{p2-hta18i}
d_{\cC^{(r,r')}_{\tY/\tX'}}\colon \cC^{(r,r')}_{\tY/\tX'} \rightarrow h'^*(\Omega'^{(r,r')})\otimes_{\co_Y}\cC^{(r,r')}_{\tY/\tX'}
\end{equation}
extending the canonical morphism $\cF^{(r,r')}_{\tY/\tX'}\rightarrow h'^*(\Omega'^{(r,r')})$ \eqref{p2-hta14f}. 
It identifies canonically with the universal $\co_Y$-derivation of $\cC^{(r,r')}_{\tY/\tX'}$. It is also a Higgs $\co_Y$-field on $\cC^{(r)}_{\tY/\tX'}$ with coefficients in $h'^*(\Omega'^{(r,r')})$. 
We denote by  
\begin{equation}\label{p2-hta18j}
\delta_{\cC^{(r,r')}_{\tY/\tX'}}\colon \cC^{(r,r')}_{\tY/\tX'} \rightarrow h'^*(\Omega')\otimes_{\co_Y}\cC^{(r,r')}_{\tY/\tX'}
\end{equation}
the $\co_Y$-derivation induced by $d_{\cC^{(r,r')}_{\tY/\tX'}}$ and the canonical injection $\Omega'^{(r,r')}\rightarrow \Omega'$. 

For any $\ur\in I$ \eqref{p2-hta15}, we denote by
\begin{equation}\label{p2-hta18k}
\delta^\ur\colon \cC^{(r_1,r_2)}_{\tY/\tX'}\otimes_{\co_Y}h'^*(\cC^{(r_3)}_{\tX'/\tX})\rightarrow h'^*(\Omega')\otimes_{\co_Y}\cC^{(r_1,r_2)}_{\tY/\tX'}\otimes_{\co_Y}h'^*(\cC^{(r_3)}_{\tX'/\tX})
\end{equation}
the $\co_Y$-derivation defined by 
\begin{equation}\label{p2-hta18l}
\delta^\ur=\delta_{\cC^{(r_1,r_2)}_{\tY/\tX'}}\otimes \id- \id \otimes h'^*(\delta'_{\cC^{(r_3)}_{\tX'/\tX}}).
\end{equation}

\begin{lem}\label{p2-hta21}
Let $r,r'$ be two rational numbers such that $r\geq r'\geq 0$, $\tgamma\in \cL_{\tX'/\tX}(\coX')$.
Then, 
\begin{itemize}
\item[{\rm (i)}] The diagram  
\begin{equation}\label{p2-hta21a}
\xymatrix{
{\cC^{(r')}_{\tY/\tX}}\ar[r]^-(0.5){\delta_{\cC^{(r')}_{\tY/\tX}}}\ar[d]_{\upphi_\tgamma^{(r,r')}}&{h^*(\Omega)\otimes_{\co_Y}\cC^{(r')}_{\tY/\tX}}\ar[d]^{h'^*(u)\otimes_{\co_Y}\upphi_\tgamma^{(r,r')}}\\
{\cC^{(r,r')}_{\tY/\tX'}}\ar[r]^-(0.5){\delta_{\cC^{(r,r')}_{\tY/\tX'}}}&{h'^*(\Omega')\otimes_{\co_Y}\cC^{(r,r')}_{\tY/\tX'},}}
\end{equation}
where $\delta_{\cC^{(r)}_{\tY/\tX}}$ is defined in \eqref{p2-hta18d}, $\delta_{\cC^{(r')}_{\tX'/\tX}}$ in \eqref{p2-hta18f},
$u$ in \eqref{p2-hta120d} and $\upphi_\tgamma^{(r,r')}$ in \eqref{p2-hta17e}, is commutative. 
\item[{\rm (ii)}] The diagram  
\begin{equation}\label{p2-hta21b}
\xymatrix{
{\cC^{(r)}_{\tY/\tX}\otimes_{\co_{Y}} h'^*(\cC^{(r')}_{\tX'/\tX})}\ar[d]_{\Lambda^{(r,r')}_\tgamma}\ar[r]^-(0.5){\partial^{(r,r')}}&
{h^*(\Omega\oplus\Omega)\otimes_{\co_Y} \cC^{(r)}_{\tY/\tX}\otimes_{\co_{Y}} h'^*(\cC^{(r')}_{\tX'/\tX})}\ar[d]^{h^*(\mu)\otimes_{\co_Y}\Lambda^{(r,r')}_\tgamma}\\
{\cC^{(r)}_{\tY/\tX}\otimes_{\co_{Y}} h'^*(\cC^{(r')}_{\tX'/\tX})}\ar[r]^-(0.5){\partial^{(r,r')}}&{h^*(\Omega\oplus\Omega)\otimes_{\co_Y}\cC^{(r)}_{\tY/\tX}\otimes_{\co_{Y}} h'^*(\cC^{(r')}_{\tX'/\tX}),}}
\end{equation}
where $\Lambda^{(r,r')}_\tgamma$ is defined in \eqref{p2-hta16d}, 
\begin{equation}\label{p2-hta21d}
\partial^{(r,r')}=(\delta_{\cC^{(r)}_{\tY/\tX}}\otimes \id) \oplus (\id \otimes h'^*(\delta_{\cC^{(r')}_{\tX'/\tX}})),
\end{equation}
and $\mu$ is the $\co_{\coX}$-linear automorphism of $\Omega\oplus \Omega$ represented by the matrix 
\begin{equation}\label{p2-hta21c} 
\begin{pmatrix}
\id_\Omega&0\\
\id_\Omega&\id_\Omega
\end{pmatrix}.
\end{equation}
\end{itemize}
\end{lem}

(i) This follows from \eqref{p2-hta17d}. 

(ii) This follows from the definition of $\Lambda^{(r,r')}_\tgamma$, especially \eqref{p2-hta16a}, and \eqref{p1-prem24j}.

\begin{lem}\label{p2-hta19}
For every $\ur=(r_1,r_2,r_3)\in I$ \eqref{p2-hta15}, the following diagram is commutative 
\begin{equation}\label{p2-hta19a}
\xymatrix{
{\cC^{(r_2)}_{\tY/\tX}\otimes_{\co_{Y}} h'^*(\cC^{(r_3)}_{\tX'/\tX})}\ar[r]^-(0.5){\partial^{(r_2,r_3)}}\ar[d]_{\uppsi^\ur}&
{h^*(\Omega\oplus \Omega)\otimes_{\co_Y}\cC^{(r_2)}_{\tY/\tX}\otimes_{\co_{Y}} h'^*(\cC^{(r_3)}_{\tX'/\tX})}\ar[d]^{h'^*(\sigma)\otimes \uppsi^\ur}\\
{\cC^{(r_1,r_2)}_{\tY/\tX'}\otimes_{\co_Y}h'^*(\cC^{(r_3)}_{\tX'/\tX})}\ar[r]^-(0.5){\updelta^\ur}&
{(h'^*(\Omega')\oplus h^*(\Omega))\otimes_{\co_Y}\cC^{(r_1,r_2)}_{\tY/\tX'}\otimes_{\co_Y}h'^*(\cC^{(r_3)}_{\tX'/\tX}),}}
\end{equation}
where $\uppsi^\ur$ is defined in \eqref{p2-hta15e}, $\partial^{(r_2,r_3)}$ in \eqref{p2-hta21d},  
\begin{equation}\label{p2-hta19b}
\updelta^\ur=(\delta_{\cC^{(r_1,r_2)}_{\tY/\tX'}}\otimes \id) \oplus (\id \otimes h'^*(\delta_{\cC^{(r_3)}_{\tX'/\tX}})),
\end{equation}
the derivations $\delta_{\cC^{(r_1,r_2)}_{\tY/\tX'}}$ and $\delta_{\cC^{(r_3)}_{\tX'/\tX}}$ being defined in \eqref{p2-hta18j} and \eqref{p2-hta18f}, and 
$\sigma\colon \cogamma^*(\Omega)\oplus \cogamma^*(\Omega)\rightarrow \Omega'\oplus \cogamma^*(\Omega)$ is the $\co_{Y}$-linear morphism represented by the matrix 
\begin{equation}\label{p2-hta19d} 
\begin{pmatrix}
u&0\\
\id_{\cogamma^*(\Omega)}&\id_{\cogamma^*(\Omega)}
\end{pmatrix},
\end{equation}
with $u$ defined in \eqref{p2-hta120d}. 

In particular, we have 
\begin{equation}\label{p2-hta19c}
\delta^\ur\circ \upphi^\ur=0,
\end{equation}
where $\delta^\ur$ is defined in \eqref{p2-hta18l} and $\upphi^\ur$ in \eqref{p2-hta15i}. 
\end{lem}

The commutativity of diagram \eqref{p2-hta19a} follows from the definition of $\uppsi^\ur$, especially \eqref{p2-hta15a}, and  \eqref{p1-prem24j}.
The second assertion is an immediate consequence of the first one.

\section{\texorpdfstring{Revisiting the local $p$-adic Simpson correspondence}{Revisiting the local p-adic Simpson correspondence}}\label{p2-rlps}

\subsection{}\label{p2-rlps1}
In this section, we let $f\colon (X,\cM_{X})\rightarrow (S,\cM_S)$ be an adequate morphism of fine logarithmic schemes, in the sense of (\cite{agt} III.4.7),
having an adequate chart $\mN\rightarrow P$ (\cite{agt} III.4.4), such that $X=\Spec(R)$ is affine and $X_s$ is non-empty. 
We denote by $X^\circ$ the maximal open subscheme of $X$ where the logarithmic structure $\cM_X$ is trivial; it is an open subscheme of $X_\eta$.
For any $X$-scheme $U$, we set
\begin{equation}\label{p2-rlps1a}
U^\circ=U\times_XX^\circ.
\end{equation}
To lighten the notation, we set 
\begin{eqnarray}
\tOmega^1_{X/S}=\Omega^1_{(X,\cM_{X})/(S,\cM_S)},\label{p2-rlps1b}\\
\tOmega^1_{R/\co_K}=\tOmega^1_{X/S}(X).\label{p2-rlps1c}
\end{eqnarray}

We endow $\coX=X\times_S\coS$ \eqref{p2-ncgt1a} with the logarithmic structure $\cM_\coX$ pullback of $\cM_X$, 
and denote by $\cof\colon (\coX,\cM_{\coX})\rightarrow (\coS,\cM_{\coS})$ the base change of $f$. 
We fix a Cartesian diagram of $\FLS$ \eqref{p1-NC1}
\begin{equation}\label{p2-rlps1d}
\xymatrix{
{(\coX,\cM_{\coX})}\ar[r]^-(0.5){i}\ar[d]_{\cof}\ar@{}[rd]|{\Box}&{(\tX,\cM_{\tX})}\ar[d]^{\tf}\\
{(\coS,\cM_{\coS})}\ar[r]^-(0.5){\iota}&{(\tS,\cM_{\tS}),}}
\end{equation}
where $\iota$ is the strict closed immersion defined in \eqref{p2-ncgt3b}, such that $\tf$ is smooth. 
Observe that this diagram is also Cartesian in the category of logarithmic schemes by \ref{p1-NC3}.
Such an $(\tS,\cM_{\tS})$-smooth deformation $\tf$ of $\cof$ exists and is unique up to isomorphism by virtue of (\cite{kato1} 3.14).

\subsection{}\label{p2-rlps2}
Let $\oy$ be a geometric point of $\oX^\circ$. The scheme $\oX$ being locally irreducible by (\cite{ag1} 4.2.7(iii)),
it is the sum of the schemes induced on its irreducible components. We denote by $\oX^\star$ the irreducible component of $\oX$ containing $\oy$.
Likewise, $\oX^\circ=\oX\times_XX^\circ$ is the sum of the schemes induced on its irreducible components
and $\oX^{\star \circ}=\oX^\star\times_{X}X^\circ$ is the irreducible component of $\oX^\circ$ containing $\oy$. 
We set $\coX^\star=\oX^\star\times_\oS\coS$ \eqref{p2-ncgt1}, 
\begin{eqnarray}
R_1&=&\Gamma(\oX^\star,\co_{\oX}),\label{p2-rlps2a}\\
R_\uptau&=&\Gamma(\coX^\star,\co_{\coX})=R_1\otimes_{\co_\oK}\co_C,\label{p2-rlps2aa}
\end{eqnarray}
and denote by $\hRun$ and $\hRtau$ their $p$-adic Hausdorff completions, which are equal.  
To lighten the notation, for any rational number $r\geq 0$, we set, with the conventions of \ref{p2-ncgt3}, 
\begin{equation}\label{p2-rlps2c}
\Omega=\txi^{-1}\tOmega^1_{R/\co_K}\otimes_RR_\uptau \ \ \ {\rm and}\ \ \ \Omega^{(r)}=p^r\Omega.
\end{equation}
For any $R_\uptau$-algebra $A$, we  consider Higgs $A$-modules with coefficients in 
$\Omega\otimes_{R_\uptau}A$ \eqref{p1-delta-con1}. We say abusively that they have coefficients in $\Omega$.
The category of these modules will be denoted by $\bHM(A,\Omega)$.

We denote by $\Delta$ the profinite group $\pi_1(\oX^{\star \circ},\oy)$ and by $(V_i)_{i\in I}$ the normalized universal cover of
$\oX^{\star \circ}$ at $\oy$ (\cite{ag1} 2.1.20). For any $i\in I$, we denote by $\oX^{V_i}$ the integral closure of $\oX$ in $V_i$.
The schemes $(\oX^{V_i})_{i\in I}$ then form a filtered  inverse system. We set
\begin{equation}\label{p2-rlps2b}
\oR=\underset{\underset{i\in I}{\longrightarrow}}{\lim}\ \Gamma(\oX^{V_i},\co_{\oX^{V_i}}).
\end{equation}
This is a normal integral domain (\cite{ag1} 4.1.10), on which $\Delta$ acts naturally by ring homomorphisms. 
We denote by $\hoR$ its $p$-adic Hausdorff completion
and by $\bRep_{\hoR}(\Delta)$ the category of $\hoR$-representations of $\Delta$ (\cite{ag2} 2.1.2). 

\subsection{}\label{p2-rlps3}
We set
\begin{equation}\label{p2-rlps3a}
\mX=\Spec(\oR)\ \ \ {\rm and} \ \ \ \hmX=\Spec(\hoR),
\end{equation}
that we endow with the logarithmic structures pullbacks of $\cM_X$,
denoted respectively by $\cM_\mX$ and $\cM_\hmX$.
The actions of $\Delta$ on $\oR$ and $\hoR$ induce left actions on
the logarithmic schemes $(\mX,\cM_\mX)$ and $(\hmX,\cM_\hmX)$; for every $\sigma\in \Delta$, the automorphism of $(\mX,\cM_\mX)$ (resp.\ $(\hmX,\cM_\hmX)$) defined by $\sigma$ is induced 
by the automorphism $\sigma^{-1}$ of $\oR$ (resp.\ $\hoR$). 
Endowing $(\coX,\cM_\coX)$ with the trivial action of $\Delta$, we have a canonical $\Delta$-equivariant morphism 
\begin{equation}\label{p2-rlps3b}
h\colon (\hmX,\cM_\hmX)\rightarrow (\coX,\cM_\coX).
\end{equation}
We set $g=\cof\circ h$. 

Following (\cite{ag2} 3.2.11), we denote by $\tmX$ one of the two schemes
\begin{equation}\label{p2-rlps3c}
\cA_2(\mX)=\Spec(\cA_2(\oR))\ \ \ {\rm or} \ \ \ \cA^{\ast}_2(\mX/S)=\Spec(\cA^{\ast}_2(\oR/\co_K)),
\end{equation}
defined in (\cite{ag2} 3.2.8 and 3.2.10), depending on whether we are in the absolute or relative case \eqref{p2-ncgt3} respectively. 
Let $\umi\colon \hmX\rightarrow \tmX$ be the closed immersion defined by the square zero ideal $\txi\co_\tmX$ of $\co_\tmX$, 
which corresponds to Fontaine's homomorphism $\theta$ (see \cite{ag2} (3.2.8.3) and (3.2.10.5)). The $\co_\hmX$-module $\txi\co_\tmX$ is invertible. 

We equip $\tmX$ with a canonical logarithmic structure $\cM'_\tmX$ defined in (\cite{ag2} 3.2.8 and 3.2.10), i.e. that depends only on $(X,\cM_X)$ and $\oy$.  
Fontaine's homomorphism $\theta$ defines a canonical morphism (\cite{ag2} (3.2.8.8) and (3.2.10.9))
\begin{equation}\label{p2-rlps3d}
\mj\colon (\hmX,\cM_\hmX)\rightarrow (\tmX,\cM'_{\tmX}),
\end{equation}
whose underlying morphism of schemes is $\umi$. 
We do not know in general if the logarithmic scheme $(\tmX,\cM'_{\tmX})$ is fine and saturated and
if $\mj$ is an exact closed immersion (see however \ref{p2-rlps4}). This is why we endow $\tmX$
with another logarithmic structure $\cM_{\tmX}$ induced by the given chart $P$ \eqref{p2-rlps1} defined in (\cite{ag2} 3.2.11). 
The logarithmic scheme $(\tmX,\cM_{\tmX})$ is fine and saturated.
There is a canonical morphism of logarithmic structures on $\tmX$
\begin{equation}\label{p2-rlps3e}
\cM_{\tmX}\rightarrow \cM'_{\tmX}.
\end{equation}
We have a canonical strict closed immersion
\begin{equation}\label{p2-rlps3f}
\mi\colon (\hmX,\cM_\hmX)\rightarrow (\tmX,\cM_{\tmX})
\end{equation}
which factors through $\mj$, and a canonical morphism $\tg\colon (\tmX,\cM_{\tmX})\rightarrow (\tS,\cM_\tS)$ that fit into a Cartesian diagram 
\begin{equation}\label{p2-rlps3g}
\xymatrix{
{(\hmX,\cM_\hmX)}\ar[r]^{\mi}\ar[d]_g&{(\tmX,\cM_{\tmX})}\ar[d]^\tg\\
{(\coS,\cM_{\coS})}\ar[r]^{\iota}&{(\tS,\cM_\tS).}}
\end{equation}

The fundamental group $\Delta$ acts on the left on the logarithmic schemes $(\tmX,\cM_{\tmX})$ and $(\tmX,\cM'_{\tmX})$, and the morphisms
$\mi$, $\mj$ and $\tg$ are $\Delta$-equivariant (see \cite{ag2} 3.2.9 and 3.2.10).

\begin{prop}[\cite{agt} II.9.13]\label{p2-rlps4}
Suppose that there exists a fine and saturated chart $M\rightarrow \Gamma(X,\cM_X)$ 
for $(X,\cM_X)$ inducing an isomorphism
\begin{equation}\label{p2-rlps4a}
M\stackrel{\sim}{\rightarrow} \Gamma(X,\cM_X)/\Gamma(X,\co^\times_X).
\end{equation}
Then, the morphism $\cM_{\tmX}\rightarrow \cM'_{\tmX}$ \eqref{p2-rlps3e} is an isomorphism.
In particular, the morphism $\mj$ \eqref{p2-rlps3d} is a strict closed immersion.
\end{prop}

Note that the condition \eqref{p2-rlps4a} is satisfied on an open affine covering of $X$ by (\cite{agt} II.5.17).

\subsection{}\label{p2-rlps5}
The conditions \ref{p2-hta1}(i)-(iv) are satisfied by the commutative diagram (see \eqref{p2-rlps1d} and \eqref{p2-rlps3f})
\begin{equation}
\xymatrix{
{(\hmX,\cM_\hmX)}\ar[r]^-(0.5){\mi}\ar[d]^h\ar@/_2pc/[dd]_g&{(\tmX,\cM_{\tmX})}\ar@/^2pc/[dd]^\tg\\
{(\coX,\cM_{\coX})}\ar[r]^-(0.5){i}\ar[d]^{\cof}\ar@{}[rd]|{\Box}&{(\tX,\cM_{\tX})}\ar[d]_{\tf}\\
{(\coS,\cM_{\coS})}\ar[r]^-(0.5){\iota}&{(\tS,\cM_{\tS}).}}
\end{equation}
For any rational number $r\geq 0$, 
we denote by $\cL^{(r)}_{\tmX/\tX}$ the $(r)$-twisted torsor of liftings of the canonical morphism $\hmX\rightarrow \coX$ to $\tmX$ over $\tX$ \eqref{p2-hta7}, 
and by $\cC^{(r)}_{\tmX/\tX}$ (resp.\ $\cF^{(r)}_{\tmX/\tX}$) the Higgs--Tate algebra (resp.\ extension) of $\tmX$ over $\tX$ of thickness $r$. 
We set  
\begin{equation}\label{p2-rlps5a}
\fF^{(r)}=\Gamma(\hmX,\cF^{(r)}_{\tmX/\tX}) \ \ \ {\rm and}\ \ \ \fC^{(r)}=\Gamma(\hmX,\cC^{(r)}_{\tmX/\tX}).
\end{equation}
We have a canonical exact sequence of $\hoR$-modules
\begin{equation}\label{p2-rlps5b}
0\rightarrow \hoR\rightarrow \fF^{(r)}\rightarrow \Omega^{(r)}\otimes_{R_\uptau}\hoR \rightarrow 0,
\end{equation}
where $\Omega^{(r)}$ is defined in \eqref{p2-rlps2c}, and a canonical isomorphism of $\hoR$-algebras
\begin{equation}\label{p2-rlps5c}
\fC^{(r)}\stackrel{\sim}{\rightarrow}\underset{\underset{n\geq 0}{\longrightarrow}}\lim\ \rS^n_{\hoR}(\fF^{(r)}). 
\end{equation}
We denote by $\hfC^{(r)}$ the $p$-adic Hausdorff completion of $\fC^{(r)}$. We set $\fC=\fC^{(0)}$, $\fF=\fF^{(0)}$ and $\hfC=\hfC^{(0)}$. 

The algebra $\fC^{(r)}$ (resp.\ extension $\fF^{(r)}$) have been introduced in (\cite{agt} II.10.5 and \cite{ag2} 3.2.16) and were called 
{\em the Higgs--Tate algebra (resp.\ extension) of thickness $r$ associated with the deformation $(\tX,\cM_\tX)$}, which is the origin of the terminology in \ref{p2-hta7}.

The action of $\Delta$ on $(\tmX,\cM_{\tmX})$ induces a canonical $\hoR$-semilinear action of $\Delta$ on $\fF^{(r)}$ 
such that the morphisms of \eqref{p2-rlps5b} are $\Delta$-equivariant (see \cite{ag2} 3.2.15). 
The latter induces an action of $\Delta$ on $\fC^{(r)}$ by ring automorphisms compatible with its action on $\hoR$.  

We denote by
\begin{equation}\label{p2-rlps5d}
d_{\fC^{(r)}}\colon \fC^{(r)}\rightarrow \Omega^{(r)}\otimes_{R_\uptau}\fC^{(r)}
\end{equation}
the universal $\hoR$-derivation of $\fC^{(r)}$ \eqref{p2-hta5g} and by 
\begin{equation}\label{p2-rlps5e}
\delta_{\fC^{(r)}}\colon \fC^{(r)}\rightarrow \Omega\otimes_{R_\uptau}\fC^{(r)}
\end{equation}
the $\hoR$-derivation induced by the canonical injection $\Omega^{(r)}\rightarrow \Omega$ \eqref{p2-hta5h}. 
These derivations are clearly $\Delta$-equivariant, and are Higgs $\hoR$-fields with coefficients in $\Omega^{(r)}$ and $\Omega$, respectively.
We denote by
\begin{equation}\label{p2-rlps5g}
\delta_{\hfC^{(r)}}\colon \hfC^{(r)}\rightarrow \Omega\otimes_{R_\uptau}\hfC^{(r)}
\end{equation}
the extension of $\delta_{\fC^{(r)}}$ to the $p$-adic completions.

\subsection{}\label{p2-rlps6}
For any rational numbers $r\geq r'\geq 0$, we have a canonical injective and $\Delta$-equivariant $\hoR$-homomorphism \eqref{p2-hta5j} 
\begin{equation}\label{p2-rlps6b}
\alpha^{r,r'}\colon \fC^{(r)}\rightarrow \fC^{(r')}.
\end{equation}
The induced homomorphism
$\halpha^{r,r'}\colon\hfC^{(r)}\rightarrow \hfC^{(r')}$ is injective by \ref{p1-thbn22}(ii). We set
\begin{equation}\label{p2-rlps6a}
\hfC^{(r+)}=\underset{\underset{t\in \mQ_{>r}}{\longrightarrow}}{\lim} \hfC^{(t)},
\end{equation}
that we identify with a sub-$\hoR$-algebra of $\hfC$. The
group $\Delta$ acts naturally on $\hfC^{(r+)}$ by  ring automorphisms in a manner
compatible with its actions on $\hoR$ and on $\hfC$.

By \eqref{p2-hta5l}, we have
\begin{equation}\label{p2-rlps6d}
(\id \otimes \alpha^{r,r'}) \circ \delta_{\fC^{(r)}}=\delta_{\fC^{(r')}}\circ  \alpha^{r,r'}.
\end{equation}
The derivations $(\delta_{\hfC^{(t)}})_{t\in \mQ_{>r}}$  induce therefore an $\hoR$-derivation
\begin{equation}\label{p2-rlps6e}
\delta_{\hfC^{(r+)}}\colon \hfC^{(r+)}\rightarrow \Omega \otimes_{R_\uptau}\hfC^{(r+)},
\end{equation}
which is none other than the restriction of $\delta_{\hfC^{(r)}}$ to $\hfC^{(r+)}$.
For simplicity, we set $\fC^\dagger=\hfC^{(0+)}$ and $\delta_{\fC^\dagger}=\delta_{\hfC^{(0+)}}$. Observe that the canonical homomorphism
$\fC^\dagger\rightarrow \hfC$ is injective and that $\delta_{\fC^\dagger}$ is induced by $d_{\hfC}=\delta_{\hfC}$. 
By \ref{p1-thbn22}(iv), we have 
\begin{equation}\label{p2-rlps6f}
\ker(\delta_{\hfC^{(r)}})=\ker(\delta_{\hfC^{(r+)}})=\hoR.
\end{equation}

We set 
\begin{equation}\label{p2-rlps6g}
\delta^\vee_{\hfC^{(r)}}=-\delta_{\hfC^{(r)}}, \ \ \ \delta^\vee_{\hfC^{(r+)}}=-\delta_{\hfC^{(r+)}}\ \ \ {\rm and} \ \ \ \delta^\vee_{\fC^\dagger}=-\delta_{\fC^\dagger}. 
\end{equation}

\subsection{}\label{p2-rlps7}
For any $\hoR$-representation $M$ of $\Delta$ (\cite{agt} II.3.1), we denote by $\mH(M)$ the $\hRun$-module \eqref{p2-rlps2a} defined by
\begin{equation}\label{p2-rlps7a}
\mH(M)=(M\otimes_{\hoR}\fC^\dagger)^\Delta.
\end{equation}
We endow it with the Higgs $\hRun$-field with coefficients in $\Omega$ induced by $\delta_{\fC^\dagger}$ \eqref{p2-rlps6e}.
We thus define  a functor \eqref{p2-rlps2}
\begin{equation}\label{p2-rlps7b}
\mH\colon \bRep_{\hoR}(\Delta) \rightarrow \bHM(\hRun,\Omega).
\end{equation}

\subsection{}\label{p2-rlps8}
For any Higgs $\hRun$--module $(N,\theta)$ with coefficients in $\Omega$ \eqref{p2-rlps2}, we denote by $\mV(N)$ the $\hoR$-module 
defined by
\begin{equation}\label{p2-rlps8a}
\mV(N)=(N\otimes_{\hRun}\fC^\dagger)^{\theta_\tot=0},
\end{equation}
where $\theta_\tot=\theta\otimes \id+\id\otimes \delta_{\fC^\dagger}$ is the total Higgs $\hRun$-field on $ N\otimes_{\hRun}\fC^\dagger$.
We endow it with the $\hoR$-semilinear action of $\Delta$ induced by its natural action on $\fC^\dagger$.
We thus define  a functor
\begin{equation}\label{p2-rlps8b}
\mV\colon \bHM(\hRun,\Omega)\rightarrow \bRep_{\hoR}(\Delta).
\end{equation}

\begin{defi}[\cite{ag2} 3.3.9] \label{p2-rlps9}
We say that an $\hoR[\frac 1 p]$-representation $M$ of $\Delta$ is {\em Dolbeault}
if the following conditions are satisfied:
\begin{itemize}
\item[(i)] $\mH(M)$ is a projective $\hRun[\frac 1 p]$-module of finite type;
\item[(ii)] the canonical morphism
\begin{equation}\label{p2-rlps9a}
\mH(M) \otimes_{\hRun}\fC^\dagger\rightarrow M\otimes_{\hoR}\fC^\dagger
\end{equation}
is an isomorphism.
\end{itemize}
\end{defi}

We denote by $\bRep_{\hoR[\frac 1 p]}^\Dol(\Delta)$ the full subcategory of $\bRep_{\hoR[\frac 1 p]}(\Delta)$ \eqref{p2-rlps2} made up of Dolbeault $\hoR[\frac 1 p]$-representations of $\Delta$.

\begin{defi}[\cite{ag2} 3.3.10]\label{p2-rlps10}
We say that a Higgs $\hRun[\frac 1 p]$-module $(N,\theta)$ with coefficients in $\Omega$
is {\em solvable} if the following conditions are satisfied:
\begin{itemize}
\item[(i)] $N$ is a projective $\hRun[\frac 1 p]$-module of finite type;
\item[(ii)] the canonical morphism
\begin{equation}\label{p2-rlps10a}
\mV(N) \otimes_{\hoR}\fC^\dagger\rightarrow N\otimes_{\hRun}\fC^\dagger
\end{equation}
is an isomorphism.
\end{itemize}
\end{defi}

We denote by $\bHM^\sol(\hRun[\frac 1 p],\Omega)$ the full subcategory of $\bHM(\hRun[\frac 1 p],\Omega)$ \eqref{p2-rlps2} 
made up of solvable Higgs $\hRun[\frac 1 p]$-modules with coefficients in $\Omega$. 

This definition inspired that of weakly twistable Higgs modules in \ref{p1-thbn30}.

\begin{remas}\label{p2-rlps100}
The functors $\mH$ \eqref{p2-rlps7b} and $\mV$ \eqref{p2-rlps8b} do not depend on the choice of the deformation $\tf$ \eqref{p2-rlps1d}, up to {\em non-canonical} isomorphism (\cite{ag2} 3.3.8).
Nevertheless, they depend a priori on the relative or absolute case considered in \ref{p2-ncgt3}.
Therefore, the notions introduced in \ref{p2-rlps9} and \ref{p2-rlps10} do not depend on the choice of the deformation $\tf$ \eqref{p2-rlps1d}, 
but they depend a priori on the relative or absolute case considered in \ref{p2-ncgt3}.
\end{remas}

\begin{prop}[\cite{ag2} 3.3.16] \label{p2-rlps11}
The functors $\mH$ \eqref{p2-rlps7b} and $\mV$ \eqref{p2-rlps8b} induce equivalences of categories quasi-inverse to each other 
between the category of Dolbeault $\hoR[\frac 1 p]$-representations of $\Delta$ and that of
solvable Higgs $\hRun[\frac 1 p]$-modules with coefficients in $\Omega$
\begin{equation}
\xymatrix{
{\bRep^\Dolb_{\hoR[\frac 1 p]}(\Delta)}\ar@<1ex>[r]^-(0.5){\mH}&{\bHM^\sol(\hRun[\frac 1 p], \Omega).}
\ar@<1ex>[l]^-(0.5){\mV}}
\end{equation}
\end{prop}

\begin{defi}\label{p2-rlps13} 
Let $M$ be an $\hoR[\frac 1 p]$-representation of $\Delta$, $(N,\theta)$ a Higgs $\hRun[\frac 1 p]$-bundle 
with coefficients in $\Omega$ \eqref{p1-delta-con6}, $r$ a rational number $>0$.  
We say that {\em $M$ and $(N,\theta)$ are $r$-associated} if there exists a $\Delta$-equivariant $\hfC^{(r)}$-linear isomorphism
of $\hfC^{(r)}$-modules with $\delta_{\hfC^{(r)}}$-connection \eqref{p2-rlps5g}, in the sense of \ref{p1-delta-con2},
\begin{equation}\label{p2-rlps13a} 
N \otimes_{\hRun}\hfC^{(r)}\stackrel{\sim}{\rightarrow} M\otimes_{\hoR}\hfC^{(r)},
\end{equation}
where $N$ (resp.\ $\hfC^{(r)}$) is endowed with the trivial (resp.\ canonical) action of $\Delta$, 
and the $\delta_{\hfC^{(r)}}$-connections are defined as in \ref{p1-delta-con5},
$N$ (resp.\ $M$) being endowed with the Higgs field $\theta$ (resp.\ zero).
\end{defi}

\begin{prop}[\cite{ag2} 3.3.12]\label{p2-rlps14} 
An $\hoR[\frac 1 p]$-representation $M$ of $\Delta$ is Dolbeault if and only if it is $r$-associated with a Higgs $\hRun[\frac 1 p]$-bundle 
$(N,\theta)$ with coefficients in $\Omega$ \eqref{p2-rlps13}, for a rational number $r>0$. 
Moreover, in this case, we have the following properties:
\begin{itemize}
\item[{\rm (i)}] The $\hoR[\frac 1 p]$-module $M$ is projective of finite type, and the action of $\Delta$ on $M$ is continuous for the $p$-adic topology {\rm (\cite{ag2} 2.1.1)}.
\item[{\rm (ii)}] We have an isomorphism of Higgs $\hRun[\frac 1 p]$-modules with coefficients in $\Omega$,
\begin{equation}\label{p2-rlps14a} 
\mH(M)\stackrel{\sim}{\rightarrow} (N,\theta).
\end{equation}
\end{itemize}
\end{prop}

\begin{prop}[\cite{ag2} 3.3.13]\label{p2-rlps15} 
A Higgs $\hRun[\frac 1 p]$-bundle $(N,\theta)$ with coefficients in $\Omega$ is solvable if and only if it is $r$-associated with an
$\hoR[\frac 1 p]$-representation $M$ of $\Delta$ \eqref{p2-rlps13}, for a rational number $r>0$. 
Moreover, in this case, we have an isomorphism of $\hoR[\frac 1 p]$-representation $M$ of $\Delta$,
\begin{equation}\label{p2-rlps15a} 
\mV(N,\theta)\stackrel{\sim}{\rightarrow} M. 
\end{equation}
\end{prop}

\subsection{}\label{p2-rlps16} 
We consider in the remaining part of this section the punctual topos $s$ ringed by the ring $\hoR$, 
and the canonical exact sequence of $\hoR$-modules \eqref{p2-rlps5b}
\begin{equation}\label{p2-rlps16a}
0\rightarrow \hoR\rightarrow \fF\rightarrow \Omega\otimes_{R_\uptau}\hoR \rightarrow 0,
\end{equation}
where $\fF=\fF^{(0)}$ is the $\hoR$-module of affine functions on the torsor $\cL_{\tmX/\tX}=\cL_{\tmX/\tX}^{(0)}$ \eqref{p2-rlps5}. 
We defined in \ref{p1-thbn9} {\em the twisting functor by the extension \eqref{p2-rlps16a}}:
\begin{equation}\label{p2-rlps16b} 
\upnu\colon 
\begin{array}[t]{clcr}
\bHM(\hoR, \Omega) &\rightarrow& \bHM(\hoR, \Omega)\\
(N,\theta)&\mapsto&(\uupnu(N,\theta), \theta_\upnu),
\end{array}
\end{equation}
where $\uupnu(N,\theta)$ denotes the $\hoR$-module underling $\upnu(N,\theta)$. We reserve the symbol $\uptau$ 
for a different twisting functor, which will be introduced later in this work. 
For a Higgs $\hoR[\frac 1 p]$-module with coefficients in $\Omega$, we introduced
the properties of being {\em weakly twistable} (resp.\ {\em twistable}) {\em by the extension \eqref{p2-rlps16a}} in \ref{p1-thbn30} (resp.\ \ref{p1-thbn14}). 
We also defined in \ref{p1-thbn40} a second twisting functor 
\begin{equation}\label{p2-rlps16c}
\upnu^\vee\colon 
\begin{array}[t]{clcr}
\bHM(\hoR,\Omega) &\rightarrow& \bHM(\hoR,\Omega)\\
(N,\theta)&\mapsto&(\uupnu^\vee(N,\theta), \theta_{\upnu^\vee}).
\end{array}
\end{equation}
Observe that the functors $\upnu$ and $\upnu^\vee$ are defined by the period ring $\fC^\dagger$ by \ref{p2-hta71}.

\begin{prop}\label{p2-rlps17} 
Let $(N,\theta)$ be a Higgs $\hRun$-module with coefficients in $\Omega$,  $(\tN,\ttheta)=(N,\theta)\otimes_{\hRun}\hoR$. 
Then, we have a canonical functorial $\hoR$-linear isomorphism 
\begin{equation}\label{p2-rlps17a} 
\mV(N,\theta)\stackrel{\sim}{\rightarrow}\uupnu(\tN,\ttheta). 
\end{equation}
In particular, the $\hoR$-representation $\mV(N,\theta)$ of $\Delta$ is canonically equipped with a $\Delta$-equivariant Higgs $\hoR$-field
with coefficients in $\Omega$, namely the one induced by $\theta\otimes \id$ on $N\otimes_{\hRun}\fC^\dagger$ \eqref{p2-rlps8a}; 
and the Higgs $\hoR$-module $\upnu(\tN,\ttheta)$ is canonically equipped with an $\hoR$-semilinear action of $\Delta$.  
\end{prop}

Indeed, the isomorphism \eqref{p2-rlps17a} is induced by the canonical identifications of the Higgs--Tate algebras in \ref{p2-rlps16}. 
The other assertions are obvious. 

\begin{prop}[\ref{p1-thbn41}]\label{p2-rlps18}
Let $(N,\theta)$ be a Higgs $\hoR[\frac 1 p]$-module with coefficients in $\Omega$, $N^\vee$ an $\hoR$-module, $r$ a rational number $>0$, 
\begin{equation}\label{p2-rlps18a}
N^\vee\otimes_{\hoR}\fC^\dagger\stackrel{\sim}{\rightarrow}N\otimes_{\hoR}\fC^\dagger,
\end{equation}
an isomorphism of $\fC^\dagger$-modules with $\delta_{\fC^\dagger}$-connection \eqref{p2-rlps6e}, in the sense of \eqref{p1-delta-con2},
where the $\delta_{\fC^\dagger}$-connections are defined as in \eqref{p1-delta-con4}, 
$N$ (resp.\ $N^\vee$) being endowed with the Higgs field $\theta$ (resp.\ $0$).  Then, 
\begin{itemize}
\item[{\rm (i)}] The isomorphism \eqref{p2-rlps18a} induces an $\hoR$-linear isomorphism $N^\vee\stackrel{\sim}{\rightarrow} \uupnu(N,\theta)$. 
We deduce a canonical Higgs $\hoR$-field $\theta^\vee$ on $N^\vee$ with coefficients in $\Omega$, so that we have an isomorphism of Higgs $\hoR$-modules
\begin{equation}\label{p2-rlps18b}
(N^\vee,\theta^\vee)\stackrel{\sim}{\rightarrow} \upnu(N,\theta),
\end{equation}
where  $\upnu$ is the functor \eqref{p2-rlps16b}. 
\item[{\rm (ii)}] The morphism \eqref{p2-rlps18a} is an isomorphism of $\fC^\dagger$-modules with $\delta^\vee_{\fC^\dagger}$-connection \eqref{p2-rlps6g}, 
where the $\delta^\vee_{\fC^\dagger}$-connections are defined as in \eqref{p1-delta-con4}, 
$N$ (resp.\ $N^\vee$) being endowed with the Higgs field $0$ (resp.\ $\theta^\vee$). 
\item[{\rm (iii)}] The isomorphism \eqref{p2-rlps18a} induces an isomorphism of Higgs $\hoR$-modules
\begin{equation}\label{p2-rlps18c}
(N,\theta)\stackrel{\sim}{\rightarrow} \upnu^\vee(N^\vee,\theta^\vee),
\end{equation}
where $\upnu^\vee$ is the functor \eqref{p2-rlps16c}. 
\end{itemize}
\end{prop}

\subsection{}\label{p2-rlps19}
Let $M$ be a Dolbeault $\hoR[\frac 1 p]$-representation of $\Delta$ \eqref{p2-rlps9}, $(N,\theta)=\mH(M)$, $(\tN,\ttheta)=(N,\theta)\otimes_{\hRun}\hoR$. 
By \ref{p2-rlps11}, we have a canonical isomorphism $M\stackrel{\sim}{\rightarrow} \mV(N,\theta)$ of $\hoR[\frac 1 p]$-representations of $\Delta$. 
Composing with \eqref{p2-rlps17a}, we deduce a $\Delta$-equivariant $\hoR$-linear isomorphism $M\stackrel{\sim}{\rightarrow} \uupnu(\tN,\ttheta)$. 
We deduce a canonical $\Delta$-equivariant Higgs $\hoR$-field $\theta_M$ on $M$  
with coefficients in $\Omega$, so that we have a $\Delta$-equivariant isomorphism of Higgs $\hoR$-modules
\begin{equation}\label{p2-rlps19a}
(M,\theta_M) \stackrel{\sim}{\rightarrow} \upnu(\tN,\ttheta). 
\end{equation}
We denote by $\bHM(\hoR[\frac 1 p],\Omega,\Delta)$ the category of $\hoR[\frac 1 p]$-representations of $\Delta$, equipped with a $\Delta$-equivariant Higgs $\hoR$-field 
with coefficients in $\Omega$. We thus define a canonical functor, depending on the deformation $\tf$ \eqref{p2-rlps1d}, 
\begin{equation}\label{p2-rlps19b}
\Theta\colon 
\begin{array}[t]{clcr}
\bRep_{\hoR[\frac 1 p]}^\Dol(\Delta)&\rightarrow&\bHM(\hoR[\frac 1 p],\Omega,\Delta),\\
M&\mapsto&(M,\theta_M),
\end{array}
\end{equation}
and a canonical functorial isomorphism
\begin{equation}\label{p2-rlps19c}
\Theta(M) \stackrel{\sim}{\rightarrow} \upnu(\mH(M)\otimes_\hRun\hoR). 
\end{equation}
By \ref{p2-rlps18}(iii), we deduce a canonical functorial $\Delta$-equivariant isomorphism of Higgs $\hoR$-modules
\begin{equation}\label{p2-rlps19d}
\mH(M)\otimes_\hRun\hoR \stackrel{\sim}{\rightarrow} \upnu^\vee(\Theta(M)).
\end{equation}
Taking $\Delta$-invariants, we obtain a canonical isomorphism of Higgs $\hRun$-modules
\begin{equation}\label{p2-rlps19e}
\mH(M)\stackrel{\sim}{\rightarrow} \upnu^\vee(\Theta(M))^\Delta.
\end{equation}

\subsection{}\label{p2-rlps20}
Consider a second Cartesian diagram of $\FLS$ 
\begin{equation}\label{p2-rlps20a}
\xymatrix{
{(\coX,\cM_{\coX})}\ar[r]^-(0.5){i'}\ar[d]_{\cof}\ar@{}[rd]|{\Box}&{(\tX',\cM_{\tX'})}\ar[d]^{\tf'}\\
{(\coS,\cM_{\coS})}\ar[r]^-(0.5){\iota}&{(\tS,\cM_{\tS}),}}
\end{equation}
where $\iota$ is the strict closed immersion defined in \eqref{p2-ncgt3b}, such that $\tf'$ is smooth. 
We associate with it objects similar to those associated with \eqref{p2-rlps1d} that we denote by the same symbols equipped with a $^\prime$ exponent.
Since $\tf$ \eqref{p2-rlps1d} is smooth and $\tX'$ is affine, there exists an isomorphism of $(\tS,\cM_{\tS})$-logarithmic schemes
\begin{equation}\label{p2-rlps20b}
(\tX',\cM_{\tX'})\stackrel{\sim}{\rightarrow} (\tX,\cM_{\tX}),
\end{equation}
extending the identity of $(\coX,\cM_{\coX})$. For every rational number $r\geq 0$, the latter induces an isomorphism of 
the $(r)$-twisted torsors of liftings of the canonical morphism $\hmX\rightarrow \coX$ to $\tmX$ over $\tX'$ and $\tX$, respectively, \eqref{p2-hta7}, 
\begin{equation}\label{p2-rlps20c}
\cL^{(r)}_{\tmX/\tX'} \stackrel{\sim}{\rightarrow} \cL^{(r)}_{\tmX/\tX}. 
\end{equation} 
We deduce by pullback a $\Delta$-equivariant $\hoR$-linear morphism $\fF^{(r)}\rightarrow \fF'^{(r)}$ that fits into a commutative diagram 
\begin{equation}\label{p2-rlps20d}
\xymatrix{
0\ar[r]&\hoR\ar@{=}[d]\ar[r]&{\fF^{(r)}}\ar[r]\ar[d]&{\Omega^{(r)}\otimes_{R_\uptau}\hoR}\ar[r]\ar[d]&0\\
0\ar[r]&\hoR\ar[r]&{\fF'^{(r)}}\ar[r]&{\Omega^{(r)}\otimes_{R_\uptau}\hoR}\ar[r]&0,}
\end{equation}
and hence a $\Delta$-equivariant isomorphism of $\hoR$-algebras
\begin{equation}\label{p2-rlps20e}
\uplambda^{(r)}\colon\fC^{(r)}\rightarrow \fC'^{(r)},
\end{equation}
which is compatible with the $\hoR$-derivations $\delta^{(r)}$ and $\delta'^{(r)}$ \eqref{p2-rlps5e}, and the homomorphisms $\alpha^{r,r'}$ and $\alpha'^{r,r'}$ \eqref{p2-rlps6b}, 
for all rational numbers $r\geq r'\geq 0$.  The latter induces a $\Delta$-equivariant homomorphism of $\hoR$-algebras
\begin{equation}\label{p2-rlps20f}
\uplambda^\dagger\colon\fC^\dagger\stackrel{\sim}{\rightarrow} \fC'^\dagger,
\end{equation}
which is compatible with the $\hoR$-derivations $\delta^\dagger$ and $\delta'^\dagger$ \eqref{p2-rlps6e}. We deduce isomorphisms of functors 
\begin{eqnarray}
\mV\stackrel{\sim}{\rightarrow} \mV', &&\mH \stackrel{\sim}{\rightarrow} \mH', \label{p2-rlps20g1}\\
\upnu\stackrel{\sim}{\rightarrow} \upnu',&& \upnu^\vee \stackrel{\sim}{\rightarrow} \upnu'^\vee.\label{p2-rlps20g2}
\end{eqnarray}

\begin{rema}\label{p2-rlps21}
The notions of Dolbeault $\hoR[\frac 1 p]$-representation of $\Delta$ \eqref{p2-rlps9} and solvable Higgs $\hRun[\frac 1 p]$-module \eqref{p2-rlps10}
do not depend on the choice of the deformation $\tf$ of $\cof$ \eqref{p2-rlps1d} by \eqref{p2-rlps20g1}.
Nevertheless, they depend a priori on the relative or absolute case considered in \ref{p2-ncgt3}.
\end{rema}

\subsection{}\label{p2-rlps22}
We keep the assumptions and notation of \ref{p2-rlps20}.
Let $M$ be a Dolbeault $\hoR[\frac 1 p]$-representation of $\Delta$ \eqref{p2-rlps9}, $(N,\theta)=\mH(M)$, $(N',\theta')=\mH'(M)$. 
The second isomorphism in \eqref{p2-rlps20g1} induces an isomorphism of Higgs $\hRun$-modules
\begin{equation}\label{p2-rlps22a}
\nu\colon (N,\theta)\stackrel{\sim}{\rightarrow} (N',\theta').
\end{equation}
By \eqref{p2-rlps19a}, $\nu$ and the first isomorphism in \eqref{p2-rlps20g2} induce an isomorphism of Higgs $\hoR$-modules
\begin{equation}\label{p2-rlps22b}
\mu \colon (M,\theta_M)\stackrel{\sim}{\rightarrow} (M,\theta'_M).
\end{equation}
The diagram 
\begin{equation}\label{p2-rlps22c}
\xymatrix{
{N\otimes_\hRun\fC^\dagger}\ar[r]\ar[d]_{\nu\otimes \uplambda^\dagger}&{M\otimes_\hoR\fC^\dagger}\ar[d]^{\id\otimes \uplambda^\dagger}\\
{N'\otimes_\hRun\fC'^\dagger}\ar[r]&{M\otimes_\hoR\fC'^\dagger,}}
\end{equation}
where the horizontal arrows are the canonical isomorphisms \eqref{p2-rlps9a}, is commutative. 
It follows that the $\hoR$-linear isomorphism $M\stackrel{\sim}{\rightarrow} M$ underlying $\mu$ \eqref{p2-rlps22b} 
is the identity of $M$. Therefore, $\theta_M=\theta'_M$.  We proved the following:

\begin{prop}\label{p2-rlps23}
For every Dolbeault $\hoR[\frac 1 p]$-representation $M$ of $\Delta$, the canonical Higgs $\hoR$-field $\theta_M$ \eqref{p2-rlps19b} does not depend 
on the choice of the deformation $\tf$ of $\cof$ \eqref{p2-rlps1d}. 
\end{prop}

Observe however that $\theta_M$ depends a priori on the relative or absolute case considered in \ref{p2-ncgt3}.

Similar canonical Higgs fields were introduced by Camargo on pro-étale vector bundles (\cite{rod} 1.0.3), and by He in his Sen theory (\cite{he} 5.28). 
They were subsequently used by Heuer \cite{heuer} and Heuer-Xu \cite{hexu} in their constructions of the $p$-adic Simpson correspondence.

\subsection{}\label{p2-rlps24} 
The Zariski topos of $\cS=\Spf(\co_C)$ is canonically equivalent to the punctual topos $s$. 
Since the $\co_C$-algebra $\hRun$ is topologically of finite presentation, the assumptions of \ref{p1-tshbn1} 
are satisfied by the triple $(\cS_\zar,\hRun,\hoR)$. We can therefore consider the notions introduced in §\ref{p1-tshbn}, 
in particular those of {\em CL-small} (resp.\ {\em locally CL-small}) Higgs $\hRun[\frac 1 p]$-modules with coefficients in $\Omega$ \ref{p1-tshbn13}, 
which are obviously equivalent. 

\begin{rema}\label{p2-rlps25} 
{\em CL-small} Higgs $\hRun[\frac 1 p]$-{\em bundles} with coefficients in $\Omega$ were called {\em small} in (\cite{agt} II.13.5) and (\cite{ag2} 3.4.23). 
Observe that any $\hRun$-module of finite type which is $\co_C$-flat is coherent by (\cite{egr1} 1.10.2(iii)). 
\end{rema}

\begin{prop}\label{p2-rlps26} 
Let $(N,\theta)$ be a Higgs $\hRun[\frac 1 p]$-bundle with coefficients in $\Omega$ \eqref{p1-delta-con6}, 
$(\tN,\ttheta)=(N,\theta)\otimes_{\hRun}\hoR$. Then, the following properties are equivalent:
\begin{itemize}
\item[{\rm (i)}]  The Higgs $\hoR[\frac 1 p]$-module $(\tN,\ttheta)$ is twistable by the extension \eqref{p2-rlps16a}.
\item[{\rm (ii)}]  The Higgs $\hoR[\frac 1 p]$-module $(\tN,\ttheta)$ is weakly twistable by the extension \eqref{p2-rlps16a}.
\item[{\rm (iii)}]  The Higgs $\hRun[\frac 1 p]$-module $(N,\theta)$ is solvable.
\item[{\rm (iv)}]  The Higgs $\hRun[\frac 1 p]$-module $(N,\theta)$ is CL-small. 
\end{itemize}
\end{prop}

Indeed, the equivalence of  (i) and (ii) was proved in \eqref{p1-thbn21}. The equivalence of (ii) and (iii) is an immediate consequence of the definitions. 
The equivalence of (iii) and (iv) was proved in (\cite{ag2} 3.4.30) as a consequence of (\cite{agt} II.13.24 and IV.5.3.10). 

\begin{remas}\label{p2-rlps27}
The implication (iv)$\Rightarrow$(i) of \ref{p2-rlps26} was proved in \ref{p1-tshbn15} in a more general context. 
\end{remas}

\begin{rema}[cf. \cite{ag2} 3.4.32]\label{p2-rlps267}
In the absolute case \eqref{p2-ncgt3}, an $\hoR[\frac 1 p]$-representation of $\Delta$ is Dolbeault \eqref{p2-rlps9}
if and only if it is {\em small} in the sense of (\cite{ag2} 3.4.19). This follows from (\cite{agt} II.14.8) and (\cite{tsuji5} 13.7).
\end{rema}

\subsection{}\label{p2-rlps28}
We set $\rT=\Hom_{R_\uptau}(\Omega,R_\uptau)$, which is a  free $R_\uptau$-module of finite type, and for any rational number $\varepsilon\geq 0$, 
\begin{equation}\label{p2-rlps28a}
\rT^{(\varepsilon)}=\Hom_{R_\uptau}(p^\varepsilon\Omega,R_\uptau). 
\end{equation}
We identify $\rT^{(\varepsilon)}$ (resp.\ $\rT^{(\varepsilon)}\otimes_{R_\uptau}\hRun$) with $p^{-\varepsilon}\rT$ (resp.\ $p^{-\varepsilon}\rT\otimes_{R_\uptau}\hRun$). 
We also set  $\cH=\rS_{\hRun}(\rT\otimes_{R_\uptau}\hRun)$ \eqref{p1-NC7},  $\tcH=\rS_{\hoR}(\rT\otimes_{R_\uptau}\hoR)$ and 
\begin{equation}\label{p2-rlps28b}
\cH^{(\varepsilon)}=\rS_{\hRun}(\rT^{(\varepsilon)}\otimes_{R_\uptau}\hRun).
\end{equation}

Let $A$ be an $\hoR[\frac 1 p]$-algebra quotient of $\tcH[\frac 1 p]$.
We denote by 
\begin{equation}\label{p2-rlps28d}
\theta_A\colon A\rightarrow A \otimes_{R_\uptau} \Omega,
\end{equation}
the {\em canonical} Higgs $\hoR$-field, induced by the homomorphism $\upmu_A\colon \tcH[\frac 1 p]\rightarrow \End_{\hoR}(A)$ 
defined by the multiplication in $A$ \eqref{p1-delta-con1j}. Following \ref{p1-thbn42}, we set 
\begin{equation}\label{p2-rlps28e}
\cL_A=\upnu(A,\theta_A),
\end{equation} 
which is naturally an $A$-module.

Let $B$ be an $\hRun[\frac 1 p]$-algebra quotient of $\cH[\frac 1 p]$. We denote by 
\begin{equation}\label{p2-rlps28c}
\theta_B\colon B\rightarrow B \otimes_{R_\uptau} \Omega
\end{equation}
the {\em canonical} Higgs $\hRun$-field, defined similarly to \eqref{p2-rlps28d}. Observe that $\theta_B$ is $B$-linear by \ref{p1-delta-con7}(i). 
We set $(\tB,\theta_\tB)=(B,\theta_B)\otimes_{\hRun}\hoR$. By \ref{p2-rlps17}, the $\tB$-module $\cL_\tB=\upnu(\tB,\theta_\tB)$ \eqref{p2-rlps28e}
is canonically equipped with a $\tB$-semilinear action of~$\Delta$.  
For any rational number $\varepsilon \geq 0$, we denote by $\cB^{(\varepsilon)}$ the image of $\cH^{(\varepsilon)}$ in $B$ \eqref{p2-rlps28b}. We clearly have
\begin{equation}\label{p2-rlps28f}
\theta_B(\cB^{(\varepsilon)})\subset p^\varepsilon \Omega \otimes_{\hRun} \cB^{(\varepsilon)}. 
\end{equation}

\begin{prop}[cf. \ref{p1-thbn44}]\label{p2-rlps29}
Let $(N,\theta)$ be a Higgs $\hoR[\frac 1 p]$-module with coefficients in $\Omega$, 
$\upmu\colon \tcH[\frac 1 p]\rightarrow \End_\hoR(N)$ 
the homomorphism of $\hoR$-algebras defined by $\theta$ \eqref{p1-delta-con1j}. 
Let $A$ be a quotient $\hoR[\frac 1 p]$-algebra of $\tcH[\frac 1 p]$ through which $\upmu$ factors, 
$\theta_A\colon A\rightarrow A\otimes_{R_\uptau}\Omega$ its canonical Higgs $\hoR$-field \eqref{p2-rlps28d}. 
Suppose that the Higgs $\hoR[\frac 1 p]$-modules $(N,\theta)$ and $(A,\theta_A)$ are weakly twistable. 
Then, the $A$-module $\cL_A=\upnu(A,\theta_A)$ is invertible \eqref{p2-rlps28e}, 
and there exists a canonical functorial $A$-linear isomorphism
\begin{equation}\label{p2-rlps29a}
N\otimes_A\cL_A\stackrel{\sim}{\rightarrow} \upnu(N,\theta). 
\end{equation}
\end{prop}

\begin{prop}[cf. \ref{p1-tshbn27}]\label{p2-rlps30}
Let $N$ be a coherent $\hRun[\frac 1 p]$-module, $\theta$ a Higgs $\hRun$-field on $N$ with coefficients in $\Omega$, 
$\upmu\colon \cH[\frac 1 p]\rightarrow \cEnd_\hRun(N)$ the homomorphism of $\hRun[\frac 1 p]$-algebras defined by $\theta$ \eqref{p1-delta-con1j}. 
We denote by $B$ the image of $\upmu$ and by $\theta_B\colon B\rightarrow B\otimes_\hRun\Omega$ 
its canonical Higgs $\hRun$-field \eqref{p2-rlps28c}. Then, the following conditions are equivalent:
\begin{itemize}
\item[{\rm (i)}] $(N,\theta)$ is CL-small in the sense of \ref{p1-tshbn13}, see \ref{p2-rlps24};
\item[{\rm (ii)}] there exists a rational number $\varepsilon >\frac{1}{p-1}$ such that 
the image $\cB^{(\varepsilon)}$ of $\cH^{(\varepsilon)}$ in $B$ \eqref{p2-rlps28b}, 
is a coherent $\hRun$-module;
\item[{\rm (iii)}] the Higgs $\hRun[\frac 1 p]$-module $(B,\theta_B)$ is CL-small. 
\end{itemize}
\end{prop}

\section{Functoriality of Higgs--Tate algebras in the local case}\label{p2-funchta}

\subsection{}\label{p2-fhtal1}
In this section, we let $f\colon (X,\cM_{X})\rightarrow (S,\cM_S)$ (resp.\ $f'\colon (X',\cM_{X'})\rightarrow (S,\cM_S)$) be an
adequate morphism of fine logarithmic schemes (\cite{agt} III.4.7) having an adequate chart $\varsigma\colon \mN\rightarrow P$
(resp.\ $\varsigma'\colon \mN\rightarrow P'$) (\cite{agt} III.4.4),
such that $X=\Spec(R)$ and $X'=\Spec(R')$ are affine, and $X_s$ and $X'_s$ are non-empty. 
We consider, moreover, a morphism of $(S,\cM_S)$-logarithmic schemes 
\begin{equation}\label{p2-fhtal1a}
\gamma\colon (X',\cM_{X'})\rightarrow (X,\cM_X),  
\end{equation} 
having a chart $\mu\colon P\rightarrow P'$ such that $\varsigma'=\mu\circ \varsigma$. 

We denote by $X^\circ$ (resp.\ $X'^\rhd$) the maximal open subscheme of $X$ (resp.\ $X'$)
where the logarithmic structure $\cM_{X}$ (resp.\ $\cM_{X'}$) is trivial.
For any $X$-scheme $U$ and $X'$-scheme $U'$, we set
\begin{eqnarray}
U^\circ=U\times_{X}X^\circ,\label{p2-fhtal1b}\\
U'^\rhd=U'\times_{X'}X'^\rhd.\label{p2-fhtal1c}
\end{eqnarray}
We have $X^\circ\subset X_\eta$ and $X'^\rhd\subset X'^\circ$. To lighten the notation, we set 
\begin{eqnarray}
\tOmega^1_{X/S}=\Omega^1_{(X,\cM_{X})/(S,\cM_S)}, && \tOmega^1_{R/\co_K}=\tOmega^1_{X/S}(X), \label{p2-fhtal1d}\\
\tOmega^1_{X'/S}=\Omega^1_{(X',\cM_{X'})/(S,\cM_S)}, &&  \tOmega^1_{R'/\co_K}=\tOmega^1_{X'/S}(X').  \label{p2-fhtale} 
\end{eqnarray}

We endow $\coX=X\times_S\coS$ \eqref{p2-ncgt1a} (resp.\ $\coX'=X'\times_S\coS$) with the logarithmic structure $\cM_{\coX}$ (resp.\ $\cM_{\coX'}$) pullback of $\cM_{X}$ (resp.\ $\cM_{X'}$), 
and denote by $\cof\colon (\coX,\cM_{\coX})\rightarrow (\coS,\cM_{\coS})$ (resp.\ $\cof'\colon (\coX',\cM_{\coX'})\rightarrow (\coS,\cM_{\coS})$) the base change of $f$ (resp.\ $f'$). 
We fix Cartesian diagrams of $\FLS$ \eqref{p1-NC1}
\begin{equation}\label{p2-fhtal1h}
\xymatrix{
{(\coX,\cM_{\coX})}\ar[r]^-(0.5){i}\ar[d]_{\cof}\ar@{}[rd]|{\Box}&{(\tX,\cM_{\tX})}\ar[d]^{\tf}\\
{(\coS,\cM_{\coS})}\ar[r]^-(0.5){\iota}&{(\tS,\cM_{\tS}),}}
\ \ \ 
\xymatrix{
{(\coX',\cM_{\coX'})}\ar[r]^-(0.5){i'}\ar[d]_{\cof'}\ar@{}[rd]|{\Box}&{(\tX',\cM_{\tX'})}\ar[d]^{\tf'}\\
{(\coS,\cM_{\coS})}\ar[r]^-(0.5){\iota}&{(\tS,\cM_{\tS}),}}
\end{equation}
where $\iota$ is the strict closed immersion defined in \eqref{p2-ncgt3b}, such that $\tf$ and $\tf'$ are smooth.

\subsection{}\label{p2-fhtal200}
Let $\oy'$ be a geometric point of $\oX'$, $\oy=\gamma(\oy')$. We take again the notation introduced in \ref{p2-rlps2} for $(f,\oy)$. 
Similarly, the scheme $\oX'$ being locally irreducible, it is the sum of the schemes induced on its irreducible components. 
We denote by $\oX'^\star$ the irreducible component of $\oX'$ containing $\oy'$. Likewise, $\oX'^\rhd$ is the sum of the schemes induced on its irreducible components
and $\oX'^{\star \rhd}=\oX'^\star\times_{X'}X'^\rhd$ is the irreducible component of $\oX'^\rhd$ containing $\oy'$. 
We set $\coX'^{\star}=\oX'^\star\times_\oS\coS$ \eqref{p2-ncgt1}, 
\begin{eqnarray}
R'_1&=&\Gamma(\oX'^\star,\co_{\oX'}),\label{p2-fhtal200a}\\
R'_\uptau&=&\Gamma(\coX'^\star,\co_{\coX'})=R'_1\otimes_{\co_\oK}\co_C,\label{p2-fhtal2b}
\end{eqnarray}
and denote by $\hRunp$ and $\hRtaup$ their $p$-adic Hausdorff completions, which are equal. 
To lighten notation, for any rational number $r\geq 0$, we set, with the conventions of \ref{p2-ncgt3}, 
\begin{equation}\label{p2-fhtal2c}
\Omega'=\txi^{-1}\tOmega^1_{R'/\co_K}\otimes_{R'}R'_\tau \ \ \ {\rm and} \ \ \ \Omega'^{(r)}=p^r\Omega'.
\end{equation}
The canonical morphism $\tOmega^1_{R/\co_K}\otimes_RR' \rightarrow \tOmega^1_{R'/\co_K}$ induces morphisms of free $R'_\uptau$-modules 
of finite type
\begin{eqnarray}
u\colon \Omega\otimes_{R_\uptau}R'_\uptau &\rightarrow& \Omega', \label{p2-fhtal200e}\\
u^{(r)}\colon \Omega^{(r)}\otimes_{R_\uptau}R'_\uptau&\rightarrow& \Omega'^{(r)}. \label{p2-fhtal200ee}
\end{eqnarray}

For any $R'_\uptau$-algebra $A$, we  consider Higgs $A$-modules with coefficients in 
$\Omega'\otimes_{R'_\uptau}A$. We say abusively  that they have coefficients in $\Omega'$.
The category of these modules will be denoted by $\bHM(A,\Omega')$.

We denote by $\Delta'$ the profinite group $\pi_1(\oX'^{\star \rhd},\oy')$ and by $(W_j)_{j\in J}$ the normalized universal cover of
$\oX'^{\star \rhd}$ at $\oy'$ (\cite{ag1} 2.1.20). For any $j\in J$, we denote by $\oX'^{W_j}$ the integral closure of $\oX'$ in $W_j$.
The schemes $(\oX'^{W_j})_{j\in J}$ then form a filtered inverse system. We set
\begin{equation}\label{p2-fhtal2g}
\oR'=\underset{\underset{j\in J}{\longrightarrow}}{\lim}\ \Gamma(\oX'^{W_j},\co_{\oX'^{W_j}}),
\end{equation}
which is naturally endowed with a discrete action of $\Delta'$ by ring automorphisms. We denote by $\hoRp$ its $p$-adic Hausdorff completion
and by $\bRep_{\hoRp}(\Delta')$ the category of $\hoRp$-representations of $\Delta'$. 

With the notation of \ref{p2-rlps2}, for any $i\in I$, we have a canonical $\oX^\circ$-morphism $\oy\rightarrow V_i$.
We deduce an $\oX'^\rhd$-morphism $\oy'\rightarrow V_i\times_{\oX^\circ}\oX'^\rhd$.
The scheme $V_i\times_{\oX^\circ}\oX'^\rhd$ being locally irreducible,
it is the sum of the schemes induced on its irreducible components.
We denote by $V'_i$ the irreducible component of $V_i\times_{\oX^\circ}\oX'^\rhd$ containing the image of $\oy'$
and by $\oX'^{V'_i}$ the integral closure of $\oX'$ in $V'_i$.
The schemes $(V'_i)_{i\in I}$ naturally form a projective system of $\oy'$-pointed connected finite étale coverings of $\oX'^{\star \rhd}$.
For any $i\in I$, we denote by $\Pi_i$ the open subgroup of $\Delta'$ corresponding to $V'_i$, in other words the kernel of the canonical action
of $\Delta'$ on the fiber of $V'_i$ above $\oy'$ (\cite{ag1} (2.1.20.2)).
We denote by $\Pi$ the closed subgroup of $\Delta'$ defined by
\begin{equation}\label{p2-fhtal2h}
\Pi=\cap_{i\in I}\Pi_i.
\end{equation}
Setting $\Delta^\intern=\Delta'/\Pi$, we have a canonical homomorphism $\Delta'\rightarrow \Delta$ which factors through an injective homomorphism
$\Delta^\intern\rightarrow \Delta$.

Consider the ring
\begin{equation}\label{p2-fhtal2i}
\oR^\intern=\underset{\underset{i\in I}{\longrightarrow}}{\lim}\ \Gamma(\oX'^{V'_i},\co_{\oX'^{V' _i}}),
\end{equation}
which is naturally endowed with a discrete action of $\Delta^\intern$ by ring automorphisms.
We have a canonical $\Delta^\intern$-equivariant homomorphism $\oR\rightarrow \oR^\intern$.
 
For all $i\in I$ and $j\in J$, there exists at most one morphism of $\oX'^{\star \rhd}$-pointed schemes $W_j\rightarrow V'_i$.
Moreover, for every $i\in I$, there exists $j\in J$ and a morphism of $\oX'^{\star \rhd}$-pointed schemes $W_j\rightarrow V'_i$.
So we have a canonical homomorphism
\begin{equation}\label{p2-fhtal2j}
\underset{\underset{i\in I}{\longrightarrow}}{\lim}\ \Gamma(\oX'^{V'_i},\co_{\oX'^{V'_i}})\rightarrow
\underset{\underset{j\in J}{\longrightarrow}}{\lim}\ \Gamma(\oX'^{W_j},\co_{\oX'^{W_j}}).
\end{equation}
We deduce a canonical $\Delta'$-equivariant homomorphism of $R'_1$-algebras
\begin{equation}\label{p2-fhtal2k}
\oR^\intern\rightarrow \oR'.
\end{equation}

\subsection{}\label{p2-fhtal3}
We take again the notation introduced in \ref{p2-rlps3}, \ref{p2-rlps5} and \ref{p2-rlps6} for the quadruple $(X,P,\oy,\tf)$. 
We associate with the quadruple $(X',P',\oy',\tf')$ analogous objects, 
that we equip with a $^\prime$ exponent: $(\hmX',\cM_{\hmX'})$, $(\tmX',\cM_{\tmX'})$, 
$\fF'^{(r)}=\Gamma(\hmX',\cF^{(r)}_{\tmX'/\tX'})$ \eqref{p2-rlps5a},  $\fC'^{(r)}=\Gamma(\hmX',\cC^{(r)}_{\tmX'/\tX'})$,
$\alpha'^{r,r'}\colon \fC'^{(r)}\rightarrow \fC'^{(r')}$...
By functoriality of the closed immersion $\mi$ \eqref{p2-rlps3f} (\cite{ag2} 5.2.3), we have a commutative diagram  
\begin{equation}\label{p2-fhtal3a}
\xymatrix{
&{(\hmX',\cM_{\hmX'})}\ar[ld]_{h'}\ar[rr]^{\mi'}\ar@{->}'[d]^{\hmj}[dd]&&{(\tmX',\cM_{\tmX'})}\ar[dd]^{\tmj}\\
{(\coX',\cM_{\coX'})}\ar[rr]^-(0.7){i'}\ar[dd]_{\cogamma}&&{(\tX',\cM_{\tX'})}\ar@/^2pc/[ddd]^-(0.45){\tf'}&\\
&{(\hmX,\cM_\hmX)}\ar[ld]_h\ar[rr]|-(0.7)\hole^{\mi}\ar[ddl]|\hole&&{(\tmX,\cM_{\tmX})}\ar[ldd]\\
{(\coX,\cM_\coX)}\ar[rr]^i\ar[d]_{\cof}&&{(\tX,\cM_\tX)}\ar[d]_{\tf}\\
{(\coS,\cM_{\coS})}\ar[rr]^\iota&&{(\tS,\cM_{\tS})}}
\end{equation}
where $\cogamma$ is the base change of $\gamma$ \eqref{p2-fhtal1a}, and $\hmj$, $\tmj$, $h$ and $h'$ are the canonical morphisms. 
The fundamental group $\Delta'$ \eqref{p2-fhtal200} acts on the left on the logarithmic scheme $(\tmX',\cM_{\tmX'})$ (see \cite{ag2} 3.2.15), 
and on the logarithmic scheme $(\tmX,\cM_{\tmX})$ via the canonical action of $\Delta$ \eqref{p2-rlps3g}. Moreover, the morphism $\tmj$ is $\Delta'$-equivariant.

\subsection{}\label{p2-fhtal4}
By \ref{p2-hta2} applied to $\tf'$, the conditions \ref{p2-hta1}(i)-(iv) are satisfied by $j=i'$ \eqref{p2-fhtal3a}. 
For any rational number $r\geq 0$, we denote by 
$\cL^{(r)}_{\tX'/\tX}$ the $(r)$-twisted torsor of liftings of $\cogamma$ to $\tX'$ over $\tX$,
and by $\cF^{(r)}_{\tX'/\tX}$ (resp.\ $\cC^{(r)}_{\tX'/\tX}$) the Higgs--Tate extension (resp.\ algebra) of $\tX'$ over $\tX$ of thickness $r$ \eqref{p2-hta7}. 
With the notation of \ref{p2-fhtal200}, we set 
\begin{equation}\label{p2-fhtal4a}
\fF^{(r)}_\uptau=\Gamma(\coX'^\star,\cF^{(r)}_{\tX'/\tX}) \ \ \ {\rm and}\ \ \ \fC^{(r)}_\uptau=\Gamma(\coX'^\star,\cC^{(r)}_{\tX'/\tX}).
\end{equation} 
We have a canonical exact sequence of $R'_\uptau$-modules
\begin{equation}\label{p2-fhtal4b}
0\rightarrow R'_\uptau\rightarrow \fF^{(r)}_\uptau\rightarrow \Omega^{(r)}\otimes_{R_\uptau}R'_\uptau \rightarrow 0
\end{equation}
and a canonical isomorphism of $R'_\uptau$-algebras
\begin{equation}\label{p2-fhtal4c}
\fC^{(r)}_\uptau\stackrel{\sim}{\rightarrow}\underset{\underset{n\geq 0}{\longrightarrow}}\lim\ \rS^n_{R'_\uptau}(\fF^{(r)}_\uptau). 
\end{equation}
We denote by $\hfC^{(r)}_\uptau$ the $p$-adic Hausdorff completion of $\fC^{(r)}_\uptau$. 
We set $\fC_\uptau=\fC^{(0)}_\uptau$, $\fF_\uptau=\fF^{(0)}_\uptau$ and $\hfC_\uptau=\hfC^{(0)}_\uptau$.

We denote by
\begin{equation}\label{p2-fhtal4d}
d_{\fC^{(r)}_\uptau}\colon \fC^{(r)}_\uptau\rightarrow \Omega^{(r)}\otimes_{R_\uptau}\fC^{(r)}_\uptau
\end{equation}
the universal $R'_\uptau$-derivation of $\fC^{(r)}_\uptau$ \eqref{p2-hta18c} and by 
\begin{equation}\label{p2-fhtal4dd}
\delta_{\fC^{(r)}_\uptau}\colon \fC^{(r)}_\uptau\rightarrow \Omega\otimes_{R_\uptau}\fC^{(r)}_\uptau
\end{equation}
the $R'_\uptau$-derivation induced by the canonical inclusion $\Omega^{(r)}\rightarrow \Omega$ \eqref{p2-hta18f}. We set 
\begin{equation}\label{p2-fhtal4e}
\delta'_{\fC^{(r)}_\uptau}=(u\otimes \id)\circ \delta_{\fC^{(r)}_\uptau}\colon \fC^{(r)}_\uptau\rightarrow \Omega'\otimes_{R'_\uptau}\fC^{(r)}_\uptau, 
\end{equation}
where $u$ is defined in \eqref{p2-fhtal200e}. We denote by
\begin{eqnarray}
\delta_{\hfC^{(r)}_\uptau}\colon \hfC^{(r)}_\uptau\rightarrow \Omega\otimes_{R_\uptau}\hfC^{(r)}_\uptau, \label{p2-fhtal4g}\\
\delta'_{\hfC^{(r)}_\uptau}\colon \hfC^{(r)}_\uptau\rightarrow \Omega'\otimes_{R'_\uptau}\hfC^{(r)}_\uptau, \label{p2-fhtal4gg}
\end{eqnarray}
the extension of $\delta_{\fC^{(r)}_\uptau}$ and $\delta'_{\fC^{(r)}_\uptau}$ to the $p$-adic completions. 

For all rational numbers $r\geq r'\geq 0$, we have a canonical $R'_\uptau$-homomorphism \eqref{p2-hta5j}
\begin{equation}\label{p2-fhtal4k}
\alpha^{r,r'}_\uptau\colon \fC^{(r)}_\uptau\rightarrow \fC^{(r')}_\uptau.
\end{equation}
The induced homomorphism $\halpha^{r,r'}_\uptau\colon\hfC^{(r)}_\uptau\rightarrow \hfC^{(r')}_\uptau$ is injective by \ref{p1-thbn22}(ii).
We set
\begin{equation}\label{p2-fhtal4m}
\hfC^{(r+)}_\uptau=\underset{\underset{t\in \mQ_{>r}}{\longrightarrow}}{\lim} \hfC^{(t)}_\uptau,
\end{equation}
that we identify with a sub-$\hRunp$-algebra of $\hfC_\uptau$. 

We have \eqref{p2-hta5l}
\begin{equation}\label{p2-fhtal4f}
(\id \otimes \alpha^{r,r'}_\uptau) \circ \delta_{\fC^{(r)}_\uptau}=\delta_{\fC^{(r')}_\uptau}\circ  \alpha^{r,r'}_\uptau.
\end{equation}
Therefore, the derivations $(\delta_{\hfC^{(t)}_\uptau})_{t\in \mQ_{>r}}$ induce an $\hRunp$-derivation
\begin{equation}\label{p2-fhtal4n}
\delta_{\hfC^{(r+)}_\uptau}\colon \hfC^{(r+)}_\uptau\rightarrow \Omega \otimes_{R_\uptau}\hfC^{(r+)}_\uptau,
\end{equation}
which is none other than the restriction of $\delta_{\hfC^{(r)}_\uptau}$ to $\hfC^{(r+)}_\uptau$.
For simplicity, we set $\fC^\dagger_\uptau=\hfC^{(0+)}_\uptau$ and $\delta_{\fC^\dagger_\uptau}=\delta_{\hfC^{(0+)}_\uptau}$. Observe that the canonical homomorphism
$\fC^\dagger_\uptau\rightarrow \hfC_\uptau$ is injective and that $\delta_{\fC^\dagger_\uptau}$ is induced by $d_{\hfC_\uptau}=\delta_{\hfC_\uptau}$.

We denote by $\mK^\bullet(\hfC^{(r)}_\uptau)$
the Dolbeault complex of $(\hfC^{(r)}_\uptau,\delta_{\hfC^{(r)}_\uptau})$ \ref{p1-delta-con1b} and by $\tmK^\bullet(\hfC^{(r)}_\uptau)$
the augmented Dolbeault complex
\begin{equation}\label{p2-fhtal4h}
\hRunp\rightarrow \mK^0(\hfC^{(r)}_\uptau)\rightarrow \mK^1(\hfC^{(r)}_\uptau)\rightarrow \dots,
\end{equation}
where $\hRunp$ is placed in degree $-1$ and the differential $\hRunp\rightarrow\hfC^{(r)}_\uptau$ is the canonical homomorphism. 
By \eqref{p2-fhtal4f}, for any rational numbers $r\geq r'\geq 0$, the homomorphism $\halpha^{r,r'}_\uptau$ induces morphisms of complexes
\begin{eqnarray}
\upiota^{r,r'}_\uptau\colon \mK^\bullet(\hfC^{(r)}_\uptau)&\rightarrow& \mK^\bullet(\hfC^{(r')}_\uptau),\label{p2-fhtal4j}\\
\tupiota^{r,r'}_\uptau\colon \tmK^\bullet(\hfC^{(r)}_\uptau)&\rightarrow& \tmK^\bullet(\hfC^{(r')}_\uptau).\label{p2-fhtal4i}
\end{eqnarray}

\begin{prop}[\ref{p1-thbn36}]\label{p2-fhtal5}
For all rational numbers $r>r'>0$, there exists a rational number $\alpha\geq 0$ depending on $r$ and $r'$ but not on the morphisms $f$ and $\gamma$ 
satisfying the conditions of \ref{p2-fhtal1}, such that
\begin{equation}\label{p2-fhtal5a}
p^\alpha\tupiota^{r,r'}_\uptau\colon \tmK^\bullet(\hfC^{(r)}_\uptau)\rightarrow \tmK^\bullet(\hfC^{(r')}_\uptau),
\end{equation}
where $\tupiota^{r,r'}_\uptau$ is the morphism \eqref{p2-fhtal4i}, is homotopic to $0$ by an $\hRunp$-linear homotopy.
\end{prop}

\begin{cor}[\ref{p1-thbn37}]\label{p2-fhtal6}
For every rational number $r\geq0$, and every $\hRunp[\frac 1 p]$-module $M$, 
the complex $M\otimes_{\hRunp}\mK^\bullet(\hfC^{(r+)}_\uptau)$ is a resolution of $M$, 
where the tensor product is defined term by term (non derived).
\end{cor}

Observe that the constructions and notation of this article are compatible with those of §\ref{p1-thbn} by \ref{p2-hta71}.

\subsection{}\label{p2-fhtal7}
By \eqref{p1-rdt4i}, for every rational number $r\geq 0$, we have a canonical isomorphism of $\co_{\hmX'}$-extensions 
\begin{equation}\label{p2-fhtal7a}
\cF^{(r)}_{\tmX'/\tX}\stackrel{\sim}{\rightarrow}\cF^{(r)}_{\tmX/\tX}\otimes_{\co_{\hmX}}\co_{\hmX'};
\end{equation}
see \eqref{p2-fhtal3a}. Setting
\begin{equation}\label{p2-fhtal7d}
\fF^{(r)}_{\hoRp}=\fF^{(r)}\otimes_{\hoR}\hoRp \ \ \ {\rm and} \ \ \ \fC^{(r)}_{\hoRp}=\fC^{(r)}\otimes_{\hoR}\hoRp,
\end{equation}
we obtain a canonical isomorphism of $\hoRp$-extensions 
\begin{equation}\label{p2-fhtal7b}
\Gamma(\hmX',\cF^{(r)}_{\tmX'/\tX})\stackrel{\sim}{\rightarrow}\fF^{(r)}_{\hoRp}, 
\end{equation}
and a canonical isomorphism of $\hoRp$-algebras 
\begin{equation}\label{p2-fhtal7c}
\Gamma(\hmX',\cC^{(r)}_{\tmX'/\tX})\stackrel{\sim}{\rightarrow}\fC^{(r)}_{\hoRp}.
\end{equation}
We will identify the sources and the targets of the above isomorphisms. 

We denote by
\begin{equation}\label{p2-fhtal7e}
d_{\fC^{(r)}_{\hoRp}}\colon \fC^{(r)}_{\hoRp}\rightarrow \Omega^{(r)}\otimes_{R_\uptau}\fC^{(r)}_{\hoRp}
\end{equation}
the universal $\hoRp$-derivation of $\fC^{(r)}_{\hoRp}$ \eqref{p2-hta18a}, 
which is none other than the base change of $d_{\fC^{(r)}}$ \eqref{p2-rlps5d} by the homomorphism $\hoR\rightarrow \hoRp$, and by 
\begin{equation}\label{p2-fhtal7g}
\delta_{\fC^{(r)}_{\hoRp}}\colon \fC^{(r)}_{\hoRp}\rightarrow \Omega\otimes_{R_\uptau}\fC^{(r)}_{\hoRp}
\end{equation}
the $\hoRp$-derivation induced by $d_{\fC^{(r)}_{\hoRp}}$ and the canonical injection $\Omega^{(r)}\rightarrow \Omega$ \eqref{p2-hta18d}. We set 
\begin{equation}\label{p2-fhtal7f}
\delta'_{\fC^{(r)}_{\hoRp}}=(u\otimes \id)\circ \delta_{\fC^{(r)}_{\hoRp}} \colon \fC^{(r)}_{\hoRp}\rightarrow \Omega'\otimes_{R_\uptau}\fC^{(r)}_{\hoRp},
\end{equation}
where $u$ is defined in \eqref{p2-fhtal200e}. 

For any rational numbers $r\geq r'\geq 0$, we denote by 
\begin{equation}\label{p2-fhtal7h}
\alpha^{r,r'}_{\hoRp}\colon \fC^{(r)}_{\hoRp}\rightarrow \fC^{(r')}_{\hoRp} 
\end{equation}
the canonical homomorphism \eqref{p2-hta5j}, 
which coincides with the base change of $\alpha^{r,r'}$ \eqref{p2-rlps6b}.

\subsection{}\label{p2-fhtal110}
For any rational numbers $r\geq r'\geq 0$, we have a canonical morphism of $\hoRp$-modules \eqref{p2-hta150h},
associated with the diagram \eqref{p2-fhtal3a},
\begin{equation}\label{p2-fhtal110a}
\upvarphi^{(r,r')}\colon \fF^{(r)}_{\hoRp}\rightarrow \fF'^{(r)} \otimes_{R'_\uptau}\fF^{(r')}_\uptau,
\end{equation}
and a canonical morphism of $\hoRp$-algebras \eqref{p2-hta150i}
\begin{equation}\label{p2-fhtal110b}
\upphi^{(r,r')}\colon \fC^{(r)}_{\hoRp}\rightarrow \fC'^{(r)} \otimes_{R'_\uptau}\fC^{(r')}_\uptau,
\end{equation}
compatible with $\upvarphi^{(r,r')}$. 
The latter induces by linearization a morphism of $(\hoRp\otimes_{R'_\uptau}\fC^{(r')}_\uptau)$-algebras \eqref{p2-hta150e}
\begin{equation}\label{p2-fhtal110c}
\uppsi^{(r,r')}\colon \fC^{(r)}_{\hoRp}\otimes_{R'_\uptau}\fC^{(r')}_\uptau\rightarrow \fC'^{(r)}\otimes_{R'_\uptau}\fC^{(r')}_\uptau.
\end{equation}

Let $t\geq t'\geq 0$ be rational numbers such that $r\geq t$ and $r'\geq t'$. By \eqref{p2-hta150k}, the diagram 
\begin{equation}\label{p2-fhtal110f}
\xymatrix{
{\fC^{(r)}_{\hoRp}}\ar[r]^-(0.5){\upphi^{(r,r')}}\ar[d]_{\alpha^{r,t}_\hoRp}
&{\fC'^{(r)} \otimes_{R'_\uptau} \fC^{(r')}_\uptau}
\ar[d]^{\alpha'^{r,t}\otimes_{R'_\uptau} \alpha_\uptau^{r',t'}}\\
{\fC^{(t)}_{\hoRp}}\ar[r]^-(0.5){\upphi^{{(t,t')}}}&{\fC'^{(t)} \otimes_{R'_\uptau} \fC^{(t')}_{\uptau},}}
\end{equation}
where the vertical arrows are defined in \eqref{p2-rlps6b}, \eqref{p2-fhtal7h} and \eqref{p2-fhtal4k}, is commutative. 
We have a similar commutative diagram for $\uppsi^{(r,r')}$. 

\subsection{}\label{p2-fhtal12}
The natural action of $\Delta'$ on $(\tmX',\cM_{\tmX'})$ induces canonical actions on the torsors $\cL_{\tmX'/\tX'}$ and $\cL_{\tmX'/\tX}$ (see \cite{agt} II.4.18). 
The morphism $\tmj$ being $\Delta'$-equivariant \eqref{p2-fhtal3a}, we easily check that the canonical isomorphism \eqref{p1-rdt4h}
\begin{equation}\label{p2-fhtal12a}
\hmj^+(\cL_{\tmX/\tX})\stackrel{\sim}{\rightarrow} \cL_{\tmX'/\tX}
\end{equation}
is $\Delta'$-equivariant, when we equip $\hmj^+(\cL_{\tmX/\tX})$ with the action induced by the canonical action of $\Delta$ on $\cL_{\tmX/\tX}$ and the canonical action of $\Delta'$ 
on the scheme $\hmX'$ \eqref{p2-fhtal3a}. 
The group $\Delta'$ acts on the torsor $h'^+(\cL_{\tX'/\tX})$ by acting trivially on the torsor $\cL_{\tX'/\tX}$ and canonically on the scheme $\hmX'$. 
We easily check that the canonical morphism \eqref{p2-hta10e}, associated with the diagram \eqref{p2-fhtal3a},
\begin{equation}\label{p2-fhtal12b}
\uppsi_\cL\colon \cL_{\tmX'/\tX'}\times h'^+(\cL_{\tX'/\tX}) \rightarrow \cL_{\tmX'/\tX} \times h'^+(\cL_{\tX'/\tX})
\end{equation}
is $\Delta'$-equivariant. 

Let $r,r'$ be rational numbers such that $r\geq r'\geq 0$. 
The action of $\Delta'$ on $\cL_{\tmX'/\tX'}$ induces an action on $\cL^{(r)}_{\tmX'/\tX'}$. 
We deduce an $\hoRp$-semilinear action of $\Delta'$ on $\fF'^{(r)}$ and an action of $\Delta'$ on $\fC'^{(r)}$ 
by ring automorphisms compatible with its action on $\hoRp$ (see \cite{ag2} 3.2.15).  

The action of $\Delta'$ on $\cL_{\tmX'/\tX}$ (resp.\ $h'^+(\cL_{\tX'/\tX})$) induces an action on $\cL^{(r)}_{\tmX'/\tX}$ (resp.\ $h'^+(\cL^{(r)}_{\tX'/\tX})$). 
By \ref{p2-fhtal7}, we deduce an $\hoRp$-semilinear action of $\Delta'$ on $\fF^{(r)}_{\hoRp}$ and 
an action of $\Delta'$ on $\fC^{(r)}_\hoRp$ by ring automorphisms compatible with its action on $\hoRp$. 
Since \eqref{p2-fhtal12a} is $\Delta'$-equivariant, these actions are induced by the canonical actions of $\Delta$ on $\fF^{(r)}$ and $\fC^{(r)}$. 

We deduce from \eqref{p2-fhtal12b} that the canonical morphism \eqref{p2-hta150d}, associated with the diagram \eqref{p2-fhtal3a},
\begin{equation}
\uppsi^{(r,r')}_\cL\colon \cL^{(r)}_{\tmX'/\tX'}\times h'^+(\cL^{(r')}_{\tX'/\tX}) \rightarrow \cL^{(r)}_{\tmX'/\tX} \times h'^+(\cL^{(r')}_{\tX'/\tX})
\end{equation}
is $\Delta'$-equivariant. Therefore, the morphism $\uppsi^{(r,r')}$ \eqref{p2-fhtal110c} is $\Delta'$-equivariant, 
where $\Delta'$ acts trivially on $\fC^{(r')}_\uptau$. 
It follows that the morphisms $\upvarphi^{(r,r')}$ \eqref{p2-fhtal110a} and $\upphi^{(r,r')}$ \eqref{p2-fhtal110b} are $\Delta'$-equivariant.

\subsection{}\label{p2-fhtal150}
For any rational numbers $r\geq r'\geq 0$, we denote by
\begin{equation}\label{p2-fhtal150a}
\delta^{(r,r')}\colon \fC'^{(r)}\otimes_{R'_\uptau}\fC^{(r')}_\uptau\rightarrow \Omega'\otimes_{R'_\uptau}\fC'^{(r)}\otimes_{R'_\uptau}\fC^{(r')}_\uptau
\end{equation}
the $\hoRp$-derivation defined by 
\begin{equation}\label{p2-fhtal150b}
\delta^{(r,r')}=\delta_{\fC'^{(r)}}\otimes \id- \id \otimes \delta'_{\fC^{(r')}_\uptau},
\end{equation}
where $\delta_{\fC'^{(r)}}$ (resp.\ $\delta'_{\fC^{(r')}_\uptau}$) is the derivation \eqref{p2-rlps5e} for $(f',\tf')$ (resp.\ \eqref{p2-fhtal4e}).  
It is a Higgs $\hoRp$-field with coefficients in $\Omega'$.  By \ref{p2-hta190}, we have a commutative diagram
\begin{equation}\label{p2-fhtal150c}
\xymatrix{
{\fC^{(r)}_{\hoRp}\otimes_{R'_\uptau}\fC^{(r')}_\uptau}\ar[rr]^-(0.5){\uppsi^{(r,r')}}\ar[d]_{\id\otimes_{R'_\uptau}\delta_{\fC^{(r')}_\uptau}}&&
{\fC'^{(r)}\otimes_{R'_\uptau}\fC^{(r')}_\uptau}\ar[d]^{-\delta^{(r,r')}}\\
{\Omega\otimes_{R_\uptau}\fC^{(r)}_{\hoRp}\otimes_{R'_\uptau}\fC^{(r')}_\uptau}\ar[rr]^-(0.5){u\otimes_{R'_\uptau} \uppsi^{(r,r')}}&&
{\Omega'\otimes_{R'_\uptau}\fC'^{(r)}\otimes_{R'_\uptau}\fC^{(r')}_\uptau,}}
\end{equation}
where $\uppsi^{(r,r')}$ is the homomorphism \eqref{p2-fhtal110c}, $u$ is the morphism defined in \eqref{p2-fhtal200e} and 
$\delta_{\fC^{(r')}_\uptau}$ is the derivation \eqref{p2-fhtal4dd}. Similarly, we have a commutative diagram
\begin{equation}\label{p2-fhtal150cc}
\xymatrix{
{\fC^{(r)}_{\hoRp}\otimes_{R'_\uptau}\fC^{(r')}_\uptau}\ar[rr]^-(0.5){\uppsi^{(r,r')}}
\ar[d]_{\delta_{\fC^{(r)}_{\hoRp}}\otimes_{R'_\uptau}\id+\id\otimes_{R'_\uptau}\delta_{\fC^{(r')}_\uptau}}&&
{\fC'^{(r)}\otimes_{R'_\uptau}\fC^{(r')}_\uptau}\ar[d]^{\id\otimes_{R'_\uptau}\delta_{\fC^{(r')}_\uptau}}\\
{\Omega\otimes_{R_\uptau}\fC^{(r)}_{\hoRp}\otimes_{R'_\uptau}\fC^{(r')}_\uptau}\ar[rr]^-(0.5){\id\otimes_{R_\uptau} \uppsi^{(r,r')}}&&
{\Omega\otimes_{R_\uptau}\fC'^{(r)}\otimes_{R'_\uptau}\fC^{(r')}_\uptau,}}
\end{equation}
where $\delta_{\fC^{(r_2)}_{\hoRp}}$ is the derivation \eqref{p2-fhtal7g}.

\section{Relative cohomologies of Higgs--Tate algebras in the local case}

\subsection{}\label{p2-fhtal2}
The assumptions and notation of §\ref{p2-funchta} remain in force throughout this section.
We assume, furthermore, that the morphism $\gamma\colon (X',\cM_{X'})\rightarrow (X,\cM_X)$ \eqref{p2-fhtal1a} 
is smooth and saturated, and that its relative chart $ \mu\colon P\rightarrow P'$ is adequate in the sense of (\cite{ag1} 5.1.11). We set 
\begin{equation}\label{p2-fhtal1f}
\tOmega^1_{X'/X}=\Omega^1_{(X',\cM_{X'})/(X,\cM_X)}\ \ \ {\rm and}\ \ \  \tOmega^1_{R'/R}=\tOmega^1_{X'/X}(X').
\end{equation}
We have a canonical split exact sequence of free $R'$-modules of finite type
\begin{equation}\label{p2-fhtal1g}
0\rightarrow \tOmega^1_{R/\co_K}\otimes_RR' \rightarrow \tOmega^1_{R'/\co_K} \rightarrow \tOmega^1_{R'/R}\rightarrow 0. 
\end{equation}

To lighten notation, for any rational number $r\geq 0$, we set, with the conventions of \ref{p2-ncgt3}, 
\begin{equation}\label{p2-fhtal2d}
\uOmega'=\txi^{-1}\tOmega^1_{R'/R}\otimes_{R'}R'_\tau \ \ \ {\rm and}\ \ \ \uOmega'^{(r)}=p^r\uOmega'.
\end{equation}
The exact sequence \eqref{p2-fhtal1g} induces an exact sequence of free $R'_\uptau$-modules of finite type 
\begin{equation}\label{p2-fhtal2e}
0\rightarrow \Omega\otimes_{R_\uptau}R'_\uptau \stackrel{u}{\rightarrow} \Omega' \rightarrow \uOmega'\rightarrow 0. 
\end{equation}

For any $R'_\uptau$-algebra $A$, we  consider Higgs $A$-modules with coefficients in 
$\uOmega'\otimes_{R'_\uptau}A$. We say abusively  that they have coefficients in $\uOmega'$.
The category of these modules will be denoted by $\bHM(A,\uOmega')$.

For any rational numbers $r\geq r'\geq 0$, we denote by $\Omega'^{(r,r')}$ the $R'_\uptau$-module deduced from this exact sequence as in \ref{p2-hta14};
it fits into a commutative diagram \eqref{p2-hta14b}
\begin{equation}\label{p2-fhtal2f}
\xymatrix{
0\ar[r]&{\Omega^{(r')}\otimes_{R_\uptau}R'_\uptau}\ar[r]^-(0.5){u^{(r')}}&{\Omega'^{(r')}}\ar[r]&{\uOmega'^{(r')}}\ar[r]&0\\
0\ar[r]&{\Omega^{(r')}\otimes_{R_\uptau}R'_\uptau}\ar[r]^-(0.5){u^{(r,r')}}\ar@{=}[u]&{\Omega'^{(r,r')}}\ar[r]\ar@{^(->}[u]&{\uOmega'^{(r)}}\ar[r]\ar[u]&0\\
0\ar[r]&{\Omega^{(r)}\otimes_{R_\uptau}R'_\uptau}\ar[r]^-(0.5){u^{(r)}}\ar[u]&{\Omega'^{(r)}}\ar[r]\ar@{^(->}[u]&{\uOmega'^{(r)}}\ar[r]\ar@{=}[u]&0,}
\end{equation}
where the morphisms $u^{(r)}$ and $u^{(r,r')}$ are induced by $u$.

\subsection{}\label{p2-fhtal8}
Let $r,r'$ be rational numbers such that $r\geq r'\geq 0$. 
We denote by $\cL^{(r,r')}_{\tmX'/\tX'}$ the twist of the torsor $\cL_{\tmX'/\tX'}$ defined in \ref{p2-hta14}, and by 
$\cF^{(r,r')}_{\tmX'/\tX'}$ (resp.\ $\cC^{(r,r')}_{\tmX'/\tX'}$) the associated $\co_{\hmX'}$-extension (resp.\ $\co_{\hmX'}$-algebra). 
We set  
\begin{equation}\label{p2-fhtal8a}
\fF'^{(r,r')}=\Gamma(\hmX',\cF^{(r,r')}_{\tmX'/\tX'}) \ \ \ {\rm and}\ \ \ \fC'^{(r,r')}=\Gamma(\hmX',\cC^{(r)}_{\tmX'/\tX'}).
\end{equation}
We have a canonical exact sequence of $\hoRp$-modules
\begin{equation}\label{p2-fhtal8b}
0\rightarrow \hoRp\rightarrow \fF'^{(r,r')}\rightarrow \Omega'^{(r,r')}\otimes_{R'_\uptau}\hoRp \rightarrow 0,
\end{equation}
where $\Omega'^{(r,r')}$ is defined in \eqref{p2-fhtal2f}, and a canonical isomorphism of $\hoRp$-algebras
\begin{equation}\label{p2-fhtal8c}
\fC'^{(r,r')}\stackrel{\sim}{\rightarrow}\underset{\underset{n\geq 0}{\longrightarrow}}\lim\ \rS^n_{\hoRp}(\fF'^{(r,r')}). 
\end{equation}
We denote by $\hfC'^{(r,r')}$ the $p$-adic Hausdorff completion of $\fC'^{(r,r')}$. 

Let
\begin{equation}\label{p2-fhtal8d}
d_{\fC'^{(r,r')}}\colon \fC'^{(r,r')} \rightarrow \Omega'^{(r,r')}\otimes_{R'_\uptau}\fC'^{(r,r')}
\end{equation}
be the universal $\hoRp$-derivation of $\fC'^{(r,r')}$ \eqref{p2-hta18i}; it is a Higgs $\hoRp$-field with coefficients in $\Omega'^{(r,r')}$.
We denote by 
\begin{equation}\label{p2-fhtal8e}
\delta_{\fC'^{(r,r')}}\colon \fC'^{(r,r')} \rightarrow \Omega'\otimes_{R'_\uptau}\fC'^{(r,r')}
\end{equation}
the $\hoRp$-derivation induced by $d_{\fC'^{(r,r')}}$ and the canonical injection $\Omega'^{(r,r')}\rightarrow \Omega'$, and by 
\begin{equation}\label{p2-fhtal8f}
\udelta_{\fC'^{(r,r')}}\colon \fC'^{(r,r')} \rightarrow \uOmega'\otimes_{R'_\uptau}\fC'^{(r,r')}
\end{equation}
the $\hoRp$-derivation induced by $\delta_{\fC'^{(r,r')}}$ and the canonical morphism $\Omega'\rightarrow \uOmega'$. Let 
\begin{eqnarray}
\delta_{\hfC'^{(r,r')}}\colon \hfC'^{(r,r')} &\rightarrow &\Omega'\otimes_{R'_\uptau}\hfC'^{(r,r')},\label{p2-fhtal8g}\\
\udelta_{\hfC'^{(r,r')}}\colon \hfC'^{(r,r')} &\rightarrow &\uOmega'\otimes_{R'_\uptau}\hfC'^{(r,r')},\label{p2-fhtal8h}
\end{eqnarray}
be the extensions of $\delta_{\fC'^{(r,r')}}$ and $\udelta_{\fC'^{(r,r')}}$, respectively, to the $p$-adic completions. 
We denote by $\mK^\bullet(\hfC'^{(r,r')})$ (resp.\ $\umK^\bullet(\hfC'^{(r,r')})$) the Dolbeault complex of $(\hfC'^{(r,r')},\delta_{\hfC'^{(r,r')}})$ (resp.\ $(\hfC'^{(r,r')},\udelta_{\hfC'^{(r,r')}})$). 

Let $t,t'$ be rational numbers such that  $t\geq t'\geq 0$, $r\geq t$ and $r'\geq t'$. We have a canonical morphism of $\hoRp$-extensions \eqref{p2-hta14i}
\begin{equation}\label{p2-fhtal8i}
\tta'^{r,r',t,t'}\colon\fF'^{(r,r')}\rightarrow \fF'^{(t,t')}
\end{equation}
that fits into a commutative diagram
\begin{equation}\label{p2-fhtal8j}
\xymatrix{
0\ar[r]&{\hoRp}\ar@{=}[d]\ar[r]&{\fF'^{(r,r')}}\ar[r]\ar[d]^-(0.5){\tta'^{r,r',t,t'}}&{\Omega'^{(r,r')}\otimes_{R'_\uptau}\hoRp}\ar[r]\ar[d]&0\\
0\ar[r]&{\hoRp}\ar[r]&{\fF'^{(t,t')}}\ar[r]&{\Omega'^{(t,t')}\otimes_{R'_\uptau}\hoRp}\ar[r]&0,}
\end{equation}
where the vertical right arrow is the canonical morphism \eqref{p2-hta14h}. 
We deduce a canonical homomorphism of $\hoRp$-algebras \eqref{p2-hta14k}
\begin{equation}\label{p2-fhtal8k}
\alpha'^{r,r',t,t'}\colon \fC'^{(r,r')}\rightarrow \fC'^{(t,t')}.
\end{equation}

It immediately follows from \eqref{p2-fhtal8j} that we have 
\begin{equation}\label{p2-fhtal8l}
(\id \otimes \alpha'^{r,r',t,t'}) \circ \delta_{\fC^{(r,r')}}=\delta_{\fC'^{(t,t')}}\circ  \alpha'^{r,r',t,t'}.
\end{equation}
Therefore, $\alpha'^{r,r',t,t'}$ induces morphisms of complexes 
\begin{eqnarray}
\upiota'^{r,r',t,t'}\colon \mK^\bullet(\hfC'^{(r,r')})&\rightarrow& \mK^\bullet(\hfC'^{(t,t')}),\label{p2-fhtal8m}\\
\uupiota'^{r,r',t,t'}\colon \umK^\bullet(\hfC'^{(r,r')})&\rightarrow& \umK^\bullet(\hfC'^{(t,t')}).\label{p2-fhtal8n}
\end{eqnarray}

\subsection{}\label{p2-fhtal9}
Let $r,r'$ be rational numbers such that $r\geq r'\geq 0$, $\tgamma\in \cL_{\tX'/\tX}(\coX')$. 
We have a canonical $\hoRp$-linear morphism \eqref{p2-hta17c}, associated with the diagram \eqref{p2-fhtal3a},
\begin{equation}\label{p2-fhtal9a}
\upvarphi_\tgamma^{(r,r')}\colon \fF^{(r')}_{\hoRp}\rightarrow \fF'^{(r,r')}
\end{equation}
that fits into a commutative diagram 
\begin{equation}\label{p2-fhtal9b}
\xymatrix{
0\ar[r]&{\hoRp}\ar@{=}[d]\ar[r]&{\fF^{(r')}_{\hoRp}}\ar[r]\ar[d]^{\upvarphi_\tgamma^{(r,r')}}&{\Omega^{(r')}\otimes_{R_\uptau}\hoRp}\ar[r]\ar[d]^{u^{(r,r')}\otimes_{R'_\uptau}\hoRp}&0\\
0\ar[r]&{\hoRp}\ar[r]&{\fF'^{(r,r')}}\ar[r]&{\Omega'^{(r,r')}\otimes_{R'_\uptau}\hoRp}\ar[r]&0,}
\end{equation}
where $u^{(r,r')}$ is defined in \eqref{p2-fhtal2f}. 
We deduce a canonical morphism of $\hoRp$-algebras \eqref{p2-hta17e} 
\begin{equation}\label{p2-fhtal9c}
\upphi_\tgamma^{(r,r')}\colon \fC^{(r')}_{\hoRp}\rightarrow \fC'^{(r,r')}. 
\end{equation}

It follows from \eqref{p2-fhtal9b} that $\upphi_\tgamma^{(r,r')}$ makes $\udelta_{\fC'^{(r,r')}}$ \eqref{p2-fhtal8f} into a $\fC^{(r')}_{\hoRp}$-derivation. 
We denote by $\utmK^\bullet_\tgamma(\hfC'^{(r,r')})$ the augmented Dolbeault complex
\begin{equation}\label{p2-fhtal9d}
\xymatrix{
{\fC^{(r')}\hotimes_{\oR}\oR'}\ar[r]^-(0.5){\upphi_\tgamma^{(r,r')}}&{\umK^0(\hfC'^{(r,r')})}\ar[r]&{\umK^1(\hfC'^{(r,r')})}\ar[r]&\dots},
\end{equation}
where $\umK^\bullet(\hfC'^{(r,r')})$ is the Dolbeault complex of $(\hfC'^{(r,r')}, \udelta_{\hfC'^{(r,r')}})$ \eqref{p2-fhtal8h}, 
$\fC^{(r')}\hotimes_{\oR}\oR'$ is placed in degree $-1$ and the tensor product $\hotimes$ is completed for the $p$-adic topology.

Let $t,t'$ be rational numbers such that $t\geq t'\geq 0$, $r\geq t$ and $r'\geq t'$. The diagram 
\begin{equation}\label{p2-fhtal9e}
\xymatrix{
{\fC^{(r')}_{\hoRp}}\ar[r]^-(0.5){\upphi_\tgamma^{(r,r')}}\ar[d]_{\alpha^{t,t'}_{\hoRp}}&{\fC'^{(r,r')}}\ar[d]^{\alpha'^{r,r',t,t'}}\\
{\fC^{(t')}_{\hoRp}}\ar[r]^-(0.5){\upphi_\tgamma^{(t,t')}}&{\fC'^{(t,t'),}}}
\end{equation}
where the vertical arrows are defined in \eqref{p2-fhtal7h} and \eqref{p2-fhtal8k}, is commutative \eqref{p2-hta17f}. 
By \eqref{p2-fhtal8l} et \eqref{p2-fhtal9e}, $\alpha'^{r,r',t,t'}$ and $\alpha^{t,t'}_{\hoRp}$ induce a morphism of complexes
\begin{equation}\label{p2-fhtal9f}
\tuupiota'^{r,r',t,t'}_\tgamma\colon \utmK^\bullet_\tgamma(\hfC'^{(r,r')})\rightarrow \utmK^\bullet_\tgamma(\hfC'^{(t,t')}),
\end{equation}
compatible with the morphism $\uupiota'^{r,r',t,t'}$ \eqref{p2-fhtal8n}

\begin{prop}\label{p2-fhtal10}
Let $\tgamma\in \cL_{\tX'/\tX}(\coX')$, $r,r',t,t'$ be rational numbers such that $r\geq r'\geq 0$, $t\geq t'\geq 0$, $r\geq t$ and $r'\geq t'$. Then,
\begin{itemize}
\item[{\rm (i)}] There exists a rational number $\alpha\geq 0$ depending on $r$ and $t$ 
but not on the morphisms $f$ and $\gamma$ satisfying the conditions of \ref{p2-fhtal1} neither on $\tgamma$, such that
\begin{equation}\label{p2-fhtal10a}
p^\alpha\tuupiota^{r,r',t,t'}_\tgamma\colon \utmK^\bullet_\tgamma(\hfC'^{(r,r')})\rightarrow \utmK^\bullet_\tgamma(\hfC'^{(t,t')}),
\end{equation}
where $\tuupiota^{r,r',t,t'}_\tgamma$ is the morphism \eqref{p2-fhtal9f}, is homotopic to $0$ by an $\hoRp$-linear homotopy.
\item[{\rm (ii)}] The canonical morphism
\begin{equation}\label{p2-fhtal10b}
\tuupiota^{r,r',t,t'}_\tgamma\otimes_{\mZ_p}\mQ_p\colon \utmK^\bullet_\tgamma(\hfC'^{(r,r')})\otimes_{\mZ_p}\mQ_p\rightarrow
\utmK^\bullet_\tgamma(\hfC'^{(t,t')})\otimes_{\mZ_p}\mQ_p
\end{equation}
is homotopic to $0$ by a continuous homotopy.
\end{itemize}
\end{prop}

Indeed, we may replace $(r,r')$ by $(r,t')$ in which case the proposition was proved in (\cite{ag2} 5.3.2). 
By \ref{p2-hta170}, the definition of the $\hoRp$-algebra $\hfC'^{(r,r')}$ \eqref{p2-fhtal8a} is equivalent to that in (\cite{ag2} (5.2.14.5)).

\subsection{}\label{p2-fhtal11}
We denote by $I$ the set of triples of rational numbers $\ur=(r_1,r_2,r_3)$ such that $r_1\geq r_2\geq r_3\geq 0$. 
For any such a triple, we have a canonical morphism of $\hoRp$-modules \eqref{p2-hta15h}, associated with the diagram \eqref{p2-fhtal3a},
\begin{equation}\label{p2-fhtal11a}
\upvarphi^\ur\colon \fF^{(r_2)}_{\hoRp}\rightarrow \fF'^{(r_1,r_2)} \otimes_{R'_\uptau}\fF^{(r_3)}_\uptau,
\end{equation}
and a canonical morphism of $\hoRp$-algebras \eqref{p2-hta15i}
\begin{equation}\label{p2-fhtal11b}
\upphi^\ur\colon \fC^{(r_2)}_{\hoRp}\rightarrow \fC'^{(r_1,r_2)} \otimes_{R'_\uptau}\fC^{(r_3)}_\uptau,
\end{equation}
compatible with $\upvarphi^{\ur}$. The latter induces by linearization a morphism of $(\hoRp\otimes_{R'_\uptau}\fC^{(r_3)}_\uptau)$-algebras \eqref{p2-hta15e}
\begin{equation}\label{p2-fhtal11c}
\uppsi^\ur\colon \fC^{(r_2)}_{\hoRp}\otimes_{R'_\uptau}\fC^{(r_3)}_\uptau\rightarrow \fC'^{(r_1,r_2)}\otimes_{R'_\uptau}\fC^{(r_3)}_\uptau.
\end{equation}

Let $\tgamma\in \cL_{\tX'/\tX}(\coX')$, $\tgamma^{(r_3)}$ the canonical image of $\tgamma$ in $\cL^{(r_3)}_{\tX'/\tX}(\coX')$, 
$\rho_{\tgamma^{(r_3)}}\colon \fF_\uptau^{(r_3)}\rightarrow R'_\uptau$ the associated splitting of the extension \eqref{p2-fhtal4b}, 
$\varrho_{\tgamma^{(r_3)}}\colon \fC_\uptau^{(r_3)}\rightarrow R'_\uptau$ the associated homomorphism of $R'_\uptau$-algebras (see \ref{p1-prem1}). 
By \ref{p2-hta160}(i), we have 
\begin{eqnarray}
\upvarphi_\tgamma^{(r_1,r_2)}&=&(\id \otimes_{R'_\uptau} \rho_{\tgamma^{(r_3)}})\circ \upvarphi^\ur,\label{p2-fhtal11d}\\
\upphi_\tgamma^{(r_1,r_2)}&=&(\id \otimes_{R'_\uptau} \varrho_{\tgamma^{(r_3)}})\circ \upphi^\ur,\label{p2-fhtal11e}
\end{eqnarray}
where $\upvarphi_\tgamma^{(r_1,r_2)}$ and $\upphi_\tgamma^{(r_1,r_2)}$ are defined in \eqref{p2-fhtal9a} and \eqref{p2-fhtal9c}.

We have a canonical isomorphism of $\hoRp$-algebras \eqref{p2-hta16d}
\begin{equation}\label{p2-fhtal11g}
\Lambda^{(r_2,r_3)}_\tgamma\colon \fC^{(r_2)}_{\hoRp}\otimes_{R'_\uptau} \fC^{(r_3)}_\uptau
\stackrel{\sim}{\rightarrow} \fC^{(r_2)}_{\hoRp}\otimes_{R'_\uptau} \fC^{(r_3)}_\uptau.
\end{equation}
By \ref{p2-hta160}(iii), we have
\begin{equation}\label{p2-fhtal11h}
\uppsi^\ur=(\upphi_\tgamma^{(r_1,r_2)}\otimes \id)\circ \Lambda^{(r_2,r_3)}_\tgamma.
\end{equation}

Let $\ur'=(r'_1,r'_2,r'_3)\in I$ be such that $r_i\geq r'_i$ for all $1\leq i\leq 3$. By \eqref{p2-hta15k}, the diagram 
\begin{equation}\label{p2-fhtal11f}
\xymatrix{
{\fC^{(r_2)}_{\hoRp}\otimes_{R'_\uptau} \fC^{(r_3)}_\uptau}\ar[r]^-(0.5){\uppsi^\ur}\ar[d]_{\alpha^{r_2,r'_2}_\hoRp\otimes_{R'_\uptau} \alpha_\uptau^{r_3,r'_3}}
&{\fC'^{(r_1,r_2)} \otimes_{R'_\uptau} \fC^{(r_3)}_\uptau}
\ar[d]^{\alpha'^{r_1,r_2,r'_1,r'_2}\otimes_{R'_\uptau} \alpha_\uptau^{r_3,r'_3}}\\
{\fC^{(r'_2)}_{\hoRp}\otimes_{R'_\uptau} \fC^{(r'_3)}_\uptau}\ar[r]^-(0.5){\uppsi^{\ur'}}&{\fC'^{(r'_1,r'_2)} \otimes_{R'_\uptau} \fC^{(r'_3)}_{\uptau},}}
\end{equation}
where the vertical arrows are defined in \eqref{p2-fhtal7h}, \eqref{p2-fhtal4k} and \eqref{p2-fhtal8k}, is commutative. We have a similar commutative diagram for $\upphi^\ur$. 

By \eqref{p2-hta16i}, the diagram 
\begin{equation}\label{p2-fhtal11i}
\xymatrix{
{\fC^{(r_2)}_{\hoRp}\otimes_{R'_\uptau} \fC^{(r_3)}_\uptau}\ar[rr]^-(0.5){\Lambda^{(r_2,r_3)}_\tgamma}\ar[d]_{\alpha^{r_2,r'_2}_{\hoRp}\otimes_{R'_\uptau} \alpha^{r_3,r'_3}_\uptau}&&
{\fC^{(r_2)}_{\hoRp}\otimes_{R'_\uptau} \fC^{(r_3)}_\uptau}\ar[d]^{\alpha^{r_2,r'_2}_{\hoRp}\otimes_{R'_\uptau} \alpha^{r_3,r'_3}_\uptau}\\
{\fC^{(r'_2)}_{\hoRp}\otimes_{R'_\uptau} \fC^{(r'_3)}_\uptau}\ar[rr]^-(0.5){\Lambda^{(r'_2,r'_3)}_\tgamma}&&{\fC^{(r'_2)}_{\hoRp}\otimes_{R'_\uptau} \fC^{(r'_3)}_{\uptau}}}
\end{equation}
is commutative. 

\begin{rema}\label{p2-fhtal111}
By \ref{p2-hta151}, for all rational numbers $r\geq r'\geq 0$, setting $\ur=(r,r,r') \in I$ \eqref{p2-fhtal11}, we have 
$\bvupphi^{(r,r')}=\bvupphi^\ur$, see \eqref{p2-fhtal110b} and \eqref{p2-fhtal11b}, and 
$\bvuppsi^{(r,r')}=\bvuppsi^\ur$, see \eqref{p2-fhtal110c} and \eqref{p2-fhtal11c}. 
\end{rema}

\subsection{}\label{p2-fhtal120}
Let $r,r'$ be rational numbers such that $r\geq r'\geq 0$. 
The action of $\Delta'$ on $\cL_{\tmX'/\tX'}$ defined in \ref{p2-fhtal12} induces an action on $\cL^{(r,r')}_{\tmX'/\tX'}$. 
We deduce an $\hoRp$-semilinear action of $\Delta'$ on $\fF'^{(r,r')}$ such that the morphisms of \eqref{p2-fhtal8b} are $\Delta'$-equivariant. 
The latter induces an action of $\Delta'$ on $\fC'^{(r,r')}$ by ring automorphisms compatible with its action on $\hoRp$. 
The Higgs $\hoRp$-field $d_{\fC'^{(r,r')}}$ \eqref{p2-fhtal8d} is $\Delta'$-equivariant. 
It follows from (\cite{ag2} 3.2.23) that the actions of $\Delta'$ on $\fF'^{(r,r')}$, $\fC'^{(r,r')}$ and  $\hfC'^{(r,r')}$ are continuous for the $p$-adic topologies.

Let $\ur=(r_1,r_2,r_3) \in I$ \eqref{p2-fhtal11}. We deduce from \eqref{p2-fhtal12b} that the canonical morphism \eqref{p2-hta15d}, 
associated with the diagram \eqref{p2-fhtal3a},
\begin{equation}
\uppsi^\ur_\cL\colon \cL^{(r_1,r_2)}_{\tmX'/\tX'}\times h'^+(\cL^{(r_3)}_{\tX'/\tX}) \rightarrow \cL^{(r_2)}_{\tmX'/\tX} \times h'^+(\cL^{(r_3)}_{\tX'/\tX})
\end{equation}
is $\Delta'$-equivariant. Therefore, the morphism $\uppsi^\ur$ \eqref{p2-fhtal11c} is $\Delta'$-equivariant, where $\Delta'$ acts trivially on $\fC^{(r_3)}_\uptau$. 
It follows that the morphisms $\upvarphi^\ur$ \eqref{p2-fhtal11a} and $\upphi^\ur$ \eqref{p2-fhtal11b} are $\Delta'$-equivariant. 

We deduce by \eqref{p2-fhtal11d} and \eqref{p2-fhtal11e} that for any $\tgamma\in \cL_{\tX'/\tX}(\coX')$, the morphisms $\upvarphi_\tgamma^{(r,r')}$ \eqref{p2-fhtal9a}
and $\upphi_\tgamma^{(r,r')}$ \eqref{p2-fhtal9c} are $\Delta'$-equivariant. We also immediately see that $\Lambda^{(r_2,r_3)}_\tgamma$ \eqref{p2-fhtal11g} is $\Delta'$-equivariant.

The following proposition generalizes (\cite{ag2} 5.2.15).

\begin{prop}\label{p2-fhtal13}
Let $\ur=(r_1,r_2,r_3)$, $\ur'=(r'_1,r'_2,r'_3)$ be two elements of $I$ \eqref{p2-fhtal11} such that $r_1> r'_1$, $r_2\geq  r'_2$ and $r_3\geq  r'_3$, 
$\oR^\intern$ the ring defined in \eqref{p2-fhtal2i}, $\Pi$ the group defined in \eqref{p2-fhtal2h}. Then,
\begin{itemize}
\item[{\rm (i)}] For every integer $n\geq 1$, the homomorphism 
\begin{equation}\label{p2-fhtal13a}
(\fC^{(r_2)}\otimes_{\oR} \oR^\intern)\otimes_{R'_\uptau}(\fC^{(r_3)}_\uptau/p^n\fC^{(r_3)}_\uptau)\rightarrow 
(\fC'^{(r_1,r_2)} \otimes_{R'_\uptau}(\fC^{(r_3)}_\uptau/p^n\fC^{(r_3)}_\uptau))^{\Pi}
\end{equation}
induced by $\uppsi^\ur$ \eqref{p2-fhtal11c} is $\alpha$-injective. Let $\cH^\ur_n$ be its cokernel. 
\item[{\rm (ii)}] There exists an integer $a\geq 0$, depending on $r_1$, $r'_1$ and $\ell=\dim(X'/X)$,
but not on the morphisms $f$ and $\gamma$ satisfying the conditions of \ref{p2-fhtal1}, such that for every integer $n\geq 1$,
the morphism $\cH^{\ur}_n\rightarrow \cH^{\ur'}_n$ induced by \eqref{p2-fhtal11f} is annihilated by $p^a$.
\item[{\rm (iii)}] There exists an integer $b\geq 0$, depending on $r_1$, $r'_1$ and $\ell=\dim(X'/X)$,
but not on the morphisms $f$ and $\gamma$ satisfying the conditions of \ref{p2-fhtal1}, such that for all integers $n,q\geq 1$,
the canonical morphism
\begin{equation}\label{p2-fhtal13b}
\rH^q(\Pi,\fC'^{(r_1,r_2)} \otimes_{R'_\uptau}(\fC^{(r_3)}_\uptau/p^n\fC^{(r_3)}_\uptau))\rightarrow \rH^q(\Pi,\fC'^{(r'_1,r'_2)} \otimes_{R'_\uptau}(\fC^{(r'_3)}_\uptau/p^n\fC^{(r'_3)}_\uptau))
\end{equation}
is annihilated by $p^b$.
\end{itemize}
\end{prop}

Indeed, since $\tf$ is smooth and $\tX'$ is affine, 
there exist $\tgamma\in \cL_{\tX'/\tX}(\coX')$ \eqref{p2-fhtal3a}. 
The proposition follows then from the analogous statement for the homomorphism $\upphi^{(r_1,r_2)}_\tgamma$ \eqref{p2-fhtal9c}
proved in (\cite{ag2} 5.2.15), in view of \eqref{p2-fhtal11h}, \eqref{p2-fhtal11i} and (\cite{agt} II.3.15).

\begin{cor}\label{p2-fhtal14}
Let $t_2,t_3$ be two rational numbers such that $t_2\geq t_3 \geq 0$, $I_{t_2,t_3}$ the subset of elements $\ur=(r_1,r_2,r_3)$ of $I$ such that 
$r_1>r_2\geq t_2$ and $r_3\geq t_3$. Then, 
\begin{itemize}
\item[{\rm (i)}] The canonical morphism
\begin{equation}\label{p2-fhtal14a}
(\fC^{(t_2)}\hotimes_{\oR}\oR^\intern\hotimes_{R'_\uptau}\fC^{(t_3)}_\uptau)[\frac 1 p]\rightarrow
\underset{\underset{\ur\in I_{t_2,t_3}}{\longrightarrow}}{\lim} ((\fC'^{(r_1,r_2)} \hotimes_{R'_\uptau}\fC^{(r_3)}_\uptau)[\frac 1 p])^\Pi,
\end{equation}
where the tensor products $\hotimes$ are completed for the $p$-adic topology, is an isomorphism.
\item[{\rm (ii)}] For every integer $i\geq 1$, we have
\begin{equation}\label{p2-fhtal14b}
\underset{\underset{\ur\in I_{t_2,t_3}}{\longrightarrow}}{\lim}\
\rH^i_\cont(\Pi,(\fC'^{(r_1,r_2)} \hotimes_{R'_\uptau}\fC^{(r_3)}_\uptau)[\frac 1 p])=0,
\end{equation}
where $(\fC'^{(r_1,r_2)} \hotimes_{R'_\uptau}\fC^{(r_3)}_\uptau)[\frac 1 p]$ is equipped with the $p$-adic topology {\rm (\cite{ag2} 2.1.1)}. 
\end{itemize}
\end{cor}

We take again the notation of \ref{p2-fhtal13}. It follows from \ref{p2-fhtal13}(i) and (\cite{gr} 2.4.2(ii)) that for any triple $\ur=(r_1,r_2,r_3)$ of $I_{t_2,t_3}$, 
the canonical sequence
\begin{equation}\label{p2-fhtal14c}
0\rightarrow \fC^{(r_2)}\hotimes_{\oR}\oR^\intern \hotimes_{R'_1}\fC^{(r_3)}_\uptau
\rightarrow (\fC'^{(r_1,r_2)} \hotimes_{R'_\uptau}\fC^{(r_3)}_\uptau)^{\Pi}\rightarrow \underset{\underset{n\geq 0}{\longleftarrow}}{\lim}\ \cH_n^{\ur}
\end{equation}
is $\alpha$-exact. Proposition (i) follows by inverting $p$ and then
passing to the inductive limit on $\ur\in I_{t_2,t_3}$, in view of \ref{p2-fhtal13}(ii).

As we have (\cite{agt} (II.3.10.2))
\begin{equation}\label{p2-fhtal14d}
\rR^1\underset{\underset{n\geq 0}{\longleftarrow}}{\lim}\ (\fC^{(r_2)}\otimes_{\oR} \oR^\intern)\otimes_{R'_1}(\fC^{(r_3)}_\uptau/p^n\fC^{(r_3)}_\uptau)=0,
\end{equation}
it follows from \ref{p2-fhtal13}(i) and (\cite{gr} 2.4.2(ii)) that the canonical morphism
\begin{equation}\label{p2-fhtal14e}
\rR^1\underset{\underset{n\geq 0}{\longleftarrow}}{\lim}\ (\fC'^{(r_1,r_2)} \otimes_{R'_\uptau}(\fC^{(r_3)}_\uptau/p^n\fC^{(r_3)}_\uptau))^{\Pi}\rightarrow \rR^1\underset{\underset{n\geq 0}{\longleftarrow}}{\lim}\ \cH_n^\ur
\end{equation}
is an $\alpha$-isomorphism.

By (\cite{agt} (II.3.10.4) and (II.3.10.5)), for every $i \geq 1$, we have a canonical exact sequence
\begin{eqnarray}
\lefteqn{
0\rightarrow \rR^1\underset{\underset{n}{\longleftarrow}}{\lim}\ \rH^{i-1}(\Pi,\fC'^{(r_1,r_2)} \otimes_{R'_\uptau}(\fC^{(r_3)}_\uptau/p^n\fC^{(r_3)}_\uptau))\rightarrow}\nonumber\\
&&
\rH^i_\cont(\Pi,\hfC'^\ur_\uptau)\rightarrow \underset{\underset{n}{\longleftarrow}}{\lim}\ \rH^{i}(\Pi,\fC'^{(r_1,r_2)} \otimes_{R'_\uptau}(\fC^{(r_3)}_\uptau/p^n\fC^{(r_3)}_\uptau))\rightarrow 0.
\end{eqnarray}
We deduce, in view of \eqref{p2-fhtal14e} and \ref{p2-fhtal13}(ii)-(iii), that for all elements $\ur=(r_1,r_2,r_3)$ and $\ur'=(r'_1,r'_2,r'_3)$ of $I$ \eqref{p2-fhtal11} 
such that $r'_1> r_1$, $r'_2\geq r_2$ and $r'_3\geq r_3$, there exists a rational number $\gamma \geq 0$ such that the canonical morphism
\begin{equation}\label{p2-fhtal14g}
\rH^i_\cont(\Pi,\fC'^{(r'_1,r'_2)} \hotimes_{R'_\uptau}\fC^{(r'_3)}_\uptau)\rightarrow \rH^i_\cont(\Pi,\fC'^{(r_1,r_2)} \hotimes_{R'_\uptau}\fC^{(r_3)}_\uptau)
\end{equation}
is annihilated by $p^\gamma$. Proposition (ii) follows by passing to the inductive limit on $\ur\in I_{t_2,t_3}$.

\subsection{}\label{p2-fhtal15}
For any $\ur=(r_1,r_2,r_3)\in I$ \eqref{p2-fhtal11}, we denote by
\begin{equation}\label{p2-fhtal15a}
\delta^\ur\colon \fC'^{(r_1,r_2)}\otimes_{R'_\uptau}\fC^{(r_3)}_\uptau\rightarrow \Omega'\otimes_{R'_\uptau}\fC'^{(r_1,r_2)}\otimes_{R'_\uptau}\fC^{(r_3)}_\uptau
\end{equation}
the $\hoRp$-derivation defined by 
\begin{equation}\label{p2-fhtal15b}
\delta^\ur=\delta_{\fC'^{(r_1,r_2)}}\otimes \id- \id \otimes \delta'_{\fC^{(r_3)}_\uptau},
\end{equation}
where $\delta_{\fC'^{(r_1,r_2)}}$ (resp.\ $\delta'_{\fC^{(r_3)}_\uptau}$) is the derivation \eqref{p2-fhtal8e} (resp.\ \eqref{p2-fhtal4e}).  
It is a Higgs $\hoRp$-field with coefficients in $\Omega'$.  By \ref{p2-hta19}, we have a commutative diagram
\begin{equation}\label{p2-fhtal15c}
\xymatrix{
{\fC^{(r_2)}_{\hoRp}\otimes_{R'_\uptau}\fC^{(r_3)}_\uptau}\ar[rr]^-(0.5){\uppsi^\ur}\ar[d]_{\id\otimes_{R'_\uptau}\delta_{\fC^{(r_3)}_\uptau}}&&
{\fC'^{(r_1,r_2)}\otimes_{R'_\uptau}\fC^{(r_3)}_\uptau}\ar[d]^{-\delta^\ur}\\
{\Omega\otimes_{R_\uptau}\fC^{(r_2)}_{\hoRp}\otimes_{R'_\uptau}\fC^{(r_3)}_\uptau}\ar[rr]^-(0.5){u\otimes_{R'_\uptau} \uppsi^\ur}&&
{\Omega'\otimes_{R'_\uptau}\fC'^{(r_1,r_2)}\otimes_{R'_\uptau}\fC^{(r_3)}_\uptau,}}
\end{equation}
where $\uppsi^\ur$ is the homomorphism \eqref{p2-fhtal11c}, $u$ is the morphism defined in \eqref{p2-fhtal200e} and $\delta_{\fC^{(r_3)}_\uptau}$ is the derivation \eqref{p2-fhtal4dd}. Similarly, we have a commutative diagram
\begin{equation}\label{p2-fhtal15cc}
\xymatrix{
{\fC^{(r_2)}_{\hoRp}\otimes_{R'_\uptau}\fC^{(r_3)}_\uptau}\ar[rr]^-(0.5){\uppsi^\ur}
\ar[d]_{\delta_{\fC^{(r_2)}_{\hoRp}}\otimes_{R'_\uptau}\id+\id\otimes_{R'_\uptau}\delta_{\fC^{(r_3)}_\uptau}}&&
{\fC'^{(r_1,r_2)}\otimes_{R'_\uptau}\fC^{(r_3)}_\uptau}\ar[d]^{\id\otimes_{R'_\uptau}\delta_{\fC^{(r_3)}_\uptau}}\\
{\Omega\otimes_{R_\uptau}\fC^{(r_2)}_{\hoRp}\otimes_{R'_\uptau}\fC^{(r_3)}_\uptau}\ar[rr]^-(0.5){\id\otimes_{R_\uptau} \uppsi^\ur}&&
{\Omega\otimes_{R_\uptau}\fC'^{(r_1,r_2)}\otimes_{R'_\uptau}\fC^{(r_3)}_\uptau,}}
\end{equation}
where $\delta_{\fC^{(r_2)}_{\hoRp}}$ is the derivation \eqref{p2-fhtal7g}.

We denote by 
\begin{equation}\label{p2-fhtal15d}
\hdelta^\ur\colon \fC'^{(r_1,r_2)}\hotimes_{R'_\uptau}\fC^{(r_3)}_\uptau\rightarrow \Omega'\otimes_{R'_\uptau}\fC'^{(r_1,r_2)}\hotimes_{R'_\uptau}\fC^{(r_3)}_\uptau
\end{equation}
the extension of $\delta^\ur$ to the $p$-adic completions and by $\mK^\bullet(\fC'^{(r_1,r_2)}\hotimes_{R'_\uptau}\fC^{(r_3)}_\uptau)$ 
the Dolbeault complex of $(\fC'^{(r_1,r_2)}\hotimes_{R'_\uptau}\fC^{(r_3)}_\uptau,-\hdelta^\ur)$. In view of \eqref{p2-fhtal15c}, $\uppsi^\ur$ induces a canonical 
morphism of complexes 
\begin{equation}\label{p2-fhtal15e}
\fC^{(r_2)}_{\hoRp}\hotimes_{R'_\uptau}\mK^\bullet(\fC^{(r_3)}_\uptau)\rightarrow  \mK^\bullet(\fC'^{(r_1,r_2)}\hotimes_{R'_\uptau}\fC^{(r_3)}_\uptau),
\end{equation}
where $\mK^\bullet(\hfC^{(r_3)}_\uptau)$ is the Dolbeault complex of $(\hfC^{(r_3)}_\uptau,\delta_{\hfC^{(r_3)}_\uptau})$ \eqref{p2-fhtal4g} and the left tensor product $\hotimes$ 
is defined term by term and is completed for the $p$-adic topology. We denote by $\tmK^\bullet(\fC'^{(r_1,r_2)}\hotimes_{R'_\uptau}\fC^{(r_3)}_\uptau)$ 
the mapping cone of the morphism \eqref{p2-fhtal15e}.

Let $\ur'=(r'_1,r'_2,r'_3)\in I$ be such that $r_i\geq r'_i$ for all $1\leq i\leq 3$.
By \eqref{p2-fhtal8l} and \eqref{p2-fhtal4f}, the homomorphism $\alpha'^{r_2,r_3,r'_2,r'_3}$ \eqref{p2-fhtal8k} and $\alpha^{r',r}_\uptau$ \eqref{p2-fhtal4k} induce a morphism of complexes
\begin{equation}\label{p2-fhtal15f}
\upiota^{\ur,\ur'}\colon \mK^\bullet(\fC'^{(r_1,r_2)}\hotimes_{R'_\uptau}\fC^{(r_3)}_\uptau)\rightarrow \mK^\bullet(\fC'^{(r'_1,r'_2)}\hotimes_{R'_\uptau}\fC^{(r'_3)}_\uptau).
\end{equation}
Moreover, the diagram 
\begin{equation}\label{p2-fhtal15g}
\xymatrix{
{\fC^{(r_2)}_{\hoRp}\hotimes_{R'_\uptau}\mK^\bullet(\fC^{(r_3)}_\uptau)}\ar[r]\ar[d]_{\alpha_{\hoRp}^{r'_2,r_2}\hotimes_{R'_\uptau}\upiota^{r_3,r'_3}_\uptau}&
{\mK^\bullet(\fC'^{(r_1,r_2)}\hotimes_{R'_\uptau}\fC^{(r_3)}_\uptau)}\ar[d]^{\upiota^{\ur',\ur}}\\
{\fC^{(r'_2)}_{\hoRp}\hotimes_{R'_\uptau}\mK^\bullet(\fC^{(r'_3)}_\uptau)}\ar[r]&{\mK^\bullet(\fC'^{(r'_1,r'_2)}\hotimes_{R'_\uptau}\fC^{(r'_3)}_\uptau),}}
\end{equation}
where $\upiota^{r_3,r'_3}_\uptau$ is the morphism \eqref{p2-fhtal4j}, is commutative.  
We deduce a morphism of complexes
\begin{equation}\label{p2-fhtal15h}
\tupiota^{\ur,\ur'}\colon \tmK^\bullet(\fC'^{(r_1,r_2)}\hotimes_{R'_\uptau}\fC^{(r_3)}_\uptau)\rightarrow \tmK^\bullet(\fC'^{(r'_1,r'_2)}\hotimes_{R'_\uptau}\fC^{(r'_3)}_\uptau).
\end{equation}

\begin{prop}\label{p2-fhtal16}
Let $\ur=(r_1,r_2,r_3)$, $\ur'=(r'_1,r'_2,r'_3)$ be two elements of $I$ \eqref{p2-fhtal11} such that $r_1> r'_1$, $r_2\geq r'_2$ and $r_3\geq r'_3$. Then, 
\begin{itemize}
\item[{\rm (i)}] There exists a rational number $\alpha\geq 0$ depending on $r_1$ and $r'_1$ 
but not on the morphisms $f$ and $\gamma$ satisfying the conditions of \ref{p2-fhtal1}, such that
\begin{equation}\label{p2-fhtal16a}
p^\alpha\tupiota^{\ur,\ur'}\colon  \tmK^\bullet(\fC'^{(r_1,r_2)}\hotimes_{R'_\uptau}\fC^{(r_3)}_\uptau)\rightarrow \tmK^\bullet(\fC'^{(r'_1,r'_2)}\hotimes_{R'_\uptau}\fC^{(r'_3)}_\uptau),
\end{equation}
where $\tupiota^{\ur,\ur'}$ is the morphism \eqref{p2-fhtal15h}, is homotopic to $0$ by an $\hoRp$-linear homotopy.
\item[{\rm (ii)}] The canonical morphism
\begin{equation}\label{p2-fhtal16b}
\tupiota^{\ur',\ur}\otimes_{\mZ_p}\mQ_p\colon  \tmK^\bullet(\fC'^{(r_1,r_2)}\hotimes_{R'_\uptau}\fC^{(r_3)}_\uptau)\otimes_{\mZ_p}\mQ_p
\rightarrow \tmK^\bullet(\fC'^{(r'_1,r'_2)}\hotimes_{R'_\uptau}\fC^{(r'_3)}_\uptau)\otimes_{\mZ_p}\mQ_p
\end{equation}
is homotopic to $0$ by a continuous homotopy.
\end{itemize}
\end{prop}

We choose a splitting $v$ of the exact sequence of $R'_\uptau$-modules \eqref{p2-fhtal2e}
\begin{equation}
\xymatrix{
0\ar[r]&{\Omega\otimes_{R_\uptau}R'_\uptau}\ar[r]^-(0.5)u&{\Omega'}\ar[r]&{\uOmega'}\ar[r]\ar@/_1pc/[l]_-(0.5)v&0.}
\end{equation}
Let $\Omega'_0$ and $\Omega'_1$ be the images of $u$ and $v$, respectively, 
so that we have $\Omega'=\Omega'_0\oplus \Omega'_1$. For any rational number $r\geq 0$, we identify the derivations $\delta_{\fC^{(r)}_\hoRp}$ \eqref{p2-fhtal7g} 
and $\delta_{\fC^{(r)}_\uptau}$ \eqref{p2-fhtal4dd} with respectively, 
\begin{eqnarray}
\delta^{(r)}\colon \fC^{(r)}_\hoRp\rightarrow \Omega'_0\otimes_{R'_\uptau}\fC^{(r)}_\hoRp,\\
\delta^{(r)}_\uptau\colon \fC^{(r)}_\uptau\rightarrow \Omega'_0\otimes_{R'_\uptau} \fC^{(r)}_\uptau.
\end{eqnarray}
For all rational numbers $r\geq r'\geq 0$, the derivation $\delta_{\fC'^{(r,r')}}$ \eqref{p2-fhtal8e} induces two derivations 
\begin{eqnarray}
\delta^{(r,r')}_0\colon \fC'^{(r,r')}\rightarrow \Omega'_0\otimes_{R'_\uptau} \fC'^{(r,r')},\\
\delta^{(r,r')}_1\colon \fC'^{(r,r')}\rightarrow \Omega'_1\otimes_{R'_\uptau} \fC'^{(r,r')}. 
\end{eqnarray}

We consider the Higgs $\hoRp$-fields 
\begin{eqnarray}
\hdelta^\ur_0=\delta^{(r_1,r_2)}_0\hotimes \id-\id\hotimes \delta^{(r_3)}_\uptau\colon \fC'^{(r_1,r_2)}\hotimes_{R'_\uptau}\fC^{(r_3)}_\uptau
&\rightarrow &\Omega'_0\otimes_{R'} \fC'^{(r_1,r_2)}\hotimes_{R'_\uptau}\fC^{(r_3)}_\uptau,\\
\hdelta^\ur_1=\delta^{(r_1,r_2)}_1\hotimes \id\colon \fC'^{(r_1,r_2)}\hotimes_{R'_\uptau}\fC^{(r_3)}_\uptau
&\rightarrow &\Omega'_1\otimes_{R'_\uptau} \fC'^{(r_1,r_2)}\hotimes_{R'_\uptau}\fC^{(r_3)}_\uptau,
\end{eqnarray}
so that $\hdelta^\ur$ \eqref{p2-fhtal15d} identifies with
\begin{equation}\label{p2-fhtal16g}
\hdelta^\ur_0\oplus \hdelta^\ur_1 \colon \fC'^{(r_1,r_2)}\hotimes_{R'_\uptau}\fC^{(r_3)}_\uptau
\rightarrow (\Omega'_0\oplus \Omega'_1)\otimes_{R'_\uptau} \fC'^{(r_1,r_2)}\hotimes_{R'_\uptau}\fC^{(r_3)}_\uptau.
\end{equation}
We then consider the double complex
\begin{equation}
\cK^{i,j}_\ur=\fC'^{(r_1,r_2)}\hotimes_{R'_\uptau}\fC^{(r_3)}_\uptau \otimes_{R'_\uptau} \wedge^i\Omega'_0\otimes_{R'_\uptau}\wedge^j\Omega'_1, 
\end{equation}
where the differentials are defined for any $\nu\in \fC'^{(r_1,r_2)}\hotimes_{R'_\uptau}\fC^{(r_3)}_\uptau$, 
$\omega_0\in \wedge^i\Omega'_0$ and  $\omega_1\in \wedge^j\Omega'_1$ by
\begin{eqnarray}
d_{\ur,1}^{i,j}(\nu\otimes \omega_0 \otimes \omega_1)&=&-\hdelta^\ur_0(\nu)\wedge \omega_0 \otimes \omega_1,\\
d_{\ur,2}^{i,j}(\nu\otimes \omega_0 \otimes \omega_1)&=&(-1)^{i+1}\hdelta^\ur_1(\nu)\wedge \omega_0 \otimes \omega_1,
\end{eqnarray}
the $\wedge$ product being taken in the left tensor $R'$-algebra $\wedge^\bullet \Omega'_0\ {^g\otimes}_{R'_\uptau} \wedge^\bullet \Omega'_1$ (see \cite{alg1-3} III § 4.7 remarks page 49).

Identifying $\Omega'_1$ and $\uOmega'$ by $v$, 
we have $\hdelta^\ur_1=\udelta_{\fC'^{(r_1,r_2)}}\hotimes \id$ \eqref{p2-fhtal8f}. 
Hence, for every integer $i\geq 0$, we have a canonical isomorphism of complexes 
\begin{equation}\label{p2-fhtal16d}
(\wedge^i\Omega'_0\otimes_{R'_\uptau} \umK^\bullet(\hfC'^{(r_1,r_2)})\hotimes_{R'_\uptau}\fC^{(r_3)}_\uptau,
- \id\otimes\udelta^\bullet_{\hfC'^{(r_1,r_2)}}\hotimes \id) 
\stackrel{\sim}{\rightarrow} (\cK^{i,\bullet}_\ur,d_{\ur,2}^{i,\bullet}), 
\end{equation}
where $\umK^\bullet(\hfC'^{(r_1,r_2)})$ is the Dolbeault complex of $(\hfC'^{(r_1,r_2)},\udelta_{\hfC'^{(r_1,r_2)}})$ \eqref{p2-fhtal8h}. 

We denote by $\fC^{(r_2)}_\hoRp\hotimes_{R'_\uptau}\fC^{(r_3)}_\uptau$ the $p$-adic completion of $\fC^{(r_2)}_\hoRp\otimes_{R'_\uptau}\fC^{(r_3)}_\uptau$, that we equip
with the Higgs $\hoRp$-field
\begin{equation}
\uplambda^{(r_2,r_3)}=\delta^{(r_2)}\hotimes \id-\id\hotimes \delta^{(r_3)}_\uptau\colon \fC^{(r_2)}_\hoRp\hotimes_{R'_\uptau}\fC^{(r_3)}_\uptau
\rightarrow \Omega'_0\otimes_{R'_\uptau} \fC^{(r_2)}_\hoRp\hotimes_{R'_\uptau}\fC^{(r_3)}_\uptau. 
\end{equation}

Since $\tf$ is smooth and $\tX'$ is affine, there exists $\tgamma\in \cL_{\tX'/\tX}(\coX')$ \eqref{p2-fhtal3a}.
In view of \ref{p2-hta21}(i) , we define a double complex by setting
\begin{eqnarray}
\tcK^{i,j}_\ur&=&\cK^{i,j}_\ur \ \ \ {\rm if} \  j\not=-1,\\
\tcK^{i,-1}_\ur&=&\fC^{(r_2)}_\hoRp\hotimes_{R'_\uptau}\fC^{(r_3)}_\uptau \otimes_{R'_\uptau} \wedge^i\Omega'_0,\label{p2-fhtal16ff}
\end{eqnarray}
where the differentials are defined for any integers $i,j$ with $j\not=-1$ and $\ell\in\{1,2\}$, by $\td_{\ur,\ell}^{i,j}=d_{\ur,\ell}^{i,j}$, and for 
any integer $i$ and any $\nu\in \fC^{(r_2)}_\hoRp\hotimes_{R'_\uptau}\fC^{(r_3)}_\uptau$ and $\omega_0\in \wedge^i\Omega'_0$, by
\begin{eqnarray}
\td_{\ur,1}^{i,-1}(\nu\otimes \omega_0)&=&\uplambda^{(r_2,r_3)}(\nu)\wedge \omega_0,\\
\td_{\ur,2}^{i,-1}(\nu\otimes \omega_0)&=&(-1)^{i}(\upphi_\tgamma^{(r_1,r_2)}\hotimes \id)(\nu)\otimes\omega_0,\label{p2-fhtal16f}
\end{eqnarray}
where $\upphi^{(r_1,r_2)}_\tgamma$ is the homomorphism \eqref{p2-fhtal9c}.

By \ref{p2-hta21}(ii), the isomorphism $\Lambda^{(r_2,r_3)}_{\tgamma}$ \eqref{p2-fhtal11g} is underlying an isomorphism of Higgs $\hoRp$-modules
\begin{equation}
(\fC^{(r_2)}_{\hoRp}\hotimes_{R'_\uptau}\fC^{(r_3)}_\uptau,\id\hotimes \delta^{(r_3)}_\uptau)\stackrel{\sim}{\rightarrow} 
(\fC^{(r_2)}_{\hoRp}\hotimes_{R'_\uptau}\fC^{(r_3)}_\uptau,-\uplambda^{r_2,r_3}). 
\end{equation}
We deduce an isomorphism of complexes 
\begin{equation}\label{p2-fhtal16c}
\fC^{(r_2)}_{\hoRp}\hotimes_{R'_\uptau} \mK^\bullet(\hfC^{(r_3)}_\uptau) \stackrel{\sim}{\rightarrow} (\tcK^{\bullet,-1}_\ur,-\td_{\ur,1}^{\bullet,-1}),
\end{equation}
where $\mK^\bullet(\hfC^{(r_3)}_\uptau)$ is the Dolbeault complex of $(\hfC^{(r_3)}_\uptau,\delta_{\hfC^{(r_3)}_\uptau})$ \eqref{p2-fhtal4g}, depending on the choice of $\tgamma$.

In view of \eqref{p2-fhtal8l} and \eqref{p2-fhtal4f}, the canonical homomorphisms $\alpha'^{r_1,r_2,r'_1,r'_2}$ \eqref{p2-fhtal8k} and $\alpha^{r_3,r'_3}_\uptau$ \eqref{p2-fhtal4k} 
induce a morphism of bicomplexes
\begin{equation}\label{p2-fhtal16e}
\jmath^{\ur',\ur}\colon \tcK^{\bullet,\bullet}(\fC'^{(r_1,r_2)}\hotimes_{R'_\uptau}\fC^{(r_3)}_\uptau)\rightarrow 
\tcK^{\bullet,\bullet}(\fC'^{(r'_1,r'_2)}\hotimes_{R'_\uptau}\fC^{(r'_3)}_\uptau).
\end{equation}

By \ref{p2-fhtal10}, there exist a rational number $\alpha$ and for every integer $j\geq -1$, an $\hoRp$-linear morphism
\begin{equation}
h^j_{\ur,\ur'}\colon \utmK^{j+1}_\tgamma(\hfC'^{(r_1,r_2)}) \rightarrow \utmK^j_\tgamma(\hfC'^{(r'_1,r'_2)})
\end{equation}
such that $(h^j_{\ur,\ur'})_{j\geq -1}$ define a homotopy from $0$ to $p^{\alpha_1} \tuupiota_{\tgamma}^{r_1,r_2,r'_1,r'_2}$ \eqref{p2-fhtal9f}. 
For any integers $i\geq 0$ and $j\geq -1$, in view of \eqref{p2-fhtal16d} and \eqref{p2-fhtal16ff}, we define an $\hoRp$-linear morphism
\begin{equation}
h^{i,j}_{\ur,\ur'}\colon \tcK^{i,j+1}_{\ur}\rightarrow \tcK^{i,j}_{\ur'}
\end{equation}
by 
\begin{eqnarray}
h^{i,j}_{\ur,\ur'}&=&(-1)^{i+1}\id_{\wedge^i\Omega'_0} \otimes_{R'_\uptau} h^j_{\ur,\ur'} \hotimes_{R'_\uptau}\id_{\fC^{(r_3)}_\uptau}, \ {\rm if}\ j\geq 0,\\
h^{i,-1}_{\ur,\ur'}&=&\id_{\wedge^i\Omega'_0} \otimes_{R'_\uptau} h^{-1}_{\ur,\ur'} \hotimes_{R'_\uptau}\id_{\fC^{(r_3)}_\uptau}.
\end{eqnarray}
We check from the proof of (\cite{ag2} 5.3.2) that for all integers $i\geq 0$ and $j\geq -1$, we have 
\begin{equation}
\td_{\ur',1}^{i,j}\circ h^{i,j}_{\ur,\ur'}+h^{i+1,j}_{\ur,\ur'}\circ\td_{\ur,1}^{i,j+1}=0.
\end{equation}

For any integer $n\geq 0$, the morphisms $h^{\bullet,\bullet}_{\ur,\ur'}$ induce an $\hoRp$-linear morphism
\begin{equation}\label{p2-fhtal16h}
t^n_{\ur,\ur'}=\oplus_{i+j=n} h^{i,j}_{\ur,\ur'} \colon \Tot^{n+1}(\tcK^{\bullet,\bullet}_{\ur})\rightarrow \Tot^{n}(\tcK^{\bullet,\bullet}_{\ur'}),
\end{equation}
where $\Tot^\bullet$ denotes the total complex (\cite{sp} \href{https://stacks.math.columbia.edu/tag/012Z}{012Z}). 
Taking into account the sign conventions in {\em loc. cit.}, we deduce from the above that $(t^n_{\ur,\ur'})_{n\geq -1}$ define a homotopy from $0$ to 
\begin{equation}
p^\alpha \Tot^\bullet(\jmath^{\ur,\ur'})\colon 
\Tot^\bullet(\tcK^{\bullet,\bullet}(\fC'^{(r_1,r_2)}\hotimes_{R'_\uptau}\fC^{(r_3)}_\uptau))\rightarrow 
\Tot^\bullet(\tcK^{\bullet,\bullet}(\fC'^{(r'_1,r'_2)}\hotimes_{R'_\uptau}\fC^{(r'_3)}_\uptau)),
\end{equation} 
where $\jmath^{\ur,\ur'}$ is the morphism \eqref{p2-fhtal16e}. 

By (\cite{alg1-3} III § 7.7 prop.~10), we have an isomorphism of bigraded algebras
\begin{equation}
\wedge^\bullet \Omega'_0\ {^g\otimes} \wedge^\bullet \Omega'_1 \stackrel{\sim}{\rightarrow}\wedge^\bullet\Omega'.
\end{equation}
Since $\hdelta^\ur=\hdelta^\ur_0\oplus \hdelta^\ur_1$ \eqref{p2-fhtal16g}, we have a canonical isomorphism 
\begin{equation}
\Tot^\bullet(\cK^{\bullet,\bullet}(\fC'^{(r_1,r_2)}\hotimes_{R'_\uptau}\fC^{(r_3)}_\uptau))\stackrel{\sim}{\rightarrow}
\mK^\bullet(\fC'^{(r_1,r_2)}\hotimes_{R'_\uptau}\fC^{(r_3)}_\uptau). 
\end{equation}
In view of \ref{p2-hta160}(iii), \eqref{p2-fhtal16c} and the sign conventions (\cite{sp} \href{https://stacks.math.columbia.edu/tag/09KA}{09KA}), we deduce an isomorphism 
\begin{equation}
\Tot^\bullet(\tcK^{\bullet,\bullet}(\fC'^{(r_1,r_2)}\hotimes_{R'_\uptau}\fC^{(r_3)}_\uptau))[1]\stackrel{\sim}{\rightarrow}
\tmK^\bullet(\fC'^{(r_1,r_2)}\hotimes_{R'_\uptau}\fC^{(r_3)}_\uptau).
\end{equation}
Moreover $\Tot^\bullet(\jmath^{\ur,\ur'})$ \eqref{p2-fhtal16e} corresponds to $\tupiota^{\ur',\ur}$ \eqref{p2-fhtal15h} under this isomorphism. 
The proposition follows then from the homotopy $(t^n_{\ur,\ur'})_{n\geq -1}$ \eqref{p2-fhtal16h}.

\begin{cor}\label{p2-fhtal17}
Let $\ut=(t_1,t_2,t_3)$ be a triple of rational numbers such that $t_1\geq t_2 > t_3 \geq 0$, $I_\ut$ the subset of elements $\ur=(r_1,r_2,r_3)$ of $I$ \eqref{p2-fhtal11} 
such that $r_1>t_1$, $r_2\geq t_2$  and $r_3>t_3$. Then, the morphism of complexes
\begin{equation}\label{p2-fhtal17a}
(\fC^{(t_2)}\hotimes_{\oR}\oR')\otimes_{\mZ_p}\mQ_p[0] \rightarrow  
\underset{\underset{(r_1,r_2,r_3)\in I_\ut}{\longrightarrow}}{\lim} \ \mK^\bullet(\fC'^{(r_1,r_2)}\hotimes_{R'_\uptau}\fC^{(r_3)}_\uptau)\otimes_{\mZ_p}\mQ_p
\end{equation}
induced by \eqref{p2-fhtal15e} is a quasi-isomorphism.
\end{cor}

Indeed, for every $\ur'=(r'_1,r'_2,r'_3)\in I_\ut$, it follows from \ref{p2-fhtal16}(ii) that the canonical morphism \eqref{p2-fhtal15e}
\begin{equation}
\fC^{(r'_2)}_{\hoRp}\hotimes_{R'_\uptau}\mK^\bullet(\fC^{(r'_3)}_\uptau)\otimes_{\mZ_p}\mQ_p
\rightarrow \underset{\underset{(r_1,r_2,r_3)\in J_{\ur'}}{\longrightarrow}}{\lim} \ \mK^\bullet(\fC'^{(r_1,r_2)}\hotimes_{R'_\uptau}\fC^{(r_3)}_\uptau)\otimes_{\mZ_p}\mQ_p,
\end{equation}
where $J_{\ur'}$ is the subset of elements $\ur=(r_1,r_2,r_3)$ of $I$ such that $r_1>r'_1$, $r_2\geq r'_2$ and $r_3\geq r'_3$, is a quasi-isomorphism. 
On the other hand, it follows from \ref{p2-fhtal5} that for all rational numbers $r_2,r'_2, r_3,r'_3$ such that $r_2= r'_2\geq 0$ and $r_3>r'_3>0$, 
the canonical morphism
\begin{equation}
\id\hotimes_{R'_\uptau} \tupiota^{r_3,r'_3}_\uptau\otimes_{\mZ_p}\mQ_p\colon 
\fC^{(r_2)}_{\hoRp}\hotimes_{R'_\uptau} \tmK^\bullet(\hfC^{(r_3)}_\uptau)\otimes_{\mZ_p}\mQ_p\rightarrow
\fC^{(r'_2)}_{\hoRp}\hotimes_{R'_\uptau} \tmK^\bullet(\hfC^{(r'_3)}_\uptau)\otimes_{\mZ_p}\mQ_p
\end{equation}
is homotopic to $0$ by a continuous homotopy. We deduce that the canonical morphism 
\begin{equation}
(\fC^{(t_2)}\hotimes_{\oR}\oR')\otimes_{\mZ_p}\mQ_p[0] \rightarrow  
\underset{\underset{r_2\in \mQ_{\geq t_2}, r_3\in \mQ_{>t_3}}{\longrightarrow}}{\lim} \ \fC^{(r_2)}_{\hoRp}\hotimes_{R'_\uptau} 
\mK^\bullet(\fC^{(r_3)}_\uptau)\otimes_{\mZ_p}\mQ_p
\end{equation}
is a quasi-isomorphism. The proposition follows.

\section{\texorpdfstring{Functoriality of the local $p$-adic Simpson correspondence by pullback}{Functoriality of the local p-adic Simpson correspondence by pullback}}\label{p2-nfspb}

\subsection{}\label{p2-fspb1}
The assumptions and notation of §\ref{p2-funchta} remain in force throughout this section.
We associate with $(f,\tf)$ \eqref{p2-fhtal1h} two functors \eqref{p2-rlps7b} and \eqref{p2-rlps8b}, 
\begin{eqnarray}
\mH\colon \bRep_{\hoR}(\Delta) \rightarrow \bHM(\hRun,\Omega),\label{p2-fspb1a}\\
\mV\colon \bHM(\hRun,\Omega)\rightarrow \bRep_{\hoR}(\Delta).\label{p2-fspb1b}
\end{eqnarray}
By \ref{p2-rlps11}, they induce equivalences of categories quasi-inverse to each other  
\begin{equation}\label{p2-fspb1c}
\xymatrix{
{\bRep^\Dolb_{\hoR[\frac 1 p]}(\Delta)}\ar@<1ex>[r]^-(0.5){\mH}&{\bHM^\sol(\hRun[\frac 1 p], \Omega).}
\ar@<1ex>[l]^-(0.5){\mV}}
\end{equation}
Similarly, we associate with $(f',\tf')$ \eqref{p2-fhtal1h} two functors 
\begin{eqnarray}
\mH'\colon \bRep_{\hoRp}(\Delta') \rightarrow \bHM(\hRunp,\Omega'),\label{p2-fspb1d}\\
\mV'\colon \bHM(\hRunp,\Omega')\rightarrow \bRep_{\hoRp}(\Delta'),\label{p2-fspb1e}
\end{eqnarray}
that induce equivalences of categories quasi-inverse to each other  
\begin{equation}\label{p2-fspb1f}
\xymatrix{
{\bRep^\Dolb_{\hoRp[\frac 1 p]}(\Delta')}\ar@<1ex>[r]^-(0.5){\mH'}&{\bHM^\sol(\hRunp[\frac 1 p], \Omega').}
\ar@<1ex>[l]^-(0.5){\mV'}}
\end{equation}

\subsection{}\label{p2-fspb2}
Consider the morphism of ringed punctual topos $(s,\hRunp)\rightarrow (s,\hRun)$ defined by the ring homomorphism
$\uplambda\colon \hRun\rightarrow \hRunp$ \eqref{p2-fhtal200a} induced by $\gamma$ \eqref{p2-fhtal1a},  
and the exact sequence of $\hRunp$-modules deduced from \eqref{p2-fhtal4b}
\begin{equation}\label{p2-fspb2b}
0\rightarrow \hRunp\rightarrow \fF_\uptau\otimes_{R'_\uptau}\hRunp\rightarrow \Omega\otimes_{R_\uptau}\hRunp \rightarrow 0. 
\end{equation}
For Higgs $\hRun$-modules with coefficients in $\Omega$ \eqref{p2-rlps2}, we defined in \eqref{p1-tphdi2b} the 
{\em base change by $\uplambda$ twisted by the extension $\fF_\uptau\otimes_{R'_\uptau}\hRunp$}:
\begin{equation}\label{p2-fspb2c}
\mT\colon \bHM(\hRun,\Omega) \rightarrow \bHM(\hRunp,\Omega).
\end{equation}
We set 
\begin{equation}\label{p2-fspb2e}
\mT'=\upmu\circ \mT\colon \bHM(\hRun,\Omega) \rightarrow \bHM(\hRunp,\Omega'), 
\end{equation} 
where
\begin{equation}\label{p2-fspb2d}
\upmu\colon \bHM(\hRunp,\Omega) \rightarrow \bHM(\hRunp,\Omega')
\end{equation} 
is the functor induced by the canonical morphism $u\colon \Omega\otimes_{R_\uptau}R'_\uptau\rightarrow \Omega'$ \eqref{p2-fhtal200e}.  

For a Higgs $\hRun[\frac 1 p]$-module with coefficients in $\Omega$, we introduced  in \ref{p1-tphdi3}
the properties of being {\em weakly twistable} (resp.\ {\em twistable}) {\em by the extension \eqref{p2-fspb2b}}.

\subsection{}\label{p2-fspb3}
We denote by $\fX$ (resp.\ $\fX'$) the formal scheme $p$-adic completion of $\coX$ (resp.\ $\coX'$) \eqref{p2-ncgt1a},
by $\fX^\star$ (resp.\ $\fX'^\star$) the open and closed formal subscheme of $\fX$ (resp.\ $\fX'$) corresponding to $\coX^\star$ \eqref{p2-rlps2} 
(resp.\ $\coX'^\star$ \eqref{p2-fhtal200}),
and by $\upgamma\colon \fX'\rightarrow \fX$ the morphism induced by the base change $\cogamma\colon \coX'\rightarrow \coX$ of $\gamma$ \eqref{p2-fhtal1a}. 
Observe that the formal schemes $\fX$ and $\fX'$ are of finite presentation over $\cS=\Spf(\co_C)$, 
globally idyllic (\cite{egr1} 2.6.12), and $\cS$-flat (\cite{egr1} 5.1.2). We set  \eqref{p2-fhtal1h}
\begin{eqnarray}
\tOmega^1_{\coX/\coS}&=&\tOmega^1_{(\coX,\cM_\coX)/(\coS,\cM_\coS)},\\
\tOmega^1_{\coX'/\coS}&=&\tOmega^1_{(\coX',\cM_{\coX'})/(\coS,\cM_\coS)},
\end{eqnarray}
and denote by $\tOmega^1_{\fX/\cS}$ and $\tOmega^1_{\fX'/\cS}$  their $p$-adic completions (\cite{egr1} 2.5.1). 
The $\co_\fX$-module $\txi^{-1}\tOmega^1_{\fX/\cS}$ \eqref{p2-ncgt3} is then canonically isomorphic to the $p$-adic completion 
of $\txi^{-1}\tOmega^1_{\coX/\coS}$; we identify these modules in the following, and similarly for $\fX'$.

We denote by $\cF_{\tX'/\tX}$ the Higgs--Tate extension of $\tX'$ over $\tX$ \eqref{p2-fhtal4}, 
and by $\hcF_{\tX'/\tX}$ its $p$-adic completion. 
We have a canonical exact sequence of $\co_{\coX'}$-modules
\begin{equation}\label{p2-fspb3a}
0\rightarrow \co_{\coX'}\rightarrow \cF_{\tX'/\tX}\rightarrow \cogamma^*(\txi^{-1}\tOmega^1_{\coX/\coS}) \rightarrow 0,
\end{equation} 
that induces an exact sequence of $\co_{\fX'}$-modules
\begin{equation}\label{p2-fspb3b}
0\rightarrow \co_{\fX'}\rightarrow \hcF_{\tX'/\tX}\rightarrow \upgamma^*(\txi^{-1}\tOmega^1_{\fX/\cS}) \rightarrow 0. 
\end{equation} 
For Higgs $\co_\fX$-modules with coefficients in $\txi^{-1}\tOmega^1_{\fX/\cS}$, 
we defined in \eqref{p1-tphdi2b} the {\em pullback by $\upgamma$ twisted by the extension $\hcF_{\tX'/\tX}$}:
\begin{equation}\label{p2-fspb3c}
\upgamma^*_\uptau\colon 
\bHM(\co_\fX,\txi^{-1}\tOmega^1_{\fX/\cS})\rightarrow \bHM(\co_{\fX'},\upgamma^*(\txi^{-1}\tOmega^1_{\fX/\cS})). 
\end{equation} 
We denote by 
\begin{equation}\label{p2-fspb3d}
\upgamma^+_\uptau\colon 
\bHM(\co_\fX,\txi^{-1}\tOmega^1_{\fX/\cS})\rightarrow \bHM(\co_{\fX'},\txi^{-1}\tOmega^1_{\fX'/\cS}) 
\end{equation} 
the functor induced by $\upgamma^*_\uptau$ and the canonical morphism $\upgamma^*(\tOmega^1_{\fX/\cS})\rightarrow \tOmega^1_{\fX'/\cS}$. 

For a Higgs $\co_\fX[\frac 1 p]$-module with coefficients in $\txi^{-1}\tOmega^1_{\fX/\cS}$, we introduced
the properties of being {\em weakly twistable} (resp.\ {\em twistable}) {\em by the extension $\hcF_{\tX'/\tX}$ \eqref{p2-fspb3b}} in \ref{p1-tphdi3}. 
These notions apply in particular to Higgs $\co_{\fX^\star}[\frac 1 p]$-modules with coefficients in $\txi^{-1}\tOmega^1_{\fX/\cS}$.

\subsection{}\label{p2-fspb4}
For a rational number $r\geq 0$, let $\cF_{\tX'/\tX}^{(r)}$ (resp.\ $\cC_{\tX'/\tX}^{(r)}$) be the Higgs--Tate extension (resp.\ algebra) of $\tX'$ over $\tX$ of thickness $r$ \eqref{p2-fhtal4}, 
$\hcF_{\tX'/\tX}^{(r)}$ (resp.\ $\hcC_{\tX'/\tX}^{(r)}$) its $p$-adic completion. 
Since the $\co_\coX$-module $\txi^{-1}\tOmega^1_{\coX/\coS}$ is $\co_C$-flat, we have a canonical exact sequence of free $\co_{\coX'}$-modules \eqref{p2-hta5d}
\begin{equation}\label{p2-fspb4a}
0\rightarrow \co_{\coX'}\rightarrow \cF^{(r)}_{\tX'/\tX}\rightarrow \cogamma^*(\txi^{-1}\tOmega^1_{\coX/\coS}) \rightarrow 0.
\end{equation} 
As explained in \ref{p2-hta71}, this extension is canonically isomorphic to the pullback of the extension $\cF_{\tX'/\tX}$ \eqref{p2-fspb3a} 
by the multiplication by $p^r$ on $\cogamma^*(\txi^{-1}\tOmega^1_{\coX/\coS})$. 
It induces by $p$-adic completion an exact sequence of $\co_{\fX'}$-modules
\begin{equation}\label{p2-fspb4b}
0\rightarrow \co_{\fX'}\rightarrow \hcF^{(r)}_{\tX'/\tX}\rightarrow \upgamma^*(\txi^{-1}\tOmega^1_{\fX/\cS}) \rightarrow 0.
\end{equation} 
By (\cite{egr1} 2.5.2(iii)), for every integer $n\geq 1$, we have a canonical isomorphism
\begin{equation}\label{p2-fspb4c}
\cF^{(r)}_{\tX'/\tX}/p^n\cF^{(r)}_{\tX'/\tX}\stackrel{\sim}{\rightarrow} \hcF^{(r)}_{\tX'/\tX}/p^n\hcF^{(r)}_{\tX'/\tX}. 
\end{equation}
We deduce that $\hcC_{\tX'/\tX}^{(r)}$ is canonically isomorphic to the $p$-adic completion of the $\co_{\fX'}$-algebra 
\begin{equation}\label{p2-fspb4d}
\underset{\underset{n\geq 0}{\longrightarrow}}\lim\ \rS^{n}_{\co_{\fX'}}(\hcF^{(r)}_{\tX'/\tX}),
\end{equation}
which is the period ring used for the definition of the functor $\upgamma^*_\uptau$ \eqref{p2-fspb3c}, see \ref{p1-thbn3}. 

We denote by  
\begin{equation}\label{p2-fspb4f}
\delta_{\cC^{(r)}_{\tX'/\tX}}\colon \cC^{(r)}_{\tX'/\tX}\rightarrow \cogamma^*(\txi^{-1}\tOmega^1_{\coX/\coS})\otimes_{\co_{\coX'}} \cC^{(r)}_{\tX'/\tX},
\end{equation}
the $\co_{\coX'}$-derivation defined in \eqref{p2-hta5h}, and by 
\begin{equation}\label{p2-fspb4g}
\delta_{\hcC^{(r)}_{\tX'/\tX}}\colon \hcC^{(r)}_{\tX'/\tX}\rightarrow \upgamma^*(\txi^{-1}\tOmega^1_{\fX/\cS})\otimes_{\co_{\fX'}} \hcC^{(r)}_{\tX'/\tX}
\end{equation}
its extension to the $p$-adic completions \eqref{p1-thbn18}.

\begin{prop}\label{p2-fspb5}
Let $(\cN,\theta)$ be a Higgs $\co_{\fX^\star}[\frac 1 p]$-bundle with coefficients in $\txi^{-1}\tOmega^1_{\fX/\cS}$ \eqref{p1-delta-con6}, 
$N=\Gamma(\fX^\star,\cN)$. We denote also by $\theta$ the Higgs $\hRun$-field with coefficients in $\Omega$ induced by $\theta$ on $N$. 
Then, the following conditions are equivalent. 
\begin{itemize}
\item[{\rm (i)}] The Higgs $\co_{\fX^\star}[\frac 1 p]$-bundle $(\cN,\theta)$ is twistable by the extension $\hcF_{\tX'/\tX}$ \eqref{p2-fspb3b}. 
\item[{\rm (ii)}] The Higgs $\co_{\fX^\star}[\frac 1 p]$-bundle $(\cN,\theta)$ is weakly twistable by the extension $\hcF_{\tX'/\tX}$. 
\item[{\rm (iii)}] The Higgs $\hRun[\frac 1 p]$-bundle $(N,\theta)$ is twistable by the extension $\fF_\uptau\otimes_{R'_\uptau}\hRunp$ \eqref{p2-fspb2b}. 
\item[{\rm (iv)}] The Higgs $\hRun[\frac 1 p]$-bundle $(N,\theta)$ is weakly twistable by the extension $\fF_\uptau\otimes_{R'_\uptau}\hRunp$. 
\end{itemize}
Moreover, if these conditions are satisfied, we have a functorial canonical isomorphism 
\begin{equation}\label{p2-fspb5j}
\mT(N,\theta)\stackrel{\sim}{\rightarrow} \Gamma(\fX'^\star,\upgamma^*_\uptau(\cN,\theta)),
\end{equation}
where $\mT$ and $\upgamma^*_\uptau$ are defined in \eqref{p2-fspb2c} and \eqref{p2-fspb3c}, respectively.
\end{prop}

Observe first that the $\hRun[\frac 1 p]$-module $N$ is projective of finite type by (\cite{agt} III.6.17). 
The equivalences (i)$\Leftrightarrow$(ii) and (iii)$\Leftrightarrow$(iv) follow from \ref{p1-thbn21}. We prove the equivalence of (i) and (iii). 
By (\cite{agt} III.6.16), we have an equivalence of monoidal abelian categories
\begin{equation}\label{p2-fspb5a}
\begin{array}[t]{clcr}
\bMod^{\coh}_\mQ(\co_{\fX^\star})&\rightarrow& \bMod^{\coh}(\co_{\fX^\star}[\frac 1 p]),\\
\cN&\mapsto& \cN[\frac 1 p], 
\end{array}
\end{equation}
where the source denotes the category of coherent $\co_{\fX^\star}$-modules up to isogeny (\cite{ag2} 2.9.1).
Recall that $\hRun$ and hence $\hRun[\frac 1 p]$ are coherent rings (\cite{egr1} 1.10.3). 
By \eqref{p2-fspb5a} and (\cite{egr1} 2.7.2, (2.7.2.4) and 2.10.5), we have an equivalence of monoidal abelian categories 
\begin{equation}\label{p2-fspb5b}
\Gamma(\fX^\star,-)\colon \bMod^{\coh}(\co_{\fX^\star}[\frac 1 p])\stackrel{\sim}{\rightarrow} \bMod^\coh(\hRun[\frac 1 p]).
\end{equation}
Similarly, we have an equivalence of monoidal abelian categories
\begin{equation}\label{p2-fspb5c}
\Gamma(\fX'^\star,-)\colon \bMod^{\coh}(\co_{\fX'^\star}[\frac 1 p])\stackrel{\sim}{\rightarrow} \bMod^\coh(\hRunp[\frac 1 p]).
\end{equation}

Let $r$ be a rational number $\geq 0$. The $\hRunp$-algebra $\hfC^{(r)}_\uptau$ \eqref{p2-fhtal4c} is topologically of finite presentation, and hence idyllic (\cite{egr1} 1.10.1). 
We set $\bL^{(r)}=\Spf(\hfC^{(r)}_\uptau)$ and denote by $\uppi^{(r)}\colon \bL^{(r)} \rightarrow \fX'^\star$ the canonical morphism. 
Similarly as in \eqref{p2-fspb5b}, we have an equivalence of monoidal abelian categories
\begin{equation}\label{p2-fspb5d}
\Gamma(\bL^{(r)},-)\colon \bMod^{\coh}(\co_{\bL^{(r)}}[\frac 1 p])\stackrel{\sim}{\rightarrow} \bMod^\coh(\hfC^{(r)}_\uptau[\frac 1 p]).
\end{equation}
On the other hand, it follows from \ref{p2-fspb4} and \ref{p1-tshbn2} that the $\co_{\fX'}$-algebra $\hcC^{(r)}_{\tX'/\tX}$ is topologically of finite presentation,
and we have $\bL^{(r)}=\Spf(\hcC^{(r)}_{\tX'/\tX}|\fX'^\star)$ (see \ref{p1-pfs4}). 
We deduce from \eqref{p2-fspb5d} and \ref{p1-pfs12}(v) that we have an equivalence of monoidal abelian categories 
\begin{equation}\label{p2-fspb5e}
\Gamma(\fX'^\star,-)\colon \bMod^\coh(\hcC^{(r)}_{\tX'/\tX}[\frac 1 p]|\fX'^\star)\stackrel{\sim}{\rightarrow} \bMod^\coh(\hfC^{(r)}_\uptau[\frac 1 p]).
\end{equation}
We have $\delta_{\hfC^{(r)}_\uptau}=\Gamma(\fX'^\star,\delta_{\hcC^{(r)}_{\tX'/\tX}})$, where these derivations are defined in \eqref{p2-fhtal4g} and \eqref{p2-fspb4g}. 

Let $\cM$ be a coherent $\co_{\fX'^\star}[\frac 1 p]$-module, $M=\Gamma(\fX'^\star,\cM)$. The canonical morphism
\begin{equation}\label{p2-fspb5f}
M\otimes_{\hRunp}\hfC^{(r)}_\uptau\rightarrow \Gamma(\bL^{(r)},\uppi^{(r)*}(\cM))
\end{equation}
is an isomorphism. It can easily been proved by taking a resolution $\co_{\fX'^\star}[\frac 1 p]^n\rightarrow \co_{\fX'^\star}[\frac 1 p]^m\rightarrow \cM\rightarrow 0$ \eqref{p2-fspb5c}. 
We deduce by \ref{p1-pfs12}(vi) that he canonical morphism
\begin{equation}\label{p2-fspb5g}
M\otimes_{\hRunp}\hfC^{(r)}_\uptau\rightarrow \Gamma(\fX'^\star,\cM\otimes_{\co_{\fX'}}\hcC^{(r)}_{\tX'/\tX})
\end{equation}
is an isomorphism. 

For every $\hcC^{(r)}_{\tX'/\tX}$-linear morphism
\begin{equation}\label{p2-fspb5h}
u\colon \cM\otimes_{\co_{\fX'}}\hcC^{(r)}_{\tX'/\tX}\rightarrow \upgamma^*(\cN)\otimes_{\co_{\fX'}}\hcC^{(r)}_{\tX'/\tX},
\end{equation}
setting 
\begin{equation}\label{p2-fspb5i}
v=\Gamma(\fX'^\star,u)\colon M\otimes_{\hRunp}\hfC^{(r)}_\uptau\rightarrow N\otimes_{\hRun}\hfC^{(r)}_\uptau,
\end{equation}
the following properties are equivalent:
\begin{itemize}
\item[(a)] $u$ is a morphism of $\hcC^{(r)}_{\tX'/\tX}$-modules with $\delta_{\hcC^{(r)}_{\tX'/\tX}}$-connection, in the sense of \ref{p1-delta-con2},
where the $\delta_{\hcC^{(r)}_{\tX'/\tX}}$-connections are defined as in \ref{p1-delta-con4},
$\upgamma^*(\cN)$ (resp.\ $\cM$) being endowed with the Higgs field $\upgamma^*(\theta)$ (resp.\ zero).
\item[(b)] $u(\cM)\subset \ker(\theta^{(r)}_\tot)$, where $\theta^{(r)}_\tot= \upgamma^*(\theta)\otimes \id+\id\otimes \delta_{\hcC^{(r)}_{\tX'/\tX}}$
is the total Higgs $\co_{\fX'^\star}$-field on $\upgamma^*(\cN)\otimes_{\co_{\fX'}}\hcC^{(r)}_{\tX'/\tX}$.
\item[(c)] $v$ is a morphism of $\hfC^{(r)}_\uptau$-modules with $\delta_{\hfC^{(r)}_\uptau}$-connection, 
where the $\delta_{\hfC^{(r)}_\uptau}$-connections are defined as in \ref{p1-delta-con4},
$N$ (resp.\ $M$) being endowed with the Higgs field $\theta$ (resp.\ zero).
\item[(d)] $v(M)\subset \ker(\theta^{(r)}_\tot)$, where $\theta^{(r)}_\tot= \theta\otimes \id+\id\otimes \delta_{\hfC^{(r)}_\uptau}$ 
is the total Higgs $\hRunp$-field on $N\otimes_{\hRun}\hfC^{(r)}_\uptau$.
\end{itemize}
Indeed, we clearly have (a)$\Leftrightarrow$(b), (c)$\Leftrightarrow$(d) and (b)$\Leftrightarrow$(d) by (\cite{egr1} 2.7.6). 
The equivalence (a)$\Leftrightarrow$(c) implies the equivalence (i)$\Leftrightarrow$(iii). 

We suppose next that the equivalent conditions (i)-(iv) are satisfied. We write $\upgamma^*_\uptau(\cN,\theta)=(\cM,\theta^\vee)$, 
which is a Higgs $\co_{\fX'^\star}[\frac 1 p]$-bundle with coefficients in $\txi^{-1}\tOmega^1_{\fX/\cS}$ by \ref{p1-thbn32}. 
By \ref{p1-thbn23}, there exists a $\hcC^{(r)}_{\tX'/\tX}$-linear isomorphism $u$ \eqref{p2-fspb5h} satisfying conditions (a). 
The Higgs field $\theta^\vee$ on $\cM$ induces a Higgs $\hRunp$-field on $M=\Gamma(\fX'^\star,\cM)$ that we denote also by $\theta^\vee$. 
It follows then from \ref{p1-thbn25}, applied first to $u$ and then to $v$ \eqref{p2-fspb5i}, that $v$ induces an isomorphism 
\begin{equation}
\mT(N,\theta)\stackrel{\sim}{\rightarrow} (M,\theta^\vee),
\end{equation}
which is clearly functorial in $(\cN,\theta)$, providing the isomorphism \eqref{p2-fspb5j}.

\begin{prop}[\ref{p1-thbn25}]\label{p2-fspb6}
Let $(N,\theta)$ be a Higgs $\hRun[\frac 1 p]$-bundle with coefficients in $\Omega$ \eqref{p1-delta-con6}, 
$N^\vee$ an $\hRunp$-module, $r$ a rational number $>0$,
\begin{equation}\label{p2-fspb6a}
N^\vee\otimes_{\hRunp}\hfC^{(r)}_\uptau\stackrel{\sim}{\rightarrow}N\otimes_{\hRun}\hfC^{(r)}_\uptau
\end{equation} 
an isomorphism of $\hfC^{(r)}_\uptau$-modules with $\delta_{\hfC^{(r)}_\uptau}$-connection, 
where the $\delta_{\hfC^{(r)}_\uptau}$-connections are defined as in \ref{p1-delta-con4}, 
$N$ (resp.\ $N^\vee$) being endowed with the Higgs field $\theta$ (resp.\ $0$).  
Then,
\begin{itemize}
\item[{\rm (i)}] The isomorphism \eqref{p2-fspb6a} induces an $\hRunp$-linear isomorphism
\begin{equation}\label{p2-fspb6b}
N^\vee \stackrel{\sim}{\rightarrow} (N\otimes_{\hRun}\hfC^{(r)}_\uptau)^{\theta^{(r)}_\tot=0},
\end{equation}
where $\theta^{(r)}_\tot=\theta\otimes \id+\id\otimes\delta_{\hfC^{(r)}_\uptau}$ is the total Higgs $\hRunp$-field on $N\otimes_{\hRun}\hfC^{(r)}_\uptau$. 
\item[{\rm (ii)}]  The Higgs field $\theta\otimes \id$ on $N\otimes_{\hRun}\hfC^{(r)}_\uptau$ induces a Higgs $\hRunp$-field $\theta^\vee$
on $N^\vee$ with coefficients in $\Omega$. The isomorphism \eqref{p2-fspb6a} induces an isomorphism of Higgs bundles
\begin{equation}\label{p2-fspb6c}
(N^\vee,\theta^\vee)\stackrel{\sim}{\rightarrow} \mT(N,\theta),
\end{equation}
where $\mT$ is defined in \eqref{p2-fspb2c}. 
\item[{\rm (iii)}] Setting $\delta^\vee_{\hfC^{(r)}_\uptau}=-\delta_{\hfC^{(r)}_\uptau}$, 
the morphism \eqref{p2-fspb6a} is an isomorphism of $\hfC^{(r)}_\uptau$-modules with $\delta^\vee_{\hfC^{(r)}_\uptau}$-connection, 
where the $\delta^\vee_{\hfC^{(r)}_\uptau}$-connections are defined as in \ref{p1-delta-con4}, 
$N$ (resp.\ $N^\vee$) being endowed with the Higgs field $0$ (resp.\ $\theta^\vee$).  
\end{itemize}
\end{prop}

\begin{prop}\label{p2-fspb7}
Every solvable Higgs $\hRun[\frac 1 p]$-bundle with coefficients in $\Omega$ in the sense of \ref{p2-rlps10}, 
is twistable by the extension $\fF_\uptau\otimes_{R'_\uptau}\hRunp$ \eqref{p2-fspb2b}. 
\end{prop}

Indeed, such a Higgs $\hRun[\frac 1 p]$-bundle is CL-small by (\cite{ag2} 3.4.30), see \ref{p2-rlps26}. 
It is therefore twistable by the extension $\fF_\uptau\otimes_{R'_\uptau}\hRunp$ by \ref{p1-tshbn15}.  

\begin{prop}\label{p2-fspb8}
Every strongly solvable Higgs $\co_{\fX}[\frac 1 p]$-bundle with coefficients in $\txi^{-1}\tOmega^1_{\fX/\cS}$ 
in the sense of {\rm (\cite{ag2} 4.6.6(iii))}, is twistable by the extension $\hcF_{\tX'/\tX}$ \eqref{p2-fspb3b}. 
\end{prop}

Indeed, such a Higgs $\co_\fX[\frac 1 p]$-bundle is CL-small by (\cite{ag2} 4.8.15).
It is therefore twistable by the extension $\hcF_{\tX'/\tX}$ by \ref{p1-tshbn15}.

\begin{teo}\label{p2-fspb9}
Let $M$ be a Dolbeault $\hoR[\frac 1 p]$-representation of $\Delta$ \eqref{p2-rlps9}. 
We set $M'=M\otimes_\hoR \hoRp$ that we equip with the 
$\hoRp$-representation of $\Delta'$ induced by the canonical homomorphism $\Delta'\rightarrow \Delta$. 
Then, $M'$ is a Dolbeault $\hoRp[\frac 1 p]$-representation of $\Delta'$, and we have a functorial canonical isomorphism 
\begin{equation}\label{p2-fspb9a}
\mH'(M')\stackrel{\sim}{\rightarrow}\mT'(\mH(M)),
\end{equation}
where the functors above are defined in \eqref{p2-fspb1a}, \eqref{p2-fspb1d} and \eqref{p2-fspb2e}. 
\end{teo}

We set $\mH(M)=(N,\theta)$, which is a Higgs $\hRun[\frac 1 p]$-bundle with coefficients in $\Omega$ \eqref{p1-delta-con6}. 
By \ref{p2-rlps14}, there exists a rational number $r>0$ such that 
the canonical morphism $N\rightarrow M\otimes_{\hoR}\fC^\dagger$ induces a $\Delta$-equivariant 
isomorphism of $\hfC^{(r)}$-modules with $\delta_{\hfC^{(r)}}$-connection \eqref{p2-rlps5g}, 
\begin{equation}\label{p2-fspb9b}
N \otimes_{\hRun}\hfC^{(r)}\stackrel{\sim}{\rightarrow} M\otimes_{\hoR}\hfC^{(r)},
\end{equation}
where $N$ (resp.\ $\hfC^{(r)}$) is endowed with the trivial (resp.\ canonical) action of $\Delta$,  
and the $\delta_{\hfC^{(r)}}$-connections are defined as in \ref{p1-delta-con4},
$N$ (resp.\ $M$) being endowed with the Higgs field $\theta$ (resp.\ zero).

For any rational number $t$ such that $0<t\leq r$, we consider the canonical homomorphism \eqref{p2-fhtal110b}
\begin{equation}\label{p2-fspb9c}
\upphi^{(r,t)}\colon \fC^{(r)}_{\hoRp}\rightarrow \fC'^{(r)}\otimes_{R'_\uptau}\fC^{(t)}_\uptau,
\end{equation}
and denote by $\hupphi^{(r,t)}$ its extension to the $p$-adic completions. It is $\Delta'$-equivariant by \ref{p2-fhtal12}. 
To lighten the notation, we denote by
\begin{eqnarray}
\delta^{(r)}=\delta_{\fC^{(r)}_\hoRp}\colon \fC^{(r)}_\hoRp\rightarrow \Omega\otimes_{R_\uptau}\fC^{(r)}_\hoRp,\label{p2-fspb9d1}\\
\delta^{(t)}_\uptau=\delta_{\fC^{(t)}_\uptau}\colon \fC^{(t)}_\uptau\rightarrow \Omega\otimes_{R_\uptau} \fC^{(t)}_\uptau,\label{p2-fspb9d2}\\
\delta'^{(r)}=\delta_{\fC'^{(r)}}\colon \fC'^{(r)}\rightarrow \Omega'\otimes_{R'_\uptau}\fC'^{(r)},\label{p2-fspb9d4}
\end{eqnarray}
the derivations defined in \eqref{p2-fhtal7g}, \eqref{p2-fhtal4dd} and the analogue of \eqref{p2-rlps5e}, respectively, and we set 
\begin{eqnarray}
\delta^{\vee(t)}_\uptau=-\delta^{(t)}_\uptau\colon \fC^{(t)}_\uptau\rightarrow \Omega\otimes_{R_\uptau} \fC^{(t)}_\uptau,\label{p2-fspb9d2a}\\
\delta'^{(t)}_\uptau=(u\otimes\id)\circ\delta^{(t)}_\uptau\colon \fC^{(t)}_\uptau\rightarrow \Omega'\otimes_{R_\uptau} \fC^{(t)}_\uptau,\label{p2-fspb9d3}\\
\delta'^{\vee(t)}_\uptau=-\delta'^{(t)}_\uptau\colon \fC^{(t)}_\uptau\rightarrow \Omega'\otimes_{R_\uptau} \fC^{(t)}_\uptau,\label{p2-fspb9d3a}
\end{eqnarray}
where $u\colon \Omega\otimes_{R_\uptau}R'_\uptau\rightarrow \Omega'$ is the canonical morphism \eqref{p2-fhtal200e}.
In the notation above, a prime exponent indicates composition with $u$, while a $\vee$ exponent denotes multiplication by $-1$.
We consider also the $\hoRp$-derivations 
\begin{eqnarray}
\delta'^{(r,t)}=\delta'^{(r)}\otimes \id+\id \otimes \delta'^{\vee(t)}_\uptau \colon \fC'^{(r)}\otimes_{R'_\uptau}\fC^{(t)}_\uptau \rightarrow 
\Omega'\otimes_{R'_\uptau}\fC'^{(r)}\otimes_{R'_\uptau}\fC^{(t)}_\uptau,\label{p2-fspb9e1}\\
\delta^{(r,t)}=\id \otimes \delta^{(t)}_\uptau \colon \fC'^{(r)}\otimes_{R'_\uptau}\fC^{(t)}_\uptau \rightarrow 
\Omega\otimes_{R_\uptau}\fC'^{(r)}\otimes_{R'_\uptau}\fC^{(t)}_\uptau.\label{p2-fspb9e2}
\end{eqnarray}
We denote the extensions to the $p$-adic completions by a $\widehat{\ }$\ .

The following diagrams 
\begin{equation}\label{p2-fspb9f}
\xymatrix{
{\fC^{(t)}_\uptau}\ar[r]^-(0.5){\delta^{(t)}_\uptau}\ar[d]&{\Omega\otimes_{R_\uptau} \fC^{(t)}_\uptau}\ar[d]\\
{\fC'^{(r)}\otimes_{R'_\uptau}\fC^{(t)}_\uptau}\ar[r]^-(0.5){\delta^{(r,t)}}&{\Omega\otimes_{R_\uptau}\fC'^{(r)}\otimes_{R'_\uptau}\fC^{(t)}_\uptau,}}
\ \ \ \ 
\xymatrix{
{\fC^{(t)}_\uptau}\ar[r]^-(0.5){\delta'^{\vee(t)}_\uptau}\ar[d]&{\Omega'\otimes_{R'_\uptau} \fC^{(t)}_\uptau}\ar[d]\\
{\fC'^{(r)}\otimes_{R'_\uptau}\fC^{(t)}_\uptau}\ar[r]^-(0.5){\delta'^{(r,t)}}&{\Omega'\otimes_{R'_\uptau}\fC'^{(r)}\otimes_{R'_\uptau}\fC^{(t)}_\uptau,}}
\end{equation}
where the vertical arrows are the canonical morphisms, are commutative. 
The following diagram 
\begin{equation}\label{p2-fspb9g}
\xymatrix{
{\fC'^{(r)}}\ar[r]^-(0.5){\delta'^{(r)}}\ar[d]&{\Omega'\otimes_{R'_\uptau} \fC'^{(r)}}\ar[d]\\
{\fC'^{(r)}\otimes_{R'_\uptau}\fC^{(t)}_\uptau}\ar[r]^-(0.5){\delta'^{(r,t)}}&{\Omega'\otimes_{R'_\uptau}\fC'^{(r)}\otimes_{R'_\uptau}\fC^{(t)}_\uptau,}}
\end{equation}
where the vertical arrows are the canonical morphisms, is commutative. Moreover, the composition 
\begin{equation}\label{p2-fspb9h}
\xymatrix{
{\fC'^{(r)}}\ar[r]&{\fC'^{(r)}\otimes_{R'_\uptau}\fC^{(t)}_\uptau}
\ar[r]^-(0.5){\delta^{(r,t)}}&{\Omega\otimes_{R_\uptau}\fC'^{(r)}\otimes_{R'_\uptau}\fC^{(t)}_\uptau}}
\end{equation}
where the first arrow is the canonical morphism, vanishes. 
By \ref{p2-hta190}, the following diagram 
\begin{equation}\label{p2-fspb9i}
\xymatrix{
{\fC^{(r)}_{\hoRp}}\ar[r]^-(0.5){\delta^{(r)}}\ar[d]_{\upphi^{(r,t)}}&{\Omega\otimes_{R_\uptau} \fC^{(r)}_{\hoRp}}\ar[d]^{\id \otimes \upphi^{(r,t)}}\\
{\fC'^{(r)}\otimes_{R'_\uptau}\fC^{(t)}_\uptau}\ar[r]^-(0.5){\delta^{(r,t)}}&{\Omega\otimes_{R_\uptau}\fC'^{(r)}\otimes_{R'_\uptau}\fC^{(t)}_\uptau}}
\end{equation}
is commutative. Moreover, the composition 
\begin{equation}\label{p2-fspb9j}
\xymatrix{
{\fC^{(r)}_\hoRp}\ar[r]^-(0.5){\upphi^{(r,t)}}&{\fC'^{(r)}\otimes_{R'_\uptau}\fC^{(t)}_\uptau}
\ar[r]^-(0.5){\delta'^{(r,t)}}&{\Omega'\otimes_{R'_\uptau}\fC'^{(r)}\otimes_{R'_\uptau}\fC^{(t)}_\uptau}}
\end{equation}
vanishes.

By \ref{p1-delta-con5} and \eqref{p2-fspb9i}, the isomorphism \eqref{p2-fspb9b} 
induces by extension of scalars by the homomorphism induced by $\upphi^{(r,t)}$ \eqref{p2-fspb9c} an isomorphism of $\fC'^{(r)}\hotimes_{R'_\uptau}\fC^{(t)}_\uptau$-modules
with $\hdelta^{(r,t)}$-connection 
\begin{equation}\label{p2-fspb9k}
N \otimes_{\hRun}\fC'^{(r)}\hotimes_{R'_\uptau}\fC^{(t)}_\uptau\stackrel{\sim}{\rightarrow} M\otimes_{\hoR}\fC'^{(r)}\hotimes_{R'_\uptau}\fC^{(t)}_\uptau, 
\end{equation}
where the tensor product $\hotimes$ is completed relatively to the $p$-adic topology and the $\hdelta^{(r,t)}$-connections are defined as in \ref{p1-delta-con4}, 
$N$ (resp.\ $M$) being endowed with the Higgs field $\theta$ (resp.\ zero).
Moreover, it follows from \eqref{p2-fspb9j} that \eqref{p2-fspb9k} is an isomorphism of $\fC'^{(r)}\hotimes_{R'_\uptau}\fC^{(t)}_\uptau$-modules
with $\hdelta'^{(r,t)}$-connection, when  $N$ and $M$ are endowed with the zero Higgs fields. 
To summarize, \eqref{p2-fspb9k} is an isomorphism
of Higgs $\hoRp$-modules with the Higgs fields indicated on the same line of this table
\begin{equation}\label{p2-fspb9l}
\begin{tabular}{|c|c|}
\hline
left hand side&right  hand  side\\
\hline
$\theta\otimes \id +\id\otimes \hdelta^{(r,t)}$ & $\id\otimes \hdelta^{(r,t)}$\\
\hline
$\id\otimes \hdelta'^{(r,t)}$ & $\id\otimes \hdelta'^{(r,t)}$\\
\hline
\end{tabular}
\end{equation}

On the other hand, by \ref{p2-fspb7}, the Higgs $\hRun[\frac 1 p]$-bundle $(N,\theta)$ is twistable by the extension $\fF_\uptau\otimes_{R'_\uptau}\hRunp$ \eqref{p2-fspb2b}. 
We set $\mT(N,\theta)=(N^\vee,\theta^\vee)$, where $\mT$ is the functor \eqref{p2-fspb2c}. 
Then, there exists a rational number $t'>0$ such that the canonical morphism $N' \rightarrow N\otimes_{\hRun}\fC^\dagger_\uptau$ 
induces an isomorphism of $\hfC^{(t')}_\uptau$-modules with $\hdelta^{(t')}_\uptau$-connection \eqref{p2-fspb9d2}
\begin{equation}\label{p2-fspb9m}
N^\vee \otimes_{\hRunp}\hfC^{(t')}_\uptau\stackrel{\sim}{\rightarrow} N\otimes_{\hRun}\hfC^{(t')}_\uptau,
\end{equation}
where the $\hdelta^{(t')}_{\uptau}$-connections are defined as in \ref{p1-delta-con4}, 
$N$ (resp.\ $N^\vee$) being endowed with the Higgs field $\theta$ (resp.\ $0$). Moreover, by \ref{p2-fspb6}, 
\eqref{p2-fspb9m} is also an isomorphism of $\hfC^{(t')}_\uptau$-modules with $\hdelta^{\vee(t')}_\uptau$-connection \eqref{p2-fspb9d2a}, 
where the $\hdelta^{\vee(t')}_{\uptau}$-connections are defined as in \ref{p1-delta-con4}, 
$N$ (resp.\ $N^\vee$) being endowed with the Higgs field $0$ (resp.\ $\theta^\vee$).
Observe that the same properties hold for every rational number $t''$ such that $0\leq t''\leq t'$. 

Let $t$ be a rational number such that $0<t\leq \inf(t',r)$. 
By \ref{p1-delta-con5} and \eqref{p2-fspb9f}, the isomorphism \eqref{p2-fspb9m} induces by extension of scalars an isomorphism of 
$\fC'^{(r)}\hotimes_{R'_\uptau}\fC^{(t)}_\uptau$-modules with $\hdelta^{(r,t)}$-connection 
\begin{equation}\label{p2-fspb9p}
N^\vee\otimes_{\hRunp}\fC'^{(r)}\hotimes_{R'_\uptau}\fC^{(t)}_\uptau\stackrel{\sim}{\rightarrow} N\otimes_{\hRun}\fC'^{(r)}\hotimes_{R'_\uptau}\fC^{(t)}_\uptau,
\end{equation}
where the $\hdelta^{(r,t)}$-connections are defined as in \ref{p1-delta-con4}, 
$N$ (resp.\ $N^\vee$) being endowed with the Higgs field $\theta$ (resp.\ $0$). 
It is also an isomorphism of $\fC'^{(r)}\hotimes_{R'_\uptau}\fC^{(t)}_\uptau$-modules with $\hdelta'^{(r,t)}$-connection,  
where the $\hdelta'^{(r,t)}$-connections are defined as in \ref{p1-delta-con4}, 
$N$ (resp.\ $N^\vee$) being endowed with the Higgs field $0$ 
(resp.\ $\theta'^\vee=(u\otimes\id)\circ\theta^\vee\colon N^\vee\rightarrow \Omega'\otimes_{R'_\uptau}N^\vee$). 
To summarize, \eqref{p2-fspb9p} is an isomorphism of Higgs $\hoRp$-modules with the Higgs fields indicated on the same line of this table
\begin{equation}\label{p2-fspb9q}
\begin{tabular}{|c|c|}
\hline
left hand side&right  hand  side\\
\hline
$\id\otimes \hdelta^{(r,t)}$ & $\theta\otimes \id+\id\otimes\hdelta^{(r,t)}$\\
\hline
$\theta'^\vee\otimes \id+\id\otimes \hdelta'^{(r,t)}$ & $\id\otimes \hdelta'^{(r,t)} $\\
\hline
\end{tabular}
\end{equation}

We set $(N',\theta')=(N^\vee,\theta'^\vee)$ and $M'=M\otimes_{\hoR}\hoRp$.
Composing \eqref{p2-fspb9k} and \eqref{p2-fspb9p}, we deduce an isomorphism of $\fC'^{(r)}\hotimes_{R'_\uptau}\fC^{(t)}_\uptau$-modules 
\begin{equation}\label{p2-fspb9r}
N' \otimes_{\hRunp}\fC'^{(r)}\hotimes_{R'_\uptau}\fC^{(t)}_\uptau\stackrel{\sim}{\rightarrow} M'\otimes_{\hoRp}\fC'^{(r)}\hotimes_{R'_\uptau}\fC^{(t)}_\uptau,
\end{equation}
which is compatible with the Higgs $\hRunp$-fields indicated on the same line of the table
\begin{equation}\label{p2-fspb9s}
\begin{tabular}{|c|c|}
\hline
left hand side&  right  hand  side\\
\hline
$\id\otimes \hdelta^{(r,t)}$ &  $\id\otimes \hdelta^{(r,t)}$\\
\hline
$\theta'\otimes \id+\id\otimes \hdelta'^{(r,t)}$ & $\id\otimes \hdelta'^{(r,t)}$\\
\hline 
\end{tabular}
\end{equation}
Moreover, these isomorphisms are compatible when $t$ varies, and they are $\Delta'$-equivariant, where $N'$ (resp.\ $\fC^{(t)}_\uptau$,
resp.\ $\hfC'^{(r)}$) is endowed with the trivial (resp.\ trivial, resp.\ canonical) action of $\Delta'$. 

We set 
\begin{equation}\label{p2-fspb9t}
\mC^{(r)}=\underset{\underset{t\in \mQ_{>0}}{\longrightarrow}}{\lim} \fC'^{(r)}\hotimes_{R'_\uptau}\fC^{(t)}_\uptau
\end{equation}
and denote by $\kappa^{(r)}$ the derivation induced by the derivations $\hdelta^{(r,t)}$ \eqref{p2-fspb9e2} (see \eqref{p2-fhtal4f}), 
for $t\in \mQ_{>0}$, and by $\mK^\bullet(\mC^{(r)})$ the Dolbeault complex 
of $(\mC^{(r)},\kappa^{(r)})$. It follows from \ref{p2-fhtal5} that $\mK^\bullet(\mC^{(r)})\otimes_{\mZ_p}\mQ_p$ is a resolution of $\hfC'^{(r)}[\frac 1 p]$. 
Thus, as $N'$ (resp.\ $M'$) is projective over $\hRunp[\frac 1 p]$ (resp.\ $\hoRp[\frac 1 p]$), taking in \eqref{p2-fspb9r} 
the inductive limits, for $t\in \mQ_{>0}$, of the kernels of the Higgs fields $\id\otimes \hdelta^{(r,t)}$, we obtain an isomorphism of $\hfC'^{(r)}$-modules 
\begin{equation}\label{p2-fspb9u}
N' \otimes_{\hRunp}\hfC'^{(r)}\stackrel{\sim}{\rightarrow} M'\otimes_{\hoRp}\hfC'^{(r)}.
\end{equation}
Moreover, by restricting the Higgs $\hRunp$-fields of the bottom line of \eqref{p2-fspb9s}, this is an isomorphism of modules with $\hdelta'^{(r)}$-connection, 
where the $\hdelta'^{(r)}$-connections are defined as in \ref{p1-delta-con4}, 
$N'$ (resp.\ $M'$) being endowed with the Higgs field $\theta'$ (resp.\ $0$), see \eqref{p2-fspb9g}.

Since \eqref{p2-fspb9u} is $\Delta'$-equivariant,  $M'$ is a Dolbeault $\hoRp[\frac 1 p]$-representation of $\Delta'$, and we have an isomorphism 
\begin{equation}
\mH'(M')\stackrel{\sim}{\rightarrow}(N',\theta').
\end{equation}
We deduce the required functorial canonical isomorphism \eqref{p2-fspb9a}.

\section{Higgs--Tate algebras in Faltings topos}\label{p2-htaft}

\subsection{}\label{p2-htaft1}
In this section, we let $f\colon (X,\cM_{X})\rightarrow (S,\cM_S)$ be an adequate morphism of fine logarithmic schemes in the sense of (\cite{agt} III.4.7). 
In particular, the $S$-scheme $X$ is assumed to be of finite type. 
We denote by $X^\circ$ the maximal open subscheme of $X$ where the logarithmic structure $\cM_X$ is trivial; it is an open subscheme of $X_\eta$.
For any $X$-scheme $U$, we set
\begin{equation}\label{p2-htaft1a}
U^\circ=U\times_XX^\circ.
\end{equation}
We denote by $\hbar\colon \oX\rightarrow X$ and $h\colon \oX^\circ\rightarrow X$ the canonical morphisms \eqref{p2-ncgt1a}. To lighten the notation, we set 
\begin{equation}\label{p2-htaft1b}
\tOmega^1_{X/S}=\Omega^1_{(X,\cM_{X})/(S,\cM_S)},
\end{equation}
that we consider as a sheaf of $X_\zar$ or $X_\et$, depending on the context \eqref{p2-ncgt6}.

We endow $\coX=X\times_S\coS$ \eqref{p2-ncgt1a} with the logarithmic structure $\cM_\coX$ pullback of $\cM_X$,
and denote by $\cof\colon (\coX,\cM_{\coX})\rightarrow (\coS,\cM_{\coS})$ the base change of $f$. 
{\em We assume in this section that there exists a Cartesian diagram of $\FLS$ \eqref{p1-NC1} that we fix}
\begin{equation}\label{p2-htaft1c}
\xymatrix{
{(\coX,\cM_{\coX})}\ar[r]^-(0.5){i}\ar[d]_{\cof}\ar@{}[rd]|{\Box}&{(\tX,\cM_{\tX})}\ar[d]^{\tf}\\
{(\coS,\cM_{\coS})}\ar[r]^-(0.5){\iota}&{(\tS,\cM_{\tS}),}}
\end{equation}
where $\iota$ is the strict closed immersion defined in \eqref{p2-ncgt3b}, such that {\em $\tf$ is smooth}.

\subsection{}\label{p2-htaft2}
For any integer $n\geq 1$, we denote by $a\colon X_s\rightarrow X$, $a_n\colon X_s\rightarrow X_n$,
$i_n\colon X_n\rightarrow X$ and $\oi_n\colon \oX_n\rightarrow \oX$ the canonical closed immersions \eqref{p2-ncgt1a}.
Since $k$ is algebraically closed, there exists a unique $S$-morphism $s\rightarrow \oS$.
It induces closed immersions $\oa\colon X_s\rightarrow \oX$ and $\oa_n\colon X_s\rightarrow \oX_n$
lifting $a$ and $a_n$, respectively.
\begin{equation}\label{p2-htaft2a}
\xymatrix{
{X_s}\ar[r]_{\oa_n}\ar@{=}[d]\ar@/^1pc/[rr]^{\oa}&{\oX_n}\ar[r]_{\oi_n}\ar[d]^{\hbar_n}&\oX\ar[d]^\hbar\\
{X_s}\ar[r]^{a_n}\ar@/_1pc/[rr]_{a}&{X_n}\ar[r]^{i_n}&X.}
\end{equation}
Since $\oa_n$ is a universal homeomorphism, we can consider $\co_{\oX_n}$
as a sheaf of rings of $X_{s,\zar}$ or $X_{s,\et}$, depending on the context \eqref{p2-ncgt6}. We set \eqref{p2-htaft1b}
\begin{equation}\label{p2-htaft2b}
\tOmega^1_{\oX_n/\oS_n}=\tOmega^1_{X/S}\otimes_{\co_X}\co_{\oX_n}, 
\end{equation}
that we also consider as a sheaf of $X_{s,\zar}$ or $X_{s,\et}$, depending on the context.

\subsection{}\label{p2-htaft45}
We denote by $\fX$ the formal $\cS$-scheme \eqref{p2-ncgt1}, $p$-adic completion of $\oX$, 
and by $\tOmega^1_{\fX/\cS}$ the $p$-adic completion of the $\co_\coX$-module $\tOmega^1_{\coX/\coS}$ (\cite{egr1} 2.5.1),
that we consider as a sheaf of $X_{s,\zar}$ or $X_{s,\et}$, depending on the context.
The $\co_\fX$-module $\txi^{-1}\tOmega^1_{\fX/\cS}$ \eqref{p2-ncgt3} is canonically 
isomorphic to the $p$-adic completion of the $\co_\coX$-module $\txi^{-1}\tOmega^1_{\coX/\coS}$. 
For any $\co_\fX$-algebra $A$ of $X_{s,\zar}$, we consider Higgs $A$-modules with coefficients in 
$\txi^{-1}\tOmega^1_{\fX/\cS}\otimes_{\co_\fX}A$ \eqref{p1-delta-con1}. We say abusively that they have coefficients in $\txi^{-1}\tOmega^1_{\fX/\cS}$.
Such a Higgs module is said to be {\em coherent} if the underlying $A$-module is coherent.
The categories of these modules will be denoted by $\bHM(A,\txi^{-1}\tOmega^1_{\fX/\cS})$ and $\bHM^\coh(A, \txi^{-1}\tOmega^1_{\fX/\cS})$, respectively. 
We set 
\begin{equation}\label{p2-htaft45b}
\cT=\cHom_{\co_{\fX}}(\tOmega^1_{\fX/\cS},\txi\co_{\fX}),
\end{equation}
which is a locally free $\co_\fX$-module of finite type, that we identify with the dual of $\txi^{-1}\tOmega^1_{\fX/\cS}$.   
We denote by $\rS_{\co_\fX}(\cT)$ the symmetric $\co_\fX$-algebra of $\cT$.
By \ref{p1-delta-con1j}, for every $\co_\fX$-algebra $A$ of $X_{s,\zar}$,  
{\em giving a Higgs $A$-field on an $A$-module $N$ with coefficients in $\txi^{-1}\tOmega^1_{\fX/\cS}$ 
is equivalent to giving an $\rS_{\co_\fX}(\cT)$-module structure on $N$ that is compatible with its $\co_\fX$-module structure}.

We denote by $\bIndHM(\co_\fX,\txi^{-1}\tOmega^1_{\fX/\cS})$
the category of Higgs ind-$\co_\fX$-modules with coefficients in $\txi^{-1}\tOmega^1_{\fX/\cS}$ \eqref{p1-indmal20},
and by $\IndSym_{\co_\fX}(\cT)$ the symmetric ind-$\co_\fX$-algebra of $\cT$ \eqref{p1-indmal6}.  
We have a canonical morphism of ind-$\co_\fX$-algebras 
\begin{equation}\label{p2-htaft45c}
\IndSym_{\co_\fX}(\cT)\rightarrow \iota_{\co_\fX}(\rS_{\co_\fX}(\cT)),
\end{equation}
where $\iota_{\co_\fX}$ is the functor \eqref{p2-cmupiso1a}, which is not an isomorphism. 
By \ref{p1-indmal20j}, {\em giving a Higgs $\co_\fX$-field on an ind-$\co_\fX$-module $\cM$ with coefficients in $\txi^{-1}\tOmega^1_{\fX/\cS}$ 
is equivalent to giving an ind-$\IndSym_{\co_\fX}(\cT)$-module structure on the ind-$\co_\fX$-module $\cM$}.

\subsection{}\label{p2-htaft3}
We denote by $E$ the category of morphisms $(V\rightarrow U)$ over $h$ \eqref{p2-htaft1}, i.e., commutative diagrams
\begin{equation}\label{p2-htaft3a}
\xymatrix{V\ar[r]\ar[d]&U\ar[d]\\
\oX^\circ\ar[r]^h&X}
\end{equation}
such that $U$ is étale over $X$ and the canonical morphism $V\rightarrow \oU^\circ$ is {\em finite étale}.
We consider it as fibered by the functor
\begin{equation}\label{p2-htaft3b}
\pi\colon E\rightarrow \Et_{/X}, \ \ \ (V\rightarrow U)\mapsto U,
\end{equation}
over the étale site of $X$ \eqref{p2-ncgt6} (\cite{ag2} 4.3.2). The fiber category over any $U\in \ob(\Et_{/X})$ is canonically equivalent to
the category $\Et_{\rf/\oU^\circ}$ of finite étale schemes over $\oU^\circ$, and for any morphism of schemes
$\mu\colon U'\rightarrow U$, the inverse image functor $\omu^{\circ+}\colon \Et_{\rf/\oU^\circ}\rightarrow \Et_{\rf/\oU'^\circ}$
is none other than the base change functor by $\omu^{\circ}$.
For any $U\in \ob(\Et_{/X})$, we denote by
\begin{equation}\label{p2-htaft3c}
\iota_U\colon \Et_{\rf/\oU^\circ}\rightarrow E
\end{equation}
the canonical functor. We consider $\pi$ as a {\em covanishing fibered site} \eqref{p2-cmt46} by endowing each fiber category with the étale topology. 
We call it {\em Faltings (covanishing) fibered site} associated with the morphism $h\colon \oX^\circ\rightarrow X$.

By \ref{p2-qfc14} and \ref{p2-cmt44}, we associate with $\pi$ three cleaved and normalized fibered categories: 
\begin{itemize}
\item[(i)] the fibered category
\begin{equation}\label{p2-htaft3f}
\cP^\vee\rightarrow (\Et_{/X})^\circ
\end{equation}
obtained by associating with any $U\in \ob(\Et_{/X})$ the category $(\Et_{\rf/\oU^\circ})^\wedge$
of presheaves of $\mU$-sets on $\Et_{\rf/\oU^\circ}$, and with any morphism $\mu\colon U'\rightarrow U$ of $\Et_{/X}$
the functor
\begin{equation}\label{p2-htaft3g}
\omu^\circ_{\rp}\colon (\Et_{\rf/\oU'^\circ})^\wedge\rightarrow (\Et_{\rf/\oU^\circ})^\wedge
\end{equation}
obtained by composing with the inverse image functor $\omu^{\circ+}\colon \Et_{\rf/\oU^\circ}\rightarrow \Et_{\rf/\oU'^\circ}$;
\item[(ii)] the fibered category 
\begin{equation}\label{p2-htaft3e}
\fF^\vee\rightarrow (\Et_{/X})^\circ
\end{equation}
obtained by associating with any $U\in \ob(\Et_{/X})$ the category $\oU^\circ_\fet$, and with any morphism
$\mu\colon U'\rightarrow U$ of $\Et_{/X}$ the functor
$\omu^\circ_*\colon \oU'^\circ_\fet\rightarrow \oU^\circ_\fet$
direct image by the morphism of topos $\omu^\circ\colon \oU'^\circ_\fet\rightarrow \oU^\circ_\fet$;
\item[(iii)] the fibered $\mU$-topos associated with $\pi$,
\begin{equation}\label{p2-htaft3d}
\fF\rightarrow \Et_{/X}.
\end{equation}
The fiber category of $\fF$ over any $U\in \ob(\Et_{/X})$ is
canonically equivalent to the finite étale topos $\oU^\circ_\fet$ of $\oU^\circ$
and the pullback functor for any morphism $\mu\colon U'\rightarrow U$ of $\Et_{/X}$ is identified with the functor
$\omu^{\circ*}\colon \oU^\circ_\fet\rightarrow \oU'^\circ_\fet$ pullback by the morphism
of topos $\omu^\circ\colon \oU'^\circ_\fet\rightarrow \oU^\circ_\fet$ (\cite{agt} VI.9.3). 
\end{itemize}

We denote by $\hE$ the category of presheaves of $\mU$-sets on $E$. 
By \ref{p2-qfc21} and \ref{p2-qfc22}, we have an equivalence of categories 
\begin{equation}\label{p2-htaft3h}
\begin{array}[t]{clcr}
\hE&\stackrel{\sim}{\rightarrow}& \bHom_{(\Et_{/X})^\circ}((\Et_{/X})^\circ,\cP^\vee)\\
F&\mapsto &\{U\mapsto F\circ \iota_U\}.
\end{array}
\end{equation}
From now on, we will identify $F$ with the section $\{U\mapsto F\circ \iota_U\}$ that is associated with it by this equivalence.
We say that $F$ is a {\em v-presheaf} if for every $U\in \ob(\Et_{/X})$, $F\circ \iota_U$ is a sheaf of $\oU^\circ_\fet$ \eqref{p2-cmt45}.

\subsection{}\label{p2-htaft4}
We endow $E$ with the {\em covanishing topology} associated with $\pi$ \eqref{p2-cmt47}, that is, with
the topology generated by the coverings $\{(V_i\rightarrow U_i)\rightarrow (V\rightarrow U)\}_{i\in I}$
of the following two types:
\begin{itemize}
\item[(v)] $U_i=U$ for every $i\in I$, and $(V_i\rightarrow V)_{i\in I}$ is an étale covering.
\item[(c)] $(U_i\rightarrow U)_{i\in I}$ is an étale covering and $V_i=U_i\times_UV$ for every $i\in I$.
\end{itemize}
The resulting covanishing site $E$  is also called the {\em Faltings  site} associated with $h$ \eqref{p2-htaft1}; it is a $\mU$-site.
We denote by $\tE$ the topos of sheaves of $\mU$-sets on $E$, called the {\em Faltings  topos} associated with $h$ (\cite{agt} VI.10.1, \cite{ag2} 4.3.2).

By (\cite{agt} VI.5.11), the functor \eqref{p2-htaft3h} induces a fully faithful functor
\begin{equation}\label{p2-htaft4a}
\tE\rightarrow \bHom_{(\Et_{/X})^\circ}((\Et_{/X})^\circ,\fF^\vee),
\end{equation}
whose essential image is made up of the sections $\{U\mapsto F_U\}$ satisfying a gluing condition.

The functors
\begin{eqnarray}
\sigma^+\colon \Et_{/X}\rightarrow E,&& U\mapsto (\oU^\circ\rightarrow U),\label{p2-htaft4b}\\
\iota_{X}\colon \Et_{\rf/\oX^\circ}\rightarrow E,&& V\mapsto (V\rightarrow X),\label{p2-htaft4c}
\end{eqnarray}
are continuous and left exact (\cite{agt} VI.10.6). Hence they define two morphisms of topos
\begin{eqnarray}
\sigma\colon \tE \rightarrow X_\et,\label{p2-htaft4d}\\
\beta\colon \tE \rightarrow \oX^\circ_\fet.\label{p2-htaft4e}
\end{eqnarray}
For every sheaf $F=\{U\mapsto F_U\}$ on $E$, we have $\beta_*(F)=F_X$.

The functor
\begin{equation}\label{p2-htaft4f}
\psi^+\colon E\rightarrow \Et_{/\oX^\circ},\ \ \ (V\rightarrow U)\mapsto V
\end{equation}
is continuous and left exact (\cite{agt} VI.10.7). It therefore defines a morphism of topos
\begin{equation}\label{p2-htaft4g}
\psi\colon \oX^\circ_\et\rightarrow \tE.
\end{equation}
We change here the notation compared to {\em loc. cit.}

We denote by $X_\et\gtimes_{X_\et}\oX^\circ_\et$ the covanishing topos of the morphism $h_\et\colon \oX^\circ_\et\rightarrow X_\et$ (\cite{agt} VI.3.12) and by
\begin{equation}\label{p2-htaft4h}
\rho\colon X_\et\gtimes_{X_\et}\oX^\circ_\et\rightarrow \tE
\end{equation}
the canonical morphism (\cite{agt} VI.10.15).

\subsection{}\label{p2-htaft5}
We denote by
\begin{equation}\label{p2-htaft5a}
\gamma\colon \tE_{/\sigma^*(X_\eta)}\rightarrow \tE
\end{equation}
the localization morphism of $\tE$ at $\sigma^*(X_\eta)$ (\cite{sga4} IV 5.2). 
We denote by $\tE_s$ the closed subtopos of $\tE$ complement of the open object $\sigma^*(X_\eta)$ (\cite{sga4} IV 8.3),
i.e., the full subcategory of $\tE$ made up  of the sheaves $F$ such that $\gamma^*(F)$
is a final object of $\tE_{/\sigma^*(X_\eta)}$ (\cite{sga4} IV 9.3.5), by
\begin{equation}\label{p2-htaft5b}
\delta\colon \tE_s\rightarrow \tE
\end{equation}
the canonical embedding and by
\begin{equation}\label{p2-htaft5c}
\sigma_s\colon \tE_s\rightarrow X_{s,\et}
\end{equation}
the morphism of topos induced by $\sigma$ (\cite{agt} (III.9.8.3)). The diagram of morphisms of topos
\begin{equation}\label{p2-htaft5d}
\xymatrix{
{\tE_s}\ar[r]^-(0.5){\sigma_s}\ar[d]_{\delta}&{X_{s,\et}}\ar[d]^{a}\\
{\tE}\ar[r]^-(0.5)\sigma&{X_\et}}
\end{equation}
is commutative up to isomorphism.

\subsection{}\label{p2-htaft6}
The scheme $\oX$ is normal and locally irreducible by (\cite{agt} III.4.2(iii)).
Moreover, the immersion $j\colon X^\circ\rightarrow X$ is quasi-compact since $X$ is Noetherian.
For any object $(V\rightarrow U)$ of $E$, we denote by $\oU^V$ the integral closure of $\oU$ in $V$.
We denote by $\ocB$ the presheaf on $E$ defined for any $(V\rightarrow U)\in \ob(E)$, by
\begin{equation}\label{p2-htaft6a}
\ocB((V\rightarrow U))=\Gamma(\oU^V,\co_{\oU^V}).
\end{equation}
It is a sheaf for the covanishing topology on $E$ by virtue of (\cite{agt} III.8.16).
For any $U\in \ob(\Et_{/X})$, we set
\begin{equation}\label{p2-htaft6b}
\ocB_{U}=\ocB\circ \iota_{U}.
\end{equation}

Unless explicitly stated otherwise, we consider $\sigma$ \eqref{p2-htaft4d}
as a morphism of ringed topos
\begin{equation}\label{p2-htaft6c}
\sigma\colon (\tE,\ocB)\rightarrow (X_\et,\hbar_*(\co_\oX)),
\end{equation}
by the canonical homomorphism $\sigma^*(\hbar_*(\co_\oX))\rightarrow \ocB$ (\cite{agt} III.8.17).

\subsection{}\label{p2-htaft11}
Let $U\in \ob(\Et_{/X})$, $\oy$ be a geometric point of $\oU^\circ$ \eqref{p2-htaft1a}.
The scheme $\oU$ being locally irreducible by (\cite{agt} III.3.3),
it is the sum of the schemes induced on its irreducible components. We denote by $\oU^\star$
the irreducible component of $\oU$ containing $\oy$.
Similarly, $\oU^\circ$ is the sum of the schemes induced on its irreducible components
and $\oU^{\star \circ}=\oU^\star\times_{X}X^\circ$ is the irreducible component of $\oU^\circ$ containing $\oy$.
We denote by $\bB_{\pi_1(\oU^{\star \circ},\oy)}$ the classifying topos of the profinite group $\pi_1(\oU^{\star \circ},\oy)$ and by
\begin{equation}\label{p2-htaft11a}
\nu_\oy\colon \oU^{\star \circ}_\fet \stackrel{\sim}{\rightarrow}\bB_{\pi_1(\oU^{\star \circ},\oy)}
\end{equation}
the fiber functor of $\oU^{\star \circ}_\fet$ at $\oy$ (\cite{agt} VI.9.8). We set
\begin{equation}\label{p2-htaft11b}
\oR^\oy_U=\nu_\oy(\ocB_U|\oU^{\star \circ}).
\end{equation}

\subsection{}\label{p2-htaft7}
For any integer $n\geq 0$, we set
\begin{equation}\label{p2-htaft7a}
\ocB_n=\ocB/p^n\ocB.
\end{equation}
By (\cite{agt} III.9.7), $\ocB_n$ is a ring of $\tE_s$. For any $U\in \ob(\Et_{/X})$, we set \eqref{p2-htaft6b}
\begin{equation}\label{p2-htaft7b}
\ocB_{U,n}=\ocB_U/p^n\ocB_U.
\end{equation}

We have a canonical homomorphism $\sigma_s^*(\co_{\oX_n})\rightarrow \ocB_n$ (\cite{agt} (III.9.9.3)).
The morphism $\sigma_s$ \eqref{p2-htaft5c} is therefore underlying a morphism of ringed topos that we denote by
\begin{equation}\label{p2-htaft7c}
\sigma_n\colon (\tE_s,\ocB_n)\rightarrow (X_{s,\et},\co_{\oX_n}).
\end{equation}
For modules, we use the notation $\sigma_n^{-1}$ (or $\sigma_s^*$)
to denote the  pullback in the sense of abelian sheaves, and we keep the notation
$\sigma_n^*$ for the pullback in the sense of  modules.

We denote by
\begin{equation}\label{p2-htaft7e}
\Sigma_n\colon (\tE_s,\ocB_n)\rightarrow (X_{s,\zar},\co_{\oX_n})
\end{equation}
the composition of the morphism $\sigma_n$ and the canonical morphism \eqref{p1-NC6e}
\begin{equation}\label{p2-htaft7f}
u_n\colon (X_{s,\et},\co_{\oX_n})\rightarrow (X_{s,\zar},\co_{\oX_n}).
\end{equation}
For modules, we use the notation $\Sigma_n^{-1}$
to denote the pullback in the sense of abelian sheaves, and we keep the notation
$\Sigma_n^*$ for the pullback in the sense of  modules.

\subsection{}\label{p2-htaft8}
For any $\mU$-topos $T$, we denote by $T^{\mN^\circ}$ the topos of inverse systems of $T$,
indexed by the ordered set $\mN$ of natural numbers \eqref{p2-ncgt5}.

We denote by $\co_{\bvoX}$ the ring $(\co_{\oX_{n+1}})_{n\in \mN}$ of $X_{s,\et}^{\mN^\circ}$ or $X_{s,\zar}^{\mN^\circ}$, depending on the context.
The reader should not confuse $\co_{\bvoX}$ and $\co_{\coX}$ \eqref{p2-ncgt1a}. We denote by
\begin{equation}\label{p2-htaft8c}
\uplambda\colon (X_{s,\zar}^{\mN^\circ},\co_{\bvoX})\rightarrow (X_{s,\zar}, \co_\fX)
\end{equation}
the canonical morphism of ringed topos \eqref{p2-cmupiso5b}, where $\fX$ is the $p$-adic completion of $\oX$ \eqref{p2-htaft45}.

We denote by $\bvocB$ the ring $(\ocB_{n+1})_{n\in \mN}$ of $\tE_s^{\mN^\circ}$ \eqref{p2-htaft7a}, and by
\begin{equation}\label{p2-htaft8a}
\bvsigma\colon (\tE_s^{\mN^\circ},\bvocB)\rightarrow(X_{s,\et}^{\mN^\circ},\co_{\bvoX})
\end{equation}
the morphism of ringed topos induced by $(\sigma_{n+1})_{n\in \mN}$ \eqref{p2-htaft7c}.
For modules, we use the notation $\bvsigma^{-1}$ to denote the pullback in the sense of abelian sheaves, and we keep the notation
$\bvsigma^*$ for the pullback in the sense of modules.

We denote by
\begin{equation}\label{p2-htaft8b}
\bvSigma\colon (\tE_s^{\mN^\circ},\bvocB)\rightarrow (X_{s,\zar}^{\mN^\circ},\co_{\bvoX})
\end{equation}
the morphism of ringed topos defined by $(\Sigma_{n+1})_{n\in \mN}$ \eqref{p2-htaft7e}. For modules, we use the notation
$\bvSigma^{-1}$ to denote the pullback in the sense of abelian sheaves, and we keep the notation
$\bvSigma^*$ for the pullback in the sense of modules.

We set
\begin{equation}\label{p2-htaft8d}
\hupsigma=\uplambda\circ \bvSigma\colon (\tE_s^{\mN^\circ},\bvocB)\rightarrow (X_{s,\zar},\co_{\fX}).
\end{equation}
For modules, we use the notation $\hupsigma^{-1}$ to denote the pullback in the sense of abelian sheaves, 
and we keep the notation $\hupsigma^*$ for the pullback in the sense of modules.

For every $\co_\fX$-module $\cF$ of $X_{s,\zar}$, we have a canonical isomorphism
\begin{equation}\label{p2-htaft8e}
\hupsigma^*(\cF)\stackrel{\sim}{\rightarrow}\bvSigma^*((\cF/p^{n+1}\cF)_{n\in \mN}).
\end{equation}
In particular, $\hupsigma^*(\cF)$ is adic (\cite{agt} III.7.18).

\subsection{}\label{p2-htaft27}
We denote by $\bMod(\bvocB)$ the category of $\bvocB$-modules of $\tE_s^{\mN^\circ}$,
by $\bIndMod(\bvocB)$ the category of ind-$\bvocB$-modules and by
\begin{equation}\label{p2-htaft27a}
\iota_{\bvocB}\colon \bMod(\bvocB)\rightarrow \bIndMod(\bvocB)
\end{equation}
the canonical fully faithful and exact functor \eqref{p1-abisoind1b}.
We denote by $\bMod_\mQ(\bvocB)$ the category of $\bvocB$-modules up to isogeny  \eqref{p1-abisoind1} and by 
\begin{equation}\label{p2-htaft27b}
\upalpha_{\bvocB}\colon \bMod_\mQ(\bvocB)\rightarrow \bIndMod(\bvocB).
\end{equation}
the canonical fully faithful and exact functor \eqref{p1-bcim5c}.
We will identify $\bMod(\bvocB)$ (resp.\ $\bMod_\mQ(\bvocB)$) with a full subcategory of $\bIndMod(\bvocB)$
by the functor $\iota_{\bvocB}$ (resp.\ $\upalpha_{\bvocB}$), which we will omit from the notation. 

We take again the notation of \ref{p2-cmupiso1}. 
By (\cite{ag2} 2.7.10), the morphism $\hupsigma$ \eqref{p2-htaft8d} induces two adjoint additive functors
\begin{eqnarray}
\rI \hupsigma^*\colon \bIndMod(\co_\fX) \rightarrow \bIndMod(\bvocB),\label{p2-htaft27c}\\
\rI \hupsigma_*\colon \bIndMod(\bvocB) \rightarrow \bIndMod(\co_\fX).\label{p2-htaft27d}
\end{eqnarray}
The functor $\rI \hupsigma^*$ (resp.\ $\rI \hupsigma_*$) is right (resp.\ left) exact.
The functor $\rI \hupsigma_*$ admits a right derived functor
\begin{equation}\label{p2-htaft27e}
\rR\rI \hupsigma_*\colon \bD^+(\bIndMod(\bvocB))\rightarrow \bD^+(\bIndMod(\co_\fX)).
\end{equation}

We set
\begin{equation}\label{p2-htaft27f}
\vupsigma_*=\kappa_{\co_\fX}\circ \rI \hupsigma_*\colon \bIndMod(\bvocB) \rightarrow \bMod(\co_\fX),
\end{equation}
where $\kappa_{\co_\fX}$ is the functor \eqref{p2-cmupiso1b}.  This functor admits a right derived functor
\begin{equation}\label{p2-htaft27g}
\rR\vupsigma_*\colon \bD^+(\bIndMod(\bvocB))\rightarrow \bD^+(\bMod(\co_\fX)),
\end{equation}
canonically isomorphic to $\kappa_{\co_\fX}\circ \rR\rI \hupsigma_*$ (\cite{ag2} 4.3.12).

\begin{rema}\label{p2-htaft270}
For every locally projective $\co_\fX[\frac 1 p]$-module of finite type $\cN$, the ind-$\bvocB$-module $\rI\hupsigma^*(\cN)$ \eqref{p2-cmupiso2c} is flat, in the sense of (\cite{ag2} 2.7.9).
Indeed, by (\cite{ag2} 2.7.17(iv)), we may assume that $\cN$ is free of finite type over $\co_\fX[\frac 1 p]$.
The desired assertion then follows from the fact that the ind-$\bvocB$-module $\rI \hupsigma^*(\co_\fX[\frac 1 p])=\bvocB_\mQ$ (\cite{ag2} (2.9.5.1)) is flat (\cite{ag2} 2.7.9.1).
\end{rema}

\subsection{}\label{p2-htaft23}
We take again the notation of \ref{p2-htaft45}. 
We denote by $\bIndHM(\bvocB,\hupsigma^*(\txi^{-1}\tOmega^1_{\fX/\cS}))$ the category of Higgs ind-$\bvocB$-modules
with coefficients in $\hupsigma^*(\txi^{-1}\tOmega^1_{\fX/\cS} )$ \eqref{p1-indmal20}.
Observe that the $\bvocB$-module $\hupsigma^*(\cT)$ \eqref{p2-htaft45b} is the dual of  $\hupsigma^*(\txi^{-1}\tOmega^1_{\fX/\cS})$.
We denote by $\IndSym_{\bvocB}(\hupsigma^*(\cT))$ the symmetric ind-$\bvocB$-algebra of $\hupsigma^*(\cT)$ \eqref{p1-indmal6}. 
It is canonically isomorphic to the ind-$\bvocB$-algebra $\rI\hupsigma^*(\IndSym_{\co_\fX}(\cT))$ \eqref{p2-htaft27c} (\cite{ag2} (2.7.10.3)). 
By \ref{p1-indmal20j}, {\em giving a Higgs $\bvocB$-field on an ind-$\bvocB$-module $\cM$ with coefficients in 
$\hupsigma^*(\txi^{-1}\tOmega^1_{\fX/\cS})$ 
is equivalent to giving an ind-$\IndSym_{\bvocB}(\hupsigma^*(\cT))$-module structure on the ind-$\bvocB$-module $\cM$}. 

The functor $\rI\hupsigma^*$ \eqref{p2-htaft27c} induces a functor that we denote again by
\begin{equation}\label{p2-htaft23a}
\rI\hupsigma^*\colon 
\begin{array}[t]{clcr}
\bIndHM(\co_\fX,\txi^{-1}\tOmega^1_{\fX/\cS})&\rightarrow& \bIndHM(\bvocB, \hupsigma^*(\txi^{-1}\tOmega^1_{\fX/\cS})),\\
(\cN,\theta)&\mapsto& (\rI\hupsigma^*(\cN),\rI\hupsigma^*(\theta)).
\end{array}
\end{equation}
In view of (\cite{ag2} 2.7.19 and 2.7.18(i)), the functor $\rI\hupsigma_*$ \eqref{p2-htaft27d} induces a functor that we denote again by
\begin{equation}\label{p2-htaft23b}
\rI\hupsigma_*\colon 
\begin{array}[t]{clcr}
\bIndHM(\bvocB, \hupsigma^*(\txi^{-1}\tOmega^1_{\fX/\cS})) &\rightarrow& \bIndHM(\co_\fX,\txi^{-1}\tOmega^1_{\fX/\cS}),\\
(\cM,\theta)&\mapsto& (\rI\hupsigma_*(\cM),\rI\hupsigma_*(\theta)).
\end{array}
\end{equation}
It follows from (\cite{ag2} 2.7.18(ii)) that the functor $\rI\hupsigma^*$ \eqref{p2-htaft23a} is a left adjoint of $\rI\hupsigma_*$ \eqref{p2-htaft23b}.
In view of (\cite{ag2} (2.7.1.5)), we set again
\begin{equation}\label{p2-htaft23c}
\vupsigma_*=\kappa_{\co_\fX}\circ \rI \hupsigma_*\colon \bIndHM(\bvocB, \hupsigma^*(\txi^{-1}\tOmega^1_{\fX/\cS})) \rightarrow \bHM(\co_\fX,\txi^{-1}\tOmega^1_{\fX/\cS}),
\end{equation}
where $\kappa_{\co_\fX}$ is the functor \eqref{p2-cmupiso1b}.

\begin{defi}\label{p2-htaft26}
We say that an ind-$\bvocB$-module $\cF$ is {\em rational} if the multiplication by $p$ on $\cF$ is an isomorphism. 
We say that a Higgs ind-$\bvocB$-module with coefficients in $\hupsigma^*(\txi^{-1}\tOmega^1_{\fX/\cS})$ is {\em rational} if the underlying ind-$\bvocB$-module is rational. 
\end{defi}

\subsection{}\label{p2-htaft9}
We denote by $\bP$ the full subcategory of $\Et_{/X}$ made up of affine schemes $U$ such that one of the following two conditions
is satisfied:
\begin{itemize}
\item[(i)] the scheme $U_s$ is empty; or
\item[(ii)] the morphism $(U,\cM_X|U)\rightarrow (S,\cM_S)$ induced by $f$ \eqref{p2-htaft1} admits an adequate chart (\cite{agt} III.4.4).
\end{itemize}

We endow $\bP$ with the topology induced by that of $\Et_{/X}$. It is a topologically generating $\mU$-small subcategory of $\Et_{/X}$, 
stable by fiber products. Since $X$ is Noetherian and therefore quasi-separated,
any object of $\bP$ is coherent over $X$. 

We denote by $u\colon \bP\rightarrow \Et_{/X}$ the canonical functor, by 
\begin{equation}\label{p2-htaft9a}
\pi_\bP\colon E_\bP\rightarrow \bP
\end{equation}
the fibered site deduced from $\pi$ \eqref{p2-htaft3b} by base change by $u$ \eqref{p2-qfc23}, 
and by $\hE_\bP$ the category of presheaves of $\mU$-sets on $E_\bP$.
For any $U\in \ob(\bP)$, we denote by $\iota_U\colon \Et_{\rf/\oU^\circ}\rightarrow E_\bP$
the canonical functor \eqref{p2-htaft3c}. By \ref{p2-qfc21}, we have an equivalence of categories 
\begin{equation}\label{p2-htaft9c}
\begin{array}[t]{clcr}
\hE_\bP&\stackrel{\sim}{\rightarrow}& \cP(E_\bP/\bP)\\
F&\mapsto &\{U\mapsto F\circ \iota_U\},
\end{array}
\end{equation}
where $\cP(E_\bP/\bP)$ is the category defined in \ref{p2-qfc20}. 
From now on, we will identify $F$ with the object $\{U\mapsto F\circ \iota_U\}$ that is associated with it by this equivalence.

Observe that $\pi_\bP$ is a covanishing fibered site in the sense of \ref{p2-cmt46}. 
We endow $E_\bP$ with the induced covanishing topology \eqref{p2-cmt47} and we denote by $\tE_\bP$ the topos of sheaves of $\mU$-sets on $E_\bP$. 
We denote by $\mu\colon E_\bP\rightarrow E$ the canonical projection and by $\mu_\rp\colon \hE\rightarrow \hE_\bP$ the functor defined by composition with $\mu$ \eqref{p2-cmt4a}. 
The functor $\mu$ is fully faithful and the category $E_\bP$ is $\mU$-small and topologically generating of the site $E$. 
By (\cite{agt} VI.5.21 and VI.5.22), the covanishing topology of $E_\bP$ is induced by that of $E$ by $\mu$. 
By (\cite{sga4} III 4.1), the functor $\mu$ is continuous and cocontinuous and $\mu_\rp$ induces an equivalence of categories
\begin{equation}\label{p2-htaft9d}
\mu_\rs\colon \tE\stackrel{\sim}{\rightarrow}\tE_\bP.
\end{equation}
Then, by (\cite{sga4} III 2.3(2)), the diagram 
\begin{equation}\label{p2-htaft9e}
\xymatrix{
{\hE}\ar[r]^-(0.5){\mu_\rp}\ar[d]&{\hE_\bP}\ar[d]\\
{\tE}\ar[r]^-(0.5){\mu_\rs}&{\tE_\bP,}}
\end{equation}
where the vertical arrows are the ``associated sheaf'' functors,
is commutative up to canonical isomorphism.

\subsection{}\label{p2-htaft10}
We denote by $\bQ$ the full subcategory of $\bP$ \eqref{p2-htaft9} made up of affine schemes $U$ such that one of the following conditions is satisfied:
\begin{itemize}
\item[(i)] the scheme $U_s$ is empty; or
\item[(ii)] there exists a fine and saturated chart $M\rightarrow \Gamma(U,\cM_X)$ for $(U,\cM_X|U)$
inducing an isomorphism
\begin{equation}\label{p2-htaft10a}
M\stackrel{\sim}{\rightarrow} \Gamma(U,\cM_X)/\Gamma(U,\co^\times_X).
\end{equation}
This chart is a priori independent of the adequate chart required in \ref{p2-htaft9}(ii).
\end{itemize}

We endow $\bQ$ with the topology induced by that of $\Et_{/X}$.
It follows from (\cite{agt} II.5.17) that $\bQ$ is a topologically generating $\mU$-small subcategory of $\Et_{/X}$.

We denote by $v\colon \bQ\rightarrow \Et_{/X}$ the canonical functor, by 
\begin{equation}\label{p2-htaft10b}
\pi_\bQ\colon E_\bQ\rightarrow \bQ
\end{equation}
the fibered site deduced from $\pi$ \eqref{p2-htaft3b}
by base change by $v$, and by $\hE_\bQ$ the category of presheaves of $\mU$-sets on $E_\bQ$.
For any $U\in \ob(\bQ)$, we denote by $\iota_U\colon \Et_{\rf/\oU^\circ}\rightarrow E_\bQ$
the canonical functor \eqref{p2-htaft3c}. 
By \ref{p2-qfc21}, we have an equivalence of categories 
\begin{equation}\label{p2-htaft10c}
\begin{array}[t]{clcr}
\hE_\bQ&\stackrel{\sim}{\rightarrow}& \cP(E_\bQ/\bQ)\\
F&\mapsto &\{U\mapsto F\circ \iota_U\},
\end{array}
\end{equation}
where $\cP(E_\bQ/\bQ)$ is the category defined in \ref{p2-qfc20}. 
From now on, we will identify $F$ with the object $\{U\mapsto F\circ \iota_U\}$ that is associated with it by this equivalence.

We denote by $\nu\colon E_\bQ\rightarrow E$ the canonical projection \eqref{p2-qfc23} and by $\nu_\rp\colon \hE\rightarrow \hE_\bQ$ the functor defined by composition with 
$\nu$ \eqref{p2-cmt4a}. The functor $\nu$ is fully faithful and the category $E_\bQ$ is $\mU$-small and topologically generating of the covanishing site $E$. 
We endow $E_\bQ$ with the topology induced by that of $E$ by $\nu$ and denote by $\tE_\bQ$ the topos of sheaves of $\mU$-sets on $E_\bQ$.
By (\cite{sga4} III 4.1), the functor $\nu$ is continuous and cocontinuous and $\nu_\rp$ induces an equivalence of categories 
\begin{equation}\label{p2-htaft10d}
\nu_\rs\colon \tE\stackrel{\sim}{\rightarrow}\tE_\bQ.
\end{equation}
Then, by (\cite{sga4} III 2.3(2)), the diagram 
\begin{equation}\label{p2-htaft10e}
\xymatrix{
{\hE}\ar[r]^-(0.5){\nu_\rp}\ar[d]&{\hE_\bQ}\ar[d]\\
{\tE}\ar[r]^-(0.5){\nu_\rs}&{\tE_\bQ,}}
\end{equation}
where the vertical arrows are the ``associated sheaf'' functors, is commutative up to an isomorphism.

Note that in general, $\bQ$  being not stable under fiber products,
we cannot speak of the covanishing topology on $E_\bQ$ associated with $\pi_\bQ$,
and even less apply (\cite{agt} VI.5.21 and VI.5.22).

\subsection{}\label{p2-htaft12}
Let $Y$ be an object of $\bP$ \eqref{p2-htaft9} such that $Y_s$ is nonempty, $((P,\gamma),(\mN,\nu),\vartheta)$
an adequate chart for the morphism $f|Y\colon (Y,\cM_X|Y)\rightarrow (S,\cM_S)$ induced by $f$,
$\oy$ a geometric point of $\oY^\circ$. Note that the assumptions of \ref{p2-rlps1} are satisfied by $f|Y$.
The scheme $\oY$ being locally irreducible by (\cite{ag1} 4.2.7(iii)),
it is the sum of the schemes induced on its irreducible components. We denote by $\oY^\star$
the irreducible component of $\oY$ containing $\oy$.
Similarly, $\oY^\circ$ is the sum of the schemes induced on its irreducible components
and $\oY^{\star \circ}=\oY^\star\times_{X}X^\circ$ is the irreducible component of $\oY^\circ$ containing $\oy$.
We denote by $\oR^\oy_Y$ the discrete representation of $\pi_1(\oY^{\star\circ},\oy)$ defined in \eqref{p2-htaft11b}
and by $\hoR^\oy_Y$ its $p$-adic Hausdorff completion. Observe that $\oR^\oy_Y$ corresponds to the algebra defined in \eqref{p2-rlps2b}. 

We set
\begin{equation}\label{p2-htaft12a}
\mY^\oy=\Spec(\oR^\oy_Y)\ \ \ {\rm and} \ \ \ \hmY^\oy=\Spec(\hoR^\oy_Y),
\end{equation}
that we endow with the logarithmic structures pullbacks of $\cM_X$,
denoted respectively by $\cM_{\mY^\oy}$ and $\cM_{\hmY^\oy}$. 
Following \ref{p2-rlps3} and (\cite{ag2} 4.4.4), we denote by $(\tmY^\oy,\cM_{\tmY^\oy})$ one or the other of the logarithmic schemes
\begin{equation}\label{p2-htaft12b}
(\cA_2(\mY^\oy),\cM_{\cA_2(\mY^\oy)})\ \ \ {\rm or} \ \ \ (\cA^{\ast}_2(\mY^\oy/S),\cM_{\cA^{\ast}_2(\mY^\oy/S)})
\end{equation}
associated with the morphism $f|Y$ and the adequate chart $((P,\gamma),(\mN,\nu),\vartheta)$
defined in (\cite{ag2} 3.2.8 and 3.2.10), depending on whether we are in the absolute or relative case \eqref{p2-ncgt3} respectively. 
We have a canonical exact closed immersion \eqref{p2-rlps3f}
\begin{equation}\label{p2-htaft12c}
\mi^\oy_Y\colon (\hmY^\oy,\cM_{\hmY^\oy})\rightarrow (\tmY^\oy,\cM_{\tmY^\oy}).
\end{equation}

The conditions \ref{p2-hta1}(i)-(iv) are satisfied by the commutative diagram \eqref{p2-htaft1c} 
\begin{equation}\label{p2-htaft12g}
\xymatrix{
{(\hmY^\oy,\cM_{\hmY^\oy})}\ar[r]^-(0.5){\mi^\oy_Y}\ar[d]&{(\tmY^\oy,\cM_{\tmY^\oy})}\ar@/^2pc/[dd]\\
{(\coX,\cM_{\coX})}\ar[r]^-(0.5){i}\ar[d]^{\cof}\ar@{}[rd]|{\Box}&{(\tX,\cM_{\tX})}\ar[d]_{\tf}\\
{(\coS,\cM_{\coS})}\ar[r]^-(0.5){\iota}&{(\tS,\cM_{\tS}).}}
\end{equation}
For any rational number $r\geq 0$, 
we denote by $\cL^{(r)}_{\tmY^\oy/\tX}$ the $(r)$-twisted torsor of liftings of the canonical morphism $\hmY^\oy\rightarrow \coX$ to $\tmY^\oy$ over $\tX$, 
and by $\cC^{(r)}_{\tmY^\oy/\tX}$ (resp.\ $\cF^{(r)}_{\tmY^\oy/\tX}$) the Higgs--Tate algebra (resp.\ extension) of $\tmY^\oy$ over $\tX$ of thickness $r$ \eqref{p2-hta7}. 
Denoting by $\tY\rightarrow \tX$ the unique étale morphism which lifts $\coY\rightarrow \coX$ (\cite{ega4} 18.1.2),
$\cL^{(r)}_{\tmY^\oy/\tX}$ identifies canonically with the Higgs--Tate torsor $\cL^{(r)}_{\tmY^\oy/\tY}$ associated with the deformation  
\begin{equation}\label{p2-htaft12h}
\xymatrix{
{(\coY,\cM_\coX|\coY)}\ar[r]\ar[d]\ar@{}[rd]|{\Box}&{(\tY,\cM_{\tX}|\tY)}\ar[d]\\
{(\coS,\cM_{\coS})}\ar[r]^-(0.5){\iota}&{(\tS,\cM_{\tS})}}
\end{equation}
induced by \eqref{p2-htaft1c}; this torsor was introduced in \ref{p2-rlps5}. 

We set  
\begin{equation}\label{p2-htaft12d}
\fF^{\oy,(r)}_Y=\Gamma(\hmY^\oy,\cF^{(r)}_{\tmY^\oy/\tX}) \ \ \ {\rm and}\ \ \ \fC^{\oy,(r)}_Y=\Gamma(\hmY^\oy,\cC^{(r)}_{\hmY^\oy/\tX}).
\end{equation}
We have a canonical exact sequence of $\hoR^\oy_Y$-modules
\begin{equation}\label{p2-htaft12e}
0\rightarrow \hoR^\oy_Y\rightarrow \fF^{\oy,(r)}_Y\rightarrow \txi^{-1}\tOmega^1_{X/S}(Y) \otimes_{\co_X(Y)}\hoR^\oy_Y \rightarrow 0,
\end{equation}
and a canonical isomorphism of $\hoR$-algebras
\begin{equation}\label{p2-htaft12f}
\fC^{\oy,(r)}_Y\stackrel{\sim}{\rightarrow}\underset{\underset{n\geq 0}{\longrightarrow}}\lim\ \rS^n_{\hoR}(\fF^{\oy,(r)}_Y). 
\end{equation}
Note that we have changed the normalization in \eqref{p2-htaft12e} compared to \ref{p2-hta7}, namely we identified $\Omega^{(r)}$ \eqref{p2-hta4d} with $\Omega$ \eqref{p2-hta4c}, 
so that the canonical morphism $\pi^{(r)}\colon \Omega^{(r)}\rightarrow \Omega$ is identified with the multiplication by $p^r$ on $\Omega$, see \ref{p2-hta71}. 

The action of $\pi_1(\oY^{\star\circ},\oy)$ on $(\tmY^\oy,\cM_{\tmY^\oy})$ induces a canonical $\hoR^\oy_Y$-semilinear action of $\pi_1(\oY^{\star\circ},\oy)$ on $\fF^{\oy,(r)}_Y$ 
that is continuous for the $p$-adic topology (see \cite{ag2} 3.2.15). The morphisms of the sequence \eqref{p2-htaft12e} are $\pi_1(\oY^{\star\circ},\oy)$-equivariant.
We deduce from this an action of $\pi_1(\oY^{\star\circ},\oy)$ on $\fC^{\oy,(r)}_Y$ by ring homomorphisms
that is continuous for the $p$-adic topology and which extends the canonical action on $\hoR^\oy_Y$. 
Note that these representations depend on
the adequate chart $((P,\gamma),(\mN,\nu),\vartheta)$. However, they do not depend on it if $Y$ is an object of $\bQ$ \eqref{p2-htaft10} by \ref{p2-rlps4}.

\subsection{}\label{p2-htaft13}
Let $Y$ be an object of $\bP$ such that $Y_s$ is nonempty, $((P,\gamma),(\mN,\nu),\vartheta)$
an adequate chart for the morphism $f|Y\colon (Y,\cM_X|Y)\rightarrow (S,\cM_S)$ induced by $f$, $r$ a rational number $\geq 0$, $n$ an integer $\geq 0$.
Recall that $\oY^\circ$ is the sum of the schemes induced on its irreducible components \eqref{p2-htaft12}.
Let $W$ be an irreducible component of $\oY^\circ$, $\Pi(W)$ its fundamental groupoid (\cite{agt} VI.9.10).
In view of \eqref{p2-htaft11b} and (\cite{agt} VI.9.11), the sheaf $\ocB_Y|W$ \eqref{p2-htaft6b} of $W_\fet$ defines a functor
\begin{equation}\label{p2-htaft13a}
\Pi(W)\rightarrow \Ens, \ \ \ \oy\mapsto \oR^\oy_Y.
\end{equation}
Let $\oy$ and $\oy'$ be geometric points of $W$. We easily check, by following the constructions in \ref{p2-htaft12},
that every morphism of $\pi_1(W,\oy,\oy')$ induces an isomorphism of algebras
\begin{equation}
\cF^{(r)}_{\tmY^\oy/\tX}\stackrel{\sim}{\rightarrow} \cF^{(r)}_{\tmY^{\oy'}/\tX},
\end{equation} 
over the isomorphism of rings $\oR^\oy_Y\stackrel{\sim}{\rightarrow} \oR^{\oy'}_Y$ defined by the functor \eqref{p2-htaft13a}. 
We deduce from this a functor
\begin{equation}\label{p2-htaft13b}
\Pi(W)\rightarrow \Ens, \ \ \ \oy\mapsto \fF^{\oy,(r)}_Y/p^n\fF^{\oy,(r)}_Y.
\end{equation}
For every geometric point $\oy$ of $W$, $\fF^{\oy,(r)}_Y/p^n\fF^{\oy,(r)}_Y$ is a continuous and discrete representation of $\pi_1(W,\oy)$
(\cite{agt} III.10.14). Consequently, by virtue of (\cite{agt} VI.9.11), the functor \eqref{p2-htaft13b} defines a $(\ocB_{Y,n}|W)$-module
$\cF^{(r)}_{W,n}$ of $W_\fet$, unique up to canonical isomorphism, where $\ocB_{Y,n}=\ocB_Y/p^n\ocB_Y$ \eqref{p2-htaft7b}.
By descent (\cite{giraud2} II 3.4.4), there exists a $\ocB_{Y,n}$-module $\cF^{(r)}_{Y,n}$ of $\oY^\circ_\fet$, unique
up to canonical isomorphism, such that for every irreducible component $W$ of $\oY^\circ$, we have $\cF^{(r)}_{Y,n}|W=\cF^{(r)}_{W,n}$.
We define similarly the $\ocB_{Y,n}$-algebra $\cC^{(r)}_{Y,n}$ from the representations $\fC^{\oy,(r)}_Y/p^n\fC^{\oy,(r)}_Y$. 

The exact sequence \eqref{p2-htaft12e} induces an exact sequence of $\ocB_{Y,n}$-modules
\begin{equation}\label{p2-htaft13c}
0\rightarrow \ocB_{Y,n}\rightarrow \cF^{(r)}_{Y,n}\rightarrow
\txi^{-1}\tOmega^1_{X/S}(Y)\otimes_{\co_X(Y)}\ocB_{Y,n} \rightarrow 0.
\end{equation}
We have a canonical isomorphism of $\ocB_{Y,n}$-algebras
\begin{equation}\label{p2-htaft13d}
\cC^{(r)}_{Y,n}\stackrel{\sim}{\rightarrow}\underset{\underset{m\geq 0}{\longrightarrow}}\lim\ \rS^m_{\ocB_{Y,n}}(\cF^{(r)}_{Y,n}).
\end{equation}
Observe that $\cF^{(r)}_{Y,n}$ and $\cC^{(r)}_{Y,n}$ depend on the choice of the adequate chart $((P,\gamma),(\mN,\nu),\vartheta)$.
However, they do not depend on it if $Y$ is an object of $\bQ$ \eqref{p2-htaft10}  by \ref{p2-rlps4}.
Observe that $\cF^{(r)}_{Y,n}$ and $\cC^{(r)}_{Y,n}$ depend on the choice of the deformation $\tf$ fixed in \eqref{p2-htaft1c}.

For all rational numbers $r\geq r'\geq 0$, we have a canonical $\ocB_{Y,n}$-linear morphism
\begin{equation}\label{p2-htaft13e}
\tta_{Y,n}^{r,r'}\colon\cF_{Y,n}^{(r)}\rightarrow \cF_{Y,n}^{(r')}
\end{equation}
which lifts the multiplication by $p^{r-r'}$ over
$\txi^{-1}\tOmega^1_{X/S}(Y)\otimes_{\co_X(Y)}\ocB_{Y,n}$ and which extends the identity over $\ocB_{Y, n}$ \eqref{p2-htaft13c}.
It induces a homomorphism of $\ocB_{Y,n}$-algebras
\begin{equation}\label{p2-htaft13f}
\alpha_{Y,n}^{r,r'}\colon \cC_{Y,n}^{(r)}\rightarrow \cC_{Y,n}^{(r')}.
\end{equation}

\subsection{}\label{p2-htaft14}
Let $g\colon Y\rightarrow Z$ be a morphism of $\bP$ \eqref{p2-htaft9} such that $Y_s$ is non-empty,
$((P,\gamma),(\mN,\nu),\vartheta)$ an adequate chart for the morphism $f|Z\colon (Z,\cM_X|Z)\rightarrow (S,\cM_S)$ induced by $f$,
$\oy$ a geometric point of $\oY^\circ$, $\oz=\ogg(\oy)$.
We endow the morphism $f|Y\colon (Y,\cM_X|Y)\rightarrow (S,\cM_S)$ induced by $f$ with the same adequate chart $((P,\gamma),(\mN, \nu),\vartheta)$.
Recall that $\oY$ and $\oZ$ are sums of the schemes induced over their irreducible components.
We denote by $\oY^\star$ the irreducible component of $\oY$ containing $\oy$ and
by $\oZ^\star$ the irreducible component of $\oZ$ containing $\oz$, so that $\ogg(\oY^\star)\subset \oZ^\star$.
The morphism $\ogg^\circ\colon \oY^\circ\rightarrow \oZ^\circ$
induces a group homomorphism $\pi_1(\oY^{\star \circ},\oy) \rightarrow\pi_1(\oZ^{\star \circ},\oz)$.
The canonical morphism $(\ogg^\circ)^*_\fet(\ocB_Z)\rightarrow \ocB_Y$ \eqref{p2-htaft9c}
determines a $\pi_1(\oY^{\star \circ},\oy)$-equivariant ring homomorphism $u\colon \oR^\oz_Z\rightarrow \oR^{\oy}_Y$. We have a Cartesian diagram 
\begin{equation}\label{p2-htaft14a}
\xymatrix{
{\hmY^{\oy}}\ar[r]\ar[d]_h&{\tmY^{\oy}}\ar[d]^\thh\\
{\hmZ^{\oz}}\ar[r]&{\tmZ^{\oz},}}
\end{equation}
where the morphisms $h$ and $\thh$ are induced by $u$. 
By \eqref{p1-rdt4h}, it induces an isomorphism of torsors under the $\hoR^{\oy}_Y$-module $\Hom_{\co_X(Y)}(\tOmega^1_{X/S}(Y),\txi^{-1}\hoR^{\oy}_Y)$,
\begin{equation}\label{p2-htaft14b}
h^+(\cL_{\tmZ^\oz/\tX})\rightarrow \cL_{\tmY^\oy/\tX},
\end{equation}
where $h^+$ is the affine pullback \eqref{p1-NC5}. For every rational number $r\geq 0$, the latter induces an 
isomorphism of torsors under the $\hoR^{\oy}_Y$-module $\Hom_{\co_X(Y)}(\tOmega^1_{X/S}(Y),\txi^{-1}\hoR^{\oy}_Y)$,
\begin{equation}\label{p2-htaft14c}
h^+(\cL^{(r)}_{\tmZ^\oz/\tX})\rightarrow \cL^{(r)}_{\tmY^\oy/\tX}.
\end{equation}
By (\cite{agt} II.4.14), we deduce an $\hoR^\oz_Y$-linear isomorphism
\begin{equation}\label{p2-htaft14d}
\fF^{\oy,(r)}_Y\stackrel{\sim}{\rightarrow} \fF^{\oz,(r)}_Z\otimes_{\hoR^\oz_Z}\hoR^\oy_Y.
\end{equation}
We easily check that it is $\pi_1(\oY^{\star \circ}, \oy)$-equivariant (\cite{agt} II.4.22). 
Its inverse determines a $\pi_1(\oY^{\star \circ},\oy)$-equivariant $\hoR^\oz_Z$-linear morphism
\begin{equation}\label{p2-htaft14e}
\fF^{\oz,(r)}_Z\rightarrow \fF^{\oy,(r)}_Y
\end{equation}
that fits into a commutative diagram
\begin{equation}\label{p2-htaft14f}
\xymatrix{
0\ar[r]&{\hoR^\oz_Z}\ar[r]\ar[d]&{\fF^{\oz,(r)}_Z}\ar[r]\ar[d]&
{\txi^{-1}\tOmega^1_{X/S}(Z)\otimes_{\co_X(Z)}\hoR^\oz_Z}\ar[r]\ar[d]&0\\
0\ar[r]&{\hoR^\oy_Y}\ar[r]&{\fF^{\oy,{(r)}}_Y}\ar[r]&
{\txi^{-1}\tOmega^1_{X/S}(Y)\otimes_{\co_X(Y)}\hoR^\oy_Y}\ar[r]&0,}
\end{equation}
where the horizontal sequences are defined in \eqref{p2-htaft12e}. 
We deduce from this a $\pi_1(\oY^{\star \circ},\oy)$-equivariant homomorphism of $\hoR^\oz_Z$-algebras
\begin{equation}\label{p2-htaft14g}
\fC^{\oz,(r)}_Z\rightarrow\fC^{\oy,(r)}_Y.
\end{equation}

By the obvious functoriality of the equivalence in (\cite{agt} VI.9.11), the morphism \eqref{p2-htaft14e} induces a $(\ogg^{\circ})_\fet^*(\ocB_{Z,n})$-linear morphism
\begin{equation}\label{p2-htaft14h}
(\ogg^{\circ})^*_\fet(\cF^{(r)}_{Z,n}) \rightarrow\cF^{(r)}_{Y,n},
\end{equation}
and therefore by adjunction, a $\ocB_{Z,n}$-linear morphism
\begin{equation}\label{p2-htaft14i}
\cF^{(r)}_{Z,n} \rightarrow \ogg^{\circ}_{\fet *}(\cF^{(r)}_{Y,n}),
\end{equation}
see (\cite{ag2} 4.4.6). It follows from \eqref{p2-htaft14f} that the diagram
\begin{equation}\label{p2-htaft14j}
\xymatrix{
0\ar[r]&{(\ogg^{\circ})^*_\fet(\ocB_{Z,n})}\ar[r]\ar[d]&{(\ogg^{\circ})^*_\fet(\cF^{(r)}_{Z,n})}\ar[r]\ar[d]&
{\txi^{-1}\tOmega^1_{X/S}(Z)\otimes_{\co_X(Z)}(\ogg^{\circ})^*_\fet(\ocB_{Z,n })}\ar[r]\ar[d]&0\\
0\ar[r]&{\ocB_{Y,n}}\ar[r]&{\cF^{(r)}_{Y,n}}\ar[r]&
{\txi^{-1}\tOmega^1_{X/S}(Y)\otimes_{\co_X(Y)}\ocB_{Y,n}}\ar[r]&0,}
\end{equation}
where the horizontal sequences are defined in \eqref{p2-htaft13c}, is commutative.
Similarly, the homomorphism \eqref{p2-htaft14g} induces a homomorphism of $(\ogg^{\circ})_\fet^*(\ocB_{Z,n})$-algebras
\begin{equation}\label{p2-htaft14k}
(\ogg^{\circ})^*_\fet(\cC^{(r)}_{Z,n}) \rightarrow\cC^{(r)}_{Y,n},
\end{equation}
and therefore by adjunction a homomorphism of $\ocB_{Z,n}$-algebras
\begin{equation}\label{p2-htaft14l}
\cC^{(r)}_{Z,n} \rightarrow\ogg^{\circ}_{\fet *}(\cC^{(r)}_{Y,n}).
\end{equation}
The morphisms \eqref{p2-htaft14h} and \eqref{p2-htaft14k} are compatible via the isomorphism \eqref{p2-htaft13d}

Note that the morphisms \eqref{p2-htaft14i} and \eqref{p2-htaft14l} depend on the choice of the adequate chart $((P,\gamma),(\mN,\nu),\vartheta)$.
However, they do not depend on it if $Y$ and $Z$ are objects of $\bQ$ \eqref{p2-htaft10} by \ref{p2-rlps4}.

\subsection{}\label{p2-htaft15}
For any rational number $r\geq 0$, any integer $n\geq 0$ and any object $Y$ of $\bP$ such that $Y_s$ is empty,
we set $\cC^{(r)}_{Y,n}=\cF^{(r)}_{Y,n}=0$. The exact sequence \eqref{p2-htaft13c} still holds in this case, since $\ocB_Y$ is a  $\oK$-algebra.
The morphisms \eqref{p2-htaft14i} and \eqref{p2-htaft14l} are then defined for any morphism
of $\bP$, and they satisfy cocycle relations of the type (\cite{egr1} (1.1.2.2)).

\subsection{}\label{p2-htaft16}
Let $r$ be a rational number $\geq 0$, $n$ an integer $\geq 0$.
The correspondences $\{U\mapsto p^n\ocB_U\}$ and $\{U\mapsto \ocB_{U,n}\}$
are naturally objects of $\cP(E/\Et_{/X})$ \eqref{p2-qfc20}. They correspond by the equivalence of categories \eqref{p2-htaft3h}
to presheaves on $E$. The canonical morphisms
\begin{eqnarray}
\{U\mapsto p^n\ocB_U\}^a&\rightarrow&p^n\ocB,\label{p2-htaft16a}\\
\{U\mapsto \ocB_{U,n}\}^a&\rightarrow&\ocB_n,\label{p2-htaft16b}
\end{eqnarray}
where the exponent $a$ means the associated sheaf, are isomorphisms
by virtue of (\cite{agt} VI.8.2 and VI.8.9). 
With the notation of  \ref{p2-htaft10}, we have $\{Y\mapsto \cB_{Y,n}\}=\nu_\rp(\{U\mapsto \cB_{U,n}\})$ in $\cP(E_\bQ/\bQ)$. 
In view of \eqref{p2-htaft10e}, we have a canonical isomorphism 
\begin{equation}\label{p2-htaft16i}
\nu_\rs(\cB_n)\stackrel{\sim}{\rightarrow}\{Y\mapsto \cB_{Y,n}\}^a. 
\end{equation}

By \ref{p2-htaft14}, we have canonical objects $\{Y\mapsto \cF^{(r)}_{Y,n} \}$ and $\{Y\mapsto \cC^{(r)}_{Y,n}\}$
of $\cP(E_\bQ/\bQ)$. They correspond by the equivalence of categories \eqref{p2-htaft10c}
to presheaves on $E_\bQ$, of modules and algebras, respectively, relative to the ring $\{Y\mapsto \cB_{Y,n}\}$. 
Hence, there exist a canonical $\ocB_n$-module $\cF^{(r)}_n$ and a canonical $\ocB_n$-algebra $\cC^{(r)}_n$ such that 
\begin{eqnarray}
\nu_\rs(\cF^{(r)}_n)&=&\{Y\mapsto \cF^{(r)}_{Y,n}\}^a,\label{p2-htaft16d}\\
\nu_\rs(\cC^{(r)}_n)&=&\{Y\mapsto \cC^{(r)}_{Y,n}\}^a.\label{p2-htaft16e}
\end{eqnarray}
We call $\cF^{(r)}_n$ (resp.\ $\cC^{(r)}_n$) the {\em Higgs--Tate $\ocB_n$-extension} (resp.\  {\em  Higgs--Tate $\ocB_n$-algebra}) 
of thickness $r$ associated with the deformation $\tf$ \eqref{p2-htaft1c}.
We set $\cF_n=\cF_n^{(0)}$ (resp.\ $\cC_n=\cC_n^{(0)}$), that we
call the {\em  Higgs--Tate $\ocB_n$-extension} (resp.\ {\em  Higgs--Tate $\ocB_n$-algebra}) associated with the deformation $\tf$.

For all rational numbers $r\geq r'\geq 0$, the morphisms \eqref{p2-htaft13e} induce a $\ocB_n$-linear morphism
\begin{equation}\label{p2-htaft16f}
\tta_n^{r,r'}\colon \cF^{(r)}_n\rightarrow \cF_n^{(r')}.
\end{equation}
The homomorphisms \eqref{p2-htaft13f} induce a homomorphism of $\ocB_n$-algebras
\begin{equation}\label{p2-htaft16g}
\alpha_n^{r,r'}\colon \cC_n^{(r)}\rightarrow \cC_n^{(r')}.
\end{equation}
For all rational numbers $r\geq r'\geq r''\geq 0$, we have
\begin{equation}\label{p2-htaft16h}
\tta_n^{r,r''}=\tta_n^{r',r''} \circ \tta_n^{r,r'} \ \ \ {\rm and}\ \ \ \alpha_n^{r, r''}=\alpha_n^{r',r''} \circ \alpha_n^{r,r'}.
\end{equation}

\begin{prop}[\cite{agt} III.10.22]\label{p2-htaft17}
Let $r$ be a rational number $\geq 0$, $n$ an integer $\geq 1$. Then,
\begin{itemize}
\item[{\rm (i)}] The sheaves $\cF^{(r)}_n$ and $\cC^{(r)}_n$ are objects of $\tE_s$.
\item[{\rm (ii)}] With the notation of \eqref{p2-htaft2b} and \eqref{p2-htaft7c}, we have a canonical locally split exact sequence of $\ocB_n$-modules
\begin{equation}\label{p2-htaft17a}
0\rightarrow \ocB_n\rightarrow \cF^{(r)}_n\rightarrow
\sigma_n^*(\txi^{-1}\tOmega^1_{\oX_n/\oS_n})\rightarrow 0.
\end{equation}
It induces for any integer $m\geq 1$ an exact sequence of $\ocB_n$-modules
\begin{equation}\label{p2-htaft17b}
0\rightarrow \rS^{m-1}_{\ocB_n}(\cF^{(r)}_n)\rightarrow \rS^m_{\ocB_n}(\cF^{(r)}_n)\rightarrow
\sigma_n^*(\rS^m_{\co_{\oX_n}}(\txi^{-1}\tOmega^1_{\oX_n/\oS_n}))\rightarrow 0.
\end{equation}
In particular, the $\ocB_n$-modules $(\rS^m_{\ocB_n}(\cF^{(r)}_n))_{m\in \mN}$ form a filtered direct system.
\item[{\rm (iii)}] We have a canonical isomorphism of $\ocB_n$-algebras
\begin{equation}\label{p2-htaft17c}
\cC^{(r)}_n \stackrel{\sim}{\rightarrow}\underset{\underset{m\geq 0}{\longrightarrow}}\lim\ \rS^m_{\ocB_n}(\cF^{(r)}_n).
\end{equation}
\item[{\rm (iv)}] For all rational numbers $r\geq r'\geq 0$, the diagram
\begin{equation}\label{p2-htaft17d}
\xymatrix{
0\ar[r]&{\ocB_n}\ar[r]\ar@{=}[d]&
{\cF^{(r)}_n}\ar[r]\ar[d]^{\tta_n^{r,r'}}&{\sigma_n^*(\txi^{-1}\tOmega^ 1_{\oX_n/\oS_n})}\ar[r]\ar[d]^{\cdot p^{r-r'}}& 0\\
0\ar[r]&{\ocB_n}\ar[r]&{\cF^{(r')}_n}\ar[r]&{\sigma_n^*(\txi^{-1}\tOmega ^1_{\oX_n/\oS_n})}\ar[r]& 0,}
\end{equation}
where the horizontal lines are the exact sequences \eqref{p2-htaft17a}
and the right vertical arrow denotes multiplication by $p^{r-r'}$, is commutative.
Moreover, the morphisms $\tta_n^{r,r'}$ and $\alpha_n^{r,r'}$ are compatible with the isomorphisms
\eqref{p2-htaft17c} for $r$ and $r'$.
\end{itemize}
\end{prop}

\subsection{}\label{p2-htaft20}
Let $r$ be a rational number $\geq 0$, $n$ an integer $\geq 1$. There is a unique $\ocB_n$-derivation of $\cC_n^{(r)}$ \eqref{p2-htaft17c}
\begin{equation}\label{p2-htaft20a}
d_{\cC_n^{(r)}}\colon \cC_n^{(r)}\rightarrow \sigma_n^*(\txi^{-1}\tOmega^1_{\oX_n/\oS_n})\otimes_{\ocB_n}\cC_n^{(r)}
\end{equation}
that extends the canonical morphism $\cF_n^{(r)}\rightarrow \sigma_n^*(\txi^{-1}\tOmega^1_{\oX_n/\oS_n})$ \eqref{p2-htaft17a}.
It canonically identifies with the universal $\ocB_n$-derivation of $\cC_n^{(r)}$ \eqref{p1-imdpa19}. 
It is also a Higgs $\ocB_n$-field on $\cC_n^{(r)}$ with coefficients in $\sigma^*_n(\txi^{-1}\tOmega^1_{\oX_n/\oS_n})$ \eqref{p1-delta-con1}.
We set 
\begin{equation}\label{p2-htaft20b}
\delta_{\cC_n^{(r)}}=p^rd_{\cC_n^{(r)}}\colon \cC_n^{(r)}\rightarrow \sigma_n^*(\txi^{-1}\tOmega^1_{\oX_n/\oS_n})\otimes_{\ocB_n}\cC_n^{(r)}.
\end{equation}
It follows from \ref{p2-htaft17}(iv) that for all rational numbers $r\geq r'\geq 0$, we have
\begin{equation}\label{p2-htaft20c}
(\id \otimes \alpha^{r,r'}_n) \circ \delta_{\cC_n^{(r)}}=\delta_{\cC_n^{(r')}} \circ \alpha^{r,r'}_n.
\end{equation}

\subsection{}\label{p2-htaft18}
Let $r$ be a rational number $\geq 0$.
For all integers $m\geq n\geq 1$, we have a canonical $\ocB_m$-linear morphism
$\cF^{(r)}_m\rightarrow \cF^{(r)}_n$ compatible with the exact sequence \eqref{p2-htaft17a}
and a canonical homomorphism of $\ocB_m$-algebras $\cC^{(r)}_m\rightarrow \cC^{(r)}_n$ such that
the induced morphisms
\begin{equation}\label{p2-htaft18a}
\cF^{(r)}_m\otimes_{\ocB_m}\ocB_n\rightarrow \cF^{(r)}_n\ \ \ {\rm and}\ \ \
\cC^{(r)}_m\otimes_{\ocB_m}\ocB_n\rightarrow \cC^{(r)}_n
\end{equation}
are isomorphisms. These morphisms form compatible systems when $m$ and $n$ vary.
With the notation of \ref{p2-htaft8}, the $\bvocB$-module $(\cF^{(r)}_{n+1})_{n\in \mN}$ of $\tE_s^{\mN^\circ}$, denote by $\bvcF^{(r)}$, is called 
the {\em  Higgs--Tate $\bvocB$-extension of thickness $r$} associated with the deformation $\tf$.
The $\bvocB$-algebra $(\cC^{(r)}_{n+1})_{n\in \mN}$ of $\tE_s^{\mN^\circ}$, denoted by $\bvcC^{(r)}$, is called 
the {\em Higgs--Tate $\bvocB$-algebra of thickness $r$} associated with the deformation $\tf$.
By (\cite{agt} III.7.3(i), (III.7.5.4) and (III.7.12.1)), the exact sequence \eqref{p2-htaft17a}
induces an exact sequence of $\bvocB$-modules
\begin{equation}\label{p2-htaft18b}
0\rightarrow \bvocB\rightarrow \bvcF^{(r)}\rightarrow
\bvsigma^*(\txi^{-1}\tOmega^1_{\bvoX/\bvoS})\rightarrow 0,
\end{equation}
where we denoted by $\txi^{-1}\tOmega^1_{\bvoX/\bvoS}$ the $\co_{\bvoX}$-module $(\txi^{-1}\tOmega^1_{\oX_{ n+1}/\oS_{n+1}})_{n\in \mN}$ of $X_{s,\et}$ \eqref{p2-htaft8}. 
Since the $\co_X$-module $\tOmega^1_{X/S}$ is locally free of finite type,
the $\bvocB$-module $\bvsigma^*(\txi^{-1}\tOmega^1_{\bvoX/\bvoS})$ is locally free of finite type
and the sequence \eqref{p2-htaft18b} is locally split.
It induces for every integer $m\geq 1$ an exact sequence of $\bvocB$-modules
\begin{equation}\label{p2-htaft18c}
0\rightarrow \rS^{m-1}_{\bvocB}(\bvcF^{(r)})\rightarrow \rS^m_{\bvocB}(\bvcF^{(r)})\rightarrow
\bvsigma^*(\rS^m_{\co_{\bvoX}}(\txi^{-1}\tOmega^1_{\bvoX/\bvoS}))\rightarrow 0.
\end{equation}
In particular, the $\bvocB$-modules $(\rS^m_{\bvocB}(\bvcF^{(r)}))_{m\in \mN}$ form a filtered direct system.
By (\cite{agt} III.7.3(i) and (III.7.12.3)), we have a canonical isomorphism of $\bvocB$-algebras
\begin{equation}\label{p2-htaft18d}
\bvcC^{(r)}\stackrel{\sim}{\rightarrow}\underset{\underset{m\geq 0}{\longrightarrow}}\lim\ \rS^m_{\bvocB}(\bvcF^{ (r)}).
\end{equation}

We set $\bvcF=\bvcF^{(0)}$ and $\bvcC=\bvcC^{(0)}$ and consider the extension \eqref{p2-htaft18b}
\begin{equation}\label{p2-htaft18h}
0\rightarrow \bvocB\rightarrow \bvcF\rightarrow
\bvsigma^*(\txi^{-1}\tOmega^1_{\bvoX/\bvoS})\rightarrow 0.
\end{equation}
We call $\bvcF$ (resp.\ $\bvcC$) the {\em  Higgs--Tate $\bvocB$-extension} (resp.\ the {\em Higgs--Tate $\bvocB$-algebra})
associated with the deformation $\tf$. 

For all rational numbers $r\geq r'\geq 0$, the morphisms $(\tta_n^{r,r'})_{n\in \mN}$ \eqref{p2-htaft16f}
induce a $\bvocB$-linear morphism
\begin{equation}\label{p2-htaft18e}
\bvtta^{r,r'}\colon \bvcF^{(r)}\rightarrow \bvcF^{(r')}.
\end{equation}
It follows from \ref{p2-htaft17}(iv) that
the diagram
\begin{equation}\label{p2-htaft17i}
\xymatrix{
0\ar[r]&{\bvocB}\ar[r]\ar@{=}[d]&
{\bvcF^{(r)}}\ar[r]\ar[d]^-(0.4){\bvtta^{r,r'}}&{\bvsigma^*(\txi^{-1}\tOmega^ 1_{\bvoX/\bvoS})}\ar[r]\ar[d]^-(0.4){\cdot p^{r-r'}}& 0\\
0\ar[r]&{\bvocB}\ar[r]&{\bvcF^{(r')}}\ar[r]&{\bvsigma^*(\txi^{-1}\tOmega ^1_{\bvoX/\bvoS})}\ar[r]& 0,}
\end{equation}
where the horizontal lines are the exact sequences \eqref{p2-htaft18b}
and the right vertical arrow denotes multiplication by $p^{r-r'}$, is commutative.
In particular, the extension \eqref{p2-htaft18b} is canonically isomorphic to the pullback of the extension \eqref{p2-htaft18h} by the morphism 
of multiplication by $p^r$ on $\bvsigma^*(\txi^{-1}\tOmega^1_{\bvoX/\bvoS})$.

The homomorphisms $(\alpha_n^{r,r'})_{n\in \mN}$ \eqref{p2-htaft16g}
induce a homomorphism of $\bvocB$-algebras
\begin{equation}\label{p2-htaft18f}
\bvalpha^{r,r'}\colon \bvcC^{(r)}\rightarrow \bvcC^{(r')}.
\end{equation}
For all rational numbers $r\geq r'\geq r''\geq 0$, we have
\begin{equation}\label{p2-htaft18g}
\bvtta^{r,r''}=\bvtta^{r',r''} \circ \bvtta^{r,r'} \ \ \ {\rm and}\ \ \ \bvalpha^{r, r''}=\bvalpha^{r',r''} \circ \bvalpha^{r,r'}.
\end{equation}

\subsection{}\label{p2-htaft19}
Let $r$ be a rational number $\geq 0$. By \eqref{p2-htaft8e}, we have a canonical $\bvocB$-linear isomorphism
\begin{equation}\label{p2-htaft19a}
\hupsigma^*(\txi^{-1}\tOmega^1_{\fX/\cS})\stackrel{\sim}{\rightarrow} \bvsigma^*(\txi^{-1}\tOmega^1_ {\bvoX/\bvoS}),
\end{equation}
where $\hupsigma$ is the morphism \eqref{p2-htaft8d}.
The $\bvocB$-derivation $(\delta_{\cC_n^{(r)}})_{n\in \mN}$ \eqref{p2-htaft20b} of $\bvcC^{(r)}$ and the isomorphism above define therefore a $\bvocB$-derivation  
\begin{equation}\label{p2-htaft19b}
\delta_{\bvcC^{(r)}}\colon \bvcC^{(r)}\rightarrow \hupsigma^*(\txi^{-1}\tOmega^1_{\fX/\cS})\otimes_{\bvocB}\bvcC^{(r)},
\end{equation}
which is also a Higgs $\bvocB$-field on $\bvcC^{(r)}$. 
We denote by $\mK^\bullet(\bvcC^{(r)})$ the Dolbeault complex of $(\bvcC^{(r)},\delta_{\bvcC^{(r)}})$ \eqref{p1-delta-con1b} 
and by $\tmK^\bullet(\bvcC^{( r)})$ the augmented Dolbeault complex 
\begin{equation}\label{p2-htaft19c}
\bvocB\rightarrow \mK^0(\bvcC^{(r)})\rightarrow \mK^1(\bvcC^{(r)})\rightarrow \dots
\rightarrow \mK^n(\bvcC^{(r)})\rightarrow \dots,
\end{equation}
where $\bvocB$ is placed in degree $-1$ and the differential $\bvocB\rightarrow\bvcC^{(r)}$ is the canonical homomorphism.

By \eqref{p2-htaft20c}, for all rational numbers $r\geq r'\geq 0$, we have
\begin{equation}\label{p2-htaft19d}
(\id \otimes \bvalpha^{r,r'}) \circ \delta_{\bvcC^{(r)}}=\delta_{\bvcC^{(r')}}\circ \bvalpha^{r,r'}.
\end{equation}
Therefore, the morphism $\bvalpha^{r,r'}$ induces a morphism of complexes
\begin{equation}\label{p2-htaft19e}
\upiota^{r,r'}\colon \tmK^\bullet(\bvcC^{(r)})\rightarrow \tmK^\bullet(\bvcC^{(r')}).
\end{equation}

We denote by $\mK^\bullet_\mQ(\bvcC^{(r)})$ and $\tmK^\bullet_\mQ(\bvcC^{(r)})$
the images of $\mK^\bullet(\bvcC^{(r)})$ and $\tmK^\bullet(\bvcC^{(r)})$ in $\bMod_{\mQ}(\bvocB) $ \eqref{p2-htaft27}.
We will also consider these complexes as complexes of $\bIndMod(\bvocB)$ via the functor $\upalpha_{\bvocB}$ \eqref{p2-htaft27b}.

\begin{prop}\label{p2-htaft41}
For every rational number $r\geq 0$, the functor
\begin{equation}\label{p2-htaft41a}
\begin{array}[t]{clcr}
\bIndMod(\bvocB)&\rightarrow &\bIndMod(\bvcC^{(r)})\\
M&\mapsto &M\otimes_{\bvocB} \bvcC^{(r)}
\end{array}
\end{equation}
is exact and faithful; in particular, $\bvcC^{(r)}$ is $\bvocB$-flat.
\end{prop}

By (\cite{ag2} 2.7.9.3), it is enough to prove that the functor 
\begin{equation}\label{p2-htaft41b}
\begin{array}[t]{clcr}
\bMod(\bvocB)&\rightarrow &\bMod(\bvcC^{(r)})\\
M&\mapsto &M\otimes_{\bvocB} \bvcC^{(r)}
\end{array}
\end{equation}
is exact and faithful. Since the $\co_X$-module $\tOmega^1_{X/S}$ is locally free of finite type,
the exact sequence \eqref{p2-htaft18b} is locally split. A local splitting of this sequence
induces, for every integer $m\geq 0$, a local splitting of the exact sequence \eqref{p2-htaft18c}.
We deduce from this that $\rS^m_{\bvocB}(\bvcF^{(r)})$ is $\bvocB$-flat and that the canonical homomorphism
$\bvocB\rightarrow \bvcC^{(r)}$ locally admits sections.
In view of \eqref{p2-htaft18d}, we deduce that the functor \eqref{p2-htaft41b} is exact and faithful.

\subsection{}\label{p2-htaft37}
We consider the ind-$\bvocB$-algebra \eqref{p1-indmal2}
\begin{equation}\label{p2-htaft37a}
\IC^\dagger=\underset{\underset{r\in \mQ_{>0}}{\longrightarrow}}{\mlq\mlq\lim \mrq\mrq}\ \bvcC^{(r)}. 
\end{equation}
By \eqref{p2-htaft19d}, we equip it with the derivation \eqref{p1-indmal21}
\begin{equation}\label{p2-htaft37b}
\delta_{\IC^\dagger}=\underset{\underset{r\in \mQ_{>0}}{\longrightarrow}}{\mlq\mlq\lim \mrq\mrq}\ \delta_{\bvcC^{(r)}}\colon 
\IC^\dagger \rightarrow \hupsigma^*(\txi^{-1}\tOmega^1_{\fX/\cS})\otimes_{\bvocB}\IC^\dagger.
\end{equation}
It is a Higgs $\bvocB$-field \eqref{p1-indmal20}. 
We denote by $\mK^\bullet(\IC^\dagger)$ the Dolbeault complex of $(\IC^\dagger,\delta_{\IC^\dagger})$, which is a complex of $\bIndMod(\bvocB)$ 
\eqref{p1-indmal20c}.

\begin{prop}\label{p2-htaft43}
The functor  \eqref{p1-indmal5c}
\begin{equation}\label{p2-htaft43a}
\begin{array}[t]{clcr}
\bIndMod(\bvocB)&\rightarrow &\bIndMod(\IC^\dagger)\\
M&\mapsto &M\otimes_{\bvocB} \IC^\dagger
\end{array}
\end{equation}
is exact and faithful; in particular, the ind-$\bvocB$-algebra $\IC^\dagger$ is flat {\rm (\cite{ag2} 2.7.9)}.
\end{prop}

Indeed, it follows from \ref{p2-htaft41} and (\cite{ag2} 2.7.9.1) that the functor \eqref{p2-htaft43a} is exact. The homomorphisms 
$(\bvalpha^{r,0})_{r\in \mQ_{>0}}$ \eqref{p2-htaft18f} induce 
a morphism of ind-$\bvocB$-algebras $\IC^\dagger\rightarrow \bvcC$ \eqref{p1-indmal2}. 
It then follows from \ref{p2-htaft41}, applied with $r=0$, that the functor \eqref{p2-htaft43a} is faithful.

\begin{cor}\label{p2-htaft44}
Let $\cM$, $\cN$ be two ind-$\bvocB$-modules,
\begin{equation}\label{p2-htaft44a}
\cM\otimes_{\bvocB}\IC^\dagger\stackrel{\sim}{\rightarrow} \cN\otimes_{\bvocB}\IC^\dagger
\end{equation}
an isomorphism of ind-$\IC^\dagger$-modules.  
\begin{itemize}
\item[{\rm (i)}]  If $\cN$ is rational \eqref{p2-htaft26}, then so is $\cM$.
\item[{\rm (i)}]  If $\cN$ is flat {\rm (\cite{ag2} 2.7.9)}, then so is $\cM$.
\end{itemize}
\end{cor}

(i)  Indeed, the ind-$\bvcC^{(r)}$-module $\cM\otimes_{\bvocB}\bvcC^{(r)}$ is rational by \eqref{p2-htaft44a}. 
Let $K$ (resp. \ $C$) be the kernel (resp.\ cokernel) of the morphism of multiplication by $p$ on $\cM$. 
By \ref{p2-htaft43}, $K\otimes_{\bvocB}\bvcC^{(r)}$ and $C\otimes_{\bvocB}\bvcC^{(r)}$ vanish and hence so do $K$ and $C$. 

(ii) Indeed, by \eqref{p1-indmal5f}--\eqref{p1-indmal5h},
for every $\IC^\dagger$-module $\cL$, we have a canonical isomorphism $\cL\otimes_{\IC^\dagger}(\IC^\dagger\otimes_{\bvocB}\cN)\stackrel{\sim}{\rightarrow}\cL\otimes_{\bvocB}\cN$ 
of $\bIndMod(\IC^\dagger)$. Since the canonical functor $\bIndMod(\IC^\dagger)\rightarrow \bIndMod(\bvocB)$ is exact,  
we deduce that the ind-$\IC^\dagger$-modules $\cN\otimes_{\bvocB}\IC^\dagger$ and $\cM\otimes_{\bvocB}\IC^\dagger$ are flat \eqref{p2-htaft44a}.
For every $\bvocB$-module $\cP$, we have a canonical isomorphism of $\bIndMod(\IC^\dagger)$
\begin{equation}
\IC^\dagger\otimes_{\bvocB}(\cP\otimes_{\bvocB}\cM) \stackrel{\sim}{\rightarrow} (\IC^\dagger\otimes_{\bvocB}\cP)\otimes_{\IC^\dagger}(\IC^\dagger\otimes_{\bvocB}\cM).
\end{equation}
We deduce by \ref{p2-htaft43} that the ind-$\bvocB$-module $\cM$ is flat.

\subsection{}\label{p2-htaft39}
For any rational number $r\geq 0$, we consider the $\co_\fX$-algebra \eqref{p1-NC7}
\begin{equation}\label{p2-htaft39a}
\cG^{(r)}=\rS_{\co_\fX}(p^r\txi^{-1}\tOmega^1_{\fX/\cS}),
\end{equation}
and let $\hcG^{(r)}$ be its $p$-adic completion. We denote by $\delta_{\cG^{(r)}}$ the $\co_\fX$-derivation of $\cG^{(r)}$ composed of 
\begin{equation}\label{p2-htaft39b}
\delta_{\cG^{(r)}}\colon \cG^{(r)}\rightarrow p^r\txi^{-1}\tOmega^1_{\fX/\cS} \otimes_{\co_\fX}\cG^{(r)}\rightarrow \txi^{-1}\tOmega^1_{\fX/\cS} \otimes_{\co_\fX}\cG^{(r)},
\end{equation}
where the first map is the universal $\co_\fX$-derivation of $\cG^{(r)}$ and the second map is the canonical injection, and by 
\begin{equation}\label{p2-htaft39c}
\delta_{\hcG^{(r)}}\colon \hcG^{(r)} \rightarrow \txi^{-1}\tOmega^1_{\fX/\cS} \otimes_{\co_\fX}\hcG^{(r)}
\end{equation}
its extension to the $p$-adic completions. Then $\delta_{\cG^{(r)}}$ and $\delta_{\hcG^{(r)}}$ are Higgs $\co_\fX$-fields. 
We denote by $\mK^\bullet(\hcG^{(r)})$ the Dolbeault complex of $(\hcG^{(r)},\delta_{\hcG^{(r)}})$ and by $\tmK^\bullet(\hcG^{( r)})$
the augmented Dolbeault complex 
\begin{equation}\label{p2-htaft39d}
\co_\fX\rightarrow \mK^0(\hcG^{(r)})\rightarrow \mK^1(\hcG^{(r)})\rightarrow \dots
\rightarrow \mK^n(\hcG^{(r)})\rightarrow \dots,
\end{equation}
where $\co_\fX$ is placed in degree $-1$ and the differential $\co_\fX\rightarrow\hcG^{(r)}$ is the canonical homomorphism.

For all rational numbers $r\geq r'\geq 0$, the canonical injection $p^r\txi^{-1}\tOmega^1_{\fX/\cS}\rightarrow p^{r'}\txi^{-1}\tOmega^1_{\fX/\cS}$ induces a homomorphism of $\co_\fX$-algebras 
$\ttc^{r,r'}\colon \cG^{(r)}\rightarrow \cG^{(r')}$. We denote by 
\begin{equation}\label{p2-htaft39e}
\httc^{r,r'}\colon \hcG^{(r)}\rightarrow \hcG^{(r')}
\end{equation}
its extension to the $p$-adic completions. We have
\begin{equation}\label{p2-htaft39f}
(\id \otimes \httc^{r,r'}) \circ \delta_{\hcG^{(r)}}=\delta_{\hcG^{(r')}}\circ \httc^{r,r'}.
\end{equation}
Therefore, $\httc^{r,r'}$ induces a morphism of complexes
\begin{equation}\label{p2-htaft39g}
\iota^{r,r'}\colon \tmK^\bullet(\hcG^{(r)})\rightarrow \tmK^\bullet(\hcG^{(r')}).
\end{equation}

\begin{prop}[\ref{p1-thbn36}]\label{p2-htaft40}
Suppose that the $\co_\fX$-module $\tOmega^1_{\fX/\cS}$ is free, and let $r$ and $r'$ be rational numbers such that $r>r'>0$. 
Then, there exists a rational number $\alpha\geq 0$ such that
\begin{equation}\label{p2-htaft40a}
p^\alpha\iota^{r,r'}\colon \tmK^\bullet(\hcG^{(r)})\rightarrow \tmK^\bullet(\hcG^{(r')}),
\end{equation}
where $\iota^{r,r'}$ is the morphism \eqref{p2-htaft39g}, is homotopic to $0$ by an $\co_\fX$-linear homotopy.
\end{prop}

Indeed, this is a special case of \ref{p1-thbn36} corresponding to the trivial extension $\co_\fX\oplus \txi^{-1}\tOmega^1_{\fX/\cS}$. 
The homotopy is in fact explicitly given in the proof of (\cite{agt} II.11.2).

\begin{cor}\label{p2-htaft33}
Suppose that the $\co_X$-module $\tOmega^1_{X/S}$ is free, and let $r$ and $r'$ be rational numbers such that $r>r'>0$. 
Then, there exists a rational number $\alpha\geq 0$ such that
\begin{equation}\label{p2-htaft33a}
p^\alpha \id_{\bvcC^{(r')}}\otimes \upiota^{r,r'}\colon \bvcC^{(r')}\otimes_{\bvocB}\tmK^\bullet(\bvcC^{(r)})\rightarrow \bvcC^{(r')}\otimes_{\bvocB} \tmK^\bullet(\bvcC^{(r')}),
\end{equation}
where $\upiota^{r,r'}$ is the morphism \eqref{p2-htaft19e}, is homotopic to $0$ by a $\bvcC^{(r')}$-linear homotopy.  
\end{cor}

Indeed, for every rational number $r\geq 0$, there is a canonical splitting of the extension of $\bvcC^{(r)}$-modules \eqref{p2-htaft18b}
\begin{equation}
0\rightarrow \bvcC^{(r)}\rightarrow \bvcC^{(r)}\otimes_{\bvocB}\bvcF^{(r)}\rightarrow
\bvcC^{(r)}\otimes_{\bvocB} \bvsigma^*(\txi^{-1}\tOmega^1_{\bvoX/\bvoS})\rightarrow 0,
\end{equation}
induced by the canonical morphism $\bvcF^{(r)}\rightarrow \bvcC^{(r)}$ \eqref{p2-htaft18d}. 
It induces a canonical isomorphism of $\bvcC^{(r)}$-algebras 
\begin{equation}\label{p2-htaft33b}
\varsigma^{(r)}\colon \bvcC^{(r)}\otimes_{\bvocB}\hupsigma^*(\hcG^{(r)})\stackrel{\sim}{\rightarrow} \bvcC^{(r)}\otimes_{\bvocB}\bvcC^{(r)},
\end{equation}
where $\hcG^{(r)}$ is the $\co_\fX$-algebra defined in \eqref{p2-htaft39a}. For all rational numbers $r\geq r'\geq 0$, the diagram 
\begin{equation}\label{p2-htaft33c}
\xymatrix{
{\bvcC^{(r)}\otimes_{\bvocB}\hupsigma^*(\hcG^{(r)})}\ar[r]^-(0.5){\varsigma^{(r)}}\ar[d]_{\bvalpha^{r,r'}\otimes \hupsigma^*(\httc^{r,r'})}
&{\bvcC^{(r)}\otimes_{\bvocB}\bvcC^{(r)}}\ar[d]^-(0.5){\bvalpha^{r,r'}\otimes \bvalpha^{r,r'}}\\
{\bvcC^{(r')}\otimes_{\bvocB}\hupsigma^*(\hcG^{(r')})}\ar[r]^-(0.5){\varsigma^{(r')}}&{\bvcC^{(r')}\otimes_{\bvocB}\bvcC^{(r')}}}
\end{equation}
is commutative. 
The derivation $\id_{\bvcC^{(r)}}\otimes_{\bvocB} \hupsigma^*(\delta_{\hcG^{(r)}})$ \eqref{p2-htaft39c} 
(resp.\ morphism $\id_{\bvcC^{(r')}}\otimes_{\bvocB}  \hupsigma^*(\iota^{r,r'})$ \eqref{p2-htaft39g}) corresponds via $\varsigma^{(r)}$ to $\id_{\bvcC^{(r)}}\otimes_{\bvocB} \delta_{\bvcC^{(r)}}$
(resp.\ $\id_{\bvcC^{(r')}}\otimes_{\bvocB} \upiota^{r,r'}$). The proposition follows then from \ref{p2-htaft40}.

\begin{cor}\label{p2-htaft42}
Let $r$ and $r'$ be rational numbers such that $r>r'>0$. 
Then, there exists a rational number $\alpha\geq 0$ such that for every ind-$\bvocB$-module $\cM$ and every integer $q$, 
the morphism of ind-$\bvocB$-modules
\begin{equation}\label{p2-htaft42a}
\rH^q(p^\alpha \id_{\cM}\otimes \upiota^{r,r'})\colon \rH^q(\cM\otimes_{\bvocB}\tmK^\bullet(\bvcC^{(r)}))\rightarrow \rH^q(\cM\otimes_{\bvocB} \tmK^\bullet(\bvcC^{(r')})),
\end{equation} 
where the tensor product is defined term by term \eqref{p1-indmal1d} (non derived), vanishes. 
\end{cor}

Indeed, the question being local on $X$ (\cite{ag2} 2.7.17(iii)), we may assume that the $\co_X$-module $\tOmega^1_{X/S}$ is free. 
It follows then from \ref{p2-htaft33} that 
\begin{equation}\label{p2-htaft42b}
p^\alpha \id_{\cM\otimes_{\bvocB}\bvcC^{(r')}}\otimes \upiota^{r,r'}\colon \cM\otimes_{\bvocB}\bvcC^{(r')}\otimes_{\bvocB}\tmK^\bullet(\bvcC^{(r)})\rightarrow 
\cM\otimes_{\bvocB}\bvcC^{(r')}\otimes_{\bvocB}\tmK^\bullet(\bvcC^{(r')})
\end{equation}
is homotopic to $0$ by a homotopy of ind-$\bvcC^{(r')}$-modules. Therefore, since $\bvcC^{(r')}$ is $\bvocB$-flat, the morphism of ind-$\bvcC^{(r')}$-modules
\begin{equation}\label{p2-htaft42c}
\rH^q(p^\alpha \id_{\cM}\otimes \upiota^{r,r'})\otimes_{\bvocB}\id_{\bvcC^{(r')}} \colon \rH^q(\cM\otimes_{\bvocB}\tmK^\bullet(\bvcC^{(r)}))\otimes_{\bvocB}\bvcC^{(r')} \rightarrow 
\rH^q(\cM\otimes_{\bvocB} \tmK^\bullet(\bvcC^{(r')}))\otimes_{\bvocB}\bvcC^{(r')}
\end{equation} 
vanishes. The proposition follows then from \ref{p2-htaft41}. 

\begin{cor}\label{p2-htaft38}
For every rational ind-$\bvocB$-module $\cM$ \eqref{p2-htaft26}, 
the canonical morphism of complexes of ind-$\bvocB$-modules
\begin{equation}\label{p2-htaft38a}
\cM[0]\rightarrow  \cM\otimes_{\bvocB}\mK^\bullet(\IC^\dagger),
\end{equation} 
where the tensor product is defined term by term \eqref{p1-indmal1d} (non derived), is a quasi-isomorphism. 
\end{cor}

Indeed, since $\cM$ is rational, it follows from \ref{p2-htaft42} and (\cite{ag2} 2.6.7.5) that 
\begin{equation}\label{p2-htaft38b}
\cM[0]\rightarrow \underset{\underset{r\in \mQ_{>0}}{\longrightarrow}}{\mlq\mlq\lim \mrq\mrq}\ \cM\otimes_{\bvocB}\mK^\bullet(\bvcC^{(r)})
\end{equation} 
is a quasi-isomorphism.

\begin{cor}[\cite{ag2} 4.4.36]\label{p2-htaft36}
The canonical morphism of complexes of ind-$\bvocB$-modules
\begin{equation}\label{p2-htaft36a}
\upalpha_{\bvocB}(\bvocB_\mQ)[0]\rightarrow \underset{\underset{r\in \mQ_{>0}}{\longrightarrow}}{\mlq\mlq\lim \mrq\mrq}\ 
\upalpha_{\bvocB}(\mK^\bullet_\mQ (\bvcC^{(r)}))
\end{equation}
is a quasi-isomorphism.
\end{cor}

\begin{prop}[\cite{ag2} 4.4.30]\label{p2-htaft47}
Let $r,r'$ be two rational numbers such that $r>r'>0$. Then,
\begin{itemize}
\item[{\rm (i)}] The canonical homomorphism
\begin{equation}\label{p2-htaft47a}
\co_{\fX}\rightarrow \hupsigma_*(\bvcC^{(r)})
\end{equation}
is injective. Let $\cL^{(r)}$ be its cokernel.
\item[{\rm (ii)}] There exists a rational number $a>0$ such that the morphism
\begin{equation}\label{p2-htaft47b}
\cL^{(r)}\rightarrow \cL^{(r')}
\end{equation}
induced by the homomorphism $\bvalpha^{r,r'}$ \eqref{p2-htaft18f} is annihilated by $p^a$.
\item[{\rm (iii)}] For every integer $q\geq 1$, there exists a rational number $b>0$ such that
the canonical morphism
\begin{equation}\label{p2-htaft47c}
\rR^q\hupsigma_*(\bvcC^{(r)})\rightarrow \rR^q\hupsigma_*(\bvcC^{(r')})
\end{equation}
is annihilated by $p^b$.
\end{itemize}
\end{prop}

\begin{cor}[\cite{ag2} 4.4.31]\label{p2-htaft48}
Let $r$, $r'$ be two rational numbers such that $r>r'>0$. Then,
\begin{itemize}
\item[{\rm (i)}] The canonical homomorphism
\begin{equation}\label{p2-htaft48a}
u^r\colon \co_{\fX}[\frac 1 p]\rightarrow \hupsigma_*(\bvcC^{(r)})[\frac 1 p]
\end{equation}
admits (as an $\co_\fX[\frac 1 p]$-linear morphism) a canonical left inverse
\begin{equation}\label{p2-htaft48b}
v^r\colon \hupsigma_*(\bvcC^{(r)})[\frac 1 p]\rightarrow \co_{\fX}[\frac 1 p].
\end{equation}
\item[{\rm (ii)}] The composed morphism
\begin{equation}\label{p2-htaft48c}
\hupsigma_*(\bvcC^{(r)})[\frac 1 p]\stackrel{v^r}{\longrightarrow} \co_{\fX}[\frac 1 p]
\stackrel{u^{r'}}{\longrightarrow} \hupsigma_*(\bvcC^{(r')})[\frac 1 p]
\end{equation}
is the canonical homomorphism.
\item[{\rm (iii)}] For every integer $q\geq 1$, the canonical morphism
\begin{equation}\label{p2-htaft48d}
\rR^q\hupsigma_*(\bvcC^{(r)})[\frac 1 p]\rightarrow \rR^q\hupsigma_*(\bvcC^{(r')})[\frac 1 p]
\end{equation}
is zero.
\end{itemize}
\end{cor}

\begin{cor}[\cite{ag2} 4.4.32]\label{p2-htaft49}
The canonical homomorphism
\begin{equation}\label{p2-htaft49a}
\co_{\fX}[\frac 1 p]\rightarrow \underset{\underset{r\in \mQ_{>0}}{\longrightarrow}}{\lim}\ \hupsigma_*(\bvcC^{(r )})[\frac 1p]
\end{equation}
is an isomorphism, and for every integer $q\geq 1$,
\begin{equation}\label{p2-htaft49b}
\underset{\underset{r\in \mQ_{>0}}{\longrightarrow}}{\lim}\ \rR^q\hupsigma_*(\bvcC^{(r)})[\frac 1 p] = 0.
\end{equation}
\end{cor}

\section{\texorpdfstring{Revisiting the global $p$-adic Simpson correspondence}{Revisiting the global p-adic Simpson correspondence}}\label{p2-rgpsc}

The assumptions and notation of §\ref{p2-htaft} remain in force throughout this section.

\subsection{}\label{p2-rgpsc1}
For any Higgs ind-$\bvocB$-module $(\cN,\theta)$ with coefficients in $\hupsigma^*(\txi^{-1}\tOmega^1_{\fX/\cS})$ \eqref{p1-indmal20}, 
we denote by $\uupnu(\cN,\theta)$ the ind-$\bvocB$-module 
\begin{equation}\label{p2-rgpsc1b}
\uupnu(\cN,\theta)=(\cN\otimes_{\bvocB}\IC^\dagger)^{\theta_\tot=0},
\end{equation}
where $\theta_\tot=\theta\otimes \id+\id\otimes \delta_{\IC^\dagger}$ is the total Higgs field on $\cN\otimes_{\bvocB} \IC^\dagger$ \eqref{p2-htaft37a} 
\eqref{p1-indmal20g}. It is clearly functorial in $(\cN,\theta)$.  
The Higgs field $\theta\otimes \id$ on $\cN\otimes_{\bvocB}\IC^\dagger$ induces a Higgs $\bvocB$-field $\theta_\upnu$ on $\uupnu(\cN,\theta)$ 
with coefficients in $\hupsigma^*(\txi^{-1}\tOmega^1_{\fX/\cS})$. 
Indeed, the question being local on $X$ (\cite{ag2} 2.7.17(i)), we may assume that the $\co_X$-module $\tOmega^1_{X/S}$ is free of finite type \eqref{p2-htaft1b}. 
Then $\theta$ (resp.\ $\delta_{\IC^\dagger}$) is determined by its components on a $\bvocB$-basis of $\hupsigma^*(\txi^{-1}\tOmega^1_{\fX/\cS})$, 
which are endomorphisms of the ind-$\bvocB$-module $\cN$ (resp.\ $\IC^\dagger$) that commute to each other. We thus define a functor \eqref{p2-htaft23}
\begin{equation}\label{p2-rgpsc1c}
\upnu\colon 
\begin{array}[t]{clcr}
\bIndHM(\bvocB, \hupsigma^*(\txi^{-1}\tOmega^1_{\fX/\cS}))&\rightarrow &\bIndHM(\bvocB, \hupsigma^*(\txi^{-1}\tOmega^1_{\fX/\cS})),\\
(\cN,\theta)&\mapsto& (\uupnu(\cN,\theta),\theta_\upnu). 
\end{array}
\end{equation}
We reserve the symbol $\uptau$ for a different twisting functor \eqref{p2-fhtft4a}, which will be introduced later. 
We also consider the functor 
\begin{equation}\label{p2-rgpsc1d}
\fh\colon 
\begin{array}[t]{clcr}
\bIndHM(\bvocB, \hupsigma^*(\txi^{-1}\tOmega^1_{\fX/\cS}))&\rightarrow &\bHM(\co_\fX,\txi^{-1}\tOmega^1_{\fX/\cS}),\\
(\cN,\theta)&\mapsto& \vupsigma_*(\cN\otimes_{\bvocB} \IC^\dagger,\theta_\tot),
\end{array}
\end{equation}
where the functor $\vupsigma_*$ is defined in \eqref{p2-htaft23c}. 

\subsection{}\label{p2-rgpsc19}
For any rational number $r\geq 0$, we set 
\begin{equation}\label{p2-rgpsc19a}
\delta^\vee_{\bvcC^{(r)}}=-\delta_{\bvcC^{(r)}}, \ \ \ {\rm and} \ \ \ \delta^\vee_{\IC^\dagger}=-\delta_{\IC^\dagger},
\end{equation} 
where $\delta_{\bvcC^{(r)}}$ is defined in \eqref{p2-htaft19b} and $\delta_{\IC^\dagger}$ in \eqref{p2-htaft37b}. 
For any Higgs ind-$\bvocB$-module $(\cN,\theta)$ with coefficients in $\hupsigma^*(\txi^{-1}\tOmega^1_{\fX/\cS})$, we denote by $\uupnu^\vee(\cN,\theta)$ the $\bvocB$-module 
\begin{equation}\label{p2-rgpsc19b}
\uupnu^\vee(\cN,\theta)=(\cN\otimes_{\bvocB}\IC^\dagger)^{\theta^\vee_\tot=0},
\end{equation}
where $\theta^\vee_\tot=\theta\otimes \id+\id\otimes \delta^\vee_{\IC^\dagger}$ is the total Higgs $\bvocB$-field on $\cN\otimes_{\bvocB}\IC^\dagger$. 
It is clear that $\uupnu^\vee(\cN,\theta)$ is functorial in $(\cN,\theta)$.  As in the construction of \eqref{p2-rgpsc1c}, 
we easily show that the Higgs field $\theta\otimes \id$ on $\cN\otimes_{\bvocB}\IC^\dagger$ induces a Higgs 
$\bvocB$-field $\theta_{\upnu^\vee}$ on $\uupnu^\vee(\cN,\theta)$ with coefficients in 
$\hupsigma^*(\txi^{-1}\tOmega^1_{\fX/\cS})$. We thus define a functor 
\begin{equation}\label{p2-rgpsc19c}
\upnu^\vee\colon 
\begin{array}[t]{clcr}
\bIndHM(\bvocB,\hupsigma^*(\txi^{-1}\tOmega^1_{\fX/\cS})) &\rightarrow& \bIndHM(\bvocB,\hupsigma^*(\txi^{-1}\tOmega^1_{\fX/\cS})),\\
(\cN,\theta)&\mapsto&(\uupnu^\vee(\cN,\theta), \theta_{\upnu^\vee}).
\end{array}
\end{equation}

\begin{defi}\label{p2-rgpsc4}
We say that a Higgs ind-$\bvocB$-module $(\cN,\theta)$ with coefficients in $\hupsigma^*(\txi^{-1}\tOmega^1_{\fX/\cS})$ is  
{\em weakly twistable by the extension \eqref{p2-htaft18h}} (or simply {\em weakly twistable} if there is no risk of ambiguity)  
if it is rational \eqref{p2-htaft26} and if the canonical morphism of ind-$\IC^\dagger$-modules \eqref{p1-indmal5c}
\begin{equation}\label{p2-rgpsc4a}
\IC^\dagger \otimes_{\bvocB}\uupnu(\cN,\theta)\rightarrow \IC^\dagger \otimes_{\bvocB}\cN
\end{equation}
is an isomorphism.  
\end{defi}

We keep the adjective {\em twistable} for a stronger notion introduced in \ref{p2-rgpsc7}. 

\begin{rema}\label{p2-rgpsc5}
We keep the assumptions of \ref{p2-rgpsc4} and let $(P_i)_{1\leq i\leq n}$ be a finite covering of the final object of $\tE_s^{\mN^\circ}$. 
By (\cite{ag2} 2.7.17(iii)), the morphism \eqref{p2-rgpsc4a} is an isomorphism if and only if its restriction over $P_i$ is an isomorphism for every $1\leq i\leq n$. 
In particular, in view of (\cite{ag2} 4.9.2), the property of being weakly twistable for a Higgs ind-$\bvocB$-module is local on the topos $X_{s,\et}$, via the composed morphism 
\begin{equation}
\xymatrix{
{\tE_s^{\mN^\circ}}\ar[r]^-(0.5){\uplambda}&{\tE_s}\ar[r]^-(0.5){\sigma_s}&{X_{s,\et}},}
\end{equation}
where the first morphism is defined in \eqref{p2-ncgt5a} and the second in \eqref{p2-htaft5c}. 
\end{rema}

\begin{rema}\label{p2-rgpsc32}
We take again the notation of \ref{p2-htaft23}. 
For every Higgs ind-$\bvocB$-module $(\cN,\theta)$ with coefficients in $\hupsigma^*(\txi^{-1}\tOmega^1_{\fX/\cS})$, the canonical morphism 
\begin{equation}\label{p2-rgpsc32a}
\uupnu(N,\theta)\rightarrow \IC^\dagger \otimes_{\bvocB}\cN
\end{equation}
is a morphism of ind-$\IndSym_{\bvocB}(\hupsigma^*(\cT))$-modules, where the ind-module structure on $\cN$ (resp.\  $\uupnu(\cN,\theta)$, resp.\ $\IC^\dagger$) 
corresponds to the Higgs field $\theta$ (resp.\ $\theta_\upnu$ \eqref{p2-rgpsc1c}, resp.\ $0$). Since the multiplications by 
$\IC^\dagger$ and $\IndSym_{\bvocB}(\hupsigma^*(\cT))$ on $\IC^\dagger \otimes_{\bvocB}\cN$ 
(resp.\ $\IC^\dagger \otimes_{\bvocB}\uupnu(\cN,\theta)$) commute, we deduce that 
the canonical morphism of ind-$\IC^\dagger$-modules
\begin{equation}\label{p2-rgpsc32b}
\IC^\dagger \otimes_{\bvocB}\uupnu(\cN,\theta)\rightarrow \IC^\dagger \otimes_{\bvocB}\cN
\end{equation}
is also a morphism of ind-$\IndSym_{\bvocB}(\hupsigma^*(\cT))$-modules. 
\end{rema}

\begin{prop}\label{p2-rgpsc18}
Let $(\cN,\theta)$ be a rational Higgs ind-$\bvocB$-module with coefficients in $\hupsigma^*(\txi^{-1}\tOmega^1_{\fX/\cS})$, $\cN^\vee$ an ind-$\bvocB$-module, 
\begin{equation}\label{p2-rgpsc18a}
\psi\colon \cN^\vee\otimes_{\bvocB}\IC^\dagger\stackrel{\sim}{\rightarrow}\cN\otimes_{\bvocB}\IC^\dagger
\end{equation}
an isomorphism of ind-$\IC^\dagger$-modules with $\delta_{\IC^\dagger}$-connection, 
where the $\delta_{\IC^\dagger}$-connections are defined as in \ref{p1-delta-con9}, 
$\cN$ (resp.\ $\cN^\vee$) being endowed with the Higgs field $\theta$ (resp.\ $0$).  Then, 
\begin{itemize}
\item[{\rm (i)}] The isomorphism $\psi$ induces an isomorphism of ind-$\bvocB$-modules 
$\cN^\vee\stackrel{\sim}{\rightarrow} \uupnu(\cN,\theta)$, where $\uupnu$ is the functor \eqref{p2-rgpsc1b}.
We deduce a canonical Higgs $\bvocB$-field $\theta^\vee$ on  $\cN^\vee$ with coefficients in $\hupsigma^*(\txi^{-1}\tOmega^1_{\fX/\cS})$, 
so that we have an isomorphism of Higgs ind-modules
\begin{equation}\label{p2-rgpsc18b}
(\cN^\vee,\theta^\vee)\stackrel{\sim}{\rightarrow} \upnu(\cN,\theta),
\end{equation}
where $\upnu$ is the functor \eqref{p2-rgpsc1c}. 
\item[{\rm (ii)}] The morphism $\psi$ is an isomorphism of ind-$\IC^\dagger$-modules with $\delta^\vee_{\IC^\dagger}$-connection \eqref{p2-rgpsc19a}, 
where the $\delta^\vee_{\IC^\dagger}$-connections are defined as in \ref{p1-delta-con9}, 
$\cN$ (resp.\ $\cN^\vee$) being endowed with the Higgs field $0$ (resp.\ $\theta^\vee$). 
\item[{\rm (iii)}] The isomorphism $\psi$ induces an isomorphism of Higgs ind-$\bvocB$-modules with coefficients in $\hupsigma^*(\txi^{-1}\tOmega^1_{\fX/\cS})$, 
\begin{equation}\label{p2-rgpsc18c}
(\cN,\theta)\stackrel{\sim}{\rightarrow} \upnu^\vee(\cN^\vee,\theta^\vee),
\end{equation}
where $\upnu^\vee$ is the functor \eqref{p2-rgpsc19c}. 
\end{itemize}
\end{prop}

(i) Indeed, the ind-$\bvocB$-module $\cN^\vee$ is rational by \ref{p2-htaft44}. The proposition follows then from \ref{p2-htaft38}. 

(ii) It follows immediately from the fact that the Higgs field $\theta^\vee$ on $\cN^\vee$ 
is induced by $\id\otimes\delta^\vee_{\IC^\dagger}$ on $\cN\otimes_{\bvocB}\IC^\dagger$, which is a $\delta^\vee_{\IC^\dagger}$-connection.

(iii) It  follows from (ii) and \ref{p2-htaft38} that $\psi^{-1}$ induces an isomorphism of the ind-$\bvocB$-modules underlying \eqref{p2-rgpsc18c}.  
 Let $u\colon \cN\rightarrow \cN\otimes_{\bvocB}\IC^\dagger$ be the canonical injection. Since $\theta^\vee_\tot\circ \psi^{-1}\circ u=0$, we have 
\begin{eqnarray}
(\theta^\vee\otimes \id)\circ (\psi^{-1}\circ u)&=&(\id\otimes\delta_{\IC^\dagger})\circ (\psi^{-1}\circ u)\nonumber \\
&=&(\theta\otimes \id+\id\otimes\delta_{\IC^\dagger})\circ u\nonumber \\
&=&\theta\otimes\id,
\end{eqnarray}
the second equality being a consequence of the assumptions. Therefore, \eqref{p2-rgpsc18c} is an isomorphism of Higgs ind-$\bvocB$-modules.

\subsection{}\label{p2-rgpsc2}
Let $r$ be a rational number $\geq 0$. We denote by $\bIndMIC(\bvcC^{(r)}/\bvocB)$ the category of ind-$\bvcC^{(r)}$-modules with integrable 
$\delta_{\bvcC^{(r)}}$-connection \eqref{p1-indmal22}.
By \ref{p1-delta-con9}, we have a functor 
\begin{equation}\label{p2-rgpsc2a}
\uppi^{(r)}\colon \begin{array}[t]{clcr}
\bIndHM(\bvocB, \hupsigma^*(\txi^{-1}\tOmega^1_{\fX/\cS}))&\rightarrow &\bIndMIC(\bvcC^{(r)}/\bvocB),\\
(\cN,\theta)&\mapsto& (\cN\otimes_{\bvocB}\bvcC^{(r)},\theta\otimes_{\bvocB}\id_{\bvcC^{(r)}}+\id_{\cN}\otimes_{\bvocB}\delta_{\bvcC^{(r)}}). 
\end{array}
\end{equation}
Restricting this functor to Higgs ind-$\bvocB$-modules with zero Higgs fields, we get a functor
\begin{equation}\label{p2-rgpsc2b}
\fS^{(r)}\colon
\begin{array}[t]{clcr}
\bIndMod(\bvocB)&\rightarrow &\bIndMIC(\bvcC^{(r)}/\bvocB),\\
\cM&\mapsto& (\bvcC^{(r)}\otimes_{\bvocB}\cM,\delta_{\bvcC^{(r)}}\otimes_{\bvocB}\id_{\cM}). 
\end{array}
\end{equation}

Let $r\geq r'\geq 0$ be rational numbers, $(\cM,\nabla)$ an ind-$\bvcC^{(r)}$-module with integrable $\delta_{\bvcC^{(r)}}$-connection. 
It follows from \ref{p1-delta-con10} that there exists a unique integrable $\delta_{\bvcC^{(r')}}$-connection
\begin{equation}\label{p2-rgpsc2c}
\nabla'\colon \cM\otimes_{\bvcC^{(r)}}\bvcC^{(r')}\rightarrow \hupsigma^*(\txi^{-1}\tOmega^1_{\fX/\cS})\otimes_{\bvocB}\cM\otimes_{\bvcC^{(r)}}\bvcC^{(r')}
\end{equation}
that fits into a commutative diagram 
\begin{equation}\label{p2-rgpsc2d}
\xymatrix{
{\cM\otimes_{\bvocB}\bvcC^{(r')}}\ar[rrr]^-(0.5){\nabla\otimes\id_{\bvcC^{(r')}}+\id_\cM\otimes \delta_{\bvcC^{(r')}}}\ar[d]
&&&{\hupsigma^*(\txi^{-1}\tOmega^1_{\fX/\cS})\otimes_{\bvocB}\cM\otimes_{\bvocB}\bvcC^{(r')}}\ar[d]\\
{\cM\otimes_{\bvcC^{(r)}}\bvcC^{(r')}}\ar[rrr]^-(0.5){\nabla'}&&&{\hupsigma^*(\txi^{-1}\tOmega^1_{\fX/\cS})\otimes_{\bvocB}\cM\otimes_{\bvcC^{(r)}}\bvcC^{(r')},}}
\end{equation}
where the vertical arrows are the canonical morphisms. We thus define a functor 
\begin{equation}\label{p2-rgpsc2e}
\varepsilon^{r,r'}\colon 
\begin{array}[t]{clcr}
\bIndMIC(\bvcC^{(r)}/\bvocB)&\rightarrow& \bIndMIC(\bvcC^{(r')}/\bvocB),\\
(\cM,\nabla)&\mapsto&(\cM\otimes_{\bvcC^{(r)}}\bvcC^{(r')},\nabla').
\end{array}
\end{equation} 
We have a canonical isomorphism of functors from the category $\bIndHM(\bvocB, \hupsigma^*(\txi^{-1}\tOmega^1_{\fX/\cS}))$ to the category $\bIndMIC(\bvcC^{(r')}/\bvocB)$, 
\begin{equation}\label{p2-rgpsc2f}
\varepsilon^{r,r'}\circ \uppi^{(r)}\stackrel{\sim}{\rightarrow} \uppi^{(r')}.
\end{equation}

We denote by 
\begin{equation}\label{p2-rgpsc2g}
\upeta^{(r)}\colon \bIndMIC(\bvcC^{(r)}/\bvocB)\rightarrow \bIndHM(\bvocB, \hupsigma^*(\txi^{-1}\tOmega^1_{\fX/\cS}))
\end{equation}
the forgetful functor. By \eqref{p2-rgpsc2d}, the homomorphism $\bvalpha^{r,r'}$ induces a canonical morphism of functors 
from $\bIndMIC(\bvcC^{(r)}/\bvocB)$ to $\bIndHM(\bvocB, \hupsigma^*(\txi^{-1}\tOmega^1_{\fX/\cS}))$,
\begin{equation}\label{p2-rgpsc2h}
\upeta^{(r)} \rightarrow \upeta^{(r')} \circ \varepsilon^{r,r'}. 
\end{equation}
The morphism \eqref{p2-rgpsc2h} and the isomorphism \eqref{p2-rgpsc2f} induce a canonical morphism of functors 
on the category $\bIndHM(\bvocB, \hupsigma^*(\txi^{-1}\tOmega^1_{\fX/\cS}))$, 
\begin{equation}\label{p2-rgpsc2i}
\upeta^{(r)} \circ \uppi^{(r)}\rightarrow  \upeta^{(r')} \circ \uppi^{(r')}.
\end{equation}

\subsection{}\label{p2-rgpsc3}
Let $r$ be rational number $\geq 0$. 
For any Higgs ind-$\bvocB$-module $(\cN,\theta)$ with coefficients in $\hupsigma^*(\txi^{-1}\tOmega^1_{\fX/\cS})$, 
we write $\uppi^{(r)}(\cN,\theta)=(\cN\otimes_{\bvocB}\bvcC^{(r)},\nabla^{(r)})$ \eqref{p2-rgpsc2a} and set 
\begin{equation}\label{p2-rgpsc3a}
\uupnu^{(r)}(\cN,\theta)=(\cN\otimes_{\bvocB}\bvcC^{(r)})^{\nabla^{(r)}=0}. 
\end{equation}
As in the construction of \eqref{p2-rgpsc1c}, 
we easily show that the Higgs field $\theta\otimes \id$ on $\cN\otimes_{\bvocB}\bvcC^{(r)}$ induces a Higgs $\bvocB$-field $\theta^{(r)}$ on $\uupnu^{(r)}(\cN,\theta)$ 
with coefficients in $\hupsigma^*(\txi^{-1}\tOmega^1_{\fX/\cS})$. We thus define a functor 
\begin{equation}\label{p2-rgpsc3b}
\upnu^{(r)}\colon \begin{array}[t]{clcr}
\bIndHM(\bvocB, \hupsigma^*(\txi^{-1}\tOmega^1_{\fX/\cS}))&\rightarrow &\bIndHM(\bvocB, \hupsigma^*(\txi^{-1}\tOmega^1_{\fX/\cS})),\\
(\cN,\theta)&\mapsto& (\uupnu^{(r)}(\cN,\theta),\theta^{(r)}). 
\end{array}
\end{equation}
Similarly, we consider the functor 
\begin{equation}\label{p2-rgpsc3e}
\fh^{(r)}=\vupsigma_*\circ \upeta^{(r)}\circ \uppi^{(r)}\colon \begin{array}[t]{clcr}
\bIndHM(\bvocB, \hupsigma^*(\txi^{-1}\tOmega^1_{\fX/\cS}))&\rightarrow &\bHM(\co_\fX, \txi^{-1}\tOmega^1_{\fX/\cS}),\\
(\cN,\theta)&\mapsto&  \vupsigma_*(\cN\otimes_{\bvocB}\bvcC^{(r)},\nabla^{(r)}),
\end{array}
\end{equation}
where the functor $\vupsigma_*$ is defined in \eqref{p2-htaft23c}.

For all rational numbers $r\geq r'\geq 0$, the morphism \eqref{p2-rgpsc2i} induces a canonical morphism of functors 
on the category $\bIndHM(\bvocB, \hupsigma^*(\txi^{-1}\tOmega^1_{\fX/\cS}))$, 
\begin{equation}\label{p2-rgpsc3c}
\upnu^{(r)}\rightarrow \upnu^{(r')},
\end{equation}
and a canonical morphism of functors 
from $\bIndHM(\bvocB, \hupsigma^*(\txi^{-1}\tOmega^1_{\fX/\cS}))$ to $\bHM(\co_\fX, \txi^{-1}\tOmega^1_{\fX/\cS})$, 
\begin{equation}\label{p2-rgpsc3f}
\fh^{(r)}\rightarrow \fh^{(r')}. 
\end{equation}
We thus obtain small filtered direct systems $(\upnu^{(r)}(\cN,\theta))_{r\geq 0}$ and $(\fh^{(r)}(\cN,\theta))_{r\geq 0}$.
By (\cite{ag2} 2.7.3 and 2.6.7.5), we have a functorial isomorphism
\begin{equation}\label{p2-rgpsc3d}
\upnu(\cN,\theta)\stackrel{\sim}{\rightarrow}
\underset{\underset{r\in \mQ_{>0}}{\longrightarrow}}{\mlq\mlq\lim \mrq\mrq}\ \upnu^{(r)}(\cN,\theta),
\end{equation}
where $\upnu$ is the functor \eqref{p2-rgpsc1c}. 
By the definition of the functor $\rI\hupsigma_*$ \eqref{p2-htaft23b} and (\cite{ag2} 2.6.4(ii)), we have a functorial isomorphism
\begin{equation}\label{p2-rgpsc3g}
\fh(\cN,\theta)\stackrel{\sim}{\rightarrow}
\underset{\underset{r\in \mQ_{>0}}{\longrightarrow}}{\lim}\ \fh^{(r)}(\cN,\theta),
\end{equation}
where $\fh$ is the functor \eqref{p2-rgpsc1d}.

\begin{defi}\label{p2-rgpsc6}
Let $\cM$ be an ind-$\bvocB$-module,
$(\cN,\theta)$ a Higgs ind-$\bvocB$-module with coefficients in $\hupsigma^*(\txi^{-1}\tOmega^1_{\fX/\cS})$.
\begin{itemize}
\item[(i)] Let $r$ be a rational number $>0$. We say that $\cM$ and $(\cN,\theta)$ are {\em $r$-associated}
if there exists an isomorphism of $\bIndMIC(\bvcC^{(r)}/\bvocB)$
\begin{equation}\label{p2-rgpsc6a}
\alpha\colon \fS^{(r)}(\cM)\stackrel{\sim}{\rightarrow} \uppi^{(r)}(\cN,\theta). 
\end{equation}
We then also say that the triple $(\cM,(\cN,\theta),\alpha)$ is {\em $r$-admissible}.
\item[(ii)] We say that $\cM$ and $(\cN,\theta)$ are {\em associated} if there exists a rational number $r>0$ such that
$\cM$ and $(\cN,\theta)$ are $r$-associated.
\end{itemize}
\end{defi}

Note that for all rational numbers  $r\geq r'>0$,
if $\cM$ and $(\cN,\theta)$ are $r$-associated, they are $r'$-associated, in view of \eqref{p2-rgpsc2f}.

\begin{defi}\label{p2-rgpsc7}
We say that a Higgs ind-$\bvocB$-module with coefficients in $\hupsigma^*(\txi^{-1}\tOmega^1_{\fX/\cS})$ is 
{\em twistable by the extension \eqref{p2-htaft18h}} (or simply {\em twistable} if there is no risk of ambiguity) 
if it is rational \eqref{p2-htaft26} and associated with an ind-$\bvocB$-module in the sense of \ref{p2-rgpsc6}(ii). 
\end{defi}

\begin{prop}\label{p2-rgpsc9}
Every Higgs ind-$\bvocB$-module $(\cN,\theta)$ with coefficients in $\hupsigma^*(\txi^{-1}\tOmega^1_{\fX/\cS})$, 
which is twistable by the extension \eqref{p2-htaft18h} in the sense of \ref{p2-rgpsc7},
is weakly twistable by the extension \eqref{p2-htaft18h} in the sense of \ref{p2-rgpsc4}. 
\end{prop}

Indeed, there exist an ind-$\bvocB$-module $\cM$, a rational number $r>0$ and an isomorphism of $\bIndMIC(\bvcC^{(r)}/\bvocB)$
\begin{equation}\label{p2-rgpsc9a}
\fS^{(r)}(\cM)\stackrel{\sim}{\rightarrow}\uppi^{(r)}(\cN,\theta). 
\end{equation}
By \ref{p1-delta-con10}, the latter induces by extension of scalars by the canonical morphism of ind-$\bvocB$-algebras 
$\bvcC^{(r)}\rightarrow \IC^\dagger$, an isomorphism of ind-$\IC^\dagger$-modules with $\delta_{\IC^\dagger}$-connection \eqref{p2-htaft37b},
\begin{equation}\label{p2-rgpsc9b}
\IC^\dagger \otimes_{\bvocB}\cM\stackrel{\sim}{\rightarrow}\IC^\dagger \otimes_{\bvocB} \cN, 
\end{equation}
where the $\delta_{\IC^\dagger}$-connections are defined as in \ref{p1-delta-con9}, 
$\cN$ (resp.\ $\cM$) being endowed with the Higgs field $\theta$ (resp.\ $0$). 
By \ref{p2-rgpsc18}(i), \eqref{p2-rgpsc9b} induces an isomorphism of ind-$\bvocB$-modules
\begin{equation}\label{p2-rgpsc9c}
\cM \stackrel{\sim}{\rightarrow} \uupnu(\cN,\theta).
\end{equation}
Moreover, the isomorphism \eqref{p2-rgpsc9b} identifies with the canonical morphism \eqref{p2-rgpsc4a}. 
Therefore, $(\cN,\theta)$ is weakly twistable.

\begin{cor}\label{p2-rgpsc10}
Every Higgs ind-$\bvocB$-module $(\cN,\theta)$ with coefficients in $\hupsigma^*(\txi^{-1}\tOmega^1_{\fX/\cS})$ which is twistable, 
is associated with the ind-$\bvocB$-module $\uupnu(\cN,\theta)$ \eqref{p2-rgpsc1b}.
\end{cor}

\begin{defi}\label{p2-rgpsc11}
Let $\cM$ be an ind-$\bvocB$-module,
$(N,\theta)$ a coherent Higgs $\co_\fX[\frac 1 p]$-module with coefficients in $\txi^{-1}\tOmega^1_{\fX/\cS}$ \eqref{p2-htaft45}, 
that we also consider as a Higgs ind-$\co_\fX$-module with coefficients in $\txi^{-1}\tOmega^1_{\fX/\cS}$ \eqref{p2-cmupiso31h}, 
$r$ a rational number $>0$. We say that $\cM$ and $(N,\theta)$ are {\em associated} (resp.\ {\em $r$-associated}) 
if so are $\cM$ and $\rI\hupsigma^*(N,\theta)$ in the sense of \ref{p2-rgpsc6}, see \eqref{p2-cmupiso31h} and \eqref{p2-htaft23a}.  
\end{defi}

\begin{defi}\label{p2-rgpsc12}\
\begin{itemize}
\item[(i)] We say that an ind-$\bvocB$-module is {\em Dolbeault} if it is associated with a Higgs $\co_\fX[\frac 1 p]$-bundle with coefficients in 
$\txi^{-1}\tOmega^1_{\fX/\cS}$ \eqref{p1-delta-con6},
in the sense of \ref{p2-rgpsc11}.
\item[(ii)] We say that a Higgs $\co_\fX[\frac 1 p]$-bundle with coefficients in $\txi^{-1}\tOmega^1_{\fX/\cS} $  is {\em solvable}
if it is associated with an ind-$\bvocB$-module.  
\end{itemize}
\end{defi}

We denote by $\bIndMod^\Dolb(\bvocB)$ the full sub-category of $\bIndMod(\bvocB)$
made up  of the Dolbeault ind-$\bvocB$-modules, and by $\bHM^\sol(\co_\fX[\frac 1 p], \txi^{-1}\tOmega^1_{\fX/\cS })$
the full subcategory of $\bHM(\co_\fX[\frac 1 p], \txi^{-1}\tOmega^1_{\fX/\cS})$
made up  of the solvable Higgs $\co_\fX[\frac 1 p]$-bundles with coefficients in $\txi^{-1}\tOmega^1_{\fX/\cS}$.

\subsection{}\label{p2-rgpsc14}
We consider the functors
\begin{equation}\label{p2-rgpsc14a}
\cH\colon 
\begin{array}[t]{clcr}
\bIndMod(\bvocB)&\rightarrow &\bHM(\co_\fX,\txi^{-1}\tOmega^1_{\fX/\cS}),\\
\cM&\mapsto& \fh(\cM,0),
\end{array}
\end{equation}
where the functor $\fh$ is defined in \eqref{p2-rgpsc1d}, and 
\begin{equation}\label{p2-rgpsc14b}
\cV\colon 
\begin{array}[t]{clcr}
\bHM^\coh(\co_\fX[\frac 1 p],\txi^{-1}\tOmega^1_{\fX/\cS})&\rightarrow &\bIndMod(\bvocB),\\
(N,\theta)&\mapsto& \uupnu(\rI\hupsigma^*(N,\theta)),
\end{array}
\end{equation}
where the functor $\rI\hupsigma^*$ (resp.\ $\uupnu$) is defined in \eqref{p2-htaft23a} (resp.\ \eqref{p2-rgpsc1b}) and 
$(N,\theta)$ is considered as a Higgs ind-$\co_\fX$-module via the functor \eqref{p2-cmupiso31h}. 

\begin{prop}[\cite{ag2} 4.5.12]\label{p2-rgpsc33}
For every Dolbeault ind-$\bvocB$-module $\cM$,
$\cH(\cM)$ \eqref{p2-rgpsc14a} is a solvable Higgs $\co_\fX[\frac 1 p]$-bundle associated with $\cM$ \eqref{p2-rgpsc11}.
\end{prop}

\begin{prop}[\cite{ag2} 4.5.18]\label{p2-rgpsc34}
For every solvable Higgs $\co_\fX[\frac 1 p]$-bundle $(N,\theta)$ with coefficients in $\txi^{-1}\tOmega^1_{\fX/\cS}$,
the ind-$\bvocB$-module $\cV(N,\theta)$ \eqref{p2-rgpsc14b} is  Dolbeault, associated with $\cN$ \eqref{p2-rgpsc11}.
\end{prop}

\begin{teo}[\cite{ag2} 4.5.20]\label{p2-rgpsc15}
The functors $\cH$ \eqref{p2-rgpsc14a} and $\cV$ \eqref{p2-rgpsc14b} induce equivalences of categories quasi-inverse to each other
\begin{equation}\label{p2-rgpsc15a}
\xymatrix{
{\bIndMod^\Dolb(\bvocB)}\ar@<1ex>[r]^-(0.5){\cH}&{\bHM^\sol(\co_\fX[\frac 1 p], \txi^{-1}\tOmega^1_{\fX/\cS}).}
\ar@<1ex>[l]^-(0.5){\cV}}
\end{equation}
\end{teo}

\begin{teo}[\cite{ag2} 4.7.4]\label{p2-rgpsc16}
Let $\cM$ be a Dolbeault ind-$\bvocB$-module.
We denote by $\mK^\bullet(\cH(\cM))$ the Dolbeault complex of the Higgs $\co_\fX[\frac 1 p]$-bundle
$\cH(\cM)$ \eqref{p2-rgpsc14a}. Then, we have a canonical functorial isomorphism of $\bD^+(\bMod(\co_\fX))$
\begin{equation}\label{p2-rgpsc16a}
\rR\vupsigma_*(\cM)\stackrel{\sim}{\rightarrow}\mK^\bullet(\cH(\cM)),
\end{equation}
where $\rR\vupsigma_*$ is the functor \eqref{p2-htaft27g}.
\end{teo}

\begin{defi}\label{p2-rgpsc20}\
\begin{itemize}
\item[(i)] A $\bvocB_\mQ$-module is said to be {\em Dolbeault} if it is so as an ind-$\bvocB$-module \eqref{p2-htaft27b}, i.e., if it is associated with
a Higgs $\co_\fX[\frac 1 p]$-bundle with coefficients in $\txi^{-1}\tOmega^1_{\fX/\cS}$ \eqref{p2-rgpsc12}.
\item[(ii)] A $\bvocB_\mQ$-module is said to be {\em strongly Dolbeault} if it is Dolbeault and adic of finite type (\cite{ag2} 4.3.14).
\item[(iii)] A Higgs $\co_\fX[\frac 1 p]$-bundle with coefficients in $\txi^{-1}\tOmega^1_{\fX/\cS}$ \eqref{p1-delta-con6}
is said to be {\em rationally solvable} (resp.\ {\em strongly solvable})
if it is associated with a $\bvocB_\mQ$-module (resp.\ an adic $\bvocB_\mQ$-module of finite type).
\end{itemize}
\end{defi}

We denote by $\bMod_\mQ^\Dolb(\bvocB)$ (resp.\ $\bMod_\mQ^\sDolb(\bvocB)$) the full subcategory of $\bMod_\mQ(\bvocB)$
made up  of the Dolbeault (resp.\ strongly Dolbeault) $\bvocB_\mQ$-modules, and by $\bHM^\qsol(\co_\fX[\frac 1 p], \txi^{-1}\tOmega^ 1_{\fX/\cS})$
(resp.\ $\bHM^\ssol(\co_\fX[\frac 1 p], \txi^{-1}\tOmega^1_{\fX/\cS})$)
the full subcategory of $\bHM(\co_\fX[\frac 1 p], \txi^{-1}\tOmega^1_{\fX/\cS})$
made up  of the rationally (resp.\ strongly) solvable Higgs $\co_\fX[\frac 1 p]$-bundles with coefficients in $\txi^{-1}\tOmega^1_{\fX/\cS}$.
These notions depend a priori on the deformation $\tf$ fixed in \eqref{p2-htaft1c}.

\subsection{}\label{p2-rgpsc21}
The composed functor $\cH\circ \upalpha_{\bvocB}$ \eqref{p2-rgpsc14a} induces clearly a functor
\begin{equation}\label{p2-rgpsc21a}
\cH_\mQ\colon \bMod_\mQ(\bvocB)\rightarrow \bHM(\co_\fX[\frac 1 p], \txi^{-1}\tOmega^1_{\fX/\cS})
\end{equation}
that fits into a strictly commutative diagram
\begin{equation}\label{p2-rgpsc21b}
\xymatrix{
{\bMod_\mQ(\bvocB)}\ar[r]^-(0.5){\cH_\mQ}\ar[d]_{\upalpha_{\bvocB}}&{\bHM(\co_\fX[ \frac 1 p], \txi^{-1}\tOmega^1_{\fX/\cS})}\ar[d]\\
{\bIndMod(\bvocB)}\ar[r]^-(0.5){\cH}&{\bHM(\co_\fX, \txi^{-1}\tOmega^1_{\fX/\cS} ),}}
\end{equation}
where the vertical arrows are the canonical functors.

\begin{prop}[\cite{ag2} 4.6.18 and 4.6.21]\label{p2-rgpsc22}
There is a canonical functor 
\begin{equation}\label{p2-rgpsc22a}
\cV_\mQ\colon \bHM^\qsol(\co_\fX[\frac 1 p], \txi^{-1}\tOmega^1_{\fX/\cS})\rightarrow \bMod^\Dolb_\mQ(\bvocB)
\end{equation}
that fits into a commutative diagram up to canonical isomorphism, 
\begin{equation}\label{p2-rgpsc22b}
\xymatrix{
{\bHM^\qsol(\co_\fX[\frac 1 p], \txi^{-1}\tOmega^1_{\fX/\cS})}\ar[r]^-(0.5){\cV_\mQ}\ar[d]&{\bMod^\Dolb_\mQ(\bvocB)}\ar[d]^{\upalpha_{\bvocB}}\\
{\bHM^\coh(\co_\fX[\frac 1 p], \txi^{-1}\tOmega^1_{\fX/\cS})}\ar[r]^-(0.5){\cV}&{\bIndMod^\Dolb(\bvocB),}}
\end{equation}
where $\cV$ is the functor \eqref{p2-rgpsc14b} and the vertical arrows are the canonical functors.
Moreover, it induces a functor that we denote again by
\begin{equation}\label{p2-rgpsc22c}
\cV_\mQ\colon \bHM^\ssol(\co_\fX[\frac 1 p], \txi^{-1}\tOmega^1_{\fX/\cS})\rightarrow \bMod^\sDolb_\mQ(\bvocB).
\end{equation}
\end{prop}

Since the functor $\upalpha_{\bvocB}$ \eqref{p2-htaft27b} is fully faithful, the functor $\cV_\mQ$ \eqref{p2-rgpsc22a} is determined by the diagram \eqref{p2-rgpsc22b}, up to isomorphism. 
We gave in (\cite{ag2} 4.6.18) an explicit construction. 

\begin{prop}[\cite{ag2} 4.6.22 and 4.6.23]\label{p2-rgpsc23}
The functors $\cH_\mQ$ \eqref{p2-rgpsc21a} and $\cV_\mQ$ \eqref{p2-rgpsc22a} induce equivalences of categories quasi-inverse to each other
\begin{equation}\label{p2-rgpsc23a}
\xymatrix{
{\bMod^\Dolb_\mQ(\bvocB)}\ar@<1ex>[r]^-(0.5){\cH_\mQ}&{\bHM^\qsol(\co_\fX[\frac 1 p], \txi^{-1}\tOmega^1_{\fX/\cS}),}
\ar@<1ex>[l]^-(0.5){\cV_\mQ}}
\end{equation}
and 
\begin{equation}\label{p2-rgpsc23b}
\xymatrix{
{\bMod^\sDolb_\mQ(\bvocB)}\ar@<1ex>[r]^-(0.5){\cH_\mQ}&{\bHM^\ssol(\co_\fX[\frac 1 p], \txi^{-1}\tOmega^1_{\fX/\cS}).}
\ar@<1ex>[l]^-(0.5){\cV_\mQ}}
\end{equation}
\end{prop}

\begin{prop}[\cite{ag2} 4.10.4 and 4.10.8]\label{p2-htaft57}
The property for an ind-$\bvocB$-module (resp.\ a $\bvocB_\mQ$-module) to be Dolbeault 
does not depend on the choice of the deformation $\tf$ fixed in \eqref{p2-htaft1c}
provided that we remain in the same setting, absolute or relative \eqref{p2-ncgt3}.
\end{prop}

The analogous proposition for the property of being solvable (resp.\ rationally solvable, resp.\ strongly solvable) is proved under a certain condition (\cite{ag2} 4.10.12).

\begin{prop}\label{p2-htaft58}
Suppose that $X$ is an object of $\bP$ \eqref{p2-htaft9} and that $X_s$ is non-empty.
Let $(N,\theta)$ be a CL-small Higgs $\co_\fX[\frac 1 p]$-module with coefficients in $\txi^{-1}\tOmega^1_{\fX/\cS}$ \eqref{p1-tshbn13},
that we consider as a Higgs ind-$\co_\fX$-module by the functor \eqref{p2-cmupiso31h}. 
Then, the Higgs ind-$\bvocB$-module $\rI\hupsigma^*(N,\theta)$ \eqref{p2-htaft23a} is twistable by the extension \eqref{p2-htaft18h}, in the sense of \ref{p2-rgpsc7}. 
\end{prop}

It follows from (\cite{ag2} 4.8.7 and 4.8.13), see (\cite{ag2} 4.6.3) for the compatibility of the different constructions. 

\begin{cor}\label{p2-htaft59}
For every locally CL-small Higgs $\co_\fX[\frac 1 p]$-module $(N,\theta)$ with coefficients in $\txi^{-1}\tOmega^1_{\fX/\cS}$ \eqref{p1-tshbn13},
the Higgs ind-$\bvocB$-module $\rI\hupsigma^*(N,\theta)$ \eqref{p2-htaft23a} is weakly twistable by the extension \eqref{p2-htaft18h}, in the sense of \ref{p2-rgpsc4}. 
\end{cor}

If follows from \ref{p2-rgpsc5}, \ref{p2-rgpsc9} and \ref{p2-htaft58}.

\begin{prop}[\cite{ag2} 4.8.15]\label{p2-rgpsc24}
We suppose that $X$ is an object of $\bP$ \eqref{p2-htaft9} and that $X_s$ is non-empty.
Then a Higgs $\co_\fX[\frac 1 p]$-bundle with coefficients in $\txi^{-1}\tOmega^1_{\fX/\cS}$ is CL-small in the sense of \eqref{p1-tshbn13},
if and only if it is strongly solvable in the sense of \ref{p2-rgpsc20}.
\end{prop}

\begin{prop}[\cite{ag2} 4.9.11]\label{p2-rgpsc240}
Every strongly solvable Higgs $\co_\fX[\frac 1 p]$-bundle with coefficients in $\txi^{-1}\tOmega^1_{\fX/\cS} $ is locally CL-small \eqref{p1-tshbn13}.
\end{prop}

\begin{cor}\label{p2-rgpsc241}
For every strongly Dolbeault $\bvocB_\mQ$-module $M$ \eqref{p2-rgpsc20}, 
the associated Higgs $\co_\fX[\frac 1 p]$-bundle $\cH_\mQ(M)$ \eqref{p2-rgpsc23a} 
is locally CL-small \eqref{p1-tshbn13}.
\end{cor}

\begin{prop}[\cite{ag2} 4.2.14]
Every nilpotent Higgs $\co_\fX[\frac 1 p]$-bundle with coefficients in $\txi^{-1}\tOmega^1_{\fX/\cS}$ {\rm (\cite{ag2} 4.2.5)} is locally CL-small. 
\end{prop}

\subsection{}\label{p2-rgpsc25}
We take again the notation of \ref{p2-htaft23}, so we denote by $\IndSym_{\bvocB}(\hupsigma^*(\cT))$ the symmetric ind-$\bvocB$-algebra of 
$\hupsigma^*(\cT)$ \eqref{p2-htaft45b} \eqref{p1-indmal6}.
Let $A$ be an ind-$\bvocB$-algebra quotient of $\IndSym_{\bvocB}(\hupsigma^*(\cT))$ \eqref{p1-indmal2},
$\upmu_A\colon \IndSym_{\bvocB}(\hupsigma^*(\cT))\rightarrow \cEnd_{\bvocB}(A)$ 
the morphism of ind-$\bvocB$-algebras induced by the multiplication of $A$ \eqref{p1-indmal15}. 
By \ref{p2-htaft23}, the morphism $\upmu_A$ defines a Higgs $\bvocB$-field
\begin{equation}\label{p2-rgpsc25a}
\theta_A\colon A\rightarrow \hupsigma^*(\txi^{-1}\tOmega^1_{\fX/\cS}) \otimes_{\co_\fX} A,
\end{equation}
that we call the {\em canonical Higgs field on $A$}. 
We set 
\begin{equation}\label{p2-rgpsc25b}
\cL_A=\upnu(A,\theta_A),
\end{equation}
where $\upnu$ is the functor \eqref{p2-rgpsc1c}, which is naturally an ind-$A$-module; indeed, its structure of ind-$\IndSym_{\bvocB}(\hupsigma^*(\cT))$-module induced by its Higgs field
factors through $A$ \eqref{p2-rgpsc32}.

\subsection{}\label{p2-rgpsc30}
Let $(\cN,\theta)$ be a Higgs ind-$\bvocB$-module with coefficients in $\hupsigma^*(\txi^{-1}\tOmega^1_{\fX/\cS})$, 
\[
\upmu\colon \IndSym_{\bvocB}(\hupsigma^*(\cT))\rightarrow \cEnd_{\bvocB}(\cN)
\]
the morphism of ind-$\bvocB$-algebras defined by $\theta$, see \ref{p2-htaft23} and \ref{p1-indmal15}. 
We assume that $\upmu$ factors 
through a quotient ind-$\bvocB$-algebra $A$ of $\IndSym_{\bvocB}(\hupsigma^*(\cT))$, and let 
\begin{equation}
\theta_A\colon A\rightarrow \hupsigma^*(\txi^{-1}\tOmega^1_{\fX/\cS}) \otimes_{\bvocB} A
\end{equation}
be the canonical Higgs $\bvocB$-field \eqref{p2-rgpsc25a}. 
Since $\theta_A$ is a morphism of ind-$A$-modules by \ref{p1-indmal23}(i), we may consider the morphism 
\begin{equation}\label{p2-rgpsc30a}
\theta_A\otimes_A\id_\cN\colon A\otimes_A\cN\rightarrow \hupsigma^*(\txi^{-1}\tOmega^1_{\fX/\cS})\otimes_{\bvocB}A\otimes_A\cN.
\end{equation}
It is the canonical Higgs $\bvocB$-field on the ind-$A$-module $A\otimes_A\cN$ by \ref{p1-indmal23}(ii). 
It fits into a commutative diagram
\begin{equation}\label{p2-rgpsc30b}
\xymatrix{
{A\otimes_A\cN}\ar[rr]^-(0.5){\theta_A\otimes\id_\cN}\ar[d]&&{\hupsigma^*(\txi^{-1}\tOmega^1_{\fX/\cS})\otimes_{\bvocB}A\otimes_A\cN}\ar[d]\\
\cN\ar[rr]^-(0.5)\theta&&{\hupsigma^*(\txi^{-1}\tOmega^1_{\fX/\cS})\otimes_{\bvocB}\cN,}}
\end{equation}
where the vertical arrows are the isomorphisms induced by $\upmu$ \eqref{p1-indmal5g}. 
Indeed, we may assume that the $\co_\fX$-module $\tOmega^1_{\fX/\cS}$ is free of finite type
by (\cite{ag2} 2.7.17(i)), in which case the assertion can be checked easily, see \ref{p1-thbn43}. 
 
The homomorphism $\upmu$ induces an isomorphism of ind-$\IC^\dagger$-modules
\begin{equation}\label{p2-rgpsc30e}
\cN\otimes_AA\otimes_{\bvocB}\IC^\dagger \rightarrow \cN\otimes_{\bvocB}\IC^\dagger.
\end{equation} 
By \eqref{p2-rgpsc30b}, it is compatible with the Higgs $\bvocB$-fields indicated on the same line of the table
\begin{equation}\label{p2-rgpsc30f}
\begin{tabular}{|c|c|}
\hline
left hand side&  right  hand  side\\
\hline
$\id\otimes \theta_A\otimes \id$ &  $\theta\otimes \id$\\
\hline
$\id\otimes \id\otimes \delta_{\IC^\dagger}$ & $\id\otimes \delta_{\IC^\dagger}$\\
\hline
\end{tabular}
\end{equation}
We deduce a morphism of $\bvocB$-modules 
\begin{equation}\label{p2-rgpsc30g}
\cN\otimes_{A}\cL_A\rightarrow \upnu(\cN,\theta),
\end{equation}
where $\cL_A=\upnu(A,\theta_A)$ \eqref{p2-rgpsc25b}. Since \eqref{p2-rgpsc30e} is a morphism of ind-$A$-modules, 
\eqref{p2-rgpsc30g} is a morphism of ind-$A$-modules by \ref{p2-rgpsc32}. 
It is therefore a morphism of Higgs ind-$\bvocB$-modules,
where $\cN\otimes_A\cL_A$ is equipped with its canonical Higgs field, see \ref{p1-indmal23}(ii).

\begin{prop}\label{p2-rgpsc31}
We keep the assumptions and notation of \ref{p2-rgpsc30}. 
Suppose moreover that the Higgs ind-$\bvocB$-modules $(\cN,\theta)$ and $(A,\theta_A)$ are weakly twistable. 
Then, the canonical morphism of ind-$A$-modules \eqref{p2-rgpsc30g}
\begin{equation}\label{p2-rgpsc31a}
\cN\otimes_A\cL_A\rightarrow \upnu(\cN,\theta)
\end{equation}
is an isomorphism. 
\end{prop}

Indeed, in view of \eqref{p1-indmal5g}, the diagram of morphisms of ind-$\IC^\dagger$-modules 
\begin{equation}\label{p2-rgpsc31b}
\xymatrix{
{\cN\otimes_A\cL_A\otimes_{\bvocB}\IC^\dagger}\ar[r]\ar[d]&{\upnu(\cN,\theta)\otimes_{\bvocB}\IC^\dagger}\ar[d]\\
{\cN\otimes_AA\otimes_{\bvocB}\IC^\dagger}\ar[r]&{\cN\otimes_{\bvocB}\IC^\dagger,}}
\end{equation}
where the upper (resp.\ lower) horizontal morphism is induced by \eqref{p2-rgpsc30g} (resp.\ is the isomorphism \eqref{p2-rgpsc30e}) and the vertical morphisms are the canonical isomorphisms 
\eqref{p2-rgpsc4a}, is commutative. The vertical and the lower horizontal morphisms being isomorphisms, we deduce that the upper horizontal morphism is an isomorphism. 
Therefore, by \ref{p2-htaft43}, the morphism \eqref{p2-rgpsc31a} is an isomorphism. 

\begin{rema}\label{p2-rgpsc39}
In the approach of Heuer \cite{heuer} and Heuer--Xu \cite{hexu}, for a smooth proper rigid variety, 
the $p$-adic Simpson correspondence is defined by twisting by a line bundle on the spectral variety, 
in a manner similar to \ref{p2-rgpsc31}, even beyond the small case.
\end{rema}

\section{Functoriality of Higgs--Tate algebras in Faltings topos}\label{p2-fhtft}

\subsection{}\label{p2-fhtft1}
In this section, we let $f\colon (X,\cM_{X})\rightarrow (S,\cM_S)$ and $f'\colon (X',\cM_{X'})\rightarrow (S,\cM_S)$ be 
adequate morphisms of fine logarithmic schemes (\cite{agt} III.4.7), and let
\begin{equation}\label{p2-fhtft1a}
g\colon (X',\cM_{X'})\rightarrow (X,\cM_X),
\end{equation} 
be a morphism of $(S,\cM_S)$-logarithmic schemes. 
In particular, the $S$-schemes $X$ and $X'$ are assumed to be of finite type. 
We denote by $X^\circ$ (resp.\ $X'^\rhd$) the maximal open subscheme of $X$ (resp.\ $X'$)
where the logarithmic structure $\cM_{X}$ (resp.\ $\cM_{X'}$) is trivial. 
For any $X$-scheme $U$ and $X'$-scheme $U'$, we set
\begin{eqnarray}
U^\circ=U\times_{X}X^\circ,\label{p2-fhtft1b}\\
U'^\rhd=U'\times_{X'}X'^\rhd.\label{p2-fhtft1c}
\end{eqnarray}
We have $X^\circ\subset X_\eta$ and $X'^\rhd\subset X'^\circ$. We denote by $\hbar\colon \oX\rightarrow X$, $h\colon \oX^\circ\rightarrow X$, 
$\hbar'\colon \oX'\rightarrow X'$ and $h'\colon \oX'^\rhd\rightarrow X'$ the canonical morphisms, and by 
\begin{equation}\label{p2-fhtft1e}
\upgamma\colon \oX'^\rhd\rightarrow \oX^\circ
\end{equation} 
the morphism induced by $g$.

We endow $\coX=X\times_S\coS$ \eqref{p2-ncgt1a} (resp.\ $\coX'=X'\times_S\coS$) with the logarithmic structure $\cM_{\coX}$ (resp.\ $\cM_{\coX'}$) pullback of $\cM_{X}$ (resp.\ $\cM_{X'}$), 
and denote by $\cof\colon (\coX,\cM_{\coX})\rightarrow (\coS,\cM_{\coS})$ (resp.\ $\cof'\colon (\coX',\cM_{\coX'})\rightarrow (\coS,\cM_{\coS})$, 
resp.\ $\cog\colon (\coX',\cM_{\coX'})\rightarrow (\coX,\cM_{\coX})$) the base change of $f$ (resp.\ $f'$, resp.\ $g$). 
We use the notation convention \eqref{p2-htaft1b} for the logarithmic differentials of $f$ and $f'$.
To lighten the notation, we set, with the convention of  \ref{p2-ncgt3}, 
\begin{eqnarray}
\tOmega^1_{\coX/\coS}=\tOmega^1_{(\coX,\cM_\coX)/(\coS,\cM_\coS)},&&\Omega=\txi^{-1}\tOmega^1_{\coX/\coS},\label{p2-fhtft9c}\\
\tOmega^1_{\coX'/\coS}=\tOmega^1_{(\coX',\cM_{\coX'})/(\coS,\cM_\coS)},&&\Omega'=\txi^{-1}\tOmega^1_{\coX'/\coS}.\label{p2-fhtft9d}
\end{eqnarray}
For any rational number $r\geq 0$, we set
\begin{equation}\label{p2-fhtft9g}
\Omega^{(r)}=p^r\Omega, \ \ \ \Omega'^{(r)}=p^r\Omega'.
\end{equation}
We have canonical morphisms of locally free $\co_{\coX'}$-modules of finite type
\begin{eqnarray}
u\colon \cog^*(\Omega)&\rightarrow& \Omega',\label{p2-fhtft90l}\\
u^{(r)}\colon \cog^*(\Omega^{(r)})&\rightarrow& \Omega'^{(r)}.\label{p2-fhtft90ll}
\end{eqnarray}

We fix Cartesian diagrams of $\FLS$ \eqref{p1-NC1}
\begin{equation}\label{p2-fhtft1d}
\xymatrix{
{(\coX,\cM_{\coX})}\ar[r]^-(0.5){i}\ar[d]_{\cof}\ar@{}[rd]|{\Box}&{(\tX,\cM_{\tX})}\ar[d]^{\tf}\\
{(\coS,\cM_{\coS})}\ar[r]^-(0.5){\iota}&{(\tS,\cM_{\tS}),}}
\ \ \ 
\xymatrix{
{(\coX',\cM_{\coX'})}\ar[r]^-(0.5){i'}\ar[d]_{\cof'}\ar@{}[rd]|{\Box}&{(\tX',\cM_{\tX'})}\ar[d]^{\tf'}\\
{(\coS,\cM_{\coS})}\ar[r]^-(0.5){\iota}&{(\tS,\cM_{\tS}),}}
\end{equation}
where $\iota$ is the strict closed immersion defined in \eqref{p2-ncgt3b}, such that $\tf$ and $\tf'$ are smooth.
We consider again the objects associated with $(f,\tf)$ introduced in §\ref{p2-htaft}, and we associate with $(f',\tf')$ similar objects that we denote 
by the same symbols equipped with a $\prime$ exponent.

\subsection{}\label{p2-fhtft2}
We denote by $\fX$ (resp.\ $\fX'$) the formal scheme $p$-adic completion of $\coX$ (resp.\ $\coX'$) \eqref{p2-ncgt1a},
which is a flat $\cS$-formal scheme of finite presentation \eqref{p2-ncgt1},
and by $\fgg\colon \fX'\rightarrow \fX$ the morphism induced by $\cog\colon \coX'\rightarrow \coX$.
We denote by $\tOmega^1_{\fX/\cS}$ (resp.\ $\tOmega^1_{\fX'/\cS}$) the $p$-adic completion of $\tOmega^1_{\coX/\coS}$
(resp.\ $\tOmega^1_{\coX'/\coS}$)
and by $\hOmega$ (resp.\ $\hOmega'$) the $p$-adic completion of $\Omega$ \eqref{p2-fhtft9c} (resp.\ $\Omega'$ \eqref{p2-fhtft9d}), 
which is canonically isomorphic to the module $\txi^{-1}\tOmega^1_{\fX/\cS}$ (resp.\ $\txi^{-1}\tOmega^1_{\fX'/\cS}$) \eqref{p2-ncgt3}; 
{\em we identify these modules in the following}. 
The morphism $u$ \eqref{p2-fhtft90l} induces a morphism of locally free $\co_{\fX'}$-modules of finite type
\begin{equation}\label{p2-fhtft2a}
\hu\colon \fgg^*(\hOmega)\rightarrow \hOmega'.
\end{equation}
For any rational numbers $r\geq r'\geq 0$, we denote by $\hOmega^{(r)}$ (resp.\ $\hOmega'^{(r)}$)
the $p$-adic completion of $\Omega^{(r)}$ (resp.\ $\Omega'^{(r)}$). 
Observe that $\hOmega^{(r)}$ (resp.\ $\hOmega'^{(r)}$) is $\cS$-flat and identifies canonically with $p^r\hOmega$ (resp.\ $p^r\hOmega'$).

\subsection{}\label{p2-fhtft3}
Consider the commutative diagram 
\begin{equation}\label{p2-fhtft3a}
\xymatrix{
{(\coX',\cM_{\coX'})}\ar[r]^-(0.5){i'}\ar[d]^{\cog}\ar@/_2pc/[dd]_{\cof'}&{(\tX',\cM_{\tX'})}\ar@/^2pc/[dd]^{\tf'}&\\
{(\coX,\cM_{\coX})}\ar[r]^-(0.5){i}\ar[d]^{\cof}\ar@{}[rd]|{\Box}&{(\tX,\cM_{\tX})}\ar[d]_{\tf}\\
{(\coS,\cM_{\coS})}\ar[r]^-(0.5){\iota}&{(\tS,\cM_{\tS}).}}
\end{equation}
By \ref{p2-hta2} applied to $\tf'$, the conditions \ref{p2-hta1}(i)-(iv) are satisfied. We denote by 
$\cL_{\tX'/\tX}$ the torsor of liftings of $\cog$ to $\tX'$ over $\tX$,
and by $\cF_{\tX'/\tX}$ (resp.\ $\cC_{\tX'/\tX}$) the Higgs--Tate extension (resp.\ algebra) of $\tX'$ over $\tX$ \eqref{p2-hta7}. 
Following the convention in \ref{p2-fhtal4}, to lighten the notation, {\em we will denote $\cF_{\tX'/\tX}$ (resp.\ $\cC_{\tX'/\tX}$) also by $\cF_\uptau$ (resp.\ $\cC_\uptau$)}. 
We have a canonical exact sequence of $\co_{\coX'}$-modules
\begin{equation}\label{p2-fhtft3b}
0\rightarrow \co_{\coX'}\rightarrow \cF_\uptau\rightarrow \cog^*(\Omega) \rightarrow 0,
\end{equation}
and a canonical isomorphism of $\co_{\coX'}$-algebras
\begin{equation}\label{p2-fhtft3c}
\cC_\uptau\stackrel{\sim}{\rightarrow}\underset{\underset{n\geq 0}{\longrightarrow}}\lim\ \rS^n_{\co_{\coX'}}(\cF_\uptau). 
\end{equation}
We denote by $\hcF_\uptau$ the $p$-adic completion of $\cF_\uptau$, so we have an exact sequence of $\co_{\fX'}$-modules
\begin{equation}\label{p2-fhtft3d}
0\rightarrow \co_{\fX'}\rightarrow \hcF_\uptau\rightarrow \fgg^*(\hOmega) \rightarrow 0.
\end{equation}

For any rational number $r\geq 0$, we denote by 
$\cL^{(r)}_{\tX'/\tX}$ the $(r)$-twisted torsor of liftings of $\cog$ to $\tX'$ over $\tX$,
and by $\cF^{(r)}_\uptau=\cF^{(r)}_{\tX'/\tX}$ (resp.\ $\cC^{(r)}_\uptau=\cC^{(r)}_{\tX'/\tX}$) the Higgs--Tate extension (resp.\ algebra) of $\tX'$ over $\tX$ of thickness $r$ \eqref{p2-hta7}. 
We have a canonical exact sequence of $\co_{\coX'}$-modules
\begin{equation}\label{p2-fhtft3e}
0\rightarrow \co_{\coX'}\rightarrow \cF^{(r)}_\uptau\rightarrow \cog^*(\Omega^{(r)}) \rightarrow 0
\end{equation}
and a canonical isomorphism of $\co_{\coX'}$-algebras
\begin{equation}\label{p2-fhtft3f}
\cC^{(r)}_\uptau\stackrel{\sim}{\rightarrow}\underset{\underset{n\geq 0}{\longrightarrow}}\lim\ \rS^n_{\co_{\coX'}}(\cF^{(r)}_\uptau). 
\end{equation}
For any integer $n\geq 0$, we set 
\begin{equation}\label{p2-fhtft3h}
\cF^{(r)}_{\uptau,n}=\cF^{(r)}_\uptau/p^n\cF^{(r)}_\uptau \ \ \ {\rm and}\ \ \ \cC^{(r)}_{\uptau,n}=\cC^{(r)}_\uptau/p^n\cC^{(r)}_\uptau.
\end{equation}
We denote by $\hcC^{(r)}_\uptau$ the $p$-adic completion of $\cC^{(r)}_\uptau$.  

For any object $U'$ of $\Et_{/X'}$, we set $\coU'=U'\times_S\coS$ \eqref{p2-ncgt1a}, 
\begin{equation}\label{p2-fhtft3g}
R'_{\uptau,U'}=\Gamma(\coU',\co_{\coX'}), \ \ \ \fF^{(r)}_{\uptau,U'}=\Gamma(\coU',\cF^{(r)}_\uptau) \ \ \ {\rm and}\ \ \ 
\fC^{(r)}_{\uptau,U'}=\Gamma(\coU',\cC^{(r)}_\uptau).
\end{equation}

For all rational numbers $r\geq r'\geq 0$, we have a canonical homomorphism of $\co_{\coX'}$-algebras \eqref{p2-hta5j}
\begin{equation}\label{p2-fhtft3i}
\alpha_{\uptau}^{r,r'}\colon \cC_\uptau^{(r)}\rightarrow \cC_\uptau^{(r')}.
\end{equation}
We denote by $\alpha_{\uptau,n}^{r,r'}\colon \cC_{\uptau,n}^{(r)}\rightarrow \cC_{\uptau,n}^{(r')}$ its reduction modulo $p^n$,  for any integer $n\geq 0$,
and by $\halpha_{\uptau}^{r,r'}\colon \hcC_{\uptau}^{(r)}\rightarrow \hcC_{\uptau}^{(r')}$ its extension to the $p$-adic completions. 

\begin{rema}\label{p2-fhtft7} 
Let $r$ be a rational number $\geq 0$, $\hcF^{(r)}_\uptau$ the $p$-adic completion of $\cF^{(r)}_\uptau$. The exact sequence \eqref{p2-fhtft3e} induces an exact sequence
\begin{equation}\label{p2-fhtft7a} 
0\rightarrow \co_{\fX'}\rightarrow \hcF^{(r)}_\uptau\rightarrow \fgg^*(\hOmega^{(r)}) \rightarrow 0,
\end{equation}
which is canonically isomorphic to the extension of $\co_{\fX'}$-modules deduced from $\hcF_\uptau$ \eqref{p2-fhtft3d} by pullback by the multiplication 
by $p^r$ on $\fgg^*(\hOmega)$. Let us temporarily denote by $C^{(r)}_\uptau$ the Higgs--Tate algebra associated with the extension $\hcF^{(r)}_\uptau$ 
in the sense of \eqref{p1-thbn1c}. The canonical morphism $\cF^{(r)}_\uptau\rightarrow \hcF^{(r)}_\uptau$ induces, for every integer $n\geq 0$, an isomorphism 
of $\co_{\fX'}$-algebras
\begin{equation}\label{p2-fhtft7b} 
\cC^{(r)}_{\uptau,n}\stackrel{\sim}{\rightarrow} C^{(r)}_\uptau/p^n C^{(r)}_\uptau. 
\end{equation}
We will therefore identify the $\co_{\fX'}$-algebra $\hcC^{(r)}_\uptau$ with the $p$-adic completion of $C^{(r)}_\uptau$. 
\end{rema}

\subsection{}\label{p2-fhtft4}
We defined in \eqref{p1-thbn9b} a twisting functor by the extension $\hcF_\uptau$ \eqref{p2-fhtft3d} 
for Higgs $\co_{\fX'}$-modules with coefficients in $\fgg^*(\hOmega)$:
\begin{equation}\label{p2-fhtft4a}
\uptau\colon \bHM(\co_{\fX'},\fgg^*(\hOmega))\rightarrow \bHM(\co_{\fX'},\fgg^*(\hOmega)).
\end{equation}
However, it only behaves well for {\em weakly twistable} Higgs bundles \eqref{p1-thbn30}. 
Composing with the pullback functor $\fgg^*$, we obtain a functor
\begin{equation}\label{p2-fhtft4b}
\fgg^*_\uptau\colon 
\begin{array}[t]{clcr}
\bHM(\co_\fX,\hOmega)&\rightarrow& \bHM(\co_{\fX'},\fgg^*(\hOmega)),\\
(N,\theta)&\mapsto&\uptau(\fgg^*(N),\fgg^*(\theta)),
\end{array}
\end{equation} 
that we call the {\em pullback functor by $\fgg$ twisted by the extension $\hcF_\uptau$} \eqref{p2-fhtft3d} 
(or simply the {\em twisted pullback functor by $\fgg$} when the extension is implicit); 
see \ref{p1-tphdi2}.

\begin{defi}\label{p2-fhtft5}
We say that a Higgs $\co_\fX[\frac 1 p]$-module $(N,\theta)$ with coefficients in $\hOmega$ is 
{\em weakly twistable} (resp.\ {\em twistable}) by the extension $\hcF_\uptau$ \eqref{p2-fhtft3d} 
if the Higgs $\co_{\fX'}[\frac 1 p]$-module $(\fgg^*(N),\fgg^*(\theta))$ with coefficients 
in $\fgg^*(\hOmega)$ is weakly twistable (resp.\ twistable) by the extension $\hcF_\uptau$ 
in the sense of \ref{p1-thbn30} (resp.\ \ref{p1-thbn14}). 
\end{defi}

\subsection{}\label{p2-fhtft6}
For any rational number $r\geq 0$, we denote by 
\begin{equation}\label{p2-fhtft6h}
d_{\cC^{(r)}_\uptau}\colon \cC^{(r)}_\uptau \rightarrow \cog^*(\Omega^{(r)})\otimes_{\co_{\coX'}}\cC^{(r)}_\uptau
\end{equation}
the universal $\co_{\coX'}$-derivation of $\cC^{(r)}_\uptau$ \eqref{p2-hta5g}, and set 
\begin{eqnarray}
\delta_{\cC^{(r)}_\uptau}=(\cog^*(\pi^{(r,0)})\otimes \id)\circ d_{\cC^{(r)}_\uptau} \colon \cC^{(r)}_\uptau &\rightarrow& 
\cog^*(\Omega)\otimes_{\co_{\coX'}}\cC^{(r)}_\uptau,\label{p2-fhtft6a}\\
\delta'_{\cC^{(r)}_\uptau}=(u\otimes \id)\circ \delta_{\cC^{(r)}_\uptau} \colon \cC^{(r)}_\uptau &\rightarrow &
\Omega'\otimes_{\co_{\coX'}}\cC^{(r)}_\uptau,\label{p2-fhtft6aa}
\end{eqnarray}
where $u$ is the morphism defined in \eqref{p2-fhtft90l}. 
They are Higgs $\co_{\coX'}$-fields. 
We denote by 
\begin{eqnarray}
\delta_{\hcC^{(r)}_\uptau}\colon \hcC^{(r)}_\uptau &\rightarrow& \fgg^*(\hOmega)\otimes_{\co_{\fX'}} \hcC^{(r)}_\uptau\label{p2-fhtft6b}\\
\delta'_{\hcC^{(r)}_\uptau}\colon \hcC^{(r)}_\uptau &\rightarrow& \hOmega'\otimes_{\co_{\fX'}} \hcC^{(r)}_\uptau\label{p2-fhtft6bb}
\end{eqnarray}
the extension of $\delta_{\cC^{(r)}_\uptau}$ and $\delta'_{\cC^{(r)}_\uptau}$, respectively, to the $p$-adic completions. 

By \eqref{p2-hta5l},  for all rational numbers $r\geq r'\geq 0$, we have 
\begin{equation}\label{p2-fhtft6d}
(\id \otimes \alpha^{r,r'}_\uptau) \circ \delta_{\cC^{(r)}_\uptau}=\delta_{\cC^{(r')}_\uptau} \circ \alpha^{r,r'}_\uptau,
\end{equation}
where $\alpha^{r,r'}_\uptau$ is defined in \eqref{p2-fhtft3i}.

We define the $\co_{\fX'}$-algebra
\begin{equation}\label{p2-fhtft6e}
\hcC^{(r+)}_\uptau=\underset{\underset{t\in \mQ_{>r}}{\longrightarrow}}{\lim}\ \hcC^{(t)}_\uptau. 
\end{equation}
By \eqref{p2-fhtft6d}, the derivations $\delta_{\hcC^{(t)}_\uptau}$ \eqref{p2-fhtft6b} induce an $\co_{\fX'}$-derivation
\begin{equation}\label{p2-fhtft6f}
\delta_{\hcC^{(r+)}_\uptau}\colon \hcC^{(r+)}_\uptau\rightarrow \fgg^*(\hOmega)\otimes_{\co_{\fX'}}\hcC^{(r+)}_\uptau,
\end{equation}
which is also a Higgs $\co_{\fX'}$-field.
For simplicity, we set $\cC^\dagger_\uptau=\hcC^{(0+)}_\uptau$, $\delta_{\uptau}=\delta_{\hcC^{(0+)}_\uptau}$ and 
\begin{equation}\label{p2-fhtft6g}
\delta'_{\uptau}=(\hu\otimes\id)\circ \delta_{\uptau}
\colon \cC^\dagger_\uptau \rightarrow \hOmega'\otimes_{\co_{\fX'}} \cC^\dagger_\uptau,
\end{equation}
where $\hu$ is the canonical morphism in \eqref{p2-fhtft2a}.

We associate with the extension \eqref{p2-fhtft3d} also an ind-$\co_{\fX'}$-algebra \eqref{p2-cmupiso32e}
\begin{equation}\label{p2-fhtft6k}
\IC^\dagger_{\uptau}=\underset{\underset{r\in \mQ_{>0}}{\longrightarrow}}{\mlq\mlq\lim \mrq\mrq}\ \hcC^{(r)}_{\uptau}. 
\end{equation}
and an $\co_{\fX'}$-derivation \eqref{p2-cmupiso32f}
\begin{equation}\label{p2-fhtft6i}
\Idelta_{\uptau}\colon \IC^\dagger_\uptau \rightarrow \fgg^*(\hOmega)\otimes_{\co_{\fX'}}\IC^\dagger_\uptau.
\end{equation}
It is a Higgs $\co_{\fX'}$-field \eqref{p1-indmal20}. We consider also the $\co_{\fX'}$-derivation
\begin{equation}\label{p2-fhtft6j}
\Idelta'_{\uptau}=(\hu\otimes \id)\circ \Idelta_{\uptau}\colon 
\IC^\dagger_\uptau\rightarrow \hOmega'\otimes_{\co_{\fX'}}\IC^\dagger_\uptau.
\end{equation}

\subsection{}\label{p2-fhtft21}
We take again the notation of \ref{p2-htaft3}--\ref{p2-htaft6}, both for $f$ and $f'$. The functor
\begin{equation}\label{p2-fhtft21a}
\Theta^+\colon E\rightarrow E', \ \ \ (V\rightarrow U)\mapsto (V\times_{X^\circ}X'^\rhd\rightarrow U\times_XX')
\end{equation}
is continuous and left exact (\cite{agt} VI.10.12). It therefore defines a morphism of topos
\begin{equation}\label{p2-fhtft21b}
\Theta\colon \tE'\rightarrow \tE.
\end{equation}
It  immediately follows  from the definitions that the squares of the diagram
\begin{equation}\label{p2-fhtft21c}
\xymatrix{
{X'_\et}\ar[d]_{g}&{\tE'}\ar[l]_-(0.5){\sigma'}\ar[d]^{\Theta}\ar[r]^-(0.5){\beta'}&
{\oX'^\rhd_\fet}\ar[d]^{\upgamma}\\
{X_\et}&{\tE}\ar[l]_{\sigma}\ar[r]^{\beta}&{\oX^\circ_\fet}}
\end{equation}
are commutative up to canonical isomorphisms.

In view of \ref{p2-htaft6} and (\cite{ag2} 6.1.11), we have a canonical ring homomorphism of $\tE$
\begin{equation}\label{p2-fhtft21d}
\ocB\rightarrow \Theta_*(\ocB').
\end{equation}
We consider in the following $\Theta$ as a morphism of ringed topos (respectively by $\ocB'$ and $\ocB$).
For modules, we use the notation $\Theta^{-1}$ to denote the pullback in the sense of abelian sheaves and we keep the notation
$\Theta^*$ for the pullback in the sense of modules.

We have a canonical isomorphism $\Theta^*(\sigma^*(X_\eta))\simeq \sigma'^*(X'_\eta)$. 
Hence, by (\cite{sga4} IV 9.4.3), there exists a morphism of topos
\begin{equation}\label{p2-fhtft21e}
\uptheta\colon \tE'_s\rightarrow \tE_s
\end{equation}
unique up to canonical isomorphism such that the diagram
\begin{equation}\label{p2-fhtft21f}
\xymatrix{
{\tE'_s}\ar[r]^{\uptheta}\ar[d]_{\delta'}&{\tE_s}\ar[d]^{\delta}\\
{\tE'}\ar[r]^\Theta&{\tE}}
\end{equation}
is commutative up to isomorphism, and even $2$-Cartesian. It follows from \eqref{p2-fhtft21c} and (\cite{sga4} IV 9.4.3)
that the diagram of morphism of topos
\begin{equation}\label{p2-fhtft21g}
\xymatrix{
{\tE'_s}\ar[r]^{\uptheta}\ar[d]_{\sigma'_s}&{\tE_s}\ar[d]^{\sigma_s}\\
{X'_{s,\et}}\ar[r]^{g_s}&{X_{s,\et}}}
\end{equation}
is commutative up to canonical isomorphism.

For every integer $n\geq 0$, the canonical homomorphism $\Theta^{-1}(\ocB)\rightarrow \ocB'$
induces a homomorphism $\uptheta^{-1}(\ocB_n)\rightarrow \ocB'_n$.
The morphism $\uptheta$ is therefore underlying a morphism of ringed topos, which we denote by
\begin{equation}\label{p2-fhtft21h}
\uptheta_n\colon (\tE'_s,\ocB'_n)\rightarrow (\tE_s,\ocB_n).
\end{equation}
In view of \eqref{p2-fhtft21g}, we immediately verify that the diagram of morphism of topos
\begin{equation}\label{p2-fhtft21i}
\xymatrix{
{(\tE'_s,\ocB'_n)}\ar[r]^{\uptheta_n}\ar[d]_{\sigma'_n}&{(\tE_s,\ocB_n)}\ar[d] ^{\sigma_n}\\
{(X'_{s,\et},\co_{\oX'_n})}\ar[r]^{\ogg_n}&{(X_{s,\et},\co_{\oX_n}) }}
\end{equation}
is commutative up to canonical isomorphism.

\subsection{}\label{p2-fhtft22}
We take again the notation of \ref{p2-htaft8}, both for $f$ and $f'$. 
The morphisms $(\uptheta_{n+1})_{n\in \mN}$ \eqref{p2-fhtft21h} induce a morphism of ringed topos 
\begin{equation}\label{p2-fhtft22a}
\bvuptheta\colon (\tE'^{\mN^\circ}_s,\bvocB') \rightarrow (\tE^{\mN^\circ}_s,\bvocB)
\end{equation}
that fits into a commutative diagram up to canonical isomorphism
\begin{equation}\label{p2-fhtft22b}
\xymatrix{
{(\tE'^{\mN^\circ}_s,\bvocB')}\ar[r]^{\bvuptheta}\ar[d]_{\hupsigma'}&{(\tE^{\mN^\circ}_s,\bvocB)}\ar[d]^{\hupsigma}\\
{(\fX',\co_{\fX'})}\ar[r]^{\fgg}&{(\fX,\co_\fX).}}
\end{equation}

\subsection{}\label{p2-fhtft39}
Let $\ox'$ be a geometric point of $X'$, $\ox=g(\ox')$, $\uX$ (resp.\ $\uX'$) the strict localization of $X$ at $\ox$ (resp.\ $X'$ at $\ox'$),
$\ug\colon \uX'\rightarrow \uX$ and $\uupgamma\colon \uoX'^\rhd\rightarrow \uoX^\circ$ the morphisms induced by $g$.
We denote by $\tuE$ (resp.\ $\tuE'$) the Faltings topos associated with the canonical morphism $\uoX^\circ \rightarrow \uX$
(resp.\ $\uoX'^\rhd \rightarrow \uX'$), by
\begin{equation}\label{p2-fhtft39a}
\uTheta\colon \tuE'\rightarrow \tuE
\end{equation}
the functoriality morphism induced by the morphism $\ug$ (see \ref{p2-fhtft21}), and by
\begin{eqnarray}
\ubeta\colon \tuE\rightarrow \uoX^\circ_\fet,\label{p2-fhtft39b}\\
\ubeta'\colon \tuE'\rightarrow \uoX'^\rhd_\fet,\label{p2-fhtft39c}
\end{eqnarray}
the canonical morphisms (\cite{agt} VI.10.6).
By (\cite{agt} (VI.10.12.6)), the diagram
\begin{equation}\label{p2-fhtft39d}
\xymatrix{
{\tuE'}\ar[r]^{\uTheta}\ar[d]_{\ubeta'}&{\tuE}\ar[d]^{\ubeta}\\
{\uoX'^\rhd_\fet}\ar[r]^{\uupgamma}&{\uoX^\circ_\fet}}
\end{equation}
is commutative up to canonical isomorphism. We denote by
\begin{eqnarray}
\theta\colon \uoX^\circ_\fet\rightarrow \tuE,\label{p2-fhtft39e}\\
\theta'\colon \uoX'^\rhd_\fet\rightarrow \tuE',\label{p2-fhtft39f}
\end{eqnarray}
the canonical sections of $\ubeta$ and $\ubeta'$ (\cite{agt} VI.10.23). The diagram
\begin{equation}\label{p2-fhtft39g}
\xymatrix{
{\uoX'^\rhd_\fet}\ar[r]^-(0.5)\uupgamma\ar[d]_{\theta'}&{\uoX^\circ_\fet}\ar[d]^{ \theta}\\
{\tuE'}\ar[r]^-(0.5){\uTheta}&{\tuE}}
\end{equation}
is commutative up to canonical isomorphism. Indeed, for every $\uX$-scheme $U$ that is étale, separated and of finite presentation,
if we denote by $U^\rf$ its $\uX$-finite part (i.e., the disjoint sum of the strict localizations of $U$ at the points of $U_\ox$), then
$U^\rf\times_\uX\uX'$ is the $\uX'$-finite part of $U\times_\uX\uX'$ (see \cite{agt} VI.10.22).

Moreover, the diagram
\begin{equation}\label{p2-fhtft39h}
\xymatrix{
{\tuE'}\ar[r]^-(0.5){\uTheta}\ar[d]_{\Phi'}&{\tuE}\ar[d]^\Phi\\
{\tE'}\ar[r]^-(0.5){\Theta}&{\tE},}
\end{equation}
where the vertical arrows are the functoriality morphisms induced by the canonical morphisms $\uX\rightarrow X$ and $\uX'\rightarrow X'$,
is commutative up to canonical isomorphism. We denote by
\begin{eqnarray}
\varphi_\ox\colon \tE\rightarrow \uoX^\circ_\fet,\label{p2-fhtft39i}\\
\varphi'_{\ox'}\colon \tE'\rightarrow \uoX'^\rhd_\fet,\label{p2-fhtft39ii}
\end{eqnarray}
the composed functors $\theta^*\circ \Phi^*$ and $\theta'^*\circ \Phi'^*$ (see \cite{ag2} 4.3.5). We therefore have a canonical isomorphism of functors
\begin{equation}\label{p2-fhtft39j}
\uupgamma^* \circ \varphi_\ox\stackrel{\sim}{\rightarrow} \varphi'_{\ox'} \circ \Theta^*.
\end{equation}

\subsection{}\label{p2-fhtft30}
We denote by $\mM$ the following category. An object of $\mM$ is a morphism of schemes $\mu\colon U'\rightarrow U$ above the morphism $g\colon X'\rightarrow X$, 
i.e., a commutative diagram 
\begin{equation}\label{p2-fhtft30a}
\xymatrix{
U'\ar[r]^-(0.5)\mu\ar[d]&U\ar[d]\\
X'\ar[r]^-(0.5)g&X,}
\end{equation}
such that the morphisms $U\rightarrow X$ and $U'\rightarrow X'$ are étale. For two objects $\mu\colon U'\rightarrow U$ and $\mu_1\colon U'_1\rightarrow U_1$ of $\mM$, 
a morphism from $\mu_1$ to $\mu$ is a pair of an $X$-morphism 
$u\colon U_1\rightarrow U$ and an $X'$-morphism $u'\colon U'_1\rightarrow U'$ such that the diagram 
\begin{equation}\label{p2-fhtft30b}
\xymatrix{
U'_1\ar[r]^-(0.45){\mu_1}\ar[d]_-(0.45){u'}&U_1\ar[d]^-(0.5)u\\
U'\ar[r]^-(0.5)\mu&U}
\end{equation}
is commutative. 
Observe that fiber products are representable in $\mM$.

We take again the notation of \ref{p2-htaft10}, both for $f$ and $f'$, and 
denote by $\mM(\bQ',\bQ)$ the full subcategory of $\mM$ of morphisms $\mu\colon U'\rightarrow U$ such that 
$U'\rightarrow X'$ (resp.\ $U\rightarrow X$) is an object of $\bQ'$ (resp.\ $\bQ$) \eqref{p2-htaft10}. We denote by $\umM(\bQ',\bQ)$ the full subcategory of $\mM(\bQ',\bQ)$ 
made of morphisms $\mu\colon U'\rightarrow U$ such that the morphism of $(S,\cM_S)$-logarithmic schemes
$(U',\cM_{X'}|U')\rightarrow (U,\cM_X|U)$ induced by $\mu$ and $g$ \eqref{p2-fhtft1a} admits a relative adequate chart in the sense of (\cite{ag1} 5.1.11).

\subsection{}\label{p2-fhtft301}
Let $\ox'$ be a geometric point of $X'$. We define the category $\mM_{\ox'}$ as follows.
The objects of $\mM_{\ox'}$ are pairs made of an object $\mu\colon U'\rightarrow U$ of $\mM$ \eqref{p2-fhtft30} and an $X'$-morphism $\iota\colon \ox'\rightarrow U'$. 
Let $(\mu\colon U'\rightarrow U, \iota\colon \ox'\rightarrow U')$ and $(\mu_1\colon U'_1\rightarrow U_1, \iota_1\colon \ox'\rightarrow U'_1)$ be two objects of
$\mM_{\ox'}$. A morphism from $(\mu_1\colon U'_1\rightarrow U_1, \iota_1\colon \ox'\rightarrow U'_1)$ to $(\mu\colon U'\rightarrow U, \iota\colon \ox'\rightarrow U')$ 
is a morphism $(u'\colon U'_1\rightarrow U',u\colon U_1\rightarrow U)$ of $\mM$ \eqref{p2-fhtft30b} such that $\iota=u'\circ \iota_1$.
Finite inverse limits are representable in $\mM_{\ox'}$ (see \cite{sga4} I 2.3).
The category $\mM_{\ox'}$ is therefore cofiltered  (\cite{sga4} I 2.7.1).

We denote by $\mM_{\ox'}(\bQ',\bQ)$ (resp. $\umM_{\ox'}(\bQ',\bQ)$) the full subcategory of $\mM_{\ox'}$ 
made of objects $(\mu\colon U'\rightarrow U, \iota\colon \ox'\rightarrow U')$ 
such that $\mu$ is an object of $\mM(\bQ',\bQ)$ (resp. $\umM(\bQ',\bQ)$) \eqref{p2-fhtft30}.
The canonical functor $\mM_{\ox'}(\bQ',\bQ)\rightarrow \mM_{\ox'}$ is
initial and the category $\mM_{\ox'}(\bQ',\bQ)$ is cofiltered by (\cite{sga4} I 8.1.3(c)) and (\cite{agt} II.5.17).
Similarly, the canonical functor $\umM_{\ox'}(\bQ',\bQ)\rightarrow \mM_{\ox'}(\bQ',\bQ)$ is
initial and the category $\umM_{\ox'}(\bQ',\bQ)$ is cofiltered by (\cite{sga4} I 8.1.3(c)) and (\cite{ag1} 5.1.12).

\subsection{}\label{p2-fhtft302}
Let $\mu\colon U'\rightarrow U$ be an object of $\mM$. We denote by 
\begin{equation}\label{p2-fhtft302a}
\tmu\colon \oU'^\rhd\rightarrow \oU^\circ
\end{equation} 
the induced morphism, see \eqref{p2-ncgt1a}, \eqref{p2-fhtft1b}  and \eqref{p2-fhtft1c}. We denote by $\tmu^+\colon \Et_{\rf/\oU^\circ}\rightarrow \Et_{\rf/\oU'^\rhd}$ 
the base change functor by $\tmu$, by $\tmu_{\rp}\colon (\Et_{\rf/\oU'^\rhd})^\wedge\rightarrow (\Et_{\rf/\oU^\circ})^\wedge$ the functor defined by 
composition with $\tmu^+$ and by $\tmu^\rp\colon  (\Et_{\rf/\oU^\circ})^\wedge\rightarrow (\Et_{\rf/\oU'^\rhd})^\wedge$ its left adjoint \eqref{p2-cmt4}. 
The functor $\tmu^+$ being left exact and continuous, it induces a morphism of topos that we denote abusively also by $\tmu\colon \oU'^\rhd_\fet\rightarrow \oU^\circ_\fet$. 

With the notation of \ref{p2-htaft6}, for every $V\in \ob(\Et_{\rf/\oU^\circ})$, setting $V'=V\times_{U^\circ}U'^\rhd$,  
$\mu$ induces a morphism $\oU'^{V'}\rightarrow \oU^V$. We deduce a ring homomorphism of $\oU^\circ_\fet$
\begin{equation}\label{p2-fhtft302b}
\ocB_U\rightarrow \tmu_*(\ocB'_{U'}). 
\end{equation}
For any integer $n\geq 0$, we denote by 
\begin{equation}\label{p2-fhtft302c}
\tmu_n\colon (\oU'^\rhd_\fet,\ocB'_{U',n})\rightarrow (\oU^\circ_\fet,\ocB_{U,n})
\end{equation}
the morphism of ringed topos induced by $\tmu$ and \eqref{p2-fhtft302b}. 

Let $\oy'$ a geometric point of $\oU'^\rhd$, $\oy=\tmu(\oy')$. The scheme $\oU$ being locally irreducible by (\cite{ag1} 4.2.7 and \cite{agt} III.3.3),
it is the sum of the schemes induced on its irreducible components. We denote by $\oU^\star$ 
the irreducible component of $\oU$ containing $\oy$.
Similarly, $\oU^\circ$ is the sum of the schemes induced on its irreducible components
and $\oU^{\star \circ}=\oU^\star\times_{X}X^\circ$ is the irreducible component of $\oU^\circ$ containing $\oy$.
We consider the discrete representation $\oR^\oy_U$ of $\pi_1(\oU^{\star\circ},\oy)$ defined in \eqref{p2-htaft11b}. 
We define similarly the scheme $\oU'^\star$ and  the discrete representation $\oR'^{\oy'}_{U'}$ of $\pi_1(\oU'^{\star\rhd},\oy')$.  
The morphism \eqref{p2-fhtft302b} induces a $\pi_1(\oU'^{\star\rhd},\oy')$-equivariant homomorphism 
\begin{equation}\label{p2-fhtft302d}
\oR^\oy_U\rightarrow \oR'^{\oy'}_{U'}.
\end{equation}

\subsection{}\label{p2-fhtft35}
Consider the strictly commutative diagram 
\begin{equation}\label{p2-fhtft35a}
\xymatrix{
E\ar[r]^{\Theta^+}\ar[d]_{\pi}&E'\ar[d]^{\pi'}\\
{\Et_{/X}}\ar[r]^{g^+}&{\Et_{/X'}}}
\end{equation}
where $\pi$ and $\pi'$ are Faltings covanishing fibered sites defined in \eqref{p2-htaft3b} relatively to $f$ and $f'$ respectively, 
$\Theta^+$ is the functor \eqref{p2-fhtft21a} and $g^+$ is the base change functor by $g$. 

We consider again the objects associated with $f$ introduced in \ref{p2-htaft10}, and we associate with $f'$ similar objects that we denote 
by the same symbols equipped with a $\prime$ exponent. We consider the composed functor 
\begin{equation}\label{p2-fhtft35b}
\Phi\colon 
\xymatrix{
{E_\bQ}\ar[r]^-(0.5){\nu}&{E}\ar[r]^-(0.5){\Theta^+}&{E'}\ar[r]^-(0.5){\tth_{E'}}&{\hE'}\ar[r]^-(0.5){\nu'_\rp}&{\hE'_{\bQ'},}}
\end{equation}
where $\tth_{E'}$ is the canonical functor and $\nu'_\rp$ is defined by composition with the canonical functor 
$\nu'\colon E'_{\bQ'}\rightarrow E'$ \eqref{p2-cmt4a}.
Recall that $E_{\bQ}$ (resp.\ $E'_{\bQ'}$) is a topologically generating $\mU$-small subcategory of the covanishing site $E$ (resp.\ $E'$) \eqref{p2-htaft10}. 
Therefore, by \ref{p2-cmt101}(a), $\Phi$ is a $\mU$-functor \eqref{p2-cmt0}. 
It then induces a functor $\Phi_\rp\colon \hE'_{\bQ'}\rightarrow \hE_\bQ$ \eqref{p2-cmt1a}. 
By \ref{p2-cmt2}, the latter admits a left adjoint
\begin{equation}\label{p2-fhtft35c}
\Phi^\rp\colon \hE_\bQ\rightarrow \hE'_{\bQ'}.
\end{equation}

By \ref{p2-cmt77}, $\Phi$ is continuous in the sense of \ref{p2-cmt3}. Let $\Phi_\rs\colon \tE'_{\bQ'}\rightarrow \tE_\bQ$ be the functor 
making strictly commutative the following diagram
\begin{equation}\label{p2-fhtft35d}
\xymatrix{
{\tE'_{\bQ'}}\ar[r]^-(0.5){\Phi_\rs}\ar[d]_{\tti'_{\bQ'}}&{\tE_\bQ}\ar[d]^{\tti_\bQ}\\
{\hE'_{\bQ'}}\ar[r]^-(0.5){\Phi_\rp}&{\hE_\bQ,}}
\end{equation}
where the vertical arrows are the canonical functors. 
By \eqref{p2-cmt77b}, the diagram 
\begin{equation}\label{p2-fhtft35e}
\xymatrix{
{\tE'}\ar[r]^-(0.5){\Theta_*}\ar[d]_{\nu'_\rs}&{\tE}\ar[d]^{\nu_\rs}\\
{\tE'_{\bQ'}}\ar[r]^-(0.5){\Phi_\rs}&{\tE_\bQ,}}
\end{equation}
where $\Theta_*$ is the direct image functor by the morphism of topos \eqref{p2-fhtft21b} 
and the vertical arrows are the equivalences of categories \eqref{p2-htaft10d}, is commutative up to isomorphism.

Let $\Phi^\rs$ be the composed functor 
\begin{equation}\label{p2-fhtft35f}
\Phi^\rs\colon 
\xymatrix{ 
{\tE_\bQ}\ar[r]^-(0.5){\tti_{\bQ}}&{\hE_\bQ}\ar[r]^-(0.5){\Phi^\rp}&{\hE'_{\bQ'}}\ar[r]^-(0.5){\tta'_{\bQ'}}&{\tE'_{\bQ'}},}
\end{equation}
where $\tta'_{\bQ'}$ is the ``associated sheaf'' functor. 
It is a left adjoint of $\Phi_\rs$ by \ref{p2-cmt5}. By \eqref{p2-cmt77d}, the diagrams  
\begin{equation}\label{p2-fhtft35g}
\xymatrix{
{\hE_\bQ}\ar[r]^-(0.5){\Phi^\rp}\ar[d]_{\tta_\bQ}&{\hE'_{\bQ'}}\ar[d]^{\tta'_{\bQ'}}\\
{\tE_\bQ}\ar[r]^-(0.5){\Phi^\rs}&{\tE'_{\bQ'},}}
\ \ \
\xymatrix{
{\tE_{\bQ}}\ar[r]^{\Phi^\rs}\ar[d]_{\nu^\rs}&{\tE'_{\bQ'}}\ar[d]^{\nu'^\rs}\\
{\tE}\ar[r]^{\Theta^*}&{\tE',}}
\end{equation}
where $\Theta^*$ is the pullback by the morphism of topos \eqref{p2-fhtft21b}, $\tta_\bQ$ is the ``associated sheaf'' functor, and
$\nu^\rs$ (resp.\ $\nu'^\rs$) is a left adjoint of $\nu_\rs$ (resp.\ $\nu'_\rs$), are commutative up to isomorphism. 

We finally consider the composed functor 
\begin{equation}\label{p2-fhtft35h}
\varphi\colon 
\xymatrix{
{\bQ}\ar[r]^-(0.5)v&{\Et_{/X}}\ar[r]^-(0.5){g^+}&{\Et_{/X'}}\ar[r]^-(0.5){\tth_{/X'}}&{(\Et_{/X'})^\wedge}\ar[r]^-(0.5){v'_\rp}&{\hbQ',}}
\end{equation} 
where $v'_\rp$ is defined by composition with the canonical functor $v'\colon \bQ'\rightarrow \Et_{/X'}$, 
and the associated category $\mI_{\varphi}$, defined in \ref{p2-cmt100},
of triples $(U,U',\mu)$, where $U\in \ob(\bQ)$, $U'\in \ob(\bQ')$ and $\mu\in\varphi(U)(U')$,
that we canonically identify with the category $\mM(\bQ,\bQ')$ defined in \ref{p2-fhtft30}.

By \ref{p2-qfc246}, the functors $(\Phi,\varphi)$ define a functor \eqref{p2-qfc246f}
\begin{equation}
\varpi\colon E'\rhd E_\bQ\rightarrow \Et_{/X'}\rhd \bQ.
\end{equation}
It is quasi-fibering in the sense of \ref{p2-qfc9}, by \ref{p2-qfc28}. 
The following notion is a special case of that introduced in \ref{p2-qfc248}. 

\begin{defi}\label{p2-fhtft36} 
For any objects $F=\{U\mapsto F_U\}$ of $\cP(E_\bQ/\bQ)$ \eqref{p2-qfc20} and $G=\{U'\mapsto G_{U'}\}$ of $\cP(E'_{\bQ'}/\bQ')$, 
an {\em $\mM(\bQ,\bQ')$-system of morphisms from $F$ to $G$} is the data for any object $\mu\colon U'\rightarrow U$ of $\mM(\bQ,\bQ')$ 
of a {\em bifunctorial} morphism
\begin{equation}\label{p2-fhtft36a}
\xi_\mu\colon \tmu^\rp(F_U)\rightarrow G_{U'},
\end{equation}
where $\tmu^\rp\colon  (\Et_{\rf/\oU^\circ})^\wedge\rightarrow (\Et_{\rf/\oU'^\rhd})^\wedge$ is the functor induced by $\tmu$ \eqref{p2-fhtft302a}.  
\end{defi}
The bifunctoriality is expressed as follows. For every morphism $(u',u)\colon (\mu_1\colon U'_1\rightarrow U_1)\rightarrow (\mu_2\colon U'_2\rightarrow U_2)$ 
of $\mM(\bQ',\bQ)$, i.e., a commutative diagram   
\begin{equation}\label{p2-fhtft36b}
\xymatrix{
U'_1\ar[r]^{\mu_1}\ar[d]_{u'}&U_1\ar[d]^u\\
U'_2\ar[r]^{\mu_2}&U_2,}
\end{equation}
the diagram 
\begin{equation}\label{p2-fhtft36c}
\xymatrix{
{\ou'^{\rhd\rp}(\tmu^\rp_2(F_{U_2}))}\ar[rr]^-(0.5){\ou'^{\rhd\rp}(\xi_{\mu_2})}\ar[d]_c&&{\ou'^{\rhd\rp}(G_{U'_2})}\ar[d]^{a'}\\
{\tmu^\rp_1(\ou^{\circ\rp}(F_{U_2}))} \ar[r]^-(0.5){\tmu^\rp_1(a)}&{\tmu_1^\rp(F_{U_1})}\ar[r]^{\xi_{\mu_1}}&{G_{U'_1},}}
\end{equation}
where $c$ is the isomorphism induced by \eqref{p2-fhtft36b}, $a\colon \ou^{\circ\rp}(F_{U_2})\rightarrow F_{U_1}$ 
(resp.\ $a'$) is the adjoint of the transition morphism of $F$ (resp.\ $G$), 
is commutative. 

\begin{rema}\label{p2-fhtft361}  
Let $F=\{U\mapsto F_U\}\in \ob(\cP(E_\bQ/\bQ))$, $G=\{U'\mapsto G_{U'}\}\in \ob(\cP(E'_{\bQ'}/\bQ'))$ 
such that for every $U\in \ob(\bQ)$ (resp.\ $U'\in \ob(\bQ')$),
$F_U$ (resp.\ $G_{U'}$) is a sheaf of $\oU^\circ_\fet$ (resp.\ $\oU'^\rhd_\fet$). 
By \ref{p2-qfc252}, giving an $\mM(\bQ,\bQ')$-system of morphisms from $F$ to $G$ is equivalent to giving, 
for any object $\mu\colon U'\rightarrow U$ of $\mM(\bQ,\bQ')$, a bifunctorial morphism
\begin{equation}\label{p2-fhtft361a} 
\zeta_\mu\colon \tmu^*(F_U)\rightarrow G_{U'},
\end{equation}
where $\tmu^*$ is the pullback functor by the morphism of topos $\tmu\colon \oU'^\rhd_\fet\rightarrow \oU^\circ_\fet$ \eqref{p2-fhtft302a}.  
\end{rema}
The bifunctoriality is expressed as follows. For every morphism $(u',u)\colon (\mu_1\colon U'_1\rightarrow U_1)\rightarrow (\mu_2\colon U'_2\rightarrow U_2)$ 
of $\mM(\bQ',\bQ)$, see \eqref{p2-fhtft36b}, the diagram 
\begin{equation}\label{p2-fhtft361b} 
\xymatrix{
{\ou'^{\rhd*}(\tmu^*_2(F_{U_2}))}\ar[rr]^-(0.5){\ou'^{\rhd*}(\xi_{\mu_2})}\ar[d]_c&&{\ou'^{\rhd*}(G_{U'_2})}\ar[d]^{a'}\\
{\tmu^*_1(\ou^{\circ*}(F_{U_2}))} \ar[r]^-(0.5){\tmu^*_1(a)}&{\tmu_1^*(F_{U_1})}\ar[r]^{\xi_{\mu_1}}&{G_{U'_1},}}
\end{equation}
where $c$ is the isomorphism induced by \eqref{p2-fhtft36b}, $a\colon \ou^{\circ*}(F_{U_2})\rightarrow F_{U_1}$ 
(resp.\ $a'$) is the adjoint of the transition morphism of $F$ (resp.\ $G$), 
is commutative.

\begin{lem}\label{p2-fhtft362} 
Let $F\in \ob(\hE_\bQ)$, $G\in \ob(\hE'_{\bQ'})$. We denote by $\{U\mapsto F_U\}$ (resp.\ $\{U'\mapsto G_{U'}\}$) the image of $F$ (resp.\ $G$)
in $\cP(E_\bQ/\bQ)$ (resp.\ $\cP(E'_{\bQ'}/\bQ')$) by the canonical equivalence of categories \eqref{p2-htaft10c},
and by $F^a$ (resp.\ $G^a$) the associated sheaf of $\tE_\bQ$ (resp.\ $\tE'_{\bQ'}$). 
Then, with the notation of \ref{p2-fhtft35}, the following data are equivalent:
\begin{itemize}
\item[(i)] a morphism $F\rightarrow \Phi_\rp(G)$ of $\hE_\bQ$;
\item[(ii)] a morphism $\Phi^\rp(F)\rightarrow G$ of $\hE'_{\bQ'}$;
\item[(iii)] an $\mM(\bQ,\bQ')$-system of morphisms from $\{U\mapsto F_U\}$ to $\{U'\mapsto G_{U'}\}$. 
\end{itemize}

Moreover, such a data determines a morphism $\Theta^*(\nu^\rs(F^a))\rightarrow \nu'^\rs(G^a)$ of $\tE'$ \eqref{p2-fhtft21b}. 
This correspondence is bifunctorial in $F$ and $G$.
\end{lem}
 
It follows from \ref{p2-qfc251}, \ref{p2-qfc253} and \eqref{p2-fhtft35g}.

\begin{rema}\label{p2-fhtft37}
We can define the notion of {\em $\mM$-system of morphisms} from an object of $\cP(E/\Et_{/X})$ to an object of $\cP(E'/\Et_{/X'})$, 
as a special case of that introduced in \ref{p2-qfc29}, and we can apply \ref{p2-qfc251}.
For instance, by \ref{p2-qfc256}, the canonical homomorphism $\ocB\rightarrow \Theta_*(\ocB')$ corresponds to the $\mM$-system of morphisms from 
$\{U \mapsto \ocB_U\}$ to $\{U'\mapsto \ocB'_{U'}\}$ defined for any object $\mu\colon U'\rightarrow U$ of $\mM$ 
by the morphism $\ocB_U\rightarrow \tmu_*(\ocB'_{U'})$ \eqref{p2-fhtft302b} of $\oU^\circ_\fet$. 
By restriction, we deduce an $\mM(\bQ,\bQ')$-system of morphisms from $\{U\in \bQ \mapsto \ocB_U\}$ to $\{U'\in \bQ'\mapsto \ocB'_{U'}\}$, 
which by \ref{p2-fhtft362} determines a homomorphism $\ocB\rightarrow \Theta_*(\ocB')$ of $\tE$. 
The latter is none other than the canonical homomorphism by \eqref{p2-cmt77b} and \ref{p2-qfc255}.
\end{rema}

\subsection{}\label{p2-fhtft32}
Recall that we use the notation introduced in \ref{p2-htaft9}--\ref{p2-htaft17} both for $(f,\tf)$ and $(f',\tf')$, 
except that for $(f',\tf')$ we equip them with a $^\prime$ exponent. 
Let $\mu\colon Y'\rightarrow Y$ be an object of $\mM(\bQ',\bQ)$ \eqref{p2-fhtft30}, $\oy'$ a geometric point of $\oY'^\rhd$, $\oy=\tmu(\oy')$. 
Since $Y\rightarrow X$ (resp.\ $Y'\rightarrow X'$) is an object of $\bQ$ (resp.\ $\bQ'$), 
the modules and algebras introduced in \ref{p2-htaft12} 
do not depend on the choice of an adequate chart for the morphism $f|Y\colon (Y,\cM_{X}|Y)\rightarrow (S,\cM_S)$ 
(resp.\ $f'|Y'\colon (Y',\cM_{X'}|Y')\rightarrow (S,\cM_S)$). 
Indeed by \ref{p2-rlps4}, with the notation of \ref{p2-rlps3}, we can obtain those introduced in \ref{p2-htaft12} 
using the logarithmic scheme $(\tmY^{\oy},\cM'_{\tmY^{\oy}})$ defined as one or the other of the logarithmic schemes
\begin{equation}\label{p2-fhtft32a}
(\cA_2(\mY^{\oy}),\cM'_{\cA_2(\mY^{\oy})})\ \ \ {\rm or} \ \ \ (\cA^{\ast}_2(\mY^{\oy}/S),\cM'_{\cA^{\ast}_2(\mY^{\oy}/S)})
\end{equation}
associated with $f|Y$ (\cite{ag2} 3.2.8 and 3.2.10), 
depending on whether we are in the absolute or relative case \eqref{p2-ncgt3} 
respectively (mind the $^\prime$ in $\cM'_{\tmY^{\oy}}$; see \cite{ag2} 3.2.13). 
The diagram
\begin{equation}\label{p2-fhtft32b}
\xymatrix{
{(\hmY'^{\oy'},\cM_{\hmY'^{\oy'}})}\ar[r]^-(0.5){\mj'^{\oy'}_{Y'}}\ar[d]_{\hupmu}&{(\tmY'^{\oy'},\cM'_{\tmY'^{\oy'}})}\ar[d]\\
{(\hmY^{\oy},\cM_{\hmY^{\oy}})}\ar[r]^-(0.5){\mj^{\oy}_Y}&{(\tmY^{\oy},\cM'_{\tmY^{\oy}}),}}
\end{equation}
where the morphism $\hupmu$ is induced by \eqref{p2-fhtft302d}
and $\mj'^{\oy'}_{Y'}$ (resp.\ $\mj^{\oy}_Y$) is the exact closed immersion \eqref{p2-rlps3d}, is commutative. 
We can then form the commutative diagram 
\begin{equation}\label{p2-fhtft32c}
\xymatrix{
&{(\hmY'^{\oy'},\cM_{\hmY'^{\oy'}})}\ar[ld]\ar[rr]^{\mj'^{\oy'}_{Y'}}\ar@{->}'[d]_-(0.5){\hupmu}[dd]&&{(\tmY'^{\oy'},\cM'_{\tmY'^{\oy'}})}\ar[dd]\ar@{.>}[ld]\\
{(\coX',\cM_{\coX'})}\ar[rr]^-(0.7){i'}\ar[dd]_{\cog}&&{(\tX',\cM_{\tX'})}\ar@/^2pc/[ddd]^{\tf'}\ar@{.>}[dd]&\\
&{(\hmY^{\oy},\cM_{\hmY^{\oy}})}\ar[ld]\ar[rr]^-(0.3){\mj^{\oy}_Y}|-(0.5)\hole|-(0.7)\hole\ar[ddl]|\hole&&{(\tmY^{\oy},\cM'_{\tmY^{\oy}})}\ar[ldd]\\
{(\coX,\cM_\coX)}\ar[rr]^i\ar[d]_{\cof}&&{(\tX,\cM_\tX)}\ar[d]_{\tf}\\
{(\coS,\cM_{\coS})}\ar[rr]^\iota&&{(\tS,\cM_{\tS}).}}
\end{equation}

\subsection{}\label{p2-fhtft330}
We keep the assumption and notation of \ref{p2-fhtft32}. 
For any rational numbers $r\geq r'\geq 0$, 
we associate with the diagram \eqref{p2-fhtft32c} a canonical morphism of $\hoRp^{\oy'}_{Y'}$-modules \eqref{p2-hta150h}
\begin{equation}\label{p2-fhtft330a}
\upvarphi^{\oy',(r,r')}_\mu\colon \fF^{\oy,(r)}_Y\otimes_{\hoR^\oy_Y}\hoRp^{\oy'}_{Y'}\rightarrow \fF'^{\oy,(r)}_{Y'} \otimes_{R'_{\uptau,Y'}}\fF^{(r')}_{\uptau,Y'},
\end{equation}
where $R'_{\uptau,Y'}$ and $\fF^{(r')}_{\uptau,Y'}$ are defined in \eqref{p2-fhtft3g}, 
and a canonical morphism of $\hoRp^{\oy'}_{Y'}$-algebras \eqref{p2-hta150i}
\begin{equation}\label{p2-fhtft330b}
\upphi^{\oy',(r,r')}_\mu\colon \fC^{\oy,(r)}_Y\otimes_{\hoR^\oy_Y}\hoRp^{\oy'}_{Y'}\rightarrow \fC'^{\oy,(r)}_{Y'} \otimes_{R'_{\uptau,Y'}}\fC^{(r')}_{\uptau,Y'},
\end{equation}
compatible with $\upvarphi^{\oy',(r,r')}_\mu$, see \ref{p2-fhtal7} and \ref{p2-fhtal110}. 
The latter induces by linearization a morphism of $(\hoRp^{\oy'}_{Y'}\otimes_{R'_{\uptau,Y'}}\fC^{(r')}_{\uptau,Y'})$-algebras \eqref{p2-hta150e}
\begin{equation}\label{p2-fhtft330c}
\uppsi^{\oy',(r,r')}_\mu\colon \fC^{\oy,(r)}_Y\otimes_{\hoR^\oy_Y}\hoRp^{\oy'}_{Y'}\otimes_{R'_{\uptau,Y'}}\fC^{(r')}_{\uptau,Y'}
\rightarrow \fC'^{\oy,(r)}_{Y'} \otimes_{R'_{\uptau,Y'}}\fC^{(r')}_{\uptau,Y'}.
\end{equation}
These morphisms are $\pi_1(\oY'^{\star\rhd},\oy')$-equivariant, when we equip $R'_{\uptau,Y'}$, $\fF^{(r')}_{\uptau,Y'}$ and $\fC^{(r')}_{\uptau,Y'}$ 
with the trivial actions; see \ref{p2-fhtal12}. 
By (\cite{agt} VI.9.11), for every integer $n\geq 0$, we deduce a canonical morphism of $\ocB'_{Y',n}$-modules of $\oY'^\rhd_\fet$
\begin{equation}\label{p2-fhtft330d}
\upvarphi^{(r,r')}_{\mu,n}\colon \tmu_n^*(\cF^{(r)}_{Y,n})\rightarrow \cF'^{(r)}_{Y',n} \otimes_{R'_{\uptau,Y'}}\fF^{(r')}_{\uptau,Y'},
\end{equation}
where $\tmu_n$ is the morphism of ringed topos \eqref{p2-fhtft302c}, 
and a canonical  morphism of $\ocB'_{Y',n}$-algebras 
\begin{equation}\label{p2-fhtft330e}
\upphi^{(r,r')}_{\mu,n}\colon \tmu_n^*(\cC^{(r)}_{Y,n})\rightarrow \cC'^{(r)}_{Y',n} \otimes_{R'_{\uptau,Y'}}\fC^{(r')}_{\uptau,Y'},
\end{equation}
compatible with $\upvarphi^{(r,r')}_{\mu,n}$. These morphisms are independent of the choice of $\oy'$. 
The morphism $\upphi^{(r,r')}_{\mu,n}$ induces by linearization a morphism of $(\ocB'_{Y',n}\otimes_{R'_{\uptau,Y'}}\fC^{(r')}_{\uptau,Y'})$-algebras 
\begin{equation}\label{p2-fhtft330f}
\uppsi^{(r,r')}_{\mu,n}\colon \tmu_n^*(\cC^{(r)}_{Y,n})\otimes_{R'_{\uptau,Y'}}\fC^{(r')}_{\uptau,Y'}
\rightarrow \cC'^{(r)}_{Y',n} \otimes_{R'_{\uptau,Y'}}\fC^{(r')}_{\uptau,Y'}.
\end{equation}

\subsection{}\label{p2-fhtft340}
Let $(u',u)\colon (\mu_1\colon Y'_1\rightarrow Y_1)\rightarrow (\mu_2\colon Y'_2\rightarrow Y_2)$ be a morphism of $\mM(\bQ',\bQ)$ \eqref{p2-fhtft30}, 
$\oy'_1$ a geometric point 
of $\oY'^\rhd_1$, $\oy'_2=\ou'^\rhd(\oy'_1)$, $\oy_1=\tmu_1(\oy'_1)$, $\oy_2=\tmu_2(\oy'_2)$ \eqref{p2-fhtft302a}, 
$r\geq r'\geq 0$ rational numbers, $n$ an integer $\geq 0$. 
\begin{equation}\label{p2-fhtft340a}
\xymatrix{
Y'_1\ar[r]^{\mu_1}\ar[d]_{u'}&Y_1\ar[d]^u\\
Y'_2\ar[r]^{\mu_2}&Y_2.}
\end{equation}
We denote by 
\begin{eqnarray}
\ou'^\rhd_n\colon (\oY'^\rhd_{1,\fet},\ocB'_{Y'_1,n})&\rightarrow&(\oY'^\rhd_{2,\fet},\ocB'_{Y'_2,n}),\\
\ou^\circ_n\colon (\oY^\circ_{1,\fet},\ocB_{Y_1,n})&\rightarrow&(\oY^\circ_{2,\fet},\ocB_{Y_2,n}),
\end{eqnarray}
the canonical morphisms of ringed topos induced by $u'$ and $u$, respectively. Diagram \eqref{p2-fhtft34a} is underlying a 
commutative diagram of morphisms of ringed topos
\begin{equation}\label{p2-fhtft340f}
\xymatrix{
{(\oY'^\rhd_{1,\fet},\ocB'_{Y'_1,n})}\ar[r]^{\tmu_{1,n}}\ar[d]_{\ou'^\rhd_n}&{(\oY^\circ_{1,\fet},\ocB_{Y_1,n})}\ar[d]^{\ou^\circ_n}\\
{(\oY'^\rhd_{2,\fet},\ocB'_{Y'_2,n})}\ar[r]^{\tmu_{2,n}}&{(\oY^\circ_{2,\fet},\ocB_{Y_2,n}),}}
\end{equation}
where the morphisms $\tmu_{i,n}$ are defined in \eqref{p2-fhtft302c}.

We check immediately that the diagram
\begin{equation}\label{p2-fhtft340b}
\xymatrix{
{(\tmY'^{\oy'_1}_1,\cM'_{\tmY'^{\oy'_1}_1})}\ar[rrr]\ar[ddd]&&&{(\tmY^{\oy_1}_1,\cM'_{\tmY^{\oy_1}_1})}\ar[ddd]\\
&{(\hmY'^{\oy'_1}_1,\cM_{\hmY'^{\oy'_1}_1})}\ar[r]\ar[d]\ar[lu]_{\mj^{\oy'_1}_{Y'_1}}&{(\hmY^{\oy_1}_1, \cM_{\hmY^{\oy_1}_1})}\ar[d]\ar[ru]^{\mj^{\oy_1}_{Y_1}}&\\
&{(\hmY'^{\oy'_2}_2,\cM_{\hmY'^{\oy'_2}_2})}\ar[r]\ar[ld]_{\mj^{\oy'_2}_{Y'_2}}&{(\hmY^{\oy_2}_2,\cM_{\hmY^{\oy_2}_2 })}\ar[rd]^{\mj^{\oy_2}_{Y_2}}&\\
{(\tmY'^{\oy'_2}_2,\cM'_{\tmY'_2})}\ar[rrr]&&&{(\tmY^{\oy_2}_2,\cM'_{\tmY_2})}}
\end{equation}
is commutative. We claim that the diagram
\begin{equation}\label{p2-fhtft340c}
\xymatrix{
{\ou'^{\rhd*}_n(\tmu_{2,n}^*(\cF^{(r)}_{Y_2,n}))}\ar[rr]^-(0.5){\ou'^{\rhd*}_n(\upvarphi^{(r,r')}_{\mu_2,n})}\ar[d]_c&&
{\ou'^{\rhd*}_n(\cF'^{(r)}_{Y' _2,n}\otimes_{R'_{\uptau,Y'_2}}\fF^{(r')}_{\uptau,Y'_2})}\ar[d]^{a'\otimes a_\uptau}\\
{\tmu_{1,n}^*(\ou^{\circ*}_n(\cF^{(r)}_{Y_2,n}))}\ar[r]^-(0.5){\tmu_{1,n}^*(a) }&{\tmu_{1,n}^*(\cF^{(r)}_{Y_1,n})}\ar[r]^-(0.5){\upvarphi^{(r,r')}_{\upmu_1,n}} &
{\cF'^{(r)}_{Y'_1,n}\otimes_{R'_{\uptau,Y'_1}}\fF^{(r')}_{\uptau,Y'_1}},}
\end{equation}
where $a\colon \ou^{\circ*}_n(\cF^{(r)}_{Y_2,n})\rightarrow \cF^{(r)}_{Y_1,n}$ 
(resp.\ $a'\colon \ou'^{\rhd*}_n(\cF'^{(r)}_{Y'_2,n})\rightarrow \cF'^{(r)}_{Y'_1,n}$) 
is the adjoint of the transition morphism of the v-presheaf
$\{Y\in \bQ^\circ\mapsto \cF^{(r)}_{Y,n}\}$ \eqref{p2-htaft16d} (resp.\ $\{Y'\in \bQ'^\circ\mapsto \cF'^{(r)}_{Y',n}\}$) \eqref{p2-cmt45},
$a_\uptau\colon \fF^{(r')}_{\uptau,Y'_2}\rightarrow \fF^{(r')}_{\uptau,Y'_1}$ is the canonical morphism \eqref{p2-fhtft3g}, 
and $c$ is the isomorphism induced by the commutative diagram \eqref{p2-fhtft340a}, is commutative. 
Indeed, we are immediately reduced to checking the commutativity of the analogous diagram for the Higgs--Tate extensions appearing in \eqref{p2-fhtft330a}, 
and furthermore to the case where $r=r'=0$ as these Higgs--Tate extensions are submodules of those with vanishing radii; 
then the statement follows from \eqref{p2-fhtft340b} by \ref{p1-rdt40} and \ref{p1-rdt9}. 
Observe that the invertibility condition required in \ref{p1-rdt4} is satisfied by \ref{p2-hta2}.

We deduce from \eqref{p2-fhtft340c} that the diagram
\begin{equation}\label{p2-fhtft340d}
\xymatrix{
{\ou'^{\rhd*}_n(\tmu_{2,n}^*(\cC^{(r)}_{Y_2,n}))}\ar[rr]^-(0.5){\ou'^{\rhd*}_n(\upphi^{(r,r')}_{\mu_2,n})}\ar[d]_d&&
{\ou'^{\rhd*}_n(\cC'^{(r)}_{Y'_2,n}\otimes_{R'_{\uptau,Y'_2}}\fC^{(r')}_{\uptau,Y'_2})}\ar[d]^{b'\otimes b_\uptau}\\
{\tmu_{1,n}^*(\ou^{\circ*}_n(\cC^{(r)}_{Y_2,n}))}\ar[r]^-(0.5){\tmu_{1,n}^*(b)}&{\tmu_{1,n}^*(\cC^{(r)}_{Y_1,n})}\ar[r]^-(0.5){\upphi^{(r,r')}_{\mu_1,n}}
&{\cC'^{(r)}_{Y'_1,n}\otimes_{R'_{\uptau,Y'_1}}\fC^{(r')}_{\uptau,Y'_1}}}
\end{equation}
is commutative, 
where $b\colon \ou^{\circ*}_n(\cC^{(r)}_{Y_2,n})\rightarrow \cC^{(r)}_{Y_1,n}$ 
(resp.\ $b'\colon \ou'^{\rhd*}_n(\cC'^{(r)}_{Y'_2,n})\rightarrow \cC'^{(r)}_{Y'_1,n}$) 
is the adjoint of the transition morphism of the v-presheaf
$\{Y\in \bQ^\circ\mapsto \cC^{(r)}_{Y,n}\}$ (resp.\ $\{Y'\in \bQ'^\circ\mapsto \cC'^{(r)}_{Y',n}\}$) \eqref{p2-htaft16e}, 
$b_\uptau\colon \fC^{(r')}_{\uptau,Y'_2}\rightarrow \fC^{(r')}_{\uptau,Y'_1}$ is the canonical morphism \eqref{p2-fhtft3g}, 
and $d$ is the isomorphism induced by the commutative diagram \eqref{p2-fhtft340a}.
Therefore, by linearization, the diagram 
\begin{equation}\label{p2-fhtft340e}
{\small 
\xymatrix{
{\ou'^{\rhd*}_n(\tmu_{2,n}^*(\cC^{(r)}_{Y_2,n})\otimes_{R'_{\uptau,Y'_2}}\fC^{(r')}_{\uptau,Y'_2})}\ar[rrr]^-(0.5){\ou'^{\rhd*}_n(\uppsi^{(r,r')}_{\mu_2,n})}\ar[d]_{d\otimes b_\uptau}&&&
{\ou'^{\rhd*}_n(\cC'^{(r)}_{Y'_2,n}\otimes_{R'_{\uptau,Y'_2}}\fC^{(r')}_{\uptau,Y'_2})}\ar[d]^{b'\otimes b_\uptau}\\
{\tmu_{1,n}^*(\ou^{\circ*}_n(\cC^{(r)}_{Y_2,n}))\otimes_{R'_{\uptau,Y'_1}}\fC^{(r')}_{\uptau,Y'_1}}\ar[rr]^-(0.5){\tmu_{1,n}^*(b)\otimes \id}&&
{\tmu_{1,n}^*(\cC^{(r)}_{Y_1,n})\otimes_{R'_{\uptau,Y'_1}}\fC^{(r')}_{\uptau,Y'_1}}\ar[r]^-(0.5){\uppsi^{(r,r')}_{\mu_1,n}}
&{\cC'^{(r)}_{Y'_1,n}\otimes_{R'_{\uptau,Y'_1}}\fC^{(r')}_{\uptau,Y'_1}}}}
\end{equation}
is commutative.

\subsection{}\label{p2-fhtft380}
Let $r\geq r'\geq 0$ be rational numbers, $n\geq 0$ an integer. We have a natural object 
\begin{equation}\label{p2-fhtft380a}
\{U'\mapsto \cF'^{(r)}_{U',n}\otimes_{R'_{\uptau,U'}}\fF^{(r')}_{\uptau,U'}\}
\end{equation}
of $\cP(E'_{\bQ'}/\bQ')$; see \ref{p2-htaft14} and \ref{p2-fhtft3}. 
By \eqref{p2-htaft10c}, it determines a presheaf of $\nu'_\rp(\ocB')$-modules on $E'_{\bQ'}$ \eqref{p2-htaft10}. 
By \eqref{p2-htaft10e} and (\cite{agt} VI.534(ii), VI.8.9 and VI.5.17), we have a canonical isomorphism
\begin{equation}\label{p2-fhtft380b}
\nu'_\rs(\cF'^{(r)}_n\otimes_{\ocB'_n}\sigma'^*_n(\cF^{(r')}_{\uptau,n}))\stackrel{\sim}{\rightarrow}
\{U'\mapsto \cF'^{(r)}_{U',n}\otimes_{R'_{\uptau,U'}}\fF^{(r')}_{\uptau,U'}\}^a,
\end{equation}
where the exponent $a$ means the associated sheaf, $\sigma'_n$ denotes the morphism of ringed topos \eqref{p2-htaft7c} and  
$\nu'_\rs$ is the equivalence of categories \eqref{p2-htaft10d} (for $f'$).

It follows from \eqref{p2-fhtft340c} that the morphisms $\upvarphi^{(r,r')}_{\mu,n}$ \eqref{p2-fhtft330d}, for $\mu\in \mM(\bQ,\bQ')$, 
induce an $\mM(\bQ,\bQ')$-system of morphisms from $\{U\mapsto \cF^{(r)}_{U,n}\}$ to 
$\{U'\mapsto \cF'^{(r)}_{U',n}\otimes_{R'_{\uptau,U'}}\fF^{(r')}_{\uptau,U'}\}$ in the sense of \ref{p2-fhtft36}; see \ref{p2-fhtft361}.
By \ref{p2-fhtft362}, the latter determines a morphism $\Theta^*(\cF^{(r)}_n)\rightarrow \cF'^{(r)}_n\otimes_{\ocB'_n}\sigma'^*_n(\cF^{(r')}_{\uptau,n})$ 
of $\tE'$, which is obviously $\Theta^*(\ocB)$-linear \eqref{p2-htaft10e}. It therefore induces a $\ocB'_n$-linear morphism 
\begin{equation}\label{p2-fhtft380c}
\upvarphi^{(r,r')}_n\colon \uptheta^*_n(\cF^{(r)}_n)\rightarrow \cF'^{(r)}_n\otimes_{\ocB'_n}\sigma'^*_n(\cF^{(r')}_{\uptau,n}),
\end{equation}
where $\uptheta_n$ is the morphism of ringed topos \eqref{p2-fhtft21h}. 

Similarly, the morphisms $\upphi^{(r,r')}_{\mu,n}$ \eqref{p2-fhtft330e}, for $\mu\in \mM(\bQ,\bQ')$, induce a morphism of $\ocB'_n$-algebras 
\begin{equation}\label{p2-fhtft380d}
\upphi^{(r,r')}_n\colon \uptheta^*_n(\cC^{(r_2)}_n)\rightarrow \cC'^{(r)}_n\otimes_{\ocB'_n}\sigma'^*_n(\cC^{(r')}_{\uptau,n}),
\end{equation}
compatible with $\upvarphi^{(r,r')}_n$. It induces a homomorphism of $\sigma'^*_n(\cC^{(r')}_{\uptau,n})$-algebras 
\begin{equation}\label{p2-fhtft380e}
\uppsi^{(r,r')}_n\colon \uptheta^*_n(\cC^{(r)}_n)\otimes_{\ocB'_n}\sigma'^*_n(\cC^{(r')}_{\uptau,n})\rightarrow 
\cC'^{(r)}_n\otimes_{\ocB'_n}\sigma'^*_n(\cC^{(r')}_{\uptau,n}).
\end{equation}

By (\cite{agt} III.7.5 and III.7.12), the homomorphisms $(\upphi^{(r,r')}_{m+1})_{m\in \mN}$ define a homomorphism of $\bvocB'$-algebras 
\begin{equation}\label{p2-fhtft380f}
\bvupphi^{(r,r')}\colon \bvuptheta^*(\bvcC^{(r)})\rightarrow \bvcC'^{(r)}\otimes_{\bvocB'}\hupsigma'^*(\hcC^{(r')}_{\uptau}),
\end{equation}
where $\bvcC^{(r)}$ (resp.\ $\bvcC'^{(r)}$, resp.\ $\hcC^{(r')}_{\uptau}$) is defined in \ref{p2-htaft18} (resp.\ \ref{p2-htaft18}, 
resp.\ \ref{p2-fhtft3}) and $\bvuptheta^*$ and $\hupsigma'$ are the morphism of ringed topos appearing in \eqref{p2-fhtft22b}. 
It induces a homomorphism of $\hupsigma'^*(\hcC^{(r_3)}_\uptau)$-algebras 
\begin{equation}\label{p2-fhtft380g}
\bvuppsi^{(r,r')}\colon \bvuptheta^*(\bvcC^{(r)})\otimes_{\bvocB'}\hupsigma'^*(\hcC^{(r')}_\uptau)\rightarrow 
\bvcC'^{(r)}\otimes_{\bvocB'}\hupsigma'^*(\hcC^{(r')}_{\uptau}).
\end{equation}

Let $t\geq t'\geq 0$ be rational numbers such that $r\geq t$ and $r'\geq t'$. 
By \ref{p2-fhtft362} and \eqref{p2-fhtal110f}, the diagram 
\begin{equation}\label{p2-fhtft380i}
\xymatrix{
{\uptheta^*_n(\cC^{(r)}_n)}\ar[r]^-(0.5){\upphi^{(r,r')}_n}\ar[d]_{\uptheta^*_n(\alpha^{r,t}_n)}
&{\cC'^{(r)}_n\otimes_{\ocB'_n}\sigma'^*_n(\cC^{(r')}_{\uptau,n})}
\ar[d]^{\alpha'^{r,t}_n\otimes \sigma'^*_n(\alpha_{\uptau,n}^{r',t'})}\\
{\uptheta^*_n(\cC^{(t)}_n)}\ar[r]^-(0.5){\upphi^{(t,t')}_n}&
{\cC'^{(t)}_n\otimes_{\ocB'_n}\sigma'^*_n(\cC^{(t')}_{\uptau,n}),}}
\end{equation}
where $\alpha^{r,t}_n$ and $\alpha'^{r,t}_n$ are defined in \eqref{p2-htaft16g} and 
$\alpha_{\uptau,n}^{r',t'}$ is induced by $\alpha_{\uptau}^{r',t'}$ \eqref{p2-fhtft3i}, is commutative. 
We deduce that the diagram 
\begin{equation}\label{p2-fhtft380h}
\xymatrix{
{\bvuptheta^*(\bvcC^{(r)})}\ar[r]^-(0.5){\bvupphi^{(r,r')}} \ar[d]_{\bvuptheta^*(\bvalpha^{r,t})}
&{\bvcC'^{(r)}\otimes_{\bvocB'}\hupsigma'^*(\bvcC^{(r')}_{\uptau})}
\ar[d]^{\bvalpha'^{r,t}\otimes_{\bvocB'} \hupsigma'^*(\halpha_\uptau^{r',t'})}\\
{\bvuptheta^*(\bvcC^{(t)})}\ar[r]^-(0.5){\bvuppsi^{(t,t')}}&
{\bvcC'^{(t)}\otimes_{\bvocB'}\hupsigma'^*(\bvcC^{(t')}_{\uptau}),}}
\end{equation}
where $\bvalpha^{r,t}$ and $\bvalpha'^{r,t}$ are defined in \eqref{p2-htaft18f} and $\halpha_\uptau^{r',t'}$  
is induced by $\alpha_{\uptau}^{r',t'}$ \eqref{p2-fhtft3i}, is commutative. 
We have similar commutative diagrams for $\uppsi^{(r,r')}_n$ and $\bvuppsi^{(r,r')}$.

\subsection{}\label{p2-rftchta150}
Let $r\geq r'\geq 0$ be rational numbers, $n\geq 1$ an integer. We denote by 
\begin{equation}\label{p2-rftchta150a}
\delta^{(r,r')}_n\colon \cC'^{(r)}_n\otimes_{\ocB'_n}\sigma'^*_n(\cC^{(r')}_{\uptau,n})\rightarrow 
\sigma'^*_n(\Omega'_n)\otimes_{\ocB'_n}\cC'^{(r)}_n\otimes_{\ocB'_n}\sigma'^*_n(\cC^{(r')}_{\uptau,n})
\end{equation}
the $\ocB'_n$-derivation defined by 
\begin{equation}\label{p2-rftchta150b}
\delta^{(r,r')}_n=\delta_{\cC'^{(r)}_n}\otimes \id- \id \otimes \sigma'^*_n(\delta'_{\cC^{(r')}_{\uptau,n}}),
\end{equation}
where $\delta_{\cC'^{(r)}_n}$ is the derivation \eqref{p2-htaft20b} (for $(f',\tf')$) 
and $\delta'_{\cC^{(r)}_{\uptau,n}}$ is the reduction modulo $p^n$ of the derivation $\delta'_{\cC^{(r')}_{\uptau}}$ \eqref{p2-fhtft6aa}.  
It is a Higgs $\ocB'_n$-field with coefficients in $\sigma'^*_n(\Omega'_n)$.  By \ref{p2-fhtft362} and \eqref{p2-fhtal150c}, we have a commutative diagram
\begin{equation}\label{p2-rftchta150c}
\xymatrix{
{\uptheta^*_n(\cC^{(r)}_n)\otimes_{\ocB'_n}\sigma'^*_n(\cC^{(r')}_{\uptau,n})}\ar[rr]^-(0.5){\id\otimes \sigma'^*_n(\delta_{\cC^{(r')}_{\uptau,n}})}
\ar[d]_{\uppsi^{(r,r')}_n}&&
{\sigma'^*_n(\cog^*(\Omega_n))\otimes_{\ocB'_n}\uptheta_n^*(\cC^{(r)}_n)\otimes_{\ocB'_n}\sigma'^*_n(\cC^{(r')}_{\uptau,n})}
\ar[d]^-(0.5){\sigma'^*_n(u_n)\otimes \uppsi^{(r,r')}_n}\\
{\cC'^{(r)}_n\otimes_{\ocB'_n}\sigma'^*_n(\cC^{(r')}_{\uptau,n})}\ar[rr]^-(0.5){-\delta^{(r,r')}_n}&&
{\sigma'^*_n(\Omega'_n)\otimes_{\ocB'_n}\cC'^{(r)}_n\otimes_{\ocB'_n}\sigma'^*_n(\cC^{(r')}_{\uptau,n}),}}
\end{equation}
where $\uppsi^{(r,r')}_n$ is the homomorphism \eqref{p2-fhtft380e} and $u_n$ (resp.\ $\delta_{\cC^{(r')}_{\uptau,n}}$) 
is the reduction modulo $p^n$ of the morphism $u$ defined in \eqref{p2-fhtft90l} (resp.\ derivation $\delta_{\cC^{(r')}_{\uptau}}$ \eqref{p2-fhtft6a}).
In particular, we have 
\begin{equation}
\delta^{(r,r')}_n\circ \upphi^{(r,r')}_n=0,
\end{equation}
where $\upphi^{(r,r')}_n$ is the homomorphism defined in \eqref{p2-fhtft380d}. 

By \ref{p2-fhtft362} and \eqref{p2-fhtal150cc}, we have a commutative diagram
\begin{equation}\label{p2-rftchta150i}
\xymatrix{
{\uptheta_n^*(\cC_n^{(r)})\otimes_{\ocB'_n}\sigma'^*_n(\cC^{(r')}_{\uptau,n})}\ar[rr]^-(0.5){b_n}
\ar[d]_{\uppsi^{(r,r')}_n}&&
{\sigma'^*_n(\cog^*(\Omega_n))\otimes_{\ocB'_n}\uptheta^*_n(\cC_n^{(r)})\otimes_{\ocB'_n}\sigma'^*_n(\cC^{(r')}_{\uptau,n})}
\ar[d]^-(0.5){\id\otimes \uppsi^{(r,r')}_n}\\
{\cC'^{(r)}_n\otimes_{\ocB'_n}\sigma'^*_n(\cC^{(r')}_{\uptau,n})}\ar[rr]^-(0.5){a_n}&&
{\sigma'^*_n(\cog^*(\Omega_n))\otimes_{\bvocB'}\cC'^{(r)}_n\otimes_{\ocB'_n}\sigma'^*_n(\cC^{(r')}_{\uptau,n}),}}
\end{equation}
where $a_n$ denotes the Higgs field 
\begin{equation}
a_n=\id\otimes \sigma'^*_n(\delta_{\cC^{(r')}_{\uptau,n}})
\end{equation}
on $\cC'^{(r)}_n\otimes_{\ocB'_n}\sigma'^*_n(\cC^{(r')}_{\uptau,n})$, $b_n$ denotes the Higgs field 
\begin{equation}
b_n=\uptheta^*_n(\delta_{\cC_n^{(r)}})\otimes\id+\id\otimes \sigma'^*_n(\delta_{\cC^{(r')}_{\uptau,n}})
\end{equation}
on $\theta^*_n(\cC_n^{(r)})\otimes_{\ocB'_n}\sigma'^*_n(\cC^{(r')}_{\uptau,n})$, and $\delta_{\cC_n^{(r)}}$ is the derivation \eqref{p2-htaft20b}.

\subsection{}\label{p2-rftchta170}
Let $r\geq r'\geq 0$ be rational numbers. We denote by 
\begin{equation}\label{p2-rftchta170a}
\bvdelta^{(r,r')}\colon \bvcC'^{(r)}\otimes_{\ocB'}\hupsigma'^*(\hcC^{(r')}_{\uptau})\rightarrow 
\hupsigma'^*(\hOmega')\otimes_{\bvocB'}\bvcC'^{(r)}\otimes_{\ocB'}\hupsigma'^*(\hcC^{(r')}_{\uptau})
\end{equation}
the $\bvocB'$-derivation defined by 
\begin{equation}\label{p2-rftchta170ab}
\bvdelta^{(r,r')}=\delta_{\bvcC'^{(r)}}\otimes \id- \id \otimes \hupsigma'^*(\delta'_{\hcC^{(r')}_{\uptau}}),
\end{equation}
where $\delta_{\bvcC'^{(r)}}$ is the derivation \eqref{p2-htaft19b} (for $(f',\tf')$) and $\delta'_{\hcC^{(r')}_{\uptau}}$ is 
the derivation \eqref{p2-fhtft6bb}. By \eqref{p2-htaft8e}, it identifies canonically with the derivation $(\delta^{(r,r')}_{n})_{\mN}$ \eqref{p2-rftchta150a}, 
and it is a Higgs $\bvocB'$-field with coefficients in $\hupsigma'^*(\hOmega')$. 

By \eqref{p2-rftchta150c}, we have a commutative diagram
\begin{equation}\label{p2-rftchta170k}
\xymatrix{
{\bvuptheta^*(\bvcC^{(r)})\otimes_{\bvocB'}\hupsigma'^*(\hcC^{(r')}_{\uptau})}\ar[rr]^-(0.5){\id\otimes \hupsigma'^*(\delta_{\hcC^{(r')}_{\uptau}})}
\ar[d]_{\bvuppsi^{(r,r')}}&&
{\hupsigma'^*(\fgg^*(\hOmega))\otimes_{\bvocB'}\bvuptheta^*(\bvcC^{(r)})\otimes_{\bvocB'}\hupsigma'^*(\hcC^{(r')}_{\uptau})}
\ar[d]^-(0.5){\hupsigma'^*(\hu)\otimes \bvuppsi^{(r,r')}}\\
{\bvcC'^{(r)}\otimes_{\ocB'}\hupsigma'^*(\hcC^{(r')}_{\uptau})}\ar[rr]^-(0.5){-\bvdelta^{(r,r')}}&&
{\hupsigma'^*(\hOmega')\otimes_{\bvocB'}\bvcC'^{(r)}\otimes_{\ocB'}\hupsigma'^*(\hcC^{(r')}_{\uptau}),}}
\end{equation}
where $\hu$ is the canonical morphism \eqref{p2-fhtft2a},
$\bvuppsi^{(r,r')}$ is the homomorphism \eqref{p2-fhtft380g} and $\delta_{\hcC^{(r')}_{\uptau}}$ is the derivation \eqref{p2-fhtft6b}. In particular, we have 
\begin{equation}\label{p2-rftchta170l}
\bvdelta^{(r,r')}\circ \bvupphi^{(r,r')}=0,
\end{equation}
where $\bvupphi^{(r,r')}$ is the homomorphism defined in \eqref{p2-fhtft380f}. 
By \eqref{p2-rftchta150i}, the diagram 
\begin{equation}\label{p2-rftchta170m}
\xymatrix{
{\bvuptheta^*(\bvcC^{(r)})\otimes_{\bvocB'}\hupsigma'^*(\hcC^{(r')}_{\uptau})}\ar[rr]^-(0.5){b}
\ar[d]_{\bvuppsi^{(r,r')}}&&
{\hupsigma'^*(\fgg^*(\hOmega))\otimes_{\bvocB'}\bvuptheta^*(\bvcC^{(r)})\otimes_{\bvocB'}\hupsigma'^*(\hcC^{(r')}_{\uptau})}
\ar[d]^-(0.5){\id\otimes \bvuppsi^{(r,r')}}\\
{\bvcC'^{(r)}\otimes_{\ocB'}\hupsigma'^*(\hcC^{(r')}_{\uptau})}\ar[rr]^-(0.5){a}&&
{\hupsigma'^*(\fgg^*(\hOmega))\otimes_{\bvocB'}\bvcC'^{(r)}\otimes_{\ocB'}\hupsigma'^*(\hcC^{(r')}_{\uptau}),}}
\end{equation}
where $a$ denotes the Higgs field 
\begin{equation}
a=\id\otimes \hupsigma'^*(\delta_{\hcC^{(r')}_{\uptau}})
\end{equation}
on $\bvcC'^{(r)}\otimes_{\ocB'}\hupsigma'^*(\hcC^{(r')}_{\uptau})$, $b$ denotes the Higgs field 
\begin{equation}
b=\bvuptheta^*(\delta_{\bvcC^{(r)}})\otimes\id+\id\otimes \hupsigma'^*(\delta_{\hcC^{(r')}_{\uptau}})
\end{equation}
on $\bvuptheta^*(\bvcC^{(r)})\otimes_{\bvocB'}\hupsigma'^*(\hcC^{(r')}_{\uptau})$, and $\delta_{\bvcC^{(r)}}$ is the derivation \eqref{p2-htaft19b},
is commutative.

\section{\texorpdfstring{Functoriality of the $p$-adic Simpson correspondence by pullback}
{Functoriality of the p-adic Simpson correspondence by pullback}}\label{p2-fgpscp}

\subsection{}\label{p2-fgpscp1}
The assumptions and notation of §\ref{p2-fhtft} remain in force throughout this section.
We consider again the objects associated with $(f,\tf)$ introduced in §\ref{p2-rgpsc}, and we associate with $(f',\tf')$ \eqref{p2-fhtft1d}
similar objects that we denote by the same symbols equipped with a $^\prime$ exponent.
In particular, we denote by $\bIndMod^\Dolb(\bvocB')$ the category of Dolbeault ind-$\bvocB'$-modules
and by $\bHM^\sol(\co_{\fX'}[\frac 1 p], \hOmega')$
the category of solvable Higgs $\co_{\fX'}[\frac 1 p]$-bundles with coefficients in $\hOmega'$, see \ref{p2-rgpsc12} and \ref{p2-fhtft2},
with respect  to the deformation $\tf'$ fixed in \eqref{p2-fhtft1d}. We denote by
\begin{equation}\label{p2-fgpscp10a}
\cH'\colon \bIndMod(\bvocB')\rightarrow \bHM(\co_{\fX'}, \hOmega')
\end{equation}
the functor defined in \eqref{p2-rgpsc14a}, associated with $(f',\tf')$.
By \ref{p2-rgpsc15}, this induces an equivalence of categories that we denote again by
\begin{equation}\label{p2-fgpscp10b}
\cH'\colon \bIndMod^\Dolb(\bvocB')\stackrel{\sim}{\rightarrow} \bHM^\sol(\co_{\fX'}[\frac 1 p], \hOmega').
\end{equation}

To lighten the notation, for any rational numbers $0\leq t\leq r$, we denote by
\begin{eqnarray}
\delta^{(r)}=\bvuptheta^*(\delta_{\bvcC^{(r)}})\colon \bvuptheta^*(\bvcC^{(r)})\rightarrow 
\hupsigma'^*(\fgg^*(\hOmega))\otimes_{\bvocB'}\bvuptheta^*(\bvcC^{(r)}),\label{p2-fgpscp1b}\\
\delta'^{(r)}=\delta_{\bvcC'^{(r)}}\colon \bvcC'^{(r)}\rightarrow \hupsigma'^*(\hOmega')\otimes_{\bvocB'}\bvcC'^{(r)},\label{p2-fgpscp1c}\\
\delta^{(t)}_\uptau=\delta_{\hcC^{(t)}_\uptau}\colon \hcC^{(t)}_\uptau\rightarrow \fgg^*(\hOmega)\otimes_{\co_{\fX'}} \hcC^{(t)}_\uptau,\label{p2-fgpscp1d}\\
\delta'^{(t)}_\uptau=\delta'_{\hcC^{(t)}_\uptau}\colon \hcC^{(t)}_\uptau\rightarrow \hOmega'\otimes_{\co_{\fX'}} \hcC^{(t)}_\uptau,\label{p2-fgpscp1e}
\end{eqnarray}
where $\delta_{\bvcC^{(r)}}$ and $\delta_{\bvcC'^{(r)}}$ are the derivations defined in \eqref{p2-htaft19b}, 
$\delta_{\hcC^{(t)}_\uptau}$ in \eqref{p2-fhtft6b} and $\delta'_{\hcC^{(t)}_\uptau}$ in \eqref{p2-fhtft6bb}, and we set 
\begin{eqnarray}
\delta^{\vee(t)}_\uptau=-\delta^{(t)}_\uptau\colon \hcC^{(t)}_\uptau\rightarrow \fgg^*(\hOmega)\otimes_{\co_{\fX'}} \hcC^{(t)}_\uptau,\label{p2-fgpscp1f}\\
\delta'^{\vee(t)}_\uptau=-\delta'^{(t)}_\uptau\colon \hcC^{(t)}_\uptau\rightarrow \hOmega'\otimes_{R_\uptau} \hcC^{(t)}_\uptau.\label{p2-fgpscp1g}
\end{eqnarray}
In the notation above, a prime exponent indicates composition with $\hu$ \eqref{p2-fhtft2a}, while a $\vee$ exponent denotes multiplication by $-1$.
We consider also the $\bvocB'$-derivations 
\begin{eqnarray}
\lefteqn{\delta'^{(r,t)}=\delta'^{(r)}\otimes \id+\id \otimes \hupsigma'^*(\delta'^{\vee(t)}_\uptau) \colon}\label{p2-fgpscp1h}\\
&&\bvcC'^{(r)}\otimes_{\bvocB'}\hupsigma'^*(\hcC^{(t)}_\uptau) 
\rightarrow \hupsigma'^*(\hOmega')\otimes_{\bvocB'}\bvcC'^{(r)}\otimes_{\bvocB'}\hupsigma'^*(\hcC^{(t)}_\uptau),\nonumber
\end{eqnarray}
\begin{eqnarray}
\lefteqn{\delta'^{\vee(r,t)}=-\delta'^{(r,t)}=-\delta'^{(r)}\otimes \id+\id \otimes \hupsigma'^*(\delta'^{(t)}_\uptau)}\label{p2-fgpscp1hv}\\
&&\bvcC'^{(r)}\otimes_{\bvocB'}\hupsigma'^*(\hcC^{(t)}_\uptau) 
\rightarrow \hupsigma'^*(\hOmega')\otimes_{\bvocB'}\bvcC'^{(r)}\otimes_{\bvocB'}\hupsigma'^*(\hcC^{(t)}_\uptau),\nonumber
\end{eqnarray}
\begin{equation}\label{p2-fgpscp1i}
\delta^{(r,t)}=\id \otimes \hupsigma'^*(\delta^{(t)}_\uptau) \colon \bvcC'^{(r)}\otimes_{\bvocB'}\hupsigma'^*(\hcC^{(t)}_\uptau) 
\rightarrow \hupsigma'^*(\fgg^*(\hOmega))\otimes_{\bvocB'}\bvcC'^{(r)}\otimes_{\bvocB'}\hupsigma'^*(\hcC^{(t)}_\uptau).
\end{equation}

The following diagrams 
\begin{equation}\label{p2-fgpscp1j}
\xymatrix{
{\hupsigma'^*(\hcC^{(t)}_\uptau)}\ar[r]^-(0.5){\delta^{(t)}_\uptau}\ar[d]&
{\hupsigma'^*(\fgg^*(\hOmega))\otimes_{\bvocB'}\hupsigma'^*(\hcC^{(t)}_\uptau)}\ar[d]\\
{\bvcC'^{(r)}\otimes_{\bvocB'}\hupsigma'^*(\hcC^{(t)}_\uptau)}\ar[r]^-(0.5){\delta^{(r,t)}}&
{\hupsigma'^*(\fgg^*(\hOmega))\otimes_{\bvocB'}\bvcC'^{(r)}\otimes_{\bvocB'}\hupsigma'^*(\hcC^{(t)}_\uptau),}}
\end{equation}
\begin{equation}\label{p2-fgpscp1k}
\xymatrix{
{\hupsigma'^*(\hcC^{(t)}_\uptau)}\ar[r]^-(0.5){\delta'^{\vee(t)}_\uptau}\ar[d]&
{\hupsigma'^*(\hOmega')\otimes_{\bvocB'}\hupsigma'^*(\hcC^{(t)}_\uptau)}\ar[d]\\
{\bvcC'^{(r)}\otimes_{\bvocB'}\hupsigma'^*(\hcC^{(t)}_\uptau)}\ar[r]^-(0.5){\delta'^{(r,t)}}&
{\hupsigma'^*(\hOmega')\otimes_{\bvocB'}\bvcC'^{(r)}\otimes_{\bvocB'}\hupsigma'^*(\hcC^{(t)}_\uptau),}}
\end{equation}
\begin{equation}\label{p2-fgpscp1l}
\xymatrix{
{\bvcC'^{(r)}}\ar[r]^-(0.5){\delta'^{(r)}}\ar[d]&
{\hupsigma'^*(\hOmega')\otimes_{\bvocB'}\bvcC'^{(r)}}\ar[d]\\
{\bvcC'^{(r)}\otimes_{\bvocB'}\hupsigma'^*(\hcC^{(t)}_\uptau)}\ar[r]^-(0.5){\delta'^{(r,t)}}&
{\hupsigma'^*(\hOmega')\otimes_{\bvocB'}\bvcC'^{(r)}\otimes_{\bvocB'}\hupsigma'^*(\hcC^{(t)}_\uptau),}}
\end{equation}
where the vertical arrows are the canonical morphisms, are commutative. 
By \eqref{p2-rftchta170m}, the following diagram 
\begin{equation}\label{p2-fgpscp1n}
\xymatrix{
{\bvuptheta^*(\bvcC^{(r)})}\ar[r]^-(0.5){\delta^{(r)}}\ar[d]_{\bvupphi^{(r,t)}}&
{\hupsigma'^*(\fgg^*(\hOmega))\otimes_{\bvocB'}\bvuptheta^*(\bvcC^{(r)})}\ar[d]^{\id \otimes \bvupphi^{(r,t)}}\\
{\bvcC'^{(r)}\otimes_{\bvocB'}\hupsigma'^*(\hcC^{(t)}_\uptau)}\ar[r]^-(0.5){\delta^{(r,t)}}&
{\hupsigma'^*(\fgg^*(\hOmega))\otimes_{\bvocB'}\bvcC'^{(r)}\otimes_{\bvocB'}\hupsigma'^*(\hcC^{(t)}_\uptau)}}
\end{equation}
is commutative. Moreover, the composition 
\begin{equation}\label{p2-fgpscp1o}
\xymatrix{
{\bvuptheta^*(\bvcC^{(r)})}\ar[r]^-(0.5){\bvupphi^{(r,t)}}&{\bvcC'^{(r)}\otimes_{\bvocB'}\hupsigma'^*(\hcC^{(t)}_\uptau)}
\ar[r]^-(0.5){\delta'^{(r,t)}}&{\hupsigma'^*(\hOmega')\otimes_{\bvocB'}\bvcC'^{(r)}\otimes_{\bvocB'}\hupsigma'^*(\hcC^{(t)}_\uptau)}}
\end{equation}
vanishes \eqref{p2-rftchta170l}.

\begin{teo}\label{p2-fgpscp2}
For every Dolbeault ind-$\bvocB$-module $\cM$ \eqref{p2-rgpsc12} such that the associated 
$\co_{\fX}[\frac 1 p]$-bundle $\cH(\cM)$ \eqref{p2-rgpsc15a} is locally CL-small \eqref{p1-tshbn13}, 
the ind-$\bvocB'$-module $\rI\bvuptheta^*(\cM)$ is Dolbeault, its associated Higgs
$\co_{\fX'}[\frac 1 p]$-bundle $\cH'(\rI\bvuptheta^*(\cM))$ \eqref{p2-fgpscp10b} is locally CL-small, 
and we have a canonical functorial isomorphism 
\begin{equation}\label{p2-fgpscp2a}
\cH'(\rI\bvuptheta^*(\cM))\stackrel{\sim}{\rightarrow}\upmu(\fgg^*_\uptau(\cH(\cM))),
\end{equation}
where $\fgg^*_\uptau$ is the twisted pullback functor by $\fgg$ defined in \eqref{p2-fhtft4b}, and 
\begin{equation}
\upmu\colon \bHM(\co_{\fX'},\fgg^*(\hOmega)) \rightarrow \bHM(\co_{\fX'},\hOmega')
\end{equation} 
is the functor induced by the canonical morphism $\hu\colon \fgg^*(\hOmega)\rightarrow \hOmega'$ \eqref{p2-fhtft2a}. 
\end{teo}

We set $(N,\theta)=\cH(\cM)$ and $(\cN,\vartheta)=\upalpha^\coh_{\co_{\fX}}(N,\theta)$ \eqref{p2-cmupiso31h}. 
There exist a rational number $r>0$ and an isomorphism of ind-$\bvcC^{(r)}$-modules with 
$\delta_{\bvcC^{(r)}}$-connection in the sense of \ref{p1-indmal22}, 
\begin{equation}\label{p2-fgpscp2b}
\cM\otimes_{\bvocB}\bvcC^{(r)} \stackrel{\sim}{\rightarrow} \rI\hupsigma^*(\cN)\otimes_{\bvocB}\bvcC^{(r)},
\end{equation}
where the $\delta_{\bvcC^{(r)}}$-connections are defined as in \ref{p1-delta-con9}, 
$\rI\hupsigma^*(\cN)$ (resp.\ $\cM$) being endowed with the Higgs field 
$\rI\hupsigma^*(\vartheta)$ (resp.\ $0$).  

Let $t$ be a rational number rational numbers such that $0<t\leq r$. 
By \ref{p1-delta-con5} and \eqref{p2-fgpscp1n}, the isomorphism \eqref{p2-fgpscp2b} 
induces by pullback by $\rI\bvuptheta$ and extension of scalars by the homomorphism $\bvupphi^{(r,t)}$ \eqref{p2-fhtft380f} an isomorphism of 
$\bvcC'^{(r)}\otimes_{\bvocB'}\hupsigma'^*(\hcC^{(t)}_\uptau)$-modules with $\delta^{(r,t)}$-connection 
\begin{equation}\label{p2-fgpscp2i1}
\rI\bvuptheta^*(\cM)\otimes_{\bvocB'}\bvcC'^{(r)}\otimes_{\bvocB'}\hupsigma'^*(\hcC^{(t)}_\uptau)
\stackrel{\sim}{\rightarrow} \rI\hupsigma'^*(\rI\fgg^*(\cN))\otimes_{\bvocB'}\bvcC'^{(r)}\otimes_{\bvocB'}\hupsigma'^*(\hcC^{(t)}_\uptau),
\end{equation}
where the $\delta^{(r,t)}$-connections are defined as in \ref{p1-delta-con4}, 
$\rI\hupsigma'^*(\rI\fgg^*(\cN))$ (resp.\ $\cM$) being endowed with the Higgs field $\rI\hupsigma'^*(\rI\fgg^*(\vartheta))$ (resp.\ zero).
Moreover, it follows from \eqref{p2-fgpscp1o} that \eqref{p2-fgpscp2i1} is an isomorphism of 
$\bvcC'^{(r)}\otimes_{\bvocB'}\hupsigma'^*(\hcC^{(t)}_\uptau)$-modules
with $\delta'^{(r,t)}$-connection, when  $\rI\hupsigma'^*(\rI\fgg^*(\cN))$ and $\cM$ are endowed with the zero Higgs fields. 
To summarize, \eqref{p2-fgpscp2i1} is an isomorphism
of Higgs $\bvocB'$-modules with the Higgs fields indicated on the same line of this table
\begin{equation}\label{p2-fgpscp2i2}
\begin{tabular}{|c|c|}
\hline
left hand side&right hand  side\\
\hline
$\id\otimes \delta^{(r,t)}$ & $\rI\hupsigma'^*(\rI\fgg^*(\vartheta))\otimes \id +\id\otimes \delta^{(r,t)}$ \\
\hline
$\id\otimes \delta'^{(r,t)}$ & $\id\otimes \delta'^{(r,t)}$ \\
\hline
\end{tabular}
\end{equation}

On the other hand, by \ref{p1-tshbn16}, the Higgs $\co_\fX[\frac 1 p]$-bundle $(N,\theta)$ being locally CL-small,
is weakly twistable by the extension $\hcF_\uptau$ \eqref{p2-fhtft3d}, in the sense of \ref{p2-fhtft5}.
Hence, by \ref{p1-thbn21}, it is twistable by the same extension.  
We set $\fgg^*_\uptau(N,\theta)=(N^\vee,\theta^\vee)$ \eqref{p2-fhtft4b}. 
Then, there exists a rational number $t'>0$ such that the canonical morphism $N^\vee \rightarrow \fgg^*(N)\otimes_{\co_{\fX'}}\cC^\dagger_\uptau$ 
induces an isomorphism of $\hcC^{(t')}_\uptau$-modules with $\delta^{(t')}_\uptau$-connection \eqref{p2-fgpscp1d}
\begin{equation}\label{p2-fgpscp2j}
N^\vee \otimes_{\co_{\fX'}}\hcC^{(t')}_\uptau\stackrel{\sim}{\rightarrow} \fgg^*(N)\otimes_{\co_{\fX'}}\hcC^{(t')}_\uptau,
\end{equation}
where the $\delta^{(t')}_{\uptau}$-connections are defined as in \ref{p1-delta-con4}, 
$N$ (resp.\ $N^\vee$) being endowed with the Higgs field $\theta$ (resp.\ $0$). Moreover, by \ref{p1-thbn25}, 
\eqref{p2-fgpscp2j} is also an isomorphism of $\hcC^{(t')}_\uptau$-modules with $\delta^{\vee(t')}_\uptau$-connection \eqref{p2-fgpscp1f}, 
where the $\delta^{\vee(t')}_{\uptau}$-connections are defined as in \ref{p1-delta-con4}, 
$N$ (resp.\ $N^\vee$) being endowed with the Higgs field $0$ (resp.\ $\theta^\vee$).
Observe that the same properties hold for every rational number $t''$ such that $0\leq t''\leq t'$.

We set $(\cN^\vee,\vartheta^\vee)=\upalpha^\coh_{\co_{\fX'}}(N^\vee,\theta^\vee)$ \eqref{p2-cmupiso31h}.
Applying the functor $\upalpha^\coh_{\delta^{(t')}_\uptau}$ \eqref{p2-cmupiso310g} to the isomorphism \eqref{p2-fgpscp2j},
we obtain an isomorphism of ind-$\hcC^{(t')}_\uptau$-modules 
with $\delta^{(t')}_\uptau$-connection \eqref{p2-fgpscp1d}
\begin{equation}\label{p2-fgpscp2j1}
\cN^\vee \otimes_{\co_{\fX'}}\hcC^{(t')}_\uptau\stackrel{\sim}{\rightarrow} \rI\fgg^*(\cN)\otimes_{\co_{\fX'}}\hcC^{(t')}_\uptau,
\end{equation}
where the $\delta^{(t')}_{\uptau}$-connections are defined as in \ref{p1-delta-con4}, 
$\cN$ (resp.\ $\cN^\vee$) being endowed with the Higgs field $\vartheta$ (resp.\ $0$). Moreover, 
\eqref{p2-fgpscp2j1} is also an isomorphism of ind-$\hcC^{(t')}_\uptau$-modules with $\delta^{\vee(t')}_\uptau$-connection \eqref{p2-fgpscp1f}, 
where the $\delta^{\vee(t')}_{\uptau}$-connections are defined as in \ref{p1-delta-con4}, 
$\cN$ (resp.\ $\cN^\vee$) being endowed with the Higgs field $0$ (resp.\ $\vartheta^\vee$).

Let $t$ be a rational number such that $0<t\leq \inf(t',r)$. 
By \ref{p1-delta-con5} and \eqref{p2-fgpscp1j}, the isomorphism \eqref{p2-fgpscp2j1} induces, by applying the functor $\rI\hupsigma'^*$
and extension of scalars, an isomorphism of 
$\bvcC'^{(r)}\otimes_{\bvocB'}\hupsigma'^*(\hcC^{(t)}_\uptau)$-modules with $\delta^{(r,t)}$-connection 
\begin{equation}\label{p2-fgpscp2k1}
\rI\hupsigma'^*(\cN^\vee)\otimes_{\bvocB'}\bvcC'^{(r)}\otimes_{\bvocB'}\hupsigma'^*(\hcC^{(t)}_\uptau)\stackrel{\sim}{\rightarrow} 
\rI \hupsigma'^*(\rI \fgg^*(\cN))\otimes_{\bvocB'}\bvcC'^{(r)}\otimes_{\bvocB'}\hupsigma'^*(\hcC^{(t)}_\uptau),
\end{equation}
where the $\delta^{(r,t)}$-connections are defined as in \ref{p1-delta-con4}, 
$\cN$ (resp.\ $\cN^\vee$) being endowed with the Higgs field $\vartheta$ (resp.\ $0$). 
By \eqref{p2-fgpscp1k}, it is also an isomorphism of $\bvcC'^{(r)}\otimes_{\bvocB'}\hupsigma'^*(\hcC^{(t)}_\uptau)$-modules with $\delta'^{(r,t)}$-connection,  
where the $\delta'^{(r,t)}$-connections are defined as in \ref{p1-delta-con4}, 
$\cN$ (resp.\ $\cN^\vee$) being endowed with the Higgs field $0$ 
(resp.\ $\vartheta'^\vee=(\hu\otimes\id)\circ\vartheta^\vee\colon \cN^\vee\rightarrow \hOmega'\otimes_{\co_{\fX'}}\cN^\vee$). 
To summarize, \eqref{p2-fgpscp2k1} is an isomorphism of Higgs $\bvocB'$-modules with the Higgs fields indicated on the same line of this table
\begin{equation}\label{p2-fgpscp2k2}
\begin{tabular}{|c|c|}
\hline
left hand side&right  hand  side\\
\hline
$\id\otimes \delta^{(r,t)}$ & $\rI \hupsigma'^*(\rI \fgg^*(\vartheta))\otimes \id+\id\otimes\delta^{(r,t)}$\\
\hline
$\rI \hupsigma'^*(\vartheta'^\vee)\otimes \id+\id\otimes \delta'^{(r,t)}$ & $\id\otimes \delta'^{(r,t)}$ \\
\hline
\end{tabular}
\end{equation}

We set $(\cN',\vartheta')=(\cN^\vee,\vartheta'^\vee)$ and $\cM'=\rI\bvuptheta^*(\cM)$.
Composing \eqref{p2-fgpscp2i1} and the inverse of \eqref{p2-fgpscp2k1}, 
we deduce an isomorphism of $\bvcC'^{(r)}\otimes_{\bvocB'}\hupsigma'^*(\hcC^{(t)}_\uptau)$-modules 
\begin{equation}\label{p2-fgpscp2l1}
\cM'\otimes_{\bvocB'}\bvcC'^{(r)}\otimes_{\bvocB'}\hupsigma'^*(\hcC^{(t)}_\uptau)
\stackrel{\sim}{\rightarrow} 
\rI\hupsigma'^*(\cN^\vee)\otimes_{\bvocB'}\bvcC'^{(r)}\otimes_{\bvocB'}\hupsigma'^*(\hcC^{(t)}_\uptau),
\end{equation}
which is compatible with the Higgs $\bvocB'$-fields indicated on the same line of the table
\begin{equation}\label{p2-fgpscp2l2}
\begin{tabular}{|c|c|}
\hline
left hand side&  right  hand  side\\
\hline
$\id\otimes \delta^{(r,t)}$ & $\id\otimes \delta^{(r,t)}$ \\
\hline
$\id\otimes \delta'^{(r,t)}$ & $\rI\hupsigma'^*(\vartheta')\otimes \id+\id\otimes \delta'^{(r,t)}$\\
\hline 
\end{tabular}
\end{equation}
These isomorphisms are compatible when $t$ varies, 
we deduce an isomorphism of $\bvcC'^{(r)}\otimes_{\bvocB'}\rI\hupsigma'^*(\IC^\dagger_\uptau)$-modules 
\begin{equation}\label{p2-fgpscp2l3}
\cM'\otimes_{\bvocB'}\bvcC'^{(r)}\otimes_{\bvocB'}\rI\hupsigma'^*(\IC^\dagger_\uptau)
\stackrel{\sim}{\rightarrow} 
\rI\hupsigma'^*(\cN')\otimes_{\bvocB'}\bvcC'^{(r)}\otimes_{\bvocB'}\rI\hupsigma'^*(\IC^\dagger_\uptau),
\end{equation}
where $\IC^\dagger_\uptau$ is the ind-algebra \eqref{p2-fhtft6k}, 
which is compatible with the Higgs $\bvocB'$-fields indicated on the same line of the table
\begin{equation}\label{p2-fgpscp2l4}
\begin{tabular}{|c|c|}
\hline
left hand side&  right  hand  side\\
\hline
$\id\otimes \kappa^{(r)}$ & $\id\otimes \kappa^{(r)}$ \\
\hline
$\id\otimes \kappa'^{(r)}$ & $\rI\hupsigma'^*(\vartheta')\otimes \id+\id\otimes \kappa'^{(r)}$\\
\hline 
\end{tabular}
\end{equation}
where
\begin{eqnarray}
\kappa^{(r)}&=&\id\otimes \hupsigma'^*(\Idelta_\uptau),\label{p2-fgpscp2l5}\\
\kappa'^{(r)}&=&\delta'^{(r)} \otimes \id-  \id\otimes \hupsigma'^*(\Idelta'_\uptau),\label{p2-fgpscp2l6}
\end{eqnarray}
and $\Idelta_\uptau$ (resp.\ $\Idelta'_\uptau$) is the derivation \eqref{p2-fhtft6i} (resp.\ \eqref{p2-fhtft6j}).

By \eqref{p1-bcim5d}, we can replace in \eqref{p2-fgpscp2l3}, 
$\IC^\dagger_\uptau$ by $\IC^\dagger_{\uptau, \mQ}$ \eqref{p2-cmupiso32g} and $\Idelta_\uptau$ by $\Idelta_{\uptau,\mQ}$ \eqref{p2-cmupiso32h}. 
Then, by \ref{p2-cmupiso33}(ii), taking the kernels of the Higgs fields $\id\otimes \kappa^{(r)}$, we obtain an isomorphism of $\bvcC'^{(r)}$-modules 
\begin{equation}
\cM'\otimes_{\bvocB'}\bvcC'^{(r)}\stackrel{\sim}{\rightarrow} \rI\hupsigma'^*(\cN')\otimes_{\bvocB'}\bvcC'^{(r)}.
\end{equation}
Moreover, by restricting the Higgs $\bvocB'$-fields of the bottom line of \eqref{p2-fgpscp2l4}, 
this is an isomorphism of ind-$\bvcC'^{(r)}$-modules with $\delta'^{(r)}$-connection, 
where the $\delta'^{(r)}$-connections are defined as in \ref{p1-delta-con4}, 
$\cN'$ (resp.\ $\cM'$) being endowed with the Higgs field $\vartheta'$ (resp.\ $0$); see \eqref{p2-fgpscp1l}.

Therefore, the $\bvocB'$-module $\cM'$ is Dolbeault, and we have an isomorphism 
\begin{equation}
\cH'(\cM')\stackrel{\sim}{\rightarrow}(N^\vee,(\hu\otimes\id)\circ\theta^\vee).
\end{equation}
The required canonical and functorial isomorphism \eqref{p2-fgpscp2a} follows. 
The Higgs $\co_{\fX'}[\frac 1 p]$-module  $(N^\vee,\theta^\vee)$ is locally CL-small by \ref{p2-cmupiso211}, 
and hence so is $\cH'(\cM')$.

\section[Refined functoriality of Higgs--Tate algebras in Faltings topos]
{Refined functoriality of Higgs--Tate algebras in Faltings topos in the smooth case}\label{p2-rfhtft}

\subsection{}\label{p2-fhtft9}
The assumptions and notation of §\ref{p2-fhtft} remain in force throughout this section.
We assume, furthermore, 
that the morphism $g\colon (X',\cM_{X'})\rightarrow (X,\cM_X)$ \eqref{p2-fhtft1a} is {\em smooth and saturated}. We set 
\begin{equation}\label{p2-fhtft9a}
\tOmega^1_{X'/X}=\tOmega^1_{(X',\cM_{X'})/(X,\cM_X)},
\end{equation}
that we consider as a sheaf of $X'_\zar$ or $X'_\et$, depending on the context \eqref{p2-ncgt6}.
We have a canonical exact sequence of locally free $\co_{X'}$-modules of finite type
\begin{equation}\label{p2-fhtft9b}
0\rightarrow g^*(\tOmega^1_{X/S})\rightarrow \tOmega^1_{X'/S}\rightarrow \tOmega^1_{X'/X}\rightarrow 0.
\end{equation}
To lighten the notation, we set, with the convention of  \ref{p2-ncgt3}, 
\begin{eqnarray}
\tOmega^1_{\coX'/\coX}=\tOmega^1_{(\coX',\cM_{\coX'})/(\coX,\cM_\coX)},&&\uOmega'=\txi^{-1}\tOmega^1_{\coX'/\coX},\label{p2-fhtft9e}
\end{eqnarray}
so we have a canonical exact sequence of locally free $\co_{\coX'}$-modules of finite type
\begin{equation}\label{p2-fhtft9l}
0\longrightarrow \cog^*(\Omega)\stackrel{u}{\longrightarrow} \Omega' \longrightarrow \uOmega' \longrightarrow 0. 
\end{equation}
For any rational number $r\geq 0$, we set
\begin{equation}\label{p2-fhtft90g}
\uOmega^{(r)}=p^r\uOmega.
\end{equation}
The exact sequence \eqref{p2-fhtft9l} induces an exact sequence of locally free $\co_{\coX'}$-modules of finite type
\begin{equation}\label{p2-fhtft9h}
0\rightarrow \cog^*(\Omega^{(r)})\rightarrow \Omega'^{(r)} \rightarrow \uOmega'^{(r)} \rightarrow 0. 
\end{equation}

Let $r,r'$ be rational numbers such that $r\geq r'\geq 0$. We denote by 
\begin{equation}\label{p2-fhtft9i}
\pi^{(r,r')}\colon \Omega^{(r)}\rightarrow \Omega^{(r')}, \ \ \ \pi'^{(r,r')}\colon \Omega'^{(r)}\rightarrow \Omega'^{(r')}  \ \ \ {\rm and} \ \ \ \upi'^{(r,r')}\colon \uOmega'^{(r)}\rightarrow \uOmega'^{(r')}
\end{equation} 
the canonical injections and by $\Omega'^{(r,r')}$ the inverse image of $\upi'^{(r,r')}(\uOmega'^{(r)})$ in $\Omega'^{(r')}$, so that we have a commutative diagram 
of canonical $\co_{\coX'}$-linear morphisms
\begin{equation}\label{p2-fhtft9j}
\xymatrix{
0\ar[r]&{\cog^*(\Omega^{(r')})}\ar[r]&{\Omega'^{(r')}}\ar[r]&{\uOmega'^{(r')}}\ar[r]&0\\
0\ar[r]&{\cog^*(\Omega^{(r')})}\ar[r]\ar@{=}[u]&{\Omega'^{(r,r')}}\ar[r]\ar@{^(->}[u]\ar@{}[ru]|{(2)}&{\uOmega'^{(r)}}\ar[r]\ar[u]_{\upi'^{(r,r')}}&0\\
0\ar[r]&{\cog^*(\Omega^{(r)})}\ar[r]\ar[u]^{\cog^*(\pi^{(r,r')})}\ar@{}[ru]|{(1)}&{\Omega'^{(r)}}\ar[r]\ar@{^(->}[u]&{\uOmega'^{(r)}}\ar[r]\ar@{=}[u]&0,}
\end{equation}
where the horizontal lines are exact, the top line and bottom line being \eqref{p2-fhtft9h}, the square $(1)$ is co-Cartesian, the square $(2)$ is Cartesian, 
and the composition of the middle vertical arrows is $\pi'^{(r,r')}$, see \ref{p2-hta14}. 

Let $t,t'$ be rational numbers such that $t\geq t'\geq 0$, $r\geq t$ and $r'\geq t'$. We have a commutative diagram of canonical $\co_{\coX'}$-linear morphisms
\begin{equation}\label{p2-fhtft9k}
\xymatrix{x
0\ar[r]&{\cog^*(\Omega^{(r')})}\ar[r]\ar[d]_{\cog^*(\pi^{(r',t')})}&{\Omega'^{(r,r')}}\ar[r]\ar[d]^{\pi'^{(r,r',t,t')}}&{\uOmega'^{(r)}}\ar[r]\ar[d]^{\upi'^{(r,t)}}&0\\
0\ar[r]&{\cog^*(\Omega^{(t')})}\ar[r]&{\Omega'^{(t,t')}}\ar[r]&{\uOmega'^{(t)}}\ar[r]&0.}
\end{equation}

\subsection{}\label{p2-fhtft200}
We denote by $\tOmega^1_{\fX'/\fX}$ the $p$-adic completion of $\tOmega^1_{\coX'/\coX}$ 
and by $\huOmega'$ the $p$-adic completion of  $\uOmega'$ \eqref{p2-fhtft9e}, 
which is canonically isomorphic to the module $\txi^{-1}\tOmega^1_{\fX'/\fX}$ \eqref{p2-ncgt3}; 
{\em we identify these modules in the following}. The exact sequence \eqref{p2-fhtft9l} 
induces an exact sequence of locally free $\co_{\fX'}$-modules of finite type
\begin{equation}\label{p2-fhtft200a}
0\longrightarrow \fgg^*(\hOmega)\stackrel{\hu}{\longrightarrow} \hOmega'\longrightarrow \huOmega'\longrightarrow 0.
\end{equation}
For any rational numbers $r\geq r'\geq 0$, we denote by $\huOmega'^{(r)}$ (resp.\ $\hOmega'^{(r,r')}$)
the $p$-adic completion of $\uOmega'^{(r)}$ (resp.\ $\Omega'^{(r,r')}$). 
Observe that $\huOmega'^{(r)}$ is $\cS$-flat and identifies canonically with  $p^r\huOmega'$.

\subsection{}\label{p2-fhtft24}
Recall that we use for $(f',\tf')$ the notation introduced in \ref{p2-htaft9}--\ref{p2-htaft17} equipped with a~$^\prime$ exponent. 
Let $r,r'$ be rational numbers such that $r\geq r'\geq 0$, $n$ an integer $\geq 0$.
Let $Y'$ be an object of $\bP'$ \eqref{p2-htaft9} 
such that $Y'_s$ is nonempty, $((P',\gamma'),(\mN,\nu'),\vartheta')$ an adequate chart for the morphism 
$f'|Y'\colon (Y',\cM_{X'}|Y')\rightarrow (S,\cM_S)$ induced by $f'$ \eqref{p2-fhtft1}, $\oy'$ a geometric point of $\oY'^\rhd$. 
Recall \eqref{p2-htaft12} that we associate with $(Y',\oy')$ the logarithmic schemes 
$(\hmY'^{\oy'},\cM_{\hmY'^{\oy'}})$ and $(\tmY'^{\oy'},\cM_{\tmY'^{\oy'}})$ \eqref{p2-htaft12g}, 
the torsor $\cL^{(r)}_{\tmY'^{\oy'}/\tX'}$, the Higgs--Tate extension of thickness $r$ \eqref{p2-htaft12e}
\begin{equation}\label{p2-fhtft24a}
0\rightarrow \hoR^{\oy'}_{Y'}\rightarrow \fF'^{\oy',(r)}_{Y'}\rightarrow \Omega'^{(r)}(\coY') \otimes_{\co_{\coX'}(\coY')}\hoR^{\oy'}_{Y'} \rightarrow 0,
\end{equation}
and the Higgs--Tate algebra $\fC'^{\oy',(r)}_{Y'}$ \eqref{p2-htaft12d}. 
Note that we have changed the normalization in \eqref{p2-htaft12e}, to align with that of \ref{p2-hta7}, namely we identified $\Omega'^{(r)}$ with $\Omega'$ \eqref{p2-fhtft9i}, 
so that the canonical morphism $\pi'^{(r,0)}\colon \Omega'^{(r)}\rightarrow \Omega'$ is identified with the multiplication by $p^r$ on $\Omega'$, see \ref{p2-hta71}.

Let $\cL^{(r,r')}_{\tmY'^{\oy'}/\tX'}$ be the twist of the torsor $\cL_{\tmY'^{\oy'}/\tX'}=\cL^{(0)}_{\tmY'^{\oy'}/\tX'}$ 
defined in \ref{p2-hta14}, $\fF'^{\oy',(r,r')}_{Y'}$ the associated $\hoR^{\oy'}_{Y'}$-extension \eqref{p2-fhtal8b}
\begin{equation}\label{p2-fhtft24b}
0\rightarrow \hoR^{\oy'}_{Y'}\rightarrow \fF'^{\oy',(r,r')}_{Y'}\rightarrow \Omega'^{(r,r')}(\coY') \otimes_{\co_{\coX'}(\coY')}\hoR^{\oy'}_{Y'} \rightarrow 0, 
\end{equation}
$\fC'^{\oy',(r,r')}_{Y'}$ the associated $\hoR^{\oy'}_{Y'}$-algebra \eqref{p2-fhtal8c}.
As in \ref{p2-htaft12}, there is a canonical $\hoR^{\oy'}_{Y'}$-semilinear action of $\pi_1(\oY'^{\star\rhd},\oy')$ on $\fF'^{\oy',(r,r')}_{Y'}$ 
that is continuous for the $p$-adic topology (see \cite{ag2} 3.2.15). The morphisms of the sequence \eqref{p2-fhtft24a} are $\pi_1(\oY'^{\star\rhd},\oy')$-equivariant.
We deduce from this an action of $\pi_1(\oY'^{\star\rhd},\oy')$ on $\fC'^{\oy',(r,r')}_{Y'}$ by ring homomorphisms
that is continuous for the $p$-adic topology and which extends the canonical action on $\hoR^{\oy'}_{Y'}$. 
Observe that these representations depend on
the adequate chart $((P',\gamma'),(\mN,\nu'),\vartheta')$. However, they do not depend on it if $Y'$ is an object of $\bQ'$ \eqref{p2-htaft10} by \ref{p2-rlps4}.

As in \ref{p2-htaft13}, we deduce from $\fF'^{\oy',(r,r')}_{Y'}$ a $\ocB'_{Y',n}$-extension of $\oY'^\rhd_\fet$, 
\begin{equation}\label{p2-fhtft24c}
0\rightarrow \ocB'_{Y',n}\rightarrow \cF'^{(r,r')}_{Y',n}\rightarrow
\Omega'^{(r,r')}(\coY')\otimes_{\co_{\coX'}(\coY')}\ocB'_{Y',n} \rightarrow 0. 
\end{equation}
Similarly, we deduce from $\fC'^{\oy',(r,r')}_{Y'}$ a $\ocB'_{Y',n}$-algebra $\cC'^{(r,r')}_{Y',n}$ of $\oY'^\rhd_\fet$. We have a canonical isomorphism of $\ocB'_{Y',n}$-algebra
\begin{equation}\label{p2-fhtft24d}
\cC'^{(r,r')}_{Y',n}\stackrel{\sim}{\rightarrow}\underset{\underset{m\geq 0}{\longrightarrow}}\lim\ \rS^m_{\ocB'_{Y',n}}(\cF'^{(r,r')}_{Y',n}).
\end{equation}
Observe that $\cF'^{(r,r')}_{Y',n}$ and $\cC'^{(r,r')}_{Y',n}$ depend on the choice of the adequate chart $((P',\gamma'),(\mN,\nu'),\vartheta')$.
However, they do not depend on it if $Y'$ is an object of $\bQ'$ \eqref{p2-htaft10} by \ref{p2-rlps4}.

Let $t,t'$ be rational numbers such that $t\geq t'\geq 0$, $r\geq t$ and $r'\geq t'$.
We easily deduce from \eqref{p2-hta14i} a canonical $\hoR^{\oy'}_{Y'}$-linear $\pi_1(\oY'^{\star\rhd},\oy')$ -equivariant morphism
\begin{equation}\label{p2-fhtft24i}
\fa'^{\oy',(r,r',t,t')}_{Y'}\colon \fF'^{\oy',(r,r')}_{Y'}\rightarrow \fF'^{\oy',(t,t')}_{Y'}
\end{equation}
that fits into a commutative diagram 
\begin{equation}\label{p2-fhtft24j}
\xymatrix{
0\ar[r]&{\hoR^{\oy'}_{Y'}}\ar@{=}[d]\ar[r]&{\fF'^{\oy',(r,r')}_{Y'}}\ar[r]\ar[d]^{\fa'^{\oy',(r,r',t,t')}_{Y'}}&
{\Omega'^{(r,r')}(\coY')\otimes_{\co_{\coX'}(\coY')}\hoR^{\oy'}_{Y'}}\ar[d]^{\pi'^{(r,r',t,t')}\otimes \id}\ar[r]&0\\
0\ar[r]&{\hoR^{\oy'}_{Y'}}\ar[r]&{\fF'^{\oy',(t,t')}_{Y'}}\ar[r]&{\Omega'^{(t,t')}(\coY')\otimes_{\co_{\coX'}(\coY')}\hoR^{\oy'}_{Y'}}\ar[r]&0,}
\end{equation}
where $\pi'^{(r,r',t,t')}$ is defined in \eqref{p2-fhtft9k}. 

Observe that $\fF'^{\oy',(r,r)}_{Y'}=\fF'^{\oy',(r)}_{Y'}$ \eqref{p2-htaft12d}. Hence, $\fF'^{\oy',(r,r')}_{Y'}$ is canonically isomorphic by \eqref{p2-fhtft24j} 
to the pullback of the extension $\fF'^{\oy',(r')}_{Y'}$ by the morphism $\pi'^{(r,r',r',r')}\otimes \id$.

We deduce from \eqref{p2-fhtft24i} a canonical $\ocB'_{Y',n}$-linear morphism
\begin{equation}\label{p2-fhtft24f}
\tta'^{(r,r',t,t')}_{Y',n}\colon\cF'^{(r,r')}_{Y',n}\rightarrow \cF'^{(t,t')}_{Y',n}
\end{equation}
that fits into a commutative diagram 
\begin{equation}\label{p2-fhtft24g}
\xymatrix{
0\ar[r]&{\ocB'_{Y', n}}\ar@{=}[d]\ar[r]&{\cF'^{(r,r')}_{Y',n}}\ar[r]\ar[d]^{\tta'^{(r,r',t,t')}_{Y',n}}&{\Omega'^{(r,r')}(\coY')\otimes_{\co_{\coX'}(\coY')}\ocB'_{Y',n}}\ar[d]^{\pi'^{(r,r',t,t')}\otimes \id}\ar[r]&0\\
0\ar[r]&{\ocB'_{Y', n}}\ar[r]&{\cF'^{(t,t')}_{Y',n}}\ar[r]&{\Omega'^{(t,t')}(\coY')\otimes_{\co_{\coX'}(\coY')}\ocB'_{Y',n}}\ar[r]&0.}
\end{equation}
It induces a homomorphism of $\ocB'_{Y',n}$-algebras
\begin{equation}\label{p2-fhtft24h}
\alpha'^{(r,r',t,t')}_{Y',n}\colon \cC'^{(r,r')}_{Y',n}\rightarrow \cC'^{(t,t')}_{Y',n}.
\end{equation}

\begin{rema}
For every object $Y'$ of $\bQ'$ \eqref{p2-htaft10}, the modules and algebras introduced in \ref{p2-fhtft24} 
do not depend on the choice of an adequate chart for the morphism $f'|Y'\colon (Y',\cM_{X'}|Y')\rightarrow (S,\cM_S)$. 
Indeed by \ref{p2-rlps4}, with the notation of \ref{p2-rlps3}, we can obtain them 
using the logarithmic scheme $(\tmY'^{\oy'},\cM'_{\tmY'^{\oy'}})$ defined as one or the other of the logarithmic schemes
\begin{equation}
(\cA_2(\mY'^{\oy'}),\cM'_{\cA_2(\mY'^{\oy'})})\ \ \ {\rm or} \ \ \ (\cA^{\ast}_2(\mY'^{\oy'}/S),\cM'_{\cA^{\ast}_2(\mY'^{\oy'}/S)})
\end{equation}
associated with $f'|Y'$ (\cite{ag2} 3.2.8 and 3.2.10), 
depending on whether we are in the absolute or relative case \eqref{p2-ncgt3} 
respectively (mind the $^\prime$ in $\cM'_{\tmY'^{\oy'}}$, see \cite{ag2} 3.2.13). 
\end{rema}

\subsection{}\label{p2-fhtft26}
For any rational numbers $r\geq r' \geq 0$, any integer $n\geq 0$ and any object $Y'$ of $\bP'$ such that $Y_s$ is empty,
we set $\cC'^{(r,r')}_{Y',n}=\cF'^{(r,r')}_{Y',n}=0$. The exact sequence \eqref{p2-fhtft24c} still holds in this case, since $\ocB'_{Y'}$ is a  $\oK$-algebra.

\subsection{}\label{p2-fhtft25}
Let $r,r'$ be rational numbers such that $r\geq r'\geq 0$, $n$ an integer $\geq 0$. We take again the notation of \ref{p2-htaft10} for $f'$. 
By \ref{p2-htaft14}, we have canonical objects $\{Y'\mapsto \cF'^{(r,r')}_{Y',n} \}$ and $\{Y'\mapsto \cC'^{(r,r')}_{Y',n}\}$ 
of $\cP(E_{\bQ'}/\bQ')$ \eqref{p2-qfc20}. They correspond by the equivalence of categories \eqref{p2-htaft10c},
to presheaves on $E'_{\bQ'}$, of modules and algebras, respectively, relative to the ring $\nu'_\rp(\ocB')$ \eqref{p2-htaft10}. 
Hence, by \eqref{p2-htaft16i}, 
there exist a canonical $\ocB'_n$-module $\cF'^{(r,r')}_n$ and a canonical $\ocB'_n$-algebra $\cC'^{(r,r')}_n$ such that 
\begin{eqnarray}
\nu'_\rs(\cF'^{(r,r')}_n)&=&\{Y'\mapsto \cF'^{(r,r')}_{Y',n}\}^a,\label{p2-fhtft25b}\\
\nu'_\rs(\cC'^{(r,r')}_n)&=&\{Y'\mapsto \cC'^{(r,r')}_{Y',n}\}^a,\label{p2-fhtft25c}
\end{eqnarray}
where the exponent $a$ means the associated sheaf and $\nu'_\rs$ is the equivalence of categories \eqref{p2-htaft10d}.

Let $t,t'$ be rational numbers such that $t\geq t'\geq 0$, $r\geq t$ and $r'\geq t'$.
The morphisms \eqref{p2-fhtft24f} induce a $\ocB'_n$-linear morphism
\begin{equation}\label{p2-fhtft25d}
\tta'^{(r,r',t,t')}_n\colon \cF'^{(r,r')}_n\rightarrow \cF'^{(t,t')}_n.
\end{equation}
The homomorphisms \eqref{p2-fhtft24h} induce a homomorphism of $\ocB'_n$-algebras
\begin{equation}\label{p2-fhtft25e}
\alpha'^{r,r',t,t'}_n\colon \cC'^{(r,r')}_n\rightarrow \cC'^{(t,t')}_n.
\end{equation}

\begin{prop}\label{p2-fhtft27}
Let $r,r'$ be rational numbers such that $r\geq r'\geq 0$, $n$ an integer $\geq 1$. Then,
\begin{itemize}
\item[{\rm (i)}] The sheaves $\cF'^{(r,r')}_n$ and $\cC'^{(r,r')}_n$ are objects of $\tE'_s$.
\item[{\rm (ii)}] With the notation of \eqref{p2-htaft7c}, setting $\Omega'^{(r,r')}_n=\Omega'^{(r,r')}/p^n\Omega'^{(r,r')}$,
we have a canonical locally split exact sequence of $\ocB_n$-modules
\begin{equation}\label{p2-fhtft27a}
0\rightarrow \ocB'_n\rightarrow \cF'^{(r,r')}_n\rightarrow
\sigma'^*_n(\Omega'^{(r,r')}_n)\rightarrow 0.
\end{equation} 
It induces for every integer $m\geq 1$, an exact sequence of $\ocB'_n$-modules
\begin{equation}\label{p2-fhtft27b}
0\rightarrow \rS^{m-1}_{\ocB'_n}(\cF'^{(r,r')}_n)\rightarrow \rS^m_{\ocB'_n}(\cF'^{(r,r')}_n)\rightarrow
\sigma'^*_n(\rS^m_{\co_{\oX'_n}}(\Omega'^{(r,r')}_n))\rightarrow 0.
\end{equation}
In particular, the $\ocB'_n$-modules $(\rS^m_{\ocB'_n}(\cF'^{(r,r')}_n))_{m\in \mN}$ form a filtered direct system.
\item[{\rm (iii)}] We have a canonical isomorphism of $\ocB'_n$-algebras
\begin{equation}\label{p2-fhtft27c}
\cC'^{(r,r')}_n \stackrel{\sim}{\rightarrow}\underset{\underset{m\geq 0}{\longrightarrow}}\lim\ \rS^m_{\ocB'_n}(\cF'^{(r,r')}_n).
\end{equation}
\item[{\rm (iv)}] We have $\cF'^{(r,r)}_n=\cF'^{(r)}_n$ \eqref{p2-htaft16d} and $\cC'^{(r,r)}_n=\cC'^{(r)}_n$ \eqref{p2-htaft16d}.
Moreover, the exact sequence \eqref{p2-fhtft27a} for $(r,r)$ corresponds to \eqref{p2-htaft17a}.
\item[{\rm (v)}] For all rational numbers $t,t'$ such that $t\geq t'\geq 0$, $r\geq t$ and $r'\geq t'$, the diagram
\begin{equation}\label{p2-fhtft27d}
\xymatrix{
0\ar[r]&{\ocB'_n}\ar[r]\ar@{=}[d]&
{\cF'^{(r,r')}_n}\ar[r]\ar[d]^{\tta'^{(r,r',t,t')}_n}&{\sigma'^*_n(\Omega'^{(r,r')}_n)}\ar[r]\ar[d]^{\sigma'^*_n(\pi'^{(r,r',t,t')})}& 0\\
0\ar[r]&{\ocB'_n}\ar[r]&{\cF'^{(t,t')}_n}\ar[r]&{\sigma'^*_n(\Omega'^{(t,t')}_n)}\ar[r]&0,}
\end{equation}
where the horizontal lines are the exact sequences \eqref{p2-fhtft27a}, is commutative.
Moreover, the morphisms $\tta'^{(r,r',t,t')}_n$ and $\alpha'^{(r,r',t,t')}_n$ are compatible with the isomorphisms
\eqref{p2-fhtft27c} for $(r,r')$ and $(t,t')$.
\end{itemize}
\end{prop}

The proof is similar to that of (\cite{agt} III.10.22). 

\subsection{}\label{p2-fhtft23}
Let $r,r'$ be rational numbers such that $r\geq r'\geq 0$, $n$ an integer $\geq 1$.  There is a unique $\ocB'_n$-derivation of $\cC'^{(r,r')}_n$ \eqref{p2-fhtft27c}
\begin{equation}\label{p2-fhtft23a}
d_{\cC'^{(r,r')}_n}\colon \cC'^{(r,r')}_n\rightarrow \sigma'^*_n(\Omega'^{(r,r')}_n)\otimes_{\ocB'_n}\cC'^{(r,r')}_n
\end{equation}
that extends the canonical morphism $\cF'^{(r,r')}_n\rightarrow \sigma'^*_n(\Omega'^{(r,r')}_n)$ \eqref{p2-fhtft27a}.
It canonically identifies with the universal $\ocB'_n$-derivation of $\cC'^{(r,r')}_n$ \eqref{p1-imdpa19}. 
It is also a Higgs $\ocB'_n$-field on $\cC'^{(r,r')}_n$ with coefficients in $\sigma'^*_n(\Omega'^{(r,r')}_n)$ \eqref{p1-delta-con1}.
We set $\Omega'_n=\Omega'/p^n\Omega'$ \eqref{p2-fhtft9d} and denote by 
\begin{equation}\label{p2-fhtft23b}
\delta_{\cC'^{(r,r')}_n}\colon \cC'^{(r,r')}_n\rightarrow \sigma'^*_n(\Omega'_n)\otimes_{\ocB'_n}\cC'^{(r,r')}_n
\end{equation}
the $\ocB'_n$-derivation induced by $d_{\cC'^{(r,r')}_n}$ and the canonical morphism $\Omega'^{(r,r')}_n\rightarrow \Omega'_n$ \eqref{p2-fhtft9j}. 

Let $t,t'$ be rational numbers such that $t\geq t'\geq 0$, $r\geq t$ and $r'\geq t'$. 
It follows from \ref{p2-fhtft27}(v) that we have
\begin{equation}\label{p2-fhtft23c}
(\id \otimes \alpha'^{(r,r',t,t')}_n) \circ \delta_{\cC'^{(r,r')}_n}=\delta_{\cC'^{(t,t')}_n} \circ \alpha'^{(r,r',t,t')}_n.
\end{equation}

\subsection{}\label{p2-fhtft28}
Let $r,r'$ be rational numbers such that $r\geq r'\geq 0$. 
For all integers $m\geq n\geq 1$, we have a canonical $\ocB'_m$-linear morphism
$\cF'^{(r,r')}_m\rightarrow \cF'^{(r,r')}_n$ compatible with the exact sequence \eqref{p2-fhtft27a}
and a canonical homomorphism of $\ocB'_m$-algebras $\cC'^{(r,r')}_m\rightarrow \cC'^{(r,r')}_n$ such that
the induced morphisms
\begin{equation}\label{p2-fhtft28a}
\cF'^{(r,r')}_m\otimes_{\ocB'_m}\ocB'_n\rightarrow \cF'^{(r,r')}_n\ \ \ {\rm and}\ \ \
\cC'^{(r,r')}_m\otimes_{\ocB'_m}\ocB'_n\rightarrow \cC'^{(r,r')}_n
\end{equation}
are isomorphisms. These morphisms form compatible systems when $m$ and $n$ vary.
With the notation of \ref{p2-fhtft22}, we consider the $\bvocB'$-module $\bvcF'^{(r,r')}=(\cF'^{(r,r')}_{n+1})_{n\in \mN}$ of $\tE'^{\mN^\circ}_s$,
and the $\bvocB'$-algebra $\bvcC'^{(r,r')}=(\cC'^{(r,r')}_{n+1})_{n\in \mN}$ of $\tE_s^{\mN^\circ}$.
By (\cite{agt} III.7.3(i), (III.7.5.4) and (III.7.12.1)), the exact sequence \eqref{p2-fhtft27a}
induces an exact sequence of $\bvocB'$-modules
\begin{equation}\label{p2-fhtft28b}
0\rightarrow \bvocB'\rightarrow \bvcF'^{(r,r')}\rightarrow \bvsigma'^*(\bvOmega'^{(r,r')})\rightarrow 0,
\end{equation}
where we denoted by $\bvOmega'^{(r,r')}$ the $\co_{\bvoX'}$-module $(\Omega'^{(r,r')}_{n+1})_{n\in \mN}$ of $X'_{s,\et}$ \eqref{p2-htaft8}. 
Since the $\co_{\coX'}$-module $\Omega'^{(r,r')}$ is locally free of finite type,
the $\bvocB'$-module $\bvsigma'^*(\bvOmega'^{(r,r')})$ is locally free of finite type
and the sequence \eqref{p2-fhtft28b} is locally split.
It induces for every integer $m\geq 1$ an exact sequence of $\bvocB'$-modules
\begin{equation}\label{p2-fhtft28c}
0\rightarrow \rS^{m-1}_{\bvocB'}(\bvcF'^{(r,r')})\rightarrow \rS^m_{\bvocB'}(\bvcF'^{(r,r')})\rightarrow
\bvsigma'^*(\rS^m_{\co_{\bvoX'}}(\bvOmega'^{(r,r')}))\rightarrow 0.
\end{equation}
In particular, the $\bvocB'$-modules $(\rS^m_{\bvocB'}(\bvcF'^{(r,r')}))_{m\in \mN}$ form a filtered direct system.
By (\cite{agt} III.7.3(i) and (III.7.12.3)), we have a canonical isomorphism of $\bvocB$-algebras
\begin{equation}\label{p2-fhtft28d}
\bvcC'^{(r,r')}\stackrel{\sim}{\rightarrow}\underset{\underset{m\geq 0}{\longrightarrow}}\lim\ \rS^m_{\bvocB'}(\bvcF'^{(r,r')}).
\end{equation}

Let $t,t'$ be rational numbers such that $t\geq t'\geq 0$, $r\geq t$ and $r'\geq t'$. The homomorphism $\alpha'^{r,r',t,t'}_n$ \eqref{p2-fhtft25e} 
induce a homomorphism of $\bvocB'$-algebras
\begin{equation}\label{p2-fhtft28e}
\bvalpha'^{r,r',t,t'}\colon \bvcC'^{(r,r')}\rightarrow \bvcC'^{(t,t')}.
\end{equation}

\subsection{}\label{p2-fhtft29}
Let $r,r'$ be rational numbers such that $r\geq r'\geq 0$. By \eqref{p2-htaft8e}, we have a canonical $\bvocB'$-linear isomorphism
\begin{equation}\label{p2-fhtft29a}
\hupsigma'^*(\hOmega')\stackrel{\sim}{\rightarrow} \bvsigma'^*(\bvOmega').
\end{equation}
The $\bvocB'$-derivation $(\delta_{\cC'^{(r,r')}_n})_{n\in \mN}$ \eqref{p2-fhtft23b} of $\bvcC'^{(r,r')}$ and the isomorphism above define therefore a $\bvocB'$-derivation  
\begin{equation}\label{p2-fhtft29b}
\delta_{\bvcC'^{(r,r')}}\colon \bvcC'^{(r,r')}\rightarrow \hupsigma'^*(\hOmega')\otimes_{\bvocB'}\bvcC'^{(r,r')},
\end{equation}
which is also a Higgs $\bvocB'$-field on $\bvcC'^{(r,r')}$.

\subsection{}\label{p2-fhtft33}
We denote by $I$ the set of triples of rational numbers $\ur=(r_1,r_2,r_3)$ such that $r_1\geq r_2\geq r_3\geq 0$. 
We take again the assumption and notation of \ref{p2-fhtft32}. 
For any element $\ur=(r_1,r_2,r_3)$ of $I$, we associate with the diagram \eqref{p2-fhtft32c} a canonical morphism of $\hoRp^{\oy'}_{Y'}$-modules \eqref{p2-hta15h}
\begin{equation}\label{p2-fhtft33a}
\upvarphi^{\oy',\ur}_\mu\colon \fF^{\oy,(r_2)}_Y\otimes_{\hoR^\oy_Y}\hoRp^{\oy'}_{Y'}\rightarrow \fF'^{\oy,(r_1,r_2)}_{Y'} \otimes_{R'_{\uptau,Y'}}\fF^{(r_3)}_{\uptau,Y'},
\end{equation}
where $R'_{\uptau,Y'}$ and $\fF^{(r_3)}_{\uptau,Y'}$ are defined in \eqref{p2-fhtft3g}, 
and a canonical morphism of $\hoRp^{\oy'}_{Y'}$-algebras \eqref{p2-hta15i}
\begin{equation}\label{p2-fhtft33b}
\upphi^{\oy',\ur}_\mu\colon \fC^{\oy,(r_2)}_Y\otimes_{\hoR^\oy_Y}\hoRp^{\oy'}_{Y'}\rightarrow \fC'^{\oy,(r_1,r_2)}_{Y'} \otimes_{R'_{\uptau,Y'}}\fC^{(r_3)}_{\uptau,Y'},
\end{equation}
compatible with $\upvarphi^{\oy',\ur}_\mu$; see \ref{p2-fhtal7} and \ref{p2-fhtal11}. 
The latter induces by linearization a morphism of $(\hoRp^{\oy'}_{Y'}\otimes_{R'_{\uptau,Y'}}\fC^{(r_3)}_{\uptau,Y'})$-algebras \eqref{p2-hta15e}
\begin{equation}\label{p2-fhtft33c}
\uppsi^{\oy',\ur}_\mu\colon \fC^{\oy,(r_2)}_Y\otimes_{\hoR^\oy_Y}\hoRp^{\oy'}_{Y'}\otimes_{R'_{\uptau,Y'}}\fC^{(r_3)}_{\uptau,Y'}
\rightarrow \fC'^{\oy,(r_1,r_2)}_{Y'} \otimes_{R'_{\uptau,Y'}}\fC^{(r_3)}_{\uptau,Y'}.
\end{equation}
These morphisms are $\pi_1(\oY'^{\star\rhd},\oy')$-equivariant, when we equip $R'_{\uptau,Y'}$, $\fF^{(r_3)}_{\uptau,Y'}$ and $\fC^{(r_3)}_{\uptau,Y'}$ 
with the trivial actions, see \ref{p2-fhtal12} and \ref{p2-fhtal120}. 
By (\cite{agt} VI.9.11), for every integer $n\geq 0$, we deduce a canonical morphism of $\ocB'_{Y',n}$-modules of $\oY'^\rhd_\fet$
\begin{equation}\label{p2-fhtft33d}
\upvarphi^{\ur}_{\mu,n}\colon \tmu_n^*(\cF^{(r_2)}_{Y,n})\rightarrow \cF'^{(r_1,r_2)}_{Y',n} \otimes_{R'_{\uptau,Y'}}\fF^{(r_3)}_{\uptau,Y'},
\end{equation}
where $\tmu_n$ is the morphism of ringed topos \eqref{p2-fhtft302c}, 
and a canonical  morphism of $\ocB'_{Y',n}$-algebras 
\begin{equation}\label{p2-fhtft33e}
\upphi^{\ur}_{\mu,n}\colon \tmu_n^*(\cC^{(r_2)}_{Y,n})\rightarrow \cC'^{(r_1,r_2)}_{Y',n} \otimes_{R'_{\uptau,Y'}}\fC^{(r_3)}_{\uptau,Y'},
\end{equation}
compatible with $\upvarphi^{\ur}_{\mu,n}$. These morphisms are independent of the choice of $\oy'$. 
The morphism $\upphi^{\ur}_{\mu,n}$ induces by linearization a morphism of $(\ocB'_{Y',n}\otimes_{R'_{\uptau,Y'}}\fC^{(r_3)}_{\uptau,Y'})$-algebras 
\begin{equation}\label{p2-fhtft33f}
\uppsi^{\ur}_{\mu,n}\colon \tmu_n^*(\cC^{(r_2)}_{Y,n})\otimes_{R'_{\uptau,Y'}}\fC^{(r_3)}_{\uptau,Y'}
\rightarrow \cC'^{(r_1,r_2)}_{Y',n} \otimes_{R'_{\uptau,Y'}}\fC^{(r_3)}_{\uptau,Y'}.
\end{equation}

\subsection{}\label{p2-fhtft34}
Let $(u',u)\colon (\mu_1\colon Y'_1\rightarrow Y_1)\rightarrow (\mu_2\colon Y'_2\rightarrow Y_2)$ 
be a morphism of $\mM(\bQ',\bQ)$ \eqref{p2-fhtft30}, $\oy'_1$ a geometric point 
of $\oY'^\rhd_1$, $\oy'_2=\ou'^\rhd(\oy'_1)$, $\oy_1=\tmu_1(\oy'_1)$, $\oy_2=\tmu_2(\oy'_2)$ \eqref{p2-fhtft302a}, 
$\ur=(r_1,r_2,r_3)$ an element of $I$, $n$ an integer $\geq 0$. 
\begin{equation}\label{p2-fhtft34a}
\xymatrix{
Y'_1\ar[r]^{\mu_1}\ar[d]_{u'}&Y_1\ar[d]^u\\
Y'_2\ar[r]^{\mu_2}&Y_2.}
\end{equation}
We take again the notation introduced in \ref{p2-fhtft340}. 
We claim that the diagram
\begin{equation}\label{p2-fhtft34c}
\xymatrix{
{\ou'^{\rhd*}_n(\tmu_{2,n}^*(\cF^{(r_2)}_{Y_2,n}))}\ar[rr]^-(0.5){\ou'^{\rhd*}_n(\upvarphi^{\ur}_{\mu_2,n})}\ar[d]_c&&
{\ou'^{\rhd*}_n(\cF'^{(r_1,r_2)}_{Y' _2,n}\otimes_{R'_{\uptau,Y'_2}}\fF^{(r_3)}_{\uptau,Y'_2})}\ar[d]^{a'\otimes a_\uptau}\\
{\tmu_{1,n}^*(\ou^{\circ*}_n(\cF^{(r_2)}_{Y_2,n}))}\ar[r]^-(0.5){\tmu_{1,n}^*(a) }&{\tmu_{1,n}^*(\cF^{(r_2)}_{Y_1,n})}\ar[r]^-(0.5){\upvarphi^{\ur}_{\upmu_1,n}} &
{\cF'^{(r_1,r_2)}_{Y'_1,n}\otimes_{R'_{\uptau,Y'_1}}\fF^{(r_3)}_{\uptau,Y'_1}},}
\end{equation}
where $a\colon \ou^{\circ*}_n(\cF^{(r_2)}_{Y_2,n})\rightarrow \cF^{(r_2)}_{Y_1,n}$ 
(resp.\ $a'\colon \ou'^{\rhd*}_n(\cF'^{(r_1,r_2)}_{Y'_2,n})\rightarrow \cF'^{(r_1,r_2)}_{Y'_1,n}$) 
is the adjoint of the transition morphism of the v-presheaf
$\{Y\in \bQ^\circ\mapsto \cF^{(r_2)}_{Y,n}\}$ \eqref{p2-htaft16d} (resp.\ $\{Y'\in \bQ'^\circ\mapsto \cF'^{(r_1,r_2)}_{Y',n}\}$ \eqref{p2-fhtft25b})  \eqref{p2-cmt45},
$a_\uptau\colon \fF^{(r_3)}_{\uptau,Y'_2}\rightarrow \fF^{(r_3)}_{\uptau,Y'_1}$ is the canonical morphism \eqref{p2-fhtft3g}, 
and $c$ is the isomorphism induced by the commutative diagram \eqref{p2-fhtft34a}, is commutative. 
Indeed, we are immediately reduced to checking the commutativity of the analogous diagram for the Higgs--Tate extensions appearing in \eqref{p2-fhtft33a}, 
and furthermore to the case where $r_1=r_2=r_3=0$ as these Higgs--Tate extensions are submodules of those with vanishing radii; 
then the statement follows from \eqref{p2-fhtft340b} by \ref{p1-rdt40} and \ref{p1-rdt9}. 
Observe that the invertibility condition required in \ref{p1-rdt4} is satisfied by \ref{p2-hta2}.

We deduce from \eqref{p2-fhtft34c} that the diagram
\begin{equation}\label{p2-fhtft34d}
\xymatrix{
{\ou'^{\rhd*}_n(\tmu_{2,n}^*(\cC^{(r_2)}_{Y_2,n}))}\ar[rr]^-(0.5){\ou'^{\rhd*}_n(\upphi^{\ur}_{\mu_2,n})}\ar[d]_d&&
{\ou'^{\rhd*}_n(\cC'^{(r_1,r_2)}_{Y'_2,n}\otimes_{R'_{\uptau,Y'_2}}\fC^{(r_3)}_{\uptau,Y'_2})}\ar[d]^{b'\otimes b_\uptau}\\
{\tmu_{1,n}^*(\ou^{\circ*}_n(\cC^{(r_2)}_{Y_2,n}))}\ar[r]^-(0.5){\tmu_{1,n}^*(b)}&{\tmu_{1,n}^*(\cC^{(r_2)}_{Y_1,n})}\ar[r]^-(0.5){\upphi^{\ur}_{\mu_1,n}}
&{\cC'^{(r_1,r_2)}_{Y'_1,n}\otimes_{R'_{\uptau,Y'_1}}\fC^{(r_3)}_{\uptau,Y'_1}}}
\end{equation}
is commutative, 
where $b\colon \ou^{\circ*}_n(\cC^{(r_2)}_{Y_2,n})\rightarrow \cC^{(r_2)}_{Y_1,n}$ 
(resp.\ $b'\colon \ou'^{\rhd*}_n(\cC'^{(r_1,r_2)}_{Y'_2,n})\rightarrow \cC'^{(r_1,r_2)}_{Y'_1,n}$) 
is the adjoint of the transition morphism of the v-presheaf
$\{Y\in \bQ^\circ\mapsto \cC^{(r_2)}_{Y,n}\}$ \eqref{p2-htaft16e} (resp.\ $\{Y'\in \bQ'^\circ\mapsto \cC'^{(r_1,r_2)}_{Y',n}\}$ \eqref{p2-fhtft25c}), 
$b_\uptau\colon \fC^{(r_3)}_{\uptau,Y'_2}\rightarrow \fC^{(r_3)}_{\uptau,Y'_1}$ is the canonical morphism \eqref{p2-fhtft3g}, 
and $d$ is the isomorphism induced by the commutative diagram \eqref{p2-fhtft34a}.
Therefore, by linearization, the diagram 
\begin{equation}\label{p2-fhtft34e}
{\small 
\xymatrix{
{\ou'^{\rhd*}_n(\tmu_{2,n}^*(\cC^{(r_2)}_{Y_2,n})\otimes_{R'_{\uptau,Y'_2}}\fC^{(r_3)}_{\uptau,Y'_2})}\ar[rrr]^-(0.5){\ou'^{\rhd*}_n(\uppsi^{\ur}_{\mu_2,n})}\ar[d]_{d\otimes b_\uptau}&&&
{\ou'^{\rhd*}_n(\cC'^{(r_1,r_2)}_{Y'_2,n}\otimes_{R'_{\uptau,Y'_2}}\fC^{(r_3)}_{\uptau,Y'_2})}\ar[d]^{b'\otimes b_\uptau}\\
{\tmu_{1,n}^*(\ou^{\circ*}_n(\cC^{(r_2)}_{Y_2,n}))\otimes_{R'_{\uptau,Y'_1}}\fC^{(r_3)}_{\uptau,Y'_1}}\ar[rr]^-(0.5){\tmu_{1,n}^*(b)\otimes \id}&&
{\tmu_{1,n}^*(\cC^{(r_2)}_{Y_1,n})\otimes_{R'_{\uptau,Y'_1}}\fC^{(r_3)}_{\uptau,Y'_1}}\ar[r]^-(0.5){\uppsi^{\ur}_{\mu_1,n}}
&{\cC'^{(r_1,r_2)}_{Y'_1,n}\otimes_{R'_{\uptau,Y'_1}}\fC^{(r_3)}_{\uptau,Y'_1}}}}
\end{equation}
is commutative.

\subsection{}\label{p2-fhtft38}
Let $\ur=(r_1,r_2,r_3)\in I$ \eqref{p2-fhtft33}, $n\geq 0$ an integer. We set 
\begin{eqnarray}
\cF'^{\ur}_n&=&\cF'^{(r_1,r_2)}_n\otimes_{\ocB'_n}\sigma'^*_n(\cF^{(r_3)}_{\uptau,n}),\label{p2-fhtft38j}\\
\cC'^{\ur}_n&=&\cC'^{(r_1,r_2)}_n\otimes_{\ocB'_n}\sigma'^*_n(\cC^{(r_3)}_{\uptau,n}),\label{p2-fhtft38k}
\end{eqnarray}
where $\sigma'_n$ denotes the morphism of ringed topos \eqref{p2-htaft7c} (for $f'$). 

We have a natural object 
\begin{equation}\label{p2-fhtft38a}
\{U'\mapsto \cF'^{(r_1,r_2)}_{U',n}\otimes_{R'_{\uptau,U'}}\fF^{(r_3)}_{\uptau,U'}\}
\end{equation}
of $\cP(E'_{\bQ'}/\bQ')$; see \ref{p2-fhtft3} and \ref{p2-fhtft25}. 
By \eqref{p2-htaft10c}, it determines a presheaf of $\nu'_\rp(\ocB')$-modules on $E'_{\bQ'}$ \eqref{p2-htaft10}. 
By \eqref{p2-htaft10e} and (\cite{agt} VI.534(ii), 8.9 and 5.17), we have a canonical isomorphism
\begin{equation}\label{p2-fhtft38b}
\nu'_\rs(\cF'^{\ur}_n)\stackrel{\sim}{\rightarrow}
\{U'\mapsto \cF'^{(r_1,r_2)}_{U',n}\otimes_{R'_{\uptau,U'}}\fF^{(r_3)}_{\uptau,U'}\}^a,
\end{equation}
where the exponent $a$ means the associated sheaf and  
$\nu'_\rs$ is the equivalence of categories \eqref{p2-htaft10d} (for $f'$).

It follows from \eqref{p2-fhtft34c} that the morphisms $\upvarphi^{\ur}_{\mu,n}$ \eqref{p2-fhtft33d}, for $\mu\in \mM(\bQ,\bQ')$, 
induce an $\mM(\bQ,\bQ')$-system of morphisms from $\{U\mapsto \cF^{(r_2)}_{U,n}\}$ to 
$\{U'\mapsto \cF'^{(r_1,r_2)}_{U',n}\otimes_{R'_{\uptau,U'}}\fF^{(r_3)}_{\uptau,U'}\}$ in the sense of \ref{p2-fhtft36}, see \ref{p2-fhtft361}.
By \ref{p2-fhtft362}, the latter determines a morphism $\Theta^*(\cF^{(r_2)}_n)\rightarrow \cF'^{\ur}_n$ 
of $\tE'$, which is obviously $\Theta^*(\ocB)$-linear \eqref{p2-htaft10e}. It therefore induces a $\ocB'_n$-linear morphism 
\begin{equation}\label{p2-fhtft38c}
\upvarphi^{\ur}_n\colon \uptheta^*_n(\cF^{(r_2)}_n)\rightarrow \cF'^{\ur}_n,
\end{equation}
where $\uptheta_n$ is the morphism of ringed topos \eqref{p2-fhtft21h}. 

Similarly, the morphisms $\upphi^{\ur}_{\mu,n}$ \eqref{p2-fhtft33e}, for $\mu\in \mM(\bQ,\bQ')$, induce a morphism of $\ocB'_n$-algebras 
\begin{equation}\label{p2-fhtft38d}
\upphi^{\ur}_n\colon \uptheta^*_n(\cC^{(r_2)}_n)\rightarrow \cC'^{\ur}_n,
\end{equation}
compatible with $\upvarphi^{\ur}_n$. It induces a homomorphism of $\sigma'^*_n(\cC^{(r_3)}_{\uptau,n})$-algebras 
\begin{equation}\label{p2-fhtft38e}
\uppsi^{\ur}_n\colon \uptheta^*_n(\cC^{(r_2)}_n)\otimes_{\ocB'_n}\sigma'^*_n(\cC^{(r_3)}_{\uptau,n})\rightarrow \cC'^{\ur}_n.
\end{equation}

We set 
\begin{equation}\label{p2-fhtft38l}
\bvcC'^{\ur}=\bvcC'^{(r_1,r_2)}\otimes_{\bvocB'}\hupsigma'^*(\hcC^{(r_3)}_{\uptau}),
\end{equation}
where $\bvcC'^{(r_1,r_2)}$ (resp.\ $\hcC^{(r_3)}_{\uptau}$) is defined in \ref{p2-fhtft28} (resp.\ \ref{p2-fhtft3}) 
and $\hupsigma'$ is the morphism of ringed topos appearing in \eqref{p2-fhtft22b}. 
By (\cite{agt} III.7.5 and III.7.12), the homomorphisms $(\upphi^{\ur}_{m+1})_{m\in \mN}$ define a homomorphism of $\bvocB'$-algebras 
\begin{equation}\label{p2-fhtft38f}
\bvupphi^{\ur}\colon \bvuptheta^*(\bvcC^{(r_2)})\rightarrow \bvcC'^{\ur},
\end{equation}
where $\bvuptheta$ is the morphism of ringed topos appearing in \eqref{p2-fhtft22b}. 
It induces a homomorphism of $\hupsigma'^*(\hcC^{(r_3)}_\uptau)$-algebras 
\begin{equation}\label{p2-fhtft38g}
\bvuppsi^{\ur}\colon \bvuptheta^*(\bvcC^{(r_2)})\otimes_{\bvocB'}\hupsigma'^*(\hcC^{(r_3)}_\uptau)\rightarrow \bvcC'^{\ur}.
\end{equation}

Let $\ur'=(r'_1,r'_2,r'_3)\in I$ be such that $r_i\geq r'_i$ for all $1\leq i\leq 3$. 
The homomorphisms $\alpha_{\uptau}^{r,r'}$ \eqref{p2-fhtft3i} and $\alpha'^{r_1,r_2,r'_1,r'_2}_n$ \eqref{p2-fhtft25e}
and $\bvalpha'^{r_1,r_2,r'_1,r'_2}$ \eqref{p2-fhtft28d} induce homomorphisms 
\begin{eqnarray}
\alpha'^{\ur,\ur'}_n=\alpha'^{r_1,r_2,r'_1,r'_2}_n\otimes \sigma'^*_n(\alpha_{\uptau,n}^{r_3,r'_3})\colon \cC'^{\ur}_n\rightarrow \cC'^{\ur'}_n,\label{p2-fhtft38m}\\
\bvalpha'^{\ur,\ur'}=\bvalpha'^{r_1,r_2,r'_1,r'_2}\otimes \hupsigma'^*(\halpha_{\uptau}^{r_3,r'_3})\colon \bvcC'^{\ur}\rightarrow \bvcC'^{\ur'}.\label{p2-fhtft38n}
\end{eqnarray}
By \ref{p2-fhtft362} and \eqref{p2-fhtal11f}, the diagram 
\begin{equation}\label{p2-fhtft38i}
\xymatrix{
{\uptheta^*_n(\cC^{(r_2)}_n)\otimes_{\ocB'_n}\sigma'^*_n(\cC^{(r_3)}_{\uptau,n})}\ar[r]^-(0.5){\uppsi^\ur_n}
\ar[d]_{\uptheta^*_n(\alpha^{r_2,r'_2}_n)\otimes_{\ocB'_n} \sigma'^*_n(\alpha_{\uptau,n}^{r_3,r'_3})}
&{\cC'^{\ur}_n}
\ar[d]^{\alpha'^{\ur,\ur'}_n}\\
{\uptheta^*_n(\cC^{(r'_2)}_n)\otimes_{\ocB'_n}\sigma'^*_n(\cC^{(r'_3)}_{\uptau,n})}\ar[r]^-(0.5){\uppsi^{\ur'}_n}&
{\cC'^{\ur'}_n,}}
\end{equation}
where $\alpha^{r_2,r'_2}_n$ is defined in \eqref{p2-htaft16g}, is commutative. 
We deduce that the diagram 
\begin{equation}\label{p2-fhtft38h}
\xymatrix{
{\bvuptheta^*(\bvcC^{(r_2)})\otimes_{\bvocB'}\hupsigma'^*(\hcC^{(r_3)}_\uptau)}\ar[r]^-(0.5){\bvuppsi^\ur}
\ar[d]_{\bvuptheta^*(\bvalpha^{r_2,r'_2})\otimes_{\bvocB'} \hupsigma'^*(\halpha_\uptau^{r_3,r'_3})}
&{\bvcC'^{\ur}}
\ar[d]^{\bvalpha'^{r_1,r'_1,r_2,r'_2}\otimes_{\bvocB'} \hupsigma'^*(\halpha_\uptau^{r_3,r'_3})}\\
{\bvuptheta^*(\bvcC^{(r'_2)})\otimes_{\bvocB'}\hupsigma'^*(\hcC^{(r'_3)}_\uptau)}\ar[r]^-(0.5){\bvuppsi^{\ur'}}&
{\bvcC'^{\ur'},}}
\end{equation}
where $\bvalpha^{r_2,r'_2}$ is defined in \eqref{p2-htaft18f}, is commutative. We have similar commutative diagrams for $\upphi^\ur_n$
and $\bvupphi^\ur$. 

\begin{rema}
By \ref{p2-fhtal111}, for all rational numbers $r\geq r'\geq 0$, setting $\ur=(r,r,r') \in I$ \eqref{p2-fhtft38}, we have 
$\bvupphi^{(r,r')}=\bvupphi^\ur$, see \eqref{p2-fhtft380f} and \eqref{p2-fhtft38f}, and 
$\bvuppsi^{(r,r')}=\bvuppsi^\ur$, see \eqref{p2-fhtft380g} and \eqref{p2-fhtft38g}. 
\end{rema}

\section{Relative Faltings topos. Cohomologies of Higgs--Tate algebras}\label{p2-rftchta}

\subsection{}\label{p2-rftchta1}
The assumptions and notation of §\ref{p2-fhtft} and §\ref{p2-rfhtft} remain in force throughout this section.
In particular, the morphism $g\colon (X',\cM_{X'})\rightarrow (X,\cM_X)$ \eqref{p2-fhtft1a} is assumed to be {\em smooth and saturated}. 
We denote, moreover, by $G$ the {\em relative Faltings $\mU$-site} associated with the pair of morphisms
$(h\colon \oX^\circ \rightarrow X,g\colon X'\rightarrow X)$ introduced in (\cite{ag1} 3.4.1). Objects of the category underlying $G$ are
the triples $(U,U'\rightarrow U, V\rightarrow U)$ made up  of an $X$-scheme $U$ and of two morphisms $U'\rightarrow U$ and $V\rightarrow U$
above $g$ and $h$ respectively, i.e., commutative diagrams of morphisms of schemes
\begin{equation}\label{p2-rftchta1a}
\xymatrix{
U'\ar[r]\ar[d]&U\ar[d]&V\ar[l]\ar[d]\\
X'\ar[r]^g&X&\oX^\circ\ar[l]_h}
\end{equation}
such that the morphisms $U\rightarrow X$ and $U'\rightarrow X'$ are étale and the morphism $V\rightarrow \oU^\circ$ is finite étale;
such an object will be denoted by $(U'\rightarrow U\leftarrow V)$. Let $(U'\rightarrow U\leftarrow V)$ and $(U'_1\rightarrow U_1\leftarrow V_1)$ be two objects of $G$.
A morphism from $(U'_1\rightarrow U_1\leftarrow V_1)$ to $(U'\rightarrow U\leftarrow V)$ consists of three morphisms
$U'_1\rightarrow U'$, $U_1\rightarrow U$ and $V_1\rightarrow V$ above $X'$, $X$ and $\oX^\circ$ respectively, that makes the diagram
\begin{equation}\label{p2-rftchta1b}
\xymatrix{
U'_1\ar[r]\ar[d]&U_1\ar[d]&V_1\ar[l]\ar[d]\\
U'\ar[r]&U&V\ar[l]}
\end{equation}
commutative.

We endow $G$ with the  {\em covanishing topology} (\cite{ag1} 3.4.1), i.e., the topology generated by the coverings
\[
\{(U'_i\rightarrow U_i\leftarrow V_i)\rightarrow (U'\rightarrow U\leftarrow V)\}_{i\in I}
\]
of the following three types:
\begin{itemize}
\item[(a)] $V_i=V$, $U_i=U$ for all $i\in I$, and $(U'_i\rightarrow U')_{i\in I}$ is a covering family.
\item[(b)] $U'_i=U'$, $U_i=U$ for all $i\in I$, and $(V_i\rightarrow V)_{i\in I}$ is a covering family.
\item[(c)] $I=\{1\}$, $U'_1=U'$ and the morphism $V_1\rightarrow V\times_{U}U_1$ is an isomorphism
(there is no condition on the morphism $U_1\rightarrow U$).
\end{itemize}

We denote by $\tG$ the {\em  relative Faltings $\mU$-topos} associated with $(h,g)$,
that is, the topos of sheaves of $\mU$-sets on $G$. We denote by
\begin{eqnarray}
\pi\colon \tG \rightarrow X'_\et,\label{p2-rftchta1c}\\
\lambda\colon \tG \rightarrow \oX^\circ_\fet,\label{p2-rftchta1d}
\end{eqnarray}
the canonical morphisms (\cite{ag1} 3.4.4).
By (\cite{ag1} 3.4.18), the canonical commutative diagram \eqref{p2-fhtft1}
\begin{equation}\label{p2-rftchta1e}
\xymatrix{
X'\ar[d]_g&{\oX'^\rhd}\ar[d]^{\upgamma}\ar[l]_{h'}\\
X&{\oX^\circ}\ar[l]_h}
\end{equation}
induces morphisms of topos
\begin{equation}\label{p2-rftchta1f}
\xymatrix{
\tE'\ar[r]^{\varpi}&\tG\ar[r]^{\lgg}&\tE}
\end{equation}
whose composition is the morphism $\Theta\colon \tE'\rightarrow \tE$ \eqref{p2-fhtft21b}.
The triangles and squares of the diagram of morphism of topos
\begin{equation}\label{p2-rftchta1g}
\xymatrix{
&\tE'\ar[d]^{\varpi}\ar[r]^{\beta'}\ar[ld]_{\sigma'} & \oX'^\rhd_\fet\ar[d]^{\upgamma_\fet}\\
X'_\et\ar[d]_{g_\et}&\tG\ar[d]_-(0.5){\lgg}\ar[r]^-(0.4){\lambda}\ar[l]_-(0.4){\pi}&\oX^\circ_\fet\\
X_\et&\tE\ar[ur]_{\beta}\ar[l]_{\sigma}&}
\end{equation}
are commutative up to canonical isomorphisms (\cite{ag1} (3.4.18.4)).

We denote by $X'_\et\gtimes_{X_\et}\oX^\circ_\et$ the oriented product of the morphisms of topos $g_\et\colon X'_\et\rightarrow X_\et$
and $f_\et\colon \oX^\circ_\et\rightarrow X_\et$ (\cite{agt} VI.3.10) and by
\begin{equation}\label{p2-rftchta1h}
\varrho\colon X'_\et\gtimes_{X_\et}\oX^\circ_\et\rightarrow \tG
\end{equation}
the canonical morphism of topos (\cite{ag1} (3.4.9.2)).
We check immediately that the squares of the diagram of morphisms of topos
\begin{equation}\label{p2-rftchta1i}
\xymatrix{
{X'_\et\gtimes_{X'_\et}\oX'^\rhd_\et}\ar[d]\ar[r]^-(0.5){\rho'}&\tE'\ar [d]^\varpi\\
{X'_\et\gtimes_{X_\et}\oX^\circ_\et}\ar[d]\ar[r]^-(0.5){\varrho}&\tG\ar[d]^\lgg\\
{X_\et\gtimes_{X_\et}\oX^\circ_\et}\ar[r]^-(0.5){\rho}&\tE,}
\end{equation}
where $\rho$ and $\rho'$ are the canonical morphisms \eqref{p2-htaft4h},
are commutative up to canonical isomorphisms (\cite{agt} VI.4.10 and \cite{ag1} 3.4.17).

\subsection{}\label{p2-rftchta12}
Let $\mu\colon U'\rightarrow U$ be an object of the category $\mM$ defined in \ref{p2-fhtft30}. 
For any presheaf $F$ on $G$, we define the presheaf $F_{U'\rightarrow U}$ on $\Et_{\rf/\oU^\circ}$ by setting for any $V\in \ob (\Et_{\rf/\oU^\circ})$,
\begin{equation}\label{p2-rftchta12a}
F_{U'\rightarrow U}(V)=F(U'\rightarrow U\leftarrow V).
\end{equation}
If $F$ is a sheaf of $\tG$, then $F_{U'\rightarrow U}$ is a sheaf of $\oU^\circ_\fet$.

\subsection{}\label{p2-rftchta9}
Let $\ox'$ be a geometric point of $X'$, $\uX'$ the strict localization of $X'$ at $\ox'$, $\uX$ the strict localization of $X$ at $g( \ox')$.
We denote by $\uG$ (resp.\ $\tuG$) the relative Faltings site (resp.\ topos) associated with the pair of morphisms
$(\uh\colon \uoX^\circ\rightarrow \uX, \ug\colon \uX'\rightarrow \uX)$ induced by $h$ and $g$ \eqref{p2-fhtft1}, by
\begin{equation}\label{p2-rftchta9a}
\Phi\colon \tuG\rightarrow \tG
\end{equation}
the functoriality morphism (\cite{ag1} (3.4.10.3)) and by
\begin{equation}\label{p2-rftchta9b}
\vartheta\colon \uoX^\circ_\fet\rightarrow \tuG
\end{equation}
the morphism defined in (\cite{ag1} (3.4.26.9)). We set
\begin{equation}\label{p2-rftchta9c}
\phi_{\ox'}=\vartheta^*\circ \Phi^*\colon \tG\rightarrow \uoX^\circ_\fet.
\end{equation}

\subsection{}\label{p2-rftchta2}
We denote by $\tG_s$ the closed subtopos of $\tG$ complement of the open object $\pi^*(X'_\eta)$ (\cite{sga4} IV 9.3.5), and by
\begin{equation}\label{p2-rftchta2a}
\kappa\colon \tG_s\rightarrow \tG
\end{equation}
the canonical embedding (\cite{sga4} IV 9.3.5) (see \cite{ag1} 6.5.2). By virtue of (\cite{sga4} IV 9.4.3), there exists a morphism of topos
\begin{equation}\label{p2-rftchta2b}
\pi_s\colon \tG_s\rightarrow X'_{s,\et}
\end{equation}
unique up to canonical isomorphism such that the diagram
\begin{equation}\label{p2-rftchta2c}
\xymatrix{
{\tG_s}\ar[r]^-(0.5){\pi_s}\ar[d]_{\kappa}&{X'_{s,\et}}\ar[d]^{a'} \\
{\tG}\ar[r]^-(0.5){\pi}&{X'_\et,}}
\end{equation}
where $a'\colon X'_s\rightarrow X$ is the canonical injection, is commutative up to isomorphism, and even $2$-Cartesian.

We have a canonical isomorphism $\varpi^*(\pi^*(X'_\eta))\simeq \sigma'^*(X'_\eta)$ \eqref{p2-rftchta1g}.
Hence, by virtue of (\cite{sga4} IV 9.4.3), there exists a morphism of topos
\begin{equation}\label{p2-rftchta2d}
\varpi_s\colon \tE'_s\rightarrow \tG_s
\end{equation}
unique up to canonical isomorphism such that the diagram
\begin{equation}\label{p2-rftchta2e}
\xymatrix{
{\tE'_s}\ar[r]^{\varpi_s}\ar[d]_{\delta'}&{\tG_s}\ar[d]^{\kappa}\\
{\tE'}\ar[r]^\varpi&{\tG}}
\end{equation}
is commutative up to isomorphism (see \ref{p2-htaft5}).

The functors $\delta'_*$ and $\kappa_*$ being exact,
for every abelian group $F$ of $\tE'_s$ and every integer $q\geq 0$, we have a canonical isomorphism
\begin{equation}\label{p2-rftchta2l}
\kappa_*(\rR^q\varpi_{s*}(F))\stackrel{\sim}{\rightarrow}\rR^q\varpi_*(\delta'_*F).
\end{equation}

It follows from \eqref{p2-rftchta1g} and (\cite{sga4} IV 9.4.3)
that the diagram of morphisms of topos
\begin{equation}\label{p2-rftchta2f}
\xymatrix{
{\tE'_s}\ar[r]^{\varpi_s}\ar[d]_{\sigma'_s}&{\tG_s}\ar[dl]^{\pi_s}\\
{X'_{s,\et}}&}
\end{equation}
is commutative up to canonical isomorphism.

We have a canonical isomorphism $\lgg^*(\sigma^*(X_\eta))\simeq \pi^*(X'_\eta)$ \eqref{p2-rftchta1g}.
Hence, by virtue of (\cite{sga4} IV 9.4.3), there exists a morphism of topos
\begin{equation}\label{p2-rftchta2h}
\lgg_s\colon \tG_s\rightarrow \tE_s
\end{equation}
unique up to canonical isomorphism such that the diagram
\begin{equation}\label{p2-rftchta2i}
\xymatrix{
{\tG_s}\ar[r]^{\lgg_s}\ar[d]_{\kappa}&{\tE_s}\ar[d]^{\delta}\\
{\tG}\ar[r]^\lgg&{\tE}}
\end{equation}
is commutative up to isomorphism.

It follows from \eqref{p2-rftchta1g} and (\cite{sga4} IV 9.4.3)
that the diagram of morphisms of topos
\begin{equation}\label{p2-rftchta2j}
\xymatrix{
{\tG_s}\ar[r]^{\lgg_s}\ar[d]_{\pi_s}&{\tE_s}\ar[d]^{\sigma_s}\\
{X'_{s,\et}}\ar[r]^{g_{s,\et}}&{X_{s,\et}}}
\end{equation}
is commutative up to canonical isomorphism.

It follows again from (\cite{sga4} IV 9.4.3) that the composition $\lgg_s\circ \varpi_s$ is the morphism \eqref{p2-fhtft21e}
\begin{equation}\label{p2-rftchta2k}
\uptheta\colon \tE'_s\rightarrow \tE_s.
\end{equation}

\subsection{}\label{p2-rftchta3}
For any object $(U'\rightarrow U\leftarrow V)$ of $G$, let $\oU'^V$ be the integral closure of $\oU'$ in $U'\times_UV$.
We denote by $\ocB^!$ the presheaf on $G$ defined for any $(U'\rightarrow U\leftarrow V)\in \ob(G)$ by
\begin{equation}\label{p2-rftchta3a}
\ocB^!(U'\rightarrow U\leftarrow V)=\Gamma(\oU'^V,\co_{\oU'^V}).
\end{equation}
Since $\oX'$ is normal and locally irreducible (\cite{agt} III.4.2(iii)),
$\ocB^!$ is a sheaf for the covanishing topology of $G$ by (\cite{ag1} 3.6.4).

By (\cite{ag1} 6.5.18), we have canonical homomorphisms
\begin{eqnarray}
\ocB^!&\rightarrow&\varpi_*(\ocB'),\label{p2-rftchta3b}\\
\ocB&\rightarrow& \lgg_*(\ocB^!),\label{p2-rftchta3c}
\end{eqnarray}
the first of which is an isomorphism by (\cite{ag1} 6.5.19).
We also have a canonical homomorphism
\begin{equation}\label{p2-rftchta3d}
\hbar'_*(\co_{\oX'})\rightarrow \pi_*(\ocB^!),
\end{equation}
where $\hbar'\colon \oX'\rightarrow X'$ is the canonical morphism. 

For any integer $n\geq 1$, we set
\begin{equation}\label{p2-rftchta3f}
\ocB^!_n=\ocB^!/p^n\ocB^!.
\end{equation}
It is a ring of $\tG_s$ (\cite{ag1} 6.5.27).
The canonical homomorphism $\pi^*(\hbar_*(\co_{\oX'}))\rightarrow \ocB^!$ \eqref{p2-rftchta3d}
induces a homomorphism $\pi_s^*(\co_{\oX'_n})\rightarrow \ocB^!_n$ of $\tG_s$ (\cite{ag1} 6.5.28).
The morphism $\pi_s$ \eqref{p2-rftchta2b} is therefore underlying a morphism of ringed topos, that we denote by
\begin{equation}\label{p2-rftchta3h}
\pi_n\colon (\tG_s,\ocB^!_n)\rightarrow (\oX'_{s,\et},\co_{\oX'_n}).
\end{equation}

The canonical homomorphism $\varpi^*(\ocB^!)\rightarrow \ocB'$ \eqref{p2-rftchta3b}
induces a homomorphism $\varpi_s^*(\ocB^!_n)\rightarrow \ocB'_n$ of $\tE'_s$.
The morphism $\varpi_s$ \eqref{p2-rftchta2d} is therefore underlying a morphism of ringed topos, that we denote by
\begin{equation}\label{p2-rftchta3i}
\varpi_n\colon (\tE'_s,\ocB'_n)\rightarrow (\tG_s,\ocB^!_n).
\end{equation}

The canonical homomorphism $\lgg^{-1}(\ocB)\rightarrow \ocB^!$ \eqref{p2-rftchta3c}
induces a homomorphism $\lgg_s^*(\ocB_n)\rightarrow \ocB^!_n$ of $\tG_s$.
The morphism $\lgg_s$ \eqref{p2-rftchta2h} is therefore underlying a morphism of ringed topos, that we denote by
\begin{equation}\label{p2-rftchta3j}
\lgg_n\colon (\tG_s,\ocB^!_n)\rightarrow (\tE_s,\ocB_n).
\end{equation}

By (\cite{ag1} (6.5.28.10)), the triangles and the square of the diagram of morphisms of ringed topos
\begin{equation}\label{p2-rftchta3k}
\xymatrix{
{(\tE'_s,\ocB'_n)}\ar[r]^-(0.5){\varpi_n}\ar[rd]_{\sigma'_n}\ar@/^2pc/[rr]^{\uptheta_n}&{(\tG_s,\ocB_n^!) }\ar[d]^-(0.5){\pi_n}\ar[r]^-(0.5){\lgg_n}&
{(\tE_s,\ocB_n)}\ar[d]^-(0.5){ \sigma_n}\\
&{(X'_{s,\et},\co_{\oX'_n})}\ar[r]^{\ogg_n}&{(X_{s,\et},\co_{\oX_n} )}}
\end{equation}
are commutative up to canonical isomorphisms (\cite{egr1} 1.2.3).

Let $\ur=(r_1,r_2,r_3)\in I$ \eqref{p2-fhtft33}. 
The homomorphism $\uppsi^{\ur}_n$ \eqref{p2-fhtft38e} induces a homomorphism of $\pi^*_n(\hcC^{(r_3)}_{\uptau,n})$-algebras 
\begin{equation}\label{p2-rftchta3l}
\lgg^*_n(\cC^{(r_2)}_n)\otimes_{\ocB^!_n}\pi^*_n(\cC^{(r_3)}_{\uptau,n})\rightarrow 
\varpi_{n*}(\cC'^{\ur}_n).
\end{equation}

\subsection{}\label{p2-rftchta4}
We take again the notation of \ref{p2-htaft8} and \ref{p2-fhtft22}. We denote by $\bvocB^!$ the ring $(\ocB^!_{n+1})_{n\in \mN}$ of $\tG_s^{\mN^\circ}$ \eqref{p2-rftchta3f} and by
\begin{equation}\label{p2-rftchta4a}
\xymatrix{
{(\tE'^{\mN^\circ}_s,\bvocB')}\ar[r]^-(0.5){\bvvarpi}\ar[rd]_{\bvsigma'}\ar@/^2pc/[rr]^{\bvuptheta}&
{(\tG^{\mN^\circ}_s,\bvocB^!)}\ar[d]^-(0.5){\bvpi}\ar[r]^-(0.5){\bvlgg}& {(\tE^{\mN^\circ}_s,\bvocB)}
\ar[d]^{\bvsigma}\\
&{(X'^ {\mN^\circ}_{s,\et},\co_{\bvoX'})}\ar[r]^-(0.5){\bvogg}&{(X^ {\mN^\circ}_{s,\et},\co_{\bvoX})}}
\end{equation}
the morphisms of ringed topos induced by $(\uptheta_{n+1})_{n\in \mN}$ \eqref{p2-fhtft21h}, $(\varpi_{n+1})_{n\in \mN}$ \eqref{p2-rftchta3i},
$(\lgg_{n+1})_{n\in \mN}$ \eqref{p2-rftchta3j}, $(\pi_{n+1})_{n\in \mN}$ \eqref{p2-rftchta3h}, 
$(\sigma'_{n+1})_{n\in \mN}$, $(\sigma_{n+1})_{n\in \mN}$ \eqref{p2-htaft7c} 
and $(\ogg_{n+1})_{n\in \mN}$ \eqref{p2-fhtft1a}.
By \eqref{p2-rftchta3k}, the diagram of morphisms of ringed topos \eqref{p2-rftchta4a} is commutative up to canonical isomorphism (\cite{egr1} 1.2.3).
We denote by
\begin{equation}\label{p2-rftchta4b}
\bvu'\colon (X'^{\mN^\circ}_{s,\et},\co_{\bvoX'})\rightarrow (X'^{\mN^\circ}_{s,\zar},\co_{\bvoX'})
\end{equation}
the canonical morphism of ringed topos and by
\begin{equation}\label{p2-rftchta4c}
\uplambda'\colon (X'^{\mN^\circ}_{s,\zar},\co_{\bvoX'})\rightarrow (X'_{s,\zar},\co_{\fX'})
\end{equation}
the morphism of ringed topos for which the functor $\uplambda'_*$ is the inverse limit functor \eqref{p2-ncgt5a}. We set
\begin{equation}\label{p2-rftchta4d}
\huppi=\uplambda'\circ \bvu' \circ \bvpi\colon (\tG_s^{\mN^\circ},\bvocB^!)\rightarrow (X'_{s,\zar},\co_{\fX'}). 
\end{equation}

By \eqref{p2-rftchta4a}, the diagram of morphisms of ringed topos
\begin{equation}\label{p2-rftchta4e}
\xymatrix{
{(\tE'^{\mN^\circ}_s,\bvocB')}\ar[r]^-(0.5){\bvvarpi}\ar[rd]_{\hupsigma'}\ar@/^2pc/[rr]^{\bvuptheta}&{(\tG^{\mN^\circ}_s,\bvocB^!)}
\ar[d]^-(0.5){\huppi}\ar[r]^-(0.5){\bvlgg}&{(\tE^{\mN^\circ}_s,\bvocB)}\ar[d]^-(0.5){\hupsigma}\\
&{(X'_{s,\zar},\co_{\fX'})}\ar[r]^{\fgg}&{(X_{s,\zar},\co_{\fX})}}
\end{equation}
is commutative up to canonical isomorphism (\cite{egr1} 1.2.3). It extends the commutative diagram \eqref{p2-fhtft22b}. 

Let $\ur=(r_1,r_2,r_3)\in I$ \eqref{p2-fhtft33}. The homomorphism $\bvuppsi^{\ur}$ \eqref{p2-fhtft38g} induces a homomorphism of $\huppi^*(\hcC^{(r_3)}_\uptau)$-algebras 
\begin{equation}\label{p2-rftchta4f}
\bvlgg^*(\bvcC^{(r_2)})\otimes_{\bvocB^!}\huppi^*(\hcC^{(r_3)}_\uptau)\rightarrow 
\bvvarpi_*(\bvcC'^{\ur}).
\end{equation}

\subsection{}\label{p2-rftchta14}
We denote by $\bMod(\bvocB^!)$ the category of $\bvocB^!$-modules of $\tG_s^{\mN^\circ}$,
by $\bIndMod(\bvocB^!)$ the category of ind-$\bvocB^!$-modules and by
\begin{equation}\label{p2-rftchta14a}
\iota_{\bvocB^!}\colon \bMod(\bvocB^!)\rightarrow \bIndMod(\bvocB^!)
\end{equation}
the canonical fully faithful and exact functor \eqref{p1-abisoind1b}.
We denote by $\bMod_\mQ(\bvocB^!)$ the category of $\bvocB^!$-modules up to isogeny \eqref{p1-abisoind1} and by 
\begin{equation}\label{p2-rftchta14b}
\upalpha_{\bvocB^!}\colon \bMod_\mQ(\bvocB^!)\rightarrow \bIndMod(\bvocB^!).
\end{equation}
the canonical fully faithful and exact functor \eqref{p1-bcim5c}.
We will identify $\bMod(\bvocB^!)$ (resp.\ $\bMod_\mQ(\bvocB^!)$) with a full subcategory of $\bIndMod(\bvocB^!)$
by the functor $\iota_{\bvocB^!}$ (resp.\ $\upalpha_{\bvocB^!}$), which we will omit from the notation. 

By \ref{p1-bcim1}, the morphism $\bvvarpi$ \eqref{p2-rftchta4a} induces two adjoint additive functors
\begin{eqnarray}
\rI \bvvarpi^*\colon \bIndMod(\bvocB^!) \rightarrow \bIndMod(\bvocB'),\label{p2-rftchta14c}\\
\rI \bvvarpi_*\colon \bIndMod(\bvocB') \rightarrow \bIndMod(\bvocB^!).\label{p2-rftchta14d}
\end{eqnarray}
The functor $\rI \bvvarpi^*$ (resp.\ $\rI \varpi_*$) is right (resp.\ left) exact.
The functor $\rI \bvvarpi_*$ admits a right derived functor
\begin{equation}\label{p2-rftchta14e}
\rR\rI \bvvarpi_*\colon \bD^+(\bIndMod(\bvocB'))\rightarrow \bD^+(\bIndMod(\bvocB^!)).
\end{equation}

\subsection{}\label{p2-rftchta10}
For local cohomological computations, it is convenient to introduce the following notation.
Let $(\oy'\rightsquigarrow \ox')$ be a point of the topos $X'_\et\gtimes_{X'_\et}\oX'^\rhd_\et$ (\cite{ag1} 3.4.23)
such that $\ox'$ is above $s$, $\uX'$ the strict localization of $X'$ at $\ox'$.
By (\cite{agt} III.3.7), $\uoX'$ is normal and strictly local (and in particular integral);
it can therefore be identified with the strict localization of $\oX'$ at $\oa'(\ox')$ \eqref{p2-htaft2}.
The $X'$-morphism $\oy'\rightarrow \uX'$
defining $(\oy'\rightsquigarrow \ox')$ lifts to an $\oX'^\rhd$-morphism $v'\colon \oy'\rightarrow \uoX'^\rhd$ and
therefore induces a geometric point of $\uoX'^\rhd$ that we (abusively) denote also by $\oy'$.
We set $\uDelta'=\pi_1(\uoX'^\rhd,\oy')$.

Let $\ox=g(\ox')$ and $\oy=\upgamma(\oy')$ \eqref{p2-fhtft1e}, which are geometric points of $X$ and
$\oX^\circ$ respectively, and let $(\oy\rightsquigarrow \ox')$ (resp.\ $(\oy\rightsquigarrow \ox)$) be the image of $(\oy'\rightsquigarrow \ox' )$
by the first (resp.\ the composition) of the canonical morphisms
\begin{equation}\label{p2-rftchta10a}
X'_\et\gtimes_{X'_\et}\oX'^\rhd_\et\rightarrow X'_\et\gtimes_{X_\et}\oX^\circ_\et
\rightarrow X_\et\gtimes_{X_\et}\oX^\circ_\et.
\end{equation}
We denote by $\uX$ the strict localization of $X$ at $\ox$.
By (\cite{agt} III.3.7), $\uoX$ is normal and strictly local (and in particular integral);
it can therefore be identified with the strict localization of $\oX$ at $\oa(\ox)$ \eqref{p2-htaft2}.
The $X$-morphism $\oy\rightarrow \uX$
defining $(\oy\rightsquigarrow \ox)$ lifts to an $\oX^\circ$-morphism $v\colon \oy\rightarrow \uoX^\circ$ and
therefore induces a geometric point of $\uoX^\circ$ which we (abusively) denote also by $\oy$.
We set $\uDelta=\pi_1(\uoX^\circ,\oy)$. We denote by
\begin{equation}\label{p2-rftchta10b}
\uupgamma\colon \uoX'^\rhd\rightarrow \uoX^\circ
\end{equation}
the morphism induced by $g$. Since $\uupgamma(\oy')=\oy$, this induces a homomorphism
\begin{equation}\label{p2-rftchta10c}
\uDelta'\rightarrow \uDelta.
\end{equation}
We denote by $\uPi$ its kernel.
We denote by $\bB_{\uDelta'}$ (resp.\ $\bB_{\uDelta}$) the classifying topos of the profinite group $\uDelta'$ (resp.\ $\uDelta$) and by
\begin{eqnarray}
\psi_{\uX',\oy'}\colon \uoX'^\rhd_\fet&\stackrel{\sim}{\rightarrow}&\bB_{\uDelta'},\label{p2-rftchta10f1}\\
\psi_{\uX,\oy}\colon \uoX^\circ_\fet&\stackrel{\sim}{\rightarrow}&\bB_{\uDelta},\label{p2-rftchta10f2}
\end{eqnarray}
the fiber functors (\cite{agt} (VI.9.8.4)).

Recall that giving a neighborhood of the point of $X_\et$ associated with $\ox$
in the site $\Et_{/X}$ is equivalent to giving an $\ox$-pointed étale $X$-scheme (\cite{sga4} IV 6.8.2).
These objects naturally form a cofiltered category, that we denote by $\fV_\ox$. 
Similarly, we denote by $\fV'_{\ox'}$ the category of $\ox'$-pointed étale $X'$-schemes.

Let $(\mu\colon U'\rightarrow U, \iota\colon \ox'\rightarrow U')$ be an object of the category $\mM_{\ox'}$ associated with $\ox'$ in \ref{p2-fhtft301}; 
we may omit $\iota$ to lighten the notation.
We have a canonical $X'$-morphism
$\uX'\rightarrow U'$ and a canonical $X$-morphism $\uX\rightarrow U$ that fit into a commutative diagram
\begin{equation}\label{p2-rftchta10k}
\xymatrix{
&\uX'\ar[r]^{\ug}\ar[d]&\uX\ar[d]\\
\ox'\ar[r]^\iota\ar[ru]&U'\ar[r]^\mu&U,}
\end{equation}
where $\ug$ is the morphism induced by $g$. We deduce from this morphisms $\uoX'\rightarrow \oU'$ and $\uoX\rightarrow \oU$.
The morphism $v'\colon \oy'\rightarrow \uoX'^\rhd$ induces a geometric point of $\oU'^\rhd$ which we also denote by $\oy'$.
Similarly, the morphism $v\colon \oy\rightarrow \uoX^\circ$
induces a geometric point of $\oU^\circ$ which we also denote by $\oy$. Observe that the diagram
\begin{equation}\label{p2-rftchta10l}
\xymatrix{
\oy'\ar[r]\ar[d]&\oy\ar[d]\\
\oU'^\rhd\ar[r]^\tmu&\oU^\circ}
\end{equation}
is commutative \eqref{p2-fhtft302a}.

The schemes $\oU$ and $\oU'$ being locally irreducible by (\cite{agt} III.3.3 and III.4.2(iii)),
they are the sums of the schemes induced on their irreducible components.
We denote by $\oU^\star$ (resp.\ $\oU'^\star$)
the irreducible component of $\oU$ (resp.\ $\oU'$) containing $\oy$ (resp.\ $\oy'$).
Similarly, $\oU^\circ$ (resp.\ $\oU'^\rhd$) is the sum of the schemes induced on its irreducible components
and $\oU^{\star \circ}=\oU^\star\times_{X}X^\circ$ (resp.\ $\oU'^{\star \rhd}=\oU'^\star\times_ {X'}X'^\rhd$)
is the irreducible component of $\oU^\circ$ (resp.\ $\oU'^\rhd$) containing $\oy$ (resp.\ $\oy'$).
We set $\Delta_U=\pi_1(\oU^{\star\circ},\oy)$ and $\Delta'_{U'}=\pi_1(\oU'^{\star\rhd},\oy')$.
By \eqref{p2-rftchta10l}, we have a canonical homomorphism $\Delta'_{U'}\rightarrow \Delta_U$. We denote by $\Pi_{U'\rightarrow U}$ its kernel.
We denote by $\bB_{\Delta'_{U'}}$ (resp.\ $\bB_{\Delta_U}$) the classifying topos of the profinite group $\Delta'_{U'}$ (resp.\ $\Delta_U$) and by
\begin{eqnarray}
\psi_{U',\oy'}\colon \oU'^{\star \rhd}_\fet&\stackrel{\sim}{\rightarrow}&\bB_{\Delta'_{U'}},\label{p2-rftchta103f}\\
\psi_{U,\oy}\colon \oU^{\star \circ}_\fet&\stackrel{\sim}{\rightarrow}&\bB_{\Delta_U},\label{p2-rftchta10f4}
\end{eqnarray}
the fiber functors.

The functors
\begin{eqnarray}
\mM_{\ox'}&\rightarrow&\fV'_{\ox'},\ \ (\mu\colon U'\rightarrow U, \iota\colon \ox'\rightarrow U') \mapsto \iota,\label{p2-rftchta10g1}\\
\mM_{\ox'}&\rightarrow&\fV_{\ox},\ \ \ (\mu\colon U'\rightarrow U, \iota\colon \ox'\rightarrow U') \mapsto \mu \circ \iota ,\label{p2-rftchta10g2}
\end{eqnarray}
where the second is defined by the fact that $\ox=g(\ox')$, are initial by virtue of (\cite{sga4} I 8.1.3(b)).
Consequently, the canonical morphisms
\begin{eqnarray}
\uoX'&\rightarrow& \underset{\underset{(U'\rightarrow U)\in \mM_{\ox'}}{\longleftarrow}}{\lim}\ \oU'^\star,\\
\uoX&\rightarrow& \underset{\underset{(U'\rightarrow U)\in \mM_{\ox'}}{\longleftarrow}}{\lim}\ \oU^\star,
\end{eqnarray}
are isomorphisms. By (\cite{agt} VI.11.8), the canonical morphisms
\begin{eqnarray}
\uDelta'&\rightarrow& \underset{\underset{(U'\rightarrow U)\in \mM_{\ox'}}{\longleftarrow}}{\lim}\ \Delta'_{U'},\\
\uDelta&\rightarrow& \underset{\underset{(U'\rightarrow U)\in \mM_{\ox'}}{\longleftarrow}}{\lim}\ \Delta_U
\end{eqnarray}
are therefore isomorphisms. We deduce from this that the canonical morphism
\begin{equation}\label{p2-rftchta10m}
\uPi\rightarrow \underset{\underset{(U'\rightarrow U)\in \mM_{\ox'}}{\longleftarrow}}{\lim}\ \Pi_{U'\rightarrow U}
\end{equation}
is an isomorphism.

We denote by
\begin{eqnarray}
\varphi'_{\ox'}\colon \tE'&\rightarrow& \uoX'^\rhd_\fet,\\
\phi_{\ox'}\colon \tG&\rightarrow& \uoX^\circ_\fet,\\
\varphi_{\ox}\colon \tE&\rightarrow& \uoX^\circ_\fet,
\end{eqnarray}
the canonical functors defined in \eqref{p2-fhtft39i}, \eqref{p2-rftchta9c} and \eqref{p2-fhtft39ii}, respectively. 
By (\cite{ag1} 6.6.8), for every abelian group $F$ of $\tE'$ and every integer $q\geq 0$, we have a canonical isomorphism
\begin{equation}\label{p2-rftchta10h}
(\rR^q\varpi_*(F))_{\varrho(\oy\rightsquigarrow \ox')}\stackrel{\sim}{\rightarrow} \rH^q(\uPi,\psi'_{\uX',\oy'}(\varphi'_{\ox'}(F))),
\end{equation}
where $\varrho$ is the morphism defined in \eqref{p2-rftchta1i}.

\subsection{}\label{p2-rftchta11}
We keep the assumptions and notation of \ref{p2-rftchta10}.
For any object $(\mu\colon U'\rightarrow U, \iota\colon \ox'\rightarrow U')$ of $\mM_{\ox'}$ \eqref{p2-fhtft301}, we set
\begin{eqnarray}
\oR'^{\oy'}_{U'}&=&\psi'_{U',\oy'}(\ocB'_{U'}|\oU'^{\star \rhd}) ,\label{p2-rftchta11a}\\
\oR^{\oy}_{U}&=&\psi_{U,\oy}(\ocB_{U}|\oU^{\star \circ}),\label{p2-rftchta11b}\\
\oR^{!\oy}_{U'\rightarrow U}&=&\psi_{U,\oy}(\ocB^!_{U'\rightarrow U}|\oU^{\star \circ} ),\label{p2-rftchta11c}
\end{eqnarray}
where $\ocB_U$ and $\ocB'_{U'}$ (resp.\ $\ocB^!_{U'\rightarrow U}$) are the sheaves defined in \eqref{p2-htaft6b} (resp.\ \eqref{p2-rftchta12a}).

Explicitly, let $(V_i)_{i\in I}$ be the normalized universal cover
 of $\oU^{\star \circ}$ at $\oy$ (\cite{agt} VI.9.8).
For every $i\in I$, $(V_i\rightarrow U)$ (resp.\ $(U'\rightarrow U\leftarrow V_i)$) is naturally an object of $E$ (resp.\ $G$). We then have
\begin{eqnarray}
\oR^{\oy}_{U}&=&\underset{\underset{i\in I}{\longrightarrow}}{\lim}\ \ocB(V_i\rightarrow U),\label{p2-rftchta11d}\\
\oR^{!\oy}_{U'\rightarrow U}&=&\underset{\underset{i\in I}{\longrightarrow}}{\lim}\ \ocB^!(U'\rightarrow U \leftarrow V_i).\label{p2-rftchta11e}
\end{eqnarray}

The canonical homomorphism $\lgg^*(\ocB)\rightarrow \ocB^!$ \eqref{p2-rftchta3c} induces for every $i\in I$ a morphism (functorial in $i$)
\begin{equation}
\ocB(V_i\rightarrow U)\rightarrow \ocB^!(U'\rightarrow U\leftarrow V_i).
\end{equation}
We deduce from this, by taking the direct limit,  a homomorphism
\begin{equation}\label{p2-rftchta11f}
\oR^{\oy}_{U}\rightarrow \oR^{!\oy}_{U'\rightarrow U}.
\end{equation}

Let $(W_j)_{j\in J}$ be the normalized universal cover of $\oU'^{\star \rhd}$ at $\oy'$.
For every $j\in J$, $(W_j\rightarrow U')$ is naturally an object of $E'$. We then have
\begin{equation}\label{p2-rftchta11g}
\oR'^{\oy'}_{U'}=\underset{\underset{j\in J}{\longrightarrow}}{\lim}\ \ocB'(W_j\rightarrow U').
\end{equation}

For every $i\in I$, we have a canonical $\oU^\circ$-morphism $\oy\rightarrow V_i$.
We deduce from this a $\oU'^\rhd$-morphism $\oy'\rightarrow V_i\times_{\oU^\circ}\oU'^\rhd$.
The scheme $V_i\times_{\oU^\circ}\oU'^\rhd$ being locally irreducible,
it is the sum of the schemes induced on its irreducible components.
We denote by $V'_i$ the irreducible component of $V_i\times_{\oU^\circ}\oU'^\rhd$ containing the image of $\oy'$.
The schemes $(V'_i)_{i\in I}$ naturally form an inverse system of connected finite étale $\oy'$-pointed  covers of $\oU'^{\star \rhd}$.
The canonical morphism
\begin{equation}
\oU'^\rhd\times_{\oU^\circ}V_i\rightarrow U'\times_{(U\times_XX')}(V_i\times_{X^\circ}X'^\rhd)
\end{equation}
is an isomorphism. We therefore have a canonical isomorphism of $\tE'$
\begin{equation}
\varpi^*((U'\rightarrow U\leftarrow V_i)^a)\stackrel{\sim}{\rightarrow} (V_i\times_{\oU^\circ}\oU'^\rhd\rightarrow U') ^a.
\end{equation}
Since the canonical homomorphism $\ocB^!\rightarrow\varpi_*(\ocB')$ \eqref{p2-rftchta3b} is an isomorphism,
we deduce from this a canonical homomorphism (functorial in $i$)
\begin{equation}
\ocB^!(U'\rightarrow U\leftarrow V_i)=\ocB'(V_i\times_{\oU^\circ}\oU'^\rhd\rightarrow U')\rightarrow \ocB'(V'_i\rightarrow U').
\end{equation}
Setting
\begin{equation}\label{p2-rftchta11h}
\oR^{\intern\oy'}_{U'\rightarrow U}=\underset{\underset{i\in I}{\longrightarrow}}{\lim}\ \ocB'(V'_i\rightarrow U '),
\end{equation}
we therefore have a canonical homomorphism
\begin{equation}
\oR^{!\oy}_{U'\rightarrow U}\rightarrow \oR^{\intern\oy'}_{U'\rightarrow U}.
\end{equation}

For all $i\in I$ and $j\in J$, there exists at most one morphism of pointed $\oU'^{\star \rhd}$-schemes $W_j\rightarrow V'_i$.
Moreover, for every $i\in I$, there exist $j\in J$ and a morphism of pointed $\oU'^{\star \rhd}$-schemes $W_j\rightarrow V'_i$.
Then, we have a homomorphism
\begin{equation}\label{p2-rftchta11i}
\oR^{\intern\oy'}_{U'\rightarrow U}=\underset{\underset{i\in I}{\longrightarrow}}{\lim}\ \ocB'(V'_i\rightarrow U ')\rightarrow
\underset{\underset{j\in J}{\longrightarrow}}{\lim}\ \ocB'(W_j\rightarrow U')=\oR'^{\oy'}_{U'}.
\end{equation}

We therefore have three canonical $\Delta'_{U'}$-equivariant homomorphisms
\begin{equation}\label{p2-rftchta11j}
\oR^\oy_U\rightarrow \oR^{!\oy}_{U'\rightarrow U} \rightarrow \oR^{\intern\oy'}_{U'\rightarrow U}\rightarrow \oR'^{ \oy'}_{U'}.
\end{equation}
These rings and these homomorphisms being functorial in $(\mu\colon U'\rightarrow U, \iota\colon \ox'\rightarrow U')\in \ob(\mM_{\ox'})$, let
{\allowdisplaybreaks
\begin{eqnarray}
\oR'^{\oy'}_{\uX'}&=&\underset{\underset{(U'\rightarrow U)\in \ob(\mM^\circ_{\ox' })}{\longrightarrow}}{\lim}\ \oR'^{\oy'}_{U'},\label{p2-rftchta11k1}\\
\oR^{\intern\oy}_{\uX'\rightarrow \uX}&=&\underset{\underset{(U'\rightarrow U)\in \ob(\mM^\circ_ {\ox'})}{\longrightarrow}}{\lim}
\oR^{\intern\oy'}_{U'\rightarrow U},\label{p2-rftchta11k2}\\
\oR^{!\oy}_{\uX'\rightarrow \uX}&=&\underset{\underset{(U'\rightarrow U)\in \ob(\mM^\circ_{ \ox'})}{\longrightarrow}}{\lim}
\oR^{!\oy}_{U'\rightarrow U},\label{p2-rftchta11k3}\\
\oR^{\oy}_{\uX}&=&\underset{\underset{(U'\rightarrow U)\in \ob(\mM^\circ_{\ox'})} {\longrightarrow}}{\lim}\ \oR^\oy_U,\label{p2-rftchta11k4}
\end{eqnarray}
}%
of which the first two are rings of $\bB_{\uDelta'}$ and the other two are rings of $\bB_{\uDelta}$.

By (\cite{agt} VI.10.37, \cite{ag2} (4.3.9.4)), as the functors \eqref{p2-rftchta10g1} and \eqref{p2-rftchta10g2} are initial, we have canonical isomorphisms
\begin{eqnarray}
\psi_{\uX',\oy'}(\varphi'_{\ox'}(\ocB'))&\stackrel{\sim}{\rightarrow} &\oR'^{\oy'}_{ \uX'},\label{p2-rftchta11l3}\\
\psi_{\uX,\oy}(\varphi_{\ox}(\ocB))&\stackrel{\sim}{\rightarrow} &\oR^{\oy}_{\uX}.\label{p2-rftchta11l4 }
\end{eqnarray}
These induce canonical isomorphisms 
\begin{eqnarray}
\ocB'_{\rho'(\oy'\rightsquigarrow \ox')}&\stackrel{\sim}{\rightarrow} &\oR'^{\oy'}_{\uX'},\label{p2- rftchta11l1}\\
\ocB_{\rho(\oy\rightsquigarrow \ox)}&\stackrel{\sim}{\rightarrow} &\oR^{\oy}_{\uX}.\label{p2-rftchta11l2}
\end{eqnarray}
Moreover, by (\cite{ag2} (6.1.17.3)), we have a canonical isomorphism
\begin{equation}\label{p2-rftchta11m}
\ocB^!_{\varrho(\oy\rightsquigarrow \ox')}\stackrel{\sim}{\rightarrow} \oR^{!\oy}_{\uX'\rightarrow \uX}.
\end{equation}

The homomorphisms \eqref{p2-rftchta11j} induce by taking the direct limit canonical $\uDelta'$-equivariant homomorphisms
\begin{equation}\label{p2-rftchta11n}
\oR^\oy_\uX\rightarrow \oR^{!\oy}_{\uX'\rightarrow \uX} \rightarrow \oR^{\intern\oy'}_{\uX'\rightarrow \uX}\rightarrow \oR'^{\oy'}_{\uX'}.
\end{equation}
By (\cite{ag1} 6.5.26), the homomorphism in the middle is an isomorphism.
By the above, the first homomorphism and the composition of the other two are identified with the homomorphisms
\begin{equation}\label{p2-rftchta11o}
\ocB_{\rho(\oy\rightsquigarrow \ox)}\rightarrow \ocB^!_{\varrho(\oy\rightsquigarrow \ox')}\rightarrow \ocB'_{\rho'(\oy'\rightsquigarrow \ox')}
\end{equation}
induced by the adjoints of the homomorphisms \eqref{p2-rftchta3c} and \eqref{p2-rftchta3b}, respectively.

\begin{rema}\label{p2-rftchta110}
Under the assumptions of \ref{p2-rftchta11}, for every object $(\mu\colon U'\rightarrow U, \iota\colon \ox'\rightarrow U')$ of $\mM_{\ox'}$ such that the restriction
$(U',\cM_{X'}|U')\rightarrow (U,\cM_X|U)$ of $g$ satisfies the assumptions of \ref{p2-fhtal1}, the rings
$\oR^\oy_U\rightarrow \oR^{\intern\oy'}_{U'\rightarrow U}\rightarrow \oR'^{\oy'}_{U'}$ coincide with the rings
$\oR\rightarrow \oR^{\intern}\rightarrow \oR'$ defined in \ref{p2-fhtal2}.
\end{rema}

\subsection{}\label{p2-rftchta13}
We keep the assumptions and notation of \ref{p2-rftchta10} and \ref{p2-rftchta11}.
Recall that $\bQ$ is a topologically generating subcategory of $\Et_{/X}$ \eqref{p2-htaft10}, that we endowed with the topology induced by the étale topology of $\Et_{/X}$.
We denote by $\bQ_\ox$ the category of neighborhoods of the point of $X_\et$ associated with $\ox$ in the site $\bQ$, or equivalently 
the category of $\ox$-pointed étale $X$-schemes of $\bQ$ (\cite{sga4} IV 6.8.2).
It is a cofiltered $\mU$-small category. Moreover, the canonical functor $\bQ\rightarrow \Et_{/X}$ induces a fully
faithful and initial functor $\bQ_\ox\rightarrow \fV_\ox$ by (\cite{sga4} I 8.1.3(c)) and (\cite{agt} II.5.17).
Similarly, for $(f',\ox')$, we denote by $\bQ'_{\ox'}$ the category of $\ox'$-pointed étale $X'$-schemes of $\bQ'$. 

With the notation of \ref{p2-fhtft301}, the functors
\begin{eqnarray}
\mM_{\ox'}(\bQ',\bQ)\rightarrow \bQ'_{\ox'},&& (\mu\colon U'\rightarrow U, \iota\colon \ox'\rightarrow U')\mapsto \iota,\label{p2-rftchta13a}\\
\mM_{\ox'}(\bQ',\bQ)\rightarrow \bQ_\ox,&& (\mu\colon U'\rightarrow U, \iota\colon \ox'\rightarrow U')\mapsto \mu\circ \iota,\label{p2-rftchta13b}
\end{eqnarray}
where the second functor is defined by the fact that $\ox=g(\ox')$, are initial by (\cite{sga4} I 8.1.3(b)) and (\cite{agt} II.5.17).

Let $r,r'$ be rational numbers such that $r\geq r'\geq 0$, $n$ an integer $\geq 0$. 
By (\cite{ag2} 4.4.20), with the notation of \ref{p2-htaft12}, we have a canonical isomorphism
\begin{equation}\label{p2-rftchta13c}
\psi_{\uX,\oy}(\varphi_{\ox}(\cC^{(r)}_n)) \stackrel{\sim}{\rightarrow} \underset{\underset{U\in \bQ^\circ_\ox}{\longrightarrow}}{\lim}\
\fC^{\oy,(r)}_{U}/p^n\fC^{\oy,(r)}_{U}.
\end{equation}
Similarly, with the notation of \ref{p2-fhtft24}, we have a canonical isomorphism
\begin{equation}\label{p2-rftchta13d}
\psi_{\uX',\oy'}(\varphi'_{\ox'}(\cC'^{(r,r')}_n)) \stackrel{\sim}{\rightarrow} \underset{\underset{U'\in \bQ'^\circ_{\ox'}}{\longrightarrow}}{\lim}\
\fC'^{\oy',(r,r')}_{U'}/p^n\fC'^{\oy',(r,r')}_{U'}.
\end{equation}

For any $\ur=(r_1,r_2,r_3)\in I$ \eqref{p2-fhtft33} and any object $(\mu\colon U'\rightarrow U, \iota\colon \ox'\rightarrow U')$ of $\mM_{\ox'}(\bQ',\bQ)$, 
we defined a canonical $\Delta'_{U'}$-equivariant homomorphism of $(\hoRp^{\oy'}_{U'}\otimes_{R'_{\uptau,U'}}\fC^{(r_3)}_{\uptau,U'})$-algebras \eqref{p2-fhtft33c}
\begin{equation}\label{p2-rftchta13e}
\uppsi^{\oy',\ur}_\mu\colon \fC^{\oy,(r_2)}_U\otimes_{\hoR^\oy_U}\hoRp^{\oy'}_{U'}\otimes_{R'_{\uptau,U'}}\fC^{(r_3)}_{\uptau,U'}
\rightarrow \fC'^{\oy,(r_1,r_2)}_{U'} \otimes_{R'_{\uptau,U'}}\fC^{(r_3)}_{\uptau,U'}.
\end{equation}

The image by the composed functor $\psi_{\uX',\oy'}\circ \varphi'_{\ox'}$ of the canonical homomorphism \eqref{p2-fhtft38e}
\begin{equation}\label{p2-rftchta13f}
\uptheta^*_n(\cC^{(r_2)}_n)\otimes_{\ocB'_n}\sigma'^*_n(\cC^{(r_3)}_{\uptau,n})\rightarrow \cC'^{\ur}_n
\end{equation}
identifies with the direct limit of the homomorphisms
$\uppsi^{\oy',\ur}_\mu$ modulo $p^n$ \eqref{p2-rftchta13e}, when $(\mu\colon U'\rightarrow U, \iota\colon \ox'\rightarrow U')$  
goes through the category $\mM^\circ_{\ox'}(\bQ',\bQ)$ (resp. $\umM^\circ_{\ox'}(\bQ',\bQ)$ \eqref{p2-fhtft301}). 
This follows from \ref{p2-qfc263}, in view of \ref{p2-cmt5}, \eqref{p2-qfc262g} and (\cite{agt} VI.10.37 and (VI.10.18.1)). 
Recall that the canonical functor $\umM^\circ_{\ox'}(\bQ',\bQ)\rightarrow \mM^\circ_{\ox'}(\bQ',\bQ)$ is cofinal \eqref{p2-fhtft301}.

\begin{prop}\label{p2-rftchta5}
Let $\ur=(r_1,r_2,r_3)$, $\ur'=(r'_1,r'_2,r'_3)$ be two elements of $I$ \eqref{p2-fhtft33} such that $r_1> r'_1$, $r_2\geq  r'_2$ and $r_3\geq  r'_3$. Then,
\begin{itemize}
\item[{\rm (i)}] For every integer $n\geq 1$, the canonical homomorphism \eqref{p2-rftchta3l}
\begin{equation}\label{p2-rftchta5a}
\lgg^*_n(\cC^{(r_2)}_n)\otimes_{\ocB^!_n}\pi^*_n(\cC^{(r_3)}_{\uptau,n})\rightarrow 
\varpi_{n*}(\cC'^{\ur}_n)
\end{equation}
is $\alpha$-injective. Let $\cH^{\ur}_n$ be its cokernel.
\item[{\rm (ii)}] There exists an integer $a\geq 0$, depending on $r_1$, $r'_1$ and $\ell=\dim(X'/X)$,
but not on the morphisms $f$ and $g$ satisfying the conditions of \ref{p2-fhtft1}, such that for every integer $n\geq 1$,
the morphism $\cH^{\ur}_n\rightarrow \cH^{\ur'}_n$ induced by \eqref{p2-fhtft38i} is annihilated by $p^a$.
\item[{\rm (iii)}] There exists an integer $b\geq 0$, depending on $r_1$, $r'_1$ and $\ell=\dim(X'/X)$,
but not on the morphisms $f$ and $g$ satisfying the conditions of \ref{p2-fhtft1}, such that for all integers $n,q\geq 1$,
the canonical morphism
\begin{equation}\label{p2-rftchta5b}
\rR^q\varpi_{n*}(\cC'^{\ur}_n)\rightarrow 
\rR^q\varpi_{n*}(\cC'^{\ur'}_n)
\end{equation}
is annihilated by $p^b$.
\end{itemize}
\end{prop}

Let $(\oy\rightsquigarrow \ox')$ be a point of $X'_\et\gtimes_{X_\et}\oX^\circ_\et$ such that $\ox'$ is above $s $,
$\uX'$ the strict localization of $X'$ at $\ox'$, $\ox=g(\ox')$, $\uX$ the strict localization of $X$ at $\ox$.
The morphism $\uupgamma\colon \uoX'^\rhd\rightarrow \uoX^\circ$ induced by $g$ being faithfully flat by (\cite{ag1} 2.4.1),
there exists a point $(\oy'\rightsquigarrow \ox')$ of $X'_\et\gtimes_{X'_\et}\oX'^\rhd_\et$
lifting the point $(\oy\rightsquigarrow \ox')$ of $X'_\et\gtimes_{X_\et}\oX^\circ_\et$ \eqref{p2-rftchta10a}.
We take again the notation introduced in \ref{p2-rftchta10}, \ref{p2-rftchta11} and \ref{p2-rftchta13}.

By \eqref{p2-rftchta2l} and \eqref{p2-rftchta10h}, for all integers $n,q\geq 0$,
we have a canonical isomorphism
\begin{equation}
\rR^q\varpi_{s*}(\cC'^{(r_1,r_2)}_n\otimes_{\ocB'_n}\sigma'^*_n(\cC^{(r_3)}_{\uptau,n}))_{\varrho(\oy\rightsquigarrow \ox')}
\stackrel{\sim}{\rightarrow}\rH^q(\uPi,\psi_{\uX',\oy'}(\varphi'_{\ox'}(\cC'^{(r_1,r_2)}_n\otimes_{\ocB'_n}\sigma'^*_n(\cC^{(r_3)}_{\uptau,n})))).
\end{equation}
In view of \eqref{p2-rftchta10m}, \eqref{p2-rftchta13d}, (\cite{ag1} (3.4.23.1))
 and (\cite{serre1} I prop.~8), the target of this isomorphism is canonically isomorphic to 
\begin{equation}
\underset{\underset{(U'\rightarrow U)\in \mM^\circ_{\ox'}(\bQ',\bQ)}{\longrightarrow}}{\lim}\
\rH^q(\Pi_{U'\rightarrow U}, \fC'^{(r_1,r_2)}_{U'}\otimes_{R'_{\uptau,U'}}(\fC^{(r_3)}_{\uptau,U'}/p^n\fC^{(r_3)}_{\uptau,U'})).
\end{equation}
In view of \eqref{p2-rftchta11m}, \eqref{p2-rftchta13c}, the description of \eqref{p2-rftchta13f} and (\cite{ag1} 6.5.15(ii)),
we deduce from this that the stalk of the morphism \eqref{p2-rftchta5a} at $\varrho(\oy\rightsquigarrow \ox')$
identifies with the direct limit of the homomorphisms
\begin{equation}
(\fC^{(r_2)}_U\otimes_{\oR_U} \oR^\intern_{U'\rightarrow U})\otimes_{R'_{\uptau,U'}}(\fC^{(r_3)}_{\uptau,U'}/p^n\fC^{(r_3)}_{\uptau,U'})\rightarrow 
(\fC'^{(r_1,r_2)}_{U'}\otimes_{R'_{\uptau,U'}}(\fC^{(r_3)}_{\uptau,U'}/p^n\fC^{(r_3)}_{\uptau,U'}))^{\Pi_{U'\rightarrow U}}
\end{equation}
deduced from $\uppsi^{\oy',\ur}_\mu$ modulo $p^n$ \eqref{p2-rftchta13e}, when $(\mu\colon U'\rightarrow U, \iota\colon \ox'\rightarrow U')$  
goes through the category $\umM^\circ_{\ox'}(\bQ',\bQ)$ \eqref{p2-fhtft301}. The proposition then follows from \ref{p2-fhtal13} in view of (\cite{ag1} 6.5.3).

\begin{cor}\label{p2-rftchta6}
Let $\ur=(r_1,r_2,r_3)$, $\ur'=(r'_1,r'_2,r'_3)$ be two elements of $I$ \eqref{p2-fhtft33} such that $r_1> r'_1$, $r_2\geq  r'_2$ and $r_3\geq  r'_3$. Then,
\begin{itemize}
\item[{\rm (i)}] The canonical homomorphism \eqref{p2-rftchta4f}
\begin{equation}\label{p2-rftchta6a}
\bvlgg^*(\bvcC^{(r_2)})\otimes_{\bvocB^!}\huppi^*(\hcC^{(r_3)}_\uptau)\rightarrow 
\bvvarpi_*(\bvcC'^{\ur})
\end{equation}
is $\alpha$-injective. Let $\bvcH^{\ur}$ be its cokernel.
\item[{\rm (ii)}] There exists a rational number $a >0$ such that the morphism
$\bvcH^{\ur}\rightarrow \bvcH^{\ur'}$ induced by \eqref{p2-fhtft38h} is annihilated by $p^a$.
\item[{\rm (iii)}] There exists a rational number $b>0$ such that for every integer $q\geq 1$,
the canonical morphism of $\tG^{\mN^\circ}_s$
\begin{equation}\label{p2-rftchta6b}
\rR^q\bvvarpi_*(\bvcC'^{\ur})\rightarrow 
\rR^q\bvvarpi_*(\bvcC'^{\ur'})
\end{equation}
is annihilated by $p^b$.
\end{itemize}
\end{cor}

This follows from \ref{p2-rftchta5} and (\cite{agt} III.7.3(i) and (III.7.5.5)).

\begin{cor}\label{p2-rftchta7}
Let $\ur=(r_1,r_2,r_3)$, $\ur'=(r'_1,r'_2,r'_3)$ be two elements of $I$ \eqref{p2-fhtft33} such that $r_1> \sup(r'_1,r_2)$, $r_2\geq  r'_2$ and $r_3\geq  r'_3$. Then,
\begin{itemize}
\item[{\rm (i)}] The canonical homomorphism \eqref{p2-rftchta4f}
\begin{equation}\label{p2-rftchta7a}
\ttu^\ur\colon (\bvlgg^*(\bvcC^{(r_2)})\otimes_{\bvocB^!}\huppi^*(\hcC^{(r_3)}_\uptau))_\mQ\rightarrow 
\bvvarpi_*(\bvcC'^{\ur})_\mQ
\end{equation}
admits a canonical left inverse
\begin{equation}\label{p2-rftchta7b}
\ttv^\ur\colon \bvvarpi_*(\bvcC'^{\ur})_\mQ\rightarrow 
(\bvlgg^*(\bvcC^{(r_2)})\otimes_{\bvocB^!}\huppi^*(\hcC^{(r_3)}_\uptau))_\mQ.
\end{equation}
\item[{\rm (ii)}] The diagram
\begin{equation}\label{p2-rftchta7c}
\xymatrix{
{\bvvarpi_*(\bvcC'^{\ur})_\mQ}\ar[d]_{\ttv^{\ur}}\ar[r]&
{\bvvarpi_*(\bvcC'^{\ur'})_\mQ}\\
{(\bvlgg^*(\bvcC^{(r_2)})\otimes_{\bvocB^!}\huppi^*(\hcC^{(r_3)}_\uptau))_\mQ}\ar[r]&
{(\bvlgg^*(\bvcC^{(r'_2)})\otimes_{\bvocB^!}\huppi^*(\hcC^{(r'_3)}_\uptau))_\mQ,}\ar[u]_{\ttu^{\ur'}}}
\end{equation}
where the horizontal morphisms are induced by \eqref{p2-htaft18f}, \eqref{p2-fhtft3i} and \eqref{p2-fhtft28e}, is commutative.
\item[{\rm (iii)}] For every integer $q\geq 1$, the canonical morphism
\begin{equation}\label{p2-rftchta7d}
\rR^q\bvvarpi_*(\bvcC'^{\ur})_\mQ\rightarrow 
\rR^q\bvvarpi_*(\bvcC'^{\ur'})_\mQ
\end{equation}
vanishes.
\end{itemize}
\end{cor}

Indeed, since $r_1>\sup(r_2,r'_1)$, there exists $\ur''=(r''_1,r''_2,r''_3)\in I$ such that $r_1>r''_1>r'_1$, $r''_2=r_2$ and $r''_3=r_3$. 
Then, by \eqref{p2-fhtft38h}, we may reduce to the case where $r'_2=r_2$ and $r'_3=r_3$. 

By \ref{p2-rftchta6}(i)-(ii), $\ttu^\ur$ is injective and there exists a unique 
$\bvocB^!_\mQ$-linear morphism $\ttv^{\ur,\ur'}$ that fits into a commutative diagram 
\begin{equation}\label{p2-rftchta7e}
\xymatrix{
{(\bvlgg^*(\bvcC^{(r_2)})\otimes_{\bvocB^!}\huppi^*(\hcC^{(r_3)}_\uptau))_\mQ}\ar[r]^-(0.5){\ttu^\ur}\ar[d]& 
{\bvvarpi_*(\bvcC'^{\ur})_\mQ}\ar[d]\ar[ld]_{\ttv^{\ur,\ur'}}\\
{(\bvlgg^*(\bvcC^{(r'_2)})\otimes_{\bvocB^!}\huppi^*(\hcC^{(r'_3)}_\uptau))_\mQ}\ar[r]^-(0.5){\ttu^{\ur'}}& 
{\bvvarpi_*(\bvcC'^{\ur'})_\mQ,}}
\end{equation}
where the vertical arrows are induced by \eqref{p2-htaft18f}, \eqref{p2-fhtft3i} and \eqref{p2-fhtft28e}. 
We deduce that the morphisms $\ttv^{\ur,\ur'}$ do not depend on $\ur'$. We denote it by $\ttv^\ur$.  
Since $\ttu^{\ur'}$ in injective, \eqref{p2-rftchta7e} implies that $\ttv^{\ur}$ is a left inverse of $\ttu^{\ur}$. Propositions (i) and (ii) follow. 
Proposition (iii) follows directly from \ref{p2-rftchta6}(iii).

\begin{cor}\label{p2-rftchta8}
Let $t_2,t_3$ be two rational numbers such that $t_2\geq t_3 \geq 0$, $I_{t_2,t_3}$ the subset of elements $\ur=(r_1,r_2,r_3)$ of $I$ \eqref{p2-fhtft33} 
such that $r_1>r_2\geq t_2$ and $r_3\geq t_3$. Then, 
\begin{itemize}
\item[{\rm (i)}] The canonical homomorphism \eqref{p2-rftchta4f} induces an isomorphism of $\bvocB^!_\mQ$-algebras
\begin{equation}\label{p2-rftchta8a}
(\bvlgg^*(\bvcC^{(t_2)})\otimes_{\bvocB^!}\huppi^*(\hcC^{(t_3)}_\uptau))_\mQ\stackrel{\sim}{\rightarrow} 
\underset{\underset{\ur\in I_{t_2,t_3}}{\longrightarrow}}{\lim}\ \bvvarpi_*(\bvcC'^{\ur})_\mQ;
\end{equation}
and for every integer $q\geq 1$, we have 
\begin{equation}\label{p2-rftchta8b}
\underset{\underset{\ur\in I_{t_2,t_3}}{\longrightarrow}}{\lim}\ \rR^q\bvvarpi_*(\bvcC'^{\ur})_\mQ= 0,
\end{equation}
where the limits are taken in $\bMod_\mQ(\bvocB^!)$. These are, in particular, representable.

\item[{\rm (ii)}] The canonical homomorphism \eqref{p2-rftchta4f} induces an isomorphism of ind-$\bvocB^!$-algebras
\begin{equation}\label{p2-rftchta8c}
\bvlgg^*(\bvcC^{(t_2)})\otimes_{\bvocB^!}\rI\huppi^*(\upalpha_{\co_{\fX'}}(\hcC^{(t_3)}_{\uptau,\mQ}))\stackrel{\sim}{\rightarrow} 
\underset{\underset{\ur\in I_{t_2,t_3}}{\longrightarrow}}{\mlq\mlq\lim\mrq\mrq}\ \rI\bvvarpi_*(\upalpha_{\bvocB'}(\bvcC'^{\ur}_\mQ)),
\end{equation}
where the functors $\upalpha_{\co_{\fX'}}$ and $\upalpha_{\bvocB'}$ are defined in \eqref{p1-bcim5c};  
and for every integer $q\geq 1$, we have
\begin{equation}\label{p2-rftchta8d}
\underset{\underset{\ur\in I_{t_2,t_3}}{\longrightarrow}}{\mlq\mlq\lim\mrq\mrq}\ \rR^q\rI\bvvarpi_*(\upalpha_{\bvocB'}(\bvcC'^{\ur}_\mQ))= 0.
\end{equation} 
\end{itemize}
\end{cor}

(i) It follows from \ref{p2-rftchta7}. 

(ii) Indeed, applying the functor $\upalpha_{\bvocB^!}$ \eqref{p2-rftchta14b}, 
we obtain an analogue of \ref{p2-rftchta7} in the category $\bIndMod(\bvocB^!)$, from which we immediately 
deduce the proposition taking into account \eqref{p1-bcim1f}.

\subsection{}\label{p2-rftchta15}
Let $\ur=(r_1,r_2,r_3)$ be an element of $I$ \eqref{p2-fhtft33}, $n$ an integer $\geq 1$. We denote by 
\begin{equation}\label{p2-rftchta15a}
\delta^\ur_n\colon \cC'^{\ur}_n\rightarrow \sigma'^*_n(\Omega'_n)\otimes_{\ocB'_n}\cC'^{\ur}_n
\end{equation}
the $\ocB'_n$-derivation defined by 
\begin{equation}\label{p2-rftchta15b}
\delta^\ur_n=\delta_{\cC'^{(r_1,r_2)}_n}\otimes \id- \id \otimes \sigma'^*_n(\delta'_{\cC^{(r_3)}_{\uptau,n}}),
\end{equation}
where $\cC'^{\ur}_n=\cC'^{(r_1,r_2)}_n\otimes_{\ocB'_n}\sigma'^*_n(\cC^{(r_3)}_{\uptau,n})$ \eqref{p2-fhtft38k}, 
$\delta_{\cC'^{(r_1,r_2)}_n}$ is the derivation \eqref{p2-fhtft23b} and $\delta'_{\cC^{(r_3)}_{\uptau,n}}$ is the reduction modulo $p^n$ of 
the derivation $\delta'_{\cC^{(r_3)}_{\uptau}}$ \eqref{p2-fhtft6aa}.  
It is a Higgs $\ocB'_n$-field with coefficients in $\sigma'^*_n(\Omega'_n)$.  By \ref{p2-fhtft362} and \eqref{p2-fhtal15c}, we have a commutative diagram
\begin{equation}\label{p2-rftchta15c}
\xymatrix{
{\uptheta^*_n(\cC^{(r_2)}_n)\otimes_{\ocB'_n}\sigma'^*_n(\cC^{(r_3)}_{\uptau,n})}\ar[rr]^-(0.5){\id\otimes \sigma'^*_n(\delta_{\cC^{(r_3)}_{\uptau,n}})}
\ar[d]_{\uppsi^\ur_n}&&
{\sigma'^*_n(\cog^*(\Omega_n))\otimes_{\ocB'_n}\uptheta_n^*(\cC^{(r_2)}_n)\otimes_{\ocB'_n}\sigma'^*_n(\cC^{(r_3)}_{\uptau,n})}
\ar[d]^-(0.5){\sigma'^*_n(u_n)\otimes \uppsi^\ur_n}\\
{\cC'^{\ur}_n}\ar[rr]^-(0.5){-\delta^\ur_n}&&
{\sigma'^*_n(\Omega'_n)\otimes_{\ocB'_n}\cC'^{\ur}_n,}}
\end{equation}
where $\uppsi^\ur_n$ is the homomorphism \eqref{p2-fhtft38e} and $u_n$ (resp.\ $\delta_{\cC^{(r_3)}_{\uptau,n}}$) 
is the reduction modulo $p^n$ of the morphism $u$ defined in \eqref{p2-fhtft90l} (resp.\ derivation $\delta_{\cC^{(r_3)}_{\uptau}}$ \eqref{p2-fhtft6a}).
In particular, we have 
\begin{equation}
\delta^\ur_n\circ \upphi^\ur_n=0,
\end{equation}
where $\upphi^\ur_n$ is the homomorphism defined in \eqref{p2-fhtft38d}. 

By \ref{p2-fhtft362} and \eqref{p2-fhtal15cc}, we have a commutative diagram
\begin{equation}\label{p2-rftchta15i}
\xymatrix{
{\uptheta_n^*(\cC_n^{(r_2)})\otimes_{\ocB'_n}\sigma'^*_n(\cC^{(r_3)}_{\uptau,n})}\ar[rr]^-(0.5){b_n}
\ar[d]_{\uppsi^\ur_n}&&
{\sigma'^*_n(\cog^*(\Omega_n))\otimes_{\ocB'_n}\uptheta^*_n(\cC_n^{(r_2)})\otimes_{\ocB'_n}\sigma'^*_n(\cC^{(r_3)}_{\uptau,n})}
\ar[d]^-(0.5){\id\otimes \uppsi^\ur_n}\\
{\cC'^{\ur}_n}\ar[rr]^-(0.5){a_n}&&
{\sigma'^*_n(\cog^*(\Omega_n))\otimes_{\bvocB'}\cC'^{\ur}_n,}}
\end{equation}
where $a_n$ denotes the Higgs field 
\begin{equation}
a_n=\id\otimes \sigma'^*_n(\delta_{\cC^{(r_3)}_{\uptau,n}})
\end{equation}
on $\cC'^{\ur}_n=\cC'^{(r_1,r_2)}_n\otimes_{\ocB'_n}\sigma'^*_n(\cC^{(r_3)}_{\uptau,n})$, $b_n$ denotes the Higgs field 
\begin{equation}
b_n=\uptheta^*_n(\delta_{\cC_n^{(r_2)}})\otimes\id+\id\otimes \sigma'^*_n(\delta_{\cC^{(r_3)}_{\uptau,n}})
\end{equation}
on $\theta^*_n(\cC_n^{(r_2)})\otimes_{\ocB'_n}\sigma'^*_n(\cC^{(r_3)}_{\uptau,n})$, and $\delta_{\cC_n^{(r_2)}}$ is the derivation \eqref{p2-htaft20b}.

We denote by $\mK^\bullet(\cC^{(r_3)}_{\uptau,n})$ the Dolbeault complex of $(\cC^{(r_3)}_{\uptau,n},\delta_{\cC^{(r_3)}_{\uptau,n}})$ 
and by $\mK^\bullet(\cC'^{\ur}_n)$ the Dolbeault complex of 
$(\cC'^{\ur}_n,-\delta_n^\ur)$.   
In view of \eqref{p2-rftchta15c}, $\uppsi^\ur_n$ and $u_n$ induce a morphism of complexes
\begin{equation}\label{p2-rftchta15d}
\uptheta^*_n(\cC^{(r_2)}_n)\otimes_{\ocB'_n}\sigma'^*_n(\mK^\bullet(\cC^{(r_3)}_{\uptau,n}))\rightarrow  
\mK^\bullet(\cC'^{\ur}_n),
\end{equation}
where on the left, the tensor product and $\sigma'^*_n$ are defined term by term. 
We denote by $\tmK^\bullet_n(\ur)$ the mapping cone of the morphism \eqref{p2-rftchta15d}. 

Let $\ur'=(r'_1,r'_2,r'_3)\in I$ be such that $r_i\geq r'_i$ for all $1\leq i\leq 3$.
By \eqref{p2-fhtft6d} and \eqref{p2-fhtft23c}, 
the homomorphisms $\alpha_{\uptau,n}^{r,r'}$ \eqref{p2-fhtft3i} and $\alpha'^{r_1,r_2,r'_1,r'_2}_n$ \eqref{p2-fhtft25e}
induce morphisms of complexes
\begin{eqnarray}
\upiota^{r_3,r'_3}_{\uptau,n}\colon \mK^\bullet(\cC^{(r_3)}_{\uptau,n})&\rightarrow &\mK^\bullet(\cC^{(r'_3)}_{\uptau,n}),\label{p2-rftchta15e}\\
\upiota^{\ur,\ur'}_n\colon \mK^\bullet(\cC'^{\ur}_n)&\rightarrow &
\mK^\bullet(\cC'^{\ur'}_n).\label{p2-rftchta15f}
\end{eqnarray}
Moreover, the diagram 
\begin{equation}\label{p2-rftchta15g}
\xymatrix{
{\uptheta^*_n(\cC^{(r_2)}_n)\otimes_{\ocB'_n}\sigma'^*_n(\mK^\bullet(\cC^{(r_3)}_{\uptau,n}))}
\ar[r]\ar[d]_{\uptheta_n^*(\alpha_n^{r'_2,r_2})\otimes \sigma'^*_n(\upiota^{r_3,r'_3}_{\uptau,n})}&
{\mK^\bullet(\cC'^{\ur}_n)}\ar[d]^{\upiota^{\ur,\ur'}_n}\\
{\uptheta^*_n(\cC^{(r'_2)}_n)\otimes_{\ocB'_n}\sigma'^*_n(\mK^\bullet(\cC^{(r'_3)}_{\uptau,n}))}\ar[r]&
{\mK^\bullet(\cC'^{\ur'}_n),}}
\end{equation}
where $\alpha^{r'_2,r_2}_n$ is the morphism \eqref{p2-htaft16g} and  
the horizontal arrows are the morphisms \eqref{p2-rftchta15d}, is commutative.  We deduce a morphism of complexes
\begin{equation}\label{p2-rftchta15h}
\tupiota^{\ur,\ur'}_n\colon \tmK^\bullet_n(\ur)\rightarrow \tmK^\bullet_n(\ur').
\end{equation}

\begin{prop}\label{p2-rftchta16}
Let $\ur=(r_1,r_2,r_3)$, $\ur'=(r'_1,r'_2,r'_3)$ be two elements of $I$ \eqref{p2-fhtft33} such that $r_1> r'_1$, $r_2\geq r'_2$ and $r_3\geq r'_3$. 
Then, there exists a rational number $\alpha\geq 0$ depending on $r_1$ and $r'_1$ 
but not on the morphisms $f$ and $g$ satisfying the conditions of \ref{p2-fhtft1}, such that for every integer $n\geq 1$, the morphism 
\begin{equation}\label{p2-rftchta16a}
p^\alpha\tupiota^{\ur,\ur'}_n\colon  \tmK^\bullet_n(\ur)\rightarrow \tmK^\bullet_n(\ur'),
\end{equation}
where $\tupiota^{\ur,\ur'}_n$ is the morphism \eqref{p2-rftchta15h}, is homotopic to $0$ by a $\ocB'_n$-linear homotopy, stalk by stalk.
\end{prop}

Let $(\oy'\rightsquigarrow \ox')$ be a point of $X'_\et\gtimes_{X'_\et}\oX'^\rhd_\et$ such that $\ox'$ is above $s$, $n$ an integer $\geq 1$.
We take again the notation introduced in \ref{p2-rftchta10}, \ref{p2-rftchta11} and \ref{p2-rftchta13}.
In view of \eqref{p2-rftchta13d} and (\cite{agt} VI.5.34(ii)), the stalk of the derivation $\delta_{\cC'^{(r_1,r_2)}_n}$ \eqref{p2-fhtft23b} 
at $\rho'(\oy' \rightsquigarrow \ox')$ identifies with the direct limit of the reductions modulo $p^n$ of the $\hoRp^{\oy'}_{U'}$-derivations
\begin{equation}
\delta_{\fC'^{\oy',(r_1,r_2)}_{U'}}\colon \fC'^{\oy',(r_1,r_2)}_{U'}\rightarrow \txi^{-1}\tOmega^1_{X'/X}(U')\otimes_{\co_{X'}(U')}\fC'^{\oy',(r_1,r_2)}_{U'}
\end{equation}
defined in \eqref{p2-fhtal8e}, when $U'$ goes through the category $\bQ'_{\ox'}$.
Moreover, the stalk of the morphism $\tupiota^{\ur,\ur'}_n$ \eqref{p2-rftchta15h} at $\rho'(\oy' \rightsquigarrow \ox')$ 
identifies with the direct limit of the reductions modulo $p^n$ of the morphisms 
\begin{equation}
\tupiota^{\ur,\ur'}_{\mu}\colon 
\tmK^\bullet(\fC'^{(r_1,r_2)}_{U'}\hotimes_{R'_{\uptau,U'}}\fC^{(r_3)}_\uptau)\rightarrow 
\tmK^\bullet(\fC'^{(r'_1,r'_2)}_{U'}\hotimes_{R'_{\uptau,U'}}\fC^{(r'_3)}_{\uptau,U'})
\end{equation}
defined in \eqref{p2-fhtal15h}, when $(\mu\colon U'\rightarrow U, \iota\colon \ox'\rightarrow U')$  
goes through the category $\umM^\circ_{\ox'}(\bQ',\bQ)$ \eqref{p2-fhtft301}.
This follows from \ref{p2-qfc263}, in view of \ref{p2-cmt5}, \eqref{p2-qfc262g} and (\cite{agt} VI.10.36 and (VI.10.18.1)). 
The proposition then follows from \ref{p2-fhtal16}(i) in view of (\cite{agt} III.9.5).

\subsection{}\label{p2-rftchta17}
Let $\ur=(r_1,r_2,r_3)$ be an element of $I$ \eqref{p2-fhtft33}. We denote by 
\begin{equation}\label{p2-rftchta17a}
\bvdelta^\ur\colon \bvcC'^{\ur}\rightarrow \hupsigma'^*(\hOmega')\otimes_{\bvocB'}\bvcC'^{\ur}
\end{equation}
the $\bvocB'$-derivation defined by 
\begin{equation}\label{p2-rftchta17ab}
\bvdelta^\ur=\delta_{\bvcC'^{(r_1,r_2)}}\otimes \id- \id \otimes \hupsigma'^*(\delta'_{\hcC^{(r_3)}_{\uptau}}),
\end{equation}
where $\bvcC'^{\ur}=\bvcC'^{(r_1,r_2)}\otimes_{\ocB'}\hupsigma'^*(\hcC^{(r_3)}_{\uptau})$ \eqref{p2-fhtft38l}, 
$\delta_{\bvcC'^{(r_1,r_2)}}$ is the derivation \eqref{p2-fhtft29b} and $\delta'_{\hcC^{(r_3)}_{\uptau}}$ is 
the derivation \eqref{p2-fhtft6bb}. By \eqref{p2-htaft8e}, it identifies canonically with the derivation $(\delta^\ur_{n})_{\mN}$ \eqref{p2-rftchta15a}, 
and it is a Higgs $\bvocB'$-field with coefficients in $\hupsigma'^*(\hOmega')$. 

By \eqref{p2-rftchta15c}, we have a commutative diagram
\begin{equation}\label{p2-rftchta17k}
\xymatrix{
{\bvuptheta^*(\bvcC^{(r_2)})\otimes_{\bvocB'}\hupsigma'^*(\hcC^{(r_3)}_{\uptau})}\ar[rr]^-(0.5){\id\otimes \hupsigma'^*(\delta_{\hcC^{(r_3)}_{\uptau}})}
\ar[d]_{\bvuppsi^\ur}&&
{\hupsigma'^*(\fgg^*(\hOmega))\otimes_{\bvocB'}\bvuptheta^*(\bvcC^{(r_2)})\otimes_{\bvocB'}\hupsigma'^*(\hcC^{(r_3)}_{\uptau})}
\ar[d]^-(0.5){\hupsigma'^*(\hu)\otimes \bvuppsi^\ur}\\
{\bvcC'^{\ur}}\ar[rr]^-(0.5){-\bvdelta^\ur}&&
{\hupsigma'^*(\hOmega')\otimes_{\bvocB'}\bvcC'^{\ur},}}
\end{equation}
where $u\colon \cog^*(\Omega)\rightarrow \Omega'$ is the canonical morphism \eqref{p2-fhtft90l},
$\bvuppsi^\ur$ is the homomorphism \eqref{p2-fhtft38g} and $\delta_{\hcC^{(r_3)}_{\uptau}}$ is the derivation \eqref{p2-fhtft6b}. In particular, we have 
\begin{equation}\label{p2-rftchta17l}
\bvdelta^\ur\circ \bvupphi^\ur=0,
\end{equation}
where $\bvupphi^\ur$ is the homomorphism defined in \eqref{p2-fhtft38f}. 
By \eqref{p2-rftchta15i}, the diagram 
\begin{equation}\label{p2-rftchta17m}
\xymatrix{
{\bvuptheta^*(\bvcC^{(r_2)})\otimes_{\bvocB'}\hupsigma'^*(\hcC^{(r_3)}_{\uptau})}\ar[rr]^-(0.5){b}
\ar[d]_{\bvuppsi^\ur}&&
{\hupsigma'^*(\fgg^*(\hOmega))\otimes_{\bvocB'}\bvuptheta^*(\bvcC^{(r_2)})\otimes_{\bvocB'}\hupsigma'^*(\hcC^{(r_3)}_{\uptau})}
\ar[d]^-(0.5){\id\otimes \bvuppsi^\ur}\\
{\bvcC'^{\ur}}\ar[rr]^-(0.5){a}&&
{\hupsigma'^*(\fgg^*(\hOmega))\otimes_{\bvocB'}\bvcC'^{\ur},}}
\end{equation}
where $a$ denotes the Higgs field 
\begin{equation}
a=\id\otimes \hupsigma'^*(\delta_{\hcC^{(r_3)}_{\uptau}})
\end{equation}
on $\bvcC'^{\ur}=\bvcC'^{(r_1,r_2)}\otimes_{\ocB'}\hupsigma'^*(\hcC^{(r_3)}_{\uptau})$, $b$ denotes the Higgs field 
\begin{equation}
b=\bvuptheta^*(\delta_{\bvcC^{(r_2)}})\otimes\id+\id\otimes \hupsigma'^*(\delta_{\hcC^{(r_3)}_{\uptau}})
\end{equation}
on $\bvuptheta^*(\bvcC^{(r_2)})\otimes_{\bvocB'}\hupsigma'^*(\hcC^{(r_3)}_{\uptau})$, and $\delta_{\bvcC^{(r_2)}}$ is the derivation \eqref{p2-htaft19b},
is commutative. 

We denote by $\mK^\bullet(\hcC^{(r_3)}_{\uptau})$
the Dolbeault complex of $(\hcC^{(r_3)}_{\uptau},\delta_{\hcC^{(r_3)}_{\uptau}})$ \eqref{p2-fhtft6b}
and by $\mK^\bullet(\bvcC'^{\ur})$ the Dolbeault complex of $(\bvcC'^{\ur},-\bvdelta^\ur)$.   
In view of \eqref{p2-rftchta17l}, we can consider $\mK^\bullet(\bvcC'^{\ur})$ 
as a complex of ind-$\bvuptheta^*(\bvcC^{(r_2)})$-modules via the homomorphism $\bvupphi^{\ur}$ \eqref{p2-fhtft38f}. 
In view of \eqref{p2-rftchta17k}, $\bvuppsi^\ur$ and $\hu$ induce a morphism of complexes of ind-$\bvuptheta^*(\bvcC^{(r_2)})$-modules 
\begin{equation}\label{p2-rftchta17b}
\bvuptheta^*(\bvcC^{(r_2)})\otimes_{\bvocB'}\hupsigma'^*(\mK^\bullet(\hcC^{(r_3)}_{\uptau}))\rightarrow  
\mK^\bullet(\bvcC'^{\ur}),
\end{equation}
where on the left, the tensor product and $\hupsigma'^*$ are defined term by term. 
We denote by $\tmK^\bullet(\ur)$ the mapping cone of the morphism \eqref{p2-rftchta17b}. 

Let $\ur'=(r'_1,r'_2,r'_3)\in I$ be such that $r_i\geq r'_i$ for all $1\leq i\leq 3$.
By \eqref{p2-fhtft6d} and \eqref{p2-fhtft23c}, we have 
\begin{equation}\label{p2-rftchta17j}
(\id \otimes \bvalpha'^{\ur,\ur'}) \circ \bvdelta^\ur=\bvdelta^{\ur'} \circ \bvalpha'^{\ur,\ur'},
\end{equation}
where $\bvalpha'^{\ur,\ur'}$ is the homomorphism defined in \eqref{p2-fhtft38n}. 
The homomorphisms $\halpha_{\uptau}^{r,r'}$ \eqref{p2-fhtft3i} and $\bvalpha'^{\ur,\ur'}$
induce then morphisms of complexes
\begin{eqnarray}
\hupiota^{r_3,r'_3}_{\uptau}\colon \mK^\bullet(\hcC^{(r_3)}_{\uptau})&\rightarrow &\mK^\bullet(\hcC^{(r'_3)}_{\uptau}),\label{p2-rftchta17c}\\
\bvupiota^{\ur,\ur'}\colon \mK^\bullet(\bvcC'^{\ur})&\rightarrow &
\mK^\bullet(\bvcC'^{\ur'}).\label{p2-rftchta17d}
\end{eqnarray}
Moreover, by \eqref{p2-rftchta15g}, the diagram 
\begin{equation}\label{p2-rftchta17e}
\xymatrix{
{\bvuptheta^*(\bvcC^{(r_2)})\otimes_{\bvocB'}\hupsigma'^*(\mK^\bullet(\hcC^{(r_3)}_{\uptau}))}
\ar[r]\ar[d]_{\bvuptheta^*(\bvalpha^{r_2,r'_2})\otimes \hupsigma'^*(\hupiota^{r_3,r'_3}_{\uptau})}&
{\mK^\bullet(\bvcC'^{\ur})}\ar[d]^{\bvupiota^{\ur',\ur}}\\
{\bvuptheta^*(\cC^{(r'_2)})\otimes_{\bvocB'}\hupsigma'^*(\mK^\bullet(\hcC^{(r'_3)}_{\uptau}))}\ar[r]&
{\mK^\bullet(\bvcC'^{\ur'}),}}
\end{equation}
where $\bvalpha^{r'_2,r_2}$ is the morphism \eqref{p2-htaft18f} and the horizontal arrows are the morphisms \eqref{p2-rftchta17b}, is commutative.  
We deduce a morphism of complexes
\begin{equation}\label{p2-rftchta17f}
\tupiota^{\ur,\ur'}\colon \tmK^\bullet(\ur)\rightarrow \tmK^\bullet(\ur').
\end{equation}

Recall that we have the canonical functors \eqref{p2-htaft27}
\begin{eqnarray}
\iota_{\bvocB'}\colon \bMod(\bvocB')&\rightarrow& \bIndMod(\bvocB'),\label{p2-rftchta17h}\\
\upalpha_{\bvocB'}\colon \bMod_\mQ(\bvocB')&\rightarrow& \bIndMod(\bvocB'), \label{p2-rftchta17g}
\end{eqnarray}
which are exact and fully faithful.
We will identify $\bMod(\bvocB')$ (resp.\ $\bMod_\mQ(\bvocB')$) with a full subcategory of $\bIndMod(\bvocB')$
by the functor $\iota_{\bvocB'}$ (resp.\ $\upalpha_{\bvocB'}$), which we will omit from the notation. 

We denote by $\tmK^\bullet_{\mQ}(\ur)$  (resp.\ $\mK^\bullet_\mQ(\bvcC'^{\ur})$) 
the image of the complex $\tmK^\bullet(\ur)$ 
(resp.\ $\mK^\bullet(\bvcC'^{\ur})$) in $\bMod_\mQ(\bvocB')$.
These complexes will also be considered as complexes of ind-$\bvocB'$-modules via $\upalpha_{\bvocB'}$.

\begin{prop}\label{p2-rftchta18}
For all elements $\ur=(r_1,r_2,r_3)$ and $\ur'=(r'_1,r'_2,r'_3)$ of $I$ \eqref{p2-fhtft33} such that $r_1> r'_1$, $r_2\geq r'_2$ and $r_3\geq r'_3$, 
and every integer $q$, the morphism of $\bvocB'_\mQ$-modules
\begin{equation}\label{p2-rftchta18a}
\rH^q(\tupiota^{\ur,\ur'}_\mQ)\colon  \rH^q(\tmK^\bullet_\mQ(\ur))\rightarrow \rH^q(\tmK^\bullet_\mQ(\ur')),
\end{equation}
where $\tupiota^{\ur,\ur'}$ is the morphism defined in \eqref{p2-rftchta17f}, vanishes.
\end{prop}

This follows from \ref{p2-rftchta16} and (\cite{agt} III.7.3(i)).

\begin{cor}\label{p2-rftchta19}
Let $\ur=(r_1,r_2,r_3)$, $\ur'=(r'_1,r'_2,r'_3)$ be two elements of $I$ \eqref{p2-fhtft33} such that $r_1> \sup(r'_1,r_2)$, $r_2\geq  r'_2$ and $r_3\geq  r'_3$,
$q$ an integer. Then,
\begin{itemize}
\item[{\rm (i)}] The canonical sequence of $\bvocB'_\mQ$-modules \eqref{p2-rftchta17b}
\begin{equation}\label{p2-rftchta19a}
0\longrightarrow{\rH^q(\bvuptheta^*(\bvcC^{(r_2)})\otimes_{\bvocB'}\hupsigma'^*(\mK^\bullet(\hcC^{(r_3)}_{\uptau})))_{\mQ}}\stackrel{u^{q,\ur}}{\longrightarrow}\rH^q(\mK^\bullet_\mQ(\bvcC'^{\ur}))\longrightarrow \rH^q(\utmK^\bullet_{\mQ}(\ur))\longrightarrow 0
\end{equation}
is exact. 
\item[{\rm (ii)}] The morphism $u^{q,\ur}$ \eqref{p2-rftchta19a} admits a canonical left inverse
\begin{equation}\label{p2-rftchta19b}
v^{q,\ur}\colon \rH^q(\mK^\bullet_\mQ(\bvcC'^{\ur})) \rightarrow 
\rH^q(\bvuptheta^*(\bvcC^{(r_2)})\otimes_{\bvocB'}\hupsigma'^*(\mK^\bullet(\hcC^{(r_3)}_{\uptau})))_{\mQ}.
\end{equation}
\item[{\rm (iii)}] The diagram
\begin{equation}\label{p2-rftchta19c}
\xymatrix{
{\rH^q(\mK^\bullet_\mQ(\bvcC'^{\ur}))}\ar[rr]^-(0.5){\rH^q(\bvupiota^{\ur,\ur'}_\mQ)}\ar[d]_{v^{q,\ur}}&&
{\rH^q(\mK^\bullet_\mQ(\bvcC'^{\ur'}))}\\
{\rH^q(\bvuptheta^*(\bvcC^{(r_2)})\otimes_{\bvocB'}\hupsigma'^*(\mK^\bullet(\hcC^{(r_3)}_{\uptau})))_{\mQ}}\ar[rr]&&{\rH^q(\bvuptheta^*(\bvcC^{(r'_2)})\otimes_{\bvocB'}\hupsigma'^*(\mK^\bullet(\hcC^{(r'_3)}_{\uptau})))_{\mQ}}\ar[u]_{u^{q,\ur'}}}
\end{equation}
where $\bvupiota^{\ur,\ur'}$ is defined in \eqref{p2-rftchta17d} and the lower horizontal morphism is induced by $\bvalpha^{r_2,r'_2}$ \eqref{p2-htaft18f} 
and $\halpha_\uptau^{r_3,r'_3}$ \eqref{p2-fhtft3i}, is commutative. 
\end{itemize}
\end{cor}

Indeed, since $r_1>\sup(r_2,r'_1)$, there exists $\ur''=(r''_1,r''_2,r''_3)\in I$ such that $r_1>r''_1>r'_1$, $r''_2=r_2$ and $r''_3=r_3$. 
Then, by \eqref{p2-fhtft38h}, we may reduce to the case where $r'_2=r_2$ and $r'_3=r_3$, which we will assume in the following proof. 

(i) By \eqref{p2-fhtft38h}, we have a commutative diagram 
\begin{equation}
\xymatrix{
{\rH^q(\utmK^\bullet_{\mQ}(\ur))}\ar[d]_{\rH^q(\tupiota^{\ur,\ur'}_\mQ)}
\ar[r]&{\rH^{q+1}(\bvuptheta^*(\bvcC^{(r_2)})\otimes_{\bvocB'}\hupsigma'^*(\mK^\bullet(\hcC^{(r_3)}_{\uptau})))_{\mQ}}\ar@{=}[d]\\
{\rH^q(\utmK^\bullet_{\mQ}(\ur'))}\ar[r]&{\rH^{q+1}(\bvuptheta^*(\bvcC^{(r_2)})\otimes_{\bvocB'}\hupsigma'^*(\mK^\bullet(\hcC^{(r_3)}_{\uptau})))_{\mQ},}}
\end{equation}
where the horizontal arrows are the canonical connecting morphisms. Hence, the upper horizontal morphism vanishes by \ref{p2-rftchta18}, which implies the proposition.  

(ii) \& (iii). The commutative diagram (without the dotted arrow)
\[
\xymatrix{
{\rH^q(\bvuptheta^*(\bvcC^{(r_2)})\otimes_{\bvocB'}\hupsigma'^*(\mK^\bullet(\hcC^{(r_3)}_{\uptau})))_{\mQ}}\ar@{=}[d]\ar@{^(->}[r]^-(0.5){u^{q,\ur}}&
{\rH^q(\mK^\bullet_\mQ(\bvcC'^{\ur}))}\ar@{.>}[ld]_-(0.5){v^{q,\ur,\ur'}}
\ar[d]^-(0.5){\rH^q(\bvupiota^{\ur,\ur'}_\mQ)}
\ar@{->>}[r]&{\rH^q(\utmK^\bullet_{\mQ}(\ur))}\ar[d]^-(0.5){\rH^q(\tupiota^{\ur,\ur'}_\mQ)}\\
{\rH^q(\bvuptheta^*(\bvcC^{(r'_2)})\otimes_{\bvocB'}\hupsigma'^*(\mK^\bullet(\hcC^{(r'_3)}_{\uptau})))_{\mQ}}\ar@{^(->}[r]^-(0.5){u^{q,\ur'}}&
{\rH^q(\mK^\bullet_\mQ(\bvcC'^{\ur'}))}\ar@{->>}[r]&
{\rH^q(\utmK^\bullet_{\mQ}(\ur'))}}
\]
shows that there exists one and only one $\bvocB'_\mQ$-linear morphism 
$v^{q,\ur,\ur'}$ as above such that $\rH^q(\bvupiota^{\ur,\ur'}_\mQ)=u^{q,\ur'}\circ v^{q,\ur,\ur'}$. Since we have
$u^{q,\ur'}\circ v^{q,\ur,\ur'}\circ u^{q,\ur}=u^{q,\ur'}$ and $u^{q,\ur'}$ is injective, we deduce that $v^{q,\ur,\ur'}$ is a left inverse of $u^{q,\ur}$.
By functoriality of $u^{q,\ur'}$ \eqref{p2-fhtft38h}, $v^{q,\ur,\ur'}$ does not depend on $\ur'$; the propositions follow.

\begin{cor}\label{p2-rftchta20}
Let $\ut=(t_1,t_2,t_3)$ be a triple of rational numbers such that $t_1\geq t_2 > t_3 \geq 0$, $I_\ut$ the subset of elements $\ur=(r_1,r_2,r_3)$ of $I$ \eqref{p2-fhtft33} 
such that $r_1>t_1$, $r_2\geq t_2$  and $r_3>t_3$. Then, the morphism of complexes of ind-$\bvocB'$-modules
\begin{equation}\label{p2-rftchta20a}
\bvuptheta^*(\bvcC^{(t_2)})_\mQ[0] \rightarrow  
\underset{\underset{\ur\in I_\ut}{\longrightarrow}}{\mlq\mlq\lim \mrq\mrq} \ 
\mK^\bullet_\mQ(\bvcC'^{\ur})
\end{equation}
induced by \eqref{p2-rftchta17b}, or equivalently by $\bvupphi^{\ur}$ \eqref{p2-fhtft38f}, is a quasi-isomorphism.
\end{cor}

Observe first that $\bIndMod(\bvocB')$ admits small direct limits and that small filtered direct limits are exact (\cite{ag2} 2.6.7.5).
The morphism \eqref{p2-rftchta17b} induces a morphism of complexes of ind-$\bvuptheta^*(\bvocB')$-modules 
\begin{equation}
\underset{\underset{r_3\in \mQ_{>t_3}}{\longrightarrow}}{\mlq\mlq\lim\mrq\mrq} \
\bvuptheta^*(\bvcC^{(t_2)})_\mQ\otimes_{\bvocB'} \hupsigma'^*(\mK^\bullet(\hcC^{(r_3)}_{\uptau}))\rightarrow  
\underset{\underset{(r_1,r_2,r_3)\in I_\ut}{\longrightarrow}}{\mlq\mlq\lim\mrq\mrq} \ 
\mK^\bullet_\mQ(\bvcC'^{\ur}),
\end{equation}
where on the left, the tensor product and $\hupsigma'^*$ are defined term by term. It is a quasi-isomorphism by \ref{p2-rftchta18}. 
By \ref{p1-thbn36}, the canonical morphism 
\begin{equation}
\bvuptheta^*(\bvcC^{(t_2)})_\mQ[0] \rightarrow  
\underset{\underset{r_3\in \mQ_{>t_3}}{\longrightarrow}}{\mlq\mlq\lim\mrq\mrq} \
\bvuptheta^*(\bvcC^{(t_2)})_\mQ\otimes_{\bvocB'} \hupsigma'^*(\mK^\bullet(\hcC^{(r_3)}_{\uptau}))
\end{equation}
is a quasi-isomorphism. The proposition follows.

\section{Base change}\label{p2-bch}

The assumptions and notation of §\ref{p2-fhtft}, §\ref{p2-rfhtft} and §\ref{p2-rftchta} remain in force throughout this section.
In particular, the morphism $g\colon (X',\cM_{X'})\rightarrow (X,\cM_X)$ \eqref{p2-fhtft1a} is assumed to be {\em smooth and saturated}. 

\subsection{}\label{p2-bch1}
For every integer $n\geq 0$, the diagram of morphisms of ringed topos \eqref{p2-rftchta3k}
\begin{equation}\label{p2-bch1a}
\xymatrix{
{(\tG_s,\ocB_n^!)}\ar[r]^-(0.5){\pi_n}\ar[d]_-(0.5){\lgg_n}&{(X'_{s,\et },\co_{\oX'_n})}\ar[d]^{\ogg_n}\\
{(\tE_s,\ocB_n)}\ar[r]^-(0.5){\sigma_n}&{(X_{s,\et},\co_{\oX_n})}}
\end{equation}
is commutative up to canonical isomorphism.
For every $\co_{\oX'_n}$-module $\cF'$ of $X'_{s,\et}$ and every integer $q\geq 0$, we have a canonical base change morphism (\cite{egr1} (1.2.3.3))
\begin{equation}\label{p2-bch1b}
\sigma_n^*(\rR^q\ogg_{n*}(\cF'))\rightarrow \rR^q\lgg_{n*}(\pi_n^*(\cF')),
\end{equation}
where $\sigma_n^*$ and $\pi_n^*$ denote the pullbacks in the sense of ringed topos.

\begin{teo}[\cite{ag1} 6.5.31]\label{p2-bch2}
Assume that the morphism $g\colon X'\rightarrow X$ is proper.
Then, there exists an integer $N\geq 0$ such that for all integers $n\geq 1$ and $q\geq 0$ and
every quasi-coherent $\co_{\oX'_n}$-module $\cF'$ of $X'_{s,\zar}$,
that we also consider as an $\co_{\oX'_n}$-module of $X'_{s,\et}$,
the kernel and the cokernel of the base change morphism \eqref{p2-bch1b}
\begin{equation}\label{p2-bch2a}
\sigma_n^*(\rR^q\ogg_{n*}(\cF'))\rightarrow \rR^q\lgg_{n*}(\pi_n^*(\cF'))
\end{equation}
are annihilated by $p^N$.
\end{teo}

\subsection{}\label{p2-bch3}
The diagram of morphisms of ringed topos \eqref{p2-rftchta4a}
\begin{equation}\label{p2-bch3a}
\xymatrix{
{(\tG^{\mN^\circ}_s,\bvocB^!)}\ar[r]^-(0.5){\bvpi}\ar[d]_-(0.5){\bvlgg}&{ (X'^{\mN^\circ}_{s,\et},\co_{\bvoX'})}\ar[d]^{\bvogg}\\
{(\tE^{\mN^\circ}_s,\bvocB)}\ar[r]^-(0.5){\bvsigma}&{(X^{\mN^\circ}_{s,\et },\co_{\bvoX})}}
\end{equation}
is commutative up to canonical isomorphism.
For every $\co_{\bvoX'}$-module $\cF'$ of $X'^{\mN^\circ}_{s,\et}$ and every integer $q\geq 0$,
we have a canonical base change morphism (\cite{egr1} (1.2.3.3))
\begin{equation}\label{p2-bch3b}
\bvsigma^*(\rR^q\bvogg_*(\cF'))\rightarrow \rR^q\bvlgg_*(\bvpi^*(\cF')),
\end{equation}
where $\bvsigma^*$ and $\bvpi^*$ denote the pullbacks in the sense of ringed topos.

\begin{prop}\label{p2-bch4}
Assume that the morphism $g\colon X'\rightarrow X$ is proper.
Then, there exists an integer $N\geq 0$ such that for every $\co_{\bvoX'}$-module $\cF'=(\cF'_n)_{n\geq 0}$ of $X'^ {\mN^\circ}_{s,\et}$,
where the $\co_{\oX'_n}$-modules $\cF'_n$ are induced by quasi-coherent $\co_{\oX'_n}$-modules of $X'_{s,\zar}$,
and every $q\geq 0$, the kernel and the cokernel of the base change morphism \eqref{p2-bch3b}
\begin{equation}\label{p2-bch4a}
\bvsigma^*(\rR^q \bvogg_*(\cF'))\rightarrow \rR^q\bvlgg_*(\bvpi^*(\cF'))
\end{equation}
are annihilated by $p^N$.
\end{prop}
This follows from \ref{p2-bch2}, (\cite{agt} III.7.3(i)) and (\cite{ag2} 6.2.10).

\subsection{}\label{p2-bch5}
The diagram of morphisms of ringed topos
\begin{equation}\label{p2-bch5a}
\xymatrix{
{(X'^{\mN^\circ}_{s,\et},\co_{\bvoX'})}\ar[r]^-(0.5){\bvu'}\ar[d]_ {\bvogg_\et}
&{(X'^{\mN^\circ}_{s,\zar},\co_{\bvoX'})}\ar[d] ^{\bvogg_\zar}\\
{(X^{\mN^\circ}_{s,\et},\co_{\bvoX})}\ar[r]^-(0.5)\bvu&
{(X^{\mN^\circ}_{s,\zar},\co_{\bvoX}),}}
\end{equation}
where $\bvu$ and $\bvu'$ are the canonical morphisms of ringed topos \eqref{p2-htaft7f}, is commutative up to canonical isomorphism.

\begin{lem}[\cite{ag1} 6.5.37]\label{p2-bch6}
Assume that the morphism $g\colon X'\rightarrow X$ is separated and quasi-compact.
Let $\cF'=(\cF'_n)_{n\in \mN}$ be an $\co_{\bvoX'}$-module of $X'^{\mN^\circ}_{s,\zar}$ such that for every integer $n\geq 0$,
the $\co_{\oX'_n}$-module $\cF'_n$ is quasi-coherent,
$q$ an integer $\geq 0$. Then, the base change morphism with respect to the diagram \eqref{p2-bch5a}
\begin{equation}\label{p2-bch6a}
\bvu^*(\rR^q\bvogg_{\zar*}(\cF'))\rightarrow \rR^q\bvogg_{\et*}(\bvu'^*(\cF'))
\end{equation}
is an isomorphism.
\end{lem}

\subsection{}\label{p2-bch7}
The diagram of morphisms of ringed topos \eqref{p2-rftchta4e}
\begin{equation}\label{p2-bch7a}
\xymatrix{
{(\tG^{\mN^\circ}_s,\bvocB^!)}\ar[r]^-(0.5){\huppi}\ar[d]_-(0.5){\bvlgg}&{ (X'_{s,\zar},\co_{\fX'})}\ar[d]^{\fgg}\\
{(\tE^{\mN^\circ}_s,\bvocB)}\ar[r]^-(0.5){\hupsigma}&{(X_{s,\zar},\co_{\fX})}}
\end{equation}
is commutative up to canonical isomorphism. We use the notation/conventions of \ref{p1-bcim6} and \ref{p1-bcim1} for these morphisms of ringed topos. 

Let $q$ be an integer. For every bounded from below complex of $\co_{\fX'}$-modules $\cM^\bullet$, 
there exists a canonical functorial base change morphism of $\bvocB$-modules \eqref{p1-bcim4d}, with respect to \eqref{p2-bch7a},
\begin{equation}\label{p2-bch7b}
\hupsigma^*(\rR^q\fgg_*(\cM^\bullet))\rightarrow \rR^q\bvlgg_*(\huppi^*(\cM^\bullet)),
\end{equation}
where the pullback $\huppi^*(\cM^\bullet)$ is defined term by term (not derived). 
By \ref{p1-bcim12}, for every $\co_{\fX'}$-module $\cM$, the base change morphism \eqref{p2-bch7b} 
for $\cM[0]$ coincides with the classical base change morphism (\cite{egr1} 1.2.3).

For every bounded from below complex of $\co_{\fX',\mQ}$-modules $M^\bullet$, 
there exists a canonical functorial base change morphism of $\bvocB_\mQ$-modules \eqref{p1-bcim4e}, with respect to \eqref{p2-bch7a},
\begin{equation}\label{p2-bch7c}
\hupsigma^*_\mQ(\rR^q\fgg_{\mQ *}(M^\bullet))\rightarrow \rR^q\bvlgg_{\mQ*}(\huppi^*_\mQ(M^\bullet)),
\end{equation}
where the pullback $\huppi^*_\mQ(M^\bullet)$ is defined term by term (not derived). 

For every bounded from below complex of ind-$\co_{\fX'}$-modules $\cF^\bullet$, 
there exists a canonical functorial base change morphism of ind-$\bvocB$-modules \eqref{p1-bcim4f}, with respect to \eqref{p2-bch7a},
\begin{equation}\label{p2-bch7d}
\rI\hupsigma^*(\rR^q\rI \fgg_*(\cF^\bullet))\rightarrow \rR^q\rI\bvlgg_*(\rI\huppi^*(\cF^\bullet)),
\end{equation}
where the pullback $\rI\huppi^*(\cF^\bullet)$ is defined term by term (not derived). 

\begin{prop}[\cite{ag2} 6.2.6]\label{p2-bch8}
The composed functors 
\begin{eqnarray}
\xymatrix{
{\bMod^\coh(\co_\fX[\frac 1 p])}\ar[r]^-(0.5){\upalpha^\coh_\fX}&{\bIndMod(\co_\fX)}\ar[r]^-(0.5){\rI\hupsigma^*}&{\bIndMod(\bvocB),}}\label{p2-bch8a}\\
\xymatrix{
{\bMod^\coh(\co_{\fX'}[\frac 1 p])}\ar[r]^-(0.5){\upalpha^\coh_{\fX'}}&{\bIndMod(\co_{\fX'})}\ar[r]^-(0.5){\rI\hupsigma'^*}&{\bIndMod(\bvocB'),}}\label{p2-bch8b}\\ 
\xymatrix{
{\bMod^\coh(\co_{\fX'}[\frac 1 p])}\ar[r]^-(0.5){\upalpha^\coh_{\fX'}}&{\bIndMod(\co_{\fX'})}\ar[r]^-(0.5){\rI\huppi^*}&{\bIndMod(\bvocB^!),}}\label{p2-bch8c}
\end{eqnarray}
where the functors $\upalpha^\coh_\fX$ and $\upalpha^\coh_{\fX'}$ are defined in \eqref{p2-cmupiso2h}, are exact.
\end{prop}

\begin{cor}\label{p2-bch13}
Let $M^\bullet$ be a bounded from below complex of ind-$\co_{\fX'}$-modules, 
$(\cU'_\lambda)_{\lambda\in \Lambda}$ a Zariski open covering of $\fX'$. 
We assume that for every $\lambda\in \Lambda$, 
there exists a bounded from below complex of coherent $\co_{\cU'_\lambda}[\frac 1 p]$-modules $N_\lambda^\bullet$ and a morphism
\begin{equation}\label{p2-bch13a}
u_\lambda\colon \upalpha^\coh_{\co_{\cU'_\lambda}}(N_\lambda^\bullet)\rightarrow M^\bullet|\cU'_\lambda,
\end{equation}
where both the source and the target are defined term by term, whose mapping cone is homotopically equivalent to zero. Then, 
for every integer $q$, $\cH^q(M^\bullet)$ is in the essential image of $\upalpha^\coh_{\co_{\fX'}}$ and the canonical morphism \eqref{p1-bcim13a}
\begin{equation}\label{p2-bch13b}
\rI \huppi^*(\cH^q(M^\bullet))\rightarrow \cH^q(\rI \huppi^*M^\bullet)
\end{equation}
is an isomorphism.  
\end{cor}

The first part is proved in \ref{p2-cmupiso20} and the second part is proved similarly 
using the exactness of the composed functor \eqref{p2-bch8c} instead of \ref{p2-cmupiso7}. 

\begin{prop}\label{p2-bch9}
Assume that the morphism $g\colon X'\rightarrow X$ is proper.
Let $\cM$ be a coherent $\co_{\fX'}$-module, $q$ an integer.
Then, there exists an integer $N\geq 0$ such that the kernel and the cokernel of the base change morphism \eqref{p2-bch7b}
\begin{equation}\label{p2-bch9a}
\hupsigma^*(\rR^q\fgg_*(\cM))\rightarrow \rR^q\bvlgg_*(\huppi^*(\cM))
\end{equation}
are annihilated by $p^N$.
\end{prop}

Indeed, we have $\huppi= \uplambda'\circ\bvu'\circ\bvpi$, $\hupsigma=\uplambda\circ\bvu\circ \bvsigma$ and the diagram 
\eqref{p2-bch7a} can be decomposed into 3 commutative squares 
\begin{equation}
\xymatrix{
{(\tG^{\mN^\circ}_s,\bvocB^!)}\ar[r]^-(0.5){\bvpi}\ar[d]_-(0.5){\bvlgg}&{(X'^{\mN^\circ}_{s,\et},\co_{\bvoX'})}\ar[r]^-(0.5){\bvu'}\ar[d]_ {\bvogg}
&{(X'^{\mN^\circ}_{s,\zar},\co_{\bvoX'})}\ar[r]^-(0.5){\uplambda'}\ar[d] _{\bvogg}&{(X'_{s,\zar},\co_{\fX'})}\ar[d]^{\fgg}\\
{(\tE^{\mN^\circ}_s,\bvocB)}\ar[r]^-(0.5){\bvsigma}&{(X^{\mN^\circ}_{s,\et},\co_{\bvoX})}\ar[r]^-(0.5)\bvu&
{(X^{\mN^\circ}_{s,\zar},\co_{\bvoX})}\ar[r]^-(0.5)\uplambda&{(X_{s,\zar},\co_{\fX}),}}
\end{equation}
where $\uplambda$ and $\uplambda'$ are the canonical morphism of ringed topos \eqref{p2-cmupiso5b}. 
The proposition then follows from \ref{p2-cmupiso15}, \ref{p2-bch4} and \ref{p2-bch6}, taking into account (\cite{egr1} 1.2.4(ii)) and 
(\cite{ag1} 2.6.3 and (2.1.18.6)).

\begin{cor}\label{p2-bch10}
Assume that the morphism $g\colon X'\rightarrow X$ is proper.
Let $M$ be a coherent $\co_{\fX'}[\frac 1 p]$-module, $q$ an integer.
Then, 
\begin{itemize}
\item[{\rm (i)}]  The $\co_\fX[\frac 1 p]$-module $\rR^q\fgg_*(M)$, is coherent. 
We consider $M$ (resp.\ $\rR^q\fgg_*(M)$) as an $\co_{\fX',\mQ}$-module (resp.\ $\co_{\fX,\mQ}$-module) 
via $\ojmath_{\co_{\fX}}$ (resp.\ $\ojmath'_{\co_{\fX'}}$) \eqref{p2-cmupiso2b}. 
\item[{\rm (ii)}] The base change morphism \eqref{p2-bch7c}
\begin{equation}\label{p2-bch10a}
\hupsigma^*_\mQ(\rR^q\fgg_{\mQ *}(M))\rightarrow \rR^q\bvlgg_{\mQ *}(\huppi^*_\mQ(M))
\end{equation}
is an isomorphism.
\end{itemize}
\end{cor}

It follows from \ref{p2-bch9} and \eqref{p1-bcim9b}.

\begin{prop}\label{p2-bch11}
Assume that the morphism $g\colon X'\rightarrow X$ is proper.
Let $\cM^\bullet$ be a bounded from below complex of ind-$\co_{\fX'}$-modules satisfying the following conditions:
\begin{itemize}
\item[{\rm (i)}] For every integer $q$, $\cH^q(\cM^\bullet)$ is in the essential image of $\upalpha_{\co_{\fX'}}^\coh$ \eqref{p2-cmupiso2c}.
\item[{\rm (ii)}] For every integer $q$, the canonical morphism \eqref{p1-bcim13a}
\begin{equation}\label{p2-bch11a}
\rI \huppi^*(\cH^q(\cM^\bullet))\rightarrow \cH^q(\rI \huppi^*(\cM^\bullet))
\end{equation}
is an isomorphism.  
\end{itemize}
Then, for every integer $q$, the base change morphism \eqref{p2-bch7d}
\begin{equation}\label{p2-bch11b}
\rI\hupsigma^*(\rR^q\rI \fgg_*(\cM^\bullet))\rightarrow \rR^q\rI\bvlgg_*(\rI\huppi^*(\cM^\bullet))
\end{equation}
is an isomorphism. 
\end{prop}

The proof is similar to that of \ref{p2-cmupiso17}.
It follows from \ref{p1-bcim14}(ii) applied to \eqref{p2-bch7a} and the thick subcategory $\bMod^{\coh}(\co_\fX[\frac 1 p])$ of 
$\bIndMod(\co_\fX)$ \eqref{p2-cmupiso3}. Recall that composed functor \eqref{p2-bch8a} is exact and observe that for all integers $q$ and $r$, 
$\rR^q\rI \fgg_*(\cH^r (\cM^\bullet))$ is in the essential image of $\upalpha_{\co_{\fX}}^\coh$
by \ref{p2-cmupiso16}(i), \eqref{p2-cmupiso18b} and \eqref{p1-bcim1f}; and the base change morphism \eqref{p2-bch7d}
\begin{equation}
\rI\hupsigma^*(\rR^q\rI \fgg_*(\cH^r(\cF^\bullet)))\rightarrow \rR^q\rI\bvlgg_*(\rI\huppi^*(\cH^r(\cF^\bullet)))
\end{equation}
is an isomorphism by \ref{p2-bch10}(ii) and \eqref{p1-bcim9f}.

\begin{teo}\label{p2-bch12} 
We assume that the morphism $g\colon X'\rightarrow X$ is proper \eqref{p2-fhtft1a}. 
Let $(N,\theta)$ be a {\em locally CL-small} Higgs $\co_{\fX'}[\frac 1 p]$-module with coefficients in $\hOmega'$ \eqref{p2-fhtft2} \eqref{p1-tshbn13}, 
such that the $\co_{\fX'}[\frac 1 p]$-module $N$ is flat, $q$ an integer. 
We set $(\cN,\vartheta)=\upalpha^\coh_{\co_{\fX'}}(N,\theta)$ \eqref{p2-cmupiso31h}  
and denote by $\mK^{\bullet}(\cN\otimes_{\co_{\fX'}}\IC^\dagger_\uptau)$ the Dolbeault complex of $\cN\otimes_{\co_{\fX'}}\IC^\dagger_\uptau$ 
equipped with the total Higgs field $\vartheta\otimes \id+\id\otimes \Idelta'_\uptau$ \eqref{p2-fhtft6j}. 
Then, the base change morphism \eqref{p2-bch7d}
\begin{equation}\label{p2-bch12a}
\rI\hupsigma^*(\rR^q\rI \fgg_*(\mK^\bullet(\cN\otimes_{\co_{\fX'}}\IC^\dagger_\uptau)))\rightarrow \rR^q\rI\bvlgg_*(\rI\huppi^*(\mK^\bullet(\cN\otimes_{\co_{\fX'}}\IC^\dagger_\uptau)))
\end{equation}
is an isomorphism of ind-$\bvocB$-modules. 
\end{teo}

It follows from \ref{p2-bch11}, applied to $\mK^\bullet(\cN\otimes_{\co_{\fX'}}\IC^\dagger)$, whose conditions are satisfied by 
the exactness of the functor \eqref{p2-bch8c}, \ref{p2-bch13} and \ref{p2-cmupiso28}(ii). 

\section{Relative cohomology of Dolbeault ind-modules}\label{p2-rcdim}

\subsection{}\label{p2-rcdim1}
The assumptions and notation of §\ref{p2-fhtft}, §\ref{p2-rfhtft} and §\ref{p2-rftchta} remain in force throughout this section.
In particular, the morphism $g\colon (X',\cM_{X'})\rightarrow (X,\cM_X)$ \eqref{p2-fhtft1a} is assumed to be {\em smooth and saturated}. 
We consider again the objects associated with $(f,\tf)$ introduced in §\ref{p2-rgpsc}, and we associate with $(f',\tf')$ \eqref{p2-fhtft1d}
similar objects that we denote by the same symbols equipped with a $^\prime$ exponent.
In particular, we denote by $\bIndMod^\Dolb(\bvocB')$ the category of Dolbeault ind-$\bvocB'$-modules
and by $\bHM^\sol(\co_{\fX'}[\frac 1 p], \hOmega')$
the category of solvable Higgs $\co_{\fX'}[\frac 1 p]$-bundles with coefficients in $\hOmega'$, see \ref{p2-rgpsc12} and \ref{p2-fhtft2},
with respect  to the deformation $\tf'$ fixed in \eqref{p2-fhtft1d}. We denote by
\begin{equation}\label{p2-rcdim1a}
\cH'\colon \bIndMod(\bvocB')\rightarrow \bHM(\co_{\fX'}, \hOmega')
\end{equation}
the functor defined in \eqref{p2-rgpsc14a}, associated with $(f',\tf')$.
By \ref{p2-rgpsc15}, this induces an equivalence of categories that we denote again by
\begin{equation}\label{p2-rcdim1b}
\cH'\colon \bIndMod^\Dolb(\bvocB')\stackrel{\sim}{\rightarrow} \bHM^\sol(\co_{\fX'}[\frac 1 p], \hOmega').
\end{equation}

\subsection{}\label{p2-rcdim5}
Let $(\cM,\theta)$ be a Higgs ind-$\bvocB'$-module with coefficients in $\hupsigma'^*(\hOmega')$ \eqref{p2-fhtft2}, 
$\ur=(r_1,r_2,r_3)\in I$ \eqref{p2-fhtft33}. We denote by $\mK^\bullet(\cM\otimes_{\bvocB'}\bvcC'^{\ur})$ the Dolbeault complex of 
$(\cM\otimes_{\bvocB'}\bvcC'^{\ur},\theta\otimes\id-\id\otimes\bvdelta^\ur)$ \eqref{p2-rftchta17a} and by 
$\mK^\bullet(\cM\otimes_{\bvocB'}\hupsigma'^*(\hcC^{(r_3)}_{\uptau}))$ 
the Dolbeault complex of $(\cM\otimes_{\bvocB'}\hupsigma'^*(\hcC^{(r_3)}_{\uptau}),
\theta\otimes\id+\id\otimes \hupsigma'^*(\delta'_{\hcC^{(r_3)}_{\uptau}}))$ \eqref{p2-fhtft6bb}. 
In view of \eqref{p2-rftchta17l}, we can consider $\mK^\bullet(\cM\otimes_{\bvocB'}\bvcC'^{\ur})$ 
as a complex of ind-$\bvuptheta^*(\bvcC^{(r_2)})$-modules via the homomorphism $\bvupphi^{\ur}$ \eqref{p2-fhtft38f}. 
To lighten the notation, we set 
\begin{equation}\label{p2-rcdim5a}
\cK^\bullet(\cM,r_2,r_3)=\bvuptheta^*(\bvcC^{(r_2)})\otimes_{\bvocB'}\mK^\bullet(\cM\otimes_{\bvocB'}\hupsigma'^*(\hcC^{(r_3)}_{\uptau})),
\end{equation}
where the tensor product is defined term by term.
By \eqref{p2-rftchta17k}, the homomorphism $\bvuppsi^\ur$ \eqref{p2-fhtft38g} induces a morphism of complexes of ind-$\bvuptheta^*(\bvcC^{(r_2)})$-modules
\begin{equation}\label{p2-rcdim5b}
\iota^\ur\colon \cK^\bullet(\cM,r_2,r_3)\rightarrow  
\mK^\bullet(\cM\otimes_{\bvocB'}\bvcC'^{\ur}). 
\end{equation}

The Higgs field $\id\otimes \hupsigma'^*(\delta_{\hcC^{(r_3)}_{\uptau}})$ \eqref{p2-fhtft6b} on 
\[
\cM\otimes_{\bvocB'}\bvcC'^{\ur}=\cM\otimes_{\bvocB'}\bvcC'^{(r_1,r_2)}\otimes_{\ocB'}\hupsigma'^*(\hcC^{(r_3)}_{\uptau})
\] 
induces a morphism of complexes of $\bvocB'$-modules
\begin{equation}\label{p2-rcdim5c}
\fd^{\ur} \colon \mK^\bullet(\cM\otimes_{\bvocB'}\bvcC'^{\ur})\rightarrow \hupsigma'^*(\fgg^*(\hOmega))\otimes_{\bvocB'} \mK^\bullet(\cM\otimes_{\bvocB'}\bvcC'^{\ur}).
\end{equation} 
Indeed, the question being local on $X'$, we may assume that the $\co_\coX$-module $\Omega$ \eqref{p2-fhtft9c} is free of finite type. 
The assertion then follows by considering  the components, with respect to a basis of $\fgg^*(\hOmega)$, 
of the derivation $\delta_{\hcC^{(r_3)}_{\uptau}}$, which are $\co_{\fX'}$-linear endomorphisms of $\hcC^{(r_3)}_{\uptau}$ 
commuting with each other since $\delta_{\hcC^{(r_3)}_{\uptau}}$ is a Higgs field. Recall here that we have 
$\bvdelta^\ur=\delta_{\bvcC'^{(r_1,r_2)}}\otimes \id- \id \otimes \hupsigma'^*(\delta'_{\hcC^{(r_3)}_{\uptau}})$ \eqref{p2-rftchta17ab}.  

Similarly, the Higgs field $\id\otimes \bvuptheta^*(\delta_{\bvcC^{(r_2)}}) \otimes \id + \id\otimes \hupsigma'^*(\delta_{\hcC^{(r_3)}_{\uptau}})$ \eqref{p2-htaft19b} on 
\[
\cM\otimes_{\bvocB'}\bvuptheta^*(\bvcC^{(r_2)})\otimes_{\bvocB'}\hupsigma'^*(\hcC^{(r_3)}_{\uptau})
\] 
induces a morphism of complexes of $\bvocB'$-modules
\begin{equation}\label{p2-rcdim5d}
\partial^{(r_2,r_3)} \colon \cK^\bullet(\cM,r_2,r_3)\rightarrow \hupsigma'^*(\fgg^*(\hOmega))\otimes_{\bvocB'} \cK^\bullet(\cM,r_2,r_3).
\end{equation} 
By \eqref{p2-rftchta17m}, the diagram 
\begin{equation}\label{p2-rcdim5e}
\xymatrix{
{\cK^\bullet(\cM,r_2,r_3)}\ar[r]^-(0.5){\partial^{(r_2,r_3)}}\ar[d]_{\iota^\ur}&
{\hupsigma'^*(\fgg^*(\hOmega))\otimes_{\bvocB'} \cK^\bullet(\cM,r_2,r_3)}\ar[d]^{\id\otimes \iota^\ur}\\
{\mK^\bullet(\cM\otimes_{\bvocB'}\bvcC'^{\ur})}\ar[r]^-(0.5){\fd^{\ur}}&
{\hupsigma'^*(\fgg^*(\hOmega))\otimes_{\bvocB'} \mK^\bullet(\cM\otimes_{\bvocB'}\bvcC'^{\ur})}}
\end{equation}
is commutative. 

\begin{lem}\label{p2-rcdim7}
We keep the assumptions and notation of \ref{p2-rcdim5} and assume moreover that $\theta=0$ and that 
the  ind-$\bvocB'$-module $\cM$ is rational \eqref{p2-htaft26} and flat {\rm (\cite{ag2} 2.7.9)}. 
Let $\ut=(t_1,t_2,t_3)$ be a triple of rational numbers such that $t_1\geq t_2 > t_3 \geq 0$, 
$I_\ut$ the subset of elements $\ur=(r_1,r_2,r_3)$ of $I$ \eqref{p2-fhtft33} 
such that $r_1>t_1$, $r_2\geq t_2$  and $r_3>t_3$. Then, the morphism of complexes of ind-$\bvocB'$-modules
\begin{equation}\label{p2-rcdim7a}
\cM\otimes_{\bvocB'}\bvuptheta^*(\bvcC^{(t_2)})[0] \rightarrow  
\underset{\underset{\ur\in I_\ut}{\longrightarrow}}{\mlq\mlq\lim \mrq\mrq} \ 
\mK^\bullet(\cM\otimes_{\bvocB'}\bvcC'^{\ur})
\end{equation}
induced by $\bvupphi^{\ur}$ \eqref{p2-fhtft38f}, is a quasi-isomorphism.
\end{lem}

Indeed, by \ref{p2-rftchta20}, the canonical morphism of complexes of ind-$\bvocB'$-modules
\begin{equation}\label{p2-rcdim7b}
\bvuptheta^*(\bvcC^{(t_2)})_\mQ[0] \rightarrow  
\underset{\underset{\ur\in I_\ut}{\longrightarrow}}{\mlq\mlq\lim \mrq\mrq} \ 
\mK^\bullet_\mQ(\bvcC'^{\ur})
\end{equation}
is a quasi-isomorphism. Moreover, $\cM$ being rational, for every $\bvocB'$-module $\cF$, the canonical morphism
$\cM\otimes_{\bvocB'}\cF\rightarrow \cM\otimes_{\bvocB'}\cF_\mQ$ (\cite{ag2} (2.6.6.7)) is an isomorphism. The proposition follows since $\cM$ is flat.

\subsection{}\label{p2-rcdim2}
Let $(N,\theta)$ be a Higgs $\co_{\fX'}[\frac 1 p]$-bundle with coefficients in $\hOmega'$ \eqref{p2-fhtft2}, 
$\ur=(r_1,r_2,r_3)\in I$ \eqref{p2-fhtft33}. 
We set $(\cN,\vartheta)=\upalpha^\coh_{\co_{\fX'}}(N,\theta)$ \eqref{p2-cmupiso31h} 
and denote by $\mK^\bullet(\cN\otimes_{\co_{\fX'}}\hcC^{(r_3)}_{\uptau})$
the Dolbeault complex of $(\cN\otimes_{\co_{\fX'}}\hcC^{(r_3)}_{\uptau},\vartheta\otimes\id+\id\otimes \delta'_{\hcC^{(r_3)}_{\uptau}})$ \eqref{p2-fhtft6bb}
and by $\mK^\bullet(\cN\otimes_{\co_{\fX'}}\IC^\dagger_{\uptau})$ 
the Dolbeault complex of $(\cN\otimes_{\co_{\fX'}}\IC^\dagger_{\uptau},\vartheta\otimes\id+\id\otimes \Idelta'_{\uptau})$ \eqref{p2-fhtft6j}.
Following \ref{p2-rcdim5}, we denote by $\mK^\bullet(\rI\hupsigma'^*(\cN)\otimes_{\bvocB'}\bvcC'^{\ur})$ the Dolbeault complex of 
$(\rI\hupsigma'^*(\cN)\otimes_{\bvocB'}\bvcC'^{\ur},\rI\hupsigma'^*(\vartheta)\otimes\id-\id\otimes\bvdelta^\ur)$, 
that we also consider as a complex of ind-$\bvuptheta^*(\bvcC^{(r_2)})$-modules via the homomorphism $\bvupphi^{\ur}$ \eqref{p2-fhtft38f}. 

Recall that for every bounded from below complex $C^\bullet$ of ind-$\bvocB^!$-modules, 
we have canonical morphisms 
\begin{equation}\label{p2-rcdim2f}
C^\bullet\rightarrow \rI\bvvarpi_*(\rI\bvvarpi^*(C^\bullet))\rightarrow \rR\rI\bvvarpi_*(\rI\bvvarpi^*(C^\bullet)),
\end{equation}
where the morphism of ringed topos $\bvvarpi$ is defined in \eqref{p2-rftchta4e}, the first arrow is induced by the adjunction morphism
$\id\rightarrow \rI\bvvarpi_*\circ \rI\bvvarpi^*$, defined term by term, and the second arrow by the canonical morphism 
$\rI\bvvarpi_*\rightarrow \rR \rI\bvvarpi_*$, 
where $\rI\bvvarpi^*$ is defined for complexes term by term \eqref{p1-bcim15e}. 

We have a canonical morphism of complexes de ind-$\bvuptheta^*(\bvcC^{(r_2)})$-modules \eqref{p2-rcdim5b}
\begin{equation}\label{p2-rcdim2b}
\iota^\ur\colon \bvuptheta^*(\bvcC^{(r_2)})\otimes_{\bvocB'}\rI\hupsigma'^*(\mK^\bullet(\cN\otimes_{\co_{\fX'}}\hcC^{(r_3)}_{\uptau}))\rightarrow  
\mK^\bullet(\rI\hupsigma'^*(\cN)\otimes_{\bvocB'}\bvcC'^{\ur}),
\end{equation}
where on the left, the tensor product and $\rI\hupsigma'^*$ are defined term by term. 
Applying the functor $\rR\rI\bvvarpi_*$ to this morphism and composing with the morphism 
\eqref{p2-rcdim2f}, we obtain a canonical morphism of $\bD^+(\bIndMod(\bvlgg^*(\bvcC^{(r_2)})))$
\begin{equation}\label{p2-rcdim2g}
\bvlgg^*(\bvcC^{(r_2)})\otimes_{\bvocB^!}\rI\huppi^*(\mK^\bullet(\cN\otimes_{\co_{\fX'}}\hcC^{(r_3)}_{\uptau}))\rightarrow  
\rR\rI\bvvarpi_*(\mK^\bullet(\rI\hupsigma'^*(\cN)\otimes_{\bvocB'}\bvcC'^{\ur})),
\end{equation}
where on the left, the tensor product and $\rI\huppi^*$ are defined term by term.

Let $\ur'=(r'_1,r'_2,r'_3)\in I$ be such that $r_i\geq r'_i$ for all $1\leq i\leq 3$.
By \eqref{p2-fhtft6d} and \eqref{p2-rftchta17j}, 
the homomorphisms $\halpha_{\uptau}^{r,r'}$ \eqref{p2-fhtft3i} and $\bvalpha'^{\ur,\ur'}$ \eqref{p2-fhtft38n}
induce morphisms of complexes
\begin{eqnarray}
\mK^\bullet(\cN\otimes_{\co_{\fX'}}\hcC^{(r_3)}_{\uptau})&\rightarrow &\mK^\bullet(\cN\otimes_{\co_{\fX'}}\hcC^{(r'_3)}_{\uptau}),\label{p2-rcdim2c}\\
\mK^\bullet(\rI\hupsigma'^*(\cN)\otimes_{\bvocB'}\bvcC'^{\ur})&\rightarrow &
\mK^\bullet(\rI\hupsigma'^*(\cN)\otimes_{\bvocB'}\bvcC'^{\ur'}).\label{p2-rcdim2d}
\end{eqnarray}
Moreover, by \eqref{p2-rftchta17e}, the diagram 
\begin{equation}\label{p2-rcdim2e}
\xymatrix{
{\bvuptheta^*(\bvcC^{(r_2)})\otimes_{\bvocB'}\rI\hupsigma'^*(\mK^\bullet(\cN\otimes_{\co_{\fX'}}\hcC^{(r_3)}_{\uptau}))}
\ar[r]^-(0.5){\iota^\ur}\ar[d]&
{\mK^\bullet(\rI\hupsigma'^*(\cN)\otimes_{\bvocB'}\bvcC'^{\ur})}\ar[d]\\
{\bvuptheta^*(\cC^{(r'_2)})\otimes_{\bvocB'}\rI\hupsigma'^*(\mK^\bullet(\cN\otimes_{\co_{\fX'}}\hcC^{(r'_3)}_{\uptau}))}\ar[r]^-(0.5){\iota^{\ur'}}&
{\mK^\bullet(\rI\hupsigma'^*(\cN)\otimes_{\bvocB'}\bvcC'^{\ur'}),}}
\end{equation}
where the vertical arrows are induced by the morphisms  \eqref{p2-rcdim2c} and \eqref{p2-rcdim2d}, is commutative.  

Let $t_2,t_3$ be two rational numbers such that $t_2\geq t_3\geq 0$, 
$I_{t_2,t_3}$ the subset of elements $\ur=(r_1,r_2,r_3)$ of $I$ such that $r_1>r_2\geq t_2$ and $r_3\geq t_3$. 
Applying the functor $\rR\rI\bvvarpi_*$ to the direct limit of the morphisms \eqref{p2-rcdim2b} over $\ur\in I_{t_2,t_3}$, and composing with the morphism 
\eqref{p2-rcdim2f}, we obtain a canonical morphism of $\bD^+(\bIndMod(\bvlgg^*(\bvcC^{(t_2)})))$
\begin{equation}\label{p2-rcdim2h}
\bvlgg^*(\bvcC^{(t_2)})\otimes_{\bvocB^!}\rI\huppi^*(\mK^\bullet(\cN\otimes_{\co_{\fX'}}\hcC^{(t_3)}_{\uptau}))\rightarrow
\rR\rI\bvvarpi_*(\underset{\underset{\ur\in I_{t_2,t_3}}{\longrightarrow}} {\mlq\mlq\lim\mrq\mrq} \ \mK^\bullet(\rI\hupsigma'^*(\cN)\otimes_{\bvocB'}\bvcC'^{\ur})). 
\end{equation}

Let $t_2$ be a rational number $> 0$, 
$I_{t_2}$ the subset of elements $\ur=(r_1,r_2,r_3)$ of $I$ such that $r_1>r_2\geq t_2$ and $r_3>0$. 
Applying the functor $\rR\rI\bvvarpi_*$ to the direct limit of the morphisms \eqref{p2-rcdim2b} and composing with the morphism 
\eqref{p2-rcdim2f}, we obtain a canonical morphism of $\bD^+(\bIndMod(\bvlgg^*(\bvcC^{(t_2)})))$
\begin{equation}\label{p2-rcdim2i}
\bvlgg^*(\bvcC^{(t_2)})\otimes_{\bvocB^!}\rI\huppi^*(\mK^\bullet(\cN\otimes_{\co_{\fX'}}\IC^\dagger_{\uptau}))\rightarrow
\rR\rI\bvvarpi_*(\underset{\underset{\ur\in I_{t_2}}{\longrightarrow}} {\mlq\mlq\lim\mrq\mrq} \ \mK^\bullet(\rI\hupsigma'^*(\cN)\otimes_{\bvocB'}\bvcC'^{\ur})). 
\end{equation}

\begin{prop}\label{p2-rcdim3}
We keep the assumptions and notation of \ref{p2-rcdim2}.
\begin{itemize}
\item[{\rm (i)}] We have a canonical spectral sequence of ind-$\bvlgg^*(\bvcC^{(r_2)})$-modules, functorial in $\ur$,
\begin{equation}\label{p2-rcdim3a}
{^\ur\rE}_1^{i,j}=\rR^j\rI\bvvarpi_*(\mK^i(\rI\hupsigma'^*(\cN)\otimes_{\bvocB'}\bvcC'^{\ur}))\Rightarrow
\rR^{i+j}\rI\bvvarpi_*(\mK^\bullet(\rI\hupsigma'^*(\cN)\otimes_{\bvocB'}\bvcC'^{\ur})),
\end{equation}
the functoriality being induced by the morphism \eqref{p2-rcdim2d}. Its edge-homomorphism 
\begin{equation}\label{p2-rcdim3f}
{^\ur\rE}_2^{q,0}=\rH^q(\rI\bvvarpi_*(\mK^\bullet(\rI\hupsigma'^*(\cN)\otimes_{\bvocB'}\bvcC'^{\ur})))\rightarrow
\rR^q\rI\bvvarpi_*(\mK^\bullet(\rI\hupsigma'^*(\cN)\otimes_{\bvocB'}\bvcC'^{\ur}))
\end{equation}
is induced by the canonical morphism $\rI\bvvarpi_*\rightarrow \rR\rI\bvvarpi_*$. 
\item[{\rm (ii)}] For every integer $j$, the derivation $\bvdelta^\ur$ \eqref{p2-rftchta17a} induces a Higgs $\bvlgg^*(\bvcC^{(r_2)})$-field, functorial in $\ur$,
\begin{equation}\label{p2-rcdim3b}
\kappa^{j,\ur}\colon \rR^j\bvvarpi_*(\bvcC'^{\ur})\rightarrow \huppi^*(\hOmega')\otimes_{\bvocB^!}\rR^j\bvvarpi_*(\bvcC'^{\ur}).
\end{equation}
We denote by $\mK^\bullet(\rI\huppi^*(\cN)\otimes_{\bvocB^!}\rR^j\bvvarpi_*(\bvcC'^{\ur}))$ the Dolbeault complex of 
\begin{equation}\label{p2-rcdim3c}
(\rI\huppi^*(\cN)\otimes_{\bvocB^!}\rR^j\bvvarpi_*(\bvcC'^{\ur}),\rI\huppi^*(\uptheta)\otimes\id-\id\otimes\kappa^{j,\ur}).
\end{equation}
\item[{\rm (iii)}] For every integer $j$, we have an isomorphism of complexes of ind-$\bvlgg^*(\bvcC^{(r_2)})$-modules
\begin{equation}\label{p2-rcdim3d}
{^\ur\rE}_1^{\bullet,j}\stackrel{\sim}{\rightarrow} \mK^\bullet(\rI\huppi^*(\cN)\otimes_{\bvocB^!}\rR^j\bvvarpi_*(\bvcC'^{\ur})). 
\end{equation}
\end{itemize}
\end{prop}

(i) Observe first that the category of ind-modules does not have enough injectives, so there are no general spectral sequences of hypercohomology 
(\cite{ega3} 0.11.4.3) and (\cite{hodge2} 1.4.5 and 1.4.6).

We set $\ovarphi(N,\theta)=(\fN,\fN',u\colon \fN\rightarrow \fN',\vartheta\colon \fN\rightarrow \fN'\otimes_{\co_{\fX'}}\hOmega')$, where $\ovarphi$ is the functor \eqref{p2-cmupiso31f}, so $\fN$ and $\fN'$ are coherent $\co_{\fX'}$-modules and $u$ and $\vartheta$ are $\co_{\fX'}$-linear morphisms 
for which there exist an integer $n\not=0$ and an $\co_{\fX'}$-linear morphism $v\colon \fN'\rightarrow \fN$ such that
$v\circ u=n\cdot \id_\fN$, $u\circ v=n\cdot \id_{\fN'}$ and such that $(\fN,(v\otimes \id)\circ \vartheta)$ and $(\fN',\vartheta\circ v)$
are Higgs $\co_{\fX'}$-modules with coefficients in $\hOmega'$. We consider the isogenies of $\bvcC'^\ur$-modules 
\begin{eqnarray}
\tu=\hupsigma'^*(u)\otimes\id\colon \hupsigma'^*(\fN)\otimes_{\bvocB'}\bvcC'^\ur\rightarrow \hupsigma'^*(\fN')\otimes_{\bvocB'}\bvcC'^\ur,\\
\tv=\hupsigma'^*(v)\otimes\id\colon \hupsigma'^*(\fN')\otimes_{\bvocB'}\bvcC'^\ur\rightarrow \hupsigma'^*(\fN)\otimes_{\bvocB'}\bvcC'^\ur,
\end{eqnarray}
and the $\bvuptheta^*(\bvcC^{(r_2)})$-linear morphism 
\begin{equation}
\tvartheta=-\hupsigma'^*(u)\otimes\bvdelta^\ur+ \hupsigma'^*(\vartheta)\otimes \id \colon 
\hupsigma'^*(\fN)\otimes_{\bvocB'}\bvcC'^\ur\rightarrow 
\hupsigma'^*(\fN')\otimes_{\bvocB'}\bvcC'^\ur \otimes_{\bvocB'}\hupsigma'^*(\hOmega').
\end{equation}
The quadruple $(\hupsigma'^*(\fN)\otimes_{\bvocB'}\bvcC'^\ur,
\hupsigma'^*(\fN')\otimes_{\bvocB'}\bvcC'^\ur,\tu,\tvartheta)$ is a Higgs $\bvuptheta^*(\bvcC^{(r_2)})$-isogeny (\cite{ag2} 2.9.9)
by \ref{p1-delta-con4}. Let $\cK^\bullet(\ur)$ be its Dolbeault complex in $\bMod_\mQ(\bvuptheta^*(\bvcC^{(r_2)}))$, see (\cite{ag2} (2.9.9.4)). 
Since the abelian category $\bMod_\mQ(\bvuptheta^*(\bvcC^{(r_2)}))$ has enough injectives, 
we have a canonical spectral sequence of $\bvlgg^*(\bvcC^{(r_2)})_\mQ$-modules, functorial in $\ur$,
\begin{equation}\label{p2-rcdim3e}
{^\ur\cE}_1^{i,j}=\rR^j\bvvarpi_{\mQ*}(\cK^i(\ur))\Rightarrow \rR^{i+j}\bvvarpi_{\mQ*}(\cK^\bullet(\ur)),
\end{equation}
where the functor $\bvvarpi_{\mQ*}\colon \bMod_\mQ(\bvuptheta^*(\bvcC^{(r_2)}))\rightarrow \bMod_\mQ(\bvlgg^*(\bvcC^{(r_2)}))$ and its derived functors are defined in \eqref{p1-bcim6}. 

In view of \eqref{p1-bcim1e} and (\cite{ag2} (2.9.7.4)), we have a canonical isomorphism 
\begin{equation}\label{p2-rcdim3g}
\upalpha_{\bvuptheta^*(\bvcC^{(r_2)})}(\cK^\bullet(\ur))\stackrel{\sim}{\rightarrow}\mK^\bullet(\rI\hupsigma'^*(\cN)\otimes_{\bvocB'}\bvcC'^{\ur}),
\end{equation}
where the functor $\upalpha_{\bvuptheta^*(\bvcC^{(r_2)})}$ is defined in \eqref{p1-bcim5c}. 
The spectral sequence \eqref{p2-rcdim3a} is then obtained from \eqref{p2-rcdim3e} 
by \eqref{p2-rcdim3g} and \eqref{p1-bcim1f}. 
The description of the edge-homomorphism \eqref{p2-rcdim3f} follows from the analogous description of the edge-homomorphism 
of \eqref{p2-rcdim3e} (\cite{ega3} 0.11.3.4).

(ii) It follows from (\cite{ag2} 6.2.2(i)), and \eqref{p2-rftchta17e} for the functoriality. 

(iii) It follows from (\cite{ag2} 6.2.2(iv)) and the fact that for every integer, 
the differentials of the complex ${^\ur\cE}_1^{\bullet,j}$ are induced by those of $\cK^\bullet(\ur)$. 

\begin{prop}\label{p2-rcdim4}
We keep the notation and assumptions of \ref{p2-rcdim2}. 
Let $t_2,t_3$ be two rational numbers such that $t_2\geq t_3 \geq 0$, $I_{t_2,t_3}$ the subset of elements $\ur=(r_1,r_2,r_3)$ of $I$ \eqref{p2-fhtft33} 
such that $r_1>r_2\geq t_2$ and $r_3\geq t_3$. Then, the canonical morphism \eqref{p2-rcdim2h} is an isomorphism of 
$\bD^+(\bIndMod(\bvlgg^*(\bvcC^{(t_2)})))$
\begin{equation}\label{p2-rcdim4a}
\bvlgg^*(\bvcC^{(t_2)})\otimes_{\bvocB^!}\rI\huppi^*(\mK^\bullet(\cN\otimes_{\co_{\fX'}}\hcC^{(t_3)}_{\uptau}))\stackrel{\sim}{\rightarrow}
\rR\rI\bvvarpi_*(\underset{\underset{\ur\in I_{t_2,t_3}}{\longrightarrow}} {\mlq\mlq\lim\mrq\mrq} \ \mK^\bullet(\rI\hupsigma'^*(\cN)\otimes_{\bvocB'}\bvcC'^{\ur})). 
\end{equation}
\end{prop}

We take again the notation of \ref{p2-rcdim3}. By \ref{p2-rftchta8}(ii) and \ref{p2-rcdim3}(iii), for every $i\geq 0$, the canonical morphism \eqref{p2-rcdim2b}
induces an isomorphism  
\begin{equation}\label{p2-rcdim4b}
\bvlgg^*(\bvcC^{(t_2)})\otimes_{\bvocB^!}\rI\huppi^*(\mK^i(\cN\otimes_{\co_{\fX'}}\hcC^{(t_3)}_{\uptau}))\stackrel{\sim}{\rightarrow}
\underset{\underset{\ur\in I_{t_2,t_3}}{\longrightarrow}}{\mlq\mlq\lim \mrq\mrq}\ {^\ur\rE}_1^{i,0},
\end{equation}
and for every $j\geq 1$, we have
\begin{equation}\label{p2-rcdim4c}
\underset{\underset{I_{t_2,t_3}}{\longrightarrow}}{\mlq\mlq\lim \mrq\mrq}\ {^\ur\rE}_1^{i,j}=0.
\end{equation}
Moreover, the isomorphisms \eqref{p2-rcdim4b} (for $i\in \mN$) form an isomorphism of complexes.
Since small filtered direct limits exist and are exact in $\bIndMod(\bvlgg^*(\bvcC^{(t_2)}))$ (\cite{ag2} (2.6.7.5)), and in view of (\cite{ag2} (2.7.10.5)),
we deduce from this that the canonical morphism \eqref{p2-rcdim2g} induces for every integer $q$ an isomorphism
\begin{equation}\label{p2-rcdim4d}
\rH^q(\bvlgg^*(\bvcC^{(t_2)})\otimes_{\bvocB^!}\rI\huppi^*(\mK^\bullet(\cN\otimes_{\co_{\fX'}}\hcC^{(t_3)}_{\uptau})))\stackrel{\sim}{\rightarrow}
\underset{\underset{\ur\in I_{t_2,t_3}}{\longrightarrow}} {\mlq\mlq\lim\mrq\mrq} \ \rR^q\rI\bvvarpi_*(\mK^\bullet(\rI\hupsigma'^*(\cN)\otimes_{\bvocB'}\bvcC'^{\ur})). 
\end{equation}
Hence, the morphism \eqref{p2-rcdim2h} is an isomorphism. 

\begin{cor}\label{p2-rcdim9}
We keep the notation and assumptions of \ref{p2-rcdim2}. 
Let $t_2$ be a rational number $>0$, $I_{t_2}$ the subset of elements $\ur=(r_1,r_2,r_3)$ of $I$ \eqref{p2-fhtft33} 
such that $r_1>r_2\geq t_2$ and $r_3>0$. Then, the canonical morphism \eqref{p2-rcdim2i} is an isomorphism of 
$\bD^+(\bIndMod(\bvlgg^*(\bvcC^{(t_2)})))$
\begin{equation}\label{p2-rcdim9a}
\bvlgg^*(\bvcC^{(t_2)})\otimes_{\bvocB^!}\rI\huppi^*(\mK^\bullet(\cN\otimes_{\co_{\fX'}}\IC^\dagger_{\uptau}))\stackrel{\sim}{\rightarrow}
\rR\rI\bvvarpi_*(\underset{\underset{\ur\in I_{t_2}}{\longrightarrow}} {\mlq\mlq\lim\mrq\mrq} \ \mK^\bullet(\rI\hupsigma'^*(\cN)\otimes_{\bvocB'}\bvcC'^{\ur})). 
\end{equation}
\end{cor}

Indeed, by \ref{p2-rcdim4}, for every integer $q$ and every rational number $t_3$ such that $t_2\geq t_3>0$,   
the canonical morphism \eqref{p2-rcdim2g} induces an isomorphism 
\begin{equation}
\rH^q(\bvlgg^*(\bvcC^{(t_2)})\otimes_{\bvocB^!}\rI\huppi^*(\mK^\bullet(\cN\otimes_{\co_{\fX'}}\hcC^{(t_3)}_{\uptau})))\stackrel{\sim}{\rightarrow}
\underset{\underset{\ur\in I_{t_2,t_3}}{\longrightarrow}} {\mlq\mlq\lim\mrq\mrq} \ \rR^q\rI\bvvarpi_*(\mK^\bullet(\rI\hupsigma'^*(\cN)\otimes_{\bvocB'}\bvcC'^{\ur})). 
\end{equation}
Taking the direct limit over $t_3\in ]0,t_2]\cap \mQ$, we deduce that the canonical morphism \eqref{p2-rcdim2g} induces an isomorphism 
\begin{equation}
\rH^q(\bvlgg^*(\bvcC^{(t_2)})\otimes_{\bvocB^!}\rI\huppi^*(\mK^\bullet(\cN\otimes_{\co_{\fX'}}\IC^\dagger_{\uptau})))\stackrel{\sim}{\rightarrow}
\underset{\underset{\ur\in I_{t_2}}{\longrightarrow}} {\mlq\mlq\lim\mrq\mrq} \ \rR^q\rI\bvvarpi_*(\mK^\bullet(\rI\hupsigma'^*(\cN)\otimes_{\bvocB'}\bvcC'^{\ur})). 
\end{equation}
Hence, the morphism \eqref{p2-rcdim2i} is an isomorphism.

\begin{teo}\label{p2-rcdim8}
Let $\cM$ be an ind-$\bvocB'$-module, $(N,\theta)$ a Higgs $\co_{\fX'}[\frac 1 p]$-bundle with coefficients in $\hOmega'$ \eqref{p2-fhtft2}. 
We assume that $\cM$ and $(N,\theta)$ are associated in the sense of \ref{p2-rgpsc11}. 
We set $(\cN,\vartheta)=\upalpha^\coh_{\co_{\fX'}}(N,\theta)$ \eqref{p2-cmupiso31h}
and denote by $\mK^{\bullet}(\cN\otimes_{\co_{\fX'}}\IC^\dagger_\uptau)$ the Dolbeault complex of $\cN\otimes_{\co_{\fX'}}\IC^\dagger_\uptau$ 
equipped with the total Higgs field $\vartheta\otimes \id+\id\otimes \Idelta'_\uptau$ \eqref{p2-fhtft6j}. 
Then, 
\begin{itemize}
\item[{\rm (i)}] There exist a rational number $r>0$ and an isomorphism of $\bD^+(\bIndMod(\bvlgg^*(\bvcC^{(r)})))$ 
\begin{equation}\label{p2-rcdim8a}
\rR\rI\bvvarpi_*(\bvuptheta^*(\bvcC^{(r)})\otimes_{\bvocB'}\cM)
\stackrel{\sim}{\rightarrow}
\bvlgg^*(\bvcC^{(r)})\otimes_{\bvocB^!}\rI\huppi^*(\mK^{\bullet}(\cN\otimes_{\co_{\fX'}}\IC^\dagger_\uptau)),
\end{equation}
where on the right, the tensor product and $\rI\huppi^*$ are defined term by term. 

\item[{\rm (ii)}] The Higgs field $\id\otimes \Idelta_\uptau$ on $\cN\otimes_{\co_{\fX'}}\IC^\dagger_\uptau$ \eqref{p2-fhtft6i}
induces a morphism of complexes of ind-$\co_{\fX'}$-modules
\begin{equation}\label{p2-rcdim8b}
\fd^\dagger \colon \mK^\bullet(\cN\otimes_{\co_{\fX'}}\IC^\dagger_\uptau)\rightarrow \fgg^*(\hOmega)\otimes_{\co_{\fX'}} 
\mK^\bullet(\cN\otimes_{\co_{\fX'}}\IC^\dagger_\uptau).
\end{equation}

\item[{\rm (iii)}] We have a commutative diagram 
\begin{equation}\label{p2-rcdim8c}
\xymatrix{
{\rR\rI\bvvarpi_*(\bvuptheta^*(\bvcC^{(r)})\otimes_{\bvocB'}\cM)}\ar[r]\ar[d]&
{\bvlgg^*(\bvcC^{(r)})\otimes_{\bvocB^!}\rI\huppi^*(\mK^{\bullet}(\cN\otimes_{\co_{\fX'}}\IC^\dagger_\uptau))}\ar[d]\\
{\rR\rI\bvvarpi_*(\bvuptheta^*(\bvcC^{(r)})\otimes_{\bvocB'}\cM)\otimes_{\bvocB^!}\huppi^*(\fgg^*(\hOmega))}\ar[r]&
{\bvlgg^*(\bvcC^{(r)})\otimes_{\bvocB^!}\rI\huppi^*(\mK^{\bullet}(\cN\otimes_{\co_{\fX'}}\IC^\dagger_\uptau))\otimes_{\bvocB^!}\huppi^*(\fgg^*(\hOmega))}} \end{equation}
where the horizontal arrows are induced by \eqref{p2-rcdim8a}, the left vertical arrow by 
$\bvuptheta^*(\delta_{\bvcC^{(r)}})\otimes \id$ \eqref{p2-htaft19b} and the right vertical arrow by 
$\bvlgg^*(\delta_{\bvcC^{(r)}})\otimes \id+\id\otimes \rI \huppi^*(\fd^\dagger)$ \eqref{p2-rcdim8b}. 
\end{itemize}
\end{teo}

(i) Indeed, there exists a rational number $r_0>0$ and an isomorphism of ind-$\bvcC'^{(r_0)}$-modules 
with integrable $\delta_{\bvcC'^{(r_0)}}$-connections, in the sense of \ref{p1-indmal22}, 
\begin{equation}\label{p2-rcdim8d}
\rho^{(r_0)}\colon \cM\otimes_{\bvocB'}\bvcC'^{(r_0)} \stackrel{\sim}{\rightarrow} \rI\hupsigma'^*(\cN)\otimes_{\bvocB'}\bvcC'^{(r_0)},
\end{equation}
where the $\delta_{\bvcC'^{(r_0)}}$-connections are defined as in \ref{p1-delta-con9}, 
$\rI\hupsigma'^*(\cN)$ (resp.\ $\cM$) being endowed with the Higgs field 
$\rI\hupsigma'^*(\vartheta)$ (resp.\ $0$).  
For every $\ur=(r_1,r_2,r_3)\in I$ \eqref{p2-fhtft33} such that $r_1\leq r_0$, we consider the composed homomorphism
\begin{equation}\label{p2-rcdim8e}
{\bvcC'^{(r_0)}}\rightarrow {\bvcC'^{(r_1)}=\bvcC'^{(r_1,r_1)}}\rightarrow{\bvcC'^{\ur}=\bvcC'^{(r_1,r_2)}\otimes_{\ocB'}\hupsigma'^*(\hcC^{(r_3)}_{\uptau}),}
\end{equation}
where the first arrow is $\alpha'^{r_0,r_1}$ \eqref{p2-htaft18f} and the second one is induced by $\alpha'^{r_1,r_1,r_1,r_2}$ \eqref{p2-fhtft28e}.
By \ref{p1-delta-con10}, the isomorphism \eqref{p2-rcdim8d} induces by extension of scalars by \eqref{p2-rcdim8e} an isomorphism of ind-$\bvcC'^{\ur}$-modules 
with integrable $\bvdelta^\ur$-connections \eqref{p2-rftchta17a}
\begin{equation}\label{p2-rcdim8f}
\rho^{\ur}\colon \cM\otimes_{\bvocB'}\bvcC'^{\ur} \stackrel{\sim}{\rightarrow} \rI\hupsigma'^*(\cN)\otimes_{\bvocB'}\bvcC'^{\ur}. 
\end{equation}
We deduce an isomorphism of complexes of ind-$\bvuptheta^*(\bvcC^{(r_2)})$-modules 
\begin{equation}\label{p2-rcdim8g}
\mK^\bullet(\cM\otimes_{\bvocB'}\bvcC'^{\ur}) \stackrel{\sim}{\rightarrow} \mK^\bullet(\rI\hupsigma'^*(\cN)\otimes_{\bvocB'}\bvcC'^{\ur}),
\end{equation}
where the Dolbeault complexes are defined as in \ref{p2-rcdim5}, $\rI\hupsigma'^*(\cN)$ (resp.\ $\cM$) being endowed with the Higgs field 
$\rI\hupsigma'^*(\vartheta)$ (resp.\ $0$). We immediately check that for every $\ur'\in I$ such that $r'_1\leq r_1$, $r'_2\leq r_2$ and $r'_3\leq r_3$, 
the diagram 
\begin{equation}\label{p2-rcdim8h}
\xymatrix{
{\mK^\bullet(\cM\otimes_{\bvocB'}\bvcC'^{\ur})}\ar[r]\ar[d]&{\mK^\bullet(\rI\hupsigma'^*(\cN)\otimes_{\bvocB'}\bvcC'^{\ur})}\ar[d]\\
{\mK^\bullet(\cM\otimes_{\bvocB'}\bvcC'^{\ur'})}\ar[r]&{\mK^\bullet(\rI\hupsigma'^*(\cN)\otimes_{\bvocB'}\bvcC'^{\ur'}),}}
\end{equation}
where the vertical arrows are induced by the homomorphism $\bvalpha'^{\ur,\ur'}$ \eqref{p2-fhtft38n}, 
is commutative. 

Let $t_2$ be a rational number such that $0<t_2<r_0$, $I_{t_2}$ be the subset of elements $\ur=(r_1,r_2,r_3)$ of $I$ \eqref{p2-fhtft33} 
such that $r_0\geq r_1>r_2\geq t_2$ and $r_3>0$. By \ref{p2-rcdim7} and \ref{p2-rcdim9}, the direct limit of the isomorphisms \eqref{p2-rcdim8g}, for $\ur\in I_{t_2}$, 
induces an isomorphism of $\bD^+(\bIndMod(\bvlgg^*(\bvcC^{(t_2)})))$ 
\begin{equation}\label{p2-rcdim8i}
\rR\rI\bvvarpi_*(\bvuptheta^*(\bvcC^{(t_2)})\otimes_{\bvocB'}\cM)
\stackrel{\sim}{\rightarrow}
\bvlgg^*(\bvcC^{(t_2)})\otimes_{\bvocB^!}\rI\huppi^*(\mK^{\bullet}(\cN\otimes_{\co_{\fX'}}\IC^\dagger_\uptau)).
\end{equation}

(ii) Indeed, the question being local on $X'$, we may assume that the $\co_\coX$-module $\Omega$ \eqref{p2-fhtft9c} is free of finite type. 
The assertion then follows by considering  the components, with respect to a basis of $\fgg^*(\hOmega)$, 
of the derivation $\Idelta_\uptau$ \eqref{p2-fhtft6i}, which are $\co_{\fX'}$-linear endomorphisms of $\IC^\dagger_\uptau$ 
commuting with each other since $\Idelta_\uptau$ is a Higgs field.

(iii) For any $\ur=(r_1,r_2,r_3)\in I$,
we denote by $\partial^\ur_\uptau$ the derivation of $\bvcC'^{\ur}=\bvcC'^{(r_1,r_2)}\otimes_{\ocB'}\hupsigma'^*(\hcC^{(r_3)}_{\uptau})$ defined by 
\begin{equation}
\partial^\ur_\uptau=\id\otimes \hupsigma'^*(\delta_{\hcC^{(r_3)}_\uptau})\colon 
\bvcC'^{\ur}\rightarrow \hupsigma'^*(\fgg^*(\hOmega))\otimes_{\bvocB'}\bvcC'^{\ur}.
\end{equation}
The isomorphism of ind-$\bvcC'^{\ur}$-modules $\rho^\ur$ \eqref{p2-rcdim8f} is then compatible with the integrable $\partial^\ur_\uptau$-connections
defined as in \ref{p1-delta-con9}, $\rI\hupsigma'^*(\cN)$ and $\cM$ being endowed with the $0$ Higgs fields. 
We deduce by \eqref{p2-rcdim5e} that diagram \eqref{p2-rcdim8c} is commutative. 
Observe that the left vertical arrow in \eqref{p2-rcdim8c} is well defined by \ref{p1-bcim22}(iii).

\begin{lem}\label{p2-rcdim11}
Let $r$ be a rational number $\geq 0$, 
$\cF$ a complex of $\bD^+(\bIndMod(\bvocB'))$, $\cG$ a complex of $\bD^+(\bIndMod(\bvocB^!))$. 
Then, we have canonical isomorphisms
\begin{eqnarray}
\rR\rI\bvvarpi_*(\cF)\otimes_{\bvocB^!}\bvlgg^*(\bvcC^{(r)})&\stackrel{\sim}{\rightarrow}& 
\rR\rI\bvvarpi_*(\cF\otimes_{\bvocB'}\bvuptheta^*(\bvcC^{(r)})),\label{p2-crindmd12a}\\
\rR\rI\bvlgg_*(\cG)\otimes_{\bvocB}\bvcC^{(r)}&\stackrel{\sim}{\rightarrow}& 
\rR\rI\bvlgg_*(\cG\otimes_{\bvocB^!}\bvlgg^*(\bvcC^{(r)})),\label{p2-crindmd12b}\\
\rR\rI\bvuptheta_*(\cF)\otimes_{\bvocB}\bvcC^{(r)}&\stackrel{\sim}{\rightarrow}& 
\rR\rI\bvuptheta_*(\cF\otimes_{\bvocB'}\bvuptheta^*(\bvcC^{(r)})),\label{p2-crindmd12c}
\end{eqnarray}
where the morphism of ringed topos $\bvvarpi$, $\bvlgg$ and $\bvuptheta$ are defined in \eqref{p2-rftchta4a}.
\end{lem}

Observe first that the tensor products are well defined in view of \ref{p2-htaft41}. 
The requested morphisms are defined as in \eqref{p1-bcim22e}. 
To see that they induce isomorphisms between their $q$th cohomology groups, 
for every integer $q$, we may, by devissage, assume that $\cF$ is an ind-$\bvocB'$-module and $\cG$ is an ind-$\bvocB^!$-module. 
We can further reduce to the case where $\cF$ is a $\bvocB'$-module and $\cG$ is a $\bvocB^!$-module by (\cite{ag2} 2.7.3 and (2.7.10.5)). 
The proposition follows then from (\cite{agt} III.7.5) and (\cite{ag2} 6.5.11).

\subsection{}\label{p2-fhtft8}
With the notation of  \ref{p2-fhtft2}, using the canonical morphism $\hu$ \eqref{p2-fhtft2a},  
we defined in \eqref{p2-cmupiso22d}, for any integer $q\geq 0$, the {\em $q$th higher direct image functor by $\fgg$
twisted by the extension $\hcF_\uptau$ \eqref{p2-fhtft3d}}
(or the {\em twisted $q$th higher direct image functor by $\fgg$} when the extension is implicit) 
\begin{equation}\label{p2-fhtft8a}
\rR^q\fgg^\uptau_*\colon 
\bHM(\co_{\fX'},\hOmega')\rightarrow \bHM(\co_\fX,\hOmega). 
\end{equation}
Similarly, we defined in \eqref{p2-cmupiso23d} the {\em $q$th higher direct image functor by $\rI \fgg$ twisted by the extension 
$\hcF_\uptau$ \eqref{p2-fhtft3d}} (or the {\em twisted $q$th higher direct image functor by $\rI \fgg$} when the extension is implicit)
\begin{equation}\label{p2-fhtft8b}
\rR^q\rI \fgg^\uptau_*\colon 
\bIndHM(\co_{\fX'},\hOmega')\rightarrow \bIndHM(\co_\fX,\hOmega). 
\end{equation}

\begin{teo}\label{p2-rcdim10}
Let $\cM$ be an ind-$\bvocB'$-module, $(N,\theta)$ a Higgs $\co_{\fX'}[\frac 1 p]$-bundle with coefficients in $\hOmega'$ \eqref{p2-fhtft2},
$q$ an integer. We set $(\cN,\vartheta)=\upalpha^\coh_{\co_{\fX'}}(N,\theta)$ \eqref{p2-cmupiso31h}. We assume 
that the following conditions are satisfied: 
\begin{itemize}
\item[{\rm (a)}] $g\colon X'\rightarrow X$ is proper \eqref{p2-fhtft1a}; 
\item[{\rm (b)}] $\cM$ and $(N,\theta)$ are associated in the sense of \ref{p2-rgpsc11};
\item[{\rm (c)}] $(N,\theta)$ is locally CL-small in the sense of \ref{p1-tshbn13}. 
\end{itemize}
Then, 
\begin{itemize}
\item[{\rm (i)}] The Higgs $\co_\fX[\frac 1 p]$-module $\rR^q\fgg^\uptau_*(N,\theta)$ \eqref{p2-fhtft8a} is (locally) CL-small; 
in particular, the underlying $\co_\fX[\frac 1 p]$-module is coherent. 
\item[{\rm (ii)}] We have a canonical functorial isomorphism of Higgs ind-$\co_{\fX}$-modules with coefficients in $\hOmega$, 
\begin{equation}\label{p2-rcdim10a}
\rR^q\rI \fgg^\uptau_*(\cN,\vartheta)\stackrel{\sim}{\rightarrow}
\upalpha^\coh_{\co_{\fX}}(\rR^q\fgg^\uptau_*(N,\theta)),
\end{equation}
where the functor $\rR^q\rI \fgg^\uptau_*$ is defined in \eqref{p2-fhtft8b}. 
\item[{\rm (iii)}] There exist a rational number $r>0$, independent of $q$, and an isomorphism of ind-$\bvcC^{(r)}$-modules with 
$\delta_{\bvcC^{(r)}}$-connection, in the sense of \ref{p1-indmal22}, 
\begin{equation}\label{p2-rcdim10c}
\bvcC^{(r)}\otimes_{\bvocB}\rR^q\rI\bvuptheta_*(\cM)\stackrel{\sim}{\rightarrow}
\bvcC^{(r)}\otimes_{\bvocB}\rI\hupsigma^*(\rR^q\rI\fgg^\uptau_*(\cN,\vartheta)),
\end{equation}
where the $\delta_{\bvcC^{(r)}}$-connections are defined as in \ref{p1-delta-con9}, 
$\rR^q\rI\bvuptheta_*(\cM)$ being endowed with the Higgs field $0$, and $\rI\hupsigma^*$ is the functor \eqref{p2-htaft23a}.  
\end{itemize}
\end{teo}

It follows from \ref{p2-cmupiso30}, \ref{p2-cmupiso47}, \ref{p2-bch12}, \ref{p2-rcdim8} and \ref{p2-rcdim11}.

\begin{rema}\label{p2-rcdim12}
We may replace the condition \ref{p2-rcdim10}(c) by the following one: 
\begin{itemize}
\item[{\rm (c')}] $\cM$ is of the form $\upalpha_{\bvocB'}(M)$, for an adic $\bvocB'_\mQ$-module of finite type $M$ (\cite{ag2} 4.3.14). 
\end{itemize}
Indeed, (b) and (c') imply that the Higgs $\co_{\fX'}[\frac 1 p]$-bundle $(N,\theta)$ is strongly solvable \eqref{p2-rgpsc20}, 
and hence locally CL-small by \ref{p2-rgpsc240}. 
\end{rema}

\begin{cor}\label{p2-rcdim13}
Assume that $g\colon X'\rightarrow X$ is proper \eqref{p2-fhtft1a}. 
Let $\cM$ be a Dolbeault ind-$\bvocB'$-module \eqref{p2-rgpsc11} such that the Higgs $\co_{\fX'}[\frac 1 p]$-module $\cH'(\cM)$ \eqref{p2-rcdim1b}
is locally CL-small \eqref{p1-tshbn13}, $q$ an integer.
Then, the Higgs $\co_\fX[\frac 1 p]$-module $\rR^q\fgg^\uptau_*(\cH'(\cM))$ \eqref{p2-fhtft8a} is (locally) CL-small, and in particular coherent, 
and there exist a rational number $r>0$, independent of $q$, and an isomorphism of ind-$\bvcC^{(r)}$-modules with 
$\delta_{\bvcC^{(r)}}$-connection, in the sense of \ref{p1-indmal22}, 
\begin{equation}\label{p2-rcdim13a}
\bvcC^{(r)}\otimes_{\bvocB}\rR^q\rI\bvuptheta_*(\cM)\stackrel{\sim}{\rightarrow}
\bvcC^{(r)}\otimes_{\bvocB}\rI\hupsigma^*(\upalpha^\coh_{\co_{\fX}}(\rR^q\fgg^\uptau_*(\cH'(\cM)))),
\end{equation}
where the $\delta_{\bvcC^{(r)}}$-connections are defined as in \ref{p1-delta-con9}, 
$\rR^q\rI\bvuptheta_*(\cM)$ being endowed with the Higgs field $0$, and 
$\rI\hupsigma^*$ (resp.\ $\upalpha^\coh_{\co_{\fX}}$) is the functor \eqref{p2-htaft23a} (resp.\ \eqref{p2-cmupiso31h}).  
\end{cor}

It follows from \ref{p2-rgpsc33} and \ref{p2-rcdim10}. 

\begin{rema}
We do not assert that the admissibility isomorphism \eqref{p2-rcdim13a} is uniquely determined by $\cM$, let alone functorial. 
However, after decreasing $r$, we can show that it depends only on $\cM$, on which it depends functorially (see the proof of \cite{ag2} 4.5.20), 
leading to the next corollary. 
\end{rema}

\begin{cor}\label{p2-rcdim14}
Let $\cM$ be a Dolbeault ind-$\bvocB'$-module \eqref{p2-rgpsc11} such that the Higgs $\co_{\fX'}[\frac 1 p]$-module $\cH'(\cM)$ \eqref{p2-rcdim1b}
is locally CL-small \eqref{p1-tshbn13}, $q$ an integer.
Suppose that the morphism $g\colon X'\rightarrow X$ is proper and the $\co_\fX[\frac 1 p]$-module underlying the 
Higgs $\co_\fX[\frac 1 p]$-module $\rR^q\fgg^\uptau_*(\cH'(\cM))$ \eqref{p2-fhtft8a} is locally projective of finite type.
Then, the ind-$\bvocB$-module $\rR^q\rI\bvuptheta_*(\cM)$ is Dolbeault, the Higgs $\co_\fX[\frac 1 p]$-module 
$\rR^q\fgg^\uptau_*(\cH'(\cM))$ is solvable \eqref{p2-rgpsc12} and (locally) CL-small, 
and we have a canonical functorial isomorphism of Higgs $\co_\fX[\frac 1 p]$-bundles
\begin{equation}\label{p2-rcdim14a}
\cH(\rR^q\rI\bvuptheta_*(\cM))\stackrel{\sim}{\rightarrow} \rR^q\fgg^\uptau_*(\cH'(\cM)),
\end{equation}
where $\cH$ is the functor \eqref{p2-rgpsc15a}.
\end{cor}

It follows from \ref{p2-rcdim13} and (\cite{ag2} 4.5.11 and 4.5.14).

\begin{rema}\label{p2-rcdim15}
By \ref{p2-rgpsc241}, we may replace in \ref{p2-rcdim13} and \ref{p2-rcdim14} the condition that $\cM$ is Dolbeault and $\cH'(\cM)$ is locally CL-small 
by the condition that $\cM$ is Dolbeault of the form $\upalpha_{\bvocB'}(M)$ for an adic $\bvocB'_\mQ$-module of finite type $M$, 
i.e. $M$ is strongly Dolbeault \eqref{p2-rgpsc20}.
\end{rema}

\section{\texorpdfstring{Dependence of the $p$-adic Simpson correspondence on the choice of lifting}
{Dependence of the p-adic Simpson correspondence on the choice of lifting}}\label{p2-dpsccd}

\subsection{}\label{p2-dpsccd1}
The assumptions and notation of §\ref{p2-fhtft} remain in force throughout this section.
We assume moreover that $(X',\cM_{X'})=(X,\cM_X)$ and $g=\id_{(X,\cM_X)}$ \eqref{p2-fhtft1a}, 
so $f=f'$; but the deformations $\tf$ and $\tf'$ may be different \eqref{p2-fhtft1d}. 
This case falls within the contexts fixed in §\ref{p2-rfhtft} and §\ref{p2-rftchta}.
Observe that we have $\Omega'=\Omega$ and $\Omega'^{(r,r')}=\Omega^{(r')}$ \eqref{p2-fhtft9j} for every rational numbers $r\geq r'\geq 0$. 

\subsection{}\label{p2-dpsccd2}
Let $\cL_{\tX'/\tX}$ be the torsor of liftings of $\id_{\coX}$ to $\tX'$ over $\tX$, 
\begin{equation}\label{p2-dpsccd2a}
0\rightarrow \co_{\coX}\rightarrow \cF_\uptau\rightarrow \Omega \rightarrow 0
\end{equation}
the associated Higgs--Tate extension \eqref{p2-fhtft3b} and 
\begin{equation}\label{p2-dpsccd2b}
\cC_\uptau\stackrel{\sim}{\rightarrow}\underset{\underset{n\geq 0}{\longrightarrow}}\lim\ \rS^n_{\co_{\coX}}(\cF_\uptau)
\end{equation}
the associated Higgs--Tate $\co_{\coX}$-algebra \eqref{p2-fhtft3c}. 
We denote by $\hcF_\uptau$ the $p$-adic completion of $\cF_\uptau$, so we have an exact sequence of $\co_{\fX}$-modules
\begin{equation}\label{p2-dpsccd2c}
0\rightarrow \co_{\fX}\rightarrow \hcF_\uptau\rightarrow \hOmega \rightarrow 0.
\end{equation}
We denote by 
\begin{equation}\label{p2-dpsccd2d}
\uptau\colon \bHM(\co_{\fX},\hOmega)\rightarrow \bHM(\co_{\fX},\hOmega)
\end{equation}
the twisting functor by the extension $\hcF_\uptau$ \eqref{p1-thbn9b}. 

\subsection{}\label{p2-dpsccd3}
We consider again the objects associated with $(f,\tf)$ introduced in §\ref{p2-rgpsc}, that we denote by the same symbols, 
equipped with a subscript $\tf$, and similarly for those associated with $(f,\tf')$. 
By \ref{p2-htaft57}, the category of Dolbeault ind-$\bvocB$-modules \eqref{p2-rgpsc12}
does not depend on the choice of $\tf$; we denote it by $\bIndMod^\Dolb(\bvocB)$. 
The analogous property for solvable Higgs $\co_{\fX}[\frac 1 p]$-bundles with coefficients in $\hOmega$ is not known (\cite{ag2} 4.10.12). 
We therefore denote by $\bHM^\sol_\tf(\co_{\fX}[\frac 1 p],\hOmega)$ 
the category  of solvable Higgs $\co_{\fX}[\frac 1 p]$-bundles with coefficients in $\hOmega$, relatively to the deformation $\tf$,
and similarly for $\tf'$. We denote by 
\begin{eqnarray}
\cH_\tf\colon \bIndMod^\Dolb(\bvocB)\stackrel{\sim}{\rightarrow} \bHM^\sol_\tf(\co_{\fX}[\frac 1 p], \hOmega),\label{p2-dpsccd3a}\\
\cH_{\tf'}\colon \bIndMod^\Dolb(\bvocB)\stackrel{\sim}{\rightarrow} \bHM^\sol_{\tf'}(\co_{\fX}[\frac 1 p], \hOmega),\label{p2-dpsccd3b}
\end{eqnarray}
the  equivalences of categories associated with $\tf$ and $\tf'$, respectively \eqref{p2-rgpsc15}. 

We take again the notation of \ref{p2-fgpscp1}, adapted to the context of \ref{p2-dpsccd1}. In particular, for any rational numbers $0\leq t\leq r$, 
we have a homomorphism of $\bvocB$-algebras \eqref{p2-fhtft380f}
\begin{equation}\label{p2-dpscc3c}
\bvupphi^{(r,t)}\colon \bvcC^{(r)}_\tf\rightarrow \bvcC^{(r)}_{\tf'}\otimes_{\bvocB}\hupsigma^*(\hcC^{(t)}_\uptau). 
\end{equation}

\begin{prop}\label{p2-dpscc4}
For every Dolbeault ind-$\bvocB$-module $\cM$, 
the $\co_{\fX}[\frac 1 p]$-bundle $\cH_\tf(\cM)$ is locally CL-small, if and only if so is the $\co_{\fX}[\frac 1 p]$-bundle $\cH_{\tf'}(\cM)$. 
Moreover, when this condition is satisfied, we have a canonical functorial isomorphism 
\begin{equation}\label{p2-dpsccd4a}
\cH_{\tf'}(\cM)\stackrel{\sim}{\rightarrow}\uptau(\cH_{\tf}(\cM)).
\end{equation}
\end{prop}

To prove the first assertion, we may assume that $X$ is affine (\cite{ag2} 4.9.4). 
Then, there exists an isomorphism of $(\tS,\cM_{\tS})$-logarithmic schemes
\begin{equation}\label{p2-dpsccd4b}
(\tX',\cM_{\tX'})\stackrel{\sim}{\rightarrow} (\tX,\cM_{\tX}),
\end{equation}
extending the identity of $(\coX,\cM_{\coX})$. We deduce that the category $\bHM^\sol_\tf(\co_{\fX}[\frac 1 p], \hOmega)$ does not depend on $\tf$, 
and an isomorphism of functors 
\begin{equation}\label{p2-dpsccd4c}
\cH_{\tf}\stackrel{\sim}{\rightarrow}\cH_{\tf'}.
\end{equation}
The first assertion follows. 

The second assertion is a special case of \ref{p2-fgpscp2}. 

\subsection{}\label{p2-dpscc5}
Let $\cM$ be a Dolbeault ind-$\bvocB$-module, $(\cN,\vartheta)=\upalpha_{\co_\fX}^\coh(\cH_\tf(\cM))$ \eqref{p2-cmupiso31h}.
By \ref{p2-rgpsc15} and the definition of the functor $\cV_{\tf}$ \eqref{p2-rgpsc14b}, we have a canonical isomorphism of ind-$\bvocB$-modules
\begin{equation}\label{p2-dpscc5a}
\cM\stackrel{\sim}{\rightarrow} \uupnu_\tf(\rI\hupsigma^*(\cN,\theta)),
\end{equation}
where the functor $\uupnu_\tf$ is defined in \eqref{p2-rgpsc1c}. We deduce a Higgs $\bvocB$-field $\theta_{\cM,\tf}$ on $\cM$ 
with coefficients in $\hupsigma^*(\hOmega)$, such that 
\eqref{p2-dpscc5a} is underlying an isomorphism of Higgs ind-$\bvocB$-modules
\begin{equation}\label{p2-dpscc5b}
(\cM,\theta_{\cM,\tf})\stackrel{\sim}{\rightarrow} \upnu_\tf(\rI\hupsigma^*(\cN,\theta)). 
\end{equation}
We thus define a functor 
\begin{equation}\label{p2-dpscc5c}
\begin{array}[t]{clcr}
\bIndMod^\Dolb(\bvocB)&\rightarrow &\bIndHM(\bvocB,\hupsigma^*(\hOmega)),\\
\cM&\mapsto & (\cM,\theta_{\cM,\tf}).
\end{array}
\end{equation}

There exist a rational number $r>0$ and an isomorphism of ind-$\bvcC^{(r)}_{\tf}$-modules with 
$\delta_{\bvcC^{(r)}_{\tf}}$-connection in the sense of \ref{p1-indmal22}, 
\begin{equation}\label{p2-dpscc5d}
\cM\otimes_{\bvocB}\bvcC^{(r)}_\tf \stackrel{\sim}{\rightarrow} \rI\hupsigma^*(\cN)\otimes_{\bvocB}\bvcC^{(r)}_\tf,
\end{equation}
where the $\delta_{\bvcC^{(r)}_{\tf}}$-connections are defined as in \ref{p1-delta-con9}, 
$\rI\hupsigma^*(\cN)$ (resp.\ $\cM$) being endowed with the Higgs field 
$\rI\hupsigma^*(\vartheta)$ (resp.\ $0$). The Higgs field $\theta_{\cM,\tf}$ on $\cM$ is then induced by the 
Higgs field $\rI\hupsigma^*(\vartheta)\otimes \id$ on $\rI\hupsigma^*(\cN)\otimes_{\bvocB}\bvcC^{(r)}_\tf$,  
or equivalently by the Higgs field $-\id\otimes \delta_{\bvcC^{(r)}_\tf}$.

For every rational number $t$ such that $0<t\leq r$, 
\eqref{p2-dpscc5d} induces by extension of scalars by the homomorphism $\bvupphi^{(r,t)}$ \eqref{p2-dpscc3c} 
an isomorphism of $\bvcC^{(r)}_{\tf'}\otimes_{\bvocB'}\hupsigma'^*(\hcC^{(t)}_\uptau)$-modules 
\begin{equation}\label{p2-dpscc5e}
\cM\otimes_{\bvocB}\bvcC^{(r)}_{\tf'}\otimes_{\bvocB}\hupsigma^*(\hcC^{(t)}_\uptau)
\stackrel{\sim}{\rightarrow} \rI\hupsigma^*(\cN)\otimes_{\bvocB}\bvcC^{(r)}_{\tf'}\otimes_{\bvocB}\hupsigma^*(\hcC^{(t)}_\uptau).
\end{equation}
It follows from \eqref{p2-fgpscp1o} that the Higgs field $\theta_{\cM,\tf}$ on $\cM$ is induced by the 
Higgs field $\rI\hupsigma^*(\vartheta)\otimes \id+\id\otimes \delta'^{\vee(r,t)}$ \eqref{p2-fgpscp1hv} 
on $\rI\hupsigma^*(\cN)\otimes_{\bvocB}\bvcC^{(r)}_{\tf'}\otimes_{\bvocB}\hupsigma^*(\hcC^{(t)}_\uptau)$.

\begin{prop}\label{p2-dpscc6}
For every Dolbeault ind-$\bvocB$-module $\cM$ such that the $\co_{\fX}[\frac 1 p]$-bundle $\cH_\tf(\cM)$ is locally CL-small, 
the canonical Higgs $\bvocB$-field $\theta_{\cM,\tf}$ \eqref{p2-dpscc5c} does not depend 
on the choice of the deformation $\tf$. 
\end{prop}

Observe however that  $\theta_{\cM,\tf}$ depends a priori on the relative or absolute case considered in \ref{p2-ncgt3}.

We keep the notation of \ref{p2-dpscc5} and set $(\cN',\vartheta')=\upalpha_{\co_\fX}^\coh(\cH_{\tf'}(\cM))$.  
By \ref{p2-dpscc4} and the proof of \eqref{p2-fgpscp2j1}, 
there exist a rational number $t>0$ and an isomorphism of ind-$\hcC^{(t)}_\uptau$-modules 
with $\delta^{(t)}_\uptau$-connection
\begin{equation}\label{p2-dpscc6a}
\cN\otimes_{\co_{\fX}}\hcC^{(t)}_\uptau\stackrel{\sim}{\rightarrow} \cN'\otimes_{\co_{\fX}}\hcC^{(t)}_\uptau,
\end{equation}
where the $\delta^{(t)}_{\uptau}$-connections are defined as in \ref{p1-delta-con4}, 
$\cN$ (resp.\ $\cN'$) being endowed with the Higgs field $\vartheta$ (resp.\ $0$). We may assume $t\leq r$. 
By \eqref{p2-fgpscp1k}, the isomorphism \eqref{p2-dpscc6a} induces by pullback by $\rI\hupsigma^*$ and extension of scalars, an isomorphism of 
$\bvcC^{(r)}_{\tf'}\otimes_{\bvocB}\hupsigma^*(\hcC^{(t)}_\uptau)$-modules  with $\delta'^{\vee(r,t)}$-connection \eqref{p2-fgpscp1hv}
\begin{equation}\label{p2-dpscc6b}
\rI\hupsigma^*(\cN)\otimes_{\bvocB}\bvcC^{(r)}_{\tf'}\otimes_{\bvocB'}\hupsigma^*(\hcC^{(t)}_\uptau)\stackrel{\sim}{\rightarrow} 
\rI \hupsigma^*(\cN')\otimes_{\bvocB}\bvcC^{(r)}_{\tf'}\otimes_{\bvocB'}\hupsigma^*(\hcC^{(t)}_\uptau),
\end{equation}
where the $\delta'^{\vee(r,t)}$-connections are defined as in \ref{p1-delta-con4}, 
$\cN$ (resp.\ $\cN'$) being endowed with the Higgs field $\vartheta$ (resp.\ $0$). 
In other words, the isomorphism \eqref{p2-dpscc6b} is compatible with the Higgs fields indicated below
\begin{equation}\label{p2-dpscc6c}
\begin{tabular}{|c|c|}
\hline
left hand side&right  hand  side\\
\hline
$\rI \hupsigma^*(\vartheta)\otimes \id+\id\otimes \delta'^{\vee(r,t)}$ & $\id\otimes \delta'^{\vee(r,t)}$ \\
\hline
\end{tabular}
\end{equation}

By the proof of \ref{p2-fgpscp2}, the composition of \eqref{p2-dpscc5e} and \eqref{p2-dpscc6b} 
induces an isomorphism of ind-$\bvcC^{(r)}_{\tf'}$-modules with 
$\delta_{\bvcC^{(r)}_{\tf'}}$-connection  
\begin{equation}\label{p2-dpscc6d}
\cM\otimes_{\bvocB}\bvcC^{(r)}_{\tf'}\stackrel{\sim}{\rightarrow} \rI\hupsigma^*(\cN')\otimes_{\bvocB}\bvcC^{(r)}_{\tf'},
\end{equation}
where the $\delta_{\bvcC^{(r)}_{\tf'}}$-connections are defined as in \ref{p1-delta-con9}, 
$\rI\hupsigma^*(\cN')$ (resp.\ $\cM$) being endowed with the Higgs field 
$\rI\hupsigma^*(\vartheta')$ (resp.\ $0$). We deduce a commutative diagram
\begin{equation}\label{p2-dpscc6e}
\xymatrix{
\cM\ar@{=}[r]\ar[d]&\cM\ar@{=}[r]\ar[d]&\cM\ar[d]\\
{\cM\otimes_{\bvocB}\bvcC^{(r)}_\tf}\ar[r]^-(0.5){\id\otimes \bvupphi^{(r,t)}}\ar[d]&
{\cM\otimes_{\bvocB}\bvcC^{(r)}_{\tf'}\otimes_{\bvocB}\hupsigma^*(\hcC^{(t)}_\uptau)}\ar[d]&{\cM\otimes_{\bvocB}\bvcC^{(r)}_{\tf'}}\ar[l]\ar[dd]\\
{\rI\hupsigma^*(\cN)\otimes_{\bvocB}\bvcC^{(r)}_\tf}\ar[r]^-(0.5){\id\otimes \bvupphi^{(r,t)}}&
{\rI\hupsigma^*(\cN)\otimes_{\bvocB}\bvcC^{(r)}_{\tf'}\otimes_{\bvocB}\hupsigma^*(\hcC^{(t)}_\uptau)}\ar[d]& \\
&{\rI\hupsigma^*(\cN')\otimes_{\bvocB}\bvcC^{(r)}_{\tf'}\otimes_{\bvocB'}\hupsigma^*(\hcC^{(t)}_\uptau)}&
{\rI\hupsigma^*(\cN')\otimes_{\bvocB}\bvcC^{(r)}_{\tf'},}\ar[l]}
\end{equation}
where the unlabelled horizontal and top vertical arrows are the canonical morphisms. 
The Higgs field $\theta_{\cM,\tf'}$ on $\cM$ is induced by the 
Higgs field $-\id\otimes \delta_{\bvcC^{(r)}_{\tf'}}$ on $\rI\hupsigma^*(\cN')\otimes_{\bvocB}\bvcC^{(r)}_{\tf'}$,  
and thus also by the Higgs field $\id\otimes \delta'^{\vee(r,t)}$ 
on $\rI\hupsigma^*(\cN')\otimes_{\bvocB}\bvcC^{(r)}_{\tf'}\otimes_{\bvocB}\hupsigma^*(\hcC^{(t)}_\uptau)$. 
In view of \eqref{p2-dpscc6c} and the remark below \eqref{p2-dpscc5e}, we deduce  
that $\theta_{\cM,\tf}=\theta_{\cM,\tf'}$, which proves the proposition. 

\begin{rema}
Canonical Higgs fields similar to those constructed in \ref{p2-dpscc6}  
were introduced by Camargo on pro-étale vector bundles (\cite{rod} 1.0.3), and by He in his Sen theory (\cite{he} 5.28). 
They were subsequently used by Heuer \cite{heuer} and Heuer-Xu \cite{hexu} in their constructions of the $p$-adic Simpson correspondence.
\end{rema}

\end{document}